\documentclass[a4paper, 10pt, openany]{book}
\usepackage{graphicx,graphics}
\usepackage[cp1251]{inputenc}
\usepackage[russian]{babel}
\usepackage{calc}

\usepackage{multicol}
\usepackage[cmtip,arrow,all,2cell]{xy}
\usepackage{pb-diagram,pb-xy}
\usepackage{amsmath,amsthm,amssymb,tabularx,stmaryrd}
\usepackage[mathscr]{eucal}
\usepackage{hyperref}
\usepackage{mathabx} 
\usepackage{skak}
\usepackage{mathdots}
\usepackage{shorttoc}

\usepackage{comment}

\usepackage{tikz}
\usepackage{tikz-3dplot}
\usetikzlibrary{arrows.meta}
\usetikzlibrary{patterns}
\usetikzlibrary{patterns.meta}
\usepackage{pgfplots}

\allowdisplaybreaks 

\makeindex

\theoremstyle{plain}

\newtheorem{tm}{Теорема}[section]

\newtheorem{lm}[tm]{Лемма}
\newtheorem{prop}[tm]{Предложение}
\newtheorem{cor}[tm]{Следствие}

\newtheorem{quest}[tm]{Вопрос}

\theoremstyle{definition}

\newtheorem{ex}[tm]{$\diamond$}
\newtheorem{exs}[tm]{$\diamond\diamond$}
\newtheorem{er}[tm]{$\triangleright$}
\newtheorem{ers}[tm]{$\triangleright\triangleright$}
\newtheorem{rem}[tm]{$\bold{!}$}
\newtheorem{df}[tm]{$\bullet$}

\newcommand{\beq}{\begin{equation}}
\newcommand{\eeq}{\end{equation}}
\newcommand{\btm}{\begin{tm}}
\newcommand{\etm}{\end{tm}}
\newcommand{\blm}{\begin{lm}}
\newcommand{\elm}{\end{lm}}
\newcommand{\bcor}{\begin{cor}}
\newcommand{\ecor}{\end{cor}}
\newcommand{\bex}{\begin{ex}}
\newcommand{\eex}{\end{ex}}
\newcommand{\bexs}{\begin{exs}}
\newcommand{\eexs}{\end{exs}}
\newcommand{\bextm}{\begin{extm}}
\newcommand{\eextm}{\end{extm}}
\newcommand{\bcx}{\begin{cex}}
\newcommand{\ecx}{\end{cex}}
\newcommand{\bers}{\begin{ers}}
\newcommand{\eers}{\end{ers}}
\newcommand{\ber}{\begin{er}}
\newcommand{\eer}{\end{er}}
\newcommand{\bdf}{\begin{df}}
\newcommand{\edf}{\end{df}}
\newcommand{\brem}{\begin{rem}}
\newcommand{\erem}{\end{rem}}
\newcommand{\bpr}{\begin{proof}}
\newcommand{\epr}{\end{proof}}

\newcommand{\bprop}{\begin{prop}}
\newcommand{\eprop}{\end{prop}}

\renewcommand{\thesection}{\S\,\arabic{section}}
\renewcommand{\thesubsection}{(\alph{subsection})}

 \makeatletter \@addtoreset {equation}{section}
 \makeatother

\input{opus.sty}

\begin{document}

\title{СТЕРЕОТИПНЫЕ ДВОЙСТВЕННОСТИ В ГЕОМЕТРИИ}

\author{С.С.Акбаров}

\setcounter{tocdepth}{4} 

\maketitle

\frontmatter

\shorttoc{Оглавление}{5} 

\chapter{ВВЕДЕНИЕ}

\section{Геометрические дисциплины как проекции теории топологических алгебр}

\subsection{Программа Клейна и геометрические дисциплины}

В 1872 году Феликс Клейн выступил в университете Эрлангена с тезисами своей знаменитой программы классификации геометрических теорий. Предложенный им взгляд на геометрию существенно повлиял на развитие этой науки, став в 20 веке одним из направлений общей работы по систематизации математических знаний. В настоящее время идеи Клейна эволюционировали в точку зрения, что геометрия рассматриваемого пространства $X$ определяется представлением его в виде подходящего типа многообразия $M$ с дополнительными структурами типа метрики, связности, кривизны и т.д., описывающими его естественные свойства.

А тип многообразия $M$ в этой схеме определяется выбором одной из четырех главных геометрических дисциплин, объектом которой $M$ должно считаться:
\bit{
\item[---] топологии,
\item[---] дифференциальной геометрии,
\item[---] комплексной геометрии,
\item[---] алгебраической геометрии.
}\eit

Философия этих дисциплин так похожа, что логично было бы думать, что они сами должны объединяться в какую-то общую схему типа клейновской, объясняющую параллели между ними и позволяющую обсуждать различия на формальном языке. Но ни в программе Клейна, ни в современных ее модификациях это сходство не объясняется. Откуда в математике появляются эти области, почему они так странно похожи, почему их четыре --- до последнего времени оставалось непонятно, и только относительно недавно появились формальные средства, позволяющие объяснить это.

Ответ пришел с неожиданной стороны: из Функционального анализа, причем лежащая в его основе интуитивная идея заимствуется из физики и описывается правилом

\medskip
\centerline{\it образ изучаемого объекта зависит от инструментов наблюдения.}
\medskip
\noindent
Подобно тому как в астрономии образ исследуемого объекта меняется при изменении инструментов наблюдения – оптического телескопа, радиотелескопа, рентгеновского телескопа, и т.д. – в математике большие геометрические дисциплины из списка выше производят впечатление образов, возникающих при изучении одной и той же реальности при помощи различных инструментов наблюдения. 

Я думаю, это уместно проиллюстрировать картинкой в таком духе:
\vglue-20pt
\beq\label{metaphysics-1}
  \begin{matrix}
  \centering
  \begin{tikzpicture}
    \begin{axis}[hide axis,xmin=-22,xmax=22,ymin=-10,ymax=10, height=9cm, width=15cm]
      \draw[line width=0.2mm] (axis cs: -15, 0) circle [x radius = 50, y radius = 20];
      \node at (axis cs:-15, 0) {\text{Реальность}};
      \draw[line width=0.2mm] (axis cs: 12, 6) circle [x radius = 90, y radius = 16];
      \node at (axis cs:12, 6) {\text{Топология}};
      \draw[line width=0.2mm] (axis cs: 12, 2) circle [x radius = 90, y radius = 16];
      \node at (axis cs:12, 2) {\text{Дифференциальная геометрия}};
      \draw[line width=0.2mm] (axis cs: 12, -2) circle [x radius = 90, y radius = 16];
      \node at (axis cs:12, -2) {\text{Комплексная геометрия}};
      \draw[line width=0.2mm] (axis cs: 12, -6) circle [x radius = 90, y radius = 16];
      \node at (axis cs:12, -6) {\text{Алгебраическая геометрия}};
      \node at (axis cs:-9, 6.5) {\text{инструменты наблюдения}};
    \draw[-](axis cs:-9.5, 1.5) parabola bend (axis cs:2, 6) (axis cs:2, 6);
    \draw[line width=1pt, ->, -stealth] (axis cs:2.1,6) -- (axis cs:2.2,6);
    \draw[-](axis cs:-9, 0.5) parabola bend (axis cs:2, 2) (axis cs:2, 2);
    \draw[line width=1pt, ->, -stealth] (axis cs:2.1,2) -- (axis cs:2.2,2);
    \draw[-](axis cs:-9, -0.5) parabola bend (axis cs:2, -2) (axis cs:2, -2);
    \draw[line width=1pt, ->, -stealth] (axis cs:2.1,-2) -- (axis cs:2.2,-2);
    \draw[-](axis cs:-9.5, -1.5) parabola bend (axis cs:2, -6) (axis cs:2, -6);
    \draw[line width=1pt, ->, -stealth] (axis cs:2.1,-6) -- (axis cs:2.2,-6);
    \end{axis}
  \end{tikzpicture}
  \end{matrix}
\eeq
\vglue-20pt
Удивительно, что за этой интуицией действительно стоит некая закономерность: по крайней мере три из четырех упомянутых геометрических дисциплин --- топология, дифференциальная геометрия и комплексная геометрия --- представляют собой картины, возникающие при исследовании различными инструментами одной и той же общей для всех реальности. В нулевом приближении такой реальностью является, как оказывается, теория топологических алгебр, а результаты, о которых идет речь,  выглядят так (см. \cite{Ak17-1,Ak17-2,Ak08}): 
\bit{\it 

\item[---] использование C*-алгебр в качестве инструментов наблюдения превращает теорию топологических алгебр в топологию,
\item[---] использование C*-алгебр с присоединенными самосопряженными нильпотентными элементами\footnote{См. определение на с.\pageref{DEF:C^*[m]}.} в качестве инструментов наблюдения превращает теорию топологических алгебр в дифференциальную геометрию, а
\item[---] использование банаховых алгебр в качестве инструментов наблюдения превращает теорию топологических алгебр в комплексную геометрию.
}\eit
\noindent
Картинка \eqref{metaphysics-1} при этом превращается в такую\footnote{Здесь под $C^*[m]$-алгебрами понимаются $C^*$-алгебры с присоединенными самосопряженными нильпотентными элементами, а под $B$-алгебрами --- банаховы алгебры, а знак вопроса -- ? -- означает, что мы пока не знаем, какими должны быть инструменты наблюдения, чтобы под их действием теория топологических алгебр превращалась в алгебраическую геометрию.}:
\vglue-20pt
\beq\label{metaphysics-2}
\begin{matrix}
  \centering
  \begin{tikzpicture}
    \begin{axis}[hide axis,xmin=-22,xmax=22,ymin=-10,ymax=10, height=10cm, width=16cm]
      \draw[line width=0.2mm] (axis cs: -15, 0) circle [x radius = 70, y radius = 20];
      \node at (axis cs:-15, .5) {\text{Теория}};
      \node at (axis cs:-15, -.5) {\text{топологических алгебр}};
      \draw[line width=0.2mm] (axis cs: 12, 6) circle [x radius = 90, y radius = 16];
      \node at (axis cs:12, 6) {\text{Топология}};
      \draw[line width=0.2mm] (axis cs: 12, 2) circle [x radius = 90, y radius = 16];
      \node at (axis cs:12, 2) {\text{Дифференциальная геометрия}};
      \draw[line width=0.2mm] (axis cs: 12, -2) circle [x radius = 90, y radius = 16];
      \node at (axis cs:12, -2) {\text{Комплексная геометрия}};
      \draw[line width=0.2mm] (axis cs: 12, -6) circle [x radius = 90, y radius = 16];
      \node at (axis cs:12, -6) {\text{Алгебраическая геометрия}};
      \node at (axis cs:-4, 6) {\text{$C^*$-алгебры}};
      \node at (axis cs:-1, 3) {\text{$C^*[m]$-алгебры}};
      \node at (axis cs:-1, -.5) {\text{$B$-алгебры}};
      \node at (axis cs:-1, -4.5) {\text{?}};
          \draw[-](axis cs:-7.5, 1.5) parabola bend (axis cs:2, 6) (axis cs:2, 6);
    \draw[line width=1pt, ->, -stealth] (axis cs:2.1,6) -- (axis cs:2.2,6);
    \draw[-](axis cs:-7, 0.5) parabola bend (axis cs:2, 2) (axis cs:2, 2);
    \draw[line width=1pt, ->, -stealth] (axis cs:2.1,2) -- (axis cs:2.2,2);
    \draw[-](axis cs:-7, -0.5) parabola bend (axis cs:2, -2) (axis cs:2, -2);
    \draw[line width=1pt, ->, -stealth] (axis cs:2.1,-2) -- (axis cs:2.2,-2);
    \draw[dashed](axis cs:-7.5, -1.5) parabola bend (axis cs:2, -6) (axis cs:2, -6);
    \draw[line width=1pt, ->, -stealth] (axis cs:2.1,-6) -- (axis cs:2.2,-6);
    \end{axis}
  \end{tikzpicture}
  \end{matrix}
\eeq
\vglue-20pt
В таком нарочито туманном, метафизическом, виде эти заявления понять, конечно, невозможно, они могут только создать некие зрительные образы. Чтобы придать им точный смысл, необходима некая система определений.

\subsection{Точные формулировки}

Языком описания этих закономерностей служит язык теории категорий, а формальной конструкцией в нем для  этих целей --- понятие оболочки \cite{Ak16,Akbarov-De-Gruyter-I}. Это вот что такое. Пусть нам дана категория $\tt K$ и два класса морфизмов в ней:
 \bit{
\item[---] класс $\varPhi$ морфизмов, называемых {\it тестовыми морфизмами} (или, собственно, {\it инструментами наблюдения}), и

\item[---] класс $\varOmega$ морфизмов, называемых {\it представлениями}.

}\eit
Тогда

 \bit{
\item[$\bullet$] Морфизм $\sigma:X\to X'$ в категории $\tt K$ называется
{\it расширением объекта $X\in\Ob({\tt K})$ в классе морфизмов $\varOmega$ относительно класса
морфизмов $\varPhi$}, если $\sigma\in\varOmega$, и для любого морфизма
$\ph:X\to B$ из класса $\varPhi$ найдется единственный морфизм $\ph':X'\to B$ в
категории ${\tt K}$, замыкающий диаграмму
 \beq\label{DEF:diagr-rasshirenie}
\begin{diagram}
\node[2]{X} \arrow{sw,t}{\varOmega\owns\sigma} \arrow{se,t}{\forall\ph\in\varPhi}\\
\node{X'}  \arrow[2]{e,b,--}{\exists!\ph'} \node[2]{B}
\end{diagram}
\eeq

\item[$\bullet$] Расширение $\rho:X\to E$ объекта
$X\in\Ob({\tt K})$ в классе морфизмов $\varOmega$
относительно класса морфизмов $\varPhi$ называется {\it оболочкой объекта
$X$ в классе морфизмов $\varOmega$
относительно класса морфизмов $\varPhi$},
если для любого другого расширения  $\sigma:X\to X'$ (объекта $X$
в классе морфизмов $\varOmega$  относительно класса
морфизмов $\varPhi$) найдется единственный морфизм $\upsilon:X'\to E$ в $\tt
K$, замыкающий диаграмму
 \beq\label{DEF:diagr-obolochka}
\begin{diagram}
\node[2]{X} \arrow{sw,t}{\forall\sigma} \arrow{se,t}{\rho}\\
\node{X'}  \arrow[2]{e,b,--}{\exists!\upsilon} \node[2]{E}
\end{diagram}
 \eeq
}\eit
Объект $E$ также называется {\it оболочкой} объекта $X$
(в классе морфизмов $\varOmega$ относительно класса
морфизмов $\varPhi$), и обозначается
    \beq\label{DEF:Env_F^L(X)}
E=\Env_{\varPhi}^\varOmega X.
 \eeq

Это определение позволяет описать упомянутые выше закономерности на формальном языке, только в этих формулировках нужно будет еще заменить понятие топологической алгебры на близкое ему понятие стереотипной алгебры. Чтобы объяснить, что это такое, нужно сначала дать определение стереотипного пространства \cite{Ak03,Akbarov-De-Gruyter-I}. Локально выпуклое пространство $X$ над полем $\C$ называется {\it стереотипным}, если оно изоморфно своему второму сопряженному пространству
\beq\label{DEF:ste}
(X^\star)^\star\cong X,
\eeq 
где под $X^\star$ понимается пространство линейных непрерывных функционалов $f:X\to\C$, наделенное топологиией равномерной сходимости на {\it  вполне ограниченных множествах}. {\it Стереотипной алгеброй} называется произвольное стереотипное пространство $A$, наделенное одновременно структурой ассоциативной алгебры над $\C$ так, чтобы операция умножения была непрерывна в следующем смысле: для всякой окрестности нуля $U\subseteq A$ и любого вполне ограниченного множества $S\subseteq A$ найдется окрестность нуля $V\subseteq A$ такая, что
\beq\label{DEF:ste-alg}
V\cdot S\subseteq U\quad \&\quad S\cdot V\subseteq U.
\eeq 
Классы стереотипных пространств и алгебр весьма широки, настолько, что можно считать, что все локально выпуклые пространства и топологические алгебры, реально испрользуемые в Анализе стереотипны, см. подробности в \cite{Akbarov-De-Gruyter-I}. 

Теперь мы можем дать аккуратную формулировку приведенных выше утверждений об оболочках топологических алгебр  \cite{Ak17-1,Ak17-2,Ak08}. 

\bit{

\item[1.] Если под $\varPhi$ понимать класс инволютивных морфизмов со значениями в $C^*$-алгебрах, а под $\varOmega$, класс морфизмов с плотным образом в области значений, то под действием соответствующего функтора оболочки $\Env_{\varPhi}^\varOmega$, называемого {\it непрерывной оболочкой}, и обозначаемого $\Env_{\mathcal C}$, категория инволютивных стереотипных алгебр превращается в категорию алгебр со свойствами, сильно напоминающими свойства алгебр ${\mathcal C}(M)$ непрерывных функций на паракомпактных локально компактных топологических пространствах $M$. Например, если $A={\mathcal E}(\R)$ --- алгебра гладких функций на вещественной прямой $\R$ (с обычными поточечными алгебраическими операциями и обычной топологией равномерной сходимости на компактах по каждой производной), то непрерывной оболочкой $A$ будет в точности алгебра ${\mathcal C}(\R)$ непрерывных функций на $\R$:
    $$
    \Env_{\mathcal C}{\mathcal E}(\R)={\mathcal C}(\R).
    $$
    В общей ситуации справедлива теорема: если $A$ --- инволютивная стереотипная подалгебра в алгебре ${\mathcal C}(M)$ гладких функций на паракомпрактном локально компактном топологическом пространстве $M$, причем спектр $A$ совпадает с $M$, то непрерывной оболочкой $A$ будет в точности алгебра ${\mathcal C}(M)$ непрерывных функций на $M$:
    $$
    \Env_{\mathcal C}A={\mathcal C}(M).
    $$

\item[2.] Если под $\varPhi$ понимать класс так называемых дифференциальных инволютивных морфизмов со значениями в $C^*$-алгебрах с присоединенными самосопряженными нильпотентнми элементами, а под $\varOmega$, как и раньше, класс морфизмов с плотным образом в области значений, то под действием соответствующего функтора оболочки $\Env_{\varPhi}^\varOmega$, называемого {\it гладкой оболочкой}, и обозначаемого $\Env_{\mathcal E}$,  категория стереотипных алгебр превращается в категорию алгебр со свойствами, странно напоминающими свойства алгебр ${\mathcal E}(M)$ гладких функций на гладких многообразиях $M$. Как иллюстрация,  если $A$ --- инволютивная стереотипная подалгебра в алгебре ${\mathcal E}(M)$ гладких функций на гладком многообразии $M$, причем спектр $A$ совпадает с $M$, а касательное пространство к $A$ в каждой точке $t\in M$ совпадает с касательным пространством к $M$ в этой точке, то гладкой оболочкой $A$ будет в точности алгебра ${\mathcal E}(M)$ гладких функций на $M$:
    $$
    \Env_{\mathcal E}A={\mathcal E}(M).
    $$

\item[3.] Если под $\varPhi$ понимать класс морфизмов со значениями в банаховых алгебрах, а под $\varOmega$, по-прежнему, класс морфизмов с плотным образом в области значений, то под действием соответствующего функтора оболочки $\Env_{\varPhi}^\varOmega$, называемого {\it голоморфной оболочкой}, и обозначаемого $\Env_{\mathcal O}$, категория стереотипных алгебр превращается в категорию алгебр по свойствам удивительно похожим на свойства алгебр ${\mathcal O}(M)$ голоморфных функций на комплексных многообразиях $M$. Например, по известной теореме А.Пирковского, голоморфной оболочкой алгебры ${\mathcal P}(M)$ многочленов на комплексном аффинном алгебраическом многообразии $M$ является в точности алгебра ${\mathcal O}(M)$ голоморфных функций на $M$:
    $$
    \Env_{\mathcal O}{\mathcal P}(M)={\mathcal O}(M)
    $$

}\eit

\section{Стереотипная двойственность в геометрии}

Глубокими (и неожиданными) результатами, оправдывающими интерес к оболочкам в функциональном анализе (независимо от описанных выше связей с геометрией), являются теоремы об обобщении теории двойственности Понтрягина на некоммутативные группы. Дело в том, что каждая оболочка из приведеннго выше списка --- $\Env_{\mathcal C}$, $\Env_{\mathcal E}$,  $\Env_{\mathcal O}$ --- порождает  некую  двойственность в соответствующей геометрической дисциплине, --- комплексной геометрии, дифференциальной геометрии, топологии, --- обобщающую обычную двойственность Понтрягина для абелевых локально компактных группп.

\subsection{Стереотипная двойственность в топологии}

Начать удобно с непрерывной оболочки $\Env_{\mathcal C}$. Оказывается, что  для любой (необязательно коммутативной) группы Мура $G$ справедливо утверждение,  иллюстрируемое следующей диаграммой \cite{Kuznetsova,Ak-Moore}:
 \beq\label{chetyrehugolnik-C-C*-0}
  \xymatrix @R=1.pc @C=2.pc
 {
 {\mathcal C}^\star(G)
 & \ar@{|->}[r]^{\Env_{\mathcal C}} & &
 \Env_{\mathcal C}{\mathcal C}^\star(G)
 \\
 & & &
 \ar@{|->}[d]^{\star}
 \\
 \ar@{|->}[u]^{\star}
 & & &
 \\
 {\mathcal C}(G)
 & &
 \ar@{|->}[l]_{\Env_{\mathcal C}}
 &
 \Env_{\mathcal C}{\mathcal C}^\star(G)^\star
 },
 \eeq
 Понимать ее нужно так: 

\bit{

\item[---] мы начинаем с алгебры ${\mathcal C}(G)$ непрерывных функций на группе $G$, и применяем к ней операцию $\star$ перехода к двойственному стереотипному пространству;
    
\item[---] в результате мы получаем алгебру ${\mathcal C}^\star(G)$ мер с компактным носителем на группе $G$; мы применяем к ней операцию $\Env_{\mathcal C}$ перехода к непрерывной оболочке;
    
\item[---] в результате мы получаем некую стереотипную алгебру $\Env_{\mathcal C}{\mathcal C}^\star(G)$; мы применяем к ней операцию $\star$ перехода к двойственному стереотипному пространству;
     
\item[---] чудесным образом оказывается, что $\Env_{\mathcal C}{\mathcal C}^\star(G)^\star$ --- тоже стереотипная алгебра относительно естественного умножения; мы применяем к ней операцию  $\Env_{\mathcal C}$ перехода к непрерывной оболочке;
    
\item[---] и вдруг обнаруживается, что результатом будет в точности исходная алгебра ${\mathcal C}(G)$ непрерывных функций на группе $G$.

}\eit

С точки зрения стереотипной теории, групповой алгеброй группы $G$ является алгебра ${\mathcal C}^\star(G)$ мер с компактным носителем на $G$. Поэтому начинать двигаться по диаграмме \eqref{chetyrehugolnik-C-C*-0} правильно с левого верхнего угла, того, где находится алгебра ${\mathcal C}^\star(G)$. И для нее круг по диаграмме будет выглядеть как равенство
\beq\label{chetyrehugolnik-C-C*-1}
\Env_{\mathcal C}\Big(\Env_{\mathcal C}{\mathcal C}^\star(G)^\star\Big)^\star\cong {\mathcal C}^\star(G)
\eeq
Если композицию операций $\Env_{\mathcal C}$ и $\star$ обозначит каким-нибудь символом, например, $\dagger$,
\beq\label{chetyrehugolnik-C-C*-2}
\dagger=\star\circ\Env_{\mathcal C},
\eeq
то равенство \eqref{chetyrehugolnik-C-C*-1} можно будет переписать в виде
\beq\label{chetyrehugolnik-C-C*-3}
\Big({\mathcal C}^\star(G)^\dagger\Big)^\dagger\cong {\mathcal C}^\star(G)
\eeq
И это логично интерпретировать как рефлексивность алгебры ${\mathcal C}^\star(G)$ относительно операции $\dagger$. Или, правильнее было бы сказать, {\it рефлексивность алгебры Хопфа ${\mathcal C}^\star(G)$ относительно операции $\dagger$}, потому что ${\mathcal C}^\star(G)$ является алгеброй Хопфа относительно естественного тензорного произведения $\circledast$ в стереотипной теории.  

Алгебры Хопфа $H$, удовлетворяющие равенству
\beq\label{DEF:neppr-reflex}
\Big(H^\dagger\Big)^\dagger\cong H
\eeq
называются {\it рефлексивными относительно непрерывной оболочки $\Env_{\mathcal C}$}.

Вторая важная деталь в этой картине состоит в том, что диаграмма \eqref{chetyrehugolnik-C-C*-0} дополняется следующей диаграммой функторов  
{\sf
 \beq\label{diagramma-kategorij-dlya-nepreryvnoi-obolochki}
 \xymatrix @R=3.pc @C=2.pc
 {
 \boxed{\begin{matrix}
  \text{алгебры Хопфа в $(\Ste,\circledast)$,}\\
  \text{рефлексивные относительно}\\
  \text{непрерывной оболочки}
 \end{matrix}}
 \ar[rr]^{H\mapsto H^\dagger} & &
 \boxed{\begin{matrix}
  \text{алгебры Хопфа в $(\Ste,\circledast)$,}\\
  \text{рефлексивные относительно}\\
  \text{непрерывной оболочки}
 \end{matrix}}
 \\
 \boxed{\text{группы Мура}} \ar[u]^(0.45){\scriptsize\begin{matrix} \mathcal{C}^\star(G)\\
 \text{\rotatebox{90}{$\mapsto$}} \\ G\end{matrix}} & &
 \boxed{\text{группы Мура}} \ar[u]_(0.45){\scriptsize\begin{matrix} \mathcal{C}^\star(G) \\
 \text{\rotatebox{90}{$\mapsto$}} \\ G\end{matrix}} \\
  \boxed{\begin{matrix}
 \text{абелевы}\\
 \text{локально компактные группы}
 \end{matrix}} \ar[u]^{e} \ar[rr]^{G\mapsto \widehat{G}} & &
  \boxed{\begin{matrix}
 \text{абелевы}\\
 \text{локально компактные группы}
 \end{matrix}}\ar[u]_{e}
 }
 \eeq }\noindent
в которой $\widehat{G}$ --- двойственная по Понтрягину группа к абелевой локально копактной группе $G$ (то есть группа характеров $\chi:G\to\T$ с поточечными алгебраическими операциями и топологией равномерной сходимости на компактах).
 
Ее нужно понимать так: 

\bit{

\item[---] если зафиксировать локально компактную абелеву группу $G$, затем перейти к ее групповой алгебре мер ${\mathcal C}^\star(G)$, а после этого к ее двойственной алгебре Хопфа ${\mathcal C}^\star(G)^\dagger$ относительно операции \eqref{chetyrehugolnik-C-C*-2}, то полученная алгебра Хопфа будет естественно изоморфной алгебре Хопфа, которая получилась бы если бы мы двигались по диаграмме \eqref{diagramma-kategorij-dlya-nepreryvnoi-obolochki} из левого нижнего угла в правый верхний другим возможным путем, сначала взяли бы двойственную по Понтрягину группу $\widehat{G}$, а потом ее групповую алгебру мер ${\mathcal C}^\star(\widehat{G})$: 
\beq\label{diagramma-kategorij-dlya-nepreryvnoi-obolochki-1}
{\mathcal C}^\star(G)^\dagger\cong {\mathcal C}^\star(\widehat{G})
\eeq    

}\eit

Тождества \eqref{chetyrehugolnik-C-C*-3} и \eqref{diagramma-kategorij-dlya-nepreryvnoi-obolochki-1} вместе интерпретируются как 

\btm\label{TH:neprer-dvoistv}  \cite{Ak-Moore}
Двойственность алгебр Хопфа \eqref{chetyrehugolnik-C-C*-3} продолжает классическую двойственность Понтрягина для абелевых локально компактных групп на класс групп Мура.
\etm

\subsection{Стереотипная двойственность в дифференциальной геометрии}

Для теорий двойственности в дифференциальной и комплексной геометриях такого же уровня элегантности в результатах пока получить не удалось, поэтому мне кажется, будет правильно сформулировать отдельно какими их желательно было бы видеть, и какими они не сегодняшний день получаются.

\paragraph{Желательный вид теории в дифференциальной геометрии.}    

Результат, аналогичный теореме \ref{TH:neprer-dvoistv}, в дифференциальной геометрии должен был бы выглядеть следующим образом. Нужно определить понятие гладкой оболочки $\Env_{\mathcal E}$ так, чтобы когда в определении \eqref{chetyrehugolnik-C-C*-2} операции $\dagger$ заменяешь $\Env_{\mathcal C}$ на $\Env_{\mathcal E}$,
\beq\label{chetyrehugolnik-E-E*-2}
\dagger=\star\circ\Env_{\mathcal E},
\eeq
тождество \eqref{chetyrehugolnik-C-C*-3} превратится в тождество
\beq\label{chetyrehugolnik-E-E*-3}
\Big({\mathcal E}^\star(G)^\dagger\Big)^\dagger\cong {\mathcal E}^\star(G)
\eeq
справедливое для всех (необязательно коммутативных) компактно порожденных групп Ли---Мура $G$. 
Алгебры Хопфа $H$, удовлетворяющие соответствующему равенству \eqref{DEF:neppr-reflex}
естественно тогда назвать {\it рефлексивными относительно гладкой оболочки $\Env_{\mathcal E}$},  
а от диаграммы \eqref{diagramma-kategorij-dlya-nepreryvnoi-obolochki} естественно ожидать, что она в этом случае превратится в диаграмму
{\sf
 \beq\label{diagramma-kategorij-dlya-gladkoi-obolochki-ideal}
 \xymatrix @R=3.pc @C=2.pc
 {
 \boxed{\begin{matrix}
  \text{алгебры Хопфа в $(\Ste,\circledast)$,}\\
  \text{рефлексивные относительно}\\
  \text{гладкой оболочки}
 \end{matrix}}
 \ar[rr]^{H\mapsto H^\dagger} & &
 \boxed{\begin{matrix}
  \text{алгебры Хопфа в $(\Ste,\circledast)$,}\\
  \text{рефлексивные относительно}\\
  \text{гладкой оболочки}
 \end{matrix}}
 \\
 \boxed{\begin{matrix}
 \text{компактно порожденные}\\
 \text{группы Ли---Мура}
 \end{matrix}} \ar[u]^(0.45){\scriptsize\begin{matrix} \mathcal{C}^\star(G)\\
 \text{\rotatebox{90}{$\mapsto$}} \\ G\end{matrix}} & &
 \boxed{\begin{matrix}
 \text{компактно порожденные}\\
 \text{группы Ли---Мура}
 \end{matrix}} \ar[u]_(0.45){\scriptsize\begin{matrix} \mathcal{C}^\star(G) \\
 \text{\rotatebox{90}{$\mapsto$}} \\ G\end{matrix}} \\
  \boxed{\begin{matrix}
 \text{компактно порожденные абелевы}\\
 \text{группы Ли}
 \end{matrix}} \ar[u]^{e} \ar[rr]^{G\mapsto \widehat{G}} & &
  \boxed{\begin{matrix}
 \text{компактно порожденные абелевы}\\
 \text{группы Ли}
 \end{matrix}}\ar[u]_{e}
 }
 \eeq }
\noindent
в которой $\widehat{G}$ --- как и раньше, двойственная по Понтрягину группа к абелевой локально копактной группе $G$ (то есть группа характеров $\chi:G\to\T$ с поточечными алгебраическими операциями и топологией равномерной сходимости на компактах). А теорема \ref{TH:neprer-dvoistv} --- в теорему 

\btm\label{TH:gladk-dvoistv-ideal} 
Двойственность алгебр Хопфа \eqref{chetyrehugolnik-E-E*-3} продолжает классическую двойственность Понтрягина для компактно порожденных абелевых групп Ли на класс компактно порожденных групп Ли---Мура.
\etm

\paragraph{Полученные на сегодняшний день результаты.}    

Выбранное на с. \pageref{DEF:C^infty-obolochka} определение гладкой оболочки не позволяет доказать теорему \ref{TH:gladk-dvoistv-ideal}. Диаграмма \eqref{diagramma-kategorij-dlya-gladkoi-obolochki-ideal} для нее превращается в диаграмму попроще
{\sf
 \beq\label{diagramma-kategorij-dlya-gladkoi-obolochki-ideal}
 \xymatrix @R=3.pc @C=2.pc
 {
 \boxed{\begin{matrix}
  \text{алгебры Хопфа в $(\Ste,\circledast)$,}\\
  \text{рефлексивные относительно}\\
  \text{гладкой оболочки}
 \end{matrix}}
 \ar[rr]^{H\mapsto H^\dagger} & &
 \boxed{\begin{matrix}
  \text{алгебры Хопфа в $(\Ste,\circledast)$,}\\
  \text{рефлексивные относительно}\\
  \text{гладкой оболочки}
 \end{matrix}}
 \\
 \boxed{\begin{matrix}
 \text{группы вида $Z\times K$,}\\
 \text{где $Z$ --- компактно порожденная}\\
 \text{абелева группа Ли, }\\
 \text{а $K$ --- компактная группа Ли}
 \end{matrix}} \ar[u]^(0.45){\scriptsize\begin{matrix} \mathcal{C}^\star(G)\\
 \text{\rotatebox{90}{$\mapsto$}} \\ G\end{matrix}} & &
 \boxed{\begin{matrix}
\text{группы вида $Z\times K$,}\\
 \text{где $Z$ --- компактно порожденная}\\
 \text{абелева группа Ли, }\\
 \text{а $K$ --- компактная группа Ли}
 \end{matrix}} \ar[u]_(0.45){\scriptsize\begin{matrix} \mathcal{C}^\star(G) \\
 \text{\rotatebox{90}{$\mapsto$}} \\ G\end{matrix}} \\
  \boxed{\begin{matrix}
 \text{конечные абелевы группы}
 \end{matrix}} \ar[u]^{e} \ar[rr]^{G\mapsto \widehat{G}} & &
  \boxed{\begin{matrix}
 \text{конечные абелевы группы}
 \end{matrix}}\ar[u]_{e}
 }
 \eeq }
а теорема \ref{TH:gladk-dvoistv-ideal} --- в теорему  

\btm\label{TH:gladk-dvoistv}  \cite{Ak17-2}
Двойственность алгебр Хопфа \eqref{chetyrehugolnik-E-E*-3} продолжает классическую двойственность Понтрягина с класса конечных абелевых групп на класс групп вида $Z\times K$, где $Z$ --- компактно порожденная абелева группа Ли,  а $K$ --- компактная группа Ли.
\etm

\subsection{Стереотипная двойственность в комплексной геометрии}

\paragraph{Желательный вид теории в комплексной геометрии.}    

Результат, аналогичный теореме \ref{TH:neprer-dvoistv}, в комплексной геометрии должен был бы выглядеть следующим образом. Нужно определить понятие голоморфной оболочки $\Env_{\mathcal O}$ так, чтобы когда в определении \eqref{chetyrehugolnik-C-C*-2} операции $\dagger$ заменяешь $\Env_{\mathcal C}$ на $\Env_{\mathcal O}$,
\beq\label{chetyrehugolnik-O-O*-2}
\dagger=\star\circ\Env_{\mathcal O},
\eeq
тождество \eqref{chetyrehugolnik-C-C*-3} превратится в тождество
\beq\label{chetyrehugolnik-O-O*-3}
\Big({\mathcal O}^\star(G)^\dagger\Big)^\dagger\cong {\mathcal O}^\star(G)
\eeq
справедливое для всех (необязательно коммутативных) конечных расширений связных комплексных групп Ли $G$. Алгебры Хопфа $H$, удовлетворяющие соответствующему равенству \eqref{DEF:neppr-reflex}
естественно тогда назвать {\it рефлексивными относительно голоморфной оболочки $\Env_{\mathcal E}$}, а от диаграммы  \eqref{diagramma-kategorij-dlya-nepreryvnoi-obolochki} естественно ожидать, что она в этом случае превратится в диаграмму
{\sf
 \beq\label{diagramma-kategorij-dlya-golomorf-obolochki-ideal}
 \xymatrix @R=3.pc @C=2.pc
 {
 \boxed{\begin{matrix}
  \text{алгебры Хопфа в $(\Ste,\circledast)$,}\\
  \text{рефлексивные относительно}\\
  \text{голоморфной оболочки}
 \end{matrix}}
 \ar[rr]^{H\mapsto H^\dagger} & &
 \boxed{\begin{matrix}
  \text{алгебры Хопфа в $(\Ste,\circledast)$,}\\
  \text{рефлексивные относительно}\\
  \text{голоморфной оболочки}
 \end{matrix}}
 \\
 \boxed{\begin{matrix}
  \text{компактно порожденные}\\
  \text{группы Ли с линейной}\\
  \text{связной компонентой единицы}
 \end{matrix}} \ar[u]^(0.45){\scriptsize\begin{matrix} \mathcal{O}^\star(G)\\
 \text{\rotatebox{90}{$\mapsto$}} \\ G\end{matrix}} & &
 \boxed{\begin{matrix}
  \text{компактно порожденные}\\
  \text{группы Ли с линейной}\\
  \text{связной компонентой единицы}
 \end{matrix}} \ar[u]_(0.45){\scriptsize\begin{matrix} \mathcal{O}^\star(G) \\
 \text{\rotatebox{90}{$\mapsto$}} \\ G\end{matrix}} \\
 \boxed{\begin{matrix}
  \text{компактно порожденные}\\
  \text{абелевы комплексные группы Ли}
 \end{matrix}} \ar[u]^{e} \ar[rr]^{G\mapsto \widehat{G}} & &
  \boxed{\begin{matrix}
  \text{компактно порожденные}\\
  \text{абелевы комплексные группы Ли}
  \end{matrix}}\ar[u]_{e}
 }
 \eeq }\noindent
в которой $\widehat{G}$ --- группа комплексных характеров $\chi:G\to\C^\times=\C\setminus\{0\}$ с поточечными алгебраическими операциями и топологией равномерной сходимости на компактах.
А теорема \ref{TH:neprer-dvoistv} --- в теорему 

\btm\label{TH:golom-dvoistv-ideal} 
Двойственность алгебр Хопфа \eqref{chetyrehugolnik-O-O*-3} продолжает классическую двойственность Понтрягина с класса компактно порожденных абелевых комплексных групп Ли на класс компактно порожденных комплексных групп Ли с линейной связной компонентой единицы.
\etm

\paragraph{Полученные на сегодняшний день результаты.}    
Для выбранного на с. \pageref{DEF:C^infty-obolochka} определения пока не получилось доказать теорему \ref{TH:golom-dvoistv-ideal}. Доказанный к настоящему времени результат иллюстрируется диаграммой
{\sf
 \beq\label{diagramma-kategorij-dlya-konechnyh-lineinyh-grupp}
 \xymatrix @R=3.pc @C=2.pc
 {
 \boxed{\begin{matrix}
  \text{алгебры Хопфа в $(\Ste,\circledast)$,}\\
  \text{рефлексивные относительно}\\
  \text{голоморфной оболочки}
 \end{matrix}}
 \ar[rr]^{H\mapsto H^\dagger} & &
 \boxed{\begin{matrix}
  \text{алгебры Хопфа в $(\Ste,\circledast)$,}\\
  \text{рефлексивные относительно}\\
  \text{голоморфной оболочки}
 \end{matrix}}
 \\
 \boxed{\begin{matrix}
  \text{конечные расширения}\\
  \text{связных линейных групп}
 \end{matrix}} \ar[u]^(0.45){\scriptsize\begin{matrix} \mathcal{O}^\star(G)\\
 \text{\rotatebox{90}{$\mapsto$}} \\ G\end{matrix}} & &
 \boxed{\begin{matrix}
  \text{конечные расширения}\\
  \text{связных линейных групп}
 \end{matrix}} \ar[u]_(0.45){\scriptsize\begin{matrix} \mathcal{O}^\star(G) \\
 \text{\rotatebox{90}{$\mapsto$}} \\ G\end{matrix}} \\
 \boxed{\begin{matrix}
 \text{абелевы конечные группы}
 \end{matrix}} \ar[u]^{e} \ar[rr]^{G\mapsto \widehat{G}} & &
  \boxed{\begin{matrix}
 \text{абелевы конечные группы}
  \end{matrix}}\ar[u]_{e}
 }
 \eeq }
и выглядит он как 

\btm\label{TH:golom-dvoistv} \cite{Ak08,ArHRC}
Двойственность алгебр Хопфа \eqref{chetyrehugolnik-O-O*-3} продолжает классическую двойственность Понтрягина для конечных абелевых групп на класс конечных расширений связных линейных комплексных групп.
\etm

\section{Перспективы}

\subsection{Расширение классов групп в теориях двойственности}

Классы групп, на которые распространяются двойственности в теоремах \ref{TH:neprer-dvoistv}, \ref{TH:gladk-dvoistv} и  \ref{TH:golom-dvoistv} (даже если включать обобщения \ref{TH:gladk-dvoistv-ideal} и \ref{TH:golom-dvoistv-ideal}), пока еще не окончательно определены. Есть некоторые ограничения на локальную структуру группы, типа теоремы Люмине---Валетта для комплексных групп Ли \cite{Luminet-Valette}, но ограничений в том, насколько широкой может быть группа пока найдено не было, и поэтому вопрос, каковы естественныые границы действия теорем \ref{TH:neprer-dvoistv}, \ref{TH:gladk-dvoistv} и  \ref{TH:golom-dvoistv}, остается открытым. В частности, остается открытым 

\begin{quest}\label{Q1}
Справедливы ли тождества \eqref{chetyrehugolnik-C-C*-3}, \eqref{chetyrehugolnik-E-E*-3} и \eqref{chetyrehugolnik-O-O*-3} для дискретных групп $G$?
\end{quest}

\subsection{Двойственность в интегральной геометрии}

Помимо описанных выше оболочек $\Env_{\mathcal C}$, $\Env_{\mathcal E}$, $\Env_{\mathcal O}$ никто не мешает рассматривать другие оболочки. В частности, можно рассмотреть оболочку, обозначим ее $\Env_{\mathcal W}$, в которой класс $\varPhi$ инструментов наблюдения состоит из всевозможных морфизмов в так называемые стереотипные алгебры фон Неймана.  

Они вот как описываются. Пусть $X$ -- произвольное гильбертово пространство. Для всякого оператора $\ph\in{\mathcal L}(X)$ символом $\ph^\bullet$ мы обозначаем сопряженный оператор, то есть удовлетворяющий условию
\beq\label{DEF:sopryazh-operator}
\langle \ph(x),y\rangle=\langle x,\ph^\bullet(y)\rangle,\qquad x,y\in X.
\eeq
Для всякого компакта $K\subseteq X$ рассмотрим полунорму на $X:X$
\beq\label{polunormy-v-W(X)-1}
\norm{\ph}_K=\max\left\{\sup_{x\in K}\norm{\ph(x)},\sup_{x\in K}\norm{\ph^\bullet(x)}\right\},\qquad \ph\in {\mathcal L}(X)
\eeq
Обозначим символом ${\mathcal W}(X)$ пространство (линейных непрерывных) операторов $\ph:X\to X$, наделенное топологией, порожденной полунормами \eqref{polunormy-v-W(X)-1}. Оказывается, что пространство ${\mathcal W}(X)$ образует стереотипную алгебру относительно обычной композиции операторов и с инволюцией $\bullet$. Ее замкнутые непосредственные подалгебры етественно считать {\it стереотипными алгебрами фон Неймана}.

Предварительный анализ показывает, что оболочка $\Env_{\mathcal W}$ относительно класса стереотипных алгебр фон Неймана обладает весьма интересными свойствами, позволяюшими надеяться, что в возникающей в результате геометрической дисциплине стереотипная двойственность окажется содержательным понятием. В частности, уже сейчас видно, что алгебры функций ${\mathcal C}(M)$, ${\mathcal E}(M)$, ${\mathcal O}(M)$, в этой геометрии заменяются на алгебру ${\mathcal I}(M)$ так называемых {\it интегралов} на (паракомпактном локально компактном) пространстве $M$, под которыми понимаются всевозможные линейные непрерывные функционалы
$$
\xi:{\mathcal C}(M)^*\to \C
$$ 
где под ${\mathcal C}(M)^*$ понимается пространство линейных непрерывных функционалов на ${\mathcal C}(M)$ с топологией равномерной сходимости на ограниченнных (а не на вполне ограниченных, как это обычно делается в стереотипной теории) множествах в ${\mathcal C}(M)$. Пространство ${\mathcal I}(M)$ --- не функциональное пространство, и факт, что оно является стереотипной алгеброй (с умножением, индуцированным из алгебры ${\mathcal C}(M)$, которая плотна  в ${\mathcal I}(M)$), глубоко неочевиден.

Структура алгебры ${\mathcal I}(M)$ оправдывает для возникающей геометрической теории название {\it интегральная геометрия}. 
А фундаментальный вопрос наших исследований в ней выглядит так: 

\begin{quest}\label{Q3}
На какие классы групп продолжается двойственность Понтрягина в интегральной геометрии?
\end{quest}

\paragraph{Соглашения и обозначения.}

Пусть $\varPhi$ -- класс морфизмов, а ${\tt L}$ -- класс объектов в категории ${\tt K}$. Мы говорим, что
\bit{
\item[---]\label{DEF:goes-from} {\it $\varPhi$ выходит из ${\tt L}$}, если для
всякого объекта $X\in{\tt L}$ найдется морфизм $\ph\in\varPhi$, выходящий из $X$:
    $$
    \forall X\in{\tt L}\qquad \exists \ph\in\varPhi\qquad \Dom\ph=X,
    $$
    в частном случае, если $\tt L$ состоит из одного объекта $X$, мы говорим, что $\varPhi$ {\it выходит из $X$},

\item[---]\label{DEF:goes-to} {\it $\varPhi$ приходит в ${\tt L}$}, если для всякого объекта $X\in{\tt L}$
найдется морфизм $\ph\in\varPhi$, приходящий в $X$:
    $$
    \forall X\in{\tt L}\qquad \exists \ph\in\varPhi\qquad \Ran\ph=X.
    $$
    в частном случае, если $\tt L$ состоит из одного объекта $X$, мы говорим, что $\varPhi$ {\it приходит в $X$}
    }\eit

Для всякого множества $I$ символ $2_I$ обозначает множество всех {\it конечных подмножеств} в $I$.\label{DEF:2_G}

В диаграммах длинная двойная стрелка $\Longrightarrow$ означает изоморфизм.

\mainmatter

\chapter{ОБОЛОЧКИ И ДЕТАЛИЗАЦИИ}\label{CH:obolochki}

\section{Оболочки и детализации}
\label{obolochka,otpechatok,uzl-razlozh}

\subsection{Оболочка}

\paragraph{Оболочка в классе морфизмов относительно класса морфизмов.}

Пусть нам даны:
 \bit{

\item[---] категория ${\tt K}$, называемая {\it объемлющей категорией},

\item[---] категория ${\tt T}$, называемая {\it притягивающей категорией},

\item[---] ковариантный функтор $F:{\tt T}\to{\tt K}$,

\item[---] два класса $\varOmega$ и $\varPhi$ морфизмов в ${\tt K}$,
принимающих значения в объектах класса $F({\tt T})$, причем $\Omega$
называется {\it классом реализующих морфизмов}, а $\varPhi$ -- {\it классом
пробных морфизмов}. }\eit

 \bit{
\item[$\bullet$] Для заданных объектов $X\in\Ob(\tt K)$ и $X'\in\Ob(\tt T)$
морфизм $\sigma:X\to F(X')$ называется {\it расширением объекта $X\in{\tt K}$
над категорией ${\tt T}$ в классе морфизмов $\varOmega$ относительно класса
морфизмов $\varPhi$}, если $\sigma\in\varOmega$, для любого объекта $B$ в $\tt T$ и для любого морфизма
$\ph:X\to F(B)$ из класса $\varPhi$ найдется единственный морфизм $\ph':X'\to B$ в
категории ${\tt T}$, для которого следующая диаграмма будет коммутативна:
 \beq\label{DEF:diagr-rasshirenie-T}
\begin{diagram}
\node[2]{X} \arrow{sw,t}{\varOmega\owns\sigma} \arrow{se,t}{\ph\in\varPhi}\\
\node{F(X')}  \arrow[2]{e,b,--}{F(\ph')} \node[2]{F(B)}
\end{diagram}
\eeq

\item[$\bullet$]\label{DEF:obolochka} Расширение $\rho:X\to F(E)$ объекта
$X\in{\tt K}$ над категорией ${\tt T}$ в классе морфизмов $\varOmega$
относительно класса морфизмов $\varPhi$ называется {\it оболочкой объекта
$X\in{\tt K}$ над категорией ${\tt T}$ в классе морфизмов $\varOmega$
относительно класса морфизмов $\varPhi$},
если для любого другого расширения  $\sigma:X\to F(X')$ (объекта $X\in{\tt K}$
над категорией ${\tt T}$ в классе морфизмов $\varOmega$  относительно класса
морфизмов $\varPhi$) найдется единственный морфизм $\upsilon:X'\to E$ в $\tt
T$, для которого будет коммутативна диаграмма
 \beq\label{DEF:diagr-obolochka-T}
\begin{diagram}
\node[2]{X} \arrow{sw,t}{\sigma} \arrow{se,t}{\rho}\\
\node{F(X')}  \arrow[2]{e,b,--}{F(\upsilon)} \node[2]{F(E)}
\end{diagram}
 \eeq
}\eit

В дальнейшем нас будет почти исключительно интересовать
случай, когда ${\tt T}={\tt K}$, а $F:{\tt K}\to{\tt K}$ -- тождественный
функтор. Определения для него полезно сформулировать отдельно.

 \bit{
\item[$\bullet$] Морфизм $\sigma:X\to X'$ в категории $\tt K$ называется
{\it расширением объекта $X\in\Ob({\tt K})$ в классе морфизмов $\varOmega$ относительно класса
морфизмов $\varPhi$}, если $\sigma\in\varOmega$, и для любого морфизма
$\ph:X\to B$ из класса $\varPhi$ найдется единственный морфизм $\ph':X'\to B$ в
категории ${\tt K}$, замыкающий диаграмму
 \beq\label{DEF:diagr-rasshirenie}
\begin{diagram}
\node[2]{X} \arrow{sw,t}{\varOmega\owns\sigma} \arrow{se,t}{\forall\ph\in\varPhi}\\
\node{X'}  \arrow[2]{e,b,--}{\exists!\ph'} \node[2]{B}
\end{diagram}
\eeq

\item[$\bullet$] Расширение $\rho:X\to E$ объекта
$X\in\Ob({\tt K})$ в классе морфизмов $\varOmega$
относительно класса морфизмов $\varPhi$ называется {\it оболочкой объекта
$X$ в классе морфизмов $\varOmega$
относительно класса морфизмов $\varPhi$},
если для любого другого расширения  $\sigma:X\to X'$ (объекта $X$
в классе морфизмов $\varOmega$  относительно класса
морфизмов $\varPhi$) найдется единственный морфизм $\upsilon:X'\to E$ в $\tt
K$, замыкающий диаграмму
 \beq\label{DEF:diagr-obolochka}
\begin{diagram}
\node[2]{X} \arrow{sw,t}{\forall\sigma} \arrow{se,t}{\rho}\\
\node{X'}  \arrow[2]{e,b,--}{\exists!\upsilon} \node[2]{E}
\end{diagram}
 \eeq
Для морфизма оболочки $\rho:X\to E$ мы будем использовать обозначение
    \beq\label{DEF:env_F^L(X)}
    \rho=\env_{\varPhi}^\varOmega X.
    \eeq
Кроме того, объект $E$ мы также будем называть {\it оболочкой} объекта $X$
(в классе морфизмов $\varOmega$ относительно класса
морфизмов $\varPhi$), и обозначать это мы будем записью
    \beq\label{DEF:Env_F^L(X)}
E=\Env_{\varPhi}^\varOmega X.
 \eeq
}\eit

\brem Понятно, что объект $\Env_{\varPhi}^\varOmega X$ (если он существует)
определен с точностью до изоморфизма. Вопрос о том, когда соответствие $X\mapsto\Env_{\varPhi}^\varOmega X$
можно определить как функтор, обсуждается ниже, начиная со с.\pageref{DIAGR:funktorialnost-env-e-E}. \erem

\brem Ясно, что если $\varOmega=\varnothing$, то ни расширений, ни оболочек в
$\varOmega$ не существует. Поэтому интерес в этой конструкции представляют
только те ситуации, когда класс $\varOmega$ непуст. Среди них мы будем выделять
следующие две: \bit{ \item[---] $\varOmega=\Epi({\tt K})$ (то есть $\varOmega$
совпадает с классом всех эпиморфизмов категории ${\tt K}$), и тогда мы будем
пользоваться обозначениями
    \begin{align}\label{env_(varPhi)^Epi}
    &     \env_{\varPhi}^{\Epi}X:=\env_{\varPhi}^{\Epi({\tt K})}X, &&     \Env_{\varPhi}^{\Epi}X:=\Env_{\varPhi}^{\Epi({\tt K})}X.
    \end{align}

\item[---] $\varOmega=\Mor({\tt K})$ (то есть $\varOmega$ совпадает с классом
всех вообще морфизмов категории ${\tt K}$), в этом случае удобно вообще
опускать упоминание о классе $\varOmega$ в формулировках и обозначениях,
поэтому мы будем говорить об {\it оболочке объекта $X\in{\tt K}$ в категории
${\tt K}$ относительно класса морфизмов $\varPhi$}, а обозначения упрощать так:
    \begin{align}\label{env_(varPhi)=env_(varPhi)^K}
    &     \env_{\varPhi}X:=\env_{\varPhi}^{\Mor({\tt K})}X, &&     \Env_{\varPhi}X:=\Env_{\varPhi}^{\Mor({\tt K})}X.
    \end{align}
 }\eit
\erem

\brem Другой вырожденный, но несмотря на это все же содержательный случай, --
когда $\varPhi=\varnothing$. Для фиксированного объекта $X$ здесь существенно
только, что $\varPhi$ не содержит ни одного морфизма, выходящего из $X$:
$$
\varPhi^X=\{\ph\in\varPhi:\ \Dom\ph=X\}=\varnothing.
$$
Тогда, очевидно, вообще любой морфизм $\sigma\in\varOmega $, выходящий из $X$,
$\sigma:X\to X'$, является расширением для $X$ (в классе морфизмов $\varOmega $
относительно класса морфизмов $\varnothing$). Если вдобавок $\varOmega=\Epi$,
то оболочкой будет терминальный объект в категории $\Epi^X$ (если он
существует). Это можно изобразить формулой
$$
\Env_\varnothing^\varOmega X=\max\Epi^X.
$$
С другой стороны, если ${\tt K}$ -- категория с нулем $0$, и $\varOmega$
содержит все морфизмы с концом в $0$, то  оболочкой всякого объекта
относительно пустого класса морфизмов является $0$:
$$
\Env_\varnothing^\varOmega X=0.
$$
\erem

\brem Еще одна крайняя ситуация -- когда $\varPhi=\Mor({\tt K})$. Для
фиксированного объекта $X$ здесь существенно только то, что $\varPhi$ содержит
единицу $X$:
$$
1_X\in\varPhi.
$$
Тогда для произвольного расширения $\sigma$ из диаграммы
$$
\xymatrix @R=2.pc @C=2.0pc 
{
X \ar[rr]^{\sigma} \ar[dr]_{1_X} & & X'\ar@{-->}[dl]\\
 & X &
}
$$
следует, что $\sigma$ должен быть коретракцией (причем такой, для которой
пунктирная стрелка единственна). Если $\varOmega\subseteq\Epi$, это
возможно только когда $\sigma$ -- изоморфизм. Как следствие, в этом случае
оболочка $X$ совпадает с $X$ (с точностью до изоморфизма):
$$
\varOmega\subseteq\Epi\quad\Longrightarrow\quad\Env_{\Mor({\tt K})}^{\varOmega} X=X.
$$
\erem

\medskip
\centerline{\bf Свойства оболочек:}

\bit{\it

\item[$1^\circ$.]\label{LM:suzhenie-verh-klassa-morfizmov}
Пусть $\varSigma\subseteq\varOmega$ и $\varPhi$ -- классы морфизмов. Тогда для всякого объекта $X$

\bit{

\item[(a)] любое расширение $\sigma:X\to X'$ в (более узком классе) $\varSigma$ относительно $\varPhi$
является расширением в (более широком классе) $\varOmega$ относительно $\varPhi$,

\item[(b)] если существуют оболочки $\env_\varPhi^\varSigma X$ и $\env_\varPhi^\varOmega X$, то существует единственный морфизм $\rho:\Env_\varPhi^\varSigma X\to\Env_\varPhi^\varOmega X$, замыкающий диаграмму
\beq\label{suzhenie-verh-klassa-morfizmov}
\begin{diagram}
\node[2]{X} \arrow{sw,t}{\env_\varPhi^\varSigma X} \arrow{se,t}{\env^\varOmega _\varPhi X}\\
\node{\Env_\varPhi^\varSigma X}\arrow[2]{e,b,--}{\rho}   \node[2]{\Env^\varOmega _\varPhi X}
\end{diagram}
\eeq

\item[(c)] если существует оболочка $\env_\varPhi^\varOmega X$ (в более широком классе), и она лежит в (более узком) классе $\varSigma$,
$$
\env_\varPhi^\varOmega X\in\varSigma,
$$
то она же будет и оболочкой $\env_\varPhi^\varSigma X$ (в более узком классе):
$$
\env_\varPhi^\varOmega X=\env_\varPhi^\varSigma X.
$$
}\eit

\item[$2^\circ$.]\label{LM:suzhenie-verh-klassa-morfizmov-2}
Пусть $\varSigma$, $\varOmega$, $\varPhi$ - классы морфизмов, и для объекта $X$ выполняется условие
 \bit{

\item[(a)] любое расширение $\sigma:X\to X'$ в $\varOmega$ относительно $\varPhi$
содержится в $\varSigma$.
 }\eit
Тогда
 \bit{

\item[(b)] оболочка $X$ относительно $\varPhi$ в классе $\varOmega$
существует тогда и только тогда, когда существует оболочка $X$
относительно $\varPhi$ в классе $\varOmega\cap\varSigma$; в этом случае эти оболочки совпадают:
$$
\env_{\varPhi}^{\varOmega}=\env_{\varPhi}^{\varOmega\cap\varSigma},
$$

\item[(c)] если $\varSigma\subseteq\varOmega$,
то оболочка $X$ относительно $\varPhi$ в (более узком классе) $\varSigma$
существует тогда и только тогда, когда существует оболочка $X$
относительно $\varPhi$ в (более широком классе) $\varOmega$, и эти оболочки совпадают:
$$
\env_{\varPhi}^{\varOmega}X=\env_{\varPhi}^{\varSigma}X.
$$
 }\eit

\item[$3^\circ$.]\label{LM:suzhenie-klassa-morfizmov} Пусть $\varPsi\subseteq\varPhi$, тогда для всякого объекта $X$ и любого класса морфизмов $\varOmega$

\bit{

\item[(a)] любое расширение $\sigma:X\to X'$ в $\varOmega$ относительно $\varPhi$ является расширением в $\varOmega$ относительно $\varPsi$,

\item[(b)] если существуют оболочки $\env_\varPsi^\varOmega X$ и
$\env_\varPhi^\varOmega X$, то существует единственный морфизм
$\alpha:\Env_\varPsi^\varOmega X\gets\Env_\varPhi^\varOmega X$, замыкающий
диаграмму \beq\label{suzhenie-klassa-morfizmov}
\begin{diagram}
\node[2]{X} \arrow{sw,t}{\env_\varPsi^\varOmega X} \arrow{se,t}{\env^\varOmega _\varPhi X}\\
\node{\Env_\varPsi^\varOmega X}   \node[2]{\Env^\varOmega _\varPhi X} \arrow[2]{w,b,--}{\alpha}
\end{diagram}
\eeq
}\eit

\item[$4^\circ$.]\label{TH:env_Psi=env_Phi}
Пусть $\varPhi\subseteq\Mor({\tt K})\circ\ \varPsi$
(то есть всякий морфизм $\ph\in\varPhi$ можно представить в виде $\ph=\chi\circ\psi$,
где $\psi\in\varPsi$), тогда для всякого объекта $X$ и любого класса морфизмов $\varOmega$
\bit{\it

\item[(a)] любое расширение $\sigma:X\to X'$ в $\varOmega$ относительно $\varPsi$, являющееся эпиморфизмом в ${\tt K}$, является расширением в $\varOmega$ относительно $\varPhi$,

\item[(b)] если существуют оболочки $\env_\varPsi^\varOmega X$ и
$\env_\varPhi^\varOmega X$, причем $\env_\varPsi^\varOmega X$ является
эпиморфизмом в ${\tt K}$, то существует единственный морфизм
$\beta:\Env_\varPsi^\varOmega X\to\Env_\varPhi^\varOmega X$, замыкающий
диаграмму
 \beq\label{env_Psi=env_Phi-0}
\begin{diagram}
\node[2]{X} \arrow{sw,t}{\env_\varPsi^\varOmega X} \arrow{se,t}{\env^\varOmega _\varPhi X}\\
\node{\Env_\varPsi^\varOmega X}\arrow[2]{e,b,--}{\beta}   \node[2]{\Env^\varOmega _\varPhi X}
\end{diagram}
\eeq
}\eit

\item[$5^\circ$.]\label{PROP:deistvie-epimorfizma-na-Env} Пусть классы морфизмов $\varOmega$, $\varPhi$ и эпиморфизм $\e:X\to Y$ в ${\tt K}$ удовлетворяют следующим условиям:
 \bit{
\item[(a)] существует оболочка $\env_{\varPhi\circ\e}^\varOmega X$ относительно
класса морфизмов $\varPhi\circ\e=\{\ph\circ\e;\ \ph\in\varPhi\}$,

\item[(b)] существует оболочка $\env_\varPhi^\varOmega Y$,

\item[(c)] композиция $\env_{\varPhi}^\varOmega Y\circ\e$ принадлежит классу $\varOmega$.

    }\eit
    Тогда существует единственный морфизм $\upsilon:\Env_{\varPhi\circ\e}^\varOmega X\gets \Env_{\varPhi}^\varOmega Y$, замыкающий диаграмму
\beq\label{deistvie-epimorfizma-na-Env}
\xymatrix @R=2.5pc @C=4.0pc
{
 X\ar@/_1ex/[rd]^{\ \env_{\varPhi}^\varOmega Y\circ\e}\ar[d]_{\env_{\varPhi\circ\e}^\varOmega X}\ar[r]^{\e} & Y\ar[d]^{\env_{\varPhi}^\varOmega Y} \\
 \Env_{\varPhi\circ\e}^\varOmega X &   \Env_{\varPhi}^\varOmega Y\ar@{-->}[l]^{\upsilon}
}
\eeq

}\eit

\bpr
1. Если морфизм $\sigma$ удовлетворяет условию \eqref{DEF:diagr-rasshirenie} с подставленным
в него классом $\varSigma$ вместо $\varOmega$, то $\sigma$ будет удовлетворять исходному
условию \eqref{DEF:diagr-rasshirenie}, потому что $\varSigma\subseteq\varOmega$. Это доказывает (a).
Отсюда  следует, в частности, что $\env_\varPsi^\varSigma X$ является расширением в $\varOmega$
относительно $\varPhi$, поэтому должна существовать единственная пунктирная стрелка в
\eqref{suzhenie-verh-klassa-morfizmov}. То есть справедливо (b). Наконец, если существует
оболочка $\env_\varPhi^\varOmega X$ (в более широком классе), и она лежит в (более узком)
классе $\varSigma$, $\env_\varPhi^\varOmega X\in\varSigma$, то из этого следует, что
$\env_\varPhi^\varOmega X$ является расширением в $\varSigma$. А с другой стороны,
любое другое расширение $\sigma:X\to X'$ в $\varSigma$, будет по только что доказанному
свойству (a) расширением в $\varOmega$, значит, существует единственный морфизм $\upsilon$ в оболочку в $\varOmega$:
$$
\begin{diagram}
\node[2]{X} \arrow{sw,t}{\sigma} \arrow{se,t}{\env_\varPhi^\varOmega X}\\
\node{X'}  \arrow[2]{e,b,--}{\upsilon} \node[2]{\Env_\varPhi^\varOmega X}
\end{diagram}
$$
Это доказывает, что $\env_\varPhi^\varOmega X$ -- оболочка в $\varSigma$. То есть справедливо (c).

2. Пусть выполнено $2^\circ$ (a). Если для какого-то объекта $X$ существует оболочка
$\env_{\varPhi}^{\varOmega}X$ в классе $\varOmega$ относительно $\varPhi$, то в сиду (a),
она должна быть расширением в более узком
классе $\varOmega\cap\varSigma$ относительно $\varPhi$. Применяя свойство $1^\circ$ (c)
на с.\pageref{LM:suzhenie-verh-klassa-morfizmov}, мы получаем, что $\env_{\varPhi}^{\varOmega}X$
является и оболочкой в более узком классе $\varOmega\cap\varSigma$, то есть
$\env_{\varPhi}^{\varOmega}X=\env_{\varPhi}^{\varOmega\cap\varSigma}X$.

Наоборот, пусть существует оболочка $\env_{\varPhi}^{\varOmega\cap\varSigma}X$
относительно $\varPhi$ в классе $\varOmega\cap\varSigma$. Тогда по свойству $1^\circ$ (a)
на с.\pageref{LM:suzhenie-verh-klassa-morfizmov}, она является расширением
относительно $\varPhi$ в более широком классе $\varOmega$. Рассмотрим какое-нибудь
другое расширение $\sigma:X\to X'$  относительно $\varPhi$ в классе $\varOmega$.
В силу (a), $\sigma$ является расширением относительно
$\varPhi$ в классе $\varOmega\cap\varSigma$. Значит, существует единственный морфизм
$\upsilon:X'\to \Env_{\varPhi}^{\varOmega\cap\varSigma} X$ в оболочку в классе
$\varOmega\cap\varSigma$, замыкающий диаграмму
$$
\begin{diagram}
\node[2]{X} \arrow{sw,t}{\sigma} \arrow{se,t}{\env_{\varPhi}^{\varOmega\cap\varSigma} X}\\
\node{X'}  \arrow[2]{e,b,--}{\upsilon} \node[2]{\Env_{\varPhi}^{\varOmega\cap\varSigma} X}
\end{diagram}
$$
Это доказывает, что $\env_{\varPhi}^{\varOmega\cap\varSigma}X$ является (не только
расширением, но и) оболочкой относительно $\varPhi$ в классе $\varOmega$. Это доказывает $2^\circ$ (b).
Утверждение $2^\circ$ (c) является его прямым следствием.

3. Пусть $\varPsi\subseteq\varPhi$. Тогда утверждение (a) очевидно:
любое расширение $\sigma:X\to X'$ относительно $\varPhi$ является расширением относительно
более узкого класса $\varPsi$. Для (b) получаем: поскольку $\env_\varPhi^\varOmega X$ есть
расширение относительно $\varPhi$, оно должно быть расширением относительно более узкого
класса $\varPsi$, поэтому существует единственный морфизм из $\Env_\varPhi^\varOmega X$ в
оболочку $\Env_\varPsi^\varOmega X$ относительно $\varPsi$, замыкающий \eqref{suzhenie-klassa-morfizmov}.

4. Пусть $\varPhi\subseteq\Mor({\tt K})\circ\ \varPsi$. Для (a) наши рассуждения будут иллюстрироваться диаграммой:
$$
\xymatrix @R=2.pc @C=5.0pc 
{
X \ar[rr]^{\sigma} \ar[dr]^{\psi}\ar@/_4ex/[ddr]_{\ph} & & X' \ar@{-->}[dl]_-{\psi'} \ar@{-->}@/^4ex/[ddl]^{\ph'}
\\
 & Y\ar[d]^{\chi} & \\
 & B &
}
$$
Если $\sigma:X\to X'$ -- какое-нибудь расширение $X$ в $\varOmega$ относительно $\varPsi$, то для произвольного морфизма $\ph\in\varPhi$, $\ph:X\to B$, мы выбираем разложение $\ph=\chi\circ\psi$, в котором $\psi\in\varPsi$. Поскольку $\sigma$ -- расширение $X$ в $\varOmega$ относительно $\varPsi$, найдется морфизм $\psi'$ такой, что $\psi=\psi'\circ\sigma$. После этого, положив $\ph'=\chi\circ\psi'$, мы получим морфизм, для которого
$$
\ph=\chi\circ\psi=\chi\circ\psi'\circ \sigma=\ph'\circ \sigma.
$$
Единственность $\ph'$ следует из эпиморфности $\sigma\in\varOmega$, и мы поэтому получаем, что $\sigma$ -- расширение $X$ в $\varOmega$ относительно $\varPhi$. После того, как (a) доказано, (b) становится его следствием: морфизм $\env^\varOmega _\varPsi X:X\to \Env^\varOmega _\varPsi X$ есть расширение $X$ в $\varOmega$ относительно $\varPsi$, значит, в силу (a), и относительно $\varPhi$ тоже. Поэтому должен существовать единственный морфизм $\beta$ из $\Env^\varOmega _\varPsi X$ в оболочку $\Env^\varOmega _\varPhi X$ относительно $\varPhi$, замыкающий \eqref{env_Psi=env_Phi-0}.

5. Для всякого морфизма $\ph:Y\to B$ из $\varPhi$ мы получаем диаграмму
$$
\xymatrix 
{
X\ar[rr]^{\env_{\varPhi}^\varOmega Y\circ\e}\ar[rd]_{\e}\ar@/_6ex/[rdd]_{\ph\circ\e} &  & \Env_{\varPhi}^\varOmega Y\ar@{-->}@/^6ex/[ldd]^{\ph'} \\
 & Y\ar[d]^{\ph}\ar[ru]_(.4){\env_{\varPhi}^\varOmega Y} & \\
 &  B  &
}
$$
Ее нужно понимать так. С одной стороны, поскольку $\env_{\varPhi}^\varOmega Y$
-- расширение относительно $\varPhi$, должен существовать морфизм $\ph'$,
замыкающий правый нижний треугольник, а значит, и периметр. С другой стороны,
если $\ph'$ -- какой-нибудь произвольный морфизм, замыкающий периметр, то есть
удовлетворяющий равенству
$$
\ph'\circ\env_{\varPhi}^\varOmega Y\circ\e=\ph\circ\e,
$$
то поскольку $\e$ -- эпиморфизм, на него можно сокращать, и поэтому
автоматически должно выполняться и равенство
$$
\ph'\circ\env_{\varPhi}^\varOmega Y=\ph,
$$
то есть должен быть коммутативен правый нижний треугольник. Отсюда следует, что
морфизм $\ph'$ должен быть единственным (поскольку, по определению оболочки, в
правом нижнем треугольнике пунктирная стрелка должна быть единственной).

Мы получаем, что у внешнего треугольника существует и притом только одна
пунктирная стрелка $\ph'$. Это верно для всякого $\ph\in\varPhi$, и вдобавок
$\env_{\varPhi}^\varOmega Y\circ\e\in\varOmega$. Поэтому можно сделать вывод,
что морфизм $\env_{\varPhi}^\varOmega Y\circ\e$ является расширением $X$ в
$\varOmega$ относительно класса морфизмов $\varPhi\circ\e$. Как следствие,
должен существовать и быть единственным морфизм $\upsilon$ из
$\Env_{\varPhi}^\varOmega Y$ в оболочку $\Env_{\varPhi\circ\e}^\varOmega X$
относительно $\varPhi\circ\e$, замыкающий диаграмму
\eqref{deistvie-epimorfizma-na-Env}. \epr

 \bit{
\item[$\bullet$] Условимся говорить, что в категории ${\tt K}$ {\it класс морфизмов $\varPhi$ порождается изнутри
 классом морфизмов $\varPsi$}, если
    \beq\label{DEF:morfizmy-porozhdayutsya-iznutri}
    \varPsi\subseteq\varPhi\subseteq\Mor({\tt K})\circ\ \varPsi.
    \eeq
  }\eit

\btm\label{TH:morfizmy-porozhdayutsya-iznutri} Пусть в категории ${\tt K}$ класс морфизмов $\varPhi$ порождается изнутри классом морфизмов $\varPsi$. Тогда для любого класса эпиморфизмов $\varOmega$ (необязательно, чтобы $\varOmega$ включал все эпиморфизмы категории ${\tt K}$) и всякого объекта $X$ существование оболочки $\env^\varOmega_\varPsi X$ эквивалентно существованию оболочки $\env^\varOmega_\varPhi X$, и эти оболочки совпадают:
 \beq\label{env_Psi=env_Phi}
\env^\varOmega _\varPsi X=\env^\varOmega _\varPhi X.
 \eeq
 \etm
\bpr 1. Пусть сначала существует оболочка $\env^\varOmega_\varPsi X$. Поскольку
она является расширением относительно $\varPsi$, и одновременно эпиморфизмом,
по свойству $4^\circ$ (a), она должна быть также расширением и относительно
$\varPhi$. Если же $\sigma:X\to X'$ -- какое-то другое расширение относительно
$\varPhi$, то по свойству $3^\circ$ (a), оно должно быть расширением
относительно $\varPsi$, и поэтому существует единственный морфизм
$\upsilon:\Env_\varPsi^\varOmega X\gets X'$, замыкающий диаграмму
$$
\begin{diagram}
\node[2]{X} \arrow{sw,t}{\env_\varPsi^\varOmega X} \arrow{se,t}{\sigma}\\
\node{\Env_\varPsi^\varOmega X}   \node[2]{X'}
\arrow[2]{w,b,--}{\exists!\upsilon}
\end{diagram}
$$
Это значит, что $\env^\varOmega _\varPsi X$ является оболочкой относительно
$\varPhi$, и справедливо \eqref{env_Psi=env_Phi}.

2. Наоборот, пусть существует оболочка $\env^\varOmega _\varPhi X$. Поскольку
она является расширением относительно $\varPhi$, по свойству $3^\circ$ (a), она
должна быть также расширением и относительно $\varPsi$. Если же $\sigma:X\to
X'$ -- какое-то другое расширение в $\varOmega$ относительно $\varPsi$, то,
поскольку $\sigma\in\Epi$, по свойству $4^\circ$ (a), оно должно быть
расширением относительно $\varPhi$, и поэтому существует морфизм
$\upsilon:X'\to \Env_\varPhi^\varOmega X$, замыкающий диаграмму
$$
\begin{diagram}
\node[2]{X} \arrow{sw,t}{\sigma} \arrow{se,t}{\env_\varPhi^\varOmega X}\\
\node{X'} \arrow[2]{e,b,--}{\exists!\upsilon} \node[2]{\Env_\varPhi^\varOmega X}
\end{diagram}
$$
Это значит, что $\env^\varOmega _\varPhi X$ является оболочкой относительно $\varPsi$, и опять справедливо \eqref{env_Psi=env_Phi}.
\epr

\bit{
\item[$\bullet$]\label{DEF:varPhi-razlich-morfizmy-snaruzhi} Условимся говорить, что класс морфизмов $\varPhi$ в категории ${\tt K}$
 {\it различает морфизмы снаружи}, если для любых двух различных параллельных морфизмов $\alpha\ne\beta:X\to Y$ найдется морфизм $\ph:Y\to M$ из класса $\varPhi$ такой, что $\ph\circ\alpha\ne\ph\circ\beta$.
}\eit

\btm\label{TH:Phi-razdel-moprfizmy}
Если класс морфизмов $\varPhi$ различает морфизмы снаружи, то для любого класса морфизмов $\varOmega$
 \bit{
\item[(i)] всякое расширение в $\varOmega$ относительно $\varPhi$ является мономорфизмом,

\item[(ii)] оболочка относительно $\varPhi$ в классе $\varOmega$ существует тогда и только тогда, когда существует оболочка относительно $\varPhi$ в классе $\varOmega\cap\Mono$; в этом случае эти оболочки совпадают:
$$
\env_{\varPhi}^{\varOmega}=\env_{\varPhi}^{\varOmega\cap\Mono},
$$

\item[(iii)] если класс $\varOmega$ содержит все мономорфизмы,
 $$
 \varOmega\supseteq\Mono,
 $$
 то существование оболочки относительно $\varPhi$ в $\Mono$ автоматически влечет за собой существование оболочки относительно $\varPhi$ в $\varOmega$ и совпадение этих оболочек:
$$
\env_{\varPhi}^{\varOmega}=\env_{\varPhi}^{\Mono}.
$$

 }\eit
\etm
\bpr
В силу свойства $2^\circ$ на с.\pageref{LM:suzhenie-verh-klassa-morfizmov-2}, утверждения (ii) и (iii) следуют из (i),
поэтому достаточно доказать (i). Предположим, что какое-то расширение $\sigma:X\to X'$ не является мономорфизмом, то есть что существуют два
различных параллельных морфизма $\alpha\ne\beta:T\to X$ такие, что
\beq\label{sigma-circ-alpha=sigma-circ-beta}
\sigma\circ\alpha=\sigma\circ\beta.
\eeq
Поскольку класс $\varPhi$ различает морфизмы снаружи, должен существовать морфизм $\ph:X\to M$, $\ph\in\varPhi$, такой, что
\beq\label{ph-circ-alpha-ne-ph-circ-beta}
\ph\circ\alpha\ne\ph\circ\beta.
\eeq
Поскольку $\sigma:X\to X'$ -- расширение относительно $\varPhi$, найдется продолжение $\ph':X'\to M$ морфизма $\ph:X\to M$: $\ph=\ph'\circ\sigma$. Теперь мы получаем:
$$
\ph\circ\alpha=\ph'\circ\sigma\circ\alpha=\eqref{sigma-circ-alpha=sigma-circ-beta}=\ph'\circ\sigma\circ\beta=\ph\circ\beta,
$$
и это противоречит \eqref{ph-circ-alpha-ne-ph-circ-beta}.
\epr

\bit{
\item[$\bullet$] Напомним, что класс морфизмов $\varPhi$ в категории ${\tt K}$ называется {\it правым идеалом}, если
$$
\varPhi\circ\Mor({\tt K})\subseteq\varPhi
$$
(то есть для любого $\ph\in\varPhi$ и любого морфизма $\mu$ категории ${\tt K}$ композиция $\ph\circ\mu$ принадлежит $\varPhi$).
}\eit

\btm\label{TH:Phi-razdel-moprfizmy-*} Если класс морфизмов $\varPhi$ различает
морфизмы снаружи и является правым идеалом в категории ${\tt K}$, то для любого класса морфизмов $\varOmega$
 \bit{
\item[(i)] всякое расширение в $\varOmega$ относительно $\varPhi$ является биморфизмом,

\item[(ii)] оболочка относительно $\varPhi$ в классе $\varOmega$ существует тогда и только тогда, когда существует оболочка относительно $\varPhi$ в классе $\varOmega\cap\Bim$; в этом случае эти оболочки совпадают:
$$
\env_{\varPhi}^{\varOmega}=\env_{\varPhi}^{\varOmega\cap\Bim}.
$$

\item[(iii)] если класс $\varOmega$ содержит все биморфизмы,
 $$
 \varOmega\supseteq\Bim,
 $$
 то существование оболочки относительно $\varPhi$ в $\Bim$ автоматически влечет за собой существование оболочки относительно $\varPhi$ в $\varOmega$ и совпадение этих оболочек:
$$
\env_{\varPhi}^{\varOmega}=\env_{\varPhi}^{\Bim}.
$$

 }\eit
\etm
\bpr
В силу свойства $2^\circ$ на с.\pageref{LM:suzhenie-verh-klassa-morfizmov-2}, утверждения (ii) и (iii) следуют из (i),
поэтому достаточно доказать (i). Пусть $\sigma:X\to X'$ -- расширение в $\varOmega$ относительно $\varPhi$.
По теореме \ref{TH:Phi-razdel-moprfizmy}(i), $\sigma$ должен быть мономорфизмом. Предположим, что он
не является эпиморфизмом. Это значит, что найдутся
два различных параллельных морфизма $\alpha\ne\beta:X'\to T$ такие, что
\beq\label{alpha-circ-sigma=beta-circ-sigma}
\alpha\circ\sigma=\beta\circ\sigma.
\eeq
Поскольку $\varPhi$ различает морфизмы снаружи, должен существовать морфизм $\ph:T\to M$, $\ph\in\varPhi$, такой, что
$$
\ph\circ\alpha\ne\ph\circ\beta.
$$
При этом, в силу \eqref{alpha-circ-sigma=beta-circ-sigma},
$$
\ph\circ\alpha\circ\sigma=\ph\circ\beta\circ\sigma.
$$
Если теперь считать, что $\varPhi$ является правым идеалом в категории ${\tt K}$,
то морфизм $\ph\circ\alpha\circ\sigma=\ph\circ\beta\circ\sigma$ будет лежать в $\varPhi$.
Тогда мы можем интерпретировать эту картину так: пробный (то есть лежащий в $\varPhi$) морфизм
$\ph\circ\alpha\circ\sigma=\ph\circ\beta\circ\sigma:X\to M$ имеет два разных продолжения
$\ph\circ\alpha\ne\ph\circ\beta:X'\to M$ вдоль $\sigma:X\to X'$. Значит $\sigma$ не может
быть расширением относительно $\varPhi$.
\epr

Важная ситуация --- когда классы $\varPhi$ и $\varOmega$ совпадают:
$$
\varPhi=\varOmega.
$$
В этом случае понятия расширения и оболочки совпадают:

\btm\label{varPhi=varOmega}
Всякое расширение объекта $X$ в классе $\varOmega$ относительно самого этого класса $\varOmega$ является оболочкой объекта $X$ в классе $\varOmega$  относительно $\varOmega$.
\etm
\bpr
Пусть $\rho:X\to E$ --- произвольное расширение объекта $X$ в классе морфизмов $\varOmega$
относительно того же самого класса морфизмов $\varOmega$. Если $\sigma:X\to X'$ --- какое-то другое расширение объекта $X$ в классе морфизмов $\varOmega$  относительно класса
морфизмов $\varOmega$, то в диаграмме \eqref{DEF:diagr-obolochka}
$$
\begin{diagram}
\node[2]{X} \arrow{sw,t}{\forall\sigma} \arrow{se,t}{\rho}\\
\node{X'}  \arrow[2]{e,b,--}{\exists!\upsilon} \node[2]{E}
\end{diagram}
$$
морфизм $\upsilon:X'\to E$ существует и единственен, потому что $\rho\in\varOmega$. Это значит, что $\rho:X\to E$ --- оболочка.
\epr

\paragraph{Оболочка в классе объектов относительно класса объектов.}

Частным случаем конструкции оболочки является ситуация, когда $\varOmega$ и/или $\varPhi$ представляют собой классы всех морфизмов в объекты из заданных подклассов класса $\Ob({\tt K})$. Точная формулировка для случая, когда оба класса $\varOmega$ и $\varPhi$ определяются таким образом, выглядит так: пусть дана категория ${\tt K}$ и два подкласса ${\tt L}$ и ${\tt M}$ в классе $\Ob({\tt K})$ объектов ${\tt K}$.

 \bit{
\item[$\bullet$] Морфизм $\sigma:X\to X'$ называется {\it расширением объекта
$X\in{\tt K}$ в классе ${\tt L}$ относительно класса ${\tt M}$}, если $X'\in{\tt L}$ и для любого объекта
$B\in{\tt M}$ и любого морфизма $\ph:X\to B$ найдется единственный морфизм $\ph':X'\to B$, замыкающий диаграмму:
$$
\put(6,-35){$\begin{matrix}\text{\rotatebox{90}{$\owns$}}\\ {\tt L}
\end{matrix}$} \put(73,-35){$\begin{matrix}\text{\rotatebox{90}{$\owns$}}\\
{\tt M} \end{matrix}$}
\begin{diagram}
\node[2]{X} \arrow{sw,t}{\sigma} \arrow{se,t}{\forall\ph}\\
\node{X'}  \arrow[2]{e,b,--}{\exists!\ph'} \node[2]{B}
\end{diagram}
$$

\item[$\bullet$]\label{DEF:obolochka-otn-klassa} Расширение $\rho:X\to E$
объекта $X\in{\tt K}$ в классе ${\tt L}$  относительно класса ${\tt M}$
называется {\it оболочкой объекта $X\in{\tt K}$  в классе ${\tt L}$
относительно класса ${\tt M}$}, и мы будем обозначать это записью
    \beq\label{DEF:env_M^L(A)}
    \rho=\env_{\tt M}^{\tt L} X,
    \eeq
если для любого другого расширения  $\sigma:X\to X'$ (объекта $X$  в
классе ${\tt L}$  относительно класса ${\tt M}$) найдется единственный морфизм
$\upsilon:X'\to E$, замыкающий диаграмму:
 \beq\label{diagr:obolochka}
\put(6,-35){$\begin{matrix}\text{\rotatebox{90}{$\owns$}}\\ {\tt L}
\end{matrix}$} \put(73,-35){$\begin{matrix}\text{\rotatebox{90}{$\owns$}}\\
{\tt L} \end{matrix}$}
\begin{diagram}
\node[2]{X} \arrow{sw,t}{\forall\sigma} \arrow{se,t}{\rho}\\
\node{X'}  \arrow[2]{e,b,--}{\exists!\upsilon} \node[2]{E}
\end{diagram}
 \eeq
Объект $E$ тоже называется {\it оболочкой} объекта $X$ (в классе объектов
${\tt L}$ относительно класса объектов ${\tt M}$), и для него мы будем
использовать обозначение
 \beq\label{DEF:Env_M^L(A)}
E=\Env_{\tt M}^{\tt L} X.
 \eeq

}\eit

Возможны также следующие две крайние ситуации в выборе подкласса ${\tt L}$:
 \bit{

\item[---] если ${\tt L} =\Ob({\tt K})$, то мы будем говорить об {\it оболочке объекта $X\in{\tt K}$ в категории ${\tt K}$ относительно класса объектов ${\tt M}$}, и обозначения будем упрощать так:
    \begin{align}\label{env_M=env_M^K}
    &     \env_{\tt M}X:=\env_{\tt M}^{\tt K}X, &&     \Env_{\tt M}X:=\Env_{\tt M}^{\tt K}X.
    \end{align}

\item[---]\label{L=M} если ${\tt L} ={\tt M}$, то (по теореме \ref{varPhi=varOmega}) понятия расширения и оболочки совпадают: {\it всякое расширение объекта $X$ в классе ${\tt L}$  относительно самого этого класса ${\tt L}$ является оболочкой объекта $X\in{\tt K}$ в классе ${\tt L}$  относительно ${\tt L}$}; чтобы не повторяться, при ${\tt L} ={\tt M}$ говорят об {\it оболочке объекта $X$ в классе ${\tt L}$}, при этом обозначения упрощаются так:
    \begin{align}
    & \env^{\tt L} X=:\env_{\tt L}^{\tt L} X, &&  \Env^{\tt L} X=:\Env_{\tt L}^{\tt L} X.
    \end{align}

}\eit

\bit{
\item[$\bullet$] Условимся говорить, что класс объектов ${\tt M}$ в категории ${\tt K}$
 {\it различает морфизмы снаружи}, если класс морфизмов со значениями в объектах из ${\tt M}$ обладает этом свойством (в смысле определения на с.\pageref{DEF:varPhi-razlich-morfizmy-snaruzhi}), то есть
 для любых двух различных параллельных морфизмов $\alpha\ne\beta:X\to Y$ найдется морфизм $\ph:Y\to M\in{\tt M}$ такой, что $\ph\circ\alpha\ne\ph\circ\beta$.
}\eit

Из теоремы \ref{TH:Phi-razdel-moprfizmy-*} сразу следует

\btm\label{TH:M-razdel-moprfizmy} Если класс объектов ${\tt M}$ различает морфизмы снаружи, то для любого класса объектов ${\tt L}$
 \bit{
\item[(i)]  всякое расширение в ${\tt L}$ относительно ${\tt M}$ является биморфизмом,

\item[(ii)] оболочка в ${\tt L}$ относительно ${\tt M}$ существует тогда и только тогда, когда существует оболочка в классе биморфизмов со значениями в ${\tt L}$ относительно ${\tt M}$; в таком случае эти оболочки совпадают:
$$
\env_{\tt M}^{\tt L}=\env_{\tt M}^{\Bim({\tt K},{\tt L})}.
$$
 }\eit
\etm

\noindent\rule{160mm}{0.1pt}\begin{multicols}{2}

\paragraph{Примеры оболочек.}

\bex {\bf Универсальная обертывающая алгебра.} Пусть ${\tt K}={\tt LieAlg}$ --
категория алгебр Ли (например, над полем $\C$), ${\tt T}={\tt Alg}$ --
категория ассоциативных алгебр (над тем же полем $\C$) с единицей и $F:{\tt
Alg}\to{\tt LieAlg}$ -- функтор, который каждую ассоциативную алгебру $A$ с
единицей представляет как алгебру Ли со скобкой Ли
$$
[x,y]=x\cdot y-y\cdot x.
$$
Тогда оболочка алгебры Ли $\frak{g}$ над категорией ${\tt Alg}$ в классе
$\Mor({\tt LieAlg},F({\tt Alg}))$ всех морфизмов из ${\tt LieAlg}$ в $F({\tt
Alg})$ относительно самого этого класса представляет собой универсальную
обертывающую алгебру $U(\frak{g})$ алгебры Ли $\frak{g}$
(см.\cite{Bourbaki-LieAlg}):
$$
\Env_{\Mor({\tt LieAlg},F({\tt Alg}))}^{\Mor({\tt LieAlg},F({\tt Alg}))}
\frak{g}=U(\frak{g}).
$$
 \eex

\bex {\bf Компактификация Стоуна-Чеха.} В категории ${\tt Tikh}$ тихоновских
топологических пространств компактификация Стоуна-Чеха $\beta:X\to\beta X$
является оболочкой пространства $X$ в классе ${\tt Com}$ компактов относительно
самого этого класса:
$$
\beta X=\Env^{\tt Com}X.
$$
\eex \bpr Здесь используется известное утверждение \cite[Теорема 3.6.1]{Engelking},
согласно которому всякое непрерывное отображение $f:X\to K$ в произвольный
компакт $K$ продолжается до непрерывного отображения $F:\beta X\to K$.
Поскольку $\beta(X)$ плотно в $\beta X$, такое продолжение $F$ единственно, и
это означает, что $\beta:X\to\beta X$ -- расширение в классе ${\tt Com}$
относительно ${\tt Com}$. В силу замечания на с.\pageref{L=M}, в случае ${\tt
L}={\tt M}$ всякое расширение автоматически является оболочкой, поэтому $\beta$
-- оболочка. \epr

\bex {\bf Пополнение} $X^\blacktriangledown$ локально выпуклого пространства
$X$ является оболочкой $X$ в категории ${\tt LCS}$ локально выпуклых
пространств относительно класса ${\tt Ban}$ пространств Банаха:
$$
X^\blacktriangledown=\Env_{\tt Ban}^{\tt LCS}X.
$$
\eex \bpr Мы будем обозначать вложение $X$ в свое пополнение символом
$\blacktriangledown_X:X\to X^\blacktriangledown$ (обозначения работы
\cite{Ak03}).

Прежде всего, всякое линейное непрерывное отображение $f:X\to B$ в произвольное
банахово пространство $B$ однозначно продолжается до линейного непрерывного
отображения $F:X^\blacktriangledown\to B$ на пополнении $X^\blacktriangledown$
пространства $X$ (здесь можно, например, сослаться на общую теорему для
равномерных пространств \cite[Теорема 8.3.10]{Engelking}). Поэтому пополнение
$\blacktriangledown_X:X\to X^\blacktriangledown$ является расширением
пространства $X$ в категории ${\tt LCS}$ локально выпуклых пространств
относительно класса ${\tt Ban}$ банаховых пространств.

Заметим далее, что класс ${\tt Ban}$ банаховых пространств различает морфизмы
снаружи в категории ${\tt LCS}$. Поэтому по теореме
\ref{TH:M-razdel-moprfizmy} любое вообще расширение $\sigma:X\to X'$
относительно ${\tt Ban}$ должно быть биморфизмом в категории ${\tt LCS}$, то
есть инъективным отображением, образ которого $\sigma(X)$ плотен в $X'$.
Покажем, что оно вдобавок является открытым отображением: для любой окрестности
нуля $U\subseteq X$ найдется окрестность нуля $V\subseteq X'$ такая, что
\beq\label{sigma(U)-supseteq-V-cap-sigma(X)}
\sigma(U)\supseteq V\cap\sigma(X)
\eeq
Здесь можно считать, что $U$ -- замкнутая выпуклая уравновешенная окрестность
нуля в $X$. Тогда множество $\Ker U=\bigcap_{\varepsilon>0}\varepsilon\cdot U$
будет замкнутым подпространством в $X$. Рассмотрим фактор-пространство $X/\Ker
U$, и наделим его топологией нормированного пространства с единичным шаром
$U+\Ker U$. Тогда пополнение $(X/\Ker U)^\blacktriangledown$ будет банаховым
пространством, которое мы обозначим $A/U$. Естественное отображение (композицию
фактор-отображения $X\to X/\Ker U$ и пополнения $X/\Ker U\to (X/\Ker
U)^\blacktriangledown$) мы будем обозначать $\pi_U:X\to X/U$. Поскольку
$\sigma:X\to X'$ -- расширение относительно ${\tt Ban}$, отображение
$\pi_U:X\to X/U$ должно продолжаться до некоторого линейного непрерывного
отображения $(\pi_U)':X'\to X/U$.
$$
\begin{diagram}
\node{X} \arrow[2]{e,t}{\sigma} \arrow{se,b}{\pi_U} \node[2]{X'}
\arrow{sw,b,--}{(\pi_U)'}
\\
\node[2]{X/U}
\end{diagram}
$$
Если обозначить теперь через $W$ единичный шар в $X/U$, то есть замыкание
множества $U+\Ker U$ в пространстве $(X/\Ker U)^\blacktriangledown=X/U$, то для
окрестности нуля $V=((\pi_U)')^{-1}(W)$ мы получим цепочку, доказывающую \eqref{sigma(U)-supseteq-V-cap-sigma(X)}:
\begin{multline*}
y\in V\cap\sigma(X)\quad\Rightarrow\\ \Rightarrow\quad  \exists x\in X:\ y=\sigma(x)\ \& \ y\in V \quad \Rightarrow \\ \Rightarrow\quad \exists x\in X:\ y=\sigma(x)\ \&\ (\pi_U)'(y)=\\= (\pi_U)'(\sigma(x))=\underbrace{\pi_U(x)\in W}_{\scriptsize\begin{matrix}\Updownarrow \\ x\in U\end{matrix}}
\quad\Rightarrow\\ \Rightarrow\quad \exists x\in U:\ y=\sigma(x)\quad\Rightarrow\quad y\in\sigma(U).
 \end{multline*}
Итак, $\sigma:X\to X'$ -- открытое и инъективное отображение, причем
$\sigma(X)$ плотно в $X'$. Это значит, что $X'$ можно представлять себе, как
подпространство в пополнении $X^\blacktriangledown$ пространства $X$ с
топологией, индуцированной из $X^\blacktriangledown$. То есть определено
единственное линейное непрерывное отображение $\upsilon:X'\to
X^\blacktriangledown$, замыкающее диаграмму
$$
\begin{diagram}
\node[2]{X} \arrow{sw,t}{\sigma} \arrow{se,t}{\blacktriangledown_X}\\
\node{X'}  \arrow[2]{e,b,--}{\upsilon} \node[2]{X^\blacktriangledown}
\end{diagram}
$$
Таким образом, $\blacktriangledown_X:X\to X^\blacktriangledown$ -- оболочка
$X$ в ${\tt LCS}$ относительно ${\tt Ban}$. \epr

\bex {\bf Оболочка Аренса---Майкла.}\label{EX:Arens-Michael-v-TopAlg}
В категории ${\TopAlg}$ топологических алгебр (с раздельно непрерывным умножением) стандартное понятие оболочки Аренса-Майкла $A^{\sf AM}$ (см. \cite{Helemskii-Polynorm-alg}) совпадает с оболочкой в описанном нами смысле относительно подкатегории ${\BanAlg}$ банаховых алгебр:
\beq\label{Arens-Michael-v-TopAlg}
A^{\sf AM}=\Env_{\BanAlg}^{\TopAlg} A.
\eeq
\eex

\end{multicols}\noindent\rule[10pt]{160mm}{0.1pt}

\subsection{Детализация}

\paragraph{Детализация в классе морфизмов посредством класса морфизмов.}

Пусть нам даны:
 \bit{

\item[---] категория ${\tt K}$, называемая {\it объемлющей категорией},

\item[---] категория ${\tt T}$, называемая {\it отталкивающей категорией},

\item[---] ковариантный функтор $F:{\tt T}\to{\tt K}$,

\item[---] два класса $\varGamma$ и $\varPhi$ морфизмов в ${\tt K}$,
выходящих из объектов класса $F({\tt T})$, причем $\Omega$
называется {\it классом реализующих морфизмов}, а $\varPhi$ -- {\it классом
пробных морфизмов}. }\eit

 \bit{
\item[$\bullet$] Для заданных объектов $X\in\Ob(\tt K)$ и $X'\in\Ob(\tt T)$
морфизм $\sigma:F(X')\to X$ называется {\it обогащением объекта $X\in{\tt K}$ в классе морфизмов
$\varGamma$ над категорией ${\tt T}$ посредством класса
морфизмов $\varPhi$}, если $\sigma\in\varGamma$, и для любого объекта $B$ в $\tt T$ и для любого морфизма
$\ph:F(B)\to X$, $\ph\in\varPhi$, найдется единственный морфизм
$\ph':B\to X'$ в категории $\tt T$, для которого следующая диаграмма будет коммутативна:
 \beq\label{DEF:suzhenie-T}
\begin{diagram}
\node[2]{X}  \\
\node{F(B)}\arrow{ne,t}{\varPhi\owns\ph}\arrow[2]{e,b,--}{F(\ph')}
\node[2]{F(X')}\arrow{nw,t}{\sigma\in\varGamma}
\end{diagram}
\eeq

\item[$\bullet$]\label{DEF:nachinka} Обогащение $\rho:F(E)\to X$ класса морфизмов $\varGamma$ в объекте $X\in{\tt K}$ над категорией ${\tt T}$ посредством класса морфизмов $\varPhi$ называется {\it детализацией над категорией ${\tt T}$ объекта $X\in{\tt K}$ в классе морфизмов $\varGamma$ посредством класса морфизмов $\varPhi$}, если для любого другого обогащения $\sigma:F(X')\to X$ (объекта $X\in{\tt K}$ в классе морфизмов $\varGamma$ над категорией ${\tt T}$ посредством класса морфизмов $\varPhi$) найдется единственный морфизм $\upsilon:E\to X'$ в $\tt T$, для которого будет коммутативна диаграмма
 \beq\label{DIAGR:otpechatok-T}
\begin{diagram}
\node[2]{X}  \\
\node{F(E)}\arrow{ne,t}{\rho}\arrow[2]{e,b,--}{F(\upsilon)}
\node[2]{F(X')}\arrow{nw,t}{\sigma}
\end{diagram}
 \eeq
}\eit

В дальнейшем нас будет почти исключительно интересовать
случай, когда ${\tt T}={\tt K}$, а $F:{\tt K}\to{\tt K}$ -- тождественный
функтор. Как и в случае с оболочками, мы сформулируем определения для этой ситуации отдельно.
 \bit{
\item[$\bullet$] Морфизм $\sigma:X'\to X$ в категории $\tt K$ называется {\it обогащением объекта $X\in\Ob({\tt K})$ в классе морфизмов $\varGamma$ посредством класса морфизмов $\varPhi$}, если $\sigma\in\varGamma$, и для любого морфизма $\ph:B\to X$, $\ph\in\varPhi$, найдется единственный морфизм $\ph':B\to X'$ в категории $\tt K$, замыкающий диаграмму
 \beq\label{DEF:suzhenie}
\begin{diagram}
\node[2]{X}  \\
\node{B}\arrow{ne,t}{\forall\ph\in\varPhi}\arrow[2]{e,b,--}{\exists!\ph'}
\node[2]{X'}\arrow{nw,t}{\sigma\in\varGamma}
\end{diagram}
 \eeq

\item[$\bullet$] Обогащение $\rho:E\to X$ класса морфизмов $\varGamma$ в объекте $X\in\Ob({\tt K})$
посредством класса морфизмов $\varPhi$ называется {\it детализацией объекта $X$ в классе морфизмов $\varGamma$ посредством класса морфизмов $\varPhi$}, если для любого другого обогащения $\sigma:X'\to X$ (объекта $X$ в классе морфизмов $\varGamma$ посредством класса морфизмов $\varPhi$) найдется единственный морфизм $\upsilon:E\to X'$ в $\tt K$, для которого будет коммутативна диаграмма
 \beq\label{DIAGR:otpechatok}
\begin{diagram}
\node[2]{X}  \\
\node{E}\arrow{ne,t}{\rho}\arrow[2]{e,b,--}{\exists!\upsilon}
\node[2]{X'}\arrow{nw,t}{\forall\sigma}
\end{diagram}
 \eeq
Для морфизма детализации $\rho:E\to X$ мы будем использовать обозначение
    \beq\label{DEF:tr_F^L(A)}
    \rho=\rf_{\varPhi}^\varGamma  X.
    \eeq
Кроме того, сам объект $E$ здесь мы также называем {\it детализацией} объекта $X$ в классе морфизмов $\varGamma$ посредством класса морфизмов $\varPhi$, и обозначаем это записью
    \beq\label{DEF:Tr_F^L(A)}
E=\Rf_{\varPhi}^\varGamma  X.
 \eeq
}\eit

\brem Как и в случае с оболочкой, детализация $\Rf_{\varPhi}^\varGamma  X$ (если она существует) определена с точностью до изоморфизма. Вопрос о том, когда соответствие $X\mapsto\Rf_{\varPhi}^\varGamma X$ можно определить как функтор,
обсуждается ниже, начиная со с.\pageref{DIAGR:funktorialnost-imp-i-I}. \erem

\brem Ясно, что если $\varGamma=\varnothing$, то ни обогащений, ни
детализаций в классе $\varGamma$ объектов категории ${\tt K}$ не существует.
Поэтому интерес в этой конструкции представляют только те ситуации, когда класс
$\varGamma$ непуст. Среди них мы будем выделять следующие две: \bit{ \item[---]
$\varGamma=\Mono({\tt K})$ (то есть $\varGamma$ совпадает с классом всех
мономорфизмов категории ${\tt K}$), и тогда мы будем пользоваться обозначениями
    \begin{align}\label{imp_(varPhi)^Mono}
 &     \rf_{\varPhi}^{\Mono}X:=\rf_{\varPhi}^{\Mono({\tt K})}X, &&
  \Rf_{\varPhi}^{\Mono}X:=\Rf_{\varPhi}^{\Mono({\tt K})}X.
    \end{align}

\item[---] $\varGamma=\Mor({\tt K})$ (то есть $\varGamma$ совпадает с классом всех вообще морфизмов категории ${\tt K}$), в этом случае удобно вообще опускать упоминание о классе $\varGamma$ в формулировках и обозначениях, поэтому мы будем говорить о {\it детализации объекта $X\in{\tt K}$ в категории ${\tt K}$ посредством класса морфизмов $\varPhi$}, и обозначения будем упрощать так:
    \begin{align}\label{imp_(varPhi)=imp_(varPhi)^K}
    &     \rf_{\varPhi}X:=\rf_{\varPhi}^{\Mor({\tt K})}X, &&     \Rf_{\varPhi}X:=\Rf_{\varPhi}^{\Mor({\tt K})}X.
    \end{align}
 }\eit
\erem

\brem Другой вырожденный, но несмотря на это все же содержательный случай, --
когда $\varPhi=\varnothing$. Для фиксированного объекта $X$ здесь существенно
только, что $\varPhi$ не содержит ни одного морфизма, приходящего в $X$:
$$
\varPhi_X=\{\ph\in\varPhi:\ \Ran\ph=X\}=\varnothing.
$$
Тогда, очевидно, вообще любой морфизм $\sigma\in\varGamma$, приходящий в $X$,
$\sigma:X\gets X'$, является обогащением $X$ (в классе морфизмов $\varGamma$ посредством класса морфизмов $\varnothing$). Если вдобавок $\varGamma=\Mono$, то детализацией будет инициальный объект в категории
$\Mono_X$ (если он существует). Это можно изобразить формулой
$$
\Rf_\varnothing^\varGamma X=\min\Mono_X.
$$
С другой стороны, если ${\tt K}$ -- категория с нулем $0$, и $\varGamma$ содержит все морфизмы с началом в $0$, то детализацией любоо объекта в $\varGamma$ посредством пустого класса морфизмов является $0$:
$$
\Rf_\varnothing^\varGamma X=0.
$$
\erem

\brem Еще одна крайняя ситуация -- когда $\varPhi=\Mor({\tt K})$. Для
фиксированного объекта $X$ здесь существенно только то, что $\varPhi$ содержит
единицу $X$:
$$
1_X\in\varPhi.
$$
Тогда для произвольного обогащения $\sigma$ из диаграммы
$$
\xymatrix @R=2.pc @C=2.0pc 
{
X   & & X'\ar[ll]_{\sigma}\\
 & B\ar[ul]^{1_X} \ar@{-->}[ur] &
}
$$
следует, что $\sigma$ должен быть коретракцией (причем такой, для которой
пунктирная стрелка единственна). В частном случае, если $\varGamma\subseteq\Mono$ это
возможно только когда $\sigma$ -- изоморфизм. Как следствие, в этом случае
детализация $X$ в $\varGamma$ совпадает с $X$ (с точностью до изоморфизма):
$$
\varGamma\subseteq\Mono\quad\Longrightarrow\quad
\Rf_{\Mor({\tt K})}^{\varGamma} X=X.
$$
\erem

\medskip
\centerline{\bf Свойства детализаций:}

\bit{\it

\item[$1^\circ$.]\label{LM:Imp,suzhenie-verh-klassa-morfizmov} Пусть
$\varSigma\subseteq\varGamma$, тогда для всякого объекта $X$ и любого класса
морфизмов $\varPhi$

\bit{

\item[(a)] любое обогащение $\sigma:X\gets X'$ в классе $\varSigma$
посредством $\varPhi$ является обогащением в классе $\varGamma$
посредством $\varPhi$,

\item[(b)] если существуют детализации $\rf_\varPhi^\varSigma X$ и
$\rf_\varPhi^\varGamma X$, то существует единственный морфизм
$\rho:\Rf_\varPhi^\varSigma X\gets\Rf_\varPhi^\varGamma X$, замыкающий
диаграмму \beq\label{Imp:suzhenie-verh-klassa-morfizmov}
\begin{diagram}
\node[2]{X} \\
\node{\Rf_\varPhi^\varSigma X} \arrow{ne,t}{\rf_\varPhi^\varSigma X}
\node[2]{\Rf^\varGamma _\varPhi X}\arrow{nw,t}{\rf^\varGamma _\varPhi X}
\arrow[2]{w,b,--}{\rho}
\end{diagram}
\eeq

\item[(c)] если существует детализация $\rf_\varPhi^\varGamma X$ (в более широком классе), и она лежит в (более узком) классе $\varSigma$,
$$
\rf_\varPhi^\varGamma X\in\varSigma,
$$
то она же будет и детализацией $\rf_\varPhi^\varSigma X$ (в более широком классе):
$$
\rf_\varPhi^\varGamma X=\rf_\varPhi^\varSigma X.
$$

}\eit

\item[$2^\circ$.]\label{LM,Imp:suzhenie-verh-klassa-morfizmov-2}
Пусть $\varSigma$, $\varGamma$, $\varPhi$ - классы морфизмов, и для объекта $X$ выполняется условие
 \bit{

\item[(a)] любое обогащение $\sigma:X\gets X'$ в классе $\varGamma$ посредством $\varPhi$
содержится в $\varSigma$.
 }\eit
Тогда
 \bit{

\item[(b)] детализация $X$ в $\varGamma$ посредством $\varPhi$ существует тогда и только тогда, когда существует детализация $X$ в $\varGamma\cap\varSigma$ посредством $\varPhi$; в этом случае эти детализации совпадают:
$$
\rf_{\varPhi}^{\varGamma}=\rf_{\varPhi}^{\varGamma\cap\varSigma},
$$

\item[(c)] если $\varSigma\subseteq\varGamma$, то существование детализации $X$ в (более узком классе) $\varSigma$ посредством $\varPhi$ автоматически влечет за собой существование детализации $X$ в (более широком классе) $\varGamma$ посредством $\varPhi$ и совпадение этих детализаций:
$$
\rf_{\varPhi}^{\varGamma}X=\rf_{\varPhi}^{\varSigma}X.
$$
 }\eit

\item[$3^\circ$.]\label{LM,Imp:suzhenie-klassa-morfizmov} Пусть
$\varPsi\subseteq\varPhi$, тогда для всякого объекта $X$ и любого класса
морфизмов $\varGamma$

\bit{\it

\item[(a)] любое обогащение $\sigma:X\gets X'$ объекта $X$ в классе $\varGamma$
посредством $\varPhi$ является обогащением объекта $X$ в классе $\varGamma$
посредством $\varPsi$,

\item[(b)] если существуют детализации $\Rf_\varPsi^\varGamma  X$ и
$\Rf_\varPhi^\varGamma  X$, то существует единственный морфизм
$\alpha:\Rf_\varPsi^\varGamma  X\to\Rf_\varPhi^\varGamma  X$, замыкающий
диаграмму \beq\label{Ipm:suzhenie-klassa-morfizmov}
\begin{diagram}
\node[2]{X}  \\
\node{\Rf_\varPsi^\varGamma  X}\arrow{ne,t}{\rf_\varPsi^\varGamma
X}\arrow[2]{e,b,--}{\alpha}   \node[2]{\Rf^\varGamma_\varPhi
X}\arrow{nw,t}{\rf^\varGamma_\varPhi X}
\end{diagram}
\eeq
}\eit

\item[$4^\circ$.]\label{TH:imp_Psi=imp_Phi}  Пусть
$\varPhi\subseteq\varPsi\circ\Mor({\tt K})$ (то есть всякий морфизм
$\ph\in\varPhi$ можно представить в виде $\ph=\psi\circ\chi$, где
$\psi\in\varPsi$), тогда для всякого объекта $X$ и любого класса морфизмов
$\varGamma$

 \bit{\it

\item[(a)]  любое обогащение $\sigma:X\gets X'$ объекта $X$ в классе $\varGamma$ посредством класса $\varPsi$, являющееся мономорфизмом в $\tt K$, является обогащением $X$ в $\varGamma$ посредством $\varPhi$,

\item[(b)]  если существуют детализации $\rf_\varPsi^\varGamma X$ и
$\rf_\varPhi^\varGamma X$, причем $\rf_\varPsi^\varGamma X$ является
мономорфизмом в $\tt K$, то существует единственный морфизм
$\beta:\Rf_\varPsi^\varGamma X\gets\Rf_\varPhi^\varGamma X$, замыкающий
диаграмму \beq\label{imp_Psi=imp_Phi-0}
\begin{diagram}
\node[2]{X}  \\
\node{\Rf_\varPsi^\varGamma  X}\arrow{ne,t}{\rf_\varPsi^\varGamma  X}   \node[2]{\Rf^\varGamma  _\varPhi X}\arrow{nw,t}{\rf^\varGamma  _\varPhi X}\arrow[2]{w,b,--}{\beta}
\end{diagram}
\eeq
}\eit

\item[$5^\circ$.]\label{PROP:deistvie-monomorfizma-na-Imp}  Пусть классы
морфизмов $\varGamma$, $\varPhi$ и мономорфизм $\mu:X\gets Y$ в ${\tt K}$
удовлетворяют следующим условиям:
 \bit{
\item[(a)]  существует детализация  $\Rf_{\mu\circ\varPhi}^\varGamma X$
посредством класса морфизмов $\mu\circ\varPhi=\{\mu\circ\ph;\ \ph\in\varPhi\}$,

\item[(b)]  существует детализация $\Rf_\varPhi^\varGamma Y$,

\item[(c)]  композиция $\mu\circ\rf_{\varPhi}^\varGamma Y$ принадлежит классу
$\varGamma$.

    }\eit
Тогда существует единственный морфизм
$\upsilon:\Rf_{\mu\circ\varPhi}^\varGamma X\to \Rf_{\varPhi}^\varGamma Y$,
замыкающий диаграмму
 \beq\label{deistvie-monomorfizma-na-Imp} \xymatrix @R=2.5pc @C=4.0pc {
 X & Y\ar[l]_{\mu} \\
 \Rf_{\mu\circ\varPhi}^\varGamma X\ar[u]^{\rf_{\mu\circ\varPhi}^\varGamma X}\ar@{-->}[r]_{\upsilon} &   \Rf_{\varPhi}^\varGamma Y\ar[u]_{\rf_{\varPhi}^\varGamma Y}\ar@/^1ex/[ul]_{\ \mu\circ\rf_{\varPhi}^\varGamma Y}
}
\eeq

}\eit

 \bit{
\item[$\bullet$]\label{DEF:morfizmy-porozhdayutsya-snaruzhi} Условимся
говорить, что в категории ${\tt K}$ {\it класс морфизмов $\varPhi$ порождается
снаружи классом морфизмов $\varPsi$}, если
    $$
    \varPsi\subseteq\varPhi\subseteq\varPsi\circ\Mor({\tt K}).
    $$
    }\eit
Следующее утверждение двойственно теореме \ref{TH:morfizmy-porozhdayutsya-iznutri} и доказывается по аналогии:

\btm\label{TH:morfizmy-porozhdayutsya-snaruzhi} Пусть в категории ${\tt K}$ класс морфизмов $\varPhi$ порождается снаружи классом морфизмов $\varPsi$. Тогда для любого класса мономорфизмов $\varGamma$ (необязательно, чтобы $\varGamma$ включал все мономорфизмы категории ${\tt K}$) и всякого объекта $X$ существование детализации $\rf^\varGamma_\varPsi X$ эквивалентно существованию детализации $\rf^\varGamma_\varPhi X$, и эти детализации совпадают:
 \beq\label{imp_Psi=imp_Phi}
\rf^\varGamma  _\varPsi X=\rf^\varGamma  _\varPhi X.
 \eeq
\etm

\bit{
\item[$\bullet$]\label{DEF:varPhi-razlich-morfizmy-iznutri} Условимся говорить, что класс морфизмов $\varPhi$ в категории ${\tt K}$
 {\it различает морфизмы изнутри}, если для любых двух различных параллельных морфизмов $\alpha\ne\beta:X\to Y$ найдется морфизм $\ph:M\to X$ из класса $\varPhi$ такой, что $\alpha\circ\ph\ne\beta\circ\ph$.
}\eit

Следующая теорема двойственна теореме \ref{TH:Phi-razdel-moprfizmy} :

\btm\label{TH:Phi-razdel-moprfizmy-iznutri}  Если класс морфизмов $\varPhi$ различает морфизмы изнутри, то для любого класса морфизмов $\varGamma$
 \bit{
\item[(i)] всякое обогащение $\varGamma$ посредством $\varPhi$ является эпиморфизмом,

\item[(ii)] детализация посредством $\varPhi$ в классе $\varGamma$ существует тогда и только тогда, когда существует детализация посредством $\varPhi$ в классе $\varGamma\cap\Mono$; в этом случае эти детализации совпадают:
$$
\rf_{\varPhi}^{\varGamma}=\rf_{\varPhi}^{\varGamma\cap\Epi},
$$

\item[(iii)] если класс $\varGamma$ содержит все эпиморфизмы,
 $$
 \varGamma\supseteq\Epi,
 $$
 то существование детализации посредством $\varPhi$ в классе $\Epi$ автоматически влечет за собой существование детализации посредством $\varPhi$ в классе $\varGamma$ и совпадение этих детализаций:
$$
\rf_{\varPhi}^{\varGamma}=\rf_{\varPhi}^{\Epi}.
$$
 }\eit
\etm

\bit{
\item[$\bullet$] Напомним, что класс морфизмов $\varPhi$ в категории ${\tt K}$ называется {\it левым идеалом}, если
$$
\Mor({\tt K})\circ\varPhi\subseteq\varPhi
$$
(то есть для любого $\ph\in\varPhi$ и любого морфизма $\mu$ категории ${\tt K}$ композиция $\mu\circ\ph$ принадлежит $\varPhi$).
}\eit

Следующее утверждение двойственно теореме \ref{TH:Phi-razdel-moprfizmy-*}

\btm\label{TH:Phi-razdel-moprfizmy-iznutri-*} Если класс морфизмов $\varPhi$ различает морфизмы изнутри и является левым идеалом в категории ${\tt K}$, то для любого класса морфизмов $\varGamma$
 \bit{
\item[(i)] всякое обогащение $\varGamma$ посредством $\varPhi$ является биморфизмом,

\item[(ii)] детализация посредством $\varPhi$ в классе $\varGamma$ существует тогда и только тогда, когда существует детализация посредством $\varPhi$ в классе $\varGamma\cap\Bim$; в этом случае эти детализации совпадают:
$$
\rf_{\varPhi}^{\varGamma}=\rf_{\varPhi}^{\varGamma\cap\Bim}.
$$

\item[(iii)] если класс $\varGamma$ содержит все биморфизмы,
 $$
 \varGamma\supseteq\Bim,
 $$
 то существование детализации посредством $\varPhi$ в классе $\Bim$ автоматически влечет за собой существование детализации посредством $\varPhi$ в классе $\varGamma$ и совпадение этих детализаций:
$$
\rf_{\varPhi}^{\varGamma}=\rf_{\varPhi}^{\Bim}.
$$
 }\eit
\etm

Важная ситуация --- когда классы $\varPhi$ и $\varGamma$ совпадают:
$$
\varPhi=\varGamma.
$$
В этом случае понятия обогащения и детализации совпадают:

\btm\label{varPhi=varGamma}
Всякое обогащение объекта $X$ в классе $\varGamma$ посредством самого этого класса $\varGamma$ является детализацией объекта $X$ в классе $\varGamma$ посредством $\varGamma$.
\etm
\bpr
Пусть $\rho:E\to X$ --- произвольное обогащение объекта $X$ в классе морфизмов $\varGamma$
посредством того же самого класса морфизмов $\varGamma$. Если $\sigma:X'\to X$ --- какое-то другое обогащение объекта $X$ в классе морфизмов $\varGamma$ посредством класса
морфизмов $\varGamma$, то в диаграмме \eqref{DIAGR:otpechatok}
$$
\begin{diagram}
\node[2]{X}  \\
\node{E}\arrow{ne,t}{\rho}\arrow[2]{e,b,--}{\exists!\upsilon}
\node[2]{X'}\arrow{nw,t}{\forall\sigma}
\end{diagram}
$$
морфизм $\upsilon:E\to X'$ существует и единственен, потому что $\rho\in\varGamma$. Это значит, что $\rho:E\to X$ --- детализация.
\epr

\paragraph{Детализация в классе объектов посредством класса объектов.}

Частным случаем конструкции детализации является ситуация, когда $\varGamma$ и/или $\varPhi$ представляют собой классы всех морфизмов из объектов из заданных подклассов класса $\Ob({\tt K})$. Точная формулировка для случая, когда оба класса $\varGamma$ и $\varPhi$ определяются таким образом, выглядит так: пусть дана категория ${\tt K}$ и два подкласса ${\tt L}$ и ${\tt M}$ в классе $\Ob({\tt K})$ объектов ${\tt K}$.

 \bit{
\item[$\bullet$] Морфизм $\sigma:X'\to X$ называется {\it обогащением объекта $X\in{\tt K}$ в классе объектов ${\tt L}$ посредством класса объектов ${\tt M}$}, если для любого объекта $B\in{\tt M}$ и любого морфизма $\ph:B\to X$ найдется единственный морфизм $\ph':B\to X'$, замыкающий диаграмму:
$$
\put(5,-35){$\begin{matrix}\text{\rotatebox{90}{$\owns$}}\\ {\tt M}
\end{matrix}$} \put(71,-35){$\begin{matrix}\text{\rotatebox{90}{$\owns$}}\\
{\tt L} \end{matrix}$}
\begin{diagram}
\node[2]{X}  \\
\node{B}\arrow{ne,t}{\forall\ph}\arrow[2]{e,b,--}{\exists!\ph'}
\node[2]{X'}\arrow{nw,t}{\sigma}
\end{diagram}
$$

\item[$\bullet$]\label{DEF:nachinka-otn-klassa} Обогащение $\rho:E\to X$
объекта $X\in{\tt K}$ в классе объектов ${\tt L}$ посредством класса объектов
${\tt M}$ называется {\it детализацией объекта $X\in{\tt K}$ в классе объектов ${\tt L}$ посредством класса объектов ${\tt M}$}, и мы будем обозначать это
записью
    \beq\label{DEF:tr_M^L(A)}
    \rho=\rf_{\tt M}^{\tt L} X,
    \eeq
если для любого другого обогащения $\sigma:X'\to X$ (объекта $X\in{\tt K}$ в классе объектов ${\tt L}$ посредством класса объектов ${\tt M}$) найдется единственный
морфизм $\upsilon:E\to X'$,
замыкающий диаграмму:
 \beq\label{diagr:sled}
\put(5,-35){$\begin{matrix}\text{\rotatebox{90}{$\owns$}}\\ {\tt L}
\end{matrix}$} \put(71,-35){$\begin{matrix}\text{\rotatebox{90}{$\owns$}}\\
{\tt L} \end{matrix}$}
\begin{diagram}
\node[2]{X}  \\
\node{E}\arrow{ne,t}{\rho}\arrow[2]{e,b,--}{\exists!\upsilon}
\node[2]{X'}\arrow{nw,t}{\forall\sigma}
\end{diagram}
 \eeq
Сам объект $E$ также называется {\it детализацией} объекта $X$ в классе объектов ${\tt L}$ (посредством класса объектов ${\tt M}$), и для него мы используем обозначение
    \beq\label{DEF:Tr_M^L(A)}
E=\Rf_{\tt M}^{\tt L} X.
 \eeq
}\eit

Возможны также следующие две крайние ситуации в выборе класса ${\tt L}$:
 \bit{

\item[---] если ${\tt L} =\Ob({\tt K})$, то мы будем говорить о {\it детализации объекта $X\in{\tt K}$ в категории ${\tt K}$ посредством класса объектов ${\tt M}$}, и обозначения будем упрощать так:
    \begin{align}\label{imp_M=imp_M^K}
    &     \rf_{\tt M}X:=\rf_{\tt M}^{\tt K}X, &&     \Rf_{\tt M}X:=\Rf_{\tt M}^{\tt K}X,
    \end{align}

\item[---] если ${\tt L} ={\tt M}$, то (по теореме \ref{varPhi=varGamma}) понятия обогащения и детализации совпадают: {\it всякое обогащение объекта $X\in{\tt K}$ в классе ${\tt L}$ посредством самого этого класса ${\tt L}$ является детализацией объекта $X\in{\tt K}$ в классе ${\tt L}$ посредством класса ${\tt L}$}; чтобы не повторяться, в этом случае говорят просто о {\it детализации объекта $X$ в классе ${\tt L}$}, при этом обозначения упрощаются так:
    \begin{align}
    & \rf_{\tt L} ^{\tt L} X=:\rf^{\tt L} X, &&  \Rf_{\tt L} ^{\tt L} X=:\Rf^{\tt L} X
    \end{align}
}\eit

\bit{
\item[$\bullet$] Условимся говорить, что класс объектов ${\tt M}$ в категории ${\tt K}$
 {\it различает морфизмы изнутри}, если класс всех морфизмов, выходящих из объектов, принадлежащих ${\tt M}$ обладает этим свойством (в смысле определения на с.\pageref{DEF:varPhi-razlich-morfizmy-iznutri}), то есть для любых двух различных параллельных морфизмов $\alpha\ne\beta:X\to Y$ найдется морфизм $\ph:M\to X$ такой, что $\alpha\circ\ph\ne\beta\circ\ph$.
}\eit

Из теоремы \ref{TH:Phi-razdel-moprfizmy-iznutri-*} сразу следует

\btm\label{TH:M-razdel-moprfizmy-iznutri} Если класс объектов ${\tt M}$ различает морфизмы изнутри, то для любого класса объектов ${\tt L}$
 \bit{
\item[(i)]  всякое обогащение в классе ${\tt L}$ посредством класса ${\tt M}$ является биморфизмом,

\item[(ii)] детализация в классе ${\tt L}$ посредством класса ${\tt M}$ существует тогда и только тогда, когда существует детализация в классе биморфизмов, выходящих из ${\tt L}$, посредством класса ${\tt M}$; в таком случае эти детализации совпадают:
$$
\rf_{\tt M}^{\tt L}=\rf_{\tt M}^{\Bim({\tt L},{\tt K})}.
$$
 }\eit
\etm

\paragraph{Примеры детализаций.}

\bex {\bf Односвязное накрытие}, используемое в топологии и теории групп Ли,
представляет собой с категорной точки зрения детализацию в классе односвязных пунктированных
накрытий посредством пустого класса морфизмов в категории связных локально
связных и полулокально односвязных пунктированных топологических пространств
(определения см. в \cite{Postnikov}). \eex

\bex {\bf Борнологизация} (``bornologification'', см. определение в
\cite{Kriegl-Michor}) $X_{\text{\rm born}}$ локально выпуклого пространства $X$
является детализацией $X$ в категории ${\tt LCS}$ локально выпуклых пространств
посредством подкатегории ${\tt Norm}$ нормированных пространств:
$$
X_{\text{\rm born}}=\Rf_{\tt Norm}^{\tt LCS}X
$$
\eex
\bpr Это следует из характеризации борнологизации, как сильнейшей локально
выпуклой топологии на $X$, при которой все вложения $X_B\to X$ будут
непрерывны, где $B$ пробегает систему всех ограниченных абсолютно выпуклых
подмножеств в $X$, а $X_B$ -- нормированное пространство с единичным шаром $B$
(см.\cite[Chapter I, Lemma 4.2]{Kriegl-Michor}).
 \epr

\bex {\bf Насыщение} $X^\blacktriangle$ псевдополного локально выпуклого
пространства $X$ является детализацией $X$ в категории ${\tt LCS}$ локально выпуклых
пространств посредством подкатегории ${\tt Smi}$
пространств Смит (см. определения в \cite{Ak03}):
$$
X^\blacktriangle=\Rf_{\tt Smi}^{\tt LCS}X
$$
 \eex

\subsection{Связь с факторизациями и узловым разложением}\label{SUBSEC:obolochki<->uzl-razlozh}

\paragraph{Факторизация категории.}

\bit{

\item[$\bullet$]
Пара морфизмов $(\mu,\e)$ называется {\it диагонализируемой}\index{диагонализируемая пара морфизмов}
\cite{General-algebra,Tsalenko-Shulgeifer,Akbarov-De-Gruyter-I}, если для любых морфизмов
$\alpha:\Dom\e\to\Dom\mu$ и $\beta:\Ran\e\to\Ran\mu$ со свойством $\mu\circ\alpha=\beta\circ\e$
существует морфизм $\delta:B\to C$, делающий коммутативной диаграмму
$$
\begin{diagram}
\node{\Dom\e}\arrow{s,l}{\alpha}\arrow{e,t}{\e}\node{\Ran\e}\arrow{s,r}{\beta}\arrow{sw,t,--}{\delta}
\\
\node{\Dom\mu}\arrow{e,b}{\mu}\node{\Ran\mu}
\end{diagram}
$$
Это обозначается записью $\mu\downarrow\e$.
}\eit

\bit{
\item[$\bullet$] Для всякого класса морфизмов $\varLambda$ в $\tt{K}$

 \bit{
\item[---] его {\it эпиморфно сопряженным классом}\index{класс!эпиморфно сопряженный} называется класс
$$
\varLambda^\downarrow=\{\e\in\Epi({\tt{K}}):\forall\lambda\in\varLambda\quad \lambda\downarrow\e\}.
$$
\item[---] его {\it мономорфно сопряженным классом}\index{класс!мономорфно сопряженный} называется класс
$$
{^\downarrow\kern-1pt\varLambda}=\{\mu\in\Mono({\tt{K}}):\forall\lambda\in\varLambda\quad \mu\downarrow\lambda\}.
$$
}\eit
}\eit
Сразу из определений видно, что для всякого класса морфизмов $\varLambda$ справедливы включения
\begin{align}
& \Iso\subseteq\varLambda^\downarrow\subseteq\Epi, && \Iso\circ\ \varLambda^\downarrow\subseteq\varLambda^\downarrow
\label{Iso-subseteq-varTheta^downarrow-subseteq-Epi}
\\
& \Iso\subseteq{^\downarrow\kern-1pt\varLambda}\subseteq\Mono, &&
{^\downarrow\kern-1pt\varLambda}\circ\Iso\subseteq{^\downarrow\kern-1pt\varLambda}
\label{Iso-subseteq-^downarrow-varTheta-subseteq-Mono}
\end{align}

\bit{

\item
Условимся говорить, что классы морфизмов $\varGamma$ и $\varOmega$
задают {\it факторизацию категории\footnote{Эта конструкция часто называется {\it бикатегорией}
\cite{General-algebra,Tsalenko-Shulgeifer}.} $\tt{K}$}\label{DEF:faktorizatsija-v-kategorii}\index{факторизация!категории}, если
 \bit{
 \item[F.1] $\varOmega$ является эпиморфно сопряженным классом для $\varGamma$:
 $$
 \varGamma^\downarrow=\varOmega
 $$

 \item[F.2] $\varGamma$ является мономорфно сопряженным классом для $\varOmega$:
 $$
 \varGamma={^\downarrow\varOmega},
 $$

 \item[F.3] композиция классов $\varGamma$ и $\varOmega$ покрывает класс всех морфизмов:
 $$
 \varGamma\circ\varOmega=\Mor({\tt{K}})
 $$
 (это означает, что всякий морфизм $\ph\in\Mor({\tt{K}})$ можно представить в виде композиции $\mu\circ\e$,
 где $\mu\in\varGamma$, $\e\in\varOmega$).
 }\eit
Если эти условия выполнены, мы будем изображать это записью
 \beq\label{K=varGamma-circledcirc-varOmega}
{\tt{K}}=\varGamma\circledcirc\varOmega.
 \eeq
}\eit

\noindent\rule{160mm}{0.1pt}\begin{multicols}{2}

\bex\label{TH:faktorizatsija-v-kategorii-s-uzlov-razlozh-1}
В категории с узловым разложением $\tt{K}$ (см. определение в \cite[p.129]{Akbarov-De-Gruyter-I})
следующие классы морфизмов задают факторизации:
$$
{\tt{K}}=\Mono\circledcirc\SEpi=\SMono\circledcirc\Epi.
$$
\eex

\end{multicols}\noindent\rule[10pt]{160mm}{0.1pt}

Следующий факт доказывается в \cite[Теорема 8.2]{Tsalenko-Shulgeifer}:

\btm\label{TH:o-faktorizatsii}
Классы $\varGamma$ и $\varOmega$ образуют факторизацию в категории $\tt{K}$
$$
{\tt{K}}=\varGamma\circledcirc\varOmega
$$
тогда и только тогда, когда выполняются условия
\bit{
 \item[(i)] $\varGamma\subseteq\Mono({\tt{K}})$ и $\varOmega\subseteq\Epi({\tt{K}})$,

 \item[(ii)] $\Iso({\tt{K}})\subseteq\varOmega\cap\varGamma$,

 \item[(iii)] для любого морфизма $\ph\in\Mor({\tt{K}})$ существует разложение
 \beq\label{faktorizatsiya-v-kat-s-faktoriz}
 \ph=\mu_{\ph}\circ\e_{\ph},\qquad \mu_{\ph}\in\varGamma,\quad \e_{\ph}\in\varOmega
 \eeq
 \item[(iv)] для любого другого разложения с теми же свойствами
$$
 \ph=\mu'\circ\e',\qquad \mu'\in\varGamma,\quad \e'\in\varOmega
$$
найдется морфизм $\theta\in\Iso({\tt{K}})$ такой что
$$
\mu'=\mu_{\ph}\circ\theta,\qquad \e'=\theta^{-1}\circ\e_{\ph}.
$$
 }\eit
\etm

 \bit{
\item[$\bullet$] Условимся говорить, что класс морфизмов $\varOmega$ в $\tt{K}$
{\it мономорфно дополняем}\label{DEF:klass-monomorfno-dopolnyaem}\index{класс!мономорфно дополняемый}, если
\beq\label{klass-monomorfno-dopolnyaem}
{\tt{K}}={^\downarrow\varOmega}\circledcirc\varOmega.
\eeq
Иными словами, для этого нужно, чтобы $\varOmega$ был эпиморно сопряженным классом
к своему мономорфно сопряженному классу
$$
\varOmega=(^\downarrow\varOmega)^\downarrow,
$$
и чтобы композиция классов ${^\downarrow\varOmega}$ и $\varOmega$ покрывала класс всех морфизмов:
$$
{^\downarrow\varOmega}\circ\varOmega=\Mor(\tt K).
$$
В этом случае класс ${^\downarrow\varOmega}$ мы будем называть {\it мономорфным дополнением}\index{дополнение!мономорфное} к $\varOmega$.
}\eit

\brem
Из \eqref{Iso-subseteq-varTheta^downarrow-subseteq-Epi} сразу следует, что если класс морфизмов
$\varOmega$ мономорфно дополняем, то
\beq\label{Iso-circ-varOmega-subseteq-varOmega}
\Iso\subseteq\varOmega\subseteq\Epi,\qquad \Iso\circ\varOmega\subseteq\varOmega
\eeq
\erem

 \bit{
\item[$\bullet$] Аналогично, мы говорим, что класс морфизмов $\varGamma$ в $\tt{K}$
{\it эпиморфно дополняем}\label{DEF:klass-epimorfno-dopolnyaem}\index{класс!эпиморфно дополняемый}, если
\beq\label{klass-epimorfno-dopolnyaem}
{\tt{K}}=\varGamma\circledcirc\varGamma^\downarrow.
\eeq
Иными словами, для этого нужно, чтобы $\varGamma$ был мономорно сопряженным классом
к своему эпиморфно сопряженному классу
$$
\varGamma={^\downarrow(\varGamma^\downarrow)},
$$
и чтобы композиция классов $\varGamma$ и $\varGamma^\downarrow$ покрывала класс всех морфизмов:
$$
\varGamma\circ\varGamma^\downarrow=\Mor(\tt K).
$$
В этом случае класс $\varGamma^\downarrow$ мы будем называть {\it эпиморфным дополнением}\index{дополнение!эпиморфное} к $\varGamma$.
}\eit

\brem
Из \eqref{Iso-subseteq-^downarrow-varTheta-subseteq-Mono} сразу следует, что если класс морфизмов
$\varGamma$ эпиморфно дополняем, то
\beq\label{varGamma-circ-Iso-subseteq-varGamma}
\Iso\subseteq\varGamma\subseteq\Mono,\qquad \varGamma\circ\Iso\subseteq\varGamma.
\eeq
\erem

\paragraph{Связь с проективным и инъективным пределами.} Сходство между понятиями оболочки и проективного предела формализуется в следующих двух утверждениях.

\blm\label{LM:obolochka-konusa}\footnote{В работе \cite{Ak16} этот результат (Lemma 3.23) приведен с ошибкой: там опущено условие \eqref{rho-in-Epi}.} Если проективный предел $\rho=\projlim \rho^i:X\to\projlim X^i$
проективного конуса $\{\rho^i:X\to X^i;\ i\in I\}$ из заданного объекта $X$ в ковариантную
(или контравариантную) систему $\{X^i,\iota^j_i\}$ является эпиморфизмом,
\beq\label{rho-in-Epi}
\rho\in\Epi,
\eeq
то он является оболочкой объекта $X$ в произвольном классе $\varOmega$, содержащем $\rho,$ относительно системы морфизмов $\{\rho^i; i\in I\}$:
 \beq\label{obolochka-konusa-Omega}
\rho=\projlim \rho^i\in\varOmega
\quad\Longrightarrow\quad
\Env_{\{\rho^i;\ i\in I\}}^\varOmega X=\projlim X^i
 \eeq
 В частности, это всегда верно при $\varOmega=\Mor({\tt K})$:
 \beq\label{obolochka-konusa}
\Env_{\{\rho^i;\ i\in I\}}^{\Mor({\tt K})}X=\projlim X^i
 \eeq
\elm
 \bpr
1. Прежде всего, морфизм $\rho$ является расширением $X$ относительно системы $\{\rho^i\}$, потому что в диаграмме
\beq\label{proof-obolochka=lim-0}
\begin{diagram}
\node{X}\arrow[2]{e,t}{\rho}\arrow{se,b}{\rho^j}\node[2]{\projlim X^i}\arrow{sw,b,--}{\pi^j}
\\
\node[2]{X^j}
\end{diagram}
\eeq
(основном свойстве проективного предела) продолжение $\pi^j$ должно быть единственным в силу условия \eqref{rho-in-Epi}.

2. Далее, пусть $\sigma:X\to X'$ -- какое-нибудь другое расширение. Тогда для всякого морфизма
$\rho^j:X\to X^j$ существует единственный морфизм $\upsilon^j:X'\to X^j$, замыкающий диаграмму
\beq\label{proof-obolochka=lim-1}
\begin{diagram}
\node{X}\arrow[2]{e,t}{\sigma}\arrow{se,b}{\rho^j}\node[2]{X'}\arrow{sw,b,--}{\upsilon^j}
\\
\node[2]{X^j}
\end{diagram}
\eeq
При этом для любых индексов $i\le j$ в диаграмме
$$
  \xymatrix @R=2.5pc @C=2.5pc
 {
 & X\ar@/_4ex/[ldd]_{\rho^i}\ar[d]^{\sigma} \ar@/^4ex/[rdd]^{\rho^j} & \\
 & X'\ar@{-->}[ld]^{\upsilon^i}\ar@{-->}[rd]_{\upsilon^j} & \\
 X^i \ar@/_2ex/[rr]^{\iota_i^j}  & & X^j
 }
$$
будут коммутативны два верхних маленьких треугольника (у каждого по одной пунктирной стороне) и большой треугольник по периметру (без пунктирных сторон). Отсюда и из единственности стрелки $\upsilon^j$ в верхнем правом маленьком треугольнике следует, что коммутативен и маленький треугольник с двумя пунктирными сторонами:
$$
\begin{cases}(\iota_i^j\circ\upsilon^i)\circ\sigma=\iota_i^j\circ(\upsilon^i\circ\sigma)=\iota_i^j\circ\rho^i=\rho^j
\\
\upsilon^j\circ\sigma=\rho^j
\end{cases}
\qquad\Longrightarrow\qquad \iota_i^j\circ\upsilon^i=\upsilon^j.
$$
Коммутативность треугольника с двумя пунктирными сторонами, в свою очередь, означает, что $X'$ с системой морфизмов $\upsilon^i$ представляет собой проективный конус ковариантной системы $\{X^i;\iota_i^j\}$. Поэтому должен существовать однозначно определенный морфизм $\upsilon$, для которого при любом $j$ в диаграмме
$$
  \xymatrix @R=2.5pc @C=2.5pc
 {
 & X\ar@/_4ex/[ldd]_{\rho}\ar[d]^{\sigma} \ar@/^4ex/[rdd]^{\rho^j} & \\
 & X'\ar@{-->}[ld]^{\upsilon}\ar[rd]_{\upsilon^j} & \\
 \projlim X^i \ar@/_2ex/[rr]^{\pi^j}  & & X^j
 }
$$
маленький нижний треугольник будет коммутативен. С другой стороны, здесь правый маленький треугольник коммутативен, потому что это повернутая диаграмма \eqref{proof-obolochka=lim-1}, а периметр коммутативен, потому что это повернутая диаграмма \eqref{proof-obolochka=lim-0}. Отсюда и из единственности морфизма $\rho$ в системе таких периметров при разных $j$ следует, что левый маленький треугольник тоже должен быть коммутативен:
$$
\left(
\forall j\qquad
\begin{cases}
\pi^j\circ\upsilon\circ\sigma=\upsilon^j\circ\sigma=\rho^j
\\
\pi^j\circ\rho=\rho^j
\end{cases}\right)\qquad\Longrightarrow\qquad \upsilon\circ\sigma=\rho
$$

Мы поняли, что существует морфизм $\upsilon$, замыкающий диаграмму \eqref{DEF:diagr-obolochka} (с $E=\projlim X^i$). Остается проверить, что такой морфизм единственен. Пусть $\upsilon'$ -- какой-то другой морфизм с тем же свойством: $\rho=\upsilon'\circ\sigma$. Рассмотрим диаграмму
$$
  \xymatrix @R=2.5pc @C=2.5pc
 {
 & X\ar@/_4ex/[ldd]_{\rho}\ar[d]^{\sigma} \ar@/^4ex/[rdd]^{\rho^j} & \\
 & X'\ar@{-->}[ld]^{\upsilon'}\ar[rd]_{\upsilon^j} & \\
 \projlim X^i \ar@/_2ex/[rr]^{\pi^j}  & & X^j
 }
$$
В ней помимо верхнего левого маленького треугольника будут коммутативны также верхний правый маленький треугольник (потому что это повернутая диаграмма \eqref{proof-obolochka=lim-1}) и периметр (потому что это повернутая диаграмма \eqref{proof-obolochka=lim-0}). Отсюда и из единственности стрелки $\upsilon^j$ в верхнем правом маленьком треугольнике следует, что коммутативен и нижний маленький треугольник:
$$
\begin{cases}
\pi^j\circ\upsilon'\circ\sigma=\pi^j\circ\rho=\rho^j
\\
\upsilon^j\circ\sigma=\rho^j
\end{cases}
\qquad\Longrightarrow\qquad
\pi^j\circ\upsilon'=\upsilon^j.
$$
Это верно для всякого
индекса $j$, поэтому морфизм $\upsilon'$ должен совпадать с построенными нами выше морфизмом $\upsilon$, то есть $\upsilon'=\upsilon$.
\epr

\blm\label{LM:obolochka-konusa-v-kat-s-uzl-razl} Пусть $\varOmega$ -- мономорфно дополняемый класс в категории
${\tt K}$, $\{X^i,\iota^j_i\}$ -- ковариантная (или контравариантная)
система и $\{\rho^i:X\to X^i;\ i\in I\}$ -- проективный конус из
заданного объекта $X$ в $\{X^i,\iota^j_i\}$. Если существует проективный предел
$\rho=\projlim\rho^i:X\to\projlim X^i$, то в
его факторизации
$$
\rho=\mu_{\rho}\circ\e_{\rho},\qquad \mu_{\rho}\in{^\downarrow\varOmega},\quad \e_{\rho}\in\varOmega
$$
эпиморфизм $\e_{\rho}$ является оболочкой объекта $X$ относительно системы морфизмов $\{\rho^i; i\in I\}$
в классе $\varOmega$:
\begin{align}\label{obolochka-konusa-v-kat-s-uzl-razl}
&
\e_{\projlim \rho^i}=
\e_{\rho}=
\env_{\{\rho^i;\ i\in I\}}^{\varOmega}X, &&
\Ran\e_{\projlim \rho^i}=\Ran\e_{\rho}=
\Env_{\{\rho^i;\ i\in I\}}^{\varOmega}X
\end{align}
\elm
 \bpr
1. По определению проективного предела всякий морфизм $\rho^j$ обладает продолжением $\pi^j$ на $\projlim X^i$.
Ограничение $\pi^j$ на $\Ran\e_{\rho}$, то есть композиция $\tau^j=\pi^j\circ \mu_{\rho}$ будет
продолжением $\rho^j$ на $\Ran\e_{\rho}$ вдоль $\e_{\rho}$:
\beq\label{proof-obolochka=lim-0-v-kat-s-uzl-razl}
\xymatrix @R=2.pc @C=5.0pc 
{
X\ar[r]_{\e_{\rho}}\ar@/_1ex/[dr]_{\rho^j}\ar@/^4ex/[rr]^{\rho}
 &  \Ran\e_{\rho}\ar[r]_{\mu_\rho}\ar@{-->}[d]^{\tau^j}   & \projlim X^i \ar@/^1ex/@{-->}[dl]^{\pi^j}
\\
& X^j &
}
\eeq
Такое продолжение $\tau^j$ будет единственно из-за эпиморфности $\e_{\rho}$, и мы можем
сделать вывод, что морфизм $\e_{\rho}$ является расширением $X$ в $\varOmega$ относительно системы $\{\rho^i\}$.

2. Далее, пусть $\sigma:X\to X'$ -- какое-нибудь другое расширение $X$ в $\varOmega$ относительно $\{\rho^i\}$.
Как в доказательстве леммы \ref{LM:obolochka-konusa}, найдем морфизм $\upsilon$ такой, что
$\upsilon\circ\sigma=\rho$.
Для него мы получим
$$
\upsilon\circ\sigma=\rho=\mu_\rho\circ\e_\rho
$$
и поскольку $\sigma\in\varOmega$, $\mu_\rho\in{^\downarrow\varOmega}$, должен существовать диагональный морфизм
$\delta$ такой что
$$
\delta\circ\sigma=\e_\rho.
$$
Этот морфизм единственен, потому что $\sigma\in\varOmega\subseteq\Epi$.
\epr

Двойственные результаты для детализаций выглядят так.

\blm\label{LM:nachinka-konusa}\footnote{В работе \cite{Ak16} этот результат (Lemma 3.25) приведен с ошибкой: там опущено условие \eqref{rho-in-Mono}.}
Если инъективный предел $\rho=\injlim \rho^i:X\gets\injlim X^i$ инъективного конуса $\{\rho^i:X\gets X^i;\ i\in I\}$ в заданный объект $X$ из ковариантной (или контравариантной) системы $\{X^i,\iota^j_i\}$ является мономорфизмом
\beq\label{rho-in-Mono}
\rho\in\Mono,
\eeq
то он является детализацией объекта $X$ в произвольном классе $\varGamma$, содержащим $\rho$, посредством системы
морфизмов $\{\rho^i; i\in I\}$:
 \beq\label{otpechatok-konusa-Omega}
\rho=\injlim \rho^i\in\varGamma
\quad\Longrightarrow\quad
\Rf_{\{\rho^i;\ i\in I\}}^\varGamma X=\injlim X^i
 \eeq
 В частности, это всегда верно при $\varGamma=\Mor({\tt K})$:
 \beq\label{otpechatok-konusa}
\Rf_{\{\rho^i;\ i\in I\}}^{\Mor({\tt K})}X=\injlim X^i
 \eeq
 \elm

\blm\label{LM:otpechatok-konusa-v-kat-s-uzl-razl}
Пусть $\varGamma$ -- эпиморфно дополняемый класс в категории ${\tt K}$,
$\{X^i,\iota^j_i\}$ -- ковариантная (или контравариантная)
система и $\{\rho^i:X\gets X^i;\ i\in I\}$ -- инъективный конус из $\{X^i,\iota^j_i\}$ в
заданный объект $X$. Если существует
инъективный предел $\rho=\injlim\rho^i:X\gets\injlim X^i$, то в его
факторизации
$$
\rho=\mu_{\rho}\circ\e_{\rho},\qquad \mu_{\rho}\in\varGamma,\quad \e_{\rho}\in\varGamma^\downarrow
$$
мономорфизм $\mu_{\rho}$ является детализацией объекта $X$ в классе $\varGamma$ посредством системы
морфизмов $\{\rho^i; i\in I\}$:
 \begin{align}\label{otpechatok-konusa-v-kat-s-uzl-razl}
& \rf_{\{\rho^i;\ i\in I\}}^{\varGamma}X=\mu_{\rho}=\mu_{\injlim \rho^i},
&& \Rf_{\{\rho^i;\ i\in I\}}^{\varGamma}X=\Dom\mu_{\rho}=\Dom\mu_{\injlim \rho^i}
 \end{align}
\elm

\paragraph{Существование оболочек и детализаций для дополняемых классов.}

\blm\label{LM:env_varPhi^varOmega_X=env_(e_ph;ph_in_varPhi)^varOmega_X}
Пусть $\varOmega$ -- мономорфно дополняемый класс в категории ${\tt K}$.
Тогда для всякого объекта $X$ и любого класса морфизмов $\varPhi$
\beq\label{env_varPhi^varOmega_X=env_(e_ph;ph_in_varPhi)^varOmega_X}
\env_{\varPhi}^{\varOmega}X=\env_{\{\e_{\ph};\ \ph\in \varPhi\}}^{\varOmega}X
\eeq
(если какая-то из этих оболочек существует, то другая тоже существует, и они совпадают).
\elm
\bpr Пусть $\ph=\mu_{\ph}\circ\e_{\ph}$ -- факторизация морфизма с $\mu_{\ph}\in{^\downarrow\varOmega}$ и
$\e_{\ph}\in\varOmega$. Нам нужно убедиться, что расширения относительно классов
$\varPhi$ и $\{\e_{\ph};\ \ph\in \varPhi\}$ одни и те же. Пусть $\sigma:X\to X'$ --
расширение $X$ в классе $\varOmega$ относительно морфизмов $\{\e_{\ph};\ \ph\in \varPhi\}$. Тогда в диаграмме
$$
\xymatrix 
{
X\ar[rr]^{\sigma}\ar[rd]_{\e_{\ph}}\ar@/_6ex/[rdd]_{\ph} &  & X'\ar@{-->}[ld]^{\e'}\ar@{-->}@/^6ex/[ldd]^{\ph'} \\
 & \Ran\e_{\ph}\ar[d]^{\mu_{\ph}} & \\
 &  Y  &
}
$$
из существования морфизма $\e'$, замыкающего верхний маленький треугольник, следует
существование морфизма $\ph'$, замыкающего правый нижний маленький треугольник, и,
поскольку оставшийся маленький треугольник, левый нижний, коммутативен, получается,
что большой треугольник (периметр) также коммутативен. Вдобавок, из того, что $\sigma$ --
эпиморфизм, следует единственность $\ph'$. Значит, $\sigma:X\to X'$ -- расширение $X$ относительно морфизмов $\varPhi$.

Наоборот, пусть $\sigma:X\to X'$ -- расширение $X$ в классе $\varOmega$ относительно морфизмов $\varPhi$. Тогда для всякого $\ph\in\varPhi$ существует морфизм $\ph'$, для которого в диаграмме
$$
\xymatrix 
{
X\ar[rr]^{\sigma}\ar[rd]_{\e_{\ph}}\ar@/_6ex/[rdd]_{\ph} &  & X'\ar@{-->}@/^6ex/[ldd]^{\ph'} \\
 & \Ran\e_{\ph}\ar[d]^{\mu_{\ph}} & \\
 &  Y  &
}
$$
будет коммутативен большой треугольник (периметр). Маленький треугольник (слева внизу) здесь также коммутативен, в силу \eqref{faktorizatsiya-v-kat-s-faktoriz}, значит коммутативен и четырехугольник
$$
\xymatrix 
{
X\ar[rr]^{\sigma}\ar[rd]_{\e_{\ph}} &  & X'\ar@{-->}@/^6ex/[ldd]^{\ph'} \\
 & \Ran\e_{\ph}\ar[d]^{\mu_{\ph}} & \\
 &  Y  &
}
$$
В нем $\sigma\in\varOmega$, а $\mu_{\ph}\in\varOmega^\downarrow$. Значит, существует диагональ $\e'$ этого четырехугольника:
$$
\xymatrix 
{
X\ar[rr]^{\sigma}\ar[rd]_{\e_{\ph}} &  & X'\ar@{-->}@/^6ex/[ldd]^{\ph'}\ar@{-->}[ld]^{\e'} \\
 & \Ran\e_{\ph}\ar[d]^{\mu_{\ph}} & \\
 &  Y  &
}
$$
Здесь, в частности, верхний треугольник коммутативен, и, поскольку это верно для любого $\ph\in\varPhi$, это означает, что  $\sigma:X\to X'$ -- расширение $X$ в классе $\varOmega$ относительно морфизмов $\{\e_{\ph};\ \ph\in\varPhi\}$.
\epr

\medskip
\centerline{\bf Свойства оболочек в мономорфно дополняемых классах:}

\bit{\it

\item[] Пусть $\varOmega$ -- мономорфно дополняемый класс морфизмов в категории $\tt K$.

\item[$1^\circ$.]\label{1^0:obolochka-otn-1-morphizma}
Для всякого морфизма $\ph:X\to Y$ в ${\tt K}$ эпиморфизм $\e_{\ph}$ в его
факторизации $\ph=\mu_{\ph}\circ\e_{\ph}$ (определяемой классами ${^\downarrow\varOmega}$ и
$\varOmega$) является оболочкой объекта $X$ в $\varOmega$ относительно морфизма $\ph$:
 \begin{align}\label{obolochka-otn-morfizma}
 & \env_{\ph}^\varOmega X=\e_{\ph}, && \Env_{\ph}^\varOmega X=\Ran\e_{\ph}
 \end{align}

\item[$2^\circ$.]\label{2^0:obolochka-otn-konechnogo-mnozhestva-morphizmov}
Если ${\tt K}$ -- категория с конечными произведениями, то в ней любой объект $X$
обладает оболочкой в $\varOmega$ относительно произвольного конечного множества морфизмов $\varPhi$,
выходящего из $X$.

\item[$3^\circ$.]\label{3^0:obolochka-otn-mnozhestva-morphizmov} Если ${\tt K}$ -- категория
с произведениями\footnote{В утверждениях $3^\circ$-$5^\circ$ предполагается, что категория
${\tt K}$ обладает произведениями над произвольныи индексным множеством, необязательно, конечным.},
то в ней любой объект $X$ обладает оболочкой в $\varOmega$ относительно произвольного множества
морфизмов $\varPhi$, выходящего из $X$.

\item[$4^\circ$.]\label{4^0:obolochka-otn-klassa-morphizmov-porozhd-mnozhestvom}
Если ${\tt K}$ -- категория с произведениями, то в ней любой объект $X$ обладает
оболочкой в $\varOmega$ относительно произвольного класса морфизмов $\varPhi$,
выходящего из  $X$, и содержащего порождающее его изнутри подмножество морфизмов.

\item[$5^\circ$.]\label{5^0:obolochka-otn-klassa-morphizmov}
Если ${\tt K}$ -- категория с произведениями, локально малая в фактор-объектах
класса $\varOmega$, то в ${\tt K}$ любой объект $X$ обладает оболочкой в классе
$\varOmega$ относительно произвольного класса морфизмов $\varPhi$,  выходящего из  $X$.
}\eit

\bpr 1. Морфизм $\e_{\ph}$ есть расширение $X$ в $\varOmega$ относительно $\ph$, как видно из диаграммы
\beq\label{e_ph=obolochka}
\xymatrix 
{
X \ar[rr]^{\e_{\ph}}  \ar[dr]_{\ph} & & \Ran\e_{\ph}\ar@{-->}[dl]^{\mu_{\ph}}\\
& Y &
}
\eeq
Пусть $\sigma:X\to N$ -- какое-то другое расширение $X$  в $\varOmega$ относительно $\ph$:
$$
\xymatrix 
{
X \ar[rr]^{\sigma} \ar[dr]_{\ph} & & N \ar@{-->}[dl]^{\exists!\nu}
\\
& Y &
}
$$
Мы получаем коммутативную диаграмму
$$
\xymatrix @R=2.5pc @C=4.0pc
{
X\ar[rr]^{\ph}\ar[rd]^{\sigma}\ar@/_2ex/[rdd]_{\e_{\ph}} & & Y \\
& N\ar[ru]^{\nu} & \\
& \Ran\e_{\ph}\ar@/_2ex/[ruu]_{\mu_{\ph}} &
}
$$
В ней $\sigma\in\varOmega$ и $\mu_{\ph}\in{^\downarrow\varOmega}$, поэтому должна существовать диагональ нижнего четырехугольника:
$$
\xymatrix @R=2.5pc @C=4.0pc
{
X\ar[rr]^{\ph}\ar[rd]^{\sigma}\ar@/_2ex/[rdd]_{\e_{\ph}} & & Y \\
& N\ar[ru]^{\nu}\ar@{-->}[d]^{\upsilon} & \\
& \Ran\e_{\ph}\ar@/_2ex/[ruu]_{\mu_{\ph}} &
}
$$
Морфизм $\upsilon$ -- тот самый морфизм из диаграммы \eqref{DEF:diagr-obolochka}, который связывает расширение $\sigma$ с оболочкой $\e_{\ph}$. Его единственность следует из эпиморфности $\sigma$.

2. Пусть дан объект $X$ и конечное множество морфизмов $\varPhi$. Понятно, что в $\varPhi$ достаточно выделить подмножество $\varPhi^X=\{\ph:X\to Y_{\ph};\ \ph\in\varPhi^X\}$, состоящее из морфизмов, выходящих из $X$,
$$
\ph\in\varPhi^X\quad\Longleftrightarrow\quad \ph\in\varPhi\quad\&\quad \Dom\ph=X.
$$
Тогда оболочка относительно $\varPhi$ -- то же самое, что оболочка относительно $\varPhi^X$. Рассмотрим произведение объектов $\prod_{\ph\in\varPhi^X}Y_{\ph}$ и соответствующее произведение морфизмов $\prod_{\ph\in\varPhi^X}\ph:X\to \prod_{\ph\in\varPhi^X}Y_{\ph}$. Оболочка $X$ относительно множества морфизмов $\varPhi^X$ будет в точности оболочкой $X$ относительно одного морфизма $\prod_{\ph\in\varPhi^X}\ph$. Далее работает свойство $1^\circ$.

3. Пусть теперь ${\tt K}$ -- категория с произведениями над произвольным (необязательно, конечным) индексным множеством. Тогда предыдущие рассуждения проходят и для случая, когда $\varPhi$ представляет собой множество (необязательно, конечное) морфизмов.

4. Пусть $\varPsi\subseteq\varPhi$ -- подмножество (а не просто класс), порождающее $\varPhi$ изнутри. По уже доказанному свойству $3^\circ$, у любого объекта $X$ существует оболочка относительно $\varPsi$. С другой стороны, в силу \eqref{env_Psi=env_Phi}, эта оболочка должна совпадать с оболочкой относительно $\varPhi$.

5. Пусть ${\tt K}$ -- категория с произведениями (над произвольным индексным множеством), $A$ -- объект ${\tt K}$ и $\varPhi$ -- класс морфизмов (необязательно, множество). Идея доказательства состоит в том, чтобы заменить {\it класс} $\varPhi$ на некое {\it множество} морфизмов $M$, относительно которого оболочка будет той же, что и относительно $\varPhi$. Как и в пункте 2, мы можем считать, что $\varPhi$ состоит из морфизмов, выходящих из $X$:
$$
\forall\ph\in\varPhi\quad \Dom\ph=X.
$$
Для каждого $\ph\in\varPhi$ рассмотрим морфизм $\e_{\ph}$.
По лемме \ref{LM:env_varPhi^varOmega_X=env_(e_ph;ph_in_varPhi)^varOmega_X}, класс $\varPhi$ можно заменить на класс  $\{\e_{\ph};\ \ph\in\varPhi\}$:
$$
\env_{\varPhi}^{\Epi}X=\env_{\{\e_{\ph};\ \ph\in\varPhi\}}^{\Epi}X.
$$
После этого остается вспомнить, что все морфизмы $\e_{\ph}$ лежат в $\varOmega$, и,
поскольку наша категория локально мала в фактор-объектах класса $\varOmega$,
из них можно выделить множество $M$, такое, что любой эпиморфизм  $\e_{\ph}$
будет изоморфен некоторому эпиморфизму $\e\in M$, то есть $\e_{\ph}=\iota\circ\e$
для некоторого изоморфизма $\iota$. Множество $M$ теперь заменяет класс
$\{\e_{\ph};\ \ph\in\varPhi\}$ (а вместе с ним, и класс $\varPhi$), и дальше работает свойство $3^\circ$.
\epr

Двойственные результаты для детализаций выглядят так.

\blm\label{LM:imp_varPhi^varGamma_X=imp_(mu_ph;ph_in_varPhi)^varGamma_X}
Пусть $\varGamma$ -- эпиморфно дополняемый класс в категории ${\tt K}$. Тогда для всякого объекта $X$ и любого класса морфизмов $\varPhi$
\beq\label{imp_varPhi^varGamma_X=imp_(mu_ph;ph_in_varPhi)^varGamma_X}
\rf_{\varPhi}^{\varGamma}X=\rf_{\{\mu_{\ph};\ \ph\in \varPhi\}}^{\varGamma}X
\eeq
(если какая-то из этих детализаций существует, то другая тоже существует, и они совпадают).
\elm

\medskip
\centerline{\bf Свойства детализаций в эпиморфно дополняемых классах:}

\bit{\it

\item[] Пусть $\varGamma$ -- эпиморфно дополняемый класс морфизмов в категории $\tt K$.

\item[$1^\circ$.]\label{1^0:otpechatok-posr-1-morphizma}
Для всякого морфизма $\ph:X\gets Y$ в ${\tt K}$ мономорфизм $\mu_{\ph}$ в его
факторизации $\ph=\mu_{\ph}\circ\e_{\ph}$ (определяемой классами $\varGamma$ и
$\varGamma^\downarrow$) является детализацией $X$ в классе $\varGamma$ посредством морфизма $\ph$:
 \begin{align}\label{otpechatok-posr-morfizma}
 & \rf_{\ph}^\varGamma X=\mu_{\ph}, && \Rf_{\ph}^\varGamma X=\Dom\mu_{\ph}
 \end{align}

\item[$2^\circ$.]\label{2^0:otpechatok-posr-konechnogo-mnozhestva-morphizmov}
Если ${\tt K}$ -- категория с конечными копроизведениями, то в ней любой объект $X$
обладает детализацией в классе $\varGamma$ посредством произвольного конечного множества
морфизмов $\varPhi$, приходящего в $X$.

\item[$3^\circ$.]\label{3^0:otpechatok-posr-mnozhestva-morphizmov}
Если ${\tt K}$ -- категория с копроизведениями\footnote{\label{FOOT:lyuboe-I}В утверждениях $3^\circ$-$5^\circ$
предполагается, что категория ${\tt K}$ обладает копроизведениями над произвольныи индексным множеством,
необязательно, конечным.}, то в ней любой объект $X$ обладает детализацией
в классе $\varGamma$ посредством произвольного множества морфизмов $\varPhi$, приходящего в $X$.

\item[$4^\circ$.]\label{4^0:otpechatok-posr-klassa-morphizmov-porozhd-mnozhestvom} Если ${\tt K}$
-- категория с копроизведениями, то в ней любой объект $X$ обладает детализацией
в классе $\varGamma$ посредством произвольного класса морфизмов $\varPhi$, приходящего в $X$, и
содержащего порождающее его снаружи подмножество морфизмов.

\item[$5^\circ$.]\label{5^0:otpechatok-posr-klassa-morphizmov} Если ${\tt K}$ --
категория с копроизведениями, локально малая в подобъектах класса $\varGamma$, то в ${\tt K}$
любой объект $X$ обладает детализацией в классе $\varGamma$ посредством произвольного класса
морфизмов $\varPhi$, приходящего в $X$.
}\eit

\paragraph{Существование оболочек и детализаций в категориях с узловым разложением.} \label{obolochki-i-nachinki-otn-1-morfizma}

Общие свойства на с.\pageref{1^0:obolochka-otn-1-morphizma}
применительно к случаям $\varOmega=\Epi$ и $\varOmega=\SEpi$ приобретают следующий вид:

\medskip
\centerline{\bf Свойства оболочек в $\Epi$ и в $\SEpi$ в категории с узловым разложением:}

\bit{\it

\item[] Пусть ${\tt K}$ -- категория с узловым разложением.

\item[$1^\circ$.]\label{1^0:obolochka-env_ph^Epi} Для всякого морфизма $\ph:X\to Y$ в ${\tt K}$
 \bit{
 \item[---] эпиморфизм $\red_\infty\ph\circ\coim_\infty\ph$ в узловом разложении $\ph$ является оболочкой объекта $X$ в классе $\Epi$ всех эпиморфизмов относительно морфизма $\ph$:
 \begin{align}\label{obolochka--env_ph^Epi}
 & \env_{\ph}^{\Epi} X=\red_\infty\ph\circ\coim_\infty\ph, && \Env_{\ph}^{\Epi}X=\Im_\infty\ph
 \end{align}

 \item[---] эпиморфизм $\coim_\infty\ph$ в узловом разложении $\ph$ является оболочкой объекта $X$ в классе $\SEpi$ строгих эпиморфизмов относительно морфизма $\ph$:
 \begin{align}\label{obolochka--env_ph^SEpi}
 & \env_{\ph}^{\SEpi} X=\coim_\infty\ph, && \Env_{\ph}^{\SEpi}X=\Coim_\infty\ph
 \end{align}

}\eit

\item[$2^\circ$.]\label{2^0:obolochka-env_Phi^Epi-otn-konechnogo-mnozhestva-morphizmov} Если ${\tt K}$ -- категория с конечными произведениями, то в ней любой объект $X$ обладает оболочками в классах $\Epi$ и $\SEpi$ относительно произвольного конечного множества морфизмов $\varPhi$, выходящего из  $X$.

\item[$3^\circ$.]\label{3^0:obolochka-env_Phi^Epi-otn-mnozhestva-morphizmov} Если ${\tt K}$ --
категория с произведениями\footnote{В этом списке свойств мы делаем те же предположения,
что в сноске \ref{FOOT:lyuboe-I}.},
то в ней любой объект $X$ обладает оболочками в классах $\Epi$ и $\SEpi$ относительно произвольного множества
морфизмов $\varPhi$, выходящего из  $X$.

\item[$4^\circ$.]\label{4^0:obolochka-env_Phi^Epi-otn-klassa-morphizmov-porozhd-mnozhestvom} Если ${\tt K}$ -- категория с произведениями, то в ней любой объект $X$ обладает оболочками в классах $\Epi$ и $\SEpi$ относительно произвольного класса морфизмов $\varPhi$, выходящего из  $X$, и содержащего порождающее его изнутри подмножество морфизмов.

\item[$5^\circ$.]\label{5^0:obolochka-env_Phi^Epi-otn-klassa-morphizmov}
Если ${\tt K}$ -- категория с произведениями, локально малая в фактор-объектах класса
$\Epi$ (соответственно, класса $\SEpi$), то в ${\tt K}$ любой объект $X$ обладает оболочкой в
классе $\Epi$  (соответственно, в классе $\SEpi$) относительно произвольного класса морфизмов
$\varPhi$,  выходящего из  $X$.
}\eit

\bprop\label{PROP:obolochka-v-Omega,Bim,Epi}
Если ${\tt K}$ -- категория с произведениями, узловым разложением и
локально малая в фактор-объектах класса $\Epi$, то в ${\tt K}$ любой объект
$X$ обладает оболочкой в любом классе $\varOmega$ содержащим биморфизмы,
    $$
    \varOmega\supseteq\Bim,
    $$
относительно произвольного правого идеала морфизмов $\varPhi$, различающего морфизмы
снаружи\footnote{См. определение на с.\pageref{DEF:varPhi-razlich-morfizmy-snaruzhi}.}
и выходящего из  $X$, причем оболочка $X$ в $\varOmega$ относительно $\varPhi$
совпадает с оболочками в $\Bim$ и в $\Epi$ относительно $\varPhi$:
$$
\env_{\varPhi}^{\varOmega}X=\env_{\varPhi}^{\Bim}X=\env_{\varPhi}^{\Epi}X.
$$
\eprop
\bpr
По свойству $5^\circ$, существует оболочка $\env_{\varPhi}^{\Epi}X$. По теореме \ref{TH:Phi-razdel-moprfizmy}(i) эта оболочка является мономорфизмом, и значит, биморфизмом: $\env_{\varPhi}^{\Epi}X\in\Bim$. Отсюда по свойству $1^\circ$(c) на с.\pageref{LM:suzhenie-verh-klassa-morfizmov} оболочка в $\Epi$ должна быть оболочкой в $\Bim$:
$\env_{\varPhi}^{\Epi}X=\env_{\varPhi}^{\Bim}X$. Теперь по теореме \ref{TH:Phi-razdel-moprfizmy-*} оболочка в $\Bim$ должна быть оболочкой в $\varOmega$: $\env_{\varPhi}^{\Bim}X=\env_{\varPhi}^{\varOmega}X$.

\epr

Двойственные результаты для детализаций выглядят так.

\medskip
\centerline{\bf Свойства детализаций в $\Mono$ и $\SMono$ в категории с узловым разложением:}

\bit{\it

\item[] Пусть ${\tt K}$ -- категория с узловым разложением.

\item[$1^\circ$.]\label{1^0:otpechatok-imp_ph^Mono} Для всякого морфизма $\ph:X\gets Y$ в ${\tt K}$
 \bit{
 \item[---] мономорфизм $\im_\infty\ph\circ\red_\infty\ph$ в узловом разложении $\ph$ является детализацией объекта $X$ в классе $\Mono$ всех мономорфизмов посредством морфизма $\ph$:
 \begin{align}\label{otpechatok-imp_ph^Mono}
 & \rf_{\ph}^{\Mono} X=\im_\infty\ph\circ\red_\infty\ph, && \Rf_{\ph}^{\Mono}X=\Coim_\infty\ph
 \end{align}

 \item[---] мономорфизм $\im_\infty\ph$ в узловом разложении $\ph$ является детализацией объекта $X$ в классе $\SMono$ строгих мономорфизмов посредством морфизма $\ph$:
 \begin{align}\label{otpechatok-imp_ph^SMono}
 & \rf_{\ph}^{\SMono} X=\im_\infty\ph, && \Rf_{\ph}^{\SMono}X=\Im_\infty\ph
 \end{align}

}\eit

\item[$2^\circ$.]\label{2^0:otpechatok-imp_Phi^Mono-posr-konechnogo-mnozhestva-morphizmov} Если ${\tt K}$ -- категория с конечными копроизведениями, то в ней любой объект $X$ обладает детализациями в классах $\Mono$ и $\SMono$ посредством  произвольного конечного множества морфизмов $\varPhi$, приходящего в $X$.

\item[$3^\circ$.]\label{3^0:otpechatok-imp_Phi^Mono-posr-mnozhestva-morphizmov} Если ${\tt K}$ -- категория
с копроизведениями\footnote{В утвреждениях $3^\circ$-$5^\circ$ предполагается, что категория
${\tt K}$ обладает копроизведениями над произвольныи индексным множеством, необязательно, конечным.},
то в ней любой объект $X$ обладает детализациями в классах $\Mono$ и $\SMono$ посредством произвольного множества морфизмов $\varPhi$, приходящего в $X$.

\item[$4^\circ$.]\label{4^0:otpechatok-imp_Phi^Mono-posr-klassa-morphizmov-porozhd-mnozhestvom} Если ${\tt K}$ -- категория с копроизведениями, то в ней любой объект $X$ обладает детализациями в классах $\Mono$ и $\SMono$ посредством произвольного класса морфизмов $\varPhi$, приходящего в $X$, и содержащего порождающее его снаружи подмножество морфизмов.

\item[$5^\circ$.]\label{5^0:otpechatok-imp_Phi^Mono-posr-klassa-morphizmov} Если ${\tt K}$ -- категория с копроизведениями, локально малая в подобъектах класса $\Mono$ (соответственно, класса $\SMono$), то в ${\tt K}$ любой объект $X$ обладает детализацией в классе $\Mono$  (соответственно, в классе $\SMono$) посредством произвольного класса морфизмов $\varPhi$, приходящего в $X$.
}\eit

\bprop\label{PROP:otpechatok-v-Omega,Bim,Epi} Если ${\tt K}$ -- категория с копроизведениями,
узловым разложением и локально малая в фактор-объектах класса $\Mono$, то в ${\tt K}$ любой
объект $X$ обладает детализацией в классе $\varGamma$, содержащим биморфизмы,
    $$
    \varGamma\supseteq\Bim,
    $$
посредством произвольного левого идеала морфизмов $\varPhi$, различающего морфизмы
изнутри\footnote{См. определение на с.\pageref{DEF:varPhi-razlich-morfizmy-iznutri}.}
и приходящего в $X$, причем детализация $X$ в $\varGamma$ посредством $\varPhi$ совпадает с детализациями в $\Bim$ и в $\Mono$ посредством $\varPhi$:
$$
\rf_{\varPhi}^{\varGamma}X=\rf_{\varPhi}^{\Bim}X=\rf_{\varPhi}^{\Mono}X.
$$
\eprop

\paragraph{Существование узлового разложения в категориях с оболочками и детализациями.}
Условимся говорить, что в категории ${\tt K}$
 \bit{
 \item[---] {\it эпиморфизмы распознают мономорфизмы}, если из того, что морфизм $\mu$ не является мономорфизмом следует, что его можно разложить в композицию $\mu=\mu'\circ\e$, в которой $\e$ -- эпиморфизм, не являющийся изоморфизмом,

 \item[---] {\it мономорфизмы распознают эпиморфизмы}, если из того, что морфизм $\e$ не является эпиморфизмом следует, что его можно разложить в композицию $\e=\mu\circ\e'$, в которой $\mu$ -- мономорфизм, не являющийся изоморфизмом,

\item[---]\label{DEF:strog-epi-razlich-mono} {\it строгие эпиморфизмы распознают мономорфизмы}, если из того, что морфизм $\mu$ не является мономорфизмом следует, что его можно разложить в композицию $\mu=\mu'\circ\e$, в которой $\e$ -- строгий эпиморфизм, не являющийся изоморфизмом,

\item[---]\label{DEF:strog-epi-razlich-mono} {\it строгие мономорфизмы распознают эпиморфизмы}, если из того, что морфизм $\e$ не является эпиморфизмом следует, что его можно разложить в композицию $\e=\mu\circ\e'$, в которой $\mu$ -- строгий мономорфизм, не являющийся изоморфизмом.

 }\eit

\btm\label{TH:env+imp=>uzl=razl} Пусть в категории ${\tt K}$
 \bit{
 \item[\rm (a)] эпиморфизмы распознают мономорфизмы, а мономорфизмы распознают эпиморфизмы,

 \item[\rm (b)]  всякий непосредственный мономорфизм является строгим мономорфизмом, и, двойственным образом, всякий непосредственный эпиморфизм является строгим эпиморфизмом,

 \item[\rm (c)]  всякий объект $X$ обладает оболочкой в классе $\Epi$ всех эпиморфизмов относительно произвольного выходящего из него морфизма, и, двойственным образом, любой объект $X$ обладает детализацией в классе $\Mono$ всех мономорфизмов посредством произвольного приходящего в него морфизма.
   }\eit
Тогда ${\tt K}$ -- категория с узловым разложением.
 \etm
\bpr Рассмотрим морфизм $\ph:X\to Y$.

1. Пусть $\e:X\to N$ -- оболочка в классе  $\Epi$ объекта $X$ относительно $\ph$, и пусть пунктирная стрелка, замыкающая диаграмму \eqref{DEF:diagr-rasshirenie}, обозначается  $\beta$:
$$
\ph=\beta\circ\e
$$
Заметим сначала, что $\beta$ -- мономорфизм. Действительно, если бы $\beta$ не был мономорфизмом, то в силу условия (a) существовало бы разложение
$\beta=\beta'\circ\pi$, в котором $\pi$ -- эпиморфизм, не являющийся изоморфизмом. Если через $N'$ обозначить область значений $\pi$, то мы получим диаграмму
\beq\label{PROOF:sushestv-uzlov-razlozhenija}
\xymatrix @R=4.0pc @C=4.0pc
{
X\ar[d]_{\e}\ar[r]^{\ph}\ar[dr]^(.7){\e'}|!{[d];[r]}\hole & Y
\\
N\ar[r]_{\pi}\ar[ur]^(.7){\beta} & N'\ar[u]_{\beta'}
}
\eeq
в которой по определению $\e'=\pi\circ\e$, причем это будет эпиморфизм, как композиция двух эпиморфизмов. Получается,  что $\e'$ -- другое расширение $X$ относительно $\ph$. Значит, существует единственный морфизм $\upsilon$, замыкающий диаграмму
$$
\begin{diagram}\dgARROWLENGTH=3em
\node[2]{X}\arrow{sw,t}{\e} \arrow{se,t}{\e'} \\
\node{N}  \node[2]{N'}\arrow[2]{w,b,--}{\upsilon}
\end{diagram}
$$
При этом получается:
$$
\pi\circ\e=\e'\quad\Longrightarrow\quad\upsilon\circ\pi\circ\e=\upsilon\circ\e'=\e=1_N\circ\e\quad\Longrightarrow\quad \upsilon\circ\pi=1_N
$$
и
$$
\upsilon\circ\e'=\e\quad\Longrightarrow\quad\pi\circ\upsilon\circ\e'=\pi\circ\e=\e'=1_{N'}\circ\e'\quad\Longrightarrow\quad \pi\circ\upsilon=1_{N'}.
$$
То есть $\pi$ должен быть изоморфизмом, а это противоречит договоренности, что $\pi$ -- не изоморфизм.

2. Точно так же доказывается, что $\beta$ -- непосредственный мономорфизм. Действительно, если рассмотреть какую-нибудь его факторизацию $\beta=\beta'\circ\pi$, то мы опять приходим к диаграмме \eqref{PROOF:sushestv-uzlov-razlozhenija}, из которой теми же рассуждениями выводится, что $\pi$ должен быть изоморфизмом.

3. Из того, что $\beta$ -- непосредственный мономорфизм и условия (b) следует, что $\beta$ -- строгий мономорфизм.

4. Обозначим далее через $\mu:M\to Y$ детализацию объекта $Y$ в классе $\Mono$ посредством морфизма $\ph$, а через $\alpha$ -- пунктирную стрелку в соответствующей диаграмме \eqref{DEF:suzhenie}, то есть
$$
\ph=\mu\circ\alpha.
$$
По аналогии с тем, как выше доказывалось, что $\beta$ -- строгий мономорфизм, в этом случае нетрудно показать, что $\alpha$ -- строгий эпиморфизм.

5. Теперь рассмотрим диаграмму
$$
\begin{diagram}\dgARROWLENGTH=3em
\node{X}\arrow[2]{r,t}{\ph}\arrow[2]{s,t}{\alpha}\arrow{se,t,-}{} \node[2]{Y} \\
 \node[2]{} \arrow{se,t}{\e} \\
\node{M}\arrow[2]{ne,t,3}{\mu}  \node[2]{N}\arrow[2]{n,r}{\beta}
\end{diagram}
$$
Поскольку, как мы поняли, $\alpha$ -- эпиморфизм, мы получаем, что $\alpha$ -- расширение в классе $\Epi$ объекта $X$ относительно $\ph$. Но с другой стороны, $\e$ -- оболочка в классе $\Epi$ объекта $X$ относительно $\ph$. Значит, должен существовать морфизм $\upsilon$, замыкающий диаграмму
$$
\begin{diagram}\dgARROWLENGTH=3em
\node{X}\arrow{s,t}{\alpha} \arrow{se,t}{\e} \\
\node{M}\arrow{e,b}{\upsilon}  \node{N}
\end{diagram}
$$
и поэтому диаграмму
\beq\label{PROOF:sushestv-uzlov-razlozhenija-2}
\begin{diagram}\dgARROWLENGTH=3em
\node{X}\arrow{e,t}{\ph}\arrow{s,t}{\alpha} \arrow{se,t}{\e}\node{Y} \\
\node{M}\arrow{e,b}{\upsilon}  \node{N}\arrow{n,r}{\beta}
\end{diagram}
\eeq
Точно так же, поскольку $\beta$ -- мономорфизм, он должен быть обогащением объекта $Y$ в классе $\Mono$ посредством $\ph$. С другой стороны, $\mu$ -- детализация объекта $Y$ в классе $\Mono$ посредством морфизма $\ph$. Значит, должен существовать морфизм $\upsilon'$, замыкающий диаграмму
$$
\begin{diagram}\dgARROWLENGTH=3em
\node[2]{Y} \\
\node{M}\arrow{ne,t}{\mu}\arrow{e,b}{\upsilon'}  \node{N}\arrow{n,r}{\beta}
\end{diagram}
$$
и поэтому диаграмму
\beq\label{PROOF:sushestv-uzlov-razlozhenija-3}
\begin{diagram}\dgARROWLENGTH=3em
\node{X}\arrow{e,t}{\ph}\arrow{s,t}{\alpha} \node{Y} \\
\node{M}\arrow{e,b}{\upsilon'}\arrow{ne,t}{\mu}  \node{N}\arrow{n,r}{\beta}
\end{diagram}
\eeq
Из \eqref{PROOF:sushestv-uzlov-razlozhenija-2} и \eqref{PROOF:sushestv-uzlov-razlozhenija-3} получаем:
$$
\underset{\scriptsize\begin{matrix}\text{\rotatebox{90}{$\owns$}}\\ \Mono\end{matrix}}{\beta}\circ\upsilon\circ\underset{\scriptsize\begin{matrix}\text{\rotatebox{90}{$\owns$}}\\ \Epi\end{matrix}}{\alpha}=\ph=\underset{\scriptsize\begin{matrix}\text{\rotatebox{90}{$\owns$}}\\ \Mono\end{matrix}}{\beta}\circ\upsilon'\circ\underset{\scriptsize\begin{matrix}\text{\rotatebox{90}{$\owns$}}\\ \Epi\end{matrix}}{\alpha}
\quad\Longrightarrow\quad
\upsilon=\upsilon'
$$
То есть должна быть коммутативна диаграмма
$$
\begin{diagram}\dgARROWLENGTH=3em
\node{X}\arrow[2]{r,t}{\ph}\arrow[2]{s,t}{\alpha}\arrow{se,t,-}{} \node[2]{Y} \\
 \node[2]{} \arrow{se,t}{\e} \\
\node{M}\arrow[2]{ne,t,3}{\mu}\arrow[2]{e,b}{\upsilon}  \node[2]{N}\arrow[2]{n,r}{\beta}
\end{diagram}
$$
В ней, $\e=\upsilon\circ\alpha$ -- эпиморфизм, поэтому $\upsilon$ тоже должен быть эпиморфизмом. С другой стороны, $\mu=\beta\circ\upsilon$ -- мономорфизм, и поэтому  $\upsilon$ должен быть также и мономорфизмом. Значит, $\upsilon$ -- биморфизм, и $\ph=\beta\circ\upsilon\circ\alpha$ -- узловое разложение для $\ph$.
\epr

\btm\label{TH:env^SEpi=>uzlo-razl} Пусть в категории ${\tt K}$
 \bit{
 \item[\rm (a)] строгие эпиморфизмы распознают мономорфизмы, а строгие мономорфизмы распознают эпиморфизмы\footnote{См.
 определения на с.\pageref{DEF:strog-epi-razlich-mono}.},

 \item[\rm (b)]  всякий объект $X$ обладает оболочкой в классе  $\SEpi$ всех строгих эпиморфизмов относительно произвольного выходящего из него морфизма, и, двойственным образом, любой объект $X$ обладает детализацией в классе $\SMono$ всех строгих мономорфизмов посредством произвольного приходящего в него морфизма.
   }\eit
Тогда ${\tt K}$ -- категория с узловым разложением.
 \etm
\bpr Рассмотрим морфизм $\ph:X\to Y$.

1. В силу условия (b), существует оболочка $\env_{\ph}^{\SEpi}X:X\to \Env_{\ph}^{\SEpi}X$ объекта $X$ в классе  $\SEpi$ всех строгих эпиморфизмов относительно морфизма $\ph$. Обозначим через $\alpha$ морфизм, продолжающий $\ph$ на $\Env_{\ph}^{\SEpi}X$:
$$
\xymatrix 
@C=4.0pc
{
X\ar[d]_{\env_{\ph}^{\SEpi}X}\ar[rd]^{\ph} &
\\
\Env_{\ph}^{\SEpi}X\ar@{-->}[r]_{\alpha} & Y
}
$$

2. Точно так же в силу (b) существует детализация $\rf_{\ph}^{\SMono}Y:\Rf_{\ph}^{\SMono}Y\to Y$ в классе $\SMono$ всех строгих мономорфизмов в объекте $Y$ посредством морфизма $\ph$. Обозначим через $\beta$ морфизм, поднимающий $\ph$ на $\Rf_{\ph}^{\SMono}X$:
$$
\xymatrix 
 @C=4.0pc
{
X\ar[rd]_{\ph}\ar@{-->}[r]^{\beta} & \Rf_{\ph}^{\SEpi}Y\ar[d]^{\rf_{\ph}^{\SMono}Y}
\\
 & Y
}
$$

3. Склеив эти треугольники по общей стороне $\ph$, и выбросив ее, мы получим четырехугольник:
$$
\xymatrix 
@C=4.0pc
{
X\ar[d]_{\env_{\ph}^{\SEpi}X}\ar[r]^{\beta} & \Rf_{\ph}^{\SEpi}Y\ar[d]^{\rf_{\ph}^{\SMono}Y}
\\
\Env_{\ph}^{\SEpi}X\ar[r]_{\alpha} & Y
}
$$
В нем $\env_{\ph}^{\SEpi}X$ -- строгий эпиморфизм, а $\rf_{\ph}^{\SMono}Y$ -- мономорфизм, поэтому существует диагональ $\delta$
\beq\label{PROOF:env^SEpi=>uzlo-razl}
\xymatrix 
@C=4.0pc
{
X\ar[d]_{\env_{\ph}^{\SEpi}X}\ar[r]^{\beta} & \Rf_{\ph}^{\SEpi}Y\ar[d]^{\rf_{\ph}^{\SMono}Y}
\\
\Env_{\ph}^{\SEpi}X\ar[r]_{\alpha}\ar@{-->}[ru]_{\delta} & Y
}
\eeq
Покажем, что $\delta$ -- биморфизм.

4. Предположим сначала, что $\delta$ -- не мономорфизм. Тогда, поскольку строгие эпиморфизмы распознают мономорфизмы (условие (a)), существует разложение $\delta=\delta'\circ\e$, в котором $\e$ -- строгий эпиморфизм, не являющийся изоморфизмом. Как следствие, будет коммутативна диаграмма:
$$
\xymatrix 
@C=4.0pc
{
X\ar[d]_{\env_{\ph}^{\SEpi}X}\ar[rd]_{\beta}\ar[r]^{\ph} & Y
\\
\Env_{\ph}^{\SEpi}X\ar[r]_{\delta}\ar@{-->}[d]_{\e} & \Rf_{\ph}^{\SEpi}Y\ar[u]_{\rf_{\ph}^{\SMono}Y} \\
M\ar@/_2ex/@{-->}[ru]_{\delta'}
}
$$
Из нее видно, что композиция $\rf_{\ph}^{\SMono}Y\circ\delta'$ является продолжением морфизма $\ph$ вдоль композиции $\e\circ\env_{\ph}^{\SEpi}X$, которая в свою очередь является строгим эпиморфизмом (как композиция двух строгих эпиморфизмов). Это означает, что $\e\circ\env_{\ph}^{\SEpi}X$ есть расширением объекта $X$ в классе $\SEpi$ относительно морфизма $\ph$. Поэтому существует морфизм $\upsilon$ из расширения $M$ в оболочку $\Env_{\ph}^{\SEpi}X$, замыкающий диаграмму \eqref{DEF:diagr-obolochka}:
$$
\xymatrix 
@C=4.0pc
{
 & X\ar@/_2ex/[dl]_{\env_{\ph}^{\SEpi}X}\ar@/^2ex/[rd]^{\e\circ\env_{\ph}^{\SEpi}X} &
\\
\Env_{\ph}^{\SEpi}X & & M\ar@{-->}[ll]_{\upsilon} \\
}
$$
Для него мы получаем $\upsilon\circ\e\circ\env_{\ph}^{\SEpi}X=\env_{\ph}^{\SEpi}X=1_M\circ\env_{\ph}^{\SEpi}X$, откуда, в силу эпиморфности $\env_{\ph}^{\SEpi}X$, следует равенство $\upsilon\circ\e=1_M$, означающее, что $\e$ -- коретракция. С другой стороны, это эпиморфизм, и вместе это означает, что $\e$ должен быть изоморфизмом. Это противоречит выбору $\e$.

5. Таким образом, $\delta$ обязан быть мономорфизмом. Аналогично доказывается, что это эпиморфизм. Теперь добавим в диаграмму \eqref{PROOF:env^SEpi=>uzlo-razl} морфизм $\ph$ и скрутим ее, представив в следующем виде:
$$
\xymatrix @R=4.0pc @C=4.0pc
{
X\ar[d]_{\env_{\ph}^{\SEpi}X}\ar[r]^{\ph}\ar[dr]^(.7){\beta}|!{[d];[r]}\hole & Y
\\
\Env_{\ph}^{\SEpi}X\ar[r]_{\delta}\ar[ur]^(.7){\alpha} & \Rf_{\ph}^{\SEpi}Y\ar[u]_{\rf_{\ph}^{\SMono}Y}
}
$$
Видно, что $\ph=\rf_{\ph}^{\SMono}Y\circ\delta\circ\env_{\ph}^{\SEpi}X$ -- узловое разложение морфизма $\ph$.
\epr

\section{Сети и функториальность}

В общем случае операции взятия оболочки и детализации, по-видимому, не обязаны быть функторами. Однако при определенных предположениях это так, и в этом параграфе мы обсудим эти условия.

\subsection{Оболочка и детализация как функторы}

Пусть $\varOmega$ и $\varPhi$ -- классы морфизмов в категории ${\tt K}$.

\bit{
    \item
    Будем говорить, что {\it оболочку $\Env^\varOmega_\varPhi$ можно определить как функтор}, если существуют
        \bit{
    \item[E.1] отображение $X\mapsto (E(X),e_X)$, сопоставляющее каждому объекту $X$ в $\tt K$ некоторый морфизм   $e_X:X\to E(X)$ в ${\tt K}$, являющийся оболочкой в $\varOmega$ относительно $\varPhi$:
$$
E(X)=\Env_{\varPhi}^\varOmega X,\qquad e_X=\env_{\varPhi}^\varOmega X
$$
    \item[E.2] отображение $\alpha\mapsto E(\alpha)$, которое всякому морфизму $\alpha:X\to Y$ в $\tt K$ ставит в соответствие морфизм $E(\alpha):E(X)\to E(Y)$ в $\tt K$, замыкающий диаграмму
    \beq\label{DIAGR:funktorialnost-env-e-E}
\xymatrix @R=2.pc @C=5.0pc 
{
X\ar[d]^{\alpha}\ar[r]^{e_X} & E(X)\ar@{-->}[d]^{E(\alpha)} \\
Y\ar[r]^{e_Y} & E(Y) \\
}
\eeq
}\eit
и при этом выполняются условия
\beq\label{tozhdestva:funktorialnost-env-e-E}
E(1_X)=1_{E(X)},\qquad E(\beta\circ\alpha)=E(\beta)\circ E(\alpha)
\eeq
Понятно, что в этом случае отображение $(X,\alpha)\mapsto(E(X),E(\alpha))$ является ковариантным функтором из ${\tt K}$ в ${\tt K}$, а отображение $X\mapsto e_X$ представляет собой естественное преобразование тождественного функтора  $(X,\alpha)\mapsto(X,\alpha)$ в функтор $(X,\alpha)\mapsto(E(X),E(\alpha))$.

    \item[$\bullet$]  Условимся также говорить, что {\it оболочку $\Env^\varOmega_\varPhi$ можно определить как идемпотентный функтор}, если в дополнение к E.1 и E.2 можно добиться условия

 \bit{
 \item[E.3] для всякого объекта $X\in\Ob({\tt K})$ морфизм $e_{E(X)}:E(X)\to E(E(X))$ представляет собой локальную единицу:
\beq\label{e_(E(X))=1_(E(X))}
E(E(X))=E(X),\qquad e_{E(X)}=1_{E(X)}\qquad X\in \Ob({\tt K}).
\eeq
 }\eit
 }\eit

\brem\label{REM:E(e_X)=1_(E(X))}
Если $\varOmega\subseteq\Epi$, то из условия \eqref{e_(E(X))=1_(E(X))} следует условие
\beq\label{E(e_X)=1_(E(X))}
E(e_X)=1_{E(X)}\qquad X\in \Ob({\tt K}).
\eeq
Действительно, подставив в \eqref{DIAGR:funktorialnost-env-e-E} морфизм $\alpha=e_X$, мы получим:
$$
\xymatrix @R=2.pc @C=8.0pc 
{
X\ar[d]^{e_X}\ar[r]^{e_X} & E(X)\ar[d]^{E(e_X)} \\
E(X)\ar[r]^{e_{E(X)}=1_{E(X)}} & E(E(X))=E(X) \\
}
$$
то есть $E(e_X)\circ e_X=1_{E(X)}\circ e_X$, и, поскольку $e_X\in\varOmega\subseteq\Epi$, на этот морфизм можно сократить: $E(e_X)=1_{E(X)}$.
\erem

   \bit{
    \item
Будем говорить, что {\it детализацию $\Rf^\varGamma_\varPhi$ можно определить как функтор}, если существуют
        \bit{
    \item[R.1] отображение $X\mapsto (I(X),i_X)$, сопоставляющее каждому объекту $X$ в $\tt K$ некоторый морфизм   $i_X:I(X)\to X$ в ${\tt K}$, являющийся детализацией в классе $\varGamma$ посредством класса $\varPhi$:
$$
I(X)=\Rf_{\varPhi}^\varGamma X,\qquad i_X=\rf_{\varPhi}^\varGamma X
$$
\item[R.2] отображение $\alpha\mapsto I(\alpha)$, которое всякому морфизму $\alpha:X\gets Y$ в $\tt K$
ставит в соответствие морфизм $I(\alpha):I(X)\gets I(Y)$ в $\tt K$, замыкающий диаграмму
    \beq\label{DIAGR:funktorialnost-imp-i-I}
\xymatrix @R=2.pc @C=5.0pc 
{
X & I(X)\ar[l]_{i_X} \\
Y\ar[u]_{\alpha} & I(Y)\ar[l]_{i_Y}\ar@{-->}[u]_{I(\alpha)} \\
}
\eeq
}\eit
и при этом выполняются условия
\beq\label{tozhdestva:funktorialnost-imp-i-I}
I(1_X)=1_{I(X)},\qquad I(\beta\circ\alpha)=I(\beta)\circ I(\alpha)
\eeq
В этом случае отображение $(X,\alpha)\mapsto(I(X),I(\alpha))$ является ковариантным функтором из ${\tt K}$ в ${\tt K}$, а отображение $X\mapsto i_X$ представляет собой естественное преобразование тождественного функтора  $(X,\alpha)\mapsto(X,\alpha)$ в функтор $(X,\alpha)\mapsto(I(X),I(\alpha))$.

     \item[$\bullet$]  Условимся также говорить, что {\it детализацию $\Rf^\varGamma_\varPhi$ можно определить как идемпотентный функтор}, если в дополнение к R.1 и R.2 можно добиться условия

 \bit{
 \item[R.3] для всякого объекта $X\in\Ob({\tt K})$ морфизм $i_{I(X)}:I(X)\gets I(I(X))$ представляет собой локальную единицу:
\beq\label{i_(I(X))=1_(I(X))}
I(I(X))=I(X),\qquad i_{I(X)}=1_{I(X)}\qquad X\in \Ob({\tt K}).
\eeq
 }\eit
 }\eit

\brem\label{REM:I(i_X)=1_(I(X))}
Если $\varGamma\subseteq\Mono$, то из условия \eqref{i_(I(X))=1_(I(X))} следует условие
\beq\label{I(i_X)=1_(I(X))}
I(i_X)=1_{I(X)}\qquad X\in \Ob({\tt K}).
\eeq
\erem

\brem Для оболочек в наиболее важных для нас случаях, когда класс
$\varOmega$ состоит из эпиморфизмов,  $\varOmega\subseteq\Epi$, тождества
\eqref{tozhdestva:funktorialnost-env-e-E} автоматически следуют из условий E.1 и E.2.
Двойственным образом, для детализаций, в случаях, когда класс $\varGamma$ состоит из мономорфизмов,
 $\varGamma\subseteq\Mono$, тождества \eqref{tozhdestva:funktorialnost-imp-i-I} автоматически
  следуют из условий R.1 и R.2.
\erem

\subsection{Сети эпиморфизмов и полурегулярная оболочка}

Пусть $\Epi^X$ --- класс всех эпиморфизмов категории ${\tt K}$, выходящих из $X$. Мы считаем его наделенным предпорядком
$\to$, определяемым правилом
 \beq\label{DEF:le-in-F_X}
\rho\to\sigma\qquad \Leftrightarrow\qquad \exists \iota\in\Mor({\tt K})\quad
\sigma=\iota\circ\rho.
 \eeq
Здесь морфизм $\iota$, если он существует, должен быть единственным и, кроме того, он будет эпиморфизмом (это
следует из того, что $\sigma$ --- эпиморфизм). Поэтому определена операция,
которая каждой паре морфизмов $\rho,\sigma\in\varOmega^X$, удовлетворяющей условию
$\rho\to\sigma$, ставит в соответствие морфизм $\iota=\iota^\sigma_\rho$ в
\eqref{DEF:le-in-F_X}:
 \beq\label{DEF:le-in-F_X-*}
\sigma=\iota^\sigma_\rho\circ\rho.
 \eeq
При этом, если $\pi\to\rho\to\sigma$, то из цепочки
$$
\iota^\sigma_\pi\circ\pi=\sigma=\iota^\sigma_\rho\circ\rho=\iota^\sigma_\rho\circ\iota^\rho_\pi\circ\pi,
$$
опять в силу эпиморфности $\pi$, следует равенство
\beq\label{iota_rho^tau=iota_rho^sigma-circ-iota_sigma^tau}
\iota^\sigma_\pi=\iota^\sigma_\rho\circ\iota^\rho_\pi.
 \eeq

\paragraph{Сети эпиморфизмов.}

\bit{ \item[$\bullet$] Пусть каждому объекту $X\in\Ob({\tt K})$ категории ${\tt
K}$ поставлено в соответствие некое подмножество ${\mathcal N}^X$ в классе
$\Epi^X$ всех эпиморфизмов категории ${\tt K}$, выходящих из $X$, причем
выполняются следующие три условия:
    \bit{
\item[(a)]\label{AX:set-Epi-a} для всякого объекта $X$ множество ${\mathcal
N}^X$ непусто и направлено влево относительно предпорядка
\eqref{DEF:le-in-F_X}, наследуемого из $\Epi^X$:
$$
\forall \sigma,\sigma'\in {\mathcal N}^X\quad \exists\rho\in{\mathcal N}^X\quad
\rho\to\sigma\ \& \ \rho\to\sigma',
$$

\item[(b)]\label{AX:set-Epi-b} для всякого объекта $X$ порождаемая множеством
${\mathcal N}^X$ ковариантная система морфизмов
 \beq\label{DEF:kategoriya-svyazyv-morpfizmov}
\Bind({\mathcal N}^X):=\{\iota_\rho^\sigma;\ \rho,\sigma\in{\mathcal N}^X,\
\rho\to\sigma\}
 \eeq
(морфизмы $\iota_\rho^\sigma$ были определены в \eqref{DEF:le-in-F_X-*};
согласно \eqref{iota_rho^tau=iota_rho^sigma-circ-iota_sigma^tau}, эта система
является ковариантным функтором из множества ${\mathcal N}^X$, рассматриваемого
как полная подкатегория в $\Epi^X$, в $\tt K$) обладает проективным пределом в
$\tt K$;

\item[(c)]\label{AX:set-Epi-c}  для всякого морфизма $\alpha:X\to Y$ и любого
элемента $\tau\in{\mathcal N}^Y$ найдется элемент $\sigma\in{\mathcal N}^X$ и
морфизм $\alpha_\sigma^\tau:\Ran\sigma\to\Ran\tau$ такие, что коммутативна
диаграмма
 \beq\label{DIAGR:set} \xymatrix @R=2.5pc @C=4.0pc {
 X\ar[r]^{\alpha}\ar@{-->}[d]_{\sigma} & Y\ar[d]^{\tau} \\
 \Ran\sigma\ar@{-->}[r]_{\alpha_\sigma^\tau} & \Ran\tau
 } \eeq
(при фиксированных $\alpha$, $\sigma$ и $\tau$ морфизм $\alpha_\sigma^\tau$,
если он существует, определен однозначно, в силу эпиморфности $\sigma$).

 }\eit
Тогда
 \bit{
\item[---] семейство множеств ${\mathcal N}=\{{\mathcal N}^X;\ X\in\Ob({\tt
K})\}$ мы будем называть {\it сетью эпиморфизмов}\label{DEF:set-epimorf} в
категории ${\tt K}$, а элементы множеств ${\mathcal N}^X$ -- {\it элементами
сети} ${\mathcal N}$,

\item[---] для всякого объекта $X$ систему морфизмов $\Bind({\mathcal N}^X)$,
определенную равенствами \eqref{DEF:kategoriya-svyazyv-morpfizmov}, мы называем
{\it системой связывающих морфизмов сети ${\mathcal N}$ над вершиной $X$}, ее
проективный предел (существование которого гарантируется условием (b))
представляет собой проективный конус, вершину которого мы будем обозначать
символом $X_{\mathcal N}$, а исходящие из нее морфизмы -- символами
$\sigma_{\mathcal N}=\projlim_{\rho\in{\mathcal N}^X} \iota^\sigma_\rho:
X_{\mathcal N}\to \Ran\sigma$:
 \beq\label{X_F-proektiv-sistema} \xymatrix @R=2.5pc @C=2.0pc {
 & X_{\mathcal N}\ar[dr]^{\sigma_{\mathcal N}}\ar[dl]_{\rho_{\mathcal N}} &  \\
 \Ran\rho\ar[rr]^{\iota^\sigma_\rho} & &  \Ran\sigma
}\qquad (\rho\to\sigma);
 \eeq
при этом, в силу \eqref{DEF:le-in-F_X-*}, сама система эпиморфизмов ${\mathcal
N}_X$ также является проективным конусом системы $\Bind({\mathcal N}^X)$:
 \beq\label{F_X-proektiv-sistema} \xymatrix @R=2.5pc @C=2.0pc {
 & X\ar[dr]^{\sigma}\ar[dl]_{\rho} &  \\
 \Ran\rho\ar[rr]^{\iota^\sigma_\rho} & &  \Ran\sigma
}\qquad (\rho\to\sigma),
 \eeq
поэтому должен существовать естественный морфизм из $X$ в вершину $X_{\mathcal N}$
проективного предела системы $\Bind({\mathcal N}^X)$. Этот морфизм мы будем
обозначать $\projlim{\mathcal N}^X$ и называть {\it локальным пределом сети
эпиморфизмов ${\mathcal N}$ на объекте $X$}:
\beq\label{DIAGR:sigma-sigma_F}
\xymatrix @R=2.5pc @C=2.0pc {
   X\ar[dr]_{\sigma}\ar@{-->}[rr]^{\projlim{\mathcal N}^X} & & X_{\mathcal N}\ar[ld]^{\sigma_{\mathcal N}}  \\
   & \Ran\sigma  &
}\qquad (\sigma\in{\mathcal N}^X). \eeq

\item[---] элемент сети $\sigma$ в диаграмме \eqref{DIAGR:set} мы называем {\it
подпоркой} элемента сети $\tau$.
 }\eit
 }\eit

Примеры сетей эпиморфизмов мы приведем на страницах \pageref{LM:A/U-set-epimorfizmov-v-InvSteAlg}, \pageref{TH:A/U-set-epimorfizmov-v-InvSteAlg-diff} и \pageref{LM:A/U-set-epimorfizmov-v-InvSteAlg-complex}.

\btm\label{TH:funktorialnost-obolochki_F}\footnote{В работе \cite{Ak16} этот результат (Theorem 3.36) приведен с ошибкой: там опущено условие \eqref{lim_N-subseteq-Epi}.} Пусть ${\mathcal N}$ -- сеть эпиморфизмов в категории ${\tt K}$, причем все локальные пределы $\projlim{\mathcal N}^X$ являются эпиморфизмами
\beq\label{lim_N-subseteq-Epi}
\{\projlim{\mathcal N}^X; \ X\in\Ob(\tt K)\}\subseteq\Epi(\tt
K).
\eeq
Тогда
 \bit{
\item[(a)] для любого объекта $X$ в ${\tt K}$ локальный предел
$\projlim{\mathcal N}^X:X\to X_{\mathcal N}$ является оболочкой $\env_{\mathcal N} X$ в категории ${\tt
K}$ относительно класса морфизмов ${\mathcal N}$:
 \beq\label{lim F_X=env_F-X}
 \projlim {\mathcal N}^X=\env_{\mathcal N} X,
 \eeq

\item[(b)] для всякого морфизма $\alpha:X\to Y$ в
${\tt K}$ при любом выборе локальных пределов $\projlim {\mathcal N}^X$ и $\projlim {\mathcal N}^Y$ формула
 \beq\label{DEF:alpha_F}
\alpha_{\mathcal N}=\projlim_{\tau\in{\mathcal
N}_Y}\projlim_{\sigma\in{\mathcal N}^X}\alpha_\sigma^\tau\circ\sigma_{\mathcal
N}
 \eeq
определяет морфизм $\alpha_{\mathcal N}:X_{\mathcal N}\to Y_{\mathcal N}$, замыкающий диаграмму
 \beq\label{DIAGR:funktorialnost-lim-F}
\xymatrix @R=2.pc @C=10.0pc 
{
X\ar[d]^{\alpha}\ar[r]^{\projlim {\mathcal N}^X=\env_{\mathcal N}X} &  X_{\mathcal N}
=\Env_{\mathcal N}X\ar@{-->}[d]^{\alpha_{\mathcal N}} \\
Y\ar[r]^{\projlim {\mathcal N}^Y=\env_{\mathcal N}Y} &  Y_{\mathcal N}=\Env_{\mathcal N}Y
},
 \eeq

\item[(c)] оболочку $\Env_{\mathcal N}$ можно определить как функтор.
}\eit
\etm
 \bpr
1. По лемме \ref{LM:obolochka-konusa} проективный предел $\projlim
{\mathcal N}^X$ является оболочкой $X$ в $\tt K$ относительно конуса морфизмов
${\mathcal N}^X$:
$$
\projlim {\mathcal N}^X=\env_{{\mathcal N}^X} X
$$
В последнем выражении можно заменить ${\mathcal N}^X$ на ${\mathcal N}$, потому
что ${\mathcal N}^X$ есть в точности подкласс в ${\mathcal N}$, состоящий из
морфизмов, выходящих из $X$:
$$
\projlim {\mathcal N}^X=\env_{{\mathcal N}^X} X=\env_{\mathcal N} X.
$$

2. Объясним далее смысл формулы \eqref{DEF:alpha_F}. Пусть $\alpha:X\to
Y$ -- какой-нибудь морфизм. Для всякого элемента сети $\tau\in{\mathcal N}^Y$
обозначим
 \beq\label{opredelenie-alpha^tau}
\alpha^\tau=\tau\circ\alpha.
 \eeq
Понятно, что семейство морфизмов $\{\alpha^\tau:X\to\Ran\tau;\
\tau\in{\mathcal N}^Y\}$ образует проективный конус системы связывающих
морфизмов $\Bind({\mathcal N}^Y)$:
 \beq\label{X-alpha^tau=konus}
 \xymatrix @R=2.5pc @C=2.0pc {
 & X\ar[dr]^{\alpha^\upsilon}\ar[dl]_{\alpha^\tau} &  \\
 \Ran\tau\ar[rr]^{\iota_\tau^\upsilon} & &  \Ran\upsilon
}\qquad (\tau\to\upsilon).
 \eeq
По свойству (c) для всякого элемента сети $\tau\in{\mathcal N}^Y$ найдутся
элемент сети $\sigma\in{\mathcal N}^X$ и морфизм
$\alpha_\sigma^\tau:\Ran\sigma\to\Ran\tau$ такие, что коммутативна
диаграмма \eqref{DIAGR:set}, причем в ней мы уже символом $\alpha^\tau$
обозначили диагональ:
\beq\label{alpha^tau=tau-circ-alpha=alpha_sigma^tau-circ-sigma}
\alpha^\tau=\tau\circ\alpha=\alpha_\sigma^\tau\circ\sigma.
\eeq
Обозначив теперь \beq\label{opredelenie-beta^tau} \alpha_{\mathcal
N}^\tau=\alpha_\sigma^\tau\circ\sigma_{\mathcal N}, \eeq мы получим диаграмму
\beq\label{DIAFR:alpha^tau} \xymatrix @R=2.5pc @C=2.0pc {
   X\ar@/_4ex/[ddr]_{\alpha^\tau}\ar[dr]_{\sigma}\ar[rr]^{\projlim{\mathcal N}^X} & & X_{\mathcal N}\ar[ld]^{\sigma_{\mathcal N}}\ar@{-->}@/^4ex/[ddl]^{\alpha_{\mathcal N}^\tau}  \\
   & \Ran\sigma\ar[d]^{\alpha_\sigma^\tau}  & \\
   & \Ran\tau &
}\qquad (\sigma\in{\mathcal N}^X). \eeq

Заметим далее, что для любого другого элемента сети $\rho\in{\mathcal N}^X$,
удовлетворяющего условию $\rho\to\sigma$, выполняется равенство, аналогичное
\eqref{opredelenie-beta^tau}: \beq\label{opredelenie-beta^tau-1}
\alpha_{\mathcal N}^\tau=\alpha_\rho^\tau\circ\rho_{\mathcal N},\qquad
\rho\to\sigma. \eeq Действительно, при $\rho\to\sigma$ диаграмма
\eqref{DIAGR:set} достраивается до диаграммы \beq\label{DIAGR:set-3} \xymatrix
@R=2.5pc @C=4.0pc {
 & X\ar[r]^{\alpha}\ar[d]_{\sigma}\ar@/_3ex/[ddl]_{\rho} & Y\ar[d]^{\tau}  \\
 & \Ran\sigma\ar[r]_{\alpha_\sigma^\tau} & \Ran\tau  \\
 \Ran\rho\ar@{-->}@/_3ex/[rru]_{\alpha_\rho^\tau}\ar[ru]_{\iota_\rho^\sigma} & &
}
\eeq
(в которой пунктирная стрелка поначалу определяется как композиция $\alpha_\sigma^\tau\circ\iota_\rho^\sigma$, и затем, поскольку такая стрелка, если она существует, определяется однозначно, делается вывод, что это и есть морфизм $\alpha_\rho^\tau$). После этого мы получаем:
$$
\alpha_{\mathcal N}^\tau=\alpha_\sigma^\tau\circ\sigma_{\mathcal
N}=\eqref{X_F-proektiv-sistema}=
\alpha_\sigma^\tau\circ\iota_\rho^\sigma\circ\rho_{\mathcal
N}=\eqref{DIAGR:set-3}= \alpha_\rho^\tau\circ\rho_{\mathcal N}.
$$

Из формулы \eqref{opredelenie-beta^tau-1} следует, что определение морфизма
$\alpha_{\mathcal N}^\tau$ формулой \eqref{opredelenie-beta^tau} не зависит от
выбора элемента сети $\sigma\in{\mathcal N}^X$, потому что если
$\sigma'\in{\mathcal N}^X$ -- какой-нибудь другой элемент, для которого тоже
существует морфизм $\alpha_{\sigma'}^\tau:\Ran\sigma'\to\Ran\tau$,
замыкающий диаграмму \eqref{DIAGR:set} (в которой $\sigma$ заменен на
$\sigma'$), то выбрав элемент сети $\rho\in{\mathcal N}^X$, отстоящий слева от
$\sigma$ и $\sigma'$,
$$
\rho\to\sigma,\qquad \rho\to\sigma',
$$
(в этот момент мы используем аксиому (a) сети эпиморфизмов) мы получим цепочку
$$
\alpha_{\mathcal N}^\tau=\alpha_\sigma^\tau\circ\sigma_{\mathcal
N}=\eqref{opredelenie-beta^tau-1}= \alpha_\rho^\tau\circ\rho_{\mathcal
N}=\eqref{opredelenie-beta^tau-1}=\alpha_{\sigma'}^\tau\circ\sigma'_{\mathcal
N}.
$$

Мы можем теперь сделать вывод, что формула \eqref{opredelenie-beta^tau}
однозначно определяет некое отображение $\tau\in{\mathcal
N}_Y\mapsto\alpha_{\mathcal N}^\tau$. Покажем, что  семейство морфизмов
$\{\alpha_{\mathcal N}^\tau:X_{\mathcal N}\to \Ran\tau;\ \tau\in{\mathcal
N}_Y\}$ является проективным конусом системы связывающих морфизмов
$\Bind({\mathcal N}^Y)$.
 \beq\label{X_F-proektiv-sistema-dllya-F_Y} \xymatrix
@R=2.5pc @C=2.0pc {
 & X_{\mathcal N}\ar[dr]^{\alpha_{\mathcal N}^\upsilon}\ar[dl]_{\alpha_{\mathcal N}^\tau} &  \\
 \Ran\tau\ar[rr]^{\iota_\tau^\upsilon} & &  \Ran\upsilon
}\qquad (\tau\to\upsilon\in{\mathcal N}^Y).
 \eeq
При $\tau\to\upsilon$ диаграмма \eqref{DIAGR:set} достраивается до диаграммы
\beq\label{DIAGR:set-4}
\xymatrix @R=2.5pc @C=4.0pc
{
  X\ar[r]^{\alpha}\ar[d]_{\sigma} & Y\ar[d]^{\tau}\ar@/^3ex/[ddr]^{\upsilon} & \\
  \Ran\sigma\ar[r]_{\alpha_\sigma^\tau}\ar@{-->}@/_3ex/[rrd]_{\alpha_\sigma^\upsilon} &
  \Ran\tau\ar[rd]_{\iota_\tau^\upsilon} & \\
  & & \Ran\upsilon
}
\eeq
(в которой пунктирная стрелка поначалу определяется как композиция $\iota_\tau^\upsilon\circ\alpha_\sigma^\tau$, и затем, поскольку такая стрелка, если она существует, определяется однозначно, делается вывод, что это и есть морфизм $\alpha_\sigma^\upsilon$). Используя эту диаграмму, мы получаем:
$$
\iota_\tau^\upsilon\circ\alpha_{\mathcal
N}^\tau=\eqref{opredelenie-beta^tau}=\iota_\tau^\upsilon\circ\alpha_\sigma^\tau\circ\sigma_{\mathcal
N}= \eqref{DIAGR:set-4}=\alpha_\sigma^\upsilon\circ\sigma_{\mathcal
N}=\eqref{opredelenie-beta^tau}=\alpha_{\mathcal N}^\upsilon.
$$
Из диаграммы \eqref{X_F-proektiv-sistema-dllya-F_Y} следует теперь, что должен
существовать естественный морфизм $\alpha_{\mathcal N}$ из $X_{\mathcal N}$ в
проективный предел $Y_{\mathcal N}$ системы связывающих морфизмов
$\Bind({\mathcal N}^Y)$:
 \beq\label{sushestvovanie-alpha_F} \xymatrix @R=2.5pc
@C=2.0pc {
   X_{\mathcal N}\ar[dr]_{\alpha_{\mathcal N}^\tau}\ar@{-->}[rr]^{\alpha_{\mathcal N}} & & Y_{\mathcal N}\ar[ld]^{\tau_{\mathcal N}}  \\
   & \Ran\tau  &
}\qquad (\tau\in{\mathcal N}^Y).
 \eeq
Вспомним теперь, что в силу аксиомы сети (b) переход от $X$ к проективному
пределу $\projlim\Bind({\mathcal N}^X)$ можно оформить как отображение.
Дальнейшие же шаги по построению $\alpha_{\mathcal N}$ (то есть выбор вершины
$X_{\mathcal N}$ конуса $\projlim\Bind({\mathcal N}^X)$, а затем стрелки
$\alpha_{\mathcal N}$, замыкающей все диаграммы \eqref{sushestvovanie-alpha_F})
также однозначны, поэтому соответствие $\alpha\mapsto\alpha_{\mathcal N}$ тоже
можно считать отображением.

3. Покажем далее, что морфизмы $\alpha_{\mathcal N}$ замыкают диаграммы вида
\eqref{DIAGR:funktorialnost-lim-F}. В диаграмме
$$
\xymatrix @R=2.pc @C=5.0pc 
{
X\ar[dd]^{\alpha}\ar[rr]^{\projlim{\mathcal N}^X}\ar[dr]_{\alpha^\tau} &  & X_{\mathcal N}\ar@{-->}[dd]^{\alpha_{\mathcal N}}\ar[dl]^{\alpha_{\mathcal N}^\tau} \\
 & \Ran\tau & \\
Y\ar[rr]^{\projlim {\mathcal N}^Y}\ar[ur]^{\tau} & &  Y_{\mathcal N}\ar[ul]_{\tau_{\mathcal N}} \\
}
$$
все внутренние треугольники коммутативны: верхний внутренний треугольник -- потому что это периметр диаграммы \eqref{DIAFR:alpha^tau},
левый внутренний треугольник -- потому что это по-другому записанная формула \eqref{opredelenie-alpha^tau},
нижний внутренний треугольник -- потому что с точностью до обозначений это диаграмма \eqref{DIAGR:sigma-sigma_F},
а правый внутренний треугольник -- потому что это повернутая диаграмма \eqref{sushestvovanie-alpha_F}. Поэтому справедливы равенства
$$
\tau_{\mathcal N}\circ\projlim{\mathcal
N}_Y\circ\alpha=\alpha^\tau=\tau_{\mathcal N}\circ\alpha_{\mathcal
N}\circ\projlim{\mathcal N}^X \qquad (\tau\in{\mathcal N}^Y)
$$
Это можно понимать так, что каждый из морфизмов $\projlim{\mathcal
N}_Y\circ\alpha$ и $\alpha_{\mathcal N}\circ\projlim{\mathcal N}^X$
представляет собой поднятие проективного конуса $\{\alpha^\tau:X\to\Ran\tau;\
\tau\in{\mathcal N}^Y\}$ для системы связывающих морфизмов $\Bind({\mathcal
N}_Y)$, о котором речь шла в диаграмме \eqref{X-alpha^tau=konus}, на
проективный предел этой системы. То есть $\projlim{\mathcal N}^Y\circ\alpha$ и
$\alpha_{\mathcal N}\circ\projlim{\mathcal N}^X$ представляют собой ту самую
пунктирную стрелку в определении проективного предела, для которой коммутативны
все диаграммы вида
$$
\xymatrix @R=2.5pc @C=2.0pc
{
   X\ar[dr]_{\alpha^\tau}\ar@{-->}[rr] & & Y_{\mathcal N}\ar[ld]^{\tau_{\mathcal N}}  \\
   & \Ran\tau  &
}\qquad (\tau\in{\mathcal N}^Y).
$$
Но такая пунктирная стрелка единственна, поэтому эти морфизмы должны совпадать:
$$
\projlim{\mathcal N}^Y\circ\alpha=\alpha_{\mathcal N}\circ\projlim{\mathcal
N}_X.
$$
Это и есть диаграмма \eqref{DIAGR:funktorialnost-lim-F}.

4. Используемая в теории множеств теорема
о полном упорядочении класса всех множеств \cite[Theorem 1.1.1]{Akbarov-De-Gruyter-I} позволяет определить операцию
перехода к локальному пределу как отображение:
$$
X\mapsto \projlim\Bind({\mathcal N}^X)
$$
(то есть существует отображение, которое каждому объекту $X\in\Ob(\tt K)$ ставит в соответствие
определенный проективный предел подкатегории $\Bind({\mathcal N}^X)$ среди всех возможных ее проективных
пределов в ${\tt K}$). Покажем, что в этом случае возникающее отображение
$(X,\alpha)\mapsto(X_{\mathcal N},\alpha_{\mathcal N})$ является функтором, то есть что справедливы тождества
 \beq\label{funktorialnost-alpha_F}
 (1_X)_{\mathcal N}= 1_{X_{\mathcal N}},\qquad
 (\beta\circ\alpha)_{\mathcal N}=\beta_{\mathcal N}\circ\alpha_{\mathcal N}.
 \eeq
Начнем с первого: пусть $\alpha=1_X:X\to X$. Тогда мы получим цепочку:
\begin{multline*}
\alpha^\tau=\eqref{opredelenie-alpha^tau}=\tau\circ\alpha=\tau\circ 1_X=\tau\quad\Longrightarrow\quad
\alpha_\sigma^\tau\circ\sigma=\eqref{alpha^tau=tau-circ-alpha=alpha_sigma^tau-circ-sigma}=\alpha^\tau=\tau=\eqref{DEF:le-in-F_X-*}=\iota_\sigma^\tau\circ\sigma
\quad\Longrightarrow\\ \Longrightarrow\quad \alpha_\sigma^\tau=\iota_\sigma^\tau \quad\Longrightarrow\quad  \alpha_{\mathcal N}^\tau=\iota_\sigma^\tau\circ\sigma_{\mathcal N}=\tau_{\mathcal N}
\end{multline*}
Мы получаем, что в диаграммах \eqref{sushestvovanie-alpha_F} можно заменить морфизм $\alpha_{\mathcal N}^\tau$ на морфизм $\tau_{\mathcal N}$:
 $$
  \xymatrix @R=2.5pc
@C=2.0pc {
   X_{\mathcal N}\ar[dr]_{\tau_{\mathcal N}}\ar@{-->}[rr]^{\alpha_{\mathcal N}} & & X_{\mathcal N}\ar[ld]^{\tau_{\mathcal N}}  \\
   & \Ran\tau  &
}\qquad (\tau\in{\mathcal N}^X).
 $$
Эти диаграммы будут коммутативны для всех $\tau\in{\mathcal N}^X$, и в них пунктирная стрелка
$\alpha_{\mathcal N}$ определяется как поднятие проективного конуса
$\{\alpha_{\mathcal N}^\tau=\tau_{\mathcal N}:X_{\mathcal N}\to\Ran\tau\}$
на проективный предел $\{\tau_{\mathcal N}:X_{\mathcal N}\to\Ran\tau\}$. Поскольку такая стрелка единственна, она должна совпадать с морфизмом $1_{X_{\mathcal N}}$, для которого все эти диаграммы коммутативны тривиальным образом: $\alpha_{\mathcal N}=1_{X_{\mathcal N}}$.

Перейдем теперь ко второму тождеству в \eqref{funktorialnost-alpha_F}. Рассмотрим последовательность морфизмов
$X\overset{\alpha}{\longrightarrow}Y\overset{\beta}{\longrightarrow}Z$.
Зафиксируем элемент сети $\upsilon\in{\mathcal N}^Z$ и, пользуясь аксиомой (c),
выберем для него сначала элемент сети $\tau\in{\mathcal N}^Y$ и морфизм
$\beta_\tau^\upsilon$ так, чтобы выполнялось равенство
$$
\upsilon\circ\beta=\beta_\tau^\upsilon\circ\tau,
$$
а затем, опять же пользуясь аксиомой сети (c), выберем элемент сети
$\sigma\in{\mathcal N}^X$ и морфизм $\alpha_\sigma^\tau$ так, чтобы выполнялось
равенство
$$
\tau\circ\alpha=\alpha_\sigma^\tau\circ\sigma.
$$
Мы получим диаграмму
$$
\xymatrix @R=2.5pc @C=2.0pc
{
   X\ar[d]^{\sigma}\ar[r]^{\alpha} & Y\ar[d]^{\tau}\ar[r]^{\beta} & Z\ar[d]^{\upsilon}  \\
   \Ran\sigma\ar[r]^{\alpha_\sigma^\tau} & \Ran\tau\ar[r]^{\beta_\tau^\upsilon}  & \Ran\upsilon
}.
$$
Если в ней удалить среднюю стрелку, то получится диаграмма
$$
\xymatrix @R=2.5pc @C=4.0pc
{
   X\ar[d]^{\sigma}\ar[r]^{\beta\circ\alpha} & Z\ar[d]^{\upsilon}  \\
   \Ran\sigma\ar[r]^{\beta_\tau^\upsilon\circ\alpha_\sigma^\tau} & \Ran\upsilon
},
$$
которую можно понимать так, что морфизм $\beta_\tau^\upsilon\circ\alpha_\sigma^\tau$ есть та самая единственно возможная пунктирная стрелка, о которой говорится в диаграмме \eqref{DIAGR:set}, только с той разницей, что в ней $Y$ поменяли на $Z$, $\alpha$ на $\beta\circ\alpha$, а $\tau$ на $\upsilon$. Мы можем поэтому сделать вывод, что существует морфизм $(\beta\circ\alpha)_\sigma^\upsilon$, совпадающий с $\beta_\tau^\upsilon\circ\alpha_\sigma^\tau$:
\beq\label{beta_tau^upsilon-circ-alpha_sigma^tau=(beta-circ-alpha)_sigma^upsilon}
\beta_\tau^\upsilon\circ\alpha_\sigma^\tau=(\beta\circ\alpha)_\sigma^\upsilon
\eeq
Это равенство используется в следующей цепочке:
$$
\underbrace{\upsilon_{\mathcal N}\circ\beta_{\mathcal N}}_{\scriptsize
\begin{matrix}\phantom{\tiny \eqref{sushestvovanie-alpha_F}} \ \text{\rotatebox{90}{$=$}} \ {\tiny \eqref{sushestvovanie-alpha_F}}\\ \beta_{\mathcal N}^\upsilon \\ \phantom{\tiny \eqref{opredelenie-beta^tau}}  \ \text{\rotatebox{90}{$=$}} \ {\tiny \eqref{opredelenie-beta^tau}} \\
\beta_\tau^\upsilon\circ\tau_{\mathcal N} \end{matrix}} \circ\ \alpha_{\mathcal
N}= \beta_\tau^\upsilon\circ\kern-3pt \underbrace{\tau_{\mathcal
N}\circ\alpha_{\mathcal N}}_{\scriptsize
\begin{matrix}\phantom{\tiny \eqref{sushestvovanie-alpha_F}} \ \text{\rotatebox{90}{$=$}} \ {\tiny \eqref{sushestvovanie-alpha_F}}\\
\alpha_{\mathcal N}^\tau \\ \phantom{\tiny \eqref{opredelenie-beta^tau}}  \ \text{\rotatebox{90}{$=$}} \ {\tiny \eqref{opredelenie-beta^tau}} \\
\alpha_\sigma^\tau\circ\sigma_{\mathcal N} \end{matrix}} =
\underbrace{\beta_\tau^\upsilon\circ\alpha_\sigma^\tau}_{\scriptsize
\begin{matrix}\phantom{\tiny \eqref{beta_tau^upsilon-circ-alpha_sigma^tau=(beta-circ-alpha)_sigma^upsilon}} \ \text{\rotatebox{90}{$=$}} \ {\tiny \eqref{beta_tau^upsilon-circ-alpha_sigma^tau=(beta-circ-alpha)_sigma^upsilon}}\\
(\beta\circ\alpha)_\sigma^\upsilon \end{matrix}} \kern-3pt\circ\
\sigma_{\mathcal N}=
\underbrace{(\beta\circ\alpha)_\sigma^\upsilon\circ\sigma_{\mathcal
N}}_{\scriptsize
\begin{matrix}\phantom{\tiny \eqref{opredelenie-beta^tau}} \ \text{\rotatebox{90}{$=$}} \ {\tiny \eqref{opredelenie-beta^tau}}\\
(\beta\circ\alpha)_{\mathcal N}^\upsilon  \end{matrix}}
=(\beta\circ\alpha)_{\mathcal
N}^\upsilon=\eqref{sushestvovanie-alpha_F}=\upsilon_{\mathcal N}\circ
(\beta\circ\alpha)_{\mathcal N}.
$$
Опустив промежуточные выкладки, мы получим такое двойное равенство:
$$
\upsilon_{\mathcal N}\circ(\beta_{\mathcal N}\circ\alpha_{\mathcal
N})=(\beta\circ\alpha)_{\mathcal N}^\upsilon=\upsilon_{\mathcal N}\circ
(\beta\circ\alpha)_{\mathcal N}.
$$
Оно выполняется для любого $\upsilon\in{\mathcal N}^Z$. Поэтому его можно
понимать так, что каждый из морфизмов $\beta_{\mathcal N}\circ\alpha_{\mathcal
N}$ и $(\beta\circ\tau)_{\mathcal N}$ представляет собой поднятие проективного
конуса $\{(\beta\circ\alpha)_{\mathcal N}^\upsilon:X_{\mathcal
N}\to\Ran\upsilon;\ \upsilon\in{\mathcal N}^Z\}$ для системы связывающих
морфизмов $\Bind({\mathcal N}^Z)$ (а это семейство действительно будет
проективным конусом в силу диаграммы \eqref{X_F-proektiv-sistema-dllya-F_Y} в
которой нужно поменять $Y$ на $Z$, а $\alpha$ на $\beta\circ\alpha$) на
проективный предел этой системы. То есть $\beta_{\mathcal
N}\circ\alpha_{\mathcal N}$ и $(\beta\circ\tau)_{\mathcal N}$ представляют
собой ту самую пунктирную стрелку в определении проективного предела, для
которой коммутативны все диаграммы вида
$$
\xymatrix @R=2.5pc @C=2.0pc
{
   X_{\mathcal N}\ar[dr]_{(\beta\circ\alpha)_{\mathcal N}^\upsilon}\ar@{-->}[rr] & & Z_{\mathcal N}\ar[ld]^{\upsilon_{\mathcal N}}  \\
   & \Ran\upsilon  &
}\qquad (\upsilon\in{\mathcal N}^Z).
$$
Но такая пунктирная стрелка единственна, поэтому эти морфизмы должны совпадать:
$$
\beta_{\mathcal N}\circ\alpha_{\mathcal N}=(\beta\circ\tau)_{\mathcal N}.
$$
Это и есть второе тождество \eqref{funktorialnost-alpha_F}.
\epr

\btm\label{TH:funktorialnost-obolochki} Пусть ${\mathcal N}$ -- сеть
эпиморфизмов в категории ${\tt K}$, порождающая класс морфизмов $\varPhi$ изнутри:
    $$
    {\mathcal N}\subseteq\varPhi\subseteq\Mor({\tt K})\circ {\mathcal N}.
    $$
Тогда для любого класса эпиморфизмов $\varOmega$ в $\tt K$, содержащего все локальные пределы
$\projlim{\mathcal N}^X$,
\beq\label{lim_N-subseteq-varOmega}
\{\projlim{\mathcal N}^X; \ X\in\Ob(\tt K)\}\subseteq\varOmega\subseteq\Epi(\tt
K),
\eeq
выполняется следующее:
 \bit{
\item[(a)] для любого объекта $X$ в ${\tt K}$ локальный предел
$\projlim{\mathcal N}^X$ является оболочкой $\env_\varPhi^\varOmega X$ в классе
$\varOmega$ относительно класса $\varPhi$: \beq\label{lim F_X=env_M^varOmega-X}
\projlim {\mathcal N}^X=\env_\varPhi^\varOmega X, \eeq

\item[(b)] оболочку $\Env_\varPhi^\varOmega$ можно определить как функтор.
 }\eit
 \etm

\bpr По теореме \ref{TH:funktorialnost-obolochki_F}
локальный предел сети $\projlim {\mathcal N}^X$ является оболочкой $X$ в классе
$\Mor({\tt K})$ всех морфизмов категории $\tt K$ относительно класса морфизмов
${\mathcal N}$:
$$
\projlim {\mathcal N}^X=\env_{\mathcal N} X:=\env_{\mathcal N}^{\Mor(\tt K)} X.
$$
С другой стороны, поскольку в силу (i) $\projlim {\mathcal N}^X$ лежит в более
узком классе $\varOmega$, по свойству $1^\circ$ (c) на
с.\pageref{LM:suzhenie-verh-klassa-morfizmov}, $\projlim {\mathcal N}^X$ должен
быть оболочкой в этом более узком классе $\varOmega$:
$$
\projlim {\mathcal N}^X=\env_{\mathcal N} X=\env_{\mathcal N}^{\Mor(\tt K)}
X=\env_{\mathcal N}^\varOmega X.
$$
Далее, поскольку класс ${\mathcal N}$ порождает класс $\varPhi$, а $\varOmega$
состоит из эпиморфизмов, в силу \eqref{env_Psi=env_Phi}, оболочка относительно
${\mathcal N}$ должна совпадать с оболочкой относительно $\varPhi$:
$$
\projlim {\mathcal N}^X=\env_{\mathcal N}^\varOmega X=\env_\varPhi^\varOmega X.
$$
Это доказывает \eqref{lim F_X=env_M^varOmega-X}. Пункт (b) следует теперь из теоремы
\ref{TH:funktorialnost-obolochki_F} (c). \epr

От левой части условия \eqref{lim_N-subseteq-varOmega} можно отказаться, если класс $\varOmega$ мономорфно дополняем:

\btm\label{TH:funktorialnost-pri-seti-Epi-i-dolonyaemosti}
Пусть ${\mathcal N}$ -- сеть эпиморфизмов в категории ${\tt K}$, порождающая класс морфизмов $\varPhi$ изнутри:
    $$
    {\mathcal N}\subseteq\varPhi\subseteq\Mor({\tt K})\circ {\mathcal N}.
    $$
Тогда для любого мономорфно
дополняемого\footnote{В смысле определения на с.\pageref{DEF:klass-monomorfno-dopolnyaem}.} в $\tt K$
класса эпиморфизмов $\varOmega$,
$$
{^\downarrow\varOmega}\circledcirc\varOmega={\tt K},
$$
выполняются следующие условия:
 \bit{
\item[(a)] для любого объекта $X$ в ${\tt K}$ морфизм $\e_{\projlim{\mathcal N}^X}$ в
факторизации \eqref{faktorizatsiya-v-kat-s-faktoriz}, заданной классами ${^\downarrow\varOmega}$ и
$\varOmega$, является оболочкой  $\env_\varPhi^{\varOmega} X$ в $\varOmega$ относительно $\varPhi$ и одновременно с оболочкой $\env_{\mathcal N}^{\varOmega} X$ в $\varOmega$ относительно $\mathcal N$:
\beq\label{im_infty-lim-F_X=env_M^Epi-X-1}
\e_{\projlim{\mathcal N}^X}=\env_{\mathcal N}^{\varOmega} X=\env_\varPhi^{\varOmega} X,
\eeq

\item[(b)] для всякого морфизма $\alpha:X\to Y$ в ${\tt K}$ при любом выборе оболочек  $\env_\varPhi^{\varOmega}X$ и $\env_\varPhi^{\varOmega}Y$ существует единственный морфизм $\Env_\varPhi^{\varOmega} \alpha:\Env_\varPhi^{\varOmega}  X\to \Env_\varPhi^{\varOmega}  Y$ в ${\tt K}$, замыкающий диаграмму
\beq\label{DIAGR:funktorialnost-env_varPhi^Epi-v-kat-s-uzl-razl-1}
\xymatrix @R=2.pc @C=5.0pc 
{
X\ar[d]^{\alpha}\ar[r]^{\env_\varPhi^{\varOmega}  X} & \Env_\varPhi^{\varOmega}  X\ar@{-->}[d]^{\Env_\varPhi^{\varOmega} \alpha} \\
Y\ar[r]^{\env_\varPhi^{\varOmega}  Y} & \Env_\varPhi^{\varOmega}  Y \\
}
\eeq

\item[(c)] если дополнительно $\tt K$ локально мала в фактор-объектах класса $\varOmega$, то
оболочку $\Env_\varPhi^{\varOmega}$ можно определить как функтор.
}\eit
 \etm \bpr
1. Поскольку класс ${\mathcal N}$ порождает класс $\varPhi$, а $\varOmega$
состоит из эпиморфизмов, в силу \eqref{env_Psi=env_Phi}, оболочка относительно
${\mathcal N}$ должна совпадать с оболочкой относительно $\varPhi$:
\beq\label{PROOF:TH:funktorialnost-pri-seti-Epi-i-dolonyaemosti-1}
\env_{\mathcal N}^\varOmega X=\env_\varPhi^\varOmega X.
\eeq
После этого равенство \eqref{obolochka-konusa-v-kat-s-uzl-razl} из леммы
\ref{LM:obolochka-konusa-v-kat-s-uzl-razl} делает очевидным \eqref{im_infty-lim-F_X=env_M^Epi-X-1}:
$$
\env_{\varPhi}^{\varOmega}X=\env_{\mathcal N}^{\varOmega}X=\e_{\projlim{\mathcal
N}_X}.
$$

2. Свойство \eqref{DIAGR:funktorialnost-env_varPhi^Epi-v-kat-s-uzl-razl-1} доказывается так. Сначала мы дополняем диаграмму \eqref{DIAGR:funktorialnost-lim-F} разложением пределов $\projlim{\mathcal
N}^X$ и $\projlim{\mathcal N}^Y$ на важные для нас составляющие:
$$
\xymatrix @R=3.pc @C=5.0pc 
{
X\ar@/^5ex/[rr]^{\projlim{\mathcal N}^X}\ar[d]^{\alpha}\ar[r]_(.3){\env_\varPhi^{\varOmega}  X} & \Env_\varPhi^{\varOmega}X=
\Dom\mu_{\projlim{\mathcal N}^X}\ar[r]_(.7){\mu_{\projlim{\mathcal N}^X}} & X_{\mathcal N}\ar[d]^{\alpha_{\mathcal N}} \\
Y\ar@/_5ex/[rr]_{\projlim{\mathcal N}^Y}\ar[r]^(.3){\env_\varPhi^{\varOmega}  Y} & \Env_\varPhi^{\varOmega}Y=
\Dom\mu_{\projlim{\mathcal N}^Y}\ar[r]^(.7){\mu_{\projlim{\mathcal N}^Y}} & Y_{\mathcal N} \\
}
$$
Представим внутренний прямоугольник в виде
$$
\xymatrix @R=3.pc @C=5.0pc 
{
X\ar[d]^{\alpha}\ar[r]^(.3){\env_\varPhi^{\varOmega}  X} & \Env_\varPhi^{\varOmega}X=\Dom\mu_{\projlim{\mathcal N}^X}
 \ar@/^2ex/[dr]^(.7){\alpha_{\mathcal N}\circ\mu_{\projlim{\mathcal N}^X}} & \\
Y\ar[r]^(.3){\env_\varPhi^{\varOmega}  Y} & \Env_\varPhi^{\varOmega}Y=\Dom\mu_{\projlim{\mathcal N}^Y} \ar[r]_(.6){\mu_{\projlim{\mathcal N}^Y}} & Y_{\mathcal N} \\
}
$$
Здесь верхняя горизонтальная стрелка, $\env_\varPhi^{\varOmega}X$, лежит в $\varOmega$, а вторая нижняя горизонтальная стрелка, $\mu_{\projlim{\mathcal N}^Y}$, лежит в $\varGamma={^\downarrow\varOmega}$. Поэтому
должен существовать морфизм $\xi$, для которого будет коммутативна диаграмма
$$
\xymatrix @R=3.pc @C=5.0pc 
{
X\ar[d]^{\alpha}\ar[r]^(.3){\env_\varPhi^{\varOmega}  X} & \Env_\varPhi^{\varOmega}X=\Dom\mu_{\projlim{\mathcal N}^X}
\ar@{-->}[d]^{\xi}\ar@/^2ex/[dr]^(.7){\alpha_{\mathcal N}\circ\mu_{\projlim{\mathcal N}^X}} & \\
Y\ar[r]^(.3){\env_\varPhi^{\varOmega}  Y} & \Env_\varPhi^{\varOmega}Y=\Dom\mu_{\projlim{\mathcal N}^Y}
 \ar[r]_(.6){\mu_{\projlim{\mathcal N}^Y}} & Y_{\mathcal N} \\
}
$$
и он как раз будет недостающей вертикальной стрелкой в \eqref{DIAGR:funktorialnost-env_varPhi^Epi-v-kat-s-uzl-razl-1}.

3. Пусть $\tt K$ локально мала в фактор-объектах класса $\varOmega$, то есть
для всякого $X$ категория $\varOmega^X=\varOmega\cap\Epi^X$ скелетно мала.
Пусть $S_X$  -- ее скелет, являющийся множеством. Воспользовавшись теоремой
\ref{TH:o-lok-malosti-v-faktor-objektah}, мы можем выбрать отображение $X\mapsto S_X$,
которое каждому объекту $X$ ставит в соответствие некий скелет $S_X$ в категории $\varOmega^X$.
Зафиксируем такое отображение.
Чтобы определить оболочку $\Env_\varPhi^{\varOmega}$ как функтор, нужно теперь
с помощью аксиомы выбора определить отображение $X\in\Ob({\tt K})\mapsto \env_\varPhi^{\varOmega}X\in S_X$,
после этого объект $\Env_\varPhi^{\varOmega}X$ определяется как область значений морфизма
$\env_\varPhi^{\varOmega}X$, а морфизм $\Env_\varPhi^{\varOmega}\alpha$
\eqref{DIAGR:funktorialnost-env_varPhi^Epi-v-kat-s-uzl-razl-1} появляется автоматически (как единственно возможный).
 \epr

\paragraph{Существование сети эпиморфизмов и полурегулярная оболочка.}

\btm[о полурегулярной оболочке]\label{TH:sushestvovanie-seti-pri-faktorizatsii} Пусть категория $\tt K$ и классы морфизмов $\varOmega$ и $\varPhi$
в ней удовлетворяют следующим условиям:
\bit{

\item[RE.1:]\label{DEF:RE.1} категория $\tt K$ проективно полна,

\item[RE.2:]\label{DEF:RE.2} класс $\varOmega$ мономорфно дополняем в $\tt K$:
${^\downarrow\varOmega}\circledcirc\varOmega={\tt K}$,

\item[RE.3:]\label{DEF:RE.3} ${\tt K}$ локально мала в фактор-объектах класса $\varOmega$, и

\item[RE.4:]\label{DEF:RE.4} класс $\varPhi$ выходит из\footnote{В смысле определения на с.\pageref{DEF:goes-from}.} $\tt K$ и является правым идеалом в ${\tt K}$:
$$
\Ob(\tt K)=\{\Dom\ph;\ \ph\in\varPhi\},\qquad \varPhi\circ\Mor({\tt K})\subseteq\varPhi.
$$
}\eit
Тогда
 \bit{
\item[(a)] в ${\tt K}$ существует сеть эпиморфизмов ${\mathcal N}$ такая, что для любого объекта $X$
в ${\tt K}$ оболочка $\env_\varPhi^{\varOmega} X$ в $\varOmega$ относительно $\varPhi$ совпадает с оболочкой $\env_{\mathcal N}^{\varOmega} X$ в $\varOmega$ относительно $\mathcal N$ и с морфизмом $\e_{\projlim{\mathcal N}^X}$ в факторизации \eqref{faktorizatsiya-v-kat-s-faktoriz}:
\beq\label{im_infty-lim F_X=env_M^Epi-X-1-*}
\e_{\projlim{\mathcal N}^X}=\env_{\mathcal N}^{\varOmega} X=\env_\varPhi^{\varOmega} X,
\eeq

\item[(b)] для всякого морфизма $\alpha:X\to Y$ в ${\tt K}$ при любом выборе оболочек  $\env_\varPhi^{\varOmega}X$ и $\env_\varPhi^{\varOmega}Y$ существует единственный морфизм $\Env_\varPhi^{\varOmega} \alpha:\Env_\varPhi^{\varOmega}  X\to \Env_\varPhi^{\varOmega}  Y$ в ${\tt K}$, замыкающий диаграмму
\beq\label{DIAGR:funktorialnost-env_varPhi^Epi-v-kat-s-uzl-razl-1-*}
\xymatrix @R=2.pc @C=5.0pc 
{
X\ar[d]^{\alpha}\ar[r]^{\env_\varPhi^{\varOmega}  X} & \Env_\varPhi^{\varOmega}  X\ar@{-->}[d]^{\Env_\varPhi^{\varOmega} \alpha} \\
Y\ar[r]^{\env_\varPhi^{\varOmega}  Y} & \Env_\varPhi^{\varOmega}  Y \\
}
\eeq

\item[(c)] оболочку $\Env_\varPhi^{\varOmega}$ можно определить как функтор.
}\eit
\etm

\bit{

\item[$\bullet$] Если выполняются условия RE.1-RE.4, то мы будем говорить, что {\it классы морфизмов $\varOmega$ и $\varPhi$ определяют в категории $\tt K$ полурегулярную оболочку}\label{DEF:polureg-obolochka}, или что {\it оболочка $\Env_\varPhi^\varOmega$ полурегулярна}.
}\eit

\brem\label{REM:env-polureg-obolochki}
Из условия RE.2 следует, что класс $\varOmega$ должен состоять из эпиморфизмов, поэтому морфизм полурегулярной оболочки $\env_\varPhi^{\varOmega} X:X\to\Env_\varPhi^{\varOmega} X$ всегда является эпиморфизмом.
\erem

\bpr
1. По условию, для всякого $X$ категория $\varOmega^X=\varOmega\cap\Epi^X$ скелетно мала.
Пусть $S_X$  -- ее скелет, являющийся множеством. Воспользовавшись теоремой
\ref{TH:o-lok-malosti-v-faktor-objektah}, мы можем и выбрать отображение $X\mapsto S_X$,
которое каждому объекту $X$ ставит в соответствие некий скелет $S_X$ в категории $\varOmega^X$.
Зафиксируем такое отображение.

Для всякого объекта $X$ в ${\tt K}$ обозначим через $\varPhi^X$ подкласс морфизмов из $\varPhi$, выходящих из $X$:
$$
\varPhi^X=\{\ph\in\varPhi:\ \Dom\ph=X\}
$$
(из RE.4 следует, что $\varPhi^X\ne\varnothing$).
Обозначим символом $2_{\varPhi^X}$ класс конечных подмножеств в $\varPhi^X$:
$$
\varPsi\in 2_{\varPhi^X}\qquad\Longleftrightarrow\qquad \varPsi\subseteq\varPhi^X\quad\&\quad \card\varPsi<\infty.
$$
Всякому объекту $X$ в ${\tt K}$ и всякому множеству морфизмов $\varPsi\in 2_{\varPhi^X}$ поставим в соответствие морфизм
$$
\overline{\varPsi}=\prod_{\psi\in\varPsi}\psi:X\to\prod_{\psi\in\varPsi}\Ran\psi,\qquad
$$
и морфизмы $\mu_{\varPsi}\in\varGamma$ и $\e_{\varPsi}\in S_X$ со свойством
\beq\label{overline(varPsi)=mu_varPsi-circ-e_(varPsi)}
\overline{\varPsi}=\mu_{\varPsi}\circ\e_{\varPsi}.
\eeq
(поскольку $S_X$ -- скелет $\varOmega^X$, такие морфизмы выбираются однозначно). Положим
\beq\label{PROOF:TH:sushestvovanie-seti-pri-faktorizatsii-0}
{\mathcal N}^X=\{\e_{\varPsi};\ \varPsi\in 2_{\varPhi^X}\}.
\eeq
Поскольку ${\mathcal N}^X\subseteq S_X$, это будет множество, а поскольку соответствие $X\mapsto S_X$ было отображением, получающееся соответствие $X\mapsto {\mathcal N}^X$ также является отображением.

2. Проверим для системы $\mathcal N$ аксиомы сети эпиморфизмов (с.\pageref{AX:set-Epi-a}). Во-первых, покажем, что ${\mathcal N}^X$ направлено влево относительно предпорядка \eqref{DEF:le-in-F_X}, наследуемого из $\Epi^X$. Для любых двух множеств $\varPsi,\varPsi'\in 2_{\varPhi^X}$ рассмотрим диаграмму
$$
\xymatrix @R=3.pc @C=5.0pc 
{
& X\ar[d]^{\overline{\varPsi\cup\varPsi'}}\ar@/_3ex/[ld]_{\overline{\varPsi}}
\ar@/^3ex/[rd]^{\overline{\varPsi'}} & \\
\Ran\overline{\varPsi}\ar@{=}[d] & \Ran\overline{\varPsi\cup\varPsi'}\ar[l]_{\pi}\ar[r]^{\pi'}\ar@{=}[d] &
\Ran\overline{\varPsi'}\ar@{=}[d]\\
\prod\limits_{\psi\in\varPsi}\Ran\psi & \prod\limits_{\psi\in\varPsi\cup\varPsi'}\Ran\psi &
\prod\limits_{\psi\in\varPsi'}\Ran\psi
}
$$
в которой $\pi$ и $\pi'$ -- естественные проекции. Пропустим стрелки, выходящие из $X$, через факторизацию \eqref{overline(varPsi)=mu_varPsi-circ-e_(varPsi)}:
\beq\label{DIAGR:sushestvovanie-seti-Epi-pri-uzlovom-razlozhenii}
\xymatrix @R=3.pc @C=5.0pc 
{
& X\ar[d]^{\e_{\varPsi\cup\varPsi'}}\ar@/_3ex/[ld]_{\e_{\varPsi}}
\ar@/^3ex/[rd]^{\e_{\varPsi'}} & \\
\Ran\e_{\varPsi}\ar[d]^{\mu_{\varPsi}} & \Ran\e_{\varPsi\cup\varPsi'}
\ar[d]^{\mu_{\varPsi\cup\varPsi'}} & \Ran\e_{\varPsi'}\ar[d]^{\mu_{\varPsi'}}
\\
\Ran\overline{\varPsi} & \Ran\overline{\varPsi\cup\varPsi'}\ar[l]_{\pi}\ar[r]^{\pi'} & \Ran\overline{\varPsi'}
}
\eeq
Левую половину диаграммы представим в виде четырехугольника:
$$
\xymatrix @R=3.pc @C=5.0pc 
{
& X\ar[d]^{\e_{\varPsi\cup\varPsi'}}\ar@/_3ex/[ld]_{\e_{\varPsi}}
\\
\Ran\e_{\varPsi}\ar[d]^{\mu_{\varPsi}} & \Ran\e_{\varPsi\cup\varPsi'}
\ar@/^3ex/[dl]^{\pi\circ\mu_{\varPsi\cup\varPsi'}}
\\
\Ran\overline{\varPsi} &
}
$$
В нем $\e_{\varPsi\cup\varPsi'}$ -- эпиморфизм, а $\mu_{\varPsi}$ -- строгий мономорфизм, поэтому должна существовать горизонтальная стрелка влево:
$$
\xymatrix @R=3.pc @C=5.0pc 
{
& X\ar[d]^{\e_{\varPsi\cup\varPsi'}}\ar@/_3ex/[ld]_{\e_{\varPsi}}
\\
\Ran\e_{\varPsi}\ar[d]^{\mu_{\varPsi}} & \Ran\e_{\varPsi\cup\varPsi'}\ar@/^3ex/[dl]^{\pi\circ\mu_{\varPsi\cup\varPsi'}}
\ar@{-->}[l]_{\delta}
\\
\Ran\overline{\varPsi} &
}
$$
По тем же причинам в \eqref{DIAGR:sushestvovanie-seti-Epi-pri-uzlovom-razlozhenii} должна существовать горизонтальная стрелка вправо, и мы получаем диаграмму
$$
\xymatrix @R=3.pc @C=5.0pc 
{
& X\ar[d]^{\e_{\varPsi\cup\varPsi'}}\ar@/_3ex/[ld]_{\e_{\varPsi}}
\ar@/^3ex/[rd]^{\e_{\varPsi'}} & \\
\Ran\e_{\varPsi} & \Ran\e_{\varPsi\cup\varPsi'}
\ar@{-->}[l]_{\delta}\ar@{-->}[r]^{\delta'} & \Ran\e_{\varPsi'}
}
$$
означающую, что в категории $\varOmega^X$ морфизм $\e_{\varPsi\cup\varPsi'}$ мажорирует морфизмы $\e_{\varPsi}$ и $\e_{\varPsi'}$:
$$
\e_{\varPsi\cup\varPsi'}\to \e_{\varPsi}\quad\&\quad \e_{\varPsi\cup\varPsi'}\to\e_{\varPsi'}.
$$

Второе условие в определении сети эпиморфизмов выполняется автоматически: поскольку категория $\tt K$ проективно полна, система связывающих морфизмов $\Bind({\mathcal N}^X)$, определенная в \eqref{DEF:kategoriya-svyazyv-morpfizmov}, всегда обладает проективным пределом.

Проверим третье условие. Пусть $\alpha:X\to Y$ -- какой-нибудь морфизм и $\varPsi\in 2_{\varPhi_Y}$.
По условию, $\varPhi$ является правым идеалом, поэтому для всякого $\psi\in\varPsi\subseteq\varPhi$ композиция $\psi\circ\alpha$ принадлежит $\varPhi$, и мы можем рассмотреть множество $\varPsi\circ\alpha\in 2_{\varPhi^X}$. Для него мы получим такую диаграмму:
 $$
  \xymatrix @R=2.5pc @C=4.0pc {
 X\ar[r]^{\alpha}\ar@/_3ex/[dr]_{\prod\limits_{\psi\in\varPsi}(\psi\circ\alpha)=\overline{\varPsi\circ\alpha}\quad} & Y\ar[d]^{\overline{\varPsi}=\prod\limits_{\psi\in\varPsi}\psi} & \\
  & \Ran\overline{\varPsi}\ar@{=}[r] & \prod\limits_{\psi\in\varPsi}\Ran\psi
 }
 $$
Пропускаем морфизмы, приходящие в $\Ran\overline{\varPsi}$, через факторизацию \eqref{overline(varPsi)=mu_varPsi-circ-e_(varPsi)}:
 $$
  \xymatrix @R=2.5pc @C=4.0pc {
 X\ar[r]^{\alpha}\ar[d]_{\e_{\varPsi\circ\alpha}} & Y\ar[d]^{\e_{\varPsi}}  \\
 \Ran\e_{\varPsi\circ\alpha}\ar@/_3ex/[dr]_{\mu_{\varPsi\circ\alpha}\quad} & \Ran\e_{\varPsi}\ar[d]^{\mu_{\varPsi}}  \\
  & \Ran\overline{\varPsi}
 }
 $$
Эту диаграмму полезно для наглядности превратить в четырехугольник
 $$
  \xymatrix @R=2.5pc @C=4.0pc {
 X\ar@/^3ex/[dr]^{\e_{\varPsi}\circ\alpha}\ar[d]_{\e_{\varPsi\circ\alpha}} &   \\
 \Ran\e_{\varPsi\circ\alpha}\ar@/_3ex/[dr]_{\mu_{\varPsi\circ\alpha}\quad} & \Ran\e_{\varPsi}\ar[d]^{\mu_{\varPsi}}  \\
  & \Ran\overline{\varPsi}
 }
 $$
после чего можно будет заметить, что в ней $\e_{\varPsi\circ\alpha}\in\varOmega$, а $\mu_{\varPsi}\in\varGamma=\varOmega^\downarrow$, поэтому должна существовать горизонтальная стрелка вправо:
 $$
  \xymatrix @R=2.5pc @C=4.0pc {
 X\ar@/^3ex/[dr]^{\e_{\varPsi}\circ\alpha}\ar[d]_{\e_{\varPsi\circ\alpha}} &   \\
 \Ran\e_{\varPsi\circ\alpha}\ar@/_3ex/[dr]_{\mu_{\varPsi\circ\alpha}\quad}\ar@{-->}[r]^{\delta} & \Ran\e_{\varPsi}\ar[d]^{\mu_{\varPsi}}  \\
  & \Ran\overline{\varPsi}
 }
 $$
Эта стрелка будет недостающей горизонтальной стрелкой в диаграмме \eqref{DIAGR:set}:
 $$
  \xymatrix @R=2.5pc @C=4.0pc {
 X\ar[r]^{\alpha}\ar[d]_{\e_{\varPsi\circ\alpha}} & Y\ar[d]^{\e_{\varPsi}}  \\
 \Ran\e_{\varPsi\circ\alpha}\ar@{-->}[r]^{\delta} & \Ran\e_{\varPsi}
 }
 $$

3. Заметим далее, что в классе $\varOmega$ оболочки относительно классов $\varPhi$, $2_{\varPhi}$, $\{\e_{\varPsi};\ \varPsi\in 2_{\varPhi}\}$ и $\mathcal N$ совпадают:
\beq\label{PROOF:TH:sushestvovanie-seti-pri-faktorizatsii-1}
\env_{\varPhi}^{\varOmega}X=\env_{2_{\varPhi}}^{\varOmega}X=\eqref{env_varPhi^varOmega_X=env_(e_ph;ph_in_varPhi)^varOmega_X}
=\env_{\{\e_{\varPsi};\ \varPsi\in 2_{\varPhi}\}}^{\varOmega}X=\eqref{PROOF:TH:sushestvovanie-seti-pri-faktorizatsii-0}=\env_{\mathcal N}^{\varOmega}X.
\eeq
Здесь второе и третье равенства следуют из \eqref{env_varPhi^varOmega_X=env_(e_ph;ph_in_varPhi)^varOmega_X} и \eqref{PROOF:TH:sushestvovanie-seti-pri-faktorizatsii-0}, и нужно только проверить первое. Для этого достаточно проверить, что понятие расширения относительно классов $\varPhi$ и $2_{\varPhi}$ не меняется. 

Действительно, пусть $\sigma:X\to X'$ --- расширение относительно $\varPhi$. Зафиксируем конечное множество $\varPsi\subseteq \varPhi^X$. Для любого $\ph\in\varPsi$, $\ph:X\to B_\ph$ существует единственная стрелка $\ph':X'\to B_\ph$, замыкающая диаграмму
\beq\label{PROOF:TH:sushestvovanie-seti-pri-faktorizatsii-2}  
\xymatrix @R=2.5pc @C=2.0pc {
 X\ar[rr]^{\sigma}\ar[dr]_{\ph} & & X'\ar@{-->}[dl]^{\ph'}  \\
 & B_\ph & 
 }
 \eeq
Поэтому для стрелки 
$$
\psi=\prod_{\ph\in\varPsi}\ph:X\to \prod_{\ph\in\varPsi} B_\ph
$$
тоже найдется стрелка 
$$
\psi'=\prod_{\ph\in\varPsi}\ph':X'\to \prod_{\ph\in\varPsi} B_\ph
$$
замыкающая диаграмму
 $$
  \xymatrix @R=2.5pc @C=2.0pc {
 X\ar[rr]^{\sigma}\ar[dr]_{\prod_{\ph\in\varPsi}\ph} & & X'\ar@{-->}[dl]^{\prod_{\ph\in\varPsi}\ph'}  \\
 & \prod_{\ph\in\varPsi} B_\ph & 
 }
 $$
И эта стрелка будет единственной, потому что $\sigma$ --- эпиморфизм.

4. После этого повторяется доказательство теоремы \ref{TH:funktorialnost-pri-seti-Epi-i-dolonyaemosti}, начиная с равенства \eqref{PROOF:TH:funktorialnost-pri-seti-Epi-i-dolonyaemosti-1}.
Сначала из равенства \eqref{obolochka-konusa-v-kat-s-uzl-razl} в лемме
\ref{LM:obolochka-konusa-v-kat-s-uzl-razl} мы выводим \eqref{im_infty-lim F_X=env_M^Epi-X-1-*}:
$$
\env_{\varPhi}^{\varOmega}X=\env_{\mathcal N}^{\varOmega}X=\e_{\projlim{\mathcal
N}_X}.
$$
А потом теми же рассуждениями доказываем свойства (b) и (c).
\epr

\paragraph{Полные объекты.}

\bprop\label{PROP:harakterizatsiya-polnoty}
Пусть $\varOmega\subseteq\Epi$, тогда для объекта $A\in\Ob({\tt K})$ следующие условия эквивалентны:
 \bit{
\item[(i)] всякое расширение $\sigma:A\to A'$ в классе $\varOmega$ относительно класса $\varPhi$ является изоморфизмом;

\item[(ii)] локальная единица $1_A:A\to A$ является оболочкой $A$ в классе $\varOmega$ относительно класса $\varPhi$;

\item[(iii)] существует оболочка объекта $A$ в классе $\varOmega$ относительно класса $\varPhi$, являющаяся изоморфизмом: $\env^\varOmega_\varPhi A\in\Iso$.
 }\eit\noindent
\eprop

 \bit{
\item[$\bullet$]\label{DEF:polnye-objekty} Объект $A$ категории $\tt K$ мы будем называть {\it полным} в классе морфизмов $\varOmega\subseteq\Epi$ относительно класса морфизмов $\varPhi$, или {\it полным относительно оболочки $\Env_\varPhi^\varOmega$}, если он удовлетворяет условиям (i)-(iii) этого предложения.
 }\eit

\bpr
1. (i)$\Longrightarrow$(ii). Пусть всякое расширение $\sigma:A\to A'$ является изоморфизмом. Тогда для локальной единицы $1_A:A\to A$ (которая тоже является расширением) мы получим коммутативную диаграмму
$$
\xymatrix @R=2.pc @C=2.0pc 
{
& A \ar[dl]_{\sigma} \ar[dr]^{1_A}\\
A'\ar[rr]_{\sigma^{-1}} & & A
}
$$
которую можно понимать как диаграмму \eqref{DEF:diagr-obolochka}, и это означает, что $1_A$ будет оболочкой.

2. Импликация (ii)$\Longrightarrow$(iii) очевидна.

3. (iii)$\Longrightarrow$(i). Пусть $\rho:A\to E$ -- какая-то оболочка, являющаяся изоморфизмом. Тогда для всякого расширения $\sigma:A\to A'$ выбрав морфизм $\upsilon$ в диаграмме \ref{DEF:diagr-obolochka} мы получим $\upsilon\circ\sigma=\rho\in\Iso$, поэтому $\sigma$ -- коретракция. С другой стороны, $\sigma\in\varOmega\subseteq\Epi$, поэтому $\sigma\in\Iso$.
\epr

\bprop\label{PROP:rasshirenie-v-polnyj-objekt}
Пусть классы морфизмов $\varOmega$ и $\varPhi$ определяют в категории $\tt K$
полурегулярную оболочку. Тогда если в расширении $\sigma:A\to S$ объекта $A$ (в классе $\varOmega$ относительно класса $\varPhi$) объект $S$ полный (в классе $\varOmega$ относительно класса $\varPhi$), то расширение $\sigma:A\to S$ является оболочкой объекта $A$ (в классе $\varOmega$ относительно класса $\varPhi$).
\eprop
\bpr
Пусть $\rho:A\to E$ --- оболочка объекта $A$. Тогда существует единственный морфизм $\upsilon:S\to E$, замыкающий диаграмму 
\eqref{DEF:diagr-obolochka}, которая выглядит в данном случае так: 
 \beq\label{rasshirenie-v-polnyj-objekt}
\xymatrix @R=2.pc @C=2.0pc 
{
& A \ar[dl]_{\sigma} \ar[dr]^{\rho} &\\
S \ar@{-->}[rr]^{\upsilon} & & E
}
 \eeq
Покажем, что $\upsilon$ --- изоморфизм. Для этого по предложению \ref{PROP:harakterizatsiya-polnoty}(i) достаточно проверить, что $\upsilon$ --- расширение объекта $S$ (в классе $\varOmega$ относительно класса $\varPhi$). Выберем произвольный морфизм $\psi:S\to B$, $\psi\in\varPhi$. По условию RE.4 на с.\pageref{DEF:RE.4}, класс $\varPhi$ является правым идеалом. Поэтому морфизм $\ph=\psi\circ\sigma$ должен тоже лежать в $\varPhi$. Поскольку $\rho:A\to E$ --- оболочка (и значит, расширение) объекта $A$, морфизм $\ph$ должен единственным образом продолжаться до некоторого морфизма $\ph':E\to B$. Мы можем теперь дополнить диаграмму \eqref{rasshirenie-v-polnyj-objekt} до диаграммы   
 \beq\label{rasshirenie-v-polnyj-objekt-1}
\xymatrix @R=3.pc @C=4.0pc 
{
& A\ar@{-->}[d]^{\ph} \ar@/_4ex/[ddl]_{\sigma} \ar@/^4ex/[ddr]^{\rho} &\\
& B & \\
S\ar[ur]^{\psi} \ar[rr]^{\upsilon} & & E \ar@{-->}[ul]_{\ph'}
}
 \eeq
Здесь периметр коммутативен, потому что это диаграмма \eqref{rasshirenie-v-polnyj-objekt}, левый внутренний треугольник коммутативен, потому что так мы выбирали $\ph$, а правый внутренний треугольник коммутативен, потому что $\ph'$ --- это продолжение $\ph$ вдоль расширения $\rho$. Поскольку $\sigma\in\varOmega$ --- эпиморфизм (в силу замечания \ref{REM:env-polureg-obolochki}), мы получаем, что нижний внутренний треугольник в \eqref{rasshirenie-v-polnyj-objekt-1} тоже коммутативен.  
\beq\label{rasshirenie-v-polnyj-objekt-2}
\xymatrix @R=3.pc @C=4.0pc 
{
& B & \\
S\ar[ur]^{\psi} \ar[rr]^{\upsilon} & & E \ar@{-->}[ul]_{\ph'}
}
 \eeq
Это значит, что морфизм $\psi$ продолжается до морфизма $\ph'$ вдоль морфизма $\upsilon$. Такое продолжение единственно, потому что 
если нам задано какое-то продолжение $\ph'$ в диаграмме \eqref{rasshirenie-v-polnyj-objekt-2}, то в диаграмме \eqref{rasshirenie-v-polnyj-objekt-1} получаются коммутативны (помимо периметра) нижний внутренний треугольник и левый внутренний треугольник. Отсюда следует, что правый внутренний треугольник тоже должен быть коммутативен:
$$
\ph'\circ\rho=\ph'\circ\upsilon\circ\sigma=\psi\circ\sigma=\ph.
$$
Мы получаем, что $\ph'$ должно быть продолжением морфизма $\ph$ вдоль $\rho$. А такое продолжение единственно.

Итак, $\upsilon$ --- расширение объекта $S$ (в классе $\varOmega$ относительно класса $\varPhi$), и по предложению \ref{PROP:harakterizatsiya-polnoty}(i), это означает, что $\upsilon$ --- изоморфизм.
\epr

\btm\label{TH:kriterij-obolochki-v-term-polnyh-objektov}
Пусть $\varPhi$ --- класс морфизмов в категории $\tt K$, являющийся правым идеалом
\beq\label{kriterij-obolochki-v-term-polnyh-objektov-1}
\varPhi\circ\Mor({\tt K})\subseteq\varPhi,
\eeq
и обладающий сетью эпиморфизмов ${\mathcal N}$, порождающей $\varPhi$ изнутри:
\beq\label{kriterij-obolochki-v-term-polnyh-objektov-2}
    {\mathcal N}\subseteq\varPhi\subseteq\Mor({\tt K})\circ {\mathcal N}.
\eeq
И пусть $\varOmega$ --- мономорфно дополняемый\footnote{В смысле определения на с.\pageref{DEF:klass-monomorfno-dopolnyaem}.} в $\tt K$ класс эпиморфизмов:
\beq\label{kriterij-obolochki-v-term-polnyh-objektov-3}
{^\downarrow\varOmega}\circledcirc\varOmega={\tt K}.
\eeq
Тогда для всякого морфизма\footnote{Здесь не предполагается, что $\eta$ лежит в $\varOmega$.} $\eta:A\to S$ в произвольный полный объект $S$ (в классе $\varOmega$ относительно класса $\varPhi$) следующие условия эквивалентны:
\bit{

\item[(i)] морфизм $\eta:A\to S$ является оболочкой (в классе $\varOmega$ относительно класса $\varPhi$);

\item[(ii)] для любого морфизма\footnote{Здесь не предполагается, что $\theta$ лежит в $\varOmega$.} $\theta:A\to E$ в произвольный полный объект $E$ (в классе $\varOmega$ относительно класса $\varPhi$) найдется единственный морфизм $\upsilon:S\to E$, замыкающий диаграмму
 \beq\label{kriterij-obolochki-v-term-polnyh-objektov}
\xymatrix @R=2.pc @C=2.0pc 
{
& A \ar[dl]_{\eta} \ar[dr]^{\theta} &\\
S \ar@{-->}[rr]^{\upsilon} & & E
}
 \eeq
}\eit
\etm
\bpr Сразу заметим, что по теореме \ref{TH:funktorialnost-pri-seti-Epi-i-dolonyaemosti} в этой ситуации классы морфизмов $\varPhi$ и $\varOmega$ порождают оболочку как функтор. 

1. (i)$\Rightarrow$(ii). Пусть сначала $\eta:A\to S$ --- оболочка, и пусть $\theta:A\to E$ --- произвольный морфизм в полный объект $E$. По предложению \ref{PROP:harakterizatsiya-polnoty}(ii), $E$ является оболочкой для самого себя:
$$
E=\Env_{\varPhi}^\varOmega E.
$$
Поэтому по теореме \ref{TH:sushestvovanie-seti-pri-faktorizatsii}, $E$ изоморфна области значений морфизма $\e_{\projlim{\mathcal N}^E}$ локального предела $\projlim{\mathcal N}^E$:
 \beq\label{PROOF:kriterij-obolochki-v-term-polnyh-objektov-1}
\xymatrix @R=2.pc @C=2.0pc 
{
& E \ar@{=>}[dl]_{\e_{\projlim{\mathcal N}^E}} \ar[dr]^{\projlim{\mathcal N}^E} &\\
\Env_{\varPhi}^\varOmega E\ar[rr]^{\mu_{\projlim{\mathcal N}^E}} & & E_{\mathcal N}
}
 \eeq
(двойная стрелка $\Longrightarrow$, по-прежнему, означает изоморфизм). 

Теперь рассмотрим какой-нибудь морфизм $\rho\in {\mathcal N}^E$. Он порождает диаграмму
 \beq\label{PROOF:kriterij-obolochki-v-term-polnyh-objektov-2}
\xymatrix @R=2.pc @C=2.0pc 
{
& A \ar[dl]_{\eta} \ar[dr]^{\theta}\ar@{-->}[dd]_{\rho\circ\theta} &\\
S\ar@{-->}[dr]_{(\rho\circ\theta)'}  & & E\ar[dl]^{\rho} \\
& \Ran(\rho) & 
}
 \eeq
Здесь поскольку $\rho\in {\mathcal N}^E$, по условию \eqref{kriterij-obolochki-v-term-polnyh-objektov-2} мы получаем, что 
$\rho\in \varPhi$. Отсюда в силу \eqref{kriterij-obolochki-v-term-polnyh-objektov-2} следует, что 
$\rho\circ\theta\in \varPhi$. Как следствие, в диаграмме \eqref{PROOF:kriterij-obolochki-v-term-polnyh-objektov-2} существует единственная стрелка $(\rho\circ\theta)'$, замыкающая левый внутренний треугольник.

Теперь если мы возьмем другой элемент сети эпиморфизмов $\sigma\in {\mathcal N}^E$, связанный с $\rho$ отношением
$$
\rho\to\sigma
$$
то есть такой, что для некоторого $\iota^\sigma_\rho\in \Bind({\mathcal N}^E)$ выполняется равенство \eqref{DEF:le-in-F_X-*}
 $$
\sigma=\iota^\sigma_\rho\circ\rho,
 $$
то мы получим диаграмму
 \beq\label{PROOF:kriterij-obolochki-v-term-polnyh-objektov-3}
\xymatrix @R=2.pc @C=4.0pc 
{
& A \ar[dl]_{\eta} \ar[dr]^{\theta}\ar@{-->}[dd]_{\rho\circ\theta}\ar@{-->}@/^5ex/[ddd]^(.4){\sigma\circ\theta} &\\
S\ar@{-->}[dr]_{(\rho\circ\theta)'}\ar@{-->}@/_4ex/[ddr]_{(\sigma\circ\theta)'}  & & E\ar[dl]^{\rho}\ar@/^4ex/[ddl]^{\sigma} \\
& \Ran(\rho)\ar[d]_{\iota^\sigma_\rho} & \\
& \Ran(\sigma) & \\}
 \eeq
В ней по построению все внутренние треугольники коммутативны, кроме левого нижнего, для которого коммутативность нужно доказывать:
 \beq\label{PROOF:kriterij-obolochki-v-term-polnyh-objektov-4}
\xymatrix @R=2.pc @C=4.0pc 
{
S\ar@{-->}[dr]_{(\rho\circ\theta)'}\ar@{-->}@/_4ex/[ddr]_{(\sigma\circ\theta)'}  & &  \\
& \Ran(\rho)\ar[d]_{\iota^\sigma_\rho} & \\
& \Ran(\sigma) & \\}
 \eeq
Для этого нужно подрисовать к нему стрелку $\eta$
 \beq\label{PROOF:kriterij-obolochki-v-term-polnyh-objektov-5}
\xymatrix @R=2.pc @C=4.0pc 
{
& A \ar[dl]_{\eta}  & \\
S\ar@{-->}[dr]_{(\rho\circ\theta)'}\ar@{-->}@/_4ex/[ddr]_{(\sigma\circ\theta)'}  & &  \\
& \Ran(\rho)\ar[d]_{\iota^\sigma_\rho} & \\
& \Ran(\sigma) & \\}
 \eeq
и заметить, что два пути, по которым можно попасть из начала в конец этой диаграммы в исходной диаграмме \eqref{PROOF:kriterij-obolochki-v-term-polnyh-objektov-3} превращаются один в другой без требования коммутативности треугольника 
\eqref{PROOF:kriterij-obolochki-v-term-polnyh-objektov-4}:
$$
(\sigma\circ\theta)'\circ\eta=\iota^\sigma_\rho\circ(\rho\circ\theta)'\circ\eta
$$
Поскольку $\eta\in\varOmega$ --- эпиморфизм, мы можем на него сократить,
$$
(\sigma\circ\theta)'=\iota^\sigma_\rho\circ(\rho\circ\theta)'
$$
--- и это равенство будет означать коммутативность треугольника \eqref{PROOF:kriterij-obolochki-v-term-polnyh-objektov-4}.

Далее, коммутативность треугольника \eqref{PROOF:kriterij-obolochki-v-term-polnyh-objektov-4} в свою очередь означает, что 
система морфизмов $\{(\sigma\circ\theta)';\ \sigma\in {\mathcal N}^E \}$ представляет собой проективный конус ковариантной системы $\Bind({\mathcal N}^E)=\{\iota^\sigma_\rho\}$. Поэтому существует единственный морфизм в проективный предел этой системы, который на с.\pageref{DIAGR:sigma-sigma_F} мы условились обозначать символом $E_{\mathcal N}$
$$
\omega:S\to E_{\mathcal N}
$$
замыкающий все диаграммы
\beq\label{PROOF:kriterij-obolochki-v-term-polnyh-objektov-6}
\xymatrix @R=2.5pc @C=2.0pc {
   S\ar[dr]_{(\sigma\circ\theta)'}\ar@{-->}[rr]^{\omega} & & E_{\mathcal N}\ar[ld]^{\sigma_{\mathcal N}}  \\
   & \Ran\sigma  &
}\qquad (\sigma\in{\mathcal N}^E). 
\eeq

Теперь рассмотрим систему морфизмов
$$
(\sigma\circ\theta)'\circ\eta=\sigma\circ\theta:A\to\Ran\sigma
$$
(равенство следует из коммутативности диаграммы \eqref{PROOF:kriterij-obolochki-v-term-polnyh-objektov-3}).
Она является проективным конусом ковариантной системы $\Bind({\mathcal N}^E)=\{\iota^\sigma_\rho\}$ с вершиной $A$. Поэтому существует единственный морфизм, 
$$
\alpha:A\to E_{\mathcal N}
$$ 
замыкающий все диаграммы 
\beq\label{PROOF:kriterij-obolochki-v-term-polnyh-objektov-7}
\xymatrix @R=2.5pc @C=2.0pc {
   A\ar[dr]_{(\sigma\circ\theta)'\circ\eta}\ar@{-->}[rr]^{\alpha} & & E_{\mathcal N}\ar[ld]^{\sigma_{\mathcal N}}  \\
   & \Ran\sigma  &
}\qquad (\sigma\in{\mathcal N}^E). 
\eeq
и одновременно все диаграммы
\beq\label{PROOF:kriterij-obolochki-v-term-polnyh-objektov-8}
\xymatrix @R=2.5pc @C=2.0pc {
   A\ar[dr]_{\sigma\circ\theta}\ar@{-->}[rr]^{\alpha} & & E_{\mathcal N}\ar[ld]^{\sigma_{\mathcal N}}  \\
   & \Ran\sigma  &
}\qquad (\sigma\in{\mathcal N}^E). 
\eeq
Диаграммы \eqref{PROOF:kriterij-obolochki-v-term-polnyh-objektov-7} и \eqref{PROOF:kriterij-obolochki-v-term-polnyh-objektov-6} вместе означают что коммутативен треугольник
 \beq\label{PROOF:kriterij-obolochki-v-term-polnyh-objektov-9}
\xymatrix @R=2.pc @C=4.0pc 
{
& A\ar[dd]_{\alpha} \ar[dl]_{\eta}  &\\
S\ar[dr]_{\omega}  & &  \\
& E_{\mathcal N} & \\}
 \eeq
А диаграмма \eqref{PROOF:kriterij-obolochki-v-term-polnyh-objektov-8} вместе с диаграммой 
\eqref{DIAGR:sigma-sigma_F} для объекта $E$
\beq\label{PROOF:kriterij-obolochki-v-term-polnyh-objektov-10}
\xymatrix @R=2.5pc @C=2.0pc {
   X\ar[dr]_{\sigma}\ar@{-->}[rr]^{\projlim{\mathcal N}^X} & & X_{\mathcal N}\ar[ld]^{\sigma_{\mathcal N}}  \\
   & \Ran\sigma  &
}\qquad (\sigma\in{\mathcal N}^X). \eeq
означает, что коммутативен треугольник
 \beq\label{PROOF:kriterij-obolochki-v-term-polnyh-objektov-11}
\xymatrix @R=2.pc @C=4.0pc 
{
& A\ar[dd]_{\alpha}  \ar[dr]^{\theta} &\\
 & & E\ar[dl]^{\projlim {\mathcal N}^X} \\
& E_{\mathcal N} & \\}
 \eeq
 Склеив \eqref{PROOF:kriterij-obolochki-v-term-polnyh-objektov-9} и \eqref{PROOF:kriterij-obolochki-v-term-polnyh-objektov-11} мы получим диаграмму
 \beq\label{PROOF:kriterij-obolochki-v-term-polnyh-objektov-12}
\xymatrix @R=2.pc @C=4.0pc 
{
& A\ar[dd]_{\alpha} \ar[dl]_{\eta} \ar[dr]^{\theta} &\\
S\ar[dr]_{\omega}  & & E\ar[dl]^{\projlim {\mathcal N}^X} \\
& E_{\mathcal N} & \\}
 \eeq
Теперь выбросим из нее стрелку $\alpha$:
 \beq\label{PROOF:kriterij-obolochki-v-term-polnyh-objektov-12}
\xymatrix @R=2.pc @C=4.0pc 
{
& A \ar[dl]_{\eta} \ar[dr]^{\theta} &\\
S\ar[dr]_{\omega}  & & E\ar[dl]^{\projlim {\mathcal N}^X} \\
& E_{\mathcal N} & \\}
 \eeq
И разложим морфизм $\projlim {\mathcal N}^X$ по составляющим лежащим в $\varOmega$ и ${^\downarrow\varOmega}$:
$$
\projlim {\mathcal N}^X=\mu\circ\e.
$$
Мы получим такую диаграмму:
 \beq\label{PROOF:kriterij-obolochki-v-term-polnyh-objektov-13}
\xymatrix @R=2.pc @C=4.0pc 
{
& A \ar[dl]_{\eta\in\varOmega} \ar[dr]^{\e\circ\theta} &\\
S\ar[dr]_{\omega}  & & \Ran\e\ar[dl]^{\mu\in {^\downarrow\varOmega}} \\
& E_{\mathcal N} & \\}
 \eeq
Поскольку правое нижнее ребро $\mu$ принадлежит классу ${^\downarrow\varOmega}$, мономорфно сопряженному классу $\varOmega$, которому принадлежит ребро $\eta$, мы получаем, что долдна существовать диагональ:
 \beq\label{PROOF:kriterij-obolochki-v-term-polnyh-objektov-14}
\xymatrix @R=2.pc @C=4.0pc 
{
& A \ar[dl]_{\eta\in\varOmega} \ar[dr]^{\e\circ\theta} &\\
S\ar[dr]_{\omega}\ar@{-->}[rr]_{\upsilon}  & & \Ran\e\ar[dl]^{\mu\in {^\downarrow\varOmega}} \\
& E_{\mathcal N} & \\}
 \eeq
Теперь выбросим нижнюю вершину:
 \beq\label{PROOF:kriterij-obolochki-v-term-polnyh-objektov-15}
\xymatrix @R=2.pc @C=4.0pc 
{
& A \ar[dl]_{\eta} \ar[dr]^{\e\circ\theta} &\\
S\ar@{-->}[rr]_{\upsilon}  & & \Ran\e \\
}
 \eeq
и заметим, что в силу \eqref{im_infty-lim-F_X=env_M^Epi-X-1}, $\Ran\e=\Env_\varPhi^\varOmega E=E$ (второе равенство выполняется, потому что $E$ --- полный объект). И одновременно $\e=1_E$. Поэтому диаграмма \eqref{PROOF:kriterij-obolochki-v-term-polnyh-objektov-15} упрощается до диаграммы
 \beq\label{PROOF:kriterij-obolochki-v-term-polnyh-objektov-16}
\xymatrix @R=2.pc @C=4.0pc 
{
& A \ar[dl]_{\eta} \ar[dr]^{\theta} &\\
S\ar@{-->}[rr]_{\upsilon}  & & E \\
}
 \eeq
Морфизм $\upsilon$ здесь единственный, потому что $\eta$ --- эпиморфизм.

2. (ii)$\Rightarrow$(i). Пусть выполняется условие (ii), то есть для любого морфизма $\theta:A\to E$ в произвольный полный объект $E$ найдется единственный морфизм $\upsilon:S\to E$, замыкающий диаграмму \eqref{kriterij-obolochki-v-term-polnyh-objektov}:
 \beq\label{PROOF:kriterij-obolochki-v-term-polnyh-objektov-17}
\xymatrix @R=2.pc @C=2.0pc 
{
& A \ar[dl]_{\eta} \ar[dr]^{\theta} &\\
S \ar@{-->}[rr]^{\upsilon} & & E
}
\eeq
Тогда, в частности, для морфизма оболочки $\env_\varOmega^\varPhi A:A\to \Env_\varOmega^\varPhi A$ мы получаем, что долден существовать единственный морфизм $\upsilon:S\to \Env_\varOmega^\varPhi A$, замыкающий диаграмму
 \beq\label{PROOF:kriterij-obolochki-v-term-polnyh-objektov-18}
\xymatrix @R=2.pc @C=2.0pc 
{
& A \ar[dl]_{\eta} \ar[dr]^{\env_\varOmega^\varPhi A} &\\
S \ar@{-->}[rr]^{\upsilon} & & \Env_\varOmega^\varPhi A
}
\eeq
А с другой стороны, по уже доказанному в пункте 1, должен существовать единственный морфизм $\xi:S\gets \Env_\varOmega^\varPhi A$, замыкающий диаграмму
 \beq\label{PROOF:kriterij-obolochki-v-term-polnyh-objektov-19}
\xymatrix @R=2.pc @C=2.0pc 
{
& A \ar[dl]_{\eta} \ar[dr]^{\env_\varOmega^\varPhi A} &\\
S  & & \Env_\varOmega^\varPhi A \ar@{-->}[ll]^{\xi}
}
\eeq
Теперь мы получим, во-первых,
$$
\upsilon\circ\xi\circ\env_\varOmega^\varPhi A=\upsilon\circ\eta=\env_\varOmega^\varPhi A=1_{\Env_\varOmega^\varPhi A}\circ\env_\varOmega^\varPhi A
$$
и, поскольку $\env_\varOmega^\varPhi A$ --- эпиморфизм, мы можем на него сократить, и у нас получится равенство
 \beq\label{PROOF:kriterij-obolochki-v-term-polnyh-objektov-20}
\upsilon\circ\xi=1_{\Env_\varOmega^\varPhi A}.
\eeq
А, во-вторых, склеив диаграммы \eqref{PROOF:kriterij-obolochki-v-term-polnyh-objektov-18} и \eqref{PROOF:kriterij-obolochki-v-term-polnyh-objektov-19} таким образом
 \beq\label{PROOF:kriterij-obolochki-v-term-polnyh-objektov-21}
\xymatrix @R=2.pc @C=4.0pc 
{
& A \ar[dl]_{\eta} \ar[d]^{\env_\varOmega^\varPhi A} \ar[dr]^{\eta} &\\
S \ar[r]^{\upsilon} & \Env_\varOmega^\varPhi A\ar[r]^{\xi} & S
}
\eeq
мы получаем, что морфизм $\xi\circ\upsilon$ замыкает диаграмму 
 \beq\label{PROOF:kriterij-obolochki-v-term-polnyh-objektov-22}
\xymatrix @R=2.pc @C=4.0pc 
{
& A \ar[dl]_{\eta}  \ar[dr]^{\eta} &\\
S \ar@{-->}[rr]^{\xi\circ\upsilon} &  & S
}
\eeq
и, поскольку по условию (ii) такой морфизм должен быть единственным, он должен совпадать с $1_S$:
 \beq\label{PROOF:kriterij-obolochki-v-term-polnyh-objektov-24}
\xi\circ\upsilon=1_S.
\eeq
Условия \eqref{PROOF:kriterij-obolochki-v-term-polnyh-objektov-20} и \eqref{PROOF:kriterij-obolochki-v-term-polnyh-objektov-24} вместе означают, что $\upsilon$ --- изоморфизм, и поэтому морфизм $\eta:A\to S$ --- оболочка.
\epr

Обозначим через $\tt L$ класс всех полных объектов в $\tt K$ (в $\varOmega$ относительно $\varPhi$).
Будем рассматривать $\tt L$ как полную подкатегорию в $\tt K$.

\bprop\label{PROP:obolochka-na-L}
В условиях теоремы \ref{TH:sushestvovanie-seti-pri-faktorizatsii} строящийся в ней
функтор оболочки $(X,\alpha)\mapsto (E(X),E(\alpha))$ на подкатегории полных объектов ${\tt L}\subseteq {\tt K}$ изоморфен тождественному функтору:
\beq\label{E:K->L=funktor}
\forall A\in {\tt L}\qquad E(A)\cong A, \qquad \forall \alpha:\underset{\scriptsize\begin{matrix}\text{\rotatebox{90}{$\owns$}}\\ {\tt L}\end{matrix}}{A}\to \underset{\scriptsize\begin{matrix}\text{\rotatebox{90}{$\owns$}}\\ {\tt L}\end{matrix}}{A'}\qquad
E(\alpha)=e_{A'}\circ\alpha\circ e_A^{-1}.
\eeq
\eprop
\bpr
Пусть $\alpha:A\to A'$ -- произвольный морфизм в $\tt L$, то есть морфизм в $\tt K$, у которого область определения и область значений лежат в $\tt L$. Тогда в диаграмме \eqref{DIAGR:funktorialnost-env-e-E} горизонтальные стрелки являются изоморфизмами, поэтому
$$
\xymatrix @R=2.pc @C=5.0pc 
{
A\ar[d]^{\alpha} & A\ar[d]^{E(\alpha)}\ar[l]_{e_A^{-1}} \\
A'\ar[r]^{e_{A'}} & A' \\
}
$$
\epr

\paragraph{Подталкивание и регулярная оболочка.}

\bit{
\item[RE.5:]\label{DEF:RE.5} Условимся говорить, что {\it класс морфизмов $\varOmega$ подталкивает класс морфизмов $\varPhi$}, если
\beq\label{Omega-podderzhivaet-Phi}
\forall\psi\in\Mor({\tt K})\qquad \forall \sigma\in\varOmega\qquad \big(\psi\circ\sigma\in\varPhi\quad\Longrightarrow\quad \psi\in\varPhi\big).
\eeq
}\eit

\brem\label{Omega-podtalkivaet-Mor(K,M)} Очевидно, для выполнения условия \eqref{Omega-podderzhivaet-Phi} достаточно (но необязательно), чтобы $\varPhi$ представлял собой класс морфизмов со значениями в некотором классе объектов $\tt M$ категории $\tt K$:
$$
\varPhi=\{\ph\in\Mor(\tt K):\ \Ran\ph\in{\tt M}\},
$$
\erem

\blm\label{LM:kompozitsiya-rasshirenij} Если класс $\varOmega$ подталкивает класс $\varPhi$, то композиция $\sigma\circ\rho:X\to X''$ любых двух последовательно взятых расширений $\rho:X\to X'$ и $\sigma:X'\to X''$ (в $\varOmega$ относительно $\varPhi$) является расширением (в $\varOmega$ относительно $\varPhi$).
\elm
\bpr
Это видно из следующей диаграммы:
$$
\xymatrix @R=2.pc @C=2.0pc 
{
X\ar@/_2ex/[dr]_{\ph}\ar[r]^{\rho} & X'\ar[r]^{\sigma}\ar@{-->}[d]^{\ph'} & X''\ar@/^2ex/@{-->}[dl]^{\ph''} \\
 & B &
}
$$
Поскольку $\rho$ -- расширение, для $\ph\in\varPhi$ существует стрелка $\ph'$, а поскольку $\varOmega$ подталкивает $\varPhi$, эта стрелка лежит в $\varPhi$. Поскольку $\sigma$ -- расширение, существует стрелка $\ph''$. При этом каждая последующая стрелка однозначно определяется предыдущей.
\epr

\bprop\label{PROP:polnota-pri-podtalkivanii}
Если в категории ${\tt K}$ класс морфимзмов $\varOmega$ подталкивает класс морфизмов $\varPhi$, то объект $A\in\Ob({\tt K})$ тогда и только тогда полон в классе $\varOmega$ относительно класса $\varPhi$, когда он изоморфен оболочке некоторого объекта $X\in\Ob({\tt K})$: 
\beq\label{A=Env^varOmega_varPhi-X}
A\cong\Env^\varOmega_\varPhi X
\eeq
 \eprop
\bpr
1. Необходимость. Если объект $A$ полный, то по предложению \ref{PROP:harakterizatsiya-polnoty}(iii), $\env^\varOmega_\varPhi A\in\Iso$. Поэтому $A\cong\Env^\varOmega_\varPhi A$.

2. Достаточность. Пусть $A\cong\Env^\varOmega_\varPhi X$ для некоторого
$X\in\Ob({\tt K})$. Тогда $A$ можно рассматривать как оболочку для $X$, то есть существует морфизм
$\rho:X\to A$, являющийся оболочкой. Пусть $\sigma:A\to A'$ -- какое-нибудь расширение для $A$.
По лемме \ref{LM:kompozitsiya-rasshirenij}, композиция $\sigma\circ\rho:X\to A'$ является расширением
для $X$, поэтому должен существовать морфизм $\upsilon$, замыкающий диаграмму
\eqref{DEF:diagr-obolochka}:
$$
\xymatrix @R=2.pc @C=2.0pc 
{
& X\ar[dl]_{\sigma\circ \rho}\ar[dr]^{\rho} & \\
A'\ar@{-->}[rr]_{\upsilon} & & A
}
$$
Теперь получаем:
$$
\upsilon\circ\sigma\circ \rho=\rho=1_A\circ \underset{\scriptsize\begin{matrix}\text{\rotatebox{90}{$\owns$}}\\ \Epi\end{matrix}}{\rho}\qquad\Longrightarrow\qquad \upsilon\circ\sigma=1_A.
$$
В последнем равенстве морфизм $\upsilon$ должен быть единственным, потому что $\sigma$ -- эпиморфизм.
Мы получаем, что расширение $\sigma:A\to A'$ подчинено расширению $1_A:A\to A$. Поскольку это верно для всякого
$\sigma$, морфизм $1_A$ должен быть  оболочкой для $A$. То есть $A$ удовлетворяет условию (ii) предложения \ref{PROP:harakterizatsiya-polnoty}, и поэтому является полным объектом.
\epr

 \bit{
\item[$\bullet$]\label{DEF:reg-obolochka} Условимся говорить, что {\it классы морфизмов $\varOmega$ и $\varPhi$ определяют в категории $\tt K$ регулярную оболочку}, или что {\it оболочка $\Env_\varPhi^\varOmega$ регулярна}, если в дополнение к условиям RE.1-RE.4 теоремы \ref{TH:sushestvovanie-seti-pri-faktorizatsii} класс $\varOmega$ подталкивает класс $\varPhi$.
}\eit

 \btm\label{TH:regulyarnaya-obolochka}
Если классы морфизмов $\varOmega$ и $\varPhi$ определяют в категории $\tt K$
регулярную оболочку, то $\Env_\varPhi^\varOmega$ можно определить как идемпотентный функтор.
 \etm
\bpr
Рассмотрим функтор оболочки $E$, строящийся в теореме \ref{TH:sushestvovanie-seti-pri-faktorizatsii}, и обозначим символом ${\tt L}_0$ класс всех объектов, являющихся образами объектов  в $\tt K$ при отображении $X\mapsto E(X)$:
\beq\label{DEF:L_0}
A\in{\tt L}_0\qquad\Longleftrightarrow\qquad \exists X\in\Ob({\tt K})\quad A=E(X).
\eeq
Определим систему изоморфизмов в $\tt K$ правилом
$$
\forall X\in\Ob({\tt K})\qquad \zeta_X=\begin{cases}1_X,& X\notin{\tt L}_0\\
e_X^{-1},& X\in{\tt L}_0 \end{cases}
$$
(это будет корректное определение в силу предложения \ref{PROP:harakterizatsiya-polnoty}). После этого рассмотрим отображения $X\mapsto F(X)$, $X\mapsto f_X$, $\alpha\mapsto F(\alpha)$, определенные правилами:
 \begin{align*}
& \forall X\in\Ob({\tt K})\qquad F(X)=\begin{cases}E(X),& X\notin{\tt L}_0\\
X,& X\in{\tt L}_0 \end{cases}, \qquad f_X=\begin{cases}e_X,& X\notin{\tt L}_0\\
1_X,& X\in{\tt L}_0 \end{cases} \\
& \forall\alpha\in\Mor({\tt K})\qquad F(\alpha)=\zeta_{\Ran E(\alpha)}\circ E(\alpha)\circ\zeta_{\Dom E(\alpha)}^{-1}
 \end{align*}
Связь с функтором $E$ выражается диаграммой
\beq\label{DIAGR:E(E(X))=E(X)}
\xymatrix @R=2.pc @C=5.0pc 
{
 &  & F(X)\ar[ddd]^{F(\alpha)} \\
X\ar[d]^{\alpha}\ar[r]^{e_X} \ar@/^3ex/[rru]^{f_X}  & E(X)\ar[d]^{E(\alpha)}\ar[ru]_{\zeta_X} & \\
Y\ar[r]^{e_Y}\ar@/_3ex/[rrd]_{f_Y} & E(Y)\ar[rd]^{\zeta_Y} & \\
 & & F(Y)
}
\eeq
Для всякого $X$ морфизм $f_X:X\to F(X)$ будет оболочкой объекта $X$, потому что $f_X$ и $e_X$ связаны изоморфизмом $\zeta_X$. Отображение $(X,\alpha)\mapsto (F(X),F(\alpha))$ будет функтором, потому что, во-первых,
 \begin{multline*}
F(\beta\circ\alpha)=\zeta_{\Ran\beta}\circ E(\beta\circ\alpha)\circ\zeta_{\Dom\alpha}^{-1}=
\zeta_{\Ran\beta}\circ E(\beta)\circ E(\alpha)\circ\zeta_{\Dom\alpha}^{-1}=\\=
\zeta_{\Ran\beta}\circ E(\beta)\circ\zeta_{\Dom\beta}^{-1}\circ \zeta_{\Ran\alpha}\circ E(\alpha)
\circ\zeta_{\Dom\alpha}^{-1}=F(\beta)\circ F(\alpha),
 \end{multline*}
и, во-вторых, при $X\notin{\tt L}_0$ диаграмма \eqref{DIAGR:E(E(X))=E(X)} приобретает вид
$$
\xymatrix @R=2.pc @C=5.0pc 
{
 &  & E(X)\ar[ddd]^{E(1_X)} \\
X\ar[d]^{\alpha}\ar[r]^{e_X} \ar@/^3ex/[rru]^{e_X}  & E(X)\ar[d]^{1_{E(X)}}\ar[ru]_{1_{E(X)}} & \\
X\ar[r]^{e_X}\ar@/_3ex/[rrd]_{e_X} & E(X)\ar[rd]^{1_{E(X)}} & \\
 & & E(X)
}
$$
поэтому
$$
F(1_X)=\zeta_X\circ E(1_X)\circ\zeta_X^{-1}=1_{E(X)}^{-1}\circ 1_{E(X)}\circ 1_{E(X)}=1_{E(X)}=1_{F(X)},
$$
а при $X\in{\tt L}_0$ диаграмма \eqref{DIAGR:E(E(X))=E(X)} превращается в диаграмму
$$
\xymatrix @R=2.pc @C=5.0pc 
{
 &  & X\ar[ddd]^{F(1_X)} \\
X\ar[d]^{1_X}\ar[r]^{e_X} \ar@/^3ex/[rru]^{1_X}  & E(X)\ar[d]^{1_{E(X)}}\ar[ru]_{e_X^{-1}} & \\
X\ar[r]^{e_X}\ar@/_3ex/[rrd]_{1_X} & E(X)\ar[rd]^{e_X^{-1}} & \\
 & & X
}
$$
и здесь видно, что если в периметр подставить вместо $F(1_X)$ морфизм $1_X$, то он остается коммутативной диаграммой. Поскольку такая стрелка единственна, мы получаем
$$
F(1_X)=1_X=1_{F(X)}.
$$

Условие \eqref{e_(E(X))=1_(E(X))} выполняется для функтора $F$ по определению: поскольку всегда
$F(X)\in{\tt L}_0$, мы получаем $f_{F(X)}=1_{F(X)}$.
\epr

\btm[описание оболочки через полные объекты]\label{TH:obolochki-i-polnye-obyekty}
Пусть классы морфизмов $\varOmega$ и $\varPhi$ определяют в категории $\tt K$
регулярную оболочку $\Env_\varPhi^\varOmega$. Тогда
 данный морфизм $\rho:X\to A$ является оболочкой (в $\varOmega$ относительно $\varPhi$) если и только если выполняются следующие условия:
  \bit{
\item[(i)] $\rho:X\to A$ -- эпиморфизм,

\item[(ii)] $A$ -- полный объект  (в $\varOmega$ относительно $\varPhi$),

\item[(iii)] для любого полного объекта $B$ (в $\varOmega$ относительно $\varPhi$) и всякого морфизма $\xi:X\to B$ найдется единственный морфизм $\xi':A\to B$, замыкающий диаграмму
 \beq\label{DIAGR:obolochki-i-polnye-obyekty}
\xymatrix @R=2.pc @C=2.0pc 
{
X\ar[rr]^{\rho}\ar[dr]_{\xi} & & A\ar@{-->}[dl]^{\xi'}\\
& B &
}
 \eeq
 }\eit
 \etm
\bpr Пусть $\rho:X\to A$ -- оболочка. Тогда, во-первых, это эпиморфизм, поскольку $\varOmega\subseteq\Epi$.
Во-вторых, в силу предложения \ref{PROP:harakterizatsiya-polnoty}, $A\cong \Env^\varOmega_\varPhi X$ --
полный объект. В-третьих, если $\xi:X\to B$ -- морфизм в полный объект $B$, то для него можно рассмотреть
диаграмму \eqref{DIAGR:funktorialnost-env_varPhi^Epi-v-kat-s-uzl-razl-1-*}, принимающую в данном случае вид
 $$
\xymatrix @R=2.pc @C=6.0pc 
{
X\ar[r]^{\rho=\env^\varOmega_\varPhi X}\ar[d]_{\xi} & A=\Env^\varOmega_\varPhi X\ar[d]^{\Env^\varOmega_\varPhi\xi}\\
B\ar[r]^{\env^\varOmega_\varPhi B} & \Env^\varOmega_\varPhi B
}
$$
В ней морфизм $\env^\varOmega_\varPhi B$ должен быть изоморфизмом, поскольку $B$ -- полный объект,
и, как следствие, определен морфизм
$$
\xi'=(\env^\varOmega_\varPhi B)^{-1}\circ \Env^\varOmega_\varPhi\xi.
$$
Он и будет пунктирной стрелкой в \eqref{DIAGR:obolochki-i-polnye-obyekty}.

Наоборот, пусть выполняются условия (i)--(iii). В наших условиях справедлива теорема
\ref{TH:sushestvovanie-seti-pri-faktorizatsii}, поэтому можно рассмотреть диаграмму
\eqref{DIAGR:funktorialnost-env_varPhi^Epi-v-kat-s-uzl-razl-1-*}:
 $$
\xymatrix @R=2.pc @C=5.0pc 
{
X\ar[d]_{\rho}\ar[r]^{\env^\varOmega_\varPhi X} & \Env^\varOmega_\varPhi X\ar[d]^{\Env^\varOmega_\varPhi\rho}\\
A\ar[r]^{\env^\varOmega_\varPhi {A}} & \Env^\varOmega_\varPhi A
}
$$
В ней морфизм $\env^\varOmega_\varPhi {A}$ должен быть изоморфизмом (потому что $A$ -- полный объект). Поэтому положив $\zeta=\env^\varOmega_\varPhi {A}^{-1}\circ \Env^\varOmega_\varPhi(\rho)$, мы получим коммутативную диаграмму
 $$
\xymatrix @R=2.pc @C=5.0pc 
{
X\ar[d]_{\rho}\ar[r]^{\env^\varOmega_\varPhi X} & \Env^\varOmega_\varPhi X\ar@{-->}[dl]^{\zeta}\\
A &
}
$$
С другой стороны, по предложению \ref{PROP:harakterizatsiya-polnoty}, $\Env^\varOmega_\varPhi X$ -- полный объект, поэтому по свойству (iii), должен существовать морфизм $\eta$, замыкающий диаграмму
 $$
\xymatrix @R=2.pc @C=5.0pc 
{
X\ar[d]_{\rho}\ar[r]^{\env^\varOmega_\varPhi X} & \Env^\varOmega_\varPhi X\\
A\ar@{-->}[ur]_{\eta} &
}
$$
Поскольку в этих диаграммах оба морфизма $\rho$ и $\env^\varOmega_\varPhi X$ -- эпиморфизмы, мы получаем, что $\zeta$ и $\eta$ -- взаимно обратные морфизмы. Таким образом, $\rho=\zeta\circ\env^\varOmega_\varPhi X$, где $\zeta\in\Iso$. В силу \eqref{Iso-circ-varOmega-subseteq-varOmega}, это означает, что $\rho\in\varOmega$, и поэтому является оболочкой.
\epr

\bcor\label{COR:rasshirenija-i-polnye-obyekty}
Пусть классы морфизмов $\varOmega$ и $\varPhi$ определяют в категории $\tt K$
регулярную оболочку. Тогда для всякого расширения $\sigma:X\to X'$ (в $\varOmega$ относительно $\varPhi$), любого полного объекта $B$ (в $\varOmega$ относительно $\varPhi$) и всякого морфизма $\xi:X\to B$ найдется единственный морфизм $\xi':X'\to B$, замыкающий диаграмму
 \beq\label{DIAGR:rasshirenija-i-polnye-obyekty}
\xymatrix @R=2.pc @C=2.0pc 
{
X\ar[rr]^{\sigma}\ar[dr]_{\xi} & & X'\ar@{-->}[dl]^{\xi'}\\
& B &
}
 \eeq
 \ecor
\bpr
Рассмотрим оболочку $\rho:X\to\Env^\varOmega_\varPhi X$ и соответствуюший морфизм $\upsilon:X'\to \Env^\varOmega_\varPhi X$, звмыкающий диаграмму \eqref{DEF:diagr-obolochka}:
$$
\xymatrix @R=2.pc @C=2.0pc 
{
& X\ar[dl]_{\sigma}\ar[dr]^{\rho} & \\
X'\ar@{-->}[rr]_{\upsilon} & & \Env^\varOmega_\varPhi X
}
$$
Теперь мы можем воспользоваться теоремой \eqref{TH:obolochki-i-polnye-obyekty} и построить сначала морфизм $\xi''$, замыкающий периметр диаграммы
\beq\label{DIAGR:rasshirenija-i-polnye-obyekty-1}
\xymatrix @R=2.pc @C=4.0pc 
{
X\ar[r]^{\sigma}\ar[dr]_{\xi}\ar@/^4ex/[rr]^{\rho} & X'\ar[r]^{\upsilon}\ar@{-->}[d]^{\xi'} &  \Env^\varOmega_\varPhi X \ar@{-->}[dl]^{\xi''} \\
& B &
}
\eeq
А затем положить $\xi'=\xi''\circ\upsilon$, и диаграмма \eqref{DIAGR:rasshirenija-i-polnye-obyekty-1} станет коммутативной. Единственность  $\xi'$ следует из того, что $\sigma\in\Epi$.
\epr

\btm\label{TH:proizvedenie-polnyh-obyektov}
Пусть классы морфизмов $\varOmega$ и $\varPhi$ определяют в категории $\tt K$
регулярную оболочку. Тогда
  \bit{
\item[(i)] произведение (если оно существует) $\prod_{i\in I}X_i$ любого семейства $\{X_i;\ i\in I\}$ полных объектов является полным объектом,

\item[(ii)] проективный предел (если он существует) $\projlim_{i\in I}X_i$ любой ковариантной (контравариантной) системы $\{X_i,\iota^j_i;\ i\in I\}$ полных объектов является полным объектом.
 }\eit
 \etm
\bpr
1. Пусть $\{X_i;\ i\in I\}$ --- семейство полных объектов. Рассмотрим какое-нибудь расширение ее произведения
$$
\sigma: \prod_{i\in I}X_i\to Y.
$$
Пусть $\pi_j: \prod_{i\in I}X_i \to X_j$ --- естественная проекция на компоненту $X_j$. Поскольку $X_j$ --- полный объект, по следствию \ref{COR:rasshirenija-i-polnye-obyekty} найдется единственный морфизм $\tau_j:Y\to X_j$, замыкающий диаграмму 
 \beq\label{proizvedenie-polnyh-obyektov-1}
\xymatrix @R=2.pc @C=2.0pc 
{
\prod_{i\in I}X_i\ar[rr]^{\sigma}\ar[dr]_{\pi_j} & & Y\ar@{-->}[dl]^{\tau_j}\\
& X_j &
}
 \eeq
Семейство морфизмов $\tau_j:Y\to X_j$ будет проективным конусом семейства объектов $\{X_i;\ i\in I\}$. Поэтому существует единственный морфизм $\tau:Y\to \prod_{i\in I}X_i$, замыкающий все диаграммы
 \beq\label{proizvedenie-polnyh-obyektov-2}
\xymatrix @R=2.pc @C=2.0pc 
{
\prod_{i\in I}X_i\ar[dr]_{\pi_j} & & Y\ar[dl]^{\tau_j}\ar@{-->}[ll]_{\tau}\\
& X_j &
}
 \eeq
Морфизмы $\sigma$ и $\tau$ будут взаимно обратными, потому что, с одной стороны, из коммутативности \eqref{proizvedenie-polnyh-obyektov-1} и  \eqref{proizvedenie-polnyh-obyektov-2} следует коммутативность всех диаграмм
$$
\xymatrix @R=2.pc @C=2.0pc 
{
\prod_{i\in I}X_i\ar@{-->}[rr]^{\tau\circ\sigma}\ar[dr]_{\pi_j} & & \prod_{i\in I}X_i\ar[dl]^{\pi_j}\\
& X_j &
}
$$
а, поскольку пунктирная стрелка здесь должна быть единственна (по определению произведения $\prod_{i\in I}X_i$), она должна совпадать с морфизмом $1_{\prod_{i\in I}X_i}$:
\beq\label{proizvedenie-polnyh-obyektov-3}
\tau\circ\sigma=1_{\prod_{i\in I}X_i}.
\eeq
Теперь умножив \eqref{proizvedenie-polnyh-obyektov-3} на морфизм $\sigma$ слева, мы получим
$$
\sigma\circ\tau\circ\sigma=\sigma\circ 1_{\prod_{i\in I}X_i}=\sigma=1_Y\circ\sigma.
$$
Поскольку $\sigma$ --- эпиморфизм, мы можем на него сократить, и получим равенство 
$$
\sigma\circ\tau=1_Y.
$$
Мы получили, что $\sigma$ --- изоморфизм, и, поскольку это верно для любого расширения $\sigma$, это значит, что $\prod_{i\in I}X_i$ --- полный объект.

2. Пусть $\{X_i,\iota^j_i;\ i\in I\}$ --- ковариантная система полных объектов.  Рассмотрим какое-нибудь расширение ее проективного предела
$$
\sigma: \projlim_{i\in I}X_i\to Y.
$$
Пусть $\pi_j: \projlim_{i\in I}X_i \to X_j$ --- естественная проекция на компоненту $X_j$. 
Поскольку $X_j$ --- полный объект, по следствию \ref{COR:rasshirenija-i-polnye-obyekty} найдется единственный морфизм $\tau_j:Y\to X_j$, замыкающий диаграмму 
 \beq\label{proizvedenie-polnyh-obyektov-4}
\xymatrix @R=2.pc @C=2.0pc 
{
\projlim_{i\in I}X_i\ar[rr]^{\sigma}\ar[dr]_{\pi_j} & & Y\ar@{-->}[dl]^{\tau_j}\\
& X_j &
}
 \eeq
Покажем, что семейство морфизмов $\tau_j:Y\to X_j$ будет проективным конусом ковариантной системы $\{X_i,\iota^j_i;\ i\in I\}$. Для любых индексоов $j\le k$ мы можем расмотреть диаграмму
 \beq\label{proizvedenie-polnyh-obyektov-4}
\xymatrix @R=2.pc @C=2.0pc 
{
\projlim_{i\in I}X_i\ar[rr]^{\sigma}\ar[dr]_{\pi_j}\ar@/_4ex/[ddr]_{\pi_k} & & Y\ar@{-->}[dl]^{\tau_j}\ar@{-->}@/^4ex/[ddl]^{\tau_k}\\
& X_j\ar[d]_{\iota_j^k} & \\
& X_k &
}
 \eeq
В ней коммутативен периметр, а также верхний внутренний треугольник и правый внутренний треугольник. Отсюда следует равенство
$$
\iota_j^k\circ\tau_j\circ\sigma=\tau_k\circ\sigma,
$$
и, поскольку $\sigma$ --- эпиморфизм, на него можно сократить:
$$
\iota_j^k\circ\tau_j=\tau_k.
$$
Это как раз и значит, что семейство $\tau_j:Y\to X_j$ --- проективный конус.
 Как следствие, существует единственный морфизм $\tau:Y\to \projlim_{i\in I}X_i$, замыкающий все диаграммы
 \beq\label{proizvedenie-polnyh-obyektov-5}
\xymatrix @R=2.pc @C=2.0pc 
{
\projlim_{i\in I}X_i\ar[dr]_{\pi_j} & & Y\ar[dl]^{\tau_j}\ar@{-->}[ll]_{\tau}\\
& X_j &
}
 \eeq
После этого повторяются все рассуждения предыдущего пункта: морфизмы $\sigma$ и $\tau$ должны быть взаимно обратными, поэтому $\sigma$ --- изоморфизм, и, поскольку расширение $\sigma$ выбиралось произвольным, это означает, что $\projlim_{i\in I}X_i$ --- полный объект.
\epr

\paragraph{Функториальность на эпиморфизмах.}

Обозначим символом ${\tt K}^{\Epi}$ подкатегорию в $\tt K$ с тем же классом объектов что и в $\tt K$, но с эпиморфизмами из $\tt K$ в качестве морфизмов:
$$
\Ob({\tt K}^{\Epi})=\Ob({\tt K}),\qquad \Mor({\tt K}^{\Epi})=\Epi({\tt K}).
$$

\btm\label{TH:funktorialnost-v-K^Epi}
 Пусть ${\tt K}$ -- категория с произведениями (над произвольным индексным множеством),
 и в ней заданы классы морфизмов $\varOmega$ и $\varPhi$, удовлетворяющие следующим условиям:
\bit{

\item[---] класс $\varOmega$ мономорфно дополняем в ${\tt K}$,

\item[---] категория ${\tt K}$ локально мала в фактор-объектах класса $\varOmega$,

\item[---] класс $\varPhi$ выходит из\footnote{В смысле определения на с.\pageref{DEF:goes-from}.} ${\tt K}$,

\item[---] классы $\varOmega$ и $\varPhi$ связаны условием: $\varPhi\circ\varOmega\subseteq\varPhi$

}\eit
Тогда
\bit{
\item[(a)] всякий объект $X$ в ${\tt K}$ обладает оболочкой $\Env_{\varPhi}^{\varOmega}X$ в классе $\varOmega$  относительно класса $\varPhi$,

\item[(b)] для всякого эпиморфизма $\pi:X\to Y$ существует единственный морфизм $\Env_{\varPhi}^{\varOmega}\pi:\Env_{\varPhi}^{\varOmega}X\to \Env_{\varPhi}^{\varOmega}Y$, замыкающий диаграмму
\beq\label{DIAGR:funktorialnost-v-K^Epi}
\xymatrix @R=2.pc @C=5.0pc 
{
X\ar[d]^{\pi}\ar[r]^{\env_{\varPhi}^{\varOmega} X} & \Env_{\varPhi}^{\varOmega} X\ar@{-->}[d]^{\Env_{\varPhi}^{\varOmega}\pi} \\
Y\ar[r]^{\env_{\varPhi}^{\varOmega} Y} & \Env_{\varPhi}^{\varOmega} Y \\
}
\eeq

\item[(c)] оболочку $\Env_{\varPhi}^\varOmega$ можно определить как функтор из ${\tt K}^{\Epi}$ в ${\tt K}^{\Epi}$.
}\eit
\etm

Для доказательства нам понадобится следующая

\blm\label{LM:Env_varPhi-circ-e^Omega-X=Env_varPhi^Omega-Y}  Если ${\tt K}$ -- категория с произведениями
(над произвольным индексным множеством), локально малая в фактор-объектах класса $\varOmega$, и
$\varOmega$ мономорфно дополняем в ${\tt K}$, то для любого класса морфизмов $\varPhi$ и всякого
эпиморфизма $\pi:X\to Y$ справедлива формула:
    \beq\label{Env_varPhi-circ-e^Omega-X=Env_varPhi^Omega-Y}
    \Env_{\varPhi\circ\pi}^\varOmega X=\Env_{\varPhi}^\varOmega Y
    \eeq
    (оболочка $X$ в ${\varOmega}$ относительно класса морфизмов $\varPhi\circ\pi=\{\ph\circ\pi;\ \ph\in\varPhi\}$ совпадает с оболочкой $Y$ в ${\varOmega}$ относительно класса морфизмов $\varPhi$).
\elm
\bpr Заметим сразу, что существование оболочек в \eqref{Env_varPhi-circ-e^Omega-X=Env_varPhi^Omega-Y} гарантируется свойством $5^\circ$ на с.\pageref{5^0:obolochka-otn-klassa-morphizmov}. Вдобавок, по свойству
$5^\circ$ на с.\pageref{PROP:deistvie-epimorfizma-na-Env}, должен существовать морфизм $\upsilon$, замыкающий диаграмму
\eqref{deistvie-epimorfizma-na-Env}:
$$
\xymatrix 
{
 & X\ar[ld]_{\env_{\varPhi}^{\varOmega} Y\circ\pi}\ar[rd]^{\env_{\varPhi\circ\pi}^{\varOmega} X} & \\
 \Env_{\varPhi}^{\varOmega} Y\ar@{-->}[rr]^{\upsilon} &   & \Env_{\varPhi\circ\pi}^{\varOmega} X
}
$$
Покажем, что существует и обратный морфизм. Для этого рассмотрим оболочку $\env_{\varPhi}^{\varOmega} Y: Y\to \Env_{\varPhi}^{\varOmega} Y$. Представим ее сначала как оболочку относительно некоторого множества морфизмов $M$, как в доказательстве свойства $5^\circ$ на с.\pageref{5^0:obolochka-otn-klassa-morphizmov}, а затем, как в доказательстве свойства $3^\circ$ на с.\pageref{3^0:obolochka-otn-mnozhestva-morphizmov}, заменим множество $M$ на единственный морфизм $\psi=\prod_{\chi\in M}\chi$. Тогда по свойству $1^\circ$ на с.\pageref{1^0:obolochka-otn-1-morphizma}, оболочка относительно $\psi$ будет описываться как эпиморфизм $\e_\psi$ в факторизации $\psi$:
$$
\env_{\varPhi}^{\varOmega} Y=\env_{M}^{\varOmega} Y=\env_{\psi}^{\varOmega} Y=\e_\psi.
$$
Мы получаем диаграмму
$$
\xymatrix 
{
X\ar@/^2ex/[rd]^(.7){\env_{\varPhi}^{\varOmega} Y\circ\pi=\e_\psi\circ\pi}
\ar@/^2ex/[rrr]^{\env_{\varPhi\circ\pi}^{\varOmega} X}\ar[dd]_{\pi} & &  & \Env_{\varPhi\circ\pi}^{\varOmega} X\ar[dd]^{(\psi\circ\pi)'}\ar@{-->}@/_2ex/[dl]_(.7){\delta} \\
 & \Env_{\varPhi}^{\varOmega} Y\ar@{=}[r] & \Ran\e_\psi\ar@/_2ex/[dr]_(.3){\mu_\psi} & \\
Y\ar@/_2ex/[rrr]_{\psi}\ar@/_2ex/[ru]_(.7){\env_{\varPhi}^{\varOmega} Y=\e_\psi}  & & & B\\
}
$$
в которой $(\psi\circ\pi)'$ -- продолжение морфизма $\psi\circ\pi\in\varPhi\circ\pi$ вдоль
оболочки $\env_{\varPhi\circ\pi}^{\varOmega} X$. Здесь существование морфизма $\delta$ следует
из того, что $\env_{\varPhi\circ\pi}^{\varOmega} X\in\varOmega$, а $\mu_\psi\in{^\downarrow\varOmega}$.
Мы получили морфизм $\delta$, замыкающий диаграмму
$$
\xymatrix 
{
 & X\ar[ld]_{\env_{\varPhi}^{\varOmega} Y\circ\pi}\ar[rd]^{\env_{\varPhi\circ\pi}^{\varOmega} X} & \\
 \Env_{\varPhi}^{\varOmega} Y &   & \Env_{\varPhi\circ\pi}^{\varOmega} X\ar@{-->}[ll]_{\delta}
}
$$

Остается убедиться, что $\upsilon$ и $\delta$ -- взаимно обратные морфизмы. Во-первых,
$$
\delta\circ\upsilon\circ\underbrace{\env_{\varPhi}^{\varOmega} Y\circ\pi}_{\scriptsize\begin{matrix}\text{\rotatebox{90}{$\owns$}}\\ \Epi\end{matrix}}=\delta\circ\env_{\varPhi\circ\pi}^{\varOmega} X=\env_{\varPhi}^{\varOmega} Y\circ\pi=1_{\Env_{\varPhi}^{\varOmega} Y}\circ\underbrace{\env_{\varPhi}^{\varOmega} Y\circ\pi}_{\scriptsize\begin{matrix}\text{\rotatebox{90}{$\owns$}}\\ \Epi\end{matrix}}
\quad\Longrightarrow\quad
\delta\circ\upsilon=1_{\Env_{\varPhi}^{\varOmega} Y}.
$$
И, во-вторых,
$$
\upsilon\circ\delta\circ\underbrace{\env_{\varPhi\circ\pi}^{\varOmega} X}_{\scriptsize\begin{matrix}\text{\rotatebox{90}{$\owns$}}\\ \Epi\end{matrix}}=
\upsilon\circ\env_{\varPhi}^{\varOmega} Y\circ\pi=\env_{\varPhi\circ\pi}^{\varOmega} X=1_{\Env_{\varPhi\circ\pi}^{\varOmega} X}\circ
\underbrace{\env_{\varPhi\circ\pi}^{\varOmega} X}_{\scriptsize\begin{matrix}\text{\rotatebox{90}{$\owns$}}\\ \Epi\end{matrix}}
\quad\Longrightarrow\quad
\upsilon\circ\delta=1_{\Env_{\varPhi\circ\pi}^{\varOmega} X}.
$$
\epr

\bpr[Доказательство теоремы \ref{TH:funktorialnost-v-K^Epi}]
Утверждение (a) следует из свойства $5^\circ$ на с.\pageref{5^0:obolochka-otn-klassa-morphizmov}. Докажем (b). По лемме \ref{LM:Env_varPhi-circ-e^Omega-X=Env_varPhi^Omega-Y}, $\Env_{\varPhi}^{\varOmega} Y=\Env_{\varPhi\circ\pi}^{\varOmega} X$, а по свойству
$3^\circ$ на с.\pageref{LM:suzhenie-klassa-morfizmov},
при переходе к более узкому классу морфизмов $\varPhi\circ\pi\subseteq\varPhi$
возникает пунктирная стрелка в верхнем треугольнике диаграммы
$$
\xymatrix @R=2.pc @C=5.0pc 
{
X\ar[dr]_{\env_{\varPhi}^{\varOmega} Y\circ\pi}\ar[dd]_{\pi}\ar[r]^{\env_{\varPhi}^{\varOmega} X} &  \Env_{\varPhi}^{\varOmega} X\ar@{-->}[d]^{\Env_{\varPhi}^{\varOmega}\pi} \\
     &             \Env_{\varPhi\circ\pi}^{\varOmega} X \ar@{=}[d] \\
Y\ar[r]^{\env_{\varPhi}^{\varOmega} Y} & \Env_{\varPhi}^{\varOmega} Y \\
}
$$
Она и будет нужной стрелкой в \eqref{DIAGR:funktorialnost-v-K^Epi}, только нужно заметить, что она
является эпиморфизмом (чтобы быть морфизмом в ${\tt K}^{\Epi}$). Это следует из свойства \cite[p.80, $3^\circ$]{Akbarov-De-Gruyter-I}: поскольку $\Env_{\varPhi}^{\varOmega}\pi\circ\env_{\varPhi}^{\varOmega}
 X=\env_{\varPhi}^{\varOmega} Y\circ\pi\in\Epi$, получаем $\Env_{\varPhi}^{\varOmega}\pi\in\Epi$.

После того, как доказаны (a) и (b), утверждение (c) становится их следствием с применением теоремы
 \ref{TH:o-lok-malosti-v-podobjektah}: поскольку $\tt K$ локально мала в фактор-объектах класса $\varOmega$,
 можно выбрать отображение $X\mapsto S_X$, которое каждому объекту ставит в соответствие скелет
 $S_X$ в категории $\varOmega\cap\Epi^X$. После этого становится возможным выбрать отображение
 $X\mapsto\env_\varPhi^\varOmega X$, а для любого эпиморфизма $\pi:X\to Y$ стрелка
 $\Env_\varPhi^\varOmega\pi$ после этого появляется автоматически из диаграммы \eqref{DIAGR:funktorialnost-v-K^Epi}.
\epr

\paragraph{Случай $\Env_{\tt L}^{\tt L}$.}

Теорема \ref{TH:sushestvovanie-seti-pri-faktorizatsii} имеет важные следствия применительно к ситуации, когда  классы пробных и реализующих морфизмов совпадают, $\varPhi=\varOmega$, и представляют собой класс всех морфизмов со
значениями в некотором фиксированном классе объектов $\tt L$ (это частный
случай ситуации, описанной на с.\pageref{DEF:obolochka-otn-klassa}, только
вдобавок мы требуем, чтобы ${\tt L}={\tt M}$).

\btm\label{TH:env^L-funktor} Пусть категория ${\tt K}$ и класс ${\tt L}$  объектов в ней
обладают следующими свойствами:
 \bit{
\item[(i)] ${\tt K}$ проективно полна,

\item[(ii)] ${\tt K}$ обладает узловым разложением,

\item[(iii)] ${\tt K}$ локально мала в фактор-объектах класса $\Epi$,

\item[(iv)] класс $\Mor({\tt K},{\tt L})$ выходит из $\tt K$:
$\forall X\in\Ob({\tt K})\quad \exists \ph\in\Mor({\tt K})\quad \Dom\ph=X\quad\&\quad\Ran\ph\in{\tt L}$,

\item[(v)] класс ${\tt L}$ различает морфизмы снаружи,

\item[(vi)] класс ${\tt L}$ замкнут относительно перехода к проективным пределам,

\item[(vii)] ${\tt K}$ --- категория с узловым разложением, а класс ${\tt L}$ замкнут в ней относительно перехода от области значений морфизма к его узловому образу: если $\Ran\alpha\in{\tt L}$, то $\Im_\infty\alpha\in{\tt L}$.
 }\eit
Тогда
 \bit{

\item[(a)] всякий объект $X$ обладает оболочкой $\env_{\tt L}^{\tt L}X$ в классе объектов ${\tt L}$ (относительно того же класса объектов ${\tt L}$)

\item[(b)] всякая оболочка $\env_{\tt L}^{\tt L}X$ является биморфизмом,

\item[(c)] оболочку $\Env_{\tt L}^{\tt L}$ можно определить как функтор.
}\eit \etm
 \bpr
Условия (i)-(v) означают, что классы $\Epi$ и $\varPhi=\Mor({\tt K},{\tt L})$ удовлетворяют посылкам теоремы
\ref{TH:sushestvovanie-seti-pri-faktorizatsii}, то есть определяют на $\tt K$ полурегулярную оболочку
$\Env_\varPhi^{\Epi}=\Env_{\tt L}^{\Epi}$.
В доказательстве теоремы \ref{TH:sushestvovanie-seti-pri-faktorizatsii} эта оболочка строится переходом
от пространств $\Ran\ph\in{\tt L}$ ($\ph\in\varPhi$) сначала к проективным пределам, которые в силу (vi) будут
лежать в $\tt L$, а потом к узловому образу, который будет лежать в $\tt L$ в силу (vii).
Значит, $\Env_{\tt L}^{\Epi}\in{\tt L}$, и поэтому по свойству $1^\circ$ на с.\pageref{LM:suzhenie-verh-klassa-morfizmov}
$\Env_{\tt L}^{\Epi}=\Env_{\tt L}^{\Epi({\tt K},{\tt L})}$. По построению, класс $\varPhi$ является правым идеалом,
а по условию (v), $\varPhi$ различает морфизмы снаружи. Поэтому по теореме \ref{TH:Phi-razdel-moprfizmy-*},
$\Env_{\tt L}^{\Epi({\tt K},{\tt L})}=\Env_{\tt L}^{\Bim({\tt K},{\tt L})}$. Далее по той же теореме
\ref{TH:Phi-razdel-moprfizmy-*} оболочка в классе $\Bim({\tt K},{\tt L})=\Mor({\tt K},{\tt L})\cap\Bim$
существует тогда и только тогда, когда существует оболочка в классе $\Mor({\tt K},{\tt L})$, и эти оболочки совпадают:
$\Env_{\tt L}^{\Bim({\tt K},{\tt L})}=\Env_{\tt L}^{\Mor({\tt K},{\tt L})}$. Мы получили цепочку
$$
\Env_{\tt L}^{\Epi}=\Env_{\tt L}^{\Epi({\tt K},{\tt L})}=
\Env_{\tt L}^{\Bim({\tt K},{\tt L})}=\Env_{\tt L}^{\Mor({\tt K},{\tt L})}=
\Env^{\tt L}_{\tt L}
$$
Это доказывает (a) и (c) и попутно оказывается доказанным (b).
 \epr

\subsection{Сети мономорфизмов и полурегулярная детализация}

Пусть $\Mono_X$ --- класс мономорфизмов в категории $\tt K$, заканчивающихся в объекте $X$. Мы считаем его наделенным предпорядком $\to$, определяемым правилом
 \beq\label{DEF:to-in-Mono(X)}
\rho\to\sigma\Longleftrightarrow\quad \exists \varkappa\in\Mor({\tt K})\quad
\rho=\sigma\circ\varkappa.
 \eeq
Здесь морфизм $\varkappa$, если он существует, должен быть единственным, и, кроме того, он будет мономорфизмом (это
следует из того, что $\rho$ -- мономорфизм). Поэтому определена операция,
которая каждой паре морфизмов $\rho,\sigma\in\varGamma_X$, удовлетворяющей условию
$\rho\to\sigma$, ставит в соответствие морфизм
$\varkappa=\varkappa^\sigma_\rho$ в \eqref{DEF:to-in-Mono(X)}:
 \beq\label{DEF:to-in-Mono(X)-*}
\rho=\sigma\circ\varkappa^\sigma_\rho.
 \eeq
При этом, если $\rho\to\sigma\to\tau$, то из цепочки
$$
\tau\circ\varkappa^\tau_\rho=\rho=\sigma\circ\varkappa^\sigma_\rho=
\tau\circ\varkappa^\tau_\sigma\circ\varkappa^\sigma_\rho,
$$
в силу мономорфности $\tau$, следует равенство
\beq\label{varkappa^gamma_alpha=varkappa^gamma_beta-circ-varkappa^beta_alpha}
\varkappa^\tau_\rho=\varkappa^\tau_\sigma\circ\varkappa^\sigma_\rho.
 \eeq

\paragraph{Сети мономорфизмов.}

\bit{ \item[$\bullet$] Пусть каждому объекту $X\in\Ob({\tt K})$ категории ${\tt
K}$ поставлено в соответствие некое подмножество ${\mathcal N}_X$ в классе
$\Mono_X$ всех мономорфизмов категории ${\tt K}$, приходящих в $X$, причем
выполняются следующие три условия:
    \bit{
\item[(a)]\label{AX:set-Mono-a} для всякого объекта $X$ множество ${\mathcal
N}_X$ непусто и направлено вправо относительно предпорядка
\eqref{DEF:to-in-Mono(X)}, наследуемого из $\Mono_X$:
$$
\forall \rho,\rho'\in {\mathcal N}_X\quad \exists\sigma\in{\mathcal N}_X\quad
\rho\to\sigma\ \& \ \rho'\to\sigma,
$$

\item[(b)]\label{AX:set-Mono-b} для всякого объекта $X$ порождаемая множеством
${\mathcal N}_X$ ковариантная система морфизмов

 \beq\label{DEF:kategoriya-svyazyv-morpfizmov-*}
\Bind({\mathcal N}_X):=\{\varkappa^\sigma_\rho;\ \rho,\sigma\in{\mathcal N}_X,\
\rho\to\sigma\}
 \eeq
(морфизмы $\varkappa^\sigma_\rho$ были определены в
\eqref{DEF:to-in-Mono(X)-*}; согласно
\eqref{varkappa^gamma_alpha=varkappa^gamma_beta-circ-varkappa^beta_alpha}, эта
система является ковариантным функтором из множества ${\mathcal N}_X$,
рассматриваемого как полная подкатегория в $\Mono_X$, в $\tt K$) обладает
инъективным пределом в $\tt K$;

\item[(c)]\label{AX:set-Mono-c}  для всякого морфизма $\alpha:X\to Y$ и любого
элемента $\sigma\in{\mathcal N}_X$ найдется элемент $\tau\in{\mathcal N}_Y$ и
морфизм $\alpha_\sigma^\tau:\Dom\sigma\to\Dom\tau$ такие, что коммутативна
диаграмма
 \beq\label{DIAGR:set-mono} \xymatrix @R=2.5pc @C=4.0pc {
 X\ar[r]^{\alpha} & Y \\
 \Dom\sigma\ar[u]^{\sigma}\ar@{-->}[r]_{\alpha_\sigma^\tau} & \Dom\tau\ar@{-->}[u]_{\tau}
 } \eeq
(при фиксированных $\alpha$, $\sigma$ и $\tau$ морфизм $\alpha_\sigma^\tau$,
если он существует, определен однозначно, в силу мономорфности $\tau$).

 }\eit
Тогда
 \bit{
\item[---] семейство множеств ${\mathcal N}=\{{\mathcal N}_X;\ X\in\Ob({\tt
K})\}$ мы будем называть {\it сетью мономорфизмов}\label{DEF:set-monomorf} в
категории ${\tt K}$, а элементы множеств ${\mathcal N}_X$ -- {\it элементами
сети} ${\mathcal N}$,

\item[---] для всякого объекта $X$ систему морфизмов $\Bind({\mathcal N}_X)$,
определенную равенствами \eqref{DEF:kategoriya-svyazyv-morpfizmov-*}, мы
называем {\it системой связывающих морфизмов сети ${\mathcal N}$ над вершиной
$X$}, ее инъективный предел (существование которого гарантируется условием (b))
представляет собой инъективный конус, вершину которого мы будем обозначать
символом $X_{\mathcal N}$, а приходящие в нее морфизмы -- символами
$\rho_{\mathcal N}=\injlim_{\sigma\in{\mathcal N}_X} \varkappa^\sigma_\rho:
X_{\mathcal N}\gets\Ran\sigma$:
 \beq\label{X_F-injektiv-sistema} \xymatrix @R=2.5pc @C=2.0pc {
 & X_{\mathcal N} &  \\
 \Dom\rho\ar[ur]^{\rho_{\mathcal N}}\ar[rr]^{\varkappa^\sigma_\rho} & &  \Dom\sigma\ar[ul]_{\sigma_{\mathcal N}}
}\qquad (\rho\to\sigma);
 \eeq
при этом, в силу \eqref{DEF:to-in-Mono(X)-*}, сама система мономорфизмов
${\mathcal N}_X$ также является инъективным конусом системы $\Bind({\mathcal
N}_X)$:
 \beq\label{F_X-injektiv-sistema} \xymatrix @R=2.5pc @C=2.0pc {
 & X &  \\
 \Dom\rho\ar[rr]^{\varkappa^\sigma_\rho}\ar[ur]^{\rho} & &  \Dom\sigma\ar[ul]_{\sigma}
}\qquad (\rho\to\sigma),
 \eeq
поэтому должен существовать естественный морфизм в $X$ из вершины $X_{\mathcal
N}$ инъективного предела системы $\Bind({\mathcal N}_X)$. Этот морфизм мы будем
обозначать $\injlim{\mathcal N}_X$ и называть {\it локальным пределом сети
мономорфизмов ${\mathcal N}$ на объекте $X$}:
 \beq\label{DIAGR:sigma-sigma_F-mono}
 \xymatrix @R=2.5pc @C=2.0pc {
   X_{\mathcal N}\ar@{-->}[rr]^{\injlim{\mathcal N}_X} & & X  \\
   & \Dom\sigma\ar[ur]_{\sigma}\ar[ul]^{\sigma_{\mathcal N}}  &
}\qquad (\sigma\in{\mathcal N}_X). \eeq

\item[---] элемент сети $\tau$ в диаграмме \eqref{DIAGR:set-mono} мы называем
{\it навесом} для элемента сети $\sigma$.
 }\eit
 }\eit

Следующие утверждения двойственны теоремам \ref{TH:funktorialnost-obolochki_F},
 \ref{TH:funktorialnost-obolochki} и \ref{TH:funktorialnost-pri-seti-Epi-i-dolonyaemosti}.

\btm\label{TH:funktorialnost-otpechatka_F}\footnote{В работе \cite{Ak16} этот результат (Theorem 3.39) приведен с ошибкой: там опущено условие \eqref{lim_N-subseteq-Mono}.} Пусть ${\mathcal N}$ -- сеть
мономорфизмов в категории ${\tt K}$, причем все локальные пределы $\injlim{\mathcal N}^X$ являются мономорфизмами
\beq\label{lim_N-subseteq-Mono}
\{\projlim{\mathcal N}^X; \ X\in\Ob(\tt K)\}\subseteq\Mono(\tt
K).
\eeq
Тогда
 \bit{
\item[(a)] для любого объекта $X$ в ${\tt K}$ локальный предел
$\injlim{\mathcal N}_X:X_{\mathcal N}\to X$ является детализацией $\rf_{\mathcal N} X$ в категории ${\tt
K}$ посредством класса морфизмов ${\mathcal N}$:
 \beq\label{lim F_X=imp_F-X}
\injlim {\mathcal N}_X=\rf_{\mathcal N} X,
 \eeq

\item[(b)] для всякого морфизма $\alpha:X\to Y$ в
${\tt K}$ при любом выборе локальных пределов $\injlim {\mathcal N}_X$ и $\injlim {\mathcal N}_Y$ формула
 \beq\label{DEF:alpha_F-mono}
\alpha_{\mathcal N}=\injlim_{\sigma\in{\mathcal
N}_X}\injlim_{\tau\in{\mathcal N}_Y}\tau_{\mathcal N}\circ\alpha_\sigma^\tau
 \eeq
определяет морфизм $\alpha_{\mathcal N}:X_{\mathcal N}\to Y_{\mathcal N}$, замыкающий диаграмму
 \beq\label{DIAGR:funktorialnost-lim-F-mono}
\xymatrix @R=2.pc @C=10.0pc 
{
X\ar[d]^{\alpha} &  X_{\mathcal N}=\Rf_{\mathcal N} X \ar@{-->}[d]^{\alpha_{\mathcal N}}\ar[l]_{\injlim {\mathcal N}_X=\rf_{\mathcal N} X} \\
Y &  Y_{\mathcal N}=\Rf_{\mathcal N} Y\ar[l]_{\injlim {\mathcal N}_Y=\rf_{\mathcal N} Y} \\
}
 \eeq

\item[(c)] детализацию $\Rf_{\mathcal N}$ можно определить как функтор.
 }\eit
\etm

\btm\label{TH:funktorialnost-otpechatka}
Пусть ${\mathcal N}$ -- сеть мономорфизмов в категории ${\tt K}$, порождающая класс морфизмов $\varPhi$ снаружи:
    $$
    {\mathcal N}\subseteq\varPhi\subseteq{\mathcal N}\circ\Mor({\tt K}).
    $$
Тогда для любого класса мономорфизмов $\varGamma$ в ${\tt K}$, содержащего все локальные пределы
$\injlim{\mathcal N}_X$,
$$
\{\injlim{\mathcal N}_X; \ X\in\Ob(\tt
K)\}\subseteq\varGamma\subseteq\Mono(\tt K),
$$
выполняются следующие условия:
\bit{

\item[(a)] для любого объекта $X$ в ${\tt K}$ локальный предел
$\injlim{\mathcal N}_X$ является детализацией $\rf_\varPhi^\varGamma X$ в классе $\varGamma$ посредством класса $\varPhi$:

 \beq\label{lim F_X=imp_M^varOmega-X}
\injlim {\mathcal N}_X=\rf_\varPhi^\varGamma X, \eeq

\item[(b)] детализацию $\Rf_\varPhi^\varGamma$ можно определить как функтор.
 }\eit
 \etm

\btm\label{TH:funktorialnost-pri-seti-Mono-i-dolonyaemosti}
Пусть ${\mathcal N}$ -- сеть мономорфизмов в категории ${\tt K}$, порождающая класс морфизмов $\varPhi$ снаружи:
    $$
    {\mathcal N}\subseteq\varPhi\subseteq{\mathcal N}\circ\Mor({\tt K}).
    $$
Тогда для любого эпиморфно дополняемого\footnote{В смысле определения на
с.\pageref{DEF:klass-epimorfno-dopolnyaem}.} в $\tt K$ класса мономорфизмов $\varGamma$,
$$
\varGamma\circledcirc\varGamma^\downarrow={\tt K},
$$
выполняются следующие условия:
 \bit{
\item[(a)] в ${\tt K}$ существует сеть мономорфизмов ${\mathcal N}$ такая, что
для любого объекта $X$ в ${\tt K}$ морфизм $\mu_{\injlim{\mathcal N}_X}$ в
факторизации \eqref{faktorizatsiya-v-kat-s-faktoriz}
является детализацией $\rf_\varPhi^{\varGamma} X$ в классе $\varGamma$ посредством класса морфизмов $\varPhi$:
\beq\label{mu_lim N_X=imp_Phi^Gamma-X}
\mu_{\injlim{\mathcal N}_X}=\rf_\varPhi^{\varGamma} X,
\eeq

\item[(b)] для всякого морфизма $\alpha:X\to Y$ в ${\tt K}$ при любом выборе детализаций $\rf_\varPhi^{\varGamma}X$ и $\rf_\varPhi^{\varGamma}Y$ существует единственный морфизм $\Rf_\varPhi^{\varGamma} \alpha:\Rf_\varPhi^{\varGamma}  X\to \Rf_\varPhi^{\varGamma}  Y$ в ${\tt K}$, замыкающий диаграмму
\beq\label{DIAGR:funktorialnost-imp_varPhi^Mono-v-kat-s-uzl-razl-1}
\xymatrix @R=2.pc @C=5.0pc 
{
X\ar[d]^{\alpha} & \Rf_\varPhi^{\varGamma}  X\ar[l]_{\rf_\varPhi^{\varGamma}X} \ar@{-->}[d]^{\Rf_\varPhi^{\varGamma} \alpha} \\
Y & \Rf_\varPhi^{\varGamma}  Y\ar[l]_{\rf_\varPhi^{\varGamma}Y} \\
}
\eeq

\item[(c)] если категория ${\tt K}$ локально мала в подобъектах класса $\varGamma$,
то детализацию $\Rf_\varPhi^{\varGamma}$ можно определить как функтор.
}\eit
\etm

\paragraph{Существование сети мономорфизмов и полурегулярная детализация.}

Двойственное утверждение для детализаций выглядит так:

\btm\label{TH:sushestvovanie-seti-Mono-pri-faktorizatsii}
Пусть категория $\tt K$ и классы морфизмов $\varGamma$ и $\varPhi$ в ней удовлетворяют следующим условиям:
\bit{

\item[RR.1:] категория $\tt K$ инъективно полна,

\item[RR.2:] класс $\varGamma$ эпиморфно дополняем в $\tt K$:
$\varGamma\circledcirc\varGamma^\downarrow={\tt K}$,

\item[RR.3:] ${\tt K}$ локально мала в подобъектах класса $\varGamma$, и

\item[RR.4:] класс $\varPhi$ приходит в\footnote{В смысле определения на с.\pageref{DEF:goes-to}.} $\tt K$ и является левым идеалом в ${\tt K}$:
$$
\Ob(\tt K)=\{\Ran\ph;\ \ph\in\varPhi\},\qquad \Mor({\tt K})\circ\varPhi\subseteq\varPhi.
$$
}\eit
Тогда
 \bit{
\item[(a)] в ${\tt K}$ существует сеть мономорфизмов ${\mathcal N}$ такая, что
для любого объекта $X$ в ${\tt K}$ морфизм $\mu_{\injlim{\mathcal N}_X}$ в
факторизации \eqref{faktorizatsiya-v-kat-s-faktoriz}
является детализацией $\rf_\varPhi^{\varGamma} X$ в классе $\varGamma$ посредством класса морфизмов $\varPhi$:
\beq\label{mu_lim N_X=imp_Phi^Gamma-X-*}
\mu_{\injlim{\mathcal N}_X}=\rf_\varPhi^{\varGamma} X,
\eeq

\item[(b)] для всякого морфизма $\alpha:X\to Y$ в ${\tt K}$ при любом выборе детализаций $\rf_\varPhi^{\varGamma}X$ и $\rf_\varPhi^{\varGamma}Y$ существует единственный морфизм $\Rf_\varPhi^{\varGamma} \alpha:\Rf_\varPhi^{\varGamma}  X\to \Rf_\varPhi^{\varGamma}  Y$ в ${\tt K}$, замыкающий диаграмму
\beq\label{DIAGR:funktorialnost-imp_varPhi^Mono-v-kat-s-uzl-razl-1-*}
\xymatrix @R=2.pc @C=5.0pc 
{
X\ar[d]^{\alpha} & \Rf_\varPhi^{\varGamma}  X\ar[l]_{\rf_\varPhi^{\varGamma}X} \ar@{-->}[d]^{\Rf_\varPhi^{\varGamma} \alpha} \\
Y & \Rf_\varPhi^{\varGamma}  Y\ar[l]_{\rf_\varPhi^{\varGamma}Y} \\
}
\eeq

\item[(c)] детализацию $\Rf_\varPhi^{\varGamma}$ можно определить как функтор.
}\eit \etm

\bit{

\item[$\bullet$] Если выполняются условия RR.1-RR.4 этой теоремы, то мы будем говорить,
что {\it классы морфизмов $\varGamma$ и $\varPhi$ определяют в категории $\tt K$
полурегулярную детализацию $\Rf_\varPhi^\varGamma$}\label{DEF:polureg-otpechatok}, или что
{\it детализация $\Rf_\varPhi^\varGamma$ полурегулярна}.
}\eit

\paragraph{Насыщенные объекты.}

\bprop\label{PROP:harakterizatsiya-nasyshennosti}
Пусть $\varGamma\subseteq\Mono$, тогда для объекта $A\in\Ob({\tt K})$ следующие условия эквивалентны:
 \bit{
\item[(i)] всякое обогащение $\sigma:A\gets A'$ в классе $\varGamma$ посредством класса $\varPhi$ является изоморфизмом;

\item[(ii)] локальная единица $1_A:A\to A$ является детализацией $A$ в классе $\varGamma$ посредством класса $\varPhi$;

\item[(iii)] существует детализация $A$ в классе $\varGamma$ посредством
класса $\varPhi$, являющаяся изоморфизмом: $\rf^\varGamma_\varPhi A\in\Iso$.
 }\eit\noindent
\eprop

 \bit{
\item[$\bullet$]\label{DEF:nasyshennye-objekty} Объект $A$ категории $\tt K$ мы будем называть {\it насыщенным} в классе
морфизмов $\varGamma\subseteq\Mono$ посредством класса морфизмов $\varPhi$, или {\it насыщенным относительно детализации $\Rf_\varPhi^\varGamma$}, если он удовлетворяет условиям (i)-(iii) этого предложения.
 }\eit

Двойственным к предложению \ref{PROP:rasshirenie-v-polnyj-objekt} является

\bprop\label{PROP:obogashenie-iz-nasyshennogo-objekta}
Пусть классы морфизмов $\varGamma$ и $\varPhi$ определяют в категории $\tt K$
полурегулярную детализацию. Тогда если в обогащении $\sigma:A\gets S$ объекта $A$ (в классе $\varGamma$ посредством класса $\varPhi$) объект $S$ насыщенный (в классе $\varGamma$ посредством класса $\varPhi$), то обогащение $\sigma:A\gets S$ является детализацией  объекта $A$ (в классе $\varGamma$ посредством класса $\varPhi$).
\eprop

\btm\label{TH:kriterij-detalizatsii-v-term-nasysh-objektov}
Пусть $\varPhi$ --- класс морфизмов в категории $\tt K$, являющийся левым идеалом
\beq\label{kriterij-detalizatsii-v-term-nasysh-objektov-1}
\Mor({\tt K})\circ\varPhi\subseteq\varPhi,
\eeq
и обладающий сетью мономорфизмов ${\mathcal N}$, порождающей $\varPhi$ снаружи:
\beq\label{kriterij-detalizatsii-v-term-nasysh-objektov-2}
    {\mathcal N}\subseteq\varPhi\subseteq{\mathcal N}\circ\Mor({\tt K}).
\eeq
И пусть $\varGamma$ --- эпиморфно дополняемый\footnote{В смысле определения на с.\pageref{DEF:klass-epimorfno-dopolnyaem}.} в $\tt K$ класс мономорфизмов:
\beq\label{kriterij-detalizatsii-v-term-nasysh-objektov-3}
\varGamma\circledcirc\varGamma^\downarrow={\tt K},
\eeq
Тогда для всякого морфизма\footnote{Здесь не предполагается, что $\eta$ лежит в $\varGamma$.} $\eta:S\to A$ из произвольного насыщенного объекта $S$ (в классе $\varGamma$ посредством класса $\varPhi$) следующие условия эквивалентны:
\bit{

\item[(i)] морфизм $\eta:S\to A$ является детализацией (в классе $\varGamma$ посредством класса $\varPhi$);

\item[(ii)] для любого морфизма\footnote{Здесь не предполагается, что $\theta$ лежит в $\varGamma$.} $\theta:E\to A$ из произвольного насыщенного объекта $E$ (в классе $\varGamma$ посредством класса $\varPhi$) найдется единственный морфизм $\upsilon:E\to S$, замыкающий диаграмму
 \beq\label{kriterij-detalizatsii-v-term-nasysh-objektov}
\xymatrix @R=2.pc @C=2.0pc 
{
& A   &\\
S\ar[ur]^{\eta}  & & E \ar[ul]_{\theta}\ar@{-->}[ll]_{\upsilon}
}
 \eeq
}\eit
\etm

Обозначим через $\tt L$ класс всех насыщенных объектов в $\tt K$ (в классе морфизмов $\varGamma\subseteq\Mono$ посредством класса морфизмов $\varPhi$). Будем рассматривать $\tt L$ как полную подкатегорию в $\tt K$.

\bprop\label{PROP:otpechatok-na-L}
В условиях теоремы \ref{TH:sushestvovanie-seti-Mono-pri-faktorizatsii} строящийся в ней
функтор детализации $(X,\alpha)\mapsto (I(X),I(\alpha))$ на подкатегории насыщенных объектов ${\tt L}\subseteq {\tt K}$ изоморфен тождественному функтору:
\beq\label{I:K->L=funktor}
\forall A\in {\tt L}\qquad I(A)\cong A, \qquad \forall \alpha:\underset{\scriptsize\begin{matrix}\text{\rotatebox{90}{$\owns$}}\\ {\tt L}\end{matrix}}{A}\gets \underset{\scriptsize\begin{matrix}\text{\rotatebox{90}{$\owns$}}\\ {\tt L}\end{matrix}}{A'}\qquad
E(\alpha)=i_A^{-1}\circ\alpha\circ i_{A'}.
\eeq
\eprop

\paragraph{Подтягивание и регулярная детализация.}

\bit{
\item[RR.5:] Условимся говорить, что {\it класс морфизмов $\varGamma$ подтягивает класс морфизмов $\varPhi$}, если
\beq\label{Omega-podtyagivaet-Phi}
\forall\psi\in\Mor({\tt K})\qquad \forall \sigma\in\varGamma\qquad \big(\sigma\circ\psi\in\varPhi\quad\Longrightarrow\quad \psi\in\varPhi\big).
\eeq
}\eit

\brem Очевидно, для выполнения условия \eqref{Omega-podtyagivaet-Phi} достаточно
(но необязательно), чтобы $\varPhi$ представлял собой класс морфизмов с областями
определения из некоторого класса объектов $\tt M$ категории $\tt K$:
$$
\varPhi=\{\ph\in\Mor(\tt K):\ \Dom\ph\in{\tt M}\},
$$
\erem

\blm\label{LM:kompozitsiya-oblastej-vliyaniya} Если класс $\varGamma$ подтягивает
класс $\varPhi$, то композиция $\sigma\circ\rho:X\gets X''$ любых двух последовательно
взятых обогащений $\sigma:X\gets X'$ и $\rho:X'\gets X''$ (в классе $\varGamma$ посредством класса $\varPhi$) является обогащением (в классе $\varGamma$ посредством класса $\varPhi$).
\elm
\bpr
Это видно из диаграммы
$$
\xymatrix @R=2.pc @C=2.0pc 
{
X & X'\ar[l]_{\sigma}  & X''\ar[l]_{\rho}  \\
 & B \ar@/^2ex/[ul]^{\ph}\ar@{-->}[u]_{\ph'} \ar@/_2ex/@{-->}[ur]_{\ph''}&
}
$$
\epr

\bprop\label{PROP:nasyshennost-pri-podtyagivanii}
Если в категории ${\tt K}$ класс морфимзмов $\varGamma$ подтягивает класс морфизмов $\varPhi$, то объект $A\in\Ob({\tt K})$ тогда и только тогда будет насыщенным в классе $\varGamma$ посредством класса $\varPhi$, когда он изоморфен детализации некоторого объекта $X\in\Ob({\tt K})$: 
\beq\label{A=Rf^varGamma_varPhi-X}
A\cong\Rf^\varGamma_\varPhi X
\eeq
 \eprop

 \bit{
\item[$\bullet$] Условимся говорить, что {\it классы морфизмов $\varGamma$ и $\varPhi$ определяют в категории $\tt K$
регулярную детализацию $\Rf_\varPhi^\varGamma$}, или что {\it детализация $\Rf_\varPhi^\varGamma$ регулярна}, если в дополнение к условиям RR.1-RR.5 теоремы \ref{TH:sushestvovanie-seti-Mono-pri-faktorizatsii} класс $\varGamma$ подтягивает класс $\varPhi$.
}\eit

 \btm\label{TH:regulyarnayi-otpechatok}
Если классы морфизмов $\varGamma$ и $\varPhi$ определяют в категории $\tt K$
регулярную детализацию, то $\Rf_\varPhi^\varGamma$ можно определить как идемпотентный функтор.
 \etm

\btm[описание детализации через насыщенные объекты]\label{TH:otpechatki-i-nasyshennye-obyekty}
Пусть классы морфизмов $\varGamma$ и $\varPhi$ определяют в категории $\tt K$
регулярную детализацию. Тогда данный морфизм $\rho:X\gets A$ является детализацией
(в классе $\varGamma$ посредством класса $\varPhi$) если и только если выполняются следующие условия:
 \bit{

 \item[(i)] $\rho:X\gets A$ -- мономорфизм,

 \item[(ii)] $A$ -- насыщенный объект (в классе $\varGamma$ посредством класса $\varPhi$),

 \item[(iii)] для любого насыщенного объекта $B$ (в классе $\varGamma$ посредством $\varPhi$) и всякого морфизма $\xi:X\gets B$ найдется единственный морфизм $\xi':A\gets B$, замыкающий диаграмму
 \beq\label{DIAGR:otpechatki-i-nasyshennye-obyekty}
\xymatrix @R=2.pc @C=2.0pc 
{
X & & A\ar[ll]_{\rho}\\
& B\ar[ul]^{\xi} \ar@{-->}[ur]_{\xi'} &
}
 \eeq
 }\eit
\etm

\bcor\label{COR:obogashenie-i-nasyshennye-obyekty}
Пусть классы морфизмов $\varGamma$ и $\varPhi$ определяют в категории $\tt K$
регулярную детализацию. Тогда для всякого обогащения $\sigma:X\gets X'$ (в $\varGamma$ посредством $\varPhi$), любого насыщенного  объекта $B$ (в $\varGamma$ посредством $\varPhi$) и всякого морфизма $\xi:X\gets B$ найдется единственный морфизм $\xi':X'\gets B$, замыкающий диаграмму
 \beq\label{DIAGR:obogashenie-i-nasyshennye-obyekty}
\xymatrix @R=2.pc @C=2.0pc 
{
X & & X'\ar[ll]_{\sigma} \\
& B\ar[ul]^{\xi}\ar@{-->}[ur]_{\xi'} &
}
 \eeq
 \ecor

\btm\label{TH:summa-nasyshennyh-obyektov}
Пусть классы морфизмов $\varGamma$ и $\varPhi$ определяют в категории $\tt K$
регулярную детализацию. Тогда
  \bit{
\item[(i)] сумма (если она существует) $\coprod_{i\in I}X_i$ любого семейства $\{X_i;\ i\in I\}$ насыщенных объектов является насыщенным объектом,

\item[(ii)] инъективный предел (если он существует) $\injlim_{i\in I}X_i$ любой ковариантной (контравариантной) системы $\{X_i,\iota^j_i;\ i\in I\}$ насыщенных объектов является насыщенным объектом.
 }\eit
 \etm

\paragraph{Функториальность на мономорфизмах.}

Двойственный результат к теореме \ref{TH:funktorialnost-v-K^Epi} выглядят так.

Обозначим символом ${\tt K}^{\Mono}$ подкатегорию в $\tt K$ с тем же классом объектов что и в $\tt K$, но с мономорфизмами из $\tt K$ в качестве морфизмов:
$$
\Ob({\tt K}^{\Mono})=\Ob({\tt K}),\qquad \Mor({\tt K}^{\Mono})=\Mono({\tt K}).
$$

\btm\label{TH:funktorialnost-v-K^Mono}
Пусть ${\tt K}$ -- категория с копроизведениями (над произвольным индексным множеством), и в
ней заданы классы морфизмов $\varGamma$ и $\varPhi$, удовлетворяющие следующим условиям:
\bit{

\item[---] $\varGamma$ эпиморфно дополняем в ${\tt K}$,

\item[---] категория ${\tt K}$ локально мала в подобъектах класса $\varGamma$,

\item[---] класс $\varPhi$ приходит в\footnote{В смысле определения на с.\pageref{DEF:goes-to}.} ${\tt K}$,

\item[---] классы $\varGamma$ и $\varPhi$ связаны условием: $\varGamma\circ\varPhi\subseteq\varPhi$

}\eit
Тогда
\bit{
\item[(a)] всякий объект $X$ в ${\tt K}$ обладает детализацией $\Rf_{\varPhi}^{\varGamma}X$ в классе $\varGamma$ посредством класса $\varPhi$,

\item[(b)] для всякого мономорфизма $\pi:X\to Y$ существует единственный морфизм $\Rf_{\varPhi}^{\varGamma}\pi:\Rf_{\varPhi}^{\varGamma}X\to \Rf_{\varPhi}^{\varGamma}Y$, замыкающий диаграмму
\beq\label{DIAGR:funktorialnost-v-K^Mono}
\xymatrix @R=2.pc @C=5.0pc 
{
X\ar[d]^{\pi} & \Rf_{\varPhi}^{\varGamma} X\ar[l]_{\Rf_{\varPhi}^{\varGamma} X}\ar@{-->}[d]^{\Rf_{\varPhi}^{\varGamma}\pi} \\
Y & \Rf_{\varPhi}^{\varGamma} Y\ar[l]_{\rf_{\varPhi}^{\varGamma} Y} \\
}
\eeq

\item[(c)] детализацию $\Rf_{\varPhi}^\varGamma$ можно определить как функтор из ${\tt K}^{\Mono}$ в ${\tt K}^{\Mono}$.
}\eit
\etm

И в доказательстве используется

\blm\label{LM:Imp_varPhi-circ-i^Gama-X=Imp_varPhi^Gamma-Y}
Если ${\tt K}$ -- категория с копроизведениями, локально малая в подобъектах класса $\varGamma$,
и $\varGamma$ мономорфно дополняем в ${\tt K}$,
то для любого класса морфизмов $\varPhi$ и всякого мономорфизма $\pi:X\gets Y$ справедлива формула:
    \beq\label{Imp_varPhi-circ-e^Omega-X=Imp_varPhi^Omega-Y}
    \Rf_{\varPhi}^\varGamma X=\Rf_{\pi\circ\varPhi}^\varGamma Y
    \eeq
(детализация $X$ в ${\varGamma}$ посредством класса морфизмов $\varPhi$ совпадает с детализацией $Y$ в ${\varGamma}$  посредством класса морфизмов $\pi\circ\varPhi=\{\pi\circ\ph;\ \ph\in\varPhi\}$).
\elm

\paragraph{Случай $\Rf_{\tt L}^{\tt L}$.}

Двойственный результат к теореме \ref{TH:env^L-funktor} выглядит так:

\btm\label{TH:imp^L-funktor} Пусть категория ${\tt K}$ и класс ${\tt L}$  объектов в ней обладают следующими свойствами:
 \bit{
\item[(i)] ${\tt K}$ инъективно полна,

\item[(ii)] ${\tt K}$ обладает узловым разложением,

\item[(iii)] ${\tt K}$ локально мала в подобъектах класса $\Mono$,

\item[(iv)] класс $\Mor({\tt L},{\tt K})$ приходит в $\tt K$:
$\forall X\in\Ob({\tt K})\quad \exists \ph\in\Mor({\tt K})\quad \Dom\ph\in{\tt L}\quad\&\quad\Ran\ph=X$,

\item[(v)] класс ${\tt L}$ различает морфизмы изнутри,

\item[(vi)] класс ${\tt L}$ замкнут относительно перехода к инъективным пределам,

\item[(vii)] ${\tt K}$ --- категория с узловым разложением, а класс ${\tt L}$ замкнут в ней   относительно перехода от области определения морфизма к его узловому кообразу: если
$\Dom\alpha\in{\tt L}$, то $\Coim_\infty\alpha\in{\tt L}$.
 }\eit
Тогда
 \bit{

\item[(a)] всякий объект $X$ обладает детализацией $\rf_{\tt L}^{\tt L}X$ в классе объектов ${\tt L}$ (посредством того же класса объектов ${\tt L}$)

\item[(b)] всякая детализация $\rf_{\tt L}^{\tt L}X$ является биморфизмом,

\item[(c)] детализацию $\Rf_{\tt L}^{\tt L}$ можно определить как функтор.
}\eit \etm

\subsection{Ядро как детализация и коядро как оболочка}

\paragraph{Ядро и коядро в произвольных категориях.}

Понятия ядра и коядра можно определить в произвольной категории, необязательно имеющей нуль или систему нулевых морфизмов.

\bit{

\item[$\bullet$]  Условимся говорить, что {\it морфизм $\psi:K\to B$ тривиален изнутри}, если для любого объекта $X$ и любой пары параллельных морфизмов $\eta,\theta:X\rightrightarrows K$ композиции $\psi\circ\eta$ и $\psi\circ\theta$ совпадают:
    $$
\forall X\quad \forall \eta,\theta:X\rightrightarrows K \quad \psi\circ\eta=\psi\circ\theta,
    $$
$$
 \xymatrix @R=2.pc @C=4.pc
 {
 X\ar@{-->}@/_2ex/[r]_{\theta} \ar@{-->}@/^2ex/[r]^{\eta}& K \ar[r]^{\psi} & B,
 }
$$

\item[$\bullet$]  Условимся говорить, что {\it морфизм $\varkappa:K\to A$ тривиализует морфизм $\ph:A\to B$ изнутри}, если композиция $\ph\circ\varkappa$ --- тривиальный изнутри морфизм:
    $$
\forall X\quad \forall \eta,\theta:K\rightrightarrows A \quad \ph\circ\varkappa\circ\eta=\ph\circ\varkappa\circ\theta,
    $$
$$
 \xymatrix @R=2.pc @C=4.pc
 {
 X\ar@{-->}@/_2ex/[r]_{\theta} \ar@{-->}@/^2ex/[r]^{\eta}& K \ar[r]^{\varkappa}& A \ar[r]^{\ph}& B,
 }
$$

\item[$\bullet$] {\it Ядром} морфизма $\ph:A\to B$ называется морфизм $\varkappa:K\to A$ такой что
\bit{

\item[K1:] морфизм $\varkappa$ тривиализует морфизм $\ph$ изнутри, и

\item[K2:] если какой-то морфизм $\lambda:L\to A$ тривиализует морфизм $\ph$ изнутри, то найдется единственный морфизм $\iota:L\to K$, замыкающий диаграмму
 \beq\label{DIAGR:ker}
 \xymatrix @R=2.pc @C=4.pc
 {
 K\ar[r]^{\varkappa}& A \ar[r]^{\ph}& B \\
 & L\ar[u]_{\lambda}\ar@{-->}[ul]^{\iota} &
 }
 \eeq
}\eit

\item[$\bullet$] Условимся говорить, что {\it категория $\tt K$ обладает ядрами}, если в ней каждый морфизм $\alpha$ обладает ядром.

}\eit

\bprop\label{PROP:monomorfnost-ker}
Ядро $\varkappa:K\to A$ всякого морфизма $\ph:A\to B$ (если оно существует) является мономорфизмом.
\eprop
\bpr
Пусть $\eta,\theta:L\rightrightarrows K$ -- два параллельных морфизма, такие что $\varkappa\circ\eta=\varkappa\circ\theta$. Обозначим $\lambda=\varkappa\circ\eta=\varkappa\circ\theta:L\to A$. Ясно, что $\lambda$ тривиализирует $\ph$ изнутри. Поэтому по условию K2 существует единственный морфизм $\iota:L\to K$ такой что $\lambda=\varkappa\circ\iota$:
 $$
 \xymatrix @R=2.pc @C=4.pc
 {
 K\ar[r]^{\varkappa}& A \ar[r]^{\ph}& B \\
 & L\ar[u]_{\lambda}\ar@{-->}[ul]^{\iota} &
 }
 $$
Из единственности $\iota$ следует $\alpha=\iota=\beta$.
\epr

Из предложения \ref{PROP:monomorfnost-ker} и условия K2 следует, что ядро $\varkappa$ (если оно существует) единственно с точностью до изоморфизма в категории $\Mono(A)$ мономорфизмов со значениями в $A$. Поэтому ему можно приписать обозначение:
$$
\ker\ph=\varkappa,\qquad \Ker\ph=K.
$$

Следующее наблюдение было подсказано автору М.Бранденбургом.

\bprop\label{PROP:neobh-usl-Ker}
Если морфизм $\ph:A\to B$ обладает ядром, то для любых двух тривиализующих его изнутри морфизмов
$$
\lambda:L\to A, \qquad \mu:M\to A
$$
их прямая сумма
$$
\lambda\sqcup\mu:L\sqcup M\to A
$$
тоже тривиализует $\ph$ изнутри.
\eprop
\bpr
Пусть $\varkappa:K\to A$ --- ядро морфизма $\ph$. Морфизмы $\lambda:L\to A$ и $\mu:M\to A$ тривиализируют $\ph$ изнутри, поэтому они пропускаются через морфизм $\varkappa$:
 $$
 \xymatrix @R=2.pc @C=4.pc
 {
 K\ar[r]^{\varkappa}&  A\ar[r]^{\ph} & B \\
 & L\ar[u]_{\lambda}\ar@{-->}[ul]^{\iota} & 
 }
 \qquad 
 \xymatrix @R=2.pc @C=4.pc
 {
 K\ar[r]^{\varkappa}&  A\ar[r]^{\ph} & B \\
 & M\ar[u]_{\mu}\ar@{-->}[ul]^{\upsilon} &
 }
 $$
Эти диаграммы дают коммутативную диаграмму 
$$
 \xymatrix @R=2.pc @C=4.pc
 {
 K\ar[r]^{\varkappa}&  A\ar[r]^{\ph} & B \\
 & L\sqcup M\ar[u]_{\lambda\sqcup \mu}\ar@{-->}[ul]^{\iota\sqcup \upsilon} &
 }
 $$ 
и, поскольку в ней морфизм $\ph\circ\varkappa$ тривиален изнутри, мы получаем что морфизм $\ph\circ\lambda\sqcup \mu=\ph\circ\varkappa\circ\iota\sqcup\upsilon$ тоже должен быть тривиален изнутри. 
\epr

\bcor
Если категория $\sf K$ обладает ядрами, то для любых двух тривиальных изнутри морфизмов
$$
\lambda:L\to B, \qquad \mu:M\to B
$$
их прямая сумма
$$
\lambda\sqcup\mu:L\sqcup M\to B
$$
тоже тривиальна изнутри.
\ecor
\bpr
Это следует из предложения \ref{PROP:neobh-usl-Ker} если положить $\ph=1_B$.
\epr

Двойственное понятие коядра определяется следующим образом.

\bit{

\item[$\bullet$] Условимся говорить, что {\it морфизм $\psi:A\to C$ тривиален снаружи}, если для любого объекта $X$ и любой пары параллельных морфизмов $\eta,\theta:C\rightrightarrows X$ композиции $\eta\circ\psi$ и $\theta\circ\psi$ совпадают:
    $$
\forall X\quad \forall \eta,\theta:C\rightrightarrows X \quad \eta\circ\psi=\theta\circ\psi,
    $$
$$
 \xymatrix @R=2.pc @C=4.pc
 {
  A \ar[r]^{\psi}& C \ar@{-->}@/_2ex/[r]_{\theta} \ar@{-->}@/^2ex/[r]^{\eta}& X.
 }
$$

\item[$\bullet$] Условимся говорить, что {\it морфизм $\chi:B\to C$ тривиализует морфизм $\ph:A\to B$ снаружи}, если композиция $\chi\circ\ph$ --- тривиальный снаружи морфизм:
    $$
\forall X\quad \forall \eta,\theta:C\rightrightarrows X \quad \eta\circ\chi\circ\ph=\theta\circ\chi\circ\ph,
    $$
$$
 \xymatrix @R=2.pc @C=4.pc
 {
  A \ar[r]^{\ph}& B \ar[r]^{\chi}& C \ar@{-->}@/_2ex/[r]_{\theta} \ar@{-->}@/^2ex/[r]^{\eta}& X.
 }
$$

\item[$\bullet$] {\it Коядром} морфизма $\ph:A\to B$ называется морфизм $\gamma:B\to C$ такой что
\bit{

\item[C1:] морфизм $\gamma$ тривиализует морфизм $\ph$ снаружи, и

\item[C2:] если какой-то морфизм $\delta:B\to D$ тривиализует $\ph$ снаружи, то найдется единственный морфизм $\iota:C\to D$, замыкающий диаграмму
 \beq\label{DIAGR:coker}
 \xymatrix @R=2.pc @C=4.pc
 {
 A\ar[r]^{\ph}& B \ar[r]^{\gamma}\ar[d]_{\delta}& C\ar@{-->}[dl]^{\iota} \\
 & D &
 }
 \eeq

}\eit

\item[$\bullet$] Условимся говорить, что {\it категория $\tt K$ обладает коядрами}, если в ней каждый морфизм $\alpha$ обладает коядром.

}\eit

Аналогично предложению \ref{PROP:monomorfnost-ker} доказывается

\bprop\label{PROP:epimorfnost-coker}
Коядро $\gamma:B\to C$ всякого морфизма $\ph:A\to B$ (если оно существует) является эпиморфизмом.
\eprop

Из этого предложения и условия C2 следует, что коядро $\gamma$ (если оно существует) единственно с точностью до изоморфизма в категории $\Epi(B)$ эпиморфизмов, выходящих из $B$. Поэтому ему можно приписать обозначение:
$$
\coker\ph=\gamma,\qquad \Coker\ph=C.
$$

Двойственным к предложению \ref{PROP:neobh-usl-Ker} является

\bprop\label{PROP:neobh-usl-Coker}
Если морфизм $\ph:A\to B$ обладает коядром, то для любых двух тривиализующих его снаружи морфизмов
$$
\delta:B\to D, \qquad \e:B\to E
$$
их прямое произведение
$$
\delta\sqcap\e:B\to D\sqcap E
$$
тоже тривиализует $\ph$ снаружи.
\eprop

\bcor
Если категория $\sf K$ обладает коядрами, то для любых двух тривиальных снаружи морфизмов
$$
\delta:A\to D, \qquad \e:A\to E
$$
их прямое произведение
$$
\delta\sqcap\e:A\to D\sqcap E
$$
тоже тривиально снаружи.
\ecor

\paragraph{Ядро и коядро в категории с нулем.}

Напомним \cite{MacLane,Akbarov-De-Gruyter-I}, что {\it нулевым объектом} или {\it нулем} в категории ${\tt K}$ называется такой объект $0$, что для любого объекта $X$ в ${\tt K}$ существует только один морфизм $X\to 0$, и только один морфизм $0\to X$. Если в категории ${\tt K}$ есть нулевой объект, то он единственен с точностью до изоморфизма. Морфизм $X\overset{\ph}{\longrightarrow}Y$ в категории ${\tt K}$ с нулевым объектом $0$ называется {\it нулевым морфизмом}, если он пропускается через нулевой объект, то есть если $\ph$ является композицией двух (единственных) морфизмов $X\overset{0_{X,0}}{\longrightarrow}0$ и $0\overset{0_{0,Y}}{\longrightarrow}Y$:
$$
\begin{diagram}
\node{X}\arrow{se,b}{0_{X,0}}\arrow[2]{e,t}{\ph}\node[2]{Y} \\
\node[2]{0}\arrow{ne,b}{0_{0,Y}}
\end{diagram}
$$
Понятно, что такой морфизм тоже единственен, поэтому он имеет специальное обозначение
$$
\ph=0_{X,Y}.
$$

\bit{

\item[$\bullet$] Условимся говорить, что {\it морфизм $\varkappa$ обнуляет морфизм $\ph$ изнутри}, если 
\beq\label{ph-circ-varkappa-*}
 \ph\circ\varkappa=0.
\eeq

\item[$\bullet$] Условимся говорить, что {\it морфизм $\gamma$ обнуляет морфизм $\ph$ снаружи}, если 
\beq\label{gamma-circ-ph-*}
 \gamma\circ\ph=0.
\eeq

}\eit

\bprop\label{PROP:harakt-ker-0}  В категории с нулем морфизм $\varkappa:K\to A$ тогда и только тогда будет ядром морфизма $\ph:A\to B$, когда выполняются следующие условия:
\bit{
\item[(i)] морфизм $\varkappa$ обнуляет морфизм $\ph$ изнутри:
\beq\label{ph-circ-varkappa=0}
 \ph\circ\varkappa=0
\eeq

\item[(ii)] для любого морфизма $\varkappa':K'\to A$, обнуляющего морфизм $\ph$ изнутри, $\ph\circ\varkappa'=0$, существует единственный морфизм $\iota: K'\to K$ такой, что $\varkappa'=\varkappa\circ\iota$
 \beq\label{DIAGR:ker-0}
 \xymatrix @R=2.pc @C=4.pc
 {
 K\ar[r]^{\varkappa}& A \ar[r]^{\ph}& B \\
 & K'\ar[u]_{\varkappa'}\ar@{-->}[ul]^{\iota}\ar[ur]_{0_{K',B}} &
 }
 \eeq
}\eit
\eprop
\bpr
1. Пусть $\varkappa:K\to A$ -- ядро морфизма $\ph:A\to B$ (в смысле определения на с.\pageref{DIAGR:ker}). Тогда $\varkappa$ тривиализует $\ph$ изнутри, и поэтому пара морфизмов $0_K,1_K:K\rightrightarrows K$ в композиции с $\ph\circ\varkappa$ должна уравниваться. Это можно интерпретировать так:
$$
\ph\circ\varkappa=\ph\circ\varkappa\circ 1_K=\ph\circ\varkappa\circ 0_K=0_{K,A}.
$$
С другой стороны, если для какого-то морфизма $\varkappa':K'\to A$ справедливо равенство $\ph\circ\varkappa'=0$, то для любых параллельных морфизмов $\eta,\theta:I\to K'$ мы получим
$$
\ph\circ\varkappa'\circ\eta=0=\ph\circ\varkappa'\circ\theta,
$$
и это будет означать, что $\varkappa'$ тривиализирует $\ph$ изнутри, поэтому по условию K2 существует единственный морфизм $\iota:K'\to K$ такой что $\varkappa'=\varkappa\circ\iota$.

2. Наоборот, пусть $\varkappa$ удовлетворяет условиям (i) и (ii). Тогда, во-первых, для любых $\eta,\theta:I\rightrightarrows  K$ мы получаем
$$
\ph\circ\varkappa\circ\eta=0=\ph\circ\varkappa\circ\theta,
$$
то есть $\varkappa$ тривиализует $\ph$ изнутри. А, во-вторых, если какой-то морфизм $\varkappa':K'\to A$ тривиализует $\ph$ изнутри, то, в частности, композиция $\ph\circ\varkappa'$ должна уравнивать морфизмы $0_{K'},1_{K'}:K'\rightrightarrows K$, и поэтому мы получим
$$
\ph\circ\varkappa'=\ph\circ\varkappa'\circ 1_{K'}=\ph\circ\varkappa'\circ 0_{K'}=0.
$$
По свойству (ii) это означает, что существует единственный морфизм $\iota:K'\to K$ такой что $\varkappa'=\varkappa\circ\iota$.
\epr

Двойственное предложение для коядра выглядит так:

\bprop\label{PROP:harakt-coker-0} В категории с нулем морфизм $\gamma:B\to C$ тогда и только тогда будет коядром морфизма $\ph:A\to B$, когда выполняются следующие условия:
\bit{
\item[(i)] морфизм $\gamma$ обнуляет морфизм $\ph$ снаружи:
\beq\label{gamma-circ-ph=0}
 \gamma\circ\ph=0
\eeq

\item[(ii)] для любого морфизма $\gamma':B\to C'$, обнуляющего морфизм $\ph$ снаружи, $\gamma'\circ\ph=0$, существует единственный морфизм $\iota: C\to C'$ такой, что $\gamma'=\iota\circ\gamma$
 \beq\label{DIAGR:coker-0}
 \xymatrix @R=2.pc @C=4.pc
 {
 A\ar[r]^{\ph}\ar[dr]_{0_{A,C'}} & B\ar[d]_{\gamma'} \ar[r]^{\gamma}& C \ar@{-->}[dl]^{\iota} \\
 & C' &
 }
 \eeq
}\eit
\eprop

\noindent\rule{160mm}{0.1pt}\begin{multicols}{2}

\paragraph{Примеры ядер и коядер.}

\bex\label{EX:Ker-v-Top}
В категории $\sf Top$ топологических пространств морфизм $\ph:A\to B$ обладает ядром только если пространство $A$ одноточечное, и тогда ядром $\ph$ будет единица $1_A:A\to A$:
$$
\ker\ph=1_A,\qquad \Ker\ph=A.
$$
Точно так же, $\ph:A\to B$ обладает коядром только если пространство $B$ одноточечное, и тогда коядром $\ph$ будет  единица $1_B:B\to B$:
$$
\coker\ph=1_B,\qquad \Coker\ph=B.
$$

\eex

\bex\label{EX:Ker-v-Vect}
В категории $_k\Vect$ векторных пространств над полем $k$ ядром морфизма $\ph:X\to Y$ является прообраз нуля, а коядром -- фактор-пространство по множеству значений:
$$
\Ker\ph=\ph^{-1}(0),\qquad \Coker\ph=Y/\ph(X).
$$
\eex

\bex\label{EX:Ker-v-LCS}
В категории $\LCS$ локально выпуклых пространств над полем $\C$ ядром морфизма $\ph:X\to Y$ является прообраз нуля (с индуцированной их $X$  топологией), а коядром -- фактор-пространство по замыканию множества значений (с топологией фактор-прространства):
$$
\Ker\ph=\ph^{-1}(0),\qquad \Coker\ph=Y/\overline{\ph(X)}.
$$
\eex

\bex\label{EX:Ker-v-Ste}
В категории $\Ste$ стереотипных пространств над полем $\C$ ядром морфизма $\ph:X\to Y$ является прообраз нуля (с топологией непосредственного подпространства в $X$), а коядром --- непосредственное стереотипное фактор-пространство по замыканию множества значений:
$$
\Ker\ph=\ph^{-1}(0)^\vartriangle,\qquad \Coker\ph=\big(Y/\overline{\ph(X)}\big)^\triangledown.
$$
\eex

\bex\label{EX:Ker-v-Alg}
Категория $_k\Alg$ унитальных алгебр над полем $k$ не обладает ни ядрами, ни коядрами. 

1. Рассмотрим например, тождественный морфизм поля $k$ в себя:
$$
\ph:k\to k,\qquad \ph(t)=t,\qquad t\in k.
$$
Предположим, что этот морфизм имеет ядро $\varkappa:K\to k$. Поскольку по предложению \ref{PROP:monomorfnost-ker}, ядро является мономорфизмом (а мономорфизмы в категории алгебр --- инъекции), мы получаем, что $K=k$ и $\varkappa$ --- тоже тождественный морфизм:
$$
\varkappa(t)=t, \qquad t\in k.
$$
Пусть теперь $X=k^2$ с покоординатным умножением. Рассмотрим два морфизма $\eta,\theta:X=k^2\rightrightarrows k=K$ 
$$
\eta(s,t)=s,\quad \theta(s,t)=t, \qquad s,t\in k. 
$$
Мы получим:
$$
\ph\circ\varkappa\circ\eta=\eta\ne\theta=\ph\circ\varkappa\circ\theta.
$$
то есть морфизм $\varkappa$ не тривиализирует морфизм $\ph$. Это значит, что $\varkappa$ не может быть ядром $\ph$.

2. Рассмотрим тождественное отображение алгебры $k^2$ (с покоординатным умножением) в себя:
$$
1_{k^2}: k^2\to k^2.
$$
Предположим, что этот морфизм обладает коядром $\gamma:k^2\to C$. Тогда по предлодению \ref{PROP:neobh-usl-Coker}, есди морфизмы $\delta:k^2\to D$ и $\e:k^2\to E$ тривиализуют $1_{k^2}$ снаружи,то их произведение $\delta\sqcap\e:k^2\to D\sqcap  E$ тоже должно тривиализовывать $1_{k^2}$ снаружи.  В частности, морфизмы $\delta,\e:k^2\to k$
$$
\delta(s,t)=s,\quad \e(s,t)=t,\qquad s,t\in k
$$
тривиализуют $1_{k^2}$ снаружи, поэтому их произведение 
$$
\delta\sqcap\e=1_{k^2}: k^2\to k^2=k\sqcap k
$$
тоже должно тривиализовывать $1_{k^2}$ снаружи. Но это не так, потому что 
$$
\delta\circ\underbrace{(\delta\sqcap\e)}_{1_{k^2}}\circ 1_{k^2}=\delta\ne\e=\e\circ\underbrace{(\delta\sqcap\e)}_{1_{k^2}}\circ 1_{k^2}.
$$

\eex

\bex\label{EX:Ker-v-Aug}
Пусть $A$ -- унитальная алгебра над полем $k$. {\it Аугментацией} на алгебре $A$ называется произвольный гомоморфизм $\e:A\to k$ алгебр над $k$. Понятно, что задать аугментацию на $A$ --- то же, что задать двусторонний идеал $I_A$ в $A$ такой, что $A$ является прямой суммой векторных пространств над $\C$
$$
A=I_A\oplus k\cdot 1_A,
$$
где $1_A$ -- единица $A$. {\it Аугментированной алгеброй} называется произвольная пара $(A,\e)$, в которой $A$ -- унитальная алгебра над $k$, а  $\e:A\to k$ -- аугментация. Класс всех аугментированных алгебр образует категорию $_k\AugAlg$, в которой морфизмом $\ph:(A,\e_A)\to (B,\e_B)$ считается произвольный гомоморфизм $\ph:A\to B$ алгебр над $k$, сохраняющий аугментацию в следующем смысле:
$$
\e_B\circ\ph=\e_A.
$$
Можно заметить, что алгебра $k$ с тождественным отображением $\id_k:k\to k$ в качестве аугментации является нулем в категории $_k\AugAlg$. Следующий факт доказывается так же как аналогичное утверждение для стереотипных алгебр с аугментапцией (см. ниже предложение \ref{PROP:AugSteAlg-imeet-ker}): {\it категория  $_k\AugAlg$ аугментированных алгебр над $k$ обладает ядрами и коядрами.}

\eex

\end{multicols}\noindent\rule[10pt]{160mm}{0.1pt}

\paragraph{Ядро и коядро как функторы.}

Пусть $\tt K$ -- произвольная категория. Для всякого морфизма $\alpha:A\to A'$ в этой категории условимся символами $\Dom(\alpha)$ и $\Ran(\alpha)$ обозначать соответственно его начало и конец:
$$
\Dom(\alpha)=A,\qquad \Ran(\alpha)=A'.
$$
Таким образом, всякий морфизм $\alpha$ ведет из $\Dom(\alpha)$ в $\Ran(\alpha)$:
$$
\alpha: \Dom(\alpha)\to \Ran(\alpha).
$$

Образуем из категории $\tt K$ новую категорию ${\tt K}^*$ по следующим правилам:

\bit{

\item[---] объектами категории ${\tt K}^*$ являются морфизмы категории ${\tt K}$:
$$
\Ob({\tt K}^*)=\Mor({\tt K}),
$$

\item[---] морфизмом категории ${\tt K}^*$ между объектами $\alpha,\beta\in\Ob({\tt K}^*)=\Mor(\tt K)$ является произвольная пара $(\ph,\ph')$ морфизмов $\ph,\ph'\in\Ob({\tt K}^*)=\Mor(\tt K)$, замыкающая диаграмму
 \beq\label{morfism-v-K^*}
 \xymatrix @R=2.pc @C=4.pc
 {
 \Dom(\alpha)\ar[d]^{\alpha}  \ar@{-->}[r]^{\ph} & \Dom(\beta) \ar[d]^{\beta}
 \\
 \Ran(\alpha) \ar@{-->}[r]^{\ph'} & \Ran(\beta)
 }
 \eeq

\item[---] композицией морфизмов $(\ph,\ph')$ и $(\chi,\chi')$ в ${\tt K}^*$ считается морфизм $(\chi\circ\ph,\chi'\circ\ph')$ (при условии, что композиции $\chi\circ\ph$ и $\chi'\circ\ph'$ определены в ${\tt K}$); наглядно это иллюстрируется коммутативной диаграммой
 \beq\label{kompoz-v-K^*}
 \xymatrix @R=2.pc @C=4.pc 
 {
 \Dom(\alpha)\ar[d]^{\alpha}  \ar@{-->}[r]^{\ph} & \Dom(\beta) \ar[d]^{\beta}\ar@{-->}[r]^{\chi} & \Dom(\gamma) \ar[d]^{\gamma}
 \\
 \Ran(\alpha) \ar@{-->}[r]^{\ph'} & \Ran(\beta)\ar@{-->}[r]^{\chi'} & \Ran(\gamma)
 }
 \eeq

}\eit

Морфизмы категории ${\tt K}^*$, по-видимому, удобно записывать в виде специальных дробей: если морфизм $(\ph,\ph'):\alpha\to\beta$, представленный диаграммой \eqref{morfism-v-K^*}, записать символом
$$
\frac{\beta,\ph}{\ph',\alpha},
$$
то закон композиции морфизмов в категории ${\tt K}^*$ будет описываться формулой
\beq
\frac{\gamma,\chi}{\chi',\beta}\circ\frac{\beta,\ph}{\ph',\alpha}=\frac{\gamma,\chi\circ\ph}{\chi'\circ\ph',\alpha}.
\eeq

Следующее утверждение очевидно.

\btm\label{TH:Ker-frac-beta,ph-ph',alpha}
Если категория $\tt K$ обладает ядрами, то
\bit{
\item[(i)] в категории ${\tt K}^*$ всякий морфизм $\frac{\beta,\ph}{\ph',\alpha}$ порождает единственный морфизм $\Ker\frac{\beta,\ph}{\ph',\alpha}: \Ker\alpha\to\Ker\beta$, замыкающий диаграмму в $\tt K$:
 \beq\label{Ker-frac-beta,ph-ph',alpha}
 \xymatrix @R=2.pc @C=4.pc
 {
 \Ker\alpha\ar[d]_{\ker\alpha}\ar@{-->}[r]^{\Ker\frac{\beta,\ph}{\ph',\alpha}} & \Ker\beta\ar[d]_{\ker\beta} \\
 \Dom(\alpha)\ar[d]^{\alpha}  \ar[r]^{\ph} & \Dom(\beta) \ar[d]^{\beta}
 \\
 \Ran(\alpha) \ar[r]^{\ph'} & \Ran(\beta)
 }
 \eeq

 \item[(ii)] отображение
 $$
\frac{\beta,\ph}{\ph',\alpha}\mapsto \Ker\frac{\beta,\ph}{\ph',\alpha}
 $$
 является ковариантным функтором из ${\tt K}^*$ в ${\tt K}$;

 \item[(iii)] отображение
 $$
\frac{\beta,\ph}{\ph',\alpha}\mapsto \frac{\ker\beta,\Ker\frac{\beta,\ph}{\ph',\alpha}}{\ph,\ker\alpha}
 $$
 является ковариантным функтором из ${\tt K}^*$ в ${\tt K}^*$.

 }\eit
\etm

\bprop\label{TH:Ker-frac-beta,ph-ph',alpha-in-Mono}
Если в диаграмме \eqref{Ker-frac-beta,ph-ph',alpha} морфизм $\ph$ является мономорфизмом, то морфизм $\Ker\frac{\beta,\ph}{\ph',\alpha}$ также является мономорфизмом.
\eprop
\bpr
Морфизм
$$
\ker\beta\circ\Ker\frac{\beta,\ph}{\ph',\alpha}=\ph\circ\ker\alpha
$$
является мономорфизмом как композиция двух мономорфизмов, $\ker\alpha$ и $\ph$, поэтому внутренний морфизм в композиции слева от знака равенства, $\Ker\frac{\beta,\ph}{\ph',\alpha}$, тоже является мономорфизмом.
\epr

Для коядер утверждение, двойственное к теореме \ref{TH:Ker-frac-beta,ph-ph',alpha} выглядит так:

\btm\label{TH:Coker-frac-beta,ph-ph',alpha}
Если категория $\tt K$ обладает коядрами, то
\bit{
\item[(i)] в категории ${\tt K}^*$ всякий морфизм $\frac{\beta,\ph}{\ph',\alpha}$ порождает единственный морфизм $\Coker\frac{\beta,\ph}{\ph',\alpha}: \Coker\alpha\to\Coker\beta$, замыкающий диаграмму в $\tt K$:
 \beq\label{Coker-frac-beta,ph-ph',alpha}
 \xymatrix @R=2.pc @C=4.pc
 {
 \Dom(\alpha)\ar[d]_{\alpha}  \ar[r]^{\ph} & \Dom(\beta) \ar[d]^{\beta}
 \\
 \Ran(\alpha) \ar[r]^{\ph'}\ar[d]_{\coker\alpha} & \Ran(\beta) \ar[d]^{\coker\beta} \\
  \Coker\alpha\ar@{-->}[r]^{\Coker\frac{\beta,\ph}{\ph',\alpha}} & \Coker\beta
 }
 \eeq

 \item[(ii)]  отображение
 $$
\frac{\beta,\ph}{\ph',\alpha}\mapsto \Coker\frac{\beta,\ph}{\ph',\alpha}
 $$
 является ковариантным функтором из ${\tt K}^*$ в ${\tt K}$;

 \item[(iii)] отображение
 $$
\frac{\beta,\ph}{\ph',\alpha}\mapsto \frac{\coker\beta,\ph'}{\Coker\frac{\beta,\ph}{\ph',\alpha},\coker\alpha}
 $$
 является ковариантным функтором из ${\tt K}^*$ в ${\tt K}^*$.

 }\eit
\etm

А утверждение, двойственное предложению \ref{TH:Ker-frac-beta,ph-ph',alpha-in-Mono}, имеет такой вид:

\bprop\label{TH:Coker-frac-beta,ph-ph',alpha-in-Epi}
Если в диаграмме \eqref{Coker-frac-beta,ph-ph',alpha} морфизм $\ph'$ является эпиморфизмом, то
  морфизм $\Coker\frac{\beta,\ph}{\ph',\alpha}$ также является эпиморфизмом.
\eprop
\bpr
Морфизм
$$
\Coker\frac{\beta,\ph}{\ph',\alpha}\circ\coker\alpha=\coker\beta\circ\ph'
$$
является эпиморфизмом как композиция двух эпиморфизмов, $\coker\beta$ и $\ph'$, поэтому внешний морфизм в композиции слева от знака равенства, $\Coker\frac{\beta,\ph}{\ph',\alpha}$, тоже является эпиморфизмом.
\epr

\paragraph{Ядро как детализация.}

Следующее утверждение мы не доказываем, потому что его доказательство получается обращением стрелок в доказательстве существенно более важной для нас теоремы \ref{TH:Coker-rasshirenie} (которую мы доказываем ниже).

\btm[о ядре как детализации]\label{TH:Ker=Ref} Пусть $\varGamma$ и $\varPhi$ --- классы морфизмов категории $\tt K$, определяющие полурегулярную детализацию $\Rf=\Rf_\varPhi^\varGamma$, и пусть морфизм $\alpha:X\to Y$ в $\tt K$ обладает следующими свойствами:
\bit{
\item[(i)] $\alpha$ обладает ядром $\ker\alpha$ в $\tt K$;

\item[(ii)] его детализация $\Rf\alpha$ обладает ядром $\ker\Rf\alpha$ в $\tt K$;

\item[(iii)] соответствующий морфизм ядер принадлежит классу $\varGamma$:
\beq\label{Ker-in-Gamma}
\Ker\frac{\alpha,\rf X}{\rf Y,\Rf\alpha}\in\varGamma,
\eeq
}\eit\noindent
Тогда ядро детализации $\alpha$ является обогащением в классе $\varGamma$ посредством класса $\varPhi$, и, как следствие, существует единственный морфизм $\upsilon$, замыкающий диаграмму
\beq\label{Ker=Ref-1}
 \xymatrix @R=2.pc @C=6.pc
 {
\Rf\Ker\alpha\ar@{-->}[dr]^{\upsilon}\ar@/^6ex/[drr]_{\rf\Ker\alpha}
\ar@/_4ex/[ddr]_{\Rf\ker\alpha} && \\
&\Ker\Rf\alpha \ar[r]_{\Ker\frac{\alpha,\rf X}{\rf Y,\Rf\alpha}}\ar[d]^{\ker\Rf\alpha}  & \Ker\alpha \ar[d]^{\ker\alpha}  \\
& \Rf X \ar[r]_{\rf X}\ar[d]^{\Rf\alpha}  & X \ar[d]^{\alpha}  \\
& \Rf Y \ar[r]_{\rf Y} & Y
 }
\eeq
Если же дополнительно к условиям (i)--(iii) выполнено условие 
\bit{
\item[(iv)] ядро $\Ker\Rf\alpha$ является насыщенным объектом\footnote{В смысле определения на с.\pageref{DEF:nasyshennye-objekty}.} в классе морфизмов $\varGamma$ посредством класса морфизмов $\varPhi$,
}\eit\noindent
то морфизм $\upsilon$ в диаграмме \eqref{Ker=Ref-1} является изоморфизмом, и поэтому ядро $\Ker\Rf\alpha$ детализации $\Rf\alpha$ морфизма $\alpha$ является детализацией $\Rf\Ker\alpha$ объекта $\Ker\alpha$ в классе $\varGamma$ посредством класса $\varPhi$:
\beq\label{Ker=Ref}
\Ker\Rf\alpha\cong \Rf\Ker\alpha.
\eeq
\etm

\paragraph{Коядро как  оболочка.}

\btm[о коядре как оболочке]\label{TH:Coker-rasshirenie} Пусть $\varOmega$ и $\varPhi$ -- классы морфизмов в категории $\tt K$, определяющие полурегулярную оболочку $\Env=\Env_\varPhi^\varOmega$, и пусть морфизм $\alpha:X\to Y$ в $\tt K$ обладает следующими свойствами:
\bit{
\item[(i)] $\alpha$ обладает коядром $\coker\alpha$ в $\tt K$;

\item[(ii)] его оболочка $\Env\alpha$ обладает коядром $\coker\Env\alpha$ в $\tt K$;

\item[(iii)] соответствующий морфизм коядер принадлежит классу $\varOmega$:
\beq\label{Coker-in-Omega}
\Coker\frac{\Env\alpha,\env X}{\env Y,\alpha}\in\varOmega
\eeq
}\eit\noindent
Тогда коядро $\coker\Env\alpha$ оболочки $\Env\alpha$ морфизма $\alpha$ является расширением объекта $\Coker\alpha$ в классе $\varOmega$ относительно класса $\varPhi$, и, как следствие, существует единственный морфизм $\upsilon$, замыкающий диаграмму
\beq\label{Coker-rasshirenie}
 \xymatrix @R=2.pc @C=6.pc
 {
  X \ar[r]^{\env X}\ar[d]_{\alpha} &  \Env X \ar[d]^{\Env\alpha}  & \\
 Y \ar[r]^{\env Y}\ar[d]^{\coker\alpha} & \Env Y \ar[d]_{\coker\Env\alpha}\ar@/^4ex/[ddr]^{\quad \Env\coker\alpha} & \\
  \Coker\alpha\ar[r]_{\Coker\frac{\Env\alpha,\env X}{\env Y,\alpha}}\ar@/_6ex/[drr]^{\qquad \env\Coker\alpha}
   & \Coker\Env\alpha\ar@{-->}[dr]^{\upsilon} & \\
  & &  \Env\Coker\alpha
  }
\eeq
Если же дополнительно к условиям (i)--(iii) выполнено условие 
\bit{
\item[(iv)] коядро  $\Coker\Env\alpha$ является полным объектом\footnote{В смысле определения на с.\pageref{DEF:polnye-objekty}.} в классе морфизмов $\varOmega$ относительно класса морфизмов $\varPhi$,
}\eit\noindent
то морфизм $\upsilon$ в диаграмме \eqref{Coker-rasshirenie} является изоморфизмом, и поэтому коядро $\Coker\Env\alpha$ оболочки $\Env\alpha$ морфизма $\alpha$ является оболочкой $\Env\Coker\alpha$ объекта $\Coker\alpha$ в классе $\varOmega$ относительно класса $\varPhi$:
\beq\label{Coker=Env}
\Coker\Env\alpha\cong \Env\Coker\alpha.
\eeq
\etm
\bpr
1. Покажем сначала, что морфизм $\Coker\frac{\Env\alpha,\env X}{\env Y,\alpha}$ является расширением объекта $\Coker\alpha$ (здесь нам понадобится условие \eqref{Coker-in-Omega}). Выберем тестовый морфизм $\ph:\Coker\alpha\to B$, $\ph\in\varPhi$.
Наши построения будут иллюстрироваться диаграммой
$$
 \xymatrix @R=3.5pc @C=2.pc
 {
 X \ar[rr]^{\env X}\ar[d]_{\alpha} &  & \Env X \ar[d]^{\Env\alpha}   \\
 Y \ar[rr]^{\env Y}\ar[d]^{\coker\alpha}\ar@/_15ex/[ddr]^{\psi} & & \Env Y \ar[d]_{\coker\Env\alpha}\ar@/^15ex/@{-->}[ddl]_{\psi'}   \\
 \Coker\alpha\ar[rr]^{\Coker\frac{\Env\alpha,\env X}{\env Y,\alpha}}\ar@/_3ex/[dr]^{\ph} & & \Coker\Env\alpha
  \ar@/^3ex/@{-->}[dl]_{\ph'} \\
  &  B\ar@/_2ex/[d]_{\eta}\ar@/^2ex/[d]^{\theta} & \\
  & C  &
 }
$$
Поскольку $\varPhi$ -- правый идеал\footnote{Это постулируется в условии RE.4 на с.\pageref{DEF:RE.4}.} в $\tt K$, композиция $\psi=\ph\circ\coker\alpha$ является тестовым морфизмом для $Y$ (то есть $\psi\in\varPhi$), поэтому существует единственное продолжение $\psi'$ морфизма $\psi$ вдоль оболочки $\env Y$. Для произвольной пары морфизмов $\eta,\theta:B\rightrightarrows C$ мы получим:
\begin{multline*}
\eta\circ\psi'\circ\Env\alpha\circ\env X=\eta\circ\psi'\circ\Env Y\circ\alpha=\eta\circ\psi\circ\alpha=\eta\circ\ph\circ\coker\alpha\circ\alpha=\\=
\theta\circ\ph\circ\coker\alpha\circ\alpha=\theta\circ\psi\circ\alpha=
\theta\circ\psi'\circ\Env Y\circ\alpha=\theta\circ\psi'\circ\Env\alpha\circ\env X,
\end{multline*}
и, поскольку $\env X\in\varOmega$ -- эпиморфизм\footnote{Класс $\varOmega$ состоит из эпиморфизмов в силу условия RE.2 на с.\pageref{DEF:RE.2}.}, его можно отбросить, и мы получим равенство
$$
\eta\circ\psi'\circ\Env\alpha=\theta\circ\psi'\circ\Env\alpha.
$$
Это верно для любых $\eta,\theta:B\rightrightarrows C$, значит, $\psi'$ тривиализирует  $\Env\alpha$ снаружи. Поэтому, $\psi'$ продолжается до некоторого морфизма $\ph'$ вдоль коядра $\coker\Env\alpha$:
$$
\psi'=\ph'\circ \coker\Env\alpha
$$
Для этого морфизма мы получим:
$$
\ph'\circ\Coker\frac{\Env\alpha,\env X}{\env Y,\alpha}\circ\coker\alpha=
\ph'\circ\coker\Env\alpha\circ\env Y=
\psi'\circ\env Y=\psi=\ph\circ\coker\alpha.
$$
По предложению \ref{PROP:epimorfnost-coker}, $\coker\alpha$ -- эпиморфизм, поэтому его можно убрать, и мы получим
$$
\ph'\circ\Coker\frac{\Env\alpha,\env X}{\env Y,\alpha}=\ph.
$$
Этот морфизм $\ph'$ единственен, потому что $\Coker\frac{\Env\alpha,\env X}{\env Y,\alpha}\in\varOmega\subseteq\Epi$.

2. Из того, что $\Coker\frac{\Env\alpha,\env X}{\env Y,\alpha}$ -- расширение, следует, что существует морфизм $\upsilon$, замыкающий нижний внутренний треугольник в диаграмме
$$
 \xymatrix @R=2.pc @C=6.pc
 {
 Y \ar[r]^{\env Y}\ar[d]^{\coker\alpha} & \Env Y \ar[d]_{\coker\Env\alpha}\ar@/^4ex/[ddr]^{\quad\Env\coker\alpha} & \\
  \Coker\alpha\ar[r]_{\Coker\frac{\Env\alpha,\env X}{\env Y,\alpha}}\ar@/_6ex/[drr]^{\qquad \env\Coker\alpha}
   & \Coker\Env\alpha\ar@{-->}[dr]^{\upsilon} & \\
  & &  \Env\Coker\alpha
  }
$$
В ней также коммутативны внутренний квадрат и периметр. Поскольку вдобавок морфизм $\env Y$ лежит в классе $\varOmega$ и значит, является эпиморфизмом, правый внутренний треугольник также должен быть коммутативен.

3. Пусть теперь выполнено условие (iv). Тогда морфизм 
$$
\Coker\frac{\Env\alpha,\env X}{\env Y,\alpha}:\Coker\alpha\to \Coker\Env\alpha
$$
является расширением в полный объект (в классе морфизмов $\varOmega$ относительно класса морфизмов $\varPhi$) и по предложению \ref{PROP:rasshirenie-v-polnyj-objekt}, этот морфизм должен быть изоморфизмом.
\epr

\section{Оболочки в моноидальных категориях}

\subsection{Оболочки, согласованные с тензорным произведением}

Пусть $\tt K$ -- моноидальная категория \cite{MacLane} с тензорным произведением $\otimes$ и единичным объектом $I$.
 \bit{\label{DEF:obolochka-soglasovana-s-tenz-proizv}
\item[$\bullet$] Условимся говорить, что {\it оболочка $\Env^\varOmega_\varPhi$ согласована с тензорным произведением $\otimes$} на $\tt K$, если выполняются следующие два условия:

 \bit{

\item[T.1]\label{DEF:T.1} Тензорное произведение $\rho\otimes\sigma:X\otimes Y\to X'\otimes Y'$ любых
двух расширений $\rho:X\to X'$ и $\sigma:Y\to Y'$ (в классе $\varOmega$
относительно класса $\varPhi$) является расширением (в $\varOmega$ относительно $\varPhi$).

\item[T.2]\label{DEF:T.2} Локальная единица $1_I:I\to I$ единичного объекта $I$
является оболочкой (в классе $\varOmega$ относительно класса $\varPhi$):
\beq\label{E(I)=I}
\env_\varPhi^\varOmega I=1_I
\eeq
 }\eit
 }\eit

Всюду ниже в этом параграфе мы будем считать, что классы $\varOmega$ и $\varPhi$ определяют регулярную оболочку
в категории $\tt K$. По теореме
\ref{TH:regulyarnaya-obolochka} это означает, что $\Env_\varPhi^\varOmega$ можно определить как идемпотентный функтор.
Мы будем обозначать его $E:{\tt K}\to{\tt K}$, а естественное преобразование тождественного функтора в $E$ мы обозначаем
$e$:
$$
E(X):=\Env_\varPhi^\varOmega X, \qquad E(\ph):=\Env_\varPhi^\varOmega\ph,\qquad e_X:=\env_\varPhi^\varOmega X.
$$
Класс всех полных объектов в $\tt K$ (в $\varOmega$ относительно $\varPhi$) мы как и раньше обозначаем $\tt L$.

\btm\label{PROP:E(e_X-otimes-e_Y)-in-Iso} Пусть $\Env_\varPhi^\varOmega$ -- регулярная оболочка,
согласованная с  тензорным произведением в ${\tt K}$. Тогда
 \bit{

\item[(i)] для любых двух объектов $A\in{\tt L}$ и $X\in\Ob({\tt K})$, оболочка $E(1_A\otimes e_X)$ морфизма $1_A\otimes e_X:A\otimes X\to A\otimes E(X)$ является изоморфизмом (в $\tt K$ и в $\tt L$):
 \beq\label{E(1_A-otimes-e_X)-in-Iso}
 E(1_A\otimes e_X)\in\Iso.
 \eeq

\item[(ii)] для любых двух объектов $X,Y\in\Ob({\tt K})$ оболочка $E(e_X\otimes e_Y)$ морфизма $e_X\otimes e_Y:X\otimes Y\to E(X)\otimes E(Y)$ является изоморфизмом (в $\tt K$ и в $\tt L$):
 \beq\label{E(e_X-otimes-e_Y)-in-in-Iso}
 E(e_X\otimes e_Y)\in\Iso,
 \eeq
 как следствие,
 \beq\label{E(X-otimes-Y)-cong-E(E(X)-otimes-E(Y))}
 E(X\otimes Y)\cong E(E(X)\otimes E(Y))
 \eeq

 }\eit
\etm
\bpr
1. Пусть $A\in{\tt L}$ и $X\in\Ob({\tt K})$. Морфизмы $1_A:A\to A$ и $e_X:X\to E(X)$ в произведении дают морфизм $1_A\otimes e_X:A\otimes X\to A\otimes E(X)$. Подставим его вместо $\alpha$ в диаграмму \eqref{DIAGR:funktorialnost-env-e-E}:
\beq\label{DIAGR:svyaz-otimes-s-env-0}
\xymatrix @R=2.pc @C=5.0pc 
{
A\otimes X \ar[d]^{1_A\otimes e_X}\ar[r]^{e_{A\otimes X}} & E(A\otimes X)\ar[d]^{E(1_A\otimes e_X)} \\
A\otimes E(X)\ar[r]^{e_{A\otimes E(X)}} & E(A\otimes E(X)) \\
}
\eeq
Из диаграммы
$$
\xymatrix @R=2.pc @C=5.0pc 
{
A\otimes X \ar@/_2ex/[dr]_{\ph}\ar[r]^{1_A\otimes e_X} & A\otimes E(X)\ar[r]^{e_{A\otimes E(X)}}\ar@{-->}[d]^{\ph'} & E(A\otimes E(X))\ar@/^2ex/@{-->}[dl]^{\ph''} \\
 & B &
}
$$
видно, что композиция $e_{A\otimes E(X)}\circ 1_A\otimes e_X$ является расширением для $A\otimes X$
(здесь в левом треугольнике используется T.1). Поэтому $e_{A\otimes E(X)}\circ 1_A\otimes e_X$
вписывается в оболочку объекта $A\otimes X$:
\beq\label{DIAGR:svyaz-otimes-s-env-00}
\xymatrix @R=2.pc @C=5.0pc 
{
A\otimes X \ar[d]^{1_A\otimes e_X}\ar[r]^{e_{A\otimes X}} & E(A\otimes X) \\
A\otimes E(X)\ar[r]^{e_{A\otimes E(X)}} & E(A\otimes E(X))\ar@{-->}[u]_{\upsilon} \\
}
\eeq
для некоторого (единственного) морфизма $\upsilon$. Вдобавок, поскольку $\varOmega\subseteq\Epi$, морфизмы
$e_{A\otimes E(X)}\circ 1_A\otimes e_X$ и $e_{A\otimes X}$, будучи расширениями, являются эпиморфизмами. Поэтому
\eqref{DIAGR:svyaz-otimes-s-env-0} и \eqref{DIAGR:svyaz-otimes-s-env-00} вместе дают
$$
\upsilon=E(1_A\otimes e_X)^{-1}.
$$

2. Для любых двух объектов $X$ и $Y$ морфизмы $e_X:X\to E(X)$ и $e_Y:Y\to E(Y)$ в произведении дают морфизм $e_X\otimes e_Y:X\otimes Y\to E(X)\otimes E(Y)$. Подставим его вместо $\alpha$ в диаграмму \eqref{DIAGR:funktorialnost-env-e-E}:
\beq\label{DIAGR:svyaz-otimes-s-env}
\xymatrix @R=2.pc @C=5.0pc 
{
X\otimes Y \ar[d]^{e_X\otimes e_Y}\ar[r]^{e_{X\otimes Y}} & E(X\otimes Y)\ar[d]^{E(e_X\otimes e_Y)} \\
E(X)\otimes E(Y)\ar[r]^{e_{E(X)\otimes E(Y)}} & E(E(X)\otimes E(Y)) \\
}
\eeq
Из диаграммы
$$
\xymatrix @R=2.pc @C=5.0pc 
{
X\otimes Y \ar@/_2ex/[dr]_{\ph}\ar[r]^{e_X\otimes e_Y} & E(X)\otimes E(Y)\ar[r]^{e_{E(X)\otimes E(Y)}}\ar@{-->}[d]^{\ph'} & E(E(X)\otimes E(Y))\ar@/^2ex/@{-->}[dl]^{\ph''} \\
 & B &
}
$$
видно, что композиция $e_{E(X)\otimes E(Y)}\circ e_X\otimes e_Y$ является расширением для $X\otimes Y$
(здесь в левом треугольнике используется T.1). Поэтому $e_{E(X)\otimes E(Y)}\circ e_X\otimes e_Y$
вписывается в оболочку объекта $X\otimes Y$:
\beq\label{DIAGR:svyaz-otimes-s-env-*}
\xymatrix @R=2.pc @C=5.0pc 
{
X\otimes Y \ar[d]^{e_X\otimes e_Y}\ar[r]^{e_{X\otimes Y}} & E(X\otimes Y) \\
E(X)\otimes E(Y)\ar[r]^{e_{E(X)\otimes E(Y)}} & E(E(X)\otimes E(Y))\ar@{-->}[u]_{\upsilon} \\
}
\eeq
для некоторого (единственного) морфизма $\upsilon$. Как и в предыдущем случае,
морфизмы $e_{E(X)\otimes E(Y)}\circ e_X\otimes e_Y$ и $e_{X\otimes Y}$, будучи расширениями,
являются эпиморфизмами, поэтому диаграммы
\eqref{DIAGR:svyaz-otimes-s-env} и \eqref{DIAGR:svyaz-otimes-s-env-*} вместе дают
$$
\upsilon=E(e_X\otimes e_Y)^{-1}.
$$
\epr

\bcor\label{COR:e_X-otimes-e_Y-in-Mono} Пусть $\Env_\varPhi^\varOmega$ -- регулярная оболочка,
согласованная с  тензорным произведением в ${\tt K}$. Тогда для любых двух объектов $X,Y\in\Ob({\tt K})$ если $e_{X\otimes Y}\in\Mono$, то $e_X\otimes e_Y\in\Mono$.
\ecor
\bpr
Из того, что $E(e_X\otimes e_Y)$ --- изоморфизм (то есть из условия \eqref{E(e_X-otimes-e_Y)-in-in-Iso}) следует, что диаграмму \eqref{DIAGR:svyaz-otimes-s-env} можно превратить в диаграмму
\beq\label{DIAGR:svyaz-otimes-s-env-1}
\xymatrix @R=2.pc @C=5.0pc 
{
X\otimes Y \ar[d]^{e_X\otimes e_Y}\ar[r]^{e_{X\otimes Y}} & E(X\otimes Y) \\
E(X)\otimes E(Y)\ar[r]^{e_{E(X)\otimes E(Y)}} & E(E(X)\otimes E(Y))\ar[u]_{E(e_X\otimes e_Y)^{-1}} \\
}
\eeq
Теперь, используя главное свойство мономорфизмов \cite[$1^\circ$, p.80]{Akbarov-De-Gruyter-I}, мы получаем импликацию
$$
E(e_X\otimes e_Y)^{-1}\circ e_{E(X)\otimes E(Y)}\circ e_X\otimes e_Y=e_{X\otimes Y}\in\Mono
\quad\Rightarrow\quad e_X\otimes e_Y\in\Mono.
$$
\epr

\paragraph{Моноидальная структура на классе полных объектов.}
Пусть $\Env_\varPhi^\varOmega$ -- регулярная оболочка, согласованная с тензорным произведением в $\tt K$,
 и $E=\Env_\varPhi^\varOmega$ -- соответствующий идемпотентный функтор,
построенный в теореме \ref{TH:regulyarnaya-obolochka}. Пусть по-прежнему ${\tt L}$ -- (полная) подкатегория
полных объектов в $\tt K$. Для любых двух объектов $A,B\in{\tt L}$ и любых двух морфизмов
$\ph,\psi\in{\tt L}$ обозначим
\beq\label{A-overset(E)(otimes)B:=E(A-otimes-B)}
A\overset{E}{\otimes}B:=E(A\otimes B),\qquad \ph\overset{E}{\otimes}\psi:=E(\ph\otimes\psi).
\eeq
Заметим тождество
\beq\label{E(E(X)-otimes-E(Y))=E(X)-overset(E)(otimes)E(Y)}
E(X)\overset{E}{\otimes}E(Y)=E(E(X)\otimes E(Y)),\qquad X,Y\in \Ob(\tt K)
\eeq
(это будет равенство объектов, потому что по предложению \ref{PROP:harakterizatsiya-polnoty}, всегда
$E(X),E(Y)\in{\tt L}$).

\btm\label{TH:sushestvovanie-tenz-proizv-v-L}
Пусть $\Env_\varPhi^\varOmega$ -- регулярная оболочка, согласованная с тензорным произведением в $\tt K$.
Тогда формулы
\eqref{A-overset(E)(otimes)B:=E(A-otimes-B)} определяют на категории $\tt L$
структуру моноидальной категории (с тензорным произведением $\overset{E}{\otimes}$ и единичным объектом $I$).
\etm
\bpr
1. Тензорное произведение локальных единиц должно быть локальной единицей. Подставим в диаграмму \eqref{DIAGR:funktorialnost-env-e-E} морфизм $1_{A\otimes B}$ вместо $\alpha$:
$$
\xymatrix @R=2.pc @C=5.0pc 
{
A\otimes B\ar[d]^{1_{A\otimes B}}\ar[r]^{e(A\otimes B)} & A\overset{E}{\otimes}B\ar[d]^{E(1_{A\otimes B})} \\
A\otimes B\ar[r]^{e(A\otimes B)} & A\overset{E}{\otimes}B \\
}
$$
Если вместо правой вертикальной стрелки $E(1_{A\otimes B})$ подставить морфизм $1_{A\overset{E}{\otimes}B}$, то диаграмма останется коммутативной. Поскольку такая вертикальная стрелка единствена (из-за того, что $e(A\otimes B)$ -- эпиморфизм), эти стрелки должны совпадать, и это используется в последнем равенстве в цепочке
$$
1_A\overset{E}{\otimes} 1_B=\eqref{A-overset(E)(otimes)B:=E(A-otimes-B)}=E(1_A\otimes 1_B)=E(1_{A\otimes B})=1_{A\overset{E}{\otimes}B}.
$$

2. Тензорное произведение коммутативных диаграмм должно быть коммутативной диаграммой. Пусть даны две коммутативные диаграммы в ${\tt L}$:
$$
\begin{diagram}
\node[2]{B}\arrow{se,l}{\chi} \\
\node{A}\arrow{ne,l}{\ph}\arrow[2]{e,r}{\psi}\node[2]{C}
\end{diagram}\qquad
\begin{diagram}
\node[2]{B'}\arrow{se,l}{\chi'} \\
\node{A'}\arrow{ne,l}{\ph'}\arrow[2]{e,r}{\psi'}\node[2]{C'}
\end{diagram}
$$
Перемножая их тензорно в ${\tt K}$, мы получим коммутативную диаграмму
$$
\begin{diagram}
\node[2]{B\otimes B'}\arrow{se,l}{\chi\otimes \chi'} \\
\node{A\otimes A'}\arrow{ne,l}{\ph\otimes
\ph'}\arrow[2]{e,r}{\psi\otimes\psi'}\node[2]{C\otimes C'}
\end{diagram}
$$
Затем применяя функтор $E$ получаем коммутативную диаграмму
$$
\begin{diagram}
\node[2]{E(B\otimes B')}\arrow{se,l}{E(\chi\otimes \chi')} \\
\node{E(A\otimes A')}\arrow{ne,l}{E(\ph\otimes
\ph')}\arrow[2]{e,r}{E(\psi\otimes\psi')}\node[2]{E(C\otimes C')}
\end{diagram}
$$
В силу \eqref{A-overset(E)(otimes)B:=E(A-otimes-B)} это будет нужная нам диаграмма
$$
\begin{diagram}
\node[2]{B\overset{E}{\otimes} B'}\arrow{se,l}{\chi\overset{E}{\otimes}\chi'} \\
\node{A\overset{E}{\otimes}A'}\arrow{ne,l}{\ph\overset{E}{\otimes}\ph'}
\arrow[2]{e,r}{\psi\overset{E}{\otimes}\psi'}\node[2]{C\overset{E}{\otimes}C'}
\end{diagram}
$$

3. Заметим, что из доказанного уже следует, что тензорное произведение изоморфизмов в $\tt L$ также является изоморфизмом:
\beq\label{Iso-E/otimes-Iso-subseteq-Iso}
\ph,\psi\in\Iso \qquad\Longrightarrow\qquad \ph\overset{E}{\otimes}\psi:=E(\ph\otimes\psi)\in\Iso
\eeq
Действительно,
$$
(\ph\otimes\psi)\circ(\ph^{-1}\otimes\psi^{-1})=(\ph\circ\ph^{-1})\otimes(\psi\circ\psi^{-1})=1\otimes 1=1,
$$
поэтому
$$
(\ph\overset{E}{\otimes}\psi)\circ (\ph^{-1}\overset{E}{\otimes}\psi^{-1})=E(\ph\otimes\psi)\circ E(\ph^{-1}\otimes\psi^{-1})=E((\ph\otimes\psi)\circ(\ph^{-1}\otimes\psi^{-1}))=E(1)=1.
$$
И точно так же
$$
(\ph^{-1}\overset{E}{\otimes}\psi^{-1})\circ (\ph\overset{E}{\otimes}\psi)=1.
$$

4. Если $\alpha_{A,B,C}: (A\otimes B)\otimes C\to A\otimes (B\otimes C)$ -- преобразование ассоциативности в $\tt K$, то преобразование ассоциативности $\alpha^E_{A,B,C}: (A\overset{E}{\otimes}B)\overset{E}{\otimes}C\to
A\overset{E}{\otimes}(B\overset{E}{\otimes}C)$ в $\tt L$ определяется диаграммой
\beq\label{DIAGR:opred-alpha-v-L}
\xymatrix @R=2.pc @C=5.0pc 
{
E((A\otimes B)\otimes C) \ar[r]^{E(\alpha_{A,B,C})} & E(A\otimes (B\otimes C))\ar[d]^{E(1_A\otimes e_{B\otimes C})} \\
E(E(A\otimes B)\otimes C)\ar[u]^{E(e_{A\otimes B}\otimes 1_C)^{-1}} & E(A\otimes E(B\otimes C)) \\
E(A\otimes B)\overset{E}{\otimes}C\ar@{=}[u] & A\overset{E}{\otimes}E(B\otimes C)\ar@{=}[u] \\
(A\overset{E}{\otimes}B)\overset{E}{\otimes}C\ar@{=}[u]\ar[r]^{\alpha^E_{A,B,C}} & A\overset{E}{\otimes}
(B\overset{E}{\otimes}C)\ar@{=}[u] \\
}
\eeq
(здесь мы воспользовались теоремой \ref{PROP:E(e_X-otimes-e_Y)-in-Iso}, из которой следует, что морфизм $E(e_{A\otimes B}\otimes 1_C)$ должен иметь обратный, и тождеством \eqref{E(E(X)-otimes-E(Y))=E(X)-overset(E)(otimes)E(Y)}).

5. Покажем, что определенное таким образом преобразование $\alpha^E$ естественно относительно тензорного произведения, то есть является морфизмом функторов:
$$
\alpha^E:\Big((A,B,C)\mapsto (A\overset{E}{\otimes}
B)\overset{E}{\otimes} C\Big)\rightarrowtail\Big((A,B,C)\mapsto A\overset{E}{\otimes} (B\overset{E}{\otimes} C)\Big)
$$
Пусть даны морфизмы $\ph:A\to A'$, $\chi:B\to B'$, $\psi: C\to C'$ в ${\tt L}$, тогда к диаграмме естественности для $\alpha$
 \beq\label{DIAGR:estestvennost-associativnosti}
\xymatrix @R=2.pc @C=5.0pc 
{
(A\otimes B)\otimes C\ar[d]_{(\ph\otimes \chi)\otimes\psi}
\ar[r]^{\alpha_{A,B,C}} & A\otimes (B\otimes C)\ar[d]_{\ph\otimes
(\chi\otimes\psi)}
\\
(A'\otimes B')\otimes C'\ar[r]_{\alpha_{A',B',C'}}& A'\otimes(B'\otimes C')
}
 \eeq
мы применяем функтор $E$:
$$
\xymatrix @R=2.pc @C=5.0pc 
{
E((A\otimes B)\otimes C)\ar[d]_{E((\ph\otimes \chi)\otimes\psi)}
\ar[r]^{E(\alpha_{A,B,C})} & E(A\otimes (B\otimes C))\ar[d]_{E(\ph\otimes(\chi\otimes\psi))}
\\
E((A'\otimes B')\otimes C')\ar[r]_{E(\alpha_{A',B',C'})}& E(A'\otimes(B'\otimes C'))
}
$$
и достраиваем полученную диаграмму так:
$$
\xymatrix @R=3.pc @C=4.0pc 
{
(A\overset{E}{\otimes}B)\overset{E}{\otimes}C \ar@{-->}[r]_{E(e_{A\otimes B}\otimes 1_C)^{-1}}
\ar@{-->}[d]_{(\ph\overset{E}{\otimes}\chi)\overset{E}{\otimes}\psi}\ar@{-->}@/^4ex/[rrr]^{\alpha^E_{A,B,C}}
 & E((A\otimes B)\otimes C)\ar[d]_{E((\ph\otimes \chi)\otimes\psi)}
\ar[r]^{E(\alpha_{A,B,C})} & E(A\otimes (B\otimes C))\ar[d]_{E(\ph\otimes(\chi\otimes\psi))}
\ar@{-->}[r]_{E(1_A\otimes e_{B\otimes C})}
& A\overset{E}{\otimes}(B\overset{E}{\otimes}C)
\ar@{-->}[d]_{\ph\overset{E}{\otimes}(\chi\overset{E}{\otimes}\psi)}
\\
(A'\overset{E}{\otimes}B')\overset{E}{\otimes}C'\ar@{-->}[r]^{E(e_{A'\otimes B'}\otimes 1_{C'})^{-1}}
\ar@{-->}@/_4ex/[rrr]_{\alpha^E_{A',B',C'}}
 & E((A'\otimes B')\otimes C')\ar[r]_{E(\alpha_{A',B',C'})}& E(A'\otimes(B'\otimes C'))
 \ar@{-->}[r]^{E(1_{A'}\otimes e_{B'\otimes C'})}
& A'\overset{E}{\otimes}(B'\overset{E}{\otimes}C')
}
$$
Выбросив из нее внутренние вершины, мы получим диаграмму естественности для $\alpha^E$:
$$
\xymatrix @R=2.pc @C=5.0pc 
{
(A\overset{E}{\otimes} B)\overset{E}{\otimes} C\ar[d]_{(\ph\overset{E}{\otimes} \chi)\overset{E}{\otimes}\psi}
\ar[r]^{\alpha^E_{A,B,C}} & A\overset{E}{\otimes} (B\overset{E}{\otimes} C)\ar[d]_{\ph\overset{E}{\otimes}
(\chi\overset{E}{\otimes}\psi)}
\\
(A'\overset{E}{\otimes} B')\overset{E}{\otimes} C'\ar[r]_{\alpha^E_{A',B',C'}}& A'\overset{E}{\otimes}(B'\overset{E}{\otimes} C')
}
$$

6. Покажем далее, что преобразовение $\alpha^E$ удовлетворяет условию ассоциативности, то есть замыкает стандартный пятиугольник. Для этого выпишем эту диаграмму для $\alpha$:
 \beq\label{5-ugolnik-dlya-otimes}
\xymatrix  @R=1pc @C=-1pc 
{
 & \text{\scriptsize $(A\otimes(B\otimes C))\otimes D$}\ar[rr]^{\alpha_{A,B\otimes C, D}}
 & & \text{\scriptsize $A\otimes((B\otimes C)\otimes D)$}\ar[ddr]!U|{1_A\otimes \alpha_{B,C,D}} &
 \\
 & & \\
\text{\scriptsize $((A\otimes B)\otimes C)\otimes D$}
 \ar[uur]!U|{\alpha_{A,B,C}\otimes 1_D}\ar[drr]_{\alpha_{A\otimes B,C,D}}
 &&&& \text{\scriptsize $A\otimes(B\otimes (C\otimes D))$} \\
  && \text{\scriptsize $(A\otimes B)\otimes (C\otimes D)$}\ar[urr]_{\alpha_{A,B,C\otimes D}}&&
}
 \eeq
Применим к ней функтор $E$ и достроим полученную диаграмму до следующей призмы:
 \beq\label{prizma-dlya-otimes^E}
\xymatrix  @R=1pc @C=0pc  
{
 & \text{\tiny $E((A\otimes(B\otimes C))\otimes D)$}\ar[d]^{\text{\tiny $E(e_{A\otimes(B\otimes C)}\otimes 1_D)$}} \ar[rr]^{\text{\tiny $E(\alpha_{A,B\otimes C, D})$}}
 & & \text{\tiny $E(A\otimes((B\otimes C)\otimes D))$}\ar@/^2ex/[ddr]^{\text{\tiny $E(1_A\otimes \alpha_{B,C,D})$}}
 \ar[d]_{\text{\tiny $E(1_A\otimes e_{(B\otimes C)\otimes D})$}} &
 \\
& \text{\tiny $E(E(A\otimes(B\otimes C))\otimes D)$}\ar@{=}[d] & & \text{\tiny $E(A\otimes E((B\otimes C)\otimes D))$}\ar@{=}[d] & \\
\text{\tiny $E(((A\otimes B)\otimes C)\otimes D)$}
\ar[d]_{\text{\tiny $E(e_{(A\otimes B)\otimes C}\otimes 1_D)$}}
 \ar@/^2ex/[uur]^{\text{\tiny $E(\alpha_{A,B,C}\otimes 1_D)$}}\ar[ddrr]^(.75){\text{\tiny $E(\alpha_{A\otimes B,C,D})$}}
 &
 \text{\tiny $E(A\otimes(B\otimes C))\overset{E}{\otimes}D$}\ar[dd]^(.3){\text{\tiny $E(1_A\otimes e_{B\otimes C})\overset{E}{\otimes}1_D$}}|!{[d]}{\hole}
 &&
 \text{\tiny $A\overset{E}{\otimes}E((B\otimes C)\otimes D)$}\ar[dd]_(.3){\text{\tiny $1_A\overset{E}{\otimes}E(e_{B\otimes C}\otimes 1_D)$}}|!{[d]}{\hole}
 & \text{\tiny $E(A\otimes(B\otimes (C\otimes D)))$}\ar[d]^{\text{\tiny $E(1_A\otimes e_{B\otimes (C\otimes D)})$}} \\
 \text{\tiny $E(E((A\otimes B)\otimes C)\otimes D)$}\ar@{=}[d]
 &
 &
 &
 &
 \text{\tiny $E(A\otimes E(B\otimes (C\otimes D)))$}\ar@{=}[d]
 \\
  \text{\tiny $E((A\otimes B)\otimes C)\overset{E}{\otimes}D$}
   \ar[d]_{\text{\tiny $E(e_{A\otimes B}\otimes 1_C)\overset{E}{\otimes} 1_D$}}
    & \text{\tiny $E(A\otimes E(B\otimes C))\overset{E}{\otimes}D$}\ar@{=}[d]
  & \text{\tiny $E((A\otimes B)\otimes (C\otimes D))$}
  \ar[ddd]_(.7){\text{\tiny $E(e_{A\otimes B}\otimes e_{C\otimes D})$}}
  \ar[uurr]^(.25){\text{\tiny $E(\alpha_{A,B,C\otimes D})$}}&
  \text{\tiny $A\overset{E}{\otimes}E(E(B\otimes C)\otimes D)$}\ar@{=}[d]
  &
   \text{\tiny $A\overset{E}{\otimes}E(B\otimes (C\otimes D))$}
   \ar[d]^{\text{\tiny $1_A\overset{E}{\otimes} E(1_B\otimes e_{C\otimes D})$}}
    \\
 \text{\tiny $E(E(A\otimes B)\otimes C)\overset{E}{\otimes}D$}\ar@{=}[dd]
  & \text{\tiny $(A\overset{E}{\otimes}(B\overset{E}{\otimes}C))\overset{E}{\otimes}D$}\ar[rr]_(.7){\text{\tiny $\alpha^E_{A,B\overset{E}{\otimes} C, D}$}}|!{[r]}{\hole}
 & & \text{\tiny $A\overset{E}{\otimes}((B\overset{E}{\otimes} C)\overset{E}{\otimes} D)$}\ar[ddr]_(.3){\text{\tiny $1_A\overset{E}{\otimes} \alpha^E_{B,C,D}$}} &
  \text{\tiny $A\overset{E}{\otimes}E(B\otimes E(C\otimes D))$}\ar@{=}[dd]
  \\
 & &
     & & \\
\text{\tiny $((A\overset{E}{\otimes} B)\overset{E}{\otimes} C)\overset{E}{\otimes} D$}
 \ar[uur]_(.7){\text{\tiny $\alpha^E_{A,B,C}\overset{E}{\otimes} 1_D$}}\ar[ddrr]_{\text{\tiny $\alpha^E_{A\overset{E}{\otimes} B,C,D}$}}
 &&
   \text{\tiny $E(E(A\otimes B)\otimes E(C\otimes D))$}\ar@{=}[dd]
 && \text{\tiny $A\overset{E}{\otimes}(B\overset{E}{\otimes} (C\overset{E}{\otimes} D))$}
  \\
 &&  &&
 \\
  && \text{\tiny $(A\overset{E}{\otimes} B)\overset{E}{\otimes} (C\overset{E}{\otimes} D)$}\ar[uurr]_{\text{\tiny $\alpha^E_{A,B,C\overset{E}{\otimes} D}$}}&&
}
 \eeq
У этой призмы верхнее основание коммутативно, потому что это результат применения функтора $E$ к диаграмме \eqref{5-ugolnik-dlya-otimes}, а коммутативность боковых граней нетрудно проверить, изменив (эквивалентным образом) вид вертикальных стрелок.

Например, коммутативность левой ближней грани станет очевидной, если ее представить как периметр следующией диаграммы:
$$
\xymatrix  @R=2pc @C=4pc 
{
 &
 \text{\tiny $E(((A\otimes B)\otimes C)\otimes D)$}\ar[r]^{\text{\tiny $E(\alpha_{A\otimes B,C,D})$}}
 \ar[d]_{\text{\tiny $E((e_{A\otimes B}\otimes 1_C)\otimes 1_D)$}}
 \ar@/_3ex/[dl]_(.7){\text{\tiny $E(e_{(A\otimes B)\otimes C}\otimes 1_D)$}\quad}
 & \text{\tiny $E((A\otimes B)\otimes (C\otimes D))$}
 \ar[d]_{\text{\tiny $E(e_{A\otimes B}\otimes (1_C\otimes 1_D))$}}  \ar@/^4ex/[dddr]^{\text{\tiny $E(e_{A\otimes B}\otimes e_{C\otimes D})$}}
 &
 \\
 \text{\tiny $E(E((A\otimes B)\otimes C)\otimes D)$}
 \ar@{=}[d] & \text{\tiny $E((E(A\otimes B)\otimes C)\otimes D)$}\ar@{=}[d] &
 \text{\tiny $E(E(A\otimes B)\otimes (C\otimes D))$}\ar@{=}[d]
 &
 \\
 \text{\tiny $E((A\otimes B)\otimes C)\overset{E}{\otimes}D$}
 \ar[d]^{\text{\tiny $E(e_{A\otimes B}\otimes 1_C)\overset{E}{\otimes}1_D$}}
 & \text{\tiny $E(((A\overset{E}{\otimes}B)\otimes C)\otimes D)$}
 \ar[r]^{\text{\tiny $E\Big(\alpha_{A\overset{E}{\otimes}B,C,D}\Big)$}}
 \ar[d]^{\text{\tiny $E(e_{(A\overset{E}{\otimes}B)\otimes C}\otimes 1_D)$}}
 & \text{\tiny $E((A\overset{E}{\otimes}B)\otimes (C\otimes D))$}
 \ar[d]^{\text{\tiny $E(1_{A\overset{E}{\otimes}B}\otimes e_{C\otimes D})$}}
 &
 \\
 \text{\tiny $E(E(A\otimes B)\otimes C)\overset{E}{\otimes}D$}
 \ar@/_3ex/@{=}[dr]
 & \text{\tiny $E(E((A\overset{E}{\otimes}B)\otimes C)\otimes D)$}\ar@{=}[d]
 & \text{\tiny $E((A\overset{E}{\otimes}B)\otimes E(C\otimes D))$}\ar@{=}[d]
 &
 \text{\tiny $E(E(A\otimes B)\otimes E(C\otimes D))$}\ar@{=}[l]  \ar@/^3ex/@{=}[dl]
 \\
 & \text{\tiny $((A\overset{E}{\otimes}B)\overset{E}{\otimes}C)\overset{E}{\otimes}D$}
\ar[r]^{\text{\tiny $\alpha^E_{A\overset{E}{\otimes}B,C,D}$}}
 & \text{\tiny $(A\overset{E}{\otimes}B)\overset{E}{\otimes}(C\overset{E}{\otimes} D)$}
 &
}
$$
Здесь верхний внутренний шестиугольник (или его можно назвать квадратом) -- результат применения функтора $E$ к диаграмме
 $$
\xymatrix  @R=2pc @C=4pc 
{
 ((A\otimes B)\otimes C)\otimes D \ar[r]^{\alpha_{A\otimes B,C,D}} \ar[d]_{(e_{A\otimes B}\otimes 1_C)\otimes 1_D} & (A\otimes B)\otimes (C\otimes D)\ar[d]^{e_{A\otimes B}\otimes (1_C\otimes 1_D)}
\\
 (E(A\otimes B)\otimes C)\otimes D\ar[r]^{\alpha_{E(A\otimes B),C,D}} &  E(A\otimes B)\otimes (C\otimes D)
}
$$
(это следствие диаграммы естественности \eqref{DIAGR:estestvennost-associativnosti} для $\alpha$). Нижний внутренний шестиугольник -- диаграмма \eqref{DIAGR:opred-alpha-v-L} для преобразования $\alpha^E$ на компонентах $A\overset{E}{\otimes}B$, $C$, $D$. Большой восьмиугольник слева можно представлять себе как ромб
 $$
\xymatrix  @R=2pc @C=2pc 
{
  & E(((A\otimes B)\otimes C)\otimes D) \ar[dl]_{E(e_{(A\otimes B)\otimes C}\otimes 1_D)\qquad}
  \ar[dr]^{\qquad E((e_{A\otimes B}\otimes 1_C)\otimes 1_D)}
  &
\\
E(E((A\otimes B)\otimes C)\otimes D)\ar[dr]_{E(E(e_{A\otimes B}\otimes 1_C)\otimes 1_D)\qquad}
 & & E((E(A\otimes B)\otimes C)\otimes D)\ar[dl]^{\qquad E(e_{E(A\otimes B)\otimes C}\otimes 1_D)}
\\
& E(E(E(A\otimes B)\otimes C)\otimes D) &
}
$$
который получается применением функтора $E$ к ромбу
 $$
\xymatrix  @R=2pc @C=2pc 
{
  &((A\otimes B)\otimes C)\otimes D \ar[dl]_{e_{(A\otimes B)\otimes C}\otimes 1_D\quad}
  \ar[dr]^{\quad (e_{A\otimes B}\otimes 1_C)\otimes 1_D}
  &
\\
E((A\otimes B)\otimes C)\otimes D\ar[dr]_{E(e_{A\otimes B}\otimes 1_C)\otimes 1_D\quad}
 & & (E(A\otimes B)\otimes C)\otimes D\ar[dl]^{\quad e_{E(A\otimes B)\otimes C}\otimes 1_D}
\\
& E(E(A\otimes B)\otimes C)\otimes D &
}
$$
который в свою очередь получается умножением справа на $D$ диаграммы
 $$
\xymatrix  @R=2pc @C=2pc 
{
  &(A\otimes B)\otimes C \ar[dl]_{e_{(A\otimes B)\otimes C}\quad}
  \ar[dr]^{\quad e_{A\otimes B}\otimes 1_C}
  &
\\
E((A\otimes B)\otimes C)\ar[dr]_{E(e_{A\otimes B}\otimes 1_C)\quad}
 & & E(A\otimes B)\otimes C\ar[dl]^{\quad e_{E(A\otimes B)\otimes C}}
\\
& E(E(A\otimes B)\otimes C) &
}
$$
а ее можно понимать как диаграмму \eqref{DIAGR:funktorialnost-env-e-E} с подставленным в нее вместо $\alpha$ морфизмом $e_{A\otimes B}\otimes 1_C$. Наконец, остающийся пятиугольник справа вверху есть треугольник
 $$
\xymatrix  @R=2pc @C=5pc 
{
  E((A\otimes B)\otimes(C\otimes D))\ar[d]_{E(e_{A\otimes B}\otimes (1_C\otimes 1_D))}
  \ar[dr]^{\qquad E(e_{A\otimes B}\otimes e_{C\otimes D})}
  &
\\
  E(E(A\otimes B)\otimes(C\otimes D))\ar[r]_{E(1_{A\otimes B}\otimes e_{C\otimes D})} &  E(E(A\otimes B)\otimes E(C\otimes D))
}
$$
получающийся применением функтора $E$ к треугольнику
 $$
\xymatrix  @R=2pc @C=5pc 
{
  (A\otimes B)\otimes(C\otimes D)\ar[d]_{e_{A\otimes B}\otimes (1_C\otimes 1_D)}
  \ar[dr]^{\qquad e_{A\otimes B}\otimes e_{C\otimes D}}
  &
\\
  E(A\otimes B)\otimes(C\otimes D)\ar[r]_{1_{A\otimes B}\otimes e_{C\otimes D}} &  E(A\otimes B)\otimes E(C\otimes D)
}
$$

Точно так же и с остальными вертикальными гранями в \eqref{prizma-dlya-otimes^E}. Поскольку вдобавок все вертикальные стрелки там -- изоморфизмы (в силу теоремы \ref{PROP:E(e_X-otimes-e_Y)-in-Iso} и свойства \eqref{Iso-E/otimes-Iso-subseteq-Iso}), мы получаем, что нижнее основание этой призмы также коммутативно, а оно и есть нужная диаграмма для $\alpha^E$.

7. Пусть $\lambda_X:I\otimes X\to X$ -- левая единица моноидальной категории $\tt K$, а $\rho_X:X\otimes I\to X$ -- правая единица. Положив для всякого $A\in\Ob({\tt L})$
\beq\label{DEF:lambda_A^E,rho_A^E}
\lambda_A^E=E(\lambda_A):I\overset{E}{\otimes}A=E(I\otimes A)\to E(A)=A,\qquad \rho_A^E=E(\rho_A):A\overset{E}{\otimes}I=E(A\otimes I)\to E(A)=A,
\eeq
мы получим левую и правую единицы для $\tt L$. Действительно, для всякого морфизма
$\ph:A\to A'$ в ${\tt L}$ диаграммы
 \beq\label{DIAGR:lambda,rho}
\xymatrix  @R=2pc @C=4pc 
{
I\otimes A \ar[d]_{1_I\otimes \ph} \ar[r]^{\lambda_A} & A\ar[d]^{\ph}
\\
I\otimes A'\ar[r]^{\lambda_{A'}} & A'
}
\qquad\qquad
\xymatrix  @R=2pc @C=4pc 
{
A\otimes I \ar[d]_{\ph\otimes 1_I} \ar[r]^{\rho_A} & A\ar[d]^{\ph}
\\
A'\otimes I\ar[r]^{\rho_{A'}} & A'
}
 \eeq
дадут нам
 $$
\xymatrix  @R=2pc @C=2.5pc 
{
I\overset{E}{\otimes}A\ar@{=}[r] & E(I\otimes A) \ar[d]_{1_I\overset{E}{\otimes}\ph=E(1_I\otimes \ph)} \ar[rr]^{\lambda_A^E=E(\lambda_A)} & & A\ar[d]^{\ph}
\\
I\overset{E}{\otimes}A'\ar@{=}[r] & E(I\otimes A')\ar[rr]^{\lambda_{A'}^E=E(\lambda_{A'})} & & A'
}
\qquad\qquad
\xymatrix  @R=2pc @C=2.5pc 
{
A\overset{E}{\otimes}I \ar@{=}[r] & E(A\otimes I)\ar[d]_{\ph\overset{E}{\otimes} 1_I=E(\ph\otimes 1_I)} \ar[rr]^{\rho_A^E=E(\rho_A)} & & A\ar[d]^{\ph}
\\
A'\overset{E}{\otimes}I \ar@{=}[r] & E(A'\otimes I)\ar[rr]^{\rho_{A'}^E=E(\rho_{A'})} & & A'
}
 $$
Кроме того, равенство $\lambda_I=\rho_I$ влечет равенство $\lambda_I^E=E(\lambda_I)=E(\rho_I)=\rho_I^E$, а диаграмма
$$
\xymatrix  @R=2pc @C=2.5pc 
{
(A\otimes I)\otimes B
 \ar[dr]_{\rho_A\otimes 1_B}
 \ar[rr]^{\alpha_{A,I,B}}
 & &
A\otimes (I\otimes B)
 \ar[dl]^{1_A\otimes\lambda_B}
\\
& A\otimes B &
}
$$
дает верхнее основание призмы
$$
\xymatrix  @R=2pc @C=5pc 
{
& E((A\otimes I)\otimes B)\ar[d]_{E(e_{A\otimes I}\otimes 1_B)}
 \ar[dr]^{E(\rho_A\otimes 1_B)}
 \ar[rr]^{E(\alpha_{A,I,B})}
 & &
E(A\otimes (I\otimes B))\ar[d]^{E(1_A\otimes e_{I\otimes B})}
 \ar[dl]_{E(1_A\otimes\lambda_B)}
 &
\\
& E(E(A\otimes I)\otimes B)\ar@{=}[d]\ar[r]_{E(\rho_A^E\otimes 1_B)}  & E(A\otimes B)\ar@{=}[dd] & E(A\otimes E(I\otimes B))\ar@{=}[d] \ar[l]^{E(1_A\otimes \lambda_B^E)} &
\\
& (A\overset{E}{\otimes}I)\overset{E}{\otimes}B\ar[dr]^{\rho_A^E\overset{E}{\otimes} 1_B}
\ar[rr]^(.7){\alpha^E_{A,I,B}}|!{[r]}{\hole} && A\overset{E}{\otimes}(I\overset{E}{\otimes}B)\ar[dl]_{1_A\overset{E}{\otimes}\lambda^E_B} &
\\
&& A\overset{E}{\otimes}B &&
}
$$
Коммутативность ее боковых граней можно считать очевидной, и вдобавок все вертикальные стрелки являются эпиморфизмами, поэтому остающееся нижнее основание также должно быть коммутативно.
\epr

\paragraph{Моноидальность функтора оболочки.}

\btm\label{TH:funktor-E-monoidalen}
Пусть $\Env_\varPhi^\varOmega$ -- регулярная оболочка, согласованная с тензорным произведением в $\tt K$.
Тогда функтор оболочки $E:{\tt K}\to {\tt L}$, построенный в теореме \ref{TH:regulyarnaya-obolochka},
является моноидальным.
\etm
\bpr
Чтобы быть моноидальным функтор $E:{\tt K}\to {\tt L}$ должен определять еще некий морфизм бифункторов
 $$
\Big((X,Y)\mapsto E(X)\overset{E}{\otimes} E(Y)\Big)\overset{E^{\otimes}}{\rightarrowtail} \Big((X,Y)\mapsto E(X\otimes Y)\Big)
 $$
в данном случае это будет семейство морфизмов
$$
E^\otimes_{X,Y}=E(e_X\otimes e_Y)^{-1}:E(X)\overset{E}{\otimes} E(Y)=E(E(X)\otimes E(Y))\to E(X\otimes Y),
$$
(по теореме \ref{PROP:E(e_X-otimes-e_Y)-in-Iso} $E(e_X\otimes e_Y)$ -- изоморфизмы,
поэтому они обладают обратными морфизмами $E(e_X\otimes e_Y)^{-1}$) и
морфизм $E^I$ в $\tt L$, переводящий единичный объект $I$ категории $\tt L$ в
образ $E(I)$ единичного объекта $I$ категории $\tt K$, и здесь это будет просто локальная единица:
$$
E^I=1_{I}:I\to I=\eqref{E(I)=I}=E(I).
$$
Проверим для этих составляющих аксиомы моноидального функтора.
Диаграмма согласованности с преобразованиями ассоциативности
\beq\label{DIAGR:assoc-v-L}
\xymatrix @R=2.5pc @C=6.0pc
{
E\big((X\otimes Y)\otimes Z\big)\ar[r]^{E(\alpha_{X,Y,Z})}& E\big(X\otimes (Y\otimes Z)\big)
\\
E(X\otimes Y)\overset{E}{\otimes}E(Z)\ar[u]^{E^{\otimes}_{X\otimes Y,Z}}
& E(X)\overset{E}{\otimes}E(Y\otimes Z)\ar[u]_{E^{\otimes}_{X,Y\otimes Z}}
 \\
\big(E(X)\overset{E}{\otimes}E(Y)\big)\overset{E}{\otimes}E(Z)
\ar[u]^{E^{\otimes}_{X,Y}\overset{E}{\otimes}1_{E(X)}}
\ar[r]^{\alpha^E_{E(X),E(Y),E(Z)}}& E(X)\overset{E}{\otimes}\big(E(Y)\overset{E}{\otimes}E(Z)\big)
\ar[u]_{1_{E(X)}\overset{E}{\otimes}E^{\otimes}_{Y,Z}}
}
\eeq
расшифровывается в данном случае так:
$$
\xymatrix @R=2.5pc @C=6.0pc
{
E\big((X\otimes Y)\otimes Z\big)\ar[d]_{E(e_{X\otimes Y}\otimes e_Z)}\ar[r]^{E(\alpha_{X,Y,Z})}& E\big(X\otimes (Y\otimes Z)\big)
\ar[d]^{E(e_X\otimes e_{Y\otimes Z})}
\\
E(E(X\otimes Y)\otimes E(Z))\ar[d]_{E(E(e_X\otimes e_Y)\otimes 1_{E(Z)})}
& E(E(X)\otimes E(Y\otimes Z))\ar[d]^{E(1_{E(X)}\otimes E(e_Y\otimes e_Z))}
 \\
E(E(E(X)\otimes E(Y))\otimes E(Z))
\ar[r]^{\alpha^E_{E(X),E(Y),E(Z)}}&
E(E(X)\otimes E(E(Y)\otimes E(Z)))
}
$$
Чтобы увидеть, что она коммутативна, представим ее как периметр следующей диаграммы:
\beq\label{DIAGR:funktor-E-monoidalen}
\xymatrix @R=3pc @C=4pc
{
&
\text{\tiny $E((X\otimes Y)\otimes Z)$}
\ar@/_3ex/[dl]_{\text{\tiny $E(e_{X\otimes Y}\otimes e_Z)$}}
\ar[d]!U|{\text{\tiny $E((e_X\otimes e_Y)\otimes e_Z)$}}
\ar[r]^{\text{\tiny $E(\alpha_{X,Y,Z})$}}&
\text{\tiny $E(X\otimes (Y\otimes Z))$}
\ar[d]!U|{\text{\tiny $E(e_X\otimes (e_Y\otimes e_Z))$}}
\ar@/^3ex/[dr]^{\text{\tiny $E(e_X\otimes e_{Y\otimes Z})$}}
&
\\
\text{\tiny $E(E(X\otimes Y)\otimes E(Z))$}
\ar@/_3ex/[dr]_(.3){\text{\tiny $E(E(e_X\otimes e_Y)\otimes 1_{E(Z)})$}\quad}
&
\text{\tiny $E((E(X)\otimes E(Y))\otimes E(Z))$}\ar[r]^{\text{\tiny $E(\alpha_{E(X),E(Y),E(Z)}$)}}
\ar[d]!U|{\text{\tiny $E(e_{E(X)\otimes E(Y)}\otimes 1_{E(Z)})$}}
&
\text{\tiny $E(E(X)\otimes (E(Y)\otimes E(Z)))$}
\ar[d]!U|{\text{\tiny $E(1_{E(X)}\otimes e_{E(Y)\otimes E(Z)})$}}
& \text{\tiny $E(E(X)\otimes E(Y\otimes Z))$}
\ar@/^3ex/[dl]^(.3){\quad\text{\tiny $E(1_{E(X)}\otimes E(e_Y\otimes e_Z))$}}
 \\
&
\text{\tiny $E(E(E(X)\otimes E(Y))\otimes E(Z))$}
\ar[r]_{\text{\tiny $\alpha^E_{E(X),E(Y),E(Z)}$}}&
\text{\tiny $E(E(X)\otimes E(E(Y)\otimes E(Z)))$} &
}
\eeq
Здесь левый внутренний четырехугольник удобно представить в виде
$$
\xymatrix @R=3pc @C=8pc
{
E((X\otimes Y)\otimes Z)
\ar[d]_{E(e_{X\otimes Y}\otimes e_Z)}
\ar[r]^{E((e_X\otimes e_Y)\otimes e_Z)}
&
E((E(X)\otimes E(Y))\otimes E(Z))
\ar[d]^{E(e_{E(X)\otimes E(Y)}\otimes 1_{E(Z)})}
 \\
E(E(X\otimes Y)\otimes E(Z))
\ar[r]_{E(E(e_X\otimes e_Y)\otimes 1_{E(Z)})}
&
E(E(E(X)\otimes E(Y))\otimes E(Z))
}
$$
Он получается примененим функтора $E$ к диаграмме
$$
\xymatrix @R=3pc @C=8pc
{
(X\otimes Y)\otimes Z
\ar[d]_{e_{X\otimes Y}\otimes e_Z}
\ar[r]^{(e_X\otimes e_Y)\otimes e_Z}
&
(E(X)\otimes E(Y))\otimes E(Z)
\ar[d]^{e_{E(X)\otimes E(Y)}\otimes 1_{E(Z)}}
 \\
E(X\otimes Y)\otimes E(Z)
\ar[r]_{E(e_X\otimes e_Y)\otimes 1_{E(Z)}}
&
E(E(X)\otimes E(Y))\otimes E(Z)
}
$$
которая в свою очередь получается умножением двух диаграмм
$$
\xymatrix @R=3pc @C=5pc
{
X\otimes Y
\ar[d]_{e_{X\otimes Y}}
\ar[r]^{(e_X\otimes e_Y)}
&
E(X)\otimes E(Y)
\ar[d]^{e_{E(X)\otimes E(Y)}}
 \\
E(X\otimes Y)
\ar[r]_{E(e_X\otimes e_Y)}
&
E(E(X)\otimes E(Y))
}
\qquad
\xymatrix @R=3pc @C=5pc
{
Z
\ar[d]_{e_Z}
\ar[r]^{e_Z}
&
E(Z)
\ar[d]^{1_{E(Z)}}
 \\
E(Z)
\ar[r]_{1_{E(Z)}}
&
E(Z)
}
$$
и из них правая тривиально коммутативна, а левая -- это транспонированная диаграмма \eqref{DIAGR:svyaz-otimes-s-env}.

Далее верхний внутренний четырехугольник в \eqref{DIAGR:funktor-E-monoidalen}
$$
\xymatrix @R=3pc @C=6pc
{
E((X\otimes Y)\otimes Z)
\ar[d]_{E((e_X\otimes e_Y)\otimes e_Z)}
\ar[r]^{E(\alpha_{X,Y,Z})}&
E(X\otimes (Y\otimes Z))
\ar[d]^{E(e_X\otimes (e_Y\otimes e_Z))}
\\
E((E(X)\otimes E(Y))\otimes E(Z))
\ar[r]^{E(\alpha_{E(X),E(Y),E(Z)})}
&
E(E(X)\otimes (E(Y)\otimes E(Z)))
}
$$
есть результат применения функтора $E$ к диаграмме
$$
\xymatrix @R=3pc @C=6pc
{
(X\otimes Y)\otimes Z
\ar[d]_{(e_X\otimes e_Y)\otimes e_Z}
\ar[r]^{\alpha_{X,Y,Z}}&
X\otimes (Y\otimes Z)
\ar[d]^{e_X\otimes (e_Y\otimes e_Z)}
\\
(E(X)\otimes E(Y))\otimes E(Z)
\ar[r]^{\alpha_{E(X),E(Y),E(Z)}}
&
E(X)\otimes (E(Y)\otimes E(Z))
}
$$
а это частный случай диаграммы \eqref{DIAGR:estestvennost-associativnosti}.

Затем, нижний внутренний четырехугольник в \eqref{DIAGR:funktor-E-monoidalen}
$$
\xymatrix @R=3pc @C=6pc
{
E((E(X)\otimes E(Y))\otimes E(Z))
\ar[r]^{E(\alpha_{E(X),E(Y),E(Z)})}
\ar[d]!U|{E(e_{E(X)\otimes E(Y)}\otimes 1_{E(Z)})}
&
E(E(X)\otimes (E(Y)\otimes E(Z)))
\ar[d]!U|{E(1_{E(X)}\otimes e_{E(Y)\otimes E(Z)})}
 \\
E(E(E(X)\otimes E(Y))\otimes E(Z))
\ar[r]_{\alpha^E_{E(X),E(Y),E(Z)}}&
E(E(X)\otimes E(E(Y)\otimes E(Z)))
}
$$
-- это частный случай диаграммы \eqref{DIAGR:opred-alpha-v-L}.

Наконец, правый внутренний четырехугольник в \eqref{DIAGR:funktor-E-monoidalen} полезно представить в виде
$$
\xymatrix @R=3pc @C=8pc
{
E(X\otimes (Y\otimes Z))
\ar[d]_{E(e_X\otimes (e_Y\otimes e_Z))}
\ar[r]^{E(e_X\otimes e_{Y\otimes Z})}
&
E(E(X)\otimes E(Y\otimes Z))
\ar[d]^{E(1_{E(X)}\otimes E(e_Y\otimes e_Z))}
\\
E(E(X)\otimes (E(Y)\otimes E(Z)))
\ar[r]_{E(1_{E(X)}\otimes e_{E(Y)\otimes E(Z)})}
&
E(E(X)\otimes E(E(Y)\otimes E(Z)))
}
$$
Он получается применением функтора $E$ к диаграмме
$$
\xymatrix @R=3pc @C=8pc
{
X\otimes (Y\otimes Z)
\ar[d]_{e_X\otimes (e_Y\otimes e_Z)}
\ar[r]^{e_X\otimes e_{Y\otimes Z}}
&
E(X)\otimes E(Y\otimes Z)
\ar[d]^{1_{E(X)}\otimes E(e_Y\otimes e_Z)}
\\
E(X)\otimes (E(Y)\otimes E(Z))
\ar[r]_{1_{E(X)}\otimes e_{E(Y)\otimes E(Z)}}
&
E(X)\otimes E(E(Y)\otimes E(Z))
}
$$
которая в свою очередь представляет собой результат умножения двух диаграмм
$$
\xymatrix @R=3pc @C=6pc
{
X
\ar[d]_{e_X}
\ar[r]^{e_X}
&
E(X)
\ar[d]^{1_{E(X)}}
\\
E(X)
\ar[r]_{1_{E(X)}}
&
E(X)
}
\qquad
\xymatrix @R=3pc @C=6pc
{
Y\otimes Z
\ar[d]_{e_Y\otimes e_Z}
\ar[r]^{e_{Y\otimes Z}}
&
E(Y\otimes Z)
\ar[d]^{E(e_Y\otimes e_Z)}
\\
E(Y)\otimes E(Z)
\ar[r]_{e_{E(Y)\otimes E(Z)}}
&
E(E(Y)\otimes E(Z))
}
$$
и из них левая тривиально коммутативна, а правая -- это видоизмененная диаграмма \eqref{DIAGR:svyaz-otimes-s-env}.

Помимо этого нужно убедиться в коммутативности диаграмм для левой и правой единицы:
$$
\xymatrix @R=2.5pc @C=4.0pc
{
E(I\otimes X)\ar[r]^{E(\lambda_X)} & E(X)
\\
I\overset{E}{\otimes} E(X)\ar[u]^{E^{\otimes}_{I,X}}
& I'\overset{E}{\otimes}E(X)\ar[u]_{\lambda'_{E(X)}}\ar[l]^{E^I\overset{E}{\otimes}1_{E(X)}}
}
\qquad
\xymatrix @R=2.5pc @C=4.0pc
{
E(X\otimes I)\ar[r]^{E(\rho_X)} & E(X)
\\
E(X)\overset{E}{\otimes}I \ar[u]^{E^{\otimes}_{X,I}}
 & E(X)\overset{E}{\otimes} I'\ar[u]_{\rho'_{E(X)}}\ar[l]^{1_{E(X)}\overset{E}{\otimes} E^I}
}
$$
В нашей ситуации они принимают следующий вид
$$
\xymatrix @R=2.5pc @C=6.0pc
{
E(I\otimes X)\ar[r]^{E(\lambda_X)}\ar[d]_{E(e_I\otimes e_X)} & E(X)
\\
E(I\otimes E(X))\ar[r]_{E(1_{E(X)}\otimes 1_{E(X)})} & E(I\otimes E(X))\ar[u]_{E(\lambda_{E(X)})}
}
\qquad
\xymatrix @R=2.5pc @C=6.0pc
{
E(X\otimes I)\ar[r]^{E(\rho_X)}\ar[d]_{E(e_X\otimes I)} & E(X)
\\
E(E(X)\otimes I)\ar[r]_{E(1_{E(X)}\otimes 1_{E(X)})} & E(E(X)\otimes I)\ar[u]_{E(\rho_{E(X)})}
}
$$
и это результат применения функтора $E$ к диаграммам \eqref{DIAGR:lambda,rho} с $X'=E(X)$ и $\ph=e_X$.
\epr

\bcor\label{COR:Env-sohranyaet-algebry}
Пусть $\Env_\varPhi^\varOmega$ -- регулярная оболочка, согласованная с тензорным произведением в $\tt K$.
Всякую алгебру (соответственно, коалгебру, биалгебру, алгебру Хопфа) $A$ в категории $\tt K$ операция $\Env_\varPhi^\varOmega$
превращает в алгебру  (соответственно, коалгебру, биалгебру, алгебру Хопфа) $\Env_\varPhi^\varOmega A$ в категории
$\tt L$.
\ecor
\bpr
Для случая алгебр и общих моноидальных функторов этот факт отмечается в \cite{Saavedra-Rivano}.
\epr

\subsection{Оболочки в категории стереотипных алгебр и двойственность относительно оболочек}

В теории стереотипных алгебр имеется несколько примеров оболочек, которые, собственно, и являются предметом изучения в этой книге. Вот их список.

\noindent\rule{160mm}{0.1pt}\begin{multicols}{2}

\bex\label{EX:Env_C-predv}
Непрерывной оболочкой $\Env_{\mathcal C}$ инволютивной стереотипной алгебры $A$ называется оболочка в классе плотных эпиморфизмов относительно класса морфизмов в $C^*$-алгебры.
\eex

\bex\label{EX:Env_E-predv}
Гладкой оболочкой $\Env_{\mathcal E}$ инволютивной стереотипной алгебры $A$ называется оболочка в классе плотных эпиморфизмов относительно класса дифференциальных гомоморфизмов в $C^*$-алгебры с присоединенными самосопряженными нильпотентными элементами.
\eex

\bex\label{EX:Env_O-predv}
Голоморфной оболочкой $\Env_{\mathcal O}$ стереотипной алгебры $A$ называется оболочка в классе плотных эпиморфизмов относительно класса морфизмов в банаховы алгебры.
\eex

\end{multicols}\noindent\rule[10pt]{160mm}{0.1pt}

Все эти оболочки регулярны, и при изучении каждой из них оказывается важным вопрос, как они действуют на алгебры Хопфа. Возникающия при этом конструкция едина для всех них, и ее удобно описать сразу как общекатегорную конструкцию.

\paragraph{Узловое разложение в $\SteAlg$.}
В моногорафии \cite{Akbarov-De-Gruyter-I} были доказаны следующие свойства категории $\SteAlg$ стереотипных алгебр (и по аналогии доказываются такие же свойства категории $\InvSteAlg$ инволютивных стереотипных алгебр).

\medskip
\centerline{\bf Свойства категорий $\SteAlg$ и $\InvSteAlg$:}

\bit{\it

\item[$1^\circ.$]\label{SteAlg-polna} Категории $\SteAlg$ и $\InvSteAlg$ полны (проективно и инъективно).

\item[$2^\circ.$]\label{SteAlg-uzl-razl} Категории $\SteAlg$ и $\InvSteAlg$ обладают узловыми разложениями.

\item[$3^\circ.$]\label{uzl-razl-v-Ste->v-SteAlg}  Для всякого морфизма $\ph:A\to B$ в категории $\SteAlg$ (или в категории $\InvSteAlg$) его узловое разложение $\ph=\im_\infty\ph\circ\red_\infty\ph\circ\coim_\infty\ph$ в категории ${\tt Ste}$ стереотипных пространств является разложением (необязательно, узловым) в категории $\SteAlg$ (соответственно, в категории $\InvSteAlg$).

\item[$4^\circ.$]\label{SteAlg-lok-mala-v-Mono} Категории $\SteAlg$ и $\InvSteAlg$ локально малы в подобъектах класса $\Mono$ (и, как следствие, в подобъектах класса $\SMono$).

\item[$5^\circ.$]\label{SteAlg-lok-mala-v-SEpi} Категории $\SteAlg$ и $\InvSteAlg$ локально малы в фактор-объектах класса $\SEpi$.

\item[$6^\circ.$]\label{SteAlg-SEpi-dissern-Mono} В категориях $\SteAlg$ и $\InvSteAlg$ строгие эпиморфизмы распознают мономорфизмы.

\item[$7^\circ.$]\label{SteAlg-SMono-dissern-Epi} В категориях $\SteAlg$ и $\InvSteAlg$ строгие мономорфизмы распознают эпиморфизмы.

}\eit

\bpr
В моногорафии \cite{Akbarov-De-Gruyter-I} это соответственно следующие результаты: Theorem 5.3.55, Theorem 5.3.64, Theorem 5.3.66, Theorem 5.3.61, Theorem 5.3.62, Theorem 5.3.57, Theorem 5.3.60.
\epr

\paragraph{Оболочки в категориях сильных стереотипных алгебр.}

В категориях $\SteAlg$ и $\InvSteAlg$ нас будут интересовать оболочки в классе так называемых плотных эпиморфизмов.

 \bit{
\item[$\bullet$]\label{DEF:DEpi} Условимся морфизм стереотипных (возможно, инволютивных) алгебр $\ph:A\to B$ называть {\it плотным}, если образ $\ph(A)$ алгебры $A$ при отображении $\ph$ плотен в алгебре $B$:
$$
\overline{\ph(A)}=B.
$$
Понятно, что плотные морфизмы являются эпиморфизмами, поэтому мы будем их также называть плотными эпиморфизмами. Класс всех плотных эпиморфизмов в категории $\InvSteAlg$ мы будем обозначать $\DEpi$. С классами $\Epi$ эпиморфизмов и $\SEpi$ строгих эпиморфизмов он связан вложениями
\beq\label{SEpi-subset-DEpi-subset-Epi}
\SEpi\subset\DEpi\subset\Epi.
\eeq
 }\eit

\btm\label{TH:SMono-circledcirc-DEpi=InvSteAlg}
Класс $\DEpi$ плотных эпиморфизмов мономорфно дополняем в категориях $\SteAlg$ и $\InvSteAlg$.
\etm
\bpr
Это следует из свойства $3^\circ$ на с.\pageref{uzl-razl-v-Ste->v-SteAlg}. Мономорфным дополнением к $\DEpi$ будет класс $\SMono_{\tt Ste}$ мономорфизмов в $\InvSteAlg$, являющихся строгими мономорфизмами в $\Ste$
\beq\label{SMono-circledcirc-DEpi=InvSteAlg}
\SMono_{\tt Ste}\circledcirc\DEpi=\InvSteAlg
\eeq
\epr

\bprop\label{TH:obolochki-i-nachinki-otn-plotnyh-epimorfizmov-v-Ste^circledast}
В категории $\SteAlg$ 
\bit{
\item[(i)] любая алгебра $A$ обладает оболочкой в классе $\DEpi$ всех плотных эпиморфизмов
относительно произвольного {\rm класса} морфизмов $\varPhi$,
выходящего из $A$:
$$
\exists \env_\varPhi^{\DEpi} A;
$$ 

\item[(ii)]  если  дополнительно к (i) $\varPhi$ различает морфизмы снаружи в $\SteAlg$,
то оболочка в классе $\DEpi$ является также оболочкой в классе $\DBim$ всех плотных биморфизмов:
    $$
    \env_\varPhi^{\DEpi} A=\env_\varPhi^{\DBim} A;
    $$

\item[(iii)] если дополнительно к (i) и (ii) $\varPhi$ является правым идеалом, то оболочка
$\env_\varPhi^{\DEpi\cap\Mono} A=\env_\varPhi^{\DBim} A$ также существует, и они совпадают:
$$
\env_\varPhi^{\DEpi} A=\env_\varPhi^{\DBim} A.
$$

    }\eit
\eprop
\bpr
Существование оболочки $\env_\varPhi^{\DEpi(\SteAlg)} A$ следует
из свойства $5^\circ$ на с.\pageref{5^0:obolochka-otn-klassa-morphizmov}.
Если $\varPhi$ различает морфизмы снаружи, то по теореме \ref{TH:Phi-razdel-moprfizmy} из существования
оболочки $\env_\varPhi^{\DEpi} A$ автоматически следует, что оболочка
$\env_\varPhi^{\DEpi\cap\Mono} A=\env_\varPhi^{\DBim} A$ также существует, и они совпадают:
$\env_\varPhi^{\DEpi} A=\env_\varPhi^{\DBim} A$. Если же вдобавок $\varPhi$ является правым идеалом, то
по теореме \ref{TH:Phi-razdel-moprfizmy-*} из существования
оболочки $\env_\varPhi^{\DEpi} A$ автоматически следует, что оболочка
$\env_\varPhi^{\DEpi\cap\Mono} A=\env_\varPhi^{\DBim} A$ также существует, и они совпадают:
$\env_\varPhi^{\DEpi} A=\env_\varPhi^{\DBim} A$.

\epr

Из теорем \ref{TH:sushestvovanie-seti-pri-faktorizatsii} и \ref{TH:regulyarnaya-obolochka} (с $\varOmega=\DEpi$) следует

\btm\label{TH:sushestvovanie-seti-v-Ste^circledast-dlya-DEpi}
Пусть в категории $\SteAlg$ задан класс морфизмов $\varPhi$,
выходящий из $\SteAlg$ и являющийся правым идеалом:
$$
\varPhi\circ\Mor(\SteAlg)\subseteq\varPhi.
$$
Тогда классы морфизмов $\DEpi$ и $\varPhi$ определяют в $\SteAlg$ полурегулярную оболочку
$\Env_\varPhi^{\DEpi}$, причем
для любого объекта $A$ в $\SteAlg$ морфизм оболочки описывается формулой
\beq\label{sushestvovanie-seti-v-Ste^circledast-dlya-DEpi-1}
\env_\varPhi^{\DEpi} A=\red_\infty \projlim{\mathcal N}^A\circ\coim_\infty \projlim{\mathcal N}^A,
\eeq
где ${\mathcal N}$ -- сеть эпиморфизмов, порожденная классами $\DEpi$ и $\varPhi$, а
$\red_\infty \projlim{\mathcal N}^A$ и $\coim_\infty \projlim{\mathcal N}^A$ -- элементы узлового разложения
\eqref{DEF:oboznacheniya-dlya-uzlov-razlozh} морфизма
$\projlim{\mathcal N}^A:A\to A_{\mathcal N}$ в категории ${\tt Ste}$ стереотипных пространств (не алгебр!).
Если вдобавок класс $\DEpi$ подталкивает класс $\varPhi$, то оболочка $\Env_\varPhi^{\DEpi}$ будет регулярна
(и поэтому ее можно определить как идемпотентный функтор).
 \etm

\paragraph{Алгебры Хопфа, рефлексивные относительно оболочки.}

Зафиксируем какую-нибудь регулярную оболочку $E$ в категории $\SteAlg$ стереотипных алгебр\footnote{Или в категории $\InvSteAlg$ инволютивных стереотипных алгебр. Нам будет удобно рассматривать этот случай параллельно, просто дополняя рассуждения подстрочными примечаниями, как например в подстрочном примечании \ref{FOOT:diagr-invol-v-reflex-alg-Hopfa} ниже. В качестве $E$ можно взять какую-нибудь из оболочек из примеров \ref{EX:Env_C-predv}-\ref{EX:Env_O-predv}.}. По замечанию \ref{REM:env-polureg-obolochki}, для всякой алгебры $A$ морфизм оболочки $e_A:A\to E(A)$ является эпиморфизмом стереотипных алгебр.
Более того, мы будем предполагать, что наша оболочка $E$ удовлетворяет условию

\bit{

\item[SE:]\label{uslovie-SE} для всякой алгебры $A$ морфизм оболочки $e_A:A\to E(A)$ является эпиморфизмом стереотипных пространств.

}\eit

\brem
Все оболочки из примеров \ref{EX:Env_C-predv}---\ref{EX:Env_O-predv} удовлетворяют условию SE. 
\erem

Пусть $H$ --- какая-нибудь алгебра Хопфа в категории $({\tt Ste},\circledast)$ (то есть алгебра Хопфа относительно тензорного произведения $\circledast$). По определению операции $E$, объект $E(H)$ является алгеброй в категории $({\tt Ste},\circledast)$ (то есть алгеброй относительно тензорного произведения $\circledast$), а морфизм $ e_H:H\to E(H)$ представляет собой морфизм в категории стереотипных алгебр $\SteAlg$ (то есть линейное непрерывное мультипликативное и сохраняющее единицу отображение).

Структура $\circledast$-алгебры на $E(H)$ наследуется из структуры алгебры в $H$. Однако при выбранном нами определении $E$, структура коалгебры на $H$ формально не определяет никакую структуру $\circledast$-коалгебры на $E(H)$. Точно так же антипод $\sigma$ в $H$ формально не обязан превращаться в какой-то морфизм в $E(H)$ (потому что функтор $E$ определен на категории стереотипных алгебр $\SteAlg$, а не на категории стереотипных пространств ${\tt Ste}$). Поэтому переход $H\mapsto E(H)$ не является операцией в категории $\circledast$-алгебр Хопфа.

Тем не менее, есть важный класс примеров, в которых операция $H\mapsto E(H)$ все-таки превращает алгебры Хопфа в алгебры Хопфа, однако с одним неожиданным уточнением: {\it алгебра Хопфа $H$ над тензорным произведением $\circledast$ при этом превращается в алгебру Хопфа $E(H)$ над другим тензорным произведением, $\odot$}. Такими алгебрами Хопфа являются, в частности, групповые алгебры ${\mathcal O}^\star(G)$ для некоторых широких классов комплексных групп Ли $G$ (как показывает пример \ref{EX:heartsuit=widehat-dlya-komp-porozhd-grupp} выше). Более того, это вообще оказывается типичной ситуацией в случаях, когда в категории $\SteAlg$ рассматривается какая-то оболочка: как было показано в \cite{Ak17-1,Ak17-2}, выделение такого класса алгебр Хопфа естественно ведет к обобщениям понтрягинской двойственности на различные классы некоммутативных групп (а в \cite{Ak17-1,Ak17-2} выделяемые таким образом алгебры Хопфа назывались {\it непрерывно рефлексивными} и {\it гладко рефлексивными}).

Интересной деталью в этих примерах является тот факт, что в них {\it морфизм оболочки $ e_H:H\to E(H)$ можно интерпретировать как гомоморфизм алгебр Хопфа}. Из-за того, что $H$ и $E(H)$ являются алгебрами Хопфа в разных моноидальных категориях (первая в категории $({\tt Ste},\circledast)$, а вторая в категории $({\tt Ste},\odot)$), понятие гомоморфизма между ними не определено. Тем не менее, термину ``гомоморфизм $\circledast$-алгебры Хопфа в $\odot$-алгебру Хопфа'' можно придать точный смысл из-за того, что бифункторы $\circledast$ и $\odot$ в категории $\Ste$ связаны естественным преобразованием, известным как преобразование Гротендика \cite[7.4]{Ak03}: каждой паре стереотипных пространств $(X,Y)$ можно поставить в соответствие морфизм стереотипных пространств $@_{X,Y}:X\circledast Y\to X\odot Y$, называемый {\it преобразованием Гротендика для пары $(X,Y)$}, так, что для любых морфизмов $\ph:X\to X'$ и $\chi:Y\to Y'$ будет коммутативна диаграмма
\beq\label{propbr-Grothendieck}
 \xymatrix @R=3.pc @C=4.pc
{
X\circledast Y\ar[r]^{@_{X,Y}}\ar[d]_{\ph\circledast\chi} &
X\odot Y\ar[d]^{\ph\odot \chi} \\
X'\circledast Y'\ar[r]^{@_{X',Y'}} &
X'\odot Y'.
}
\eeq

Как следствие, для любых двух морфизмов стереотипных пространств $\ph:X\to X'$ и $\chi:Y\to Y'$ естественным образом определен морфизм из тензорного произведения $X\circledast Y$ в тензорное произведение $X'\odot Y'$: таким морфизмом будет диагональ диаграммы \eqref{propbr-Grothendieck}:
\beq\label{propbr-Grothendieck-1}
 \xymatrix @R=3.pc @C=4.pc
{
X\circledast Y\ar[r]^{@_{X,Y}}\ar@{-->}[dr]\ar[d]_{\ph\circledast\chi} &
X\odot Y\ar[d]^{\ph\odot \chi} \\
X'\circledast Y'\ar[r]^{@_{X',Y'}} &
X'\odot Y'
}
\eeq

Из этого, в свою очередь, следует, что если $A$ -- алгебра в категории $(\Ste,\circledast)$, а $B$ -- алгебра в категории $(\Ste,\odot)$, то гомоморфизм между ними можно определить как морфизм стереотипненых пространств $\ph:A\to B$, для которого будут коммутативны диаграммы
$$
 \xymatrix @R=2.pc @C=2.pc
{
& A\odot A\ar[dr]^{\quad  \ph \odot  \ph } & \\
A\circledast A\ar[ur]^{@}\ar[dr]^{\quad  \ph \circledast  \ph }\ar[dd]_{\mu_A} & & B\odot B\ar[dd]_{\mu_B} \\
& B\circledast B\ar[ur]^{@} & \\
A\ar[rr]^{ \ph } && B
}\qquad
\xymatrix @R=2.pc @C=2.pc
{
A\ar[rr]^{ \ph } & & B \\
& \C\ar[ul]^{\iota_A}\ar[ur]_{\iota_B} &
}
$$
(здесь $\mu_A$ и $\mu_B$ --- умножения, а $\iota_A$ и $\iota_B$ --- единицы в $A$ и $B$).

Аналогично, если $A$ -- коалгебра в категории $(\Ste,\circledast)$, а $B$ -- коалгебра в категории $(\Ste,\odot)$, то гомоморфизм между ними можно определить как морфизм стереотипненых пространств $\ph:A\to B$, для которого будут коммутативны диаграммы
$$
 \xymatrix @R=2.pc @C=2.pc
{
& A\odot A\ar[dr]^{\quad  \ph \odot  \ph } & \\
A\circledast A\ar[ur]^{@}\ar[dr]^{\quad  \ph \circledast  \ph } & & B\odot B \\
& B\circledast B\ar[ur]^{@} & \\
A\ar[rr]^{ \ph }\ar[uu]^{\varkappa_A} && B\ar[uu]^{\varkappa_B}
}\qquad
\xymatrix @R=2.pc @C=2.pc
{
A\ar[rr]^{ \ph }\ar[dr]_{\e_A} & & B\ar[dl]^{\e_B} \\
& \C &
}
$$
(здесь $\varkappa_A$ и $\varkappa_B$ --- коумножения, а $\e_A$ и $\e_B$ --- коединицы в $A$ и $B$).

Эти предварительные замечания оправдывают следующее определение.

Условимся говорить, что стереотипная алгебра Хопфа $H$  в категории $({\tt Ste},\circledast)$  {\it рефлексивна относительно регулярной оболочки $E$}\label{DEF:algebra-Hopfa-reflex-otn-obolochki}, если на оболочке $E(H)$ определена структура алгебры Хопфа в категории $({\tt Ste},\odot)$ так, что выполняются следующие четыре условия:\footnote{В этом определении диаграмма \eqref{DIAG:reflex-otn-obolochki-1} с левой диаграммой в \eqref{DIAG:reflex-otn-obolochki-3} означают, что оператор $ e_H$ является гомоморфизмом $\circledast$-алгебры $H$ в $\odot$-алгебру $E(H)$, а диаграмма \eqref{DIAG:reflex-otn-obolochki-2} с правой диаграммой в \eqref{DIAG:reflex-otn-obolochki-3} означают, что оператор $ e_H$ является гомоморфизмом $\circledast$-коалгебры $H$ в $\odot$-коалгебру $E(H)$.}
\bit{
\item[(i)] {\it оболочка является гомоморфизмом алгебр Хопфа}: морфизм оболочки $ e_H:H\to E(H)$ является гомоморфизмом алгебр Хопфа в том смысле, что коммутативны следующие диаграммы:\footnote{\label{FOOT:diagr-invol-v-reflex-alg-Hopfa} В случае, когда рассматриваемые алгебры считаются инволютивными, к этому списку еще добавляется диаграмма $$\xymatrix @R=3.pc @C=4.pc
{
H\ar[r]^{e_H}\ar[d]_{\bullet} & E(H)\ar[d]^{E(\bullet)} \\
H\ar[r]^{e_H} & E(H)
}
$$
в которой $\bullet$ и $E(\bullet)$ --- инволюции.}
\beq\label{DIAG:reflex-otn-obolochki-1}
 \xymatrix @R=2.pc @C=2.pc
{
& H\odot H\ar[dr]^{\quad   e_H \odot   e_H } & \\
H\circledast H\ar[ur]^{@}\ar[dr]^{\quad   e_H \circledast   e_H }\ar[dd]_{\mu} & & E(H)\odot E(H)\ar[dd]_{E(\mu)} \\
& E(H)\circledast E(H)\ar[ur]^{@} & \\
H\ar[rr]^{  e_H } && E(H)
}
\eeq
\beq\label{DIAG:reflex-otn-obolochki-2}
 \xymatrix @R=2.pc @C=2.pc
{
& H\odot H\ar[dr]^{\quad  e_H \odot   e_H } & \\
H\circledast H\ar[ur]^{@}\ar[dr]^{\quad   e_H \circledast   e_H } & & E(H) \odot E(H) \\
& E(H)\circledast E(H)\ar[ur]^{@} & \\
H\ar[rr]^{  e_H }\ar[uu]^{\varkappa} && E(H) \ar[uu]^{E(\varkappa)}
}
\eeq
\beq\label{DIAG:reflex-otn-obolochki-3}
 \xymatrix @R=2.pc @C=2.pc
{
H\ar[rr]^{  e_H } & & E(H) \\
& \C\ar[ul]^{\iota}\ar[ur]_{E(\iota)} &
}\qquad
 \xymatrix @R=2.pc @C=2.pc
{
H\ar[rr]^{  e_H }\ar[dr]_{\e} & & E(H) \ar[dl]^{E(\e)} \\
& \C &
}
\eeq
\beq\label{DIAG:reflex-otn-obolochki-4}
 \xymatrix @R=3.pc @C=4.pc
{
H\ar[r]^{  e_H }\ar[d]_{\sigma} & E(H) \ar[d]^{E(\sigma)} \\
H\ar[r]^{  e_H } & E(H)
}
\eeq
-- здесь $@$ -- преобразование Гротендика из \eqref{propbr-Grothendieck}, $\mu$, $\iota$, $\varkappa$, $\e$, $\sigma$ -- структурные морфизмы (умножение, единица, коумножение, коединица, антипод) в $H$, а $E(\mu)$, $E(\iota)$, $E(\varkappa)$, $E(\e)$, $E(\sigma)$ -- структурные морфизмы в $E(H) $.

\item[(ii)]\label{(env-H)^star:H^star-gets-(Env-H)^star} {\it сопряженное отображение к оболочке является оболочкой}: отображение $(  e_H )^\star:H^\star\gets E(H)^\star$, сопряженное к морфизму оболочки $  e_H :H\to E(H) $, является оболочкой в том же смысле:
\beq\label{(env-H)^star:H^star-gets-(Env-H)^star-0}
(  e_H )^\star= e_{E(H)^\star}
\eeq

\item[(iii)]\label{@:H-circledast-H->H-odot-H-Mono} {\it мономорфность преобразования Гротендика исходной алгебры}: преобразование Гротендика $@:H\circledast H\to H\odot H$ является мономорфизмом стереотипных пространств (то есть инъективным отображением),

\item[(iv)]\label{@:E(H)-circledast-E(H)->E(H)-odot-E(H)-Epi} {\it эпиморфность преобразования Гротендика оболочки}: преобразование Гротендика $@:E(H)\circledast E(H)\to E(H)\odot E(H)$ является эпиморфизмом стереотипных пространств (то есть его образ плотен в области значений),

}\eit

\brem\label{REM:DEF-reflex-Hopf-1} Условия (iii) и (iv) в этом списке требуют комментария. 
\bit{
\item[1)] Нужно сказать, что (iii) и (iv) двойственны друг другу в том смысле, что если требовать выполнение условия (iii) --- мономорфности морфизма $@:H\circledast H\to H\odot H$ --- для алгебры $H$, то для симметрии нужно требовать такого же условия  --- мономорфности морфизма $@:E(H)^\star\circledast E(H)^\star\to E(H)^\star\odot E(H)^\star$ ---  для двойственной алгебры $H^\dagger=E(H)^\star$ (в смысле определения \eqref{DEF:H^dagger:=E(H)^star}). А это то же самое, что условие (iv) --- эпиморфность морфизма $@:E(H)\circledast E(H)\to E(H)\odot E(H)$. 

\item[2)] Эти условия вводятся потому, что (iii) используется ниже в теореме \ref{TH:H-cong-(H^*)^*}, а (iv) появляется в  теореме \ref{TH:H-holom-refl=>H^+-holom-refl} как необходимавя посылка для симметрии с (iii), о которой мы уже сказали. 

\item[3)] Полезно добавить, что условия (iii) и (iv) будут выполнены автоматически, если потребовать более сильного условия 
\bit{
\item[(v)]\label{holom-reflex=>approximation} пространства $H$ и $E(H)$ обладают стереотипной аппроксимацией.
}\eit
Это условие, однако, существенно труднее проверить (и именно поэтому мы ограничиваемся требованиями (iii) и (iv)).
}\eit
\erem

Класс всех рефлексивных алгебр Хопфа относительно оболочки $E$ желательно как-нибудь обозначить, мы для этого выберем символ $\Hopf^E$. Этот класс образует категорию, в которой морфизмами являются обычные морфизмы алгебр Хопфа в категории $(\Ste,\circledast)$.

\btm\label{TH:Hopf-na-H^heartsuit}
Из условия (iv) автоматически следует, что на $E(H)$ существует не более одной структуры алгебры Хопфа в $({\tt Ste},\odot)$, для которой выполняются условия (i) и (ii).
\etm
\bpr
Здесь используется изначально обговоренное условие SE на с.\pageref{uslovie-SE}: отображение $ e_H: H\to  E(H)$ должно быть эпиморфизмом стереотипных пространств. Отсюда следует, что отображение $ @\circ e_H \circledast  e_H :H\circledast H\to  E(H) \odot E(H) $ тоже является  эпиморфизмом, как композиция двух эпиморфизмов
$ e_H \circledast e_H :H\circledast H\to  E(H) \circledast E(H) $ и $@:E(H) \circledast E(H)  \to  E(H) \odot E(H)$.
\epr

Условия (i) и (ii) удобно изображать в виде диаграммы
 \beq\label{diagramma-golomorfnoi-refleksivnosti}
 \xymatrix @R=1.pc @C=1.pc
 {
 H
 & \ar@{|->}[r]^{E} & &
 E(H)
 \\
 & & &
 \ar@{|->}[d]^{\star}
 \\
 \ar@{|->}[u]^{\star}
 & & &
 \\
 H^\star
 & &
 \ar@{|->}[l]_{E}
 &
 E(H)^\star
 }
 \eeq
которую мы называем {\it диаграммой рефлексивности относительно оболочки $E$}, и в которую вкладываем следующий смысл:
 \bit{
\item[1)] в углах квадрата стоят алгебры Хопфа, причем $H$ -- алгебра Хопфа в $({\tt Ste},\circledast)$, затем следует алгебра Хопфа $ E(H) $ в $({\tt Ste},\odot)$, и далее категории $({\tt Ste},\circledast)$ и $({\tt Ste},\odot)$ чередуются,

\item[2)] чередование операций $E$ и $\star$ (с какого места ни начинай) на четвертом шаге возвращает к исходной алгебре Хопфа (конечно, с точностью до изоморфизма функторов).
 }\eit

Чтобы объяснить смысл термина ``рефлексивность'', обозначим однократное последовательное применение операций $E$ и $\star$ каким-нибудь символом, например $\dagger$:
\beq\label{DEF:H^dagger:=E(H)^star}
H^\dagger:=E(H)^\star
\eeq
Поскольку на $E(H) $ определена структура алгебры Хопфа относительно $\odot$, на сопряженном пространстве $H^\dagger=E(H)^\star$ определена структура алгебры Хопфа относительно $\circledast$.

Тогда диаграмму \eqref{diagramma-golomorfnoi-refleksivnosti} можно понимать как утверждение, что $H$ будет изоморфна своей второй двойственной в смысле операции $\dagger$ алгебре Хопфа:
 \beq\label{H-cong-(H^*)^*}
H\cong (H^\dagger)^\dagger.
 \eeq
Замыкание пути в диаграмме \eqref{diagramma-golomorfnoi-refleksivnosti} интерпретируется как изоморфизм стереотипных пространств, однако теорема \ref{TH:H-cong-(H^*)^*} ниже показывает, что этот изоморфизм автоматически будет изоморфизмом алгебр Хопфа в $(\Ste,\circledast)$, и более того, такое семейство изоморфизмов $I_H:H\to  (H^\dagger)^\dagger$ будет естественным преобразованием тождественного функтора $H\mapsto H$ в категории $\Hopf^E$ в функтор $H\mapsto (H^\dagger)^\dagger$ (и, поскольку каждое такое преобразование есть изоморфизм в $\Hopf^E$, это будет изоморфизм функторов).

Однако, прежде чем это обсуждать, заметим вот что:

\btm\label{TH:H-holom-refl=>H^+-holom-refl}
Пусть $E$ --- регулярная оболочка в категории $\SteAlg$. Тогда
\bit{
\item[(i)] на классе $\Hopf^E$ алгебр Хопфа, рефлексивных относительно оболочки $E$, операция $H\mapsto H^\dagger$ может быть определена как отображение, и

\item[(ii)] это отображение $H\mapsto H^\dagger$ переводит класс $\Hopf^E$ в себя.
}\eit
\etm
\bpr
Операцию $H\mapsto H^\dagger$ можно определить как отображение, потому что каждую из операций, входящих в ее определение, $H\mapsto E(H)$ и $X\mapsto X^\star$ можно определить как отображение. Поэтому здесь важно доказать утверждение  (ii).

Применив функтор $\star$ к диаграммам \eqref{DIAG:reflex-otn-obolochki-1}-\eqref{DIAG:reflex-otn-obolochki-4}, мы получим диаграммы
\beq\label{DIAG:reflex-otn-obolochki-1*}
 \xymatrix @R=2.pc @C=2.pc
{
& H^\star\odot H^\star\ar[dl]_{@} & \\
H^\star\circledast H^\star & & E(H)^{\star}\odot E(H)^{\star}\ar[ul]_{\quad  ( e_H)^\star \odot  ( e_H)^\star }\ar[dl]_{@}  \\
& E(H)^{\star}\circledast E(H)^{\star}\ar[ul]_{\quad  ( e_H)^\star \circledast  ( e_H)^\star } & \\
H^\star\ar[uu]^{\mu^\star}  && E(H)^{\star}\ar[uu]^{E(\mu)^\star}\ar[ll]_{ ( e_H)^\star }
}
\eeq
\beq\label{DIAG:reflex-otn-obolochki-2*}
 \xymatrix @R=2.pc @C=2.pc
{
& H^\star\odot H^\star\ar[dl]_{@} & \\
H^\star\circledast H^\star\ar[dd]_{\varkappa^\star} & & E(H)^{\star} \odot E(H)^{\star}\ar[ul]_{\quad ( e_H)^\star \odot  ( e_H)^\star }
\ar[dl]_{@} \ar[dd]_{E(\varkappa)^\star} \\
& E(H)^{\star}\circledast E(H)^{\star}\ar[ul]_{\quad  ( e_H)^\star \circledast  ( e_H)^\star }  & \\
H^\star && E(H)^{\star} \ar[ll]_{ ( e_H)^\star }
}
\eeq
\beq\label{DIAG:reflex-otn-obolochki-3*}
 \xymatrix @R=2.pc @C=2.pc
{
H^\star\ar[dr]_{\iota^\star} & & E(H)^{\star}\ar[ll]_{ ( e_H)^\star }
\ar[dl]^{E(\iota)^\star} \\
& \C &
}\qquad
 \xymatrix @R=2.pc @C=2.pc
{
H^\star & & E(H)^{\star}\ar[ll]_{( e_H)^\star} \\
& \C\ar[ul]_{\e^\star}\ar[ur]_{E(\e)^\star} &
}
\eeq
\beq\label{DIAG:reflex-otn-obolochki-4*}
 \xymatrix @R=3.pc @C=4.pc
{
H^\star & E(H)^{\star}\ar[l]_{ ( e_H)^\star } \\
H^\star\ar[u]^{\sigma^\star} & E(H)^{\star}\ar[l]_{ ( e_H)^\star }\ar[u]_{E(\sigma)^\star}
}
\eeq

Заметим после этого, что согласно \eqref{(env-H)^star:H^star-gets-(Env-H)^star-0},
отображение $(  e_H )^\star:H^\star\gets E(H)^{\star}$ является оболочкой алгебры $E(H)^{\star}$:
\beq\label{H^star=E(E(H)^star)}
H^\star=E(E(H)^\star)
\eeq
Отсюда мы можем сделать вывод, что оболочка $E(E(H)^\star)$, как ее ни выбирай (она ведь определяется неоднозначно, а только с точностью до изоморфизма $\circledast$-алгебр), обладает структурой $\odot$-алгебры Хопфа, у которой  структурные морфизмы $E(E(\mu)^\star)$ (коумножение), $E(E(\varkappa)^\star)$ (умножение), $E(E(\iota)^\star)$ (коединица), $E(E(\e)^\star)$ (единица), $E(E(\sigma)^\star)$ (антипод) делают коммутативными следующие диаграммы:
\beq\label{DIAG:reflex-otn-obolochki-1**}
 \xymatrix @R=2.pc @C=2.pc
{
& E(E(H)^\star)\odot E(E(H)^\star)\ar[dl]_{@} & \\
E(E(H)^\star)\circledast E(E(H)^\star) & & E(H)^{\star}\odot E(H)^{\star}\ar[ul]_{\quad   e_{E(H)^{\star}} \odot   e_{E(H)^{\star}} }\ar[dl]_{@}  \\
& E(H)^{\star}\circledast E(H)^{\star}\ar[ul]_{\quad   e_{E(H)^{\star}} \circledast   e_{E(H)^{\star}} } & \\
E(E(H)^\star)\ar[uu]^{E(E(\varkappa)^\star)}  && E(H)^{\star}\ar[uu]^{E(\mu)^\star}\ar[ll]_{  e_{E(H)^{\star}} }
}
\eeq
\beq\label{DIAG:reflex-otn-obolochki-2**}
 \xymatrix @R=2.pc @C=2.pc
{
& E(E(H)^\star)\odot E(E(H)^\star)\ar[dl]_{@} & \\
E(E(H)^\star)\circledast E(E(H)^\star)\ar[dd]_{E(E(\mu)^\star)} & & E(H)^{\star} \odot E(H)^{\star}\ar[ul]_{\quad  e_{E(H)^{\star}} \odot   e_{E(H)^{\star}} }
\ar[dl]_{@} \ar[dd]_{E(\varkappa)^\star} \\
& E(H)^{\star}\circledast E(H)^{\star}\ar[ul]_{\quad   e_{E(H)^{\star}} \circledast   e_{E(H)^{\star}} }  & \\
E(E(H)^\star) && E(H)^{\star} \ar[ll]_{  e_{E(H)^{\star}} }
}
\eeq
\beq\label{DIAG:reflex-otn-obolochki-3**}
 \xymatrix @R=2.pc @C=2.pc
{
E(E(H)^\star)\ar[dr]_{E(E(\iota)^\star)} & & E(H)^{\star}\ar[ll]_{  e_{E(H)^{\star}} }
\ar[dl]^{E(\iota)^\star} \\
& \C &
}\qquad
 \xymatrix @R=2.pc @C=2.pc
{
E(E(H)^\star) & & E(H)^{\star}\ar[ll]_{ e_{E(H)^{\star}}} \\
& \C\ar[ul]^{E(E(\e)^\star)}\ar[ur]_{E(\e)^\star} &
}
\eeq
\beq\label{DIAG:reflex-otn-obolochki-4**}
 \xymatrix @R=3.pc @C=4.pc
{
E(E(H)^\star) & E(H)^{\star}\ar[l]_{  e_{E(H)^{\star}} } \\
E(E(H)^\star)\ar[u]^{E(E(\sigma)^\star)} & E(H)^{\star}\ar[l]_{ e_{E(H)^{\star}} }\ar[u]_{E(\sigma)^\star}
}
\eeq
(это так, потому что конкретная оболочка $H^\star$ алгебры $E(H)^{\star}$ обладает этими свойствами).

Но это в точности диаграммы \eqref{DIAG:reflex-otn-obolochki-1}-\eqref{DIAG:reflex-otn-obolochki-4}, только с подставленным $E(H)^{\star}=H^\dagger$ вместо $H$. Мы поняли, что $E(H)^{\star}=H^\dagger$ обладает свойством (i) на с.\pageref{DIAG:reflex-otn-obolochki-1}.

Далее, рассмотрим отображение $( e_H)^\star:E(H)^{\star}\to H^\star$. Согласно \eqref{(env-H)^star:H^star-gets-(Env-H)^star-0}, оно является оболочкой алгебры $E(H)^{\star}$. Его сопряженное отображение $( e_H)^{\star\star}:H^{\star\star}\to E(H)^{\star\star}$ можно вписать в коммутативную диаграмму
$$
 \xymatrix @R=3.pc @C=4.pc
{
H^{\star\star}\ar[r]^{( e_H)^{\star\star}} &  E(H)^{\star\star} \\
H\ar[r]^{ e_H}\ar[u]^{i_H} & E(H) \ar[u]^{i_{E(H)}}
}
$$
в которой $i_H$ и $i_{E(H)}$ --- изоморфизмы $\circledast$-алгебр, а $ e_H$ --- оболочка алгебры $H$. Отсюда можно заключить, что $( e_H)^{\star\star}$ --- оболочка алгебры $H^{\star\star}$.

Мы получаем, что сопряженное отображение к оболочке алгебры $E(H)^{\star}$ также является оболочкой, то есть \eqref{(env-H)^star:H^star-gets-(Env-H)^star-0} выполняется для $E(H)^{\star}=H^\dagger$ подставленного вместо $H$.

Наконец, условие (iii) и (iv) на с.\pageref{@:H-circledast-H->H-odot-H-Mono} для алгебры $H^\dagger=E(H)^\star$ также выполнены, потому что 
\bit{

\item[---] из эпиморфности отображения $@:E(H)\circledast E(H)\to E(H)\odot E(H)$ следует мономорфность двойственного отображения $@:E(H)^\star\circledast E(H)^\star\to E(H)^\star\odot E(H)^\star$, то есть отображения $@:H^\dagger\circledast H^\dagger\to H^\dagger\odot H^\dagger$, а 

\item[---] из мономорфности отображения $@:H\circledast H\to H\odot H$ следует эпиморфность двойственного отображения $@:H^\star\circledast H^\star\to H^\star\odot H^\star$, то есть в силу \eqref{H^star=E(E(H)^star)}, отображения $@:H^\dagger\circledast H^\dagger\to H^\dagger\odot H^\dagger$.

}\eit
\epr

\paragraph{Функториальность перехода к алгебре Хопфа, двойственной относительно оболочки.}

Теперь покажем, что операция $H\mapsto H^\dagger$, описанная в теореме \ref{TH:H-holom-refl=>H^+-holom-refl}, является функтором и более того, двойственностью на категории $\Hopf^E$: 

\btm\label{TH:H-cong-(H^*)^*}
Пусть $E$ --- регулярная оболочка в категории $\SteAlg$. Тогда
\bit{
\item[(i)] на категории $\Hopf^E$ алгебр Хопфа, рефлексивных относительно оболочки $E$, операция $H\mapsto H^\dagger$ может быть определена как контравариантный функтор, и

\item[(ii)] этот функтор $H\mapsto H^\dagger$ является двойственностью на $\Hopf^E$ (то есть его квадрат $H\mapsto (H^\dagger)^\dagger$  изоморфен тождественному функтору $H\mapsto H$).
}\eit
\etm
\bpr
Нам нужно доказать тождество функторов
\beq\label{H-cong-(H^*)^*-0}
H\cong E(E(H)^\star)^\star
\eeq
Если добавить еще одну звездочку, мы получим
$$
H^\star\cong E(E(H)^\star)^{\star\star}
$$
Две подряд идущие звездочки справа можно убрать, потому что функтор $\star\star$ изоморфен тождественному функтору
\beq\label{H-cong-(H^*)^*-1}
X\cong X^{\star\star}
\eeq
(это изоморфизм функторов не только в $\Ste$, но и в $\Hopf^E$), и тогда мы получим тождество
\beq\label{H-cong-(H^*)^*-2}
H^\star\cong E(E(H)^\star)
\eeq
Теперь заметим, что если нам удастся доказать тождество \eqref{H-cong-(H^*)^*-2} как изоморфизм функторов, действующих на $\Hopf^E$, то тождество \eqref{H-cong-(H^*)^*-0} тоже превратится в изоморфизм функторов, потому что
$$
H\cong\eqref{H-cong-(H^*)^*-1}\cong H^{\star\star}\cong\eqref{H-cong-(H^*)^*-2}\cong E(E(H)^\star)^\star.
$$
После того, как мы это поняли, мы сосредоточимся на доказательстве  \eqref{H-cong-(H^*)^*-2}.

1. Сначала определим систему морфизмов
\beq\label{H-cong-(H^*)^*-3}
\alpha_H:H^\star\to E(E(H)^\star),\qquad H\in {\Hopf^E},
\eeq
про которую мы дальше будем доказывать, что это изоморфизм функторов $\star$ и $E\circ\star\circ E$. Здесь используется свойство \eqref{(env-H)^star:H^star-gets-(Env-H)^star-0}: по нему отображение $(  e_H )^\star:E(H)^\star\to H^\star$, сопряженное к морфизму оболочки $ e_H :H\to E(H)$, также является оболочкой. В частности, оно является расширением, поэтому существует морфизм \eqref{H-cong-(H^*)^*-3}, замыкающий диаграмму
\beq\label{H-cong-(H^*)^*-4}
 \xymatrix @R=2.pc @C=1.5pc
{
& E(H)^{\star} \ar[dl]_{( e_H)^\star}\ar[dr]^{ e_{E(H)^{\star}}} &  \\
H^\star\ar@{-->}[rr]_{\alpha_H} & & E(E(H)^\star)
}
\eeq
Из того, что обе наклонные стрелки в этом треугольнике являются оболочками, следует, что $\alpha_H$ является изоморфизмом в категории $\Ste^\circledast$ (стереотипных алгебр относительно тензорного произведения $\circledast$). Наша первая задача --- убедиться, что $\alpha_H$ --- изоморфизм не только в категории $\Ste^\circledast$, но и в категории стереотипных алгебр Хопфа относительно тензорного произведения $\odot$.

2. Покажем, что $\alpha_H$ сохраняет умножение, то есть, что коммутативна диаграмма
\beq\label{H-cong-(H^*)^*-4-1}
 \xymatrix @R=2.pc @C=6pc
{
H^\star\odot H^\star\ar[r]^{\alpha_H\odot\alpha_H}\ar[d]_{\varkappa^\star} &  E(E(H)^\star)\odot E(E(H)^\star)
\ar[d]^{E(E(\varkappa)^\star)} \\
H^\star\ar[r]_{\alpha_H} & E(E(H)^\star)
}
\eeq
(мы обозначаем буквой $\varkappa$ коумножение в $H$, а $\varkappa^\star$ и $E(E(\varkappa)^\star)$ --- соответствующие умножения в $H^\star$ и в $E(E(H)^\star)$). Чтобы это доказать, дополним диаграмму \eqref{H-cong-(H^*)^*-4-1} до диаграммы
\beq\label{H-cong-(H^*)^*-5}
 \xymatrix @R=3.pc @C=8.pc
{
H^\star\odot H^\star\ar[rr]^{\alpha_H\odot\alpha_H}\ar[ddd]_{\varkappa^\star} &  & E(E(H)^\star)\odot E(E(H)^\star)
\ar[ddd]^{E(E(\varkappa)^\star)} \\
& E(H)^{\star}\circledast E(H)^{\star}
\ar[ul]^{@\circ ( e_H)^\star\circledast ( e_H)^\star\qquad}
\ar[ur]_{\quad\quad @\circ  e_{E(H)^{\star}}\circledast  e_{E(H)^{\star}}}
\ar[d]^{E(\varkappa)^\star} & \\
& E(H)^{\star} \ar[dl]_{( e_H)^\star}\ar[dr]^{ e_{E(H)^{\star}}} & \\
H^\star\ar[rr]_{\alpha_H} & & E(E(H)^\star)
}
\eeq
Здесь нужно заметить, что внутренние фигуры в этой диаграмме коммутативны. Например, нижний внутренний треугольник
$$
 \xymatrix @R=2.pc @C=1.5pc
{
& E(H)^{\star} \ar[dl]_{( e_H)^\star}\ar[dr]^{ e_{E(H)^{\star}}} &  \\
H^\star\ar[rr]_{\alpha_H} & & E(E(H)^\star)
}
$$
--- представляет собой просто диаграмму \eqref{H-cong-(H^*)^*-4}. А верхний внутренний треугольник
$$
 \xymatrix @R=2.pc @C=1.5pc
{
H^\star\odot H^\star\ar[rr]^{\alpha_H\odot\alpha_H} &  & E(E(H)^\star)\odot E(E(H)^\star)
 \\
& E(H)^{\star}\circledast E(H)^{\star}
\ar[ul]^{@\circ ( e_H)^\star\circledast ( e_H)^\star\quad\quad}
\ar[ur]_{\quad\quad @\circ  e_{E(H)^{\star}}\circledast  e_{E(H)^{\star}}} &
}
$$
--- можно представить как периметр диаграммы
$$
 \xymatrix @R=2.pc @C=1.5pc
{
H^\star\odot H^\star\ar[rr]^{\alpha_H\odot\alpha_H} &  & E(E(H)^\star)\odot E(E(H)^\star)
 \\
H^\star\circledast H^\star\ar[rr]^{\alpha_H\circledast\alpha_H}\ar[u]_{@} &  & E(E(H)^\star)\circledast E(E(H)^\star) \ar[u]_{@}
 \\
& E(H)^{\star}\circledast E(H)^{\star}
\ar[ul]^{@\circ ( e_H)^\star\circledast ( e_H)^\star\quad\quad}
\ar[ur]_{\quad\quad @\circ  e_{E(H)^{\star}}\circledast  e_{E(H)^{\star}}} &
}
$$
(в которой нижний внутренний треугольник --- диаграмма \eqref{H-cong-(H^*)^*-4}, помноженная на себя тензорным произведением $\circledast$).

Далее, левый внутренний четырехугольник в \eqref{H-cong-(H^*)^*-5}
$$
 \xymatrix @R=4.pc @C=8.pc
{
H^\star\odot H^\star\ar[d]_{\varkappa^\star} &  E(H)^{\star}\circledast E(H)^{\star}
\ar[l]_{@\ \circ ( e_H)^\star\circledast ( e_H)^\star\qquad}
\ar[d]^{E(\varkappa)^\star} \\
H^\star & E(H)^{\star} \ar[l]_{( e_H)^\star}
}
$$
--- коммутативен, потому что он получается  из диаграммы \eqref{DIAG:reflex-otn-obolochki-2}
$$
 \xymatrix @R=4.pc @C=8.pc
{
H\circledast H\ar[r]^{ e_H\odot e_H\circ @} & E(H) \odot E(H) \\
H\ar[r]^{  e_H }\ar[u]^{\varkappa} & E(H) \ar[u]^{E(\varkappa)}
}
$$
действием операции $\star$.

Наконец, правый внутренний четырехугольник в \eqref{H-cong-(H^*)^*-5}
$$
 \xymatrix @R=4.pc @C=8pc
{
E(H)^{\star}\circledast E(H)^{\star}
\ar[r]^{@\ \circ\  e_{E(H)^{\star}}\circledast  e_{E(H)^{\star}}}
\ar[d]^{E(\varkappa)^\star} & E(E(H)^\star)\odot E(E(H)^\star)
\ar[d]^{E(E(\varkappa)^\star)} \\
E(H)^{\star}
\ar[r]^{ e_{E(H)^{\star}}} & E(E(H)^\star)
}
$$
--- коммутативен, потому что представляет собой диаграмму \eqref{DIAG:reflex-otn-obolochki-1}  с подставленными в нее $E(H)^{\star}$ вместо $H$ и $E(\varkappa)^\star$ вместо $\mu$.

После того, как мы поняли, что все внутренние диаграммы в \eqref{H-cong-(H^*)^*-5} коммутативны, мы можем заметить, что если в ней двигаться из вершины $E(H)^{\star}\circledast E(H)^{\star}$ в вершину $E(E(H)^\star)$, то получающиеся морфизмы будут совпадать. В частности если двигаться из $E(H)^{\star}\circledast E(H)^{\star}$ сначала в $H^\star\odot H^\star$, а потом (двумя возможными путями) в $E(E(H)^\star)$, то получающиеся морфизмы тоже будут совпадать:
\beq\label{H-cong-(H^*)^*-6}
 \xymatrix @R=3.pc @C=2.pc
{
H^\star\odot H^\star\ar[rr]^{\alpha_H\odot\alpha_H}\ar[dd]_{\varkappa^\star} &  & E(E(H)^\star)\odot E(E(H)^\star)
\ar[dd]^{E(E(\varkappa)^\star)} \\
& E(H)^{\star}\circledast E(H)^{\star}
\ar[ul]_{\quad @\circ ( e_H)^\star\circledast ( e_H)^\star}
 & \\
H^\star\ar[rr]_{\alpha_H} & & E(E(H)^\star)
}
\eeq
Теперь заметим, что стрелка из центра в левую верхнюю вершину --- $@\ \circ \ ( e_H)^\star\circledast ( e_H)^\star$ --- эпиморфизм в категории $\Ste$. Действительно, $( e_H)^\star: E(H)^{\star}\to H^\star$ --- эпиморфизм стереотипных пространств, потому что по условию \eqref{(env-H)^star:H^star-gets-(Env-H)^star-0} он представляет собой оболочку $e_{E(H)^\star}$, а по условию SE на с.\pageref{uslovie-SE}, всякая такая оболочка всегда эпиморфизм стереотипных пространств. Отсюда следует, что произведение $(e_H)^\star\circledast (e_H)^\star:E(H)^{\star}\circledast E(H)^{\star}\to H^\star\circledast H^\star$ --- тоже эпиморфизм стереотипных пространств (здесь используется \cite[Theorem 4.4.72]{Akbarov-De-Gruyter-I}). С другой стороны, по условию (iii) на с.\pageref{@:H-circledast-H->H-odot-H-Mono}, преобразование Гротендика $@:H\circledast H\to H\odot H$ --- мономорфизм. Поэтому двойственное отображение $@:H^\star\circledast H^\star\to H^\star\odot H^\star$ --- эпиморфизм. Мы получаем, что композиция этих эпиморфизмов $@\circ ( e_H)^\star\circledast ( e_H)^\star$ --- тоже эпиморфизм.

Отсюда можно, наконец, заключить, что в диаграмме \eqref{H-cong-(H^*)^*-6} центральную вершину можно отбросить, и диаграмма станет коммутативной. А это и будет диаграмма \eqref{H-cong-(H^*)^*-4-1}.

3. Покажем, что $\alpha_H$ сохраняет единицу, то есть, что коммутативна диаграмма
\beq\label{H-cong-(H^*)^*-7}
 \xymatrix @R=2.pc @C=1.pc
{
& \C\ar[ld]_{\e^\star}\ar[rd]^{E(E(\e)^\star)} &  \\
H^\star\ar[rr]_{\alpha_H} & & E(E(H)^\star)
}
\eeq
(мы обозначаем буквой $\e$ коединицу в $H$, а $\e^\star$ и $E(E(\e)^\star)$ --- соответствующие единицы в $H^\star$ и в $E(E(H)^\star)$). Для этого дополним диаграмму \eqref{H-cong-(H^*)^*-7} до диаграммы
\beq\label{H-cong-(H^*)^*-8}
 \xymatrix @R=2.5pc @C=5pc
{
& \C\ar[d]^{E(\e)^\star}\ar@/_4ex/[ldd]_{\e^\star}\ar@/^4ex/[rdd]^{E(E(\e)^\star)} &  \\
& E(H)^{\star}\ar[ld]_{( e_H)^\star}\ar[rd]^{ e_{E(H)^{\star}}} &  \\
H^\star\ar[rr]_{\alpha_H} & & E(E(H)^\star)
}
\eeq
В этой диаграмме нижний внутренний треугольник коммутативен, потому что это диаграмма \eqref{H-cong-(H^*)^*-4}. Верхний левый внутренний треугольник
$$
 \xymatrix @R=2.5pc @C=2.5pc
{
& \C\ar[dr]^{E(\e)^\star}\ar[ld]_{\e^\star}
 &  \\
H^\star & & E(H)^{\star}\ar[ll]_{( e_H)^\star}
 }
$$
--- коммутативен, потому что это результат применения функтора $\star$ к правой диаграмме в \eqref{DIAG:reflex-otn-obolochki-3}:
$$
 \xymatrix @R=2.pc @C=2.pc
{
H\ar[rr]^{  e_H }\ar[dr]_{\e} & & E(H) \ar[dl]^{E(\e)} \\
& \C &
}
$$
А верхний правый внутренний треугольник
$$
 \xymatrix @R=2.5pc @C=2.5pc
{
& \C\ar[dr]^{E(E(\e)^\star)}\ar[ld]_{E(\e)^\star}
 &  \\
E(H)^{\star}\ar[rr]_{ e_{E(H)^{\star}}} & & E(E(H)^\star)
 }
$$
--- коммутативен, потому что это правая диаграмма в \eqref{DIAG:reflex-otn-obolochki-3**}.

Мы получаем, что все внутренние треугольники в \eqref{H-cong-(H^*)^*-8} коммутативны, значит периметр тоже коммутативен, а это и есть диаграмма \eqref{H-cong-(H^*)^*-7}.

4. Покажем, что $\alpha_H$ сохраняет коумножение, то есть, что коммутативна диаграмма
\beq\label{H-cong-(H^*)^*-9}
 \xymatrix @R=2.pc @C=6pc
{
H^\star\odot H^\star\ar[r]^{\alpha_H\odot\alpha_H} &  E(E(H)^\star)\odot E(E(H)^\star)
 \\
H^\star\ar[r]_{\alpha_H}\ar[u]^{\mu^\star} & E(E(H)^\star)
\ar[u]_{E(E(\mu)^\star)}
}
\eeq
(мы обозначаем буквой $\mu$ умножение в $H$, а $\mu^\star$ и $E(E(\mu)^\star)$ --- соответствующие коумножения в $H^\star$ и в $E(E(H)^\star)$). Чтобы это доказать, дополним диаграмму \eqref{H-cong-(H^*)^*-9} до диаграммы
\beq\label{H-cong-(H^*)^*-10}
 \xymatrix @R=3.pc @C=8.pc
{
H^\star\odot H^\star\ar[rr]^{\alpha_H\odot\alpha_H} &  & E(E(H)^\star)\odot E(E(H)^\star)
 \\
& E(H)^{\star}\circledast E(H)^{\star}
\ar[ul]^{@\circ ( e_H)^\star\circledast ( e_H)^\star\qquad}
\ar[ur]_{\quad\quad @\circ  e_{E(H)^{\star}}\circledast  e_{E(H)^{\star}}}
 & \\
& E(H)^{\star}\ar[u]_{E(\mu)^\star} \ar[dl]_{( e_H)^\star}\ar[dr]^{ e_{E(H)^{\star}}} & \\
H^\star\ar[rr]_{\alpha_H} \ar[uuu]^{\mu^\star} & & E(E(H)^\star)
\ar[uuu]_{E(E(\mu)^\star)}
}
\eeq
Здесь внутренние треугольники, верхний и нижний, коммутативны, потому что мы это уже проверили, когда рассматривали диаграмму \eqref{H-cong-(H^*)^*-5}. Левый внутренний четырехугольник в \eqref{H-cong-(H^*)^*-10}
$$
 \xymatrix @R=4.pc @C=8.pc
{
H^\star\odot H^\star &  E(H)^{\star}\circledast E(H)^{\star}
\ar[l]_{@\ \circ ( e_H)^\star\circledast ( e_H)^\star\qquad}
 \\
H^\star\ar[u]^{\mu^\star} & E(H)^{\star} \ar[l]_{( e_H)^\star}
\ar[u]_{E(\mu)^\star}
}
$$
--- коммутативен, потому что он получается  из диаграммы \eqref{DIAG:reflex-otn-obolochki-1}
$$
 \xymatrix @R=4.pc @C=8.pc
{
H\circledast H\ar[r]^{ e_H\odot e_H\circ @}
\ar[d]_{\mu} & E(H) \odot E(H) \ar[d]_{E(\mu)} \\
H\ar[r]^{  e_H } & E(H)
}
$$
действием операции $\star$.

А правый внутренний четырехугольник в \eqref{H-cong-(H^*)^*-10}
$$
 \xymatrix @R=4.pc @C=8pc
{
E(H)^{\star}\circledast E(H)^{\star}
\ar[r]^{@\ \circ\  e_{E(H)^{\star}}\circledast  e_{E(H)^{\star}}}
 & E(E(H)^\star)\odot E(E(H)^\star)
 \\
E(H)^{\star}
\ar[r]^{ e_{E(H)^{\star}}}
\ar[u]_{E(\mu)^\star} & E(E(H)^\star)
\ar[u]_{E(E(\mu)^\star)}
}
$$
--- коммутативен, потому что представляет собой диаграмму \eqref{DIAG:reflex-otn-obolochki-2}  с подставленными в нее $E(H)^{\star}$ вместо $H$ и $E(\mu)^\star$ вместо $\varkappa$.

После того, как мы поняли, что все внутренние диаграммы в \eqref{H-cong-(H^*)^*-10} коммутативны, мы можем заметить, что если в ней двигаться из вершины $E(H)^{\star}$ в вершину $E(E(H)^\star)\odot E(E(H)^\star)$, то получающиеся морфизмы будут совпадать. В частности если двигаться из $E(H)^{\star}$ сначала в $H^\star$, а потом (двумя возможными путями) в $E(E(H)^\star)\odot E(E(H)^\star)$, то получающиеся морфизмы тоже будут совпадать:
\beq\label{H-cong-(H^*)^*-11}
 \xymatrix @R=3.pc @C=2.pc
{
H^\star\odot H^\star\ar[rr]^{\alpha_H\odot\alpha_H} &  & E(E(H)^\star)\odot E(E(H)^\star)
 \\
& E(H)^{\star} \ar[dl]^{( e_H)^\star} & \\
H^\star\ar[rr]_{\alpha_H} \ar[uu]^{\mu^\star} & & E(E(H)^\star)
\ar[uu]_{E(E(\mu)^\star)}
}
\eeq
Теперь заметим, что стрелка из центра в левую нижнюю вершину --- $( e_H)^\star: E(H)^{\star}\to H^\star$ --- эпиморфизм в категории $\Ste$, потому что по условию \eqref{(env-H)^star:H^star-gets-(Env-H)^star-0} он представляет собой оболочку $e_{E(H)^\star}$, а по условию SE на с.\pageref{uslovie-SE}, всякая такая оболочка всегда эпиморфизм стереотипных пространств. Отсюда следует, что если в диаграмме \eqref{H-cong-(H^*)^*-11} центральную вершину отбросить, то она станет коммутативной. А это и будет диаграмма \eqref{H-cong-(H^*)^*-9}.

5. Теперь покажем, что $\alpha_H$ сохраняет коединицу, то есть, что коммутативна диаграмма
\beq\label{H-cong-(H^*)^*-12}
 \xymatrix @R=2.pc @C=1.pc
{
& \C &  \\
H^\star\ar[rr]_{\alpha_H}\ar[ru]^{\iota^\star} & & E(E(H)^\star)
\ar[lu]_{E(E(\iota)^\star)}
}
\eeq
(мы обозначаем буквой $\iota$ единицу в $H$, а $\iota^\star$ и $E(E(\iota)^\star)$ --- соответствующие коединицы в $H^\star$ и в $E(E(H)^\star)$). Для этого дополним диаграмму \eqref{H-cong-(H^*)^*-12} до диаграммы
\beq\label{H-cong-(H^*)^*-13}
 \xymatrix @R=2.5pc @C=5pc
{
& \C &  \\
& E(H)^{\star}\ar[u]_{E(\iota)^\star}
\ar[ld]_{( e_H)^\star}\ar[rd]^{ e_{E(H)^{\star}}} &  \\
H^\star\ar[rr]_{\alpha_H}\ar@/^4ex/[ruu]^{\iota^\star} & & E(E(H)^\star)\ar@/_4ex/[luu]_{E(E(\iota)^\star)}
}
\eeq
В этой диаграмме нижний внутренний треугольник коммутативен, потому что это диаграмма \eqref{H-cong-(H^*)^*-4}. Верхний левый внутренний треугольник
$$
 \xymatrix @R=2.5pc @C=2.5pc
{
& \C
 &  \\
H^\star\ar[ru]^{\iota^\star} & & E(H)^{\star}\ar[ll]_{( e_H)^\star}
\ar[ul]_{E(\iota)^\star}
 }
$$
--- коммутативен, потому что это результат применения функтора $\star$ к левой диаграмме в \eqref{DIAG:reflex-otn-obolochki-3}:
$$
 \xymatrix @R=2.pc @C=2.pc
{
H\ar[rr]^{  e_H } & & E(H)  \\
& \C\ar[ul]^{\iota}\ar[ur]_{E(\iota)} &
}
$$
А верхний правый внутренний треугольник
$$
 \xymatrix @R=2.5pc @C=2.5pc
{
& \C
 &  \\
E(H)^{\star}\ar[rr]_{ e_{E(H)^{\star}}}
\ar[ru]^{E(\iota)^\star} & & E(E(H)^\star)
\ar[ul]_{E(E(\iota)^\star)}
 }
$$
--- коммутативен, потому что это левая диаграмма в \eqref{DIAG:reflex-otn-obolochki-3**}.

Мы получаем, что все внутренние треугольники в \eqref{H-cong-(H^*)^*-13} коммутативны. Отсюда следует, что если в диаграмме
\beq\label{H-cong-(H^*)^*-14}
 \xymatrix @R=2.5pc @C=5pc
{
& \C &  \\
& E(H)^{\star}
\ar[ld]_{( e_H)^\star} &  \\
H^\star\ar[rr]_{\alpha_H}\ar@/^4ex/[ruu]^{\iota^\star} & & E(E(H)^\star)\ar@/_4ex/[luu]_{E(E(\iota)^\star)}
}
\eeq
двигаться из центра в верхнюю вершину (двумя возможными путями), то получающиеся морфизмы будут совпадать. Поскольку первый морфизм при таком движении --- $( e_H)^\star: E(H)^{\star} \to H^\star$ --- является оболочкой $e_{E(H)^\star}$ (в силу условия \eqref{(env-H)^star:H^star-gets-(Env-H)^star-0}), и значит, по условию SE на с.\pageref{uslovie-SE}, эпиморфизмом стереотипных пространств, мы получаем, что его можно отбросить, и это значит, что должна быть коммутативна диаграмма \eqref{H-cong-(H^*)^*-12}.

6. Теперь нужно проверить, что $\alpha_H$ сохраняет антипод:
\beq\label{H-cong-(H^*)^*-15}
 \xymatrix @R=3.pc @C=4.pc
{
H^\star\ar[r]^{ \alpha_H }\ar[d]_{\sigma^\star} & E(E(H)^\star) \ar[d]^{E(E(\sigma)^\star)} \\
H^\star\ar[r]^{ \alpha_H } & E(E(H)^\star)
}
\eeq
(мы обозначаем буквой $\sigma$ антипод в $H$, а $\sigma^\star$ и $E(E(\sigma)^\star)$ --- обозначения для антиподов в $H^\star$ и в $E(E(H)^\star)$). Мы дополняем диаграмму \eqref{H-cong-(H^*)^*-15} до диаграммы
\beq\label{H-cong-(H^*)^*-16}
 \xymatrix @R=3.pc @C=4.pc
{
H^\star\ar[rr]^{ \alpha_H }\ar[ddd]_{\sigma^\star} & & E(E(H)^\star) \ar[ddd]^{E(E(\sigma)^\star)} \\
& E(H)^{\star}\ar[ul]^{( e_H)^\star}\ar[ur]_{ e_{E(H)^{\star}}}
\ar[d]^{E(\sigma)^\star} & \\
& E(H)^{\star}\ar[dl]_{( e_H)^\star}\ar[dr]^{ e_{E(H)^{\star}}} & \\
H^\star\ar[rr]^{ \alpha_H } & & E(E(H)^\star)
}
\eeq
Здесь верхний и нижний внутренние треугольники коммутативны, потому что это диаграммы \eqref{H-cong-(H^*)^*-4}. Левый внутренний четырехугольник
$$
 \xymatrix @R=3.pc @C=4.pc
{
H^\star\ar[d]_{\sigma^\star} & E(H)^{\star}\ar[l]_{( e_H)^\star}
\ar[d]^{E(\sigma)^\star}  \\
H^\star & E(H)^{\star}\ar[l]_{( e_H)^\star}
}
$$
--- коммутативен, потому что это диаграмма \eqref{DIAG:reflex-otn-obolochki-4*}. А правый внутренний четырехугольник
$$
 \xymatrix @R=3.pc @C=4.pc
{
 E(H)^{\star}\ar[r]^{ e_{E(H)^{\star}}}
\ar[d]^{E(\sigma)^\star} & E(E(H)^\star) \ar[d]^{E(E(\sigma)^\star)} \\
 E(H)^{\star}\ar[r]^{ e_{E(H)^{\star}}} & E(E(H)^\star)
}
$$
--- коммутативен, потому что это диаграмма \eqref{DIAG:reflex-otn-obolochki-4**}.
Из коммутативности всех внутренних фигур в \eqref{H-cong-(H^*)^*-16} следует, что если в диаграмме
$$
 \xymatrix @R=3.pc @C=4.pc
{
H^\star\ar[rr]^{ \alpha_H }\ar[dd]_{\sigma^\star} & & E(E(H)^\star) \ar[dd]^{E(E(\sigma)^\star)} \\
& E(H)^{\star}\ar[ul]^{( e_H)^\star}
 & \\
H^\star\ar[rr]^{ \alpha_H } & & E(E(H)^\star)
}
$$
двигаться из центра в правый нижний угол, то два возможных пути при этом будут равны (как морфизмы). Поскольку первый морфизм в этих путях --- $( e_H)^\star$ --- представляет собой оболочку (по свойству \eqref{(env-H)^star:H^star-gets-(Env-H)^star-0}), и значит, по свойству SE на с.\pageref{uslovie-SE}, эпиморфизм стереотипных пространств, его можно отбросить, и мы получим коммутативную диаграмму \eqref{H-cong-(H^*)^*-16}.

7. Все сказанное доказывает, что $\alpha_H$ является системой морфизмов $\odot$-алгебр Хопфа. Теперь проверим, что она будет естественным преобразованием функтора $\star$ в функтор $E\circ\star\circ E$, то есть что для любого морфизма $\ph:H\to J$ в категории $\Hopf^E$ коммутативна диаграмма
\beq\label{H-cong-(H^*)^*-17}
 \xymatrix @R=3.pc @C=4.pc
{
J^\star\ar[r]^{ \alpha_J }\ar[d]_{\ph^\star} & E(E(J)^\star) \ar[d]^{E(E(\ph)^\star)} \\
H^\star\ar[r]^{ \alpha_H } & E(E(H)^\star)
}
\eeq
Чтобы это проверить, мы дополняем ее до диаграммы
\beq\label{H-cong-(H^*)^*-18}
 \xymatrix @R=3.pc @C=4.pc
{
J^\star\ar[rr]^{ \alpha_J }\ar[ddd]_{\ph^\star} & & E(E(J)^\star) \ar[ddd]^{E(E(\ph)^\star)}
 \\
& E(J)^{\star}\ar[d]^{E(\ph)^\star}
\ar[ul]^{( e_J)^\star}
\ar[ur]_{ e_{E(J)^{\star}}} & \\
& E(H)^{\star}\ar[dl]_{( e_H)^\star}
\ar[dr]^{ e_{E(H)^{\star}}} & \\
H^\star\ar[rr]^{ \alpha_H } & & E(E(H)^\star)
}
\eeq
Здесь верхний и нижний внутренние треугольники коммутативны, потому что это диаграммы \eqref{H-cong-(H^*)^*-4}. Левый внутренний четырехугольник
$$
 \xymatrix @R=3.pc @C=4.pc
{
J^\star\ar[d]_{\ph^\star} & E(J)^{\star}\ar[d]^{E(\ph)^\star}
\ar[l]_{( e_J)^\star}  \\
H^\star & E(H)^{\star}
\ar[l]_{( e_H)^\star}
}
$$
--- коммутативен, потому что это результат применения функтора $\star$ к диаграмме
$$
 \xymatrix @R=3.pc @C=4.pc
{
J\ar[r]^{ e_J} & E(J)
  \\
H\ar[u]^{\ph}\ar[r]^{ e_H} & E(H)\ar[u]_{E(\ph)}
}
$$
(которая в свою очередь коммутативна потому что оболочка $E$ была выбрана как функтор). А правый внутренний четырехугольник
$$
 \xymatrix @R=3.pc @C=4.pc
{
 E(J)^{\star}\ar[d]^{E(\ph)^\star}
\ar[r]^{ e_{E(J)^{\star}}} & E(E(J)^\star) \ar[d]^{E(E(\ph)^\star)} \\
 E(H)^{\star}
\ar[r]^{ e_{E(H)^{\star}}} & E(E(H)^\star)
}
$$
--- коммутативен, опять же, потому что оболочка $E$ была выбрана как функтор.

У нас получается, что в \eqref{H-cong-(H^*)^*-18} все внутренние диаграммы коммутативны. Вдобавок морфизм $( e_J)^\star: E(J)^{\star}\to J^\star$ является оболочкой в силу \eqref{(env-H)^star:H^star-gets-(Env-H)^star-0}, и поэтому эпиморфизм. Вместе это означает, что периметр диаграммы \eqref{H-cong-(H^*)^*-18} коммутативен. А это и есть диаграмма \eqref{H-cong-(H^*)^*-17}.

8. Мы уже говорили в самом начале, что $\alpha_H$ являются изоморфизмами. Из этого следует, что семейство $\{\alpha_H;\ H\in \Hopf^E\}$ является не просто естественным преобразованием функторов, но и их изоморфизмом.
\epr

\chapter{ПОСТРОЙКИ}

\section{Постройки}

\subsection{Постройки в категориях}

\paragraph{Фундамент.}

Пусть нам дано следующее:
\bit{

\item[1)] категория $\sf K$,

\item[2)] частично упорядоченное множество $( I,\le)$,

\item[3)] ковариантная система $\{\iota_i^j:X_i\to X_j,\ i\le j\}$ в $\sf K$ над $ I$

\item[4)] контравариантная система $\{\varkappa_j^i:X_i\gets X_j,\ i\le j\}$ в $\sf K$ над $ I$

\item[5)] для любой пары индексов $i\le j$ морфизм $\iota_i^j$ является коретракцией для морфизма $\varkappa_j^i$ (а морфизм $\varkappa_j^i$ --- ретракцией для $\iota_i^j$):
\beq\label{varkappa_j^i-circ-iota_i^j=1}
\varkappa_j^i\circ\iota_i^j=1_{X_i},\qquad i\in I
\eeq
}\eit
Эти данные удобно представлять картинкой
\beq\label{DEF:diagr-building}
 \xymatrix  @R=3pc @C=3pc
 {
X_k\ar@{-->}[r]^{1_{X_k}}
 &  X_k\ar@/^3ex/[d]_{\varkappa_k^j} \ar@/^12ex/[dd]^{\varkappa_k^i} \\
X_j\ar@{-->}[r]^{1_{X_j}}\ar@/^3ex/[u]_{\iota_j^k}  &  X_j \ar@/^3ex/[d]_{\varkappa_j^i} \\
X_i\ar@{-->}[r]^{1_{X_i}} \ar@/^12ex/[uu]^{\iota_i^k} \ar@/^3ex/[u]_{\iota_i^j}  & X_i
 }\put(40,-50){$(i\le j\le k)$}
\eeq
\noindent
Условимся такую пару систем ${\varPhi}=(\{\iota_i^j\},\{\varkappa^i_j\})$ называть {\it фундаментом} (foundation) в категории $\sf K$ над направленным множеством $I$.

\paragraph{Постройки.}

Пусть нам дан какой-то фундамент ${\varPhi}=(\{\iota_i^j\},\{\varkappa_j^i\})$ и тройка $(X,\{\sigma_i\},\{\pi^i\})$, в которой $X$ --- некий объект в категории $\sf K$, $\{\sigma_i;\ i\in I\}$ --- семейство морфизмов
$$
\sigma_i: X_i\to X,\qquad i\in I,
$$
с которым $X$ образует инъективный конус над ковариантной системой $\{\iota_i^j\}$
\beq\label{building-1}
 \xymatrix  @R=3pc @C=3pc
 {
 & X & \\
X_i\ar@/^3ex/[ur]^{\sigma_i}\ar[rr]_{\iota_i^j} &  & X_j\ar@/_3ex/[ul]_{\sigma_j} \\
 }\put(40,-25){$(i\le j)$}
\eeq
а $\{\pi^i;\ i\in I\}$ --- семейство морфизмов
$$
\pi^i: X\to X_i,\qquad i\in I,
$$
с которым $X$ образует проективный конус над контравариантной системой $\{\varkappa^i_j\}$
\beq\label{building-2}
 \xymatrix  @R=3pc @C=3pc
 {
 & X\ar@/_3ex/[dl]_{\pi^i}\ar@/^3ex/[dr]^{\pi^j} & \\
X_i &  & X_j\ar[ll]_{\varkappa^i_j} \\
 }\put(40,-25){$(i\le j)$}
\eeq
И вдобавок пусть для всякого $i\in I$ будет коммутативна диаграмма
\beq\label{building-3}
 \xymatrix  @R=3pc @C=3pc
 {
 & X\ar@/^3ex/[dr]^{\pi^i} & \\
X_i\ar@/^3ex/[ur]^{\sigma_i}\ar[rr]_{1_{X_i}} &  & X_i \\
 }\put(40,-25){$(i\in I)$}
\eeq
Эти данные удобно представлять картинкой
\beq\label{DEF:buillding}
 \xymatrix  @R=3pc @C=3pc
 {
& X \ar@/^3ex/[dr]_{\pi^j} \ar@/^12ex/[ddr]^{\pi^i} & \\
X_j\ar@{-->}[rr]^{1_{X_j}}\ar@/^3ex/[ur]_{\sigma_j}  &  & X_j \ar@/^1ex/[d]_{\varkappa_j^i} \\
X_i\ar@{-->}[rr]^{1_{X_i}} \ar@/^12ex/[uur]^{\sigma_i} \ar@/^1ex/[u]_{\iota_i^j} & & X_i
 }\put(40,-50){$(i\le j)$}
\eeq
\noindent
Тогда такую тройку $(X,\{\sigma_i\},\{\pi^i\})$ мы будем называть {\it постройкой} (building) на фундаменте ${\varPhi}$ в категории $\sf K$ над направленным множеством $I$.

\brem
Условия \eqref{building-1}---\eqref{building-2}---\eqref{building-3} автоматически влекут за собой коммутативность диаграмм
\beq\label{building-4}
 \xymatrix  @R=3pc @C=3pc
 {
 & X\ar@/^3ex/[dr]^{\pi^j} & \\
X_i\ar@/^3ex/[ur]^{\sigma_i}\ar[rr]_{\iota_i^j} &  & X_j \\
 }
\qquad 
\xymatrix  @R=3pc @C=3pc
 {
 & X\ar@/^3ex/[dr]^{\pi^i} & \\
X_j\ar@/^3ex/[ur]^{\sigma_j}\ar[rr]_{\varkappa_j^i} &  & X_i \\
 } \put(40,-25){$(i\le j)$}
\eeq
\erem
\bpr
1. Дополним левую диаграмму до диаграммы
$$
 \xymatrix  @R=3pc @C=3pc
 {
 & X\ar@/^3ex/[dr]^{\pi^j} & \\
X_i\ar@/^3ex/[ur]^{\sigma_i}\ar[r]_{\iota_i^j} & X_j\ar[u]^{\sigma_j}\ar[r]_{1_{X_j}}  & X_j \\
 }
$$
Здесь левый внутренний треугольник коммутативен, потому что это диаграмма \eqref{building-1}, а правый внутренний треугольник коммутативен потому что это диаграмма \eqref{building-3}. Значит периметр, то есть левая диаграмма в \eqref{building-4} тоже коммутативен.

2. Дополним правую диаграмму в \eqref{building-4} до диаграммы
$$
\xymatrix  @R=3pc @C=3pc
 {
 & X\ar[d]^{\pi^j}\ar@/^3ex/[dr]^{\pi^i} & \\
X_j\ar@/^3ex/[ur]^{\sigma_j}\ar[r]_{1_{X_j}} & X_j\ar[r]_{\varkappa_j^i} & X_i \\
 }
$$
Здесь левый внутренний треугольник коммутативен, потому что это диаграмма \eqref{building-3}, а правый внутренний треугольник коммутативен потому что это диаграмма \eqref{building-2}. Значит периметр, то есть правая диаграмма в \eqref{building-4} тоже коммутативен.
\epr

\paragraph{Морфизмы построек.}

Пусть далее $(X,\{\sigma_i\},\{\pi^i\})$ и $(Y,\{\varsigma_i\},\{\varpi^i\})$ --- две постройки на фундаменте ${\varPhi}$ в категории $\sf K$ над направленным множеством $I$. Морфизм 
$$
\ph:X\to Y
$$
называется {\it морфизмом построек} $(X,\{\sigma_i\},\{\pi^i\})$ и $(Y,\{\varsigma_i\},\{\varpi^i\})$, если он является одновременно морфизмом инъективных конусов, то есть коммутативны все диаграммы
\beq\label{building-morphism-1}
 \xymatrix  @R=3pc @C=3pc
 {
 X\ar[rr]^{\ph} &  & Y \\
 & X_i\ar@/^3ex/[ul]^{\sigma_i} \ar@/_3ex/[ur]_{\varsigma_i} &
 }\put(40,-25){$(i\in I)$}
\eeq
и проективных конусов, то есть коммутативны все диаграммы
\beq\label{building-morphism-2}
 \xymatrix  @R=3pc @C=3pc
 {
 X\ar[rr]^{\ph}\ar@/_3ex/[dr]_{\pi_i} &  & Y\ar@/^3ex/[dl]^{\varpi_i} \\
 & X_i  &
 }\put(40,-25){$(i\in I)$}
\eeq
Морфизмы построек превращают класс ${\sf K}({\varPhi})$ всех построек на фундаменте ${\varPhi}$ в категорию.

\brem\label{REM:building-morphism-3}
Если $\ph:X\to Y$ --- морфизм построек, то коммутативны диаграммы
\beq\label{building-morphism-3}
 \xymatrix  @R=3pc @C=3pc
 {
X\ar@{-->}[r]^{\ph} &  Y\ar@/^3ex/[d]^{\varpi^i} \\
X_i\ar@/^3ex/[u]^{\sigma_i}\ar[r]_{1_{X_i}} & X_i 
 }\put(40,-25){$(i\in I)$}
\eeq 
\beq\label{building-morphism-3-1}
 \xymatrix  @R=3pc @C=3pc
 {
X\ar@{-->}[r]^{\ph} &  Y\ar@/^3ex/[d]^{\varpi^j} \\
X_i\ar@/^3ex/[u]^{\sigma_i}\ar[r]_{\iota_i^j} & X_j 
 }\put(40,-25){$(i\le j)$}
\eeq 
\beq\label{building-morphism-3-2}
 \xymatrix  @R=3pc @C=3pc
 {
X\ar@{-->}[r]^{\ph} &  Y\ar@/^3ex/[d]^{\varpi^i} \\
X_j\ar@/^3ex/[u]^{\sigma_j}\ar[r]_{\varkappa_j^i} & X_i 
 }\put(40,-25){$(i\le j)$}
\eeq

\erem
\bpr
1. Если коммутативны диаграммы \eqref{building-morphism-1} и \eqref{building-morphism-2}, то можно рассмотреть диаграмму 
$$
\xymatrix  @R=3pc @C=3pc
 {
X\ar@{-->}[r]^{\ph} &  Y\ar@/^3ex/[d]^{\varpi^i} \\
X_i\ar@/^3ex/[u]^{\sigma_i}\ar@/^3ex/[ur]_{\varsigma_i}\ar[r]_{1_{X_i}} & X_i 
 }\put(40,-25){$(i\in I)$}
$$
в которой внутренние  треугольники коммутативны, будучи диаграммами \eqref{building-morphism-1} и \eqref{building-3}, поэтому периметр тоже коммутативен.

2. После того, как доказана коммутативность \eqref{building-morphism-3}, мы заменяем в ней индекс $i$ на $j$
\beq\label{building-morphism-3-3}
 \xymatrix  @R=3pc @C=3pc
 {
X\ar@{-->}[r]^{\ph} &  Y\ar@/^3ex/[d]^{\varpi^j} \\
X_j\ar@/^3ex/[u]^{\sigma_j}\ar[r]_{1_{X_j}} & X_j 
 }\put(40,-25){$(j\in I)$}
\eeq 
и дополняем ее до диаграммы
$$
 \xymatrix  @R=3pc @C=3pc
 {
& X\ar@{-->}[r]^{\ph} &  Y\ar@/^3ex/[d]^{\varpi^j} \\
X_i\ar[r]_{\iota_i^j}\ar@/^3ex/[ur]^{\sigma_i} & X_j\ar@/^3ex/[u]^{\sigma_j}\ar[r]_{1_{X_j}} & X_j 
 }\put(40,-25){$(i\le j)$}
$$ 
(в которой левый внутренний треугольник --- диаграмма \eqref{building-1}). Теперь мы выбрасываем стрелку $\sigma_j$ и склеиваем вершины $X_j$ по стрелке $1_{X_j}$, и у нас получается диаграмма \eqref{building-morphism-3-1}.

3. Наконец, чтобы доказать \eqref{building-morphism-3-2}, мы дополняем \eqref{building-morphism-3-3} до диаграммы
$$
 \xymatrix  @R=3pc @C=3pc
 {
 X\ar@{-->}[r]^{\ph} &  Y\ar@/^3ex/[d]^{\varpi^j}\ar@/^3ex/[dr]^{\varpi^i} & \\
 X_j\ar@/^3ex/[u]^{\sigma_j}\ar[r]_{1_{X_j}} & X_j\ar[r]_{\varkappa_j^i} & X_i
 }\put(40,-25){$(i\le j)$}
$$ 
(в которой правый внутренний треугольник --- диаграмма \eqref{building-2} для постройки $Y$). Теперь мы выбрасываем стрелку $\pi^j$ и склеиваем вершины $X_j$ по стрелке $1_{X_j}$, и у нас получается диаграмма \eqref{building-morphism-3-2}.
\epr

\subsection{Каркас и здание}

\paragraph{Теорема о каркасе и здании.}

Категорию $\sf K$ мы называем {\it полной}, если в ней существуют инъективные и проективные пределы у любого малого функтора.

\btm\label{TH:frame->edifice}
Если индексное множество $I$ направлено по возрастанию, а исходная категория $\sf K$ полна, то категория ${\sf K}({\varPhi})$ всех построек в $\sf K$ на фундаменте ${\varPhi}$ над индексным множеством $I$ обладает 
\bit{

\item[---] инициальным объектом, которым в данном случае будет инъективный предел $\injlim_{k\in I} X_k$ в категории $\sf K$ с системой инъекций $\sigma_i$, с которой он образует инъективный конус над ковариантной системой $\iota_i^j$, и подходящим образом выбранной системой проекций $\pi^i$, и   

\item[---] финальным объектом, которым в данном случае будет проективный предел $\projlim_{k\in I} X_k$ в категории $\sf K$ с системой проекций $\pi^i$, с которой он образует проективный конус над контравариантной системой $\varkappa_j^i$, и подходящим образом выбранной системой инъекций $\sigma_i$.
}\eit

\etm

\bpr
Мы покажем, что в ${\sf K}({\varPhi})$ есть инициальный объект, а доказательство существования финального объекта проводится по аналогии. Им будет инъективный предел $\injlim_{k\in I} X_k$ в $\sf K$ системы коретракций $\{\iota_i^j\}$. У него инъекциями
$$
X_i\overset{\sigma_i}{\longrightarrow}\injlim_{k\in I} X_k, \qquad i\in I,
$$
будут естественные инъекции в инъективный предел (образующие с ним инъективный конус, и поэтому для них диаграмма \eqref{building-1} коммутативна автоматически). Нам нужно объяснить, как строятся морфизмы 
$$
\injlim_{k\in I} X_k\overset{\pi_i}{\longrightarrow} X_i,\qquad i\in I,
$$
удовлетворяющие условиям \eqref{building-2} и \eqref{building-3}, и почему получающаяся постройка
$$
(\injlim_{k\in I} X_k,\ \{\sigma_i\},\ \{\pi^i\})
$$
является инициальным объектом в категории ${\sf K}({\varPhi})$. Это довольно громоздкая часть работы, мы ее опишем в несколько этапов.

1. Заметим, что для любых индексов $i,j,k\in I$ со свойством $i\le j\le k$ коммутативна диаграмма
\beq\label{frame->edifice-1}
 \xymatrix  @R=1pc @C=3pc
 {
X_k\ar@/^2ex/[dr]^{\varkappa_k^i} & \\
& X_i \\
X_j\ar@/^3ex/[uu]^{\iota_j^k}\ar@/_2ex/[ur]_{\varkappa_j^i} & 
 }\put(40,-30){$(i\le j\le k)$}
\eeq 
Действительно,
$$
\varkappa_k^i=\varkappa_j^i\circ\varkappa_k^j \quad\Rightarrow\quad 
\varkappa_k^i\circ\iota_j^k=\varkappa_j^i\circ\underbrace{\varkappa_k^j\circ\iota_j^k}_{1_{X_j}}=\varkappa_j^i
$$

2. Из \eqref{frame->edifice-1} следует, что при фиксированном $i\in I$ система морфизмов $\{\varkappa_j^i;\ j\ge i\}$ является инъективным конусом для ковариантной системы $\{\iota_j^k:X_j\to X_k,\ i\le j\le k\}$. Поэтому существует единственный морфизм 
$$
\tau^i:\injlim_{k\ge i}X_k\to X_i,
$$
замыкающий все диаграммы
\beq\label{frame->edifice-2}
 \xymatrix  @R=1pc @C=3pc
 {
\injlim_{k\ge i}X_k\ar@{-->}@/^2ex/[dr]^{\tau^i} & \\
& X_i \\
X_j\ar@/^3ex/[uu]^{\sigma_j}\ar@/_2ex/[ur]_{\varkappa_j^i} & 
 }\put(40,-30){$(i\le j)$}
\eeq

3. Вспомним теперь, что множество $I$ направлено по возрастанию. Из этого следует, что для всякого $i\in I$ естественный морфизм
$$
\upsilon_i:\injlim_{k\ge i}X_k\to \injlim_{k\in I}X_k
$$
является изоморфизмом. Как следствие, мы можем определить морфизм
$$
\pi^i=\tau^i\circ\upsilon_i^{-1}:\injlim_{k\in I}X_k\to X_i
$$
И у нас получится диаграмма
\beq\label{frame->edifice-3}
 \xymatrix  @R=1pc @C=3pc
 {
\injlim_{k\in I}X_k\ar@{-->}@/^2ex/[dr]^{\pi^i} & \\
& X_i \\
\injlim_{k\ge i}X_k\ar@/^3ex/[uu]^{\upsilon_i}\ar@/_2ex/[ur]_{\tau^i} & 
 }\put(40,-30){$(i\in I)$}
\eeq 

4. Заметим, что для любых $i\le j$ будет коммутативна диаграмма
\beq\label{frame->edifice-2-1}
 \xymatrix  @R=1pc @C=3pc
 {
\injlim_{k\in I}X_k\ar@{-->}@/^2ex/[dr]^{\pi^i} & \\
& X_i \\
X_j\ar@/^3ex/[uu]^{\sigma_j}\ar@/_2ex/[ur]_{\varkappa_j^i} & 
 }\put(40,-30){$(i\le j)$}
\eeq 
Действительно, мы ее можем вписать в диаграмму
\beq\label{frame->edifice-2-2}
 \xymatrix  @R=2pc @C=3pc
 {
\injlim_{k\in I}X_k\ar@{-->}@/^2ex/[drr]^{\pi^i} & & \\
& \injlim_{k\ge i}X_k\ar[r]^{\tau^i}\ar[ul]^{\upsilon_i} & X_i \\
X_j\ar@/^3ex/[uu]^{\sigma_j}\ar[ur]^{\sigma_j}\ar@/_2ex/[urr]_{\varkappa_j^i} & &
 }\put(40,-30){$(i\le j)$}
\eeq 
В ней все внутренние треугольники коммутативны, поэтому периметр -- то есть диаграмма \eqref{frame->edifice-2-1} -- тоже должен быть коммутативен. Можно заметить еще, что из \eqref{frame->edifice-2-1} следует коммутативность диаграммы \eqref{building-3} для наших морфизмов : если положить $j=i$, то мы получим $\varkappa_i^i=1_{X_i}$, и тогда
\beq\label{frame->edifice-2-1-1}
 \xymatrix  @R=1pc @C=3pc
 {
\injlim_{k\in I}X_k\ar@{-->}@/^2ex/[dr]^{\pi^i} & \\
& X_i \\
X_i\ar@/^3ex/[uu]^{\sigma_i}\ar@/_2ex/[ur]_{1_{X_i}} & 
 }\put(40,-30){$(i\in I)$}
\eeq

5. Покажем теперь, что семейство морфизмов $\{\pi^i;\ i\in I\}$ образует проективный конус контравариантной системы $\{\varkappa_j^i:X_i\gets X_j,\ i\le j\}$, то есть что для любых $i\le j$ коммутативна диаграмма \eqref{building-2}:
\beq\label{frame->edifice-4}
 \xymatrix  @R=1pc @C=3pc
 {
 & X^j \ar[dd]^{\varkappa_j^i}\\
\injlim_{k\in I}X_k\ar@/^2ex/[ur]^{\pi^j}\ar@/_2ex/[dr]_{\pi^i} & \\
& X_i \\
 }\put(40,-10){$(i\le j)$}
\eeq 
Для этого возьмем произвольный индекс $k\ge j$ и рассмотрим диаграмму
\beq\label{frame->edifice-5}
 \xymatrix  @R=1pc @C=3pc
 {
 & & X^j \ar[dd]^{\varkappa_j^i}\\
X_k\ar@/_3ex/[drr]_{\varkappa_k^i}\ar[r]^{\sigma_k} & \injlim_{k\in I}X_k\ar@/^2ex/[ur]_{\pi^j}\ar@/_2ex/[dr]^{\pi^i} & \\
& & X_i \\
 }\put(40,-10){$(i\le j\le k)$}
\eeq 
Нам важно понять, что в ней движение из $X_k$ в $X_i$ тремя возможными путями дает один результат:
\beq\label{frame->edifice-6}
\varkappa_j^i\circ\pi^j\circ\sigma_k=\varkappa_k^i=\pi^i\circ\sigma_k,\qquad i\le j\le k
\eeq 
Для этого дополним \eqref{frame->edifice-5} до диаграммы
\beq\label{frame->edifice-7}
 \xymatrix  @R=1pc @C=3pc
 {
 & & X^j \ar[dd]^{\varkappa_j^i}\\
X_k\ar@/^3ex/[urr]^{\varkappa_k^j}\ar@/_3ex/[drr]_{\varkappa_k^i}\ar[r]^{\sigma_k} & \injlim_{k\in I}X_k\ar@/^2ex/[ur]_{\pi^j}\ar@/_2ex/[dr]^{\pi^i} & \\
& & X_i \\
 }\put(40,-10){$(i\le j\le k)$}
\eeq 
В ней левый нижний треугольник 
\beq\label{frame->edifice-7-1}
 \xymatrix  @R=1pc @C=3pc
 {
X_k\ar@/_3ex/[drr]_{\varkappa_k^i}\ar[r]^{\sigma_k} & \injlim_{k\in I}X_k\ar@/_2ex/[dr]^{\pi^i} & \\
& & X_i \\
 }
\eeq
--- коммутативен потому что это диаграмма \eqref{frame->edifice-2-1}. Левый верхний треугольник  
\beq\label{frame->edifice-7-2}
 \xymatrix  @R=1pc @C=3pc
 {
 & & X^j \\
X_k\ar@/^3ex/[urr]^{\varkappa_k^j}\ar[r]^{\sigma_k} & \injlim_{k\in I}X_k\ar@/^2ex/[ur]_{\pi^j} & \\
 }
\eeq
--- коммутативен потому что это диаграмма \eqref{frame->edifice-2-1} с подставленным $j$ вместо $i$.
А периметр 
\beq\label{frame->edifice-7-3}
 \xymatrix  @R=1pc @C=3pc
 {
 & & X^j \ar[dd]^{\varkappa_j^i}\\
X_k\ar@/^3ex/[urr]^{\varkappa_k^j}\ar@/_3ex/[drr]_{\varkappa_k^i} &  & \\
& & X_i \\
 }
\eeq
--- коммутативен, потому что это правый треугольник в диаграмме \eqref{DEF:diagr-building} (с измененными индексами).

В результате мы получаем, что второе равенство в \eqref{frame->edifice-6} справедливо из-за \eqref{frame->edifice-7-1}:
$$
\varkappa_k^i=\eqref{frame->edifice-7-1}=\pi^i\circ\sigma_k
$$
а первое доказывается цепочкой
$$
\varkappa_j^i\circ\pi^j\circ\sigma_k=\eqref{frame->edifice-7-2}=\varkappa_j^i\circ\varkappa_k^j=\eqref{frame->edifice-7-3}=\varkappa_k^i
$$

Теперь заметим, что равенства \eqref{frame->edifice-6} можно понимать, как заявление, что для любого $k\ge j$ коммутативны диаграммы
$$
 \xymatrix  @R=1pc @C=3pc
 {
\injlim_{k\in I}X_k\ar@{-->}@/^2ex/[dr]^{\varkappa_j^i\circ\pi^j} & \\
& X_i \\
X_k\ar@/^3ex/[uu]^{\sigma_k}\ar@/_2ex/[ur]_{\varkappa_k^i} & 
 }\qquad 
 \xymatrix  @R=1pc @C=3pc
 {
\injlim_{k\in I}X_k\ar@{-->}@/^2ex/[dr]^{\pi^i} & \\
& X_i \\
X_k\ar@/^3ex/[uu]^{\sigma_k}\ar@/_2ex/[ur]_{\varkappa_k^i} & 
 }
$$ 
Поскольку $I$ --- направленное по возрастанию множество, мы можем заменить их на диаграммы
$$
 \xymatrix  @R=1pc @C=3pc
 {
\injlim_{k\ge j}X_k\ar@{-->}@/^2ex/[dr]^{\varkappa_j^i\circ\pi^j} & \\
& X_i \\
X_k\ar@/^3ex/[uu]^{\sigma_k}\ar@/_2ex/[ur]_{\varkappa_k^i} & 
 }\qquad 
 \xymatrix  @R=1pc @C=3pc
 {
\injlim_{k\ge j}X_k\ar@{-->}@/^2ex/[dr]^{\pi^i} & \\
& X_i \\
X_k\ar@/^3ex/[uu]^{\sigma_k}\ar@/_2ex/[ur]_{\varkappa_k^i} & 
 }\put(20,-25){$(j\le k)$}
$$ 
и их коммутативность означает, что пунктирные стрелки в них --- естественные морфизмы из инъективного предела в инъективный конус $\{\varkappa_k^i:X_k\to X_i;\ k\ge j\}$. Поскольку такой морфизм единственен, пунктирная стрелка в диаграмме слева --- та же самая, что в диаграмме справа:
$$
\varkappa_j^i\circ\pi^j=\pi^i.
$$
Это и есть коммутативность диаграммы \eqref{frame->edifice-4}.

6. Мы поняли, что инъективный предел $\injlim_{k\in I}X_k$ с построенными морфизмами $\sigma_i$ и $\pi^i$ удовлетворяет условиям \eqref{building-1}---\eqref{building-2}---\eqref{building-3} (диаграмма \eqref{building-2} здесь --- диаграмма \eqref{frame->edifice-4}, а диаграмма \eqref{building-3} --- диаграмма \eqref{frame->edifice-2-1-1}), то есть $\injlim_{k\in I}X_k$ является постройкой на фундаменте ${\varPhi}$. Нам остается проверить, что эта постройка --- инициальный объект в категории всех построек. Пусть $(Y,\{\varsigma_i\},\{\varpi^i\})$ --- какая-нибудь другая постройка на фундаменте ${\varPhi}$. Поскольку $(Y,\{\varsigma_i\})$ --- инъективный конус в категории ${\sf K}$, существует единственный морфизм $\ph:\injlim_{k\in I}X_k\to Y$ в категории ${\sf K}$, замыкающий все диаграммы \eqref{building-morphism-1} (с подставленным $X=\injlim_{k\in I}X_k$ и нам сейчас удобно будет еще заменить в них $i$ на $j$):
\beq\label{frame->edifice-8}
 \xymatrix  @R=3pc @C=3pc
 {
 \injlim_{k\in I}X_k\ar[rr]^{\ph} &  & Y \\
 & X_j\ar@/^3ex/[ul]^{\sigma_j} \ar@/_3ex/[ur]_{\varsigma_j} &
 }\put(40,-25){$(j\in I)$}
\eeq
 Нам нужно проверить, что он будет замыкать и диаграммы \eqref{building-morphism-2} (с подставленным $X=\injlim_{k\in I}X_k$).
Дополним \eqref{frame->edifice-8} до диаграммы
\beq\label{frame->edifice-8-1}
 \xymatrix  @R=3pc @C=3pc
 {
 \injlim_{k\in I}X_k\ar[r]^{\ph} &  Y\ar[d]^{\varpi^i} \\
 X_j\ar[u]^{\sigma_j} \ar[ur]_{\varsigma_j}\ar[r]_{\varkappa_j^i} & X_i
 }\put(40,-25){$(i\le j)$}
\eeq
Здесь правый внутренний треугольник --- диаграмма \eqref{building-4} с подставленной вместо $X$ постройкой $Y$. Поскольку внутренние треугольники коммутативны, периметр тоже должен быть комммутативен:
\beq\label{frame->edifice-8-2}
 \xymatrix  @R=3pc @C=3pc
 {
 \injlim_{k\in I}X_k\ar[r]^{\ph} &  Y\ar[d]^{\varpi^i} \\
 X_j\ar[u]^{\sigma_j} \ar[r]_{\varkappa_j^i} & X_i
 }\put(40,-25){$(i\le j)$}
\eeq
Добавим сюда морфизм $\pi^i$ из \eqref{frame->edifice-2-1}:
\beq\label{frame->edifice-8-2}
 \xymatrix  @R=3pc @C=3pc
 {
 \injlim_{k\in I}X_k\ar[r]^{\ph}\ar@{-->}[rd]^{\pi^i} &  Y\ar[d]^{\varpi^i} \\
 X_j\ar[u]^{\sigma_j} \ar[r]_{\varkappa_j^i} & X_i
 }\put(40,-25){$(i\le j)$}
\eeq
И заменим инъективный предел $\injlim_{k\in I}X_k$ на его хвост $\injlim_{k\ge i}X_k$, изоморфный ему поскольку $I$ направлено по возрастанию:
\beq\label{frame->edifice-8-3}
 \xymatrix  @R=3pc @C=3pc
 {
 \injlim_{k\ge i}X_k\ar[r]^{\ph}\ar@{-->}[rd]^{\pi^i} &  Y\ar[d]^{\varpi^i} \\
 X_j\ar[u]^{\sigma_j} \ar[r]_{\varkappa_j^i} & X_i
 }\put(40,-25){$(i\le j)$}
\eeq
Теперь зафиксируем $i\in I$ и задумаемся, что у нас получилось: периметр и левый внутренний треугольник коммутативны при любом $j\ge i$. Это значит, что морфизмы $\pi^i$ и $\varpi^i\circ\ph$ --- морфизмы инъективного конуса $\{\injlim_{k\ge i}X_k,\ \sigma_j,\ j\ge i\}$ в инъективный конус $\{X_i,\ \varkappa_j,\ j\ge i\}$ над ковариантной системой $\iota_j^k$ (что второе семейство --- инъективный конус над $\iota_j^k$ доказано в \eqref{frame->edifice-1}). Поскольку таких морфизмов из проективного предела $\injlim_{k\ge i}X_k$ существует всего один, морфизмы $\pi^i$ и $\varpi^i\circ\ph$ должны совпадать:
$$
\pi^i=\varpi^i\circ\ph
$$
То есть должна быть коммутативна диаграмма  
\beq\label{frame->edifice-8-4}
 \xymatrix  @R=3pc @C=3pc
 {
 \injlim_{k\ge i}X_k\ar[r]^{\ph}\ar[rd]^{\pi^i} &  Y\ar[d]^{\varpi^i} \\
  & X_i
 }\put(40,-25){$(i\in I)$}
\eeq
Мы заменяем в ней хвост $\injlim_{k\ge i}X_k$ на исходный предел $\injlim_{k\in I}X_k$ и получаем, что должна быть коммутативна диаграмма
\beq\label{frame->edifice-8-5}
 \xymatrix  @R=3pc @C=3pc
 {
 \injlim_{k\in I}X_k\ar[r]^{\ph}\ar[rd]^{\pi^i} &  Y\ar[d]^{\varpi^i} \\
  & X_i
 }\put(40,-25){$(i\in I)$}
\eeq
И это верно для всякого $i\in I$

Мы получили что хотели: морфизм $\ph$ замыкает не только диаграммы \eqref{building-morphism-1}, но и диаграммы \eqref{building-morphism-2} (с подставленным $X=\injlim_{k\in I}X_k$). Значит, $\ph$ --- морфизм построек. Его единственность следует из единственности морфизма из инъективного предела в диаграммах \eqref{frame->edifice-8}.
\epr

Теорема \ref{TH:frame->edifice} делает осмысленными следующие определения.

\bit{

\item[$\bullet$] Инициальный объект в категории ${\sf K}({\varPhi})$ всех построек на фундаменте ${\varPhi}$ мы будем называть {\it каркасом} (frame) на фундаменте ${\varPhi}$ и обозначать $\fram {\varPhi}$; по теореме \ref{TH:frame->edifice}, как объект в ${\sf K}$ каркас совпадает с инъективным пределом ковариантной системы $\iota_i^j$:
\beq\label{DEF:frame}
\fram {\varPhi}=\injlim_{k\in I} X_k
\eeq 

\item[$\bullet$] Финальный объект в категории ${\sf K}({\varPhi})$ всех построек на фундаменте ${\varPhi}$ мы будем называть {\it зданием} (edifice) на фундаменте ${\varPhi}$ и обозначать $\edif {\varPhi}$; по теореме \ref{TH:frame->edifice}, как объект в ${\sf K}$ здание совпадает с проективным пределом контравариантной системы $\varkappa_j^i$:
\beq\label{DEF:edifice}
\edif {\varPhi}=\projlim_{k\in I} X_k
\eeq 

\item[$\bullet$] Между каркасом и зданием (как между любыми инициальным и финальным объектами) имеется единственный морфизм построек 
\beq\label{DEF:theta:frame-body}
\theta_{\varPhi}:\fram {\varPhi}\to \edif {\varPhi}
\eeq
}\eit

\btm\label{TH:frame-X-edifice}
Если $X$ --- постройка на фундаменте ${\varPhi}$, то существуют два морфизма построек $\alpha:\fram {\varPhi} \to Y$ и $\beta:Y\to \edif{\varPhi}$, замыкающие диаграмму
\beq\label{DEF:frame-X-edifice}
 \xymatrix  @R=3pc @C=3pc
 {
 & X\ar[rd]^{\beta} & \\
 \fram{\varPhi}\ar[rr]^{\theta_{\varPhi}}\ar[ru]^{\alpha} & & \edif {\varPhi}
 }
 \eeq
 \etm
 \bpr
 Поскольку $\fram{\varPhi}$ --- инициальный объект в ${\sf K}(\varPhi)$, существует единственный морфизм $\alpha:\fram {\varPhi} \to X$. С другой стороны, поскольку $\edif{\varPhi}$ --- финальный объект в ${\sf K}(\varPhi)$, существует единственный морфизм  $\beta:X\to \edif{\varPhi}$. Наконец, их композиция $\beta\circ\alpha:\fram{\varPhi} \to \edif{\varPhi}$ --- морфизм из $\fram{\varPhi}$ в $\edif{\varPhi}$, а такой морфизм в этой категории единственный и мы его обозначили символом $\theta_{\varPhi}$.
 \epr

Результаты теоремы \ref{TH:frame->edifice} можно проиллюстрировать так: если $(X,\{\varsigma_i\},\{\varpi^i\})$ --- постройка на фундаменте $\varPhi$, то коммутативна диаграмма
\beq\label{DEF:frame-body}
 \xymatrix  @R=3pc @C=1pc
 {
& Y\ar[rd]^{\beta}\ar@/^30ex/[ddddr]^{\varpi^i} & \\
\fram {\varPhi}\ar@{=}[d]\ar[rr]^{\theta_{\varPhi}}\ar[ru]^{\alpha} & & \edif {\varPhi}\ar@{=}[d] \\ 
\injlim_{k\in I} X_k\ar[rr]^{\theta_{\varPhi}}
 & & \projlim_{k\in I} X_k\ar@/^2ex/[d]_{\pi^j} \ar@/^6ex/[dd]^{\pi^i} \\
X_j\ar@{-->}[rr]^{1_{X_j}}\ar@/^2ex/[u]_{\sigma_j}  & &  X_j \ar@/^2ex/[d]_{\varkappa_j^i} \\
X_i\ar@{-->}[rr]^{1_{X_i}} \ar@/^6ex/[uu]^{\sigma_i}\ar@/^30ex/[uuuur]^{\varsigma_i} \ar@/^2ex/[u]_{\iota_i^j} & & X_i
 }\put(40,-50){$(i\le j)$}
\eeq
Здесь $\sigma_i$ --- естественные инъекции в инъективный предел $\injlim_{k\in I} X_k$, а  $\pi^i$ --- естественные проекции из проективного предела $\projlim_{k\in I} X_k$. Коммутативность диаграмм 
\beq\label{frame->body}
 \xymatrix  @R=3pc @C=3pc
 {
\fram\varPhi\ar@{-->}[r]^{\theta_\varPhi} & \edif\varPhi\ar@/^3ex/[d]^{\pi^i} \\
X_i\ar@/^3ex/[u]^{\sigma_i}\ar[r]_{1_{X_i}} & X_i 
 }\put(40,-25){$(i\in I)$}
\eeq 
--- следствие замечания \ref{REM:building-morphism-3}.

Полезно отметить следующее

\bprop\label{frame->body-ph}
Пусть $\ph:\fram\varPhi\to \edif\varPhi$ -- какой-нибудь морфизм в $\sf K$, замыкающий все диаграммы
\beq\label{frame->body-ph}
 \xymatrix  @R=3pc @C=3pc
 {
\fram\varPhi\ar@{-->}[r]^{\ph} & \edif\varPhi\ar@/^3ex/[d]^{\pi^j} \\
X_j\ar@/^3ex/[u]^{\sigma_j}\ar[r]_{1_{X_j}} & X_j 
 }\put(40,-25){$(j\in I)$}
\eeq 
Тогда $\ph=\theta_\varPhi$.
\eprop
\bpr
Дополним \eqref{frame->body-ph} до диаграммы 
\beq\label{frame->body-ph-1}
 \xymatrix  @R=3pc @C=3pc
 {
\fram\varPhi\ar[r]^{\ph} & \edif\varPhi\ar@/^3ex/[d]^{\pi^j}\ar@/^3ex/[dr]^{\pi^i} &  \\
X_j\ar@/^3ex/[u]^{\sigma_j}\ar[r]_{1_{X_j}} & X_j\ar[r]_{\varkappa_j^i} & X_j
 }\put(40,-25){$(i\le j)$}
\eeq 
(здесь правый внутренний треугольник --- диаграмма \eqref{building-2} для постройки $\edif\varPhi$). Выбросив стрелку $\pi^j$ и склеив вершины $X_j$ по стрелке $1_{X_j}$, мы получим диаграмму
\beq\label{frame->body-ph-2}
 \xymatrix  @R=3pc @C=3pc
 {
\fram\varPhi\ar[r]^{\ph} & \edif\varPhi\ar@/^3ex/[d]^{\pi^i}  \\
X_j\ar@/^3ex/[u]^{\sigma_j}\ar[r]_{\varkappa_j^i} & X_j
 }\put(40,-25){$(i\le j)$}
\eeq
Дальше повторяются рассуждения из доказательства теоремы \ref{TH:frame->edifice}. Мы обозначаем
$$
\ph^i=\pi^i\circ\ph,
$$
и делаем замену $\fram\varPhi=\injlim_{k\in I}X_k$ и $\edif\varPhi=\projlim_{k\in I}X_k$:
\beq\label{frame->body-ph-3}
 \xymatrix  @R=2pc @C=3pc
 {
 \injlim_{k\in I}X_k\ar[r]^{\ph}\ar@{-->}[rd]^{\ph^i} &  \projlim_{k\in I}X_k\ar[d]^{\pi^i}  \\
 X_j\ar[u]^{\sigma_j} \ar[r]_{\varkappa_j^i} & X_i\\
 }\put(40,-25){$(i\le j)$}
\eeq
В этой диаграмме, в частности, нижний левый треугольник коммутативен для любых $i\le j$, поэтому в соответствии с рассуждениями пунктов 2 и 4 доказательства теоремы \ref{TH:frame->edifice}, морфизм $\ph^i$ в ней определяется однозначно как морфизм из инъективного предела ковариантной системы $\{\iota_j^k:X_i\to X_j,\ i\le j\le k\}$. С другой стороны, в соответствии с рассуждениями пункта 6, морфизм $\ph$ тоже определяется однозначно семейством морфизмов $\{\ph^i;\ i\in I\}$, как морфизм в проективный предел контравариантной системы $\{\varkappa_j^i:X_i\gets X_j,\ i\le j\}$. Значит, $\ph$  однозначно определяется системой морфизмов $\{\iota_i^j;\varkappa_j^i;\ i\le j\}$ и совпадает с $\theta_\varPhi$.
\epr

\bprop\label{PROP:pi^i-circ-theta_varPhi-retraktsija}
Для всякого индекса $i\in I$ 
\bit{
\item[---] морфизмы $\sigma_i$ и  $\theta_\varPhi\circ\sigma_i$ являются коретркациями (и поэтому мономорфизмами), а

\item[---] морфизмы $\pi^i$ и $\pi^i\circ\theta_\varPhi$  являются ретракциями (и поэтому эпиморфизмами).
}\eit
\eprop
\bpr
Это следует из равенства $\pi^i\circ\theta_\varPhi\circ\sigma_i=1_{X_i}$.
\epr

\paragraph{Ограничение и расширение фундамента.}

Пусть теперь ${\varPhi}=(\{\iota_i^j\},\{\varkappa_j^i\})$ --- фундамент в полной категории $\sf K$ над направленным множеством $I$, а $J$ --- некое подмножество в $I$, тоже направленное относительно бинарного  отношения $\le$. Мы можем рассмотреть морфизмы $\iota_i^j$ и $\varkappa_j^i$ с индексами $i,j$, бегающими по множеству $J$, и это тоже будет некий фундамент в $\sf K$, но уже над направленным множеством $J$. Обозначим его символом
$$
{\varPhi}_J=(\{\iota_i^j\},\{\varkappa_j^i\}),\qquad i,j\in J.
$$
и условимся называть его {\it ограничением фундамента} ${\varPhi}$, и обозначать это записью
$$
{\varPhi}_J\subseteq {\varPhi}.
$$
А фундамент ${\varPhi}$ при этом мы будем называть {\it расширением фундамента} ${\varPhi}_J$.

Каркас и здание на ограничении и расширении фундамента связаны следующей теоремой

\btm
Пусть ${\varPhi}_J$ --- ограничение фундамента ${\varPhi}$ в полной категории $\sf K$. Тогда естественные морфизмы между каркасами этих фундаментов 
$$
\fram\varPhi_J=\injlim_{i\in J} X_i\overset{\sigma_J}{\longrightarrow}\injlim_{i\in I} X_i=\fram\varPhi
$$
и между их зданиями
$$
\edif\varPhi=\injlim_{i\in I} X_i\overset{\pi^J}{\longrightarrow}\injlim_{i\in J} X_i=\edif\varPhi_J
$$
замыкают диаграмму
\beq\label{frame->body-for-inner-building}
 \xymatrix  @R=3pc @C=3pc
 {
\fram\varPhi\ar[r]^{\theta_{\varPhi}} &  \edif\varPhi\ar@{-->}[d]^{\pi^J} \\
\fram\varPhi_J \ar@{-->}[u]^{\sigma_J}\ar[r]^{\theta_{{\varPhi}_J}} &  \edif\varPhi_J
 }
\eeq 

\etm
\bpr
Пусть, как и раньше, инъекции в инъективный предел $\injlim_{i\in I} X_i$ обозначаются $\sigma_i$ 
$$
\sigma_i:X_i\to \injlim_{i\in J} X_i
$$
а проекции из проективного предела  $\projlim_{i\in I} X_i$ обозначаются $\varkappa^i$
$$
\varkappa^i:\projlim_{i\in J} X_i\to X_i.
$$
В дополнение к этому мы условимся инъекции в инъективный предел  $\injlim_{i\in J} X_i$ обозначать $\varsigma_j$
$$
\varsigma_j:X_j\to \injlim_{i\in J} X_i
$$
а проекции из проективного предела  $\projlim_{i\in J} X_i$ обозначать $\varpi^j$
$$
\varpi^j:\projlim_{i\in J} X_i\to X_j.
$$
Тогда для всякого $j\in J$ будут коммутативны диаграммы
\beq\label{frame->body-for-inner-building-1}
 \xymatrix  @R=3pc @C=3pc
 {
\injlim_{i\in I} X_i  \\
\injlim_{i\in J} X_i\ar[u]_{\sigma_J} \\
X_j \ar[u]_{\varsigma_j} \ar@/^10ex/[uu]^{\sigma_j}
 }\qquad 
 \xymatrix  @R=3pc @C=3pc
 {
\projlim_{i\in I} X_i\ar[d]_{\pi^J}\ar@/^10ex/[dd]^{\pi^j}  \\
\projlim_{i\in J} X_i \ar[d]_{\varpi_j} \\
X_j  
 }
\eeq
Мы их можем вписать в такую диаграмму:
\beq\label{frame->body-for-inner-building-2}
 \xymatrix  @R=3pc @C=3pc
 {
\injlim_{i\in I} X_i\ar[r]^{\theta_{\varPhi}} & \projlim_{i\in I} X_i\ar[d]_{\pi^J}\ar@/^10ex/[dd]^{\pi^j} \\
\injlim_{i\in J} X_i\ar[u]_{\sigma_J} & \projlim_{i\in J} X_i \ar[d]_{\varpi_j} \\
X_j\ar[r]^{1_{X_j}} \ar[u]_{\varsigma_j} \ar@/^10ex/[uu]^{\sigma_j} & X_j  
 }
\eeq
В ней периметр будет коммутативен, потому что это диаграмма \eqref{frame->body}. Поэтому внутренний шестиугольник тоже должен быть коммутативен:
\beq\label{frame->body-for-inner-building-3}
 \xymatrix  @R=3pc @C=3pc
 {
\injlim_{i\in I} X_i\ar[r]^{\theta_{\varPhi}} & \projlim_{i\in I} X_i\ar[d]_{\pi^J} \\
\injlim_{i\in J} X_i\ar[u]_{\sigma_J} & \projlim_{i\in J} X_i \ar[d]_{\varpi_j} \\
X_j\ar[r]^{1_{X_j}} \ar[u]_{\varsigma_j}  & X_j  
 }
\eeq
Это можно интерпретировать как коммутативность прямоугольника
\beq\label{frame->body-for-inner-building-3}
 \xymatrix  @R=3pc @C=5pc
 {
\injlim_{i\in J} X_i\ar[r]^{\pi^J\circ\theta_{\varPhi}\circ\sigma_J} & \projlim_{i\in J} X_i \ar[d]_{\varpi_j} \\
X_j\ar[r]^{1_{X_j}} \ar[u]_{\varsigma_j}  & X_j  
 }
\eeq
И это верно для любого $j\in J$. По предложению \ref{frame->body-ph} это значит, что 
$$
\pi^J\circ\theta_{\varPhi}\circ\sigma_J=\theta_{{\varPhi}_J}
$$ 
Это нам и нужно было доказать.
\epr

Фундамент ${\varPhi}$ мы будем называть {\it индуктивным}, если его индексное множество $I$ изоморфно как частично упорядоченное множество множеству натуральных чисел $\N$. В частности, ограничение ${\varPhi}_J$ фундамента ${\varPhi}$ будет называтьсяя {\it индуктивным}, если его индексное множество $J$ изоморфно $\N$ (как частично упорядоченное множество). 

Из теоремы \ref{TH:strogij-projlim} следует

\blm\label{LM:ob-indukt-postr-v-LCS}
Пусть ${\varPhi}=(\{\iota_i^j\},\{\varkappa_j^i\})$ --- (индуктивный) фундамент в категории $\LCS$ над направленным множеством $\N$, причем все пространства $X_i$ стереотипны и насыщены. Тогда проективный предел $\projlim_{i\in I} X_i$ в категории $\LCS$ является стереотипным и насыщенным пространством.
\elm

Заметим, что в отличие от диаграммы \eqref{frame->body}, в диаграмме \eqref{frame->body-for-inner-building} вертикальные морфизмы слева --- $\sigma_J$ --- не обязаны быть коретракциями, а морфизмы справа --- $\pi^J$ --- не обязаны быть ретракциями. Однако в важных случаях они будут обладать этими свойствами, причем внешний обратный для морфизма $\sigma_J$ и внутренний обратный для морфизма $\pi^J$ согласованы с морфизмами $\theta_{{\varPhi}_J}$ и  $\theta_{\varPhi}$. Речь идет вот о каких ситуациях:

\bit{

\item[---] существует морфизм $\tau$, замыкающий диаграмму
\beq\label{complementable-inner-building-1}
 \xymatrix  @R=3pc @C=3pc
 {
\fram\varPhi\ar[rr]^{\theta_{\varPhi}}\ar@{-->}[rd]_{\tau}  & & \edif\varPhi\ar[dd]^{\pi^J} \\
 &  \fram\varPhi_J\ar[rd]^{\theta_{{\varPhi}_J}} & \\
\fram\varPhi_J\ar[ru]^{1_{\fram\varPhi_J}} \ar[uu]^{\sigma_J}\ar[rr]^{\theta_{{\varPhi}_J}} & & \edif\varPhi_J
 }
\eeq 
и тогда мы говорим, что {\it каркас на фундаменте ${\varPhi}_J$ дополняется в каркасе на расширенном фундаменте ${\varPhi}$},

\item[---] существует морфизм $\rho$, замыкающий диаграмму
\beq\label{complementable-inner-building-2}
 \xymatrix  @R=3pc @C=3pc
 {
\fram\varPhi\ar[rr]^{\theta_{\varPhi}}  & & \edif\varPhi\ar[dd]^{\pi^J} \\
 &  \edif\varPhi_J\ar@{-->}[ru]_{\rho}\ar[rd]^{1_{\edif\varPhi_J}} & \\
\fram\varPhi_J \ar[uu]^{\sigma_J}\ar[rr]^{\theta_{{\varPhi}_J}}\ar[ru]^{\theta_{{\varPhi}_J}} & & \edif\varPhi_J
 }
\eeq 
и тогда мы говорим, что {\it здание на фундаменте ${\varPhi}_J$ дополняется в здании на расширенном фундаменте ${\varPhi}$},

\item[---] существуют оба морфизма $\tau$ и $\rho$, замыкающие диаграммы \eqref{complementable-inner-building-1} и \eqref{complementable-inner-building-2}, и тогда для них будет коммутативна диаграмма 
\beq\label{complementable-inner-building}
 \xymatrix  @R=3pc @C=3pc
 {
\fram\varPhi\ar[rrr]^{\theta_{\varPhi}}\ar@{-->}[rd]_{\tau} & & & \edif\varPhi\ar[dd]^{\pi^J} \\
& \fram\varPhi\ar[r]^{\theta_{\varPhi}} &  \edif\varPhi_J\ar@{-->}[ru]_{\rho}\ar[rd]^{1_{\edif\varPhi_J}} & \\
\fram\varPhi_J \ar[uu]^{\sigma_J}\ar[rrr]^{\theta_{{\varPhi}_J}}\ar[ru]^{1_{\fram\varPhi_J}} & & & \edif\varPhi_J
 }
\eeq 
и в этом случае мы говорим, что {\it фундамент ${\varPhi}_J$ дополняется в расширенном фундаменте ${\varPhi}$}.

}\eit

\paragraph{Проекторы здания.}

Для данного фундамента ${\varPhi}$ определим семейство морфизмов $P^i:\edif\varPhi\to \edif\varPhi$ формулой 
\beq\label{DEF:proektory-zdanija}
P^i=\theta_{\varPhi}\circ\sigma_i\circ\pi^i, \qquad 
 \xymatrix  @R=2pc @C=3pc
 {
 \fram\varPhi\ar[r]^{\theta_{\varPhi}} &  \edif\varPhi  \\
 X_i \ar[u]^{\sigma_i} & \edif\varPhi \ar[l]^{\pi^i} \ar@{-->}[u]_{P^i}\\
 }\put(40,-25){$(i\in I)$}
\eeq
Такие морфизмы мы будем называть {\it проекторами здания} $\edif\varPhi$ на фундаменте ${\varPhi}$.

\bigskip

\centerline{\bf Свойства проекторов здания:}

\bit{\it

\item[$1^\circ$.] Каждый морфизм $P^i$ является проектором (или, в другой терминологии, идемпотентом):
$$
P^i\circ P^i=P^i.
$$

\item[$2^\circ$.] В композиции двух проекторов проектор с большим индексом исчезает:   
$$
P^i\circ P^j=P^i=P^j\circ P^i\qquad (i\le j).
$$

\item[$3^\circ$.] Для любых индексов $i\le j$ коммутативна диаграмма 
\beq\label{P^i-theta_B-sigmaj}
 \xymatrix  @R=2pc @C=3pc
 {
 X_i\ar[r]^{\iota_i^j} & X_j \\ 
 \edif\varPhi\ar[u]_{\pi^i}\ar@{-->}[r]^{P^i} &  \edif\varPhi\ar[u]^{\pi^j}  \\
 X_j\ar@/^7ex/[uu]^{\varkappa_j^i}\ar[u]_{\theta_{\varPhi}\circ\sigma_j} \ar[r]_{\varkappa_j^i} & X_i\ar[u]^{\theta_{\varPhi}\circ\sigma_i}\ar@/_7ex/[uu]_{\iota_i^j}\\
 }\put(40,-25){$(i\le j)$}
\eeq

\item[$4^\circ$.] Если ${\varPhi}$ --- фундамент в категории с базисным разложением, то  
\beq\label{Im P^i=Im (theta_varPhi-circ-sigma_i)}
\Im P^i=\Im (\theta_{\varPhi}\circ\sigma_i),\qquad \Ker P^i=\Ker\pi^i
\eeq

}\eit

\bpr

1.
$$
(P^i)^2=\theta_{\varPhi}\circ\sigma_i\circ\underbrace{\pi^i\circ\theta_{\varPhi}\circ\sigma_i}_{\scriptsize\begin{matrix}\| \\ 1_{X_i}\end{matrix}}\circ\pi^i=\theta_{\varPhi}\circ\sigma_i\circ 1_{X_i}\circ \pi^i=\theta_{\varPhi}\circ\sigma_i\circ\pi^i=P^i
$$

2. Пусть $i\le j$. Тогда, во-первых,
$$
P^i\circ P^j=\theta_{\varPhi}\circ\sigma_i\circ\underbrace{\pi^i\circ\theta_{\varPhi}\circ\sigma_j}_{\scriptsize\begin{matrix} \|\put(3,0){\eqref{building-morphism-3-2}} \\ \varkappa_j^i\end{matrix}}\circ\pi^j=\theta_{\varPhi}\circ\sigma_i\circ \underbrace{\varkappa_j^i\circ \pi^j}_{\scriptsize\begin{matrix}\|\put(3,0){\eqref{DEF:frame-body}} \\ \pi^i\end{matrix}}=\theta_{\varPhi}\circ\sigma_i\circ\pi^i=P^i
$$ 
и, во-вторых,
$$
P^j\circ P^i=\theta_{\varPhi}\circ\sigma_j\circ\underbrace{\pi^j\circ\theta_{\varPhi}\circ\sigma_i}_{\scriptsize\begin{matrix} \|\put(3,0){\eqref{building-morphism-3-1}} \\ \iota_i^j\end{matrix}}\circ\pi^i=\theta_{\varPhi}\circ\underbrace{\sigma_j\circ \iota_i^j}_{\scriptsize\begin{matrix}\|\put(3,0){\eqref{DEF:frame-body}} \\ \sigma_i\end{matrix}}\circ \pi^i=\theta_{\varPhi}\circ\sigma_i\circ\pi^i=P^i
$$

3. В диаграмме \eqref{P^i-theta_B-sigmaj} внутренние треугольники по бокам --- это диаграммы \eqref{building-morphism-3-2} и \eqref{building-morphism-3-1}. Коммутативность верхнего внутреннего прямоугольника в \eqref{P^i-theta_B-sigmaj} доказывается цепочкой
$$
\pi^j\circ P^i=\underbrace{\pi^j\circ \theta_{\varPhi}\circ\sigma_i}_{\scriptsize \begin{matrix}\|\put(3,0){\eqref{building-morphism-3-1}} \\ \iota_i^j \end{matrix}}\circ\pi^i
=\iota_i^j \circ\pi^i
$$
а нижнего --- цепочкой
$$
P^i\circ\theta_{\varPhi}\circ\sigma_j=\theta_{\varPhi}\circ\sigma_i\circ\underbrace{\pi^i\circ\theta_{\varPhi}\circ\sigma_j}_{\scriptsize \begin{matrix}\|\put(3,0){\eqref{building-morphism-3-2}} \\ \varkappa_j^i \end{matrix}}
=
\theta_{\varPhi}\circ\sigma_i\circ\varkappa_j^i
$$

4. Равенства \eqref{Im P^i=Im (theta_varPhi-circ-sigma_i)} следуют из предложения \ref{PROP:pi^i-circ-theta_varPhi-retraktsija}, по которому $\pi^i$ --- ретракция, а $\theta_{\varPhi}\circ\sigma_i$ --- коретракция.
\epr

\paragraph{Проекторы каркаса.}

Для данного фундамента ${\varPhi}$ определим семейство морфизмов $Q^i:\fram\varPhi\to \fram\varPhi$ формулой 
\beq\label{DEF:proektory-karkasa}
Q^i=\sigma_i\circ\pi^i\circ \theta_{\varPhi}, \qquad 
 \xymatrix  @R=2pc @C=3pc
 {
 \fram\varPhi\ar@{-->}[d]_{Q^i}\ar[r]^{\theta_{\varPhi}} &  \edif\varPhi\ar[d]^{\pi^i}  \\
 \fram\varPhi  & X_i \ar[l]^{\sigma_i} \\
 }\put(40,-25){$(i\in I)$}
\eeq
Такие морфизмы мы будем называть {\it проекторами каркаса} $\fram\varPhi$ на фундаменте ${\varPhi}$.

\bigskip

\centerline{\bf Свойства проекторов каркаса:}

\bit{\it

\item[$1^\circ$.] Каждый морфизм $Q^i$ является проектором (или, в другой терминологии, идемпотентом):
$$
Q^i\circ Q^i=Q^i.
$$

\item[$2^\circ$.] В композиции двух проекторов проектор с большим индексом исчезает:   
$$
Q^i\circ Q^j=Q^i=Q^j\circ Q^i\qquad (i\le j).
$$

\item[$3^\circ$.] Для любых индексов $i\le j$ коммутативна диаграмма 
\beq\label{pi^j-theta_B-Q^i}
 \xymatrix  @R=2pc @C=3pc
 {
 X_i\ar[r]^{\iota_i^j} & X_j \\ 
 \fram\varPhi\ar[u]_{\pi^i\circ\theta_{\varPhi}}\ar@{-->}[r]^{Q^i} &  \fram\varPhi\ar[u]^{\pi^j\circ\theta_{\varPhi}}  \\
 X_j\ar@/^7ex/[uu]^{\varkappa_j^i}\ar[u]_{\sigma_j} \ar[r]_{\varkappa_j^i} & X_i\ar[u]^{\sigma_i}\ar@/_7ex/[uu]_{\iota_i^j}\\
 }\put(40,-25){$(i\le j)$}
\eeq

\item[$4^\circ$.] Если ${\varPhi}$ --- фундамент в категории с базисным разложением, то  
\beq\label{Im Q^i=Im-sigma_i}
\Im Q^i=\Im \sigma_i,\qquad \Ker Q^i=\Ker(\pi^i\circ\theta_{\varPhi})
\eeq

}\eit

\bpr

1.
$$
(Q^i)^2=\sigma_i\circ\underbrace{\pi^i\circ\theta_{\varPhi}\circ\sigma_i}_{\scriptsize\begin{matrix}\| \\ 1_{X_i}\end{matrix}}\circ\pi^i\circ\theta_{\varPhi}=\sigma_i\circ 1_{X_i}\circ \pi^i\circ\theta_{\varPhi}=\sigma_i\circ\pi^i\circ\theta_{\varPhi}=Q^i
$$

2. Пусть $i\le j$. Тогда, во-первых,
$$
Q^i\circ Q^j=\sigma_i\circ\underbrace{\pi^i\circ\theta_{\varPhi}\circ\sigma_j}_{\scriptsize\begin{matrix} \|\put(3,0){\eqref{building-morphism-3-2}} \\ \varkappa_j^i\end{matrix}}\circ\pi^j\circ\theta_{\varPhi}=\sigma_i\circ \underbrace{\varkappa_j^i\circ \pi^j}_{\scriptsize\begin{matrix}\|\put(3,0){\eqref{DEF:frame-body}} \\ \pi^i\end{matrix}}\circ\theta_{\varPhi}=\sigma_i\circ\pi^i\circ\theta_{\varPhi}=Q^i
$$ 
и, во-вторых,
$$
Q^j\circ Q^i=\sigma_j\circ\underbrace{\pi^j\circ\theta_{\varPhi}\circ\sigma_i}_{\scriptsize\begin{matrix} \|\put(3,0){\eqref{building-morphism-3-1}} \\ \iota_i^j\end{matrix}}\circ\pi^i\circ\theta_{\varPhi}=\underbrace{\sigma_j\circ \iota_i^j}_{\scriptsize\begin{matrix}\|\put(3,0){\eqref{DEF:frame-body}} \\ \sigma_i\end{matrix}}\circ \pi^i\circ\theta_{\varPhi}=\sigma_i\circ\pi^i\circ\theta_{\varPhi}=Q^i
$$

3. В диаграмме \eqref{pi^j-theta_B-Q^i} внутренние треугольники по бокам --- это диаграммы \eqref{building-morphism-3-2} и \eqref{building-morphism-3-1}. Коммутативность верхнего внутреннего прямоугольника в \eqref{pi^j-theta_B-Q^i} доказывается цепочкой
$$
\pi^j\circ \theta_{\varPhi}\circ Q^i=\underbrace{\pi^j\circ \theta_{\varPhi}\circ\sigma_i}_{\scriptsize \begin{matrix}\|\put(3,0){\eqref{building-morphism-3-1}} \\ \iota_i^j \end{matrix}}\circ\pi^i\circ \theta_{\varPhi}
=\iota_i^j \circ\pi^i\circ \theta_{\varPhi}
$$
а нижнего --- цепочкой
$$
Q^i\circ\sigma_j=\sigma_i\circ\underbrace{\pi^i\circ\theta_{\varPhi}\circ\sigma_j}_{\scriptsize \begin{matrix}\|\put(3,0){\eqref{building-morphism-3-2}} \\ \varkappa_j^i \end{matrix}}
=
\sigma_i\circ\varkappa_j^i
$$

4. Равенства \eqref{Im Q^i=Im-sigma_i} следуют из предложения \ref{PROP:pi^i-circ-theta_varPhi-retraktsija}, по которому $\pi^i\circ\theta_{\varPhi}$ --- ретракция, а $\sigma_i$ --- коретракция.
\epr

\subsection{Постройки в категориях $\LCS$ и $\Ste$}

\paragraph{Теоремы о наследовании стереотипности.}

\btm[о наследовании стереотипности зданием]\label{TH:o-postr-v-LCS}
Пусть ${\varPhi}=(\{\iota_i^j\},\{\varkappa_j^i\})$ --- фундамент в категории $\LCS$ над направленным множеством $I$, причем 
\bit{

\item[(i)] все пространства $X_i$ полны и насыщены (и поэтому стереотипны),

\item[(ii)] здание $\edif\varPhi_J$ на всяком индуктивном ограничении ${\varPhi}_J\subseteq {\varPhi}$ фундамента ${\varPhi}$ дополняется\footnote{В смысле определения на с.\pageref{complementable-inner-building-2}.} в здании $\edif\varPhi$ на фундаменте ${\varPhi}$ в категории $\LCS$.

}\eit\noindent
Тогда 
\bit{

\item[(a)] здание $\edif\varPhi$ на фундаменте ${\varPhi}$ в категории $\LCS$ является полным и насыщенным (и поэтому стереотипным) пространством, и

\item[(b)] морфизм $\theta_\varPhi$ плотно отображает каркас $\fram\varPhi$ в здание $\edif\varPhi$. 
}\eit
\etm

\brem
В формулировке этой теоремы можно всюду заменить требование полноты требованием псевдополноты.
\erem

\bpr
Пусть $X=\edif\varPhi=\LCS\text{-}\projlim_{i\in I} X_i$ --- здание на фундаменте ${\varPhi}$, и $P^i$ --- проекторы, определенные формулой 
\eqref{DEF:proektory-zdanija}.

1. Сначала нужно заметить, что для всякого $x\in X$ направленность проекций $\{ P^k(x)\}$ стремится к $x$ в $X$:
\beq\label{(theta-circ-sigma_k-circ-pi_k)(x)->x}
 P^k(x)\overset{X}{\underset{k\to\infty}{\longrightarrow}}x
\eeq
Зафискируем окрестность нуля $U$. Поскольку $X$ --- подпространство в произведении $\prod_{i\in  I} X_i$, найдутся
индексы $i_1,...,i_n$ такие что
$$
\bigcap_{s=1}^n\Ker\pi_{i_s}\subseteq U
$$
Подберем индекс $k\in I$ так чтобы
$$
\forall s=1,...,n\quad k\ge i_s.
$$
Тогда если $y\in\Ker\pi_k$ то для любого $s=1,...,n$ мы получим
$$
\pi_{i_s}(y)=\varkappa_{i_s}^k(\pi_k(y))=0
$$
и это значит, что 
$$
\Ker\pi_k\subseteq \bigcap_{s=1}^n\Ker\pi_{i_s}\subseteq U
$$
Более того, для всякого $l\ge k$ мы получим
$$
\Ker\pi_l\subseteq \Ker\pi_k\subseteq \bigcap_{s=1}^n\Ker\pi_{i_s}\subseteq U
$$
Теперь при $l\ge k$ мы получаем цепочку:
$$
\pi_l( P^l(x)-x)=(\pi_l\circ\theta\circ\sigma_l)(\pi_l(x))-\pi_l(x)=\pi_l(x)-\pi_l(x)=0
$$
$$
\Downarrow
$$
$$
 P^l(x)-x\in \Ker\pi_l\subseteq U
$$
$$
\Downarrow
$$
$$
 P^l(x)\in x+U
$$
Поскольку $U$ выбиралась произвольной, это доказывает \eqref{(theta-circ-sigma_k-circ-pi_k)(x)->x}.

2. Из \eqref{(theta-circ-sigma_k-circ-pi_k)(x)->x} сразу следует условие (b).

3. Далее нужно заметить, что для всякого фиксированного $x\in X$ направленность $\{ P^k(x)\}$ вполне ограничена в $X$. Это следует из того, что в объемлющем пространстве $\prod_{i\in  I} X_i$ проекции на конечные произведения образуют вполне ограниченное семейство.

4. Пусть $F$ --- множество линейных функционалов на $X$ (не предполагаемых априори непрерывными на $X$, но) непрерывных на каждом вполне ограниченном множестве $S\subseteq X$. Нам нужно показать, что $F$ равностепенно непрерывно на $X$.

a) Сначала обозначим
\beq\label{R_k=theta_varPhi(sigma_k(X_k))}
R_k=\theta_{\varPhi}(\sigma_k(X_k)),\qquad N_k=\Ker\pi^k
\eeq
и
\beq\label{R_infty=bigcup_{k-in-I}R_k}
R_\infty=\bigcup_{k\in I}R_k
\eeq
и заметим, что существует индекс $l\in I$ такой, что
\beq\label{f|_R_infty-cap-N_l=0}
\forall f\in F\quad f\Big|_{R_\infty\cap N_l}=0
\eeq
Предположим, что это неверно:
\beq\label{f|_R_infty-cap-N_l=0}
\forall l\in I\quad \exists f\in F\quad \exists x\in R_\infty\cap N_l \quad f(x)=1.
\eeq
Тогда существуют последовательности $k_i\in I$, $f_i\in F$, $x_i\in R_\infty$ со свойствами
\beq\label{f(x_i)=1*}
k_i\le k_{i+1},\qquad f_i\in F,\qquad x_i\in R_{k_{i+1}}\cap N_{k_i},\qquad f_i(x_i)=1.
\eeq
Положим
$$
J=\{k_i;\ i\in\N\}
$$
и рассмотрим ограничение ${\varPhi}_J\subseteq {\varPhi}$ фундамента $\varPhi$. Здание на нем 
$$
Y=\LCS\text{-}\projlim_{i\to\infty} R_{k_i}
$$
по условию (b) нашей теоремы, можно считать подпространством в здании $X$ на фундаменте ${\varPhi}$,
$$
Y\subseteq X,
$$
причем с индуцированной из $X$ топологией, и, что важнее, формально более слабой топологией, порожденной всевозможными отображениями $$ 
\pi_{k_i}|_Y:Y\to X_{k_i}.
$$

С другой стороны, по лемме \ref{LM:ob-indukt-postr-v-LCS} $Y$ --- насыщенное пространство. А система функционалов $F|_Y$ равностепенно непрерывна на вполне ограниченных множествах в $Y$. Значит, функционалы $F|_Y$ равностепенно непрерывны на $Y$ и поэтому должны обнуляться на ядре $\Ker U$ какой-то окрестности нуля $U$ в пространстве $Y$: 
$$
F|_Y|_{\Ker U}=0.
$$
Но топология пространства $Y=\LCS\text{-}\projlim_{i\to\infty} R_{k_i}$ порождается проекциями $\pi_{k_i}|_Y:Y\to X_{k_i}$, значит, функционалы $F|_Y$ должны обнуляться на ядре $\Ker\pi_{k_i}|_Y$ какого-то оператора $\pi_{k_i}|_Y$:
$$
F|_Y|_{\Ker(\pi_{k_i}|_Y)}=0.
$$
Упростим эту запись,
$$
F|_{Y\cap \Ker(\pi_{k_i}|_Y)}=0,
$$
и заметим, что множество, на которое ограничиваются функционалы $F$, содержит множество $R_{k_{i+1}}\cap N_{k_i}$: из коммутативности верхнего внутреннего четырехугольника в \eqref{complementable-inner-building-2} следует
$$
Y\supseteq (\theta_{\varPhi}\circ\sigma_{k_{i+1}})(X_{k_{i+1}})=R_{k_{i+1}}, 
$$
откуда 
$$
Y\cap \Ker(\pi_{k_i}|_Y)=Y\cap\{x\in Y:\ \pi_{k_i}(x)=0\}=Y\cap Y\cap \Ker\pi_{k_i}=Y\cap N_{k_i}\supseteq R_{k_{i+1}}\cap N_{k_i}.
$$
В результате, мы получаем, что функционалы $F$ должны обнуляться и на множестве $R_{k_{i+1}}\cap N_{k_i}$:
$$
F|_{R_{k_{i+1}}\cap N_{k_i}}=0,
$$ 
А это противоречит \eqref{f(x_i)=1*}. Значит, \eqref{f|_R_infty-cap-N_l=0} не может быть верно.

b) Теперь покажем, что функционалы $F$ обнуляются на некотором пространстве $N_l$:
\beq\label{F|_N_l=0*}
\exists l\in I\quad \forall f\in F\quad f\Big|_{N_l}=0
\eeq
Пусть $l$ --- индекс из \eqref{f|_R_infty-cap-N_l=0}, и $x\in N_l$. Тогда в силу \eqref{(theta-circ-sigma_k-circ-pi_k)(x)->x}, 
$$
 P^k(x)\overset{X}{\underset{k\to 0}{\longrightarrow}} x
$$ 
причем в силу пункта 2, множество $\{x\}\cup \{ P^k(x);\ k\in I\}$ вполне ограничено в $X$. С другой стороны, в силу \eqref{P_k(N_l)-subseteq-N_l}, условие $x\in N_l$ влечет условие $ P^k(x)\in N_l$, и поэтому
$$
\forall f\in F\quad f( P^k(x))=0.
$$
Теперь, поскольку функционалы $f\in F$ непрерывны на вполне ограниченном множестве $\{x\}\cup \{ P^k(x);\ k\in I\}$, мы получаем:
$$
0=f( P^k(x))\overset{X}{\underset{k\to 0}{\longrightarrow}} f(x)\quad\Rightarrow\quad f(x)=0.
$$  
И это доказывает \eqref{F|_N_l=0*}.

c) Из условия \eqref{F|_N_l=0*} следует, что при разложении в прямую сумму
$$
X=R_l\oplus N_l
$$
каждый функционал $f\in F$ обнуляется на подпространстве $N_l=\Ker P^l=\Ker\pi_l$ (определяясь своим ограничением $f\big|_{R_l}$ на подпространство $R_l= P^l(X)$). Поэтому он имеет вид 
\beq
f=f\big|_{R_l}\circ  P^l
\eeq
Поскольку $R_l$ насыщено, система ограничений $\{f\big|_{R_l};\ f\in F\}$ (по наследству от $F$, равностепенно непрерывная на вполне ограниченных множествах в $R_l$) равностепенно непрерывна на $R_l$. Значит, система функционалов
$$
F=\{f\big|_{R_l}\circ  P^l;\ f\in F\}
$$
равностепенно непрерывна на $X$.  
\epr

\btm[о наследовании стереотипности каркасом]\label{TH:o-nasl-ster-karkasom}
Пусть ${\varPhi}=(\{\iota_i^j\},\{\varkappa_j^i\})$ --- фундамент в категории $\LCS$ над направленным множеством $I$, причем 
\bit{

\item[(i)] все пространства $X_i$ полны и насыщены (и поэтому стереотипны),

\item[(ii)] каркас $\fram\varPhi_J$ на всяком индуктивном ограничении ${\varPhi}_J\subseteq {\varPhi}$ фундамента ${\varPhi}$ дополняется\footnote{В смысле определения на с.\pageref{complementable-inner-building-1}.} в каркасе $\fram\varPhi$ на фундаменте ${\varPhi}$ в категории $\LCS$.

}\eit\noindent
Тогда 
\bit{

\item[(a)] каркас $\fram\varPhi$ на фундаменте ${\varPhi}$ в категории $\LCS$ является полным и насыщенным (и поэтому стереотипным) пространством, и

\item[(b)] морфизм $\theta_\varPhi$ инъективно отображает каркас $\fram\varPhi$ в здание $\edif\varPhi$. 
}\eit
\etm

\brem
В формулировке этой теоремы можно всюду заменить требование полноты требованием псевдополноты.
\erem
\bpr
1. Сначала заметим, что всякое ограниченное множество $B\subseteq \fram\varPhi=\injlim_{i\in I}X_i$ содержится в некотором $X_j$. Предположим, что это не так, то есть для всякого $X_i$ найдется $x\in B$ такой, что $x\notin X_i$. Тогда можно выбрать последовательности $\{x_k\}$ и $\{X_{i_k}\}$ такие что 
$$
i_k\le i_{k+1},\quad x_k\in B\cap X_{i_{k+1}},\quad x_k\notin X_{i_k},\qquad k\in \N.
$$ 
Это можно понимать так, что у нас имеется ограниченная последовательность $\{x_k\}$ в строгом инъективном пределе
$$
\LCS\text{-}\injlim_{k\in\N}X_{i_k}
$$
не содержащаяся ни в каком $X_{i_k}$. Это противоречит теореме \ref{TH:strogij-injlim-ogr}.

2. Из пункта 1 следует, что пространство $\fram\varPhi=\injlim_{i\in I}X_i$ псевдополно (потому что все $X_i$ псевдополны). С другой стороны, инъективный предел наследует свойство насыщенности, поэтому пространство $\fram\varPhi=\injlim_{i\in I}X_i$ насыщено. Для нас важно, что вместе это означает стереотипность пространства $\fram\varPhi=\injlim_{i\in I}X_i$.

3. Теперь если перейти к сопряженному пространству, то мы получим, что оно должно быть (тоже стереотипно и) зданием сопряженного фундамента ${\varPhi}^\star=(\{(\varkappa_j^i)^\star\},\{(\iota_i^j)^\star\})$:
$$
(\fram\varPhi)^\star=\l\LCS\text{-}\injlim_{i\in I}X_i\r^\star=\l\Ste\text{-}\injlim_{i\in I}X_i\r^\star=
\Ste\text{-}\projlim_{i\in I}X_i^\star=\LCS\text{-}\projlim_{i\in I}X_i^\star=\edif(\varPhi^\star)
$$
А по теореме \ref{TH:o-postr-v-LCS}, $\edif(\varPhi^\star)$ полно и насыщено. Значит, сопряженное ему пространство $\fram\varPhi$ должно быть насыщено (это мы уже отмечали) и полно.

4.  Наконец, по теореме \ref{TH:o-postr-v-LCS}, морфизм
$$
\theta_{\varPhi^\star}:\fram(\varPhi^\star) \to\edif(\varPhi^\star)
$$
является эпиморфизмом. Отсюда следует, что сопряженный ему морфизм
$$
\theta_{\varPhi}:\fram\varPhi\to\edif\varPhi
$$
должен быть мономорфизмом, то есть инъекцией.
\epr

\paragraph{Теоремы о тензорных произведениях фундаментов.}

Пусть ${\varPhi}=(\{\iota_i^j\},\{\varkappa_j^i\})$ и ${\varPsi}=(\{\lambda_p^q\},\{\mu_q^p\})$ --- два фундамента в категории $\Ste$ над направленными множествами $I$ и  $P$ соответственно. Рассмотрим фнудамент 
$$
\varPhi\odot\varPsi,
$$ 
определенный следующим образом:
\bit{
\item[1)] $\varPhi\odot\varPsi$ считается фундаментом над частично упорядоченным множеством
$$
I\times P,
$$
в котором порядок задается правилом
$$
(i,p)\le (j,q)\quad\Leftrightarrow\quad i\le j\ \&\ p\le q. 
$$ 

\item[2)] ковариантная система в $\varPhi\odot\varPsi$ определяется как система морфизмов
$$
\iota_i^j\odot\lambda_p^q: X_i\odot Y_k\to X_j\odot Y_l,\qquad (i,p)\le (j,q),
$$

\item[3)] контравариантная система в $\varPhi\odot\varPsi$ определяется как система морфизмов
$$
\varkappa_j^i\odot\mu_q^p: X_i\odot Y_p\gets X_j\odot Y_q,\qquad (i,p)\le (j,q),
$$

\item[4)] тогда для всякой пары  $(i,p)\le (j,q)$ морфизм $\iota_i^j\odot\lambda_p^q$ будет коретракцией для морфизма 
$\varkappa_j^i\odot\mu_q^p$, потому что
$$
(\varkappa_j^i\odot\mu_q^p) \circ (\iota_i^j\odot\lambda_p^q)=
(\varkappa_j^i\circ\iota_i^j)\odot(\mu_q^p\circ\lambda_p^q)=1_{X_i}\odot 1_{Y_p}=1_{X_i\odot Y_p}.
$$

}\eit

Диаграмма \eqref{DEF:diagr-building} при этом преобразуется в диаграмму
\beq\label{DEF:diagr-building-Phi-odot-Psi}
 \xymatrix  @R=3pc @C=3pc
 {
X_k\odot Y_r\ar@{-->}[r]^{1_{X_k\odot Y_r}}
 &  X_k\odot Y_r\ar@/^3ex/[d]_{\varkappa_k^j\odot\mu_r^q} \ar@/^12ex/[dd]^{\varkappa_k^i\odot\mu_r^p} \\
X_j\odot Y_q\ar@{-->}[r]^{1_{X_j\odot Y_q}}\ar@/^3ex/[u]_{\iota_j^k\odot\lambda_q^r}  &  X_j\odot Y_q \ar@/^3ex/[d]_{\varkappa_j^i\odot\mu_q^p} \\
X_i\odot Y_p\ar@{-->}[r]^{1_{X_i\odot Y_p}} \ar@/^12ex/[uu]^{\iota_i^k\odot\lambda_p^r} \ar@/^3ex/[u]_{\iota_i^j\odot\lambda_p^q}  & X_i\odot Y_p
 }\put(20,-50){$((i,p)\le (j,q)\le (k,r))$}
\eeq

\btm[о наследовании стереотипности слабым тензорным произведением]\label{TH:o-nasl-ster-odot}
Пусть ${\varPhi}=(\{\iota_i^j\},\{\varkappa_j^i\})$ и ${\varPsi}=(\{\lambda_p^q\},\{\mu_q^p\})$ --- два фундамента в категории $\Ste$ над направленными множествами $I$ и  $P$ соответственно, причем 
\bit{

\item[(i)] все пространства $X_i\odot Y_p$ полны и насыщены,

\item[(ii)] здание $\edif\varPhi_J$ на всяком индуктивном ограничении ${\varPhi}_J\subseteq {\varPhi}$ фундамента ${\varPhi}$ дополняется\footnote{В смысле определения на с.\pageref{complementable-inner-building-2}.} в здании $\edif\varPhi$ на фундаменте ${\varPhi}$ в категории $\Ste$.

\item[(iii)] здание $\edif\varPsi_Q$ на всяком индуктивном ограничении ${\varPsi}_Q\subseteq {\varPsi}$ фундамента ${\varPsi}$ дополняется в здании $\edif\varPsi$ на фундаменте ${\varPsi}$ в категории $\Ste$.

}\eit\noindent
Тогда 
\bit{

\item[(a)] здание $\edif^{\Ste}(\varPhi\odot\varPsi)$ на фундаменте $\varPhi\odot\varPsi$ в категории $\Ste$ совпадает со зданием $\edif^{\LCS}(\varPhi\odot\varPsi)$ на этом фундаменте в категории $\LCS$ 
\beq\label{TH:o-nasl-ster-odot}
\edif^{\Ste}(\varPhi\odot\varPsi)=\Ste\text{-}\projlim_{(i,p)\in I\times P}X_i\odot Y_p=
\LCS\text{-}\projlim_{(i,p)\in I\times P}X_i\otimes Y_p=\edif^{\LCS}(\varPhi\odot\varPsi)
\eeq    
    и при этом является полным и насыщенным пространством, и

\item[(b)] морфизм $\theta_{\varPhi\odot\varPsi}$ плотно отображает каркас $\fram(\varPhi\odot\varPsi)$ в здание $\edif(\varPhi\odot\varPsi)$. 
}\eit
\etm
\bpr
Здесь действует теорема \ref{TH:o-postr-v-LCS}, только нужно проверить условие (ii) этой теоремы. Пусть $S\subseteq I\times P$ --- индуктивное подмножество (то есть изомофрное $\N$ как частично упорядоченное множество). Обозначим
$$
J=\{j\in I: \quad \exists q\in P\quad (j,q)\in S\},\qquad Q=\{q\in P: \quad \exists j\in I\quad (j,q)\in S\}.
$$
Эти множества тоже будут изоморфны $\N$. Поэтому здания на их фундаментах  $\edif^{\Ste}\varPhi_J$ и $\edif^{\Ste}\varPsi_Q$ дополняемы в зданиях $\edif^{\Ste}\varPhi$ и $\edif^{\Ste}\varPsi$, то есть найдутся морфизмы $\rho_\varPhi$ и $\rho_\varPsi$ такие, что коммутативны диаграммы \eqref{complementable-inner-building-2}:
\beq\label{PROOF:o-nasl-ster-odot-1}
 \xymatrix  @R=3pc @C=2pc
 {
\fram^{\Ste}\varPhi\ar[rr]^{\theta_{\varPhi}}  & & \edif^{\Ste}\varPhi\ar[dd]^{\pi^J} \\
 &  \edif^{\Ste}\varPhi_J\ar@{-->}[ru]_{\rho_\varPhi}\ar[rd]^{1_{\edif^{\Ste}\varPhi_J}} & \\
\fram^{\Ste}\varPhi_J \ar[uu]^{\sigma_J}\ar[rr]^{\theta_{{\varPhi}_J}}\ar[ru]^{\theta_{{\varPhi}_J}} & & \edif^{\Ste}\varPhi_J
 }
\qquad
 \xymatrix  @R=3pc @C=2pc
 {
\fram^{\Ste}\varPsi\ar[rr]^{\theta_{\varPsi}}  & & \edif^{\Ste}\varPsi\ar[dd]^{\pi^Q} \\
 &  \edif^{\Ste}\varPsi_Q\ar@{-->}[ru]_{\rho_\varPsi}\ar[rd]^{1_{\edif^{\Ste}\varPsi_Q}} & \\
\fram^{\Ste}\varPsi_Q \ar[uu]^{\sigma_Q}\ar[rr]^{\theta_{{\varPsi}_Q}}\ar[ru]^{\theta_{{\varPsi}_Q}} & & \edif^{\Ste}\varPsi_Q
 }
\eeq

1. Выделим в диаграммах \eqref{PROOF:o-nasl-ster-odot-1} правые внутренние треугольники:
\beq\label{PROOF:o-nasl-ster-odot-2}
 \xymatrix  @R=3pc @C=2pc
 {
  & & \edif^{\Ste}\varPhi\ar[dd]^{\pi^J} \\
 &  \edif^{\Ste}\varPhi_J\ar@{-->}[ru]^{\rho_\varPhi}\ar[rd]_{1_{\edif^{\Ste}\varPhi_J}} & \\
 & & \edif^{\Ste}\varPhi_J
 }
\qquad
 \xymatrix  @R=3pc @C=2pc
 {
  & & \edif^{\Ste}\varPsi\ar[dd]^{\pi^Q} \\
 &  \edif^{\Ste}\varPsi_Q\ar@{-->}[ru]^{\rho_\varPsi}\ar[rd]_{1_{\edif^{\Ste}\varPsi_Q}} & \\
 & & \edif^{\Ste}\varPsi_Q
 }
\eeq
Перемножив их тензорным произведением $\odot$, мы получим диаграмму
\beq\label{PROOF:o-nasl-ster-odot-3}
 \xymatrix  @R=3pc @C=2pc
 {
  & \edif^{\Ste}\varPhi\odot\edif^{\Ste}\varPsi\ar[dd]^{\pi^J\odot\pi^Q} \\
   \edif^{\Ste}\varPhi_J\odot\edif^{\Ste}\varPsi_Q\ar@{-->}[ru]^{\rho_\varPhi\odot\rho_\varPsi} \ar[rd]_{1_{\edif^{\Ste}\varPhi_J}\odot 1_{\edif^{\Ste}\varPsi_Q}} & \\
 & \edif^{\Ste}\varPhi_J\odot\edif^{\Ste}\varPsi_Q
 }
\eeq
Теперь заметим, что 
\beq\label{PROOF:o-nasl-ster-odot-4}
\edif^{\Ste}\varPhi\odot\edif^{\Ste}\varPsi=\Ste\text{-}\projlim_{i\in I}X_i\odot \Ste\text{-}\projlim_{p\in P}Y_p=
\Ste\text{-}\projlim_{(i,p)\in I\times P}X_i\odot Y_p=\edif^{\Ste}(\varPhi\odot\varPsi),
\eeq
\begin{multline}\label{PROOF:o-nasl-ster-odot-5}
\edif^{\Ste}\varPhi_J\odot\edif^{\Ste}\varPsi_Q=\Ste\text{-}\projlim_{j\in J}X_j\odot \Ste\text{-}\projlim_{q\in Q}Y_q=
\Ste\text{-}\projlim_{(j,q)\in J\times Q}X_j\odot Y_q=\\=
(\text{множество $S\subseteq J\times Q$ конфинально в $J\times Q$})=
\Ste\text{-}\projlim_{(j,q)\in S}X_j\odot Y_q=\edif^{\Ste}(\varPhi\odot\varPsi)_S,
\end{multline}
и
\beq\label{PROOF:o-nasl-ster-odot-6}
1_{\edif^{\Ste}\varPhi_J}\odot 1_{\edif^{\Ste}\varPsi_Q}=
1_{\edif^{\Ste}\varPhi_J\odot \edif^{\Ste}\varPsi_Q}=1_{\edif^{\Ste}(\varPhi\odot\varPsi)_S}
\eeq
Поэтому диаграмму \eqref{PROOF:o-nasl-ster-odot-3} можно интерпретировать как диаграмму
\beq\label{PROOF:o-nasl-ster-odot-7}
 \xymatrix  @R=3pc @C=2pc
 {
   & \edif^{\Ste}(\varPhi\odot\varPsi)\ar[dd]^{\pi^J\odot\pi^Q} \\
   \edif^{\Ste}(\varPhi\odot\varPsi)_S \ar@{-->}[ru]^{\rho_\varPhi\odot\rho_\varPsi} \ar[rd]_{1_{\edif^{\Ste}(\varPhi\odot\varPsi)_S}} & \\
  & \edif^{\Ste}(\varPhi\odot\varPsi)_S
 }
\eeq

2. Выделим в диаграммах \eqref{PROOF:o-nasl-ster-odot-1} левые внутренние четырехугольники:
\beq\label{PROOF:o-nasl-ster-odot-8}
 \xymatrix  @R=3pc @C=2pc
 {
\fram^{\Ste}\varPhi\ar[rr]^{\theta_{\varPhi}}  & & \edif^{\Ste}\varPhi \\
\fram^{\Ste}\varPhi_J \ar[u]^{\sigma_J}\ar[rr]^{\theta_{{\varPhi}_J}} & & \edif^{\Ste}\varPhi_J\ar@{-->}[u]_{\rho_\varPhi}
 }
\qquad
 \xymatrix  @R=3pc @C=2pc
 {
\fram^{\Ste}\varPsi\ar[rr]^{\theta_{\varPsi}}  & & \edif^{\Ste}\varPsi \\
\fram^{\Ste}\varPsi_Q \ar[u]^{\sigma_Q}\ar[rr]^{\theta_{{\varPsi}_Q}} & & \edif^{\Ste}\varPsi_Q\ar@{-->}[u]_{\rho_\varPsi}
 }
\eeq
Зафискируем $j\in J$ и $i\in  I$ так, чтобы $j\le i$, а также индексы  $q\in Q$ и $p\in P$ так, чтобы $q\le p$, а затем дополним диаграммы \eqref{PROOF:o-nasl-ster-odot-8} до диаграмм
\beq\label{PROOF:o-nasl-ster-odot-9}
 \xymatrix  @R=3pc @C=2pc
 {
X_i\ar@{-->}[r]^{\sigma_i} & \fram^{\Ste}\varPhi\ar[r]^{\theta_{\varPhi}}   & \edif^{\Ste}\varPhi \\
X_j\ar@{-->}[r]^{\sigma_i^J}\ar@{-->}[u]^{\iota_j^i} & \fram^{\Ste}\varPhi_J \ar[u]^{\sigma_J}\ar[r]^{\theta_{{\varPhi}_J}}  & \edif^{\Ste}\varPhi_J\ar[u]_{\rho_\varPhi}
 }
\qquad
 \xymatrix  @R=3pc @C=2pc
 {
Y_p\ar@{-->}[r]^{\sigma_p} &\fram^{\Ste}\varPsi\ar[r]^{\theta_{\varPsi}} & \edif^{\Ste}\varPsi \\
Y_q\ar@{-->}[r]^{\sigma_q^J}\ar@{-->}[u]^{\lambda_q^p} &\fram^{\Ste}\varPsi_Q \ar[u]^{\sigma_Q}\ar[r]^{\theta_{{\varPsi}_Q}} & \edif^{\Ste}\varPsi_Q\ar[u]_{\rho_\varPsi}
 }
\eeq
и выбросим средние столбцы:
\beq\label{PROOF:o-nasl-ster-odot-10}
 \xymatrix  @R=3pc @C=2pc
 {
X_i\ar[rr]^{\theta_{\varPhi}\circ\sigma_i}  & & \edif^{\Ste}\varPhi \\
X_j\ar[rr]^{\theta_{{\varPhi}_J}\circ\sigma_i^J}\ar[u]^{\iota_j^i} & & \edif^{\Ste}\varPhi_J\ar[u]_{\rho_\varPhi}
 }
\qquad
 \xymatrix  @R=3pc @C=2pc
 {
Y_p\ar[rr]^{\theta_{\varPsi}\circ\sigma_p}  & & \edif^{\Ste}\varPsi \\
Y_q\ar[rr]^{\theta_{{\varPsi}_Q}\circ\sigma_q^J}\ar[u]^{\lambda_q^p} & & \edif^{\Ste}\varPsi_Q\ar[u]_{\rho_\varPsi}
 }
\eeq
Перемножив эти диаграммы тензорным произведением $\odot$, мы получим диаграмму
\beq\label{PROOF:o-nasl-ster-odot-11}
 \xymatrix  @R=3pc @C=4pc
 {
X_i\odot Y_p\ar[rr]^{\theta_{\varPhi}\circ\sigma_i\odot \theta_{\varPsi}\circ\sigma_p}  & & \edif^{\Ste}\varPhi \odot \edif^{\Ste}\varPsi \\
X_j\odot Y_q\ar[rr]^{\theta_{{\varPhi}_J}\circ\sigma_i^J\odot \theta_{{\varPsi}_Q}\circ\sigma_q^J}\ar[u]^{\iota_j^i\odot \lambda_q^p} & & \edif^{\Ste}\varPhi_J
\odot \edif^{\Ste}\varPsi_Q \ar[u]_{\rho_\varPhi\odot \rho_\varPsi}
 }
\eeq
Применив \eqref{PROOF:o-nasl-ster-odot-4} и \eqref{PROOF:o-nasl-ster-odot-5}, мы получим
\beq\label{PROOF:o-nasl-ster-odot-12}
 \xymatrix  @R=3pc @C=4pc
 {
X_i\odot Y_p\ar[rr]^{\theta_{\varPhi}\circ\sigma_i\odot \theta_{\varPsi}\circ\sigma_p}  & & \edif^{\Ste}(\varPhi\odot\varPsi) \\
X_j\odot Y_q\ar[rr]^{\theta_{{\varPhi}_J}\circ\sigma_i^J\odot \theta_{{\varPsi}_Q}\circ\sigma_q^J}\ar[u]^{\iota_j^i\odot \lambda_q^p} & & \edif^{\Ste}(\varPhi\odot\varPsi)_S\ar[u]_{\rho_\varPhi\odot \rho_\varPsi}
 }
\eeq
Теперь заметим, что морфизмы
\beq\label{PROOF:o-nasl-ster-odot-12-1}
 \xymatrix  @R=3pc @C=4pc
 {
X_i\odot Y_p\ar[rr]^{\theta_{\varPhi}\circ\sigma_i\odot \theta_{\varPsi}\circ\sigma_p}  & & \edif^{\Ste}(\varPhi\odot\varPsi) \\
 }
\eeq
образуют инъективный конус: если $i\le k$ и $p\le r$, то 
\beq\label{PROOF:o-nasl-ster-odot-12-2}
 \xymatrix  @R=1pc @C=3pc
 {
X_k\odot Y_r\ar@/^4ex/[rrd]^{\quad\theta_{\varPhi}\circ\sigma_k\odot \theta_{\varPsi}\circ\sigma_r}  & & \\
 & &  \edif^{\Ste}(\varPhi\odot\varPsi) \\
X_i\odot Y_p\ar[uu]^{\iota_i^k\odot \lambda_p^r} \ar@/_4ex/[rru]_{\quad\theta_{\varPhi}\circ\sigma_i\odot \theta_{\varPsi}\circ\sigma_p}  & & \\
 }
\eeq
--- это следует из того, что морфизмы 
\beq\label{PROOF:o-nasl-ster-odot-12-3}
 \xymatrix  @R=1pc @C=3pc
 {
X_i\ar[rr]^{\theta_{\varPhi}\circ\sigma_i}  & & \edif^{\Ste}\varPhi \\
 },
\qquad
 \xymatrix  @R=1pc @C=3pc
 {
Y_p\ar[rr]^{\theta_{\varPsi}\circ\sigma_p}  & & \edif^{\Ste}\varPsi \\
 }
\eeq
образуют инъективные конусы.

Свойство \eqref{PROOF:o-nasl-ster-odot-12-2} в свою очередь означает, что в диаграмме \eqref{PROOF:o-nasl-ster-odot-12} 
мы можем перейти к инъективноум пределу в верхнем левом углу:
\beq\label{PROOF:o-nasl-ster-odot-13}
 \xymatrix  @R=3pc @C=6pc
 {
\injlim_{(i,p)\in I\times P} X_i\odot Y_p\ar[rr]^{\injlim_{(i,p)\in I\times P}\theta_{\varPhi}\circ\sigma_i\odot \theta_{\varPsi}\circ\sigma_p}  & & \edif^{\Ste}(\varPhi\odot\varPsi) \\
X_j\odot Y_q\ar[rr]^{\theta_{{\varPhi}_J}\circ\sigma_i^J\odot \theta_{{\varPsi}_Q}\circ\sigma_q^J}\ar[u]^{\iota_j\odot \lambda_q} & & \edif^{\Ste}(\varPhi\odot\varPsi)_S\ar[u]_{\rho_\varPhi\odot \rho_\varPsi}
 }
\eeq
а потом в нижнгем леврм углу
\beq\label{PROOF:o-nasl-ster-odot-14}
 \xymatrix  @R=3pc @C=6pc
 {
\injlim_{(i,p)\in I\times P} X_i\odot Y_p\ar[rr]^{\injlim_{(i,p)\in I\times P}\theta_{\varPhi}\circ\sigma_i\odot \theta_{\varPsi}\circ\sigma_p}  & & \edif^{\Ste}(\varPhi\odot\varPsi) \\
\injlim_{(j,q)\in J\times Q}X_j\odot Y_q\ar[rr]^{\injlim_{(j,q)\in J\times Q} \theta_{{\varPhi}_J}\circ\sigma_i^J\odot \theta_{{\varPsi}_Q}\circ\sigma_q^J}\ar[u]^{\iota\odot \lambda} & & \edif^{\Ste}(\varPhi\odot\varPsi)_S\ar[u]_{\rho_\varPhi\odot \rho_\varPsi}
 }
\eeq
Теперь можно заметить, что в левом столбце стоят каркасы: 
\beq\label{PROOF:o-nasl-ster-odot-15}
 \xymatrix  @R=3pc @C=6pc
 {
\fram^{\Ste} (\varPhi\odot\varPsi)\ar[rr]^{\theta_{\varPhi\odot\varPsi}}  & & \edif^{\Ste}(\varPhi\odot\varPsi) \\
\fram^{\Ste} (\varPhi\odot\varPsi)_S\ar[rr]^{\theta_{(\varPhi\odot\varPsi)_S}}\ar[u]^{\sigma_S} & & \edif^{\Ste}(\varPhi\odot\varPsi)_S\ar[u]_{\rho_\varPhi\odot \rho_\varPsi}
 }
\eeq

3. Теперь мы можем объединить \eqref{PROOF:o-nasl-ster-odot-7} и \eqref{PROOF:o-nasl-ster-odot-15} в одну диаграмму:
\beq\label{PROOF:o-nasl-ster-odot-16}
 \xymatrix  @R=3pc @C=2pc
 {
\fram^{\Ste}(\varPhi\odot\varPsi)\ar[rr]^{\theta_{\varPhi\odot\varPsi}^{\Ste}}  & & \edif^{\Ste}(\varPhi\odot\varPsi)\ar[dd]^{\pi^S} \\
 &  \edif^{\Ste}(\varPhi\odot\varPsi)_S\ar@{-->}[ru]_{\rho_{\varPhi\odot\varPsi}}\ar[rd]^{1_{\edif^{\Ste}(\varPhi\odot\varPsi)_S}} & \\
\fram^{\Ste}(\varPhi\odot\varPsi)_S \ar[uu]^{\sigma_S^{\Ste}}\ar[rr]^{\theta_{(\varPhi\odot\varPsi)_S}}\ar[ru]^{\theta_{(\varPhi\odot\varPsi)_S}} & & \edif^{\Ste}(\varPhi\odot\varPsi)_S
 }
\eeq
Теперь заметим, что 
\beq\label{PROOF:o-nasl-ster-odot-17}
\edif^{\Ste}(\varPhi\odot\varPsi)_S=
\Ste\text{-}\projlim_{(i,p)\in S}X_i\odot Y_p=(\text{теорема \ref{TH:strogij-projlim}})=
\LCS\text{-}\projlim_{(i,p)\in S}X_i\odot Y_p=\edif^{\LCS}(\varPhi\odot\varPsi)_S
\eeq
и
\beq\label{PROOF:o-nasl-ster-odot-18}
\fram^{\Ste}(\varPhi\odot\varPsi)_S=
\Ste\text{-}\injlim_{(i,p)\in S}X_i\odot Y_p=(\text{теорема \ref{TH:strogij-injlim}})=
\LCS\text{-}\injlim_{(i,p)\in S}X_i\odot Y_p=\fram^{\LCS}(\varPhi\odot\varPsi)_S
\eeq
Поэтому \eqref{PROOF:o-nasl-ster-odot-16} можно заменить на диаграмму
\beq\label{PROOF:o-nasl-ster-odot-19}
 \xymatrix  @R=3pc @C=2pc
 {
\fram^{\Ste}(\varPhi\odot\varPsi)\ar[rr]^{\theta_{\varPhi\odot\varPsi}^{\Ste}}  & & \edif^{\Ste}(\varPhi\odot\varPsi)\ar[dd]^{\pi^S} \\
 &  \edif^{\LCS}(\varPhi\odot\varPsi)_S\ar@{-->}[ru]_{\rho_{\varPhi\odot\varPsi}}\ar[rd]^{1_{\edif^{\LCS}(\varPhi\odot\varPsi)_S}} & \\
\fram^{\LCS}(\varPhi\odot\varPsi)_S \ar[uu]^{\sigma_S^{\Ste}}\ar[rr]^{\theta_{(\varPhi\odot\varPsi)_S}}\ar[ru]^{\theta_{(\varPhi\odot\varPsi)_S}} & & \edif^{\LCS}(\varPhi\odot\varPsi)_S
 }
\eeq
Теперь мы можем заметить, что 
\beq\label{PROOF:o-nasl-ster-odot-20}
\fram^{\Ste}(\varPhi\odot\varPsi)=\Big(\fram^{\LCS}(\varPhi\odot\varPsi)\Big)^\triangledown
\eeq
а морфизм $\sigma_S$ представляет собой композицию
\beq\label{PROOF:o-nasl-ster-odot-21}
\sigma_S^{\Ste}=\triangledown\circ \sigma_S^{\LCS}
\eeq
Из этого следует, что в \eqref{PROOF:o-nasl-ster-odot-19} можно заменить левую верхниюю вершину на $\fram^{\LCS}(\varPhi\odot\varPsi)$:
\beq\label{PROOF:o-nasl-ster-odot-22}
 \xymatrix  @R=3pc @C=2pc
 {
\fram^{\LCS}(\varPhi\odot\varPsi)\ar[rr]^{\theta_{\varPhi\odot\varPsi}^{\LCS}}  & & \edif^{\Ste}(\varPhi\odot\varPsi)\ar[dd]^{\pi^S} \\
 &  \edif^{\LCS}(\varPhi\odot\varPsi)_S\ar@{-->}[ru]_{\rho_{\varPhi\odot\varPsi}}\ar[rd]^{1_{\edif^{\LCS}(\varPhi\odot\varPsi)_S}} & \\
\fram^{\LCS}(\varPhi\odot\varPsi)_S \ar[uu]^{\sigma_S^{\LCS}}\ar[rr]^{\theta_{(\varPhi\odot\varPsi)_S}}\ar[ru]^{\theta_{(\varPhi\odot\varPsi)_S}} & & \edif^{\LCS}(\varPhi\odot\varPsi)_S
 }
\eeq
Наконец, заметим, что 
\beq\label{PROOF:o-nasl-ster-odot-23}
\edif^{\Ste}(\varPhi\odot\varPsi)=\Big(\edif^{\LCS}(\varPhi\odot\varPsi)\Big)^\vartriangle
\eeq
и при этом 
\bit{
\item[---] отображение $\theta_{\varPhi\odot\varPsi}^{\LCS}$ непрерывно, если его понимать как отображение со значениями в  
$\edif^{\LCS}(\varPhi\odot\varPsi)$, 
\beq\label{PROOF:o-nasl-ster-odot-24}
\theta_{\varPhi\odot\varPsi}^{\LCS}:\fram^{\LCS}(\varPhi\odot\varPsi)\to \edif^{\LCS}(\varPhi\odot\varPsi)
\eeq
потому что при переходе от $\edif^{\Ste}(\varPhi\odot\varPsi)$ к $\edif^{\LCS}(\varPhi\odot\varPsi)$ топология ослабляется,

\item[---] по той же причине отображение $\rho_{\varPhi\odot\varPsi}$ непрерывно, если его понимать как отображение со значениями в  
$\edif^{\LCS}(\varPhi\odot\varPsi)$, 
\beq\label{PROOF:o-nasl-ster-odot-25}
\rho_{\varPhi\odot\varPsi}:\edif^{\LCS}(\varPhi\odot\varPsi)_S\to \edif^{\LCS}(\varPhi\odot\varPsi)
\eeq

\item[---] наконец, отображение $\pi^S$ должно быть непрерывно, если его понимать как отображение из $\edif^{\LCS}(\varPhi\odot\varPsi)$, потому что оно представляет собой естественный морфизм из проективного предела над $I\times P$ в проективный предел над подмножеством $S\subseteq I\times P$:  
\beq\label{PROOF:o-nasl-ster-odot-26}
\pi^S:\edif^{\LCS}(\varPhi\odot\varPsi)\to \edif^{\LCS}(\varPhi\odot\varPsi)_S.
\eeq

}\eit
Вместе все это означает , что мы можем переписать диаграмму \eqref{PROOF:o-nasl-ster-odot-22} в виде
\beq\label{PROOF:o-nasl-ster-odot-27}
 \xymatrix  @R=3pc @C=2pc
 {
\fram^{\LCS}(\varPhi\odot\varPsi)\ar[rr]^{\theta_{\varPhi\odot\varPsi}^{\LCS}}  & & \edif^{\LCS}(\varPhi\odot\varPsi)\ar[dd]^{\pi^S} \\
 &  \edif^{\LCS}(\varPhi\odot\varPsi)_S\ar@{-->}[ru]_{\rho_{\varPhi\odot\varPsi}}\ar[rd]^{1_{\edif^{\LCS}(\varPhi\odot\varPsi)_S}} & \\
\fram^{\LCS}(\varPhi\odot\varPsi)_S \ar[uu]^{\sigma_S^{\LCS}}\ar[rr]^{\theta_{(\varPhi\odot\varPsi)_S}}\ar[ru]^{\theta_{(\varPhi\odot\varPsi)_S}} & & \edif^{\LCS}(\varPhi\odot\varPsi)_S
 }
\eeq
Это и есть условие (ii) теоремы \ref{TH:o-postr-v-LCS}.
\epr

\paragraph{Согласованные проекторы на стереотипном пространстве.}

Рассмотрим следующую ситуацию. Пусть $X$ --- стереотипное пространство, а $( I,\le)$ --- решетка. Условимся семейство  $\{P^k; \ k\in I\}$ (линейных непрерывных) операторов 
$$
P^k:X\to X
$$ 
называть {\it системой согласованных проекторов} на пространстве $X$ над решеткой $( I,\le)$, если для любых индексов $k,l\in I$ выполняется равенство
\beq\label{P_k-circ-P_l=P_k-vee-l}
P^k\circ P^l=P^{k\wedge l},
\eeq
--- здесь $k\wedge l$ --- точная нижняя грань элементов $k$ и $l$ в $I$.

\medskip
\centerline{\bf Свойства системы согласованных проекторов:}

\bit{\it

\item[$1^\circ$.] Для любых двух индексов $k,l\in I$ операторы $P^k$ и $P^l$ коммутируют:
\beq\label{P_k-circ-P_l=P_k-vee-l-1}
P^k\circ P^l=P^{k\wedge l}=P^{l\wedge k}=P^l\circ P^k
\eeq

\item[$2^\circ$.] Для всякого индекса $k\in I$ 
\bit{

\item[(a)] оператор $P^k$ является проектором
$$
(P^k)^2=P^k,
$$

\item[(b)] его ядро 
\beq\label{N_k=Ker-P_k}
N_k=\Ker P^k
\eeq
(замкнуто в $X$ и) является множеством значений дополняющего проектора $\id_X-P^k$:
\beq\label{N_k=(id_X-P_k)(X)}
N_k=(\id_X-P^k)(X)
\eeq

\item[(c)] а его множество значений 
\beq\label{X_k=P_k(X)}
X_k=P^k(X)
\eeq
замкнуто в $X$ и, будучи наделено топологией, индуцированной из $X$, является одновременно кообразом и образом оператора $P^k$ в категории $\LCS$:
$$
 \xymatrix  @R=3pc @C=3pc
 {
X & & X \ar[ll]_{P^k}\ar@/^3ex/[dl]^{\coim P^k} \\
& X_k\ar@/^3ex/[ul]^{\im P^k} &
 }
$$

\item[(d)] пространство $X$ распадается в прямую сумму локально выпуклых пространств $X_k$ и $N_k$:
\beq\label{X=X_k-oplus-N_k}
X=X_k\oplus N_k
\eeq 

}\eit

\item[$3^\circ$.] Для любых двух индексов $k,l\in I$ оператор $P^k$ переводит ядро $N_l$ и образ $X_l$ оператора $P^l$ в себя:
\beq\label{P_k(N_l)-subseteq-N_l}
P^k(N_l)\subseteq N_l,\qquad P^k(X_l)\subseteq X_l 
\eeq

\item[$4^\circ$.] При увеличении индекса $k$ пространства $X_k$ расширяются:
\beq\label{X_l-supseteq-X_k}
k\le l \quad\Rightarrow\quad X_k\subseteq X_l
\eeq

\item[$5^\circ$.] Для любых индексов $k\le l$ существует единственный оператор $P^k_l:X_l\to X_l$, замыкающий диаграмму  в категории $\LCS$:
$$
 \xymatrix  @R=3pc @C=3pc
 {
X & & X\ar[ll]_{P^k}\ar[d]^{\coim P^l} \\
X_l\ar[u]^{\im P^l} & & X_l  \ar@{-->}[ll]^{P^k_l}
 }
$$
причем
\bit{

\item[(a)] оператор $P^k_l$ является проектором:
$$
(P^k_l)^2=P^k_l,
$$

\item[(b)] а его множество значений $P^k_l(X_l)$ совпадает с пространством $X_k$:
$$
P^k_l(X_l)=X_k
$$
и является одновременно кообразом и образом этого оператора в категории $\LCS$:
$$
 \xymatrix  @R=3pc @C=3pc
 {
X_l & & X_l \ar[ll]_{P^k_l}\ar@/^3ex/[dl]^{\coim P^k_l} \\
& X_k\ar@/^3ex/[ul]^{\im P^k_l} &
 }
$$

}\eit

}\eit

Вместе свойства $2^\circ$ и $5^\circ$ можно изобразить диаграммой:
$$
 \xymatrix  @R=3pc @C=3pc
 {
X & & X\ar@/^12ex/[ddl]^{\coim P^k}
\ar[d]_{\coim P^l}\ar[ll]_{P^k} \\
X_l\ar[u]_{\im P^l} & & X_l\ar@/^3ex/[dl]_{\coim P^k_l} \ar@{-->}[ll]_{P^k_l}  \\
& X_k\ar@/^3ex/[ul]_{\im P^k_l}\ar@/^12ex/[uul]^{\im P^k} &
 }\put(40,-50){$(k\le l)$}
$$

\btm\label{TH:building-generated-by-P^k}
Система операторов $\varPhi=(\{\im P_i^j\},\{\coim P_i^j\})$ образует фундамент в категории $\Ste$ над упорядоченным множеством $I$,
\beq\label{DEF:diagr-building-generated-by-P^k}
 \xymatrix  @R=3pc @C=3pc
 {
X_k\ar@{-->}[r]^{1_{X_k}}
 &  X_k\ar@/^3ex/[d]_{\coim P_k^j} \ar@/^12ex/[dd]^{\coim P_k^i} \\
X_j\ar@{-->}[r]^{1_{X_j}}\ar@/^3ex/[u]_{\im P_j^k}  &  X_j \ar@/^3ex/[d]_{\coim P_j^i} \\
X_i\ar@{-->}[r]^{1_{X_i}} \ar@/^12ex/[uu]^{\im P_i^k} \ar@/^3ex/[u]_{\im P_i^j}  & X_i
 }\put(40,-50){$(i\le j\le k)$}
\eeq
а пространство $X$ c системой операторов $(\{\im P^j\},\{\coim P^j\})$ --- постройку на этом фундаменте.
\etm

Из теоремы \ref{TH:frame-X-edifice} мы получаем

\btm\label{TH:frame->edifice-generated-by-P^k}
Существуют и единственны морфизмы построек 
\beq
\alpha:\fram\varPhi\to X,\qquad 
\beta:X\to\edif\varPhi
\eeq
замыкающие диаграмму
\beq\label{frame->edifice-generated-by-P^k}
 \xymatrix  @R=3pc @C=3pc
 {
 & X\ar[dr]^{\beta} & \\
\fram\varPhi\ar[ur]^{\alpha}\ar[rr]_{\theta_B} & & \edif\varPhi
 }
\eeq 
\etm

\btm\label{TH:X=edifice-in-Ste}
Пусть $\{P^k; \ k\in I\}$ --- система согласованных проекторов на стереотипном пространстве $X$ над направленной по возрастанию решеткой $I$, причем выполняются следующие условия:
\bit{

\item[(i)] все пространства $X_k=P^k(X)$ полны и насыщены (и поэтому стереотипны), 

\item[(ii)] здание $\edif\varPhi_J$ на всяком индуктивном ограничении ${\varPhi}_J\subseteq {\varPhi}$ фундамента ${\varPhi}$ дополняется\footnote{В смысле определения на с.\pageref{complementable-inner-building-2}.} в здании $\edif\varPhi$ на фундаменте ${\varPhi}$ в категории $\Ste$,

\item[(iii)] морфизм $\beta:X\to\edif\varPhi$ является непосредственным мономорфизмом в категории $\Ste$. 

}\eit
Тогда 
\bit{
\item[(a)] здание $\edif\varPhi$ (в категории $\Ste$) представляет собой (не только проективный предел в категории $\Ste$, но и) проективный предел контравариантной системы  $\{\coim P^i_j\}$  в категории $\LCS$:
\beq\label{X=LCS-projlim_(k->infty)X_k}
\edif\varPhi=\Ste\text{-}\projlim_{k\in I} X_k=\LCS\text{-}\projlim_{k\in I} X_k
\eeq

\item[(b)] пространство $\edif\varPhi$ полно и насыщено (и поэтому стереотипно), 

\item[(c)] морфизм $\beta:X\to\edif\varPhi$ является изоморфизмом стереотипных пространств,
\beq\label{X=edif-varPhi}
X\cong\edif\varPhi
\eeq

\item[(d)] морфизм $\alpha:\fram\varPhi\to X$ является плотным отображением:
\beq\label{overline(alpha(fram-varPhi))=X}
\overline{\alpha(\fram\varPhi)}=X
\eeq

}\eit

\etm

\brem
В формулировке этой теоремы можно всюду заменить требование полноты требованием псевдополноты.
\erem

\bpr
1. По теореме \ref{TH:o-postr-v-LCS}, проективный предел $\LCS\text{-}\projlim_{k\to\infty} X_k$ --- стереотипное пространство, поэтому он совпадает с проективным пределом $\Ste\text{-}\projlim_{k\to\infty} X_k$, и отсюда мы получаем цепочку  \eqref{X=LCS-projlim_(k->infty)X_k}. 

2. В силу условия (a) теоремы \ref{TH:o-postr-v-LCS} пространство $\LCS\text{-}\projlim_{k\to\infty} X_k$ насыщено. Поэтому по формуле \eqref{X=LCS-projlim_(k->infty)X_k} пространство $\edif\varPhi$ тоже насыщено.

3. Из коммутативности диаграммы \eqref{frame->edifice-generated-by-P^k} следует, что образ $X$ в $\edif\varPhi$ содержит образ $\theta_\varPhi$:
$$
\theta_\varPhi(\fram\varPhi)\subseteq \beta(X).
$$
Отсюда в силу условия (b) теоремы \ref{TH:o-postr-v-LCS} мы получаем, что $\beta(X)$ плотно в $\edif\varPhi$:
$$
\overline{\beta(X)}=\edif\varPhi.
$$
Поэтому если $\beta$ --- непосредственный мономорфизм, то он является изоморфизмом. То есть справедливо \eqref{X=edif-varPhi}.

4. По условию (b) теоремы \ref{TH:o-postr-v-LCS} морфизм $\theta_\varPhi$ имеет плотный образ в $\edif\varPhi$.  С другой стороны, как мы уже поняли, $\beta$ --- изоморфизм. Значит, морфизм 
$$
\alpha=\beta^{-1}\circ\theta_\varPhi
$$
тоже имеет плотный образ в своей области значений, то есть в $X$.
\epr

Двойственное утверждение выглядит так:

\btm\label{TH:X=frame-in-Ste}
Пусть $\{P^k; \ k\in I\}$ --- система согласованных проекторов на стереотипном пространстве $X$ над направленной по возрастанию решеткой $I$, причем выполняются следующие условия:
\bit{

\item[(i)] все пространства $X_k=P^k(X)$ полны и насыщены (и поэтому стереотипны), 

\item[(ii)] каркас $\fram\varPhi_J$ на всяком индуктивном ограничении ${\varPhi}_J\subseteq {\varPhi}$ фундамента ${\varPhi}$ дополняется\footnote{В смысле определения на с.\pageref{complementable-inner-building-1}.} в каркасе $\fram\varPhi$ на фундаменте ${\varPhi}$ в категории $\Ste$,

\item[(iii)] морфизм $\alpha:\fram\varPhi\to X$ является непосредственным эпиморфизмом в категории $\Ste$. 

}\eit
Тогда 
\bit{
\item[(a)] каркас $\fram\varPhi$ (в категории $\Ste$) представляет собой (не только инъективный предел в категории $\Ste$, но и) инъективный предел ковариантной системы  $\{\im P^i_j\}$  в категории $\LCS$:
\beq\label{X=LCS-injlim_(k->infty)X_k}
\fram\varPhi=\Ste\text{-}\injlim_{k\in I} X_k=\LCS\text{-}\injlim_{k\in I} X_k
\eeq

\item[(b)] пространство $\fram\varPhi$ полно и насыщено (и поэтому стереотипно), 

\item[(c)] морфизм $\alpha:\fram\varPhi\to X$ является изоморфизмом стереотипных пространств,
\beq\label{X=fram-varPhi}
X\cong\fram\varPhi
\eeq

\item[(d)] морфизм $\beta:X\to\edif\varPhi$ является инъективным отображением:
\beq\label{Ker(beta)=0}
\Ker\beta=0.
\eeq

}\eit

\etm

\brem
В формулировке этой теоремы можно всюду заменить требование насыщенности требованием псевдонасыщенности.
\erem

\bpr
1. По теореме \ref{TH:o-nasl-ster-karkasom}, инъективный предел $\LCS\text{-}\injlim_{k\to\infty} X_k$ --- стереотипное пространство, поэтому он совпадает с инъективным пределом $\Ste\text{-}\injlim_{k\to\infty} X_k$, и отсюда мы получаем цепочку  \eqref{X=LCS-injlim_(k->infty)X_k}. 

2. В силу условия (a) теоремы \ref{TH:o-nasl-ster-karkasom} пространство $\LCS\text{-}\injlim_{k\to\infty} X_k$ полно и насыщено. Поэтому по формуле \eqref{X=LCS-injlim_(k->infty)X_k} пространство $\fram\varPhi$ тоже полно и насыщено.

3. В силу коммутативности диаграммы \eqref{frame->edifice-generated-by-P^k},
\beq\label{PROOF:X=frame-in-Ste-1}
\theta_\varPhi=\beta\circ\alpha.
\eeq
При этом в силу условия (b) теоремы \ref{TH:o-nasl-ster-karkasom}, $\theta_\varPhi$ --- мономорфизм. Значит, $\alpha$ --- тоже мономорфизм. С другой стороны, по условию (iii), $\alpha$ --- непосредственный эпиморфизм. Значит, $\alpha$ --- изоморфизм, то есть справедливо \eqref{X=fram-varPhi}.

4. Поскольку $\alpha$ --- изоморфизм, мы можем переписать \eqref{PROOF:X=frame-in-Ste-1} в виде
$$
\beta=\theta_\varPhi\circ\alpha^{-1}.
$$
Здесь $\alpha^{-1}$ --- изоморфизм, а $\theta_\varPhi$ --- мономорфизм (по теореме \ref{TH:o-nasl-ster-karkasom}(b)). Значит, $\beta$ --- тоже мономорфизм.
\epr

\chapter{ЛОКАЛЬНО ВЫПУКЛЫЕ РАССЛОЕНИЯ}

\section{Локально выпуклые расслоения и дифференциально-геомет\-ри\-ческие конструкции}

Нам понадобятся некоторые факты из теории расслоений топологических векторных
пространств. В этом изложении мы следуем идеологии монографии М.~Ж.~Дюпре и
Р.~М.~Жилетта \cite{Dupre-Gillette}, обобщая результаты на случай локально выпуклых пространств.

Нам будет удобно сразу ввести следующее определение.

\bit{

\item[$\bullet$] {\it Векторным расслоением} над полем $\C$ мы будем называть семерку
$(\Xi,M,\pi,\cdot,+,\{0_t;\ t\in M\},-)$, в которой
\bit{

\item[1)] $\Xi$ -- множество, называемое {\it пространством расслоения},

\item[2)] $M$ -- множество, называемое {\it базой расслоения},

\item[3)] $\pi:\Xi\to M$ -- их сюръективное отображение, называемое {\it проекцией расслоения},

\item[4)] $\cdot:\C\times\Xi\to \Xi$ -- отображение, называемое {\it послойным умножением на скаляры},

\item[5)] $+:\Xi\underset{M}{\sqcap} \Xi\to \Xi$ -- отображение, называемое {\it послойным сложением}\footnote{\label{FOOTNOTE:Xi-_M-Xi} Здесь под $\Xi\underset{M}{\sqcap} \Xi$ понимается расслоенное произведение $\Xi$ на себя над $M$, то есть подмножество декартова произведения $\Xi\times \Xi$, состоящее из пар $(\xi,\zeta)$ таких, что $\pi(\xi)=\pi(\zeta)$.}

\item[6)] $\{0_t;\ t\in M\}$ -- семейство элементов пространства $\Xi$, называемых {\it нулями}.

\item[7)] $-:\Xi\to \Xi$ -- отображение, называемое {\it послойным обращением},

}\eit\noindent
причем  для каждой точки $t\in M$ операции $\cdot,+,-$ с элементом $0_t$ задают на прообразе $\Xi_t:=\pi^{-1}(t)\subseteq \Xi$ (называемом {\it слоем над точкой $t$}) структуру векторного пространства над $\C$ с нулем $0_t$.

\item[$\bullet$]  {\it Сечением} векторного расслоения $(\Xi,M,\pi,\cdot,+,\{0_t;\ t\in M\},-)$ мы называем произвольное отображение $x:M\to \Xi$ такое, что
$$
\pi\circ x=\id_M.
$$

\item[$\bullet$] {\it Подрасслоением} векторного расслоения $(\Xi,M,\pi,\cdot,+,\{0_t;\ t\in M\},-)$ мы назваем всякое подмножество $\Psi$ в $\Xi$ пересечение которого с каждым слоем $\pi^{-1}(t)$ образует в нем (непустое) векторное подпространство. Понятно, что $\Psi$ будет расслоением над той же базой относительно тех же структурных элементов.

}\eit

\subsection{Локально выпуклые расслоения}

\bit{

\item[$\bullet$] Пусть векторное расслоение $(\Xi,M,\pi,\cdot,+,\{0_t;\ t\in M\},-)$ над $\C$ наделено следующей дополнительной структурой:
\bit{

\item[1)] пространство расслоения $\Xi$ и база расслоения $M$ наделены топологиями, относительно которых проекция $\pi:\Xi\to M$ является (не только сюръективным, но и) непрерывным и открытым отображением,

\item[2)] задано множество $\mathcal P$ функций $p:\Xi\to\R_+$, называемых {\it полунормами},

}\eit\noindent
причем выполняются следующие условия:
\bit{

\item[(a)] на каждом слое $\pi^{-1}(t)$ ограничения $p|_{\pi^{-1}(t)}:\Xi_t\to\R_+$ функций $p\in{\mathcal P}$ представляют собой систему полунорм, определяющих структуру (отделимого) локально выпуклого пространства на $\pi^{-1}(t)$, топология которого совпадает с индуцированной из $\varXi$,

\item[(b)] отображение послойного умножения на скаляры $\C\times \Xi\to \Xi$
непрерывно:
$$
\Big(\lambda_i\overset{\C}{\underset{i\to\infty}{\longrightarrow}}\lambda,\quad
\xi_i\overset{\varXi}{\underset{i\to\infty}{\longrightarrow}} \xi\Big)
\quad\Longrightarrow\quad
\lambda_i\cdot\xi_i\overset{\varXi}{\underset{i\to\infty}{\longrightarrow}}
\lambda\cdot\xi.
$$

\item[(c)] отображение послойного сложения\footnote{См. подстрочное примечание \ref{FOOTNOTE:Xi-_M-Xi}.} $\Xi\underset{M}{\sqcap} \Xi\to \Xi$ непрерывно:
$$
\Big(\xi_i\overset{\varXi}{\underset{i\to\infty}{\longrightarrow}} \xi,\quad \zeta_i\overset{\varXi}{\underset{i\to\infty}{\longrightarrow}} \zeta,\quad \pi(\xi_i)=\pi(\zeta_i),\quad \pi(\xi)=\pi(\zeta)\Big)
\quad\Longrightarrow\quad
\xi_i+\zeta_i\overset{\varXi}{\underset{i\to\infty}{\longrightarrow}} \xi+\zeta.
$$

\item[(d)] всякая полунорма $p\in{\mathcal P}$ представляет собой полунепрерывное сверху отображение $p:\Xi\to\R_+$, то есть для любого $\e>0$ множество точек $\{\xi\in \Xi: \ p(\xi)<\e\}$ открыто, или, что то же самое,
$$
\Big(\xi_i\overset{\varXi}{\underset{i\to\infty}{\longrightarrow}} \xi,\quad p(\xi)<\e\Big)
\quad\Longrightarrow\quad
\underbrace{p(\xi_i)<\e}_{\text{для почти всех $i$}},
$$

\item[(e)]\label{uslovie-e-DEF-rassloeniya} для любой точки $t\in M$ и всякой окрестности $V$ точки $0_t$ в $\Xi$ найдутся полунорма $p\in{\mathcal P}$, число $\e>0$ и открытое множество $U$ в $M$, содержащее $t$, такие, что
$$
\{ \xi\in \pi^{-1}(U): \ p(\xi)<\e\}\subseteq V;
$$
иными словами, выполняется импликация
$$
\Big(\pi(\xi_i)\overset{M}{\underset{i\to\infty}{\longrightarrow}} t\quad\&\quad \forall p\in{\mathcal P}\quad p(\xi_i)\underset{i\to\infty}{\longrightarrow} 0\Big)
\quad\Longrightarrow\quad
\xi_i\overset{\varXi}{\underset{i\to\infty}{\longrightarrow}} 0_t.
$$
}\eit
Тогда систему $(\Xi,M,\pi,\cdot,+,\{0_t;\ t\in M\},-)$ с описанными топологиями на $M$ и $\Xi$ и системой полунорм $\mathcal P$ мы называем {\it локально выпуклым расслоением}.
}\eit

\brem Условие (b) в этом списке можно заменить на формально более слабое условие, что для всякого $\lambda\in\C$ отображение умножения $\xi\mapsto\lambda\cdot\xi$ непрерывно из $\varXi$ в $\varXi$:
$$
\xi_i\overset{\varXi}{\underset{i\to\infty}{\longrightarrow}} \xi
\quad\Longrightarrow\quad \forall\lambda\in\C\quad
\lambda\cdot\xi_i\overset{\varXi}{\underset{i\to\infty}{\longrightarrow}} \lambda\cdot\xi.
$$
Действительно, при выполнении этого условия из $\lambda_i\to \lambda$ и $\xi_i\to\xi$ будет следовать с одной стороны,
$$
\forall p\in{\mathcal P}\qquad p(\lambda_i\cdot\xi_i-\lambda\cdot\xi_i)\le \abs{\lambda_i-\lambda}\cdot p(\xi_i)\overset{\R}{\underset{i\to\infty}{\longrightarrow}} 0
$$
а с другой,
$$
\pi(\xi_i)\to\pi(\xi)\quad\Longrightarrow\quad \pi(\lambda_i\cdot\xi_i-\lambda\cdot\xi_i)=\pi(\xi_i)\overset{M}{\underset{i\to\infty}{\longrightarrow}}\pi(\xi)
$$
Вместе в силу (e) это дает
$$
\lambda_i\cdot\xi_i-\lambda\cdot\xi_i\overset{\varXi}{\underset{i\to\infty}{\longrightarrow}} 0_{\pi(\xi)},
$$
а это в свою очередь в силу (c) дает
$$
\lambda_i\cdot\xi_i=\underbrace{\lambda_i\cdot\xi_i-\lambda\cdot\xi_i}_{\scriptsize\begin{matrix}\downarrow\\ 0_{\pi(\xi)} \end{matrix}}+\underbrace{\lambda\cdot\xi_i}_{\scriptsize\begin{matrix}\downarrow\\ \lambda\cdot\xi \end{matrix}}\overset{\varXi}{\underset{i\to\infty}{\longrightarrow}} 0_{\pi(\xi)}+\lambda\cdot\xi=\lambda\cdot\xi.
$$
\erem

\brem Из условий (b) и (c) следует, что в (c) можно заменить послойное сложение послойным вычитанием:
\beq\label{xi_i->xi&zeta_i->zeta=>xi_i-zeta_i->xi-zeta}
\Big(\xi_i\overset{\varXi}{\underset{i\to\infty}{\longrightarrow}} \xi,\quad \zeta_i\overset{\varXi}{\underset{i\to\infty}{\longrightarrow}} \zeta,\quad \pi(\xi_i)=\pi(\zeta_i),\quad \pi(\xi)=\pi(\zeta)\Big)
\quad\Longrightarrow\quad
\xi_i-\zeta_i\overset{\varXi}{\underset{i\to\infty}{\longrightarrow}} \xi-\zeta.
\eeq
\erem

\bprop\label{PROP:harakt-shod-v-Xi}
В локально выпуклом расслоении $(\Xi,M,\pi)$ соотношение
$$
\xi_i\overset{\varXi}{\underset{i\to\infty}{\longrightarrow}} \xi
$$
эквивалентно следующим условиям:
\bit{

\item[(i)] $\pi(\xi_i)\overset{M}{\underset{i\to\infty}{\longrightarrow}} \pi(\xi)$,

\item[(ii)] для всякой полунормы $p\in{\mathcal P}$ и любого $\e>0$ существует направленность $\zeta_i\in \varXi$ и элемент $\zeta\in\varXi$ такие, что
\beq\label{trebovanie-dlya-shod-napravl}
\zeta_i\overset{\varXi}{\underset{i\to\infty}{\longrightarrow}} \zeta
\quad
\underbrace{\pi(\zeta_i)=\pi(\xi_i)}_{\text{для почти всех $i$}},
\quad
\pi(\zeta)=\pi(\xi),
\quad
p(\zeta-\xi)<\e,
\quad
\underbrace{p(\zeta_i-\xi_i)<\e}_{\text{для почти всех $i$}}.
\eeq
}\eit
\eprop
 \bpr
Понятно, что здесь нужно доказать достаточность. Пусть выполнены (i) и (ii). Зафиксируем $p\in{\mathcal P}$ и $\e>0$ и подберем направленность $\zeta_i$, описанную в (ii). Пусть $V$ -- произвольная окрестность точки $\xi$ в $\varXi$. Поскольку отображение $\pi$ открыто, образ $\pi(V)$ множества $V$ должен быть окрестностью точки $\pi(\xi)$ в $M$. Поэтому из условия $\pi(\xi_i)\overset{M}{\underset{i\to\infty}{\longrightarrow}} \pi(\xi)$ следует, что почти все точки $\pi(\xi_i)$ лежат в $\pi(V)$:
$$
\exists i_V:\quad \forall i\ge i_V \quad \pi(\xi_i)\in \pi(V).
$$
Значит, существуют точки $\{\zeta^V_i;\ i\ge i_V\}$ такие, что
$$
\zeta^V_i\in V,\qquad \pi(\zeta^V_i)=\pi(\xi_i).
$$
Мы получили двойную направленность $\{\zeta^V_i;\ V\in{\mathcal U}(\xi),\ i\ge i_V\}$, в которой верхний индекс $V$ пробегает систему ${\mathcal U}(\xi)$ всех окрестностей точки $\xi$ в $\varXi$, упорядоченную по включению в сторону сужения, со следующими свойствами:
$$
\pi(\xi_i)=\pi(\zeta^V_i),\qquad \zeta^V_i\overset{\varXi}{\underset{\scriptsize\begin{matrix}i\to\infty\\ V\to \{\xi\}\end{matrix}}{\longrightarrow}}\xi
$$
Вместе с условиями $\pi(\zeta)=\pi(\xi)$, $\pi(\zeta_i)=\pi(\xi_i)$ (для почти всех $i$) и $\zeta_i\overset{\varXi}{\underset{i\to\infty}{\longrightarrow}} \zeta$ из \eqref{trebovanie-dlya-shod-napravl} это в силу \eqref{xi_i->xi&zeta_i->zeta=>xi_i-zeta_i->xi-zeta} дает соотношение
$$
\zeta_i-\zeta^V_i\overset{\varXi}{\underset{\scriptsize\begin{matrix}i\to\infty\\ V\to \{\xi\}\end{matrix}}{\longrightarrow}}\zeta-\xi
$$
Оно, в свою очередь, вместе с неравенством $p(\zeta-\xi)<\e$ из \eqref{trebovanie-dlya-shod-napravl}, в силу (d), дает неравенство
$$
p(\zeta_i-\zeta^V_i)<\e
$$
(верное для почти всех $i$), из которого затем выводится неравенство
$$
p(\xi_i-\zeta^V_i)=p(\xi_i-\zeta_i+\zeta_i-\zeta^V_i)\le \overbrace{p(\xi_i-\zeta_i)}^{\scriptsize\begin{matrix}\e \\ \phantom{\tiny{\eqref{trebovanie-dlya-shod-napravl}}}\ \text{\rotatebox{90}{$<$}}\ \tiny{\eqref{trebovanie-dlya-shod-napravl}}\end{matrix}}+p(\zeta_i-\zeta^V_i)<\e+\e=2\e.
$$
(также верное для почти всех $i$). Добавив к нему очевидное соотношение
$$
\pi(\xi_i-\zeta^V_i)=\pi(\xi_i)\overset{M}{\underset{\scriptsize\begin{matrix}i\to\infty\\ V\to \{\xi\}\end{matrix}}{\longrightarrow}}\pi(\xi)
$$
мы, в силу (e) получим:
$$
\xi_i-\zeta^V_i\overset{\varXi}{\underset{\scriptsize\begin{matrix}i\to\infty\\ V\to \{\xi\}\end{matrix}}{\longrightarrow}}0_{\pi(\xi)}.
$$
Теперь, применяя (c), мы получаем:
$$
\xi_i=\underbrace{\xi_i-\zeta^V_i}_{\scriptsize\begin{matrix}\downarrow\\ 0_{\pi(\xi)} \end{matrix}}+\underbrace{\zeta^V_i}_{\scriptsize\begin{matrix}\downarrow\\ \xi \end{matrix}}\overset{\varXi}{\underset{\scriptsize\begin{matrix}i\to\infty\\ V\to \{\xi\}\end{matrix}}{\longrightarrow}} 0_{\pi(\xi)}+\xi=\xi.
$$
 \epr

\paragraph{Непрерывные сечения локально выпуклого расслоения.}
Мы рассмотрим здесь {\it непрерывные сечения локально выпуклых расслоений} $\pi:\Xi\to M$, то есть непрерывные отображения $x:M\to \Xi$ такие, что
$$
\pi\circ x=\id_M.
$$
Множество всех непрерывных сечений обозначается $\Sec(\pi)$ и наделяется структурой левого
${\mathcal C}(M)$-модуля и топологией равномерной сходимости на компактах в $M$.

\vglue10pt \centerline{\bf Свойства непрерывных сечений:} \vglue10pt
{\it
\begin{itemize}
\item[$1^\circ$.] Пространство $\Sec(\pi)$ непрерывных сечений всякого локально выпуклого расслоения $\pi:\varXi\to M$ над произвольным паракомпактным локально компактным пространством $M$ является локально выпуклым ${\mathcal C}(M)$-модулем (с совместно непрерывным умножением) относительно послойного умножения
    $$
    (a\cdot x)(t)=a(t)\cdot x(t),\qquad a\in {\mathcal C}(M),\ x\in \Sec(\pi)
    $$
    и полунорм
$$
p_T(x)=\sup_{t\in T}p(x(t)),\qquad p\in{\mathcal P}.
$$
где $T$ -- всевозможные компакты в $M$.

\item[$2^\circ$.] \label{LM:poltnost-v-sloyah=>plotnost-v-secheniyah} Пусть $M$ -- паракомпактное локально компактное пространство, $\pi:\varXi\to M$ -- локально выпуклое расслоение, и $X$ --  ${\mathcal C}(M)$-подмодуль в ${\mathcal C}(M)$-модуле непрерывных сечений $\pi$,
$$
X\subseteq\Sec(\pi)
$$
плотный в каждом слое:
\beq\label{plotnost-v-sloyah}
\forall t\in M\qquad \overline{\{x(t),\ x\in X\}}=\pi^{-1}(t).
\eeq
Тогда $X$ плотен в $\Sec(\pi)$:
$$
\overline{X}=\Sec(\pi).
$$
\end{itemize}
} \vglue10pt

\bpr
1. Если $x\in\Sec(\pi)$ и $a\in {\mathcal C}(M)$, то послойное произведение $a\cdot x$ будет непрерывным сечением из-за свойства (b):
$$
t_i\overset{M}{\underset{i\to\infty}{\longrightarrow}}t \quad\Longrightarrow\quad
\Big(
a(t_i)\overset{\C}{\underset{i\to\infty}{\longrightarrow}}a(t),\quad x(t_i)\overset{\varXi}{\underset{i\to\infty}{\longrightarrow}} x(t)
\Big)
\quad\Longrightarrow\quad
a(t_i)\cdot x(t_i)\overset{\varXi}{\underset{i\to\infty}{\longrightarrow}} a(t)\cdot x(t).
$$
Если $T$ -- компакт в $M$, то для всякого непрерывного сечения $x\in\Sec(\pi)$ образ $x(T)$
является компактным пространством в $\Xi$. Как следствие, образ $p(x(T))$ при
непрерывном отображении $p\in{\mathcal P}$ из $\Xi$ в $\R$ с топологией,
порожденной базой открытых множеств вида $(-\infty,\e)$, $\e\in\R$, тоже должен быть
компактным подмножеством в $\R$. Поэтому каждое его покрытие множествами вида
$(-\infty,\e)$ содержит конечное подпокрытие, и это означает, что $p(x(T))$
ограничено в $\R$ в обычном смысле. То есть конечна величина
$$
p_T(x)=\sup_{t\in T}p(x(t)).
$$
Она, очевидно, будет полунормой на $\Sec(\pi)$, и из цепочки
$$
p_T(a\cdot x)=\sup_{t\in T}p(a(t)\cdot x(t))\le \sup_{t\in
T}\Big(\abs{a(t)}\cdot p(x(t))\Big)\le \sup_{t\in T}\abs{a(t)}\cdot\sup_{t\in
T}p(x(t))
$$
следует, что такие полунормы превращают $\Sec(\pi)$ в локально выпуклый
${\mathcal C}(M)$-модуль (с совместно непрерывным умножением).

2. Пусть $y\in \Sec(\pi)$ и $\e>0$. Из \eqref{plotnost-v-sloyah} следует, что
для любой полунормы $p\in{\mathcal P}$ и всякой точки $t\in M$ существует непрерывное
сечение $x_t\in\Sec(\pi)$ со свойством
$$
p(x_t(t)-y(t))<\e.
$$
Из того, что полунорма $p$ полунепрерывна сверху (а отображения $x_t$ и $y$ непрерывны), следует, что множество
$$
U_t=\{s\in M: \ p(x_t(s)-y(s))<\e\}
$$
открыто, и поэтому является окрестностью точки $t$. Как следствие, семейство $\{U_t;\ t\in M\}$ является открытым покрытием пространства $M$.

Зафиксируем теперь какой-нибудь компакт $T\subseteq M$. Его покрытие $\{U_t;\ t\in T\}$ содержит некоторое конечное подпокрытие $\{U_{t_1},...,U_{t_n}\}$. Подберем подчиненное ему разбиение единицы
$$
0\le a_i\le 1,\qquad \supp a_i\subseteq U_{t_i},\qquad \sum_{i=1}^n a_i(t)=1\quad (t\in T)
$$
и положим
$$
x=\sum_{i=1}^n a_i\cdot x_{t_i}
$$
(поскольку $\Sec(\pi)$ есть ${\mathcal C}(M)$-модуль, $x\in \Sec(\pi)$). Тогда
\begin{multline*}
p_T(x-y)=\sup_{t\in T}p(x(t)-y(t))=\sup_{t\in T}p\l \sum_{i=1}^n a_i(t)\cdot
x_{t_i}(t)-\sum_{i=1}^n a_i(t)\cdot y(t)\r=\\= \sup_{t\in T}p\l \sum_{i=1}^n
a_i(t)\cdot \big(x_{t_i}(t)-y(t)\big)\r\le \sup_{t\in T}\sum_{i=1}^n a_i(t)\cdot
p\l \big(x_{t_i}(t)-y(t)\big)\r<\sup_{t\in T}\sum_{i=1}^n a_i(t)\cdot \e=\e.
\end{multline*}
\epr

\paragraph{Задание локально выпуклого расслоения системами сечений и полунорм.}

\bprop\label{PROP:sushestv-topologii-v-Xi} Пусть даны
\bit{

\item[1)] расслоение векторных пространств $(\Xi,M,\pi,\cdot,+,\{0_t;\ t\in M\},-)$ над $\C$,

\item[2)] векторное пространство $X$ его сечений,

\item[3)] система ${\mathcal P}$ функций на $\varXi$,

\item[4)] топология на базе $M$,
}\eit\noindent
причем выполняются следующие условия:
\bit{

\item[(i)] на каждом слое ограничения $p|_{\pi^{-1}(t)}$ функций $p\in{\mathcal P}$ образуют систему полунорм, превращающую $\pi^{-1}(t)$ в (отделимое) локально выпуклое пространство;

\item[(ii)] система ${\mathcal P}$ направлена по возрастанию: для любых двух функций $p,q\in{\mathcal P}$ найдется функция $r\in{\mathcal P}$, мажорирующая $p$ и $q$:
 \beq\label{p<r&q<r}
p(\upsilon)\le r(\upsilon),\qquad q(\upsilon)\le r(\upsilon),\qquad \upsilon\in \varXi
 \eeq

\item[(iii)] для любого сечения $x\in X$ и любой полунормы $p\in{\mathcal P}$
функция $t\in M\mapsto p(x(t))$ полунепрерывна сверху на $M$,

\item[(iv)] для всякой точки $t\in M$ множество $\{x(t);\ x\in X\}$ плотно в
локально выпуклом пространстве $\pi^{-1}(t)$.

}\eit
Тогда существует единственная топология на $\varXi$, необходимая для того, чтобы система 
$$
(\Xi,M,\pi,\cdot,+,\{0_t;\ t\in M\},-)
$$
с заданной топологией на $M$ и системой полунорм $\mathcal P$ превратилась в локально выпуклое расслоение, множество непрерывных сечений которого содержит $X$:
$$
X\subseteq\Sec(\pi).
$$
При этом базу такой топологии в $\varXi$ образуют множества вида
 \beq\label{baza-topologii-v-varXi}
W(x,U,p,\e)=\{\xi\in\varXi: \ \pi(\xi)\in U\ \& \ p(\xi-x(\pi(\xi)))<\e\},
 \eeq
где $x\in X$, $p\in{\mathcal P}$, $\e>0$ и $U$ -- открытое множество в $M$.
\eprop

\bpr 1. Покажем сначала, что множества \eqref{baza-topologii-v-varXi} действительно образуют
базу некоторой топологии в $\varXi$. Прежде всего, они покрывают $\varXi$, потому что если $\xi\in\varXi$, то из условия (iii) следует, что для любых $\e>0$ и $p\in{\mathcal P}$ найдется $x\in X$ такой что
$$
p(\xi-x(\pi(\xi)))<\e,
$$
и если теперь выбрать какую-нибудь открытую окрестность $U$ точки $\pi(\xi)$, то точка $\xi$ будет лежать в множестве $W(x,U,p,\e)$.

Проверим далее вторую аксиому базы: рассмотрим точку $\xi$, какие-нибудь ее базисные окрестности $W(x,U,p,\e)$ и $W(y,V,q,\delta)$,
 \beq\label{sushestv-topologii-v-Xi-8}
\xi\in W(x,U,p,\e)\cap W(y,V,q,\delta),
 \eeq
и покажем, что существует ее базисная окрестность $W(z,O,r,\sigma)$ такая, что
 \beq\label{sushestv-topologii-v-Xi-7}
\xi\in W(z,O,r,\sigma)\subseteq W(x,U,p,\e)\cap W(y,V,q,\delta).
 \eeq
Включение \eqref{sushestv-topologii-v-Xi-8} означает выполнение условий
 \beq\label{sushestv-topologii-v-Xi-6}
\pi(\xi)\in U,\qquad p(\xi-x(\pi(\xi)))<\e,\qquad  \pi(\xi)\in V,\qquad
q(\xi-y(\pi(\xi)))<\delta.
 \eeq
Рассмотрим слой $\pi^{-1}(\pi(\xi))$. Условия
 \beq\label{sushestv-topologii-v-Xi-0}
p(\xi-x(\pi(\xi)))<\e,\qquad  q(\xi-y(\pi(\xi)))<\delta
 \eeq
можно понимать так, что точка $\xi$ лежит в пересечении окрестностей точек $x(\pi(\xi))$ и $y(\pi(\xi))$, определяемых полунормами $p$ и $q$ с радиусами $\e$ и $\delta$. Поэтому (здесь мы используем (ii)) существует некая полунорма $r\in{\mathcal P}$ и число $\sigma>0$ такие, что $r$-окрестность точки $\xi$ радиуса $2\sigma$ содержится в описанных $p$- и $q$-окрестностях:
 \beq\label{sushestv-topologii-v-Xi-1}
\forall \zeta\in\pi^{-1}(\pi(\xi))\qquad r(\zeta-\xi)<2\sigma\quad\Longrightarrow\quad
\Big(p(\zeta-x(\pi(\xi)))<\e\quad\&\quad  q(\zeta-y(\pi(\xi)))<\delta\Big).
 \eeq
Немного уменьшим $\sigma$, если это необходимо, так, чтобы условия \eqref{sushestv-topologii-v-Xi-0} можно было заменить на
 \beq\label{sushestv-topologii-v-Xi-4}
p(\xi-x(\pi(\xi)))<\e-2\sigma,\qquad  q(\xi-y(\pi(\xi)))<\delta-2\sigma.
 \eeq
Воспользуемся затем условием (iv) и выберем сечение $z\in X$ так, чтобы
 \beq\label{sushestv-topologii-v-Xi-5}
r(z(\pi(\xi))-\xi)<\sigma.
 \eeq
Тогда мы получим цепочку импликаций
 \begin{multline*}
r(\zeta-z(\pi(\xi)))<\sigma\quad\Longrightarrow\quad
r(\zeta-\xi)\le r(\zeta-z(\pi(\xi)))+r(z(\pi(\xi))-\xi)<\sigma+\sigma=2\sigma
\quad\Longrightarrow
\\ \overset{\eqref{sushestv-topologii-v-Xi-1}}{\Longrightarrow}\quad
\Big(p(\zeta-x(\pi(\xi)))<\e\quad\&\quad  q(\zeta-y(\pi(\xi)))<\delta\Big).
 \end{multline*}
Иными словами, в слое $\pi^{-1}(\pi(\xi))$ $r$-окрестность радиуса $\sigma$ точки $z(\pi(\xi))$ также содержится в описанных $p$- и $q$-окрестностях:
 \beq\label{sushestv-topologii-v-Xi-2}
\forall \zeta\in\pi^{-1}(\pi(\xi))\qquad r(\zeta-z(\pi(\xi)))<\sigma\quad\Longrightarrow\quad
\Big(p(\zeta-x(\pi(\xi)))<\e\quad\&\quad  q(\zeta-y(\pi(\xi)))<\delta\Big).
 \eeq
Заметим, что в силу условия (ii), можно считать, что полунорма $r$ мажорирует полунормы $p$ и $q$, то есть выполняется \eqref{p<r&q<r}.
Примем это и рассмотрим множество
 \beq\label{sushestv-topologii-v-Xi-3}
O=\left\{s\in U\cap V:\quad p(z(s)-x(s))<\e-\sigma \quad \&\quad q(z(s)-y(s))<\delta-\sigma \right\}.
 \eeq
В силу (iii) оно открыто в $M$. Кроме того, оно содержит точку $\pi(\xi)$, потому что, во-первых, $\pi(\xi)\in U\cap V$ в силу \eqref{sushestv-topologii-v-Xi-6}, во-вторых,
$$
p(z(\pi(\xi))-x(\pi(\xi)))\le \underbrace{p(z(\pi(\xi))-\xi)}_{\scriptsize\begin{matrix}
\phantom{\tiny{\eqref{p<r&q<r}}}\quad\text{\rotatebox{90}{$\ge$}}\quad\tiny{\eqref{p<r&q<r}}  \\
r(z(\pi(\xi))-\xi) \\
\phantom{\tiny{\eqref{sushestv-topologii-v-Xi-5}}}\quad\text{\rotatebox{90}{$>$}}\quad\tiny{\eqref{sushestv-topologii-v-Xi-5}} \\
\sigma
\end{matrix}}+\underbrace{p(\xi-x(\pi(\xi)))}_{\scriptsize\begin{matrix}
\phantom{\tiny{\eqref{sushestv-topologii-v-Xi-5}}}\quad\text{\rotatebox{90}{$>$}}\quad\tiny{\eqref{sushestv-topologii-v-Xi-4}} \\
\e-2\sigma
\end{matrix}}<\e-\sigma,
$$
и, в-третьих,
$$
p(z(\pi(\xi))-y(\pi(\xi)))\le \underbrace{p(z(\pi(\xi))-\xi)}_{\scriptsize\begin{matrix}
\phantom{\tiny{\eqref{p<r&q<r}}}\quad\text{\rotatebox{90}{$\ge$}}\quad\tiny{\eqref{p<r&q<r}}  \\
r(z(\pi(\xi))-\xi) \\
\phantom{\tiny{\eqref{sushestv-topologii-v-Xi-5}}}\quad\text{\rotatebox{90}{$>$}}\quad\tiny{\eqref{sushestv-topologii-v-Xi-5}} \\
\sigma
\end{matrix}}+\underbrace{p(\xi-y(\pi(\xi)))}_{\scriptsize\begin{matrix}
\phantom{\tiny{\eqref{sushestv-topologii-v-Xi-5}}}\quad\text{\rotatebox{90}{$>$}}\quad\tiny{\eqref{sushestv-topologii-v-Xi-4}} \\
\delta-2\sigma
\end{matrix}}<\delta-\sigma.
$$
Включение $\pi(\xi)\in O$ вместе с условием \eqref{sushestv-topologii-v-Xi-5} означают
 \beq\label{sushestv-topologii-v-Xi-9}
\xi\in W(z,O,r,\sigma).
 \eeq
Далее, для всякой точки $\zeta\in\varXi$ мы получим
$$
\pi(\zeta)\in O\quad\&\quad r(\zeta-z(\pi(\zeta)))<\sigma
\quad\Longrightarrow\quad
p(\zeta-x(\pi(\zeta)))\le \underbrace{p(\zeta-z(\pi(\zeta)))}_{\scriptsize\begin{matrix}
\phantom{\tiny{\eqref{p<r&q<r}}}\quad\text{\rotatebox{90}{$\ge$}}\quad\tiny{\eqref{p<r&q<r}}  \\
r(\zeta-z(\pi(\zeta))) \\
\text{\rotatebox{90}{$>$}} \\
\sigma
\end{matrix}}
+\underbrace{p(z(\pi(\zeta))-x(\pi(\zeta)))}_{\scriptsize\begin{matrix}
\phantom{\tiny{\eqref{sushestv-topologii-v-Xi-3}}}\quad\text{\rotatebox{90}{$>$}}\quad\tiny{\eqref{sushestv-topologii-v-Xi-3}} \\
\e-\sigma
\end{matrix}}<\e
$$
и поэтому
$$
W(z,O,r,\sigma)\subseteq W(x,U,p,\e).
$$
Аналогично,
$$
W(z,O,r,\sigma)\subseteq W(y,V,q,\delta).
$$
Вместе с \eqref{sushestv-topologii-v-Xi-9} это дает \eqref{sushestv-topologii-v-Xi-7}.

2. Итак, мы поняли, что множества \eqref{baza-topologii-v-varXi} образуют базу некоторой топологии на $\Xi$. Заметим теперь, что относительно этой топологии отображение $\pi:\Xi\to M$ непрерывно и открыто. Пусть
$$
\xi_i\overset{\Xi}{\underset{i\to\infty}{\longrightarrow}}\xi.
$$
Рассмотрим произвольную окрестность $U$ точки $\pi(\xi)$. Выберем какую-нибудь полунорму $p\in{\mathcal P}$, и, пользуясь условием (iv), подберем сечение $x\in X$ так, чтобы $p(\xi-x(\pi(\xi)))<1$. Тогда множество $W(x,U,p,1)$ будет окрестностью точки $\xi$, поэтому $\xi\in W(x,U,p,1)$ для почти всех $i$, а отсюда следует, что $\pi(\xi_i)\in U$ для почти всех $i$. Это доказывает соотношение
$$
\pi(\xi_i)\overset{M}{\underset{i\to\infty}{\longrightarrow}}\pi(\xi),
$$
то есть непрерывность $\pi$.

Теперь рассмотрим произвольную базисную окрестность $W(x,U,p,\e)$. При проекции $\pi$ она отображается в открытое множество $U$, причем сюръективно, потому что у каждой точки $t\in U$ найдется прообраз $x(t)\in W(x,U,p,\e)$. Это значит, что отображение $\pi$ открыто.

3. Покажем далее, что относительно полученной топологии в $\Xi$ всякое сечение $x\in X$ будет непрерывным отображением. Пусть
$$
t_i\overset{M}{\underset{i\to\infty}{\longrightarrow}}t.
$$
Рассмотрим какую-нибудь базисную окрестность $W(y,U,p,\e)$ точки $x(t)$. Условие $x(t)\in W(y,U,p,\e)$ означает
 \beq\label{sushestv-topologii-v-Xi-10}
t\in U,\qquad p(x(t)-y(t))<\e.
 \eeq
Положим
$$
V=\{s\in U:\quad p(x(s)-y(s))<\e\}.
$$
В силу условия (iii) это множество открыто в $M$. В силу \eqref{sushestv-topologii-v-Xi-10} оно содержит точку $t$, то есть является окрестностью точки $t$. Поэтому для почти всех индексов $i$ мы последовательно получаем
$$
t_i\in V\quad\Longrightarrow\quad p(x(t_i)-y(t_i))<\e\quad\Longrightarrow\quad x(t_i)\in W(y,V,p,\e)\subseteq W(y,U,p,\e).
$$
Это доказывает соотношение
 \beq\label{sushestv-topologii-v-Xi-12}
x(t_i)\overset{\Xi}{\underset{i\to\infty}{\longrightarrow}}x(t),
 \eeq
то есть непрерывность $x$.

4. Теперь начнем проверять, что относительно заданной топологии на $\Xi$ тройка $(\Xi,M,\pi)$ становится
локально выпуклым расслоением. Проверим прежде всего непрерывность послойного умножения на скаляры. Пусть
$$
\lambda_i\overset{\C}{\underset{i\to\infty}{\longrightarrow}}\lambda,\quad
\xi_i\overset{\varXi}{\underset{i\to\infty}{\longrightarrow}} \xi
$$
Рассмотрим сначала случай $\lambda\ne 0$. Выберем какую-нибудь окрестность $W(x,U,p,\e)$ точки $\lambda\cdot\xi$:
$$
\pi(\lambda\cdot\xi)\in U,\qquad p(\lambda\cdot\xi-x(\pi(\xi)))<\e.
$$
Тогда
$$
\pi(\xi)\in U,\qquad p\l \xi-\frac{1}{\lambda}\cdot x(\pi(\xi))\r <\frac{\e}{\abs{\lambda}}.
$$
то есть множество $W\l \frac{1}{\lambda}\cdot x,U,p,\frac{\e}{\abs{\lambda}}\r$ является окрестностью для точки $\xi$, и поэтому $\xi_i\in W\l \frac{1}{\lambda}\cdot x,U,p,\frac{\e}{\abs{\lambda}}\r$ для почти всех $i$:
 \beq\label{sushestv-topologii-v-Xi-11}
\exists i_0\qquad\forall i\ge i_0\qquad
\pi(\xi_i)\in U,\qquad p\l \xi_i-\frac{1}{\lambda}\cdot x(\pi(\xi_i))\r <\frac{\e}{\abs{\lambda}}.
 \eeq
Теперь для почти всех $i$ получаем
\begin{multline*}
p(\lambda_i\cdot\xi_i-x(\pi(\lambda_i\cdot\xi_i)))=p(\lambda_i\cdot\xi_i-x(\pi(\xi_i)))=
\abs{\lambda_i}\cdot p\l \xi_i-\frac{1}{\lambda_i}\cdot x(\pi(\xi_i))\r\le \\ \le
\abs{\lambda_i}\cdot p\l \xi_i-\frac{1}{\lambda}\cdot x(\pi(\xi))\r+\abs{\lambda_i}\cdot p\l \frac{1}{\lambda}\cdot x(\pi(\xi))-\frac{1}{\lambda_i}\cdot x(\pi(\xi_i))\r=\\=
\abs{\lambda_i}\cdot p\l \xi_i-\frac{1}{\lambda}\cdot x(\pi(\xi))\r+p\l \frac{\lambda_i}{\lambda}\cdot x(\pi(\xi))-x(\pi(\xi_i))\r=
\\=
\abs{\lambda_i}\cdot p\l \xi_i-\frac{1}{\lambda}\cdot x(\pi(\xi))\r+p\l \frac{\lambda_i}{\lambda}\cdot x(\pi(\xi))-x(\pi(\xi))\r
+p\Big( x(\pi(\xi))-x(\pi(\xi_i))\Big)=
\\=
\underbrace{\abs{\lambda_i}}_{\scriptsize\begin{matrix}
\text{\rotatebox{90}{$>$}} \\
2\abs{\lambda}
\end{matrix}}\cdot \underbrace{p\l \xi_i-\frac{1}{\lambda}\cdot x(\pi(\xi))\r}_{\scriptsize\begin{matrix}
\phantom{\tiny{\eqref{sushestv-topologii-v-Xi-11}}}\quad\text{\rotatebox{90}{$>$}}\quad\tiny{\eqref{sushestv-topologii-v-Xi-11}} \\
\frac{\e}{\abs{\lambda}}
\end{matrix}}
+\underbrace{\abs{\frac{\lambda_i}{\lambda}-1}}_{\scriptsize\begin{matrix}
\text{\rotatebox{90}{$>$}} \\
\e
\end{matrix}}\cdot p\Big(x(\pi(\xi))\Big)
+\underbrace{p\Big( x(\pi(\xi))-x(\pi(\xi_i))\Big)}_{\scriptsize\begin{matrix}
\phantom{\tiny{\eqref{sushestv-topologii-v-Xi-12}}}\quad\text{\rotatebox{90}{$>$}}\quad\tiny{\eqref{sushestv-topologii-v-Xi-12}} \\
\e
\end{matrix}}<3\e+\e\cdot p\big(x(\pi(\xi))\big)
\end{multline*}
Добавляя к этому условие $\pi(\xi_i)\in U$ из \eqref{sushestv-topologii-v-Xi-11}, мы получаем
$$
\lambda_i\cdot\xi_i\in W(x,U,p,3\e+\e\cdot p(x(\pi(\xi))))
$$
для почти всех $i$. Поскольку $\e$ выбирался произвольным,
 \beq\label{lambda_i-xi_i->lambda-xi}
\lambda_i\cdot\xi_i\overset{\varXi}{\underset{i\to\infty}{\longrightarrow}}
\lambda\cdot\xi.
 \eeq

Теперь рассмотрим случай $\lambda=0$. Выберем какую-нибудь окрестность $W(x,U,p,\e)$ точки $0_{\pi(\xi)}$:
$$
\pi(\xi)\in U,\qquad p(x(\pi(\xi)))<\e.
$$
Подберем $\delta>0$ так, чтобы
$$
p(x(\pi(\xi)))<\e-\delta.
$$
Поскольку $\xi_i\to\xi$, для почти всех $i$ выполняется
$$
\pi(\xi_i)\in U,\qquad p(x(\pi(\xi_i)))<\e-\delta.
$$
Теперь если $p(\xi)\ne 0$, то для почти всех $i$
$$
p(\lambda_i\cdot\xi_i-x(\pi(\xi_i)))\le \underbrace{\abs{\lambda_i}}_{\scriptsize\begin{matrix}
\text{\rotatebox{90}{$>$}}\\ \frac{\delta}{2p(\xi)} \end{matrix}}\cdot \underbrace{p(\xi_i)}_{\scriptsize\begin{matrix}
\text{\rotatebox{90}{$>$}} \\ 2p(\xi)\end{matrix}} +\underbrace{p(x(\pi(\xi_i)))}_{\scriptsize\begin{matrix}
\text{\rotatebox{90}{$>$}} \\ \e-\delta\end{matrix}}<\e.
$$
Если же $p(\xi)=0$, то из $\xi_i\to\xi$ следует $p(\xi_i)<1$ для почти всех $i$, поэтому
$$
p(\lambda_i\cdot\xi_i-x(\pi(\xi_i)))\le \underbrace{\abs{\lambda_i}}_{\scriptsize\begin{matrix}
\text{\rotatebox{90}{$>$}}\\ \delta \end{matrix}}\cdot \underbrace{p(\xi_i)}_{\scriptsize\begin{matrix}
\text{\rotatebox{90}{$>$}} \\ 1\end{matrix}} +\underbrace{p(x(\pi(\xi_i)))}_{\scriptsize\begin{matrix}
\text{\rotatebox{90}{$>$}} \\ \e-\delta\end{matrix}}<\e.
$$
В любом случае для почти всех $i$ получаем
$$
\pi(\xi_i)\in U,\qquad p(\lambda_i\cdot\xi_i-x(\pi(\xi_i)))<\e.
$$
то есть $\lambda_i\cdot\xi_i\in W(x,U,p,\e)$. Это тоже означает \eqref{lambda_i-xi_i->lambda-xi}.

5. Докажем непрерывность послойного сложения. Пусть
$$
\xi_i\overset{\varXi}{\underset{i\to\infty}{\longrightarrow}} \xi,\quad \upsilon_i\overset{\varXi}{\underset{i\to\infty}{\longrightarrow}} \upsilon,\quad \pi(\xi_i)=\pi(\upsilon_i),\quad \pi(\xi)=\pi(\upsilon)
$$
Выберем какую-нибудь базисную окрестность $W(z,U,p,\e)$ точки $\xi+\upsilon$:
$$
\pi(\xi+\upsilon)\in U,\qquad p(\xi+\upsilon-z(\pi(\xi+\upsilon)))<\e.
$$
Подберем число $\sigma>0$ так, чтобы
\beq\label{sushestv-topologii-v-Xi-18}
p(\xi+\upsilon-z(\pi(\xi+\upsilon)))<\e-2\sigma.
\eeq
В слое $\pi^{-1}(\pi(\xi+\upsilon))=\pi^{-1}(\pi(\xi))=\pi^{-1}(\pi(\upsilon))$ операция сложения непрерывна, поэтому существуют базисные окрестности $W(x,V_x,q,\delta)$ и $W(y,V_y,r,\delta)$ точек $\xi$ и $\upsilon$ такие, что
$$
\forall \xi'\in W(x,V_x,q,\delta)\cap \pi^{-1}(\pi(\xi))\quad \forall \upsilon'\in W(y,V_y,r,\delta)\cap \pi^{-1}(\pi(\upsilon))\quad
\xi'+\upsilon'\in W(z,U,p,\e-2\sigma).
$$
То есть
\begin{multline} \label{sushestv-topologii-v-Xi-13}
\forall \xi',\upsilon'\in \pi^{-1}(\pi(\xi+\upsilon))\quad
q(\xi'-x(\pi(\xi)))<\delta\quad\&\quad r(\upsilon'-y(\pi(\upsilon)))<\delta \quad\Longrightarrow \\ \Longrightarrow\quad
p(\xi'+\upsilon'-z(\pi(\xi+\upsilon)))<\e-2\sigma.
\end{multline}
При этом из включений $\xi\in W(x,V_x,q,\delta)$ и $\upsilon\in W(y,V_y,r,\delta)$ следует во-первых, что
 \beq\label{sushestv-topologii-v-Xi-17}
q(\xi-x(\pi(\xi)))<\delta\quad\&\quad r(\upsilon-y(\pi(\upsilon)))<\delta,
 \eeq
и во-вторых, что
$$
\pi(\xi)\in V_x,\qquad \pi(\upsilon)\in V_y,
$$
то есть $V_x$ и $V_y$ являются окрестностями для $\pi(\xi)$ и $\pi(\upsilon)$ соответственно. Уменьшив, если необходимо, эти окрестности $V_x$ и $V_y$, можно считать, что они совпадают и лежат в $U$:
$$
V_x=V_y=V\subseteq U.
$$

Кроме того, очевидно, можно выбрать полунормы $q$ и $r$ так, чтобы они мажорировали $p$ (всюду на $\varXi$)
\beq\label{sushestv-topologii-v-Xi-14}
p\le q, \qquad p\le r,\qquad
\eeq
а число $\delta$ так, чтобы
\beq\label{sushestv-topologii-v-Xi-15}
\delta<\frac{\sigma}{2}.
\eeq

Теперь рассмотрим множество
 \beq\label{sushestv-topologii-v-Xi-16}
O=\left\{s\in V:\quad p(x(s)+y(s)-z(s))<\e-\sigma \right\}.
 \eeq
(оно открыто в $M$ в силу условия (iii)).

Открытое множество $W(x,O,q,\delta)$ является окрестностью точки $\xi$, потому что, во-первых,
 $$
p(x(\pi(\xi))+y(\underbrace{\pi(\xi)}_{\scriptsize\begin{matrix}
\|
 \\
\pi(\upsilon)
 \end{matrix}})-z(\underbrace{\pi(\xi)}_{\scriptsize\begin{matrix}
\|
 \\
\pi(\xi+\upsilon)
 \end{matrix}})\le \underbrace{p(x(\pi(\xi))-\xi)}_{\scriptsize\begin{matrix}
\phantom{\tiny{\eqref{sushestv-topologii-v-Xi-14}}}\quad\text{\rotatebox{90}{$\ge$}}\quad\tiny{\eqref{sushestv-topologii-v-Xi-14}}
 \\
q(x(\pi(\xi))-\xi)
 \\
\phantom{\tiny{\eqref{sushestv-topologii-v-Xi-17}}}\quad\text{\rotatebox{90}{$>$}}\quad\tiny{\eqref{sushestv-topologii-v-Xi-17}}
 \\
 \delta
 \\
\phantom{\tiny{\eqref{sushestv-topologii-v-Xi-15}}}\quad\text{\rotatebox{90}{$>$}}\quad\tiny{\eqref{sushestv-topologii-v-Xi-15}}
 \\
\frac{\sigma}{2}
 \end{matrix}}
+\underbrace{p(y(\pi(\upsilon))-\upsilon)}_{\scriptsize\begin{matrix}
\phantom{\tiny{\eqref{sushestv-topologii-v-Xi-14}}}\quad\text{\rotatebox{90}{$\ge$}}\quad\tiny{\eqref{sushestv-topologii-v-Xi-14}}
 \\
r(x(\pi(\upsilon))-\upsilon)
 \\
\phantom{\tiny{\eqref{sushestv-topologii-v-Xi-17}}}\quad\text{\rotatebox{90}{$>$}}\quad\tiny{\eqref{sushestv-topologii-v-Xi-17}}
 \\
\delta \\
\phantom{\tiny{\eqref{sushestv-topologii-v-Xi-15}}}\quad\text{\rotatebox{90}{$>$}}\quad\tiny{\eqref{sushestv-topologii-v-Xi-15}}
 \\
\frac{\sigma}{2}
\end{matrix}}
+\underbrace{p(\xi+\upsilon-z(\pi(\xi+\upsilon)))}_{\scriptsize\begin{matrix}
\phantom{\tiny{\eqref{sushestv-topologii-v-Xi-18}}}\quad\text{\rotatebox{90}{$>$}}\quad\tiny{\eqref{sushestv-topologii-v-Xi-18}}
 \\
\e-2\sigma
 \end{matrix}} <\e-\sigma
 $$
и значит, $\pi(\xi)\in O$. А, во-вторых, в силу \eqref{sushestv-topologii-v-Xi-17}, $q(\xi-x(\pi(\xi)))<\delta$.

Точно так же $\upsilon\in W(y,O,r,\delta)$.

Покажем теперь, что окрестности $W(x,O,q,\delta)$ и $W(y,O,r,\delta)$ точек $\xi$ и $\upsilon$ удовлетворяют условию
 \beq\label{sushestv-topologii-v-Xi-19}
\forall \xi'\in W(x,O,q,\delta)\quad \forall \upsilon'\in W(y,O,r,\delta)\quad \Big( \pi(\xi')=\pi(\upsilon')\quad\Longrightarrow\quad
\xi'+\upsilon'\in W(z,U,p,\e)\Big).
 \eeq
Действительно, из $\xi'\in W(x,O,q,\delta)$, $\upsilon'\in W(y,O,r,\delta)$ и $\pi(\xi')=\pi(\upsilon')$ следует, во-первых, что
$$
\pi(\xi')=\pi(\upsilon')\in O\subseteq V\subseteq U.
$$
И, во-вторых,
$$
p(\xi'+\upsilon'-z(\pi(\xi'+\upsilon')))\le
\underbrace{p(\xi'-x(\pi(\xi')))}_{\scriptsize\begin{matrix}
\phantom{\tiny{\eqref{sushestv-topologii-v-Xi-14}}}\quad\text{\rotatebox{90}{$\ge$}}\quad\tiny{\eqref{sushestv-topologii-v-Xi-14}}
 \\
q(\xi'-x(\pi(\xi')))
 \\
\text{\rotatebox{90}{$>$}}
 \\
 \delta
 \\
\phantom{\tiny{\eqref{sushestv-topologii-v-Xi-15}}}\quad\text{\rotatebox{90}{$>$}}\quad\tiny{\eqref{sushestv-topologii-v-Xi-15}}
 \\
\frac{\sigma}{2}
 \end{matrix}}
+\underbrace{p(\upsilon'-y(\pi(\upsilon')))}_{\scriptsize\begin{matrix}
\phantom{\tiny{\eqref{sushestv-topologii-v-Xi-14}}}\quad\text{\rotatebox{90}{$\ge$}}\quad\tiny{\eqref{sushestv-topologii-v-Xi-14}}
 \\
r(\upsilon'-y(\pi(\upsilon')))
 \\
\text{\rotatebox{90}{$>$}}
 \\
\delta \\
\phantom{\tiny{\eqref{sushestv-topologii-v-Xi-15}}}\quad\text{\rotatebox{90}{$>$}}\quad\tiny{\eqref{sushestv-topologii-v-Xi-15}}
 \\
\frac{\sigma}{2}
\end{matrix}}
+\underbrace{p(x(\xi')+y(\upsilon')-z(\pi(\xi'+\upsilon')))}_{\scriptsize\begin{matrix}
\phantom{\tiny{\eqref{sushestv-topologii-v-Xi-16}}}\quad\text{\rotatebox{90}{$>$}}\quad\tiny{\eqref{sushestv-topologii-v-Xi-16}}
 \\
\e-\sigma
 \end{matrix}} <\e
$$
Теперь мы получаем цепочку
$$
\xi_i\overset{\varXi}{\underset{i\to\infty}{\longrightarrow}} \xi,\quad \upsilon_i\overset{\varXi}{\underset{i\to\infty}{\longrightarrow}} \upsilon,\quad \pi(\xi_i)=\pi(\upsilon_i),\quad \pi(\xi)=\pi(\upsilon)
$$
$$
\Downarrow
$$
$$
\xi_i\in W(x,O,q,\delta),\quad \upsilon_i\in W(y,O,r,\delta)\qquad \text{для почти всех $i$}
$$
$$
\phantom{\scriptsize\text{\eqref{sushestv-topologii-v-Xi-19}}}\quad\Downarrow\quad{\scriptsize\text{\eqref{sushestv-topologii-v-Xi-19}}}
$$
$$
\xi_i+\upsilon_i\in W(z,U,p,\e)\qquad \text{для почти всех $i$}
$$
Здесь $W(z,U,p,\e)$ выбиралась как произвольная базисная окрестность точки $\xi+\upsilon$. Поэтому
$$
\xi_i+\upsilon_i\overset{\varXi}{\underset{i\to\infty}{\longrightarrow}} \xi+\upsilon.
$$

6. Пусть $p$ -- произвольная полунорма из $\mathcal P$ и $\e>0$. Покажем, что множество $W=\{\xi\in\Xi: \ p(\xi)<\e\}$ открыто в $\Xi$. Зафиксируем какую-нибудь точку $\xi\in W$. Из условия $p(\xi)<\e$ следует, что найдется $\sigma>0$ такое, что $p(\xi)<\e-2\sigma$. В силу (iv), найдется $x\in X$ такой, что $p(\xi-x(\pi(\xi)))<\sigma$. Положим
$$
O=\{t\in M:\ p(x(t))<\e-\sigma\}.
$$
Тогда базисная окрестность $W(x,O,p,\sigma)$ содержит точку $\xi$, потому что, во-первых,
$$
p(x(\pi(\xi)))\le \underbrace{p(x(\pi(\xi))-\xi)}_{\scriptsize\begin{matrix}
\text{\rotatebox{90}{$>$}} \\
\sigma
 \end{matrix}}+\underbrace{p(\xi)}_{\scriptsize\begin{matrix}
\text{\rotatebox{90}{$>$}}
 \\
\e-2\sigma
 \end{matrix}}<\e-\sigma
 \quad\Longrightarrow\quad \pi(\xi)\in O,
$$
и, во-вторых, в силу выбора $x$, выполняется $p(\xi-x(\pi(\xi)))<\sigma$. То есть $W(x,O,p,\sigma)$ -- окрестность точки $\xi$.

С другой стороны, $W(x,O,p,\sigma)$ содержится в множестве $W$, потому что если $\upsilon\in W(x,O,p,\sigma)$, то
$$
p(\upsilon)\le \underbrace{p(\upsilon-x(\pi(\upsilon)))}_{\scriptsize\begin{matrix}
\text{\rotatebox{90}{$>$}} \\
\sigma
 \end{matrix}}+\underbrace{x(\pi(\upsilon))}_{\scriptsize\begin{matrix}
\text{\rotatebox{90}{$>$}}
 \\
\e-\sigma
 \end{matrix}}<\e.
$$

7. Покажем, что выполняется условие (e) на с.\pageref{uslovie-e-DEF-rassloeniya}:
для любой точки $t\in M$ и всякой окрестности $V$ точки $0_t$ в $\Xi$ найдутся полунорма $p\in{\mathcal P}$, число $\sigma>0$ и открытое множество $O$ в $M$, содержащее $t$, такие, что
$$
\{ \xi\in \pi^{-1}(O): \ p(\xi)<\sigma\}\subseteq V.
$$
Поскольку топология в $\Xi$ порождена окрестностями \eqref{baza-topologii-v-varXi}, существует некая базисная окрестность $W(x,U,p,\e)$ точки $0_t$, содержащаяся в $V$:
$$
0_t\in W(x,U,p,\e)\subseteq V.
$$
Это значит выполнение двух условий:
$$
t\in U\quad\&\quad p(0_t-x(t))<\e.
$$
Подберем $\sigma>0$ так, чтобы
$$
p(0_t-x(t))<\e-\sigma
$$
и положим
$$
O=\{s\in U:\ p(x(s))<\e-\sigma\}.
$$
В силу условия (iii) это будет открытое множество в $M$. Оно содержит точку $t$, потому что
$$
p(x(t))\le \underbrace{p(x(t)-0_t)}_{\scriptsize\begin{matrix}
\text{\rotatebox{90}{$>$}}
 \\
\e-\sigma
 \end{matrix}}+\underbrace{p(0_t)}_{\scriptsize\begin{matrix}
\|
 \\
0
 \end{matrix}}<\e-\sigma
$$
Заметим, что если $\xi\in \pi^{-1}(O)$ и $p(\xi)<\sigma$, то
$$
p(\xi-x(\pi(\xi)))\le \underbrace{p(\xi)}_{\scriptsize\begin{matrix}
\text{\rotatebox{90}{$>$}}
 \\
\sigma
 \end{matrix}}+\kern-25pt\underbrace{p(x(\pi(\xi)))}_{\scriptsize\begin{matrix}
\phantom{\tiny{(\pi(\xi)\in O)}}\quad\text{\rotatebox{90}{$>$}}\quad\tiny{(\pi(\xi)\in O)}
 \\
\e-\sigma
 \end{matrix}}\kern-25pt<\e
$$
Теперь мы получаем цепочку
$$
\{ \xi\in \pi^{-1}(O): \ p(\xi)<\sigma\}\subseteq W(x,O,p,\e)\subseteq W(x,U,p,\e)\subseteq V.
$$

8. Докажем единственность такой топологии. Она следует из того, что в ней
сходимость направленности
\beq\label{PROOF:edinstvennost-topologii-na-rassloenii}
\xi_i\overset{\varXi}{\underset{i\to\infty}{\longrightarrow}} \xi \eeq
однозначно определяется поведением $\pi$, ${\mathcal P}$ и $X$ на точках
$\xi_i$ и $\xi$, точнее, следующими двумя условиями:

\bit{

\item[(a)] $\pi(\xi_i)\overset{M}{\underset{i\to\infty}{\longrightarrow}} \pi(\xi)$,

\item[(b)] для любого сечения $x\in X$, любой полунормы $p\in{\mathcal P}$ и любого $\e>0$ условие
$$
p(\xi_i-x(\pi(\xi_i)))<p(\xi-x(\pi(\xi)))+\e
$$
выполняется для почти всех индексов $i$.
}\eit
Действительно, если выполняется \eqref{PROOF:edinstvennost-topologii-na-rassloenii}, то условие (a) будет просто следствием непрерывности отображения $\pi:\varXi\to M$. А условие (b) доказывается так: из $\xi_i\overset{\varXi}{\underset{i\to\infty}{\longrightarrow}} \xi$ следует $x(\pi(\xi_i))\overset{\varXi}{\underset{i\to\infty}{\longrightarrow}} x(\pi(\xi))$, и вместе, в силу полунепрерывности сверху полунормы $p$, это дает цепочку неравенств, справедливую для почти всех $i$:
$$
p(\xi_i-x(\pi(\xi_i)))\le \kern-10pt\underbrace{p(\xi_i-\xi)}_{\scriptsize\begin{matrix}\text{\rotatebox{90}{$>$}}\\ \frac{\e}{2}\\ \text{для почти всех $i$}\end{matrix}}\kern-10pt+p(\xi-x(\pi(\xi)))+\underbrace{p(x(\pi(\xi))-x(\pi(\xi_i)))}_{\scriptsize\begin{matrix}\text{\rotatebox{90}{$>$}}\\ \frac{\e}{2}\\ \text{для почти всех $i$}\end{matrix}}<\frac{\e}{2}+p(\xi-x(\pi(\xi)))+\frac{\e}{2}.
$$
Наоборот, пусть выполняются условия (a) и (b). Тогда по условию (iii), для любой полунормы $p\in{\mathcal P}$ и любого $\e>0$ можно выбрать сечение $x\in X$ такое, что
\beq\label{PROOF:edinstvennost-topologii-na-rassloenii-1}
p(\xi-x(\pi(\xi)))<\e.
\eeq
После этого, мы получим во-первых,
$$
x(\pi(\xi_i))\overset{\varXi}{\underset{i\to\infty}{\longrightarrow}} x(\pi(\xi))
$$
(в силу условия (a) и непрерывности отображения $x$), во-вторых,  для всех $i$
$$
\pi(x(\pi(\xi)))=\pi(\xi),\qquad \pi(x(\pi(\xi_i)))=\pi(\xi_i)
$$
(потому что $x$ -- сечение для $\pi$), и, в-третьих, для почти всех $i$
$$
p(\xi_i-x(\pi(\xi_i)))\le \underbrace{p(\xi-x(\pi(\xi)))}_{\scriptsize\begin{matrix}\phantom{\tiny{\eqref{PROOF:edinstvennost-topologii-na-rassloenii-1}}}\ \text{\rotatebox{90}{$>$}}\ \tiny{\eqref{PROOF:edinstvennost-topologii-na-rassloenii-1}}\\ \e\end{matrix}}+\e<2\e.
$$
(в силу условия (b)).
Вместе все это означает, что точки $\zeta_i=x(\pi(\xi_i))$, $\zeta=x(\pi(\xi))$  удовлетворяют условию \eqref{trebovanie-dlya-shod-napravl} (в котором $\e$ заменено на $2\e$), и поскольку $p\in{\mathcal P}$ и $\e>0$ выбирались произвольными, по предложению \ref{PROP:harakt-shod-v-Xi} это влечет соотношение \eqref{PROOF:edinstvennost-topologii-na-rassloenii}.
\epr

\paragraph{Морфизмы расслоений.}

Пусть даны два локально выпуклых расслоения $\pi:\varXi\to M$, $\rho:\varOmega\to N$ и два непрерывных
отображения $\mu:\varXi\to\varOmega$ и $\sigma:M\to N$, замыкающие диаграмму
\beq\label{predmorfizm-rassloenij}
 \xymatrix @R=2pc @C=3pc
 {
 \varXi\ar[r]^{\mu}\ar[d]_{\pi} & \varOmega\ar[d]^{\rho}\\
   M\ar[r]^{\sigma} & N
 }
\eeq
Тогда для всякой точки $t\in M$ и любого $\xi\in\pi^{-1}(t)$ мы получим
$$
\rho\big(\mu(\xi)\big)=\sigma(\pi(\xi))=\sigma(t),
$$
и это означает, что отображение $\mu$ переводит слой $\pi^{-1}(t)$ в слой $\rho^{-1}(\sigma(t))$:
$$
\mu\Big(\pi^{-1}(t)\Big)\subseteq \rho^{-1}(\sigma(t)).
$$

Условимся называть {\it морфизмом} локально выпуклого расслоения $\pi:\varXi\to M$ в локально
выпуклое расслоение $\rho:\varOmega\to N$ всякую пару непрерывных
отображений $\mu:\varXi\to\varOmega$ и $\sigma:M\to N$, замыкающее диаграмму \eqref{predmorfizm-rassloenij}
и линейное на каждом слое, то есть такое, что для всякой точки $t\in M$
порождаемое отображение слоев
$$
\mu:\pi^{-1}(t)\to\rho^{-1}(\sigma(t))
$$
линейно (и непрерывно в силу непрерывности $\mu$).

\paragraph{Двойственное расслоение.}

Всякому локально выпуклому расслоению $\pi:\varXi\to M$ над $\C$ можно поставить в соответствие расслоение векторных пространств над $\C$
$$
\pi^\star:\varXi^\star\to M\qquad\Big|\qquad \varXi^\star=\bigsqcup_{t\in M}\pi^{-1}(t)^\star,\qquad \forall u\in \pi^{-1}(t)^\star\quad \pi^\star(u)=t.
$$
Мы будем называть его {\it сопряженным расслоением векторных пространств} к расслоению $\pi:\varXi\to M$.
Для всякого компакта $K\subset \varXi$ и любой точки $u\in \varXi^\star$ положим
$$
p_K(u)=\sup_{\xi\in K\cap\pi^{-1}\big(\pi^\star(u)\big)}\abs{u(\xi)}.
$$

Рассмотрим локально выпуклое расслоение, называемое {\it тривиальным} со слоем $\C$:
$$
\pi_M:\C\times M\to M\qquad\Big|\qquad \pi_M(\lambda,t)=t,\qquad \lambda\in\C,\ t\in M.
$$
Обозначим символом $\pi_{\C}$ проекцию на первую компоненту:
$$
\pi_{\C}:\C\times M\to\C\qquad\Big|\qquad \pi_{\C}(\lambda,t)=\lambda,\qquad  \lambda\in\C,\ t\in M.
$$
Слоем $\pi_M^{-1}(t)$ всякой точки $t\in M$ в этом расслоении является множество $\C\times\{t\}$, которое при отображении $\pi_{\C}$ отождествляется с полем $\C$. Как следствие, если $\mu:\varXi\to \C\times M$ -- морфизм расслоений, то на каждом слое определена композиция
$$
\pi^{-1}(t)\overset{\mu}{\longrightarrow}\pi_M^{-1}(t)=\C\times\{t\}\overset{\pi_{\C}}{\longrightarrow}\C,
$$
являющаяся линейным непрерывным функционалом на $\pi^{-1}(t)$. Мы можем сделать вывод, что формула
\beq\label{x(t)=pi_C-circ-mu}
x(t)=\pi_{\C}\circ\mu\Big|_{\pi^{-1}(t)},\qquad t\in M,
\eeq
определяет некое сечение $x:M\to \varXi^\star$ сопряженного расслоения векторных пространств $\pi^\star:\varXi^\star\to M$. Обозначим через $X$ множество всех таких сечений.

Рассмотрим теперь в каждом слое $\pi^{-1}(t)^\star$ множество $X_t$ функционалов $u$, представимых в виде $u=\pi_{\C}\circ\mu\Big|_{\pi^{-1}(t)}$ для некоторого морфизма расслоений $\mu:\varXi\to\C\times M$, и пусть $\overline{X_t}$ обозначает замыкание этого пространства в пространстве $\pi^{-1}(t)^\star$ (относительно топологии, порожденной полунормами ${\mathcal P}=\{p_K;\ K\subseteq \varXi\}$).

\bprop\label{PROP:o-dvoistvennom-rasloenii}
Если база $M$ является хаусдорфовым пространством, то на подрасслоении
$$
\pi_\star:\varXi_\star=\bigsqcup_{t\in M}\overline{X_t}\to M
$$
множество сечений $X$ и множество полунорм $\mathcal P$ задают топологию, превращающую $\varXi_\star$ в локально выпуклое расслоение в соответствии с предложением \ref{PROP:sushestv-topologii-v-Xi}.
\eprop

\bit{

\item[$\bullet$] Расслоение $\pi_\star:\varXi_\star\to M$ мы будем называть {\it двойственным расслоением} к расслоению $\pi:\varXi\to M$.
}\eit

\bpr
Здесь нужно просто проверить условия (i)-(iv) предложения \ref{PROP:sushestv-topologii-v-Xi}.

1. Полунормы $p_K$ на каждом слое $\pi^{-1}(t)^\star$ задают локально выпуклую топологию, которая будет отделима, например потому что сильнее топологии поточечной сходимости, которая задается полунормами вида $p_{\xi}$, где $\xi$ пробегает слой $\pi^{-1}(t)$.

2. Система $\mathcal P$ направлена по возрастанию: любые две полунормы $p_K$ и $p_L$ мажорируются полунормой $p_{K\cup L}$.

3. Покажем, что для любого сечения $x\in X$ и любой полунормы $p_K$ функция $t\in T\mapsto p_K(x(t))$ полунепрерывна сверху. Пусть $\e>0$, рассмотрим множество $O=\{t\in M:\ p_K(x(t))<\e\}$. Если оно пусто, то оно автоматически открыто, поэтому важно рассмотреть случай, когда оно непусто. Пусть $t_0\in O$, то есть
$$
p_K(x(t_0))=\sup_{\xi\in K\cap\pi^{-1}(t_0)}\abs{\pi_{\C}\big(\mu(\xi)\big)}<\e
$$
Рассмотрим множество $U=\{\xi\in\varXi:\ \abs{\pi_{\C}\big(\mu(\xi)\big)}<\e\}$.
Оно открыто и содержит $K\cap\pi^{-1}(t_0)$:
\beq\label{K-cap-pi^(-1)(t_0)-subseteq-U}
K\cap\pi^{-1}(t_0)\subseteq U.
\eeq
Наша задача показать, что существует окрестность $V$ точки $t_0$ такая, что
$$
K\cap\pi^{-1}(V)\subseteq U.
$$
Предположим, что это не так. Тогда для всякой окрестности $V$ точки $t_0$ найдется точка
$$
\xi_V\in K\cap\pi^{-1}(V)\setminus U.
$$
Поскольку направленность $\{\xi_V;\ V\to\{t_0\}\}$ лежит в компакте $K$, она в этом компакте должна иметь некоторую предельную точку $\xi\in K$ \cite[Теорема 3.1.23]{Engelking}. Под действием проекции $\pi$ направленность $\{\xi_V;\ V\to\{t_0\}\}$ превращается в направленность $\{\pi(\xi_V);\ V\to\{t_0\}\}$ с предельной точкой $\pi(\xi)$. При этом $\{\pi(\xi_V);\ V\to\{t_0\}\}$ сходится к $t_0$. Поскольку пространство $M$ хаусдорфово, предельная точка $\pi(\xi)$ должна совпадать с пределом $t_0$. Мы получаем, что
$$
\xi\in K\cap\pi^{-1}(t_0)\setminus U.
$$
Это противоречит \eqref{K-cap-pi^(-1)(t_0)-subseteq-U}.

4. Для всякой точки $t\in M$ множество $\{x(t);\ x\in X\}$ плотно в пространстве $\overline{X_t}$ по определению самого $\overline{X_t}$.
\epr

\subsection{Расслоение значений и морфизмы модулей}

\paragraph{Расслоение значений модуля над коммутативной инволютивной алгеброй.}

Пусть $A$ -- коммутативная инволютивная стереотипная алгебра. Для всякой точки $t\in\Spec(A)$ обозначим через $I_t$ ядро $t$:
$$
I_t=\{a\in A:\ t(a)=0\}.
$$
Пусть далее $X$ -- левый стереотипный модуль над $A$. В соответствии с \eqref{DEF:M-cdot-N}, мы обозначаем символом $I_t\cdot X$ подмодуль в $X$, состоящий из всевозможных сумм элементов вида $a\cdot x$, где $a\in I_t$ и $x\in X$:
$$
I_t\cdot X=\left\{ \sum_{i=1}^k a_i\cdot x_i ;\ a_i\in I_t, \ x_i\in X, k\in \Bbb N\right\}
$$
а $\overline{I_t\cdot X}$ -- его замыкание.
Положив
$$
\Jet^0_A X=\bigsqcup_{t\in M} (X/\overline{I_t\cdot X})^\triangledown,\qquad
\pi^0_{A,X}(x+\overline{I_t\cdot X})=t,
$$
мы получим сюръекцию $\pi^0_{A,X}: \Jet^0_A X\to\Spec(A)$. Для всякого элемента
$x\in X$ обозначим через $\jet^0(x)$ сечение сюръекции $\pi^0_{A,X}$, действующее
по формуле
$$
\jet^0(x)(t)=x+\overline{I_t\cdot X},\qquad t\in M.
$$
Обозначим через ${\mathcal P}(X)$ множество всех непрерывных полунорм $p:X\to
\R_+$ локально выпуклого пространства $X$. Всякая полунорма $p\in{\mathcal
P}(X)$ порождает полунорму $p^0$ на стереотипном фактор-пространстве
$(X/\overline{I_t\cdot X})^\triangledown$ по формуле \beq\label{norma-v-X/K_t}
p^0(x+\overline{I_t\cdot X}):=\inf_{y\in \overline{I_t\cdot
X}}p(x+y)=\inf_{y\in I_t\cdot X}p(x+y),\qquad x\in X. \eeq Тогда $p^0$ можно
считать функцией на $ \Jet^0_AX$, действие которой на каждом слое описывается
формулой \eqref{norma-v-X/K_t}.

Обозначим далее через $\sigma:A\to {\mathcal C}(M)$ естественное отображение
алгебры $A$ в алгебру непрерывных функций на своем спектре $M=\Spec(M)$:
\beq\label{A->C(M)} \sigma(a)(t)=t(a),\qquad a\in A,\ t\in M \eeq и заметим
тождество: \beq\label{j^0(a-cdot-x)=sigma(a)-cdot-j^0(x)} \jet^0(a\cdot
x)=\sigma(a)\cdot \jet^0(x),\qquad a\in A,\ x\in X. \eeq Оно доказывается
переносом в левую часть и подстановкой агрумента $t\in M$:
$$
\jet^0(a\cdot x)(t)-(\sigma(a)\cdot \jet^0(x))(t)=\underbrace{a\cdot
x+\overline{I_t\cdot X}}_{\scriptsize\begin{matrix}\text{\rotatebox{90}{$=$}}
\\ \jet^0(a\cdot
x)(t)\end{matrix}}-\underbrace{t(a)}_{\scriptsize\begin{matrix}\text{\rotatebox{90}{$=$}}
\\ \sigma(a)(t)\end{matrix}}\cdot \underbrace{(x+\overline{I_t\cdot
X})}_{\scriptsize\begin{matrix}\text{\rotatebox{90}{$=$}} \\
\jet^0(x)(t)\end{matrix}}\subseteq
\underbrace{(a-t(a))}_{\scriptsize\begin{matrix}\text{\rotatebox{90}{$\owns$}}
\\ I_t\end{matrix}}\cdot x+\overline{I_t\cdot X}\subseteq \overline{I_t\cdot X}
$$

Чтобы применить предложение \ref{PROP:sushestv-topologii-v-Xi} нам не хватает следующей леммы.

\blm\label{LM:poluneprer-sverhu-p^0(j^0(x)(t))} Для всякого элемента $x\in X$ и
любой непрерывной полунормы $p:X\to\R_+$ отображение $t\in\Spec(A)\mapsto
p^0(\jet^0(x)(t))\in\R_+$ полунепрерывно сверху. \elm \bpr Пусть $t\in\Spec(A)$ и
$\e>0$. Условие
$$
p^0(\jet^0(x)(t))=p^0(x+\overline{I_t\cdot X})=\inf_{y\in \overline{I_t\cdot
X}}p(x+y)=\inf_{y\in I_t\cdot X}p(x+y)<\e
$$
означает, что для некоторого $y\in I_t\cdot X$ выполняется неравенство
$$
p(x+y)<\e.
$$
Это в свою очередь означает, что существуют $m\in\N$, векторы $y_1,...,y_m\in X$ и векторы $a^1,...,a^m\in I_t$ такие, что
\beq\label{p-l-x+sum_(k=1)^m-l-a^k-cdot-y_k-r<e}
p\l x+\sum_{k=1}^m a^k\cdot y_k\r<\e.
\eeq
Для всякой последовательности чисел $\lambda=(\lambda^1,...,\lambda^m)$, $\lambda^k\in\C$, обозначим
$$
f(\lambda)=p\l x+\sum_{k=1}^m (a^k-\lambda^k\cdot 1_A)\cdot y_k\r.
$$
Функция $\lambda\mapsto f(\lambda)$ в точке $\lambda=0$ совпадает с левой частью \eqref{p-l-x+sum_(k=1)^m-l-a^k-cdot-y_k-r<e} и поэтому удовлетворяет неравенству
$$
f(0)<\e.
$$
С другой стороны, она непрерывна, как композиция аффинного отображения из $\C^m$ в локально выпуклое пространство $X$ и непрерывной функции $p$ на $X$. Значит, должно существовать число $\delta>0$ такое, что
\beq\label{max_(k)|lambda^k|<delta=>f(lambda)<e}
\forall \lambda \qquad \max_{1\le k\le m}\abs{\lambda^k}<\delta\quad\Longrightarrow\quad f(\lambda)<\e.
\eeq

Рассмотрим далее множество
$$
U=\{s\in\Spec(A):\ \forall k\in\{1,...,m\}\quad \abs{s(a^k)}<\delta\}.
$$
Оно открыто и содержит точку $t$ (потому что включения $a^k\in I_t$ означают систему равенств $t(a^k)=0$, $1\le k\le m$).
С другой стороны, для всякой точки $s\in U$, рассмотрев последовательность
$$
\lambda^k=s(a^k),
$$
мы получим, во-первых, $s(a^k-\lambda^k\cdot 1_A)=s(a^k)-\lambda^k\cdot s(1_A)=s(a^k)-s(a^k)=0$, то есть
$$
a^k-\lambda^k\cdot 1_A\in I_s
$$
и, во-вторых, $\max_{k}\abs{\lambda^k}=\max_{k}\abs{s(a^k)}<\delta$, то есть, в силу \eqref{max_(k)|lambda^k|<delta=>f(lambda)<e},
$$
p\bigg( x+\sum_{k=1}^m \underbrace{(a^k-\lambda^k\cdot 1_A)}_{\scriptsize\begin{matrix}\text{\rotatebox{90}{$\owns$}}\\ I_s\end{matrix}}\cdot y_k\bigg)=f(\lambda)<\e
$$
Это можно понимать так, что для некоторой точки $z\in I_s\cdot X$ выполняется неравенство
$$
p(x+z)<\e,
$$
и поэтому
$$
p^0(\jet^0(x)(s))=p^0(x+\overline{I_s\cdot X})=\inf_{z\in \overline{I_s\cdot
X}}p(x+z)=\inf_{z\in I_s\cdot X}p(x+z)<\e
$$
Это верно для всякой точки $s$ из окрестности $U$ точки $t$, и это то, что нам нужно было доказать.
\epr

\btm\label{TH:J^0_AX}
Для всякого стереотипного модуля $X$ над коммутативной инволютивной стереотипной алгеброй $A$ прямая сумма стереотипных фактор-модулей
$$
 \Jet^0_AX=\bigsqcup_{t\in \Spec(A)} (X/\overline{I_t\cdot X})^\triangledown
$$
обладает единственной топологией, превращающей проекцию
$$
\pi^0_{A,X}: \Jet^0_A X\to\Spec(A),\qquad \pi^0_{A,X}(x+\overline{I_t\cdot
X})=t,\qquad t\in\Spec(A),\ x\in X
$$
в локально выпуклое расслоение с системой полунорм $\{p^0;\ p\in{\mathcal P}(X)\}$, для которого отображение
$$
x\in X\mapsto \jet^0(x)\in\Sec(\pi^0_{A,X})\quad\Big|\quad
\jet^0(x)(t)=x+\overline{I_t\cdot X},\quad t\in\Spec(A),
$$
непрерывно переводит $X$ в стереотипный $A$-модуль $\Sec(\pi^0_{A,X})$ сечений
$\pi^0_{A,X}$. При этом

 \bit{
\item[(i)] базу топологии $\Jet^0_AX$ образуют множества
$$
W(x,U,p,\e)=\set{\xi\in  \Jet^0_AX: \ \pi^0_{A,X}(\xi)\in U\ \& \
p^0\Big(\xi-\jet^0(x)\big(\pi^0_{A,X}(\xi)\big)\Big)<\e }
$$
где $x\in X$, $p\in{\mathcal P}(X)$, $\e>0$, $U$ -- открытое множество в $\Spec(A)$;

\item[(ii)] при фиксированных $x\in X$ и $s\in M$ локальную базу топологии $\Jet^0_AX$ в точке $\zeta=\jet^0(x)(s)$ образуют множества
$$
W(x,U,p,\e)=\set{\xi\in  \Jet^0_AX: \ \pi^0_{A,X}(\xi)\in U\ \& \
p^0\Big(\xi-\jet^0(x)\big(\pi^0_{A,X}(\xi)\big)\Big)<\e }
$$
где $p\in{\mathcal P}(X)$, $\e>0$, $U$ -- окрестность точки $s=\pi^0(\zeta)$ в $\Spec(A)$;

\item[(iii)] если вдобавок спектр $\Spec(A)$ является паракомпактным локально
компактным пространством, то отображение $\jet^0:X\to\Sec(\pi^0_{A,X})$ имеет
плотный образ. }\eit \etm

\bit{

\item[$\bullet$]\label{DEF:rassloenie-znachenij} Расслоение $\pi^0_A: \Jet^0_AX\to \Spec(A)$ называется {\it
расслоением значений} модуля $X$ над алгеброй $A$.

}\eit

\bpr
(i). Обозначим $M=\Spec(A)$. Из леммы \ref{LM:poluneprer-sverhu-p^0(j^0(x)(t))} и предложения
\ref{PROP:sushestv-topologii-v-Xi} следует существование и единственность
топологии на $ \Jet^0_AX$, для которой проекция $\pi^0: \Jet^0_A X\to\Spec(A)$
представляет собой локально выпуклое расслоение с полунормами $p^0$, а сечения
вида $\jet^0(x)$, $x\in X$, будут непрерывными. Непрерывность отображения $x\in
X\mapsto \jet^0(x)\in\Sec(\pi^0_{A,X})$ доказывается импликацией
$$
p^0(\jet^0(x)(t))=\inf_{y\in I_t\cdot X}p(x+y)\le p(x)\quad \Longrightarrow \quad
p^0_T(\jet^0(x)(t))=\sup_{t\in T}p^0(\jet^0(x)(t))\le p(x)
$$
для всякого компакта $T\subseteq M$. Свойство (i) также следует из предложения
\ref{PROP:sushestv-topologii-v-Xi}.

Докажем (ii). Пусть $x\in X$, $s\in M$ и $\zeta=\jet^0(x)(s)$. Рассмотрим какую-нибудь базисную окрестность $W(y,V,p,\delta)$ точки $\zeta$, то есть $s\in V$ и $p^0(\zeta-\jet^0(y)(s))<\delta$. Выберем $\e>0$ так, чтобы
$$
p^0(\zeta-\jet^0(y)(s))<\delta-\e.
$$
По предложению \ref{PROP:sushestv-topologii-v-Xi} (iii), множество
$$
U=\{t\in V: p^0(\jet^0(x)(t)-\jet^0(y)(t))<\delta-\e\}
$$
является окрестностью точки $s$ в $M$. Теперь мы получим, что $\zeta\in W(x,U,p,\e)\subseteq W(y,V,p,\delta)$. Действительно, если
$\xi\in W(x,U,p,\e)$, то, во-первых, $\pi^0(\xi)=t\in U\subseteq V$, и, во-вторых, $p^0(\xi-\jet^0(x)(t))<\e$, поэтому
$$
p^0(\xi-\jet^0(y)(t))\le p^0(\xi-\jet^0(x)(t))+p^0(\jet^0(x)(t)-\jet^0(y)(t))<\e+\delta-\e=\delta.
$$

Докажем (iii). Заметим сначала, что множество $\jet^0(X)=\{\jet^0(x);\ x\in X\}$
плотно в порожденном им ${\mathcal C}(\Spec(A))$-модуле
$$
{\mathcal C}(M)\cdot \jet^0(X)=\set{\sum_{i=1}^k b_i\cdot \jet^0(x_i);\ b_i\in {\mathcal C}(M),\ x_i\in X}.
$$
Для этого зафиксируем $b\in {\mathcal C}(M)$ и $x\in X$. Поскольку $A$ -- инволютивная
алгебра, отображение $\sigma:A\to {\mathcal C}(M)$ должно иметь плотный образ в
${\mathcal C}(M)$. Поэтому найдется направленность $a_i\in A$ такая, что
$$
\sigma(a_i)\overset{{\mathcal C}(M)}{\underset{i\to\infty}{\longrightarrow}}b.
$$
Как следствие,
$$
\jet^0(a_i\cdot x)=\eqref{j^0(a-cdot-x)=sigma(a)-cdot-j^0(x)}=\sigma(a_i)\cdot
\jet^0(x)\overset{{\mathcal C}(M)}{\underset{i\to\infty}{\longrightarrow}}b\cdot \jet^0(x).
$$
Далее заметим, что поскольку образ $X$ при каждой проекции $\pi^0_t$ плотен в
фактор-пространстве $(X/\overline{I_t\cdot X})^\triangledown$, множество
$\jet^0(X)(\Spec(A))$ будет плотно в каждом слое расслоения $ \Jet^0_AX$. Поэтому его
надмножество $({\mathcal C}(M)\cdot \jet^0(X))(\Spec(A))$ также  будет плотно в каждом слое
расслоения $ \Jet^0_AX$. Теперь мы можем применить свойство $2^\circ$ на
с.\pageref{LM:poltnost-v-sloyah=>plotnost-v-secheniyah}: если $M$ паракомпактно
и локально компактно, то ${\mathcal C}(M)$-модуль ${\mathcal C}(M)\cdot \jet^0(X)$ плотен в
$\Sec(\pi^0_{A,X})$ (а, как мы уже поняли $\jet^0(X)$ плотен в ${\mathcal C}(M)\cdot
\jet^0(X)$). \epr

\btm\label{TH:Sec(J^0X)-kak-A-modul}
Пусть $X$ -- левый стереотипный модуль над коммутативной инволютивной стереотипной алгеброй $A$ с паракомпактным локально компактным инволютивным спекром $\Spec(A)$. Тогда формула
\beq\label{DEF:Sec(pi^0_(A,X))-kak-A-modul}
(a\cdot x)(t)=t(a)\cdot x(t),\qquad x\in \Sec(\pi^0_{A,X}),\quad a\in A,\quad t\in\Spec(A),
\eeq
наделяет пространство $\Sec(\pi^0_{A,X})$ сечений расслоения значений структурой левого стереотипного $A$-модуля, для которой отображение $\jet^0:X\to\Sec(\pi^0_{A,X})$ является морфизмом стереотипных $A$-модулей.
\etm
\bpr
1. Покажем сначала, что формула \eqref{DEF:Sec(pi^0_(A,X))-kak-A-modul} наделяет $\Sec(\pi^0_{A,X})$ структурой  стереотипного $A$-модуля. Пусть
$$
a_i\overset{A}{\underset{i\to\infty}{\longrightarrow}}0.
$$
Тогда для любого компакта $T\subseteq\Spec(A)$ мы получаем равномерное по $t\in T$ стремление к нулю в $\C$,
$$
t(a_i)\overset{\C}{\underset{i\to\infty, t\in T}{\rightrightarrows}}0,
$$
и поэтому для любого компакта $K\subseteq\Sec(\pi^0_{A,X})$
$$
(a_i\cdot x)(t)=t(a_i)\cdot x(t)\overset{X}{\underset{i\to\infty, t\in T, x\in K}{\rightrightarrows}}0.
$$
Наоборот, если
$$
x_i\overset{\Sec(\pi^0_{A,X})}{\underset{i\to\infty}{\longrightarrow}}0,
$$
то для любых компактов $K\subseteq A$ и $T\subseteq\Spec(A)$ множество $\{t(a);\ t\in T, a\in K\}$ будет компактом в $\C$, поэтому
$$
(a\cdot x_i)(t)=t(a)\cdot x_i(t)\overset{X}{\underset{i\to\infty, t\in T, x\in K}{\rightrightarrows}}0.
$$

2. Теперь убедимся, что отображение $\jet^0:X\to\Sec(\pi^0_{A,X})$ является морфизмом $A$-модулей. Для $a\in A$ и $x\in X$ получаем:
\begin{multline*}
\jet^0(a\cdot x)(t)-\big(a\cdot \jet^0(x)\big)(t)=a\cdot x+\overline{I_t\cdot X}-t(a)\cdot \Big(x+\overline{I_t\cdot X}\Big)=a\cdot x-t(a)\cdot x+\overline{I_t\cdot X}-t(a)\cdot\overline{I_t\cdot X}=\\=
\underbrace{\big(a-t(a)\big)}_{\scriptsize\begin{matrix}\text{\rotatebox{90}{$\owns$}}\\ I_t \end{matrix}}\cdot x+\overline{I_t\cdot X}-t(a)\cdot\overline{I_t\cdot X}\subseteq \overline{I_t\cdot X}+
\overline{I_t\cdot X}-t(a)\cdot\overline{I_t\cdot X}\subseteq \overline{I_t\cdot X}
\end{multline*}
и поэтому
$$
\jet^0(a\cdot x)(t)=\big(a\cdot \jet^0(x)\big)(t).
$$
\epr

\paragraph{Морфизмы модулей и их связь с морфизмами расслоений значений.}

\btm\label{TH:morf-modul->morf-rassl-znach} Всякий морфизм стереотипных модулей
$D:X\to Y$ над коммутативной инволютивной алгеброй $A$ определяет единственный
морфизм расслоений значений $\jet^0(D): \Jet^0_AX\to \Jet^0_A(Y)$,
$$
 \xymatrix @R=2pc @C=1.2pc
 {
  \Jet^0_AX\ar[rr]^{\jet^0(D)}\ar[dr]_{\pi^0_{A,X}} & &  \Jet^0_AY\ar[dl]^{\pi^0_{A,Y}}\\
  & \Spec(A) &
 }
$$
удовлетворяющий тождеству
 \beq\label{j^0(Dx)=j^0(D)-circ-j^0(x)}
\jet^0(Dx)=\jet^0(D)\circ \jet^0(x),\qquad x\in X.
 \eeq
$$
 \xymatrix @R=2pc @C=1.2pc
 {
  \Jet^0_AX\ar[rr]^{\jet^0(D)} & &  \Jet^0_AY\\
  & \Spec(A)\ar[ul]^{\jet^0(x)}\ar[ur]_{\jet^0(Dx)} &
 }
$$
\etm
\bpr
Из очевидного вложения
\beq\label{D((I_t-X))-subseteq-(I_t-Y)}
D\l \overline{I_t\cdot X}\r\subseteq \overline{I_t\cdot Y}
\eeq
следует существование естественного отображения фактор-пространств:
$$
X/\ \overline{I_t\cdot X}\owns x+\overline{I_t\cdot X}\mapsto Dx+\overline{I_t\cdot Y}\in Y/\ \overline{I_t\cdot Y}
$$
Оно непрерывно, поскольку исходное отображение $D$ непрерывно, и значит, существует естественное (также, непрерывное) отображение стереотипных фактор-пространств (то есть псевдопополнений обычных фактор-пространств):
$$
\jet^0(D):\Big( X/\ \overline{I_t\cdot X}\Big)^\triangledown\to \Big( Y/\
\overline{I_t\cdot Y}\Big)^\triangledown
$$
Это верно для всякого $t\in\Spec(A)$, поэтому возникает отображение прямых сумм:
$$
\jet^0(D):\bigsqcup_{t\in\Spec(A)}\Big( X/\ \overline{I_t\cdot
X}\Big)^\triangledown\to \bigsqcup_{t\in\Spec(A)}\Big( Y/\ \overline{I_t\cdot
Y}\Big)^\triangledown
$$

Тождество \eqref{j^0(Dx)=j^0(D)-circ-j^0(x)} проверяется прямым вычислением:
для любых $t\in\Spec(A)$ и $x\in X$
$$
\big(\jet^0(D)\circ \jet^0(x)\big)(t)=\jet^0(D)\big(\jet^0(x)(t)\big)=\jet^0(D)\l
x+\overline{I_t\cdot X}\r=Dx+\overline{I_t\cdot Y}=\jet^0(Dx)(t).
$$

Остается проверить непрерывность отображения $\jet^0(D)$. Зафиксируем точку
$\zeta=\jet^0(x)(t)\in  \Jet^0_AX$, $x\in X$, $t\in T$ и рассмотрим какую-нибудь
базисную окрестность ее образа
$\jet^0(D)(\zeta)=\jet^0(D)\big(\jet^0(x)(t)\big)=\jet^0(Dx)(t)$ при отображении $\jet^0(D)$:
$$
W(y,V,q,\e)=\set{\upsilon\in  \Jet^0_A(Y): \ \pi^0_{A,Y}(\upsilon)\in U\ \& \
q^0\Big(\upsilon-\jet^0(y)\big(\pi^0_{A,Y}(\upsilon)\big)\Big)<\e }.
$$
где $q:Y\to\R_+$ -- произвольная непрерывная полунорма на $Y$, $V$ --
окрестность точки $t$ в $\Spec(A)$, $y\in Y$, $\e>0$. Поскольку $\jet^0(Dx)(t)\in
W(y,U,q,\e)$, должно выполняться условие
$$
q^0\Big(\jet^0(Dx)(t)-\jet^0(y)(t)\Big)<\e.
$$
Из него следует, что существуют число $\delta>0$ и окрестность $U$ точки $t$, лежащая в $V$ такие, что
$$
\pi^0_{A,Y}(\upsilon)\in U\quad \& \quad
q^0\Big(\upsilon-\jet^0(Dx)\big(\pi^0_{A,Y}(\upsilon)\big)\Big)<\delta
\qquad\Longrightarrow\qquad
q^0\Big(\jet^0(Dx)\big(\pi^0_{A,Y}(\upsilon)\big)-\jet^0(y)\big(\pi^0_{A,Y}(\upsilon)\big)\Big)<\e
$$
Отсюда следует, что окрестность
$$
W(Dx,U,q,\delta)=\set{\upsilon\in  \Jet^0_A(Y): \ \pi^0_{A,Y}(\upsilon)\in U\quad \& \quad
q^0\Big(\upsilon-\jet^0(Dx)\big(\pi^0_{A,Y}(\upsilon)\big)\Big)<\delta }.
$$
содержится в окрестности $W(y,V,q,\e)$ точки $\jet^0(Dx)(t)$:
$$
W(Dx,U,q,\delta)\subseteq W(y,V,q,\e).
$$

Далее вспомним, что $D:X\to Y$ -- непрерывное линейное отображение. Как
следствие, существует непрерывная полунорма $p:X\to\R_+$, удовлетворяющая
условию
\beq\label{q(Dx)-le-p(x)} q(Dx)\le p(x),\qquad x\in X.
\eeq
Заметим, что из этого неравенства следует неравенство для полунорм на расслоениях значений:
\beq\label{q^0(j^0(D)(xi))-le-p^0(xi)} q^0\Big(\jet^0(D)(\xi)\Big)\le
p^0(\xi),\qquad t\in\Spec(A),\ \xi\in  \Jet^0_AX.
\eeq
Это достаточно доказать для точек $\xi=\jet^0(x)(t)$, $x\in X$, поскольку они плотны в каждом слое:
\begin{multline*}
q^0\Big(\jet^0(D)\big(\jet^0(x)(t)\big)\Big)=\eqref{j^0(Dx)=j^0(D)-circ-j^0(x)}=q^0\Big(\jet^0(Dx)(t)\Big)=q^0\Big(Dx+\overline{I_t\cdot
Y}\Big)= \inf_{v\in I_t\cdot Y}q(Dx+v)
\overset{\scriptsize\eqref{D((I_t-X))-subseteq-(I_t-Y)}}{\le}
 \\ \le \inf_{u\in I_t\cdot X}q(Dx+Du)=
\inf_{u\in I_t\cdot X}q\Big(D(x+u)\Big)\le\eqref{q(Dx)-le-p(x)}\le \inf_{u\in
I_t\cdot X}p(x+u)=p^0\Big(x+\overline{I_t\cdot X}\Big)=p^0\Big(\jet^0(x)(t)\Big)
\end{multline*}

Теперь можно показать, что окрестность
$$
W(x,U,p,\delta)=\set{\xi\in  \Jet^0_AX: \ \pi^0_{A,X}(\xi)\in U\ \& \
p^0\Big(\xi-\jet^0(x)\big(\pi^0_{A,X}(\xi)\big)\Big)<\delta }.
$$
точки $\zeta=\jet^0(x)(t)$ при отображении $\jet^0(D)$ переходит в окрестность
$W(Dx,U,q,\delta)$ точки $\jet^0(D)(\zeta)=\jet^0(Dx)(t)$:
\beq\label{j^0(D)(W(x,U,p,delta))-subseteq-W(Dx,U,q,delta)} \jet^0(D)\Big(
W(x,U,p,\delta) \Big)\subseteq W(Dx,U,q,\delta) \eeq Действительно, для всякой
точки $\xi\in W(x,U,p,\delta)$ условие $\pi^0_{A,X}(\xi)\in U$ оказывается
полезным в конце цепочки
$$
\pi^0_{A,Y}(\jet^0(D)(\xi))=\pi^0_{A,X}(\xi)\in U
$$
а условие $p^0\Big(\xi-\jet^0(x)\big(\pi^0_{A,X}(\xi)\big)\Big)<\delta$ в конце
цепочки
\begin{multline*}
q_{\jet^0(D)(\xi)}\Big(\jet^0(D)(\xi)-\jet^0(Dx)\big(\pi^0_{A,Y}(\jet^0(D)(\xi))\big)\Big)=
q_{\pi^0_{A,X}(\xi)}\Big(\jet^0(D)(\xi)-\jet^0(Dx)\big(\pi^0_{A,X}(\xi)\big)\Big)=\eqref{j^0(Dx)=j^0(D)-circ-j^0(x)}=\\=
q_{\pi^0_{A,X}(\xi)}\bigg(\jet^0(D)(\xi)-\jet^0(D)\Big(\jet^0(x)\big(\pi^0_{A,X}(\xi)\big)\Big)\bigg)=
q_{\pi^0_{A,X}(\xi)}\bigg(\jet^0(D)\Big(\xi-\jet^0(x)\big(\pi^0_{A,X}(\xi)\big)\Big)\bigg)\le \eqref{q^0(j^0(D)(xi))-le-p^0(xi)}\le \\
\le p^0\Big(\xi-\jet^0(x)\big(\pi^0_{A,X}(\xi)\big)\Big)<\delta.
\end{multline*}
Вместе то и другое означает, что $\jet^0(D)(\xi)\in W(Dx,U,q,\delta)$, а это и
требовалось. \epr

\paragraph{Действие гомоморфизмов в $C^*$-алгебры на спектры.}

Ниже нам понадобится следующий технический результат.

\btm\label{LM:overline(ph(I_t)-cdot-B)=B} Пусть $\ph :A\to B$ --- инволютивный гомоморфизм стереотипных алгебр, причем $A$ коммутативна, $B$ --- $C^*$-алгебра и $\ph(A)$ плотно в $B$. Тогда

\bit{

\item[(i)] для всякой точки $t\in\Spec(B)$ выполняется равенство
\beq\label{overline(ph(Ker(t-circ-ph)))=Ker_t}
\overline{\ph\big(\Ker (t\circ\ph)\big)}=\Ker t
\eeq

\item[(ii)] для всякой точки
$s\in\Spec(A)\setminus\big(\Spec(B)\circ\ph \big)$ выполняется равенство
\beq\label{overline(ph(Ker_s))=B}
\overline{\ph(\Ker s)}=B
\eeq
}\eit
\etm

Для доказательства нам понадобится следующее обозначение. Пусть $X$ -- левый стереотипный модуль над стереотипной алгеброй $A$. Если $M$
-- подпространство в $A$ и $N$ -- подпространство в $X$, то мы обозначаем
\beq\label{DEF:M-cdot-N} M\cdot N=\left\{ \sum_{i=1}^k m_i\cdot n_i ;\ m_i\in
M, n_i\in N, k\in \Bbb N\right\} \eeq

Очевидно, для любых $L,M,N\subseteq A$
\beq\label{(LM)N=L(MN)}
(L\cdot M)\cdot N=L\cdot (M\cdot N)
\eeq
и
\beq\label{(M^N)^=(MN)^}
\overline{M\cdot N}=\overline{\overline{M} \cdot N}= \overline{M\cdot
\overline{N}}=\overline{\overline{M}\cdot \overline{N}}
\eeq

\blm\label{LM:I-ideal=>ph(I)-ideal}
Пусть  $\ph :A\to B$ -- морфизм стереотипных алгебр. Тогда
\bit{
\item[(i)] для любых подмножеств $M,N\subseteq A$ выполняется включение
$$
\ph(\overline{M\cdot N})\subseteq \overline{\ph (M)\cdot \ph (N)}
$$
\item[(ii)] для всякого левого (правого) идеала $I$ в $A$ множество $\overline{\ph(I)}$ есть левый (правый) идеал в $\overline{\ph(A)}$.
}\eit
\elm
\bpr Часть (i) отмечалась в \cite[Lemma 13.10]{Ak03}, поэтому остается доказать (ii).
Пусть $b\in \overline{\ph(A)}$ и $y\in\overline{\ph(I)}$. Тогда
$$
b\underset{\infty\gets i}{\longleftarrow} \ph(a_i),\qquad
y\underset{\infty\gets j}{\longleftarrow} \ph(x_j),
$$
для некоторых направленностей $a_i\in A$ и $x_j\in I$. Отсюда
$$
b\cdot y\underset{\infty\gets i}{\longleftarrow} \ph(a_i)\cdot y\underset{\infty\gets j}{\longleftarrow} \ph(a_i)\cdot \ph(x_j)=\ph(a_i\cdot x_j).
$$
и поэтому $b\cdot y\in \overline{\ph(I)}$.
\epr

\bpr[Доказетельство теоремы \ref{LM:overline(ph(I_t)-cdot-B)=B}]
1. Пусть $t\in\Spec(B)$. Прежде всего,
$$
a\in \Ker (t\circ\ph)\quad\Rightarrow\quad t(\ph(a))=(t\circ\ph)(a)=0\quad\Rightarrow\quad \ph(a)\in\Ker t.
$$
Отсюда $\ph\big(\Ker (t\circ\ph)\big)\subseteq\Ker t$, и поэтому $\overline{\ph\big(\Ker (t\circ\ph)\big)}\subseteq\Ker t$.

Пусть наоборот, $b\in \Ker t$. Поскольку $\ph(A)$ плотно в $B$, для всякого $\e>0$ можно подобрать элемент $a_{\e}\in A$ такой, что
$$
\norm{b-\ph(a_{\e})}<\e.
$$
Положим
$$
a'_{\e}=a_{\e}-t\big(\ph(a_{\e})\big)\cdot 1_A.
$$
Тогда, во-первых,
$$
t\big(\ph(a'_{\e})\big)=t\big(\ph(a_{\e})\big)-t\big(\ph(a_{\e})\big)\cdot t\big(\ph(1_A)\big)=0.
$$
То есть $a'_{\e}\in\Ker (t\circ\ph)$. А, во-вторых,
$$
\norm{\ph(a_{\e})-\ph(a'_{\e})}=\norm{\ph(a_{\e})-\ph(a_{\e})+t\big(\ph(a_{\e})\big)\cdot\ph(1_A)}=\abs{t\big(\ph(a_{\e})\big)}=
\Big|\underbrace{t(b)}_{\scriptsize\begin{matrix}\text{\rotatebox{90}{$=$}}\\ 0\end{matrix}}-t\big(\ph(a_{\e})\big)\Big|\le \norm{t}\cdot \norm{b-\ph(a_{\e})}<\e.
$$
Отсюда
$$
\norm{b-\ph(a'_{\e})}\le \norm{b-\ph(a_{\e})}+\norm{\ph(a_{\e})-\ph(a'_{\e})}<2\e.
$$
Мы получили, что для любых $b\in \Ker t$ и $\e>0$ найдется элемент $a'_{\e}\in\Ker(t\circ\ph)$ такой, что $\norm{b-\ph(a'_{\e})}<2\e$. Это доказывает включение $\overline{\ph\big(\Ker (t\circ\ph)\big)}\supseteq\Ker t$.

2. Зафиксируем $s\in\Spec(A)\setminus\big(\Spec(B)\circ\ph \big)$. Для всякой точки $t\in \Spec(B)$ найдется элемент $x_t\in A$, отделяющий $s$ от $t\circ\ph$:
$$
s(x_t)\ne t(\ph(x_t)).
$$
Как следствие, элемент $a_t=x_t-s(x_t)\cdot 1_A$ должен обладать свойством
\beq\label{s(a_t)=0}
s(a_t)=s(x_t)-s(x_t)\cdot s(1_A)=s(x_t)-s(x_t)=0
\eeq
и свойством
\beq\label{t(ph(a_t))-ne-0}
t(\ph(a_t))=t(\ph(x_t))-s(x_t)\cdot t(\ph(1_A))=t(\ph(x_t))-s(x_t)\ne 0.
\eeq
Из \eqref{t(ph(a_t))-ne-0} следует, что множества
$$
U_t=\{r\in\Spec(A):\ r(a_t)\ne 0 \}
$$
покрывают компакт $\Spec(B)\circ\ph$. Значит, среди них имеется конечное подпокрытие $U_{t_1},...,U_{t_n}$:
$$
\bigcup_{i=1}^n U_{t_i}\supseteq \Spec(B)\circ\ph.
$$
Положим
$$
a=\sum_{i=1}^n a_{t_i}\cdot a_{t_i}^\bullet.
$$
Этот элемент в точке $s$ равен нулю, потому что в силу \eqref{s(a_t)=0}, все $a_{t_i}$ в ней равны нулю,
$$
s(a)=\sum_{i=1}^n s(a_{t_i})\cdot \overline{s(a_{t_i})}=0.
$$
Как следствие, $a\in \Ker s$, и значит,
$$
\ph(a)\in \overline{\ph(\Ker s)}.
$$
А, с другой стороны, элемент $\ph(a)$ отличен от нуля всюду на $\Spec(B)$, потому что на каждой точке $t\in\Spec(B)$ какой-нибудь элемент $\ph(a_{t_i})$ отличен от нуля:
$$
t(\ph(a))=\sum_{i=1}^n t(\ph(a_{t_i}))\cdot \overline{t(\ph(a_{t_i}))}=\sum_{i=1}^n \abs{t(\ph(a_{t_i}))}^2>0.
$$
То есть $\ph(a)$ отличен от нуля всюду на спектре $\Spec(B)$ коммутативной $C^*$-алгебры $B\cong C(\Spec(B))$, и как следствие, он обратим в $B$:
$$
b\cdot\ph(a)=1_B
$$
для некоторого $b\in B$. Поскольку здесь $\ph(a)$ лежит в левом идеале $\overline{\ph(\Ker s)}$ алгебры $B=\overline{\ph(A)}$ (по лемме \ref{LM:I-ideal=>ph(I)-ideal}(ii)), мы получаем, что $1_B$ тоже лежит в $\overline{\ph(\Ker s)}$. Поэтому
$$
\overline{\ph(\Ker s)}\supseteq B\cdot\overline{\ph(\Ker s)}\supseteq B\cdot 1_B=B.
$$
\epr

\bcor\label{COR:overline(ph(I_t)-cdot-B)=B} Пусть $\ph :A\to B$ --- инволютивный гомоморфизм инволютивных стереотипных алгебр, причем $A$ коммутативна, а $B$ --- $C^*$-алгебра. Тогда для всякой точки
$$
s\in\Spec(A)\setminus\big(\Spec\big(\overline{\ph(A)}\big)\circ\ph\big)
$$
выполняется равенство
\beq
\overline{\ph(\Ker s)\cdot B}=B
\eeq
\ecor
\bpr По теореме \ref{LM:overline(ph(I_t)-cdot-B)=B}(ii), множество $\overline{\ph(\Ker s)}$ содержит единицу алгебры $\overline{\ph(A)}$, которая совпадает с единицей алгебры $B$:
$$
\overline{\ph(\Ker s)}\owns 1_B.
$$
Поэтому
$$
\overline{\ph(\Ker s)\cdot B}\supseteq \overline{\ph(\Ker s)}\cdot B\supseteq 1_B\cdot B=B.
$$
\epr

\paragraph{Морфизмы со значениями в $C^*$-алгебре и теорема Даунса-Хофманна.}

Пусть $B$ -- подалгебра в центре стереотипной алгебры $F$, причем спектр $\Spec(B)$ -- паракомпактное локально компактное пространство. Алгебру $F$ можно считать (формально, левым) модулем над $B$. Рассмотрим расслоение значений\footnote{Расслоение значений было определено на с.\pageref{DEF:rassloenie-znachenij}.}
$\pi^0_B: \Jet^0_B F\to \Spec(B)$. Для всякой точки спектра $t\in\Spec(B)$ справедливо равенство
$$
\overline{I_t\cdot F}=\overline{F\cdot I_t},
$$
означающее, что модули $\overline{I_t\cdot F}$ являются двусторонними идеалами в $F$. Поэтому каждый слой
$$
\Big(F/\overline{I_t\cdot F}\Big)^\triangledown
$$
является стереотипной алгеброй, а проекция  $F\to \Big(F/\overline{I_t\cdot F}\Big)^\triangledown$ --
гомоморфизмом стереотипных алгебр. Отсюда следует, что пространство непрерывных сечений
$\Sec(\pi^0_{B,F})$ тоже наделено структурой стереотипной алгебры, а отображение
$v:F\to\Sec(\pi^0_{B,F})$ является гомоморфизмом стереотипных алгебр.

В частном случае, когда $F$ -- $C^*$-алгебра, слои $\Big(F/\overline{I_t\cdot F}\Big)^\triangledown$ и алгебра сечений $\Sec(\pi^0_{B,F})$ также являются $C^*$-алгебрами.

Следующий вариант теоремы Даунса-Хофманна \cite{Dauns-Hofmann} отмечался в монографии М.~Дюпре и Р.~Жилетта \cite[Theorem 2.4]{Dupre-Gillette} (а для случая $B=Z(F)$ -- в работе Т.~Бекера \cite{Becker}):

\btm[Дж.~Даунс, К.~Х.~Хофманн]\label{TH:Dauns-Hofmann} Пусть $F$ -- $C^*$-алгебра и $B$ -- ее замкнутая инволютивная подалгебра, лежащая в центре $F$:
$$
B\subseteq Z(F).
$$
Тогда отображение $v:F\to\Sec(\pi^0_{B,F})$, переводящее $F$ в алгебру непрерывных сечений
расслоения значений $\pi^0_{B,F}: \Jet^0_BF\to\Spec(B)$ над алгеброй $B$, является
изоморфизмом $C^*$-алгебр:
$$
F\cong\Sec(\pi^0_{B,F}).
$$
\etm \bpr Алгебра $B$ является коммутативной $C^*$-алгеброй, поэтому ее спектр
должен быть компактом. Отсюда следует по теореме
\ref{TH:J^0_AX}, что отображение
$v:F\to\Sec(\pi^0_{B,F})$ не только непрерывно, но и имеет плотный образ в
$\Sec(\pi^0_{B,F})$.

Если теперь $\pi:F\to {\mathcal B}(X)$ -- какое-то неприводимое представление $F$, то
центр $Z(F)$ оно переводит в скалярные кратные единицы. Иными словами, $\pi$
отображает $Z(F)$, а значит и $B$, в подалгебру $\C\cdot 1_{{\mathcal B}(X)}$ алгебры
${\mathcal B}(X)$. Как следствие, должен существовать характер $t\in\Spec(B)$ такой, что
$$
\pi(a)=t(a)\cdot 1_{{\mathcal B}(X)},\qquad a\in B.
$$
Отсюда, в свою очередь, следует, что $\pi$ обнуляется на $I_t\cdot F$, потому что
$$
\pi(a\cdot
x)=\pi(a)\cdot\pi(x)=\underbrace{t(a)}_{\scriptsize\begin{matrix}\text{\rotatebox{90}{$=$}}\\
0\end{matrix}}\cdot 1_{{\mathcal B}(X)}\cdot\pi(x)=0,\qquad a\in I_t,\ x\in F.
$$
Значит, $\pi$ обнуляется и на $\overline{I_t\cdot F}$:
$$
\pi\big|_{\overline{I_t\cdot F}}=0.
$$
Из этого можно сделать вот какой вывод: если $x$ -- какой-то ненулевой вектор
из $F$, то, поскольку найдется неприводимое представление $\pi:F\to {\mathcal B}(X)$,
отличное от нуля на $x$, найдется и какая-то точка $t\in\Spec(B)$ для которой
$x\notin\overline{I_t\cdot F}$. Это означает, что $\jet^0(x)(t)\ne 0$.

Мы поэтому можем заключить, что отображение $v:F\to\Sec(\pi^0_BF)$ инъективно.
С другой стороны, как мы уже заметили, оно имеет плотный образ. Поcкольку
непрерывный гомоморфизм $C^*$-алгебр всегда имеет замкнутый образ \cite[Теоремы
3.1.5 и 3.1.6]{Murphy}, вложение $F\to \Sec(\pi)$ обязано быть изоморфизмом
$C^*$-алгебр. \epr

Теорема Даунса-Хофманна \ref{TH:Dauns-Hofmann} позволяет усилить теорему \ref{TH:Sec(J^0X)-kak-A-modul} в важном частном случае, когда $A$-модуль $X$ представляет из себя $C^*$-алгебру, в которую алгебра $A$ отображается некоторым гомоморфизмом.

\btm\label{TH:B-cong-Sec(val_AB)} Пусть $\ph :A\to F$ -- гомоморфизм инволютивных стереотипных алгебр, причем $A$ коммутативна, $F$ -- $C^*$-алгебра, и $\ph(A)$ лежит в центре $F$:
$$
\ph(A)\subseteq Z(F).
$$
Тогда отображение $v:F\to\Sec(\pi^0_{A,F})$, переводящее $F$ в алгебру непрерывных сечений
расслоения значений $\pi^0_{A,F}: \Jet^0_AF\to\Spec(A)$ над алгеброй $A$, является
изоморфизмом $C^*$-алгебр:
\beq\label{B-cong-Sec(val_AB)}
 F\cong\Sec(\pi^0_{A,F})
\eeq
\etm

\bpr Обозначим $B=\overline{\ph(A)}$. По теореме Даунса-Хофманна \ref{TH:Dauns-Hofmann},
$$
F\cong \Sec(\pi^0_{B,F}).
$$
Поэтому нам достаточно проверить равенство
$$
\Sec(\pi^0_{B,F})\cong\Sec(\pi^0_{A,F}).
$$
Из условия (i) теоремы \ref{LM:overline(ph(I_t)-cdot-B)=B} следует, что у каждой точки $t\in\Spec(B)$ слои $\Jet^0_BF(t)$ и $\Jet^0_AF(t\circ\ph)$ совпадают:
$$
(\pi^0_{A,F})^{-1}(t\circ\ph)=F\Big/\
\overline{\Ker (t\circ\ph)\kern-15pt\underset{\scriptsize\begin{matrix}\uparrow\\ \text{действие}\\ \text{$A$ на $F$}\end{matrix}}{\cdot}\kern-15pt F}=F\Big/\ \overline{\ph(\Ker (t\circ\ph))\kern-15pt\underset{\scriptsize\begin{matrix}\uparrow\\ \text{действие}\\ \text{$F$ на $F$}\end{matrix}}{\cdot}\kern-15pt F}=F\Big/\
\overline{\underbrace{\overline{\ph(\Ker (t\circ\ph))}}_{\scriptsize\begin{matrix}\phantom{\tiny{\eqref{overline(ph(Ker(t-circ-ph)))=Ker_t}}}\ \|\ \tiny{\eqref{overline(ph(Ker(t-circ-ph)))=Ker_t}}\\ \Ker t\end{matrix}}\cdot F}=
F\Big/\ \overline{\Ker t\cdot F}=(\pi^0_{B,F})^{-1}(t).
$$
А из условия (ii) теоремы \ref{LM:overline(ph(I_t)-cdot-B)=B} -- что над всеми точками $s\in\Spec(A)$, лежащими вне $\Spec(B)\circ\ph$, слои у $\pi^0_{A,F}$ нулевые:
$$
\big(\pi^0_{A,F}\big)^{-1}(s)=F/\overline{\Ker
s\kern-15pt\underset{\scriptsize\begin{matrix}\uparrow\\ \text{действие}\\ \text{$A$ на $F$}\end{matrix}}{\cdot}\kern-15pt F}=F\Big/\ \overline{\ph(\Ker
s)\kern-15pt\underset{\scriptsize\begin{matrix}\uparrow\\ \text{действие}\\ \text{$F$ на $F$}\end{matrix}}{\cdot}\kern-15pt F}=F\Big/\ \overline{\underbrace{\overline{\ph(\Ker
s)}}_{\scriptsize\begin{matrix}\phantom{\tiny{\eqref{overline(ph(Ker_s))=B}}}\ \|\ \tiny{\eqref{overline(ph(Ker_s))=B}}\\ F\end{matrix}}\cdot F}=F/\overline{F\cdot F}=F/F=\{0\}.
$$
Отсюда следует, что, во-первых, сечение $x\in\Sec(\pi^0_{A,F})$ полностью
определяется своими значениями на компакте $\Spec(B)\circ\ph$, то есть своим
ограничением на $\Spec(B)\circ\ph$. И, во-вторых, норма $x$ совпадает с нормой
его ограничения на $\Spec(B)\circ\ph$. Мы получаем, что $\Sec(\pi^0_{B,F})$ и
$\Sec(\pi^0_{A,F})$ одинаковы по набору элементов и по норме, значит, они
изоморфны как банаховы пространства. \epr

\subsection{Расслоение струй и дифференциальные операторы}

Если $I$ -- левый идеал в алгебре $A$, то в духе обозначений на с.\pageref{DEF:M-cdot-N},  для всякого $n\in\N$ мы определяем степень $I^n$ как линейное подпространство, порожденное всевозможными произведениями элементов из $I$ длины $n$:
$$
I^n=\sp\{a_1\cdot...\cdot a_n;\ a_1,...,a_n\in I\}.
$$
А замкнутая степень $\overline{I^n}$ -- как замыкание $I^n$:
$$
\overline{I^n}=\csp\{a_1\cdot...\cdot a_n;\ a_1,...,a_n\in I\}.
$$
Понятно, что это будет замкнутый левый идеал в $A$.

\paragraph{Расслоение струй.}

Для всякого $t\in \Spec(A)$ пусть, как и раньше, $I_t=\{a\in A:\ t(a)=0\}$ --
идеал в $A$, состоящий из элементов, обращающихся в нуль в точке $t$. Если
теперь $X$ -- левый модуль над $A$, то для всякого номера $n\in\Z_+$ рассмотрим
идеал $I_t^{n+1}$, порожденный им подмодуль $\overline{I_t^{n+1}\cdot X}$ в $X$
и фактор-модуль
 \beq\label{DEF:prostr-struj}
 \Jet_t^n(X)=\big( X/\overline{I_t^{n+1}\cdot X}\big)^\triangledown.
 \eeq
Он называется {\it модулем струй} порядка $n$ модуля $X$ в точке $t$.

Как обычно, всякая непрерывная полунорма $p:X\to\R_+$ определяет полунорму на
фактор-пространстве $\Jet_t^n(X)=\big( X/\overline{I_t^{n+1}\cdot
X}\big)^\triangledown$ по формуле
 \beq\label{polunorma-na-X/I_t^n-X}
p_t^n(x+\overline{I_t^{n+1}\cdot X}):=\inf_{y\in \overline{I_t^{n+1}\cdot
X}}p(x+y),\qquad x\in X.
 \eeq
Рассмотрим теперь прямую сумму множеств
$$
\Jet_A^n X=\bigsqcup_{t\in\Spec(A)}  \Jet_t^n(X)=\bigsqcup_{t\in\Spec(A)} \Big( X/\
\overline{I_t^{n+1}\cdot X}\ \Big)^\triangledown
$$
и обозначим через $\pi^n_{A,X}$ естественную проекцию $\Jet_A^n X$ на $\Spec(A)$:
$$
\pi^n_{A,X}:\Jet_A^n X\to \Spec(A),\qquad \pi^n_{A,X}\Big(
x+\overline{I_t^{n+1}\cdot X}\Big)=t,\qquad t\in \Spec(A),\ x\in X.
$$
Кроме того, для всякого вектора $x\in X$ мы рассмотрим отображение
$$
\jet^n_{A,X}(x):\Spec(A)\to \Jet_A^n X\quad\Big|\quad
\jet_{A,X}^n(x)(t)=x+\overline{I_t^{n+1}\cdot X}.
$$
Понятно, что для всякого $x\in X$
$$
\pi^n_{A,X}\circ \jet_A^n(x)=\id_{\Spec(A)}.
$$

\blm\label{LM:poluneprer-sverhu-p^n_t(j(x)(t))} Для всякого элемента $x\in X$ и
любой непрерывной полунормы $p:X\to\R_+$ отображение $t\in\Spec(A)\mapsto
p^n_t(\jet_A^n(x)(t))\in\R_+$ полунепрерывно сверху. \elm \bpr Пусть
$t\in\Spec(A)$ и $\e>0$. Условие
$$
p^n_t(\jet_A^n(x)(t))=p^n_t(x+\overline{I_t^{n+1}\cdot X})=\inf_{y\in
\overline{I_t^{n+1}\cdot X}}p(x+y)=\inf_{y\in I_t^{n+1}\cdot X}p(x+y)<\e
$$
означает, что для некоторого $y\in I_t^{n+1}\cdot X$ выполняется неравенство
$$
p(x+y)<\e.
$$
Это в свою очередь означает, что существуют $m\in\N$, векторы $y_1,...,y_m\in X$ и матрица векторов $\{a^k_i;\ 1\le i\le n+1,\ 1\le k\le m\}\subseteq I_t$ такие, что
\beq\label{p-l-x+sum_(k=1)^m-l-prod_(i=1)^(n+1) a^k_i-r-cdot-y_k-r<e}
p\l x+\sum_{k=1}^m \l \prod_{i=1}^{n+1} a^k_i\r\cdot y_k\r<\e.
\eeq
Для всякой матрицы чисел
$$
\lambda=\{\lambda^k_i;\ 1\le i\le n+1,\ 1\le k\le m\},\qquad \lambda^k_i\in\C
$$
обозначим
$$
f(\lambda)=p\l x+\sum_{k=1}^m\l\prod_{i=1}^{n+1} (a^k_i-\lambda^k_i\cdot 1_A)\r\cdot y_k\r.
$$
Функция $\lambda\mapsto f(\lambda)$ в точке $\lambda=0$ совпадает с левой частью \eqref{p-l-x+sum_(k=1)^m-l-prod_(i=1)^(n+1) a^k_i-r-cdot-y_k-r<e} и поэтому удовлетворяет неравенству
$$
f(0)<\e.
$$
С другой стороны, она непрерывна, как композиция многочлена от $m\cdot (n+1)$ комплексных переменных со значениями в локально выпуклом пространстве $X$ и непрерывной функции $p$ на $X$. Значит, должно существовать число $\delta>0$ такое, что
\beq\label{max_(i,k)|lambda^k_i|<delta=>f(lambda)<e}
\forall \lambda \qquad \max_{i,k}\abs{\lambda^k_i}<\delta\quad\Longrightarrow\quad f(\lambda)<\e.
\eeq

Рассмотрим далее множество
$$
U=\{s\in\Spec(A):\ \forall i,k\quad \abs{s(a^k_i)}<\delta\}.
$$
Оно открыто и содержит точку $t$ (потому что включения $a^k_i\in I_t$ означают систему равенств $t(a^k_i)=0$, $1\le i\le n+1,\ 1\le k\le m$).
С другой стороны, для всякой точки $s\in U$, рассмотрев матрицу
$$
\lambda^k_i=s(a^k_i),
$$
мы получим, во-первых, $s(a^k_i-\lambda^k_i\cdot 1_A)=s(a^k_i)-\lambda^k_i\cdot s(1_A)=s(a^k_i)-s(a^k_i)=0$, то есть
$$
a^k_i-\lambda^k_i\cdot 1_A\in I_s
$$
и, во-вторых, $\max_{i,k}\abs{\lambda^k_i}=\max_{i,k}\abs{s(a^k_i)}<\delta$, то есть, в силу \eqref{max_(i,k)|lambda^k_i|<delta=>f(lambda)<e},
$$
p\bigg( x+\sum_{k=1}^m \bigg(\prod_{i=1}^{n+1} \underbrace{(a^k_i-\lambda^k_i\cdot 1_A)}_{\scriptsize\begin{matrix}\text{\rotatebox{90}{$\owns$}}\\ I_s\end{matrix}}\bigg)\cdot y_k\bigg)=f(\lambda)<\e
$$
Это можно понимать так, что для некоторой точки $z\in I_s^{n+1}\cdot X$ выполняется неравенство
$$
p(x+z)<\e,
$$
которое, в свою очередь, влечет за собой неравенство
$$
p^n_s(\jet_A^n(x)(s))=p^n_s(x+\overline{I_s^{n+1}\cdot X})=\inf_{z\in
\overline{I_s\cdot X}}p(x+z)=\inf_{z\in I_s^{n+1}\cdot X}p(x+z)<\e
$$
Оно верно для всякой точки $s$ из окрестности $U$ точки $t$, и как раз это нам и нужно было доказать.
\epr

\btm\label{TH:rassloenie-struj} Для всякого стереотипного модуля $X$ над
инволютивной стереотипной алгеброй $A$ прямая сумма стереотипных фактор-модулей
$$
\Jet_A^nX=\bigsqcup_{t\in\Spec(A)}  \Jet_t^n(X)=\bigsqcup_{t\in\Spec(A)}
(X/\overline{I_t^{n+1}\cdot X})^\triangledown
$$
обладает единственной топологией, превращающей проекцию
$$
\pi^n_{A,X}:\Jet_A^n X\to\Spec(A),\qquad \pi^n_{A,X}(x+\overline{I_t^{n+1}\cdot
X})=t,\qquad t\in\Spec(A),\ x\in X
$$
в локально выпуклое расслоение с системой полунорм $\{p^n;\ p\in{\mathcal P}(X)\}$, для которого отображение
$$
x\in X\mapsto \jet^n(x)\in\Sec(\pi^n)\quad\Big|\quad
\jet^n(x)(t)=x+\overline{I_t^{n+1}\cdot X},\quad t\in\Spec(A),
$$
непрерывно переводит $X$ в стереотипный $A$-модуль $\Sec(\pi^n_{A,X})$ непрерывных сечений
$\jet^n_AX$. При этом базу топологии $\Jet_A^nX$ образуют множества
$$
W(x,U,p,\e)=\set{\xi\in \Jet_A^n X: \ \pi^n(\xi)\in U\ \& \
p^n_{\pi^n(\xi)}\Big(\xi-\jet^n(x)\big(\pi^n(\xi)\big)\Big)<\e }
$$
где $x\in X$, $p\in{\mathcal P}(X)$, $\e>0$, $U$ -- открытое множество в $M$;
 \etm

\bit{

\item[$\bullet$] Расслоение $\pi^n_{A,X}:\Jet_A^n X\to \Spec(A)$ называется {\it
расслоением струй порядка $n$} модуля $X$ над алгеброй $A$.

}\eit

\bpr Из леммы \ref{LM:poluneprer-sverhu-p^n_t(j(x)(t))} и предложения
\ref{PROP:sushestv-topologii-v-Xi} следует существование и единственность
топологии на $\Jet_A^nX$, для которой проекция $\jet^n_A X:\Jet_A^n X\to\Spec(A)$
представляет собой локально выпуклое расслоение с полунормами $p^n_t$, а
сечения вида $\jet^n(x)$, $x\in X$, будут непрерывными. Непрерывность отображения
$x\in X\mapsto \jet^n(x)\in\Sec(\jet^n_AX)$ доказывается импликацией
$$
p^n_t(\jet^n(x)(t))=\inf_{y\in I_t^{n+1}\cdot X}p(x+y)\le p(x)\quad
\Longrightarrow \quad \sup_{t\in T}p^n_t(\jet^n(x)(t))\le p(x)
$$
для всякого компакта $T\subseteq M$. Свойство (i) также следует из предложения
\ref{PROP:sushestv-topologii-v-Xi}.

Докажем (ii). Заметим сначала, что множество $\jet^n(X)=\{\jet^n(x);\ x\in X\}$
плотно в порожденном им ${\mathcal C}(M)$-модуле
$$
{\mathcal C}(M)\cdot \jet^n(X)=\set{\sum_{i=1}^k b_i\cdot \jet^n(x_i);\ b_i\in {\mathcal C}(M),\ x_i\in X}.
$$
Для этого зафиксируем $b\in {\mathcal C}(M)$ и $x\in X$. Поскольку $A$ -- инволютивная
алгебра, отображение $\sigma:A\to {\mathcal C}(M)$ должно иметь плотный образ в
${\mathcal C}(M)$. Поэтому найдется направленность $a_i\in A$ такая, что
$$
\sigma(a_i)\overset{{\mathcal C}(M)}{\underset{i\to\infty}{\longrightarrow}}b.
$$
Как следствие,
$$
\jet^n(a_i\cdot x)=\eqref{j^0(a-cdot-x)=sigma(a)-cdot-j^0(x)}=\sigma(a_i)\cdot
\jet^n(x)\overset{{\mathcal C}(M)}{\underset{i\to\infty}{\longrightarrow}}b\cdot \jet^n(x).
$$
Далее заметим, что поскольку образ $X$ при каждой проекции $\pi^n_t$ плотен в
фактор-пространстве $(X/\overline{I_t\cdot X})^\triangledown$, множество
$\jet^n(X)(\Spec(A))$ будет плотно в каждом слое расслоения $\Jet^n_AX$. Поэтому его
надмножество $({\mathcal C}(M)\cdot \jet^n(X))(\Spec(A))$ также  будет плотно в каждом слое
расслоения $\Jet^n_AX$. Теперь мы можем применить свойство $2^\circ$ на
с.\pageref{LM:poltnost-v-sloyah=>plotnost-v-secheniyah}: если $M$ паракомпактно
и локально компактно, то ${\mathcal C}(M)$-модуль ${\mathcal C}(M)\cdot \jet^n(X)$ плотен в
$\Sec(\pi^n)$ (а, как мы уже поняли $\jet^n(X)$ плотен в ${\mathcal C}(M)\cdot \jet^n(X)$). \epr

\paragraph{Дифференциальные операторы и их связь с морфизмами расслоений струй.}

Пусть $X$ и $Y$ -- два левых стереотипных модуля над стереотипной алгеброй $A$.
Для всякого линейного (над $\C$) отображения $D:X\to Y$ и любого элемента $a\in
A$ отображение $[D,a]:X\to Y$, действующее по формуле 
\beq\label{DEF:[P,a]}
[D,a](x)=D(a\cdot x)-a\cdot D(x),\qquad x\in X, 
\eeq 
называется {\it коммутатором} отображения $D$ с элементом $a$. Если даны два элемента $a,b\in
A$, то коммутатор отображения $D$ с парой элементов $(a,b)$ определяется как
коммутатор $[[D,a],b]:X\to Y$ отображения $[D,a]$ с элементом $b$. Точно так же
по индукции определяется коммутатор с произвольной конечной последовательностью
элементов $(a_0,...,a_n)$:
$$
[...[D,a_0],... a_n]
$$

Линейное (над $\C$) непрерывное отображение $D:X\to Y$ называется {\it
дифференциальным оператором} из $X$ в $Y$, если существует число $n\in\Z_+$
такое, что для любых $a_0,...,a_n\in A$ выполняется равенство
\beq\label{DEF:diff-oper} [...[D,a_0],...a_n]=0. 
\eeq 
Наименьшее из чисел
$n\in\Z_+$, для которых справедливо \eqref{DEF:diff-oper}, называется {\it
порядком} дифференциального оператора $D$ и обозначается $\ord D$.

Множество всех дифференциальных операторов из $X$ в $Y$ порядка не больше $n$
мы будем обозначать символом $\Diff^n(X,Y)$. Если дополнительно ввести
обозначение
$$
\Diff^{-1}(X,Y)=\{D\in Y\oslash X: \ D=0\},
$$
то эту последовательность пространств можно определить следующими индуктивным правилом:
\begin{align}
& \Diff^{n+1}(X,Y)=\{D\in Y\oslash X:\ \forall a\in A\quad [D,a]\in \Diff^n(X,Y)\} \label{D^n(X,Y)}
\end{align}
Очевидно, что
\beq\label{[D^n,A]-subseteq-D^(n-1)}
[\Diff^n,A]\subseteq \Diff^{n-1}
\eeq

\blm
Если $D:X\to Y$ -- дифференциальный оператор порядка $n\in\Z_+$, то подмодуль $\overline{I^{n+1}\cdot X}$ он переводит в подмодуль $\overline{I\cdot Y}$:
\beq\label{D((I^(n+1)-X))-subseteq-(I-Y)}
D\in \Diff^n(X,Y)\qquad\Longrightarrow\quad D\l \overline{I^{n+1}\cdot X}\r\subseteq \overline{I\cdot Y}.
\eeq
\elm
\bpr
Поскольку оператор $D$ непрерывен, достаточно доказать включение
\beq\label{D(I^(n+1)-X)-subseteq-I-Y}
D\in \Diff^n(X,Y)\qquad\Longrightarrow\quad D(I^{n+1}\cdot X)\subseteq I\cdot Y.
\eeq
Это делается индукцией по $n\in\Z_+$. При $n=0$ утверждение принимает вид
$$
D\in \Diff^0(X,Y)\qquad\Longrightarrow\quad D(I\cdot X)\subseteq I\cdot Y,
$$
и это очевидно, потому что дифференциальный оператор порядка $n=0$ представляет собой просто однородное отображение
$$
D(a\cdot x)=a\cdot D(x),\qquad a\in A,\ x\in X.
$$

Предположим, что мы доказали вложение \eqref{D(I^(n+1)-X)-subseteq-I-Y} для $n=k$:
$$
D\in \Diff^k(X,Y)\qquad\Longrightarrow\quad D(I^{k+1}\cdot X)\subseteq I\cdot Y.
$$
Тогда для $n=k+1$ получаем: если $D\in \Diff^{k+1}(X,Y)$, то для всякого $a\in A$ должно выполняться включение $[D,a]\in \Diff^k(X,Y)$, поэтому, в силу предположения индукции, $[D,a](I^{k+1}\cdot X)\subseteq I\cdot Y$. Это означает, что для всякого вектора $x\in X$ и любой последовательности $a_0,...,a_k\in I$ выполняется включение
$$
\underbrace{[D,a](a_0\cdot ...\cdot a_k\cdot x)}_{\scriptsize\begin{matrix}\text{\rotatebox{90}{$=$}}
\\ D(a\cdot a_0\cdot ...\cdot a_k\cdot x)-a\cdot D(a_0\cdot ...\cdot a_k\cdot x)\end{matrix}}\kern-30pt\in I\cdot Y
$$
из которого следует
$$
D(a\cdot a_0\cdot ...\cdot a_k\cdot x)\in I\cdot Y+\underbrace{a\cdot D(a_0\cdot ...\cdot a_k\cdot x)}_{\scriptsize\begin{matrix}\text{\rotatebox{90}{$\owns$}}
\\ I\cdot Y\end{matrix}}\subseteq I\cdot Y
$$
Поскольку это верно для любых $a,a_0,...,a_k\in I$, мы получаем нужное нам включение $D(I^{k+2}\cdot X)\subseteq I\cdot Y$.
\epr

\btm\label{TH:diff-oper->rassl-struuj} Всякий дифференциальный оператор $D:X\to
Y$ порядка $n$ определяет некий морфизм расслоений струй $\jet_n[D]:\Jet_A^n X\to
\Jet_A^0(Y)$,
$$
 \xymatrix @R=2pc @C=1.2pc
 {
  \Jet^n_AX\ar[rr]^{\jet^n[D]}\ar[dr]_{\pi^n_{A,X}} & &  \Jet^0_AY\ar[dl]^{\pi^0_{A,Y}}\\
  & \Spec(A) &
 }
$$
удовлетворяющий тождеству
 \beq\label{j^0(Dx)=j_n[D]-circ-j^n(x)}
\jet^0(Dx)=\jet_n[D]\circ \jet^n(x),\qquad x\in X.
 \eeq
$$
 \xymatrix @R=2pc @C=1.2pc
 {
  \Jet^n_AX\ar[rr]^{\jet^n[D]} & &  \Jet^0_AY\\
  & \Spec(A)\ar[ul]^{\jet^n(x)}\ar[ur]_{\jet^0(Dx)} &
 }
$$
\etm \bpr Из вложения
\eqref{D((I^(n+1)-X))-subseteq-(I-Y)}, примененного к идеалу $I_t$, где
$t\in\Spec(A)$ \beq\label{D((I_t^(n+1)-X))-subseteq-(I_t-Y)} D\l
\overline{I_t^{n+1}\cdot X}\r\subseteq \overline{I_t\cdot Y} \eeq следует
существование естественного отображения фактор-пространств:
$$
X/\ \overline{I_t^{n+1}\cdot X}\owns x+\overline{I_t^{n+1}\cdot X}\mapsto Dx+\overline{I_t\cdot Y}\in Y/\ \overline{I_t\cdot Y}
$$
Оно непрерывно, поскольку исходное отображение $D$ непрерывно, и значит, существует естественное (также непрерывное) отображение стереотипных фактор-пространств (то есть псевдопополнений обычных фактор-пространств):
$$
\jet_n[D]_t:\Big( X/\ \overline{I_t^{n+1}\cdot X}\Big)^\triangledown\to \Big( Y/\ \overline{I_t\cdot Y}\Big)^\triangledown
$$
Это верно для всякого $t\in\Spec(A)$, поэтому возникает отображение прямых сумм:
$$
\jet_n[D]:\bigsqcup_{t\in\Spec(A)}\Big( X/\ \overline{I_t^{n+1}\cdot X}\Big)^\triangledown\to \bigsqcup_{t\in\Spec(A)}\Big( Y/\ \overline{I_t\cdot Y}\Big)^\triangledown
$$

Тождество \eqref{j^0(Dx)=j_n[D]-circ-j^n(x)} проверяется прямым вычислением: для любых $t\in\Spec(A)$ и $x\in X$
$$
\big(\jet_n[D]\circ \jet^n(x)\big)(t)=\jet_n[D]\big(\jet^n(x)(t)\big)=\jet_n[D]\l x+\overline{I_t^{n+1}\cdot X}\r=Dx+\overline{I_t\cdot Y}=\jet^0(Dx)(t).
$$

Остается проверить непрерывность отображения $\jet_n[D]$. Зафиксируем точку $\zeta=\jet^n(x)(t)\in \Jet_A^n X$, $x\in X$, $t\in T$ и рассмотрим какую-нибудь базисную окрестность ее образа $\jet_n[D](\zeta)=\jet_n[D]\big(\jet^n(x)(t)\big)=\jet^0(Dx)(t)$ при отображении $\jet_n[D]$:
$$
W(y,V,q,\e)=\set{\upsilon\in \Jet_A^n(Y): \ \pi_Y(\upsilon)\in U\ \& \ q_{\pi_Y(\upsilon)}^n\Big(\upsilon-\jet^0(y)\big(\pi_Y(\upsilon)\big)\Big)<\e }.
$$
где $q:Y\to\R_+$ -- произвольная непрерывная полунорма на $Y$, $V$ -- окрестность точки $t$ в $\Spec(A)$, $y\in Y$, $\e>0$. Поскольку $\jet^0(Dx)(t)\in W(y,U,q,\e)$, должно выполняться условие
$$
q_t^0\Big(\jet^0(Dx)(t)-\jet^0(y)(t)\Big)<\e.
$$
Из него следует, что существуют число $\delta>0$ и окрестность $U$ точки $t$, лежащая в $V$ такие, что
$$
\pi_Y(\upsilon)\in U\& \ q_{\pi_Y(\upsilon)}^n\Big(\upsilon-\jet^0(Dx)\big(\pi_Y(\upsilon)\big)\Big)<\delta \qquad\Longrightarrow\qquad q_{\pi_Y(\upsilon)}^0\Big(\jet^0(Dx)\big(\pi_Y(\upsilon)\big)-\jet^0(y)\big(\pi_Y(\upsilon)\big)\Big)<\e
$$
Отсюда следует, что окрестность
$$
W(Dx,U,q,\delta)=\set{\upsilon\in \Jet_A^n(Y): \ \pi_Y(\upsilon)\in U\ \& \ q_{\pi_Y(\upsilon)}^0\Big(\upsilon-\jet^0(Dx)\big(\pi_Y(\upsilon)\big)\Big)<\delta }.
$$
содержится в окрестности $W(y,V,q,\e)$ точки $\jet^0(Dx)(t)$:
$$
W(Dx,U,q,\delta)\subseteq W(y,V,q,\e).
$$

Далее вспомним, что $D:X\to Y$ -- непрерывное линейное отображение. Как следствие, существует непрерывная полунорма $p:X\to\R_+$, удовлетворяющая условию
\beq\label{q(Dx)-le-p(x)}
q(Dx)\le p(x),\qquad x\in X.
\eeq
Заметим, что из этого неравенства следует неравенство для полунорм на расслоениях струй:
\beq\label{q_t^0(j_n[D](xi))-le-p_t^n(xi)}
q_t^0\Big(\jet_n[D](\xi)\Big)\le p_t^n(\xi),\qquad t\in\Spec(A),\ \xi\in \Jet_A^n X.
\eeq
Это достаточно доказать для точек $\xi=\jet^n(x)(t)$, $x\in X$, поскольку они плотны в каждом слое:
\begin{multline*}
q_t^0\Big(\jet_n[D]\big(\jet^n(x)(t)\big)\Big)=\eqref{j^0(Dx)=j_n[D]-circ-j^n(x)}=q_t^0\Big(\jet^0(Dx)(t)\Big)=q_t^0\Big(Dx+\overline{I_t\cdot Y}\Big)=
\inf_{v\in I_t\cdot Y}q(Dx+v)
\kern-30pt
\overset{\scriptsize\begin{matrix}\eqref{D(I^(n+1)-X)-subseteq-I-Y}
\\
\Downarrow\\
I_t\cdot Y\supseteq D(I_t^{n+1}\cdot Y)
\\
\Downarrow
\end{matrix}}{\le}
\kern-20pt
 \\ \le \inf_{u\in I_t^{n+1}\cdot X}q(Dx+Du)=
\inf_{u\in I_t^{n+1}\cdot X}q\Big(D(x+u)\Big)\le\eqref{q(Dx)-le-p(x)}\le \inf_{u\in I_t^{n+1}\cdot X}p(x+u)=p_t^n\Big(x+\overline{I_t^{n+1}\cdot X}\Big)=p_t^n\Big(\jet^n(x)(t)\Big)
\end{multline*}

Теперь можно показать, что окрестность
$$
W(x,U,p,\delta)=\set{\xi\in \Jet_A^n X: \ \pi_X(\xi)\in U\ \& \ p_{\pi_X(\xi)}^n\Big(\xi-\jet^n(x)\big(\pi_X(\xi)\big)\Big)<\delta }.
$$
точки $\zeta=\jet^n(x)(t)$ при отображении $\jet_n[D]$ переходит в окрестность $W(Dx,U,q,\delta)$ точки $\jet_n[D](\zeta)=\jet^0(Dx)(t)$:
\beq\label{j_n[D](W(x,U,p,delta))-subseteq-W(Dx,U,q,delta)}
\jet_n[D]\Big( W(x,U,p,\delta) \Big)\subseteq W(Dx,U,q,\delta)
\eeq
Действительно, для всякой точки $\xi\in W(x,U,p,\delta)$ условие $\pi_X(\xi)\in U$ оказывается полезным в конце цепочки
$$
\pi_Y(\jet_n[D](\xi))=\pi_X(\xi)\in U
$$
а условие $p_{\pi_X(\xi)}^n\Big(\xi-\jet^n(x)\big(\pi_X(\xi)\big)\Big)<\delta$ в конце цепочки
\begin{multline*}
q_{\jet_n[D](\xi)}^0\Big(\jet_n[D](\xi)-\jet^0(Dx)\big(\pi_Y(\jet_n[D](\xi))\big)\Big)=
q_{\pi_X(\xi)}^0\Big(\jet_n[D](\xi)-\jet^0(Dx)\big(\pi_X(\xi)\big)\Big)=\eqref{j^0(Dx)=j_n[D]-circ-j^n(x)}=\\=
q_{\pi_X(\xi)}^0\bigg(\jet_n[D](\xi)-\jet_n[D]\Big(\jet^n(x)\big(\pi_X(\xi)\big)\Big)\bigg)=
q_{\pi_X(\xi)}^0\bigg(\jet_n[D]\Big(\xi-\jet^n(x)\big(\pi_X(\xi)\big)\Big)\bigg)\le \eqref{q_t^0(j_n[D](xi))-le-p_t^n(xi)}\le \\
\le
p_{\pi_X(\xi)}^n\Big(\xi-\jet^n(x)\big(\pi_X(\xi)\big)\Big)<\delta.
\end{multline*}
Вместе то и другое означает, что $\jet_n[D](\xi)\in W(Dx,U,q,\delta)$, а это и требовалось.
\epr

\paragraph{Дифференциальные операторы на алгебрах.}

Всякий гомоморфизм алгебр $\ph:A\to B$ задает на $B$ структуру левого $A$-модуля по формуле:
$$
a\cdot y=\ph(a)\cdot y,\qquad a\in A,\ y\in Y=B.
$$
Если теперь $D:A\to B$ -- произвольное линейное (над $\C$) отображение, то формула \eqref{DEF:[P,a]} его коммутатора с элементом $a\in A$ в данной ситуации имеет вид
\beq\label{DEF:[D,a]}
[D,a](x)=D(a\cdot x)-\ph(a)\cdot D(x),\qquad a,x\in A,
\eeq
Этот оператор мы будем называть {\it коммутатором оператора $D:A\to B$ и элемента $a\in A$ относительно гомоморфизма} $\ph:A\to B$. Кроме того, для любого элемента $b\in B$ мы будем рассматривать линейное (над $\C$) отображение $b\cdot D:A\to B$, определенное формулой
$$
(b\cdot D)(x)=b\cdot D(x),\qquad x\in A.
$$

\bprop Справедливы тождества:
\begin{align}
& [\ph,a]=0, && a\in A \label{[ph,a]=0} \\
& [b\cdot D, a]=b\cdot [D,a]+[b,\ph(a)]\cdot D, && a\in A,\ b\in B,\ D\in B\oslash A, \label{[b-cdot-P, a]} \\
& [b\cdot\ph,a]=[b,\ph(a)]\cdot\ph, && a\in A, \ b\in B. \label{[b-cdot-ph,a_0]}
\end{align}
\eprop
\bpr
Первое и второе тождества доказываются прямым вычислением: при $x\in A$ мы получаем
$$
[\ph,a](x)=\ph(a\cdot x)-\ph(a)\cdot\ph(x)=0,
$$
и
\begin{multline*}
[b\cdot D, a](x)=(b\cdot D)(a\cdot x)-\ph(x)\cdot(b\cdot D)(x)=b\cdot D(a\cdot x)-\ph(x)\cdot b\cdot D(x)=\\=
b\cdot D(a\cdot x)-b\cdot \ph(x)\cdot D(x)+b\cdot \ph(x)\cdot D(x)-\ph(x)\cdot b\cdot D(x)=
b\cdot \Big( D(a\cdot x)-\ph(x)\cdot D(x)\Big)+\Big(b\cdot \ph(x)-\ph(x)\cdot b\Big)\cdot D(x)=\\=
b\cdot [D,a](x)+[b,\ph(x)]\cdot D(x)=\Big(b\cdot [D,a]+[b,\ph(a)]\cdot D\Big)(x).
\end{multline*}
А третье после этого становится следствием первого и второго:
$$
[b\cdot\ph,a]=\eqref{[b-cdot-P, a]}=b\cdot\kern-8pt\underbrace{[\ph,a_0]}_{\tiny\begin{matrix} \phantom{\eqref{[ph,a]=0}} \ \text{\rotatebox{90}{$=$}} \ \eqref{[ph,a]=0}
\\  0  \end{matrix}}\kern-8pt +[b,\ph(a)]\cdot \ph
$$
\epr

Заданный гомоморфизм алгебр $\ph:A\to B$ порождает последовательность пространств дифференциальных операторов из $A$ в $B$, которую мы будем обозначать $\Diff^n(\ph)$ или просто $\Diff^n$. Она индуктивно определяется правилами
\begin{align}
& \Diff^{-1}(\ph)=\{D\in B\oslash A: \ D=0\}, \label{D^(-1)(ph)} \\
& \Diff^{n+1}(\ph)=\{D\in B\oslash A:\ \forall a\in A\quad [D,a]\in \Diff^n(\ph)\} \label{D^n(ph)}
\end{align}
Понятно, что последовательность пространств $\Diff^n(\ph)$ расширяется:
$$
0=\Diff^{-1}(\ph)\subseteq \Diff^0(\ph)\subseteq \Diff^1(\ph)\subseteq ...\subseteq \Diff^n(\ph)\subseteq \Diff^{n+1}(\ph)\subseteq...
$$

\bit{
\item[$\bullet$] Линейное (над $\C$) непрерывное отображение $D:A\to B$ называется {\it
дифференциальным оператором} из $A$ в $B$ относительно гомоморфизма $\ph:A\to B$, если существует число $n\in\Z_+$
такое, что 
\beq\label{DEF:diff-oper-A->B} 
D\in \Diff^n(\ph). 
\eeq 
Это эквивалентно тождеству
\beq\label{DEF:diff-oper-A->B-1} [...[D,a_0],...a_n]=0. 
\eeq 
для любых $a_0,...,a_n\in A$.

\item[$\bullet$] Наименьшее из чисел $n\in\Z_+$, для которых справедливо \eqref{DEF:diff-oper-A->B}, называется {\it
порядком} дифференциального оператора $D$ и обозначается $\ord D$.

}\eit

\noindent\rule{160mm}{0.1pt}\begin{multicols}{2}

\bex\label{EX:ph-diff-oper-por-0}
Сам гомоморфизм $\ph:A\to B$ является дифференциальным оператором порядка 0 относительно самого себя, потому что 
$$
[\ph,a_0](x)=\eqref{DEF:[D,a]}=\ph(a\cdot x)-\ph(a)\cdot\ph(x)=0
$$
\eex

\brem[не всякое дифференцирование является дифференциальным оператором!]\label{EX:differentsirovanie}
{\it Дифференцирование} алгебры $A$, то есть оператор $D:A\to A$, удовлетворяющий тождеству
\beq\label{DEF:differentsirovanie}
D(x\cdot y)=D(x)\cdot y + x\cdot D(y)
\eeq
является дифференциальным оператором порядка 1 относительно тождественного оператора $\ph=\id_A:A\to A$ тогда и только тогда, когда его множество значений $D(A)$ лежит в центре алгебры $A$:
$$
D(A)\subseteq Z(A).
$$
\erem
\bpr
Сначала заметим, что 
$$
[D,a]=a\cdot \id_A.
$$
Отсюда следует, что 
$$
[[D,a],b]=(D(a)\cdot b-b\cdot D(a))\cdot \id_A
$$
Чтобы $D$ был дифференциальным оператором порядка 1, нужно, чтобы оператор $[[D,a],b]$ был равен нулю при всех $a,b\in A$. Это будет только если 
$$
D(a)\cdot b=b\cdot D(a), \quad a,b\in A.
$$ 
\epr

\end{multicols}\noindent\rule[10pt]{160mm}{0.1pt}

Кроме того нам будет важна последовательность подпространств в $B$, обозначаемых $Z^n(\ph)$, или просто $Z^n$, и определяемых индуктивными правилами
\begin{align}
& Z^0(\ph)=0 \label{Z^0}\\
& Z^{n+1}(\ph)=\{b\in B:\ \forall a\in A\quad [b,\ph(a)]\in Z^n(\ph)\} \label{Z^n}
\end{align}
Очевидно, что пространства $Z^n(\ph)$ также образуют расширяющуюся последовательность:
\beq\label{Z^n(ph)}
0=Z^0(\ph)\subseteq Z^1(\ph)\subseteq ...\subseteq Z^n(\ph)\subseteq Z^{n+1}(\ph)\subseteq...
\eeq

\bprop Справедливы включения:
\begin{align}
& [Z^q\cdot \Diff^p,A]\subseteq Z^q\cdot \Diff^{p-1}+Z^{q-1}\cdot \Diff^p, \label{[Z^q-cdot-D^p,A]}\\
& Z^0\cdot \Diff^p\subseteq \Diff^{-1}=0,\qquad p\ge 0 \label{Z.D^p=0} \\
& Z^q\cdot \Diff^0\subseteq \Diff^{q-1},\qquad q\ge 0 \label{Z^q.D^0->D^(q-1)} \\
& Z^q\cdot \Diff^p\subseteq \Diff^{q+p-1},\qquad q\ge 0 \label{Z^q.D^p->D^(q+p-1)}
\end{align}
\eprop
\bpr

1. Для доказательства \eqref{[Z^q-cdot-D^p,A]} берем произвольные $b\in Z^q$, $D\in \Diff^p$ и $a\in A$, тогда
$$
[b\cdot D, a]=\overset{\tiny\begin{matrix}Z^q \\ \text{\rotatebox{90}{$\in$}} \end{matrix}}{b}\cdot \underbrace{[\overset{\tiny\begin{matrix}\Diff^p \\ \text{\rotatebox{90}{$\in$}} \end{matrix}}{D},a]}_{\tiny\begin{matrix}\text{\rotatebox{90}{$\owns$}} \\ \Diff^{p-1}\end{matrix}}+\underbrace{[\overset{\tiny\begin{matrix}Z^q \\ \text{\rotatebox{90}{$\in$}} \end{matrix}}{b},\ph(a)]}_{\tiny\begin{matrix}\text{\rotatebox{90}{$\owns$}} \\ Z^{q-1}\end{matrix}}\cdot \overset{\tiny\begin{matrix}\Diff^p \\ \text{\rotatebox{90}{$\in$}} \end{matrix}}{D}\in Z^q\cdot \Diff^{p-1}+Z^{q-1}\cdot \Diff^p
$$

2. Утверждение \eqref{Z.D^p=0} очевидно, потому что при умножении на нуль всегда получается нуль.

3. Формула \eqref{Z^q.D^0->D^(q-1)} доказывается индукцией. При $q=0$ она верна, потому что превращается в частный случай формулы \eqref{Z.D^p=0} с $p=0$:
$$
Z^0\cdot \Diff^0\subseteq \Diff^{-1}=0.
$$
Предположим, мы ее доказали для какого-то $q=n$:
\beq\label{Z^n-cdot-D^0-subseteq-D^(n-1)}
Z^n\cdot \Diff^0\subseteq \Diff^{n-1}
\eeq
Тогда для $q=n+1$ мы получим: если $b\in Z^q=Z^{n+1}$ и $D\in \Diff^0$, то
$$
\Big(\forall a\in A\quad
[b\cdot D, a]=b\cdot \underbrace{[\overset{\tiny\begin{matrix}\Diff^0 \\ \text{\rotatebox{90}{$\in$}} \end{matrix}}{D},a]}_{\tiny\begin{matrix}\text{\rotatebox{90}{$=$}} \\ 0\end{matrix}}+\underbrace{[b,\ph(a)]}_{\tiny\begin{matrix}\text{\rotatebox{90}{$\owns$}} \\ Z^n\end{matrix}}\cdot \overset{\tiny\begin{matrix}\Diff^0 \\ \text{\rotatebox{90}{$\in$}} \end{matrix}}{D}\overset{\eqref{Z^n-cdot-D^0-subseteq-D^(n-1)}}{\in} \Diff^{n-1}\Big)
\quad\Longrightarrow\quad b\cdot D\in \Diff^n=\Diff^{q-1}
$$

4. Формула \eqref{Z^q.D^p->D^(q+p-1)} доказывается индукцией по $p$. При $p=0$ она верна, потому что превращается в уже доказанную формулу  \eqref{Z^q.D^0->D^(q-1)}. Предположим, что мы ее доказали для какого-то $p=n$:
\beq\label{Z^q.D^n->D^(q+n-1)}
Z^q\cdot \Diff^n\subseteq \Diff^{q+n-1},\qquad q\ge 0
\eeq
Тогда для $p=n+1$ и $a_1,...,a_q\in A$ мы получим
$$
[Z^q\cdot \Diff^{n+1},a_1]\overset{\eqref{[Z^q-cdot-D^p,A]}}{\subseteq}
\underbrace{Z^q\cdot \Diff^n}_{\tiny\begin{matrix}
\phantom{\eqref{Z^q.D^n->D^(q+n-1)}}
\ \text{\rotatebox{90}{$\supseteq$}}\ \eqref{Z^q.D^n->D^(q+n-1)} \\ \Diff^{q+n-1}\end{matrix}}
+Z^{q-1}\cdot \Diff^{n+1}\subseteq \Diff^{q+n-1}+Z^{q-1}\cdot \Diff^{n+1}
$$
$$
\Downarrow
$$
\begin{multline*}
[[Z^q\cdot \Diff^{n+1},a_1],a_2]\subseteq \underbrace{[\Diff^{q+n-1},a_2]}_{\tiny\begin{matrix}
\phantom{\eqref{[D^n,A]-subseteq-D^(n-1)}}
\ \text{\rotatebox{90}{$\supseteq$}}\ \eqref{[D^n,A]-subseteq-D^(n-1)} \\ \Diff^{q+n-2}\end{matrix}}+[Z^{q-1}\cdot \Diff^{n+1},a_2]\subseteq \\ \subseteq
\Diff^{q+n-2}+\underbrace{Z^{q-1}\cdot \Diff^n}_{\tiny\begin{matrix}
\phantom{\eqref{Z^q.D^n->D^(q+n-1)}}
\ \text{\rotatebox{90}{$\supseteq$}}\ \eqref{Z^q.D^n->D^(q+n-1)} \\ \Diff^{q+n-2}\end{matrix}}
+Z^{q-1}\cdot \Diff^{n+1}\subseteq \Diff^{q+n-2}+Z^{q-2}\cdot \Diff^{n+1}
\end{multline*}
$$
\Downarrow
$$
$$
\cdots
$$
$$
\Downarrow
$$
$$
[...[[Z^q\cdot \Diff^{n+1},a_1],a_2],... a_q]\subseteq \Diff^{q+n-q}+Z^{q-q}\cdot \Diff^{n+1}=
\Diff^n+Z^0\cdot \Diff^{n+1}=\Diff^n
$$
Это верно для любых $a_1,...,a_q\in A$, поэтому
$$
Z^q\cdot \Diff^{n+1}\subseteq \Diff^{q+n}=\Diff^{q+p-1}
$$
\epr

\bprop
Для всякого элемента $b\in B$
\beq\label{b.ph-in-D^n(ph)<=>b-in-Z^(n+1)(ph)}
b\cdot\ph\in\Diff^n(\ph)\quad\Longleftrightarrow\quad b\in Z^{n+1}(\ph)
\eeq
\eprop
\bpr
Здесь используется формула \eqref{[b-cdot-ph,a_0]}, которую нужно обобщить до тождества
$$
[...[b\cdot\ph,a_0],...a_n]=[...[b,\ph(a_0)],...\ph(a_n)]\cdot\ph, \qquad a_0,...,a_n\in A.
$$
Из него видно, что если $b\in Z^{n+1}(\ph)$, то коэффициент при $\ph$ справа будет нулевым, и поэтому вся правая часть будет равна нулю. Поскольку это верно для любых $a_0,...,a_n\in A$, мы получаем, что $b\cdot\ph\in\Diff^n(\ph)$. И наоборот, если $b\cdot\ph\in\Diff^n(\ph)$, то левая часть будет нулевой, поэтому, в частности, при подстановке единицы, мы получаем
$$
0=[...[b\cdot\ph,a_0],...a_n](1_A)=[...[b,\ph(a_0)],...\ph(a_n)]\cdot\ph(1_A)=[...[b,\ph(a_0)],...\ph(a_n)]\cdot 1_B=[...[b,\ph(a_0)],...\ph(a_n)].
$$
Опять это верно для любых $a_0,...,a_n\in A$, поэтому $b\in Z^{n+1}(\ph)$.
\epr

\paragraph{Дифференциальные операторы со значениями в $C^*$-алгебре.}

Пусть $A$ -- инволютивная замкнутая подалгебра в $C^*$-алгебре $B$. Рассмотрим последовательность $Z^n(A)$ подпространств в $B$, определенных индуктивными правилами:
\begin{align}
& Z^0(A)=0 \label{Z^0(A)}\\
& Z^{n+1}(A)=\{b\in B:\ \forall a\in A\quad [b,a]\in Z^n(A)\} \label{Z^n(A)}
\end{align}

\btm\label{TH:stabilizatsiya-Z^n} Для всякой $C^*$-алгебры $B$ и любой ее замкнутой инволютивной подалгебры $A$ последовательность подпространств $Z^n(A)$ стабилизируется начиная с номера $n=1$:
\beq\label{Z^n(A)=Z^(n+1)(A)}
Z^n(A)=Z^{n+1}(A),\qquad n\ge 1.
\eeq
\etm

Идея доказательства принадлежит Ю.~Н.~Кузнецовой и опирается на две леммы.

Условимся {\it проектором} в локально выпуклом пространстве $X$ называть всякий линейный непрерывный оператор $P:X\to X$ со свойством идемпотентности:
$$
P^2=P.
$$
Всякий такой оператор действует как тождественный на своем образе \cite[с.37]{Bourbaki-tvs}:
\beq\label{y-in-Im(P)=>P(y)=y}
\forall y\in \Im P\qquad P(y)=y.
\eeq

\blm\label{LM:TP=PT} Для всякого локально выпуклого пространства $X$, любого проектора $P$ в $X$ и любого линейного непрерывного оператора $T$ в $X$ условие коммутирования
$$
[T,P]=0
$$
эквивалентно инвариантности образа $\Im P$ и ядра $\Ker P$ оператора $P$ относительно оператора $T$:
\beq\label{TP=PT-2}
T(\Im P)\subseteq \Im P\quad\&\quad T(\Ker P)\subseteq \Ker P.
\eeq
\elm
\bpr
1. Пусть $TP=PT$. Тогда если $x\in \Im P$, то есть $x=Px$, то $Tx=TPx=PTx\in \Im P$. Если же $x\in \Ker P$, то есть $Px=0$, то $PTx=TPx=0$ и поэтому $Tx\in \Ker P$.

2. Наоборот, пусть выполнено \eqref{TP=PT-2}. Тогда для всякого вектора $x\in X$ обозначив
$$
x_\parallel=Px,\qquad x_\perp=x-x_\parallel,
$$
мы получим:
$$
\begin{cases}x_\parallel\in\Im P\quad\Longrightarrow\quad Tx_\parallel\in\Im P\\
x_\perp\in\Ker P\quad\Longrightarrow\quad Tx_\perp\in\Ker P\end{cases},
$$
поэтому
$$
PTx=PT(x_\parallel+x_\perp)=P\Big(\underbrace{Tx_\parallel}_{\tiny\begin{matrix}\text{\rotatebox{90}{$\owns$}}\\
\Im P\end{matrix}}+\underbrace{Tx_\perp}_{\tiny\begin{matrix}\text{\rotatebox{90}{$\owns$}}\\
\Ker P\end{matrix}}\Big)=P\underbrace{Tx_\parallel}_{\tiny\begin{matrix}\text{\rotatebox{90}{$\owns$}}\\
\Im P\end{matrix}}=\eqref{y-in-Im(P)=>P(y)=y}=Tx_\parallel=TPx.
$$
\epr

\blm\label{LM:TP=PT-(2)} Для всякого локально выпуклого пространства $X$, любого проектора $P$ в $X$ и любого линейного непрерывного оператора $T$ в $X$ условие
$$
[[T,P],P]=0
$$
влечет за собой условие
$$
[T,P]=0.
$$
\elm
\bpr Из леммы \ref{LM:TP=PT} следует, что пространства $\Im P$ и $\Ker P$ инвариантны для оператора $[T,P]$:
$$
[T,P](\Im P)\subseteq \Im P\quad\&\quad [T,P](\Ker P)\subseteq \Ker P.
$$
Первое из этих условий означает, что для всякого $x\in\Im P$ справедлива цепочка
$$
TPx-\underbrace{PTx}_{\tiny\begin{matrix}\text{\rotatebox{90}{$\owns$}}\\
\Im P\end{matrix}}\in \Im P\quad\Longrightarrow\quad T\underbrace{Px}_{\tiny\begin{matrix}\text{\rotatebox{90}{$=$}}\\
x\end{matrix}}\in \Im P\quad\Longrightarrow\quad Tx\in \Im P.
$$
А второе -- что для всякого $x\in\Ker P$ справедлива цепочка
$$
T\underbrace{Px}_{\tiny\begin{matrix}\text{\rotatebox{90}{$=$}}\\ 0\end{matrix}}-PTx\in \Ker P\quad\Longrightarrow\quad PTx\in \Ker P\quad\Longrightarrow\quad \underbrace{P^2Tx}_{\tiny\begin{matrix}\text{\rotatebox{90}{$=$}}\\ PTx\end{matrix}}=0
\quad\Longrightarrow\quad Tx\in \Ker P.
$$
Вместе эти две цепочки означают выполнение условия \eqref{TP=PT-2}, которое опять по лемме \ref{LM:TP=PT} эквивалентно равенству $[T,P]=0$.
\epr

\bpr[Доказательство теоремы \ref{TH:stabilizatsiya-Z^n}] 1. Пусть сначала $B$ -- алгебра ограниченных операторов на каком-нибудь гильбертовом пространстве $X$: $B={\mathcal B}(X)$. Заметим, что при замене $A$ на его бикоммутант $A^{!!}$ в ${\mathcal B}(X)$ последовательность подпространств $Z^n(A)$ не меняется:
$$
Z^n(A^{!!})=Z^n(A).
$$
Поэтому можно с самого начала считать, что $A$ совпадает со своим бикомутантом, то есть является алгеброй фон Неймана в ${\mathcal B}(X)$. Тогда $A$ порождается системой своих ортогональных проекторов $P\in A$: $A=\{P\}^{!!}$. Для всякого такого проектора $P$ и любого элемента $T\in Z^2(A)$ мы получаем равенство $[[T,P],P]=0$, которое по лемме \ref{LM:TP=PT-(2)} влечет $[T,P]=0$. Поскольку проекторы $P$ порождают $A$, это влечет равенство
$$
[T,a]=0,\qquad a\in A.
$$
То есть $T\in Z^1(A)$, и мы получаем включение $Z^1(A)\supseteq Z^2(A)$. Оно эквивалентно равенству $Z^1(A)=Z^2(A)$, которое в свою очередь влечет \eqref{Z^n(A)=Z^(n+1)(A)}.

2. Рассмотрим теперь общий случай, когда $B$ -- произвольная $C^*$-алгебра. Тогда $B$ вкладывается в некоторую алгебру ${\mathcal B}(X)$, и мы получаем цепочку
$$
A\subseteq B\subseteq {\mathcal B}(X).
$$
Пусть, как и раньше $Z^n(A)$ -- последовательность подпространств в $B$, определенных равенствами \eqref{Z^n(A)}, а $Z^n_X(A)$ -- такая же последовательность в ${\mathcal B}(X)$ (то есть подпространства, определяемые равенствами \eqref{Z^n(A)}, в которых $B$ заменено на ${\mathcal B}(X)$). Покажем, что
 \beq\label{Z^n(A)=B-cap-Z^n_H(A)}
Z^n(A)=B\cap Z^n_X(A),\qquad n\in \Z_+.
 \eeq
Для $n=0$ это верно по определению. Предположим, мы доказали это для всех индексов, не превышающих какого-то $n\in\Z_+$. Тогда для $n+1$ мы получим цепочку:
\begin{multline*}
c\in B\cap Z^{n+1}_X(A)\quad\Leftrightarrow\quad \begin{cases}c\in B \\ \forall a\in A\quad [c,a]\in Z^n_X(A) \end{cases}
\quad\Leftrightarrow\quad \begin{cases} c\in B \\ \forall a\in A\quad  [c,a]=c\cdot a-a\cdot c\in B \\ \forall a\in A\quad  [c,a]\in Z^n_X(A) \end{cases}
\quad\Leftrightarrow\\ \Leftrightarrow\quad \forall a\in A\quad [c,a]\in B\cap Z^n_X(A)\kern-20pt\underset{\tiny\begin{matrix}\uparrow\\ \text{предположение}\\ \text{индукции}\end{matrix}}{=}\kern-20pt Z^n(A) \quad\Leftrightarrow\quad c\in Z^{n+1}(A).
\end{multline*}

Для последовательности $Z^n_X(A)$ мы уже доказали равенства
$$
Z^1_X(A)=Z^n_X(A),\qquad n\ge 1.
$$
Вместе с формулой \eqref{Z^n(A)=B-cap-Z^n_H(A)} они влекут
$$
Z^n(A)=B\cap Z^n_X(A)=B\cap Z^1_X(A)=Z^1(A).
$$
\epr

Из теоремы \ref{TH:stabilizatsiya-Z^n} следует

\btm\label{TH:stabilizatsiya-Z^n(ph)} Если $A$ -- инволютивная стереотипная алгебра, а $B$ -- $C^*$-алгебра, то для всякого инволютивного гомоморфизма стереотипных алгебр $\ph:A\to B$ последовательность подпространств $Z^n(\ph)$, определенных формулами \eqref{Z^0} и \eqref{Z^n}, стабилизируется, начиная с номера 1:
\beq\label{Z^n(ph)=Z^(n+1)(ph)}
Z^n(\ph)=Z^{n+1}(\ph),\qquad n\ge 1.
\eeq
\etm
\bpr
Замкнутый образ $\overline{\ph(A)}$ алгебры $A$ при отображении $\ph$ является инволютивной подалгеброй в $B$, причем
$$
Z^n(\ph)=Z^n\Big(\overline{\ph(A)}\Big).
$$
Поэтому \eqref{Z^n(ph)=Z^(n+1)(ph)} является просто следствием \eqref{Z^n(A)=Z^(n+1)(A)}.
\epr

\subsection{Касательное и кокасательное расслоения}

 \bit{
 \item[$\bullet$] {\it Комплексным характером}\label{DEF:Spec_C(G)}\index{характер!комплексный} стереотипной алгебры $A$ над $\C$ мы называем произвольный (непрерывный и сохраняющий единицу) гомоморфизм $s:A\to \C$. Множество  всех комплексных характеров алгебры $A$ мы называем ее {\it комплексным спектром} и обозначаем  ${\C}\Spec(A)$. Это множество мы наделяем топологией равномерной сходимости на вполне ограниченных множествах в $A$.

 \item[$\bullet$] {\it Вещественным характером}\label{DEF:RSpec(G)}\index{характер!вещественный} стереотипной алгебры $A$ над $\R$ мы называем произвольный (непрерывный и сохраняющий единицу) гомоморфизм $s:A\to \R$. Множество  всех вещественных характеров алгебры $A$ мы называем ее {\it вещественным спектром} и обозначаем  ${\R}\Spec(A)$. Это множество мы также наделяем топологией равномерной сходимости на вполне ограниченных множествах в $A$.

 \item[$\bullet$] {\it Инволютивным характером}\label{DEF:InvSpec(G)}\index{характер!инволютивный} стереотипной алгебры $A$ над $\C$ с инволюцией $\bullet$ мы называем произвольный (непрерывный и сохраняющих единицу) гомоморфизм $s:A\to \C$ 
    \beq\label{DEF:invol-harakter}
    s(x^\bullet)=\overline{s(x)},\qquad x\in A.
    \eeq      
    Множество всех инволютивных характеров алгебры $A$ мы называем ее {\it инволютивным спектром} и обозначчаем $\Spec(A)$. Это множество мы также наделяем топологией равномерной сходимости на вполне ограниченных множествах в $A$.  Можно заметить, что инволютивный спектр изоморфен вещественной части комплексного спектра и вещественному спектру от вещественной части алгебры $A$:
\beq\label{DEF:Spec[A]}
\Spec[A]=\Real({\C}\Spec[A])\cong {\R}\Spec[\Real A].
\eeq

 \item[$\bullet$]
{\it Комплексным} {\it касательным вектором}\index{касательный вектор!комплексный} к стереотипной алгебре $A$ над $\C$ в точке $s\in\C\Spec(A)$ называется всякий линейный непрерывный функционал $\tau:A\to\C$, удовлетворяющий условию Лейбница:
 \beq\label{DEF:kasatelnyj-vektor}
 \tau(a\cdot b)=s(a)\cdot\tau(b)+\tau(a)\cdot s(b), \qquad
 a,b\in A
 \eeq
Множество всех комплексных касательных векторов к $A$ в точке $s$
называется {\it комплексным} {\it касательным пространством} к $A$ в точке
$s\in\C\Spec(A)$ и обозначается ${\C}T_s[A]$. Оно наделяется топологией,
представляющей собой псевдонасыщение топологии равномерной сходимости на вполне
ограниченных множествах в $A$. Как следствие, ${\C}T_s[A]$ --- стереотипное
пространство и непосредственное подпространство в $A^\star$.

 \item[$\bullet$]
{\it Вещественным} {\it касательным вектором}\index{касательный вектор!вещественный} к стереотипной алгебре $A$ над $\R$ в точке $s\in\R\Spec(A)$ называется всякий линейный непрерывный функционал $\tau:A\to\R$, удовлетворяющий условию Лейбница:
 \beq\label{DEF:kasatelnyj-vektor-R}
 \tau(a\cdot b)=s(a)\cdot\tau(b)+\tau(a)\cdot s(b), \qquad
 a,b\in A
 \eeq
Множество всех вещественных касательных векторов к $A$ в точке $s$
называется ({\it вещественным}) {\it касательным пространством} к $A$ в точке
$s\in\R\Spec(A)$ и обозначается ${\R}T_s[A]$. Оно наделяется топологией,
представляющей собой псевдонасыщение топологии равномерной сходимости на вполне
ограниченных множествах в $A$. Как следствие, ${\R}T_s[A]$ --- стереотипное
пространство и непосредственное подпространство в $A^\star$.

 \item[$\bullet$]
{\it Инволютивным} {\it касательным вектором}\index{касательный вектор!инволютивный} к стереотипной алгебре $A$ над $\C$ с инволюцией $\bullet$ в точке $s\in\Spec(A)$ называется всякий линейный непрерывный функционал $\tau:A\to\C$, удовлетворяющий условию Лейбница:
 \beq\label{DEF:kasatelnyj-vektor}
 \tau(a\cdot b)=s(a)\cdot\tau(b)+\tau(a)\cdot s(b), \qquad
 a,b\in A
 \eeq
и дополнительному условию сохранения инволюции: 
 \beq\label{DEF:kasatelnyj-vektor-inv}
 \tau(a^\bullet)=\tau(a)^\bullet, \qquad a\in A.
 \eeq
Множество всех инволютивных касательных векторов к $A$ в точке $s$
называется {\it инволютивным} {\it касательным пространством} к $A$ в точке
$s\in\Spec(A)$ и обозначается $T_s[A]$. Оно наделяется топологией,
представляющей собой псевдонасыщение топологии равномерной сходимости на вполне
ограниченных множествах в $A$. Как следствие, $T_s[A]$ --- стереотипное
пространство и непосредственное подпространство в $A^\star$. Можно заметить, что инволютивное касательное пространство изоморфно  вещественной части комплексного касательного пространства и вещественному касательному пространству от вещественной части алгебры $A$: 
\beq\label{DEF:T_s[A]}
T_s[A]=\Real({\C}T_s[A])\cong {\R}T_s[\Real A].
\eeq

 }\eit

\bit{
 \item[$\bullet$]
{\it Комплексным кокасательным пространством}\index{касательное пространство!комплексное} инволютивной стереотипной алгебры
$A$ в точке $s\in\Spec(A)$ называется стереотипное фактор-пространство 
 \beq\label{DEF:T_s^(C-star)[A]}
{\C}T_s^\star[A]:=\Big(I_s[A]\Big/\overline{I_s^2}[A]\Big)^\smalltriangledown.
 \eeq
идеала
\beq
I_s[A]=\{a\in A:\ s(a)=0\}
\eeq 
по его замкнутому квадрату (подидеалу в $I_s[A]$)
\beq
\overline{I_s^2}[A]=\overline{\sp\{a\cdot b;\ a,b\in I_s[A]\}} 
\eeq
Элементы этого пространства называются {\it комплексными кокасательными
векторами}\index{кокасательный вектор!комплексный} (алгебры $A$ в точке $s\in\Spec(A)$).

 \item[$\bullet$]
Поскольку оба идеала $I_s[A]$ и $\overline{I_s^2}[A]$ замкнуты относительно
инволюции, формула
$$
(a+\overline{I_s^2}[A])^\bullet=a^\bullet+\overline{I_s^2}[A],\qquad a\in I_s[A]
$$
определяет некую инволюцию на фактор-пространстве
${\C}T_s^\star[A]=\Big(I_s[A]/\overline{I_s^2}[A]\Big)^\smalltriangledown$ (как
стереотипном $A$-бимодуле).\footnote{Понятно, что это определение выбрано так, чтобы
фактор-отображение $\pi:I_s\to (I_s/\overline{I_s^2})^\smalltriangledown$ переводило
инволюцию в $I_s$, индуцированную из $A$, в определенную таким образом
инволюцию на $(I_s/\overline{I_s^2})^\smalltriangledown$:
$$
\pi(a^\bullet)=\pi(a)^\bullet.
$$
}

\item[$\bullet$] {\it Кокасательным вектором} (или {\it вещественным
кокасательным вектором}) алгебры $A$ в точке $s\in\Spec(A)$ называется
произвольный вещественный вектор в
пространстве ${\C}T_s^\star[A]$, то есть произвольный комплексный кокасательный
вектор $\xi\in {\C}T_s^\star[A]$, инвариантный относительно инволюции $\bullet$ в
${\C}T_s^\star[A]$:
 \beq\label{DEF:kokasat-vektor}
\xi^\bullet=\xi.
 \eeq
Множество всех кокасательных векторов алгебры $A$ в точке $s\in\Spec(A)$
называется {\it кокасательным пространством} (алгебры $A$ в точке
$s\in\Spec(A)$) и обозначается $T_s^\star[A]$. Очевидно, $T_s^\star[A]$
является вещественной частью ${\C}T_s^\star[A]$:
$$
T_s^\star[A]=\Real{\C}T_s^\star[A].
$$

 }\eit

\paragraph{Кокасательное расслоение $T^\star[A]$.}

Пусть $A$ -- инволютивная стереотипная алгебра, и $t\in\Spec(A)$. Заметим, что всякая непрерывная полунорма $p:A\to\R_+$ определяет полунорму $p'_t$ на фактор-пространстве $\C T^\star_t[A]=\big(
I_t/\overline{I_t^2}\big)^\triangledown$ по формуле
 \beq\label{polunorma-na-CT^*(X)}
p'_t(x+\overline{I_t^2}):=\inf_{y\in
\overline{I_t^2}}p(x+y),\qquad x\in I_t.
 \eeq
Рассмотрим теперь прямую сумму множеств $\C T^\star_t(A)$
$$
\C T^\star(A)=\bigsqcup_{t\in\Spec(A)} \C T^\star_t(A)=\bigsqcup_{t\in\Spec(A)} \Big(
I_t/\ \overline{I_t^2}\ \Big)^\triangledown
$$
и обозначим через $\pi$ естественную проекцию $\C T^\star[A]$ на $\Spec(A)$:
$$
\pi:\C T^\star[A]\to \Spec(A),\qquad \pi\Big( x+\overline{I_t^2}\Big)=t,\qquad t\in \Spec(A),\ x\in I_t[A].
$$
Кроме того, для всякого вектора $x\in A$ мы рассмотрим отображение
$$
T^\star(x):\Spec(A)\to\C T^\star(X)\quad\Big|\quad T^\star(x)(t)=x-t(x)\cdot
1_A+\overline{I_t^2}.
$$
Понятно, что для всякого $x\in A$
$$
\pi\circ T^\star(x)=\id_{\Spec(A)},
$$
т.е. $T^\star(x)$ является сечением векторного расслоения $\pi$.

\blm\label{LM:poluneprer-sverhu-p^*_t(T^*(x)(t))} Для всякого элемента $x\in A$
и любой непрерывной полунормы $p:A\to\R_+$ отображение $t\in\Spec(A)\mapsto
p'_t(T^\star(x)(t))\in\R_+$ полунепрерывно сверху. \elm \bpr Пусть $t\in\Spec(A)$ и
$\e>0$. Условие
$$
p'_t(T^\star(x)(t))=\inf_{y\in
\overline{I_t^2}}p(x-t(x)\cdot 1_A+y)=\inf_{y\in I_t^2}p(x-t(x)\cdot 1_A+y)<\e
$$
означает, что для некоторого $y\in I_t^2$ выполняется неравенство
$$
p(x-t(x)\cdot 1_A+y)<\e.
$$
Это в свою очередь означает, что существует число $m\in\N$ и векторы
$a^1,...,a^m,b^1,...,b^m\in I_t$ такие, что
 \beq\label{p-l-x+sum_(k=1)^m-a^k_i-cdot-b^k-r<e}
p\l x-t(x)\cdot 1_A+\sum_{k=1}^m a^k\cdot b^k\r<\e.
 \eeq
Для любых чисел
$$
\lambda\in\C,\qquad \alpha=\{\alpha^k;\ 1\le k\le m\}\subseteq\C,
\beta=\{\beta^k;\ 1\le k\le m\}\subseteq\C
$$
обозначим
$$
f(\lambda,\alpha,\beta)=p\l x-\lambda\cdot 1_A+\sum_{k=1}^m (a^k-\alpha^k\cdot
1_A)\cdot (b^k-\beta^k\cdot 1_A)\r.
$$
Функция $(\lambda,\alpha,\beta)\mapsto f(\lambda,\alpha,\beta)$ в точке
$(\lambda,\alpha,\beta)=(t(x),0,0)$ совпадает с левой частью
\eqref{p-l-x+sum_(k=1)^m-a^k_i-cdot-b^k-r<e} и поэтому удовлетворяет
неравенству
$$
f(t(x),0,0)<\e.
$$
С другой стороны, она непрерывна, как композиция многочлена от $2m+1$
комплексных переменных со значениями в локально выпуклом пространстве $X$ и
непрерывной функции $p$ на $X$. Значит, должно существовать число $\delta>0$
такое, что
 \beq\label{max_k|lambda^k|<delta=>f(lambda)<e} \forall \lambda,\alpha,\beta \qquad
 |\lambda-t(x)|<\delta\quad\&\quad\max_k\abs{\alpha^k}<\delta
 \quad\&\quad\max_k\abs{\beta^k}<\delta\quad\Longrightarrow\quad f(\lambda,\alpha,\beta)<\e.
 \eeq

Рассмотрим далее множество
$$
U=\{s\in\Spec(A):\ \abs{s(x)-t(x)}<\delta\quad \&\quad
\max_k\abs{s(a^k)}<\delta\quad \&\quad \max_k\abs{s(b^k)}<\delta\}.
$$
Оно открыто и содержит точку $t$ (потому что $\abs{t(x)-t(x)}=0<\delta$, и
включения $a^k,b^k\in I_t$ означают систему равенств $t(a^k)=t(b^k)=0$, $1\le
k\le m$). С другой стороны, для всякой точки $s\in U$, рассмотрев числа
$$
\lambda=s(x),\qquad \alpha^k=s(a^k),\qquad \beta^k=s(b^k),
$$
мы получим, во-первых, $s(a^k-\alpha^k\cdot 1_A)=s(a^k)-\alpha^k\cdot
s(1_A)=s(a^k)-s(a^k)=0$, то есть
$$
a^k-\alpha^k\cdot 1_A\in I_s
$$
во-вторых, по тем же причинам,
$$
b^k-\beta^k\cdot 1_A\in I_s
$$
и, в-третьих,
$$
\abs{\lambda-t(x)}<\delta,\qquad
\max_k\abs{\alpha^k}=\max_k\abs{s(a^k)}<\delta, \qquad
\max_k\abs{\beta^k}=\max_k\abs{s(b^k)}<\delta,
$$
то есть, в силу \eqref{max_k|lambda^k|<delta=>f(lambda)<e},
$$
p\bigg( x-\lambda\cdot 1_A+\sum_{k=1}^m  \underbrace{(a^k-\alpha^k\cdot
1_A)}_{\scriptsize\begin{matrix}\text{\rotatebox{90}{$\owns$}}\\
I_s\end{matrix}}\cdot \underbrace{(b^k-\beta^k\cdot
1_A)}_{\scriptsize\begin{matrix}\text{\rotatebox{90}{$\owns$}}\\
I_s\end{matrix}}\bigg)=f(\lambda,\alpha,\beta)<\e
$$
Это можно понимать так, что для некоторой точки $z\in I_s^2$ выполняется
неравенство
$$
p(x-s(x)\cdot 1_A+z)<\e,
$$
которое, в свою очередь, влечет за собой неравенство
$$
p_s(T^\star(x)(s))=p_s(x-s(x)\cdot 1_A+\overline{I_s^2})=\inf_{z\in
\overline{I_s^2}}p(x-s(x)\cdot 1_A+z)=\inf_{z\in I_s^2}p(x-s(x)\cdot 1_A+z)<\e
$$
Оно верно для всякой точки $s$ из окрестности $U$ точки $t$, и как раз это нам
и нужно было доказать. \epr

В следующей теореме мы описываем топологию на кокасательном расслоении
$T^\star(A)=\bigsqcup_{t\in\Spec(A)} T^\star_t(A)$, однако точно теми же рассуждениями вводится топология на комплексном касательном расслоении $\C T^\star(A)=\bigsqcup_{t\in\Spec(A)} \C T^\star_t(A)$.

\btm\label{TH:kokasatelnoe-rassloenie} Для всякой инволютивной стереотипной
алгебры $A$ прямая сумма стереотипных фактор-модулей
$$
T^\star(A)=\bigsqcup_{t\in\Spec(A)} T^\star_t(A)=\bigsqcup_{t\in\Spec(A)} \Real\Big(
I_t/\ \overline{I_t^2}\ \Big)^\triangledown
$$
обладает единственной топологией, превращающей проекцию
$$
\pi:T^\star(A)\to\Spec(A),\qquad \pi(x+\overline{I_t^2})=t,\qquad
t\in\Spec(A),\ x\in A
$$
в локально выпуклое расслоение с системой полунорм $\{p';\ p\in{\mathcal P}(A)\}$, для которого отображение
$$
x\in A\mapsto T^\star(x)\in\Sec(\pi)\quad\Big|\quad
T^\star(x)(t)=x+\overline{I_t^2},\quad t\in\Spec(A),
$$
непрерывно переводит $A$ в стереотипное пространство $\Sec(\pi)$ непрерывных сечений
$\pi$. При этом
 базу топологии $T^\star(A)$ образуют множества
$$
W(x,U,p,\e)=\set{\xi\in T^\star(A): \ \pi(\xi)\in U\ \& \
p_{\pi(\xi)}\Big(\xi-T^\star(x)\big(\pi(\xi)\big)\Big)<\e }
$$
где $x\in A$, $p\in{\mathcal P}(A)$, $\e>0$, $U$ -- открытое множество в $M$;
 \etm

\bit{

\item[$\bullet$] Расслоение $\pi: T^\star(A)\to \Spec(A)$ называется {\it
кокасательным расслоением} алгебры $A$.

}\eit

\bpr Из леммы \ref{LM:poluneprer-sverhu-p^*_t(T^*(x)(t))} и предложения
\ref{PROP:sushestv-topologii-v-Xi} следует существование и единственность
топологии на $T^\star(A)$, для которой проекция $\pi:T^\star(A)\to\Spec(A)$
представляет собой локально выпуклое расслоение с полунормами $p'$, а сечения
вида $T^\star(x)$, $x\in X$, будут непрерывными. Непрерывность отображения $x\in
A\mapsto T^\star(x)\in\Sec(\pi)$ доказывается импликацией
$$
p'_t(\pi(x)(t))=\inf_{y\in I_t^2}p(x+y)\le p(x)\quad \Longrightarrow \quad
\sup_{t\in T}p'_t(\pi^n(x)(t))\le p(x)
$$
для всякого компакта $T\subseteq M$. Структура базы топологии в $T^\star(A)$
также следует из предложения \ref{PROP:sushestv-topologii-v-Xi}. \epr

\paragraph{Касательное расслоение $T[A]$.}

Пусть $\ph:A\to B$ -- гомоморфизм инволютивных стереотипных алгебр. Линейное
(над $\C$) непрерывное отображение $D:A\to B$ называется {\it
дифференцированием} из $A$ в $B$ относительно $\ph$, если справедливо тождество
 \beq\label{DEF:differentsirovanie}
D(x^\bullet)=D(x)^\bullet,\qquad D(x\cdot y)=D(x)\cdot\ph(y)+\ph(x)\cdot D(y),\qquad x,y\in
A.
 \eeq
Множество всех дифференцирований из $A$ в $B$ относительно $\ph$ мы обозначаем
$\Der_{\ph}(A,B)$. В частном случае, если $A=B$ и $\ph=\id_A$, принято говорить
просто о дифференцированиях в $A$, и при этом используется обозначание
$$
\Der(A):=\Der_{\id_A}(A,A).
$$
Пространства $\Der_{\ph}(A,B)$ и $\Der(A)$ наделяются топологиями
непосредственных подпространств в пространствах операторов $\Mor(A,B)$ и
$\Mor(A,A)$.

\btm\label{TH:differentsirovanie->rassl-struuj} Если множество значений
дифференцирования $D:A\to B$ лежит в центре алгебры $B$,
 \beq\label{D(A)-subseteq-Z(A)}
D(A)\subseteq Z(B),
 \eeq
то $D$ является дифференциальным оператором порядка 1. Как следствие, в этом
случае $D$ определяет некий морфизм расслоений струй $\jet_1[D]:\Jet_A^1(A)\to
\Jet_A^0(B)$, удовлетворяющий тождеству
 \beq\label{j^0(Dx)=j_1[D]-circ-j^1(x)}
\jet^0(Dx)=\jet_1[D]\circ \jet^1(x),\qquad x\in A.
 \eeq
 \etm
\bpr Из тождества \eqref{DEF:differentsirovanie} следует
$$
[D,a](x)=D(a\cdot x)-a\cdot D(x)=D(a)\cdot x+a\cdot D(x)-a\cdot D(x)=D(a)\cdot
x,
$$
из которого, в свою очередь,
$$
[[D,a],b](x)=[D,a](b\cdot x)-b\cdot [D,a](x)=D(a)\cdot b\cdot x-b\cdot
D(a)\cdot x=\eqref{D(A)-subseteq-Z(A)}=0.
$$
То есть $D$ является дифференциальным оператором порядка 1. После этого
применяется теорема \ref{TH:diff-oper->rassl-struuj}.
 \epr

\btm
Всякое дифференцирование $D:A\to A$ алгебры $A$ в каждой точке спектра $t\in\Spec(A)$ определяет касательный вектор по формуле
$$
D^T(t)(x)=t(D(x)),\qquad x\in A.
$$
\etm
\bpr
Действительно,
$$
D^T(t)(a\cdot b)=t(D(a\cdot b))=t(Da\cdot b+a\cdot Db)=t(Da)\cdot
t(b)+t(a)\cdot t(Db)=D^T(t)(a)\cdot t(b)+t(a)\cdot D^T(t)(b).
$$
\epr

Из предложения \ref{PROP:o-dvoistvennom-rasloenii} следует

\btm\label{kasatelnoe-rassloenie} Пусть у инволютивной стереотипной алгебры $A$ спектр $\Spec(A)$ является хаусдорфовым пространством, а множество $\Der(A)$ дифференцирований имеет плотный след в каждом касательном пространстве:
$$
\overline{\{D_t;\ D\in\Der(A)\}}=T_t[A].
$$
Тогда прямая сумма стереотипных пространств
$$
T[A]=\bigsqcup_{t\in\Spec(A)}T_t[A]
$$
обладает топологией, превращающей проекцию
$$
T_A:T[A]\to\Spec(A),\qquad T_A(\tau)=t,\ \tau\in T_t[A]
$$
в сопряженное расслоение к кокасательному расслоению $T^\star[A]$. При этом отображение
$$
D\in\Der(A)\mapsto D^T\in\Sec(T[A])\quad\Big|\quad D^T(s)(a)=s(Da),\qquad a\in
A,\quad s\in\Spec(A),
$$
непрерывно переводит $\Der(A)$ в стереотипное пространство $\Sec(T[A])$ непрерывных сечений
расслоения $T[A]$.  \etm

\bit{

\item[$\bullet$] Расслоение $T_A:T[A]\to\Spec(A)$ называется {\it касательным
расслоением} инволютивной стереотипной алгебры $A$.

}\eit

\chapter{АНАЛИЗ НА ГРУППАХ}

\section{Групповые алгебры на локально компактных группах}

\subsection{Важные классы групп}

\paragraph{LP-группы и решетка $\lambda(G)$.}

Локально компактная группа $G$ называется {\it LP-группой} (или {\it Ли-проективной группой}), если она представима как проективный предел групп Ли $G_i$:
\beq\label{DEF:LP-gruppa}
G=\projlim_{i\to\infty}G_i.
\eeq 
Пусть $G$ --- LP-группа. Зафиксируем какую-нибудь компактную нормальную подгруппу $H_0\subseteq G$, фактор-группа по которой $G/H_0$ является группой Ли. Пусть далее $\lambda(G)$ обозначает систему компактных нормальных подгрупп $K$ в локально компактной группе $G$, фактор-группы по которым $G/K$ являются группами Ли, и которые содержатся в $H_0$.

\btm\label{TH:K-wedge-L=K-cap-L}
Множество $\lambda(G)$, будучи упорядочено по включению, 
$$
K\le L\quad\Leftrightarrow\quad K\subseteq L,\qquad K,L\in\lambda(G)
$$
образует решетку: для любых двух групп $K,L\in\lambda(G)$
\bit{

\item[---] их точной нижней гранью является их пересечение:
\beq\label{K-wedge-L=K-cap-L}
K\wedge L=K\cap L
\eeq

\item[---] а точной верхней гранью является их поэлементное произведение:
\beq\label{K-vee-L=K-cdot-L}
K\vee L=K\cdot L
\eeq
}\eit 
\etm
Здесь для доказательства используется следующая

\blm\label{LM:K-cdot-L-komp-norm-podgr} Если $K,L\subseteq G$ --- две компактные нормальные подгруппы в локально компактной группе $G$, то их поэлементное произведение 
$$
K\cdot L=\{a\cdot x;\ a\in K,\ x\in L\}
$$
--- тоже компактная нормальная подгруппа в $G$. 
\elm
\bpr
Пусть $a,b\in K$ и $x,y\in L$. Тогда
$$
\underbrace{a\cdot x}_{\scriptsize\begin{matrix}\text{\rotatebox{90}{$\owns$}}\\ K\cdot L\end{matrix}}\cdot 
\underbrace{b\cdot y}_{\scriptsize\begin{matrix}\text{\rotatebox{90}{$\owns$}}\\ K\cdot L\end{matrix}}=
\underbrace{a\cdot b}_{\scriptsize\begin{matrix}\text{\rotatebox{90}{$\owns$}}\\ K \end{matrix}}\cdot 
\underbrace{\underbrace{b^{-1}\cdot x\cdot b}_{\scriptsize\begin{matrix}\text{\rotatebox{90}{$\owns$}}\\ L\end{matrix}} \cdot y}_{\scriptsize\begin{matrix}\text{\rotatebox{90}{$\owns$}}\\ L\end{matrix}}\in K\cdot L
$$
--- и это доказывает,  что $K\cdot L$ --- подгруппа в $G$. А с другой стороны, для любого $g\in G$
$$
g^{-1}\cdot K\cdot L\cdot g=g^{-1}\cdot K\cdot g\cdot g^{-1}\cdot L\cdot g\subseteq K\cdot L,
$$
--- и это доказывает, что $K\cdot L$ --- нормальная подгруппа в $G$.
\epr

\btm\label{TH:LP-gruppa-lambda}
Если $G$ --- LP-группа, то она является проективным пределом своих фактор-групп (Ли) по подгруппам системы $\lambda(G)$:
\beq\label{DEF:LP-gruppa-lambda}
G=\projlim_{K\in\lambda(G)}G/K.
\eeq 

\etm

\btm\label{TH:LP-gruppa-v-LC-gruppe}
Всякая локально компактная группа $G$ содержит некоторую открытую LP-подгруппу $G_0$.
\etm

\paragraph{Центральные группы.}
Напомним, что центром $Z(G)$ группы $G$ называется множество элементов $a\in G$, коммутирующих со всем остальными элементами:
$$
a\cdot x=x\cdot a,\qquad x\in G.
$$
Локально  компактная группа $G$ называется {\it центральной}\label{DEF:zentralnaya-gruppa}, если ее фактор-группа по центру $G/Z(G)$ компактна.

\btm[Х.Фрейденталь \cite{Freudenthal}, см. тж.\cite{Palmer}]\label{TH:stroenie-zentralnoi-gruppy}
Связная локально компактная группа $G$ центральна тогда и только тогда, когда она является прямым произведением векторной группы $\R^n$ и компактной группы $K$:
 \beq\label{stroenie-zentralnoi-gruppy}
 G=\R^n\times K.
 \eeq
\etm

\paragraph{Компактные надстройки абелевых групп.}

Условимся называть локально компактную группу $G$ {\it компактной надстройкой абелевой группы}\label{DEF:komp-nadstr-abelevoi-gruppy}, если в $G$ существуют замкнутые подгруппы $Z$ и $K$ со следующими свойствами:
\bit{

\item[1)] $Z$ -- абелева группа,

\item[2)] $K$ -- компактная группа,

\item[3)] $Z$ и $K$ коммутируют:
$$
\forall a\in Z,\quad \forall y\in K\qquad a\cdot y=y\cdot a,
$$

\item[4)] $Z$ и $K$ в произведении дают $G$:
$$
\forall x\in G\qquad \exists a\in Z\quad\exists y\in K\qquad x=a\cdot y.
$$
}\eit
Если нужно уточнить, какие группы в этой конструкции используются, то про группу $G$ мы будем говорить, что она является {\it  надстройкой абелевой группы $Z$ с помощью компактной группы $K$}. Можно заметить, что всякая такая группа $G$ центральна.

\blm\label{LM:Z-cdot-K=(Z-times-K)/C}
Пусть $G=Z\cdot K$ --- надстройка абелевой локально компактной группы $Z$ с помощью компактной группы $K$. Рассмотрим подгруппу
$$
H=Z\cap K,
$$
ее вложение в декартово произведение
$$
\iota:H\to Z\times K\quad\Big|\quad \iota(x)=(x,x^{-1}),\quad x\in H,
$$
и образ этого вложения
$$
\iota(H)\subseteq Z\times K.
$$
Тогда
$$
G\cong(Z\times K)/\iota(H).
$$
\elm
\bpr
Рассмотрим отображение $\varkappa:Z\times K\to G$, действующее по формуле
$$
\varkappa(x,y)=x\cdot y.
$$
Из условия 3) следует, что оно является гомоморфизмом групп, а из условия 4) --- что оно сюръективно. Его ядром является группа
\begin{multline*}
\Ker\varkappa=\{(x,y)\in Z\times K: \ x\cdot y=e\} =\{(x,y)\in Z\times K: \ y=x^{-1}\}=\\=
 \{(x,x^{-1}); \ x\in Z,\ x^{-1}\in K\}= \{(x,x^{-1}); \ x\in Z\cap K=H\}=\iota(H).
\end{multline*}
Это означает, что у нас имеется точная последовательность групп
$$
     \xymatrix 
 {
 1\ar[r] & H\ar[r]^{\iota} & Z\times K\ar[r]^{\varkappa} & Z\cdot K\ar[r] & 1.
 }
$$
Как следствие, гомоморфизм $\varkappa$ поднимается до некоторого гомоморфизма $\sigma$ на фактор-группу:
$$
     \xymatrix 
 {
 1\ar[r] & H\ar[r]^{\iota} & Z\times K\ar[d]_{\pi}\ar[r]^{\varkappa} & Z\cdot K\ar[r] & 1 \\
  & & (Z\times K)/\iota(H)\ar@{-->}[ru]_{\sigma} & &
 }
$$
причем этот гомоморфизм $\sigma$ должен быть изоморфизмом групп. Остается только проверить, что $\sigma$ непрерывен в обе стороны.

Непрерывность $\sigma$ в прямую сторону следует из непрерывности $\varkappa$: для всякой окрестности единицы $W\subseteq Z\cdot K$
$$
x\in \varkappa^{-1}(W)\quad\Leftrightarrow\quad \varkappa(x)=\sigma(\pi(x))\in W \quad\Leftrightarrow\quad
\pi(x)\in \sigma^{-1}(W)
$$
То есть
$$
\sigma^{-1}(W)=\pi\big(\varkappa^{-1}(W)\big),
$$
и это будет открытое множество, как образ открытого множества $\varkappa^{-1}(W)$ при открытом отображении $\pi$.

А непрерывность $\sigma$ в обратную сторону, то есть открытость $\sigma$, следует из открытости $\varkappa$: если $U\subseteq Z$ и $V\subseteq K$ --- окрестности единицы, то
$$
\varkappa(U,V)=U\cdot V
$$
--- открытое множество, поэтому $\varkappa$ --- открытое отображение, поэтому для любой окрестности единицы $W\subseteq (Z\times K)/\iota(H)$ множество
$$
\sigma(W)=\varkappa(\pi^{-1}(W))
$$
должно быть открыто.
\epr

\btm[Кузнецова \cite{Kuznetsova}]\label{TH:zentral-Lie}
Всякая центральная группа Ли $G$ является конечным расширением некоторой компактной надстройки абелевой группы.
\etm
\bpr
В силу \cite[Theorem 4.4]{Grosser-Moskowitz-1}, группа $G$ представима в виде $G=\R^n\times G_1$, где $G_1$ содержит открытую компактную нормальную подгруппу $K$. Пусть $C$ -- центр $G$. Отождествим $K$ с $\{0\}\times K$ и рассмотрим группу $H=C\cdot K$. Понятно, она будет компактной надстройкой абелевой группы. Кроме того, она содержит открытую группу $\R^n\times K$, поэтому она открыта. Из нормальности $K$ в $G_1$ следует, что $K=\{0\}\times K$ нормальна в $G$, поэтому $H$ также нормальна в $G$:
$$
x\cdot H=x\cdot C\cdot K=C\cdot x\cdot K=C \cdot K\cdot x=H\cdot x.
$$
Поскольку группа $G$ центральна, ее фактор-группа по центру $G/C$ должна быть компактна. Как следствие, фактор-группа $G/H$ тоже должна быть компактна (потому что она является непрерывным образом компактной группы $G/C$). С другой стороны, поскольку $H=C\cdot K$ открыта, $G/H$ должна быть дискретна. Значит, $G/H$ конечна.
\epr

\paragraph{SIN-группы.}

Множество $U$ в группе $G$ называется {\it нормальным}, если оно инвариантно относительно сопряжений:
$$
a\cdot U\cdot a^{-1}\subseteq U,\qquad a\in G.
$$
Локально компактная группа $G$ называется {\it SIN-группой}\label{DEF:SIN-gruppa}, если в ней нормальные окрестности единицы образуют локальную базу. Эквивалентное условие: левая и правая равномерные структуры на $G$ эквивалентны. В частности, все такие группы унимодулярны.

Класс SIN-групп включает абелевы, компактные и дискретные группы. Следующее утверждение принадлежит С.Гроссер и М.Московиц \cite[2.13]{Grosser-Moskowitz-2}:

\btm\label{TH:stroenie-SIN}
Всякая SIN-группа $G$ является дискретным расширением группы $\R^n\times K$, где $n\in\Z_+$, а $K$ -- компактная группа:
\beq\label{SIN-kak-rasshirenie}
1\to \R^n\times K=N\to G\to D\to 1
\eeq
($D$ -- дискретная группа).
\etm

\paragraph{Группы Мура.}

Локально компактная группа $G$ называется {\it группой Мура}\label{DEF:Moore-gruppa}, если у нее все непрерывные (в обычном смысле) унитарные неприводимые представления конечномерны.

\btm\label{TH:Moore->SIN}
Всякая группа Мура является SIN-группой.\footnote{См. \cite[p.1452]{Palmer}.}
\etm

\btm\label{TH:Moore->AM}
Всякая группа Мура аменабельна.\footnote{См. \cite[p.1486]{Palmer}.}
\etm

\btm\label{TH:Moore->quotient}
Всякая (отделимая) фактор-группа $G/H$ произвольной группы Мура $G$ является группой Мура.\footnote{Это утверждение очевидно.}
\etm

\bcor\label{TH:G=Moore->D=Moore}
Если $G$ -- группа Мура, то в ее представлении \eqref{SIN-kak-rasshirenie} группа $D$ -- тоже группа Мура (и поэтому, в частности, $D$ аменабельна).
\ecor

\btm\label{TH:D=Moore->konech-rassh-Abelevoi}
Всякая дискретная группа Мура является конечным расширением абелевой группы.\footnote{См. \cite[Theorem 12.4.26 and p.1397]{Palmer}.}
\etm

\paragraph{Группы Ли---Мура.}

{\it Группой Ли---Мура} мы называем произвольную группу Ли $G$, которая является одновременно группой Мура.

\btm\footnote{Ю.~Н.~Кузнецова \cite{Kuznetsova}.}\label{TH:Lie-Moore}
Всякая группа Ли---Мура $G$ является конечным расширением некоторой компактной надстройки абелевой группы.
\etm
\bpr
В силу \cite[Theorem 12.4.27]{Palmer}, $G$ является конечным расширением некоторой центральной группы $G_1$. По теореме \ref{TH:zentral-Lie}, $G_1$ является конечным расширением некоторой компактной надстройки абелевой группы.
\epr

\btm\label{TH:Moore=lim-Lie-Moore}
Всякая группа Мура является проективным пределом групп Ли---Мура, или точнее, своих фактор-групп Ли, которые автоматически являются группами Ли---Мура.
\etm

\brem
Это, в частности, означает, что всякая группа Мура является LP-группой.
\erem

\bpr
Согласно \cite[12.6.5]{Palmer}, группа Мура $G$ является проективным пределом своих фактор-групп Ли $G/H$, которые по теореме \ref{TH:Moore->quotient} должны быть группами Ли---Мура.
\epr

\subsection{Групповые алгебры ${\mathcal C}^\star(G)$ и ${\mathcal E}^\star(G)$}

Пусть $G$ -- локально компактная группа и $A$ -- стереотипная алгебра. Отображение $\pi:G\to A$ мы называем {\it представлением $G$ в $A$}, если оно непрерывно и мультипликативно:
$$
\pi(x\cdot y)=\pi(x)\cdot\pi(y),\qquad \pi(1)=1.
$$
В этом случае (см. \cite[Theorem 10.12]{Ak03}) существует единственный гомоморфизм стереотипных алгебр $\dot{\pi}:{\mathcal C}^\star(G)\to A$, замыкающий диаграмму
\begin{gather}\label{C*(G)=grup-alg-A}
\begin{diagram}
\node{G} \arrow[2]{e,t}{\delta} \arrow{se,b}{\pi}
\node[2]{\mathcal{C}^\star(G)} \arrow{sw,b}{\dot{\pi}}
\\
\node[2]{A}
\end{diagram}
\end{gather}
(здесь $\delta$ -- вложение в качестве дельта-функций).

Если $G$ -- вещественная группа Ли, то можно рассматривать также групповую алгебру ${\mathcal E}^\star(G)$ распределений с компактным носителем на $G$. Определяется она по аналогии: сначала рассматривается алгебра ${\mathcal E}(G)$ гладких функций на $G$ с обычной топологией равномерной сходимости на компактах по каждой частной производной в каждой локальной карте. Затем ${\mathcal E}^\star(G)$ определяется как стереотипное двойственное пространство к ${\mathcal E}(G)$. Те же формулы определяют на ${\mathcal E}^\star(G)$ структуру алгебры и инволюцию и точно так же между представлениями группы $G$ и алгебры ${\mathcal E}^\star(G)$ имеется связь в виде модифицированной диаграммы \eqref{C*(G)=grup-alg-A}
\begin{gather}\label{E*(G)=grup-alg-A}
\begin{diagram}
\node{G} \arrow[2]{e,t}{\delta} \arrow{se,b}{\pi}
\node[2]{\mathcal{E}^\star(G)} \arrow{sw,b}{\ddot{\pi}}
\\
\node[2]{A}
\end{diagram}
\end{gather}
с той только разницей, что под $\pi$ здесь понимается гладкое в (обычном смысле) отображение \cite[Theorem 10.12]{Ak03}.

\subsection{Интегральное исчисление в групповых алгебрах}

Пусть $G$ --- локально компактная группа, ${\mathcal C}_{00}(G)$ --- пространство непрерывных финитных функций на $G$.

\paragraph{Мера как интеграл.}
Действие меры $\alpha\in \mathcal{C}^\star (G)$ на функцию $u\in \mathcal{C}(G)$ обычно записывают в виде интеграла:
\beq\label{mu(u)=int_Gu(t)mu(dt)}
\alpha(u)=\int_G u(t)\ \alpha(\d t),\qquad u\in \mathcal{C}(G).
\eeq
Для функции $f\in \mathcal{C}(G)$ и меры $\alpha\in
\mathcal{C}^\star (G)$ мы обозначаем символом $f\cdot\mu$ меру, действующую по формуле
\beq\label{(f-cdot-alpha)(u)}
(f\cdot\alpha)(u)=\alpha(u\cdot f)=\int_G u(t)\cdot f(t)\ \alpha(\d t),\qquad u\in \mathcal{C}(G).
\eeq

\paragraph{Антипод и сдвиг.}
Для функций $u\in \mathcal{C}(G)$ и мер $\alpha\in
\mathcal{C}^\star (G)$ обозначим через $\widetilde{u}$ и $\widetilde{\alpha}$
их {\it антиподы}
\beq\label{widetilde(u)(t)=u(t^(-1))}
\widetilde{u}(t)=u(t^{-1}), \qquad \widetilde{\alpha}(u)=\alpha(\widetilde{u})
\eeq
а для всякого $a\in G$ символами $a\cdot u$, $u\cdot a$ и $a\cdot \alpha$, $\alpha\cdot a$ обозначаются  {\it сдвиги}:
\beq\label{DEF:u-cdot-a}
(a\cdot u)(t)=u(t\cdot a), \qquad (u\cdot a)(t)=u(a\cdot t) 
\eeq
\beq\label{eq10.6}
(a\cdot \alpha)(u)=\alpha(u\cdot a), \qquad (\alpha\cdot a)(u)=\alpha(a\cdot u)
\eeq
Очевидно, будут справедливы тождества
\beq\label{eq10.7}
\widetilde{a\cdot u}=\widetilde{u}\cdot a^{-1}, \qquad \widetilde{u\cdot
a}=a^{-1}\cdot \widetilde{u} 
\eeq
\beq\label{eq10.8}
\widetilde{a\cdot \alpha}=\widetilde{\alpha }\cdot a^{-1}, \qquad
\widetilde{\alpha \cdot a}=a^{-1}\cdot \widetilde{\alpha} 
\eeq
Если обозначить символом $\delta^a$ дельта-функционал
\beq\label{delta^a(u)=u(a)}
\delta^a(u)=u(a),
\eeq
то
\beq\label{eq10.9}
\widetilde{\delta^a}=\delta^{a^{-1}}, \qquad a\cdot \delta^b=\delta^{a\cdot b},
\qquad \delta^b\cdot a=\delta^{b\cdot a} 
\eeq

\paragraph{Свертка мер.}

Свертка функции с мерой определяется тождеством
\beq\label{eq10.10}
(\alpha * u) (t)=\alpha\l \widetilde{t\cdot u}\r, \qquad (u * \alpha) (t)=\alpha\l
\widetilde{u\cdot t}\r
\eeq
Последовательно доказываются тождества
\beq
  \delta^a * u=u\cdot a^{-1}, \qquad u * \delta^a=a^{-1}\cdot u
\label{eq10.11}
\eeq
\beq
\widetilde{\alpha * u}=\widetilde{u} * \widetilde{\alpha},
 \qquad
\widetilde{u * \alpha}=\widetilde{\alpha} * \widetilde{u} \label{eq10.12}
\eeq
\beq
\alpha (\beta * u)=\beta(\alpha * \widetilde{u}),
 \qquad
\alpha (u * \beta)=\beta(\widetilde{u} * \alpha) \label{eq10.13}
\eeq
Из них следует цепочка тождеств
\beq
\alpha \l u * \widetilde{\beta} \r= \alpha \l \widetilde{\beta * \widetilde{u}}
\r= \widetilde{\alpha} \l \beta * \widetilde{u} \r= \beta \l \widetilde{\alpha}
* u \r \label{eq10.14}
\eeq
которая и объявляется {\it сверткой мер}:\index{свертка!мер}
\beq\label{DEF:svertka}
(\alpha * \beta) (u)=\alpha \l u * \widetilde{\beta} \r= \alpha \l
\widetilde{\beta * \widetilde{u}} \r= \widetilde{\alpha} \l \beta *
\widetilde{u} \r= \beta \l \widetilde{\alpha}
* u \r
\eeq
Прямым вычислением проверяется, что определенная таким образом свертка мер удовлетворяет тождествам
\beq\label{DEF:svertka-s-delta^a}
\delta^a * \beta=a\cdot \beta, \qquad \beta * \delta^a=\beta \cdot a,
\eeq
\beq\label{delta^a*delta^b=delta^(a-cdot-b)}
\delta^a * \delta^b=\delta^{a\cdot b}
\eeq
\beq\label{eq10.18}
\widetilde{\alpha * \beta}=\widetilde{\beta} * \widetilde{\alpha}
\eeq

\btm
Свертка мер, определенная формулой \eqref{DEF:svertka} совпадает со сверткой, определенной формулой
\beq\label{sevrtka-kak-dvoinoi-integral}
\alpha * \beta (u)=(\alpha\boxasterisk \beta)(w)\Big|_{w(s,t)=u(s\cdot t)}=
\int_G \l \int_G u(s\cdot t) \alpha(\d s) \r \beta(\d t)= \int_G \l \int_G
u(s\cdot t) \beta(\d t)\r \alpha(\d s)
\eeq
\etm

\paragraph{Мера Хаара.}
Пусть $\mu$ --- левоинвариантная мера Хаара на локально компактной группе $G$, то есть функционал на ${\mathscr C}_{00}(G)$ со свойством
\beq\label{t-cdot-mu=mu-1}
\mu(u\cdot t)=\mu(u), \qquad u\in{\mathcal C}_{00}(G),\ t\in G.
\eeq
или, в эквивалентной записи,
\beq\label{t-cdot-mu=mu}
t\cdot \mu =\mu,\qquad t\in G.
\eeq
Напомним \cite[15.11]{Hewitt-Ross-2}, что модулярная функция $\Delta$ на $G$ определяется равенством:
\beq\label{DEF:modulyarnaya-funktsiya}
\Delta(t)=\frac{\mu(t^{-1}\cdot u)}{\mu(u)},\qquad s\in G,
\eeq
в котором $u$ --- произвольная ненулевая неотрицательная финитная непрерывная функция на $G$ (значение $\Delta(t)$ не зависит от выбора $u$). Модулярная функция $\Delta$ не зависит от выбора левоинвариантной меры Хаара $\mu$ и является (непрерывным) характером на $G$:
\beq\label{svoista-modulyarnoi-funktsii}
\Delta(1_G)=1,\quad \Delta(s\cdot t)=\Delta(s)\cdot\Delta(t),\quad \Delta(s^{-1})=\Delta(s)^{-1},\qquad s,t\in G.
\eeq
Свойства меры Хаара $\mu$ удобно собрать в одном месте:
\beq\label{svoistva-mery-Haara-1}
\mu\big(u\cdot t\big)= \mu(u),\qquad \mu\big(t\cdot u\big)= \frac{\mu(u)}{\varDelta(t)},\qquad  \mu\big(\widetilde{u}\big)= \mu \l \frac{u}{\varDelta}\r,\qquad u\in {\mathcal C}_{00}(G), \ t\in G.
\eeq
и в эквивалентной записи это выглядит так:
\beq\label{svoistva-mery-Haara-2}
t\cdot \mu=\mu,\qquad \mu\cdot t= \frac{\mu}{\varDelta(t)},\qquad \widetilde{\mu}=  \frac{1}{\varDelta}\cdot\mu,\qquad t\in G.
\eeq
(в последней формуле произведение функции и меры понимается в смысле \eqref{(f-cdot-alpha)(u)}).

Тождества \eqref{svoistva-mery-Haara-2} удобно записывать в виде свойств дифференциала:
\beq\label{svoistva-mery-Haara-3-diff}
\mu(\d(a\cdot t))=\mu(\d t),\qquad
\mu(\d (t\cdot a))= \frac{\mu(\d t)}{\varDelta(a)},\qquad
\mu(\d(t^{-1}))= \frac{\mu(\d t)}{\varDelta(t)},\qquad a\in G.
\eeq

Символ $\mu$ часто удобно опускать в самом интеграле
\beq\label{int_G-u(t)-dt=int_G-u(t)-mu(dt)=mu(u)}
\int_G u(t)\ \d t:=\int_G u(t)\ \mu(\d t)=\mu(u),\qquad u\in{\mathcal C}_{00}(G),
\eeq
и в тождествах \eqref{svoistva-mery-Haara-3-diff}:
\beq\label{svoistva-mery-Haara-3}
\d(a\cdot t)= \d t,\qquad
\d (t\cdot a)= \frac{\d t}{\varDelta(a)},\qquad
\d(t^{-1})= \frac{\d t}{\varDelta(t)},\qquad a\in G.
\eeq
\brem
Для унимодулярных групп, то есть таких, у которых модулярная функция тождественна единице, в частности, для компактных групп, тождества \eqref{svoistva-mery-Haara-2} и \eqref{svoistva-mery-Haara-3} принимают вид
\beq\label{svoistva-mery-Haara-2-unimod}
t\cdot \mu=\mu,\qquad \mu\cdot t= \mu,\qquad \widetilde{\mu}=  \mu,\qquad t\in G.
\eeq
\beq\label{svoistva-mery-Haara-4}
\d(a\cdot t)= \d t,\qquad
\d (t\cdot a)= \d t,\qquad
\d(t^{-1})= \d t,\qquad a\in G.
\eeq
\erem

\bpr
Первое из равенств \eqref{svoistva-mery-Haara-3} доказывается цепочкой
\begin{multline*}
\int u(t)\ \d(a\cdot t)=\begin{vmatrix}a\cdot t=s\\ t=a^{-1}\cdot s\end{vmatrix}=
\int u(a^{-1}\cdot s)\ \d s=\int (u\cdot a^{-1})(s)\ \d s=\\=\mu(u\cdot a^{-1})=\eqref{svoistva-mery-Haara-2}=\mu(u)=\int u(t)\ \d t.
\end{multline*}
Второе --- цепочкой
\begin{multline*}
\int u(t)\ \d(t\cdot a)=\begin{vmatrix}t\cdot a=s\\ t=s\cdot a^{-1}\end{vmatrix}=
\int u(s\cdot a^{-1})\ \d s=\int (a^{-1}\cdot u)(s)\ \d s=\mu(a^{-1}\cdot u)=\\=\eqref{svoistva-mery-Haara-2}=\Delta(a)\cdot\mu(u)=\mu(u\cdot\Delta(a))=\int u(t)\cdot \Delta(a)\cdot\d t.
\end{multline*}
А третье --- цепочкой
$$
\int u(t)\ \d(t^{-1})=\begin{vmatrix}t^{-1}=s\\ t=s^{-1}\end{vmatrix}=
\int u(s^{-1})\ \d s=\int \widetilde{u}(s)\ \d s=\mu(\widetilde{u}) =\eqref{svoistva-mery-Haara-2}=\mu\left(\frac{u}{\Delta}\right)
=\int \frac{u(t)}{\Delta(t)}\ \d t.
$$
\epr

\paragraph{Автоморфизмы компактной группы.}

Нам понадобится следующее утверждение, относящееся к математическому фольклору:

\btm\label{TH:mu(ph(X))=mu(X)}
Всякий (непрерывный) автоморфизм $\ph:K\to K$ компактной группы $K$ сохраняет меру Хаара на $K$:
\beq\label{mu(ph(X))=mu(X)}
\mu(\ph(X))=\mu(X),\qquad X\subseteq K.
\eeq
\etm
\bpr
Поскольку $\ph$ -- гомеоморфизм, он сохраняет борелевские множества, и поэтому мы можем рассмотреть меру
$$
m(X)=\mu(\ph(X))
$$
на борелевских множествах $X\subseteq K$. Для всякой точки $t\in K$ мы получаем
$$
m(t\cdot X)=\mu(\ph(t\cdot X))=\mu(\ph(t)\cdot\ph(X))=\mu(\ph(X))=m(X).
$$
То есть, $m$ инвариантна на $K$, как и $\mu$. Значит, они отличаются на константу:
$$
m(X)=C\cdot\mu(X).
$$
С другой стороны, $\ph(K)=K$, поэтому
$$
m(K)=\mu(\ph(K))=\mu(K),
$$
и значит, $C=1$.
\epr

\paragraph{Свертка функций.}
Свертка функций $f,g\in{\mathscr C}_{00}(G)$ определяется двумя эквивалентными формулами:
\beq\label{DEF:svertka-funktsij}
(f*g)(t)=\int f(s)\cdot g(s^{-1}\cdot t)\cdot\mu (\d s)=\int \frac{f(t\cdot s^{-1})\cdot g(s)}{\varDelta(s)}\cdot \mu(\d s)
\eeq

\bpr[Доказательство эквивалентности:]
\begin{multline*}
(f*g)(t)=\int f(s)\cdot g(s^{-1}\cdot t)\cdot\mu (\d s)=\begin{vmatrix} s^{-1}\cdot t=\tau\\
s=t\cdot \tau^{-1} \\
\d s=\d \l t\cdot \tau^{-1} \r=\frac{\d\tau}{\varDelta(\tau)}
\end{vmatrix}=\\=
\int f(t\cdot\tau^{-1})\cdot g(\tau)\cdot \frac{\mu(\d\tau)}{\varDelta(\tau)}=
\int \frac{f(t\cdot s^{-1})\cdot g(s)}{\varDelta(s)}\cdot \mu(\d s)
\end{multline*}
\epr

\medskip
\centerline{\bf Свойства свертки функций:}

\bit{\it

\item[$1^\circ$.] Антипод свертки функций:
\beq\label{widetilde(f*g)}
\widetilde{f*g}=\widetilde{g}*\widetilde{f},
\eeq

\item[$2^\circ$.] Связь свертки функций со сверткой функции и меры:
\beq\label{f*g=(f-cdot-mu)*g}
f*g=(f\cdot\mu)*g
\eeq

\item[$3^\circ$.] Связь свертки функций со сверткой мер:
\beq\label{(f-cdot-mu)*(g-cdot-mu)=(f*g)-cdot-mu}
(f\cdot\mu) * (g\cdot\mu)=(f*g)\cdot\mu
\eeq
}\eit

\bpr
1. Тождество \eqref{widetilde(f*g)} доказывается цепочкой
\begin{multline*}
\widetilde{f*g}(s)=\eqref{widetilde(u)(t)=u(t^(-1))}=f*g(s^{-1})=\eqref{DEF:svertka-funktsij}=
\int f(t)\cdot g(t^{-1}\cdot s^{-1})\cdot \d t=\int f(t)\cdot g((s\cdot t)^{-1})\cdot \d t=\\=
\int \widetilde{f}(t^{-1})\cdot \widetilde{g}(s\cdot t)\cdot \d t=
\begin{vmatrix}s\cdot t=r, \ t=s^{-1}\cdot r, \ t^{-1}=r^{-1}\cdot s
\\ \d t=\d(s^{-1}\cdot r)=\eqref{svoistva-mery-Haara-3}=\d r\end{vmatrix}=
\int \widetilde{f}(r^{-1}\cdot s)\cdot \widetilde{g}(r)\cdot \d r =\eqref{DEF:svertka-funktsij}=
\big(\widetilde{g}*\widetilde{f}\big)(s).
\end{multline*}

2. Тождество \eqref{f*g=(f-cdot-mu)*g} доказывается цепочкой
\begin{multline*}
(f*g)(t)=\int f(s)\cdot g(s^{-1}\cdot t)\cdot\mu (\d s)=\int f(s)\cdot (t\cdot g)(s^{-1})\cdot\mu (\d s)=\\= \int f(s)\cdot (\widetilde{t\cdot g})(s)\cdot\mu (\d s)=(f\cdot\mu)(\widetilde{t\cdot g})=\eqref{eq10.10}= ((f\cdot\mu)*g)(t)
\end{multline*}

3. Для доказательства \eqref{(f-cdot-mu)*(g-cdot-mu)=(f*g)-cdot-mu} зафиксируем $u\in{\mathscr C}(G)$. Тогда, с одной стороны,
\begin{multline*}
\big((f\cdot\mu) * (g\cdot\mu)\big)(u)=(f\cdot\mu)(u*(\widetilde{g\cdot\mu}))=\int
(u*(\widetilde{g\cdot\mu}))(s)\cdot f(s)\cdot \mu(\d s)=\eqref{eq10.10}=\\=
\int(\widetilde{g\cdot\mu})(\widetilde{u\cdot s})\cdot f(s)\cdot \mu(\d s)=
\int(g\cdot\mu)(u\cdot s)\cdot f(s)\cdot \mu(\d s)=\\=
\int\int (u\cdot s)(t)\cdot g(t)\cdot\mu(\d t)\cdot f(s)\cdot \mu(\d s)=
\int\int u(s\cdot t)\cdot f(s)\cdot g(t)\cdot \mu(\d s)\cdot \mu(\d t)
\end{multline*}
А с другой,
\begin{multline*}
\big((f * g)\cdot\mu)\big)(u)=\int u(t)\cdot (f*g)(t)\cdot \mu(\d t)=
\int\int  u(t)\cdot f(s)\cdot g(s^{-1}\cdot t)\cdot \mu(\d s)\cdot \mu(\d t)=\\=
\int f(s)\cdot \l \int  u(t)\cdot g(s^{-1}\cdot t)\cdot \mu(\d t)\r \cdot \mu(\d s)=\\=
\begin{vmatrix} s^{-1}\cdot t=\tau\\
t=s\cdot \tau \\
\d t=\d \l s\cdot\tau \r=\d\tau
\end{vmatrix}=
\int f(s)\cdot \l\int  u(s\cdot\tau)\cdot g(\tau)\cdot \mu(\d \tau)\r\cdot \mu(\d s)=
\int\int  u(s\cdot\tau)\cdot f(s)\cdot g(\tau)\cdot \mu(\d s)\cdot \mu(\d \tau)
\end{multline*}
\epr

\btm
Сопряженным оператором к оператору свертки с функцией $f$ слева
$$
M(u)=f*u
$$
является оператор свертки с мерой $\widetilde{f\cdot\mu}$ слева:
$$
M^\star(\beta)=\widetilde{f\cdot\mu}*\beta
$$
\etm

\bpr
$$
M^\star(\beta)(u)=\beta(M(u))=\beta(f*u)=\eqref{f*g=(f-cdot-mu)*g} =\beta((f\cdot\mu)*u)=\eqref{DEF:svertka}=((\widetilde{f\cdot\mu})*\beta)(u)
$$
\epr

\paragraph{Инволюции в $\mathcal{C}(G)$, $\mathcal{C}^\star(G)$ и $L_1(G)$.} 
На алгебрах $\mathcal{C}(G)$ и $\mathcal{C}^\star (G)$ удобно ввести два вида инволюции:
\bit{

\item[---] {\it плоская инволюция} (или {\it поточечная инволюция}) на $\mathcal{C}(G)$ определяется формулой
\beq\label{ploskaya-involutsija-v-C(G)}
\overline{u}(t)=\overline{u(t)},\qquad u\in\mathcal{C}(G),\ t\in G.
 \eeq

\item[---] {\it объемная инволюция} (или {\it инволюция с подкруткой}) на $\mathcal{C}(G)$ определяется формулой
\beq\label{involutsija-v-C(G)}
u^\bullet(t)=\frac{\overline{u(t^{-1})}}{\Delta(t)},\qquad t\in G.
 \eeq
 или, что эквивалентно, формулой
\beq\label{involutsija-v-C(G)-1}
u^\bullet=\frac{\overline{\widetilde{u}}}{\Delta},\qquad u\in\mathcal{C}(G).
 \eeq
(где $\Delta$ --- модулярная функция на $G$, определенная формулой \eqref{DEF:modulyarnaya-funktsiya}),

\item[---] {\it плоская инволюция} на $\mathcal{C}^\star (G)$ определяется формулой
 \beq\label{ploskaya-involutsija-v-C^star(G)}
\overline{\alpha}(u)=\overline{\alpha(\overline{u})},\qquad \alpha\in\mathcal{C}^\star (G).
 \eeq
\item[---] {\it объемная инволюция} (или {\it инволюция с подкруткой}) на $\mathcal{C}^\star (G)$ определяется формулой
 \beq\label{involutsija-v-C^star(G)}
\alpha^\bullet(u)=\overline{\alpha(\overline{\widetilde{u}})},\qquad \alpha\in\mathcal{C}^\star (G).
 \eeq
 
 }\eit

\noindent\rule{160mm}{0.1pt}\begin{multicols}{2}

\bex
Плоская и объемная инволюции дельта-функции описываются формулами
$$
\overline{\delta^a}=\delta^a,\qquad (\delta^a)^\bullet=\delta^{a^{-1}}.
$$
\eex

\bex
Мера Хаара на компактной группе $G$ вещественна: 
\beq\label{mu(overline(u))-G-komp}
\overline{\mu}=\mu.
\eeq
Это эквивалентно условию что $\mu$ перестановочна с плоской инволюцией функции:
\beq\label{mu(overline(u))-G-komp-1}
\mu(\overline{u})=\overline{\mu(u)},\qquad u\in \mathcal{C}(G)
\eeq
Можно заметить, что это справедливо также и для левоинвариантной меры Хаара на произвольной локально компактной группе $G$ (определенной выше формулой \eqref{t-cdot-mu=mu-1}):
\beq\label{mu(overline(u))}
\mu(\overline{u})=\overline{\mu(u)},\qquad u\in \mathcal{C}_{00}(G)
\eeq

\eex

\bex
На банаховой алгебре $L_1(G)$ инволюция обычно определяется как объемная, то есть заданная формулой \eqref{involutsija-v-C(G)-1}:
\beq\label{involutsija-na-L_1(G)}
f^\bullet=\frac{\overline{\widetilde{f}}}{\Delta}, \qquad f\in L_1(G).
\eeq
\eex

\end{multicols}\noindent\rule[10pt]{160mm}{0.1pt}

\btm\label{TH:(f-cdot-mu)^bullet=f^bullet-cdot-mu}
Для всякой левоинвариантной меры Хаара $\mu$ на $G$ и любой финитной функции $f\in{\mathscr C}_{00}(G)$ справедливо 
\bit{
\item[---] равенство, связывающее плоскую инволюцию на $\mathcal{C}(G)$ с плоской инволюцией на $\mathcal{C}^\star (G)$:
\beq\label{overline(f-cdot-mu)=overline(f)-cdot-mu}
\overline{f\cdot\mu}=\overline{f}\cdot\mu,
\eeq
и

\item[---] равенство, связывающее объемную инволюцию на $\mathcal{C}(G)$ с объемной инволюцией на $\mathcal{C}^\star (G)$:
\beq\label{(f-cdot-mu)^bullet=f^bullet-cdot-mu}
(f\cdot\mu)^\bullet=f^\bullet\cdot\mu
\eeq

}\eit
\etm 
\bpr Для всякой функции $u\in\mathcal{C}(G)$ мы получим, во-первых,
\begin{multline*}
\overline{f\cdot\mu}(u)=\eqref{ploskaya-involutsija-v-C^star(G)}=
\overline{(f\cdot\mu)(\overline{u})}=\eqref{(f-cdot-alpha)(u)}=
\overline{\mu(f\cdot\overline{u})}=\eqref{mu(overline(u))}=\\=
\mu\Big(\overline{f\cdot\overline{u}}\Big)=
\mu\Big(\overline{f}\cdot u\Big)=\eqref{(f-cdot-alpha)(u)}=(\overline{f}\cdot\mu)(u),
\end{multline*}
и, во-вторых,
\begin{multline*}
(f\cdot\mu)^\bullet(u)=\eqref{involutsija-v-C^star(G)}=
\overline{(f\cdot\mu)(\overline{\widetilde{u}})}=\eqref{(f-cdot-alpha)(u)}=
\overline{\mu(f\cdot\overline{\widetilde{u}})}=\eqref{mu(overline(u))}=
\mu\Big(\overline{f\cdot\overline{\widetilde{u}}}\Big)=
\mu\Big(\overline{f}\cdot\widetilde{u}\Big)=\\=
\mu\left(\Big(\widetilde{\overline{f}}\cdot u\Big)\widetilde{\phantom{\cdot}} \right)=\eqref{svoistva-mery-Haara-1}=
\mu\left(\frac{\widetilde{\overline{f}}\cdot u}{\Delta}\right)=\eqref{(f-cdot-alpha)(u)}=
\left(\frac{\widetilde{\overline{f}}}{\Delta}\cdot\mu\right)(u)=\eqref{involutsija-v-C(G)-1}=
\left(f^\bullet\cdot\mu\right)(u).
\end{multline*}

\epr

\paragraph{Действие $\mathcal{C}^\star (G)$ на $\mathcal{C}(G)$.}
Если $X$   -- стереотипный левый (правый) модуль над
стереотипной алгеброй  $A$, то сопряженное пространство $X^\star$ есть
стереотипный правый (левый) модуль над $A$ относительно операции
\begin{equation}
(f\cdot a)(x)=f(a\cdot x) \qquad ((a\cdot f)(x)=f(x\cdot a)) \label{eq11.1}
\end{equation}
Если в качестве $X$ взята алгебра $\mathcal{C}^\star (G)$, то умножение
функции $u\in \mathcal{C}(G)$ слева и справа на меру $\alpha\in
\mathcal{C}^\star (G)$ определяется формулами
$$
\beta(\alpha\cdot u)=(\beta *\alpha)(u), \qquad \beta(u\cdot \alpha)=(\alpha
*\beta)(u), \qquad \beta\in \mathcal{C}^\star (G)
$$
В силу \eqref{DEF:svertka}, это означает, что
$$
\beta(\alpha\cdot u)=(\beta
*\alpha)(u)=\beta\l\widetilde{\alpha*\widetilde{u}}\r, \qquad \beta(u\cdot
\alpha)=(\alpha *\beta)(u)=\beta\l \widetilde{\alpha}* u\r,
$$
То есть, в силу \eqref{eq10.12}, справедливы явные формулы:
\begin{equation}
\alpha\cdot u=\widetilde{\alpha*\widetilde{u}}=u*\widetilde{\alpha}, \qquad
u\cdot \alpha=\big(\widetilde{u} * \alpha\big)\widetilde{\phantom{.}}=\widetilde{\alpha}* u
\label{eq11.2}
\end{equation}
В частности,
\begin{equation} \label{eq11.3}
\delta_a\cdot u=u*\delta_{a^{-1}}=a\cdot u, \qquad u\cdot
\delta_a=\delta_{a^{-1}}*u=u\cdot a \qquad (a\in G),
\end{equation}
где $a\cdot u$ и $u\cdot a$ определены равенствами \eqref{DEF:u-cdot-a}.

\subsection{Групповые алгебры как аугментированные стереотипные алгебры}

Групповая алгебра мер ${\mathcal C}^\star(G)$ превращается в объект категории $\AugInvSteAlg$ аугментированных инволютивных стереотипных алгебр (см. определение на с.\pageref{DEF:augmented-ster-alghebras}), будучи наделена объемной инволюцией \eqref{involutsija-v-C^star(G)} 
 \beq\label{involution-in-C*(G)}
\alpha^\bullet(u)=\overline{\alpha(\overline{\widetilde{u}})},\qquad \alpha\in {\mathcal C}^\star(G), \ u\in {\mathcal C}(G)
 \eeq
и аугментацией
 \beq\label{augmentation-in-C*(G)}
 \e(\alpha)=\alpha(1),\qquad \alpha\in {\mathcal C}^\star(G)
 \eeq
(здесь 1 -- тождественная единица в ${\mathcal C}(G)$).

\btm\label{TH:coker-groups=>coker-group-algebras}
Пусть в цепочке локально компактных групп
$$
 \xymatrix  
{
H\ar[r]^{\lambda} & G\ar[r]^{\varkappa} &F
}
$$
второй морфизм является коядром первого в категории локально компактных групп\footnote{Это можно понимать так, что $\varkappa$ является фактор-отображением группы $G$ по замыканию $\overline{\lambda(H)}$ подгруппы $\lambda(H)$.}:
$$
\coker\lambda=\varkappa.
$$
Тогда в соответствующей цепочке групповых алгебр
$$
 \xymatrix 
{
{\mathcal C}^\star(H)\ar[r]^{{\mathcal C}^\star(\lambda)} & {\mathcal C}^\star(G)\ar[r]^{{\mathcal C}^\star(\varkappa)} &
{\mathcal C}^\star(F)
}
$$
второй морфизм является коядром первого в категории $\AugInvSteAlg$ аугментированных инволютивных стереотипных алгебр:
\beq\label{coker-C^star(lambda)=C^star(varkappa)}
\coker{\mathcal C}^\star(\lambda)={\mathcal C}^\star(\varkappa).
\eeq
\etm
\bpr
Введем более простые обозначения:
$$
\lambda'={\mathcal C}^\star(\lambda),\quad \varkappa'={\mathcal C}^\star(\varkappa).
$$
Кроме того, пусть $\e=\e_{{\mathcal C}^\star(H)}$ --- аугментация на ${\mathcal C}^\star(H)$, а $\iota=\iota_{{\mathcal C}^\star(F)}$ --- вложение $\C$ в ${\mathcal C}^\star(F)$.

1. Сначала отметим следующее тождество:
\beq\label{varkappa'-circ-lambda'=iota-circ-e}
\varkappa'\circ\lambda'=\iota\circ\e.
\eeq
$$
 \xymatrix @R=2.pc @C=4.pc
{
H\ar[r]^{\lambda}\ar[d]_{\delta_H} & G\ar[r]^{\varkappa}\ar[d]_{\delta_G} &F\ar[d]_{\delta_F} \\
{\mathcal C}^\star(H)\ar[r]^{\lambda'}
\ar@{-->}@/_2ex/[rd]^{\e} & {\mathcal C}^\star(G)\ar[r]^{\varkappa'} &
{\mathcal C}^\star(F)\\
& \C \ar@{-->}@/_2ex/[ru]^{\iota} &
}
$$
Его удобно проверить на элементах $\delta^t$, $t\in H$:
$$
(\varkappa'\circ\lambda')(\delta^t)=\varkappa'(\lambda'(\delta^t))=\varkappa'(\delta^{\lambda(t)})=\delta^{\varkappa(\lambda(t))}=
\delta^{1_F}=1_{{\mathcal C}^\star(F)}=\iota(1)=\iota(\e(\delta^t))=(\iota\circ\e)(\delta^t).
$$
Поскольку элементы $\delta^t$, $t\in H$, полны в ${\mathcal C}^\star(H)$ \cite[Lemma 8.2]{Ak03}, это доказывает  \eqref{varkappa'-circ-lambda'=iota-circ-e}.

2. Покажем далее, что $\varkappa'$ тривиализует $\lambda'$ снаружи. Пусть $\eta,\theta:{\mathcal C}^\star(F)\rightrightarrows B$ --- два параллельных морфизма. Тогда
$$
\eta\circ\varkappa'\circ\lambda'=\eqref{varkappa'-circ-lambda'=iota-circ-e}=
\eta\circ\iota\circ\e=\theta\circ\iota\circ\e=\eqref{varkappa'-circ-lambda'=iota-circ-e}=
\theta\circ\varkappa'\circ\lambda'.
$$

3. Пусть теперь $\gamma:{\mathcal C}^\star(G)\to A$ --- какой-то морфизм аугментированных стереотипных алгебр, тривиализующий $\lambda'$ снаружи, то есть
\beq\label{eta-circ-gamma-circ-lambda'=theta-circ-gamma-circ-lambda'}
\eta\circ\gamma\circ\lambda'=\theta\circ\gamma\circ\lambda'
\eeq
для любых $\eta,\theta:A\rightrightarrows B$. Пусть $\eta=\id_A$ (тождественное отображение $A$ в себя) и $\theta=\iota_A\circ\e_A$ (где $\e_A:A\to\C$ -- аугментация в $A$, $\iota_A:\C\to A$ --- вложение $\C$ в $A$). Рассмотрим диаграмму:
$$
 \xymatrix @R=2.pc @C=4.pc
{
H\ar[r]^{\lambda}\ar[d]_{\delta_H} & G\ar[rr]^{\varkappa}\ar[d]_{\delta_G} & & F\ar[dl]_{\delta_F}\ar@{-->}@/^6ex/[ddll]^{\ph} \\
{\mathcal C}^\star(H)\ar[r]^{\lambda'}
 & {\mathcal C}^\star(G)\ar[r]^{\varkappa'}\ar[d]_{\gamma} &
{\mathcal C}^\star(F)\ar@{-->}[dl]_{\gamma'}&  \\
& A\ar@/_3ex/[d]_{\iota_A\circ\e_A}\ar@/^3ex/[d]^{\id_A} & & \\
& A & &
}
$$
Для всякого $t\in H$ мы получим:
\begin{multline*}
\gamma(\delta_G^{\lambda(t)})=
\gamma\big(\lambda'(\delta_H^t)\big)=\id_A\Big(\gamma\big(\lambda'(\delta_H^t)\big)\Big)=
\eqref{eta-circ-gamma-circ-lambda'=theta-circ-gamma-circ-lambda'}=
(\iota_A\circ\e_A)\Big(\gamma\big(\lambda'(\delta_H^t)\big)\Big)= (\iota_A\circ\e_A)\Big(\gamma\big(\delta_H^{\lambda(t)}\big)\Big)=\\=
\iota_A\bigg(\e_A\Big(\gamma\big(\delta_H^{\lambda(t)}\big)\Big)\bigg)=
\iota_A\bigg(\e_{{\mathcal C}^\star(G)}\big(\delta_H^{\lambda(t)}\big)\bigg)=
\iota_A(1)=1.
\end{multline*}
То есть гомоморфизм $t\in G\mapsto \gamma(\delta^t)$ (группы $G$ в алгебру $A$) тождественен на подгруппе $\lambda(H)\subseteq G$. Это значит, что он продолжается до некоторого гомоморфизма $\ph:F=\coker\lambda=G/\overline{\lambda(H)}\to A$. Этот гомоморфизм $\ph$ в свою очередь продолжается до гомоморфизма на групповую алгебру $\gamma':{\mathcal C}^\star(F)\to A$ \cite[Theorem 10.12]{Ak03}, который и будет продолжением гомоморфизма $\gamma$, поскольку на элементах $G$ они совпадают (здесь мы опять используем \cite[Lemma 8.2]{Ak03}).
\epr

По аналогии доказывается

\btm\label{TH:coker-groups=>coker-group-algebras-E}
Пусть в цепочке групп Ли
$$
 \xymatrix  
{
H\ar[r]^{\lambda} & G\ar[r]^{\varkappa} &F
}
$$
второй морфизм является коядром первого в категории локально компактных групп\footnote{Это можно понимать так, что $\varkappa$ является фактор-отображением группы $G$ по замыканию $\overline{\lambda(H)}$ подгруппы $\lambda(H)$.}:
$$
\coker\lambda=\varkappa.
$$
Тогда в соответствующей цепочке групповых алгебр распределений
$$
 \xymatrix 
{
{\mathcal E}^\star(H)\ar[r]^{{\mathcal E}^\star(\lambda)} & {\mathcal E}^\star(G)\ar[r]^{{\mathcal E}^\star(\varkappa)} &
{\mathcal E}^\star(F)
}
$$
второй морфизм является коядром первого в категории $\AugInvSteAlg$ аугментированных инволютивных стереотипных алгебр:
\beq\label{coker-E^star(lambda)=E^star(varkappa)}
\coker{\mathcal E}^\star(\lambda)={\mathcal E}^\star(\varkappa).
\eeq
\etm

\section{Представления групп и отражения}

\subsection{Представления локально компактной группы}

\paragraph{Представления в гильбертовом пространстве.}

\btm\label{TH:nepr-predstavleniya} Пусть локально компактная группа $G$ действует унитарными операторами в гильбертовом пространстве $X$, то есть задано отображение $\pi:G\to{\mathcal L}(X)$, удовлетворяющее условиям
\beq\label{DEF:representation-of-G-in-X }
\pi(1_G)=1_X,\qquad \pi(s\cdot t)=\pi(s)\cdot\pi(t),\qquad \pi(s^{-1})=\pi(s)^\bullet,\qquad s,t\in G.
\eeq
Тогда следующие условия эквивалентны:
 \bit{
\item[(i)] отображение $\pi:G\to{\mathcal L}(X)$ непрерывно относительно сильной операторной топологии на ${\mathcal L}(X)$,

\item[(ii)] отображение $(s,x)\in G\times X\mapsto \pi(s)x\in X$ непрерывно,

\item[(iii)] отображение $\pi:G\to{\mathcal L}(X)$ непрерывно относительно стереотипной топологии пространства ${\mathcal L}(X)$.

\item[(iv)] существует (необходимо, единственный) инволютивный непрерывный гомоморфизм стереотипных алгебр $\dot{\pi}:\mathcal{C}^\star(G)\to {\mathcal L}(X)$, замыкающий диаграмму
\begin{gather}\label{C*(G)=grup-alg}
\begin{diagram}
\node{G} \arrow[2]{e,t}{\delta} \arrow{se,b}{\pi}
\node[2]{\mathcal{C}^\star(G)} \arrow{sw,b}{\dot{\pi}}
\\
\node[2]{{\mathcal L}(X)}
\end{diagram},
\end{gather}
 }\eit
\etm

\bit{
\item[$\bullet$] В дальнейшем под {\it унитарным представлением локально компактной группы $G$} (в гильбертовом пространстве $X$) мы понимаем отображение $\pi:G\to{\mathcal L}(X)$, удовлетворяющее \eqref{DEF:representation-of-G-in-X } и условиям (i)-(iv) этой теоремы.
}\eit

\bpr
Эквивалентность (i)$\Leftrightarrow$(ii) -- стандартный факт (см. напр. \cite[глава 5, \S 1]{Barut-Raczka}), а
(iii)$\Leftrightarrow$(iv) доказывается в \cite[Theorem 10.12]{Ak03}. Докажем (ii)$\Rightarrow$(iii). Пусть отображение $(s,x)\in G\times X\mapsto \pi(s)x\in X$ непрерывно. Тогда для всякого компакта $K\subseteq X$
$$
\pi(g)x\underset{g\to g_0}{\longrightarrow}\pi(g_0)x
$$
равномерно по $x\in K$. Это означает, что отображение $\pi:G\to X:X$ непрерывно (здесь $X:X$ обозначает пространство линейных непрерывных отображений из $X$ в $X$, наделенное топологией равномерной сходимости на компактах). Для всякого компакта $T\subseteq G$ его образ $\pi(T)$ является компактом в $X:X$, поэтому при псевдонасыщении пространства $X:X$ (то есть при переходе от $X:X$ к ${\mathcal L}(X)$) топология на $\pi(T)$ не меняется. Отсюда следует, что ограничение $\pi|_T:T\to {\mathcal L}(X)$ непрерывно. Поскольку это верно для любого компакта $T\subseteq G$, отображение $\pi:G\to {\mathcal L}(X)$ тоже непрерывно. Мы доказали (ii)$\Rightarrow$(iii). Импликация $(ii)\Leftarrow(iii)$ доказывается обратными рассуждениями.
\epr

\paragraph{Разложение по характерам нормальной компактной подгруппы.}

Напомним некоторые определения из теории представлений \cite{Kowalski}. Пусть $\{X_i;\ i\in I\}$ -- семейство гильбертовых пространств. Рассмотрим множество $X=\dot{\bigoplus}_{i\in I}X_i$, состоящее из семейств $x=\{x_i;\ x_i\in X_i\}$, квадратично суммируемых по норме:
$$
x=\{x_i;\ i\in I\}\in \dot{\bigoplus_{i\in I}}X_i\qquad\Longleftrightarrow\qquad  \forall i\in I\quad x_i\in X_i\quad
\&
\quad \sum_{i\in I} \norm{x_i}_{X_i}^2<\infty.
$$
Оно образует гильбертово пространство относительно послойных алгебраических операций и скалярного произведения
$$
\langle x,y \rangle=\sum_{i\in I}\langle x_i,y_i\rangle_i
$$
(где $\langle x_i,y_i\rangle_i$ -- скалярное произведение в $X_i$). Это пространство называется {\it гильбертовой прямой суммой пространств $X_i$}.

Если теперь $G$ -- локально компактная группа и $\{\pi_i:G\to {\mathcal L}(X_i)\}$ -- семейство ее унитарных представлений в пространствах $X_i$, то отображение
$$
\pi:G\to {\mathcal L}\left(\dot{\bigoplus_{i\in I}}X_i\right)\qquad\Big|\qquad \pi(g)\{x_i;\ x_i\in X_i\}=\{\pi_i(g)x_i;\ x_i\in X_i\}
$$
является унитарным представлением группы $G$ и называется {\it гильбертовой прямой суммой представлений $\pi_i$}.

Линейное непрерывное отображение гильбертовых пространств $\alpha:X\to Y$ называется
 \bit{
\item[---] {\it морфизмом представления $\pi:G\to{\mathcal L}(X)$ в представление $\rho:G\to{\mathcal L}(Y)$}, если
$$
\rho(t)\alpha(x)=\alpha(\pi(t)x),\qquad t\in G,\ x\in X,
$$

\item[---] {\it подпредставлением} представления $\rho:G\to{\mathcal L}(Y)$, если $\alpha$ является морфизмом представлений, оно инъективно, и $\alpha(X)$ замкнуто в $Y$ и инвариантно относительно действия $\rho$:
    $$
    \rho(t)\alpha(X)\subseteq\alpha(X),\qquad t\in G,
    $$

\item[---] {\it изоморфизмом представлений} $\pi:G\to{\mathcal L}(X)$ и $\rho:G\to{\mathcal L}(Y)$, если $\alpha$ является морфизмом этих представлений и изоморфизмом гильбертовых пространств $X$ и $Y$.

}\eit

Представление $\pi:G\to{\mathcal L}(X)$ в гильбертовом пространстве $X$ называется
 \bit{
\item[---] {\it неприводимым}, если у него нет нетривиальных (то есть отличных от $X$ и 0) подпредставлений,

\item[---] {\it полупростым}, если оно изоморфно некоторой гильбертовой прямой сумме неприводимых представлений:
 \beq\label{DEF:semisimple}
\pi=\dot{\bigoplus_{i\in I}}\sigma_i,
 \eeq

\item[---] {\it изотипическим}, если оно полупростое, причем в разложении \eqref{DEF:semisimple} все представления $\sigma_i$ изоморфны.

 }\eit

Для всякой локально компактной группы $G$ мы обозначаем символом $\widehat{G}$ ее {\it двойственный объект}, то есть множество неприводимых унитарных представлений $G$ такое, что
 \bit{
\item[---] любые два представления $\rho,\sigma\in\widehat{G}$ изоморфны тогда и только тогда, когда они совпадают:
$$
\rho\cong\sigma\qquad\Longleftrightarrow\qquad \rho=\sigma,
$$

\item[---] любое неприводимое унитарное представление $\pi:G\to{\mathcal L}(X)$ изоморфно некоторому $\sigma\in\widehat{G}$.
}\eit
Известно, что двойственный объект $\widehat{G}$ существует для всех локально компактных групп $G$.

Пусть $K$ -- нормальная компактная подгруппа в локально компактной группе $G$, $\mu_K$ -- нормированная мера Хаара на $K$. Для всякого представления $\sigma\in\widehat{K}$, $\sigma:K\to {\mathcal B}(X_\sigma)$  определим меру $\nu_\sigma\in{\mathcal C}^\star(G)$ формулой
\beq\label{DEF:nu_sigma-1}
\nu_\sigma(u)=\dim X_\sigma\cdot \int_K \overline{\tr\sigma(s)}\cdot u(s)\cdot \mu_K(\d s)=\dim X_\sigma\cdot \int_K \tr\sigma(s^{-1})\cdot u(s)\cdot \mu_K(\d s),\qquad u\in{\mathcal C}(G).
\eeq
или, что эквивалентно, формулой
\beq\label{DEF:nu_sigma}
\nu_\sigma=\dim X_\sigma\cdot \int_K \overline{\tr\sigma(s)}\cdot\delta_s\cdot \mu_K(\d s)=\dim X_\sigma\cdot \int_K \tr\sigma(s^{-1})\cdot\delta_s\cdot \mu_K(\d s)
\eeq

 \vglue10pt
\centerline{\bf Свойства мер $\nu_\sigma$:}
 \vglue10pt

\bit{\it

\item[$1^\circ$.] Для всякого $\sigma\in\widehat{K}$ мера $\nu_\sigma$ центральная в ${\mathcal C}^\star(G)$:
 \beq\label{chi_pi-centr}
 \nu_\sigma*\alpha=\alpha*\nu_\sigma,\qquad \alpha\in {\mathcal C}^\star(G),
 \eeq

\item[$2^\circ$.] Меры $\nu_\sigma$ образуют систему ортогональных проекторов в алгебре ${\mathcal C}^\star(G)$:
 \beq\label{chi_pi-ortog-proj}
 \nu_\sigma*\nu_\sigma=\nu_\sigma,\qquad \nu_\sigma*\nu_\rho=0,\qquad \sigma\ne\rho\in \widehat{K}.
 \eeq

}\eit

\bpr
1. Поскольку дельта-функционалы полны в ${\mathcal C}^\star(G)$, тождество \eqref{chi_pi-centr} эквивалентно тождеству
$$
\delta^t*\nu_\sigma*\delta^{t^{-1}}=\nu_\sigma,\qquad t\in G.
$$
или, в силу \eqref{DEF:svertka-s-delta^a}, тождеству
$$
t\cdot\nu_\sigma\cdot t^{-1}=\nu_\sigma,\qquad t\in G.
$$
Для его  доказательства используется теорема \ref{TH:mu(ph(X))=mu(X)}: поскольку отображение $x\mapsto t\cdot x\cdot t^{-1}$ является автоморфизмом группы $K$, оно сохраняет меру Хаара $\mu_K$ на $K$. Поэтому в следующей цепочке справедлива замена переменной:
\begin{multline*}
t\cdot\nu_\sigma\cdot t^{-1}(u)=\nu_\sigma(t^{-1}\cdot u\cdot t)=\dim X_\sigma\cdot \int_K\tr\sigma(s^{-1})
(t^{-1}\cdot u\cdot t)(s)\cdot \mu_K(\d s)=\\=
\dim X_\sigma\cdot \int_K\tr\sigma(s^{-1})
u(t\cdot s\cdot t^{-1})\cdot \mu_K(\d s)=\begin{vmatrix}t\cdot s\cdot t^{-1}=r \\ s=t^{-1}\cdot r\cdot t\\ \mu_K(\d s)=\mu_K(\d r)\end{vmatrix}=
\dim X_\sigma\cdot \int_K\tr\sigma(t^{-1}\cdot r^{-1}\cdot t) u(r)\cdot \mu_K(\d r)=\\=
\dim X_\sigma\cdot \int_K\tr\sigma(r^{-1}) u(r)\cdot \mu_K(\d r)=\nu_\sigma(u).
\end{multline*}

2. Обозначим $\chi_\sigma$ характеры представлений $\sigma\in\widehat{K}$:
$$
\chi_\sigma(s)=\tr\sigma(s).
$$
Как известно \cite[(27.24)]{Hewitt-Ross-2},
\beq\label{e_sigma*e_tau}
\chi_\sigma*\chi_\tau=\begin{cases}0, &\sigma\ne\tau\\ \frac{1}{\dim X_\sigma}\cdot\chi_\sigma,& \sigma=\tau\end{cases}.
\eeq
С другой стороны, понятно,
\beq\label{chi_sigma=dim X_sigma-widetilde(e_sigma)-mu_K}
\nu_\sigma=\dim X_\sigma\cdot\widetilde{\chi_\sigma}\cdot\mu_K.
\eeq
Поэтому
\begin{multline*}
\nu_\sigma*\nu_\tau=(\dim X_\sigma\cdot\widetilde{\chi_\sigma}\cdot\mu_K)*(\dim X_\tau\cdot\widetilde{\chi_\tau}\cdot\mu_K)=\eqref{(f-cdot-mu)*(g-cdot-mu)=(f*g)-cdot-mu}=
\dim X_\sigma\cdot\dim X_\tau\cdot(\widetilde{\chi_\sigma}*\widetilde{\chi_\tau})\cdot\mu_K=\eqref{widetilde(f*g)}=\\=
\dim X_\sigma\cdot\dim X_\tau\cdot(\widetilde{\chi_\tau*\chi_\sigma})\cdot\mu_K=\eqref{e_sigma*e_tau}=
\begin{Bmatrix}
0, &\sigma\ne\tau\\ \dim X_\sigma\cdot\dim X_\sigma\cdot\frac{1}{\dim X_\sigma}\cdot\widetilde{\chi_\sigma}
\cdot\mu_K,& \sigma=\tau
\end{Bmatrix}=\\=
\begin{Bmatrix}
0, &\sigma\ne\tau\\ \dim X_\sigma\cdot\widetilde{\chi_\sigma}
\cdot\mu_K,& \sigma=\tau
\end{Bmatrix}=\eqref{chi_sigma=dim X_sigma-widetilde(e_sigma)-mu_K}=
\begin{Bmatrix}
0, &\sigma\ne\tau\\ \nu_\sigma,& \sigma=\tau
\end{Bmatrix}.
\end{multline*}
\epr

Пусть $\pi:G\to {\mathcal L}(X)$ -- унитарное представление. Для всякого $\sigma\in\widehat{K}$ обозначим
 \beq\label{DEF:T_pi}
 \pi_\sigma(t)=\dot{\pi}(\nu_\sigma*\delta_t),\qquad  \dot{\pi}_\sigma(\alpha)=\dot{\pi}(\nu_\sigma*\alpha),\qquad t\in G,\quad \alpha\in {\mathcal C}^\star(G),\qquad \sigma\in \widehat{K}.
 \eeq
где $\dot{\pi}$ -- гомоморфизм алгебр из \eqref{C*(G)=grup-alg}.

 \vglue10pt
\centerline{\bf Свойства $\pi_\sigma$:}
 \vglue10pt

\bit{\it

\item[$1^\circ$.] Для всякого $\sigma\in\widehat{K}$ отображение $\pi_\sigma:G\to {\mathcal L}(X)$ является представлением группы $G$, а $\dot{\pi}_\sigma:{\mathcal C}^\star(G)\to {\mathcal L}(X)$ -- его продолжением на групповую алгебру:
\beq
\dot{\pi}_\sigma=(\pi_\sigma)^\cdot
\eeq

\item[$2^\circ$.] Для всякого $\sigma\in\widehat{K}$ пространство
\beq
X_\sigma=\overline{\dot{\pi}_\sigma({\mathcal C}^\star(G))X}
\eeq
инвариантно относительно $\pi$ (и поэтому определяет подпредставление представления $\pi$).

\item[$3^\circ$.] Пространство $X$ раскладывается в гильбертову прямую сумму подпространств $X_\sigma$:
\beq\label{X=oplus_pi-X_pi}
X=\dot{\bigoplus_{\sigma\in\text{$\widehat{K}$}}}\ \ X_\sigma.
\eeq
а представление $\pi$ в гильбертову прямую сумму представлений $\pi_\sigma$:
\beq\label{T=oplus_pi-T_pi}
\pi=\dot{\bigoplus_{\sigma\in\text{$\widehat{K}$}}}\ \ \pi_\sigma.
\eeq

\item[$4^\circ$.] При ограничении на подгруппу $K$ разложение \eqref{T=oplus_pi-T_pi} превращается в разложение
 $\pi|_K$ на изотипические компоненты, кратные представлениям $\sigma$:
\beq\label{T=oplus_pi-T_pi|_K}
\pi|_K^\cdot=\dot{\bigoplus_{\sigma\in\text{$\widehat{K}$}}}\ \ \pi_\sigma|_K^\cdot, \qquad \pi_\sigma|_K^\cdot=\dot{\bigoplus_{i\in I_\sigma}}\sigma_i^\cdot,\qquad \sigma_i\cong\sigma\qquad (i\in I_\sigma).
\eeq
}\eit
\bpr
1. Прежде всего,
\begin{multline*}
\pi_\sigma(s\cdot t)=\dot{\pi}(\nu_\sigma*\delta_{s\cdot t})=\dot{\pi}(\nu_\sigma*\delta_s*\delta_t)=\eqref{chi_pi-ortog-proj}=
\dot{\pi}(\nu_\sigma*\nu_\sigma*\delta_s*\delta_t)=\\=\eqref{chi_pi-centr}=\dot{\pi}(\nu_\sigma*\delta_s *\nu_\sigma*\delta_t)=
\dot{\pi}(\nu_\sigma*\delta_s)\cdot \dot{\pi}(\nu_\sigma*\delta_t)=\pi_\sigma(s)\cdot \pi_\sigma(t)
\end{multline*}

2. Пусть $x\in \dot{\pi}_\sigma({\mathcal C}^\star(G))X$, то есть
$$
x=\sum_{i=1}^n \dot{\pi}_\sigma(\beta_i)x_i,\qquad \beta_i\in{\mathcal C}^\star(G),\ x_i\in X.
$$
Тогда для любого $\alpha\in{\mathcal C}(G)$
\begin{multline*}
\dot{\pi}(\alpha)x=\sum_{i=1}^n \dot{\pi}(\alpha)\dot{\pi}_\sigma(\beta_i)x_i=\sum_{i=1}^n \dot{\pi}(\alpha)\dot{\pi}(\nu_\sigma*\beta_i)x_i=\sum_{i=1}^n \dot{\pi}(\alpha*\nu_\sigma*\beta_i)x_i=\\=
\eqref{chi_pi-centr}=\sum_{i=1}^n \dot{\pi}(\nu_\sigma*\alpha*\beta_i)x_i=\sum_{i=1}^n \dot{\pi}_\sigma(\alpha*\beta_i)x_i
\in \dot{\pi}_\sigma(G)X.
\end{multline*}
Мы получаем
$$
\dot{\pi}({\mathcal C}^\star(G))\Big(\dot{\pi}_\sigma({\mathcal C}^\star(G))X\Big)\subseteq \dot{\pi}_\sigma({\mathcal C}^\star(G))X,
$$
откуда
$$
\dot{\pi}({\mathcal C}^\star(G))\Big(\overline{\dot{\pi}_\sigma({\mathcal C}^\star(G))X}\Big)\subseteq \overline{\dot{\pi}_\sigma({\mathcal C}^\star(G))X}.
$$

3. Свойства $3^\circ$ и $4^\circ$ мы докажем вместе. Рассмотрим сначала ограничение $\pi|_K:K\to{\mathcal L}(X)$. Как известно, оно раскладывается в гильбертову прямую сумму неприводимых представлений:
$$
\pi|_K=\dot{\bigoplus_{\sigma\in\text{$\widehat{K}$}}}\ \ \dot{\bigoplus_{i\in I_\sigma}}\sigma_i,\qquad \sigma_i\cong\sigma\qquad (i\in I_\sigma).
$$
При этом проекция на изотипическую компоненту
$$
\varPhi_\sigma:X=\dot{\bigoplus_{\sigma\in\text{$\widehat{K}$}}}\ \ \dot{\bigoplus_{i\in I_\sigma}}X_i\to \dot{\bigoplus_{i\in I_\sigma}}X_i
$$
и, одновременно, сплетающий оператор между представлениями $\pi$ и $\pi_\sigma$ (где $\pi_\sigma$ определено в \eqref{DEF:T_pi}),
$$
\varPhi_\sigma\circ\pi=\pi_\sigma\circ\varPhi_\sigma,
$$
описывается формулой \cite[Theorem 5.5.1(2)]{Kowalski}
\beq\label{DEF:varPhi_sigma}
\varPhi_\sigma=\dim X_\sigma\cdot\int_K\overline{\tr\sigma(s)}\cdot \pi(s)\cdot\mu_K(\d s)=\dim X_\sigma\cdot\int_K \tr\sigma(s^{-1})\cdot \pi(s)\cdot\mu_K(\d s)=\dot{\pi}(\nu_\sigma).
\eeq
Заметим, что
$$
\dot{\pi}_\sigma(\alpha)=\dot{\pi}(\nu_\sigma*\alpha)=
\dot{\pi}(\nu_\sigma)\cdot\dot{\pi}(\alpha)=\varPhi_\sigma\cdot \dot{\pi}(\alpha),\qquad \alpha\in{\mathcal C}^\star(G).
$$
Отсюда
$$
X_\sigma=\overline{\dot{\pi}_\sigma({\mathcal C}^\star(G))X}=\overline{\varPhi_\sigma \dot{\pi}({\mathcal C}^\star(G))X}=
\overline{\varPhi_\sigma X}=\overline{\dot{\bigoplus_{i\in I_\sigma}}X_i}=\dot{\bigoplus_{i\in I_\sigma}}X_i.
$$
Мы получаем разложение в гильбертову сумму
$$
X=\dot{\bigoplus_{\sigma\in\text{$\widehat{K}$}}}\ \ \dot{\bigoplus_{i\in I_\sigma}}X_i=\dot{\bigoplus_{\sigma\in\text{$\widehat{K}$}}}\ \ X_\sigma.
$$
Это доказывает и \eqref{T=oplus_pi-T_pi|_K} и \eqref{T=oplus_pi-T_pi}.
\epr

\subsection{Непрерывные по норме представления}

\paragraph{Непрерывные по норме представления.}

Вернемся снова к ситуации, когда задано представление $\pi:G\to A$ группы $G$ в стереотипной алгебре $A$ (см. диаграмму \eqref{C*(G)=grup-alg-A}). Пусть теперь $A$ -- банахова алгебра. Тогда непрерывность $\pi$ будет непрерывностью по норме:
$$
x_i\to x\quad\Longrightarrow\quad \norm{\pi(x_i)-\pi(x)}\to 0.
$$

\bex\label{EX:gladkost-nepr-po-norme-predstavleniya}
В случае, если $G$ -- вещественная группа Ли, любое ее непрерывное по норме представление $\pi:G\to A$ в произвольной банаховой алгебре $A$ обязательно будет гладким отображением, и, как следствие, будет (однозначно) продолжаться (в смысле диаграммы \eqref{E*(G)=grup-alg-A}) до гомоморфизма $\ddot{\pi}:{\mathcal E}^\star(G)\to A$ групповой алгебры ${\mathcal E}^\star(G)$ распределений с компактным носителем (описанной выше на с.\pageref{E*(G)=grup-alg-A}).
\eex
\bpr
Гладкость $\pi$ доказывается по аналогии с классическим результатом, когда $A$ -- конечномерная алгебра, см. \cite[4.21]{Michor}. После этого существование $\ddot{\pi}$ вытекает из \cite[10.12]{Ak03}.
\epr

Пусть теперь ${\mathcal B}(X)$, как обычно обозначает банахово пространство линейных непрерывных операторов на гильбертовом пространстве $X$ (то есть то же самое пространство ${\mathcal L}(X)$, но с обычной топологией, порожденной нормой оператора). Обозначим $\iota:{\mathcal B}(X)\to{\mathcal L}(X)$ естественное вложение (очевидно, оно непрерывно).

Следующее утверждение очевидно.

\btm\label{TH:nepr-po-norme-predstavleniya} Для представления $\pi:G\to{\mathcal L}(X)$ локально компактной группы $G$ в гильбертовом пространстве $X$ следующие условия эквивалентны:
 \bit{
\item[(i)] отображение $\pi:G\to{\mathcal B}(X)$ непрерывно,

\item[(ii)] морфизм стереотипных алгебр $\dot{\pi}:\mathcal{C}^\star(G)\to {\mathcal L}(X)$ из \eqref{C*(G)=grup-alg} поднимается до морфизма стереотипных алгебр $\ph:\mathcal{C}^\star(G)\to {\mathcal B}(X)$:
\begin{gather}\label{C*(G)=grup-alg-B(X)}
\begin{diagram}
\node[2]{\mathcal{C}^\star(G)} \arrow{sw,t,--}{\ph}\arrow{se,t}{\dot{\pi}}
\\
\node{{\mathcal B}(X)}\arrow[2]{e,b}{\iota} \node[2]{{\mathcal L}(X)}
\end{diagram}.
\end{gather}
 }\eit
\etm

\bit{
\item[$\bullet$] Представление $\pi:G\to{\mathcal L}(X)$, удовлетворяющее этим условиям, мы назывем {\it непрерывным по норме}.\label{DEF:nepr-po-norme-predst}

\item[$\bullet$] В частном случае, когда все операторы $\pi(t)$, $t\in G$, унитарны, представление $\pi:G\to{\mathcal L}(X)$, удовлетворяющее условиям этой теоремы, мы назывем {\it непрерывным по норме унитарным представлением}.\label{DEF:nepr-po-norme--unit-predst}
}\eit

\btm \cite[Corollary 2]{Shtern}\label{TH:nepr-po-norme-perdst-K} Унитарное представление $\pi:K\to{\mathcal L}(X)$ компактной группы $K$ непрерывно по норме в том и только в том случае, если в его разложении \eqref{T=oplus_pi-T_pi|_K}
$$
\dot{\pi}=\dot{\bigoplus_{\sigma\in\text{$\widehat{K}$}}}\ \ \dot{\pi}_\sigma,
$$
только конечное число изотипических компонент $\dot{\pi}_\sigma$ нетривиально.
\etm

\btm \cite[Theorem 4.8]{Kuznetsova}\label{TH:Kuz-2}
Если $K$ -- нормальная компактная подгруппа в $G$ и унитарное представление $\pi:G\to {\mathcal L}(X)$ непрерывно по норме, то в разложении \eqref{X=oplus_pi-X_pi} $X_\sigma\ne 0$ только для конечного набора индексов $\sigma\in\widehat{G}$, и поэтому сумма в \eqref{T=oplus_pi-T_pi} в этом случае конечна.
\etm

\btm\label{TH:nepr-po-norme-perdst-G} Если $\pi:G\to{\mathcal L}(X)$ --- унитарное непрерывное по норме представление локально компактной группы $G$, то его ядро 
$$
\Ker\pi=\{g\in G:\ \pi(g)=1_X\}
$$
является замкнутой подгруппой в $G$, фактор-группа по которой $G/\Ker\pi$ является группой Ли.
\etm
\bpr
1. Предположим, что это неверно то есть существует унитарное непрерывное по норме представление $\pi:G\to{\mathcal L}(X)$ локально компактной группы $G$, такое что фактор-группа по ядру $G/\Ker\pi$ не является группой Ли. Представим $\pi$ в виде композиции
$$
\pi=\rho\circ\tau
$$
где $\tau:G\to G/\Ker\pi$ --- фактор-отображение, а $\rho:G/\Ker\pi\to{\mathcal L}(X)$ --- инъективное унитарное непрерывное по норме представление локально компактной группы $G/\Ker\pi$. 

2. Мы получили, что из нашего предположения следует, что существует {\it инъективное} унитарное непрерывное по норме представление $\pi:G\to{\mathcal L}(X)$ локально компактной группы $G$, не являющейся группой Ли. Выберем теперь открытую LP-подгруппу $G'$ в $G$, то есть такую, которая представляется в виде проективного предела своих фактор-групп Ли:
$$
G'=\projlim_{K\in\lambda(G')} G'/K
$$
--- здесь $\lambda(G')$ --- система компактных нормальных подгрупп $K$ в $G'$, фактор-группы по которым $G'/K$ являются группами Ли (это всегда можно сделать, см. \cite[5.6]{Glushkov} или \cite[4.0]{Montgomery-Zippin}). Ограничение 
$$
\pi\Big|_{G'}:G'\to{\mathcal L}(X)
$$
будет по-прежнему инъективным унитарным непрерывным про норме представлением.

3. Теперь мы получаем, что существует инъективное унитарное непрерывное по норме представление $\pi:G\to{\mathcal L}(X)$ некоей LP-группы $G$, не являющейся группой Ли. Теперь выберем какую-нибудь группу $K\in\lambda(G)$ (то есть компактную нормальную подгруппу $K$ в $G$, фактор-группа по которой $G/K$ является группой Ли).  Ограничение 
$$
\pi\Big|_{K}:K\to{\mathcal L}(X)
$$
будет по-прежнему инъективным унитарным непрерывным про норме представлением. 
По теореме \ref{TH:nepr-po-norme-perdst-K}, в его разложении \eqref{T=oplus_pi-T_pi|_K}
$$
\dot{\pi}=\dot{\bigoplus_{\sigma\in\text{$\widehat{K}$}}}\ \ \dot{\pi}_\sigma,
$$
только конечное число изотипических компонент $\dot{\pi}_\sigma$ нетривиально. Выберем в каждом множестве $I_\sigma$ из \eqref{T=oplus_pi-T_pi|_K} по одному представителю $i_{I_\sigma}\in I_\sigma$. Тогда отображение
$$
x\in K\mapsto \ph(x)=\dot{\bigoplus_{\sigma\in\text{$\widehat{K}$}:\ \pi_\sigma\ne 1}}\ \  \sigma_{i_{I_\sigma}}
$$
будет тоже инъективным унитарным непрерывным по норме представлением компактной группы $K$, но в конечномерном гильбертовом пространстве 
$$
\dot{\bigoplus_{\sigma\in\text{$\widehat{K}$}:\ \pi_\sigma\ne 1}}\ \  X_\sigma.
$$
Значит, $K$ --- группа Ли.

4. Мы получили, что $G/K$ и $K$ --- группы Ли. Значит, $G$ --- тоже группа Ли. А это противоречит выбору $G$. 
\epr

\bcor\label{COR:G->B-propusk-cherez-Lie}
Пусть $G$ --- локально компактная группа. Всякий ее непрерывный унитарный гомоморфизм $\pi:G\to B$ в произвольную $C^*$-алгебру $B$ пропускается через фактор-отображение 
$$
\rho_K:G\to G/K
$$
на некоторую фактор-группу $G/K$, являющуюся группой Ли:
\beq\label{G->B-propusk-cherez-Lie}
\xymatrix @R=2.pc @C=3.0pc 
{
G\ar@{-->}[r]^{\rho_K}\ar[dr]_{\pi} &  G/K\ar@{-->}[d]^{\pi_K} \\
 & B
}
\eeq
(где $\pi_K:G/K\to B$ --- непрерывный гомоморфизм).
\ecor
\bpr
Вложим $C^*$-алгебру $B$ изометрически в алгебру операторов ${\mathcal B}(X)$ на подходящем гильбертовом пространстве $X$:
$$
\sigma:B\to {\mathcal B}(X).
$$ 
Композиция  
$$
\sigma\circ\pi:G\to B\to {\mathcal B}(X).
$$ 
будет иметь то же ядро, что и гомоморфизм $\pi$:
$$
K=\Ker (\sigma\circ\pi)=\Ker\pi.
$$ 
По теореме \ref{TH:nepr-po-norme-perdst-G}, $K$ --- нормальная подгруппа в $G$, фактор-группа про которой $G/K$ --- группа Ли. Поэтому $\pi$ раскладывается в композицию \eqref{G->B-propusk-cherez-Lie}.
\epr

\paragraph{Индуцированные представления.}

Напомним конструкцию индуцированного представления.
Пусть $N$ -- открытая нормальная подгруппа в локально компактной группе $G$. Обозначим $F=G/N$ и будем считать $G$ расширением группы $N$ по группе $F$:
\beq\label{G-kak-rasshirenie}
\xymatrix 
{
1\ar[r] & N \ar[r]^{\eta} & G\ar[r]^{\ph} & F\ar[r]& 1
}
\eeq
(здесь  $\eta$ -- естественное вложение, а $\ph$ -- фактор отображение).  Зафиксируем какую-нибудь функцию $\sigma:F\to G$, являющуюся коретракцией для фактор-отображения $\ph$,
\beq\label{ph(sigma(t))=t}
\ph(\sigma(t))=t,\qquad t\in F,
\eeq
и сохраняющую единицу:
\beq\label{sigma(1_D)=1_G}
\sigma(1_F)=1_G.
\eeq
Тогда для всякого $g\in G$ элемент $\sigma(\ph(g))$ будет принадлежать тому же классу смежности относительно $N$, что и $g$,
$$
g\in \sigma(\ph(g))\cdot N
$$
то есть
\beq\label{g-sigma(ph(g))^(-1)-in-N}
g\cdot \sigma(\ph(g))^{-1}\in N,\qquad g\in G.
\eeq

Условимся называть {\it непрерывным (по норме) представлением} группы $G$ в гильбертовом пространстве $X$ всякое непрерывное отображение $\pi:G\to{\mathcal B}(X)$, удовлетворяющее условиям гомоморфности
\beq\label{pi(g-cdot-h)=pi(g)-cdot-pi(h)-nepr}
\pi(g\cdot h)=\pi(g)\cdot\pi(h),\qquad \pi(1_G)=1_B,\qquad g,h\in G,
\eeq

\btm\footnote{\cite[Lemma 3.6]{Kuznetsova}.} \label{TH:Ind-Preds}
Пусть SIN-группа $G$ представлена как расширение \eqref{G-kak-rasshirenie} некоторой своей открытой нормальной подгруппы $N$, и  $\pi:N\to{\mathcal B}(X)$ -- непрерывное по норме представление. Рассмотрим пространство $L_2(F,X)$ суммируемых с квадратом функций $\xi:F\to X$ (относительно считающей меры $\card$ на $F$). Тогда формула
\begin{multline}\label{ind-representation}
\pi'(g)(\xi)(t)=\pi\Big(\underbrace{\sigma(t)\cdot g\cdot \sigma\big(\ph(\sigma(t)\cdot g)\big)^{-1}}_{\scriptsize\begin{matrix}\phantom{\quad\tiny \eqref{g-sigma(ph(g))^(-1)-in-N}}
\text{\rotatebox{90}{$\owns$}}{\quad\tiny \eqref{g-sigma(ph(g))^(-1)-in-N}}\\ N\end{matrix}}\Big)\Big(\xi\big(\underbrace{\ph(\sigma(t)\cdot g)}_{\scriptsize\begin{matrix}
\text{\rotatebox{90}{$\owns$}}\\ F\end{matrix}}\big)\Big)=\\=
\pi\Big(\sigma(t)\cdot g\cdot \sigma\big(t\cdot\ph(g))\Big)^{-1}\Big)\Big(\xi\big(t\cdot\ph(g))\Big),
 \qquad \xi\in L_2(F,X),\quad t\in F,\quad g\in G,
\end{multline}
определяет непрерывное по норме представление $\pi':G\to{\mathcal B}(L_2(F,X))$.
\etm

\bit{

\item[$\bullet$] Определенное таким образом представление $\pi':G\to {\mathcal B}(L_2(F,X))$ называется {\it представлением, индуцированным представлением $\pi:N\to {\mathcal B}(X)$}.

}\eit

\bpr
Зафиксируем $\e>0$. Поскольку $\pi$ непрерывно по норме, найдется окрестность единицы $U\subseteq N$ такая что
$$
\norm{\pi(h)-1_X}<\e,\qquad h\in U.
$$
Поскольку $N$ открыто в $G$, $U$ является окрестностью единицы и в $G$. При этом $G$ -- SIN-группа, поэтому найдется окрестность единицы $V\subseteq U$, инвариантная относительно сопряжений: $g\cdot V\cdot g^{-1}\subseteq V$ для любых $g\in G$. Теперь для всякого $h\in V\subseteq U\subseteq N$ и любого $t\in F$ мы получим
$$
\ph(\sigma(t)\cdot h)=\ph(\sigma(t))\cdot\ph(h)=t\cdot 1=t,
$$
и отсюда при $h\in V$ и $\xi\in L_2(F,X)$
\begin{multline*}
\norm{\pi'(h)(\xi)-\xi}^2=\sum_{t\in D}\norm{\pi'(h)(\xi)(t)-\xi(t)}^2=
\sum_{t\in D}\norm{\pi\Big(\sigma(t)\cdot h\cdot\sigma\big(\ph(\sigma(t)\cdot h)\big)^{-1}\Big)\Big(\xi\big(\ph(\sigma(t)\cdot h)\big)\Big)-\xi(t)}^2=\\=
\sum_{t\in D}\norm{\pi\big(\sigma(t)\cdot h\cdot\sigma(t)^{-1}\big)\big(\xi(t)\big)-\xi(t)}^2\le
\sum_{t\in D}\|\pi\big(\underbrace{\sigma(t)\cdot h\cdot\sigma(t)^{-1}}_{\scriptsize\begin{matrix}\text{\rotatebox{90}{$\owns$}}\\ V\end{matrix}}\big)-1_X\|^2\cdot\norm{\xi(t)}^2<
\e^2\cdot \sum_{t\in D}\norm{\xi(t)}^2=\e^2\cdot\norm{\xi}^2
\end{multline*}
\epr

\paragraph{Пространство $\Trig(G)$ непрерывных по норме тригонометрических многочленов.}
Условимся {\it непрерывным по норме тригонометрическим многочленом} на $G$ называть произвольную линейную комбинацию матричных элементов непрерывных по норме унитарных неприводимых представлений группы $G$
\beq\label{DEF:Trig(G)}
 u(t)=\sum_{i=1}^n\lambda_i\cdot\langle \sigma_i(t)x_i,y_i\rangle,\quad t\in G,
\eeq
($\sigma_i:G\to{\mathcal B}(X_i)$ -- непрерывные по норме унитарные неприводимые представления, $x_i,y_i\in X_i$, $\lambda_i\in\C$).

\bex\label{EX:trig-monom}
Пусть $G$ -- локально компактная группа, $U$ -- произвольная $C^*$-окрестность нуля в ${\mathcal C}^\star(G)$ и $v$ -- произвольное чистое состояние $C^*$-алгебры ${\mathcal C}^\star(G)/U$. Тогда функция
$$
G\overset{\delta}{\longrightarrow}{\mathcal C}^\star(G)\overset{\pi_U}{\longrightarrow}{\mathcal C}^\star(G)/U\overset{v}{\longrightarrow}\C
$$
является непрерывным по норме тригонометрическим многочленом на $G$.
\eex
\bpr
Рассмотрим ГНС-представление $\sigma:{\mathcal C}^\star(G)/U\to {\mathcal B}(X)$, порожденное функционалом $v$:
$$
v(a)=\langle \sigma(a)x,x\rangle,\qquad a\in {\mathcal C}^\star(G)/U,
$$
($x\in X$). Поскольку $v$ -- чистое состояние, представление $\sigma$ будет неприводимым \cite[10.2.3]{Kad-Ring-2}. Поскольку $G$ полно отображается в алгебру ${\mathcal C}^\star(G)$ (это значит, что линейные комбинации дельта-функций полны в ${\mathcal C}^\star(G)$) и значит, и в фактор-алгебру ${\mathcal C}^\star(G)/U$, мы получаем, что композиция
$$
G\overset{\delta}{\longrightarrow}{\mathcal C}^\star(G)\overset{\pi_U}{\longrightarrow}{\mathcal C}^\star(G)/U\overset{\sigma}{\longrightarrow}{\mathcal B}(X)
$$
является неприводимым представлением группы $G$. Значит, функция
$$
u(t)=v(\pi_U(\delta(t)))=\langle \sigma(\pi_U(\delta(t)))x,x\rangle
$$
является непрерывным по норме тригонометрическим многочленом.
\epr

\brem
Если группа $G$ компактна или абелева, то пространство $\Trig(G)$ является алгеброй относительно операции поточечного умножения. Однако, по-видимому, существуют (некомпактные и неабелевы) группы, для которых это не так. Если не требовать непрерывности тригонометрического многочлена относительно $C^*$-полунорм на ${\mathcal C}^\star(G)$, то нужным контрпримером становится группа $G=H_3(\R)$ верхнетреугольних матриц с вещественными коэффициентами и единицами на диагонали: на ней пространство $\Trig(G)$ (без требования непрерывности относительно $C^*$-полунорм) не замкнуто относительно поточечного умножения -- этот пример был сообщен автору И.~Чои.
\erem

\paragraph{Алгебра $k(G)$ непрерывных по норме матричных элементов.}
{\it Непрерывным по норме матричным элементом} на группе $G$ мы называем функцию на $G$ вида
\beq\label{DEF:k(G)}
 u(t)=\langle \pi(t)x,y\rangle,\quad t\in G,
\eeq
где $\pi:G\to{\mathcal B}(X)$ -- какое-нибудь непрерывное по норме унитарное представление в некотором гильбертовом пространстве $X$, и $x,y\in X$. Множество всех таких функций мы будем обозначать $k(G)$.

\btm
Для всякой локально компактной группы $G$ пространство $k(G)$ образует инволютивную алгебру над $\C$ относительно поточечных алгебраических операций.
\etm
\bit{
\item[$\bullet$] Мы будем наделять алгебру $k(G)$ сильнейшей локально выпуклой топологией над $\C$. Это превращает $k(G)$ в инволютивную стереотипную алгебру (в силу \cite[Example 5.1.4]{Akbarov-De-Gruyter-I}).
}\eit
\bpr
Умножение на скаляр в $k(G)$ эквивалентно умножению вектора $x$ на этот скаляр в формуле \eqref{DEF:k(G)}. Рассмотрим  две функции вида \eqref{DEF:k(G)}
$$
u(t)=\langle \pi(t)x,y\rangle,\quad v(t)=\langle \pi'(t)x',y'\rangle,\qquad t\in G,
$$
где $\pi:G\to{\mathcal B}(X)$ и $\pi':G\to{\mathcal B}(X')$ -- непрерывные по норме унитарные представления. Если рассмотреть гильбертову прямую сумму $X\oplus X'$, то есть пространство пар $(\xi,\xi')$, $\xi\in X$, $\xi'\in X'$, со скалярным произведением
$$
\langle (\xi,\xi'),(\upsilon,\upsilon')\rangle=\langle \xi,\upsilon\rangle+\langle \xi',\upsilon'\rangle,
$$
то формула
$$
(\pi\oplus\pi')(t)(\xi,\xi')=(\pi(t)\xi,\pi'(t)\xi'),\qquad t\in G,
$$
определяет непрерывное по норме представление $\pi\oplus\pi':G\to{\mathcal B}(X\oplus X')$, матричным элементном которого будет сумма функций $u$ и $v$:
$$
\langle (\pi\oplus\pi')(t)(x,x'),(y,y')\rangle=\langle (\pi(t)x,\pi'(t)x'),(y,y')\rangle=
\langle \pi(t)x,y\rangle+\langle \pi'(t)x',y'\rangle=u(t)+v(t).
$$
С другой стороны, если рассмотреть гильбертово тензорное произведение $X\dot{\otimes} X'$ \cite{Kad-Ring}, то есть пополнение алгебраического тензорного произведения $X\otimes X'$ относительно скалярного произведения, определяемого равенством
$$
\langle \xi\otimes\xi',\upsilon\otimes\upsilon'\rangle=\langle \xi,\upsilon\rangle\cdot \langle \xi',\upsilon'\rangle, \qquad \xi,\upsilon\in X,\quad \xi',\upsilon'\in X',
$$
то формула
$$
(\pi\otimes\pi')(t)\xi\otimes\xi'=\pi(t)\xi\otimes\pi'(t)\xi',\qquad t\in G,
$$
определяет непрерывное по норме представление $\pi\otimes\pi':G\to{\mathcal B}(X\dot{\otimes} X')$, матричным элементном которого будет произведение функций $u$ и $v$:
$$
\langle (\pi\otimes\pi')(t)x\otimes x',y\otimes y'\rangle=\langle \pi(t)x\otimes\pi'(t)x',y\otimes y'\rangle=
\langle \pi(t)x,y\rangle\cdot\langle \pi'(t)x',y'\rangle=u(t)\cdot v(t).
$$
Замкнутость $k(G)$ относительно поточечной инволюции следует из существования сопряжения в $X$ \cite[27.25,27.26]{Hewitt-Ross-2}. Тождественная единица лежит в $k(G)$ потому что ее можно представить в виде \eqref{DEF:k(G)}, положив $\pi(t)=1$ и взяв произвельные $x$ и $y$ так, чтобы $\langle x,y\rangle=1$.
\epr

Условимся говорить, что функция $u\in{\mathcal C}(G)$ {\it подчинена $C^*$-полунорме},\label{DEF:positive-def-func-sub-C*} если для некоторой непрерывной $C^*$-полунормы $p$ на ${\mathcal C}^\star(G)$ выполняется неравенство
$$
|\alpha(u)|\le p(\alpha),\qquad \alpha\in{\mathcal C}^\star(G).
$$
Понятно, что все функции из алгебры $k(G)$ подчинены $C^*$-полунормам.

Напомним еще, что функция $u\in{\mathcal C}(G)$ называется {\it положительно определенной} \cite[\S 32]{Hewitt-Ross-2}, если выполняется неравенство
\beq\label{(alpha^bullet*alpha)(u)-ge-0}
(\alpha^\bullet*\alpha)(u)\ge 0,\qquad \alpha\in {\mathcal C}^\star(G).
\eeq

\btm\label{TH:u(t)=langle-pi(t)x,x-rangle} Для функции $u\in{\mathcal C}(G)$ следующие условия эквивалентны:
\bit{
\item[(i)] $u$ представима в виде
\beq\label{u(t)=langle-pi(t)x,x-rangle}
u(t)=\langle\pi(t)x,x\rangle,\qquad t\in G,
\eeq
для некоторого непрерывного по норме унитарного представления $\pi:G\to{\mathcal B}(X)$ и некоторого $x\in X$ (и, как следствие, $u\in k(G)$);

\item[(ii)] $u$ положительно определена и подчинена некоторой $C^*$-полунорме.
}\eit
\etm
\bpr Здесь нужно только проверить импликацию (i)$\Longleftarrow$(ii). Пусть выполняется (ii).
Поскольку $u$ подчинена $C^*$-полунорме $p$, эту функцию можно продолжить как функционал на $C^*$-фактор-алгебре ${\mathcal C}^\star(G)/p$:
$$
u=v\circ\rho,
$$
где $\rho:{\mathcal C}^\star(G)\to {\mathcal C}^\star(G)/p$ -- фактор-отображение, а $v:{\mathcal C}^\star(G)/p\to\C$ -- некоторый непрерывный функционал на $C^*$-алгебре ${\mathcal C}^\star(G)/p$. При этом, поскольку $u$ положительно определена, а инволюция и умножение в ${\mathcal C}^\star(G)/p$ наследуются из ${\mathcal C}^\star(G)$, $v$ будет положительным функционалом. Рассмотрим ГНС-представление $\sigma:{\mathcal C}^\star(G)/p\to {\mathcal B}(X)$, порожденное функционалом $v$. Тогда
$$
v(a)=\langle \sigma(a)x,x\rangle,\qquad a\in {\mathcal C}^\star(G)/p,
$$
для некоторого $x\in X$, и если положить $\pi=\sigma\circ\rho$, то
$$
\alpha(u)=v(\rho(\alpha))=\langle \sigma(\rho(\alpha))x,x\rangle=\langle \pi(\alpha)x,x\rangle,
\qquad \alpha\in {\mathcal C}^\star(G).
$$
\epr

Очевидно,
\beq\label{Trig(G)-subseteq-k(G)}
\Trig(G)\subseteq k(G).
\eeq

\btm\label{TH:Trig(G)=k(G)-na-kompaktah}
Если группа $G$ компактна, то
\bit{
\item[(a)] вложение \eqref{Trig(G)-subseteq-k(G)} превращается в равенство
$$
\Trig(G)=k(G);
$$
\item[(b)] алгебра $k(G)$ совпадает с множеством матричных элементов \eqref{DEF:k(G)} всевозможных {\rm конечномерных} непрерывных унитарных представлений $\pi:G\to{\mathcal B}(X)$;

\item[(c)] алгебра $k(G)$ совпадает с множеством функций вида
 \beq\label{u(t)=f(pi(t))}
u(t)=f(\pi(t)),\qquad t\in G,
\eeq
где $\pi:G\to{\mathcal B}(X)$ -- произвольное непрерывное по норме унитарное представление, а $f:{\mathcal B}(X)\to\C$ -- произвольный линейный непрерывный функционал.
}\eit
\etm
\bpr Здесь достаточно убедиться, что всякая функция вида \eqref{u(t)=f(pi(t))} лежит в $\Trig(G)$. Пусть $\pi:G\to{\mathcal B}(X)$ -- непрерывное по норме унитарное представление, и $f:{\mathcal B}(X)\to\C$ -- линейный непрерывный функционал.

1. Рассмотрим сначала случай, когда гильбертово пространство $X$ конечномерно. Тогда представление $\pi:G\to{\mathcal B}(X)$ раскладывается в сумму неприводимых,
$$
\pi=\bigoplus_{i=1}^n\pi_i:G\to{\mathcal B}\left(\bigoplus_{i=1}^nX_i\right),
$$
Выберем в каждом $X_i$ ортогональный нормированный базис $e_{ij}$. Тогда функционалы вида
$$
A\in{\mathcal B}(X)\mapsto \lambda_{kl}^{ij}\cdot\langle Ae_{ij},e_{kl}\rangle\in C
$$
образуют базис а пространстве ${\mathcal B}^\star(X)$, поэтому $f\in{\mathcal B}^\star(X)$ раскладывается в сумму
$$
u(t)=f(\pi(t))=\sum_{i,j,k,l}\lambda_{kl}^{ij}\cdot\langle\pi(t)e_{ij},e_{kl}\rangle= \sum_{i,j,k,l}\lambda_{kl}^{ij}\cdot\langle\underbrace{\pi_i(t)e_{ij}}_{ \scriptsize\begin{matrix}\text{\rotatebox{90}{$\owns$}}\\ X_i\end{matrix}}, \underbrace{e_{kl}}_{ \scriptsize\begin{matrix}\text{\rotatebox{90}{$\owns$}}\\ X_k\end{matrix}}\rangle=
\sum_{i,j,l}\lambda_{kl}^{ij}\cdot\langle\underbrace{\pi_i(t)e_{ij}}_{ \scriptsize\begin{matrix}\text{\rotatebox{90}{$\owns$}}\\ X_i\end{matrix}}, \underbrace{e_{il}}_{ \scriptsize\begin{matrix}\text{\rotatebox{90}{$\owns$}}\\ X_i\end{matrix}}\rangle
$$
Последняя сумма уже лежит в $\Trig(G)$, потому что каждое $\pi_i$ является неприводимым представлением.

2. Если $X$ бесконечномерно, то по теореме А.И.Штерна \ref{TH:nepr-po-norme-perdst-K} непрерывное по норме представление $\pi:G\to{\mathcal B}(X)$ все равно пропускается через некоторое конечномерное унитарное неприводимое представление $\sigma$ (см. диаграмму ниже, в которой $\ph$ -- гомоморфизм инволютивных алгебр). Положив $g=f\circ\ph$, мы получим функционал на конечномерной алгебре ${\mathcal B}(Y)$, порождающий ту же самую функцию $u$
$$
\xymatrix @R=2.pc @C=3.0pc 
{
G\ar@/_4ex/[ddr]_{u}\ar[rr]^{\sigma}\ar[dr]_{\pi} &  & {\mathcal B}(Y)\ar@{-->}[dl]^{\ph}\ar@/^4ex/@{-->}[ddl]^{g} \\
& {\mathcal B}(X)\ar[d]^{f} & \\
& \C &
}
$$
и по уже доказанному, $u\in\Trig(G)$.
\epr

\brem
Если группа $G$ некомпактна, то \eqref{Trig(G)-subseteq-k(G)} не обязано быть равенством. Например, этого не будет для группы $G=\Z$. Этот пример был сообщен автору А.~И.~Дегтяревым: достаточно взять регулярное представление $\pi:\Z\to L_2(\Z)$, $\pi(k)u(l)=u(l+k)$, тогда обозначив
$$
\chi_0(k)=\begin{cases}0,& k\ne 0\\ 1,& k=1 \end{cases}
$$
мы получим матричный элемент
$$
f(k)=\langle \pi(k)\chi_0,\chi_0\rangle=\chi_0(k),\qquad k\in\Z,
$$
не представимый в виде тригонометрического многочлена. Действительно, предположим, что $f$ -- тригонометрический многочлен на $\Z$, то есть функция вида
$$
f(k)=\sum_{n=1}^N\lambda_n\cdot\e_n^k,
$$
где $\e_n\in\C$, $|\e_n|=1$, $\lambda_n\in\C$. Тогда для вычисления коэффициентов $\lambda_k$ можно составить систему
$$
\sum_{n=1}^N\lambda_n\cdot\e_n^k=f(k)=\chi_0(k)=0,\qquad k=1,...,N
$$
которая будет линейно независима, потому что ее определитель есть определитель Вандермонда, и он ненулевой. Поскольку столбец свободных членов у этой системы нулевой, это значит, что коэффициенты $\lambda_n$ должны быть нулевыми. То есть $f=0$, а это невозможно, потому что $f(0)=\chi_0(0)=1$.
\erem

Напомним, что в теории алгебраических групп (см. \cite[Глава VI]{Chevalley} и \cite[45.6]{Zhelobenko-2}) всякой компактной группе Ли $G$ сопоставляется алгебра функций на ней, называемая {\it представляющим кольцом} и состоящая из всевозможных матричных элементов произвольных непрерывных конечномерных представлений группы $G$. Это в точности множество функций, упоминаемое в утверждении (b) теоремы \ref{TH:Trig(G)=k(G)-na-kompaktah}. Как следствие, справедлива

\btm\label{TH:k(G)=predst-koltso}
Если $G$ -- компактная группа Ли, то алгебра $k(G)=\Trig(G)$ совпадает с представляющим кольцом группы $G$.
\etm

Напомним, что {\it многочленом} на полной линейной группе $GL_n(\R)$ (или $GL_n(\C)$) называется функция вида
\beq\label{DEF:mnogochlen-na-GL_n(R)}
u(g)=\frac{f(g)}{(\det g)^k},\qquad g\in GL_n(\R),
\eeq
где $f$ -- многочлен от коэффициентов матрицы $g$, а $k\in\Z_+$. Множество всех многочленов на $GL_n(\C)$ мы будем обозначать ${\mathcal P}(GL_n(\C))$. Оно образует инволютивную алгебру относительно поточечных операций.

Следующий результат считается классическим (см. \cite[Глава VI]{Chevalley} и \cite[45.6]{Zhelobenko-2}):

\btm
 Для всякой компактной вещественной группы Ли $G$

\bit{

\item[--] cуществует некоторое $n\in\N$ и непрерывный инъективный гомоморфизм (вложение) групп $\sigma:G\to GL_n(\C)$;

\item[--] при любом таком вложении $\sigma:G\to GL_n(\C)$ алгебра $k(G)$ изоморфна фактор-алгебре алгебры многочленов ${\mathcal P}(GL_n(\C))$ по идеалу $I$ многочленов, обнуляющихся на $G$,
    \beq\label{k(G)=P(GL_n(C))/I}
    k(G)={\mathcal P}(GL_n(\C))/I,
    \eeq
    и такой же изоморфизм справедлив для вещественных частей:
    \beq\label{Re-k(G)=Re-P(GL_n(C))/Re-I}
    \Real k(G)=\Real{\mathcal P}(GL_n(\C))/\Real I;
    \eeq

\item[--] как следствие, группа $G$ обладает естественной струкрутой вещественной алгебраической группы, для которой $\Real k(G)$ является алгеброй (вещественных) многочленов,
    $$
    {\mathcal P}(G)=\real k(G)
    $$
    а множество $\C G$ общих нулей идеала $I$
    \beq\label{DEF:CG-dlya-komp-grupp}
    \C G=\{g\in GL_n(\C):\ \forall u\in I\ u(g)=0 \}
    \eeq
    (или, что эквивалентно, комплексный спектр $\C\Spec k(G)$ алгебры $k(G)$) обладает естественной струкрутой комплексной алгебраической группы, для которой $k(G)$ является алгеброй (комплексных) многочленов:
    $$
    {\mathcal P}(\C G)=k(G).
    $$

\item[--] при естественном вложении $G\subseteq\C G$ группа $G$ является вещественной формой группы $\C G$, то есть
\bit{
\item[(a)] $G$ имеет непустое пересечение с любой связной компонентой $\C G$, и

\item[(b)] касательная алгебра $L(G)$ есть вещественная часть касательной алгебры $L(\C G)$:
\beq\label{L(G)=Real-L(CG)}
L(G)=\Real L(\C G)
\eeq

}\eit
}\eit
\etm

\bcor\label{COR:kasat-prostr-k-Trig(G)}
Пусть $G$ -- компактная вещественная группа Ли. Тогда
\bit{

\item[(i)] каждый инволютивный характер $s\in\Spec(k(G))$ является значением в некоторой точке $a\in G$
\beq\label{Spec(k(G))}
s(u)=u(a),\qquad u\in k(G),
\eeq
поэтому спектр алгебры $k(G)$ совпадает с $G$:
\beq\label{Spec-k(G)=R-Spec-Real-k(G)=G}
\Spec k(G)=\R\Spec\Real k(G)=G
\eeq

\item[(ii)] всякий касательный вектор $\tau\in T_a[k(G)]$ к алгебре $k(G)$ в произвольной точке $a\in G$ представляет собой дифференцирование вдоль некоторой однопараметрической подгруппы $x:\R\to G$:
\beq\label{kasat-prostr-k-Trig(G)}
\tau(u)=\lim_{t\to 0}\frac{u(a\cdot x(t))-u(a)}{t},\qquad u\in k(G),
\eeq
и однозначно продолжается до касательного вектора алгебры ${\mathcal E}(G)$ в точке $a\in G$; поэтому касательные пространства у $k(G)$ и ${\mathcal E}(G)$ совпадают:
\beq\label{T_a[k(G)]=T_a(G)}
T_a[k(G)]=T_a[{\mathcal E}(G)]=T_a(G).
\eeq
}\eit
\ecor
\bpr
Оба утверждения следуют из \cite[Example5.1.37]{Akbarov-De-Gruyter-I}: формула \eqref{Spec-k(G)=R-Spec-Real-k(G)=G} следует из \cite[(5.99)]{Akbarov-De-Gruyter-I}\footnote{Формула \eqref{Spec-k(G)=R-Spec-Real-k(G)=G} справедлива для произвольных компактных групп, необязательно групп Ли, см.\cite[(30.30)]{Hewitt-Ross-2}.}, поскольку $G$ -- вещественное алгебраическое многообразие, а $k(G)$ -- комплексификация $\Real k(G)$. А с другой стороны, из \cite[(5.100)]{Akbarov-De-Gruyter-I} следует, что в \eqref{kasat-prostr-k-Trig(G)} $\tau$ можно считать дифференцированием вдоль некоторой гладкой кривой в $G$, а это эквивалентно дифференцированию вдоль некоторой однопараметрической подгруппы.
\epr

\subsection{Отражения в стереотипных алгебрах}

Пусть $A$ -- стереотипная алгебра и $A^\star$ -- ее сопряженное
стереотипное пространство, рассматриваемое как двусторонний
стереотипный $A$-модуль с операциями
$$
(a\cdot u\cdot b)(x)=u(b\cdot x\cdot a), \qquad a,b,x\in A, u\in
A^\star
$$
Линейное непрерывное отображение $@ :A^\star \to A $ называется {\it
отражением} на $A$, если выполнены следующие два тождества:
 \begin{align}
& @(a\cdot u\cdot b)=a\cdot @ u\cdot b, && a,b\in A,
u\in A^\star \label{DEF:reflection-1}\\
& u(@ v)=v(@ u), && u,v \in A^\star
\label{DEF:reflection-2}
 \end{align}
Первое означает, что оператор $@:A^\star\to A$ является морфизмом $A$-бимодулей, а второе -- что он совпадает со своим сопряженным при отождествлении $i_A:A\to A^{\star\star}$:
\beq\label{DEF:reflection-2-star}
@^\star=@.
\eeq
Если обозначить $@(u,v)=v(@ u)$, то мы получим билинейную форму
$@:A^\star\times A^\star\to \Bbb C$ со свойствами
 \begin{align}
& @(a\cdot u\cdot b,v)=@(u, b\cdot v\cdot a), &&
a,b\in A, u,v\in A^\star \label{DEF:reflection-1*}
 \\
&
@(u,v)=@(v,u),&& u,v \in A^\star \label{DEF:reflection-2*}
 \end{align}
Поэтому отражение можно представлять себе как билинейную форму на
сопряженном к $A$ прост\-ранстве $A^\star$, удовлетворяющую
тождествам (i)$'$-(ii)$'$.

Отражение $@$ мы назовем {\it невырожденным}, если оно имеет
тривиальное ядро ($\Ker @=0$), или, что равносильно, если образ $@$
плотен в $A^\star$. На языке билинейных форм это означает, что
$@(\cdot,\cdot)$ должна быть невырождена:
$$
\forall u\ne 0 \quad \exists v\quad @(u,v)\ne 0
$$

\noindent\rule{160mm}{0.1pt}\begin{multicols}{2}

\paragraph{Отражения в алгебре ${\mathcal L}(X)$.}

\bex{\it Отражения в стереотипной алгебре операторов $\mathcal{L}(X)$.}\label{ex-14.2} Для всякого стереотипного пространства $X$
преобразование Гротендика $@_X:\mathcal{L}^\star(X)\to \mathcal{L}(X)$
\beq\label{eq14.3}
f\circledast x\in X^\star \circledast X=\mathcal  L^\star (X) \mapsto f\odot x
\in X^\star \odot X=\mathcal{L}(X) 
\eeq
является единственным с точностью до скалярного множителя отражением на алгебре
операторов $\mathcal{L}(X)$. Мы называем это отражение $@_X$ нормальным
отражением на $\mathcal{L}(X)$. Отражение $@_X$ тогда и только тогда будет
точным, когда пространство $X$ обладает свойством стереотипной аппроксимации.
\eex

\bpr  Покажем, что преобразование Гротендика действительно является отражением.
Понятно, что каждое тождество (i), (ii) достаточно проверять для одномерных
тензоров $u,v\in \mathcal  L^\star (X)$, поскольку $X^\star \otimes X$  плотно
в $X^\star \circledast X=\mathcal{L}^\star (X)$ \cite[Proposition 4.4.23]{Akbarov-De-Gruyter-I}.

1. Проверим (i). Пусть $u=f\circledast x, v=g\circledast y, \, f,g\in X^\star,
x,y\in X$. Тогда $@_X u=f\odot x, @_X v=g\odot y$, и 
\begin{multline*}
u(@_X v)= f\circledast x
\, (g\odot y)=\\=\cite[(4.185)]{Akbarov-De-Gruyter-I}= g\circledast y \, (f\odot x)= v(@_X u)
\end{multline*}

2. Для тождества (ii) также выберем одномерный тензор $u=f\circledast x$. Тогда
$@_X u=f\odot x$ и для любых $\ph, \psi \in \mathcal{L}(X)$ получаем 
\begin{multline*}
@_X(\ph
\cdot u \cdot \psi)= @_X(\ph \cdot f\circledast x \cdot \psi)= \cite[(4.181)]{Akbarov-De-Gruyter-I}=\\=
@_X((f\circ \psi)\circledast \ph(x))=\eqref{eq14.3}=\\= (f\circ \psi)\odot
\ph(x)=\cite[(4.182)]{Akbarov-De-Gruyter-I}= \ph \circ f\odot x \circ  \psi=\\= \ph \cdot @_X u \cdot
\psi
\end{multline*}

Итак, доказано, что преобразование Гротендика $@_X :\mathcal  L^\star (X) \to
\mathcal{L}(X)$ является отражением. Пусть теперь $@ :\mathcal  L^\star (X) \to
\mathcal{L}(X)$ --- какое-нибудь другое отражение. Зафиксируем $f\in X^\star$ и
$x\in X$ такие, что
$$f(x)=1.$$
Обозначим $u =f\circledast x, \quad \pi=f\odot x$ и заметим, что
\beq\label{eq14.4}
\pi \cdot u= u= u\cdot \pi 
\eeq
Действительно, с одной стороны, 
\begin{multline*}
(\pi \cdot u)(\ph)= u(\ph \circ \pi)=
f\circledast x(\ph \circ f\odot x)=\\= f\circledast x(f\odot \ph(x))= f(x)\cdot
f(\ph(x))= f\circledast x(\ph)= u(\ph)
\end{multline*}
 а с другой, 
\begin{multline*}
(u\cdot \pi)(\ph)= u(\pi
\circ \ph)= f\circledast x(f\odot x \circ \ph)=\\= f\circledast x((f\circ
\ph)\odot x)= f(x)\cdot f(\ph(x))= (f\circledast x)(\ph)= u(\ph)
\end{multline*}

Из равенств \eqref{eq14.4} и тождества (ii) следует
$$
\pi \circ @ u=@ u=@ u \circ \pi
$$
поэтому если обозначить $@ u=\psi$, то мы получим
$$
(f\odot x) \circ \psi=\psi=\psi \circ (f\odot  x)
$$
или
$$
(f\circ \psi)\odot x =\psi=f\odot \psi (x)
$$

Таким образом, оператор $\psi$ одномерен и имеет вид
$$
\psi = \lambda \cdot f\odot x
$$
где число $\lambda \in \C$ удовлетворяет равенствам
$$
\lambda f=f\circ \psi, \quad\quad \lambda x=\psi(x)
$$

Итак, для фиксированного тензора $f\circledast x$
\beq
@ (f\circledast x)=\lambda \cdot f\odot x \label{eq14.5}
\eeq
Покажем, что $\forall g\in X^\star, \quad \forall y\in X$
\beq
@ (g\circledast y)=\lambda \cdot g\odot y \label{eq14.6}
\eeq
Действительно, из формул
\begin{align*}
& f\odot y \cdot f\circledast x \cdot g\odot x= f(x)^2 \cdot g\circledast y=
g\circledast y \\
& f\odot y \cdot f\odot x \cdot g\odot x= f(x)^2
\cdot g\odot y= g\odot y
\end{align*}
следует $@(g\circledast y)= f\odot y \cdot @(f\circledast x) \cdot g\odot
x=\eqref{eq14.5}= f\odot y \cdot \lambda \cdot f\odot x \cdot g\odot x= \lambda
\cdot g\odot y$. Мы доказали \eqref{eq14.6}, откуда $@=\lambda @_X$. Вторая
часть предложения следует из \cite[Theorem 4.4.65]{Akbarov-De-Gruyter-I}. \epr

\paragraph{Отражения в алгебрах ${\mathcal C}^\star(G)$ и ${\mathcal E}^\star(G)$.}

На алгебрах $\mathcal C^\star(G)$ мер на компактной группе
$G$ и $\mathcal E^\star(G)$ распределений на компактной группе Ли $G$
стандартное отражение $@_G$ определяется эквивалентными формулами:
 \beq\label{DEF:@_G-C^*(G)}
@_G u=\int_G u(t^{-1})\cdot \delta_t \;\d t, \quad u,v \in \mathcal C(G),\mathcal
E(G)
 \eeq
 \beq\label{DEF:@_G-C^*(G)-1}
 @_G(u,v)=\int_G
u(t^{-1})\cdot  v (t) \;\d t,\quad u,v \in \mathcal C(G),\mathcal
E(G)
 \eeq
где $\delta_t(u)=u(t)$ -- $\delta$-функционал, а интегрирование
ведется по нормированной мере Хаара: $\mu(G)=1$.

Это не единственный пример отражения в $\mathcal C^\star(G)$ (и в $\mathcal E^\star(G)$), но типичный. Такое отражение называется {\it нормальным отражением} в $\mathcal C^\star(G)$ (и в $\mathcal E^\star(G)$).

 \vglue10pt
\centerline{\bf Свойства нормальных отражений} 
\centerline{в $\mathcal C^\star(G)$ (и в $\mathcal E^\star(G)$):}
 \vglue10pt

\bit{\it
\item[$1^\circ$.] Перестановочность с антиподом:
\beq\label{@_G-widetilde(u)}
@_G\widetilde{u} =\widetilde{@_G u},\quad  u\in {\mathcal C}(G)  
\eeq
\item[$2^\circ$.] Замена на антипод при свертке:
\beq\label{(@_Gu)*v}
(@_G u) * v  =\widetilde{u}* v, \quad u,v\in {\mathcal C}(G)  
\eeq
\item[$3^\circ$.] Перестановочность со сверткой функций:
\beq \label{@_G(u*v)}
@_G (u * v)  =@_G v* @_Gu, \quad u,v\in {\mathcal C}(G) 
\eeq

\item[$4^\circ$.] Замена на антипод при внесении меры под знак отражения:
\beq\label{alpha*(@_Gu)}
\alpha * @_G u   =@_G(u*\widetilde{\alpha}), \quad u\in {\mathcal C}(G),\ \alpha\in {\mathcal C}^\star(G)
\eeq

}\eit
\bpr
1.
\begin{multline*}
@_G\widetilde{u}(v)=\int_G\widetilde{u}(t^{-1})\cdot v(t)\d t=\int_Gu(t)\cdot v(t)\d t =\\= \int_Gu(t)\cdot\widetilde{v}(t^{-1})\d t=@u(\widetilde{v})=\widetilde{@u}(v).
\end{multline*}

2.
\begin{multline*}
\big((@_Gu)*v\big)(s)=\eqref{eq10.10}=@_Gu(\widetilde{s\cdot v})=\\=\int_G u(t^{-1})\cdot \widetilde{s\cdot v}(t)\d t=\eqref{DEF:@_G-C^*(G)}=\\=
\int_G u(t^{-1})\cdot (s\cdot v)(t^{-1})\d t=\int_G \widetilde{u}(t)\cdot v(t^{-1}\cdot s)\d t=\\=
\widetilde{u}*v(s)
\end{multline*}

3.
\begin{multline*}
@(u*v)=\eqref{(@_Gu)*v}=@(@\widetilde{u}*v)=\eqref{@_G-widetilde(u)}=\\=
@(\widetilde{@u}*v)=\eqref{eq11.2}=@(v\cdot @u)=\eqref{DEF:reflection-1}=\\=@v*@u
\end{multline*}

4. Тождество \eqref{alpha*(@_Gu)} достаточно проверить для случая $\alpha=\delta^t$:
\begin{multline*}
\delta^t*@u=@(\delta^t\cdot u)=@(t\cdot u)=@(u*\delta^{t^{-1}})=@(u*\widetilde{\delta^t})
\end{multline*}

\epr

\paragraph{Отражения в алгебрах ${\mathcal O}^\star(G)$ и ${\mathcal P}^\star(G)$.}

На алгебрах $\mathcal O^\star(G)$ аналитических функционалов и $\mathcal P^\star(G)$ многочленов на редуктивной комплексной алгебраической группе Ли $G$ стандартное отражение $@_G$ определяется формулами похожими на \eqref{DEF:@_G-C^*(G)} и \eqref{DEF:@_G-C^*(G)-1}, но только интегрирование в них должно вестись по вещественной форме $G_{\Bbb
R}$ группы $G$ \cite{Vinberg}:
:
 \beq\label{DEF:@_G-O^*(G)}
@_G u=\int_{G_{\R}} u(t^{-1})\cdot \delta_t \;\d t, \quad u,v \in \mathcal C(G),\mathcal
E(G)
 \eeq
 \beq\label{DEF:@_G-O^*(G)-1}
 @_G(u,v)=\int_{G_{\R}}
u(t^{-1})\cdot  v (t) \;\d t,\quad u,v \in \mathcal C(G),\mathcal
E(G)
 \eeq
где $\delta_t(u)=u(t)$ -- $\delta$-функционал, а интегрирование
ведется по нормированной мере Хаара: $\mu(G_{\R})=1$.

Это не единственный пример отражения в $\mathcal O^\star(G)$ (и в $\mathcal P^\star(G)$), но типичный. Такое отражение называется {\it нормальным отражением} в $\mathcal O^\star(G)$ (и в $\mathcal P^\star(G)$).

По аналогии со свойствами на с.\pageref{@_G-widetilde(u)} доказываются

 \vglue10pt
\centerline{\bf Свойства нормальных отражений} 
\centerline{в $\mathcal O^\star(G)$ (и в $\mathcal P^\star(G)$):}
 \vglue10pt

\bit{\it
\item[$1^\circ$.] Перестановочность с антиподом:
\beq\label{@_G-widetilde(u)-O}
@_G\widetilde{u} =\widetilde{@_G u},\quad  u\in {\mathcal O}(G)  
\eeq
\item[$2^\circ$.] Замена на антипод при свертке:
\beq\label{(@_Gu)*v-O}
(@_G u) * v  =\widetilde{u}* v, \quad u,v\in {\mathcal O}(G)  
\eeq
\item[$3^\circ$.] Перестановочность со сверткой функций:
\beq \label{@_G(u*v)-O}
@_G (u * v)  =@_G v* @_Gu, \quad u,v\in {\mathcal O}(G) 
\eeq

\item[$4^\circ$.] Замена на антипод при внесении функционала под знак отражения:
\beq\label{alpha*(@_Gu)-O}
\alpha * @_G u   =@_G(u*\widetilde{\alpha}), \quad u\in {\mathcal O}(G),\ \alpha\in {\mathcal O}^\star(G)
\eeq

}\eit

\end{multicols}\noindent\rule[10pt]{160mm}{0.1pt}

\paragraph{След.}

Напомним, что след на алгебре операторов ${\mathcal L}(X)$ на конечномерном пространстве $X$ над полем $\C$ определяется как линейный функционал $\tr:{\mathcal L}(X)\to\C$ со свойством
\beq\label{tr(f-odot-x)=f(x)}
\tr(f\odot x)=f(x),\qquad x\in X,\ f\in X^\star.
\eeq
Известно, что он удовлетворяет тождеству
\beq\label{tr(A-circ-B)=tr(B-circ-A)}
\tr (A\circ B)=\tr(B\circ A),\qquad A,B\in{\mathcal L}(X).
\eeq

Следующее тождество связывает след с отражением $@_X$ на алгебре ${\mathcal L}(X)$:
\beq\label{tr(@_Xu-circ-A)=u(A)}
\tr (@_Xu\circ A)=u(A)=\tr(A\circ @_Xu),\qquad A\in{\mathcal L}(X),\ u\in {\mathcal L}^\star(X).
\eeq
\bpr
Это достаточно проверить для функционалов вида $u=f\circledast x$:
\begin{multline*}
\tr (@_Xu\circ A)=\tr (@_X(f\circledast x)\circ A)=
\tr ((f\odot x)\circ A)= \tr ((f\circ A)\odot x)=\\=\eqref{tr(f-odot-x)=f(x)}=(f\circ A)(x)=f(A(x))=
(f\circledast x)(A)=u(A)
\end{multline*}
Это доказывает первое равенство в \eqref{tr(@_Xu-circ-A)=u(A)}. Второе доказывается аналогично.
\epr

 \blm\label{LM:pi-A-sig} Справедливы тождества:
 \begin{align}\label{pi-A-sig}
 \int_K \pi(t)\circ A \circ \sigma(t^{-1})\;\d t &=
 \left\{\begin{array}{cl}
 \frac{\tr A}{\dim X_\sigma} \cdot \id_{X_\sigma}, & \pi=\sigma \\
 0, & \pi\ne \sigma \end{array}\right\},
 &&
 \begin{matrix}
 \pi,\sigma\in \Sigma(K),\\ A:X_\sigma\to X_\pi
 \end{matrix}
\\ \label{u-pi-A-v-sig}
 \int_K u\left(\pi(t)\right)\cdot v\left(\sigma(t^{-1})\right)\;\d t &=
 \left\{\begin{array}{cl}
 \frac{1}{\dim X_\sigma} \cdot @_{X_\sigma}(u,v), & \pi=\sigma \\
 0, & \pi\ne \sigma \end{array}\right\},
 &&
 \begin{matrix}
\pi,\sigma\in \Sigma(K),\\ u\in \mathcal L^\star(X_\pi),\ v\in \mathcal L^\star(X_\sigma)
\end{matrix}
 \\ \label{u-pi-sig}
 \int_K u\left(\pi(t)\right)\cdot \sigma(t^{-1})\;\d t &=
 \left\{\begin{array}{cl}
 \frac{1}{\dim X_\sigma} \cdot @_{X_\sigma}u, & \pi=\sigma \\
 0, & \pi\ne \sigma \end{array}\right\},
 &&
 \begin{matrix}
\pi,\sigma\in \Sigma(K),\\ u\in \mathcal L^\star (X_\pi)
\end{matrix}
 \end{align}
\elm
\bpr
1. Обозначим $B=\int_K \pi(t)\circ A \circ \sigma(t^{-1})\;\d t$.
Тогда для всякого $s\in K$ получим
\begin{multline*}
\pi(s)\circ B=\pi(s)\circ \int_K \pi(t)\circ A \circ \sigma(t^{-1})\;\d t=\int_K
\pi(s)\circ \pi(t)\circ A \circ \sigma(t^{-1})\;\d t=\\=\int_K \pi(s\cdot t)\circ A \circ \sigma(t^{-1})\;\d t=
\begin{pmatrix}st=\tau \\ t=s^{-1}\tau \\
t^{-1}=\tau^{-1}s\end{pmatrix}=\int_K \pi(\tau)\circ A \circ
\sigma(\tau^{-1}\cdot s)\;\d\tau=\\=\int_K \pi(\tau)\circ A \circ
\sigma(\tau^{-1})\circ\sigma(s)\;\d\tau=\int_K \pi(\tau)\circ A \circ
\sigma(\tau^{-1})\;\d\tau\circ\sigma(s)=B\circ\sigma(s)
\end{multline*}
то есть, $B$ -- сплетающий оператор между двумя неприводимыми
представлениями. Значит, если $\pi\ne\sigma$, то есть, $\pi$ и
$\sigma$ не эквивалентны, то, по лемме Шура, $B=0$. Если же
$\pi=\sigma$, то $B$ должен быть кратен тождественному оператору, то
есть
\beq\label{PROOF:pi-A-sig}
\int_K \sigma(t)\circ A \circ \sigma(t^{-1})\;\d t=\lambda\cdot\id_{X_\sigma}
\eeq
при некотором $\lambda\in\Bbb C$. Чтобы вычислить этот скаляр,
применим функционал следа к левой и правой частям:
\begin{multline*}
\tr A=\tr \int_K A\;\d t=\int_K \tr A\;\d t= \int_K \tr\left[A \circ
\sigma(t^{-1})\circ\sigma(t)\right]\;\d t= \int_K
\tr\left[\sigma(t)\circ A \circ \sigma(t^{-1})\right]\;\d t=\\=\tr \int_K\sigma(t)\circ A \circ \sigma(t^{-1})\;\d t=\eqref{PROOF:pi-A-sig}=\lambda\cdot\tr\left[\id_{X_\sigma}\right] =\lambda\cdot\dim
X_\sigma
\end{multline*}
Отсюда $\lambda=\frac{\tr A}{\dim X_\sigma}$, и мы доказали
\eqref{pi-A-sig}.

2. Рассмотрим в качестве $u,v$ одномерные функционалы:
$u=f\circledast x$, $v=g\circledast y$, где $x\in X_\pi$, $f\in
X_\pi^\star$, $y\in X_\sigma$, $g\in X_\sigma^\star$. Тогда
 \begin{multline}\nonumber
 \int_K u\left(\pi(t)\right)\cdot v\left(\sigma(t^{-1})\right)\;\d t =
 \int_K f\left(\pi(t)x\right)\cdot g\left(\sigma(t^{-1})y\right)\;\d t
 =\\=
  \int_K f\left(g\left(\sigma(t^{-1})y\right)\cdot\pi(t)x\right)\;\d t
  =
    \int_K f\left(\Big(g\odot\pi(t)x\Big)\left(\sigma(t^{-1})y\right)\right)\;\d t
 =\\=
    \int_K f\left(\Big(\pi(t)\circ g\odot x\Big)\left(\sigma(t^{-1})y\right)\right)\;\d t=
    \int_K f\left(\Big(\pi(t)\circ g\odot x\circ \sigma(t^{-1})\Big) y\right)\;\d t=
\eqref{tr(f-odot-x)=f(x)}=\\=
    \int_K \tr \left(f\odot \Big(\pi(t)\circ g\odot x\circ \sigma(t^{-1})\Big) y\right)\;\d t=
    \int_K \tr \Big(\pi(t)\circ g\odot x\circ \sigma(t^{-1})\circ f\odot y\Big)\;\d t
 =\\=
\tr\left[\int_K \left\{\pi(t)\circ g\odot x \circ \sigma(t^{-1})\circ
f\odot y \right\}\;\d t \right]
 =
\tr\left[\left\{\int_K \pi(t)\circ g\odot x \circ \sigma(t^{-1}) \;\d t \right\} \circ f\odot y  \right] =\\=\eqref{pi-A-sig}=
 \left\{\begin{array}{cl}
 \tr\left[\frac{\tr g\odot x}{\dim X_\sigma} \cdot f\odot y\right], & \pi=\sigma \\
 0, & \pi\ne \sigma \end{array}\right\}
 =
 \left\{ \begin{array}{cl}
 \frac{1}{\dim X_\sigma} \cdot g(x)\cdot f(y), & \pi=\sigma \\
 0, & \pi\ne \sigma \end{array}\right\}
 =\\=
 \left\{ \begin{array}{cl}
\frac{1}{\dim X_\sigma} \cdot @_{X_\sigma} (f\circledast x,
g\circledast y), & \pi=\sigma \\ 0, & \pi\ne \sigma
\end{array}\right\}=
 \left\{ \begin{array}{cl}
\frac{1}{\dim X_\sigma} \cdot @_{X_\sigma} (u, v), & \pi=\sigma
\\ 0, & \pi\ne \sigma
\end{array}\right\}
 \end{multline}
Таким образом, формула \eqref{u-pi-A-v-sig} справед\-лива для
одномерных
 функционалов $u,v$. Значит, она справедлива и для всех остальных.

3. Выбросив из \eqref{u-pi-A-v-sig} аргумент $v$, сразу же получаем
\eqref{u-pi-sig}.
\epr

\paragraph{Связь между отражениями и неприводимыми представлениями.}

\btm[о неприводимых представлениях]\label{ex-14.22}
Пусть $\pi :G\to \mathcal{L}(X)$ --- непрерывное представление компактной группы
$G$ в конечномерном пространстве $X$, и $\ph : \mathcal{C}^\star (G)\to
\mathcal{L}(X)$ --- соответствующий морфизм стереотипных алгебр \cite[Theorem 5.1.22]{Akbarov-De-Gruyter-I} 
$$
\begin{diagram}
\node[2]{\mathcal{L}(X)}
\\
\node{G} \arrow[2]{e,t}{\delta} \arrow{ne,t}{\pi}
\node[2]{\mathcal{C}^\star(G)} \arrow{nw,t}{\varphi}
\end{diagram}
$$
Пусть $@_G : \mathcal{C}(G)\to \mathcal{C}^\star (G)$ --- нормальное отражение
на алгебре $\mathcal{C}^\star (G)$ (определенное формулой \eqref{DEF:@_G-C^*(G)}), и пусть $@_X : \mathcal{L}^\star (X)\to \mathcal{L}(X)$ --- нормальное отражение на алгебре $\mathcal{L}(X)$ (то есть преобразование Гротендика из \eqref{eq14.3}). Тогда следующие
условия эквивалентны:
\begin{itemize}
\item[(i)] $\pi$ --- неприводимое представление;

\item[(ii)] для какого-нибудь ненулевого отражения $@$ на $\mathcal{L}(X)$ отображение $\ph$ является
морфизмом алгебр с отражением $\ph :(\mathcal{C}^\star (G), @_G) \to (\mathcal{L}(X), @)$;

\item[(iii)] для некоторой ненулевой константы $\lambda \in \C$ отображение $\ph$ является морфизмом алгебр с отражением $\ph:(\mathcal{C}^\star (G), @_G) \to (\mathcal{L}(X), \lambda @_X)$.
\end{itemize}
При этом,
\beq
\ph \circ @_G \circ \ph^\star = \frac{1}{\dim X} \cdot @_X \label{@_G<->@_X}
\eeq
где $\dim X$ --- размерность пространства $X$. В частности,
\begin{itemize}
\item[(a)] $\pi$ --- одномерное представление $\Leftrightarrow$ $\ph$ --- морфизм
алгебр с отражением $\ph :(\mathcal{C}^\star (G), @_G) \to (\mathcal{L}(X),
@_X)$.
\end{itemize}
\etm
\bpr
Здесь используется следующая формула
\beq\label{eq14.24}
(\ph \circ @_G \circ \ph^\star)(u)= \int_G u(\pi(g)) \cdot \pi (g^{-1}) \cdot
\d g \quad\quad\quad u \in \mathcal{L}^\star(X)
\eeq
которая доказывается такой цепочкой: если $u \in \mathcal  L^\star(X)$, то
$$
\ph^\star(u)\in \mathcal{C}(G) :\quad \ph^\star(u)(g)=(u\circ
\ph)(g)=u(\ph(g))=u(\pi(g))
$$
$$
\Downarrow
$$
$$
@_G \ph^\star (u)\in \mathcal{C}^\star (G) :\quad @_G \ph^\star (u)= \int_G
\ph^\star (u)(g) \cdot \delta^{g^{-1}} \cdot \d g= \int_G u(\pi(g)) \cdot
\delta^{g^{-1}} \cdot \d g
$$
$$
\Downarrow
$$
$$
\ph @_G \ph^\star (u)\in \mathcal{L}(X) :\quad \ph @_G \ph^\star (u)= \int_G
u(\pi(g)) \cdot \ph(\delta^{g^{-1}}) \cdot \d g= \int_G u(\pi(g)) \cdot
\pi(g^{-1}) \cdot \d g
$$
Из \eqref{eq14.24} мы получаем \eqref{@_G<->@_X}:
$$
(\ph \circ @_G \circ \ph^\star)(u)=\eqref{eq14.24}= \int_G u(\pi(g)) \cdot \pi (g^{-1}) \cdot
\d g=\eqref{u-pi-sig}=\frac{1}{\dim X} \cdot @_X u
$$
Мы доказали, импликацию $(i)\Rightarrow (iii)$. Импликация $(iii)\Rightarrow (ii)$
очевидна. Докажем $(ii)\Rightarrow (i)$. Пусть $\mathcal{L}(X)$ наделено
некоторым ненулевым отражением $@$, и $\ph$ является морфизмом алгебр с
отражением $\ph :(\mathcal{C}^\star (G), @_G) \to (\mathcal{L}(X), @)$, то есть
коммутативна диаграмма
$$
\begin{diagram}
\node{\mathcal{C}(G)} \arrow{e,t}{\AT_G} \node{\mathcal{C}^\star(G)}
\arrow{s,r}{\varphi}
\\
\node{\mathcal{L}^\star(X)} \arrow{e,t}{\AT} \arrow{n,l}{\varphi^\star}
\node{\mathcal{L}(X)}
\end{diagram}.
$$
В силу примера \ref{ex-14.2}, $@= \lambda \cdot @_X$, $\lambda \ne 0$.
Поэтому $\mathcal{L}(X)=\Im @ \subseteq \Im \ph$, откуда следует, что $\pi$ --- неприводимое представление. \epr

\paragraph{Нормированные характеры $\chi_\sigma$.}

Пусть $\Sigma(K)$ -- двойственный объект к компактной группе $K$, то
есть полная система унитарных неприводимых представлений $\sigma :K\to
\mathcal L(X_\sigma)$ $(\dim X_\sigma<\infty)$. Для любого $\sigma
\in \Sigma(K)$ обозначим через $\chi_\sigma$ ``нормированный'' характер представления
$\sigma$:
\beq\label{DEF:chi_sigma}
\chi_\sigma (t)=\dim X_\sigma\cdot \tr \sigma (t).
\eeq

\blm\label{LM:chi-sig} Справедливы тождества
 \begin{align}
 \label{pi-chi}
\int_K \chi_\pi(t)\cdot
 \sigma(t^{-1})\;\d t &=
 \left\{\begin{array}{cl}
 \id_{X_\sigma}, & \pi=\sigma \\
 0, & \pi\ne \sigma \end{array}\right\},
 &&
\\ \label{chi-chi}
 \int_K \chi_\pi(t)\cdot\chi_\sigma(t^{-1})\;\d t &=
 \left\{\begin{array}{cl}
 \big(\dim X_\sigma\big)^2, & \pi=\sigma \\
 0, & \pi\ne \sigma \end{array}\right\}.
 &&
 \end{align}
 \elm
 \bpr 1. Подставив $u=\dim_{X_\sigma}\cdot\tr$ в \eqref{u-pi-sig}, получим \eqref{pi-chi}:
 \begin{multline*}
 \int_K \chi_\pi(t)\cdot
 \sigma(t^{-1})\;\d t=
 \int_K \dim_{X_\sigma}\cdot\tr\left[\pi(t)\right]\cdot \sigma(t^{-1})\;\d t =\\=
 \left\{\begin{array}{cl}
 \frac{\dim_{X_\sigma}}{\dim X_\sigma} \cdot @_{X_\sigma}\tr, & \pi=\sigma \\
 0, & \pi\ne \sigma \end{array}\right\}=
 \left\{\begin{array}{cl}
 \id_{X_\sigma}, & \pi=\sigma \\
 0, & \pi\ne \sigma \end{array}\right\}
 \end{multline*}

2. Формула \eqref{chi-chi} (вариант соотношений ортогональности для характеров) получается применением к обеим частям \eqref{pi-chi} функционала $u=\dim_{X_\sigma}\cdot\tr$.
\epr

 \vglue10pt
\centerline{\bf Свойства характеров $\chi_\sigma$:}
 \vglue10pt

\bit{\it

\item[$1^\circ$.] Унитарность:
\beq\label{widehat(chi_sigma)=overline(chi_sigma)}
\chi_\sigma(t^{-1})=\overline{\chi_\sigma(t)},\qquad t\in G.
\eeq

\item[$2^\circ$.] Центральность относительно свертки:
\beq\label{chi_sigma*u=u*chi_sigma}
u*\chi_\sigma=\chi_\sigma*u,\qquad u\in {\mathcal C}(K).
\eeq

\item[$3^\circ$.] Ортогональность относительно свертки:
\beq\label{chi_pi*chi_sigma}
\chi_\pi*\chi_\sigma =\begin{cases} \chi_\sigma,& \pi = \sigma
\\ 0, &\pi\ne\sigma \end{cases},\qquad \pi,\sigma\in \widehat{K}.
\eeq

}\eit

\bpr

1. Тождество \eqref{widehat(chi_sigma)=overline(chi_sigma)} доказывается цепочкой
\begin{multline*}
\chi_\sigma (t^{-1})=\dim X_\sigma\cdot\tr \sigma (t^{-1})=\dim X_\sigma\cdot\tr\big( \sigma (t)^{-1}\big)=
\dim X_\sigma\cdot\tr \big(\sigma(t)^\bullet\big)=\\=
\dim X_\sigma\cdot\tr \big(\overline{\sigma(t)}^T\big)=\dim X_\sigma\cdot\tr \overline{\sigma(t)} =\overline{\dim X_\sigma\cdot\tr\sigma(t)}=\overline{\chi_\sigma (t)}.
\end{multline*}

2. Для любого $u\in {\mathcal C}(K)$ мы получим
\begin{multline*}
u*\chi_\sigma(s)=\eqref{DEF:svertka-funktsij}=\int_K u(t)\cdot\chi_\sigma(t^{-1}\cdot s)\d t=
\int_K u(t)\cdot\dim X_\sigma\cdot\tr \Big(\sigma(t^{-1}\cdot s)\Big)\d t=\\=
\dim X_\sigma\cdot\tr \Big(\int_K u(t)\cdot\sigma(t^{-1}\cdot s)\d t\Big)=
\dim X_\sigma\cdot\tr \Big(\int_K u(t)\cdot\sigma(t^{-1})\circ\sigma(s)\d t\Big)=\\=
\dim X_\sigma\cdot\tr \Big(\int_K u(t)\cdot\sigma(t^{-1})\d t\circ\sigma(s)\Big)=\eqref{tr(A-circ-B)=tr(B-circ-A)}=
\dim X_\sigma\cdot\tr \Big(\sigma(s)\circ\int_K u(t)\cdot\sigma(t^{-1})\d t\Big)=\\=
\dim X_\sigma\cdot\tr \Big(\int_K\sigma(s)\circ u(t)\cdot\sigma(t^{-1})\d t\Big)=
\dim X_\sigma\cdot\tr \Big(\int_K u(t)\cdot\sigma(s\cdot t^{-1})\d t\Big)=\\=
\int_K u(t)\cdot\dim X_\sigma\cdot\tr \Big(s\cdot\sigma(t^{-1})\Big)\d t=
\int_K u(t)\cdot\chi_\sigma(s\cdot t^{-1})\d t=\eqref{DEF:svertka-funktsij}=\chi_\sigma*u(s).
\end{multline*}

3. Формула \eqref{chi_pi*chi_sigma} следует из \eqref{pi-chi}:
\begin{multline*}
\chi_\pi*\chi_\sigma(s)=\int_K \chi_\pi(t)\cdot\chi_\sigma(t^{-1}\cdot s)\d t=
\int_K \chi_\pi(t)\cdot \dim_{X_\sigma}\cdot \tr \sigma(t^{-1}\cdot s)\d t=\\=
\dim_{X_\sigma}\cdot\tr \bigg[\int_K \chi_\pi(t)\cdot \sigma(t^{-1}\cdot s)\d t\bigg]=
\dim_{X_\sigma}\cdot\tr \bigg[\int_K \chi_\pi(t)\cdot \sigma(t^{-1})\circ \sigma(s)\d t\bigg]=\\=
\dim_{X_\sigma}\cdot\tr \bigg[\int_K \chi_\pi(t)\cdot \sigma(t^{-1})\d t \circ\ \sigma(s)\bigg]=\eqref{pi-chi}=\\=
\left\{\begin{array}{cl}
\dim_{X_\sigma}\cdot \tr\Big[\id_{X_\sigma}\circ \ \sigma(s)\Big], & \pi=\sigma \\
 0, & \pi\ne \sigma \end{array}\right\}
 =
  \left\{\begin{array}{cl}
 \dim X_\sigma \cdot\tr \sigma(s), & \pi=\sigma \\
 0, & \pi\ne \sigma \end{array}\right\}
 =
 \left\{\begin{array}{cl}
 \chi_\sigma(s), & \pi=\sigma \\
 0, & \pi\ne \sigma \end{array}\right\}
\end{multline*}
 \epr

Обозначим символом ${\mathcal K}_\sigma(K)$ пространство функций из ${\mathcal C}(K)$, не меняющихся при свертке с $\chi_\sigma$:
\beq\label{DEF:K_sigma(K)}
{\mathcal K}_\sigma(K)=\{u\in{\mathcal C}(K): u*\chi_\sigma=u\}.
\eeq
Очевидно, это замкнутое подпространство в ${\mathcal C}(K)$.

 \vglue10pt
\centerline{\bf Свойства пространства ${\mathcal K}_\sigma(K)$:}
 \vglue10pt

\bit{\it

\item[$1^\circ$.] ${\mathcal K}_\sigma(K)$ является образом ${\mathcal C}(K)$ при отображении $u\mapsto u*\chi_\sigma$:
\beq\label{C(K)*chi_sigma-subseteq-K_sigma(K)}
{\mathcal C}(K)*\chi_\sigma={\mathcal K}_\sigma(K)
\eeq

\item[$2^\circ$.] ${\mathcal K}_\sigma(K)$ является образом ${\mathcal C}^\star(K)$ при отображении $\alpha\mapsto \alpha*\chi_\sigma$:
\beq\label{C(K)^star*chi_sigma-subseteq-K_sigma(K)}
{\mathcal C}^\star(K)*\chi_\sigma={\mathcal K}_\sigma(K)
\eeq

\item[$3^\circ$.] ${\mathcal K}_\sigma(K)$ является двусторонним идеалом в ${\mathcal C}(K)$ относительно свертки и
алгеброй относительно свертки с единицей $\chi_\sigma$.

\item[$4^\circ$.] ${\mathcal K}_\sigma(K)$ состоит в точности из функций вида $w\circ \sigma$, где $w\in {\mathcal L}^\star(X_\sigma)$, а отображение $w\mapsto w\circ \sigma$ устанавливает изоморфизм между пространствами ${\mathcal L}^\star(X_\sigma)$ и ${\mathcal K}_\sigma(K)$:
    \beq\label{K_sigma(K)=L^star(X_sigma)-circ-sigma}
    {\mathcal K}_\sigma(K)={\mathcal L}^\star(X_\sigma)\circ\sigma\cong {\mathcal L}^\star(X_\sigma).
    \eeq
    Как следствие, алгебра ${\mathcal K}_\sigma(K)$ конечномерна, и ее размерность равна
\beq\label{dim-K_sigma(K)}
\dim {\mathcal K}_\sigma(K)=\dim {\mathcal L}^\star(X_\sigma)=(\dim X_\sigma)^2.
\eeq

}\eit

\bpr
1. Если  $u\in {\mathcal C}(K)*\chi_\sigma$, то есть $u=v*\chi_\sigma$ для некоторой функции $v\in {\mathcal C}(K)$, то
$$
u=v*\chi_\sigma=\eqref{chi_pi*chi_sigma}=v*\chi_\sigma*\chi_\sigma=u*\chi_\sigma
$$
и поэтому $u\in {\mathcal K}_\sigma(K)$. Наоборот, если $u\in {\mathcal K}_\sigma(K)$, то
$$
u=u*\chi_\sigma,
$$
и поэтому $u\in {\mathcal C}(K)*\chi_\sigma$.

2. Пусть $\alpha\in{\mathcal C}^\star(K)$. Выберем направленность $u_i\in {\mathcal C}(K)$ так, чтобы
$$
@_Gu_i\overset{{\mathcal C}^\star(K)}{\underset{i\to\infty}{\longrightarrow}}\alpha
$$
(это всегда можно сделать, потому что отражение $@_G$ точное). Поскольку
отображение $(\alpha,u)\in{\mathcal C}^\star(K)\times {\mathcal C}(K) \mapsto\alpha* u\in {\mathcal C}(K)$ является непрерывной блиненейной формой, при свертке с $\chi_\sigma$ мы получаем:
$$
\underbrace{@_Gu_i*\chi_\sigma}_{\scriptsize\begin{matrix}\|\put(2,0){\eqref{(@_Gu)*v}}\\
\widetilde{u_i}*\chi_\sigma \\ \text{\rotatebox{90}{$\owns$}}\\
{\mathcal K}_\sigma(K)
\end{matrix}}\overset{{\mathcal C}(K)}{\underset{i\to\infty}{\longrightarrow}}\alpha*\chi_\sigma
$$
Поскольку ${\mathcal K}_\sigma(K)$ -- замкнутое подпространство в ${\mathcal C}^\star(K)$, мы получаем, что
$\alpha*\chi_\sigma\in {\mathcal K}_\sigma(K)$.

3. Если $u\in {\mathcal C}(K)$ и $v\in{\mathcal K}_\sigma(K)$, то
$$
(u*v)*\chi_\sigma=u*(v*\chi_\sigma)=\eqref{DEF:K_sigma(K)}=u*v,
$$
и поэтому $u*v\in {\mathcal K}_\sigma(K)$. И, в силу \eqref{chi_sigma*u=u*chi_sigma}, наоборот, $v*u\in {\mathcal K}_\sigma(K)$. Равенство \eqref{chi_pi*chi_sigma} при $\pi=\sigma$ имеет вид $\chi_\sigma*\chi_\sigma=\chi_\sigma$, и это влечет включение $\chi_\sigma\in{\mathcal K}_\sigma(K)$. С другой стороны, тождество $u*\chi_\sigma=u$ для элементов $u\in {\mathcal K}_\sigma(K)$ (которым определяется ${\mathcal K}_\sigma(K)$ в \eqref{DEF:K_sigma(K)}) означает, что $\chi_\sigma$ является единицей в ${\mathcal K}_\sigma(K)$ относительно свертки.

4. Пусть $u\in{\mathcal L}^\star(X_\sigma)\circ\sigma$, то есть $u=w\circ\sigma$ для некоторого $w\in {\mathcal L}^\star(X_\sigma)$. Тогда
\begin{multline*}
u*\chi_\sigma(s)=\int_K u(t)\cdot\chi_\sigma(t^{-1}\cdot s)\d t=\int_K w(\sigma(t))\cdot\dim X_\sigma\cdot\tr \Big(\sigma(t^{-1}\cdot s)\Big)\d t=\\=
\dim X_\sigma\cdot\tr \Big(\int_K w(\sigma(t))\cdot\sigma(t^{-1}\cdot s)\d t\Big)=
\dim X_\sigma\cdot\tr \Big(\int_K w(\sigma(t))\cdot\sigma(t^{-1})\circ\sigma(s)\d t\Big)=\\=
\dim X_\sigma\cdot\tr \Big(\int_K w(\sigma(t))\cdot\sigma(t^{-1})\d t\circ\sigma(s)\Big)=\eqref{u-pi-sig}=
\dim X_\sigma\cdot\tr \Big(\frac{1}{\dim X_\sigma}\cdot @_{X_\sigma}w\circ\sigma(s)\Big)=\\=\tr \Big(@_{X_\sigma}w\circ\sigma(s)\Big)=\eqref{tr(@_Xu-circ-A)=u(A)}=w(\sigma(s))=u(s),
\end{multline*}
поэтому $u\in {\mathcal K}_\sigma(K)$.

Наоборот, пусть $u\in {\mathcal K}_\sigma(K)$. Обозначим
$$
A=\int_K u(t)\cdot\sigma(t^{-1})\d t.
$$
Это будет оператор в конечномерном пространстве $X_\sigma$, поэтому у него должен быть прообраз при отражении $@_{X_\sigma}$:
$$
w=@_{X_\sigma}^{-1}A.
$$
Мы получим:
\begin{multline*}
u*\chi_\sigma(s)=\int_K u(t)\cdot\chi_\sigma(t^{-1}\cdot s)\d t=
\int_K u(t)\cdot\tr \Big(\sigma(t^{-1}\cdot s)\Big)\d t=
\tr \Big(\int_K u(t)\cdot\sigma(t^{-1}\cdot s)\d t\Big)=\\=
\tr \Big(\int_K u(t)\cdot\sigma(t^{-1})\cdot\sigma(s)\d t\Big)=
\tr \Big(\int_K u(t)\cdot\sigma(t^{-1})\d t\cdot\sigma(s)\Big)=\tr(A\circ\sigma(s))=\\=
\eqref{tr(@_Xu-circ-A)=u(A)}=@_{X_\sigma}^{-1}A(\sigma(s))=w(\sigma(s)).
\end{multline*}
То есть $u=w\circ\sigma$.

Мы доказали первое равенство в \eqref{K_sigma(K)=L^star(X_sigma)-circ-sigma}:
$$
{\mathcal K}_\sigma(K)={\mathcal L}^\star(X_\sigma)\circ\sigma.
$$
Его можно понимать так, что образ пространства ${\mathcal L}^\star(X_\sigma)$ при отображении $\dot{\sigma}^\star: {\mathcal C}(K)\gets {\mathcal L}^\star(X_\sigma)$, сопряженном к гомоморфизму алгебр $\dot{\sigma}: {\mathcal C}^\star(K)\to {\mathcal L}(X_\sigma)$, продолжающему представление $\sigma: K\to {\mathcal L}(X_\sigma)$, совпадает с ${\mathcal K}_\sigma(K)$:
$$
{\mathcal K}_\sigma(K)=\dot{\sigma}^\star({\mathcal L}^\star(X_\sigma)).
$$
Поскольку представление $\sigma$ неприводимо, и значит, отображение $\dot{\sigma}$ сюръективно, сопряженное отображение $\dot{\sigma}^\star$ будет инъективно. Поэтому оно является биекцией между ${\mathcal L}^\star(X_\sigma)$ и ${\mathcal K}_\sigma(K)$. С другой стороны, пространство ${\mathcal L}(X_\sigma)$ конечномерно, поэтому эта биекция является одновременно изоморфизмом стеретипных пространств. Одновременно это доказывает \eqref{dim-K_sigma(K)}.
\epr

\paragraph{Меры $\nu_\sigma$ и алгебры ${\mathcal C}_\sigma^\star(K)$.}

Пусть $K$ -- компактная группа, $\mu_K$ -- нормированная мера Хаара на $K$. Для всякого представления $\sigma\in\widehat{K}$, $\sigma:K\to {\mathcal B}(X_\sigma)$, рассмотрим нормированный характер $\chi_\sigma$ из \eqref{DEF:chi_sigma} и  определим меру $\nu_\sigma\in{\mathcal C}^\star(K)$ формулой
\beq
\nu_\sigma(u)=\int_K \overline{\chi_\sigma(s)}\cdot u(s)\cdot \mu_K(\d s)=\int_K \chi_\sigma(s^{-1})\cdot u(s)\cdot \mu_K(\d s),\qquad u\in{\mathcal C}(G),
\eeq
или, эквивалентно, формулой
\beq
\nu_\sigma=\int_K \overline{\chi_\sigma(s)}\cdot\delta_s\cdot \mu_K(\d s)= \int_K \chi_\sigma(s^{-1})\cdot\delta_s\cdot \mu_K(\d s).
\eeq
или формулой
\beq\label{nu_sigma=@_G-chi_sigma-1}
\nu_\sigma=\widetilde{\chi_\sigma}\cdot\mu_K
\eeq
или формулой
\beq\label{nu_sigma=@_G-chi_sigma}
\nu_\sigma=@_G\chi_\sigma
\eeq

 \vglue10pt
\centerline{\bf Свойства мер $\nu_\sigma$:}
 \vglue10pt

\bit{\it

\item[$1^\circ$.] Для всякого $\sigma\in\widehat{K}$ мера $\nu_\sigma$ -- центральная в ${\mathcal C}^\star(K)$:
 \beq\label{chi_pi-centr}
 \nu_\sigma*\alpha=\alpha*\nu_\sigma,\qquad \alpha\in {\mathcal C}^\star(K),
 \eeq

\item[$2^\circ$.] Меры $\nu_\sigma$ образуют систему ортогональных проекторов в алгебре ${\mathcal C}^\star(K)$:
 \beq\label{chi_pi-ortog-proj}
 \nu_\pi*\nu_\sigma=\begin{cases}\nu_\sigma,& \pi=\sigma\\ 0,& \pi\ne\sigma\end{cases}\qquad \pi,\sigma\in \widehat{K}.
 \eeq

}\eit

Рассмотрим пространство
\beq\label{DEF:C_sigma^star(K)}
{\mathcal C}^\star_\sigma(K)=\{\alpha\in {\mathcal C}^\star(K):\ \nu_\sigma*\alpha=\alpha \}
\eeq
Очевидно, оно замкнуто в ${\mathcal C}^\star(K)$.

 \vglue10pt
\centerline{\bf Свойства пространства ${\mathcal C}_\sigma^\star(K)$:}
 \vglue10pt

\bit{\it

\item[$1^\circ$.] ${\mathcal C}_\sigma^\star(K)$ является образом отображения $\alpha\mapsto\nu_\sigma*\alpha$:
\beq\label{nu_sigma*C^star(K)=C^star_sigma(K)}
\nu_\sigma* {\mathcal C}^\star(K)={\mathcal C}^\star_\sigma(K).
\eeq

\item[$2^\circ$.] ${\mathcal C}_\sigma^\star(K)$ является двусторонним идеалом в ${\mathcal C}^\star(K)$ и алгеброй относительно свертки с единицей $\nu_\sigma$.

\item[$3^\circ$.] Отражение $@_G$ биективно отображает пространство ${\mathcal K}_\sigma(K)$ на пространство ${\mathcal C}_\sigma^\star(K)$:
\beq\label{@_G(K_sigma)=C_sigma^star(K)}
@_G{\mathcal K}_\sigma(K)={\mathcal C}_\sigma^\star(K)
\eeq
и является изоморфизмом алгебр
$$
{\mathcal K}_\sigma(K)\cong{\mathcal C}_\sigma^\star(G).
$$

\item[$4^\circ$.] Гомоморфизм $\dot{\sigma}:{\mathcal C}^\star(G)\to{\mathcal L}(X_\sigma)$ биективно отображает подалгебру ${\mathcal C}_\sigma^\star(G)$ на ${\mathcal L}(X_\sigma)$:
    \beq\label{dot-sigma-C_sigma^star(G)=L(X_sigma)}
    \dot{\sigma}\Big({\mathcal C}_\sigma^\star(G)\Big)={\mathcal L}(X_\sigma).
    \eeq

\item[$5^\circ$.] Если $\pi,\sigma\in\widehat{K}$ и $\pi\ne\sigma$, то
\beq\label{dot-pi-C_sigma^star(G)=0}
\dot{\sigma}\Big({\mathcal C}_\pi^\star(G)\Big)=0
\eeq

}\eit

\bpr
1. Если $\alpha\in {\mathcal C}^\star_\sigma(K)$, то $\alpha=\nu_\sigma*\alpha$, и поэтому $\alpha\in \nu_\sigma*{\mathcal C}^\star(K)$. Наоборот, если $\alpha\in \nu_\sigma*{\mathcal C}^\star(K)$, то есть $\alpha=\nu_\sigma*\beta$ для всякой меры $\beta\in {\mathcal C}^\star(K)$, то
$$
\nu_\sigma*\alpha=\nu_\sigma*\nu_\sigma*\beta=\nu_\sigma*\beta=\alpha
$$
поэтому $\alpha\in {\mathcal C}_\sigma^\star(K)$.

2. Если $\alpha\in {\mathcal C}^\star_\sigma(K)$ и $\beta\in {\mathcal C}^\star(K)$, то
$$
\alpha*\beta=\eqref{DEF:C_sigma^star(K)}=(\nu_\sigma*\alpha)*\beta=\nu_\sigma*(\alpha*\beta),
$$
и поэтому $\alpha*\beta\in {\mathcal C}^\star_\sigma(K)$. Точно так же в силу \eqref{chi_pi-centr},
$\beta*\alpha\in {\mathcal C}^\star_\sigma(K)$. С другой стороны, в силу \eqref{chi_pi-ortog-proj}, $\nu_\sigma*\nu_\sigma=\nu_\sigma$, поэтому $\nu_\sigma\in {\mathcal C}^\star_\sigma(K)$, а для любого $\alpha\in {\mathcal C}^\star_\sigma(K)$ мы получаем
$$
\nu_\sigma*\alpha=\eqref{DEF:C_sigma^star(K)}=\alpha
$$
то есть $\nu_\sigma$ -- единица в ${\mathcal C}^\star_\sigma(K)$.

3. Пусть $u\in {\mathcal K}_\sigma(K)$, то есть $u=f\circ\sigma$ для некоторого $f\in{\mathcal L}^\star(X_\sigma)$. Тогда
\begin{multline*}
\nu_\sigma*@u=\eqref{alpha*(@_Gu)}=
@(u*\widetilde{\nu_\sigma})=\eqref{nu_sigma=@_G-chi_sigma}=@(u*\widetilde{@\chi_\sigma})=\\=
\eqref{@_G-widetilde(u)}=
@(u*@\widetilde{\chi_\sigma})=\eqref{(@_Gu)*v}=@(u*\chi_\sigma)=\eqref{DEF:K_sigma(K)}=@u,
\end{multline*}
и значит, $@u\in {\mathcal C}_\sigma^\star(K)$. Пусть наоборот, $\alpha\in {\mathcal C}_\sigma^\star(K)$. Тогда
$$
\alpha=\nu_\sigma*\alpha=@\chi_\sigma*\alpha=@(\widetilde{\alpha}*\chi_\sigma)
$$
то есть $\alpha=@u$, где $u=\widetilde{\alpha}*\chi_\sigma\in {\mathcal K}_\sigma(K)$ в силу \eqref{C(K)^star*chi_sigma-subseteq-K_sigma(K)}. Остается отметить, что в силу \eqref{@_G(u*v)}, отображение $@_G$ мультипликативно,
$$
@_G(u*v)=@_Gu*@_Gv,
$$
а, в силу \eqref{nu_sigma=@_G-chi_sigma}, сохраняет единицу:
$$
@_G\chi_\sigma=\nu_\sigma.
$$

4. Пусть $A\in {\mathcal L}(X_\sigma)$. Положим
$$
w=@_G^{-1}(A),\quad u=w\circ\sigma,\quad \alpha=\dim X_\sigma\cdot @_G(u).
$$
Тогда
$$
w\in{\mathcal L}^\star(X_\sigma)\quad\overset{\eqref{K_sigma(K)=L^star(X_sigma)-circ-sigma}}{\Longrightarrow}\quad
u=w\circ\sigma\in {\mathcal K}_\sigma(K)
\quad\overset{\eqref{@_G(K_sigma)=C_sigma^star(K)}}{\Longrightarrow}\quad
\alpha=@_G(u)\in {\mathcal C}_\sigma^\star(K).
$$
При этом,
$$
\dot{\sigma}(\alpha)=\dim X_\sigma\cdot \dot{\sigma}(@_G(\dot{\sigma}^\star(w)))=
\dim X_\sigma\cdot (\dot{\sigma}\circ @_G\circ \dot{\sigma}^\star)(w)=\eqref{@_G<->@_X}=
\dim X_\sigma\cdot \frac{1}{\dim X_\sigma}@_Xw=@_Xw=A.
$$

5. Пусть $\pi,\sigma\in\widehat{K}$ и $\pi\ne\sigma$. Рассмотрим оператор $\dot{\sigma}(\nu_\pi)\in{\mathcal L}(X_\sigma)$. Его можно описать действием на нем произвольного функционала $w\in{\mathcal L}^\star(X_\sigma)$:
\begin{multline*}
w\Big(\dot{\sigma}(\nu_\pi)\Big)=\nu_\pi(w\circ\sigma)=@\chi_\pi(w\circ\sigma)=\int_K\chi_\pi(t)\cdot (w\circ\sigma)(t^{-1})\d t=\\=\int_K\chi_\pi(t)\cdot w\Big(\sigma(t^{-1})\Big)\d t=
w\Big(\int_K\chi_\pi(t)\cdot \sigma(t^{-1})\d t\Big)=\eqref{pi-chi}=0.
\end{multline*}
\epr

\paragraph{Гомоморфизмы ${\mathcal C}^\star(K)\to\prod_{i=1}^{k} {\mathcal L}(X_{\sigma_i})$.}

Пусть $\sigma_1,...,\sigma_k\in\widetilde{K}$ -- различные неприводимые представления. Рассмотрим их прямую сумму
$$
\prod_{i=1}^{k}\sigma_i: K\to \prod_{i=1}^{k} {\mathcal L}(X_{\sigma_i})
$$
и порожденный им гомоморфизм алгебр
\beq\label{(sigma_1-oplus...oplus-sigma_k)^cdot}
\l\prod_{i=1}^{k}\sigma_i\r^\cdot:{\mathcal C}^\star(K)\to \prod_{i=1}^{k} {\mathcal L}(X_{\sigma_i})
\eeq

\btm\label{TH:plotnost-C*(K)-v-prod-B(X_sigma)} Гомоморфизм \eqref{(sigma_1-oplus...oplus-sigma_k)^cdot} является сюръективным отображением.
\etm
\bpr
Это можно доказать с помощью теоремы Петера---Вейля, но нам это удобнее будет вывести из уже известных нам свойств пространств ${\mathcal C}_\sigma^\star(K)$.
 Зафиксируем операторы $A_i\in {\mathcal L}(X_{\sigma_i})$, $i=1,...,k$.
 Воспользуемся формулой \eqref{dot-sigma-C_sigma^star(G)=L(X_sigma)}, и подберем меры $\alpha_i\in{\mathcal C}_{\sigma_i}^\star(K)$ так, что
$$
\dot{\sigma_i}(\alpha_i)=A_i.
$$
Положим
$$
\alpha=\alpha_1+...+\alpha_k.
$$
В силу \eqref{dot-pi-C_sigma^star(G)=0},
$$
\dot{\sigma_i}(\alpha_j)=\begin{cases}A_i,& i=j\\ 0,& i\ne j\end{cases}.
$$
Поэтому
$$
\dot{\sigma_i}(\alpha)=\alpha_i,
$$
и отсюда
$$
\l\prod_{i=1}^{k}\sigma_i\r^\cdot(\alpha)=
\l\dot{\sigma_1}(\alpha),...,\dot{\sigma_k}(\alpha)\r=(A_1,...,A_k).
$$
\epr

\section{Усреднение по малым подгруппам в ${\mathcal C}^\star(G)$}

\subsection{Меры $\mu_K$ и операторы $\pi_K^\star$.}

\paragraph{Меры $\mu_K$.}

Пусть $G$ --- локально компактная группа.

Пусть $K\subseteq G$ --- компактная нормальная подгруппа в локально компактной группе $G$, и $\mu_K$ -- нормированная мера Хаара на $K$:
$$
\mu_K(K)=1.
$$
Помимо того, что $\mu_K$ есть мера на $K$, мы ее также будем рассматривать как меру на $G$, 
доопределив нулем вне $K$. Это можно понимать как определение в виде функционала на пространстве ${\mathcal C}(G)$, действующего по формуле
\beq\label{mu_K(u)=int_K-u(t)mu_K(dt)}
\mu_K(u)=\int_Ku(t)\ \mu_K(\d t),\qquad u\in {\mathcal C}(G).
\eeq

\medskip
\centerline{\bf Свойства мер $\mu_K$:}

\bit{\it

\item[$1^\circ$.] Носителем меры $\mu_K$ является множество $K$:
\beq\label{supp-mu_K=K}
\supp \mu_K=K
\eeq

\item[$2^\circ$.] Сдвиги на элементы $K$ не меняют $\mu_K$:
\beq\label{x-cdot-mu_K=mu_K}
x\cdot \mu_K=\mu_K=\mu_K\cdot x,\qquad x\in K
\eeq

\item[$3^\circ$.] Антипод не меняет $\mu_K$:
\beq\label{widetilde(mu_K)=mu_K}
\widetilde{\mu_K}=\mu_K
\eeq

\item[$4^\circ$.] Плоская инволюция не меняет $\mu_K$:
\beq\label{overline(mu_K)=mu_K}
\overline{\mu_K}=\mu_K
\eeq

\item[$5^\circ$.] Сопряжения в группе $G$ не меняют $\mu_K$:
\beq\label{t^(-1)-cdot-mu_K-cdot-t=mu_K}
t^{-1}\cdot \mu_K\cdot t=\mu_K,\qquad t\in G.
\eeq

\item[$6^\circ$.] Мера $\mu_K$ перестановочна со сдвигами
\beq\label{mu_K-cdot-t=t-cdot-mu_K}
\mu_K\cdot t=t\cdot \mu_K,\qquad t\in G
\eeq
и поэтому центральная:
\beq\label{mu_K*alpha=alpha*mu_K}
\mu_K*\alpha=\alpha*\mu_K,\qquad \alpha\in{\mathcal C}^\star(G)
\eeq

\item[$7^\circ$.] Свертка двух мер $\mu_K$ и $\mu_L$ дает меру $\mu_{K\cdot L}$\footnote{Здесь используется лемма \ref{LM:K-cdot-L-komp-norm-podgr}.}:
\beq\label{mu_K*mu_L=mu_(KL)}
\mu_K*\mu_L=\mu_{K\cdot L}
\eeq
в частности,
\beq\label{mu_K*mu_K=mu_K}
\mu_K*\mu_K=\mu_K
\eeq

}\eit

\bpr
1. Свойство \eqref{supp-mu_K=K} очевидно.

2. Свойство \eqref{x-cdot-mu_K=mu_K} --- следствие двусторонней инвариантности меры Хаара на компактных группах.

3. Тождество \eqref{widetilde(mu_K)=mu_K} --- следствие унимодулярности компактных групп:
\begin{multline*}
\widetilde{\mu_K}(u)=\mu_K(\widetilde{u})=\int_K\widetilde{u}(t)\ \mu_K(\d t)=\int_K u(t^{-1})\ \mu_K(\d t)=\\=\begin{vmatrix}t^{-1}=s,\ t=s^{-1}, \\ \mu_K(\d s^{-1})=\eqref{svoistva-mery-Haara-3-diff}=\mu_K(\d s)\end{vmatrix}=
\int_K u(s)\ \mu_K(\d s)=\mu_K(u)
\end{multline*}

3. Тождество \eqref{overline(mu_K)=mu_K} --- следствие вещественности меры Хаара:
$$
\overline{\mu_K}(u)=\overline{\mu_K(\overline{u})}=\overline{\int_K\overline{u}(t)\ \mu_K(\d t)}=
\int_K\overline{\overline{u}(t)}\ \mu_K(\d t)=\int_Ku(t)\ \mu_K(\d t)=\mu_K(u)
$$

4. Тождество \eqref{t^(-1)-cdot-mu_K-cdot-t=mu_K} --- следствие теоремы \ref{TH:mu(ph(X))=mu(X)}: поскольку $K$ --- нормальная подгруппа в $G$, для всякого $t\in G$ отображение $\tau\mapsto t\cdot \tau\cdot t^{-1}$ --- автоморфизм $K$, поэтому оно не меняет меру $\mu_K$,  и мы получаем
\begin{multline*}
(t^{-1}\cdot \mu_K\cdot t)(u)=\mu_K(t\cdot u\cdot t^{-1})=\int_K(t\cdot u\cdot t^{-1})(s)\ \mu_K(\d s)=
\int_K u(t^{-1}\cdot s\cdot t)\ \mu_K(\d s)=\\=
\begin{vmatrix}t^{-1}\cdot s\cdot t=\tau,\ s=t\cdot \tau\cdot t^{-1}, \\ \mu_K(\d (t\cdot \tau\cdot t^{-1}))=\eqref{mu(ph(X))=mu(X)}=\mu_K(\d \tau)\end{vmatrix}= \int_K u(\tau)\ \mu_K(\d\tau)=\mu_K(u).
\end{multline*}

5. Свойство \eqref{mu_K-cdot-t=t-cdot-mu_K} следует сразу из \eqref{t^(-1)-cdot-mu_K-cdot-t=mu_K}, а из него следует \eqref{mu_K*alpha=alpha*mu_K}.

6. Для доказательства \eqref{mu_K*mu_L=mu_(KL)} нужно сначала заметить, что на подгруппе $K\cdot L$ мера $\mu_K*\mu_L$ есть нормированная мера Хаара. Действительно, во-первых, $\mu_K*\mu_L$ --- правоинвариантная мера на $K\cdot L$, потому что для любых $x\in K$ и $y\in L$ мы получим
\begin{multline*}
(\mu_K*\mu_L)\cdot(x\cdot y)=\eqref{DEF:svertka-s-delta^a}=(\mu_K*\mu_L)*\delta^{x\cdot y}=\mu_K*\mu_L*\delta^x*\delta^y=\eqref{mu_K*alpha=alpha*mu_K}=\\=
\mu_K*\delta^x*\mu_L*\delta^y=\eqref{DEF:svertka-s-delta^a}=(\mu_K\cdot x)*(\mu_L\cdot y)=\eqref{x-cdot-mu_K=mu_K}=\mu_K*\mu_L
\end{multline*}
Во-вторых, $\mu_K*\mu_L$ --- положительная мера на $K\cdot L$, потому что если $u\ge 0$, то для всякого $t\in G$
$$
(u*\widetilde{\mu_L})(t)=\eqref{eq10.10}=\widetilde{\mu_L}(\widetilde{u\cdot t})=
\mu_L(u\cdot t)\ge 0
$$
откуда
$$
u*\widetilde{\mu_L}\ge 0
$$
и это дает
$$
(\mu_K*\mu_L)(u)=\eqref{DEF:svertka}=\mu_K(u*\widetilde{\mu_L})\ge 0
$$
И, в-третьих, $\mu_K*\mu_L$ --- нормированная мера на $K\cdot L$, потому что
$$
(\mu_K*\mu_L)(1)=\eqref{DEF:svertka}=\mu_K(1*\widetilde{\mu_L})=\mu_K(1)=1.
$$
Теперь, поскольку $\mu_K*\mu_L$ --- нормированная мера Хаара на $K\cdot L$, она на этой подгруппе совпадает с 
$\mu_{K\cdot L}$:
$$
\mu_K*\mu_L\Big|_{K\cdot L}=\mu_{K\cdot L}\Big|_{K\cdot L}.
$$
Остается заметить, что вне этой подгруппы мера $\mu_K*\mu_L$ должна быть нулевой: 
$$
\supp u\subseteq G\setminus (K\cdot L) \quad\Rightarrow\quad 
\supp (u*\widetilde{\mu_L})\subseteq G\setminus K \quad\Rightarrow\quad 
(\mu_K*\mu_L)(u)=\eqref{DEF:svertka}=\mu_K(u*\widetilde{\mu_L})=0
$$
Можно заметить, что формула \eqref{mu_K*mu_K=mu_K} удобно доказывается цепочкой 
\begin{multline*}
(\mu_K*\mu_K)(u)=\eqref{DEF:svertka}=\mu_K(\widetilde{\mu_K}*u)=\eqref{mu_K(u)=int_K-u(t)mu_K(dt)}=
\int_K(\widetilde{\mu_K}*u)(t)\ \mu_K(\d t)=\eqref{eq10.10}=\\=\int_K\widetilde{\mu_K}(\widetilde{t\cdot u})\ \mu_K(\d t)=
\eqref{widetilde(u)(t)=u(t^(-1))}=\int_K\mu_K(\widetilde{\widetilde{t\cdot u}})\ \mu_K(\d t)=
\int_K\mu_K(t\cdot u)\ \mu_K(\d t)=\eqref{eq10.6}=\\=
\int_K(\mu_K\cdot t)(u)\ \mu_K(\d t)=\eqref{svoistva-mery-Haara-2-unimod}=
\int_K \underbrace{\mu_K(u)}_{\scriptsize\begin{matrix}\text{не зависит}\\ \text{от $t$}\end{matrix}}\ \mu_K(\d t)=\mu_K(u)
\end{multline*}
\epr

\btm\label{TH:mu_K->delta^1}
Направленность мер $\{\mu_K:\ K\in\lambda(G)\}$ вполне ограничена\footnote{В доказательстве этой теоремы, а также в теоремах \ref{TH:mu_K^C->delta^1}, \ref{TH:pi_K-vpolne-ogr-v-L(Env-C*(G))} и \ref{TH:Env_C-C(G)-approx} ниже, где этот результат используется, условие полной ограниченности семейства $\{\mu_K; \ K\in\lambda(G)\}$ несущественно.} а алгебре ${\mathcal C}^\star(G)$ и стремится к ее единице:
\beq\label{mu_K->delta^1}
\mu_K\overset{{\mathcal C}^\star(G)}{\underset{K\to 0}{\longrightarrow}}\delta^{1_G}=1_{{\mathcal C}^\star(G)}
\eeq
\etm
\bpr
1. Пространство ${\mathcal C}^\star(G)$ обладает свойством Гейне---Бореля, поэтому чтобы доказать полную ограниченность множества $\{\mu_K:\ K\in\lambda(G)\}$ достаточно просто доказать его ограниченность. Пусть $B\subset {\mathcal C}(G)$ --- компакт. На компакте $H_0\subseteq G$ все функции $u\in B$ ограничены некоторой общей константой $C$:
$$
\sup_{u\in B}\sup_{t\in H_0}\abs{u(t)}=C<\infty.
$$
Поэтому
\begin{multline*}
\sup_{K\in\lambda(G)}\sup_{u\in B}\abs{\mu_K(u)}=
\sup_{K\in\lambda(G)}\sup_{u\in B}\abs{\int_K u(t)\ \mu_K(\d t)}\le \sup_{K\in\lambda(G)}\sup_{u\in B}\int_K \abs{u(t)}\ \mu_K(\d t)\le\\ \le \sup_{K\in\lambda(G)} \sup_{u\in B} \sup_{t\in K}\abs{u(t)}\cdot \mu_K(K)
\le \sup_{K\in\lambda(G)} \sup_{u\in B} \sup_{t\in H_0}\abs{u(t)}\cdot \mu_K(K)
\le \sup_{K\in\lambda(G)} C\cdot \mu_K(K)
\le C
\end{multline*}
То есть множество функционалов $\{\mu_K:\ K\in\lambda(G)\}$ ограничено на компакте $B\subset {\mathcal C}(G)$. Поскольку это верно для всякого компакта $B\subset {\mathcal C}(G)$, множество $\{\mu_K:\ K\in\lambda(G)\}$ ограничено в ${\mathcal C}(G)$.

2. Теперь докажем соотношение \eqref{mu_K->delta^1}. Пусть снова $B\subset {\mathcal C}(G)$ --- компакт. Зафиксируем $\e>0$. Поскольку на множестве $H_0$ функции $u\in B$ равностиепенно непрерывны, найдется окрестность единицы $U\subseteq G$ такая, что 
$$
\sup_{u\in B}\sup_{s\in K,\ t\in U}\abs{u(s\cdot t)-u(s)}<\e
$$
В частности, при $s=1_G$, получаем
\beq\label{TH:mu_K->delta^1-1}
\sup_{u\in B}\sup_{t\in U}\abs{u(t)-u(1_G)}=\sup_{u\in B}\sup_{t\in U}\abs{u(1_G\cdot t)-u(1_G)}<\e
\eeq
Поэтому при $K\subseteq U$,
\begin{multline*}
\sup_{u\in B}\abs{\mu_K(u)-\delta^{1_G}(u)}=
\sup_{u\in B}\abs{\int_Ku(t)\ \mu_K(\d t)-u(1_G)}=
\sup_{u\in B}\abs{\int_Ku(t)\ \mu_K(\d t)-\int_Ku(1_G)\ \mu_K(\d t)}=\\
\sup_{u\in B}\abs{\int_K(u(t)-u(1_G))\ \mu_K(\d t)}\le 
\sup_{u\in B}\sup_{t\in K}\abs{u(t)-u(1_G)}\cdot \mu_K(K)<\eqref{TH:mu_K->delta^1-1}<
\sup_{u\in B}\e\cdot 1=\e.
\end{multline*}
\epr

\paragraph{Операторы $\pi_K$ и $\pi_K^\star$.}

Для всякой компактной нормальной подгруппы $K\subseteq G$ рассмотрим операторы свертки с мерой $\mu_K$:
\bit{

\item[---] оператор $\pi_K:{\mathcal C}(G)\to {\mathcal C}(G)$ определяется тождеством 
\beq\label{pi_K(u)=mu_K*u}
\pi_K(u)=\mu_K*u=u*\mu_K, \qquad u\in {\mathcal C}(G),
\eeq

\item[---] а оператор $\pi_K^\star:{\mathcal C}^\star(G)\to {\mathcal C}^\star(G)$ тождеством
\beq\label{pi_K^star(alpha)=mu_K*alpha}
\pi_K^\star(\alpha)=\mu_K*\alpha=\alpha*\mu_K, \qquad \alpha\in {\mathcal C}^\star(G).
\eeq
}\eit

\bpr
Здесь нужно проверить, что свертка с мерой $\mu_K$ справа и слева дает один и тот же результат. В \eqref{pi_K^star(alpha)=mu_K*alpha} это уже доказано в \eqref{mu_K*alpha=alpha*mu_K}, а в \eqref{pi_K(u)=mu_K*u} это доказывается цепочкой
\begin{multline*}
(\mu_K*u)(t)=\eqref{eq10.10}=\mu_K(\widetilde{t\cdot u})=\eqref{widetilde(u)(t)=u(t^(-1))}=
\widetilde{\mu_K}(t\cdot u)=\eqref{widetilde(mu_K)=mu_K}=\mu_K(t\cdot u)=\eqref{eq10.6}=
(\mu_K\cdot t)(u)=\\=\eqref{mu_K-cdot-t=t-cdot-mu_K}=(t\cdot \mu_K)(u)=\eqref{eq10.6}=\mu_K(u\cdot t)=\eqref{widetilde(mu_K)=mu_K}=
\widetilde{\mu_K}(u\cdot t)=\mu_K(\widetilde{u\cdot t})=\eqref{eq10.10}=(u*\mu_K)(t)
\end{multline*}
\epr

\btm
Операторы $\pi_K:{\mathcal C}(G)\to {\mathcal C}(G)$ и $\pi_K^\star:{\mathcal C}^\star(G)\to {\mathcal C}^\star(G)$ сопряжены друг другу
\beq\label{(pi_K)^star=pi_K^star}
(\pi_K)^\star=\pi_K^\star
\eeq
и действуют как операторы усреднения по подгруппе $K$:
\beq\label{pi_K(u)=int_K-s-u-mu_K(ds)}
\pi_K(u)=\int_K u\cdot s\ \mu_K(\d s)=\int_K s\cdot u\ \mu_K(\d s).
\eeq
\beq\label{pi_K^star(alpha)=int_K-alpha-s-mu_K(ds)}
\pi_K^\star(\alpha)=\int_K\alpha\cdot s\ \mu_K(\d s)=\int_K s\cdot\alpha\ \mu_K(\d s)
\eeq

\etm

\bpr
1. Сначала доказывается равенство \eqref{(pi_K)^star=pi_K^star}:
\begin{multline*}
(\pi_K^\star)(u)=(\mu_K*\alpha)(u)=\eqref{DEF:svertka}=\alpha(\widetilde{\mu_K}*u)=\alpha(\mu_K*u)=
\alpha(\pi_K(u))=(\alpha\circ\pi_K)(u)=\pi_K^\star(\alpha)(u)
\end{multline*}

2. Затем равенство \eqref{pi_K(u)=int_K-s-u-mu_K(ds)}: 
\begin{multline*}
\pi_K(u)(t)=(\mu_K*u)(t)=\eqref{eq10.10}=\mu_K(\widetilde{t\cdot u})=\eqref{widetilde(u)(t)=u(t^(-1))}=
\widetilde{\mu_K}(t\cdot u)=\eqref{widetilde(mu_K)=mu_K}=\mu_K(t\cdot u)=\eqref{eq10.6}=\\=
(\mu_K\cdot t)(u)=\mu_K(t\cdot u)=\int_K (t\cdot u)(s)\ \mu_K(\d s)=
\int_K (u\cdot s)(t)\ \mu_K(\d s)=\l \int_K (u\cdot s)\ \mu_K(\d s)\r (t).
\end{multline*}

3. Наконец, \eqref{pi_K^star(alpha)=int_K-alpha-s-mu_K(ds)}:
\begin{multline*}
\pi_K^\star(\alpha)(u)=\eqref{(pi_K)^star=pi_K^star}=\alpha(\pi_K(u))=\eqref{pi_K(u)=int_K-s-u-mu_K(ds)}=
\alpha\l \int_K s\cdot u\ \mu_K(\d s)\r=
 \int_K \alpha(s\cdot u)\ \mu_K(\d s)=\\= \int_K (\alpha\cdot s)(u)\ \mu_K(\d s)=
 \int_K \alpha\cdot s\ \mu_K(\d s)(u)
\end{multline*}

\epr

Из теоремы \ref{TH:mu_K->delta^1} следует

\btm\label{LM:pi_K-vpolne-ogr-v-L(C(G))}
Пусть $G$ --- LP-группа. Тогда направленность операторов $\{\pi_K; \ K\in\lambda(G)\}$ из \eqref{pi_K(u)=mu_K*u} вполне ограничена в ${\mathcal L}({\mathcal C}(G))$ и стремится в этом пространстве к тождественному оператору
\beq\label{pi_K->id-v-L(C(G))}
\pi_K\overset{{\mathcal L}({\mathcal С}(G))}{\underset{K\to 0}{\longrightarrow}}\id_{{\mathcal L}({\mathcal C}(G))}
\eeq
Точно так же направленность операторов $\{\pi_K^\star; \ K\in\lambda(G)\}$ из \eqref{pi_K^star(alpha)=mu_K*alpha} вполне ограничена в ${\mathcal L}({\mathcal C}^\star(G))$  и стремится в этом пространстве к тождественному оператору
\beq\label{pi_K'->id-v-L(C'(G))}
\pi_K^\star\overset{{\mathcal L}({\mathcal C}^\star(G))}{\underset{K\to 0}{\longrightarrow}}\id_{{\mathcal L}({\mathcal C}^\star(G))}
\eeq
\etm

\medskip
\centerline{\bf Свойства операторов $\pi_K$ и $\pi_K^\star$:}

\bit{\it

\item[$1^\circ$.] 
Для любых двух компактных нормальных подгрупп $K,L\subseteq G$ справедливы равенства
\beq\label{pi_K-circ-pi_L=pi_(K-cdot-L)}
\pi_K\circ\pi_L=\pi_{K\cdot L},
\eeq
\beq \label{pi_K*-circ-pi_L*=pi_(K-cdot-L)*}
\pi_K^\star\circ\pi_L^\star=\pi_{K\cdot L}^\star
\eeq

\item[$2^\circ$.] Операторы $\pi_K$ коммутируют:
\beq\label{pi_K-circ-pi_L=pi_L-circ-pi_K}
\pi_K\circ\pi_L=\pi_L\circ\pi_K
\eeq
и таким же образом, операторы $\pi_K^\star$ коммутируют:
\beq\label{pi_K*-circ-pi_L*=pi_L*-circ-pi_K*}
\pi_K^\star\circ\pi_L^\star=\pi_L^\star\circ\pi_K^\star.
\eeq

\item[$3^\circ$.] Операторы $\pi_K$ и $\pi_K^\star$ являются проекторами:
\beq\label{pi_K^2=pi_K}
\pi_K^2=\pi_K,
\eeq
\beq\label{pi_K*^2=pi_K*}
(\pi_K^\star)^2=\pi_K^\star.
\eeq

\item[$4^\circ$.] Операторы $\pi_K$ и $\pi_K^\star$ сохраняют антипод и инволюцию в ${\mathcal C}(G)$ и ${\mathcal C}^\star(G)$:
\beq\label{pi_K(overline(u))=overline(pi_K(u))}
\pi_K(\widetilde{u})=\widetilde{\pi_K(u)},\quad \pi_K(\overline{u})=\overline{\pi_K(u)},\qquad u\in {\mathcal C}(G),
\eeq
\beq\label{pi_K^star(alpha^bullet)=pi_K^star(alpha)^bullet}
\pi_K^\star(\widetilde{\alpha})=\widetilde{\pi_K^\star(\alpha)},\quad \pi_K^\star(\alpha^\bullet)=\pi_K^\star(\alpha)^\bullet,\quad \qquad \alpha\in {\mathcal C}^\star(G).
\eeq

\item[$5^\circ$.] Операторы $\pi_K^\star$ сохраняют умножение\footnote{Свойства $4^\circ$ и $5^\circ$ в этом списке означают, в частности, что операторы $\pi_K^\star$ являются морфизмами в категории $\InvSteAlg^0$, определенной ниже на с.\pageref{DEF:InvSteAlg^0}.} в ${\mathcal C}^\star(G)$ 
\beq\label{pi_K^star(alpha*beta)=pi_K^star(alpha)*pi_K^star(beta)}
\pi_K^\star(\alpha*\beta)=\pi_K^\star(\alpha)*\pi_K^\star(\beta),\qquad \alpha,\beta\in {\mathcal C}^\star(G).
\eeq
но не единицу, если $K\ne 1_G$:  
\beq
\pi_K^\star(\delta^{1_G})=\mu_K*\delta^{1_G}=\mu_K\ne \delta^{1_G}
\eeq

}\eit

\bpr
1. Равенство \eqref{pi_K-circ-pi_L=pi_(K-cdot-L)} доказывается цепочкой
\begin{multline*}
(\pi_K\circ\pi_L)(u)=\pi_K(\pi_L(u))=\mu_K*(\mu_L*u)=(\mu_K*\mu_L)*u=(\mu_L*\mu_K)*u=\\=
\mu_L*(\mu_K*u)=\pi_L(\pi_K(u))=(\pi_L\circ\pi_K)(u)
\end{multline*}

4. Первое равенство в \eqref{pi_K(overline(u))=overline(pi_K(u))} доказывается цепочкой
\begin{multline*}
\pi_K(\widetilde{u})(t)=\eqref{pi_K(u)=mu_K*u}=(\mu_K*\widetilde{u})(t)=\eqref{eq10.10}=\mu_K\l \widetilde{t\cdot \widetilde{u}}\r
=\eqref{eq10.7}=\mu_K\l u\cdot t^{-1}\r=\eqref{widetilde(mu_K)=mu_K}=\\=\widetilde{\mu_K}\l u\cdot t^{-1}\r=\mu_K\l\widetilde{ u\cdot t^{-1}}\r=\eqref{eq10.10}=
(u*\mu_K)(t^{-1})=\eqref{pi_K(u)=mu_K*u}=\pi_K(u)(t^{-1})=\widetilde{\pi_K(u)}(t)
\end{multline*}
а второе --- цепочкой
\begin{multline*}
\pi_K(\overline{u})(t)=\eqref{pi_K(u)=mu_K*u}=(\mu_K*\overline{u})(t)=\eqref{eq10.10}=\mu_K\l \widetilde{t\cdot \overline{u}}\r
=\mu_K\l \overline{\widetilde{t\cdot u}}\r=\\=\overline{\overline{\mu_K}\l \widetilde{t\cdot u}\r}=\eqref{overline(mu_K)=mu_K}=
\overline{\mu_K\l \widetilde{t\cdot u}\r}=\eqref{eq10.10}=\overline{(\mu_K*u)(t)}=\eqref{pi_K(u)=mu_K*u}=\overline{\pi_K(u)(t)}.
\end{multline*}
Теперь  мы получаем:
\begin{multline*}
\pi_K^\star(\widetilde{\alpha})(u)=\eqref{(pi_K)^star=pi_K^star}=\widetilde{\alpha}(\pi_K(u))=\eqref{widetilde(u)(t)=u(t^(-1))}=
\alpha(\widetilde{\pi_K(u)})=\eqref{pi_K(overline(u))=overline(pi_K(u))}=\\=
\alpha(\pi_K(\widetilde{u}))=\eqref{(pi_K)^star=pi_K^star}=\pi_K^\star(\alpha)(\widetilde{u})=\eqref{widetilde(u)(t)=u(t^(-1))}=
\widetilde{\pi_K^\star(\alpha)}(u)
\end{multline*}
и
\begin{multline*}
\pi_K^\star(\widetilde{\alpha})(u)=\eqref{(pi_K)^star=pi_K^star}=\widetilde{\alpha}(\pi_K(u))=\eqref{widetilde(u)(t)=u(t^(-1))}=
\alpha(\widetilde{\pi_K(u)})=\eqref{pi_K(overline(u))=overline(pi_K(u))}=\\=
\alpha(\pi_K(\widetilde{u}))=\eqref{(pi_K)^star=pi_K^star}=\pi_K^\star(\alpha)(\widetilde{u})=\eqref{widetilde(u)(t)=u(t^(-1))}=
\widetilde{\pi_K^\star(\alpha)}(u)
\end{multline*}

\epr

\subsection{Стереотипная аппроксимация в ${\mathcal C}(M)$ и ${\mathcal C}(G)$}

\paragraph{Аппроксимация в ${\mathcal C}(M)$.}

Пусть $M$ --- полное метрическое пространство, и $d$ --- расстояние на нем. Обозначим символом ${\mathcal C}_\natural(M)$ алгебру непрерывных функций $f:M\to\C$ с обычным поточечным умножением и топологией равномерной сходимости на компактах в $M$. Иначе говоря, топология на ${\mathcal C}_\natural(M)$ задается системой полунорм
\beq\label{norm(f)_T=max_t-abs(f(t))}
\norm{f}_T=\max_{t\in T}\abs{f(t)},
\eeq
в которых $T$ пробегает множество всех компактов в $M$. Операция умножения в ${\mathcal C}_\natural(M)$ будет непрерывной билинейной формой в смысле теории стереотипных пространств: для любой окрестности нуля $U\subseteq {\mathcal C}_\natural(M)$ и любого вполне ограниченного множества $F\subseteq {\mathcal C}_\natural(M)$ найдется окрестность нуля $V\subseteq {\mathcal C}_\natural(M)$ такая, что $V\cdot F=F\cdot V\subseteq U.$ Отсюда в силу \cite[Theorem 3.6.6]{Akbarov-De-Gruyter-I} следует, что при псевдонасыщении $\vartriangle$ операция умножения остается непрерывной билинейной формой в этом же смысле. Поэтому псевдонасыщение ${\mathcal C}_\natural(M)^\vartriangle$ пространства ${\mathcal C}_\natural(M)$ должно быть стереотипной алгеброй. Далее везде мы символом ${\mathcal C}(M)$ будем обозначать эту алгебру ${\mathcal C}_\natural(M)^\vartriangle$:
$$
{\mathcal C}(M):={\mathcal C}_\natural(M)^\vartriangle.
$$

\btm\label{TH:C(M)-obl-approx-kogda-M-metricheskoe}
Для всякого полного метрического пространства $M$ алгебра ${\mathcal C}(M)$ обладает стереотипной аппроксимацией.
\etm
\bpr
1. Для произвольного $k\in\N$ рассмотрим систему $\{U_{\frac{1}{k}}(t); t\in M\}$ открытых шаров радиуса $\frac{1}{k}$ в $M$. Это будет открытое покрытие $M$, поэтому ему подчинено некоторое локально конечное разбиение единицы $\{\eta^k_t;t\in M\}$ \cite[Теорема 5.1.3, Теорема 5.1.9]{Engelking}:
$$
\eta^k_t\in{\mathcal C}(M),\quad 0\le \eta^k_t\le 1,\quad \Supp\eta^k_t\subseteq U_{\frac{1}{k}}(t),\quad \sum_{t\in M}\eta^k_t=1.
$$
Пусть $M^k\subseteq M$ --- множество, состоящее из точек $t\in M$ у которых соответствующая функция $\eta^k_t$ ненулевая: $t\in M^k$ $\Leftrightarrow$ $\eta^k_t\ne 0.$

2. Определим далее последовательность операторов $P^k:{\mathcal C}(M)\to {\mathcal C}(M)$
\beq\label{C(M)-obl-approx-kogda-M-metricheskoe-1}
P^kf=\sum_{t\in M^k} f(t)\cdot \eta^k_t,\qquad f\in{\mathcal C}(M)
\eeq
Поскольку ряд справа локально конечен (то есть в окрестности каждой точки $x\in M$ только конечное число его слагаемых отлично от нуля), функция $P^kf\in{\mathcal C}(M)$ корректно определена, и мы получаем линейное отображение $P^k:{\mathcal C}(M)\to {\mathcal C}(M)$.

Заметим, что это отображение будет непрерывно в топологии ${\mathcal C}_\natural(M)$. Действительно, если $\{f_i;i\to\infty\}$ --- направленность, стремящаяся к нулю в ${\mathcal C}_\natural(M)$, то есть равномерно на компактах $T\subseteq M$, $\max_{x\in T}\abs{f_i(x)}\underset{i\to\infty}{\longrightarrow}0,$
то, поскольку множество
\beq\label{M^k_T=(t-in-T:eta^k_t|_T-ne-0)}
M^k_T=\{t\in T: \eta^k_t\Big|_T\ne 0\}
\eeq
конечно, мы получаем
\begin{multline*}
\max_{x\in T}\abs{P^kf_i(x)}=
\max_{x\in T}\abs{\sum_{t\in M^k} f_i(t)\cdot \eta^k_t(x)}=
\max_{x\in T}\abs{\sum_{t\in M^k_T} f_i(t)\cdot \eta^k_t(x)}\le\\ \le
\max_{x\in T}\sum_{t\in M^k_T}\abs{f_i(t)\cdot \eta^k_t(x)}\le
\max_{x\in T}\max_{t\in M^k_T}\abs{f_i(t)}\cdot \underbrace{\sum_{t\in M^k_T}\eta^k_t(x)}_{\scriptsize\begin{matrix}\| \\1 \end{matrix}}\le
\max_{t\in M^k_T} \abs{f_i(t)}
\underset{i\to\infty}{\longrightarrow}0.
\end{multline*}

Из того, что операторы  $P^k:{\mathcal C}_\natural(M)\to {\mathcal C}_\natural(M)$ непрерывны, следует, что их псевдонасыщения $(P^k)^\vartriangle:{\mathcal C}(M)^\vartriangle\to {\mathcal C}(M)^\vartriangle$, то есть операторы $P^k:{\mathcal C}(M)\to {\mathcal C}(M)$, тоже непрерывны \cite[Theorem 3.3.15]{Akbarov-De-Gruyter-I}.

3. Покажем далее, что представление  \eqref{C(M)-obl-approx-kogda-M-metricheskoe-1} оператора $P^k$ можно понимать как разложение $P^k$ в сходящийся ряд из одномерных операторов в ${\mathcal L}({\mathcal C}(M))$:
\beq\label{C(M)-obl-approx-kogda-M-metricheskoe-2}
P^k=\sum_{t\in M^k} \delta_t\odot \eta^k_t
\eeq
($\delta_t$ --- $\delta$-функционал в точке $t$). Это делается в несколько этапов. Сначала покажем, что формулу \eqref{C(M)-obl-approx-kogda-M-metricheskoe-2} можно понимать как сходимость ряда в пространстве ${\mathcal C}_\natural(M):{\mathcal C}(M)$ (мы используем обозначения \cite[3.5.1]{Akbarov-De-Gruyter-I}: по ним ${\mathcal C}_\natural(M):{\mathcal C}(M)$ есть пространство линейных непрерывных отображений $\ph:{\mathcal C}(M)\to{\mathcal C}_\natural(M)$ с топологией равномерной сходимости на компактах). Для этого зафиксируем вполне ограниченное множество $F\subseteq {\mathcal C}(M)$ и базисную окрестность нуля в ${\mathcal C}_\natural(M)$, то есть множество вида $B_{\e}(T)=\{g\in {\mathcal C}_\natural(M): \norm{g}_T<\e\},$ где $T$ --- компакт в $M$. Вспомним о множестве $M^k_T$, которое мы определили формулой \eqref{M^k_T=(t-in-T:eta^k_t|_T-ne-0)}. Оно конечно (поскольку семейство $\eta^k_t$ локально конечно), поэтому для всякого конечного множества $N\supseteq M^k_T$ мы получаем цепочку
\begin{multline*}
\sum_{t\in N} (\delta_t\odot\eta^k_t)(f)\Big|_T=
\sum_{t\in N} f(t)\cdot \eta^k_t\Big|_T=\sum_{t\in N\setminus M^k_T} f(t)\cdot \underbrace{\eta^k_t\Big|_T}_{\scriptsize\begin{matrix}\|\\ 0\end{matrix}}+\sum_{t\in M^k_T} f(t)\cdot \eta^k_t\Big|_T=\\=
\sum_{t\in M^k_T} f(t)\cdot \eta^k_t\Big|_T=P^kf\Big|_T
\end{multline*}
Отсюда мы получаем цепочку
$
\Big(\forall N\supseteq M^k_T\quad \forall f\in F\quad \|\sum_{t\in N} (\delta_t\odot\eta^k_t)(f)-P^kf\|_T=0\Big)
$
$\Longrightarrow$
$\Big(
\forall N\supseteq M^k_T\quad \forall f\in F\quad \sum_{t\in N} (\delta_t\odot\eta^k_t)(f)-P^kf\in B_{\e}(T)\Big)
$
$\Longrightarrow$
$\Big(
\forall N\supseteq M^k_T\quad \sum_{t\in N}\delta_t\odot\eta^k_t-P^k\in B_{\e}(T):F\Big)
$
(здесь $B_{\e}(T):F$ обозначает множество операторов $\ph:{\mathcal C}(M)\to {\mathcal C}_\natural(M)$ со свойством $\ph(B_{\e}(T))\subseteq F$, см. \cite[p.421]{Akbarov-De-Gruyter-I}; по определению топологии в ${\mathcal C}_\natural(M):{\mathcal C}(M)$, это будет базисная окрестность нуля в ${\mathcal C}_\natural(M):{\mathcal C}(M)$). И это доказывает, что \eqref{C(M)-obl-approx-kogda-M-metricheskoe-2} действительно можно понимать как сходимость ряда в пространстве ${\mathcal C}_\natural(M):{\mathcal C}(M)$.

Далее, из того, что ряд справа в \eqref{C(M)-obl-approx-kogda-M-metricheskoe-2}  сходится в ${\mathcal C}_\natural(M):{\mathcal C}(M)$ следует, что его частичные суммы
$P^k_N=\sum_{t\in N} \delta_t\odot \eta^k_t$, $N\in 2_{M^k}$
($2_{M^k}$ обозначает множество всех конечных подмножеств в $M^k$)
образуют вполне ограниченное множество в ${\mathcal C}_\natural(M):{\mathcal C}(M)$ \cite[Proposition 4.4.78]{Akbarov-De-Gruyter-I}. В силу \cite[Theorem 3.5.1]{Akbarov-De-Gruyter-I}, это означает, что семейство операторов $\{P^k_N;N\in 2_{M^k}\}$ равностепенно непрерывно на любом вполне ограниченном множестве $F\subseteq {\mathcal C}(M)$ и имеет вполне ограниченный образ на нем $\bigcup_{N\in 2_{M^k}}P^k_N(F)\subseteq {\mathcal C}_\natural(M)$.
Поскольку система вполне ограниченных множеств в ${\mathcal C}_\natural(M)$ и ${\mathcal C}(M)$ одинакова, и топология на них тоже не меняется \cite[Theorem 3.3.17]{Akbarov-De-Gruyter-I}, мы можем заключить, что $\bigcup_{N\in 2_{M^k}}P^k_N(F)$ вполне ограничено в ${\mathcal C}(M)$, и $\{P^k_N;N\in 2_{M^k}\}$ равностепенно непрерывна как система отображений из $F$ (с индуцированной из ${\mathcal C}(M)$ равномерной структурой) в $\bigcup_{N\in 2_{M^k}}P^k_N(F)\subseteq {\mathcal C}(M)$, или, иными словами, как система отображений из $F$ в ${\mathcal C}(M)$. Поскольку это верно для любого вполне ограниченного множества $F\subseteq {\mathcal C}(M)$, мы, в силу \cite[Theorem 3.5.1]{Akbarov-De-Gruyter-I}, заключаем, что система операторов $\{P^k_N;N\in 2_{M^k}\}$ является вполне ограниченным множеством в ${\mathcal C}(M):{\mathcal C}(M)$ (с убранным символом $\natural$ в числителе дроби).

Замыкание $\overline{\{P^k_N;N\in 2_{M^k}\}}$ этого множества в ${\mathcal C}(M):{\mathcal C}(M)$ будет компактом в ${\mathcal C}(M):{\mathcal C}(M)$. Формально более слабая топология пространства ${\mathcal C}_\natural(M):{\mathcal C}(M)$ все равно хаусдорфова, и поэтому разделяет точки $\overline{\{P^k_N;N\in 2_{M^k}\}}$. Отсюда следует, что на множестве $\overline{\{P^k_N;N\in 2_{M^k}\}}$ эти топологии совпадают. Поэтому из сходимости направленности $\{P^k_N;N\in 2_{M^k}\}$ к элементу $P^k\in \overline{\{P^k_N;N\in 2_{M^k}\}}$ относительно топологии ${\mathcal C}_\natural(M):{\mathcal C}(M)$ следует сходимость $\{P^k_N;N\in 2_{M^k}\}$ к элементу $P^k\in \overline{\{P^k_N;N\in 2_{M^k}\}}$ относительно топологии ${\mathcal C}(M):{\mathcal C}(M)$:
$P^k_N\overset{{\mathcal C}(M):{\mathcal C}(M)}{\underset{N\to M}{\longrightarrow}}P^k.$

Поскольку множество $\{P^k_N;N\in 2_{M^k}\}\cup\{P^k\}$ вполне ограничено в ${\mathcal C}(M):{\mathcal C}(M)$, при переходе к псевдонасыщению $\vartriangle$ его топология не меняется. Отсюда следует, что направленность
$\{P^k_N;N\in 2_{M^k}\}$ сходится к элементу $P^k$ в топологии пространства $({\mathcal C}(M):{\mathcal C}(M))^\vartriangle={\mathcal L}({\mathcal C}(M))$:
$P^k_N\overset{{\mathcal L}({\mathcal C}(M))}{\underset{N\to M}{\longrightarrow}}P^k.$
Это как раз и означает, что равенство  \eqref{C(M)-obl-approx-kogda-M-metricheskoe-2} справедливо в пространстве ${\mathcal L}({\mathcal C}(M))$.

4. Покажем, что операторы $P^k$ аппроксимируют тождественный оператор $I$ в пространстве ${\mathcal L}({\mathcal C}(M))$:
\beq\label{C(M)-obl-approx-kogda-M-metricheskoe-5}
P^k\overset{{\mathcal L}({\mathcal C}(M))}{\underset{N\to M}{\longrightarrow}}I.
\eeq
Это тоже проделывается в несколько этапов. Опять зафиксируем вполне ограниченное множество $F\subseteq {\mathcal C}(M)$ и компакт $T\subseteq M$. Тогда
\begin{multline*}
  \norm{P^kf-f}_T=\max_{x\in T}\abs{P^kf(x)-f(x)}=
  \max_{x\in T}\Big|\sum_{t\in M^k} f(t)\cdot \eta^k_t(x)-f(x)\Big|=\\=
  \max_{x\in T}\Big|\sum_{t\in M^k} f(t)\cdot \eta^k_t(x)-f(x)\cdot \underbrace{\sum_{t\in M^k}\eta^k_t(x)}_{\scriptsize\begin{matrix}\|\\ 1\end{matrix}}\Big|=
  \max_{x\in T}\Big|\sum_{t\in M^k} \big(f(t)-f(x)\big)\cdot \eta^k_t(x)\Big|\le \\ \le
  \max_{x\in T}\sum_{t\in M^k} \Big|f(t)-f(x)\Big|\cdot \eta^k_t(x)
  \end{multline*}
В последней сумме если $\eta^k_t(x)\ne 0$, то $x\in \Supp\eta^k_t\subseteq U_{\frac{1}{k}}(t)$, и значит $d(x,t)<\frac{1}{k}$. Поэтому
\begin{multline}\label{C(M)-obl-approx-kogda-M-metricheskoe-3}
\norm{P^kf-f}_T\le \max_{x\in T}\sum_{t\in M^k} \kern-22pt \underbrace{\Big|f(t)-f(x)\Big|}_{ \scriptsize\begin{matrix}\text{\rotatebox{90}{$\ge$}} \\
\sup_{x,s\in T,\ d(s,x)< \frac{1}{k}}\Big|f(s)-f(x)\Big| \end{matrix}}\kern-22pt \cdot \eta^k_t(x)\le\\ \le
\sup_{x,s\in T,\ d(s,x)< \frac{1}{k}}\Big|f(s)-f(x)\Big|\cdot \underbrace{\sum_{t\in M^k}
\eta^k_t(x)}_{\scriptsize\begin{matrix}\|\\ 1\end{matrix}}=
\sup_{x,s\in T,\ d(s,x)< \frac{1}{k}}\Big|f(s)-f(x)\Big|
\end{multline}
Покажем, что последняя величина стремится к нулю равномерно по $f\in F$:
\beq\label{C(M)-obl-approx-kogda-M-metricheskoe-4}
\sup_{f\in F}\sup_{x,s\in T,\ d(s,x)< \frac{1}{k}}\Big|f(s)-f(x)\Big|\underset{k\to\infty}{\longrightarrow}0
\eeq
Предположим, что это не так:
$
\sup_{f\in F}\sup_{x,s\in T,\ d(s,x)< \frac{1}{k}}\Big|f(s)-f(x)\Big|\underset{k\to\infty}{\not\longrightarrow}0.
$
Это значит, что существует некое число $\e>0$ и последовательности $k_n\to \infty$, $f_n\in F$, $x_n\in T$, $s_n\in U_{\frac{1}{k_n}}(x_n)$ такие что
$
\forall n\in\N\quad \Big|f_n(s_n)-f_n(x_n)\Big|>\e.
$
Поскольку $T$ --- компакт, из последовательности $x_n\in T$ можно выбрать сходящуюся подпоследовательность $x_{n_i}$:
$
x_{n_i}\underset{i\to\infty}{\longrightarrow}x\in T.
$
Тогда
$
s_{n_i}\underset{i\to\infty}{\longrightarrow}x\in T,
$
и мы получим неравенство
$
\forall i\in\N\quad \Big|f_{n_i}(s_{n_i})-f_{n_i}(x_{n_i})\Big|>\e,
$
в котором $f_{n_i}\in F$, а $x_{n_i}\to x$ и $s_{n_i}\to x$. Это значит, что множество функций $F$ не равностепенно непрерывно на компакте $\{x_{n_i}\}\cup \{s_{n_i}\}\cup\{x\}$, и поэтому $F$ не будет вполне ограничено в ${\mathcal C}(M)$ \cite[8.2.10]{Engelking}. А это противоречит выбору $F\subseteq {\mathcal C}(M)$.

Мы доказали \eqref{C(M)-obl-approx-kogda-M-metricheskoe-4} и вместе с \eqref{C(M)-obl-approx-kogda-M-metricheskoe-3} это дает
$0\le\sup_{f\in F}\norm{P^kf-f}_T$ $\le$\break $\sup_{f\in F} \sup_{x\in T,\ d(s,x) < \frac{1}{k}}\Big|f(s)-f(x)\Big|\underset{k\to\infty}{\longrightarrow}0$,
и значит
$\sup_{f\in F}\norm{P^kf-f}_T\underset{k\to\infty}{\longrightarrow}0.$
Это верно для любого компакта $T\subseteq M$, поэтому мы можем сказать, что для любого вполне ограниченного множества $F\subseteq {\mathcal C}(M)$ значения $P^kf-f$ стремятся к нулю в пространстве ${\mathcal C}_\natural(M)$ равномерно по $f\in F$:
\beq\label{C(M)-obl-approx-kogda-M-metricheskoe-6}
P^kf-f\overset{{\mathcal C}_\natural(M)}{\underset{f\in F}{\underset{k\to\infty}{\rightrightarrows}}}0.
\eeq
Вспомним, что $P^k$ --- не просто направленность, а {\it последовательность}. Вместе с \eqref{C(M)-obl-approx-kogda-M-metricheskoe-6} это дает, что множество $\bigcup_{k\in \N}(P^k-I)(F)$ должно быть вполне ограничено в ${\mathcal C}_\natural(M)$. Поэтому при псевдонасыщении $\vartriangle$ топология на $\bigcup_{k\in \N}(P^k-I)(F)\cup\{0\}$ не поменяется \cite[Theorem 3.3.17]{Akbarov-De-Gruyter-I}, и мы можем заключить, что
$P^kf-f$ стремится к нулю в пространстве ${\mathcal C}_\natural(M)^\vartriangle={\mathcal C}(M)$ равномерно по $f\in F$:
\beq\label{C(M)-obl-approx-kogda-M-metricheskoe-7}
P^kf-f\overset{{\mathcal C}(M)}{\underset{f\in F}{\underset{k\to\infty}{\rightrightarrows}}}0.
\eeq
И это верно для всякого вполне ограниченного множества $F\subseteq {\mathcal C}(M)$. Значит, справедливо соотношение
$
P^k\overset{{\mathcal C}(M):{\mathcal C}(M)}{\underset{N\to M}{\longrightarrow}}I.
$
Вспомним опять, что $P^k$ --- последовательность. Поскольку она сходится к $I$ в пространстве ${\mathcal C}(M):{\mathcal C}(M)$, множество $\{P^k\}\cup\{I\}$ компактно. Поэтому мы можем считать, что $P^k$ сходится к $I$ в компакте $\{P^k\}\cup\{I\}\subseteq {\mathcal C}(M):{\mathcal C}(M)$. Когда мы применяем к пространству ${\mathcal C}(M):{\mathcal C}(M)$ операцию псевдонасыщения $\vartriangle$, топология на компакте $\{P^k\}\cup\{I\}$ не меняется \cite[Theorem 3.3.17]{Akbarov-De-Gruyter-I}. Поэтому мы можем сказать, что $P^k$ сходится к $I$ в компакте $\{P^k\}\cup\{I\}\subseteq ({\mathcal C}(M):{\mathcal C}(M))^\vartriangle={\mathcal L}({\mathcal C}(M))$. Значит, $P^k$ сходится к $I$ в пространстве ${\mathcal L}({\mathcal C}(M))$. То есть справедливо \eqref{C(M)-obl-approx-kogda-M-metricheskoe-5}.

5. Итак, мы получили, что единичный оператор $I$ аппроксимируется в пространстве ${\mathcal L}({\mathcal C}(M))$ операторами $P^k$, в силу \eqref{C(M)-obl-approx-kogda-M-metricheskoe-5}, а сами операторы $P^k$ аппроксимируются в ${\mathcal L}({\mathcal C}(M))$ конечномерными (частичными суммами ряда в \eqref{C(M)-obl-approx-kogda-M-metricheskoe-2}). Значит, $I$ аппроксимируется конечномерными операторами в ${\mathcal L}({\mathcal C}(M))$, а это и есть свойство стереотипной аппроксимации для ${\mathcal C}(M)$.
\epr

\paragraph{Аппроксимация в ${\mathcal C}(G)$ и ${\mathcal C}^\star(G)$.}

\btm\label{TH:C(G)-approx}
Для всякой локально компакной группы $G$ пространства ${\mathcal C}(G)$ и ${\mathcal C}^\star(G)$ обладают стереотипной аппроксимацией.
\etm
\bpr
Это достаточно доказать для ${\mathcal C}(G)$.

1. Если $G$ --- группа Ли, то это верно по теореме \ref{TH:C(M)-obl-approx-kogda-M-metricheskoe}.

2. Пусть далее $G$ --- LP-группа. Для фиксированного $K\in\lambda(G)$ мы получим, что $G/K$ --- группа Ли, и поэтому по пункту 1, пространство ${\mathcal C}(G:K)$ обладает стереотипной аппроксимацией. Поэтому можно подобрать направленность $\{\ph_K^i; \ i\to\infty\}$ конечномерных операторов, аппроксимирующих тождественный оператор $\id_{{\mathcal L}({\mathcal C}(G:K))}$:
\beq\label{TH:C(G)-approx-2}
\ph_K^i\overset{{\mathcal L}({\mathcal C}(G:K))}{\underset{i\to\infty}{\longrightarrow}}\id_{{\mathcal L}({\mathcal C}(G:K))}
\eeq
Теперь мы получаем, что (конечномерные) операторы $\im\pi_K\circ\ph_K^i\circ\coim\pi_K$ аппроксимируют тождественный оператор $\id_{{\mathcal C}(G)}$ в силу \eqref{TH:C(G)-approx-2} и \eqref{pi_K->id-v-L(C(G))}:
$$
\im\pi_K\circ\ph_K^i\circ\coim\pi_K\overset{{\mathcal L}({\mathcal C}(G))}{\underset{i\to\infty}{\longrightarrow}}\im\pi_K\circ\id_{{\mathcal L}({\mathcal C}(G:K))}\circ\coim\pi_K=\pi_K\overset{{\mathcal L}({\mathcal C}(G))}{\underset{K\to 0}{\longrightarrow}}\id_{{\mathcal L}({\mathcal C}(G))}
$$
$$
 \xymatrix 
 {
 {\mathcal C}(G)\ar[r]^{\id_{{\mathcal L}({\mathcal C}(G))}}\ar@{-->}[d]_{\coim\pi_K}& {\mathcal C}(G)\\
 {\mathcal C}(G:K)\ar@{-->}[r]_{\ph_K^i} & {\mathcal C}(G:K)\ar@{-->}[u]_{\im\pi_K}
 }
$$
(эта диаграмма не коммутативна, она только показывает, на каких пространствах действуют рассматриваемые операторы).

3. Наконец, если $G$ --- произвольная локально компактная группа, то мы можем выбрать в ней открытую LP-подгруппу $G_0$, и тогда как топологическое пространство $G$ представляется в виде прямой суммы топологических пространств вида $G_0$:
$$
G=\coprod_{i\in I}G_i,\qquad G_i\cong G_0.
$$
А пространство ${\mathcal C}(G)$ тогда будет прямым произведением пространств ${\mathcal C}(G_i)$:
$$
{\mathcal C}(G)=\prod_{i\in I}{\mathcal C}(G_i).
$$
При этом в пункте 2 мы уже показали, что пространства ${\mathcal C}(G_i)$ обладают аппроксимацией. Значит, ${\mathcal C}(G)$ тоже обладает аппроксимацией.
\epr

\section{Пространства Брюа}

В работе \cite{Bruhat} 1961 года Франсуа Брюа предложил определение гладкой функции на произвольной локально компактной группе $G$. Мы упомянем его конструкцию среди примеров функциональных пространств на группах. Это, с одной стороны, будет полезно как иллюстрация применения теоремы о наследовании стереотипности \ref{TH:o-nasl-ster-karkasom}, а с другой --- просто как напоминание коллегам о незаслуженно забытой ветви гармонического анализа.

\subsection{Пространство ${\mathcal D}(G)$.}

Построения Брюа начинаются с пространства ${\mathcal D}(G)$ финитных гладких функций на локально компактной группе $G$.

\paragraph{Случай LP-группы $G$.}

Пусть $G_0$ --- LP-группа. Для всякой подгруппы $K\in\lambda(G_0)$ фактор-группа $G_0/K$ является группой Ли. Поэтому определено  пространство ${\mathcal D}(G_0/K)$ финитных гладких функций на $G_0/K$. Функции, получающиеся как композиции функций из ${\mathcal D}(G_0/K)$ с фактор-отображением $\rho_K:G_0\to G_0/K$, образуют некое пространство функций на $G_0$
\beq\label{DEF:D(G_0:K)}
{\mathcal D}(G_0:K)=\{v\circ \rho_K;\ v\in {\mathcal D}(G_0/K)\}
\eeq 
называемое {\it пространством финитных гладких функций на $G$, инвариантных относительно $K$}. Изоморфизм векторных пространств
\beq\label{DEF:D(G_0:K)-topologiya}
{\mathcal D}(G_0:K)\cong {\mathcal D}(G_0/K)
\eeq 
наделяет каждое пространство ${\mathcal D}(G_0:K)$ естественной топологией.

\blm\label{PROP:D(G_0:K)-nasysheno} Пространство ${\mathcal D}(G_0:K)$ полно и насыщено (и поэтому стереотипно).
\elm

Объединение этих пространств
\beq\label{DEF:D(G_0)}
{\mathcal D}(G_0)=\bigcup_{K\in\lambda(G_0)}{\mathcal D}(G_0:K)
\eeq 
называется {\it пространством финитных гладких функций на $G$}. Оно наделяется сильнейшей локально выпуклой топологией, отнсительно которой все вложения
\beq\label{DEF:D(G_0:K)->D(G_0)}
{\mathcal D}(G_0:K)\subseteq {\mathcal D}(G_0)
\eeq
непрерывны\footnote{В учебниках по топологическим векторным пространствам (например, в \cite{Schaefer} и \cite{Jarchow}) такая топология называется индуктивной топологией порожденной семейством отображений (в данном случае, отображениями \eqref{DEF:D(G_0:K)->D(G_0)}).}. По-другому эту топологию можно описать как топологию локально выпуклого инъективного предела пространств  ${\mathcal D}(G_0:K)$:
\beq\label{DEF:D(G_0)-topologiya}
{\mathcal D}(G_0)=\LCS\text{-}\injlim_{K\in\lambda(G_0)}{\mathcal D}(G_0:K)
\eeq 

\btm\label{PROP:D(G_0)-nasysheno} Для всякой LP-группы $G_0$ пространство ${\mathcal D}(G_0)$ полно и насыщено (и поэтому стереотипно).
\etm
\bpr
Для любых двух групп $K,L\in\lambda(G)$ с условием $K\supseteq L$, справедливо с ожной стороны, вложение   
$$
{\mathcal D}(G_0:K)\subseteq {\mathcal D}(G_0:L),
$$
а, с другой стороны, оператор $\pi_K$ будет  переводить ${\mathcal D}(G_0:L)$ в ${\mathcal D}(G_0:K)$:
$$
\pi_K({\mathcal D}(G_0:L))\subseteq {\mathcal D}(G_0:K).
$$
Если обозначить первое вложение
$$
\iota_K^L:{\mathcal D}(G_0:K)\to {\mathcal D}(G_0:L)
$$ 
а вторую проекцию
$$
\varkappa_K^L:{\mathcal D}(G_0:L)\to {\mathcal D}(G_0:K),
$$
то мы получиим, что система морфизмов
$$
\varPhi=(\{\iota_K^L\},\{\varkappa_K^L\})
$$
образует фундамент в категории $\LCS$ в смысле определения на с.\pageref{DEF:diagr-building}. 
 
Этот фундамент будет удовлетворять условиям теоремы о наследовании стереотипности \ref{TH:o-nasl-ster-karkasom}. Значит, его каркас
$$
\fram\varPhi=\LCS\text{-}\injlim_{K\in\lambda(G_0)}{\mathcal D}(G_0:K)={\mathcal D}(G_0)
$$
является полным и насыщенным пространством.
\epr

\paragraph{Случай произвольной локально компактной группы $G$.}

Пусть $G$ --- произвольная локально компактная группа. Тогда по теореме \ref{TH:LP-gruppa-lambda} в $G$ можно выбрать открытую LP-подгруппу $G_0$. Для нее уже определено пространство ${\mathcal D}(G_0)$. Его левые сдвиги будут образовывать функциональные пространства, которые определяются как пространства ${\mathcal D}(x\cdot G_0)$ гладких финитных функций на левых смежных классах $G_0$:
\beq
{\mathcal D}(x\cdot G_0)=\{u\circ\tau_x;\ u\in {\mathcal D}(G_0)\}
\eeq
где
\beq
\tau_x(y)=x^{-1}\cdot y,\qquad y\in x\cdot G_0
\eeq
После этого пространство ${\mathcal D}(G)$ определяется как объединение пространств ${\mathcal D}(x\cdot G_0)$,
\beq
{\mathcal D}(G)=\bigcup_{x\cdot G_0\in G/G_0}{\mathcal D}(x\cdot G_0)
\eeq
Оно наделяется сильнейшей локально выпуклой топологией, отнсительно которой все вложения
\beq\label{DEF:D(x-cdot-G_0)->D(G)}
{\mathcal D}(x\cdot G_0)\subseteq {\mathcal D}(G)
\eeq
непрерывны\footnote{В терминологии \cite{Schaefer} и \cite{Jarchow} это индуктивная топология, порожденная семейством отображений \eqref{DEF:D(x-cdot-G_0)->D(G)}.}. По-другому эту топологию можно описать как топологию локально выпуклой прямой суммы пространств  ${\mathcal D}(G_0:K)$:
\beq\label{DEF:D(G)-topologiya}
{\mathcal D}(G)=\LCS\text{-}\oplus_{x\cdot G_0\in G/G_0}{\mathcal D}(x\cdot G_0)
\eeq 
Можно заметить, что если в этом определении поменять левые сдвиги на правые, результат не изменится.

\btm\label{PROP:D(G)-polno-i-nasysheno} Для всякой локально компактной группы $G$ пространство ${\mathcal D}(G)$ полно и насыщено (и поэтому стереотипно).
\etm
\bpr
Это следует из теоремы \ref{PROP:D(G_0)-nasysheno} и формулы \eqref{DEF:D(G)-topologiya}.
\epr

\subsection{Пространство ${\mathcal E}(G)$.}

После того,как построено пространство ${\mathcal D}(G)$ финитных гладких функций на локально компактной группе $G$, пространство ${\mathcal E}(G)$ произвольных гладких функций строится следующим образом.

Условимся говорить, что две функции $u,v:G\to\C$ на локально  компактной группе $G$ {\it имеют одинаковый росток в точке} $a\in G$, и изображать это записью
\beq\label{DEF:rostok}
u\equiv v\ (\kern-8pt\mod a)
\eeq 
если существует окрестность $U$ точки $a$, на которой они совпадают:
\beq\label{DEF:rostok-1}
u(t)=v(t),\quad t\in U.
\eeq 
Пространство ${\mathcal E}(G)$ гладких функций на $G$ определяется следующим образом: функция $u:G\to\C$ считается лежащей в ${\mathcal E}(G)$ (то есть гладкой функцией), если для любой точке $a\in G$ существует функция $v\in{\mathcal D}(G)$, имеющая с $u$ одинаковый росток: 
\beq\label{DEF:E(G)}
\forall a\in G\quad\exists v\in{\mathcal D}(G)\quad u\equiv v\ (\kern-8pt\mod a)
\eeq 
Эквивалентное определение выглядит так: функция $u:G\to\C$ лежит в ${\mathcal E}(G)$, если для любой функции $v\in{\mathcal D}(G)$ произведение $u\cdot v$ принадлежит ${\mathcal D}(G)$:
\beq\label{DEF:E(G)-1}
\forall v\in{\mathcal D}(G)\quad u\cdot v\in {\mathcal D}(G).
\eeq 

Пространство ${\mathcal E}(G)$ наделяется инициальной локально выпуклой топологией относительно ссемейства отображений
\beq\label{DEF:E(G)-topologiya}
u\in {\mathcal E}(G)\quad\mapsto\quad  u\cdot v\in {\mathcal D}(G),\qquad v\in {\mathcal D}(G).
\eeq

Для функции $u\in {\mathcal E}(G)$ ее {\it  носитель} $\supp u$ определяется как множество точек $a\in G$, в которых росток функции $u$ отличен от нуля:
\beq\label{DEF:supp-u}
\supp u=\{a\in G: \ u\not\equiv 0\ (\kern-8pt\mod a)\}
\eeq 
Иначе это множество определяется как наименьшее замкнутое множество $M$ в  $G$, вне которого функция $u$ равна нулю. Можно заметить, что пространство ${\mathcal D}(G)$ состоит из функций $u\in {\mathcal E}(G)$ с компактным носителем:
$$
u\in {\mathcal D}(G)\quad\Leftrightarrow\quad u\in {\mathcal E}(G) \quad \&\quad \text{$\supp u$ --- компакт.}
$$

Для функционала $f\in {\mathcal D}(G)^\star$ его {\it  носитель} $\supp f$ определяется как множество точек $a\in G$ у которых каждая окрестность $U$ содержит носитель некоторой функции $u\in{\mathcal D}(G)$, на которой функционал $f$ отличен от нуля:
\beq\label{DEF:supp-f}
\supp f=\{a\in G: \quad \forall U\owns a \quad \exists u\in{\mathcal D}(G)\quad \supp u\subseteq U\quad \& \quad f(u)\ne 0 \}.
\eeq 
Иначе это множество определяется как наименьшее замкнутое множество $M$ в  $G$, вне которого функционал $f$ равен нулю (то есть равен нулю на всякой функции $u$, носитель которой не пересекается с $M$). Можно заметить, что пространство ${\mathcal E}(G)^\star$ состоит из функционалов $f\in {\mathcal D}(G)^\star$ с компактным носителем:
$$
f\in {\mathcal E}(G)^\star\quad\Leftrightarrow\quad f\in {\mathcal D}(G)^\star \quad \&\quad \text{$\supp u$ --- компакт.}
$$
Если $F\subseteq {\mathcal D}(G)^\star$ --- какое-то множество функционалов, то его носитель определяется как замыкание объединения носителей его элементов:
\beq\label{DEF:supp-F}
\supp F=\overline{\bigcup_{f\in F}\supp f}
\eeq

Пространство ${\mathcal E}(G)$ образует алгебру относительно поточечного умножения. По свойствам эта алгебра очень сильно похожа на обычную алгебру гладких функций на группе Ли. В частности, в ${\mathcal E}(G)$  справедливы теорема о разбиении единицы, о существовании базиса в касательном пространстве, о разложении дифференциальных операторов и т.д. (см. \cite{Bruhat,Akbarov-glad-strukt,Akbarov-diff-geom,Akbarov-cotang-bundle}). Мы здесь отметим только два важных результата: о стереотипности пространства ${\mathcal E}(G)$ и о том, что оно образует стереотипную алгебру относительно поточечного умножения.

\btm\label{PROP:E(G)-polno-i-nasysheno} Для всякой локально компактной группы $G$ пространство ${\mathcal E}(G)$ полно и насыщено (и поэтому стереотипно).
\etm
\bpr
1. Полнота пространства ${\mathcal E}(G)$ следует из способа задания топологии на нем: если $\{u_i\}$ --- направленность Коши в ${\mathcal E}(G)$, то для всякого $v\in {\mathcal D}(G)$ направленность $\{u_i\cdot v\}$ будет направленностью Коши в ${\mathcal D}(G)$, и, поскольку по теореме \ref{PROP:D(G)-polno-i-nasysheno}, ${\mathcal D}(G)$ полно, существует предел
$$
u_i\cdot v\overset{{\mathcal D}(G)}{\underset{i\to\infty}{\longrightarrow}} w.
$$
Выбрав в качестве $v\in {\mathcal D}(G)$ функцию, тождественно равную единице в окрестности фиксированной точки $a\in G$, мы можем заметить, что этот предел $w\in {\mathcal D}(G)$ не зависит от такого выбора $v$ в окрестности точки $a$. Таким образом мы можем определить $w$ локально. А этого достаточно для определения элемента пространства ${\mathcal E}(G)$. После этого нужно только проверить, что направленность $\{u_i\}$ стремится к $w$ в пространстве ${\mathcal E}(G)$.

2. Нам нужно показать, что пространство ${\mathcal E}(G)$ насыщено. Пусть $F$ --- множество линейных функционалов на ${\mathcal E}(G)$, равностепенно непрерывных на каждом вполне ограниченном множестве $S\subseteq{\mathcal E}(G)$. Поскольку ${\mathcal D}(G)$ непрерывно вкладывается в ${\mathcal E}(G)$, система функционалов $F$ будет равностепенно непрерывна и на вполне ограниченных множествах $T\subseteq{\mathcal D}(G)$. Поскольку по теореме \ref{PROP:D(G)-polno-i-nasysheno}, пространство ${\mathcal D}(G)$  насыщено, система функционалов $F$ должна быть равностепенно непрерывна на всем пространстве ${\mathcal D}(G)$. Запомним это.

3. Покажем далее, что система функционалов $F$ имеет компактный носитель $\supp F$. Предположим, что это не так, то есть что для всякого компакта $K\subseteq G$ найдется функция $u\in {\mathcal D}(G)$  с носителем, не пересекающимся с $K$,
$$
K\cap\supp u=\varnothing,
$$
и какой-то функционал $f\in F$ такие, что
$$
f(u)\ne 0.
$$
Тогда можно подобрать последовательность функций $u_k\in {\mathcal D}(G)$ и функционалов $f_k\in F$ таких, что 
\bit{

\item[1)] последовательность функций $u_k\in {\mathcal D}(G)$ локально конечна (то есть на каждом компакте $K\subseteq G$ только конечный набор функций $u_k$ отличен от нуля), и
    
\item[2)] $f_k(u_k)\ne 0$.     
}\eit

Подберем последовательность чисел $\lambda_k\in\C$ удовлетворяющих условию
\beq\label{PROOF:E(G)-polno-i-nasysheno-1}
\abs{f_n\l\sum_{k=1}^n\lambda_k\cdot u_k\r}\ge n,\qquad n\in\N.
\eeq
Это всегда можно сделать, выбрав сначала $\lambda_1$ так, чтобы
$$
\Big|f_1\l\lambda_1\cdot u_1\r\Big|=\Big|\lambda_1\cdot \underbrace{f_1(u_1)}_{\scriptsize\begin{matrix}\text{\rotatebox{90}{$\ne$}}\\ 0\end{matrix}}\Big|\ge 1,
$$
затем выбрав $\lambda_2$ так, чтобы
$$
\Big|f_2\l\lambda_1\cdot u_1+\lambda_2\cdot u_2\r\Big|=\Big|\lambda_1\cdot f_2(u_1)+\lambda_2\cdot \underbrace{f_2(u_2)}_{\scriptsize\begin{matrix}\text{\rotatebox{90}{$\ne$}}\\ 0\end{matrix}}\Big|\ge 2,
$$
и так далее. Заметим теперь, что из локальной конечности семейства функций $u_k$ следует, что ряд
$$
\sum_{k=1}^\infty\lambda_k\cdot u_k
$$
сходится в пространстве ${\mathcal E}(G)$ (просто потому что если его умножить на произвольную функцию $v\in{\mathcal D}(G)$, то у него только конечный набор слагаемых будет отличен от нуля).

Поэтому частичные суммы этого ряда 
\beq\label{PROOF:E(G)-polno-i-nasysheno-2}
\sum_{k=1}^n\lambda_k\cdot u_k,\qquad n\in\N
\eeq
вполне ограничены в пространстве ${\mathcal E}(G)$. Отсюда, в частности, следует, что функционалы $f\in F$ должны быть равномерно ограничены на множестве \eqref{PROOF:E(G)-polno-i-nasysheno-2}. Но из \eqref{PROOF:E(G)-polno-i-nasysheno-1} следует, что это не так:
$$
\abs{f_n\l\sum_{k=1}^n\lambda_k\cdot u_k\r}\ge n\underset{n\to\infty}{\longrightarrow}\infty.
$$
Это противоречие означает, что носитель $\supp F$ все-таки компактен. 

4. Теперь подберем функцию $v\in{\mathcal D}(G)$, равную единице в окрестности носителя $\supp F$:
$$
v\equiv 1\ (\kern-8pt\mod a),\qquad a\in \supp F.
$$
Тогда действие функционалов $f\in F$ на функциях $u\in{\mathcal E}(G)$ можно будет представить как композицию умножения на функцию 
$v\in{\mathcal D}(G)$, а потом как действие $f$ на $u\cdot v$:
$$
u\in{\mathcal E}(G) \ \mapsto \ u\cdot v\in{\mathcal D}(G) \ \mapsto \ f(u\cdot v)=f(u)\in\C
$$
В этой цепочке первая операция непрерывна (в силу выбора топологии на ${\mathcal E}(G)$), а вторая равностепенно непрерывна по $f\in F$ (в силу пункта 2 выше). Мы получаем, что функционалы $f\in F$ равностепенно непрервыны на пространстве ${\mathcal E}(G)$.
\epr

\btm\label{TH:E(G)-ster-algebra} Для всякой локально компактной группы $G$ пространство ${\mathcal E}(G)$ образует стереотипную алгебру относительно поточечного умножения.
\etm

Для доказательства нам понадобится две леммы. Первую из них мы доказывать не будем, просто сошлемся на статью Брюа \cite[Th\'eor\'eme 2]{Bruhat}:

\blm\label{LM:ogr-mnozh-v-D(G)} Если $G_0$ --- LP-группа, то всякое ограниченное множество $S$ в ${\mathcal D}(G_0)$ содержится в некотором подпространстве вида ${\mathcal D}(G_0:K)$, где $K\in\lambda(G_0)$ (и является ограниченным множеством в ${\mathcal D}(G_0:K)$). 
\elm

А вторую докажем:

\blm\label{LM:umnozhenie-v-D(G)} Для всякой локально компактной группы $G$ поточечное умножение а ${\mathcal D}(G_0)$ стереотипно непрерывно: для всякого вполне ограниченного множества $S\subseteq {\mathcal D}(G)$ и любой окрестности нуля $U\subseteq {\mathcal D}(G)$ найдется окрестность нуля $W\subseteq {\mathcal D}(G)$ такая что
$$
W\cdot S \subseteq U.
$$
\elm
\bpr
1. Сначала рассмотрим случай,  когда $G$ --- LP-группа. Тогда если $S\subseteq {\mathcal D}(G)$ --- вполне ограниченное множество,  то по лемме \ref{LM:ogr-mnozh-v-D(G)} $S$ содержится в некотором ${\mathcal D}(G:K)$, где $K\in\lambda(G)$. Если $U\subseteq {\mathcal D}(G)$ --- окрестность нуля, то для всякого $L\in\lambda(G)$ множество
$$
\{u\in {\mathcal D}(G:L):\ \forall v\in S \quad u\cdot v\in U \}
$$    
должно быть окрестностью нуля в ${\mathcal D}(G:L)$ (потому что $G/L$ --- группа Ли). Поэтому множество
$$
W=\{u\in {\mathcal D}(G):\ \forall v\in S \quad u\cdot v\in U \}=\bigcup_{L\in\lambda(G)} \{u\in {\mathcal D}(G:L):\ \forall v\in S \quad u\cdot v\in U \}
$$    
будет окрестностью нуля в ${\mathcal D}(G)$.

2. Если $G$ --- произвольная локально компактная группа,  то выбрав в ней открытую LP-подгруппу $G_0$, мы можем разложить пространство ${\mathcal D}(G)$ в прямую сумму \eqref{DEF:D(G)-topologiya}:
$$
{\mathcal D}(G)=\LCS\text{-}\oplus_{x\cdot G_0\in G/G_0}{\mathcal D}(x\cdot G_0)
$$
в которой поточечное умножение будет покомпонентным, и поэтому его стереотипная непрерывность  наследуется из компонент.
\epr

\bpr[Доказательство теоремы \ref{TH:E(G)-ster-algebra}.]
Здесь нужно убедиться, что поточечное умножение в ${\mathcal E}(G)$ непрерывно в стереотипном смысле.
Пусть $S\subseteq {\mathcal E}(G)$ --- вполне ограниченное множество и $\{u_i\}\subseteq {\mathcal E}(G)$ --- направленность, стремящаяся к нулю:
\beq\label{PROOF:E(G)-ster-algebra-1}
u_i\overset{{\mathcal E}(G)}{\underset{i\to\infty}{\longrightarrow}}0.
\eeq
Зафиксируем $w\in {\mathcal D}(G)$ и подберем функцию $\eta\in {\mathcal D}(G)$ так, чтобы 
$$
\eta\equiv 1\ (\kern-8pt\mod a),\qquad a\in \supp w.
$$
Тогда 
\beq\label{PROOF:E(G)-ster-algebra-2}
u_i\cdot v\cdot w=u_i\cdot\eta \cdot v\cdot w
\eeq
причем из \eqref{PROOF:E(G)-ster-algebra-1} будет следовать (по определению топологии в ${\mathcal E}(G)$), что 
\beq\label{PROOF:E(G)-ster-algebra-3}
u_i\cdot\eta \overset{{\mathcal D}(G)}{\underset{i\to\infty}{\longrightarrow}}0.
\eeq
а с другой стороны (опять по определению топологии в ${\mathcal E}(G)$) мы получим, что множество
$$
\{v\cdot w;\ v\in S\}
$$
вполне ограничено в ${\mathcal D}(G)$. 

По лемме \ref{LM:umnozhenie-v-D(G)} вместе это означает, что направленность $u_i\cdot v\cdot w$ должна равномерно по  $v\in S$ стремиться к нулю в пространстве ${\mathcal D}(G)$:
\beq\label{PROOF:E(G)-ster-algebra-4}
u_i\cdot v\cdot w=\eqref{PROOF:E(G)-ster-algebra-3}=u_i\cdot\eta \cdot v\cdot w
\overset{{\mathcal D}(G)}{\underset{i\to\infty}{\underset{v\in S}{\rightrightarrows}}}0
\eeq
Это верно для всякого $w\in {\mathcal D}(G)$ поэтому по выбору топологии в ${\mathcal E}(G)$, мы получаем, что 
направленность $u_i\cdot v$ должна равномерно по  $v\in S$ стремиться к нулю в пространстве ${\mathcal E}(G)$:
\beq\label{PROOF:E(G)-ster-algebra-5}
u_i\cdot v
\overset{{\mathcal E}(G)}{\underset{i\to\infty}{\underset{v\in S}{\rightrightarrows}}}0
\eeq
Это нам и нужно было доказать.
\epr 
\chapter{ДВОЙСТВЕННОСТЬ В ТОПОЛОГИИ}

\section{Непрерывные оболочки}

\subsection{Определение и функториальные свойства непрерывной оболочки}

\paragraph{$C^*$-полунормы.} Пусть $A$ -- инволютивная стереотипная алгебра. Непрерывная полунорма $p:A\to\R_+$ называется {\it $C^*$-полунормой}, если она удовлетворяет тождеству
\beq\label{DEF:C^*-polunorma}
p(a^\bullet\cdot a)=p(a)^2,\qquad a\in A.
\eeq
Согласно знаменитой теореме Себестьена \cite{Sebestyen}, всякая такая полунорма субмультипликативна:
$$
p(a\cdot b)\le p(a)\cdot p(b),\qquad a,b\in A.
$$
Множество всех $C^*$-полунорм на $A$ мы обозначаем $\Sn(A)$. Очевидно, что $p\in\Sn(A)$ тогда и только тогда, когда $p$ можно представить как композицию некоторого (непрерывного, инволютивного и унитального) гомоморфизма $\ph:A\to B$ в некоторую $C^*$-алгебру $B$ и нормы $\norm{\cdot}$ на $B$:
$$
p(x)=\norm{\ph(x)},\qquad x\in A.
$$

\btm\label{TH:p(a)=norm(pi(a))}
Всякая непрерывная $C^*$-полунорма $p$ на $A$ представима в виде нормы некоторого унитарного (непрерывного по норме\footnote{См. определение на с.\pageref{DEF:nepr-po-norme--unit-predst}.}) представления $\pi:A\to{\mathcal B}(X)$:
\beq\label{p(a)=norm(pi(a))}
p(a)=\norm{\pi(a)},\qquad a\in A.
\eeq
\etm
\bpr
Пространство $\Ker p=\{x\in A:\ p(x)=0\}$ является замкнутым идеалом в $A$. Полунорма $p$ пропускается через некоторую $C^*$-норму $p'$ на фактор-алгебре $A/\Ker p$. Пусть $B$ -- пополнение $A/\Ker p$ относительно $p'$. Понятно, что $B$ является $C^*$-алгеброй, поэтому ее можно изометрически вложить в некоторую алгебру вида ${\mathcal B}(X)$, где $X$ -- гильбертово пространство \cite[Теорема 3.4.1]{Murphy}. Композиция отображений $A\to A/\Ker p\to B\to {\mathcal B}(X)$ будет искомым представлением $\pi$ алгебры $A$.
\epr

\paragraph{Определение непрерывной оболочки и функториальность.}

Условимся символом ${\tt C}^*$ обозначать класс $C^*$-алгебр.

 \bit{
\item[$\bullet$] {\it Непрерывной оболочкой}\label{DEF:nepr-obolochka} $\env_{\mathcal C} A:A\to\Env_{\mathcal C} A$ инволютивной стереотипной алгебры $A$ называется ее оболочка в классе $\DEpi$ плотных эпиморфизмов категории $\InvSteAlg$ инволютивных стереотипеных алгебр относительно класса $\Mor(\InvSteAlg,{\tt C}^*)$  морфизмов в $C^*$-алгебры:
\beq\label{DEF::Env_C}
\Env_{\mathcal C} A=\Env_{{\tt C}^*}^{\DEpi}A
\eeq
 }\eit

Более подробно, под {\it непрерывным расширением} инволютивной стереотипной алгебры $A$ мы понимаем плотный эпиморфизм $\sigma:A\to A'$ инволютивных стереотипных алгебр такой, что для любой $C^*$-алгебры $B$ и инволютивного гомоморфизма $\ph:A\to B$ найдется (необходимо, единственный) гомоморфизм инволютивных стереотипных алгебр $\ph':A'\to B$, замыкающий диаграмму
\beq\label{DEF:diagr-nepr-rasshirenie}
 \xymatrix @R=2pc @C=1.2pc
 {
  A\ar[rr]^{\sigma}\ar[dr]_{\ph} & & A'\ar@{-->}[dl]^{\ph'} \\
  & B &
 }
\eeq
А {\it непрерывная оболочка} инволютивной стереотипной алгебры $A$ определяется как непрерывное расширение $\rho:A\to \Env_{\mathcal C} A$ такое, что для любого непрерывного расширения $\sigma:A\to A'$ найдется (необходимо, единственный) гомоморфизм  инволютивных стереотипных алгебр $\upsilon:A'\to \Env_{\mathcal C} A$, замыкающий диаграмму
$$
 \xymatrix @R=2pc @C=1.2pc
 {
  & A\ar[ld]_{\sigma}\ar[rd]^{\rho} &   \\
  A'\ar@{-->}[rr]_{\upsilon} &  & \Env_{\mathcal C} A
 }
$$

\btm\label{TH:Env_C-reg-obolochka}
Непрерывная оболочка $\Env_{\mathcal C}$ регулярна и согласована с проективным тензорным произведением\footnote{В смысле определений на страницах \pageref{DEF:reg-obolochka} и \pageref{DEF:obolochka-soglasovana-s-tenz-proizv}.} $\circledast$ в $\InvSteAlg$.
\etm
\bpr
1. Сначала докажем регулярность, то есть выполнение условий RE.1 - RE.5 на с.\pageref{DEF:RE.1} и \pageref{DEF:RE.5}. Пусть $\varPhi=\Mor(\InvSteAlg,{\tt C}^*)$ -- класс гомоморфизмов в $C^*$-алгебры.
\bit{
\item[RE.1:] Категория $\InvSteAlg$ инволютивных стереотипных алгебр полна, в частности, проективно полна (любой функтор из малой категории в $\InvSteAlg$ имеет проективный предел).
\item[RE.2:] По теореме \ref{TH:SMono-circledcirc-DEpi=InvSteAlg}, класс $\DEpi$ всех плотных эпиморфизмов, в которых ищется оболочка, мономорфно дополняем в категории $\InvSteAlg$.
\item[RE.3:] Категория $\InvSteAlg$ также локально мала в фактор-объектах класса $\DEpi$ (потому что ${\tt Ste}$ локально мала в фактор-объектах класса $\Epi$).
\item[RE.4:] Какой бы ни была инволютивная стереотипная алгебра $A$, всегда найдется выходящий из нее морфизм $\ph:A\to B$ в некоторую $C^*$-алгебру $B$ (например, в качестве $B$ можно взять нулевую алгебру $B=0$ и положить $\ph=0$). По определению на с.\pageref{DEF:goes-from}, это означает, что класс $\varPhi$ морфизмов в $C^*$-алгебры выходит из категории $\InvSteAlg$. Кроме того, этот класс является правым идеалом в $\InvSteAlg$ (композиция $\ph\circ\psi$ любого морфизма $\ph:A\to B$ в $C^*$-алгебру $B$ с любым другим морфизмом $\psi:A'\to A$ -- снова морфизм со значениями в $C^*$-алгебре).

\item[RE.5:] Если $\psi\circ\sigma\in\varPhi$ и $\sigma\in\DEpi$, то у композиции $\psi\circ\sigma:A\to C$ область значений $C$ есть $C^*$-алгебра, и поэтому у множителя $\psi$ область значений та же же самая, и значит $\psi\in\varPhi$. Это значит, что класс $\DEpi$ подталкивает класс $\varPhi$.
 }\eit

2. Теперь проверим согласованность с тензорным произведением $\circledast$, то есть условия T.1 и T.2 на с.\pageref{DEF:obolochka-soglasovana-s-tenz-proizv}.
\bit{
\item[T.1:] Пусть $\rho:A\to A'$ и $\sigma:B\to B'$ -- непрерывные расширения и пусть $\ph:A\circledast B\to C$ -- гомоморфизм в $C^*$-алгебру $C$. По лемме \ref{LM:ph:A-circledast-B->C} он представим в виде
$$
\ph(a\circledast b)=\alpha(a)\cdot\beta(b),\qquad a\in A,\ b\in B,
$$
где $\alpha:A\to C$ и $\beta:B\to C$ -- некоторые морфизмы стереотипных алгебр с коммутирующими образами:
$$
\alpha(a)\cdot\beta(b)=\beta(b)\cdot\alpha(a),\qquad a\in A,\ b\in B.
$$
Продолжим $\alpha$ и $\beta$ до морфизмов $\alpha'$ и $\beta'$ так, чтобы замыкались диаграммы
$$
\xymatrix @R=2.pc @C=2.0pc 
{
A\ar[dr]_{\alpha}\ar[rr]^{\rho} & & A'\ar@{-->}[dl]^{\alpha'}\\
& C &
}
\qquad
\xymatrix @R=2.pc @C=2.0pc 
{
B\ar[dr]_{\beta}\ar[rr]^{\sigma} & & B'\ar@{-->}[dl]^{\beta'}\\
& C &
}
$$
Поскольку $\rho$ и $\sigma$ -- плотные эпиморфизмы, образы $\alpha'$ и $\beta'$ также должны коммутировать:
$$
\alpha'(a)\cdot\beta'(b)=\beta'(b)\cdot\alpha'(a),\qquad a\in A',\ b\in B'.
$$
Поэтому мы можем снова воспользоваться леммой \ref{LM:ph:A-circledast-B->C} и определить морфизм
$$
\ph'(a\circledast b)=\alpha'(a)\cdot\beta'(b),\qquad a\in A',\ b\in B'.
$$
Очевидно, он будет продолжением $\ph$ (единственно возможным, потому что $\rho$ и $\sigma$ -- плотные эпиморфизмы).

\item[T.2:] Тождественное отображение $1_{\C}:\C\to\C$ является непрерывным расширением (например, просто потому что это изоморфизм алгебр). С другой стороны, если $\rho:\C\to A'$ -- какое-то другое непрерывное расширение, то, поскольку оно должно быть плотно, оно будет сюръективно. При этом, $A'$ не может быть нулем, потому что иначе достраивалась бы диаграмма
$$
\xymatrix @R=2.pc @C=2.0pc 
{
\C\ar[rr]^{\rho}\ar[dr]_{1_{\C}} & & 0\ar@{-->}[dl] \\
& \C &
},
$$
что невозможно. Мы получаем, что $\rho:\C\to A'$ должен быть изоморфизмом алгебр, и поэтому имеет смысл диаграмма
$$
\xymatrix @R=2.pc @C=2.0pc 
{
& \C\ar[dl]_{\rho}\ar[dr]^{1_{\C}} &  \\
A'\ar[rr]^{\rho^{-1}} & & \C
},
$$
означающая, что расширение $\rho$ вкладывается в расширение $1_{\C}$.
}\eit
\epr

\bcor\label{COR:Env_C-idemp-funktor}
Непрерывную оболочку можно определить как идемпотентный ковариантный функтор из $\InvSteAlg$ в $\InvSteAlg$: существуют
\bit{
\item[1)] отображение $A\mapsto (\Env_{\mathcal C}A,\env_{\mathcal C}A)$, сопоставляющее каждой инволютивной стереотипной алгебре $A$ инволютивную стереотипную алгебру $\Env_{\mathcal C}A$ и морфизм стереотипных алгебр $\env_{\mathcal C}A:A\to\Env_{\mathcal C}A$, являющийся непрерывной оболочкой алгебры $A$, и

\item[2)] отображение $\ph\mapsto \Env_{\mathcal C}(\ph)$, сопоставляющее каждому морфизму инволютивных стереотипных алгебр $\ph:A\to B$ морфизм инволютивных стереотипных алгебр $\Env_{\mathcal C}(\ph):\Env_{\mathcal C}A\to \Env_{\mathcal C}B$, замыкающий диаграмму
\beq\label{DIAGR:funktorialnost-env_C}
\xymatrix @R=2.pc @C=5.0pc 
{
A\ar[d]^{\ph}\ar[r]^{\env_{\mathcal C}A} & \Env_{\mathcal C}A\ar@{-->}[d]^{\Env_{\mathcal C}(\ph)} \\
B\ar[r]^{\env_{\mathcal C}B} & \Env_{\mathcal C}B \\
}
\eeq
}\eit
причем выполняются тождества
\beq\label{funktorialnost-env_C-v-Ste^circledast}
\Env_{\mathcal C}(1_A)=1_{\Env_{\mathcal C}A},\quad \Env_{\mathcal C}(\beta\circ\alpha)=\Env_{\mathcal C}(\beta)\circ \Env_{\mathcal C}(\alpha),
\eeq
\beq\label{funktorialnost-env_C-v-Ste^circledast-2}
\Env_{\mathcal C}(\Env_{\mathcal C}A)=\Env_{\mathcal C}A,\quad \env_{\mathcal C}\Env_{\mathcal C}A=1_{\Env_{\mathcal C}A},
\eeq
\beq\label{Env_C(C)=C}
\Env_{\mathcal C}\C=\C
\eeq
\ecor

\bit{
\item[$\bullet$] Морфизм $\Env_{\mathcal C}(\ph)$ называется {\it непрерывной оболочкой морфизма} $\ph$.
}\eit

\noindent\rule{160mm}{0.1pt}\begin{multicols}{2}

\bex\label{PROP:nepr-obolochka-C*-algebry}
Если $A$ --- $C^*$-алгебра, то $\Env_{\mathcal C} A=A$.
\eex
\bpr
Понятно, что тождественное отображение $\id_A:A\to A$ является непрерывным расширением алгебры $A$. Если $\sigma:A\to A'$ --- какое-то другое непрерывное расширение алгебры $A$, то поскольку $A$ --- $C^*$-алгебра, в диаграмме
$$
\xymatrix @R=2.pc @C=3.0pc 
{
A\ar[rr]^{\sigma}\ar[dr]_{\id_A}& & A'\ar@{-->}[dl]^{\ph'} \\
& A &
}
$$
существует единственная пунктирная стрелка $\ph'$. Это означает, что $\id_A$ --- непрерывная оболочка.
\epr

\bex\label{EX:Env_C-C_F=C_F}
Если $F$ -- конечная группа, то ее групповую алгебру удобно обозначать $\C_F$ и представлять, как сопряженное пространство к алгебре $\C^F$ всех функций на $G$:
$$
\C_F=(\C^F)^\star.
$$
Группа $F$ действует сдвигами на пространстве функций $L_2(F)$, поэтому алгебру $\C_F$ тоже можно считать действующей на $L_2(F)$. Это действие можно понимать как вложение $\C_F$ в алгебру операторов ${\mathcal B}(L_2(F))$. При этом, ${\mathcal B}(L_2(F))$ является $C^*$-алгеброй, поэтому $\C_F$ саму можно считать $C^*$-алгеброй. Как следствие, в силу примера  \ref{PROP:nepr-obolochka-C*-algebry}, ее непрерывная оболочка должна совпадать с ней самой:
\beq\label{Env_C-C_F=C_F}
\Env_{\mathcal C} \C_F=\C_F
\eeq
\eex

\end{multicols}\noindent\rule[10pt]{160mm}{0.1pt}

\paragraph{Сеть $C^*$-фактор-отображений в категории $\InvSteAlg$.}
Конструкцию непрерывной оболочки можно описать несколько более наглядно следующим образом. Условимся окрестность нуля $U$ в инволютивной стереотипной алгебре $A$ называть {\it $C^*$-окрестностью нуля}\label{DEF:C*-neighbourhood-of-zero}, если она является прообразом единичного шара при некотором (непрерывном) гомоморфизме $D:A\to B$ в какую-нибудь $C^*$ алгебру $B$:
$$
U=\{x\in A:\ \norm{D(x)}\le 1\}
$$
Это эквивалентно тому, что $U$ является единичным шаром некоторой (непрерывной) $C^*$-полунормы $p$ на $A$:
\beq\label{DEF:C^*-okrestnost-kak-shar-C^*-polunormy}
U=\{x\in A:\ p(x)\le 1\}.
\eeq
У каждой $C^*$-окрестности нуля $U$ в $A$ ядро
$$
\Ker U=\bigcap_{\lambda>0}\lambda\cdot U
$$
совпадает с ядром гомоморфизма $D$ и поэтому является замкнутым идеалом в $A$. Рассмотрим фактор-алгебру $A/\Ker U$ и наделим ее нормой, в которой класс $U+\Ker U$ является единичным шаром. Эту алгебру $A/\Ker U$ можно считать подалгеброй в $B$ с индуцированной из этого пространства нормой, или, что то же самое, с нормой, порожденной полунормой $p$ из \eqref{DEF:C^*-okrestnost-kak-shar-C^*-polunormy}. Ее пополнение мы будем обозначать $A/U$ или $A/p$
\beq\label{A/U=(A/Ker U)^blacktriangledown}
A/U=A/p=(A/\Ker U)^\blacktriangledown
\eeq
и называть {\it фактор-алгеброй алгебры $A$ по $C^*$-окрестности нуля $U$} или {\it по $C^*$-полунорме $p$}. Понятно, что $A/U=A/p$ является $C^*$-алгеброй (и ее можно считать замкнутой подалгеброй в $B$). Соответствующее отображение
\beq\label{DEF:pi_U}
\pi_U=\pi_p:A\to A/U=A/p
\eeq
мы будем называть {\it фактор-отображением} алгебры $A$ по $C^*$-окрестности нуля $U$, или по $C^*$-полунорме $p$, или {\it $C^*$-фактор-отображением} алгебры $A$.  Норму $C^*$-алгебры $A/p$ мы обозначаем $\norm{\cdot}_{A/p}$, и тогда полунорма $p$ будет раскладываться в композицию отображений:
\beq\label{p=norm(cdot)_A/p-circ-ph_U}
\xymatrix @R=1.pc @C=4.0pc
{
A\ar@/_6ex/@{-->}[rr]_{p}\ar[r]^(.4){\pi_p} & A/p\ar[r]^(.6){\norm{\cdot}_{A/p}}& \R_+
}
\eeq

Множество всех $C^*$-окрестностей нуля в алгебре $A$ мы будем обозначать символом ${\mathcal U}_{C^*}^A$:
\beq\label{DEF:U_C^*^A}
U\in {\mathcal U}_{C^*}^A\quad\Leftrightarrow\quad U=\{x\in A:\ p(x)\le 1\},\qquad p(x^*\cdot x)=p(x)^2.
\eeq
а множество всех $C^*$-фактор-отображений алгебры $A$ мы будем обозначать символом ${\mathcal N}_{C^*}^A$:
\beq\label{DEF:N_C^*^A}
{\mathcal N}_{C^*}^A=\{\pi_U:A\to A/U;\ U\in {\mathcal U}_{C^*}^A\}
\eeq

\medskip
\centerline{\bf Свойства $C^*$-окрестностей нуля в категории $\InvSteAlg$:}

\bit{\it

\item[$1^\circ$.]\label{LM:ph=ph_U-circ-pi_U-0}
Для всякого гомоморфизма $\ph:A\to B$ в $C^*$-алгебру $B$ найдется $C^*$-окрестность нуля $U\subseteq A$ и гомоморфизм $\ph_U:A/U\to B$, замыкающий диаграмму
\beq\label{ph=ph_U-circ-pi_U-0}
 \xymatrix @R=2pc @C=1.2pc
 {
  A\ar[rr]^{\ph}\ar[dr]_{\pi_U} & & B  \\
  & A/U\ar@{-->}[ur]_{\ph_U} &
  }
\eeq

\item[$2^\circ$.]\label{LM:varkappa^U'_U-0}
Если $U$ и $U'$ -- две $C^*$-окрестности нуля в $A$, причем $U\supseteq U'$, то найдется единственный морфизм $\varkappa^{U'}_U:A/U\gets A/U'$, замыкающий диаграмму
\beq\label{varkappa^U'_U-0}
 \xymatrix @R=2pc @C=1.2pc
 {
  & A\ar[ld]_{\pi_U}\ar[rd]^{\pi_{U'}} &   \\
  A/U &  & A/ U'\ar@{-->}[ll]^{\varkappa^{U'}_U}
 }
\eeq

\item[$3^\circ$.]\label{LM:C*-okr-nulya-uporyadocheny} Пересечение $U\cap U'$ любых двух $C^*$-окрестностей нуля $U$ и $U'$ в $A$ является $C^*$-окрестностью нуля.

}\eit 

\brem\label{REM:ph=ph_U-circ-pi_U-0-*}
Свойство $1^\circ$  в этом списке означает, что система $C^*$-фактор-отображений $\pi_U:A\to A/U$ порождает класс $\varPhi$ морфизмов в $C^*$-алгебры изнутри:
\beq\label{ph=ph_U-circ-pi_U-0-*}
    {\mathcal N}\subseteq\varPhi\subseteq\Mor(\InvSteAlg)\circ {\mathcal N}.
\eeq
\erem

\bpr Здесь не вполне очевидно только свойство $3^\circ$. Вот его доказательство. Пусть $D:A\to B$ и $D':A\to B'$ -- гомоморфизмы, порождающие $U$ и $U'$.
$$
U=\{x\in A:\ \norm{D(x)}\le 1\},\qquad U'=\{x\in A:\ \norm{D'(x)}\le 1\}.
$$
Рассмотрим $C^*$-алгебру $B\oplus B'$ с нормой
$$
\norm{b\oplus b'}=\max\{\norm{b},\norm{b'}\},\qquad b\in B,\quad b'\in B'.
$$
Отображение
$$
D'':A\to B\oplus B'\quad\Big|\quad D''(x)=D(x)\oplus D'(x), \quad x\in A,
$$
будет гомоморфизмом, и для всякого $x\in A$ мы получим
$$
x\in U\cap U'\quad\Longleftrightarrow\quad \norm{D(x)}\le 1\ \&\ \norm{D''(x)}\le 1\quad\Longleftrightarrow\quad
\norm{D''(x)}=\max\{\norm{D(x)},\norm{D"(x)}\}\le 1.
$$
\epr

\btm\label{LM:A/U-set-epimorfizmov-v-InvSteAlg}
Система ${\mathcal N}_{C^*}^A$ $C^*$-фактор-отображений образует сеть эпиморфизмов\footnote{См. определение на с.\pageref{DEF:set-epimorf}.} в категории $\InvSteAlg$ инволютивных стереотипных алгебр, то есть обладает следующими свойствами:
  \bit{
\item[(a)] у всякой алгебры $A$ есть хотя бы одна $C^*$-окрестность нуля $U$, и множество всех $C^*$-окрестностей нуля в $A$ направлено относительно предпорядка
$$
U\le U'\quad\Longleftrightarrow\quad U\supseteq U',
$$

\item[(b)] для всякой алгебры $A$ система морфизмов $\varkappa_U^{U'}$ из \eqref{varkappa^U'_U-0} ковариантна, то есть для любых трех $C^*$-окрестностей нуля $U\supseteq U'\supseteq U''$ коммутативна диаграмма
$$
 \xymatrix @R=2pc @C=1.2pc
 {
  A/U &  & A/ U''\ar[ll]_{\varkappa^{U''}_U}\ar[dl]^{\varkappa^{U''}_{U'}}\\
  & A/U \ar[ul]^{\varkappa^{U'}_U} &
 }
$$
и эта система $\varkappa_U^{U'}$ обладает проективным пределом в $\InvSteAlg$;

\item[(c)] для всякого гомоморфизма $\alpha:A\gets A'$ в $\InvSteAlg$ и любой $C^*$-окрестности нуля $U$ в $A$ найдется $C^*$-окрестность нуля $U'$ в $A'$ и гомоморфизм $\alpha_U^{U'}:A/U\gets A'/U'$ такие, что коммутативна
диаграмма
 \beq\label{DIAGR:set-C^*} \xymatrix @R=2.5pc @C=4.0pc {
 A\ar[d]_{\pi_U} & A'\ar@{-->}[d]^{\pi_{U'}}\ar[l]_{\alpha} \\
A/U & A'/U'\ar@{-->}[l]^{\alpha_U^{U'}}
 } \eeq
 }\eit
\etm

В соответствии с пунктом (b) этой теоремы, существует проективный предел $\projlim_{U'\in {\mathcal U}_{C^*}^A} A/ U'$ системы $\varkappa_U^{U'}$. Как следствие, существует единственная стрелка $\pi:A\to \projlim_{U'\in {\mathcal U}_{C^*}^A} A/ U'$ в  $\InvSteAlg$, замыкающая все диаграммы
\beq\label{A-to-leftlim-A/U-0}
 \xymatrix @R=2pc @C=1.2pc
 {
  & A\ar[ld]_{\pi_U}\ar@{-->}[rd]^{\pi} &   \\
  A/U &  & \projlim_{U'\in {\mathcal U}_{C^*}^A} A/ U'\ar[ll]^{\varkappa_U}
 }
\eeq
Образ $\pi(A)$ отображения $\pi$ является (инволютивной подалгеброй и) подпространством в стереотипном пространстве $\projlim_{U'\in {\mathcal U}_{C^*}^A} A/ U'$. Поэтому оно порождает некое непосредственное подпространство в $\projlim_{U'\in {\mathcal U}_{C^*}^A} A/ U'$, или оболочку $\Env\pi(A)$ \cite[(4.68)]{Akbarov-De-Gruyter-I}, то есть наибольшее стереотипное пространство содержащееся в  $\projlim_{U'\in {\mathcal U}_{C^*}^A} A/ U'$ и имеющее $\pi(A)$ плотным подпространством. 
Иными словами, $\Env\pi(A)$ --- узловой образ морфизма $\pi$:
$$
\Env\pi(A)=\Im_\infty \pi
$$
Обозначим через $\rho:A\to \Env\pi(A)$ поднятие морфизма $\pi$ в $\Env\pi(A)$.
\beq\label{rho:A-to-Env(pi(A))}
 \xymatrix @R=2pc @C=1.2pc
 {
  A\ar[rr]^{\pi}\ar@{-->}[rd]_{\rho} & & \projlim_{U'\in {\mathcal U}_{C^*}^A} A/ U'   \\
   & \Env\pi(A)=\Im_\infty \pi\ar[ru]_{\iota} & 
 }
\eeq
(здесь $\iota$ --- естественное погружение подпространства $\Env\pi(A)$ в объемлющее пространство $\projlim_{U'\in {\mathcal U}_{C^*}^A} A/ U'$).

\btm\label{TH:opisanie-nepr-obolochki} Морфизм $\rho:A\to \Env\pi(A)$ является непрерывной оболочкой алгебры $A$:
  \beq\label{Env_EA=Env-pi(A)}
\Env_{\mathcal C} A=\Env\pi(A)=\Im_\infty \pi.
  \eeq
(а морфизм $\iota:\Env\pi(A)\to \projlim_{U'\in {\mathcal U}_{C^*}^A} A/ U'$ --- естественным морфизмом оболочки $\Env_{\mathcal C}A$ в проективный предел $\projlim_{0\gets U'} A/ U'$).
\etm
\bpr
В силу \eqref{ph=ph_U-circ-pi_U-0-*}, система $C^*$-фактор-отображений $\pi_U:A\to A/U$ порождает класс $\varPhi$ морфизмов в $C^*$-алгебры изнутри. С другой стороны, по теореме \ref{TH:SMono-circledcirc-DEpi=InvSteAlg}, класс $\DEpi$ всех плотных эпиморфизмов, в которых ищется оболочка, мономорфно дополняем в категории $\InvSteAlg$. Отсюда по теореме \ref{TH:funktorialnost-pri-seti-Epi-i-dolonyaemosti}, оболочки алгебры $A$ относительно классов $\varPhi$ и ${\mathcal N}$ совпадают между собой и представляют собой морфизм $\rho$.
\epr

\brem
Непрерывную оболочку $\rho:A\to\Env_{\mathcal E}A$ можно представляеть себе как композицию элементов $\red_\infty$ и $\coim_\infty$ узлового разложения морфизма $\pi:A\to\projlim_{U'\in {\mathcal U}_{C^*}^A} A/ U'$ в категории $\tt Ste$ стереотипных пространств (не алгебр!):
  \beq\label{C-envelope=im_infty-lim N_X}
\env_{\mathcal C} A=\red_\infty\pi\circ\coim_\infty\pi.
  \eeq
Наглядно это изображается диаграммой
  \beq\label{DIAGR:C-envelope=im_infty-lim N_X}
\xymatrix @R=3.pc @C=4.0pc 
{
A\ar[d]_{\coim_\infty\pi}\ar@{-->}[drrr]^{\env_{\mathcal C} A}\ar[rrr]^{\pi=\projlim_{U'\in {\mathcal U}_{C^*}^A}\pi_{U'}} &&& \projlim_{U'\in {\mathcal U}_{C^*}^A} A/U' &   \\
\Coim_\infty\pi\ar[rrr]_{\red_\infty\pi} &&&  \Im_\infty \pi \ar[u]_{\im_\infty\pi}\ar@{=}[r] & \Env_{\mathcal C}A
}
 \eeq
\erem

\paragraph{Непрерывная оболочка как функтор на $\InvSteAlg^0$.}

Пусть $\InvSteAlg^0$ обозначает категорию инволютивных стереотипных алгебр (с единицей, но) в которой морфизмами считаются непрерывные инволютивные линейные мультипликативные отображения (без требования, чтобы они сохраняли единицу):
\beq\label{DEF:InvSteAlg^0}
\ph:A\to B\quad \Big|\quad \ph(x\cdot y)=\ph(x)\cdot\ph(y) \quad \&\quad \ph(x^\bullet)=\ph(x)^\bullet,\quad x,y\in A.
\eeq
Как и в случае категории $\InvSteAlg$, под {\it плотным эпиморфизмом} в $\InvSteAlg^0$ мы понимаем такой морфизм (непрерывное инволютивное линейное мультипликативное отображение) $\ph:A\to B$, у которого множество значений $\ph(A)$ плотно в области значений:
$$
\overline{\ph(A)}=B.
$$
Класс всех таких морфизмов мы будем обозначать $\DEpi^0$.

\blm\label{LM:dense-epi-v-invstealg^0}
Классы плотных эпиморфизмов $\DEpi^0$ и $\DEpi$ совпадают:
\beq\label{DEpi^0=DEpi}
\DEpi^0=\DEpi
\eeq
(иными словами, всякий плотный эпиморфизм $\ph:A\to B$ в категории $\InvSteAlg^0$ автоматически сохраняет единицу).
\elm
\bpr
Пусть $1_A$ --- единица в алгебре $A$. Покажем, что $\ph(1_A)$ --- единица в алгебре $B$. Зафиксируем элемент $b\in B$. Поскольку $\ph$ --- плотное отображение, найдется направленность $a_i\in A$ такая, что
$$
\ph(a_i)\underset{i\to\infty}{\longrightarrow} b.
$$
Теперь мы получаем:
$$
\ph(1_A)\cdot b=\ph(1_A)\cdot\lim_{i\to\infty}\ph(a_i)=\lim_{i\to\infty}\Big(\ph(1_A)\cdot\ph(a_i)\Big)=
\lim_{i\to\infty}\ph(1_A\cdot a_i)=\lim_{i\to\infty}\ph(a_i)=b.
$$
И точно так же
$$
b\cdot\ph(1_A)=\lim_{i\to\infty}\ph(a_i)\cdot\ph(1_A)=\lim_{i\to\infty}\Big(\ph(a_i)\cdot\ph(1_A)\Big)=
\lim_{i\to\infty}\ph(a_i\cdot 1_A)=\lim_{i\to\infty}\ph(a_i)=b.
$$
\epr

\blm\label{LM:downarrow-DEpi}
Класс плотных эпиморфизмов $\DEpi^0=\DEpi$ мономорфно дополняем в категории $\InvSteAlg^0$:
\beq\label{downarrow-DEpi}
{^\downarrow\DEpi^0}\circledcirc\DEpi^0=\InvSteAlg^0
\eeq
\elm
\bpr
Это доказывается так же, как и для $\InvSteAlg$. Берем морфизм $\ph:A\to B$ в $\InvSteAlg^0$ (то есть непрерывное линейное мультипликативное отображение). По \cite[Theorem 4.2.26]{Akbarov-De-Gruyter-I}, $\ph$, как морфизм в категории $\Ste$ стереотипных пространств, имеет узловое разложение
\beq\label{downarrow-DEpi-1}
\xymatrix @R=2.pc @C=5.0pc 
{
A\ar[r]^{\ph}\ar[d]_{\coim_\infty\ph} & B \\
\Coim_\infty\ph\ar[r]_{\red_\infty\ph} & \Im_\infty\ph\ar[u]_{\im_\infty\ph}
}
\eeq

1. Покажем сначала, что каждая вершина в этой диаграмме --- стереотипная алгебра, а каждая стрелка --- мультипликативное отображение. 

a) По построению, пространство $\Coim_\infty\ph$ должно быть стереотипной алгеброй (как инъективный предел трансфинитной последовательности фактор-алгебр), а морфизм $\coim_\infty\ph$ должен быть мультипликативен (как предел трансфинитной последовательности мультипликативных морфизмов), поэтому  $\coim_\infty\ph:A\to\Coim_\infty\ph$ --- морфизм в категории $\InvSteAlg^0$.

b) Опять по построению, пространство $\Im_\infty\ph$ должно быть стереотипной алгеброй, возможно без единицы (как проективный предел трансфинитной последовательности подалгебр), а морфизм $\im_\infty\ph$ должен быть мультипликативен (как предел трансфинитной последовательности мультипликативных морфизмов), поэтому  $\im_\infty\ph:\Im_\infty\ph\to B$ --- морфизм в категории $\InvSteAlg^0$. 

c) Дальше для вычислений желательно ввести какие-нибудь обозначения для стрелок. Пусть
$$
\pi=\coim_\infty\ph,\qquad \rho=\red_\infty\ph,\qquad \sigma=\im_\infty\ph.
$$
Про $\pi$ и $\sigma$ мы уже доказали, что это мультипликативные отображения. Поэтому для любых $x,y\in A$ мы получим:
$$
\sigma(\rho(\pi(x)\cdot\pi(y)))=\sigma(\rho(\pi(x\cdot y)))=\ph(x\cdot y)=\ph(x)\cdot\ph(y)=
\sigma(\rho(\pi(x)))\cdot \sigma(\rho(\pi(y)))=\sigma(\rho(\pi(x))\cdot \rho(\pi(y)))
$$
То есть
$$
\sigma(\rho(\pi(x)\cdot\pi(y)))=\sigma(\rho(\pi(x))\cdot \rho(\pi(y))),\qquad x,y\in A
$$
Поскольку $\sigma$ --- инъекция, мы можем ее выбросить:
$$
\rho(\pi(x)\cdot\pi(y))=\rho(\pi(x))\cdot \rho(\pi(y)),\qquad x,y\in A.
$$
Поскольку $\pi$ --- плотный эпиморфизм, мы можем заменить $\pi(x)$ и $\pi(y)$ на произвольные объекты, лежащие в $\Coim_\infty\ph$:
$$
\rho(x\cdot y)=\rho(x)\cdot \rho(y),\qquad x,y\in \Coim_\infty\ph.
$$

2. Итак, что $\Coim_\infty\ph$ и $\Im_\infty\ph$ --- стереотипные алгебры, возможно без единицы, а $\coim_\infty\ph$, $\red_\infty\ph$ и $\im_\infty\ph$ --- мультипликативные отображения. Покажем, что алгебры $\Coim_\infty\ph$ и $\Im_\infty\ph$ содержат единицу. Это следует из доказательства леммы \ref{LM:dense-epi-v-invstealg^0}: поскольку  
$\coim_\infty\ph:A\to\Coim_\infty\ph$ --- плотный эпиморфизм, в его области значений $\Coim_\infty\ph$ единицей будет элемент $\coim_\infty\ph(1_A)$. Точно так же, поскольку  
$\red_\infty\ph:\Coim_\infty\ph\to\Im_\infty\ph$ --- плотный эпиморфизм, в его области значений $\Im_\infty\ph$ единицей будет элемент $\coim_\infty\ph(1_{\Coim_\infty\ph})$.

3. Мы поняли, что \eqref{downarrow-DEpi-1} можно считать диаграммой в категории $\InvSteAlg^0$. Теперь если положить
$$
\e=\red_\infty\ph\circ\coim_\infty\ph, \qquad \sigma=\im_\infty\ph
$$
то мы получим
$$
\ph=\sigma\circ\e,
$$
где $\e\in\DEpi^0$, а $\sigma\in\SMono$. То есть, мономорфным дополнением к классу $\DEpi^0$ в категории $\InvSteAlg^0$ является класс морфизмов, являющихся строгими мономорфизмами в $\Ste$.
\epr

Как и в категории $\InvSteAlg$, в $\InvSteAlg^0$ мы можем определить оболочку в классе всех плотных эпиморфизмов $\DEpi^0$ относительно класса $\varPhi^0$ всех морфизмов (необязательно, сохраняющих единицу) в $C^*$-алгебры. Эту оболочку естественно называть {\it непрерывной оболочкой в} $\InvSteAlg^0$, и ее желательно как-нибудь обозначить, например, так:
\beq\label{DEF::Env^0_C}
\Env^0_{\mathcal C} A=\Env_{{\tt C}^*}^{\DEpi^0}A
\eeq
Важное для нас наблюдение заключается в том, что такая оболочка будет совпадать с определенной в \eqref{DEF::Env^0_C} непрерывной оболочкой в  $\InvSteAlg$:

\btm\label{TH:Env^0_C-A=Env_C-A}
Для всякой стереотипной алгебры $A$ ее непрерывная оболочка в $\InvSteAlg^0$ совпадает с ее непрерывной оболочкой в  $\InvSteAlg$  
\beq\label{Env^0_C-A=Env_C-A}
\Env^0_{\mathcal C}A=\Env_{\mathcal C}A
\eeq
\etm

\bpr 1. Сначала покажем, что соответствующие расширения совпадают. 

a) Пусть $\sigma:A\to A'$ --- расширение в категории $\InvSteAlg^0$ в классе $\DEpi^0$ относительно морфизмов в $C^*$-алгебры. Тогда, во-первых, по лемме \ref{LM:dense-epi-v-invstealg^0}, $\sigma$ --- морфизм и в категории $\InvSteAlg$. А, во-вторых, если $\ph:A\to B$ --- морфизм в $\InvSteAlg$ в $C^*$-алгебру, то $\ph$ --- морфизм и в категории $\InvSteAlg^0$, поэтому существует морфизм $\ph':A'\to B$ в $\InvSteAlg^0$, замыкающий диаграмму
\beq\label{TH:Env^0_C-A=Env_C-A-1}
 \xymatrix @R=2pc @C=1.2pc
 {
  A\ar[rr]^{\sigma}\ar[dr]_{\ph} & & A'\ar@{-->}[dl]^{\ph'} \\
  & B &
 }
\eeq
Опять по лемме \ref{LM:dense-epi-v-invstealg^0}, $\ph'$ должен быть морфизмом в $\InvSteAlg$. А поскольку $\sigma$ --- эпиморфизм, такой морфизм $\ph'$ будет единственным.

b) Наоборот, пусть $\sigma:A\to A'$ --- расширение в категории $\InvSteAlg$ в классе $\DEpi$ относительно морфизмов в $C^*$-алгебры. Тогда, во-первых, $\sigma$ --- морфизм и в категории $\InvSteAlg^0$. А, во-вторых, если $\ph:A\to B$ --- морфизм в $\InvSteAlg^0$ в $C^*$-алгебру, то по лемме \ref{LM:dense-epi-v-invstealg^0}, $\ph$ --- морфизм и в категории $\InvSteAlg$, поэтому существует морфизм $\ph':A'\to B$ в $\InvSteAlg$, замыкающий диаграмму \eqref{TH:Env^0_C-A=Env_C-A-1}.
Опять очевидно, $\ph'$ должен быть морфизмом в $\InvSteAlg^0$. А поскольку $\sigma$ --- эпиморфизм, такой морфизм $\ph'$ будет единственным.

2. Теперь покажем, что оболочки совпадают. 

a) Пусть $\rho:A\to E$ --- оболочка в категории $\InvSteAlg^0$  в классе $\DEpi^0$ относительно морфизмов в $C^*$-алгебры. Тогда, во-первых, по лемме \ref{LM:dense-epi-v-invstealg^0}, $\rho$ --- морфизм и в категории $\InvSteAlg$. А, во-вторых, если $\sigma:A\to A'$ --- расширение в категории $\InvSteAlg$, то $\sigma$ --- морфизм и в категории $\InvSteAlg^0$, поэтому существует морфизм $\upsilon:A'\to E$ в $\InvSteAlg^0$, замыкающий диаграмму
\beq\label{TH:Env^0_C-A=Env_C-A-2}
 \xymatrix @R=2pc @C=1.2pc
 {
  & A\ar[ld]_{\sigma}\ar[rd]^{\rho} &   \\
  A'\ar@{-->}[rr]_{\upsilon} &  & E
 }
\eeq
Опять по лемме \ref{LM:dense-epi-v-invstealg^0}, $\upsilon$ должен быть морфизмом в $\InvSteAlg$. А поскольку $\sigma$ --- эпиморфизм, такой морфизм $\upsilon$ будет единственным.

b) Наоборот, пусть $\rho:A\to E$ --- оболочка в категории $\InvSteAlg$ в классе $\DEpi$ относительно морфизмов в $C^*$-алгебры. Тогда, во-первых, $\rho$ --- морфизм и в категории $\InvSteAlg^0$. А, во-вторых, если $\sigma:A\to A'$ --- морфизм в $\InvSteAlg^0$ в $C^*$-алгебру, то по лемме \ref{LM:dense-epi-v-invstealg^0}, $\sigma$ --- морфизм и в категории $\InvSteAlg$, поэтому существует морфизм $\upsilon:A'\to E$ в $\InvSteAlg$, замыкающий диаграмму \eqref{TH:Env^0_C-A=Env_C-A-2}.
Опять очевидно, $\upsilon$ должен быть морфизмом в $\InvSteAlg^0$. А поскольку $\sigma$ --- эпиморфизм, такой морфизм $\upsilon$ будет единственным.
\epr

\btm\label{TH:Env^0_C-idemp-funktor}
Непрерывную оболочку можно доопределить как идемпотентный ковариантный функтор из $\InvSteAlg^0$ в $\InvSteAlg^0$: для всякого морфизма $\ph:A\to B$ в $\InvSteAlg^0$ существует единственный морфизм $\Env^0_{\mathcal C}(\ph):\Env_{\mathcal C}A\to \Env_{\mathcal C}B$ в $\InvSteAlg^0$, замыкающий диаграмму
\beq\label{DIAGR:funktorialnost-env^0_C}
\xymatrix @R=2.pc @C=5.0pc 
{
A\ar[d]^{\ph}\ar[r]^{\env_{\mathcal C}A} & \Env_{\mathcal C}A\ar@{-->}[d]^{\Env^0_{\mathcal C}(\ph)} \\
B\ar[r]^{\env_{\mathcal C}B} & \Env_{\mathcal C}B \\
}
\eeq
причем отображение $\ph\mapsto \Env^0_{\mathcal C}(\ph)$ обладает свойствами
\beq\label{funktorialnost-env^0_C-v-Ste^circledast}
\Env^0_{\mathcal C}(1_A)=1_{\Env_{\mathcal C}A},\quad \Env^0_{\mathcal C}(\beta\circ\alpha)=\Env^0_{\mathcal C}(\beta)\circ \Env^0_{\mathcal C}(\alpha),
\eeq
\etm

Доказательство этого утверждения опирается на понятие сети эпиморфизмов определенное в \eqref{DEF:pi_U}:
\beq\label{DEF:pi_U-v-InvSteAlg^0}
\pi_U=\pi_p:A\to A/U=A/p
\eeq
Здесь важное наблюдение состоит в том, что эта система морфизмов обладает в категории $\InvSteAlg^0$ теми же свойствами, что и в категории $\InvSteAlg$. В частности, свойства $1^\circ$ на с.\pageref{LM:ph=ph_U-circ-pi_U-0} не меняются (и в них везде пунктирные стрелки существуют, потому что они существуют в категории $\Ste$).

\medskip
\centerline{\bf Свойства $C^*$-окрестностей нуля в категории $\InvSteAlg^0$:}

\bit{\it

\item[$1^\circ$.]\label{LM:ph=ph_U-circ-pi_U-0-v-InvSteAlg^0}
Для всякого морфизма $\ph:A\to B$ в категории $\InvSteAlg^0$ в $C^*$-алгебру $B$ найдется $C^*$-окрестность нуля $U\subseteq A$ и морфизм $\ph_U:A/U\to B$ в категории $\InvSteAlg^0$, замыкающий диаграмму
\beq\label{ph=ph_U-circ-pi_U-0-v-InvSteAlg^0}
 \xymatrix @R=2pc @C=1.2pc
 {
  A\ar[rr]^{\ph}\ar[dr]_{\pi_U} & & B  \\
  & A/U\ar@{-->}[ur]_{\ph_U} &
  }
\eeq

\item[$2^\circ$.]\label{LM:varkappa^U'_U-0-v-InvSteAlg^0}
Если $U$ и $U'$ -- две $C^*$-окрестности нуля в $A$, причем $U\supseteq U'$, то найдется единственный морфизм $\varkappa^{U'}_U:A/U\gets A/U'$ в категории $\InvSteAlg^0$, замыкающий диаграмму
\beq\label{varkappa^U'_U-0-v-InvSteAlg^0}
 \xymatrix @R=2pc @C=1.2pc
 {
  & A\ar[ld]_{\pi_U}\ar[rd]^{\pi_{U'}} &   \\
  A/U &  & A/ U'\ar@{-->}[ll]^{\varkappa^{U'}_U}
 }
\eeq

\item[$3^\circ$.]\label{LM:C*-okr-nulya-uporyadocheny-v-InvSteAlg^0} Пересечение $U\cap U'$ любых двух $C^*$-окрестностей нуля $U$ и $U'$ в $A$ является $C^*$-окрестностью нуля.

}\eit

А теорема \ref{LM:A/U-set-epimorfizmov-v-InvSteAlg} приобретает вид

\btm\label{LM:A/U-set-epimorfizmov-v-InvSteAlg^0}
Система $\pi_U:A\to A/U$ $C^*$-фактор-отображений образует сеть эпиморфизмов\footnote{См. определение на с.\pageref{DEF:set-epimorf}.} в категории $\InvSteAlg^0$, то есть обладает следующими свойствами:
  \bit{
\item[(a)] у всякой алгебры $A$ есть хотя бы одна $C^*$-окрестность нуля $U$, и множество всех $C^*$-окрестностей нуля в $A$ направлено относительно предпорядка
$$
U\le U'\quad\Longleftrightarrow\quad U\supseteq U',
$$

\item[(b)] для всякой алгебры $A$ система морфизмов $\varkappa_U^{U'}$ из \eqref{varkappa^U'_U-0} ковариантна, то есть для любых трех окрестностей нуля $U\supseteq U'\supseteq U''$ коммутативна диаграмма
$$
 \xymatrix @R=2pc @C=1.2pc
 {
  A/U &  & A/ U''\ar[ll]_{\varkappa^{U''}_U}\ar[dl]^{\varkappa^{U''}_{U'}}\\
  & A/U \ar[ul]^{\varkappa^{U'}_U} &
 }
$$
и эта система $\varkappa_U^{U'}$ обладает проективным пределом в $\InvSteAlg^0$;

\item[(c)] для всякого морфизма $\alpha:A\gets A'$ в $\InvSteAlg^0$ и любой $C^*$-окрестности нуля $U$ в $A$ найдется $C^*$-окрестность нуля $U'$ в $A'$ и морфизм $\alpha_U^{U'}:A/U\gets A'/U'$ в $\InvSteAlg^0$ такие, что коммутативна
диаграмма
 \beq\label{DIAGR:set-C^*} \xymatrix @R=2.5pc @C=4.0pc {
 A\ar[d]_{\pi_U} & A'\ar@{-->}[d]^{\pi_{U'}}\ar[l]_{\alpha} \\
A/U & A'/U'\ar@{-->}[l]^{\alpha_U^{U'}}
 } \eeq
 }\eit
\etm

\bpr[Доказательство теоремы \ref{TH:Env^0_C-idemp-funktor}.]
Система $C^*$-фактор-отображений $\pi_U:A\to A/U$ порождает класс $\varPhi^0$ морфизмов в $C^*$-алгебры изнутри в силу свойства $1^\circ$ на с.\ref{LM:ph=ph_U-circ-pi_U-0-v-InvSteAlg^0}:
    $$
    {\mathcal N}\subseteq\varPhi^0\subseteq\Mor(\InvSteAlg^0)\circ {\mathcal N}.
    $$
С другой стороны, по лемме \ref{LM:downarrow-DEpi}, класс $\DEpi^0$ всех плотных эпиморфизмов, в которых ищется оболочка, мономорфно дополняем в категории $\InvSteAlg^0$. Отсюда по теореме \ref{TH:funktorialnost-pri-seti-Epi-i-dolonyaemosti}, оболочки алгебры $A$ относительно классов $\varPhi^0$ и ${\mathcal N}$ совпадают между собой и представляют собой морфизм $\rho$ из \eqref{rho:A-to-Env(pi(A))} (и это дает еще одно доказательство теоремы \ref{TH:Env^0_C-A=Env_C-A}). Но, что для нас важнее здесь, по следствиям (b) и (c) теоремы \ref{TH:funktorialnost-pri-seti-Epi-i-dolonyaemosti}, для всякого морфизма $\ph:A\to B$ в $\InvSteAlg^0$ существует единственный морфизм $\Env^0_{\mathcal C}(\ph):\Env_{\mathcal C}A\to \Env_{\mathcal C}B$ в $\InvSteAlg^0$, замыкающий диаграмму \eqref{DIAGR:funktorialnost-env^0_C}, и отображение $\ph\mapsto \Env^0_{\mathcal C}(\ph)$ обладает свойствами
\eqref{funktorialnost-env^0_C-v-Ste^circledast}.  
\epr

\subsection{Непрерывные алгебры и непрерывное тензорное произведение}

\paragraph{Непрерывные алгебры.}

Инволютивную стереотипную алгебру $A$ мы называем {\it непрерывной алгеброй}, если она является полным объектом в категории $\InvSteAlg$ инволютивных стереотипных алгебр относительно непрерывной оболочки $\Env_{\mathcal C}$ в смысле определения на с.\pageref{DEF:polnye-objekty}, то есть удовлетворяет следующим равносильным условиям\footnote{Условие (iv) в этом списке --- действие предложения \ref{PROP:polnota-pri-podtalkivanii}.}:
\bit{
\item[(i)] всякое непрерывное расширение $\sigma:A\to A'$ алгебры $A$ является изоморфизмом в $\InvSteAlg$;

\item[(ii)] локальная единица $1_A:A\to A$ является непрерывной оболочкой $A$;

\item[(iii)] непрерывная оболочка алгебры $A$ является изоморфизмом в $\InvSteAlg$: $\env_{\mathcal C} A\in\Iso$,

\item[(iv)] алгебра $A$ изоморфна в $\InvSteAlg$ непрерывной оболочке какой-нибудь алгебры $B$: 
$$
A\cong\Env_{\mathcal C}B.
$$
}\eit\noindent 
Класс всех непрерывных алгебр мы обозначаем ${\mathcal C}\text{-}{\tt Alg}$. Он образует полную подкатегорию в категории $\InvSteAlg$.

\noindent\rule{160mm}{0.1pt}\begin{multicols}{2}

\bex\label{EX:C(M)-nepr-alg}
Пусть $M$ --- паракомпактное локально компактное топологическое пространство. Алгебра ${\mathcal C}(M)$ непрерывных функций $u:M\to \C$ с поточечными алгебраическими операциями и топологией равномерной сходимости на компактах в $M$, непрерывна.
\eex
\bpr
Мы сошлемся здесь на доказательство следующего примера \ref{EX:C(M,B(X))-nepr-alg}, в котором описывается более общий случай. 
\epr

\bex\label{EX:C(M,B(X))-nepr-alg}
Пусть $M$ --- паракомпактное локально компактное топологическое пространство, $X$ --- конечномерное гильбертово пространство над $\C$, ${\mathcal B}(X)$ --- алгебра операторов на $X$, ${\mathcal C}(M,{\mathcal B}(X))$ --- алгебра непрерывных отображений $u:M\to {\mathcal B}(X)$ с поточечными алгебраическими операциями и топологией равномерной сходимости на компактах в $M$.
Покажем, что алгебра ${\mathcal C}(M,{\mathcal B}(X))$ непрерывна.
\eex
\bpr
Пусть 
$$
\sigma:{\mathcal C}(M,{\mathcal B}(X))\to A'
$$
--- произвольное непрерывное расширение алгебры ${\mathcal C}(M,{\mathcal B}(X))$. Нам нужно показать, что $\sigma$ --- изоморфизм.

1. Для всякого компакта $S\subseteq M$ алгебра ${\mathcal C}(S,{\mathcal B}(X))$ является банаховой. Поэтому оператор ограничения
\begin{multline*}
\ph_S:{\mathcal C}(M,{\mathcal B}(X))\to {\mathcal C}(S,{\mathcal B}(X))\quad\Big|\\
 \ph(u)=u\big|_S,\qquad u\in {\mathcal C}(M,{\mathcal B}(X)),
\end{multline*}
является морфизмом в банахову алгебру. Знначит, он должен единственным образом продолжаться до некоторого морфизма $\ph_S'$ вдоль расширения $\sigma$:
\beq\label{PROOF:EX:C(M,B(X))-nepr-alg-1}
 \xymatrix @R=2pc @C=1.2pc
 {
  {\mathcal C}(M,{\mathcal B}(X))\ar[rr]^{\sigma}\ar[dr]_{\ph_S} & & A'\ar@{-->}[dl]^{\ph_S'} \\
  & {\mathcal C}(S,{\mathcal B}(X)) &
 }
\eeq

2. Покажем, что возникающая система морфизмов $\ph_S':A'\to {\mathcal C}(S,{\mathcal B}(X))$ является проективным конусом контравариантной системы морфизмов, действующих как ограничения с большего компакта на меньший:
$$
\ph_S^T:{\mathcal C}(T,{\mathcal B}(X))\to {\mathcal C}(S,{\mathcal B}(X)),\qquad S\subseteq T.
$$
Действительно, для любых компактов $S\subseteq T$ мы получим диаграмму
\beq\label{PROOF:EX:C(M,B(X))-nepr-alg-2}
 \xymatrix @R=2.5pc @C=2pc
 {
 {\mathcal C}(M,{\mathcal B}(X))\ar[rr]^{\sigma}\ar@/_5ex/[ddr]_{\ph_S}\ar[dr]_{\ph_T} & & A' \ar@/^5ex/[ddl]^{\ph_S'}\ar[dl]^{\ph_T'}\\
  & {\mathcal C}(T,{\mathcal B}(X))\ar[d]_{\ph_S^T} & \\
  & {\mathcal C}(S,{\mathcal B}(X)) &
 }
 \eeq
в которой периметр и верхний внутренний треугольник коммутативны, потому что это варианты диаграммы \eqref{PROOF:EX:C(M,B(X))-nepr-alg-1}, а левый внутренний треугольник коммутативен, потому в нем все морфизмы --- операторы ограничения. Поскольку вдобавок $\sigma$ --- эпиморфизм, отсюда следует, что правый внутренний треугольник 
\beq\label{PROOF:EX:C(M,B(X))-nepr-alg-3}
 \xymatrix @R=2.5pc @C=2pc
 {
 & & A' \ar@/^5ex/[ddl]^{\ph_S'}\ar[dl]^{\ph_T'}\\
  & {\mathcal C}(T,{\mathcal B}(X))\ar[d]_{\ph_S^T} & \\
  & {\mathcal C}(S,{\mathcal B}(X)) &
 }
 \eeq
--- тоже коммутативен.

3. Из коммутативности \eqref{PROOF:EX:C(M,B(X))-nepr-alg-3} следует, что должен существовать единственный морфизм
$$
\upsilon: A'\to {\mathcal C}(M,{\mathcal B}(X)),
$$
замыкающий все диаграммы
\beq\label{PROOF:EX:C(M,B(X))-nepr-alg-4}
 \xymatrix @R=2pc @C=1.2pc
 {
  {\mathcal C}(M,{\mathcal B}(X))\ar[dr]_{\ph_S} & & A'\ar@{-->}[ll]_{\upsilon}\ar[dl]^{\ph_S'} \\
  & {\mathcal C}(S,{\mathcal B}(X)) &
 }
\eeq
Композиция $\upsilon\circ\sigma$ должна замыкать все диаграммы
\beq\label{PROOF:EX:C(M,B(X))-nepr-alg-5}
 \xymatrix @R=2pc @C=1.2pc
 {
  {\mathcal C}(M,{\mathcal B}(X))\ar@{-->}[rr]^{\upsilon\circ\sigma}\ar[dr]_{\ph_S} & & {\mathcal C}(M,{\mathcal B}(X))\ar[dl]^{\ph_S} \\
  & {\mathcal C}(S,{\mathcal B}(X)) &
 }
\eeq
и, поскольку таким свойством обладает только тождественный морфизм $1_{{\mathcal C}(M,{\mathcal B}(X))}$, мы получаем, что 
\beq\label{PROOF:EX:C(M,B(X))-nepr-alg-6}
\upsilon\circ\sigma=1_{{\mathcal C}(M,{\mathcal B}(X))}
\eeq
Умножив это равенство на $\sigma$ слева, мы получим
$$
\sigma\circ\upsilon\circ\sigma=\sigma\circ 1_{{\mathcal C}(M,{\mathcal B}(X))}=\sigma=1_{A'}\circ\sigma.
$$
Поскольку здесь $\sigma$ --- эпиморфизм, мы можем на него сократить, и получим
\beq\label{PROOF:EX:C(M,B(X))-nepr-alg-7}
\sigma\circ\upsilon=1_{A'}
\eeq

4. Равенства \eqref{PROOF:EX:C(M,B(X))-nepr-alg-6} и \eqref{PROOF:EX:C(M,B(X))-nepr-alg-7} вместе означают, что $\sigma$ --- изоморфизм, а это нам и нужно было доказать. 
\epr
\end{multicols}\noindent\rule[10pt]{160mm}{0.1pt}

Из \eqref{ph=ph_U-circ-pi_U-0-*} и теоремы \ref{TH:kriterij-obolochki-v-term-polnyh-objektov} следует

\btm\label{TH:kriterij-nepr-obolochki-v-term-nepr-algebr}
Для всякого морфизма инволютивных стереотипных алгебр\footnote{Здесь не предполагается, что $\eta$ лежит в $\DEpi$.} $\eta:A\to S$ в произвольную непрерывную алгебру $S$ следующие условия эквивалентны:
\bit{

\item[(i)] морфизм $\eta:A\to S$ является непрерывной оболочкой;

\item[(ii)] для любого морфизма инволютивных стереотипных алгебр\footnote{Здесь не предполагается, что $\theta$ лежит в $\DEpi$.} $\theta:A\to E$ в произвольную непрерывную алгебру $E$ найдется единственный морфизм инволютивных стереотипных алгебр $\upsilon:S\to E$, замыкающий диаграмму
 \beq\label{kriterij-nepr-obolochki-v-term-nepr-algebr}
\xymatrix @R=2.pc @C=2.0pc 
{
& A \ar[dl]_{\eta} \ar[dr]^{\theta} &\\
S \ar@{-->}[rr]^{\upsilon} & & E
}
 \eeq
}\eit
\etm

Из теоремы \ref{TH:proizvedenie-polnyh-obyektov} следует

\btm\label{TH:proizvedenie-neprer-algebr}\phantom{.}
  \bit{
\item[(i)] Произведение $\prod_{i\in I}A_i$ любого семейства $\{A_i;\ i\in I\}$ непрерывных алгебр является непрерывной алгеброй.

\item[(ii)] Проективный предел $\projlim_{i\in I}A_i$ любой ковариантной (контравариантной) системы $\{A_i,\iota^j_i;\ i\in I\}$ непрерывных алгебр является непрерывной алгеброй.
 }\eit
 \etm

\paragraph{Непрерывное тензорное произведение инволютивных стереотипных алгебр.}

Пусть $\Env_{\mathcal C}$ -- функтор непрерывной оболочки, определенный в следствии \ref{COR:Env_C-idemp-funktor}. Для любых двух инволютивных стереотипных алгебр $A$ и $B$ определим их {\it непрерывное тензорное произведение} равенством
\beq\label{C/circledast}
A\overset{\mathcal C}{\circledast} B=\Env_{\mathcal C}(A\circledast B)
\eeq
Каждой паре $\alpha:A\to A'$ и $\beta:B\to B'$ морфизмов алгебр можно поставить в соответствие морфизм
\beq\label{DEF:alpha-C/circledast-beta}
\alpha\overset{\mathcal C}{\circledast}\beta=\Env_{\mathcal C}(\alpha\circledast\beta):
A\overset{\mathcal C}{\circledast} B\to A'\overset{\mathcal C}{\circledast} B'.
\eeq
Наконец, каждой паре элементов $a\in A$, $b\in B$ можно поставить в соответствие элементарный тензор
\beq\label{DEF:a-overset-C-circledast-b}
a\overset{\mathcal C}{\circledast}b=\env_{\mathcal C}(a\circledast b)
\eeq

\blm\label{LM:polnota-a-overset-C-circledast-b}
Элементарные тензоры $a\overset{\mathcal C}{\circledast}b$, $a\in A$, $b\in B$, полны в $A\overset{\mathcal C}{\circledast} B$.
\elm
\bpr
Тензоры $a\circledast b$ полны в $A\circledast B$, а образ $\env_{\mathcal C}$ плотен в $A\overset{\mathcal C}{\circledast} B$.
\epr

Ниже нам понадобится следующая конструкция. Для любых двух полунорм $q\in \Sn(A)$ и $r\in \Sn(B)$ рассмотрим полунорму $q\otimes_{\max}r$ на $A\circledast B$, определенную как композиция отображений
\beq\label{DEF:q-otimes_max-r}
\xymatrix @R=1.pc @C=4.0pc
{
A\circledast B\ar@/_6ex/@{-->}[rrr]_{q\otimes_{\max}r}\ar[r]^(.4){\pi_q\circledast\pi_r} & A/q\circledast B/r=A/q\widehat{\otimes}B/r\ar[r]^{\tau} & A/q\underset{\max}{\otimes}B/r\ar[r]^(.6){\norm{\cdot}_{\max}}& \R_+
}
\eeq
где $\pi_q:A\to A/q$ и $\pi_r:B\to B/r$ -- $C^*$-фактор-отображения, определенные в \eqref{DEF:pi_U}, $\pi_q\circledast\pi_r$ -- их сильное стереотипное тензорное произведение, $\tau$ -- естественное отображение тензорных произведений, и $\norm{\cdot}_{\max}$ -- норма максимального тензорного произведения $C^*$-алгебр.

\blm\label{LM:p(z)-le-q-otimes_max-r(z)}
Всякая $C^*$-полунорма $p$ на $A\circledast B$ подчинена некоторой $C^*$-полунорме $q\otimes_{\max}r$ на $A\circledast B$:
\beq\label{p(z)-le-q-otimes_max-r(z)}
p(x)\le (q\otimes_{\max}r)(x),\qquad x\in A\circledast B
\eeq
\elm
\bpr
Обозначим $C=(A\circledast B)/p$. Тогда $p$ представляет собой композицию отображений
$$
\xymatrix @R=1.pc @C=4.0pc
{
A\circledast B\ar@/_6ex/@{-->}[rr]_{p}\ar[r]^(.4){\pi_p} & (A\circledast B)/p=C\ar[r]^(.6){\norm{\cdot}_C}& \R_+
}
$$
где $\pi_p$ -- проекция из \eqref{DEF:pi_U}, а $\norm{\cdot}_C$ -- норма на $C$. Рассмотрим полунормы
$$
q(a)=p(a\circledast 1_B),\qquad r(b)=p(1_A\circledast b),\qquad a\in A, \ b\in B.
$$
Покажем, что гомоморфизм $\pi_p:A\circledast B\to C$ продолжается до некоторого гомоморфизма $\sigma:A/q\otimes_{\max} B/r\to C$:
\beq\label{PROOF:p(z)-le-q-otimes_max-r(z)}
\xymatrix @R=3.pc @C=4.0pc
{
A\circledast B\ar[dr]_{\pi_p}\ar[r]^(.4){\pi_q\circledast\pi_r} & A/q\circledast B/r\ar@{-->}[d]^{\rho}\ar[r]^{\tau}& A/q\otimes_{\max} B/r\ar@{-->}[dl]^{\sigma} \\
& C &
}
\eeq
Для этого рассмотрим гомоморфизмы
$$
\alpha:A\to C\quad\Big|\quad \alpha(a)=\pi_p(a\circledast 1_B),\qquad a\in A,
$$
$$
\beta:B\to C\quad\Big|\quad \beta(b)=\pi_p(1_A\circledast b),\qquad b\in B.
$$
Из того, что полунормы $q$ и $r$ являются ограничениями полунормы $p$ при отображениях $a\mapsto a\circledast 1_B$ и $b\mapsto 1_A\circledast b$, следует, что $\alpha$ и $\beta$ продолжаются до гомоморфизмов $A/q\to C$ и $B/r\to C$:
$$
\xymatrix @R=2.pc @C=2.0pc
{
A\ar[dr]_{\alpha}\ar[rr]^{\pi_q} & & A/q\ar@{-->}[dl]^{\alpha_q} \\
& C &
}
\qquad
\xymatrix @R=2.pc @C=2.0pc
{
B\ar[dr]_{\beta}\ar[rr]^{\pi_r} & & B/r\ar@{-->}[dl]^{\beta_r} \\
& C &
}
$$

С другой стороны, по лемме \ref{LM:ph:A-circledast-B->C}, образы отображений $\alpha$ и $\beta$ коммутируют,
$$
\alpha(a)\cdot\beta(b)=\beta(b)\cdot\alpha(a),\qquad a\in A,\ b\in B,
$$
а $\pi_q$ и $\pi_r$ -- плотные эпимофризмы. Отсюда следует, что образы отображений $\alpha_q$ и $\beta_r$ также коммутируют:
$$
\alpha_q(a')\cdot\beta_r(b')=\beta_r(b')\cdot\alpha_q(a'),\qquad a'\in A/q,\ b'\in B/r.
$$
Как следствие, опять по лемме \ref{LM:ph:A-circledast-B->C}, найдется гомоморфизм $\rho$, замыкающий левый внутренний треугольник в \eqref{PROOF:p(z)-le-q-otimes_max-r(z)}. Затем, поскольку $A/q$ и $B/r$ -- $C^*$-алгебры, гомоморфизм $\rho$ продолжается до гомоморфизма $\sigma$ на $A/q\otimes_{\max}B/r$.

Наконец, после того, как построено отображение $\sigma$ в \eqref{PROOF:p(z)-le-q-otimes_max-r(z)}, мы получаем, что
полунорма $z\mapsto \norm{\sigma(z)}_C$ должна быть подчинена полунорме $z\mapsto\norm{z}_{\max}$ (потому что гомоморфизм $C^*$-алгебр не увеличивает норму, \cite[Теорема 2.1.7]{Murphy}):
$$
\norm{\sigma(z)}_C\le \norm{z}_{\max},\qquad z\in A/q\otimes_{\max}B/r.
$$
Отсюда уже следует \eqref{p(z)-le-q-otimes_max-r(z)}:
$$
p(x)=\norm{\pi_p(x)}_C=\norm{\sigma(\tau((\pi_q\circledast\pi_r)(x)))}_C
\le \norm{\tau((\pi_q\circledast\pi_r)(x))}_{\max}=(q\otimes_{\max}r)(x),\qquad x\in A\circledast B.
$$
\epr

\btm\label{TH:C-circledast->odot} Пусть $A$ и $B$ --- непрерывные алгебры.
Тогда существует единственное линейное непрерывное отображение $\eta_{A,B}:A\overset{\mathcal C}{\circledast}B\to A\odot B$, замыкающее диаграмму
    \beq\label{DIAGR:eta_(A,B)}
\xymatrix @R=2.pc @C=5.0pc 
{
A\circledast B\ar[dr]_{\env_{\mathcal C}{A\circledast B}\quad}\ar[rr]^{@_{A,B}} & & A\odot B \\
& A\overset{\mathcal C}{\circledast}B\ar[ur]_{\eta_{A,B}} & \\
}
\eeq
а система отображений  $\eta_{A,B}:A\overset{\mathcal C}{\circledast}B\to A\odot B$ является естественным преобразованием функтора $(A,B)\mapsto A\overset{\mathcal C}{\circledast}B\in{\mathcal C}\text{-}{\tt Alg}$ в функтор $(A,B)\mapsto A\odot B\in {\tt Ste}$.
\etm
\bpr
1. Сначала построим систему отображений $\eta_{A,B}$. 
Пусть для любых полунорм $q\in\Sn(A)$ и $r\in\Sn(B)$ символы $A/q$ и $B/r$ обозначают $C^*$-фактор-алгебры алгебр $A$ и $B$ по этим полунормам. Рассмотрим преобразование Гротендика
\beq\label{PROOF:C-circledast->odot-0}
\xymatrix @R=1.pc @C=5.0pc
{
A/q\circledast B/r\ar[r]^{@_{A/q,B/r}}  & A/q\odot B/r.
}
\eeq
Переходя к проективному пределу, мы получаем цепочку морфизмов
\begin{multline*}
\projlim\limits_{p\in\Sn(A\circledast B)}\kern-4pt(A\circledast B)/p =(\text{лемма \ref{LM:p(z)-le-q-otimes_max-r(z)}})= \kern-10pt\projlim\limits_{q\in\Sn(A),\ r\in \Sn(B)}\kern-10pt A\circledast B/(q\underset{\max}{\otimes}r)=\\=
\kern-10pt\projlim\limits_{q\in\Sn(A),\ r\in \Sn(B)}\kern-10pt A/q\circledast B/r\stackrel{\projlim\limits @_{A/q,B/r}}{\longrightarrow}  \projlim\limits_{q\in\Sn(A),\ r\in \Sn(B)}\kern-10pt A/q\odot B/r=\\ =\cite[(4.153)]{Akbarov-De-Gruyter-I}= 
 \projlim\limits_{q\in\Sn(A)} A/q\odot \projlim\limits_{r\in \Sn(B)} B/r 
\end{multline*}
Дополним ее до цепочки
\begin{multline}\label{PROOF:C-circledast->odot-1}
A\overset{\mathcal C}{\circledast} B= \Env_{\mathcal C}(A\circledast B)
\stackrel{\projlim\limits_{p\in\Sn(A\circledast B)}\rho_p}{\longrightarrow} \\
\stackrel{\projlim\limits_{p\in\Sn(A\circledast B)}\rho_p}{\longrightarrow}
\projlim\limits_{p\in\Sn(A\circledast B)}\kern-4pt(A\circledast B)/p =(\text{лемма \ref{LM:p(z)-le-q-otimes_max-r(z)}})= \kern-10pt\projlim\limits_{q\in\Sn(A),\ r\in \Sn(B)}\kern-10pt A\circledast B/(q\underset{\max}{\otimes}r)=\\=
\kern-10pt\projlim\limits_{q\in\Sn(A),\ r\in \Sn(B)}\kern-10pt A/q\circledast B/r\stackrel{\projlim\limits @_{A/q,B/r}}{\longrightarrow}  \projlim\limits_{q\in\Sn(A),\ r\in \Sn(B)}\kern-10pt A/q\odot B/r=\\ =\cite[(4.153)]{Akbarov-De-Gruyter-I}= 
 \projlim\limits_{q\in\Sn(A)} A/q\odot \projlim\limits_{r\in \Sn(B)} B/r 
\end{multline}
и заметим, что элементарный тензор $a\overset{\mathcal C}{\circledast} b$ по ней переходит в тензор $a\odot b$, лежащий в пространстве $A\odot  B$, являющимся подпространством в $\projlim\limits_{q\in\Sn(A)} A/q\odot \projlim\limits_{r\in \Sn(B)} B/r$:
$$
a\overset{\mathcal C}{\circledast} b\in A\overset{\mathcal C}{\circledast} B\mapsto  a\odot b\in A\odot B\subseteq \projlim\limits_{q\in\Sn(A)} A/q\odot \projlim\limits_{r\in \Sn(B)} B/r.
$$
Это следует из того, что сам тензор $a\overset{\mathcal C}{\circledast} b\in A\overset{\mathcal C}{\circledast} B$ есть образ тензора $a\circledast b\in A\circledast B$, а для любых $q\in\Sn(A)$ и $r\in\Sn(B)$ образ $a_q\circledast b_r$ тензора $a\circledast b$ в $A/q\circledast B/r$ превращается под действием отображения \eqref{PROOF:C-circledast->odot-0} в тензор   
$a_q\odot b_r$ (а такая нитка $\{a_q\odot b_r;\ q,r\}$ в пространстве $\projlim\limits_{q\in\Sn(A)} A/q\odot \projlim\limits_{r\in \Sn(B)} B/r$ соответствует тензору $a\odot b\in A\odot B\subseteq \projlim\limits_{q\in\Sn(A)} A/q\odot \projlim\limits_{r\in \Sn(B)} B/r$).

С другой стороны, поскольку $A$ --- непрерывная алгебра, она должна быть непосредственным подпространством в 
$\projlim\limits_{q\in\Sn(A)} A/q$, а поскольку $B$ --- непрерывная алгебра, она должна быть непосредственным подпространством в $\projlim\limits_{r\in\Sn(B)} B/r$:
$$
A\osubarr \projlim\limits_{q\in\Sn(A)} A/q, \quad B\osubarr \projlim\limits_{r\in\Sn(B)} B/r
$$
По теореме \ref{TH:odot-sohranyaet-osubarr} это означает, что $A\odot B$ --- непосредственное подпространство в  
$\projlim\limits_{q\in\Sn(A)} A/q\odot \projlim\limits_{r\in\Sn(B)} B/r$:
$$
A\odot B\osubarr \projlim\limits_{q\in\Sn(A)} A/q\odot \projlim\limits_{r\in \Sn(B)} B/r.
$$
Мы получаем, что образ пространства $A\overset{\mathcal C}{\circledast} B$ под действием отображения  \eqref{PROOF:C-circledast->odot-1} лежит в $A\odot B$, которое является непосредственным подпространством в $\projlim\limits_{q\in\Sn(A)} A/q\odot \projlim\limits_{r\in\Sn(B)} B/r$. По теореме \ref{TH:Z-subarr-X,Y-osubarr-X,Z-subseteq-Y} это значит, что отображение  \eqref{PROOF:C-circledast->odot-1} поднимается до некоторого оператора 
$$
A\overset{\mathcal C}{\circledast} B \longrightarrow A\odot B
$$
Этот оператор и есть $\eta_{A,B}$.

2. Теперь покажем, что отображения $\eta_{A,B}$ представляют собой морфизм функторов. Пусть $\alpha:A\to A'$ и $\beta:B\to B'$ -- морфизмы непрерывных алгебр. Рассмотрим диаграмму
\beq\label{PROOF:C-circledast->odot}
\xymatrix @R=2.pc @C=5.0pc 
{
A\circledast B\ar[dd]_{\alpha\circledast\beta}\ar[dr]_{\env_{\mathcal C}{A\circledast B}\quad}\ar[rr]^{@_{A,B}} & & A\odot B\ar[dd]_{\alpha\odot\beta} \\
& A\overset{\mathcal C}{\circledast}B\ar[dd]_(.2){\alpha\overset{\mathcal C}{\circledast}\beta}\ar[ur]_{\eta_{A,B}} & \\
A'\circledast B'\ar[dr]_{\env_{\mathcal C}{A'\circledast B'}\quad}\ar[rr]^(.7){@_{A',B'}}|!{[r]}\hole & & A'\odot B' \\
& A'\overset{\mathcal C}{\circledast}B'\ar[ur]_{\eta_{A',B'}} & \\
}
\eeq
В ней верхнее и нижнее основания коммутативны, потому что это диаграммы \eqref{DIAGR:eta_(A,B)}, дальняя боковая грань коммутативна, потому что это диаграмма, выражающая естественность преобразования Гротендика $@$, а левая ближняя боковая грань коммутативна, потому что это диаграмма функториальности оболочки \eqref{DIAGR:funktorialnost-env_varPhi^Epi-v-kat-s-uzl-razl-1-*}. С другой стороны, из леммы \ref{LM:polnota-a-overset-C-circledast-b} следует, что отображение $\env_{\mathcal C}{A\circledast B}$ является  эпиморфизмом стереотипных пространств. Вместе все это означает, что оставшаяся правая боковая грань в диаграмме
\eqref{PROOF:C-circledast->odot} также коммутативна, что и нужно было доказать.
\epr

\brem\label{REM:eta(a-overset-C-circledast-b)=a-odot-b}
Из диаграммы \eqref{DIAGR:eta_(A,B)} следует, что при отображении $\eta_{A,B}$ элементарные тензоры $a\overset{\mathcal C}{\circledast}b$ переходят в элементарные тензоры $a\odot b$:
\beq\label{eta(a-overset-C-circledast-b)=a-odot-b}
\eta_{A,B}(a\overset{\mathcal C}{\circledast}b)=a\odot b.
\eeq
\erem

Из теорем \ref{TH:Env_C-reg-obolochka} и \ref{TH:sushestvovanie-tenz-proizv-v-L} следует

\btm\label{TH:C-obolochka=monoidalnyi-funktor} Формула
\eqref{C/circledast}
 определяет в ${\mathcal C}\text{-}{\tt Alg}$ тензорное произведение, превращающее ${\mathcal C}\text{-}{\tt Alg}$ в моноидальную категорию, а функтор непрерывной оболочки $\Env_{\mathcal C}$ является моноидальным функтором из моноидальной категории $(\InvSteAlg,\circledast)$ инволютивных стереотипных алгебр в моноидальную категорию $({\mathcal C}\text{-}{\tt Alg},\overset{\mathcal C}{\circledast})$ непрерывных алгебр. Соответствующий морфизм бифункторов
 $$
\Big((A,B)\mapsto \Env_{\mathcal C}(A)\overset{\mathcal C}{\circledast} \Env_{\mathcal C}(B)\Big)\overset{E^{\circledast}}{\rightarrowtail} \Big((A,B)\mapsto \Env_{\mathcal C}(A\circledast B)\Big)
 $$
определяется формулой
$$
E^\circledast_{A,B}=\Env_{\mathcal C}(\env_{\mathcal C}A\circledast \env_{\mathcal C}B)^{-1}:\Env_{\mathcal C}(A)\overset{\mathcal C}{\circledast} \Env_{\mathcal C}(B)=\Env_{\mathcal C}(\Env_{\mathcal C}(A)\circledast \Env_{\mathcal C}(B))\to \Env_{\mathcal C}(A\circledast B),
$$
а морфизмом $E^{\C}$ в ${\mathcal C}\text{-}{\tt Alg}$, переводящий единичный объект $\C$ категории ${\mathcal C}\text{-}{\tt Alg}$ в образ $\Env_{\mathcal C}(\C)$ единичного объекта $\C$ категории $\InvSteAlg$, будет локальная единица:
$$
E^{\C}=1_{\C}:\C\to\C=\Env_{\mathcal C}(\C).
$$
\etm

\paragraph{Действие непрерывной оболочки на биалгебры.}

\btm\label{LM:koalg-v-CAlg->koalg-v-odot}
Если $A$ --- непрерывная алгебра и одновременно коалгебра в моноидальной категории $({\mathcal C}\text{-}{\tt Alg},\overset{\mathcal C}{\circledast})$  непрерывных алгебр со структурными морфизмами
$$
\varkappa:A\to A\overset{\mathcal C}{\circledast} A,\qquad \e:A\to \C,
$$
то $A$ является коалгеброй в моноидальной категории $({\tt Ste},\odot)$ стереотипных пространств со структурными морфизмами
$$
\lambda=\eta_{A,A}\circ\varkappa:A\to A\odot A,\qquad \e:A\to\C.
$$
При этом всякий морфизм $\ph:A\to B$ коалгебр в моноидальной категории $({\mathcal C}\text{-}{\tt Alg},\overset{\mathcal C}{\circledast})$, будет морфизмом $A$ в $B$ как коалгебр в моноидальной категории $({\tt Ste},\odot)$.
\etm
\bpr
1. Рассмотрим диаграмму ассоциативности для $\varkappa$
$$
\xymatrix @R=2.pc @C=4.0pc 
{
&& A\ar[dll]_{\varkappa}\ar[drr]^{\varkappa} && \\
A\overset{\mathcal C}{\circledast}A\ar[dr]_{\varkappa\overset{\mathcal C}{\circledast}1_A} && && A\overset{\mathcal C}{\circledast}A \ar[dl]^{1_A\overset{\mathcal C}{\circledast}\varkappa}\\
& (A\overset{\mathcal C}{\circledast}A)\overset{\mathcal C}{\circledast}A\ar[rr]_{\alpha_{A,A,A}} && A\overset{\mathcal C}{\circledast}(A\overset{\mathcal C}{\circledast}A)
}
$$
и дополним ее до диаграммы
$$
\xymatrix @R=1.5pc @C=4.0pc 
{
&& A\ar[dll]_{\varkappa}\ar[drr]^{\varkappa}\ar[ddddd]_(.7){1_A}|!{[dd]}\hole && \\
A\overset{\mathcal C}{\circledast}A\ar[dr]_{\varkappa\overset{\mathcal C}{\circledast}1_A}\ar[ddddd]_{\eta_{A,A}} && && A\overset{\mathcal C}{\circledast}A \ar[dl]^{1_A\overset{\mathcal C}{\circledast}\varkappa}\ar[ddddd]_{\eta_{A,A}}\\
& (A\overset{\mathcal C}{\circledast}A)\overset{\mathcal C}{\circledast}A\ar[rr]_(.7){\alpha^{\tt C}_{A,A,A}}\ar[dd]^{\eta_{A\overset{\mathcal C}{\circledast}A,A}} && A\overset{\mathcal C}{\circledast}(A\overset{\mathcal C}{\circledast}A)\ar[dd]_{\eta_{A,A\overset{\mathcal C}{\circledast}A}} & \\
&&&& \\
& (A\overset{\mathcal C}{\circledast}A)\odot A\ar[ddd]^(.7){\eta_{A,A}\odot 1_A} && A\odot (A\overset{\mathcal C}{\circledast}A)\ar[ddd]_(.7){1_A\odot\eta_{A,A}} & \\
&& \quad A\quad \ar[dll]_(.3){\lambda}|!{[d];[ll]}\hole\ar[drr]^(.3){\lambda}|!{[d];[rr]}\hole && \\
A\odot A\ar[dr]_{\lambda\odot 1_A} && && A\odot A \ar[dl]^{1_A\odot \lambda}\\
& (A\odot A)\odot A\ar[rr]_{\alpha^{\odot}_{A,A,A}} && A\odot (A\odot A) &
}
$$
В ней верхнее основание коммутативно, и нужно доказать коммутативность нижнего. Для этого нужно просто проверить коммутативность боковых граней. Две дальние боковые грани коммутативны просто потому что представляют собой определение морфизма $\lambda$. Коммутативность левой ближней боковой грани
\beq\label{DIAGR-1:koalg-v-CAlg->koalg-v-odot}
\xymatrix @R=2.pc @C=5.0pc 
{
A\overset{\mathcal C}{\circledast}A\ar[dd]_{\eta_{A,A}}\ar[r]^{\varkappa\overset{\mathcal C}{\circledast}1_A} & (A\overset{\mathcal C}{\circledast}A)\overset{\mathcal C}{\circledast}A
\ar[d]^{\eta_{A\overset{\mathcal C}{\circledast}A,A}} \\
 & (A\overset{\mathcal C}{\circledast}A)\odot A\ar[d]^{\eta_{A,A}\odot 1_A} \\
A\odot A\ar[r]_{\lambda\odot 1_A} & (A\odot A)\odot A
}
\eeq
достаточно проверить на элементарных тензорах. Зафиксируем $a,b\in A$, и представим $\varkappa(a)$ как предел направленности сумм элементарных тензоров (здесь в первый раз используется лемма \ref{LM:polnota-a-overset-C-circledast-b}):
$$
\sum_{i\in I_s}x_i^s\overset{\mathcal C}{\circledast}y_i^s\underset{s\to\infty}{\longrightarrow}\varkappa(a)
$$
Тогда двигаясь по диаграмме \eqref{DIAGR-1:koalg-v-CAlg->koalg-v-odot} элементарный тензор $a\overset{\mathcal C}{\circledast}b$ даст следующие элементы:
$$
\xymatrix @R=2.pc @C=2.0pc 
{
a\overset{\mathcal C}{\circledast}b\ar@{|->}[dd]_{\eta_{A,A}}\ar@{|->}[rrr]^{\varkappa\overset{\mathcal C}{\circledast}1_A} & & & \varkappa(a)\overset{\mathcal C}{\circledast}b\ar@{=}[r]
 & \lim\limits_{s\to\infty}(\sum\limits_{i\in I_s}x_i^s\overset{\mathcal C}{\circledast}y_i^s)\overset{\mathcal C}{\circledast}b\ar@{|->}[d]^{\eta_{A\overset{\mathcal C}{\circledast}A,A}} \\
& & & & \lim\limits_{s\to\infty}(\sum\limits_{i\in I_s}x_i^s\overset{\mathcal C}{\circledast}y_i^s)\odot b\ar@{|->}[d]^{\eta_{A,A}\odot 1_A} \\
a\odot b\ar@{|->}[rr]^(.4){\lambda\odot 1_A} & &\eta_{A,A}(\varkappa(a))\odot b\ar@{=}[r] & \eta_{A,A}\left(\lim\limits_{s\to\infty}(\sum\limits_{i\in I_s}x_i^s\overset{\mathcal C}{\circledast}y_i^s)\right)\odot b
\ar@{=}[r]& \lim\limits_{s\to\infty}(\sum\limits_{i\in I_s}x_i^s\odot y_i^s)\odot b
}
$$
Поскольку тензоры $a\overset{\mathcal C}{\circledast}b$ полны в $A\overset{\mathcal C}{\circledast}A$ (здесь второй раз используется лемма \ref{LM:polnota-a-overset-C-circledast-b}), это доказывает коммутативность \eqref{DIAGR-1:koalg-v-CAlg->koalg-v-odot}.

Тем же приемом доказывается коммутативность правой ближней боковой грани
$$
\xymatrix @R=2.pc @C=5.0pc 
{
 A\overset{\mathcal C}{\circledast}(A\overset{\mathcal C}{\circledast}A)
\ar[d]_{\eta_{A,A\overset{\mathcal C}{\circledast}A}} & A\overset{\mathcal C}{\circledast}A\ar[dd]_{\eta_{A,A}}\ar[l]_{\varkappa\overset{\mathcal C}{\circledast}1_A} \\
A\odot(A\overset{\mathcal C}{\circledast}A)\ar[d]_{1_A\odot\eta_{A,A}} & \\
A\odot (A\odot A) & A\odot A\ar[l]_{\lambda\odot 1_A}
}
$$

Для центральной ближней боковой грани
$$
\xymatrix @R=2.pc @C=5.0pc 
{
(A\overset{\mathcal C}{\circledast}A)\overset{\mathcal C}{\circledast}A\ar[r]^{\alpha^{\tt C}_{A,A,A}}
 \ar[d]_{\eta_{A\overset{\mathcal C}{\circledast}A,A}} & A\overset{\mathcal C}{\circledast}(A\overset{\mathcal C}{\circledast}A)\ar[d]_{\eta_{A,A\overset{\mathcal C}{\circledast}A}} \\
(A\overset{\mathcal C}{\circledast}A)\odot A\ar[d]_{\eta_{A,A}\odot1_A} & A\odot(A\overset{\mathcal C}{\circledast}A)\ar[d]_{1_A\odot\eta_{A,A}} \\
(A\odot A)\odot A\ar[r]^{\alpha^{\odot}_{A,A,A}} & A\odot (A\odot A)
}
$$
нужно рассмотреть тройку элементов $a,b,c\in A$. Двигаясь по ней они дадут очевидную картину:
$$
\xymatrix @R=2.pc @C=5.0pc 
{
(a\overset{\mathcal C}{\circledast}b)\overset{\mathcal C}{\circledast}c\ar@{|->}[r]^{\alpha^{\tt C}_{A,A,A}}
 \ar@{|->}[d]_{\eta_{A\overset{\mathcal C}{\circledast}A,A}} & a\overset{\mathcal C}{\circledast}(b\overset{\mathcal C}{\circledast}c)\ar@{|->}[d]_{\eta_{A,A\overset{\mathcal C}{\circledast}A}} \\
(a\overset{\mathcal C}{\circledast}b)\odot c\ar@{|->}[d]_{\eta_{A,A}\odot1_A} & a\odot(b\overset{\mathcal C}{\circledast}c)\ar@{|->}[d]_{1_A\odot\eta_{A,A}} \\
(a\odot b)\odot c\ar@{|->}[r]^{\alpha^{\odot}_{A,A,A}} & a\odot (b\odot c)
}
$$
Точно так же проверяются диаграммы для коединицы.

2. Пусть $\ph:A\to B$ -- морфизм коалгебр в моноидальной категории $({\mathcal C}\text{-}{\tt Alg},\overset{\mathcal C}{\circledast})$, то есть морфизм $A$ в $B$ как стереотипных пространств, замыкающий диаграммы
\beq\label{PROOF:LM:koalg-v-CAlg->koalg-v-odot}
\xymatrix @R=2.pc @C=5.0pc 
{
A\ar[d]_{\varkappa_A}\ar[r]^{\ph}  & B\ar[d]^{\varkappa_B} \\
A\overset{\mathcal C}{\circledast}A\ar[r]^{\ph\overset{\mathcal C}{\circledast}\ph} & B\overset{\mathcal C}{\circledast}B
}
\qquad
\xymatrix @R=2.pc @C=2.0pc 
{
A\ar[dr]_{\e_A}\ar[rr]^{\ph}  & & B\ar[dl]^{\e_B} \\
 & \C &
}
\eeq
где $\varkappa_A$, $\varkappa_B$, $\e_A$, $\e_B$ -- структурные морфизмы. Тогда левую из этих диаграмм можно будет дополнить до диаграммы
$$
\xymatrix @R=2.pc @C=5.0pc 
{
A\ar[d]_{\varkappa_A}\ar[r]^{\ph}  & B\ar[d]^{\varkappa_B} \\
A\overset{\mathcal C}{\circledast}A
\ar[d]_{\eta_{A,A}}\ar[r]^{\ph\overset{\mathcal C}{\circledast}\ph} & B\overset{\mathcal C}{\circledast}B
\ar[d]_{\eta_{B,B}} \\
A\odot A\ar[r]_{\ph\odot\ph} & B\odot B,
}
$$
которая будет коммутативна, потому что по теореме \ref{TH:C-circledast->odot} морфизмы $\eta$ представляют собой естественное преобразование функторов. Вместе с правой диаграммой в \eqref{PROOF:LM:koalg-v-CAlg->koalg-v-odot} это будет означать, что $\ph:A\to B$ -- морфизм коалгебр в $(\Ste,\odot)$.
\epr

\btm\label{TH:C-obolochka-sohranyaet-Hopfov}
Пусть $H$ --- инволютивная биалгебра в категории $({\tt Ste},\circledast)$ стереотипных пространств с коумножением $\varkappa$ и коединицей $\e$. Тогда
 \bit{
\item[(i)] непрерывная оболочка $\Env_{\mathcal C}H$ является коалгеброй в моноидальной категории $({\mathcal C}\text{-}{\tt Alg},\overset{\mathcal C}{\circledast})$ непрерывных алгебр с коумножением и коединицей
    \beq\label{varkappa^E,e^E}
    \varkappa_E=\Env_{\mathcal C}(\env_{\mathcal C}H\circledast \env_{\mathcal C}H)\circ \Env_{\mathcal C}(\varkappa),\qquad \e_E=\Env_{\mathcal C}(\e),
    \eeq

\item[(ii)] непрерывная оболочка $\Env_{\mathcal C}H$ является коалгеброй в моноидальной категории $({\tt Ste},\odot)$ стереотипных пространств с коумножением и коединицей
    \beq\label{varkappa^odot,e^odot}
    \varkappa_\odot=\eta_{\Env_{\mathcal C}H,\Env_{\mathcal C}H}\circ \Env_{\mathcal C}(\env_{\mathcal C}H\circledast \env_{\mathcal C}H)\circ \Env_{\mathcal C}(\varkappa)=
    \eta_{\Env_{\mathcal C}H,\Env_{\mathcal C}H}\circ \varkappa_E,\qquad \e_\odot=\Env_{\mathcal C}(\e),
    \eeq

\item[(iii)] морфизм $\env_{\mathcal C}H^\star:H^\star\gets (\Env_{\mathcal C}H)^\star$, сопряженный к морфизму оболочки $\env_{\mathcal C}H:H\to \Env_{\mathcal C}H$, является морфизмом стереотипных алгебр, если $(\Env_{\mathcal C}H)^\star$ рассматривается как алгебра с умножением и единицей, сопряженными к \eqref{varkappa^odot,e^odot}, а $H^\star$ -- как алгебра с умножением и единицей
    $$
    \varkappa^\star\circ @_{H^\star,H^\star},\qquad \e^\star.
    $$

 }\eit
\etm
\bpr
Инволютивная биалгебра в категории $({\tt Ste},\circledast)$ -- то же самое, что коалгебра в категории $\InvSteAlg$ инволютивных стереотипных алгебр. Поэтому по теореме \ref{TH:C-obolochka=monoidalnyi-funktor} $\Env_{\mathcal C}H$ -- коалгебра в категории $({\mathcal C}\text{-}{\tt Alg},\overset{\mathcal C}{\circledast})$ с коумножением и коединицей \eqref{varkappa^E,e^E}. После этого применяется теорема \ref{LM:koalg-v-CAlg->koalg-v-odot}, и мы получаем, что $H$ -- коалгебра в категории $({\tt Ste},\odot)$ с коумножением и коединицей \eqref{varkappa^odot,e^odot}. Остается проверить (iii). Для этого нужно заметить, что коммутативна диаграмма
$$
\xymatrix @R=1.5pc @C=4.0pc 
{
\Env_{\mathcal C}H\ar[dr]_{\Env_{\mathcal C}(\varkappa)}\ar@/^8ex/[ddrr]^{\varkappa_E}\ar@/_30ex/[dddddr]_{\varkappa_\odot} & & \\
 H\ar[u]_{\env_{\mathcal C}H}\ar[d]^{\varkappa}  & \Env_{\mathcal C}(H\circledast H)\ar[dr]_{\Env_{\mathcal C}(\env_{\mathcal C}H\circledast \env_{\mathcal C}H)\quad\qquad}  & \\
 H\circledast H\ar[ur]_{\quad\env_{\mathcal C}{H\circledast H}}\ar[dr]^{\qquad\env_{\mathcal C}H\circledast \env_{\mathcal C}H}\ar[dd]_{@_{H,H}} & &  \Env_{\mathcal C}(\Env_{\mathcal C}H\circledast \Env_{\mathcal C}H)\ar@{=}[dd] \\
 & \Env_{\mathcal C}H\circledast \Env_{\mathcal C}H\ar[ur]_{\qquad \env_{\mathcal C}{\Env_{\mathcal C}H\circledast \Env_{\mathcal C}H}}\ar[dd]^{@_{\Env_{\mathcal C}H,\Env_{\mathcal C}H}} & \\
 H\odot H\ar[dr]^(.4){\qquad\env_{\mathcal C}H\odot \env_{\mathcal C}H} & & \Env_{\mathcal C}H\overset{\mathcal C}{\circledast} \Env_{\mathcal C}H\ar[dl]^{\qquad \eta_{\Env_{\mathcal C}H,\Env_{\mathcal C}H}} \\
 & \Env_{\mathcal C}H\odot \Env_{\mathcal C}H & \\
}
$$
Переходя к сопряженным пространствам, мы получим диаграмму
$$
\xymatrix @R=1.5pc @C=4.0pc 
{
(\Env_{\mathcal C}H)^\star\ar[d]^{\env_{\mathcal C}H^\star} & & \\
 H^\star  & \Env_{\mathcal C}(H\circledast H)^\star\ar[ul]^{\Env_{\mathcal C}(\varkappa)^\star}\ar[dl]^(.4){\env_{\mathcal C}{H\circledast H}^\star}  & \\
 H^\star\odot H^\star\ar[u]_{\varkappa^\star}  & &  \Env_{\mathcal C}(\Env_{\mathcal C}H\circledast \Env_{\mathcal C}H)^\star\ar[ul]^{\Env_{\mathcal C}(\env_{\mathcal C}H\circledast \env_{\mathcal C}H)^\star\quad\qquad} \ar@/_8ex/[uull]_{\varkappa_E^\star}\ar@{=}[dd]\ar[dl]^{\qquad \env_{\mathcal C}{\Env_{\mathcal C}H\circledast \Env_{\mathcal C}H}^\star} \\
 & (\Env_{\mathcal C}H)^\star\odot (\Env_{\mathcal C}H)^\star\ar[ul]_{\qquad\env_{\mathcal C}H^\star\odot \env_{\mathcal C}H^\star}  & \\
 H^\star\circledast H^\star\ar[uu]^{@_{H^\star,H^\star}} & & (\Env_{\mathcal C}H\overset{\mathcal C}{\circledast} \Env_{\mathcal C}H)^\star \\
 & (\Env_{\mathcal C}H)^\star\circledast (\Env_{\mathcal C}H)^\star\ar@/^30ex/[uuuuul]^{\varkappa_\odot^\star}\ar[ul]_{\qquad\env_{\mathcal C}H^\star\circledast \env_{\mathcal C}H^\star} \ar[uu]_{@_{(\Env_{\mathcal C}H)^\star,(\Env_{\mathcal C}H)^\star}}\ar[ur]_{\qquad \eta_{\Env_{\mathcal C}H,\Env_{\mathcal C}H}^\star} & \\
}
$$
От нее нужно оставить левый нижний фрагмент, который и представляет собой диаграмму согласованности операций умножения:
$$
\xymatrix @R=2.pc @C=6.0pc 
{
H^\star & (\Env_{\mathcal C}H)^\star\ar[l]_{\env_{\mathcal C}H^\star} \\
H^\star\odot H^\star\ar[u]^{\varkappa^\star} & \\
H^\star\circledast H^\star\ar[u]^{@_{H^\star,H^\star}} & (\Env_{\mathcal C}H)^\star\circledast (\Env_{\mathcal C}H)^\star
\ar[l]_{\env_{\mathcal C}H^\star\circledast \env_{\mathcal C}H^\star}\ar[uu]_{\varkappa_\odot^\star}
}
$$
\epr

\btm\label{TH:C-obolochka-sohranyaet-inv-Hopfov}
Пусть $H$ -- инволютивная алгебра Хопфа в категории $({\tt Ste},\circledast)$ стереотипных пространств. Тогда
 \bit{
\item[(i)] непрерывная оболочка $\Env_{\mathcal C}H$, как коалгебра в моноидальных категориях $({\mathcal C}\text{-}{\tt Alg},\overset{\mathcal C}{\circledast})$ и $({\tt Ste},\odot)$, обладает согласованными между собой антиподом $\Env_{\mathcal C}(\sigma)$ и инволюцией $\Env_{\mathcal C}(\bullet)$, однозначно определяемыми диаграммами в категории $\tt{Ste}$
    \beq\label{C-obolochka-sohranyaet-inv-Hopfov}
    \xymatrix @R=2.pc @C=4.0pc 
{
H\ar[d]_{\sigma}\ar[r]^{\env_{\mathcal C}H} & \Env_{\mathcal C}H\ar@{-->}[d]^{\Env_{\mathcal C}(\sigma)} \\
H\ar[r]^{\env_{\mathcal C}H} & \Env_{\mathcal C}H
}
\qquad
    \xymatrix @R=2.pc @C=4.0pc 
{
H\ar[d]_{\bullet}\ar[r]^{\env_{\mathcal C}H} & \Env_{\mathcal C}H\ar@{-->}[d]^{\Env_{\mathcal C}(\bullet)} \\
H\ar[r]^{\env_{\mathcal C}H} & \Env_{\mathcal C}H
}
\eeq

\item[(ii)] морфизм $\env_{\mathcal C}H^\star:H^\star\gets (\Env_{\mathcal C}H)^\star$, сопряженный к морфизму оболочки $\env_{\mathcal C}H:H\to \Env_{\mathcal C}H$, является инволютивным гомоморфизмом стереотипных алгебр над $\circledast$, если $H^\star$ и $(\Env_{\mathcal C}H)^\star$ наделяются структурой сопряженных инволютивных алгебр к инволютивным коалгебрам с антиподом $H$ и $\Env_{\mathcal C}H$ по \cite[p.657, $4^\circ$]{Akbarov-De-Gruyter-I}.
 }\eit
\etm
\bpr
1. Обозначим через $H^{\op}$ алгебру $H$ с противоположным умножением:
$$
\mu^{\op}=\mu\circ\diamond.
$$
(здесь $\diamond$ -- преобразование симметрии в категории $\Ste$). Пусть $\op_H:H\to H^{\op}$ обозначает тождественное отображение $H$ на себя (но так, чтобы область значений считалась алгеброй с противоположным умножением). Это будет антигомоморфизм алгебр. Непрерывные оболочки алгебр $H$ и $H^{\op}$ также связаны естественным антигомоморфизмом, который мы условимся обозначать $\Env_{\mathcal C}(\op)$:
\beq\label{E(op_H)}
    \xymatrix @R=2.pc @C=4.0pc 
{
H\ar[d]_{\op_H}\ar[r]^{\env_{\mathcal C}H} & \Env_{\mathcal C}H\ar@{-->}[d]^{\Env_{\mathcal C}(\op_H)} \\
H^{\op}\ar[r]^{\env_{\mathcal C}H} & \Env_{\mathcal C}(H^{\op})
}
\eeq
Это можно доказать, заметив, что при факторизации по любой $C^*$-полунорме $p$ коммутативна диаграмма
$$
    \xymatrix @R=2.pc @C=4.0pc 
{
H\ar[d]_{\op_H}\ar[r]^{\pi_p} & H/p\ar@{-->}[d]^{\op_{H/p}} \\
H^{\op}\ar[r]^{\pi_p} & (H/p)^{\op}
}
$$
Отсюда следует, что существует естественный антигомоморфизм между проективными пределами:
$$
    \xymatrix @R=2.pc @C=4.0pc 
{
H\ar[d]_{\op_H}\ar[r]^{\projlim_p\pi_p} & \projlim_p H/p\ar@{-->}[d] \\
H^{\op}\ar[r]^{\projlim_p\pi_p} & \projlim_p (H/p)^{\op}
}
$$
А затем переходя к непосредственным подпространствам, порожденным образами $H$ и $H^{\op}$, мы получаем искомую пунктирную стрелку в \eqref{E(op_H)}.

После того, как $\Env_{\mathcal C}(\op_H)$ определено для любой стереотипной алгебры $H$, отображение $\Env_{\mathcal C}(\sigma)$ можно задать формулой
$$
\Env_{\mathcal C}(\sigma)=\Env_{\mathcal C}(\op_{H^{\op}})\circ \Env_{\mathcal C}(\op_H\circ\sigma)
$$
или, нагляднее, диаграммой
$$
    \xymatrix @R=2.pc @C=8.0pc 
{
H\ar[d]^{\sigma}\ar[r]^{\env_{\mathcal C}H}\ar@/_7ex/[ddd]_{\sigma} & \Env_{\mathcal C}H\ar[dd]_{\Env_{\mathcal C}(\op_H\circ\sigma)}\ar@{-->}@/^7ex/[ddd]^{\Env_{\mathcal C}(\sigma)} \\
H\ar[d]^{\op_H} &  \\
H^{\op}\ar[d]^{\op_{H^{\op}}}\ar[r]^{\env_{\mathcal C}{H^{\op}}} & \Env_{\mathcal C}(H^{\op})\ar[d]_{\Env_{\mathcal C}(\op_{H^{\op}})}\\
H\ar[r]^{\env_{\mathcal C}{H}} & \Env_{\mathcal C}H
}
$$

2. Существование $\Env_{\mathcal C}(\bullet)$ доказывается тем же приемом, что и существование $\Env_{\mathcal C}(\op_H)$. Сначала для произвольной $C^*$-полунормы $p$ отмечается диаграмма
$$
    \xymatrix @R=2.pc @C=4.0pc 
{
H\ar[d]_{\bullet}\ar[r]^{\pi_p} & H/p\ar@{-->}[d]^{\bullet} \\
H\ar[r]^{\pi_p} & H/p
}
$$
А затем переходом к проективному пределу и непосредственному подпространству, порожденному образом $H$, получаем правую диаграмму в \eqref{C-obolochka-sohranyaet-inv-Hopfov}.

3. В свойстве (ii) нужно только проверить, что $\env_{\mathcal C}H^\star$ сохраняет инволюцию. Для $a\in H$ и $f\in (\Env_{\mathcal C}H)^\star$ мы получаем:
\begin{multline*}
\env_{\mathcal C}H^\star(f^\bullet)(a)=f^\bullet(\env_{\mathcal C}H(a))=\overline{f(\Env_{\mathcal C}(\sigma)(\env_{\mathcal C}H(a))^\bullet)}=
\overline{f((\bullet\circ \Env_{\mathcal C}(\sigma)\circ \env_{\mathcal C}H)(a))}=\\=
\overline{f((\bullet\circ \env_{\mathcal C}H\circ\sigma)(a))}=
\overline{f((\env_{\mathcal C}H\circ\bullet\circ\sigma)(a))}=\overline{f(\env_{\mathcal C}H(\sigma(a)^\bullet))}=
\overline{\env_{\mathcal C}H^\star(f)(\sigma(a)^\bullet)}=\env_{\mathcal C}H^\star(f)^\bullet(a).
\end{multline*}
\epr

\paragraph{Непрерывное тензорное произведение с ${\mathcal C}(M)$.}

Пусть $X$ -- стереотипное пространство, и $M$ -- паракомпактное локально компактное топологическое пространство. Рассмотрим алгебру ${\mathcal C}(M)$ непрерывных функций на $M$ и пространство ${\mathcal C}(M,X)$ непрерывных  функций на $M$ со значениями в $X$. Мы наделяем ${\mathcal C}(M)$ и ${\mathcal C}(M,X)$ стандартной топологией равномерной сходимости на компактах в $M$
$$
u_i\overset{{\mathcal C}(M,X)}{\longrightarrow} 0\quad\Longleftrightarrow\quad \forall \text{компакта}\ K\subseteq M\quad  u_i|_K\overset{{\mathcal C}(K,X)}{\longrightarrow} 0,\qquad u\in {\mathcal C}(M,X),\quad t\in M.
$$
и поточечными алгебраическими операциями:
$$
(\lambda\cdot u)(t)=\lambda\cdot u(t)\qquad (u+v)(t)=u(t)+v(t),\qquad u,v\in {\mathcal C}(M,X),\quad \lambda\in\C,\quad t\in M.
$$

Из \cite[Theorem 8.1]{Ak03} следует

\bprop Пространство ${\mathcal C}(M,X)$ является стереотипным модулем над ${\mathcal C}(M)$, и
\beq\label{C(M,X)-cong-C(M)-odot-X}
{\mathcal C}(M,X)\cong {\mathcal C}(M)\odot X
\eeq
\eprop

В дальнейшем нас будет интересовать случай, когда $X=A$ -- непрерывная (и поэтому стереотипная) алгебра. Пространство ${\mathcal C}(M,A)$ при этом мы также наделяем структурой стереотипной алгебры с поточечным умножением
$$
(u\cdot v)(t)=u(t)\cdot v(t),\qquad u,v\in {\mathcal C}(M,A),\quad t\in M.
$$
Из \eqref{C(M,X)-cong-C(M)-odot-X} следует, что ${\mathcal C}(M,A)$ является стереотипным $A$-модулем.

\btm\label{TH:C(M)-circledast-A->C(M,A)} Для всякой непрерывной алгебры $A$ и любого паракомпактного локально компактного пространства $M$ естественное отображение
\beq\label{C(M)-circledast-A->C(M,A)}
\iota:{\mathcal C}(M)\circledast A\to {\mathcal C}(M,A) \quad\Big|\quad \iota(u\circledast a)(t)=u(t)\cdot a,\quad u\in {\mathcal C}(M),\ a\in A,\ t\in M,
\eeq
является непрерывной оболочкой:
\beq\label{C(M)-circledast-A=C(M,A)}
{\mathcal C}(M)\overset{\mathcal C}{\circledast} A\cong{\mathcal C}(M,A).
\eeq
\etm

\blm\label{LM:iota-in-DEpi-C(M)}
Отображение $\iota:{\mathcal C}(M)\circledast A\to {\mathcal C}(M,A)$ имеет образ, плотный в ${\mathcal C}(M,A)$.
\elm
\bpr
Пусть $f\in {\mathcal C}(M,A)$. Зафиксируем компакт $K\subseteq M$ и выпуклую окрестность нуля $U\subseteq A$. Поскольку $f$ равномерно непрерывно на $K$, найдется окружение диагонали $V$ в $K$ такое, что
$$
|s-t|<V\quad\Longrightarrow\quad f(s)-f(t)\in U.
$$
Выберем конечный набор точек $t_1,...,t_n\in K$ такой что окрестности $t_k+V$ образуют покрытие $K$. Подберем для них разбиение единицы на компакте $K$, то есть функции $\eta_1,...,\eta_n\in{\mathcal C}(M)$ такие, что
$$
\eta_k\big|_{K\setminus(t_k+V)}=0,\qquad 0\le \eta_k\le 1,\qquad\sum_{k=1}^n\eta_k\big|_K=1.
$$
Положим
$$
x=\sum_{k=1}^n\eta_k\circledast f(t_k)\in{\mathcal C}(M)\circledast A.
$$
Тогда для $s\in K$ мы получим
$$
f(s)-\iota(x)(s)=\sum_{k=1}^n\eta_k(s)\cdot f(s)-\sum_{k=1}^n\eta_k(s)\cdot f(t_k)=\sum_{k=1}^n\eta_k(s)\cdot \big(f(s)-f(t_k)\big)\in U.
$$
\epr

\blm\label{LM:J^0_C(M)C(M)-circledast-A-cong-J^0_C(M)C(M,A)}
Модули ${\mathcal C}(M)\circledast A$ и ${\mathcal C}(M,A)$ над алгеброй ${\mathcal C}(M)$ имеют изоморфные расслоения значений:
\beq\label{J^0_C(M)C(M)-circledast-A-cong-J^0_C(M)C(M,A)}
\Jet^0_{{\mathcal C}(M)}{\mathcal C}(M)\circledast A\cong \Jet^0_{{\mathcal C}(M)}{\mathcal C}(M,A),\qquad n\in\N
\eeq
\elm
\bpr
Для каждой точки $t\in M$ идеал $I_t$ имеет коразмерность 1 в ${\mathcal C}(M)$, поэтому можно воспользоваться леммой \ref{PROP:[(X-circledast-Z)/(Y-circledast-Z)]^triangledown-cong-[(X-odot-Z)/(Y-odot-Z)]^triangledown}:
\begin{multline*}
\Jet^0_{{\mathcal C}(M)}{\mathcal C}(M)\circledast A=[({\mathcal C}(M)\circledast A)/ (I_t\circledast A)]^\vartriangle=\eqref{[(X-circledast-Z)/(Y-circledast-Z)]^triangledown-cong-[(X-odot-Z)/(Y-odot-Z)]^triangledown}=\\=
[({\mathcal C}(M)\odot A)/ (I_t\odot A)]^\vartriangle=\Jet^0_{{\mathcal C}(M)}{\mathcal C}(M)\odot A=\eqref{C(M,X)-cong-C(M)-odot-X}=
\Jet^0_{{\mathcal C}(M)}{\mathcal C}(M,A)
\end{multline*}
\epr

\blm\label{TH:diff-oper<-morfizm-rassl-znachenij} Пусть $M$ -- паракомпактное локально компактное пространство, $F$ -- $C^*$-алгебра и $\ph:{\mathcal C}(M)\to F$ -- гомоморфизм инволютивных стереотипных алгебр, причем $\ph({\mathcal C}(M))$ лежит в центре $F$:
$$
\ph({\mathcal C}(M))\subseteq Z(F).
$$
Тогда для любого стереотипного пространства $X$ всякий морфизм расслоений значений
$$
\mu:\Jet_{{\mathcal C}(M)}^0({\mathcal C}(M,X))\to \Jet_{{\mathcal C}(M)}^0(F)
$$
определяет единственный морфизм стереотипных ${\mathcal C}(M)$-модулей $D:{\mathcal C}(M,X)\to F$, удовлетворяющий тождеству
\beq\label{j^0(Dx)=mu-circ-j^0(x)} \jet^0(\ph(x))=\mu\circ \jet^0(x),\qquad
x\in {{\mathcal C}(M,X)}.
 \eeq \elm
\bpr По теореме \ref{TH:B-cong-Sec(val_AB)}, отображение
$v:F\to\Sec(\pi^0_{{{\mathcal C}(M)},F})$, переводящее $F$ в алгебру непрерывных сечений расслоения
значений $\pi^0_{{{\mathcal C}(M)},F}: \Jet^0_{{\mathcal C}(M)}F\to\Spec({{\mathcal C}(M)})$ над алгеброй ${{\mathcal C}(M)}$, является
изоморфизмом $C^*$-алгебр:
$$
F\cong\Sec(\pi^0_{{{\mathcal C}(M)},F}).
$$
Рассмотрим обратный изоморфизм $v^{-1}:\Sec(\pi^0_{{{\mathcal C}(M)},F})\to F$:
\beq\label{lambda(j^0(b))=b} v^{-1}(\jet^0(b))=b,\qquad b\in F.
\eeq
Тогда всякому морфизму расслоений значений $\mu:\Jet_{{\mathcal C}(M)}^0({{\mathcal C}(M,X)})\to \Jet_{{\mathcal C}(M)}^0(F)$ можно поставить в
соответствие оператор $\ph:{{\mathcal C}(M,X)}\to F$ по формуле
\beq\label{Da=lambda(mu-circ-j^0(a))} \ph(x)=v^{-1}\Big(\mu\circ \jet^0(x)\Big),\qquad
x\in {{\mathcal C}(M,X)}.
\eeq
Он, очевидно, будет удовлетворять тождеству
\eqref{j^0(Dx)=mu-circ-j^0(x)}.
\epr

\blm Отображение $\iota:{\mathcal C}(M)\circledast A\to {\mathcal C}(M,A)$ является непрерывным расширением.
\elm
\bpr
Пусть $\ph:{\mathcal C}(M)\circledast A\to B$ -- морфизм в какую-нибудь $C^*$-алгебру $B$. По лемме \ref{LM:ph:A-circledast-B->C}, он представим в виде
\beq\label{ph(u-circledast-a)=eta(u)-cdot-alpha(a)}
\ph(u\circledast a)=\eta(u)\cdot\alpha(a)=\alpha(a)\cdot\eta(u),\qquad u\in {\mathcal C}(M), \ a\in A,
\eeq
где $\eta:{\mathcal C}(M)\to B$, $\alpha:A\to B$ -- некоторые морфизмы инволютивных стереотипных алгебр. Рассмотрим оператор $\eta$ и обозначим буквой $C$ его образ в $B$:
$$
C=\overline{\eta({\mathcal C}(M))}.
$$
Пусть $F$ -- коммутант алгебры $C$ в $B$:
$$
F=C^!=\{x\in B:\quad \forall c\in C\quad x\cdot c=c\cdot x\}.
$$
Поскольку алгебра $C$ коммутативна, она тоже лежит в $F$, и более того, в центре $F$:
$$
C\subseteq Z(F).
$$
Заметим, что образ оператора $\ph$ лежит в $F$:
$$
\ph({\mathcal C}(M)\circledast A)\subseteq F,
$$
потому что
\begin{multline*}
\ph(v\circledast a)\cdot \eta(u)=\ph(v\circledast a)\cdot \ph(u\circledast 1)=\ph\big((v\circledast a)\cdot (u\circledast 1)\big)=\ph\big((v\cdot u)\circledast 1\big)=\ph\big((u\cdot v)\circledast 1\big)=\\=\ph\big((u\circledast 1)\cdot (v\circledast a)\big)=\ph(u\circledast 1)\cdot \ph(v\circledast a)=\eta(u)\cdot \ph(v\circledast a)
\end{multline*}

Чтобы убедиться, что $\iota$ -- непрерывное расширение, нам надо показать, что существует (единственный) гомоморфизм $\ph':{\mathcal C}(M,A)\to F$, продолжающий $\ph$:
 \beq\label{prodolzhenie-D-na-C(M,A)}
 \xymatrix @R=2pc @C=1.2pc
 {
 {\mathcal C}(M)\circledast A\ar[rr]^{\iota}\ar[dr]_{\ph} & & {\mathcal C}(M,A)\ar@{-->}[dl]^{\ph'}\\
  & F &
 }
 \eeq
Гомоморфизм $\ph:{\mathcal C}(M)\circledast A\to F$ является ${\mathcal C}(M)$-морфизмом, и значит, по теореме
\ref{TH:morf-modul->morf-rassl-znach} ему соответствует некий морфизм расслоений значений
$\jet_0[\ph]:\Jet_{{\mathcal C}(M)}^0[{\mathcal C}(M)\circledast A]\to \Jet_{{\mathcal C}(M)}^0(F)=\pi^0_A F$, удовлетворяющий тождеству
$$
\jet^0(\ph(x))=\jet_n[\ph]\circ \jet^n(x),\qquad x\in {\mathcal C}(M)\circledast A.
$$
По лемме \ref{LM:J^0_C(M)C(M)-circledast-A-cong-J^0_C(M)C(M,A)}, расслоения значений алгебр ${\mathcal C}(M)\circledast A$ и ${\mathcal C}(M,A)$ изоморфны. Обозначим этот изоморфизм $\mu:\Jet_{{\mathcal C}(M)}^0[{\mathcal C}(M)\circledast A]\gets \Jet^0_{{\mathcal C}(M)}[{\mathcal C}(M,A)]$. Рассмотрим композицию $\nu=\jet_0[\ph]\circ\mu:\Jet_{{\mathcal C}(M)}[{\mathcal C}(M,A)]\to \Jet_{{\mathcal C}(M)}^0(F)=\pi^0_A F$:
$$
 \xymatrix @R=2pc @C=1.2pc
 {
 \Jet_{{\mathcal C}(M)}^0[{\mathcal C}(M)\circledast A]\ar[dr]_{\jet_0[\ph]} & & \Jet^0_{{\mathcal C}(M)}[{\mathcal C}(M,A)]\ar@{-->}[dl]^{\quad\nu=\jet_0[\ph]\circ\mu}\ar[ll]_{\mu}\\
  & \Jet^0_{{\mathcal C}(M)}[F] &
 }
$$
По лемме \ref{TH:diff-oper<-morfizm-rassl-znachenij} этой пунктирной стрелке $\nu$ соответствует морфизм  $\ph':{\mathcal C}(M,A)\to F$ (над алгеброй ${\mathcal C}(M)$), удовлетворяющий тождеству
$$
\jet^0(\ph'f)=\jet_0[\ph]\circ \jet^0(f),\qquad  f\in {\mathcal C}(M,A).
$$
Для всякого $x\in {\mathcal C}(M)\circledast A$ мы получим
\beq\label{PROOF:TH:C(M,A)-1}
\jet^0\big(\ph'(\iota(x))\big)=\jet_0[\ph]\circ \jet^0(\iota(x))=\jet_0[\ph]\circ
\jet^0(x)=\jet^0(\ph(x)).
\eeq
Заметим далее, что по теореме \ref{TH:B-cong-Sec(val_AB)}, отображение
$\jet^0=v:F\to\Sec(\pi^0_{{\mathcal E}(M)}F)=\Sec(\Jet^0_{{\mathcal E}(M)}F)$, переводящее $F$ в алгебру непрерывных сечений
расслоения значений $\pi^0_{{\mathcal E}(M)}F: \Jet^0_{{\mathcal E}(M)}F\to\Spec({{\mathcal E}(M)})$ над алгеброй ${{\mathcal E}(M)}$, является
изоморфизмом $C^*$-алгебр:
$$
F\cong\Sec(\pi^0_{{\mathcal E}(M)}F).
$$
Поэтому к \eqref{PROOF:TH:C(M,A)-1} можно применить оператор, обратный $\jet^0$, и мы получим равенство
$$
\ph'(\iota(x))=\ph(x).
$$
То есть $\ph'$ продолжает $\ph$ в диаграмме \eqref{prodolzhenie-D-na-C(M,A)}. По лемме \ref{LM:iota-in-DEpi-C(M)} элементы вида $\iota(u\circledast a)$ полны в ${\mathcal C}(M,A)$, по этой причине такое продолжение $\ph'$ будет единственно.

По построению, оператор $\ph'$ будет морфизмом относительно алгебры ${\mathcal C}(M)$, однако нам этого недостаточно: нужно чтобы $\ph'$ был гомоморфизмом алгебр. Чтобы это проверить, выберем $u,v\in {\mathcal C}(M)$ и $a,b\in A$. Тогда
$$
\ph'(\iota(u\circledast a)\cdot\iota(v\circledast b))=\ph'(\iota(u\circledast a\cdot v\circledast b))=\ph(u\circledast a\cdot v\circledast b)=\ph(u\circledast a)\cdot\ph(v\circledast b)=\ph'(\iota(u\circledast a))\cdot\ph'(\iota(v\circledast b))
$$
При этом по лемме \ref{LM:iota-in-DEpi-C(M)} элементы вида $\iota(u\circledast a)$ полны в ${\mathcal C}(M,A)$. Отсюда следует, что их можно заменить произвольными векторами из ${\mathcal C}(M,A)$, и значит $\ph'$ должен быть  гомоморфизмом.
\epr

\blm Отображение $\iota:{\mathcal C}(M)\circledast A\to {\mathcal C}(M,A)$ является непрерывной оболочкой.
\elm
\bpr

Пусть $\sigma:{\mathcal C}(M)\circledast A\to C$ -- какое-то другое непрерывное расширение. Нам нужно убедиться, что существует морфизм $\upsilon$, замыкающий диаграмму
$$
 \xymatrix @R=2pc @C=1.2pc
 {
{\mathcal C}(M)\circledast A\ar[rr]^{\sigma}\ar[dr]_{\iota} & & C \ar@{-->}[dl]^{\upsilon}\\
  & {\mathcal C}(M,A) &
 }
$$

Зафиксируем компакт $K\subseteq M$ и гомоморфизм $\eta:A\to B$ в какую-нибудь $C^*$-алгебру, и положим
$$
D(u\circledast a)(t)=u(t)\circledast \eta(a),\qquad u\in {\mathcal C}(M),\quad a\in A,\quad t\in K.
$$
Отображение $D$ будет гомоморфизмом из ${\mathcal C}(M)\circledast A$ в алгебру ${\mathcal C}(K)\circledast B$. Та, в свою очередь естественно отображается в максимальное тензорное произведение $C^*$-алгебр ${\mathcal C}(K)\underset{\max}{\otimes}B$, а оно изоморфно ${\mathcal C}(K)\odot B$ и ${\mathcal C}(K,B)$:
$$
{\mathcal C}(K)\circledast B\to {\mathcal C}(K)\underset{\max}{\otimes}B\cong \cite[(5.25)]{Akbarov-De-Gruyter-I} \cong
{\mathcal C}(K)\odot B\cong \eqref{C(M,X)-cong-C(M)-odot-X} \cong {\mathcal C}(K,B).
$$

Поэтому мы можем считать $D$ морфизмом в $C^*$-алгебру ${\mathcal C}(K,B)$:
$$
D:{\mathcal C}(M)\circledast A\to {\mathcal C}(K)\circledast B\to {\mathcal C}(K)\underset{\max}{\otimes}B \cong\cite[(5.25)]{Akbarov-De-Gruyter-I}\cong{\mathcal C}(K)\odot B\cong\eqref{C(M,X)-cong-C(M)-odot-X}\cong {\mathcal C}(K,B)
$$
Поскольку $\sigma:{\mathcal C}(M)\circledast A\to C$ есть непрерывное расширение, гомоморфизм $D:{\mathcal C}(M)\circledast A\to {\mathcal C}(K,B)$ должен однозначно продолжаться до некоторого
гомоморфизма $D':C\to {\mathcal C}(K,B)$:
 \beq\label{DIAGR:env-C(M)-circledast-A-1}
 \xymatrix @R=2pc @C=1.2pc
 {
{\mathcal C}(M)\circledast A\ar[rr]^{\sigma}\ar[dr]_{D} & & C \ar@{-->}[dl]^{D'}\\
  & {\mathcal C}(K,B) &
 }
 \eeq
Если теперь зафиксировать $c\in C$ и менять компакт $K\subset M$, то
возникающие при этом непрерывные функции $D'(c)$ на $K$ будут согласованы между
собой тем, что на пересечении своих областей определения они совпадают. Поэтому
определена некая общая непрерывная функция $\iota_B'(c):M\to B$, обладающая тем
свойством, что ее ограничение на каждый компакт $K$ будет совпадать с
соответствующей функцией $D'(c)$:
$$
\iota'(c)\big|_K=D'(c),\qquad K\subset M.
$$
Иными словами, определено некое отображение
$\iota'_B:C\to{\mathcal C}(M,B)$ (по построению это будет гомоморфизм алгебр),
для которого будет коммутативна следующая диаграмма, уточняющая
\eqref{DIAGR:env-C(M)-circledast-A-1}:
 \beq\label{DIAGR:env-C(M)-circledast-A-*-M}
 \xymatrix @R=2.5pc @C=2pc
 {
 {\mathcal C}(M)\circledast A\ar[rr]^{\sigma}\ar@/_5ex/[ddr]_{D}\ar[dr]_{\iota_B} & & C \ar@/^5ex/[ddl]^{D'}\ar@{-->}[dl]^{\iota'_B}\\
  & {\mathcal C}(M,B)\ar[d]_{\rho_K} & \\
  & {\mathcal C}(K,B) &
 }
 \eeq
(здесь $\rho_K$ -- отображение ограничения на компакт $K$).

Пусть теперь $U$ -- $C^*$-окрестность нуля\footnote{$C^*$-окрестности нуля были определены на с.\pageref{DEF:C*-neighbourhood-of-zero}.} в $A$, соответствующая гомоморфизму $\eta:A\to B$. Из того, что $\sigma$ -- плотный эпиморфизм, следует, что верхний внутренний треугольник в \eqref{DIAGR:env-C(M)-circledast-A-*-M} можно достроить до диаграммы
 \beq\label{DIAGR:env-E(M)-circledast-A-*-M-1}
 \xymatrix @R=2.5pc @C=2pc
 {
 {\mathcal C}(M)\circledast A\ar[rr]^{\sigma}\ar@/_5ex/[ddr]_{\iota_B}\ar[dr]_{\vartheta_U} & & C \ar@/^5ex/[ddl]^{\iota'_B}\ar@{-->}[dl]^{\vartheta'_U}\\
  & {\mathcal C}(M,A/U)\ar[d]_{\eta_U\oslash 1_M} & \\
  & {\mathcal C}(M,B) &
 }
 \eeq
где $\eta_U:A/U\to B$ -- морфизм из \eqref{ph=ph_U-circ-pi_U-0}, и
$$
\vartheta_U(u\circledast a)(t)=u(t)\cdot\pi_U(a),\qquad u\in{\mathcal E}(M),\quad a\in A,\quad t\in M,
$$
$$
(\eta_U\oslash 1_M)(h)(t)=\eta_U(h(t)),\qquad h\in {\mathcal E}(M,A/U),\quad t\in M.
$$
Из определения $\vartheta_U$ сразу следует, что если $U'\subseteq U$ -- какая-то другая $C^*$-окрестность нуля, то
\beq\label{vartheta_U=(varkappa^U'_U-oslash-1_M)-cdot-vartheta_U'-0}
\vartheta_U=(\varkappa^{U'}_U\oslash 1_M)\cdot\vartheta_{U'},\qquad U\supseteq U',
\eeq
где $\varkappa^{U'}_U$ -- морфизм из \eqref{varkappa^U'_U-0}, и
$$
(\varkappa^{U'}_U\oslash 1_M)(h)(t)=\varkappa^{U'}_U(h(t)),\qquad h\in {\mathcal E}(M,A/U'),\quad t\in M.
$$
Равенство \eqref{vartheta_U=(varkappa^U'_U-oslash-1_M)-cdot-vartheta_U'-0} будет левым нижним внутренним треугольником в диаграмме $$
 \xymatrix @R=2.5pc @C=3pc
 {
{\mathcal C}(M)\circledast A\ar[rr]^{\sigma}\ar[dr]_(.6){\vartheta_{U'}}\ar@/_5ex/[ddr]_{\vartheta_U} & & C \ar[dl]^(.6){\vartheta'_{U'}} \ar@/^5ex/[ddl]^{\vartheta'_U}\\
 & {\mathcal C}(M,A/U')\ar@{-->}[d]^{\varkappa^{U'}_U\oslash 1_M} & \\
  & {\mathcal C}(M,A/U) &
 }
$$
При этом периметр и верхний внутренний треугольник здесь будут вариантами верхнего внутреннего треугольника в \eqref{DIAGR:env-C(M)-circledast-A-*-M}, и вдобавок $\sigma$ -- эпиморфизм. Как следствие, оставшийся правый нижний внутренний треугольник тоже должен быть коммутативен.

Это означает, что морфизмы $\vartheta'_U:C\to {\mathcal C}(M,A/U)$ образуют проективный конус системы $\varkappa^{U'}_U\oslash 1_M$, и поэтому существует морфизм $\vartheta'$ в проективный предел:
$$
 \xymatrix @R=2.5pc @C=3pc
 {
{\mathcal C}(M)\circledast A\ar[rr]^{\sigma}\ar[dr]_(.6){\vartheta}\ar@/_5ex/[ddr]_{\vartheta_U} & & C \ar@{-->}[dl]^(.6){\vartheta'} \ar@/^5ex/[ddl]^{\vartheta'_U}\\
 & \projlim\limits_{0\gets U'}{\mathcal C}(M,A/U')\ar[d]^{\varkappa_U\oslash 1_M} & \\
  & {\mathcal C}(M,A/U) &
 }
$$
Теперь заметим цепочку
$$
\projlim\limits_{0\gets U'}{\mathcal C}(M,A/U')=\projlim\limits_{0\gets U'}({\mathcal C}(M)\odot A/U')=\cite[(4.152)]{Akbarov-De-Gruyter-I}=
{\mathcal C}(M)\odot \projlim\limits_{0\gets U'}A/U'={\mathcal C}(M,\projlim\limits_{0\gets U'}A/U')
$$
и подставим последнее пространство в нашу диаграмму:
$$
 \xymatrix @R=2.5pc @C=3pc
 {
{\mathcal C}(M)\circledast A\ar[rr]^{\sigma}\ar[dr]_(.6){\vartheta}\ar@/_5ex/[ddr]_{\vartheta_U} & & C \ar@{-->}[dl]^(.6){\vartheta'} \ar@/^5ex/[ddl]^{\vartheta'_U}\\
 & {\mathcal C}(M,\projlim\limits_{0\gets U'}A/U')\ar[d]^{\varkappa_U\oslash 1_M} & \\
  & {\mathcal C}(M,A/U) &
 }
$$

Еще раз вспомним, что $\sigma$ -- плотный эпиморфизм. Из этого следует, что стрелка $\vartheta'$ поднимается до некоторой стрелки $\upsilon$ со значениями в пространстве ${\mathcal C}(M,\Im\pi)$ функций, принимающих значения в образе отображения $\pi:A\to\projlim\limits_{0\gets U'}A/U'$, или, что то же самое, в непосредственном подпространстве, порожденном множеством значений отображения $\pi$, а это пространство как раз совпадает с $A$, поскольку $A$ -- непрерывная алгебра:
$$
\Im\pi\cong\Env_{\mathcal C}A\cong A
$$
Мы получаем диаграмму
$$
 \xymatrix @R=2.5pc @C=3pc
 {
{\mathcal C}(M)\circledast A\ar[rr]^{\sigma}\ar[dr]_(.6){\iota}\ar@/_5ex/[ddr]_{\vartheta} & & C \ar@{-->}[dl]^(.6){\upsilon} \ar@/^5ex/[ddl]^{\vartheta'}\\
 & {\mathcal C}(M,A)\ar[d]^{\im\pi\oslash 1_M} & \\
  & {\mathcal C}(M,\projlim\limits_{0\gets U'}A/U') &
 }
$$
где $\pi$ -- морфизм из \eqref{A-to-leftlim-A/U-0}.
\epr

\subsection{Непрерывные оболочки коммутативных алгебр}

Условимся непрерывное отображение топологических пространств $\e:X\to Y$ называть {\it наложением}, если всякий компакт $T\subseteq Y$ содержится в образе некоторого компакта $S\subseteq X$. Если пространство $Y$ хаусдорфово, то это автоматически означает, что $\e$ должно быть сюръективно. Если вдобавок отображение $\e$ инъективно, то мы называем его {\it точным наложением}.\label{DEF:nalozhenie} В точном наложении $\e:X\to Y$ пространство $Y$ можно представлять себе как новую, более слабую топологизацию пространства $X$, которая не меняет систему компактов и топологию на каждом компакте.

\paragraph{${\mathcal C}(M)$, как непрерывная оболочка своих подалгебр.}

\bit{
 \item[$\bullet$]
{\it Инволютивным спектром} $\Spec(A)$ инволютивной стереотипной алгебры $A$ над $\C$ мы называем множество ее инволютивных характеров, то есть гомоморфизмов $\chi:A\to \C$ (также непрерывных, инволютивных и сохраняющих единицу). Это множество мы наделяем топологией равномерной сходимости на компактах в $A$.

 \item[$\bullet$] Если $\ph:A\to B$ --- морфизм инволютивных стереотипных алгебр, то определено непрерывное отображение 
 \beq
 \ph^{\Spec}:\Spec(B)\to \Spec(A) \quad \Big| \quad \ph^{\Spec}(\chi)=\chi\circ\ph,\quad \chi\in \Spec(B)
 \eeq
называемое {\it отображением спектров, сопряженным морфизму} $\ph$.

}\eit

\btm\label{C-obolochka-podalgebry-v-C(M)} Пусть $A$ --- коммутативная инволютивная стереотипная алгебра, и $M$ --- паракомпактное локально компактное пространство. Морфизм инволютивных стереотипных алгебр
$$
\iota:A\to \mathcal{C}(M).
$$
тогда и только тогда будет непрерывной оболочкой алгебры $A$\beq\label{Env_C_A=C(M)}
\Env_{\mathcal C} A=\mathcal{C}(M)
\eeq
когда сопряженное отображение спектров $\iota^{\Spec}:\Spec(A)\gets M$ является точным наложением. 
\etm

\bpr
1. Докажем сначала необходимость. Пусть выполняется \eqref{Env_C_A=C(M)}. Зафиксируем компакт $T\subseteq\Spec(A)$ и рассмотрим отображение $\iota_T:A\to{\mathcal C}(T)$, определенное правилом
$$
\iota_T(a)(t)=t(a),\qquad t\in T,\ a\in A.
$$
или, что то же самое, правилом
\beq\label{delta^t-circ-ph_T=t}
\delta^t\circ \iota_T=t,\qquad t\in T
\eeq
Это будет инволютивный гомоморфизм в $C^*$-алгебру, поэтому он продолжается на оболочку $\Env_{\mathcal C} A=\mathcal{C}(M)$:
$$
 \xymatrix @R=2pc @C=1.2pc
 {
 A\ar[rr]^{\iota}\ar[dr]_{\iota_T} & & {\mathcal C}(M)\ar@{-->}[dl]^{\pi}\\
  & {\mathcal C}(T) &
 }
$$
где $\pi$ -- некий инволютивный гомоморфизм. Сопряженное отображение спектров $\pi^{\Spec}:\Spec({\mathcal C}(T))\to\Spec({\mathcal C}(M))$ переводит $T=\Spec({\mathcal C}(T))$ в некий компакт $S=\pi^{\Spec}(T)\subseteq \Spec({\mathcal C}(M))$. При отображении спектров $\iota^{\Spec}:\Spec(A)\gets\Spec({\mathcal C}(M))$ компакт $S$ сюръективно отобразится на компакт $T$, потому что любой
элемент $t\in T$ встраивается в диаграмму 
$$
 \xymatrix @R=2pc @C=1.2pc
 {
 A\ar[rr]^{\iota}\ar[dr]_{\iota_T}\ar@/_6ex/[ddr]_{t=\delta^t\circ \iota_T} & & {\mathcal C}(M)\ar@{-->}[dl]^{\pi}\ar@/^6ex/[ddl]^{\delta^t\circ\pi}\\
  & {\mathcal C}(T)\ar[d]_{\delta^t} & \\
  & \C & 
 }
$$
означающую, что для $t\in T$ существует элемент $s=\delta^t\circ\pi\in S$ такой что
$$
\iota^{\Spec}(s)=s\circ\iota=\delta^t\circ\pi\circ\iota=\delta^t\circ\iota_T=\eqref{delta^t-circ-ph_T=t}=t.
$$
Это доказывает, что $\iota^{\Spec}:\Spec(A)\gets\Spec({\mathcal C}(M))$ является наложением. Покажем далее, что оно инъективно. Если бы это было не так, мы получили бы, что какие-то две точки $s\ne s'\in M$ при отображении $\iota^{\Spec}$ слипаются:
$$
s\circ\sigma=s'\circ\sigma=t\in\Spec(A)
$$
Это можно понимать так, что характер $t:A\to\C$ имеет два разных продолжения на ${\mathcal C}(M)$:
$$
 \xymatrix @R=2pc @C=1.2pc
 {
 A\ar[rr]^{\iota}\ar[dr]_{t} & & {\mathcal C}(M)\ar@/_1ex/[dl]_{s}\ar@/^1ex/[dl]^{s'}\\
  & \C &
 }
$$
Но $\iota$ -- непрерывная оболочка, и поэтому характер $t:A\to\C$, будучи инволютивным гомоморфизмом в $C^*$-алгебру $\C$, должен продолжаться единственным образом.

2. Теперь докажем достаточность. Пусть $\iota^{\Spec}$ -- точное наложение. Тогда, в частности, $\iota^{\Spec}$ ---  сюръективное отображение, и поэтому алгебра $\iota(A)$ разделяет точки $M$. Вдобавок $\iota(A)$ содержит единицу, а значит и все константы. Поэтому, по теореме Стоуна-Вейерштрасса алгебра $\iota(A)$ должна быть плотна в $\mathcal{C}(M)$.

Покажем, что инъекция $\iota:A\to \mathcal{C}(M)$ является непрерывным расширением. Пусть $\ph:A\to B$ -- морфизм $A$
в какую-нибудь $C^*$-алгебру $B$. Чтобы построить пунктирную стрелку $\ph'$,
замыкающую диаграмму \eqref{DEF:diagr-nepr-rasshirenie},
$$
 \xymatrix @R=2pc @C=1.2pc
 {
 A\ar[rr]^{\iota}\ar[dr]_{\ph} & & {\mathcal C}(M)\ar@{-->}[dl]^{\ph'}\\
  & B &
 }
$$
достаточно считать, что $B$ коммутативна и что $\ph(A)$ плотно в $B$ (потому что иначе можно будет заменить $B$ на замыкание $\overline{\ph(A)}$ в $B$, которое будет коммутативной подалгеброй в $B$). Тогда из коммутативности $B$
будет следовать, что $B$ имеет вид ${\mathcal C}(T)$, а из плотности $\ph(A)$  в $B$ -- что компакт $T$ инъективно вкладывается в $\Spec(A)$. Поскольку по условию теоремы, $\iota^{\Spec}:\Spec(A)\gets M$ является наложением, $T\subseteq\Spec(A)$ является образом некоторого компакта $S\subseteq M$, $S\cong T$. Тогда отображение $\ph$ представляется как композиция инъекции $\iota:A\to{\mathcal C}(M)$ с отображением $\pi_T:{\mathcal C}(M)\to {\mathcal C}(T)$ ограничения на компакт $T$:
$$
 \xymatrix @R=2pc @C=1.2pc
 {
 A\ar[rr]^{\iota}\ar[dr]_{\ph} & & {\mathcal C}(M)\ar@{-->}[dl]^{\pi_T}\\
  & {\mathcal C}(T) &
 }
$$
Мы получаем, что $\ph'=\pi_T$. При этом пунктирная стрелка будет единственной из-за того, что $A$ плотно в ${\mathcal C}(M)$.

Проверим после этого, что $\iota:A\to {\mathcal C}(M)$ является максимальным
расширением, то есть если взять какое-то другое расширение $\sigma:A\to C$, то
найдется единственный морфизм $\upsilon:C\to {\mathcal C}(M)$, замыкающий диаграмму
 \beq\label{TOP:C->C(M)}
 \xymatrix @R=2pc @C=1.2pc
 {
 & A\ar[rd]^{\iota}\ar[ld]_{\sigma} & \\
 C\ar@{-->}[rr]_{\upsilon} & & {\mathcal C}(M)
 }
 \eeq
Для всякого компакта $T\subseteq M$ гомоморфизм
$$
\iota_T:A\to{\mathcal C}(T)\quad\big|\quad \iota_T(a)(t)=t(a),\qquad t\in T\subseteq\Spec(A)
$$
однозначно продолжается до некоторого гомоморфизма $\iota'_T:A'\to{\mathcal C}(T)$
$$
 \xymatrix @R=2pc @C=1.2pc
 {
 A\ar[rr]^{\sigma}\ar[dr]_{\iota_T} & & A'\ar@{-->}[dl]^{\iota'_T}\\
  & {\mathcal C}(T) &
 }
$$
Если $T\subseteq S$ -- два компакта в $M$, то эта диаграмма дополняется до диаграммы
$$
 \xymatrix @R=2.5pc @C=3pc
 {
 A\ar[rr]^{\sigma}\ar@/_3ex/[ddr]_{\iota_T}\ar[dr]_{\iota_S} & & A'\ar@{-->}[dl]^{\iota'_S}\ar@/^3ex/@{-->}[ddl]^{\iota'_T}\\
  & {\mathcal C}(S)\ar[d]^{\pi_T^S} & \\
  & {\mathcal C}(T) &
 }
$$
где $\pi_T^S$ -- отображение ограничения на компакт $T$. Правый нижний треугольник в этой диаграмме означает, что морфизмы $\iota'_T:A'\to {\mathcal C}(T)$ образуют проективный конус в контравариантной системе $\pi_T^S$. Поэтому существует стрелка $\iota'$, замыкающая правый нижний треугольник во всех диаграммах
$$
 \xymatrix @R=2.5pc @C=3pc
 {
 A\ar[rr]^{\sigma}\ar@/_3ex/[ddr]_{\iota_T}\ar[dr]_{\iota} & & A'\ar@{-->}[dl]^{\iota'}\ar@/^3ex/[ddl]^{\iota'_T}\\
  & {\mathcal C}(M)\ar[d]^{\pi_T} & \\
  & {\mathcal C}(T) &
 }
$$
Поскольку периметр и левый нижний треугольник в каждой такой диаграмме также коммутативны, мы получаем
$$
\pi_T\circ\iota'\circ\sigma=\iota'_T\circ\sigma=\iota_T=\pi_T\circ\iota.
$$
Поскольку это верно для всех $T$, мы получаем, равенство
$$
\iota'\circ\sigma=\iota.
$$
\epr

\noindent\rule{160mm}{0.1pt}\begin{multicols}{2}

\paragraph{Контрпример.}

\bex\label{EX:C(8)} {\it
Существует плотная инволютивная стереотипная подалгебра $A$ в ${\mathcal C}(\R)$, у которой сопряженное отображение спектров $\iota^{\Spec}:\Spec(A)\gets M$ является биекцией, но не наложением, и, как следствие, непрерывная оболочка $A$ не совпадает с ${\mathcal C}(\R)$:}
$$
\Spec(A)=M,\qquad \Env_{\mathcal C}A\ne {\mathcal C}(\R)
$$
\eex
\bpr
Такой алгеброй будет алгебра $A$ непрерывных функций на $\R$, у которых предел на бесконечности равен значению в какой-нибудь фиксированной точке, например, в нуле:
$$
u\in A\quad\Longleftrightarrow\quad u\in {\mathcal C}(\R)\quad\&\quad \lim_{t\to\infty}u(t)=u(0).
$$
Алгебра $A$ инволютивна, содержит константы и разделяет точки $\R$, поэтому она плотна в ${\mathcal C}(\R)$. Мы наделяем $A$ топологией равномерной сходимости на всей прямой $\R$. Спектр $A$, как легко понять, представляет собой прямую $\R$, только с новой топологией, в которой базисными окрестностями точки 0 становятся множества вида
$$
(-\infty, A)\cup(a,b)\cup(B,+\infty)
$$
где $A<a<0<b<B$ -- числа на $\R$ (а базисные окрестности остальных точек не меняются). Эту топологию удобно представлять как индуцированную на $\R$ при вложении $\R$ в фигуру, внешне похожую на цифру 8, и которую поэтому удобно обозначать символом 8 (при этом вложении концы прямой $\R$ загибаются и подходят к точке $0\in\R$). Спектр $A$, понятное дело, будет гомеоморфен этой восьмерке 8,
$$
\Spec A\cong 8,
$$
а сама алгебра $A$ изоморфна (как стереотипная алгебра) алгебре ${\mathcal C}(8)$ (алгебре непрерывных функций на компакте 8). Поэтому непрерывная оболочка $A$ совпадает с $A$, и не изоморфна ${\mathcal C}(\R)$:
$$
\Env_{\mathcal C}A\cong A\cong{\mathcal C}(8)\not\cong {\mathcal C}(\R).
$$
\epr

\end{multicols}\noindent\rule[10pt]{160mm}{0.1pt}

\paragraph{Преобразование Гельфанда как непрерывная оболочка.}

 \bit{

 \item[$\bullet$]
{\it Преобразованием Гельфанда} инволютивной стереотипной алгебры $A$ мы называем естественное
отображение ${\mathcal G}_A:A\to C(M)$ алгебры $A$ в алгебру $C(M)$ функций на инволютивном спектре $M=\Spec(A)$, непрерывных на каждом компакте $K\subseteq M$:
 \beq\label{TOP:vlozh-A-v-C(Spec(A))}
{\mathcal G}_A(x)(t)=t(x),\qquad t\in M=\Spec(A),\ x\in A.
 \eeq

 }\eit

\bprop\label{LM:rho-in-DEpi} Если спектр $M=\Spec(A)$ алгебры $A$ является паракомпактным локально компактным пространством,
морфизм ${\mathcal G}_A:A\to C(M)$ является плотным эпиморфизмом.\footnote{Предложение \ref{LM:rho-in-DEpi} останется частично верным для случая, когда спектр $M=\Spec(A)$ --- произвольное $k$-пространство. В этом случае нужно определить топологию на алгебре $C(M)$ как псевдонасыщение топологии равномерной сходимости на компактах в $M$. Это превращает $C(M)$ в стереотипную алгебру, и тогда преобразование Гельфанда ${\mathcal G}_A:A\to C(M)$ будет морфизмом стереотипных алгебр. Остается, однако, сохранится ли при таком ослаблении посылки свойство эпиморфности у ${\mathcal G}_A:A\to C(M)$.}
\eprop
\bpr  В первой части этого утверждения неочевидной является только непрерывность отображения ${\mathcal G}_A$. Пусть $U$ -- базисная окрестность нуля в $C(M)$, то есть $U=\{f\in C(M):\ \sup_{t\in T}|f(t)|\le \e\}$ для некоторого компакта $T\subseteq M$ и числа $\e>0$. Ее прообразом при отображении ${\mathcal G}_A:A\to C(M)$ будет множество $\{x\in A:\ \sup_{t\in T}|t(x)|\le \e\}=\e\cdot {^\circ T}$, то есть гомотетия поляры ${^\circ T}$ компакта $T$. Поскольку пространство $A$ стереотипно, ${^\circ T}$ является в нем окрестностью нуля. Это доказывает, что отображение ${\mathcal G}_A:A\to C(M)$ непрерывно, если пространство $C(M)$ наделено топологией равномерной сходимости на компактах в $M$.

Пусть далее $M=\Spec(A)$ -- паракомпактное локально компактное пространство. Для всякого компакта $K\subseteq M$ образ ${\mathcal G}_K(A)$ алгебры $A$ в $C(K)$ при отображении ${\mathcal G}_K$ будет инволютивной подалгеброй в $C(K)$, содержащей единицу (а значит и
все константы) и разделяющей точки $t\in K$. Поэтому, по теореме
Стоуна-Вейерштрасса, ${\mathcal G}_K(A)$ плотно в $C(K)$. Это верно для всякого
отображения ${\mathcal G}_K=\pi_K\circ\gamma$, где $K$ -- компакт в $M$.
Поскольку топология $C(M)$ есть проективная топология относительно отображений
$\pi_K$, мы получаем, что образ ${\mathcal G}_A(A)$ алгебры $A$ в $C(M)$ плотен
в $C(M)$. \epr

\btm\label{TH:Gelfand-Env_C_A=C(M)} Если $A$ -- коммутативная инволютивная
стереотипная алгебра с паракомпактным локально компактным инволютивным спектром
$M=\Spec(A)$, то ее преобразование Гельфанда ${\mathcal G}_A:A\to C(M)$ является ее непрерывной
оболочкой
\beq\label{Gelfand-Env_C_A=C(M)}
\Env_{\mathcal C} A=\mathcal{C}(M)
\eeq
 \etm
\bpr 
В данном случае сопряженное отображение спектров ${\mathcal G}_A^{\Spec}:\Spec(A)\gets M$ будет просто гомеоморфизмом, и тем более, точным наложением. Поэтому по теореме о ключевом представлении \ref{C-obolochka-podalgebry-v-C(M)} должно выполняться \eqref{Gelfand-Env_C_A=C(M)}.
 \epr

\noindent\rule{160mm}{0.1pt}\begin{multicols}{2}

\bex Непрерывная оболочка алгебры $\Trig(G)=k(G)$ тригонометрических многочленов на компактной группе $G$ совпадает с алгеброй ${\mathcal C}(G)$ непрерывных функций на $G$:
\beq
\Env_{\mathcal C}\Trig(G)={\mathcal C}(G).
\eeq
\eex
\bpr
Согласно \cite[(30.30)]{Hewitt-Ross-2}, $\Spec(\Trig(G))=G$. Теперь можно применить теорему \ref{TH:Gelfand-Env_C_A=C(M)}.
\epr

\end{multicols}\noindent\rule[10pt]{160mm}{0.1pt}

\subsection{Предварительные результаты}

\paragraph{Непрерывные оболочки аугментированных стереотипных алгебр.}

\btm\label{LM:A-augste=>Env_C-A-augste}
Пусть $(A,\e)$ --- инволютивная стереотипная алгебра с аугментацией. Тогда
 \bit{

 \item[(i)] непрерывная оболочка $\Env_{\mathcal C} \e:\Env_{\mathcal C} A\to\Env_{\mathcal C} \C=\C$ аугментации $\e$ на $A$ является аугментацией на непрерывной оболочке $\Env_{\mathcal C} A$ алгебры $A$;

 \item[(ii)] оболочка  $\env_{\mathcal C} A:A\to\Env_{\mathcal C} A$ является морфизмом аугментированых инволютивных стереотипных алгебр.
    \beq\label{env_C-e}
     \xymatrix @R=2.pc @C=4.pc
 {
 A\ar[r]^{\env_{\mathcal C} A}\ar[d]_{\e} &  \Env_{\mathcal C} A \ar@{-->}[d]^{\Env_{\mathcal C} \e} \\
 \C\ar[r]_{\id_{\C}=\env_{\mathcal C} \C} & \C
 }
    \eeq
  }\eit
\etm
\bpr
Это сразу следует из диаграммы \eqref{env_C-e}, которая коммутативна в силу \eqref{DIAGR:funktorialnost-env_C} и \eqref{Env_C(C)=C}.
\epr

\btm\label{TH:Env_C-functor-v-AugSteAlg}
Если $\ph:(A,\e_A)\to(B,\e_B)$ -- морфизм в категории аугментированных инволютивных стереотипных алгебр, то его непрерывная оболочка $\Env_{\mathcal C} \ph:(\Env_{\mathcal C} A,\Env_{\mathcal C} \e_A)\to(\Env_{\mathcal C} B,\Env_{\mathcal C} \e_B)$ -- тоже морфизм в категории аугментированных инволютивных стереотипных алгебр в силу коммутативности диаграммы
    \beq\label{Env_C-functor-v-AugSteAlg}
     \xymatrix @R=2.pc @C=4.pc
 {
 A\ar[rr]^{\env_{\mathcal C} A}\ar[dd]_{\ph}\ar[dr]_{\e_A} & &  \Env_{\mathcal C} A \ar@{-->}[dd]^{\Env_{\mathcal C} \ph}\ar@{-->}[dl]^{\ \Env_{\mathcal C} \e_A} \\
 & \C & \\
 B\ar[rr]_{\env_{\mathcal C} B}\ar[ur]^{\e_B} & & \Env_{\mathcal C} B\ar@{-->}[ul]_{\ \Env_{\mathcal C} \e_B}
 }
    \eeq
\etm
\bpr
Здесь коммутативны периметр и все внутренние треугольники, кроме правого (у которого все стороны пунктирны), и его коммутативность мы должны доказать. Но он будет коммутативен из-за того, что $\env_{\mathcal C} A$ --- эпиморфизм.
\epr

\bprop\label{LM:exten_C=>exten_C^Aug}
Если $A$ --- инволютивная стереотипная алгебра с аугментацией $\e:A\to\C$, и $\sigma:A\to A'$ --- ее непрерывное расширение, то алгебра $A'$ обладает единственной аугментацией $\e':A'\to\C$, для которой морфизм $\sigma$ становится морфизмом аугментированных стереотипных алгебр, и более того, расширением в категории $\AugInvSteAlg$ аугментированных инволютивных стереотипных алгебр в классе $\DEpi$ плотных эпиморфизмов относительно класса $\AugC*$ аугментированных $C^*$-алгебр.
\eprop
\bpr
1. Рассмотрим диаграмму:
    \beq\label{exten_C=>exten_C^Aug}
     \xymatrix 
 {
 A\ar[rr]^{\sigma}\ar[rd]_{\e} & & A' \ar@{-->}[dl]^{\e'} \\
 & \C &
 }
    \eeq
Поскольку $\C$ --- $C^*$-алгебра, морфизм $\e'$ существует и однозначно определен. Коммутативность этой диаграммы означает, что $\sigma$ будет морфизмом аугментированных стереотипных алгебр $\sigma:(A,\e)\to(A',\e')$.

2. Покажем, что этот морфизм $\sigma:(A,\e)\to(A',\e')$ является расширением в $\DEpi$ относительно $\AugC*$. Пусть $\ph:(A,\e)\to(B,\delta)$ --- морфизм в аугментированную $C^*$-алгебру. Рассмотрим диаграмму в категории стеретипных алгебр:
$$
     \xymatrix 
 {
 A\ar@/_5ex/[rdd]_{\ph}\ar[rr]^{\sigma}\ar[rd]_{\e} & & A' \ar[dl]^{\e'}\ar@{-->}@/^5ex/[ldd]^{\ph'} \\
 & \C & \\
 & B\ar[u]^{\delta} &
 }
$$
В ней пунктирная стрелка, $\ph'$, существует, единственна и замыкает периметр, поскольку $B$ --- $C^*$-алгебра, а $\sigma$ --- непрерывная оболочка. Одновременно верхний внутренний треугольник коммутативен, поскольку это просто диаграмма \eqref{exten_C=>exten_C^Aug}, а левый внутренний треугольник коммутативен, поскольку $\ph:(A,\e)\to(B,\delta)$ --- морфизм аугментированных алгебр. Вдобавок $\sigma$ --- эпиморфизм, поэтому правый внутренний треугольник тоже должен быть коммутативен:
$$
\delta\circ\ph'=\e'.
$$
Это означает, что $\ph'$ --- морфизм в категории $\AugInvSteAlg$, и поскольку он единственен, $\sigma$ --- расширение в $\AugInvSteAlg$ (в классе $\DEpi$ относительно класса $\AugC*$ аугментированных $C^*$-алгебр).
\epr

\bcor
Для всякой стереотипной алгебры $A$ с аугментацией $\e:A\to\C$ существует единственный морфизм стереотипных алгебр с аугментацией $\upsilon:(A,\e)\to\Env_{\AugC*}^{\DEpi}(A,\e)$, замыкающий диаграмму
    \beq\label{DEF:Env_C-e}
     \xymatrix 
 {
 & (A,\e)\ar[dl]_{\env_{\mathcal C} A}\ar[dr]^{\env_{\AugC*}^{\DEpi}(A,\e)} & \\
 (\Env_{\mathcal C} A,\Env_{\mathcal C} \e) \ar@{-->}[rr]_{\upsilon}  & & \Env_{\AugC*}^{\DEpi}(A,\e)
 }
    \eeq
\ecor
\bpr
В обозначениях предложения \ref{LM:exten_C=>exten_C^Aug}, здесь $\sigma$ --- не просто расширение, а оболочка $\env_{\mathcal C} A$. И в этом случае аугментация $\e'$ на $A'=\Env_{\mathcal C} A$ --- это морфизм $\Env_{\mathcal C} \e$ из \eqref{env_C-e}. Поскольку по предложению \ref{LM:exten_C=>exten_C^Aug} $\sigma=\env_{\mathcal C} A$ --- расширение в категории аугментированных стереотипных алгебр, должен существовать и быть единственным морфизм $\upsilon$ в диаграмме \eqref{DEF:Env_C-e}.
\epr

\section{Непрерывные оболочки групповых алгебр}

\subsection{Алгебра ${\mathcal K}(G)$}

Для всякой локально компактной группы $G$ ее групповая алгебра мер ${\mathcal C}^\star(G)$ является инволютивной алгеброй Хопфа относительно проективного стереотипного тензорного произведения $\circledast$. Поэтому, по теоремам  \ref{TH:C-obolochka-sohranyaet-Hopfov} и \ref{TH:C-obolochka-sohranyaet-inv-Hopfov}, ее непрерывная оболочка $\Env_{\mathcal C}{\mathcal C}^\star(G)$ должна быть коалгеброй с согласованными антиподом и инволюцией в категориях  ${\mathcal C}\text{-}{\tt Alg}$ непрерывных алгебр и $({\tt Ste},\odot)$ стереотипных пространств. Обозначим символом ${\mathcal K}(G)$ пространство, стереотипно сопряженное к $\Env_{\mathcal C}{\mathcal C}^\star(G)$:
 \beq\label{DEF:K(G)}
{\mathcal K}(G):=\Big(\Env_{\mathcal C}{\mathcal C}^\star(G)\Big)^\star.
 \eeq

\paragraph{Структурное представление ${\mathcal K}(G)$.}

\btm\label{TH:K(G)=lim}
Алгебра ${\mathcal K}(G)$ как стереотипное пространство представима в виде узлового кообраза (в категории стереотипных пространств)
 \beq\label{K(G)=lim}
{\mathcal K}(G)=\Coim_{\infty}\ph^\star
 \eeq
отображения $\ph^\star$, сопряженного к естественному морфизму стереотипных пространств
$$
\ph:{\mathcal C}^\star(G)\to \projlim_{p\in\Sn({\mathcal C}^\star(G))}{\mathcal C}^\star(G)/p.
$$
\etm
\bpr
Рассмотрим диаграмму \eqref{DIAGR:C-envelope=im_infty-lim N_X} для $A={\mathcal C}^\star(G)$
$$
\xymatrix @R=3.pc @C=4.0pc 
{
{\mathcal C}^\star(G)\ar[d]_{\coim_\infty\ph}\ar[rrr]^{\ph=\projlim_{p\in\Sn({\mathcal C}^\star(G))}\rho_p} & & & \projlim_{p\in\Sn({\mathcal C}^\star(G))} {\mathcal C}^\star(G)/p &   \\
\Coim_\infty\ph\ar[rrr]_{\red_\infty\ph} & & & \Im_\infty \ph \ar[u]_{\im_\infty\ph}\ar@{=}[r] & \Env_{\mathcal C}{\mathcal C}^\star(G)
}
$$
Сопряженная диаграмма имеет вид
\beq\label{DIAGR:K(G)}
\xymatrix @R=2.pc @C=6.0pc 
{
{\mathcal C}(G) & \Big(\projlim_{p\in\Sn({\mathcal C}^\star(G))} {\mathcal C}^\star(G)/p\Big)^\star \ar[l]_{\ph^\star} \ar[d]^{\coim_\infty\ph^\star} & \injlim_{p\in\Sn({\mathcal C}^\star(G))}\Big({\mathcal C}^\star(G)/p\Big)^\star \ar@{=}[l]\\
\Im_\infty\ph^\star\ar[u]^{\im_\infty\ph^\star}
 &  \Coim_\infty \ph^\star\ar[l]^{\red_\infty\ph^\star} & {\mathcal K}(G)=\Big(\Env_{\mathcal C}{\mathcal C}^\star(G)\Big)^\star\ar@{=}[l]
}
\eeq
\epr

\paragraph{Иммерсия ${\mathcal K}(G)\subarr {\mathcal C}(G)$.}

По теореме \ref{TH:C-obolochka-sohranyaet-Hopfov}, $\Env_{\mathcal C}{\mathcal C}^\star(G)$ является коалгеброй в $({\tt Ste},\odot)$, а по теореме \ref{TH:C-obolochka-sohranyaet-inv-Hopfov}, в ней имеются согласованные между собой антипод и инволюция. Отсюда в силу \cite[$4^\circ$ на с. 655 и $4^\circ$ на с.657]{Akbarov-De-Gruyter-I} справедлива

\btm Для всякой локально компактной группы $G$ пространство ${\mathcal K}(G)$ является алгеброй в категории $({\tt Ste},\circledast)$ (то есть стереотипной алгеброй) с согласованными антиподом и инволюцией.
\etm

По теореме \ref{TH:C-obolochka-sohranyaet-inv-Hopfov}(ii) морфизм
\beq\label{DEF:env_C^star}
\env_{\mathcal C}^\star=\big(\env_{\mathcal C}{{\mathcal C}^\star(G)}\big)^\star:{\mathcal K}(G)=\Big(\Env_{\mathcal C}{\mathcal C}^\star(G)\Big)^\star\to {\mathcal C}^\star(G)^\star={\mathcal C}(G),
\eeq
сопряженный к морфизму оболочки, является инволютивным гомоморфизмом алгебр:
\beq\label{env_C^star-homomorphism}
\env_{\mathcal C}^\star1=1,\qquad  \env_{\mathcal C}^\star(u\cdot v)=\env_{\mathcal C}^\star(u)\cdot\env_{\mathcal C}^\star(v),\qquad\env_{\mathcal C}^\star \overline{u}=\overline{\env_{\mathcal C}^\star u},
\qquad u,v\in {\mathcal K}(G)
\eeq

Покажем, что у него нулевое ядро.\label{Ker-e_C-C^star(G)^star=0} Действительно, группа $G$ вкладывается в алгебру $\Env_{\mathcal C}{\mathcal C}^\star(G)$ как композиция дельта-отображения и оболочки
$$
G\overset{\delta}{\longrightarrow}{\mathcal C}^\star(G)\overset{\env_{\mathcal C}{{\mathcal C}^\star(G)}}{\longrightarrow}\Env_{\mathcal C}{\mathcal C}^\star(G).
$$
При этом образ $G$ полон (то есть линейные комбинации его элементов плотны) в ${\mathcal C}^\star(G)$ (в силу \cite[Lemma 8.2]{Ak03}), а образ ${\mathcal C}^\star(G)$ плотен в $\Env_{\mathcal C}{\mathcal C}^\star(G)$. Поэтому образ $G$ полон в $\Env_{\mathcal C}{\mathcal C}^\star(G)$. Отсюда следует, что всякий элемент $u\in {\mathcal K}(G):=\Big(\Env_{\mathcal C}{\mathcal C}^\star(G)\Big)^\star$ однозначно определяется композицией
$$
u\circ \env_{\mathcal C}{{\mathcal C}^\star(G)}\circ\delta:G\to\C,
$$
которую можно понимать как ограничение $u$ на группу $G$. В частности, если эта композиция равна нулю в ${\mathcal C}(G)$, то $u=0$ в ${\mathcal K}(G)$.

Важный для нас вывод состоит в том, что ${\mathcal K}(G)$ можно понимать как некую инволютивную подалгебру в ${\mathcal C}(G)$:

\btm\label{TH:K(G)->C(G)} Для всякой локально компактной группы $G$  отображение $u\mapsto u\circ \env_{\mathcal C}{{\mathcal C}^\star(G)}\circ\delta$
совпадает с отображением $\env_{\mathcal C}{{\mathcal C}^\star(G)}^\star$, сопряженным к $\env_{\mathcal C}{{\mathcal C}^\star(G)}$:
\beq\label{K(G)->C(G)}
\env_{\mathcal C}{{\mathcal C}^\star(G)}^\star(u)=u\circ \env_{\mathcal C}{{\mathcal C}^\star(G)}\circ\delta
\eeq
и инъективно и гомоморфно вкладывает ${\mathcal K}(G)$ в ${\mathcal C}(G)$ в качестве инволютивной подалгебры (и поэтому операции сложения, умножения и инволюции в ${\mathcal K}(G)$ являются поточечными).
\etm

\paragraph{Цепочка $\Trig(G)\subseteq k(G)\subseteq {\mathcal K}(G)\subseteq {\mathcal C}(G)$.}

Вспомним о пространствах функций $\Trig(G)$ и $k(G)$, определенных на с.\pageref{DEF:k(G)} и \pageref{DEF:Trig(G)}.

\btm\label{TH:V-k-K-C} Для всякой локально компактной группы $G$ справедлива цепочка теоретико-множественных включений,
 \beq\label{V-k-K-C}
\Trig(G)\subseteq k(G)\subseteq {\mathcal K}(G)\subseteq {\mathcal C}(G),
 \eeq
причем
\bit{
\item[(i)] всегда
 \beq\label{Trig-plotno-v-K}
\overline{\Trig(G)}={\mathcal K}(G),
 \eeq

\item[(ii)] если $G$ -- компактная группа, то
 \beq\label{Trig=K}
\Trig(G)=k(G)={\mathcal K}(G)
 \eeq
 }\eit
\etm
\bpr
1. В цепочке \eqref{V-k-K-C} первое включение очевидно, а третье уже отмечалось в теореме \ref{TH:K(G)->C(G)}. Докажем второе: $k(G)\subseteq {\mathcal K}(G)$. Пусть $u\in k(G)$, то есть выполняется \eqref{DEF:k(G)}, где $\pi:G\to{\mathcal B}(X)$ -- непрерывное по норме унитарное представление. По теореме \ref{TH:nepr-po-norme-predstavleniya} $\pi$ порождает некий (непрерывный) гомоморфизм инволютивных алгебр $\psi:{\mathcal C}^\star(G)\to{\mathcal B}(X)$. Поскольку здесь ${\mathcal B}(X)$ -- $C^*$-алгебра, он должен продолжаться до (непрерывного) гомоморфизма инволютивных алгебр $\psi':\Env_{\mathcal C}{\mathcal C}^\star(G)\to{\mathcal B}(X)$. Рассмотрим функционал
$$
f(\beta)=\langle \psi'(\beta)x,y\rangle,\qquad \beta\in \Env_{\mathcal C}{\mathcal C}^\star(G)
$$
на ${\mathcal K}^\star(G)$. Он будет порождать функцию на $G$, совпадающую с $u$:
$$
\xymatrix @R=3.pc @C=6.0pc 
{
G\ar@/_4ex/[ddr]_{u}\ar[r]^{\delta}\ar[dr]_{\pi} & {\mathcal C}^\star(G)\ar[r]^{\env_{\mathcal C}{{\mathcal C}^\star(G)}}
\ar@{-->}[d]_{\psi} & \Env_{\mathcal C}{\mathcal C}^\star(G)\ar@{-->}[dl]^{\psi'}\ar@/^4ex/@{-->}[ddl]^{f} \\
& {\mathcal B}(X)\ar[d] & \\
& \C &
}
$$
Это как раз и означает, что $u\in {\mathcal K}(G)$.

2. Покажем, что $\Trig(G)$ плотно в ${\mathcal K}(G)$. Каждую алгебру ${\mathcal C}^\star(G)/p$ можно изометрически вложить в алгебру вида ${\mathcal B}(X)$:
$$
{\mathcal C}^\star(G)/p\to {\mathcal B}(X).
$$
При этом, во-первых, по теореме Хана-Банаха всякий функционал $f\in ({\mathcal C}^\star(G)/p)^\star$ продолжается до некоторого функционала $g\in{\mathcal B}(X)^\star$, а, во-вторых, всякий функционал $g\in{\mathcal B}(X)^\star$ приближается в ${\mathcal B}(X)^\star$ линейными комбинациями чистых состояний, то есть (пример \ref{EX:trig-monom}) функционалами, порождающими функции из $\Trig(G)$. Это значит, что при сопряженном отображении
$$
({\mathcal C}^\star(G)/p)^\star\gets{\mathcal B}(X)^\star
$$
функционалы, порождающие функции из $\Trig(G)$, переходят в плотное подмножество в $({\mathcal C}^\star(G)/p)^\star$. Поскольку это верно для каждой полунормы $p\in\Sn({\mathcal C}^\star(G))$, мы получаем, что функционалы из 
$$
\Big(\Ste\text{-}\kern-5pt\projlim_{p\in\Sn({\mathcal C}^\star(G))} {\mathcal C}^\star(G)/p\Big)^\star=\Ste\text{-}\kern-8pt\injlim_{p\in\Sn({\mathcal C}^\star(G))}\Big({\mathcal C}^\star(G)/p\Big)^\star,
$$
порождающие функции из $\Trig(G)$, плотны в
$$
\Big(\Ste\text{-}\kern-5pt\projlim_{p\in\Sn({\mathcal C}^\star(G))} {\mathcal C}^\star(G)/p\Big)^\star=\Ste\text{-}\kern-8pt\injlim_{p\in\Sn({\mathcal C}^\star(G))}\Big({\mathcal C}^\star(G)/p\Big)^\star.
$$
С другой стороны, в диаграмме \eqref{DIAGR:K(G)} видно, что
$$
\Big(\Ste\text{-}\kern-5pt\projlim_{p\in\Sn({\mathcal C}^\star(G))} {\mathcal C}^\star(G)/p\Big)^\star=\Ste\text{-}\kern-8pt\injlim_{p\in\Sn({\mathcal C}^\star(G))}\Big({\mathcal C}^\star(G)/p\Big)^\star
$$
плотно отображается в ${\mathcal K}(G)$. Вместе это значит, что $\Trig(G)$ плотно в ${\mathcal K}(G)$.
\epr

\paragraph{Сдвиг в ${\mathcal K}(G)$.}

\btm\label{LM:sdvig-v-K(G)}
Для всякой локально компактной группы $G$ сдвиг (правый и левый) на произвольный элемент $a\in G$ является изоморфизмом стереотипной алгебры ${\mathcal K}(G)$.
\etm
\bpr
Докажем это для оператора левого сдвига: пусть
$$
M_a:{\mathcal C}^\star(G)\to{\mathcal C}^\star(G)\quad\Big|\quad  M_a(\alpha)=\delta^a*\alpha,\qquad \alpha\in {\mathcal C}^\star(G).
$$
Обозначим $\eta_a=\env_{\mathcal C}(\delta^a)$ и положим
$$
N_a:\Env_{\mathcal C}{\mathcal C}^\star(G)\to
\Env_{\mathcal C}{\mathcal C}^\star(G)\quad\Big|\quad N_a(\omega)=\eta_a*\omega,\qquad \omega\in \Env_{\mathcal C}{\mathcal C}^\star(G)
$$
(здесь $*$ -- умножение в алгебре $\Env_{\mathcal C}{\mathcal C}^\star(G)$). Тогда
$$
\env_{\mathcal C}(M_a(\alpha))=\env_{\mathcal C}(\delta^a*\alpha)=\env_{\mathcal C}(\delta^a)*\env_{\mathcal C}(\alpha)=N_a(\env_{\mathcal C}(\alpha)),
$$
то есть коммутативна диаграмма
\beq\label{PROOF:LM:sdvig-v-K(G)}
\xymatrix @R=2.pc @C=8.0pc 
{
{\mathcal C}^\star(G)\ar[r]^{\env_{\mathcal C}}\ar[d]_{M_a} & \Env_{\mathcal C}{\mathcal C}^\star(G)\ar[d]^{N_a}\\
{\mathcal C}^\star(G)\ar[r]^{\env_{\mathcal C}} & \Env_{\mathcal C}{\mathcal C}^\star(G)
}
\eeq
Как следствие, коммутативна двойственная диаграмма
$$
\xymatrix @R=2.pc @C=8.0pc 
{
{\mathcal C}(G) & {\mathcal K}(G)\ar[l]_{\env_{\mathcal C}^\star}\\
{\mathcal C}(G)\ar[u]^{M_a^\star} & {\mathcal K}(G)\ar[l]_{\env_{\mathcal C}^\star}\ar[u]_{N_a^\star}
}
$$
Из нее следует, что, во-первых, $N_a^\star$ -- оператор левого сдвига на элемент $a$ в функциональной алгебре ${\mathcal K}(G)$, потому что при гомоморфной инъекции $\env_{\mathcal C}^\star$ (описанной в теореме \ref{TH:K(G)->C(G)}) он превращается в оператор $M_a^\star$ сдвига на элемент $a$ на алгебре ${\mathcal C}(G)$. И, во-вторых, что $N_a^\star$ -- гомоморфизм алгебр, потому что
\begin{multline*}
\env_{\mathcal C}^\star N_a^\star(u\cdot v)=M_a^\star\Big(\env_{\mathcal C}^\star(u\cdot v)\Big)=\eqref{env_C^star-homomorphism}=
M_a^\star\Big(\env_{\mathcal C}^\star(u)\cdot\env_{\mathcal C}^\star(v)\Big)=\\=
M_a^\star\env_{\mathcal C}^\star(u)\cdot M_a^\star\env_{\mathcal C}^\star(v)=
\env_{\mathcal C}^\star(N_a^\star u)\cdot\env_{\mathcal C}^\star(N_a^\star v)=\eqref{env_C^star-homomorphism}=
\env_{\mathcal C}^\star(N_a^\star u \cdot N_a^\star v)
\end{multline*}
и, поскольку, по теореме \ref{TH:K(G)->C(G)}, $\env_{\mathcal C}^\star$ -- инъективное отображение,
$$
N_a^\star(u\cdot v)=N_a^\star u \cdot N_a^\star v.
$$
\epr

\subsection{Непрерывные оболочки групповых алгебр}

\paragraph{Дискретная группа.}

Для дискретной группы $D$ ее групповая алгебра представляет собой алгебру функций на $D$ с конечным носителем:
$$
{\mathcal C}^\star(D)=\C_D=\{\alpha=\{\alpha_x,x\in D\}:\quad \card\{x\in D:\alpha_x\ne 0\}<\infty\}.
$$
Свертка в $\C_D$ определяется своим действием на дельта-функционалах \eqref{delta^a*delta^b=delta^(a-cdot-b)}.

Заметим, что всякая $C^*$-полунорма $p$ на $\C_D$ переводит единицу либо в нуль, либо в единицу:
$$
p(\delta^e)=p(\delta^e*\delta^e)=p(\delta^e*(\delta^e)^\bullet)=p(\delta^e)^2\quad\Longrightarrow\quad p(\delta^e)=0\quad\vee\quad p(\delta^e)=1.
$$
В первом случае $p$ будет вообще любой элемент переводить в нуль (потому что $p$ всегда субмультипликативна). Поэтому если $p\ne 0$, то $p(\delta^e)=1$. Более того, в этом случае вообще любой дельта-функционал переводится в единицу:
$$
1=p(\delta^e)=p(\delta^a*\delta^{a^{-1}})=p(\delta^a*(\delta^a)^\bullet)=p(\delta^a)^2\quad\Longrightarrow\quad p(\delta^a)=1.
$$
Из этого, в свою очередь следует, что всякая $C^*$-полунорма $p$ на $\C_D$ подчинена $\ell_1$-норме:
\beq\label{||pi(alpha)||-le-||alpha||_1}
p(\alpha)\le \norm{\alpha}_1,\qquad \alpha\in\C_D,
\eeq
потому что
$$
p(\alpha)=p\Big(\sum_{x\in D}\alpha_x\cdot\delta^x\Big)\le
\sum_{x\in D}|\alpha_x|\cdot p(\delta^x)\le\sum_{x\in D}|\alpha_x|\cdot 1=\norm{\alpha}_1.
$$
Из \eqref{||pi(alpha)||-le-||alpha||_1} следует, что для всякого $\alpha\in\C_D$ существует точная грань по всем $C^*$-полунормам
\beq\label{||alpha||^*}
\norm{\alpha}_\bullet=\sup_{p\in \Sn(\C_D)}p(\alpha)\le \norm{\alpha}_1.
\eeq
Это будет $C^*$-полунорма на $\C_D$, потому что
$$
\norm{\alpha*\alpha^\bullet}_\bullet=\sup_{p\in \Sn(\C_D)}p(\alpha*\alpha^\bullet) =\sup_{p\in \Sn(\C_D)}p(\alpha)^2=
\Big(\sup_{p\in \Sn(\C_D)}p(\alpha)\Big)^2=
\norm{\alpha}_\bullet^2.
$$
Более того, это будет норма на $\C_D$, потому что если $\alpha\ne0$, то при левом регулярном представлении $\pi:D\to{\mathcal B}(L_2(D))$ оно переходит в ненулевой элемент, который отделяется от нуля нормой в ${\mathcal L}(L_2(D))$, а она определяет $C^*$-полунорму на $\C_D$. Пополнение алгебры $\C_D$ относительно этой нормы совпадает с пополнением $\ell_1(D)$ относительно нее, называется {\it групповой $C^*$-алгеброй} группы $D$ и обозначается $C^*(D)$ \cite{Dixmier}.

\bprop\label{PROP:nepr-obol-diskr-gruppy}
Для дискретной группы $D$ непрерывной оболочкой ее групповой алгебры ${\mathcal C}^\star(D)=\C_D$ является групповая $C^*$-алгебра $C^*(D)$:
$$
\Env_{\mathcal C}\C_D=C^*(D).
$$
\eprop
\bpr
Пусть $\rho:\C_D\to C^*(D)$ -- отображение пополнения относительно нормы $\norm{\cdot}_\bullet$. Покажем, что оно является непрерывным расширением. Пусть $\ph:\C_D\to B$ -- инволютивный гомоморфизм в какую-то $C^*$-алгебру $B$. Чтобы достроить диаграмму
\beq\label{diagr-ph-ph'}
 \xymatrix @R=2pc @C=1.2pc
 {
 \C_D\ar[rr]^{\rho}\ar[dr]_{\ph} & & C^*(D)\ar@{-->}[dl]^{\ph'} \\
  & B &
 }
\eeq
достаточно считать, что $\ph$ имеет плотный образ в $B$. Тогда $B$ можно считать пополнением алгебры $\C_D$ относительно какой-то $C^*$-полунормы $p$ (с факторизацией по ядру этой нормы). Но $p$, будучи $C^*$-полунормой, должна быть подчинена норме $\norm{\cdot}_\bullet$. Поэтому $\Ker\norm{\cdot}_\bullet\subseteq\Ker p$. Отсюда следует, что $\ph$ можно представить как композицию $\C_D\to \C_D/\Ker\norm{\cdot}_\bullet\to B$. А потом такое представление достраивается до \eqref{diagr-ph-ph'}.

Теперь покажем, что $\rho:\C_D\to C^*(D)$ -- непрерывная оболочка. Пусть $\sigma:\C_D\to A$ -- какое-то другое непрерывное расширение. Тогда, поскольку $\rho:\C_D\to C^*(D)$ -- гомоморфизм в $C^*$-алгебру, он должен пропускаться через $\sigma$:
$$
 \xymatrix @R=2pc @C=1.2pc
 {
 \C_D\ar[rr]^{\sigma}\ar[dr]_{\rho} & & A\ar@{-->}[dl]^{\upsilon} \\
  & C^*(D) &
 }
$$
\epr

\paragraph{Абелевы локально компактные группы.}

Пусть $C$ -- коммутативная локально компактная группа. Вспомним алгебру $\mathcal{C}(C)$ непрерывных функций и алгебру  $\mathcal{C}^\star(C)$ мер с компактным носителем на $C$.
Формула
\beq\label{Fourier-transform}
\overbrace{{\mathcal F}_C(\alpha)(\chi)}^{\scriptsize \begin{matrix}
\text{значение функции ${\mathcal F}_C(\alpha)\in {\mathcal C}(\widehat{C})$}\\
\text{в точке $\chi\in \widehat{C}$} \\ \downarrow \end{matrix}}\kern-35pt=\kern-50pt\underbrace{\alpha(\chi)}_{\scriptsize \begin{matrix}\uparrow \\
\text{действие функционала $\alpha\in{\mathcal C}^\star(C)$}\\
\text{на функцию $\chi\in \widehat{C}\subseteq {\mathcal C}(C)$ }\end{matrix}}
\kern-50pt \qquad (\chi\in \widehat{C},\quad \alpha\in {\mathcal C}^\star(C))
\eeq
определяет отображение
$$
{\mathcal F}_C:{\mathcal C}^\star(C)\to {\mathcal C}(\widehat{C})
$$
являющееся гомоморфизмом инволютивных стереотипных алгебр, и называемое {\it преобразованием Фурье} на группе $C$.

Следующий факт был доказан в \cite[Theorem 5.53]{Ak16} (для оболочек Кузнецовой в \cite[Theorem 2.11]{Kuznetsova}).

\bprop\label{PROP:Env_C-C^star(G)=C(widehat(G))}
Преобразование Фурье на коммутативной локально компактной группе $C$ является непрерывной оболочкой групповой алгебры ${\mathcal C}^\star(C)$. Как следствие,
\beq\label{Env_C-C^star(G)=C(widehat(G))}
\Env_{\mathcal C}{\mathcal C}^\star(C)={\mathcal C}(\widehat{C}).
\eeq
\eprop
\bpr Спектр алгебры ${\mathcal C}^\star(C)$ гомеоморфен
двойственной группе $\widehat{C}$, поэтому все следует из теоремы \ref{C-obolochka-podalgebry-v-C(M)}. \epr

\paragraph{Компактные группы.}

\bprop\label{PROP:nepr-obol-komp-gruppy}
Для компактной группы $K$ непрерывной оболочкой ее групповой алгебры ${\mathcal C}^\star(K)$ является декартово произведение алгебр ${\mathcal B}(X_\pi)$, где $\pi$ пробегает двойственный объект $\widehat{K}$, а $X_\pi$ представляет собой пространство представления $\pi$:
\beq\label{nepr-obol-komp-gruppy}
\Env_{\mathcal C}{\mathcal C}^\star(K)=\prod_{\pi\in\widehat{K}}{\mathcal B}(X_\pi).
\eeq
\eprop
\bpr
1. Сначала покажем, что отображение $P=\prod_{\pi\in\widehat{K}}\pi:{\mathcal C}^\star(K)\to\prod_{\pi\in\widehat{K}}{\mathcal B}(X_\pi)$ является непрерывным расширением. Заметим сразу, что по теореме \ref{TH:plotnost-C*(K)-v-prod-B(X_sigma)}, композиция отображения $P$ с каждой конечной проекцией $\prod_{\pi\in\widehat{K}}{\mathcal B}(X_\pi)\to \prod_{\pi\in S}{\mathcal B}(X_\pi)$, $S\in 2_{\widehat{K}}$, является сюръекцией. Отсюда следует, что $P$ -- плотный эпиморфизм.  Далее, пусть $\psi:{\mathcal C}^\star(K)\to B$ -- инволютивный гомоморфизм в какую-то $C^*$-алгебру $B$. Рассмотрим какое-нибудь вложение $C^*$-алгебр $\eta:B\to {\mathcal B}(X)$. Композиция $\psi=\eta\circ\ph$ порождает непрерывное по норме представление $\rho=\psi\circ\delta:K\to {\mathcal B}(X)$, которое по теореме \ref{TH:nepr-po-norme-perdst-K} раскладывается в прямую сумму унитарных неприводимых представлений, среди которых только конечное число не эквивалентно друг другу. Это означает, в частности, что существует конечное множество $M\subseteq\widehat{K}$ такое, что $\psi$ представимо в виде композиции $\psi'\circ P_M$, где $P_M$ -- естественная проекция ${\mathcal C}^\star(K)$ в прямое произведение $\prod_{\pi\in M}{\mathcal B}(X_\pi)$. Отсюда в свою очередь следует, что гомоморфизм $\ph$ обнуляется на ядре $P_M$: $\Ker P_M\subseteq\Ker\ph$. Вдобавок, алгебра $\prod_{\pi\in M}{\mathcal B}(X_\pi)$ конечномерна и изоморфна фактор-алгебре ${\mathcal C}^\star(K)/\Ker P_M$. Поэтому $\ph$ представимо в виде компрозиции $\ph'\circ P_M$, и мы получаем диаграмму
$$
 \xymatrix @R=2pc @C=1.2pc
 {
{\mathcal C}^\star(K)\ar[rr]^{P_M}\ar[dr]_{\ph}\ar@/_4ex/[ddr]_{\psi} & &
\prod\limits_{\pi\in M}{\mathcal B}(X_\pi)\ar@{-->}[dl]^{\ph'} \ar@{-->}@/^4ex/[ddl]^{\psi'}
\\
  & B\ar[d]^{\eta} & \\
  & {\mathcal B}(X) &
 }
$$
Из нее в свою очередь строится диаграмма
$$
 \xymatrix @R=2pc @C=1.2pc
 {
{\mathcal C}^\star(K)\ar[rr]^{P}\ar[dr]_{P_M}\ar@/_4ex/[ddr]_{\ph} & & \prod\limits_{\pi\in\widehat{K}}{\mathcal B}(X_\pi)\ar@{-->}[dl]^{Q_M} \ar@{-->}@/^4ex/[ddl]^{\ph''}
\\
  & \prod\limits_{\pi\in M}{\mathcal B}(X_\pi) \ar[d]^{\ph'} & \\
  & B &
 }
$$
в которой $Q_M$ -- естественная проекция прямого произведения в свое подпроизведение.

2. Докажем далее, что  $P:{\mathcal C}^\star(K)\to\prod_{\pi\in\widehat{K}}{\mathcal B}(X_\pi)$ -- непрерывная оболочка. Пусть $Q:{\mathcal C}^\star(K)\to A$ -- какое-то другое непрерывное расширение. Тогда для всякого представления $\sigma\in\widehat{K}$ найдется (единственный) морфизм $\alpha_\sigma:A\to {\mathcal B}(X_\sigma)$, такой что $\ph_\sigma=\alpha_\sigma\circ Q$. Семейству морфизмов $\{\alpha_\sigma:A\to {\mathcal B}(X_\sigma);\ \sigma\in\widetilde{K}\}$ соответствует морфизм $\upsilon:A\to \prod_{\pi\in\widehat{K}}{\mathcal B}(X_\pi)$ такой что $\alpha_\sigma=\iota_\sigma\circ\upsilon$ для всякого $\sigma\in\widehat{K}$. Мы получаем диаграмму
$$
 \xymatrix @R=2pc @C=1.2pc
 {
& {\mathcal C}^\star(K)\ar@/_4ex/[ddl]_{Q}\ar@/^4ex/[ddr]^{P}\ar[d]^{\ph_\sigma} &
\\
& {\mathcal B}(X_\sigma) &
\\
A\ar[rr]_{\upsilon}\ar[ur]_{\alpha_\sigma} & & \prod\limits_{\pi\in\widehat{K}}{\mathcal B}(X_\pi)\ar[ul]_{\iota_\sigma}
 }
$$
в которой внутренние маленькие треугольники коммутативны в силу свойств отображений $P$, $Q$, $\upsilon$. Как следствие,
$$
\iota_\sigma\circ\upsilon\circ Q=\alpha_\sigma\circ Q=\ph_\sigma=\iota_\sigma\circ P,
$$
и поскольку это верно для любого $\sigma$, получаем
$$
\upsilon\circ Q=P,
$$
то есть в этой диаграмме коммутативен также и периметр. Единственность морфизма $\upsilon$ следует из единственности морфизмов $\alpha_\sigma$.
\epr

\paragraph{Произведение $Z\times K$ абелевой и компактной группы.}

Пусть $Z$ -- абелева локально компактная группа, а $K$ -- компактная группа (необязательно, абелева).

\bprop\label{TH:env_C^star(R^n-times-K)}
Непрерывной оболочкой групповой алгебры ${\mathcal C}^\star(Z\times K)$ является алгебра ${\mathcal C}(\widehat{Z},\prod_{\sigma\in\widehat{K}}{\mathcal B}(X_\sigma))$ непрерывных отображений из двойственной по Понтрягину группы $\widehat{Z}$ в
декартово произведение алгебр ${\mathcal B}(X_\sigma)$, где $\sigma$ пробегает двойственный объект $\widehat{K}$, а $X_\sigma$  представляет собой пространство представления $\sigma$:
\beq\label{env_C^star(R^n-times-K)}
\Env_{\mathcal C}{\mathcal C}^\star(Z\times K)={\mathcal C}\Big(\widehat{Z},\prod_{\sigma\in\widehat{K}}{\mathcal B}(X_\sigma)\Big)=
\prod_{\sigma\in\widehat{K}}{\mathcal C}\big(\widehat{Z},{\mathcal B}(X_\sigma)\big)={\mathcal C}(\widehat{Z})\odot \prod_{\sigma\in\widehat{K}}{\mathcal B}(X_\sigma)=\Env_{\mathcal C}{\mathcal C}^\star(Z)\odot \Env_{\mathcal C}{\mathcal C}^\star(K).
\eeq
\eprop
\bpr Первое равенство доказывается цепочкой
\begin{multline*}
\Env_{\mathcal C}{\mathcal C}^\star(Z\times K)=\Env_{\mathcal C}\Big({\mathcal C}^\star(Z)\circledast{\mathcal C}^\star(K)\Big)=
\eqref{E(X-otimes-Y)-cong-E(E(X)-otimes-E(Y))}=
\Env_{\mathcal C}\Big(\Env_{\mathcal C}{\mathcal C}^\star(Z)\circledast\Env_{\mathcal C}{\mathcal C}^\star(K)\Big)=\eqref{Env_C-C^star(G)=C(widehat(G))}=\\=
\Env_{\mathcal C}\Big({\mathcal C}(\widehat{Z})\circledast\prod_{\sigma\in\widehat{K}}{\mathcal B}(X_\sigma)\Big)=\eqref{C/circledast}=
{\mathcal C}(\widehat{Z})\overset{\mathcal C}{\circledast}\prod_{\sigma\in\widehat{K}}{\mathcal B}(X_\sigma)=\eqref{C(M)-circledast-A=C(M,A)}=
{\mathcal C}\Big(\widehat{Z},\prod_{\sigma\in\widehat{K}}{\mathcal B}(X_\sigma)\Big)
\end{multline*}
Второе равенство в \eqref{env_C^star(R^n-times-K)} очевидно, третье является следствием \eqref{C(M,X)-cong-C(M)-odot-X} и \cite[(4.151)]{Akbarov-De-Gruyter-I}, а последнее -- следствием \eqref{Env_C-C^star(G)=C(widehat(G))} и \eqref{nepr-obol-komp-gruppy}.
\epr

\paragraph{Надстройка $Z\cdot K$ абелевой группы $Z$ с помощью компактной группы $K$.}

Напомним, что понятие компактной надстройки абелевой группы (или надстройки абелевой группы с помощью компактной группы) было введено на с.\pageref{DEF:komp-nadstr-abelevoi-gruppy}. По лемме \ref{LM:Z-cdot-K=(Z-times-K)/C} всякая надстройка $G=Z\cdot K$ абелевой локально компактной группы $Z$ с помощью компактной группы $K$ представима в виде фактор-группы
$$
G\cong(Z\times K)/\iota(H),
$$
в которой $H=Z\cap K$ --- подгруппа в $G$, а $\iota$ --- вложение $H$ в декартово произведение
$$
\iota:H\to Z\times K\quad\Big|\quad \iota(x)=(x,x^{-1}),\quad x\in H.
$$

Гомоморфизм групп $\iota:H\to Z\times K$ порождает морфизм групповых алгебр ${\mathcal C}^\star(\iota):{\mathcal C}^\star(H)\to {\mathcal C}^\star(Z\times K)$, который в свою очередь порождает морфизм оболочек $\Env_{\mathcal C}{\mathcal C}^\star(\iota):\Env_{\mathcal C}{\mathcal C}^\star(H)\to \Env_{\mathcal C}{\mathcal C}^\star(Z\times K)$. В результате возникает диаграмма
\beq\label{DEF:xi:Env-C*(H)->Env-C*(Z-times-K)}
 \xymatrix  @R=2.pc @C=6.pc
{
 H\ar[r]^{\iota}\ar[d]_{\delta_H} & Z\times K \ar[d]^{\delta_{Z\times K}}
 \\
{\mathcal C}^\star(H)\ar[r]^{{\mathcal C}^\star(\iota)}
 \ar[d]_{\env_{\mathcal C}{\mathcal C}^\star(H)} & {\mathcal C}^\star(Z\times K)
 \ar[d]^{\env_{\mathcal C}{\mathcal C}^\star(Z\times K)} \\
\Env_{\mathcal C}{\mathcal C}^\star(H)\ar[r]^{\Env_{\mathcal C}{\mathcal C}^\star(\iota)}\ar@{=}[d] & \Env_{\mathcal C}{\mathcal C}^\star(Z\times K)\ar@{=}[d]\\
{\mathcal C}(\widehat{H})\ar[r]^{\xi} & \prod_{\sigma\in\widehat{K}}{\mathcal C}\big(\widehat{Z},{\mathcal B}(X_\sigma)\big)
}
\eeq
в которой $\xi$ -- морфизм, соответствующий $\Env_{\mathcal C}{\mathcal C}^\star(\iota)$ при отождествлении, описываемом вертикальными равенствами. Заметим, что по теореме \ref{TH:Env_C-functor-v-AugSteAlg} последние три горизонтальные стрелки --- морфизмы аугментированных стереотипных алгебр.

Ниже нам понадобится описание действия морфизма $\xi$. Зафиксируем $\sigma\in\widehat{K}$ и рассмотрим отображение
$$
\xi_\sigma:{\mathcal C}(\widehat{H})\to {\mathcal C}(\widehat{Z},{\mathcal B}(X_\sigma))
$$
которое каждой функции $f\in {\mathcal C}(\widehat{H})$ ставит в соответствие проекцию на компоненту ${\mathcal C}(\widehat{Z},{\mathcal B}(X_\sigma))$ ее образа $\xi(f)\in \prod_{\sigma\in\widehat{K}}{\mathcal C}\big(M_\sigma,{\mathcal B}(X_\sigma)\big)$:
$$
\xi_\sigma(f)=\xi(f)_\sigma, \qquad f\in {\mathcal C}(\widehat{H}).
$$
На элементы $h\in H$, рассматриваемые как функции $h:\widehat{H}\to\C$, это отображение действует по формуле
\beq\label{xi_sigma(h)(chi)}
\xi_\sigma(h)(\chi)=\underbrace{h(\chi\big|_H)}_{\scriptsize \begin{matrix}\|\\ \chi(h) \end{matrix}}\cdot\sigma(h^{-1})=h(\chi\big|_H)\cdot\overline{\sigma(h)},\qquad \chi\in \widehat{Z},\ h\in H.
\eeq
А действие $\xi_\sigma$ на произвольной функции $f\in {\mathcal C}(\widehat{H})$ легче всего описать так: нужно функцию $f$ приблизить в пространстве ${\mathcal C}(\widehat{H})$ линейными комбинациями характеров на $\widehat{H}$, то есть элементами группы $H$,
\beq\label{f=sum_i-lambda_i-cdot-h_i}
f\longleftarrow\sum_{i}\lambda_i\cdot h_i,\qquad h_i\in H,\ \lambda_i\in \C.
\eeq
и тогда получится, что  
\beq\label{xi_sigma(f)(chi)}
\xi_\sigma(f)(\chi)\longleftarrow\sum_{i}\lambda_i\cdot\xi_\sigma(h_i)(\chi)=\eqref{xi_sigma(h)(chi)}=
\sum_{i}\lambda_i\cdot h_i(\chi\big|_H)\cdot \overline{\sigma(h_i)}.
\eeq
Здесь предел существует и единственен, потому что, как видно в диаграмме \eqref{DEF:xi:Env-C*(H)->Env-C*(Z-times-K)}, предел   
\eqref{f=sum_i-lambda_i-cdot-h_i} в ${\mathcal C}(\widehat{H})$ должен под действием непрерывного отображения $\xi$ превращаться в некий предел в пространстве $\prod_{\sigma\in\widehat{K}}{\mathcal C}\big(\widehat{Z},{\mathcal B}(X_\sigma)\big)$, и для фиксированного $\sigma\in\widehat{K}$ в некий предел в пространстве ${\mathcal C}\big(\widehat{Z},{\mathcal B}(X_\sigma)\big)$.

\bit{

\item[$\bullet$] Для всякого $\sigma\in\widehat{K}$ рассмотрим множество $M_\sigma\subseteq \widehat{Z}$, состоящее из характеров $\chi\in\widehat{Z}$, для которых справедливо тождество  
\beq\label{charakterizatsija-M_sigma}
\sigma(h)=h(\chi\big|_H)\cdot 1_{{\mathcal B}(X_\sigma)}, \qquad h\in H.
\eeq
То есть 
\beq\label{1^0:stroenie-M_sigma}
M_\sigma=\{\chi\in\widehat{Z}:\quad \forall h\in H\quad \sigma(h)=h(\chi\big|_H)\cdot 1_{{\mathcal B}(X_\sigma)} \}.
\eeq
Очевидно, $M_\sigma$ --- замкнутое подмножество в $\widehat{Z}$.

\item[$\bullet$] Для всякого $\sigma\in\widehat{K}$ пусть $I_\sigma$ обозначает замкнутый двусторонний идеал алгебры  ${\mathcal C}\big(M_\sigma,{\mathcal B}(X_\sigma)\big)$, порожденный образом $\xi_\sigma(\e_{{\mathcal C}(\widehat{H})}^{-1}(0))$ идеала 
\beq\label{e_(C(hat(H)))^(-1)(0)}
\e_{{\mathcal C}(\widehat{H})}^{-1}(0)=\{f\in {\mathcal C}(\widehat{H}):\ f(1_{\widehat{H}})=0\}
\eeq
в ${\mathcal C}(\widehat{H})$ (здесь $1_{\widehat{H}}$ --- единица группы $\widehat{H}$).

}\eit

\btm\label{TH:charakterizatsija-M_sigma}
Множество $M_\sigma\subseteq \widehat{Z}$ представляет собой в точности множество нулей в пространстве $\widehat{Z}$ замкнутого двустороннего идеала $I_\sigma$ алгебры  ${\mathcal C}\big(M_\sigma,{\mathcal B}(X_\sigma)\big)$:
\beq\label{DEF:M_sigma}
M_\sigma=\{\chi\in \widehat{Z}: \ \forall f\in{\mathcal C}(\widehat{H})\quad f(1_{\widehat{H}})=0\ \Rightarrow \ \xi_\sigma(f)(\chi)=0   \}.
\eeq
\etm

\bpr
1. Пусть $\chi\in M_\sigma$, то есть $\chi\in \widehat{Z}$ --- характер, для которого выполняется тождество \eqref{charakterizatsija-M_sigma}. Заметим, что тогда (для этого фиксированного $\chi\in \widehat{Z}$) будет справедливо тождество 
\beq\label{PROOF:charakterizatsija-M_sigma-1}
\xi_\sigma(f)(\chi)=f(1_{\widehat{H}}) \cdot  1_{{\mathcal B}(X_\sigma)},\qquad f\in {\mathcal C}(\widehat{H}).
\eeq
Действительно, если функцию $f\in {\mathcal C}(\widehat{H})$ приблизить линейными комбинациями характеров \eqref{f=sum_i-lambda_i-cdot-h_i}, то мы получим:  
\begin{multline*}
\xi_\sigma(f)(\chi)\longleftarrow\sum_{i}\lambda_i\cdot\xi_\sigma(h_i)(\chi)=\eqref{xi_sigma(h)(chi)}=
\sum_{i}\lambda_i\cdot h_i(\chi\big|_H)\cdot \sigma(h_i)^{-1}=\eqref{charakterizatsija-M_sigma}=\\=
\sum_{i}\lambda_i\cdot \underbrace{h_i(\chi\big|_H) \cdot  h_i(\chi\big|_H)^{-1}}_{\scriptsize \begin{matrix}\|\\
1\\ \|\\  h_i(1_{\widehat{H}}) \end{matrix}}
\cdot 1_{{\mathcal B}(X_\sigma)}
= \sum_{i}\lambda_i\cdot h_i(1_{\widehat{H}}) \cdot 1_{{\mathcal B}(X_\sigma)}
 \longrightarrow
f(1_{\widehat{H}}) \cdot  1_{{\mathcal B}(X_\sigma)}
\end{multline*}
 Теперь если функция $f\in {\mathcal C}(\widehat{H})$ обнуляется в точке $1_{\widehat{H}}\in\widehat{H}$,
\beq\label{PROOF:charakterizatsija-M_sigma-2}
 f(1_{\widehat{H}})=0,
\eeq
то мы получаем 
$$
\xi_\sigma(f)(\chi)=\eqref{PROOF:charakterizatsija-M_sigma-1}=f(1_{\widehat{H}}) \cdot  1_{{\mathcal B}(X_\sigma)}=\eqref{PROOF:charakterizatsija-M_sigma-2}=0
$$
Мы доказали включение
$$
M_\sigma\subseteq \{\chi\in \widehat{Z}: \ \forall f\in{\mathcal C}(\widehat{H})\quad f(1_{\widehat{H}})=0\ \Rightarrow \ \xi_\sigma(f)(\chi)=0   \}.
$$

2. Наоборот, пусть $\chi\in \widehat{Z}$ --- характер со свойством
\beq\label{PROOF:charakterizatsija-M_sigma-3}
\forall f\in{\mathcal C}(\widehat{H})\quad f(1_{\widehat{H}})=0\ \Rightarrow \ \xi_\sigma(f)(\chi)=0. 
\eeq
Зафиксируем произвольный элемент $h\in H$ и рассмотрим функцию $f\in{\mathcal C}(\widehat{H})$, заданную формулой
$$
f(\psi)=\psi(h)-1,\qquad \psi\in\widehat{H}.
$$
Эта функция в точке $1_{\widehat{H}}$ будет равна нулю (потому что $1_{\widehat{H}}$ --- характер на $H$, тождественно равный единице, и значит на элементе $h$ он тоже будет равен единице):
\beq\label{f(1|_hat-H)=0}
f(1_{\widehat{H}})=1_{\widehat{H}}(h)-1=1-1=0
\eeq
Поэтому, в силу \eqref{PROOF:charakterizatsija-M_sigma-3}, в точке $\chi$ должна обнуляться функция $\xi_\sigma(f)$:
\begin{multline*}
0=\xi_\sigma(f)(\chi)=\xi_\sigma(h)(\chi)-\xi_\sigma(1)(\chi)=\eqref{xi_sigma(h)(chi)}=
h(\chi\big|_H) \cdot \sigma(h^{-1})-1(\chi\big|_H) \cdot \sigma(1)=\\= 
h(\chi\big|_H) \cdot \sigma(h)^{-1}-1_{{\mathcal B}(X_\sigma)}
\end{multline*}
Отсюда мы получаем
$$
\sigma(h)=h(\chi\big|_H)\cdot 1_{{\mathcal B}(X_\sigma)}.
$$
И это верно для всякого $h\in H$. Мы доказали \eqref{charakterizatsija-M_sigma}.
То есть $\chi\in M_\sigma$.  И это доказывает включение
$$
M_\sigma\supseteq \{\chi\in \widehat{Z}: \ \forall f\in{\mathcal C}(\widehat{H})\quad f(1_{\widehat{H}})=0\ \Rightarrow \ \xi_\sigma(f)(\chi)=0   \}.
$$

\epr

\medskip
\centerline{\bf Свойства множеств $M_\sigma$:}

\bit{\it

\item[$1^\circ$.]\label{2^0:stroenie-M_sigma} Если представление $\sigma:K\to{\mathcal B}(X_\sigma)$ постоянно на подгруппе $H\subseteq K$
\beq\label{sigma|_H=1*}
\sigma\big|_H=1_{{\mathcal B}(X_\sigma)}
\eeq
то множество $M_\sigma$ совпадает с аннулятором $H^\bot$ подгруппы $H$ в группе $Z$:
\beq\label{M_sigma=H^bot}
M_\sigma=H^\bot=\{\chi\in\widehat{Z}:\quad \chi\big|_H=1\}
\eeq

\item[$2^\circ$.]\label{3^0:stroenie-M_sigma} Если представление $\sigma:K\to{\mathcal B}(X_\sigma)$ непостоянно на подгруппе $H\subseteq K$
\beq\label{sigma|_H-ne-1*}
\sigma\big|_H\ne 1_{{\mathcal B}(X_\sigma)}
\eeq
то множество $M_\sigma$ не совпадает с $H^\bot$:
\beq\label{M_sigma-ne-H^bot}
M_\sigma\ne H^\bot
\eeq
}\eit

\bpr
1. Пусть выполняется \eqref{sigma|_H=1*}. 
Тогда при заданном $\chi\in\widehat{Z}$ тождество \eqref{charakterizatsija-M_sigma}
$$
\sigma(h)=h(\chi\big|_H)\cdot 1_{{\mathcal B}(X_\sigma)}, \qquad h\in H
$$
эквивалентно тождеству  
$$
h(\chi\big|_H)=1,\qquad h\in H.
$$
И мы получаем 
$$
M_\sigma=\eqref{1^0:stroenie-M_sigma}=
\{\chi\in\widehat{Z}:\quad \forall h\in H\quad \sigma(h)=h(\chi\big|_H)\cdot 1_{{\mathcal B}(X_\sigma)} \}=
\{\chi\in\widehat{Z}:\quad \forall h\in H\quad h(\chi\big|_H)=1 \}=H^\bot.
$$

2. Пусть наоборот, выполняется \eqref{sigma|_H-ne-1*}. Тогда найдется $h\in H$ такой что
\beq\label{sigma(h)-ne-1_B(X_sigma)}
\sigma(h)\ne 1_{{\mathcal B}(X_\sigma)}
\eeq
Рассмотрим функцию
$$
f=h-1_H
$$
Она будет принадлежать идеалу \eqref{e_(C(hat(H)))^(-1)(0)}, потому что
$$
f(1_{\widehat{H}})=h(1_{\widehat{H}})-1_H(1_{\widehat{H}})=1-1=0
$$
Но при этом,
\begin{multline*}
\xi_\sigma(f)(1_{\widehat{Z}})=\xi_\sigma(h)(1_{\widehat{Z}})-\xi_\sigma(1_H)(1_{\widehat{Z}})=\eqref{xi_sigma(h)(chi)}=
h(1_{\widehat{Z}}\big|_H)\cdot\overline{\sigma(h)}-1_H(1_{\widehat{Z}}\big|_H)\cdot\overline{\sigma(1_H)}=\\=
h(1_{\widehat{H}})\cdot\overline{\sigma(h)}-1_H(1_{\widehat{H}})\cdot\overline{\sigma(1_H)}=
\overline{\sigma(h)}-\overline{1_{{\mathcal B}(X_\sigma)}}\ne\eqref{sigma(h)-ne-1_B(X_sigma)}\ne 0
\end{multline*}
То есть, элемент $\xi_\sigma(f)\in I_\sigma$ идеала $I_\sigma$ не обращается в нуль в точке $1_{\widehat{Z}}\in \widehat{Z}$. Это значит, что точка $1_{\widehat{Z}}\in \widehat{Z}$ не может лежать в множестве $M_\sigma$, которое по теореме \ref{TH:charakterizatsija-M_sigma} является множеством общих нулей идеала $I_\sigma$:
$$
1_{\widehat{Z}}\notin M_\sigma.
$$
Но, с другой стороны, она лежит в подгруппе $H^\bot$ (потому что это единица группы $\widehat{Z}$)
$$
1_{\widehat{Z}}\in H^\bot.
$$
Поэтому 
$$
M_\sigma\ne H^\bot.
$$

\epr

\blm\label{LM:stroenie-Coker(xi)}
Коядро отображения $\xi$ в \eqref{DEF:xi:Env-C*(H)->Env-C*(Z-times-K)} в категории $\AugInvSteAlg$ аугментированных инволютивных стереотипных алгебр имеет вид
\beq\label{stroenie-Coker(xi)}
\Coker(\xi)\cong
\prod_{\sigma\in\widehat{K}}{\mathcal C}\big(M_\sigma,{\mathcal B}(X_\sigma)\big)
\eeq
где $\{M_\sigma,\sigma\in\widehat{K}\}$  --- множества из \eqref{DEF:M_sigma}.
\elm
\bpr 
По теореме \ref{PROP:AugSteAlg-imeet-ker}, коядро $\Coker(\xi)$ есть непосредственная фактор-алгебра алгебры $\prod_{\sigma\in\widehat{K}}{\mathcal C}(\widehat{Z},{\mathcal B}(X_\sigma))$ 
\beq\label{coker-xi-prod-sigma}
\Coker\ph=\Big(\prod_{\sigma\in\widehat{K}}{\mathcal C}(\widehat{Z},{\mathcal B}(X_\sigma))/I\Big)^\triangledown
\eeq
по замкнутому двустороннему идеалу $I$ в $\prod_{\sigma\in\widehat{K}}{\mathcal C}(\widehat{Z},{\mathcal B}(X_\sigma))$, порожденному образом $\xi(\e_{{\mathcal C}(\widehat{H})}^{-1}(0))$ идеала \eqref{e_(C(hat(H)))^(-1)(0)}.

По лемме \ref{LM:prod-A/prod-I} такая фактор-алгебра изоморфна призведению фактор-алгебр алгебр ${\mathcal C}\big(\widehat{Z},{\mathcal B}(X_\sigma)\big)$ по идеалам $I_\sigma$, которые являются проекциями идеала $I$ на алгебры ${\mathcal C}(\widehat{Z},{\mathcal B}(X_\sigma))$. А всякая такая проекция --- это замкнутый двусторонний идеал в ${\mathcal C}(\widehat{Z},{\mathcal B}(X_\sigma))$, порожденный образом $\xi_\sigma(\e_{{\mathcal C}(\widehat{H})}^{-1}(0))$ идеала \eqref{e_(C(hat(H)))^(-1)(0)}, то есть идеал $I_\sigma$, опреленный на с.\pageref{e_(C(hat(H)))^(-1)(0)}.

Наконец, по лемме \ref{LM:idealy-v-C(M,B(X))} такие фактор-алгебры изоморфны алгебрам ${\mathcal C}\big(M_\sigma,{\mathcal B}(X_\sigma)\big)$, где каждое $M_\sigma$ --- множество общих нулей идеала $I_\sigma$, то есть в точности множество 
\eqref{DEF:M_sigma}:
\begin{multline*}
\Coker(\xi)\cong \bigg(\Big(\prod_{\sigma\in\widehat{K}}{\mathcal C}\big(\widehat{Z},{\mathcal B}(X_\sigma)\big)\Big)/I\bigg)^\vartriangle\cong \bigg(\prod_{\sigma\in\widehat{K}}\Big({\mathcal C}\big(\widehat{Z},{\mathcal B}(X_\sigma)\big)/I_\sigma\Big)\bigg)^\vartriangle\cong\\ \cong
\bigg(\prod_{\sigma\in\widehat{K}}{\mathcal C}\big(M_\sigma,{\mathcal B}(X_\sigma)\big)\bigg)^\vartriangle
\cong
\prod_{\sigma\in\widehat{K}}{\mathcal C}\big(M_\sigma,{\mathcal B}(X_\sigma)\big).
\end{multline*} 
\epr

\bcor\label{LM:Coker(xi)-nepr-alg}
Коядро отображения $\xi$ в \eqref{DEF:xi:Env-C*(H)->Env-C*(Z-times-K)} в категории $\AugInvSteAlg$ аугментированных инволютивных стереотипных алгебр является непрерывной алгеброй.
\ecor
\bpr
В силу примера \ref{EX:C(M,B(X))-nepr-alg}, каждая алгебра ${\mathcal C}\big(M_\sigma,{\mathcal B}(X_\sigma)\big)$ непрерывна.
Поэтому по теореме \ref{TH:proizvedenie-neprer-algebr} их произведение 
$$
\prod_{\sigma\in\widehat{K}}{\mathcal C}\big(M_\sigma,{\mathcal B}(X_\sigma)\big)\cong\Coker\xi
$$
тоже должно быть непрерывно.
\epr

Давайте теперь в диаграмме \eqref{DEF:xi:Env-C*(H)->Env-C*(Z-times-K)} уберем нижнюю строчку, но дополним ее еще одним столбцом.
Пусть, по-прежнему, $Z\cdot K$ -- надстройка абелевой группы $Z$ с помощью компактной группы $K$.
Рассмотрим цепочку гомоморфизмов групп
\beq\label{H->Z-times-K->Z-cdot-K}
 \xymatrix  
{
H\ar[r]^{\iota} & Z\times K\ar[r]^{\varkappa} & Z\cdot K
}
\eeq
в которой $H$, $\iota$, $\varkappa$ определены в лемме \ref{LM:Z-cdot-K=(Z-times-K)/C}, причем
$$
\varkappa\cong\coker\iota.
$$
Отсюда по теореме \ref{TH:coker-groups=>coker-group-algebras} в соответствующей цепочке морфизмов групповых алгебр  в категории $\AugInvSteAlg$ аугментированных инволютивных стереотипных алгебр
\beq\label{C*(H)->C*(Z-times-K)->C*(Z-cdot-K)}
 \xymatrix  @R=2.pc @C=4.pc
{
{\mathcal C}^\star(H)\ar[r]^{{\mathcal C}^\star(\iota)} & {\mathcal C}^\star(Z\times K)\ar[r]^{{\mathcal C}^\star(\varkappa)} & {\mathcal C}^\star(Z\cdot K)
}
\eeq
второй морфизм является коядром первого:
$$
{\mathcal C}^\star(\varkappa)=\coker{\mathcal C}^\star(\iota).
$$
Иными словами,
\beq\label{C^star(Z-cdot-K)=Coker-C^star(iota)}
{\mathcal C}^\star(Z\cdot K)=\Coker{\mathcal C}^\star(\iota).
\eeq

Рассмотрим далее цепочку, полученную из \eqref{C*(H)->C*(Z-times-K)->C*(Z-cdot-K)} применением к ней функтора непрерывной оболочки:
\beq\label{Env-C*(H)->Env-C*(Z-times-K)->Env-C*(Z-cdot-K)}
 \xymatrix  @R=2.pc @C=4.pc
{
\Env_{\mathcal C} {\mathcal C}^\star(H)\ar[r]^{\Env_{\mathcal C} {\mathcal C}^\star(\iota)} & \Env_{\mathcal C} {\mathcal C}^\star(Z\times K)\ar[r]^{\Env_{\mathcal C} {\mathcal C}^\star(\varkappa)} & \Env_{\mathcal C} {\mathcal C}^\star(Z\cdot K)
}
\eeq
По теореме \ref{TH:Env_C-functor-v-AugSteAlg} это будет цепочка морфизмов аугментированных стереотипных алгебр.
Теперь диаграмма \eqref{DEF:xi:Env-C*(H)->Env-C*(Z-times-K)} превращается в диаграмму
\beq\label{ph:Coker(Env_C(iota^*))->B-1}
 \xymatrix  @R=2.pc @C=3.pc
{
H\ar[d]_{\delta_H}\ar[r]^{\iota} & Z\times K\ar[d]_{\delta_{Z\times K}}\ar[rr]^{\varkappa} & &  Z\cdot K\ar[d]_{\delta_{Z\cdot K}} \\
{\mathcal C}^\star(H)\ar[ddd]_{\env_{\mathcal C}{\mathcal C}^\star(H)}\ar[r]^{{\mathcal C}^\star(\iota)} & {\mathcal C}^\star(Z\times K) \ar[rr]^{{\mathcal C}^\star(\varkappa)} \ar[ddd]_{\env_{\mathcal C}{\mathcal C}^\star(Z\times K)}\ar[rd]_{\coker{\mathcal C}^\star(\iota)\quad} &  & {\mathcal C}^\star(Z\cdot K)
\ar[ddd]_{\env_{\mathcal C} {\mathcal C}^\star(Z\cdot K)} 
\\
& & \Coker{\mathcal C}^\star(\iota)\ar@{=}[ru]_{\quad\eqref{C^star(Z-cdot-K)=Coker-C^star(iota)}} \ar[d]_{\Coker\frac{\Env_{\mathcal C}{\mathcal C}^\star(\iota),\env_{\mathcal C}{\mathcal C}^\star(H)}{\env_{\mathcal C}{\mathcal C}^\star(Z\times K),{\mathcal C}^\star(\iota)}} &
\\
& & \Coker\Env_{\mathcal C}{\mathcal C}^\star(\iota)\ar@{-->}[rd]^{\upsilon}\ar@{=}[dd]|\hole & 
\\
\Env_{\mathcal C}{\mathcal C}^\star(H)\ar[r]_{\Env_{\mathcal C}{\mathcal C}^\star(\iota)}\ar@{=}[d] & \Env_{\mathcal C}{\mathcal C}^\star(Z\times K) \ar[rr]_(.7){\Env_{\mathcal C} {\mathcal C}^\star(\varkappa)} \ar[ru]^{\coker\Env_{\mathcal C}{\mathcal C}^\star(\iota)\quad}\ar@{=}[d] &  &
\Env_{\mathcal C} {\mathcal C}^\star(Z\cdot K)\\
{\mathcal C}(\widehat{H})\ar[r]_{\xi} & \prod_{\sigma\in\widehat{K}}{\mathcal C}\big(\widehat{Z},{\mathcal B}(X_\sigma)\big)\ar[r]_{\coker\xi}
& \prod_{\sigma\in\widehat{K}}{\mathcal C}\big(M_\sigma,{\mathcal B}(X_\sigma)\big) &
}
\eeq
в которой $\upsilon$ --- естественный морфизм, связывающий $\coker\Env_{\mathcal C}{\mathcal C}^\star(\iota)$ и 
$\Env_{\mathcal C} {\mathcal C}^\star(\varkappa)$.

\blm\label{LM:Coker(Env_C(iota^*))->Env_C-C^*(Z-cdot-K)}
Морфизм
$$
\upsilon:\Coker\Env_{\mathcal C} {\mathcal C}^\star(\iota)\to \Env_{\mathcal C} {\mathcal C}^\star(Z\cdot K)
$$
в диаграмме \eqref{ph:Coker(Env_C(iota^*))->B-1} является изоморфизмом стереотипных алгебр.
\elm
\bpr
По следствию \ref{LM:Coker(xi)-nepr-alg}, коядро  $\Coker\Env_{\mathcal C} {\mathcal C}^\star(\iota)$ является непрерывной алгеброй (то есть полным объектом относительно непрерывной оболочки). Поэтому по теореме \ref{TH:Coker-rasshirenie}, морфизм $\upsilon$ должен быть изоморфизмом.  
\epr

\btm\label{TH:Env_C-C*(Z-cdot-K)}
Пусть $Z\cdot K$ --- надстройка абелевой локально компактной группы $Z$ с помощью компактной группы $K$. Тогда
\bit{

\item[(i)] непрерывная оболочка $\Env_{\mathcal C} {\mathcal C}^\star(Z\cdot K)$ ее групповой алгебры имеет вид
\beq\label{Env_C-C*(Z-cdot-K)}
\Env_{\mathcal C} {\mathcal C}^\star(Z\cdot K)\cong
\prod_{\sigma\in\widehat{K}}{\mathcal C}\big(M_\sigma,{\mathcal B}(X_\sigma)\big),
\eeq
где $\{M_\sigma;\ \sigma\in\widehat{K}\}$ --- семейство замкнутых подмножеств в двойственной группе $\widehat{Z}$ к группе $Z$, определенное формулой \eqref{DEF:M_sigma},

\item[(ii)] пространство $\Env_{\mathcal C} {\mathcal C}^\star(Z\cdot K)$ насыщено (как стереотипное пространство).
}\eit
\etm
\bpr
Первая часть этого утверждения следует из лемм \ref{LM:Coker(Env_C(iota^*))->Env_C-C^*(Z-cdot-K)} и \ref{LM:stroenie-Coker(xi)}, а вторая --- из того факта, что каждое пространство ${\mathcal C}\big(M_\sigma,{\mathcal B}(X_\sigma)\big)$, будучи прямым произведением пространств Фреше, насыщено.
\epr

\bcor\label{COR:Env_C-Z.K=LCS-lim} Пусть $G=Z\cdot K$ --- надстройка абелевой локально компактной группы $Z$ с помощью компактной группы $K$. Тогда непрерывная оболочка групповой алгебры ${\mathcal C}^\star(G)$ совпадает с локально выпуклой оболочкой Кузнецовой, то есть с проективным пределом ее $C^*$-фактор-алгебр в категории локально выпуклых пространств (и в категории топологических алгебр):
\beq\label{Env_C-Z.K=LCS-lim}
\Env_{\mathcal C}{\mathcal C}^\star(G)={\tt LCS}\text{-}\kern-3pt\projlim_{p\in\Sn({\mathcal C}^\star(G))}{\mathcal C}^\star(G)/p
\eeq
\ecor
\bpr
Пусть для всякого $\tau\in\widehat{K}$
$$
\pi_\tau:\prod_{\sigma\in\widehat{K}}{\mathcal C}\big(M_\sigma,{\mathcal B}(X_\sigma)\big)\to {\mathcal C}\Big(M_\tau,{\mathcal B}(X_\tau)\Big),
$$
--- естественная проекция, а $\norm{\cdot}_\tau$ --- норма в алгебре ${\mathcal B}(X_\tau)$. Всякому конечному множеству $S\subseteq \widehat{K}$ и любому семейству компактов
$$
T_\tau\subseteq M_\tau,\qquad \tau\in S
$$
поставим в соответствие полунорму 
$$
p_{\{T_\tau\} }(\alpha)=\sum_{\tau\in S}\sup_{t\in T_\tau}\norm{\pi_\tau(\env_{\mathcal C}(\alpha))(t)}_\tau,\qquad\alpha\in {\mathcal C}^\star(G).
$$
Это будет $C^*$-полунорма, и более того, такие полунормы образуют конфинальную подсистему в множестве ${\Sn}({\mathcal C}^\star(G))$ всех $C^*$-полунорм на алгебре ${\mathcal C}^\star(G)$. По этой причине, проективный предел системы фактор-алгебр по этим полунормам совпадает с проективным пределом системы фактор-алгебр по всем $C^*$-полунормам:
\beq\label{Env_C-Z.K=LCS-lim-1}
{\tt LCS}\text{-}\kern-3pt\projlim_{p_{\{T_\tau\} }}{\mathcal C}^\star(G)/p_{\{T_\tau\}}=
{\tt LCS}\text{-}\kern-3pt\projlim_{p\in\Sn({\mathcal C}^\star(G))}{\mathcal C}^\star(G)/p
\eeq
Но с другой стороны, из самого выбора этих полунорм видно, что этот проективный предел совпадает с алгеброй 
$\prod_{\sigma\in\widehat{K}}{\mathcal C}\big(M_\sigma,{\mathcal B}(X_\sigma)\big)$:
\beq\label{Env_C-Z.K=LCS-lim-2}
{\tt LCS}\text{-}\kern-3pt\projlim_{p_{\{T_\tau\} }}{\mathcal C}^\star(G)/p_{\{T_\tau\}}=
\prod_{\sigma\in\widehat{K}}{\mathcal C}\big(M_\sigma,{\mathcal B}(X_\sigma)\big)\cong \Env_{\mathcal C} {\mathcal C}^\star(G)
\eeq
Равенства \eqref{Env_C-Z.K=LCS-lim-1} и \eqref{Env_C-Z.K=LCS-lim-2} вместе дают \eqref{Env_C-Z.K=LCS-lim}.
\epr

\paragraph{Группы Ли---Мура.}

Пусть $G$ -- группа Ли---Мура. Тогда по теореме \ref{TH:Lie-Moore} $G$ --- конечное расширение некоторой компактной надстройки абелевой группы (Ли):
\beq\label{G-Lie-Moore}
1\to Z\cdot K=N\to G\to F\to 1
\eeq
($Z$ -- абелева группа Ли, $K$ -- компактная группа Ли, $F$ -- конечная группа). Эта цепочка порождает цепочку гомоморфизмов групповых алгебр
$$
\C\to {\mathcal C}^\star(Z\cdot K)={\mathcal C}^\star(N)\to {\mathcal C}^\star(G)\to {\mathcal C}^\star(F)=\C_F\to \C.
$$
(второе равенство в этой цепочке следует из примера \ref{EX:Env_C-C_F=C_F}) и цепочку гомоморфизмов непрерывных оболочек
\beq\label{Env_C-C*(G)-Lie-Moore-0}
\C\to \Env_{\mathcal C} {\mathcal C}^\star(Z\cdot K)=\Env_{\mathcal C} {\mathcal C}^\star(N)\to \Env_{\mathcal C} {\mathcal C}^\star(G)\to \Env_{\mathcal C} {\mathcal C}^\star(F)=\eqref{Env_C-C_F=C_F}=\C_F\to \C.
\eeq

Напомним, что множество всех $C^*$-полунорм на $A$ мы условились обозначать $\Sn(A)$. В частном случае, когда $A={\mathcal C}^\star(G)$ мы будем употреблять обозначение
$$
\Sn(G)=\Sn\big({\mathcal C}^\star(G)\big).
$$
Рассмотрим ситуацию, описанную в \eqref{G-kak-rasshirenie}: пусть $N$ -- открытая нормальная подгруппа в локально компактной группе $G$, и $F=G/N$ --- фактор-группа. Будем считать $G$ расширением группы $N$ по группе $F$:
\beq\label{G-kak-rasshirenie-1}
\xymatrix 
{
1\ar[r] & N \ar[r]^{\eta} & G\ar[r]^{\ph} & F\ar[r]& 1
}
\eeq
(здесь  $\eta$ -- естественное вложение, а $\ph$ -- фактор отображение).
Пусть далее $\theta:{\mathcal C}^\star(N)\to {\mathcal C}^\star(G)$ обозначает отображение, индуцированное вложением $\eta:N\subseteq G$ в цепочке \eqref{G-kak-rasshirenie-1}, 
\beq\label{DIAGR:theta+Env-theta-0}
 \xymatrix  @R=2.pc @C=6.pc
{
 N\ar[r]^{\eta}\ar[d]_{\delta_N} & G \ar[d]^{\delta_{G}}
 \\
{\mathcal C}^\star(N)\ar[r]^{\theta}
  & {\mathcal C}^\star(G)
  \\
}
\eeq

\bprop\label{PROP:otkrytost-Env-C^*(N)->Env-C^*(G)}
Для всякой полунормы $p\in\Sn(N)$ найдется полунорма $p'\in\Sn(G)$ со следующими свойствами:

\bit{

\item[(i)] $p$ мажорируется ограничением $p'\in\Sn(G)$ на ${\mathcal C}^\star(N)$:
\beq\label{otkrytost-Env-C^*(N)->Env-C^*(G)}
p(\alpha)\le p'(\theta(\alpha)),\qquad\alpha\in {\mathcal C}^\star(N).
\eeq

\item[(ii)] для всякого класса смежности $L\in G/N$ найдется элемент $x_L\in L$ для которого выполняется неравенство:
\beq\label{otkrytost-Env-C^*(N)->Env-C^*(G)-1}
p(x_L^{-1}\cdot\beta)\le p'(\beta),\qquad\beta\in {\mathcal C}^\star(L).
\eeq

}\eit

\beq\label{DIAGR:theta+Env-theta-1}
 \xymatrix  @R=2.pc @C=6.pc
{
 N\ar[rr]^{\eta}\ar[d]_{\delta_N} & & G \ar[d]^{\delta_{G}}
 \\
{\mathcal C}^\star(N)\ar[rr]^{\theta}\ar[dr]_{p}
  & & {\mathcal C}^\star(G)\ar@{-->}[dl]^{p'}
  \\
  & \R_+ & 
}
\eeq
(треугольник внизу не обязан быть коммутативным).

\eprop
\bpr
Нужно $p$ представить как полунорму, порожденную некоторым представлением $\dot{\pi}:{\mathcal C}^\star(N)\to{\mathcal B}(X)$,
$$
p(\alpha)=\norm{\dot{\pi}(\alpha)},\qquad\alpha\in {\mathcal C}^\star(N),
$$
потом рассмотреть индуцированное представление $\dot{\pi}':{\mathcal C}^\star(G)\to{\mathcal B}(L_2(D,X))$ и положить
$$
p'(\beta)=\norm{\dot{\pi}'(\beta)},\qquad\beta\in {\mathcal C}^\star(G).
$$

1. Для доказательства \eqref{otkrytost-Env-C^*(N)->Env-C^*(G)} рассмотрим сначала случай, когда $\alpha$ --- линейная комбинация дельта-функций из $N$:
\beq\label{alpha=sum_i-lambda_i-cdot-delta^(g_i)}
\alpha=\sum_i\lambda_i\cdot\delta^{g_i},\qquad g_i\in N,
\eeq
и любого $\xi\in L_2(F,X)$ мы получим:
\begin{multline*}
\norm{\dot{\pi}'(\alpha)(\xi)(1)}= \norm{\dot{\pi}'\l\sum_i\lambda_i\cdot\delta^{g_i}\r(\xi)(1)}=
\norm{\sum_i\lambda_i\cdot\dot{\pi}'(\delta^{g_i})(\xi)(1)}=
\norm{\sum_i\lambda_i\cdot\pi'(g_i)(\xi)(1)}=\eqref{ind-representation}=\\=
\norm{\sum_i\lambda_i\cdot\pi(\sigma(1)\cdot g_i\cdot \sigma(1)^{-1})(\xi(1))}=\eqref{sigma(1_D)=1_G}=
\norm{\sum_i\lambda_i\cdot\pi(g_i)(\xi(1))}=
\norm{\sum_i\lambda_i\cdot\dot{\pi}(\delta^{g_i})(\xi(1))}=\\=
\norm{\dot{\pi}\l\sum_i\lambda_i\cdot\delta^{g_i}\r(\xi(1))}=
\norm{\dot{\pi}(\alpha)(\xi(1))}
\end{multline*}
$$
\Downarrow
$$
$$
\norm{\dot{\pi}'(\alpha)(\xi)}=\sqrt{\sum_{t\in F}\norm{\dot{\pi}'(\alpha)(\xi)(t)}^2}\ge \norm{\dot{\pi}'(\alpha)(\xi)(1)}=\norm{\dot{\pi}(\alpha)(\xi(1))}
$$
$$
\Downarrow
$$
$$
p'(\alpha)=\norm{\dot{\pi}'(\alpha)}=\sup_{\norm{\xi}\le 1}\norm{\dot{\pi}'(\alpha)(\xi)}\ge \sup_{\norm{\xi}\le 1}\norm{\dot{\pi}(\alpha)(\xi(1))}=
\sup_{\norm{\zeta}\le 1}\norm{\dot{\pi}(\alpha)(\zeta)}=
\norm{\dot{\pi}(\alpha)}=p(\alpha).
$$
Меры вида \eqref{alpha=sum_i-lambda_i-cdot-delta^(g_i)} плотны в ${\mathcal C}^\star(N)$, поэтому неравенство \eqref{otkrytost-Env-C^*(N)->Env-C^*(G)} справедливо для всех $\alpha\in{\mathcal C}^\star(N)$.

2. Теперь докажем \eqref{otkrytost-Env-C^*(N)->Env-C^*(G)-1}. Зафиксируем $L\in G/N=F$ и положим $x_L=\sigma(L)^{-1}$. Тогда
\beq\label{ph(x_L)=L^(-1)}
\ph(x_L)=\ph(\sigma(L)^{-1})=\ph(\sigma(L))^{-1}=L^{-1}
\eeq
Опять для мер $\alpha\in{\mathcal C}^\star(N)$ вида \eqref{alpha=sum_i-lambda_i-cdot-delta^(g_i)} и любого $\xi\in L_2(F,X)$ мы получим:
\begin{multline*}
\norm{\dot{\pi}'(x_L\cdot \alpha)(\xi)(L)}= \norm{\dot{\pi}'\l\sum_i\lambda_i\cdot\delta^{x_L\cdot g_i}\r(\xi)(L)}=
\norm{\sum_i\lambda_i\cdot\dot{\pi}'(\delta^{x_L\cdot g_i})(\xi)(L)}=
\norm{\sum_i\lambda_i\cdot\pi'(x_L\cdot g_i)(\xi)(L)}=\\=\eqref{ind-representation}=
\norm{\sum_i\lambda_i\cdot\pi(\sigma(L)\cdot x_L\cdot g_i\cdot \sigma(L\cdot \ph(x_L\cdot g_i))^{-1})(\xi(L\cdot \ph(x_L\cdot g_i)))}=\\=
\bigg\|\sum_i\lambda_i\cdot
\pi(\sigma(L)\cdot \underbrace{x_L}_{\scriptsize\begin{matrix} \|\\ \sigma(L)^{-1}\end{matrix}}\cdot g_i\cdot \sigma(L\cdot 
\underbrace{\ph(x_L)}_{\scriptsize\begin{matrix}\put(-28,0){\eqref{ph(x_L)=L^(-1)}} \|\\ L^{-1}\end{matrix}}\cdot \underbrace{\ph(g_i)}_{\scriptsize\begin{matrix}\|\\ 1\end{matrix}})^{-1})(\xi(L\cdot \underbrace{\ph(x_L)}_{\scriptsize\begin{matrix}\put(-28,0){\eqref{ph(x_L)=L^(-1)}} \|\\ L^{-1}\end{matrix}}\cdot \underbrace{\ph(g_i)}_{\scriptsize\begin{matrix}\|\\ 1\end{matrix}}))\bigg\|=\\=
\bigg\|\sum_i\lambda_i\cdot
\pi(\sigma(L)\cdot \sigma(L)^{-1}\cdot g_i\cdot \sigma(L\cdot 
L^{-1}\cdot  1)^{-1})(\xi(L\cdot L^{-1}\cdot 1))\bigg\|=\\=
\bigg\|\sum_i\lambda_i\cdot
\pi( g_i\cdot \sigma(1)^{-1})(\xi(1))\bigg\|=
\eqref{sigma(1_D)=1_G}=
\norm{\sum_i\lambda_i\cdot\pi(g_i)(\xi(1))}=
\norm{\sum_i\lambda_i\cdot\dot{\pi}(\delta^{g_i})(\xi(1))}=\\=
\norm{\dot{\pi}\l\sum_i\lambda_i\cdot\delta^{g_i}\r(\xi(1))}=
\norm{\dot{\pi}(\alpha)(\xi(1))}
\end{multline*}
$$
\Downarrow
$$
$$
\norm{\dot{\pi}'(x_L\cdot \alpha)(\xi)}=\sqrt{\sum_{t\in F}\norm{\dot{\pi}'(x_L\cdot \alpha)(\xi)(t)}^2}\ge \norm{\dot{\pi}'(x_L\cdot \alpha)(\xi)(L)}=\norm{\dot{\pi}(\alpha)(\xi(1))}
$$
$$
\Downarrow
$$
$$
p'(x_L\cdot \alpha)=\norm{\dot{\pi}'(x_L\cdot \alpha)}=\sup_{\norm{\xi}\le 1}\norm{\dot{\pi}'(x_L\cdot \alpha)(\xi)}\ge \sup_{\norm{\xi}\le 1}\norm{\dot{\pi}(\alpha)(\xi(1))}=
\sup_{\norm{\zeta}\le 1}\norm{\dot{\pi}(\alpha)(\zeta)}=
\norm{\dot{\pi}(\alpha)}=p(\alpha).
$$
$$
\Downarrow
$$
\beq\label{otkrytost-Env-C^*(N)->Env-C^*(G)-2}
p'(x_L\cdot \alpha)\ge p(\alpha).
\eeq
Меры вида \eqref{alpha=sum_i-lambda_i-cdot-delta^(g_i)} плотны в ${\mathcal C}^\star(N)$, поэтому неравенство \eqref{otkrytost-Env-C^*(N)->Env-C^*(G)-2} справедливо для всех $\alpha\in{\mathcal C}^\star(N)$. Заменой
$$
\beta=x_L\cdot \alpha
$$
оно превращается в неравенство \eqref{otkrytost-Env-C^*(N)->Env-C^*(G)-1}.
\epr

Пусть далее $\Env_{\mathcal C}\theta:\Env_{\mathcal C} {\mathcal C}^\star(N)\to \Env_{\mathcal C} {\mathcal C}^\star(G)$ --- морфизм оболочек, порожденный морфизмом $\theta:{\mathcal C}^\star(N)\to {\mathcal C}^\star(G)$ из \eqref{DIAGR:theta+Env-theta-0}:
\beq\label{DIAGR:theta+Env-theta}
 \xymatrix  @R=2.pc @C=6.pc
{
 N\ar[r]^{\eta}\ar[d]_{\delta_N} & G \ar[d]^{\delta_{G}}
 \\
{\mathcal C}^\star(N)\ar[r]^{\theta}
 \ar[d]_{\env_{\mathcal C}} & {\mathcal C}^\star(G)
 \ar[d]^{\env_{\mathcal C}} \\
\Env_{\mathcal C}{\mathcal C}^\star(N)\ar[r]^{\Env_{\mathcal C}\theta} & \Env_{\mathcal C}{\mathcal C}^\star(G)
}
\eeq
Каждая полунорма $p\in{\Sn}(N)$ однозначно продолжается до некоторой полунормы на $\Env_{\mathcal C}{\mathcal C}^\star(N)$. Мы будем обозначать это продолжение той же буквой $p$. Точно так же, каждая полунорма $q\in{\Sn}(G)$ однозначно продолжается до некоторой полунормы на $\Env_{\mathcal C}{\mathcal C}^\star(G)$, и мы будем обозначать это продолжение той же буквой $q$.
Из предложения \ref{PROP:otkrytost-Env-C^*(N)->Env-C^*(G)} следует

\bprop\label{PROP:otkrytost-Env-C^*(N)->Env-C^*(G)-1}
Для всякой полунормы $p\in\Sn(N)$ найдется полунорма $p'\in\Sn(G)$ такая, что
\beq\label{p(x)-le-q(Env-theta(x))}
p(x)\le p'(\Env_{\mathcal C}\theta(x)),\qquad x\in\Env_{\mathcal C}{\mathcal C}^\star(N).
\eeq
\eprop
\bpr
Пусть $x=\env_{\mathcal C}\alpha$, где $\alpha\in {\mathcal C}^\star(N)$. Тогда для полунормы $p'$ из предложения \ref{PROP:otkrytost-Env-C^*(N)->Env-C^*(G)} мы получим:
$$
p(x)=p(\env_{\mathcal C}\alpha)=p(\alpha)\le p'\big(\theta(\alpha)\big)= p'\Big(\env_{\mathcal C}\big(\theta(\alpha)\big)\Big)= p'(\Env_{\mathcal C}\theta(\env_{\mathcal C}\alpha))=p'(\Env_{\mathcal C}\theta(x))
$$
Элементы вида $x=\env_{\mathcal C}\alpha$, где $\alpha\in {\mathcal C}^\star(N)$, плотны в
$\Env_{\mathcal C}{\mathcal C}^\star(N)$, поэтому неравенство \eqref{p(x)-le-q(Env-theta(x))} справедливо для всех $x\in\Env_{\mathcal C}{\mathcal C}^\star(N)$.
\epr

\bprop\label{PROP:Env_C(theta)} Пусть в цепочке \eqref{G-kak-rasshirenie-1} $N=Z\cdot K$ --- компактная надстройка абелевой группы (и по-прежнему открытая нормальная подгруппа в $G$). Тогда непрерывная оболочка $\Env_{\mathcal C}\theta:\Env_{\mathcal C} {\mathcal C}^\star(N)\to \Env_{\mathcal C} {\mathcal C}^\star(G)$ морфизма $\theta:{\mathcal C}^\star(N)\to {\mathcal C}^\star(G)$ является инъективным и открытым\footnote{См. определение на с.\pageref{DEF:open-map}.} отображением стереотипных пространств.
\eprop
\bpr
Здесь используется следствие \ref{COR:Env_C-Z.K=LCS-lim}, согласно которому $\Env_{\mathcal C} {\mathcal C}^\star(N)$ представляет собой локально выпуклый проективный предел фактор-пространств ${\mathcal C}^\star(N)/p$, где  $p\in\Sn(N)$:
\beq\label{E(C*(G))=LCS-leftlim}
\Env_{\mathcal C}{\mathcal C}^\star(N)={\tt LCS}\text{-}\kern-3pt\projlim_{p\in\Sn({\mathcal C}^\star(N))}{\mathcal C}^\star(N)/p
\eeq

1. Из \eqref{E(C*(G))=LCS-leftlim} следует, что полунормы $p\in\Sn(N)$ отделяют элементы $\Env_{\mathcal C} {\mathcal C}^\star(N)$. То есть если $0\ne x\in \Env_{\mathcal C} {\mathcal C}^\star(N)$, то найдется полунорма $p\in\Sn(N)$ такая, что $p(x)>0$. По предложению \ref{PROP:otkrytost-Env-C^*(N)->Env-C^*(G)} мы можем подобрать полунорму $p'\in\Sn(G)$ так, чтобы выполнялось \eqref{otkrytost-Env-C^*(N)->Env-C^*(G)}, и тогда мы в силу \eqref{p(x)-le-q(Env-theta(x))} получим
$$
0<p(x)\le p'(\Env_{\mathcal C}\theta(x)),
$$
и значит, $\Env_{\mathcal C}\theta(x)\ne 0$. Это доказывает инъективность $\Env_{\mathcal C}\theta$.

2. Из \eqref{E(C*(G))=LCS-leftlim} следует, кроме того, что топология $\Env_{\mathcal C} {\mathcal C}^\star(N)$ порождается полунормами $p\in\Sn(N)$ (без псевдонасыщения). Поэтому базисные окрестности нуля в $\Env_{\mathcal C} {\mathcal C}^\star(N)$ можно считать порожденными полунормами $p\in\Sn(N)$. Для всякой такой окрестности нуля
$$
U=\{x\in \Env_{\mathcal C} {\mathcal C}^\star(N):\ p(x)\le\e \}
$$
мы по предложению \ref{PROP:otkrytost-Env-C^*(N)->Env-C^*(G)} подбираем полунорму $p'\in\Sn(G)$ так, чтобы выполнялось \eqref{p(x)-le-q(Env-theta(x))}, и тогда для окрестности нуля
$$
V=\{y\in \Env_{\mathcal C}{\mathcal C}^\star(G):\ p'(y)\le\e \}
$$
мы получим
\begin{multline*}
y\in \Env_{\mathcal C}\theta\big(\Env_{\mathcal C}{\mathcal C}^\star(N)\big)\cap V\quad\Longrightarrow\quad
\exists x\in \Env_{\mathcal C}{\mathcal C}^\star(N)\quad y=\Env_{\mathcal C}\theta(x)\in V\quad\Longrightarrow \\ \Longrightarrow\quad p(x)\le p'\big(\Env_{\mathcal C}\theta(x)\big)\le\e\quad\Longrightarrow\quad x\in U \quad\Longrightarrow\quad y=\Env_{\mathcal C}\theta(x)\in\Env_{\mathcal C}\theta(U).
\end{multline*}
Это доказывает открытость $\Env_{\mathcal C}\theta$.
\epr

\blm\label{LM:1_L-in-K(G)}
Пусть $N$ --- открытая нормальная подгруппа в локально компактной группе $G$. Тогда для всякого класса смежности $L\in G/N$ его характеристическая функция 
\beq\label{DEF:1_L}
1_L(x)=\begin{cases}1, & x\in L\\ 0, & x\notin L \end{cases}
\eeq 
является элементом пространства ${\mathcal K}(G)$:
\beq\label{1_L-in-K(G)}
1_L\in {\mathcal K}(G)
\eeq
\elm
\bpr
Рассмотрим тривиальное представление $\pi:N\to\C$, $\pi(t)=1$. Его индуцированное представление
$\pi':G\to {\mathcal L}(L_2(D))$ определяется формулой \eqref{ind-representation}, которая в данном случае принимает вид
$$
\pi'(g)(\xi)(t)=\xi\big(\ph(\sigma(t)\cdot g)\big)=
\xi\big(t\cdot\ph(g)\big),\qquad \xi\in L_2(D),\quad t\in D,\quad g\in G.
$$
Рассмотрим в качестве $\xi$ характеристическую функцию единицы $1_D$ группы $D$:
$$
\xi(t)=\begin{cases}1, & t=1_D\\ 0, & t\ne 1_D\end{cases}
$$
Тогда
$$
f(g)=\langle\pi'(g)(\xi),\xi\rangle=\sum_{t\in D}\pi'(g)(\xi)(t)\cdot\overline{\xi(t)}=
\sum_{t\in D}\xi\big(t\cdot\ph(g)\big)\cdot \xi(t)=\begin{Bmatrix}1, & \ph(g)=1\\ 0, & \ph(g)\ne 1\end{Bmatrix}=
\begin{Bmatrix}1, & g\in N\\ 0, & g\notin N\end{Bmatrix}.
$$
То есть характеристическая функция $1_N$ лежит в $k(G)$, и значит, по теореме \ref{TH:V-k-K-C}, в  ${\mathcal K}(G)$. Из теоремы \ref{LM:sdvig-v-K(G)} следует, что все ее сдвиги тоже лежат в ${\mathcal K}(G)$.
\epr

Для всякого класса смежности $L\in G/N$ обозначим
\beq\label{DEF:K_L(G)}
{\mathcal K}_L(G)=1_L\cdot{\mathcal K}(G),\qquad {\mathcal K}_{G\setminus L}(G)=(1-1_L)\cdot{\mathcal K}(G)
\eeq
(здесь 1 -- единица алгебры ${\mathcal K}(G)$). Мы наделяем эти пространства структурой непосредственных подпространств в ${\mathcal K}(G)$ (см.\cite{Ak16}).

\blm\label{LM:K_L(G)-K_G-L(G)} Пусть $N$ --- открытая нормальная подгруппа в локально компактной группе $G$. Тогда 
\bit{
\item[(i)] для классов смежности $L\in G/N$ пространства ${\mathcal K}_L(G)$ и ${\mathcal K}_{G\setminus L}(G)$ дополняют друг друга в ${\mathcal K}(G)$:
\beq\label{K_N+K_(G-N)=K(G)}
{\mathcal K}_L(G)\oplus {\mathcal K}_{G\setminus L}(G)={\mathcal K}(G)
\eeq
(то есть ${\mathcal K}(G)$ есть прямая сумма в категории стереотипных пространств).

\item[(ii)] пространства ${\mathcal K}_L(G)$ изоморфны друг другу и пространству ${\mathcal K}_N(G)$:
\beq\label{K_L(G)=K_N(G)}
{\mathcal K}_L(G)\cong {\mathcal K}_N(G),\qquad L\in G/N.
\eeq

}\eit
\elm
\bpr
\eqref{K_N+K_(G-N)=K(G)} следует из \eqref{1_L-in-K(G)}, а \eqref{K_L(G)=K_N(G)} --- из теоремы \ref{LM:sdvig-v-K(G)},  по которой сдвиг является изоморфизмом пространства ${\mathcal K}(G)$, и поэтому подходящим сдвигом можно изоморфно отобразить ${\mathcal K}_L(G)$ в ${\mathcal K}_N(G)$.
\epr

Рассмотрим в диаграмме \eqref{DIAGR:theta+Env-theta} морфизм $\Env_{\mathcal C}\theta:\Env_{\mathcal C} ({\mathcal C}^\star(N))\to \Env_{\mathcal C} {\mathcal C}^\star(G)$  и обозначим его сопряженное отображение буквой $\psi=(\Env_{\mathcal C}\theta)^\star:{\mathcal K}(N)\gets {\mathcal K}(G)$.

\blm\label{LM:K(N)=K_N(G)} Пусть в цепочке \eqref{G-kak-rasshirenie} $N=Z\cdot K$ --- компактная надстройка абелевой группы (и по-прежнему открытая нормальная подгруппа в локально компактной группе $G$).  Тогда морфизм стереотипных пространств $\psi=(\Env_{\mathcal C}\theta)^\star:{\mathcal K}(N)\gets {\mathcal K}(G)$ обладает следующими свойствами:
\bit{
\item[(i)] его ядром является вторая компонента в разложении \eqref{K_N+K_(G-N)=K(G)} (с $L=N$):
\beq\label{ph(K_(G-N)(G))=0}
\Ker\psi={\mathcal K}_{G\setminus N}(G).
\eeq
\item[(ii)] ограничение $\psi|_{{\mathcal K}_N(G)}:{\mathcal K}_N(G)\to {\mathcal K}(N)$ является изоморфизмом стереотипных алгебр.
\beq\label{K_N(G)=K(N)}
{\mathcal K}_N(G)\cong {\mathcal K}(N).
\eeq

}\eit
\elm
\bpr
1. Для доказательства (i) рассмотрим диаграмму
$$
\xymatrix @R=2.pc @C=8.0pc 
{
{\mathcal C}^\star(N)\ar[r]^{\env_{\mathcal C} {{\mathcal C}^\star(N)}}\ar[d]_{\theta} & \Env_{\mathcal C} ({\mathcal C}^\star(N))\ar[d]^{\Env_{\mathcal C}\theta}\\
{\mathcal C}^\star(G)\ar[r]^{\env_{\mathcal C} {{\mathcal C}^\star(G)}} & \Env_{\mathcal C} {\mathcal C}^\star(G)
}
$$
Если $u\in \Ker\psi$, то мы получаем цепочку
\begin{multline*}
0=\psi(u)=u\circ \Env_{\mathcal C}\theta\quad\Longrightarrow\quad 0=u\circ \Env_{\mathcal C}\theta\circ \env_{\mathcal C} {{\mathcal C}^\star(N)}=
u\circ \env_{\mathcal C} {{\mathcal C}^\star(G)}\circ \theta \quad\Longrightarrow\\ \Longrightarrow\quad 0=u\circ \env_{\mathcal C} {{\mathcal C}^\star(G)}\circ \theta\circ\delta^N \quad\Longrightarrow\quad 0=u\Big|_N \quad\Longrightarrow\\ \Longrightarrow\quad 0=u\cdot 1_N
\quad\Longrightarrow\quad u=u\cdot (1-1_N) \quad\Longrightarrow\quad u\in {\mathcal K}_{G\setminus N}(G).
\end{multline*}
Наоборот, если $u\in {\mathcal K}_{G\setminus N}(G)$, то возникает обратная цепочка
\begin{multline*}
0=\psi(u)=u\circ \Env_{\mathcal C}\theta\quad\overset{\env_{\mathcal C} {{\mathcal C}^\star(N)}\in\Epi}{\Longleftarrow}\quad 0=u\circ \Env_{\mathcal C}\theta\circ \env_{\mathcal C} {{\mathcal C}^\star(N)}=
u\circ \env_{\mathcal C} {{\mathcal C}^\star(G)}\circ \theta \quad\overset{\overline{\sp\delta^N}={\mathcal C}^\star(N)}{\Longleftarrow}\\ \Longleftarrow\quad 0=u\circ \env_{\mathcal C} {{\mathcal C}^\star(G)}\circ \theta\circ\delta^N \quad\Longleftarrow\quad 0=u\Big|_N \quad\Longleftarrow\\ \Longleftarrow\quad 0=u\cdot 1_N
\quad\Longleftarrow\quad \exists v\in {\mathcal K}(G)\quad u=v\cdot (1-1_N) \quad\Longleftarrow\quad u\in {\mathcal K}_{G\setminus N}(G).
\end{multline*}

2. Докажем (ii). Вспомним, что по предложению \ref{PROP:Env_C(theta)} отображение $\Env_{\mathcal C} \theta:\Env_{\mathcal C} {\mathcal C}^\star(N)\to \Env_{\mathcal C} {\mathcal C}^\star(G)$ инъективно и открыто. Отсюда следует \cite[Theorem 4.2.29]{Akbarov-De-Gruyter-I}, что сопряженное отображение $\psi:{\mathcal K}(G)\to {\mathcal K}(N)$ замкнуто. С другой стороны, из того, что $\Env_{\mathcal C}\theta$ инъективно, следует, что $\psi=(\Env_{\mathcal C}\theta)^\star$ --- эпиморфизм стереотипных пространств. Вместе с его замкнутостью это значит, что $\psi$ сюръективно. 

Отсюда следует, что его кообраз $\psi|_{{\mathcal K}_N(G)}:{\mathcal K}(G)/\Ker\psi={\mathcal K}_N(G)\to {\mathcal K}(N)$ --- тоже замкнутое сюръективное отображение. Заметим, что $\psi|_{{\mathcal K}_N(G)}$ вдобавок инъективно:
$$
u\in {\mathcal K}_N(G)\ \&\ \psi(u)=0\quad\Longrightarrow\quad u\in {\mathcal K}_N(G)\cap\Ker\psi={\mathcal K}_N(G)\cap{\mathcal K}_{G\setminus N}(G)=0\quad\Longrightarrow\quad u=0.
$$
Наконец, еще одно важное замечание: пространство ${\mathcal K}(N)$ кополно (в смысле определения на с.\pageref{DEF:kopolno}), потому что его сопряженое пространство
$$
 {\mathcal K}(N)^\star =\Env_{\mathcal C} {\mathcal C}^\star(N)=\eqref{Env_C-C*(Z-cdot-K)}=\prod_{\sigma\in\widehat{K}}{\mathcal C}\big(M_\sigma,{\mathcal B}(X_\sigma)\big),
$$
-- полно. Таким образом, $\psi|_{{\mathcal K}_N(G)}:{\mathcal K}_N(G)\to {\mathcal K}(N)$ -- замкнутое биективное линейное непрерывное отображение стереотипных пространств, причем его область значений -- кополное пространство. По теореме \ref{TH:zamk-biektsija-v-kopolnoe-prostrancsvo}, $\psi|_{{\mathcal K}_N(G)}:{\mathcal K}_N(G)\to {\mathcal K}(N)$ -- изоморфизм стереотипных пространств.
\epr

Заметим теперь вот что. Характеристические функции $1_L$ в сумме дают единицу алгебры ${\mathcal K}(G)$,
\beq\label{DEF:1=sums_L-1_L}
1=\sum_{L\in G/N}1_L,
\eeq
(здесь конечное число слагаемых) поэтому пространство ${\mathcal K}(G)$ раскладывается в прямую сумму пространств ${\mathcal K}_L(G)$:
\beq\label{K(G)=oplus_L-K_L(G)}
{\mathcal K}(G)=\bigoplus_{L\in G/N}{\mathcal K}_L(G).
\eeq
Обозначим $E_N$ образ пространства $\Env_{\mathcal C} {\mathcal C}^\star(N)$ в $\Env_{\mathcal C} {\mathcal C}^\star(G)$ при отображении $\Env_{\mathcal C} {\mathcal C}^\star(N)\to \Env_{\mathcal C} {\mathcal C}^\star(G)$. 
По предложению \ref{PROP:Env_C(theta)}, это отображение открыто, поэтому
\beq\label{E_N=Env_C^star(N)}
E_N\cong \Env_{\mathcal C} {\mathcal C}^\star(N)={\mathcal K}(N)^\star
\eeq 

Пусть для всякого $L\in G/N$ символ $E_L$ обозначает сдвиг пространства $E_N$ на произвольный элемент $x_L\in L$:
\beq\label{DEF:E_L=x_L-cdot-E_N}
E_L=x_L\cdot E_N.
\eeq
Поскольку по теореме \ref{LM:sdvig-v-K(G)} сдвиг является автоморфизмом пространства ${\mathcal K}(G)$, он является автоморфизмом и сопряженного пространства $\Env_{\mathcal C} {\mathcal C}^\star(G)$, и поэтому пространство $E_L$ корректно определено. С другой стороны, его определение не зависит от выбора элемента $x_L\in L$, потому что если $y_l\in L$ --- какой-то другой элемент, то 
$$
y_L=x_L\cdot h
$$
для некоторого элемента $h\in N$, и мы получим
$$
y_L\cdot E_N=x_L\cdot h\cdot E_N=x_L\cdot E_N.
$$
С другой стороны, из того, что сдвиг является автоморфизмом пространства $\Env_{\mathcal C} {\mathcal C}^\star(G)$ сразу следует, что пространства $E_L$ изоморфны $E_N$:
\beq\label{E_L-cong-E_N}
E_L\cong E_N,\qquad L\in G/N.
\eeq

При таких обозначениях сопряженное пространство к каждому ${\mathcal K}_L(G)$ можно отождествить с пространством $E_L$, 
\beq\label{K_L(G)^star=E_L}
{\mathcal K}_L(G)^\star=E_L,
\eeq
поскольку
$$
{\mathcal K}_L(G)^\star= \Big( {\mathcal K}(G)\cdot 1_L \Big)^\star=
 \Big( {\mathcal K}(G)\cdot 1_N\cdot x_L \Big)^\star=\Big( {\mathcal K}_N(G)\cdot x_L \Big)^\star=
 x_L\cdot {\mathcal K}_N(G)^\star=x_L\cdot {\mathcal K}(N)^\star=x_L\cdot E_N=E_L.
$$
Как следствие, 
$$
\Env_{\mathcal C}  {\mathcal C}^\star(G)={\mathcal K}(G)^\star=\eqref{K(G)=oplus_L-K_L(G)}=
\l\bigoplus_{L\in G/N}{\mathcal K}_L(G)\r^\star=
\bigoplus_{L\in G/N}{\mathcal K}_L(G)^\star=\eqref{K_L(G)^star=E_L}=\bigoplus_{L\in G/N} E_L,
$$
то есть справедливо разложение 
\beq\label{Env_C-C*(konech-rasshirenie-ZK)}
\Env_{\mathcal C}  {\mathcal C}^\star(G)=\bigoplus_{L\in G/N} E_L.
\eeq

\btm\label{TM:Env_C-C*(G)-Lie-Moore}
Пусть $G$ --- локально компактная группа, являющаяся конечным расширением надстройки $N=Z\cdot K$ некоторой абелевой группы $Z$ с помощью некоторой компактной группы $K$, то есть такая, для которой существует цепочка непрерывных гомоморфизмов  
\beq\label{G-kak-rasshirenie-2}
\xymatrix 
{
1\ar[r] & Z\cdot K=N \ar[r]^{\eta} & G\ar[r]^{\ph} & F\ar[r]& 1
}
\eeq
в которой $F$ --- конечная группа, и
\beq
\eta(N)=\Ker\ph.
\eeq
Тогда 
\bit{

\item[(i)] непрерывная оболочка алгебры мер ${\mathcal C}^\star(G)$ как стереотипное пространство совпадает с непрерывной оболочкой алгебры мер ${\mathcal C}^\star(N)$, возведенной в степень $\card F$: 
\beq\label{E(C*(G))=E(C*(N))^(card-F)}
\Env_{\mathcal C} {\mathcal C}^\star(G)\cong \Big(\Env_{\mathcal C} {\mathcal C}^\star(N)\Big)^{\card F} 
\eeq

\item[(ii)] алгебра ${\mathcal K}(G)$ как стереотипное пространство совпадает с алгеброй ${\mathcal K}(N)$, возведенной в степень $\card F$: 
\beq\label{K(G)=K(N)^(card-F)}
{\mathcal K}(G)\cong \Big({\mathcal K}(N)\Big)^{\card F} 
\eeq

\item[(iii)] пространство $\Env_{\mathcal C} {\mathcal C}^\star(G)$ насыщено.

}\eit

\etm
\bpr
Формулы \eqref{E(C*(G))=E(C*(N))^(card-F)} и \eqref{K(G)=K(N)^(card-F)} переходят одна в другую под действием операции сопряжения $\star$, поэтому достаточно доказать одну из них, например, вторую.
Это делается так:
\begin{multline*}
{\mathcal K}(G)=\eqref{K(G)=oplus_L-K_L(G)}=\oplus_{L\in G/N}{\mathcal K}_L(G)
\cong \eqref{K_L(G)=K_N(G)} \cong \\ \cong \oplus_{L\in G/N}{\mathcal K}_N(G)\cong {\mathcal K}_N(G)^{\card G/N}\cong {\mathcal K}_N(G)^{\card F}
\cong \eqref{K_N(G)=K(N)} \cong {\mathcal K}(N)^{\card F}
\end{multline*}

Утверждение (iii) следует из теоремы \ref{TH:Env_C-C*(Z-cdot-K)} (ii).
\epr

\bcor\label{COR:stroenie-nepr-Env-Lie-Moore} Пусть $G$ --- группа Ли---Мура. Тогда
\beq\label{stroenie-nepr-Env-Lie-Moore}
\Env_{\mathcal C} {\mathcal C}^\star(G)\cong 
\prod_{s\in S}{\mathcal C}(M_s),
\eeq
где $\{M_s;\ s\in S\}$ --- некое семейство паракомпактных локально компактных топологических пространств.
\ecor
\bpr
Здесь нужно воспользоваться теоремой \ref{TH:Lie-Moore} и представить $G$ как конечное расширение компактной надстройки абелевой группы
$$
\xymatrix 
{
1\ar[r] & Z\cdot K=N \ar[r]^{\eta} & G\ar[r]^{\ph} & F\ar[r]& 1
}
$$
и после этого заметить, что в формуле \eqref{E(C*(G))=E(C*(N))^(card-F)} пространство $\Env_{\mathcal C} {\mathcal C}^\star(N)$ описано равенством \eqref{Env_C-C*(Z-cdot-K)}:
$$
\Env_{\mathcal C} {\mathcal C}^\star(N)= 
\Env_{\mathcal C} {\mathcal C}^\star(Z\cdot K)\cong\eqref{Env_C-C*(Z-cdot-K)}\cong
\prod_{\sigma\in\widehat{K}}{\mathcal C}\big(M_\sigma,{\mathcal B}(X_\sigma)\big)=
\prod_{\sigma\in\widehat{K}}{\mathcal C}(M_\sigma)^{(\dim X_\sigma)^2}.
$$
\epr

\bcor\label{COR:nepr-Env-Lie-Moore=LCS-lim} Пусть $G$ --- группа Ли---Мура. Тогда непрерывная оболочка групповой алгебры ${\mathcal C}^\star(G)$ совпадает с локально выпуклой оболочкой Кузнецовой, то есть с проективным пределом ее $C^*$-фактор-алгебр в категории локально выпуклых пространств (и в категории топологических алгебр):
\beq\label{nepr-Env-Lie-Moore=LCS-lim}
\Env_{\mathcal C}{\mathcal C}^\star(G)={\tt LCS}\text{-}\kern-3pt\projlim_{q\in\Sn({\mathcal C}^\star(G))}{\mathcal C}^\star(G)/q
\eeq
\ecor
\bpr
По предложению \ref{PROP:otkrytost-Env-C^*(N)->Env-C^*(G)}, всякая полунорма $p\in\Sn(N)$ мажорируется некоторой полунормой $p'\in\Sn(G)$. Докажем цепочку равенств
\beq\label{nepr-Env-Lie-Moore=LCS-lim-1}
\Env_{\mathcal C} {\mathcal C}^\star(G)=
\LCS\text{-}\kern-3pt\projlim_{p\in \Sn(N)}{\mathcal C}^\star(G)/p'=
\LCS\text{-}\kern-3pt\projlim_{q\in\Sn(G)}{\mathcal C}^\star(G)/q
\eeq

1. Чтобы доказать первое равенство в \eqref{nepr-Env-Lie-Moore=LCS-lim-1}\, достаточно убедиться, что $\Env_{\mathcal C} {\mathcal C}^\star(G)$ полно и что полунормы $\{p';\ p\in\Sn(N)\}$ задают топологию на этом пространстве. 
Вспомним разложение \eqref{Env_C-C*(konech-rasshirenie-ZK)}:
\beq\label{nepr-Env-Lie-Moore=LCS-lim-3}
\Env_{\mathcal C} {\mathcal C}^\star(G)=\oplus_{L\in G/N} E_L
\eeq
При этом 
\beq\label{nepr-Env-Lie-Moore=LCS-lim-4}
E_L\cong\eqref{E_L-cong-E_N}\cong E_N\cong\eqref{E_N=Env_C^star(N)}\cong {\mathcal K}(N)^\star
\cong \Env_{\mathcal C} {\mathcal C}^\star(N) \cong 
\prod_{\sigma\in\widehat{K}}{\mathcal C}\big(M_\sigma,{\mathcal B}(X_\sigma)\big),\qquad L\in G/N.
\eeq
(изоморфизмы в $\LCS$). Из этой цепочки сразу видно, что все пространства $E_L$ полны, и поэтому их прямая сумма \eqref{nepr-Env-Lie-Moore=LCS-lim-3}, то есть  пространство $\Env_{\mathcal C} {\mathcal C}^\star(G)$ тоже полно.

Теперь покажем, что $\{p';\ p\in\Sn(N)\}$ задают топологию на этом пространстве. Для этого достаточно убедиться, что они задают топологию на каждой компоненте $E_L$. Это следует из формулы \eqref{otkrytost-Env-C^*(N)->Env-C^*(G)-1}: по ней полунормы $\{p';\ p\in\Sn(G)\}$ мажорируют на $E_L$ полунормы вида 
\beq\label{nepr-Env-Lie-Moore=LCS-lim-5}
\beta\in E_L \mapsto p(x_L^{-1}\cdot\beta)
\eeq
Поэтому порождаемая ими топология должна быть не слабее топологии, порождаемой на $E_L$ полунормами \eqref{nepr-Env-Lie-Moore=LCS-lim-5}. Но с другой стороны, полунормы \eqref{nepr-Env-Lie-Moore=LCS-lim-5} порождают на $E_L$ топологию этого пространства, потому что при изоморфизме (здесь применяется теорема \ref{LM:sdvig-v-K(G)})
$$
\beta\in E_L \mapsto x_L^{-1}\cdot\beta\in E_N
$$ 
они превращаются в полунормы $\{p;\ p\in\Sn(N)\}$, которые по уже доказанной формуле \eqref{Env_C-Z.K=LCS-lim} порождают топологию этого пространства.

Мы получаем, что полунормы $\{p';\ p\in\Sn(N)\}$ порождают на $E_L$ топологию, не слабее, чем собственная топология этого пространства. С другой стороны, она не может быть сильнее, потому что все эти полунормы непрерывны на $E_L$. Значит, полунормы $\{p';\ p\in\Sn(N)\}$ порождают на $E_L$ в точности топологию этого пространства. 

И это доказывает первое равенство в \eqref{nepr-Env-Lie-Moore=LCS-lim-1}.

2. Чтобы доказать первое равенство в \eqref{nepr-Env-Lie-Moore=LCS-lim-1}, теперь достаточно (вспомнить, что $\Env_{\mathcal C} {\mathcal C}^\star(G)$ полно и) заметить, что семейство полунорм $\{q;\ q\in\Sn(G)\}$ задает топологию на этом пространстве. Это следует из того, что содержащееся в нем подсемейство $\{p';\ p\in\Sn(N)\}$, как мы уже поняли, задает топологию на этом пространстве.

\epr

\paragraph{Метризуемые группы Мура.}

\btm\label{TH:metrizuemaya-gruppa-Mura}
Группа Мура $G$ тогда и только тогда метризуема, когда она является проективным пределом некоей последовательности своих фактор-групп Ли---Мура:
\beq\label{metrizuemaya-gruppa-Mura}
G=\projlim_{n\in\N}G/K_n.
\eeq 
\etm

Пусть $G$ --- произвольная метризуемая группа Мура. Представим ее в виде предела \eqref{metrizuemaya-gruppa-Mura}. Каждая проекция $$
G\to G/K_n
$$
индуцирует гомоморфизм алгебр непрерывных функций
\beq\label{Env_C-C*(G)->projlim-Env_C-C*(G/K_n)-metriz-3}
\ph_n:{\mathcal C}(G)\gets{\mathcal C}(G/K_n)
\eeq
и сопряженный ему гомоморфизм алгебр мер
\beq\label{Env_C-C*(G)->projlim-Env_C-C*(G/K_n)-metriz-2}
\ph_n^\star:{\mathcal C}^\star(G)\to{\mathcal C}^\star(G/K_n)
\eeq
А тот в свою очередь индуцирует гомоморфизм оболочек
\beq\label{Env_C-C*(G)->projlim-Env_C-C*(G/K_n)-metriz-1}
\Env_{\mathcal C}\ph_n^\star:\Env_{\mathcal C}{\mathcal C}^\star(G)\to \Env_{\mathcal C}{\mathcal C}^\star(G/K_n)
\eeq
и, как следствие, гомоморфизм оболочки в проективный предел
\beq\label{Env_C-C*(G)->projlim-Env_C-C*(G/K_n)-metriz}
\Ste\text{-}\projlim_{n\in\N}\Env_{\mathcal C}\ph_n^\star: \Env_{\mathcal C} {\mathcal C}^\star(G)\to \Ste\text{-}\projlim_{n\in\N}
\Env_{\mathcal C} {\mathcal C}^\star (G/K_n).
\eeq

\btm\label{TM:Env_C-C*(G)-metriz-gruppa}
Если $G$ ---  метризуемая группа Мура, то 
\bit{
\item[(a)] морфизм \eqref{Env_C-C*(G)->projlim-Env_C-C*(G/K_n)-metriz} является изоморфизмом в категории $\InvSteAlg$, и

\item[(b)] проективный предел в \eqref{Env_C-C*(G)->projlim-Env_C-C*(G/K_n)-metriz} совпадает с проективным пределом в категории локально выпуклых пространств $\LCS$:
\beq\label{Env_C-C*(G)=projlim-Env_C-C*(G/K)-metriz}
\Env_{\mathcal C} {\mathcal C}^\star(G)\cong \Ste\text{-}\projlim_{n\in\N}
\Env_{\mathcal C} {\mathcal C}^\star (G/K_n)\cong \LCS\text{-}\projlim_{n\in\N}
\Env_{\mathcal C} {\mathcal C}^\star (G/K_n).
\eeq

\item[(c)] пространство $\Env_{\mathcal C} {\mathcal C}^\star(G)$ насыщено.

}\eit
\etm

Для доказательства нам понадобится

\blm\label{LM:Env_C-C*(G)-osubarr-metriz-gruppa}
Гомоморфизм \eqref{Env_C-C*(G)->projlim-Env_C-C*(G/K_n)-metriz} является непосредственным мономорфизмом стереотипных пространств:
\beq\label{Env_C-C*(G)-osubarr-metriz-gruppa}
\Ste\text{-}\projlim_{n\in\N}\Env_{\mathcal C}\ph_n^\star: \Env_{\mathcal C} {\mathcal C}^\star(G)\osubarr
 \Ste\text{-}\projlim_{n\in\N}\Env_{\mathcal C} {\mathcal C}^\star (G/K_n).
\eeq
\elm
\bpr 1. Обозначим $A={\mathcal C}^\star(G)$ и рассмотрим множество ${\mathcal U}={\mathcal U}_{C^*}^A$ $C^*$-окрестностей нуля в алгебре $A$. Оно частично упорядочено отношением
$$
U\le V\quad\Leftrightarrow\quad U\supseteq V.
$$
По теореме \ref{TH:opisanie-nepr-obolochki}, непрерывная оболочка $\Env_{\mathcal C}A$ алгебры $A$ представляет собой  оболочку $\Env\pi(A)$ в проективном пределе $\Ste\text{-}\projlim_{U\in {\mathcal U}_{C^*}^A}A/U$ множества значений $\pi(A)$ морфизма
$$
\pi:A\to \Ste\text{-}\projlim_{U\in {\mathcal U}_{C^*}^A}A/U
$$
Или, что то же самое, узловой образ этого морфизма:
\beq\label{PROOF:Env_C-C*(G)-osubarr-metriz-gruppa-0}
\Env_{\mathcal C}A=\Env\pi(A)=\Im_\infty\pi
\eeq

2. Каждому индексу $n\in\N$ поставим в соответствие 
\bit{

\item[---] алгебру $A_n={\mathcal C}^\star(G/K_n)$,

\item[---]  множество ${\mathcal U}_{C^*}^{A_n}$ $C^*$-окрестностей нуля в алгебре $A_n$,

\item[---]  множество ${\mathcal U}_n$ прообразов окрестностей из ${\mathcal U}_{C^*}^{A_n}$ при отображении $\ph_n^\star:A\to A_n$:
\beq\label{PROOF:Env_C-C*(G)-osubarr-metriz-gruppa-1}
{\mathcal U}_n=\{(\ph_n^\star)^{-1}(V);\ V\in {\mathcal U}_{C^*}^{A_n}\}
\eeq

}\eit

По следствию \ref{COR:G->B-propusk-cherez-Lie}, у нас получится, что множества ${\mathcal U}_n$ покрывают множество ${\mathcal U}$:  
\beq\label{PROOF:Env_C-C*(G)-osubarr-metriz-gruppa-2}
{\mathcal U}=\bigcup_{n\in\N} {\mathcal U}_n
\eeq
А, с другой стороны, очевидно, они образуют возрастающую последовательность:
\beq\label{PROOF:Env_C-C*(G)-osubarr-metriz-gruppa-3}
{\mathcal U}_n\subseteq {\mathcal U}_{n+1},\qquad n\in\N.
\eeq
Вместе \eqref{PROOF:Env_C-C*(G)-osubarr-metriz-gruppa-2} и \eqref{PROOF:Env_C-C*(G)-osubarr-metriz-gruppa-3} означают, что система множеств 
$\{{\mathcal U}_n\}$ является покрытием частично упорядоченного множества ${\mathcal U}$ (в смысле определения на с.\pageref{DEF:pokrytie-POSet}).

3. Заметим далее, что, снова по теореме \ref{TH:opisanie-nepr-obolochki}, непрерывная оболочка $\Env_{\mathcal C}A_n$ алгебры $A_n$ представляет собой  оболочку $\Env\pi_n(A_n)$ в проективном пределе $\Ste\text{-}\projlim_{V\in {\mathcal U}_{C^*}^{A_n}}A_n/V$ множества значений $\pi_n(A_n)$ морфизма
$$
\pi_n:A_n\to \Ste\text{-}\projlim_{V\in {\mathcal U}_{C^*}^{A_n}}A_n/V
$$
То есть справедливо равенство   
\beq\label{PROOF:Env_C-C*(G)-osubarr-metriz-gruppa-4}
\Env_{\mathcal C}A_n=\Env\pi_n(A_n)=\Im_\infty\pi_n
\eeq
Из того, что всякий морфизм $\ph_n^\star:A\to A_n$ сюръективен, можно заключить, что множество окрестностей ${\mathcal U}_{C^*}^{A_n}$ здесь можно заменить множеством окрестностей ${\mathcal U}_n$, а алгебру $A_n$ на алгебру $A$: 
$$
\pi_n:A\to \Ste\text{-}\projlim_{U\in {\mathcal U}_n}A/U
$$
и формула \eqref{PROOF:Env_C-C*(G)-osubarr-metriz-gruppa-4} при этом не  изменится:
\beq\label{PROOF:Env_C-C*(G)-osubarr-metriz-gruppa-5}
\Env_{\mathcal C}A_n=\Env\pi_n(A_n)=\Im_\infty\pi_n
\eeq

4. Рассмотрим  морфизмы $\e(\pi_n)=\red_\infty\pi_n\circ\coim_\infty\pi_n$:
\beq\label{PROOF:Env_C-C*(G)-osubarr-metriz-gruppa-6}
\e(\pi_n):A\to \Im_\infty\pi_n=\Env_{\mathcal C}A_n
\eeq
и проективный предел этих морфизмов:
\beq\label{PROOF:Env_C-C*(G)-osubarr-metriz-gruppa-7}
\Ste\text{-}\projlim_{n\in\N}\e(\pi_n):A\to \Ste\text{-}\projlim_{n\in\N}\Env_{\mathcal C}A_n
\eeq
Если рассмотреть повторный проективный предел
\beq\label{PROOF:Env_C-C*(G)-osubarr-metriz-gruppa-8}
\Ste\text{-}\projlim_{n\in\N}\Big(\Ste\text{-}\projlim_{U\in {\mathcal U}_n}A/U\Big)
=\eqref{DEF:povtornyj-proj-limit} =\Ste\text{-}\projlim_{U\in {\mathcal U}}A/U
\eeq
то по теореме \ref{TH:ob-uzlov-obraze-v-povtornom-proj-predele}, мы получим, что узловой образ $\pi$ совпадает с узловым образом $\Ste\text{-}\projlim_{n\in\N}\e(\pi_n)$ 
\beq\label{PROOF:Env_C-C*(G)-osubarr-metriz-gruppa-9}
\Im_\infty\pi=\Im_\infty\Ste\text{-}\projlim_{n\in\N}\e(\pi_n)
\eeq

Это означает, что пространство $\Env_{\mathcal C}A$ являтся непосредственным подпространством в проективном пределе пространств   
$\Env_{\mathcal C}A_n$
\beq\label{PROOF:Env_C-C*(G)-osubarr-metriz-gruppa-10}
\Env_{\mathcal C}A\osubarr\Ste\text{-}\projlim_{n\in\N}\Env_{\mathcal C}A_n
\eeq
\epr

\bpr[Доказательство теоремы \ref{TM:Env_C-C*(G)-metriz-gruppa}.]
Вспомним проекторы $\pi_{K_n}^\star:{\mathcal C}^\star(G)\to {\mathcal C}^\star(G)$, определенные в \eqref{pi_K^star(alpha)=mu_K*alpha}.
По свойству \eqref{pi_K^star(alpha*beta)=pi_K^star(alpha)*pi_K^star(beta)} эти отображения мультипликативны, то есть $\pi_{K_n}^\star$ --- морфизмы в категории $\InvSteAlg^0$. Поэтому по теореме \ref{TH:Env^0_C-idemp-funktor} можно к ним применить функтор $\Env_{\mathcal C}^0$, и они превратятся в систему $\Env_{\mathcal C}^0\pi_{K_n}^\star$ эндоморфизмов в категории $\InvSteAlg^0$ алгебры $\Env_{\mathcal C}{\mathcal C}^\star(G)$. С другой стороны, по свойству \eqref{pi_K*^2=pi_K*} операторы $\pi_{K_n}^\star$ представляют собой систему проекторов для алгебры ${\mathcal C}^\star(G)$:
$$
(\pi_{K_n}^\star)^2=\pi_{K_n}^\star
$$
Поскольку $\Env_{\mathcal C}^0$ --- функтор, мы отсюда получаем равенство 
$$
\Big(\Env_{\mathcal C}^0\pi_{K_n}^\star\Big)^2=\Env_{\mathcal C}^0\pi_{K_n}^\star,
$$
означающее, что операторы $\Env_{\mathcal C}^0\pi_{K_n}^\star$ образуют систему проекторов на стереотипном 
пространстве $\Env_{\mathcal C}{\mathcal C}^\star(G)$. По теореме \ref{TH:building-generated-by-P^k} эта система 
проекторов порождает некий фундамент $\varPhi$ в категории $\Ste$ (и в категории $\InvSteAlg$) над упорядоченным множеством $\N$.

Покажем, что система проекторов $\Env_{\mathcal C}^0\pi_{K_n}^\star$ для алгебры $\Env_{\mathcal C}{\mathcal C}^\star(G)$ удовлетворяет посылкам теоремы \ref{TH:X=edifice-in-Ste}. 

1) Прежде всего, для всякой группы ${K_n}\in\lambda(G)$  оператор $\pi_{K_n}^\star$ раскладывается в композицию морфизмов в категории 
$\InvSteAlg^0$
$$
{\mathcal C}^\star(G)\overset{\rho_{K_n}}{\longrightarrow}{\mathcal C}^\star(G:{K_n})\overset{\sigma_{K_n}}{\longrightarrow}{\mathcal C}^\star(G)
$$
причем
$$
\rho_{K_n}\circ\sigma_{K_n}=1_{{\mathcal C}^\star(G:{K_n})}
$$
Поскольку $\Env_{\mathcal C}^0$ --- функтор, мы отсюда получаем, что оболочка $\Env_{\mathcal C}^0\pi_{K_n}^\star$ тоже раскладывается в композицию
$$
\Env_{\mathcal C}^0{\mathcal C}^\star(G)\overset{\Env_{\mathcal C}^0\rho_{K_n}}{\longrightarrow} \Env_{\mathcal C}^0{\mathcal C}^\star(G:{K_n})\overset{\Env_{\mathcal C}^0\sigma_{K_n}}{\longrightarrow} \Env_{\mathcal C}^0{\mathcal C}^\star(G)
$$
причем 
$$
\Env_{\mathcal C}^0\rho_{K_n}\circ\Env_{\mathcal C}^0\sigma_{K_n}=1_{\Env_{\mathcal C}^0{\mathcal C}^\star(G:{K_n})}
$$
Это означает, что образом проектора $\Env_{\mathcal C}^0\pi_{K_n}^\star$ должно быть пространство $\Env_{\mathcal C}{\mathcal C}^\star(G:{K_n})$:
$$
\Env_{\mathcal C}^0\pi_{K_n}^\star(\Env_{\mathcal C}{\mathcal C}^\star(G))=\Env_{\mathcal C}{\mathcal C}^\star(G:{K_n})
=\Env_{\mathcal C}{\mathcal C}^\star(G/{K_n})
$$ 
А по теореме \ref{TM:Env_C-C*(G)-Lie-Moore} последнее пространство стереотипно и насыщено.

2) Далее, пусть $n_i$ --- неубывающая последовательность индексов
$$
n_1\le n_2\le n_3 \le...
$$
Соответствующая им система подгрупп $K_{n_i}$ будет невозрастать по включению,
$$
K_{n_1}\supseteq K_{n_2}\supseteq K_{n_3}\supseteq...
$$
Положим 
$$
H=\bigcap_{i=1}^\infty K_{n_i}.
$$

Рассмотрим две возможные ситуации.

2.1) Последовательность индексов $n_i$ стабилизируется:
$$
\exists k\in\N\qquad \forall i\ge k \qquad n_k=n_i.
$$
Тогда 
$$
H=\bigcap_{i=1}^\infty K_{n_i}=K_{n_k},
$$
а ограничение  $\varPhi_J$ фундамента $\varPhi$, порождаемое системой индексов $J=\{n_i\}$, имеет здание
$$
\edif\varPhi_J=\Env_{\mathcal C} {\mathcal C}^\star(G/K_{n_k})=\Ste\text{-}\projlim_{i\in\N}
\Env_{\mathcal C} {\mathcal C}^\star (G/K_{n_i})\cong \LCS\text{-}\projlim_{i\in\N}
\Env_{\mathcal C} {\mathcal C}^\star (G/K_{n_i}).
$$
и оно дополняется в здании 
$$
\edif\varPhi=\Ste\text{-}\projlim_{n\in\N}
\Env_{\mathcal C} {\mathcal C}^\star (G/{K_n}),
$$
потому что коретракцией к морфизму 
$$
\edif\varPhi_J=\Env_{\mathcal C} {\mathcal C}^\star(G/K_{n_k})=\Ste\text{-}\projlim_{i\in\N}
\Env_{\mathcal C} {\mathcal C}^\star (G/K_{n_i})
\gets \Ste\text{-}\projlim_{n\in\N}\Env_{\mathcal C} {\mathcal C}^\star (G/{K_n})=\edif\varPhi
$$
будет морфизм
$$
\edif\varPhi_J=\Env_{\mathcal C} {\mathcal C}^\star(G/K_{n_k})
\stackrel{\Env_{\mathcal C}\im\pi_{K_{n_k}}^\star}{\longrightarrow}\Env_{\mathcal C} {\mathcal C}^\star(G)
\to \Ste\text{-}\projlim_{n\in\N}\Env_{\mathcal C} {\mathcal C}^\star (G/{K_n})=\edif\varPhi
$$

2.2) Последовательность индексов $n_i$ уходит в бесконечность:
$$
n_i\underset{i\to\infty}{\longrightarrow}\infty
$$
Тогда 
$$
H=\bigcap_{i=1}^\infty K_{n_i}=1,
$$
а ограничение  $\varPhi_J$ фундамента $\varPhi$, порождаемое системой индексов $J=\{n_i\}$, имеет то же здание, что и у $\varPhi$:
$$
\edif\varPhi_J=\edif\varPhi.
$$
Понятно, что $\edif\varPhi_J$ тривиально дополняемо в $\edif\varPhi$.

3) По лемме \ref{LM:Env_C-C*(G)-osubarr-metriz-gruppa}, морфизм \eqref{Env_C-C*(G)->projlim-Env_C-C*(G/K_n)-metriz}
 является непосредственным мономорфизмом в категории $\InvSteAlg$.

Теперь мы можем применить теорему \ref{TH:X=edifice-in-Ste}: формула \eqref{Env_C-C*(G)=projlim-Env_C-C*(G/K)-metriz}
 следует из формул 
\eqref{X=edif-varPhi} и \eqref{X=LCS-projlim_(k->infty)X_k}, а утверждение (с) --- из утверждения (b) 
теоремы \ref{TH:X=edifice-in-Ste}.
\epr

\paragraph{Общие группы Мура.}

Пусть $G$ -- произвольная группа Мура. Тогда по теореме \ref{TH:Moore=lim-Lie-Moore} ее можно представить как проективный предел групп Ли---Мура:
$$
G=\projlim_{K\in\lambda(G)} G/K
$$
Каждая проекция $G\to G/K$ индуцирует гомоморфизм алгебр непрерывных функций
$$
{\mathcal C}(G)\gets{\mathcal C}(G/K)
$$
и сопряженный ему гомоморфизм алгебр мер
$$
{\mathcal C}^\star(G)\to{\mathcal C}^\star(G/K)
$$
А тот в свою очередь индуцирует гомоморфизм оболочек
$$
\Env_{\mathcal C}{\mathcal C}^\star(G)\to \Env_{\mathcal C}{\mathcal C}^\star(G/K)
$$
и, как следствие, гомоморфизм оболочки в проективный предел
\beq\label{Env_C-C*(G)->projlim-Env_C-C*(G/K)}
\Env_{\mathcal C} {\mathcal C}^\star(G)\to \Ste\text{-}\projlim_{K\in\lambda(G)}
\Env_{\mathcal C} {\mathcal C}^\star (G/K).
\eeq

\btm\label{TM:Env_C-C*(G)=projlim-Env_C-C*(G/K)}
Если $G$ --- группа Мура, то 
\bit{
\item[(a)] морфизм \eqref{Env_C-C*(G)->projlim-Env_C-C*(G/K)} является изоморфизмом в категории $\InvSteAlg$, и

\item[(b)] проективный предел в \eqref{Env_C-C*(G)->projlim-Env_C-C*(G/K)} совпадает с проективным пределом в категории локально выпуклых пространств $\LCS$:
\beq\label{Env_C-C*(G)=projlim-Env_C-C*(G/K)}
\Env_{\mathcal C} {\mathcal C}^\star(G)\cong \Ste\text{-}\projlim_{K\in\lambda(G)}
\Env_{\mathcal C} {\mathcal C}^\star (G/K)\cong \LCS\text{-}\projlim_{K\in\lambda(G)}
\Env_{\mathcal C} {\mathcal C}^\star (G/K).
\eeq

\item[(c)] пространство $\Env_{\mathcal C} {\mathcal C}^\star(G)$ насыщено.

}\eit
\etm

Для доказательства нам понадобится

\blm\label{LM:Env_C-C*(G)-osubarr}
Гомоморфизм \eqref{Env_C-C*(G)->projlim-Env_C-C*(G/K)} является непосредственным мономорфизмом стереотипных пространств:
\beq\label{Env_C-C*(G)-osubarr}
\Ste\text{-}\projlim_{n\in\N}\Env_{\mathcal C}\ph_n^\star: \Env_{\mathcal C} {\mathcal C}^\star(G)\osubarr
 \Ste\text{-}\projlim_{n\in\N}\Env_{\mathcal C} {\mathcal C}^\star (G/K_n).
\eeq
\elm
\bpr 1. Обозначим $A={\mathcal C}^\star(G)$ и рассмотрим множество ${\mathcal U}={\mathcal U}_{C^*}^A$ $C^*$-окрестностей нуля в алгебре $A$. Оно частично упорядочено отношением 
$$
U\le V\quad\Leftrightarrow\quad U\supseteq V.
$$
По теореме \ref{TH:opisanie-nepr-obolochki}, непрерывная оболочка $\Env_{\mathcal C}A$ алгебры $A$ представляет собой  оболочку $\Env\pi(A)$ в проективном пределе $\Ste\text{-}\projlim_{U\in {\mathcal U}_{C^*}^A}A/U$ множества значений $\pi(A)$ морфизма 
$$
\pi:A\to \Ste\text{-}\projlim_{U\in {\mathcal U}_{C^*}^A}A/U
$$
Или, что то же самое, узловой образ этого морфизма:
\beq\label{PROOF:Env_C-C*(G)-osubarr-0}
\Env_{\mathcal C}A=\Env\pi(A)=\Im_\infty\pi
\eeq

2. Каждой группе $K\in\lambda(G)$ поставим в соответствие 
\bit{

\item[---] алгебру $A_K={\mathcal C}^\star(G/K)$,

\item[---]  множество ${\mathcal U}_{C^*}^{A_K}$ $C^*$-окрестностей нуля в алгебре $A_K$,

\item[---]  множество ${\mathcal U}_K$ прообразов окрестностей из ${\mathcal U}_{C^*}^{A_K}$ при отображении $\ph_K^\star:A\to A_K$:
\beq\label{PROOF:Env_C-C*(G)-osubarr-1}
{\mathcal U}_K=\{(\ph_K^\star)^{-1}(V);\ V\in {\mathcal U}_{C^*}^{A_K}\}
\eeq

}\eit

По следствию \ref{COR:G->B-propusk-cherez-Lie}, у нас получится, что множества ${\mathcal U}_K$ покрывают множество ${\mathcal U}$:  
\beq\label{PROOF:Env_C-C*(G)-osubarr-2}
{\mathcal U}=\bigcup_{K\in\lambda(G)} {\mathcal U}_K
\eeq
А, с другой стороны, очевидно, они направлены по включению:
\beq\label{PROOF:Env_C-C*(G)-osubarr-3}
\forall K,L\in\lambda(G) \qquad \exists M=K\cap L\in\lambda(G)\qquad {\mathcal U}_K\cup{\mathcal U}_L \subseteq {\mathcal U}_M.
\eeq
Вместе \eqref{PROOF:Env_C-C*(G)-osubarr-2} и \eqref{PROOF:Env_C-C*(G)-osubarr-3} означают, что система множеств 
$\{{\mathcal U}_K\}$ является покрытием частично упорядоченного множества ${\mathcal U}$ (в смысле определения на с.\pageref{DEF:pokrytie-POSet}).

3. Заметим далее, что, снова по теореме \ref{TH:opisanie-nepr-obolochki}, непрерывная оболочка $\Env_{\mathcal C}A_K$ алгебры $A_K$ представляет собой  оболочку $\Env\pi_K(A_K)$ в проективном пределе $\Ste\text{-}\projlim_{V\in {\mathcal U}_{C^*}^{A_K}}A_K/V$ множества значений $\pi_K(A_K)$ морфизма
$$
\pi_K:A_K\to \Ste\text{-}\projlim_{V\in {\mathcal U}_{C^*}^{A_K}}A_K/V
$$
То есть справедливо равенство   
\beq\label{PROOF:Env_C-C*(G)-osubarr-4}
\Env_{\mathcal C}A_K=\Env\pi_K(A_K)=\Im_\infty\pi_K
\eeq
Из того, что всякий морфизм $\ph_K^\star:A\to A_K$ сюръективен, можно заключить, что множество окрестностей ${\mathcal U}_{C^*}^{A_K}$ здесь можно заменить множеством окрестностей ${\mathcal U}_K$, а алгебру $A_K$ на алгебру $A$: 
$$
\pi_K:A\to \Ste\text{-}\projlim_{U\in {\mathcal U}_K}A/U
$$
и формула \eqref{PROOF:Env_C-C*(G)-osubarr-4} при этом не  изменится:
\beq\label{PROOF:Env_C-C*(G)-osubarr-5}
\Env_{\mathcal C}A_K=\Env\pi_K(A_K)=\Im_\infty\pi_K
\eeq

4. Рассмотрим  морфизмы $\e(\pi_K)=\red_\infty\pi_K\circ\coim_\infty\pi_K$:
\beq\label{PROOF:Env_C-C*(G)-osubarr-6}
\e(\pi_K):A\to \Im_\infty\pi_K=\Env_{\mathcal C}A_K
\eeq
и проективный предел этих морфизмов:
\beq\label{PROOF:Env_C-C*(G)-osubarr-7}
\Ste\text{-}\projlim_{K\in\lambda(G)}\e(\pi_K):A\to \Ste\text{-}\projlim_{K\in\lambda(G)}\Env_{\mathcal C}A_K
\eeq
Если рассмотреть повторный проективный предел
\beq\label{PROOF:Env_C-C*(G)-osubarr-8}
\Ste\text{-}\projlim_{K\in\lambda(G)}\Big(\Ste\text{-}\projlim_{U\in {\mathcal U}_K}A/U\Big)
=\eqref{DEF:povtornyj-proj-limit} =\Ste\text{-}\projlim_{U\in {\mathcal U}}A/U
\eeq
то по теореме \ref{TH:ob-uzlov-obraze-v-povtornom-proj-predele}, мы получим, что узловой образ $\pi$ совпадает с узловым образом $\Ste\text{-}\projlim_{K\in\lambda(G)}\e(\pi_K)$ 
\beq\label{PROOF:Env_C-C*(G)-osubarr-9}
\Im_\infty\pi=\Im_\infty\Ste\text{-}\projlim_{K\in\lambda(G)}\e(\pi_K)
\eeq

Это означает, что пространство $\Env_{\mathcal C}A$ являтся непосредственным подпространством в проективном пределе пространств   
$\Env_{\mathcal C}A_K$
\beq\label{PROOF:Env_C-C*(G)-osubarr-10}
\Env_{\mathcal C}A\osubarr\Ste\text{-}\projlim_{K\in\lambda(G)}\Env_{\mathcal C}A_K
\eeq
\epr

\bpr[Доказательство теоремы \ref{TM:Env_C-C*(G)=projlim-Env_C-C*(G/K)}.]
Вспомним проекторы $\pi_K^\star:{\mathcal C}^\star(G)\to {\mathcal C}^\star(G)$, определенные в \eqref{pi_K^star(alpha)=mu_K*alpha}.
По свойству \eqref{pi_K^star(alpha*beta)=pi_K^star(alpha)*pi_K^star(beta)} эти отображения мультипликативны, то есть $\pi_K^\star$ --- морфизмы в категории $\InvSteAlg^0$. Поэтому по теореме \ref{TH:Env^0_C-idemp-funktor} можно к ним применить функтор $\Env_{\mathcal C}^0$, и они превратятся в систему $\Env_{\mathcal C}^0\pi_K^\star$ эндоморфизмов в категории $\InvSteAlg^0$ алгебры $\Env_{\mathcal C}{\mathcal C}^\star(G)$. С другой стороны, по свойству \eqref{pi_K*^2=pi_K*} операторы $\pi_K^\star$ представляют собой систему проекторов для алгебры ${\mathcal C}^\star(G)$:
$$
(\pi_K^\star)^2=\pi_K^\star
$$
Поскольку $\Env_{\mathcal C}^0$ --- функтор, мы отсюда получаем равенство 
$$
\Big(\Env_{\mathcal C}^0\pi_K^\star\Big)^2=\Env_{\mathcal C}^0\pi_K^\star,
$$
означающее, что операторы $\Env_{\mathcal C}^0\pi_K^\star$ образуют систему согласованных проекторов на стереотипном пространстве $\Env_{\mathcal C}{\mathcal C}^\star(G)$. По теореме \ref{TH:building-generated-by-P^k} эта система проекторов порождает некий фундамент $\varPhi$ в категории $\Ste$ (и в категории $\InvSteAlg$) над упорядоченным множеством $\lambda(G)$.

Покажем, что система проекторов $\Env_{\mathcal C}^0\pi_K^\star$ для алгебры $\Env_{\mathcal C}{\mathcal C}^\star(G)$ удовлетворяет посылкам теоремы \ref{TH:X=edifice-in-Ste}. 

1) Прежде всего, для всякой группы $K\in\lambda(G)$  оператор $\pi_K^\star$ раскладывается в композицию морфизмов в категории 
$\InvSteAlg^0$
$$
{\mathcal C}^\star(G)\overset{\rho_K}{\longrightarrow}{\mathcal C}^\star(G:K)\overset{\sigma_K}{\longrightarrow}{\mathcal C}^\star(G)
$$
причем
$$
\rho_K\circ\sigma_K=1_{{\mathcal C}^\star(G:K)}
$$
Поскольку $\Env_{\mathcal C}^0$ --- функтор, мы отсюда получаем, что оболочка $\Env_{\mathcal C}^0\pi_K^\star$ тоже раскладывается в композицию
$$
\Env_{\mathcal C}^0{\mathcal C}^\star(G)\overset{\Env_{\mathcal C}^0\rho_K}{\longrightarrow} \Env_{\mathcal C}^0{\mathcal C}^\star(G:K)\overset{\Env_{\mathcal C}^0\sigma_K}{\longrightarrow} \Env_{\mathcal C}^0{\mathcal C}^\star(G)
$$
причем 
$$
\Env_{\mathcal C}^0\rho_K\circ\Env_{\mathcal C}^0\sigma_K=1_{\Env_{\mathcal C}^0{\mathcal C}^\star(G:K)}
$$
Это означает, что образом проектора $\Env_{\mathcal C}^0\pi_K^\star$ должно быть пространство $\Env_{\mathcal C}{\mathcal C}^\star(G:K)$:
$$
\Env_{\mathcal C}^0\pi_K^\star(\Env_{\mathcal C}{\mathcal C}^\star(G))=\Env_{\mathcal C}{\mathcal C}^\star(G:K)
=\Env_{\mathcal C}{\mathcal C}^\star(G/K)
$$ 
А по теореме \ref{TM:Env_C-C*(G)-metriz-gruppa} последнее пространство стереотипно и насыщено.

2) Далее, пусть $K_i$ -- убывающая по включению система компактных подгрупп,
$$
K_1\supseteq K_2\supseteq K_3\supseteq...
$$
и пусть 
$$
H=\bigcap_{i=1}^\infty K_i.
$$
Система  $K_i$ порождает некое ограничение  $\varPhi_J$ фундамента $\varPhi$. По теореме \ref{TM:Env_C-C*(G)-metriz-gruppa}, здание над фундаментом $\varPhi_J$, то есть проективный предел  системы  $\Env_{\mathcal C}{\mathcal C}^\star(G/K_i)$, совпадает с алгеброй $\Env_{\mathcal C}{\mathcal C}^\star(G/H)$:
$$
\Env_{\mathcal C} {\mathcal C}^\star(G/H)=\edif\varPhi_J=\Ste\text{-}\projlim_{i\in\N}
\Env_{\mathcal C} {\mathcal C}^\star (G/K_i)\cong \LCS\text{-}\projlim_{i\in\N}
\Env_{\mathcal C} {\mathcal C}^\star (G/K_i).
$$
С другой стороны, оно дополняется в здании $\Ste\text{-}\projlim_{K\in\lambda(G)}
\Env_{\mathcal C} {\mathcal C}^\star (G)$, потому что коретракцией к морфизму 
$$
\edif\varPhi_J=\Env_{\mathcal C} {\mathcal C}^\star(G/H)=\Ste\text{-}\projlim_{i\in\N}
\Env_{\mathcal C} {\mathcal C}^\star (G/K_i)
\gets \Ste\text{-}\projlim_{K\in\lambda(G)}\Env_{\mathcal C} {\mathcal C}^\star (G/K)=\edif\varPhi
$$
будет морфизм
$$
\edif\varPhi_J=\Env_{\mathcal C} {\mathcal C}^\star(G/H)\stackrel{\Env_{\mathcal C}\im\pi_H^\star}{\longrightarrow}\Env_{\mathcal C} {\mathcal C}^\star(G)
\to \Ste\text{-}\projlim_{K\in\lambda(G)}\Env_{\mathcal C} {\mathcal C}^\star (G/K)=\edif\varPhi
$$

3) По лемме \ref{LM:Env_C-C*(G)-osubarr}, морфизм \eqref{Env_C-C*(G)->projlim-Env_C-C*(G/K)} является непосредственным мономорфизмом в категории $\InvSteAlg$.

Теперь мы можем применить теорему \ref{TH:X=edifice-in-Ste}: формула \eqref{Env_C-C*(G)=projlim-Env_C-C*(G/K)} следует из формул 
\eqref{X=edif-varPhi} и \eqref{X=LCS-projlim_(k->infty)X_k}, а утверждение (с) --- из утверждения (b) теоремы \ref{TH:X=edifice-in-Ste}.
\epr

\paragraph{Дополнение: непрерывная оболочка групповой алгебры распределений ${\mathcal E}^\star(G)$.}

Если $G$ -- вещественная группа Ли, то помимо групповой алгебры ${\mathcal C}^\star(G)$ мер с компактным носителем на $G$ можно рассмотреть групповую алгебру ${\mathcal E}^\star(G)$ распределений с компактным носителем на $G$ (см.\cite{Ak03}). Если обозначить естественное включение ${\mathcal E}(G)\subseteq{\mathcal C}(G)$ каким-нибудь символом, например, $\lambda$, то мы получим диаграмму
\beq\label{Env_C-C^star(G)-Env_C-E^star(G)}
 \xymatrix @R=2.5pc @C=4pc
 {
{\mathcal C}^\star(G)\ar[r]^{\lambda^\star}\ar[d]_{\env_{\mathcal C}{\mathcal C}^\star(G)} & {\mathcal E}^\star(G)\ar[d]^{\env_{\mathcal C}{\mathcal E}^\star(G)}
\\
\Env_{\mathcal C}{\mathcal C}^\star(G)\ar[r]^{\Env_{\mathcal C}(\lambda^\star)} & \Env_{\mathcal C}{\mathcal E}^\star(G)
 }
\eeq

\btm\label{TH:Env_C-C^star(G)=Env_C-E^star(G)}
Для всякой вещественной группы Ли $G$ непрерывные оболочки групповых алгебр ${\mathcal C}^\star(G)$ и ${\mathcal E}^\star(G)$ совпадают:
\beq\label{Env_C-C^star(G)=Env_C-E^star(G)}
\Env_{\mathcal C}{\mathcal C}^\star(G)=\Env_{\mathcal C}{\mathcal E}^\star(G).
\eeq
\etm
\bpr
Рассмотрим композицию
$$
 \xymatrix @R=2.5pc @C=4pc
 {
{\mathcal C}^\star(G)\ar[r]^{\lambda^\star} & {\mathcal E}^\star(G)\ar[r]^{\env_{\mathcal C}{\mathcal E}^\star(G)}&
\Env_{\mathcal C}{\mathcal E}^\star(G)
 }
$$
и убедимся, что она является непрерывным расширением алгебры ${\mathcal C}^\star(G)$. Действительно, если $\ph:{\mathcal C}^\star(G)\to B$ -- морфизм в какую-нибудь $C^*$-алгебру $B$, то ему соответствует непрерывное (по норме) представление $\pi=\ph\circ\delta:G\to B$, которое, в силу примера \ref{EX:gladkost-nepr-po-norme-predstavleniya}, будет гладким, и значит, будет продолжаться до некоторого гомоморфизма $\ph'=\ddot{\pi}:{\mathcal E}^\star(G)\to B$. Этот гомоморфизм $\ddot{\pi}$ будет продолжать $\ph$, потому что они совпадают на элементах вида $\delta^x$, $x\in G$, которые полны в ${\mathcal C}^\star(G)$ и ${\mathcal E}^\star(G)$. После того, как построен гомоморфизм $\ph'$, он однозначно продолжается до некоторого гомоморфизма $\ph''$ на алгебру $\Env_{\mathcal C}{\mathcal E}^\star(G)$, как на непрерывное расширение алгебры ${\mathcal E}^\star(G)$:
$$
 \xymatrix @R=3pc @C=4pc
 {
{\mathcal C}^\star(G)\ar@/_3ex/[dr]_{\ph}\ar[r]^{\lambda^\star} & {\mathcal E}^\star(G)\ar[r]^{\env_{\mathcal C}{\mathcal E}^\star(G)}\ar@{-->}[d]^{\ph'=\ddot{\pi}} &
\Env_{\mathcal C}{\mathcal E}^\star(G)\ar@/^3ex/@{-->}[dl]^{\ph''}\\
& B & \\
 }
$$
Это доказывает, что  $\Env_{\mathcal C}{\mathcal E}^\star(G)$ -- непрерывное расширение алгебры ${\mathcal C}^\star(G)$. Значит, существует единственный морфизм $\upsilon$ из $\Env_{\mathcal C}{\mathcal E}^\star(G)$ в непрерывную обоолочку ${\mathcal C}^\star(G)$, замыкающий диаграмму
$$
 \xymatrix @R=2.5pc @C=4pc
 {
& {\mathcal C}^\star(G)\ar@/_3ex/[dl]_{\lambda^\star}\ar@/^3ex/[rd]^{\quad\env_{\mathcal C}{\mathcal C}^\star(G)\circ\lambda^\star} &
\\
\Env_{\mathcal C}{\mathcal C}^\star(G) & & \Env_{\mathcal C}{\mathcal E}^\star(G)\ar@{-->}[ll]_{\upsilon}
 }
$$
Из того, что ${\mathcal C}^\star(G)$ плотно в $\Env_{\mathcal C}{\mathcal C}^\star(G)$ и $\Env_{\mathcal C}{\mathcal E}^\star(G)$, следует, что $\upsilon$ -- обратный морфизм к $\Env_{\mathcal C}(\lambda^\star)$ из \eqref{Env_C-C^star(G)-Env_C-E^star(G)}.
\epr

\subsection{Свойства алгебры ${\mathcal K}(G)$ для групп Мура}

\paragraph{Отображение $\omega_{G.H}^\star:{\mathcal K}(G)\circledast{\mathcal K}(H)\to{\mathcal K}(G\times H)$.}

Пусть $G$ и $H$ -- две локально компактные группы. Каждой паре функций $u\in{\mathcal C}(G)$ и $v\in{\mathcal C}(H)$ поставим в соответствие функцию на декартовом произведении $G\times H$:
\beq\label{u-boxdot-v}
(u\boxdot v)(s,t)=u(s)\cdot v(t),\qquad s\in G,\ t\in H.
\eeq
Как известно \cite[Theorem 8.4]{Ak16}, существует единственное линейное непрерывное отображение
$$
\iota:{\mathcal C}(G)\odot {\mathcal C}(H)\to {\mathcal C}(G\times H)
$$
удовлетворяющее тождеству
$$
\iota(u\odot v)=u\boxdot v, \qquad u\in{\mathcal C}(G),\quad v\in{\mathcal C}(H).
$$
($\iota$ является изоморфизмом стереотипных алгебр). Обозначим
\beq\label{omega_G,H}
\omega_{G,H}=\eta_{{\mathcal K}^\star(G),{\mathcal K}^\star(H)}\circ \Env_{\mathcal C}(\env_{\mathcal C}{{\mathcal C}^\star(G)}\circledast \env_{\mathcal C}{{\mathcal C}^\star(H)})\circ \Env_{\mathcal C}(\iota^\star).
\eeq
Действие этого отображения иллюстрируется следующей диаграммой:
{\scriptsize
$$
\xymatrix @R=2.pc @C=4.0pc 
{
{\mathcal C}^\star(G\times H)\ar[rr]^{\env_{\mathcal C}{{\mathcal C}^\star(G\times H)}}\ar[dd]_{\iota^\star} & &
\Env_{\mathcal C}{\mathcal C}^\star(G\times H)\ar[dd]^{\Env_{\mathcal C}(\iota^\star)} & {\mathcal K}^\star(G\times H)\ar@{=}[l]
\ar@{-->}[dddd]^{\omega_{G,H}} \\ &&& \\
{\mathcal C}^\star(G)\circledast{\mathcal C}^\star(H)\ar[rr]^{\env_{\mathcal C}\big({\mathcal C}^\star(G)\circledast {\mathcal C}^\star(G)\big)}\ar[dd]|{\env_{\mathcal C}{{\mathcal C}^\star(G)}\circledast \env_{\mathcal C}{{\mathcal C}^\star(H)}} & & \Env_{\mathcal C}\big({\mathcal C}^\star(G)\circledast{\mathcal C}^\star(H)\big)\ar[dd]|{\Env_{\mathcal C}\big(\env_{\mathcal C}{{\mathcal C}^\star(G)}\circledast \env_{\mathcal C}{{\mathcal C}^\star(H)}\big)} & \\
&&& \\
\Env_{\mathcal C}{\mathcal C}^\star(G)\circledast \Env_{\mathcal C}{\mathcal C}^\star(H) \ar[rr]^(.45){\env_{\mathcal C}\big(\Env_{\mathcal C}{\mathcal C}^\star(G)\circledast \Env_{\mathcal C}{\mathcal C}^\star(H)\big)}
& & \Env_{\mathcal C}(\Env_{\mathcal C}{\mathcal C}^\star(G)\circledast \Env_{\mathcal C}{\mathcal C}^\star(H))\ar[r]^(.6){\eta_{{\mathcal K}^\star(G),{\mathcal K}^\star(H)}} & {\mathcal K}^\star(G)\odot {\mathcal K}^\star(H) \\
&& {\mathcal K}^\star(G)\overset{\mathcal C}{\circledast} {\mathcal K}^\star(H)\ar@{=}[u] & \\
{\mathcal K}^\star(G)\circledast {\mathcal K}^\star(H)\ar@{=}[uu]\ar@/_18ex/[uurrr]_(.8){\quad @_{{\mathcal K}^\star(G),{\mathcal K}^\star(H)}} & &  & \\
 &  &
}
$$
}
Рассмотрим сопряженное отображение:
$$
\xymatrix @R=2.pc @C=8.0pc 
{
{\mathcal K}(G\times H) & {\mathcal K}(G)\circledast{\mathcal K}(H)
\ar[l]_{\omega_{G,H}^\star}
}
$$

\btm\label{TH:omega^star(u-circledast-v)=u-boxdot-v}
Для любых двух локально компактных групп $G$ и $H$ справедливо тождество
\beq\label{omega^star(u-circledast-v)=u-boxdot-v}
\omega_{G,H}^\star(u\circledast v)=u\boxdot v,\qquad u\in{\mathcal K}(G),\quad v\in{\mathcal K}(H).
\eeq
\etm
\bpr
Рассмотрим сопряженную диаграмму
{\scriptsize
$$
\xymatrix @R=4.pc @C=4.0pc 
{
{\mathcal C}(G\times H) & &
\big(\Env_{\mathcal C}{\mathcal C}^\star(G\times H)\big)^\star\ar[ll]_{\big(\env_{\mathcal C}{\mathcal C}^\star(G\times H)\big)^\star} & {\mathcal K}(G\times H)\ar@{=}[l]
 \\
{\mathcal C}(G)\odot{\mathcal C}(H)\ar[u]^{\iota}
& & \Env_{\mathcal C}\big({\mathcal C}^\star(G)\circledast{\mathcal C}^\star(H)\big)^\star\ar[ll]_{\Big(\env_{\mathcal C}\big({\mathcal C}^\star(G)\circledast {\mathcal C}^\star(G)\big)\Big)^\star}\ar[u]_{\Env_{\mathcal C}(\iota^\star)^\star}
 & \\
(\Env_{\mathcal C}{\mathcal C}^\star(G))^\star\odot (\Env_{\mathcal C}{\mathcal C}^\star(H))^\star
 \ar[u]|{\big(\env_{\mathcal C}{{\mathcal C}^\star(G)}\circledast \env_{\mathcal C}{{\mathcal C}^\star(H)}\big)^\star}
& & \Env_{\mathcal C}(\Env_{\mathcal C}{\mathcal C}^\star(G)\circledast \Env_{\mathcal C}{\mathcal C}^\star(H))^\star\ar[u]|{\Env_{\mathcal C}\big(\env_{\mathcal C}{{\mathcal C}^\star(G)}\circledast \env_{\mathcal C}{{\mathcal C}^\star(H)}\big)^\star}\ar[ll]_{\env_{\mathcal C}{\big(\Env_{\mathcal C}{\mathcal C}^\star(G)\circledast \Env_{\mathcal C}{\mathcal C}^\star(H)\big)}^\star}
 & {\mathcal K}(G)\circledast {\mathcal K}(H)\ar[l]_(.35){\eta_{{\mathcal K}^\star(G),{\mathcal K}^\star(H)}^\star} \ar@{-->}[uu]_{\omega_{G,H}^\star}\ar@/^10ex/[dlll]^{@_{{\mathcal K}(G),{\mathcal K}(H)}} \\
{\mathcal K}(G)\odot {\mathcal K}(H)\ar@{=}[u] & &  & \\
 &  &
}
$$
}
Двигаясь из правого нижнего угла в левый верхний элементарный тензор $u\circledast v$, с одной стороны, проделает путь
$$
\xymatrix @R=2.pc @C=8.0pc 
 {
 u\circledast v\ar[r]^{@} & u\odot v\ar[r]^{(\env_{\mathcal C}{{\mathcal C}^\star(G)}\circledast \env_{\mathcal C}{{\mathcal C}^\star(H)})^\star} & u\odot v\ar[r]^{\iota}& u\boxdot v
 }
$$
а с другой - путь
$$
\xymatrix @R=2.pc @C=8.0pc 
 {
 u\circledast v\ar[r]^{\omega_{G,H}^\star} & \omega_{G,H}^\star(u\odot v)\ar[r]^(.35){\big(\env_{\mathcal C}{\mathcal C}^\star(G\times H)\big)^\star} & \big(\env_{\mathcal C}{\mathcal C}^\star(G\times H)\big)^\star(\omega_{G,H}^\star(u\odot v))
 }
$$
Отсюда
$$
\big(\env_{\mathcal C}{\mathcal C}^\star(G\times H)\big)^\star(\omega_{G,H}^\star(u\odot v))=u\boxdot v,
$$
и поскольку по теореме \ref{TH:K(G)->C(G)}, $\big(\env_{\mathcal C}{\mathcal C}^\star(G\times H)\big)^\star$ можно считать теоретико-множественным вложением,
$$
\omega_{G,H}^\star(u\odot v)=u\boxdot v.
$$
\epr

\btm\label{TH:K(A-times-K)-cong-K(A)-circledast-K(K)}
Если $C$ -- абелева локально компактная группа, а $K$ -- компактная группа, то отображение $\omega_{C,K}:{\mathcal K}(C)\circledast {\mathcal K}(K)\to {\mathcal K}(C\times K)$ является изоморфизмом:
\beq\label{K(A-times-K)-cong-K(A)-circledast-K(K)}
{\mathcal K}(C\times K)\cong{\mathcal K}(C)\circledast {\mathcal K}(K)
\eeq
\etm
\bpr
Это следует из \eqref{env_C^star(R^n-times-K)}:
\begin{multline*}
{\mathcal K}(C\times K)=\Big(\Env_{\mathcal C}{\mathcal C}^\star(C\times K)\Big)^\star\cong \eqref{env_C^star(R^n-times-K)}\cong \Big(\Env_{\mathcal C}{\mathcal C}^\star(C)\odot\Env_{\mathcal C}{\mathcal C}^\star(K) \Big)^\star\cong\\ \cong \Big(\Env_{\mathcal C}{\mathcal C}^\star(C)\Big)^\star\circledast\Big(\Env_{\mathcal C}{\mathcal C}^\star(K)\Big)^\star\cong{\mathcal K}(C)\circledast {\mathcal K}(K)
\end{multline*}
\epr

\btm\label{TH:u-boxdot-v-polno-v-K(G-times-R^n)}
Для любых двух групп Мура $G$ и $H$ функции вида $u\boxdot v$  (определенные в \eqref{u-boxdot-v}), где $u\in{\mathcal K}(G)$ и $v\in{\mathcal K}(H)$, полны в ${\mathcal K}(G\times H)$.
\etm
\bpr
Рассмотрим функцию $w$ вида \eqref{DEF:k(G)} на группе $G\times H$. Поскольку $G\times H$ -- тоже группа Мура, неприводимое представление $\pi:G\times H\to{\mathcal B}(X)$ должно быть конечномерным. Пусть $e_1,...,e_n$ -- ортонормированный базис в $X$, и пусть $\rho:G\to{\mathcal B}(X)$ и $\sigma:H\to{\mathcal B}(X)$ -- представления, действующие по формулам
$$
\rho(s)=\pi(s,1_H),\qquad \sigma(t)=\pi(1_G,t),\qquad s\in G,\quad t\in H.
$$
Тогда
\begin{multline*}
w(s,t)=\Big\langle\pi(s,t)x,y\Big\rangle=\Big\langle\pi\Big((s,1_H)\cdot(1_G,t)\Big)x,y\Big\rangle=
\Big\langle\Big(\pi(s,1_H)\cdot\pi(1_G,t)\Big)x,y\Big\rangle=\Big\langle\Big(\rho(s)\cdot\sigma(t)\Big)x,y\Big\rangle=\\=
\langle\sigma(t)x,\rho(s)^\bullet y\rangle=\sum_{i=1}^n \langle\sigma(t)x,e_i\rangle\cdot\langle e_i,\rho(s)^\bullet y\rangle=
\sum_{i=1}^n \langle\sigma(t)x,e_i\rangle\cdot\langle \rho(s)e_i,y\rangle=\sum_{i=1}^n (u_i\boxdot v_i)(s,t)
\end{multline*}
где
$$
u_i(s)=\langle \rho(s)e_i,y\rangle,\qquad v_i(t)=\langle\sigma(t)x,e_i\rangle.
$$
Из этого следует, что пространство $\Trig(G\times H)$ непрерывных по норме тригонометрических многочленов на $G$ содержится в линейной оболочке множества функций $u\boxdot v$, где $u\in{\mathcal K}(G)$ и $v\in{\mathcal K}(H)$. Значит,
$$
\overline{\sp\{u\boxdot v; \ u\in{\mathcal K}(G), \ v\in{\mathcal K}(H)\}}\supseteq\overline{\Trig(G\times H)}=\eqref{Trig-plotno-v-K}={\mathcal K}(G\times H).
$$
\epr

\paragraph{Спектр алгебры ${\mathcal K}(G)$.}

Вспомним, что мы определили группы Мура на с.\pageref{DEF:Moore-gruppa}, как те, у которых все унитарные неприводимые представления $\pi:G\to{\mathcal L}(X)$ конечномерны. По теореме \ref{TH:Moore->SIN}, всякая такая группа является SIN-группой, поэтому для нее справедливо представление \eqref{SIN-kak-rasshirenie},
$$
1\to \R^n\times K=N\to G\to D\to 1
$$
в котором $K$ -- компактная группа, а $D$ -- дискретная. По следствию \ref{TH:G=Moore->D=Moore} $D$ здесь тоже является группой Мура.

\bit{
 \item[$\bullet$] {\it Инволютивным характером} на инволютивной стереотипенгой алгебре $A$ над $\C$ мы называем произвольный (непрерывный, инволютивный и сохраняющий единицу) гомоморфизм $s:A\to \C$. Пространство всех инволютивных характеров на $A$ с топологией равномерной сходимости на вполне ограниченных множествах в $A$ мы называем {\it инволютивным спектром}\label{DEF:Spec_R(G)} (или просто {\it спектром}) алгебры $A$, и обозначаем $\Spec(A)$.
 }\eit

Перед следующей важной теоремой \ref{TH:Spec-K(G)=G} нам понадобится несколько лемм.

\blm\label{LM:Spec-K(G)-dlya-diskretnoi-gruppy}
Если $G$ -- аменабельная дискретная группа, то отображение спектров $G\to \Spec{\mathcal K}(G)$ является биекцией.
\elm
\bpr Здесь используются идеи \cite[Theorem 6.3]{Kuznetsova}.
По предложению \ref{PROP:nepr-obol-diskr-gruppy}, непрерывной оболочкой групповой алгебры для $G$ в данном случае будет $C^*$-алгебра этой группы:
$$
\Env_{\mathcal C}{\mathcal C}^\star(G)=C^*(G)
$$
Ее сопряженное пространство (в банаховом смысле) представляет собой классическую алгебру Фурье-Стилтьеса $B(G)$, поэтому алгебра ${\mathcal K}(G)$ совпадает как множество с $B(G)$ \cite{Eymard,Renault}:
$$
{\mathcal K}(G)=\Env_{\mathcal C}{\mathcal C}^\star(G)^\star=C^*(G)^\star=B(G).
$$
(это равенство векторных пространств, но топология на $B(G)$ сильнее топологии на ${\mathcal K}(G)$). Алгебра $B(G)$ содержит идеал $A(G)$, называемый алгеброй Фурье \cite{Eymard,Renault}:
$$
{\mathcal K}(G)=B(G)\supseteq A(G).
$$
Из аменабельности группы $G$ следует, что аннулятор алгебры $A(G)$ в пространстве  $C^*(G)$ нулевой \cite{Hulanicki}, \cite[Theorem 4.21]{Paterson}:
$$
A(G)^\perp=0
$$
Как следствие $A(G)$ плотно в ${\mathcal K}(G)$:
$$
\overline{A(G)}={\mathcal K}(G).
$$
Если $\chi:{\mathcal K}(G)\to\C$ -- инволютивный характер на ${\mathcal K}(G)$, то он будет инволютивным характером и на $A(G)$, и по теореме Эймара \cite{Eymard}, на функциях $u\in A(G)$ он представляет собой дельта-функцию в некоторой точке $a\in G$:
$$
\chi(u)=u(a),\qquad A(G).
$$
Поскольку $A(G)$ плотно в ${\mathcal K}(G)$, это верно и для функций $u\in{\mathcal K}(G)$.
\epr

\blm\label{LM:Spec-K(G)-dlya-kompaktnoi-gruppy}
Если $G$ -- компактная группа, то отображение спектров $G\to \Spec{\mathcal K}(G)$ является гомеомеорфизмом.
\elm
\bpr
Всякий инволютивный характер $\chi:{\mathcal K}(G)\to\C$ является инволютивным характером на алгебре $\Trig(G)$ тригонометрических многочленов на $G$, поэтому \cite[30.5]{Hewitt-Ross-2} существует точка $a\in G$ такая, что
$$
\chi(u)=u(a),\qquad u\in\Trig(G).
$$
При этом в силу \eqref{Trig-plotno-v-K} алгебра $\Trig(G)$ плотна в алгебре ${\mathcal K}(G)$. Значит,
$$
\chi(u)=u(a),\qquad u\in{\mathcal K}(G).
$$
Мы получаем, что отображение $G\to \Spec{\mathcal K}(G)$ сюръективно. С другой стороны, любые две точки $s\ne t$ на компактной группе $G$ разделяются тригонометрическим многочленом $u\in\Trig(G)$, и это означает, что отображение $G\to \Spec{\mathcal K}(G)$ инъективно. Таким образом, оно биективно и непрерывно. При этом $G$ -- компакт. Значит, оно -- гомеоморфизм.
\epr

\blm\label{LM:Spec-K(G)-dlya-kommut-gruppy}
Если $G$ -- абелева локально компактная группа, то отображение спектров $G\to \Spec{\mathcal K}(G)$ является гомеомеорфизмом.
\elm
\bpr
По предложению \ref{PROP:Env_C-C^star(G)=C(widehat(G))}, $\Env_{\mathcal C}{\mathcal C}^\star(G)={\mathcal C}(\widehat{G})$.
Поэтому
$$
{\mathcal K}(G)={\mathcal C}^\star(\widehat{G}).
$$
Всякий характер $f:{\mathcal C}^\star(\widehat{G})\to\C$ в композиции с отображением $\delta:\widehat{G}\to {\mathcal C}^\star(\widehat{G})$ дает комплексный характер $f\circ\delta:\widehat{G}\to\C^\times$, а если $f$ -- инволютивный характер, то $f\circ\delta$ действует в окружность $\T$. То есть $f\circ\delta$ -- характер на группе $\widehat{G}$, и поэтому
$$
(f\circ\delta)(\chi)=f(\delta^\chi)=\chi(a),\qquad \chi\in \widehat{G}
$$
для некоторой точки $a\in G$. Поскольку линейные комбинации дельта-функций $\delta^\chi$ полны в ${\mathcal C}^\star(\widehat{G})$, это равенство продолжается на все элементы ${\mathcal C}^\star(\widehat{G})$:
$$
f(u)=u(a),\qquad u\in{\mathcal K}(G)={\mathcal C}^\star(\widehat{G})
$$
Мы получаем, что отображение $G\to \Spec{\mathcal K}(G)=\Spec{\mathcal C}^\star(\widehat{G})$ является сюръекцией. С другой стороны, любые две точки $s\ne t$ на $G$, то есть два характера на $\widehat{G}$ различаются каким-то дельта-функционалом $\delta_\chi\in{\mathcal C}^\star(\widehat{G})$, $\chi\in \widehat{G}$,
$$
\delta^\chi(s)=\chi(s)\ne\chi(t)=\delta^\chi(t)
$$
и это означает, что отображение спектров $G\to \Spec{\mathcal K}(G)=\Spec{\mathcal C}^\star(\widehat{G})$ является инъекцией. Итак, оно биективно и непрерывно. Нам остается убедиться, что оно открыто (непрерывно в обратную сторону). Пусть $f_i\to f$ в $\Spec{\mathcal C}^\star(\widehat{G})$, и $a_i$, $a$ -- соответствующие точки в $G$. Если $K$  -- компакт в $\widehat{G}$, то множество $\{\delta^\chi;\ \chi\in K\}$ будет компактом в ${\mathcal C}^\star(\widehat{G})$, поэтому $f_i(\delta^\chi)$ будет стремиться к $f(\delta^\chi)$ равномерно по $\chi\in K$:
$$
\chi(a_i)=f_i(\delta^\chi)\underset{i\to\infty}{\underset{\chi\in K}{\rightrightarrows}}f(\delta^\chi)=\chi(a).
$$
Это верно для всякого компакта $K\subseteq \widehat{G}$, поэтому $a_i\to a$ в $G$.
\epr

\blm\label{LM:Spec(R^n-times-K)}
Для всякой компактной группы $K$ и любого $n\in\N$  отображение спектров $\R^n\times K\to \Spec{\mathcal K}(\R^n\times K)$ является гомеомеорфизмом.
\elm
\bpr
Пусть $\chi:{\mathcal K}(\R^n\times K)\to\C$ -- инволютивный характер. Рассмотрим вложения
$$
\rho:{\mathcal K}(\R^n)\to {\mathcal K}(\R^n\times K),\qquad \rho(u)=u\boxdot 1_K
$$
и
$$
\sigma:{\mathcal K}(K)\to {\mathcal K}(\R^n\times K),\qquad \rho(v)=1_{\R^n}\boxdot v.
$$
Композиция $\chi\circ\rho$ является инволютивным характером на алгебре ${\mathcal K}(\R^n)$, спектр которой по лемме \ref{LM:Spec-K(G)-dlya-kommut-gruppy} совпадает с $\R^n$. Поэтому существует точка $a\in\R^n$ такая, что
$$
\chi(u\boxdot 1_K)=(\chi\circ\rho)(u)=u(a), \qquad u\in {\mathcal K}(\R^n)
$$
С другой стороны, композиция $\chi\circ\sigma$ является инволютивным характером на алгебре ${\mathcal K}(K)$, содержащей алгебру $\Trig(G)$, поэтому по лемме \ref{LM:Spec-K(G)-dlya-kompaktnoi-gruppy} существует точка $b\in K$ такая, что
$$
\chi(1_{\R^n}\boxdot v)=(\chi\circ\sigma)(v)=v(b), \qquad v\in {\mathcal K}(K).
$$
Теперь для функций вида $u\boxdot v$ мы получим
$$
\chi(u\boxdot v)=\chi(u\boxdot 1_K\cdot 1_{\R^n}\boxdot v)=\chi(u\boxdot 1_K)\cdot\chi(1_{\R^n}\boxdot v)=u(a)\cdot v(b)=(u\boxdot v)(a,b).
$$
Это тождество продолжается на линейные комбинации функций вида $u\boxdot v$, а затем по теореме \ref{TH:u-boxdot-v-polno-v-K(G-times-R^n)}, на все пространство ${\mathcal K}(\R^n\times K)$. Это доказывает сюръективность отображения $\R^n\times K\to \Spec{\mathcal K}(\R^n\times K)$.

Докажем инъективность: пусть $(s,a)\ne(t,b)$ в $\R^n\times K$. Тогда либо $s\ne t$, либо $a\ne b$. В первом случае можно подобрать функцию $u\in{\mathcal K}(\R^n)$ так, чтобы $u(s)\ne u(t)$ (в этот момент мы используем лемму \ref{LM:Spec-K(G)-dlya-kommut-gruppy}), и тогда мы получим, что функция $u\boxdot 1$ (лежащая в ${\mathcal K}(\R^n\times K)$ по теореме \ref{TH:omega^star(u-circledast-v)=u-boxdot-v}) разделяет точки $(s,a)$ и $(t,b)$:
$$
(u\boxdot 1)(s,a)=u(s)\cdot 1=u(s)\ne u(t)=u(t)\cdot 1=(u\boxdot 1)(t,b).
$$
А во втором случае, когда $a\ne b$, то же самое можно сделать, подобрав (с помощью леммы  \ref{LM:Spec-K(G)-dlya-kompaktnoi-gruppy}) функцию из ${\mathcal K}(\R^n)$, разделяющую эти точки.

Мы получаем, что  отображение спектров $\R^n\times K\to \Spec{\mathcal K}(\R^n\times K)$ биективно (и непрерывно).

Докажем его открытость (непрерывность в обратную сторону). Пусть $(s_i,a_i)\to (s,a)$ в $\Spec{\mathcal K}(\R^n\times K)$. Тогда выбрав компакт $T\subseteq {\mathcal K}(\R^n)$ мы можем рассмотреть множество $\rho(T)\subseteq{\mathcal K}(\R^n\times K)$, также являющееся компактом, и для него мы получим равномерную сходимость по $u\in T$:
$$
u(s_i)=u(s_i)\cdot 1= (u\boxdot 1)(s_i,a_i)=\rho(u)(s_i,a_i)\underset{i\to\infty}{\underset{u\in T}{\rightrightarrows}}\rho(u)(s,a)=(u\boxdot 1)(s,a)=u(s)\cdot 1=u(s).
$$
Поскольку это верно для всякого компакта $T\subseteq {\mathcal K}(\R^n)$, мы получаем, что $s_i\to s$ в $\Spec{\mathcal K}(\R^n)$, а по лемме \ref{LM:Spec-K(G)-dlya-kommut-gruppy} это означает, что $s_i\to s$ в $\R^n$.

И точно так же (с помощью леммы  \ref{LM:Spec-K(G)-dlya-kompaktnoi-gruppy}) доказывается, что $a_i\to a$ в $K$.
\epr

\btm\label{TH:Spec-K(G)=G}
Если $G$ -- группа Мура, то инволютивный спектр алгебры ${\mathcal K}(G)$ топологически изоморфен $G$:
\beq\label{Spec-K(G)=G}
\Spec{\mathcal K}(G)=G
\eeq
\etm
\bpr 1. Сначала покажем, что отображение спектров $G\to \Spec{\mathcal K}(G)$ является сюръекцией. Пусть $\chi:{\mathcal K}(G)\to\C$ -- инволютивный характер. Гомоморфизм $G\to D$ из \eqref{SIN-kak-rasshirenie} порождает гомоморфизм ${\mathcal C}^\star(G)\to {\mathcal C}^\star(D)$, который затем порождает гомоморфизм $\Env_{\mathcal C} {\mathcal C}^\star(G)\to \Env_{\mathcal C} {\mathcal C}^\star(D)$, а тот, в свою очередь,
гомоморфизм ${\mathcal K}(G)\gets {\mathcal K}(D)$. Обозначим его $\ph:{\mathcal K}(D)\to {\mathcal K}(G)$. Композиция $\chi\circ\ph:{\mathcal K}(D)\to \C$ есть инволютивный непрерывный характер на ${\mathcal K}(D)$, причем $D$ здесь является группой Мура (в силу следствия  \ref{TH:G=Moore->D=Moore}), и значит, аменабельной группой (по теореме \ref{TH:Moore->AM}). Поэтому по лемме \ref{LM:Spec-K(G)-dlya-diskretnoi-gruppy} $\chi$ является дельта-функцией:
$$
(\chi\circ\ph)(u)=u(L),\qquad u\in {\mathcal K}(D),
$$
для некоторого $L\in G/N$. Рассмотрим пространство ${\mathcal K}_L(G)$ из \eqref{DEF:K_L(G)} и обозначим $\rho_L$ его вложение в ${\mathcal K}(G)$. Кроме того обозначим символом $\sigma$ вложение ${\mathcal K}(N)\to {\mathcal K}_N(G)$, то есть изоморфизм, определяемый леммой \ref{LM:K(N)=K_N(G)}(ii). Пусть далее $b\in L$, то есть $L=N\cdot b$, и $\tau_b:{\mathcal K}(G)\to{\mathcal K}(G)$ -- оператор сдвига на элемент $b^{-1}$ (действующий на ${\mathcal K}(G)$ по теореме \ref{LM:sdvig-v-K(G)}):
$$
\tau_bu=b^{-1}\cdot u,\qquad {\mathcal K}(G).
$$
Он переводит пространство ${\mathcal K}(N)$ в пространство ${\mathcal K}_L(G)$, поэтому определено отображение $\sigma_L=\tau_b\circ\sigma:{\mathcal K}(N)\to {\mathcal K}_L(G)$. Положив теперь
$$
\chi_L=\chi\circ\rho_L,\qquad \chi_N=\chi_L\circ\sigma_L
$$
и обозначив символом $\rho_L:{\mathcal K}_L(G)\to {\mathcal K}(G)$ естественное вложение, мы получим коммутативную диаграмму
$$
\xymatrix @R=2.pc @C=3.0pc 
{
{\mathcal K}(N)\ar[r]^{\sigma_L}\ar[dr]_{\chi_N} & {\mathcal K}_L(G)\ar[d]^{\chi_L}\ar[r]^{\rho_L} & {\mathcal K}(G)\ar[dl]^{\chi}\\
 & \C &
}
$$
Поскольку $\chi_N$ -- характер на ${\mathcal K}(N)$ по лемме \ref{LM:Spec(R^n-times-K)} он должен быть дельта-функцией:
\beq\label{chi_N(u)=u(a)}
\chi_N(u)=u(a),\qquad u\in{\mathcal K}(N)
\eeq
для некоторого $a\in N$. Тогда
\begin{multline*}
\chi(u)=\chi(1_L)\cdot\chi(u)=\chi(1_L\cdot u)=\chi_L(1_L\cdot u)=\chi_N(\sigma_L^{-1}(1_L\cdot u))=\eqref{chi_N(u)=u(a)}=
\sigma_L^{-1}(1_L\cdot u)(a)=\sigma(\sigma_L^{-1}(1_L\cdot u))(a)=\\=
(\sigma\circ\sigma_L^{-1})(1_L\cdot u))(a)=(\sigma\circ(\tau_b\circ\sigma)^{-1})(1_L\cdot u))(a)=
(\sigma\circ\sigma^{-1}\circ\tau_b^{-1})(1_L\cdot u))(a)=(\sigma\circ\sigma^{-1}\circ\tau_b^{-1})(1_L\cdot u))(a)=\\=
\tau_{b^{-1}}(1_L\cdot u))(a)=(b\cdot(1_L\cdot u))(a)=(1_L\cdot u)(\underbrace{a\cdot b}_{\scriptsize\begin{matrix}\text{\rotatebox{90}{$\owns$}}\\ L\end{matrix}})=u(a\cdot b)=\delta^{a\cdot b}(u).
\end{multline*}

2. Теперь убедимся, что отображение спектров $G\to \Spec{\mathcal K}(G)$ является инъекцией. Пусть $a\ne b\in G$. Если $a\cdot b^{-1}\notin N$, то есть $a\notin b\cdot N$, то характеристическая функция $1_L\in{\mathcal K}(G)$ класса $L=b\cdot N$ из леммы \ref{LM:1_L-in-K(G)} будет различать $a$ и $b$:
$$
1_L(a)=0\ne 1=1_L(b).
$$
Пусть $a\in b\cdot N$, то есть $a\cdot b^{-1}\in N$. Тогда по лемме \ref{LM:Spec(R^n-times-K)} мы можем подобрать функцию $u\in{\mathcal K}(N)$ так, чтобы
$$
u(a\cdot b^{-1})\ne u(1)
$$
Затем по лемме \ref{LM:K(N)=K_N(G)}, найдется функция $v\in{\mathcal K}_N(G)$ такая, что $u\big|_N=v\big|_N$, и поэтому
$$
v(a\cdot b^{-1})\ne v(1)
$$
Далее по теореме \ref{LM:sdvig-v-K(G)} сдвиг $b^{-1}\cdot v$ будет снова лежать в ${\mathcal K}(G)$, и для этой функции мы получим
$$
(b^{-1}\cdot v)(a)=v(a\cdot b^{-1})\ne v(1)=v(b\cdot b^{-1})=(b^{-1}\cdot v)(b).
$$

3. Остается проверить открытость отображения $G\to \Spec{\mathcal K}(G)$. Пусть $a_i\to a$ в $\Spec{\mathcal K}(G)$. Из теоремы \ref{LM:sdvig-v-K(G)} сразу следует, что $a_i\cdot a^{-1}\to 1$ в $\Spec{\mathcal K}(G)$. Для характеристической функции $1_N\in{\mathcal K}(G)$ подгруппы $N$ мы получаем $1_N(a_i\cdot a^{-1})\to 1_N(1)=1$, поэтому начиная с некоторого индекса все $a_i\cdot a^{-1}$ должны лежать в $N$. Возьмем компакт $S\subseteq {\mathcal K}(N)$. По лемме \ref{LM:K(N)=K_N(G)} для него можно подобрать компакт $T\subseteq {\mathcal K}(G)$, состоящий из функций, у которых ограничения на $N$ лежат в $S$, и при этом получится биекция между $T$ и $S$. Поскольку $a_i\cdot a^{-1}\to 1$ в $\Spec{\mathcal K}(G)$, мы получаем
$$
v(a_i\cdot a^{-1})\underset{i\to\infty}{\underset{v\in T}{\rightrightarrows}}v(1)
$$
и это эквивалентно
$$
u(a_i\cdot a^{-1})\underset{i\to\infty}{\underset{u\in S}{\rightrightarrows}}u(1).
$$
Это верно для всякого компакта $S\subseteq {\mathcal K}(N)$, значит $a_i\cdot a^{-1}\to 1$ в $\Spec{\mathcal K}(N)$. Но в лемме \ref{LM:Spec(R^n-times-K)} мы уже доказали, что $\Spec{\mathcal K}(N)=N$, поэтому мы получаем, что $a_i\cdot a^{-1}\to 1$ в $N$, и значит, в $G$.
\epr

\paragraph{Непрерывная оболочка алгебры ${\mathcal K}(G)$.}

Из теорем \ref{TH:Spec-K(G)=G} и \ref{C-obolochka-podalgebry-v-C(M)} сразу следует

\btm\label{TH:E(K(G))=C(G)}
Если $G$ -- группа Мура, то непрерывной оболочкой алгебры ${\mathcal K}(G)$ является алгебра ${\mathcal C}(G)$:
\beq\label{E(K(G))=C(G)}
\Env_{\mathcal C} {\mathcal K}(G)={\mathcal C}(G)
\eeq
\etm

\bcor\label{COR:K(G)-plotno-v-C(G)}
Если $G$ -- группа Мура, то алгебра ${\mathcal K}(G)$ плотна в алгебре ${\mathcal C}(G)$:
\beq\label{K(G)-plotno-v-C(G)}
\overline{{\mathcal K}(G)}={\mathcal C}(G)
\eeq
\ecor
\bpr
Это формально следует уже из теоремы \ref{TH:Spec-K(G)=G} (более того, из ее ослабленного варианта, в котором утверждается только, что спектр ${\mathcal K}(G)$ биективен $G$) и теоремы Стоуна---Вейерштрасса, но можно также это вывести из теоремы \ref{TH:E(K(G))=C(G)}. По определению, непрерывная оболочка $\env_{\mathcal C}$ является плотным эпиморфизмом, поэтому из равенства \eqref{E(K(G))=C(G)} следует, что ${\mathcal K}(G)$ плотно в ${\mathcal C}(G)$.
\epr

\subsection{Усреднение по малым подгруппам в $\Env_{\mathcal C}{\mathcal C}^\star(G)$}

\paragraph{Элементы $\mu_K^{\mathcal C}$ в $\Env_{\mathcal C}{\mathcal C}^\star(G)$.}

Вспомним меры $\mu_K\in{\mathcal C}^\star(G)$, определенные в \eqref{mu_K(u)=int_K-u(t)mu_K(dt)}. Для всякой подгруппы $K\in\lambda(G)$ мера $\mu_K\in{\mathcal C}^\star(G)$ под действием отображения оболочки
$$
\env_{\mathcal C}:{\mathcal C}^\star(G)\to \Env_{\mathcal C}{\mathcal C}^\star(G)
$$
превращается в некий элемент
$$
\mu_K^{\mathcal C}=\env_{\mathcal C}(\mu_K)\in \Env_{\mathcal C}{\mathcal C}^\star(G).
$$

\btm\label{TH:mu_K^C->delta^1}
Направленность $\{\mu_K^{\mathcal C}:\ K\in\lambda(G)\}$ вполне ограничена\footnote{В доказательстве этой теоремы, а также в теоремах \ref{TH:pi_K-vpolne-ogr-v-L(Env-C*(G))} и \ref{TH:Env_C-C(G)-approx}, где этот результат используется, условие полной ограниченности семейства $\{\mu_K^{\mathcal C}; \ K\in\lambda(G)\}$ несущественно.} в алгебре $\Env_{\mathcal C}{\mathcal C}^\star(G)$ и стремится к ее единице:
\beq\label{mu_K^C->delta^1}
\mu_K^{\mathcal C}\overset{\Env_{\mathcal C}{\mathcal C}^\star(G)}{\underset{K\to 0}{\longrightarrow}} 
1_{\Env_{\mathcal C}{\mathcal C}^\star(G)}
\eeq
\etm
\bpr
По теореме \ref{TH:mu_K->delta^1}, направленность мер $\{\mu_K;\ K\in\lambda(G)\}$ вполне ограничена в ${\mathcal C}^\star(G)$ и стремится к единице этой алгебры:
$$
\mu_K\overset{{\mathcal C}^\star(G)}{\underset{K\to 0}{\longrightarrow}}\delta^1=1_{{\mathcal C}^\star(G)}
$$
Поэтому под действием гомоморфизма $\env_{\mathcal C}:{\mathcal C}^\star(G)\to \Env_{\mathcal C}{\mathcal C}^\star(G)$
она превращается в направленность элементов $\{\mu_K^{\mathcal C};\ K\in\lambda(G)\}$, которая вполне ограничена в $\Env_{\mathcal C}{\mathcal C}^\star(G)$ и стремится к единице этой алгебры:
$$
\mu_K^{\mathcal C}\overset{\Env_{\mathcal C}{\mathcal C}^\star(G)}{\underset{K\to 0}{\longrightarrow}}\env_{\mathcal C}(\delta^1)=1_{\Env_{\mathcal C}{\mathcal C}^\star(G)}
$$
\epr

\paragraph{Операторы $\Env_{\mathcal C}\pi_K^\star$ в $\Env_{\mathcal C}{\mathcal C}^\star(G)$.}

Вспомним проекторы $\pi_K^\star:{\mathcal C}^\star(G)\to {\mathcal C}^\star(G)$, определенные в \eqref{pi_K^star(alpha)=mu_K*alpha}.
Будучи гомоморфизмами стереотипных алгебр, они порождают гомоморфизмы оболочек
$$
\Env_{\mathcal C}\pi_K^\star:\Env_{\mathcal C}{\mathcal C}^\star(G)\to \Env_{\mathcal C}{\mathcal C}^\star(G)
$$

\btm\label{TH:pi_K-mu_K^C}
Справедливо тождество
\beq\label{pi_K-mu_K^C}
\Env_{\mathcal C}\pi_K^\star(\alpha)=\mu_K^{\mathcal C}*\alpha,\qquad \alpha\in\Env_{\mathcal C}{\mathcal C}^\star(G) 
\eeq
(в котором $*$ --- операция умножения в алгебре $\Env_{\mathcal C}{\mathcal C}^\star(G)$).
\etm
\bpr
Тождество \eqref{pi_K^star(alpha)=mu_K*alpha} под действием гомоморфизма $\env_{\mathcal C}:{\mathcal C}^\star(G)\to \Env_{\mathcal C}{\mathcal C}^\star(G)$ дает цепочку
$$
(\Env_{\mathcal C}\pi_K^\star)(\env_{\mathcal C}\alpha)=\env_{\mathcal C}(\pi_K^\star(\alpha))=\env_{\mathcal C}\mu_K*\env_{\mathcal C}\alpha=\mu_K^{\mathcal C}*\env_{\mathcal C}\alpha,\qquad \alpha\in {\mathcal C}^\star(G) 
$$
и поскольку элементы $\env_{\mathcal C}\alpha$ плотны в $\Env_{\mathcal C}{\mathcal C}^\star(G)$, мы можем заменить их на произвольные элементы $\alpha\in\Env_{\mathcal C}{\mathcal C}^\star(G)$
$$
(\Env_{\mathcal C}\pi_K^\star)(\alpha)=\mu_K^{\mathcal C}*\alpha,\qquad \alpha\in \Env_{\mathcal C}{\mathcal C}^\star(G) 
$$
\epr

\btm\label{TH:pi_K-vpolne-ogr-v-L(Env-C*(G))} 
Пусть $G$ --- локально компактная группа. Тогда направленность операторов $\{\Env_{\mathcal C}\pi_K^\star; \ K\in\lambda(G)\}$ вполне ограничена\footnote{В доказательстве этой теоремы, а также в теореме \ref{TH:Env_C-C(G)-approx}, где этот результат используется, условие полной ограниченности семейства $\{\Env_{\mathcal C}\pi_K^\star; \ K\in\lambda(G)\}$ несущественно.} в ${\mathcal L}(\Env_{\mathcal C}{\mathcal C}^\star(G))$ и стремится в этом пространстве к тождественному оператору
\beq\label{pi_K->id-v-L(Env-C*(G))}
\Env_{\mathcal C}\pi_K^\star\overset{{\mathcal L}(\Env_{\mathcal C}{\mathcal C}^\star(G))}{\underset{K\to 0}{\longrightarrow}}\id_{{\mathcal L}(\Env_{\mathcal C}{\mathcal C}^\star(G))}
\eeq
\etm
\bpr
Здесь используется теорема о представлении \cite[Theorem 5.3.2]{Akbarov-De-Gruyter-I}. Рассмотрим стереотипную алгебру $A=\Env_{\mathcal C}{\mathcal C}^\star(G)$ (в которой умножение мы обозначаем звездочкой $*$). Она является левым стереотипным модулем над самой собой. Поэтому в силу \cite[Theorem 5.3.2]{Akbarov-De-Gruyter-I} ее действие умножением на самой себе можно представить как гомоморфизм стереотипных алгебр
$$
\ph:A\to{\mathcal L}(A)\quad\Big|\quad \ph(a)(b)=a*b,\quad a,b\in A.
$$
Отсюда следует, что направленность $a_K=\mu_K^{\mathcal C}$ в $A$, стремящаяся по теореме \ref{TH:mu_K^C->delta^1} к единице в $A$, 
$$
a_K=\mu_K^{\mathcal C}\overset{A}{\underset{K\to 0}{\longrightarrow}} 
1_A
$$
под действием гомоморфизма $\ph$ превратится в направленность в ${\mathcal L}(A)$, стремящуюся к единице в ${\mathcal L}(A)$: 
$$
\ph(a_K)=\ph(\mu_K^{\mathcal C})\overset{A}{\underset{K\to 0}{\longrightarrow}} 
1_{{\mathcal L}(A)}
$$
Но операторы $\ph(a_K)=\ph(\mu_K^{\mathcal C})$ --- это в точности операторы $\Env_{\mathcal C}\pi_K^\star$:
$$
\ph(\mu_K^{\mathcal C})(b)=\mu_K^{\mathcal C}*b=\eqref{pi_K-mu_K^C}=\Env_{\mathcal C}\pi_K^\star(b).
$$
Значит, 
$$
\Env_{\mathcal C}\pi_K^\star \overset{A}{\underset{K\to 0}{\longrightarrow}} 1_{{\mathcal L}(A)}.
$$
\epr

\paragraph{Аппроксимация в $\Env_{\mathcal C}{\mathcal C}^\star(G)$.}

Группа Мура не обязана быть компактно порожденной, например, существуют абелевы не компактно порожденные группы. Но если все же группа обладает этим свойством, для нее справедлива следующая теорема.

\btm\label{TH:Env_C-C(G)-approx} 
Для всякой компактно порожденной группы Мура $G$ пространство $\Env_{\mathcal C}{\mathcal C}^\star(G)$ обладает стереотипной аппроксимацией.
\etm 

\bpr
1. Сначала рассмотрим случай, когда $G$ --- (компактно порожденная) группа Ли, причем компактная надстройка абелевой группы. Тогда в ее представлении $G=Z\cdot K$ абелева группа $Z$ компактно порождена, как замкнутая подгруппа компактно порожденной группы $G$. Как следствие, двойственная группа $\widehat{Z}$ -- группа Ли. Поэтому в представлении \eqref{Env_C-C*(Z-cdot-K)} пространства $M_\sigma$ -- метризуемые. Отсюда по теореме \ref{TH:C(M)-obl-approx-kogda-M-metricheskoe} все пространства 
$$
{\mathcal C}\big(M_\sigma,{\mathcal B}(X_\sigma)\big)\cong {\mathcal C}(M_\sigma)^{(\dim X_\sigma)^2}
$$ 
обладают стереотипной аппроксимацией. Значит, их произведение
$$
\Env_{\mathcal C} {\mathcal C}^\star(Z\cdot K)\cong
\prod_{\sigma\in\widehat{K}}{\mathcal C}\big(M_\sigma,{\mathcal B}(X_\sigma)\big),
$$
--- тоже обладает стереотипной аппроксимацией.

2. Теперь пусть $G$ --- произвольная (компактно порожденная) группа Ли---Мура. Тогда по теореме \ref{TH:Lie-Moore}, $G$ --- конечное расширение некоторой компактной надстройки абелевой группы $N=Z\cdot K$:
$$
1\to N\to G\to F\to 1.
$$
Согласно \eqref{Env_C-C*(konech-rasshirenie-ZK)}, пространство $\Env_{\mathcal C}  {\mathcal C}^\star(G)$ раскладывается в (конечную) прямую сумму:
\beq\label{Env_C-C*(konech-rasshirenie-ZK)-1}
\Env_{\mathcal C}  {\mathcal C}^\star(G)=\bigoplus_{L\in G/N} E_L.
\eeq
В нем все пространства $E_L$ изоморфны $E_N$ и значит, в силу \eqref{E_N=Env_C^star(N)}, пространству $\Env_{\mathcal C}  {\mathcal C}^\star(N)$:
$$
E_L\cong \Env_{\mathcal C}  {\mathcal C}^\star(N).
$$
Поэтому, по уже доказанному в предыдущем пункте, они обладают стереотипной аппроксимацией. Значит, пространство $\Env_{\mathcal C}  {\mathcal C}^\star(G)$ также обладает стереотипной аппроксимацией.

3. Наконец, пусть $G$ --- произвольная (компактно порожденная) группа Мура. Для всякой подгруппы $K\in\lambda(G)$ фактор-группа $G/K$ является (компактно порожденной) группой Ли---Мура, поэтому по уже доказанному, пространство $\Env_{\mathcal C}{\mathcal C}^\star(G:K)\cong \Env_{\mathcal C}{\mathcal C}^\star(G/K)$ обладает стереотипной аппроксимацией. Пусть 
$$
\ph_K^i\in {\mathcal L}\l\Env_{\mathcal C}  {\mathcal C}^\star(G:K)\r
$$ 
--- направленность конечномерных операторов, аппроксимирующих тождественный оператор:
\beq\label{ph_K^i->id-v-L(C'(G:K))}
\ph_K^i\overset{{\mathcal L}(\Env_{\mathcal C}  {\mathcal C}^\star (G))}{\underset{i\to\infty}{\longrightarrow}}\id_{{\mathcal L}(\Env_{\mathcal C}  {\mathcal C}^\star (G:K))}
\eeq
Тогда из \eqref{ph_K^i->id-v-L(C'(G:K))} и \eqref{pi_K->id-v-L(Env-C*(G))} мы получим:
\begin{multline*}
\im\Env_{\mathcal C}\pi_K^\star\circ\ph_K^i\circ\coim\Env_{\mathcal C}\pi_K^\star=\\
=\im\rho_K^\star\circ\ph_K^i\circ\coim\rho_K^\star\overset{{\mathcal L}(\Env_{\mathcal C}  {\mathcal C}^\star (G))}{\underset{i\to\infty}{\longrightarrow}}\im\rho_K^\star\circ\id_{{\mathcal L}(\Env_{\mathcal C}  {\mathcal C}^\star (G:K))}\circ\coim\rho_K^\star=\\= \rho_K^\star\overset{{\mathcal L}(\Env_{\mathcal C}  {\mathcal C}^\star(G))}{\underset{K\to 0}{\longrightarrow}}\id_{{\mathcal L}({\mathcal K}^\star(G))}=\id_{{\mathcal L}(\Env_{\mathcal C}{\mathcal C}^\star(G))}
\end{multline*}
$$
 \xymatrix  @R=2.5pc @C=2.5pc
 {
 \Env_{\mathcal C}{\mathcal C}^\star(G)\ar[r]^{\id_{{\mathcal L}(\Env_{\mathcal C}{\mathcal C}^\star(G))}}\ar@{-->}[d]_{\coim\Env_{\mathcal C}\pi_K^\star}& \Env_{\mathcal C} {\mathcal C}^\star(G)\\
 \Env_{\mathcal C}{\mathcal C}^\star(G:K)\ar@{-->}[r]_{\ph_K^i} & \Env_{\mathcal C}{\mathcal C}^\star(G:K)\ar@{-->}[u]_{\im\Env_{\mathcal C}\pi_K^\star}
 }
$$
(эта диаграмма не коммутативна, она только показывает, на каких пространствах действуют рассматриваемые операторы).
\epr

\subsection{Преобразования Гротендика алгебр ${\mathcal C}^\star(G)$ и $\Env_{\mathcal C}{\mathcal C}^\star(G)$}

Покажем, что алгебры ${\mathcal C}^\star(G)$ и $\Env_{\mathcal C}{\mathcal C}^\star(G)$ обладают свойствами (iii) и (iv) из определения рефлексивной алгебры Хопфа на с.\pageref{@:H-circledast-H->H-odot-H-Mono}.

\paragraph{Преобразование Гротендика алгебры ${\mathcal C}^\star(G)$.}

\btm\label{TH:@:C^star(G)-circledast-C^star(G)->C^star(G)-odot-C^star(G)}
Преобразование Гротендика $@:{\mathcal C}^\star(G)\circledast {\mathcal C}^\star(G)\to {\mathcal C}^\star(G)\odot {\mathcal C}^\star(G)$ является биморфизмом (и, в частности, мономорфизмом).
\etm
\bpr
Это следует из теоремы \ref{TH:C(G)-approx}.
\epr

\paragraph{Преобразование Гротендика алгебры $\Env_{\mathcal C}{\mathcal C}^\star(G)$.}

\btm\label{TH:@:Env_C-C^star(G)-circledast-Env_C-C^star(G)->Env_C-C^star(G)-odot-Env_C-C^star(G)}
Для всякой группы Мура $G$ преобразование Гротендика $@:\Env_{\mathcal C}{\mathcal C}^\star(G)\circledast \Env_{\mathcal C}{\mathcal C}^\star(G)\to \Env_{\mathcal C}{\mathcal C}^\star(G)\odot \Env_{\mathcal C}{\mathcal C}^\star(G)$ является эпиморфизмом.
\etm
\bpr
Здесь используется теорема \ref{TH:o-nasl-ster-odot} о наследовании стереотипности слабым тензорным произведением. Чтобы ее применить, вспомним доказательство теоремы \ref{TM:Env_C-C*(G)=projlim-Env_C-C*(G/K)}: операторы $\Env_{\mathcal C}^0\pi_K^\star$, $K\in\lambda(G)$, образуют систему согласованных проекторов на стереотипном пространстве $\Env_{\mathcal C}{\mathcal C}^\star(G)$. По теореме \ref{TH:building-generated-by-P^k} эта система проекторов порождает некий фундамент $\varPhi$ в категории $\Ste$ (и в категории $\InvSteAlg$) над упорядоченным множеством $\lambda(G)$. Эта система проекторов $\Env_{\mathcal C}^0\pi_K^\star$ для алгебры $\Env_{\mathcal C}{\mathcal C}^\star(G)$ удовлетворяет посылкам теоремы \ref{TH:X=edifice-in-Ste}. Это, в частности, значит, что фундамент $\varPhi$ удовлетворяет условию (ii) теоремы \ref{TH:o-nasl-ster-odot} (а также условию (iii), потому что в качестве $\varPsi$ мы собираемся использовать тот же $\varPhi$). 

Чтобы применить теорему \ref{TH:o-nasl-ster-odot}, нужно проверить условие (i) в ней: все пространства $X_i\odot Y_p$ должны быть полны и насыщены. В данном случае компонента фундамента $X_i$ представляет собой пространство вида $\Env_{\mathcal C}{\mathcal C}^\star(G:K)$, $K\in\lambda(G)$. Поэтому
\beq\label{PROOF:@:Env_C-C^star(G)-circledast-Env_C-C^star(G)->Env_C-C^star(G)-odot-Env_C-C^star(G)-1}
X_i=\Env_{\mathcal C}{\mathcal C}^\star(G:K)=\eqref{stroenie-nepr-Env-Lie-Moore}=\prod_{s\in S}{\mathcal C}(M_s)
\eeq
И, точно так же, $Y_p$ --- некоторое пространство $\Env_{\mathcal C}{\mathcal C}^\star(G:L)$, $L\in\lambda(G)$, поэтому
\beq\label{PROOF:@:Env_C-C^star(G)-circledast-Env_C-C^star(G)->Env_C-C^star(G)-odot-Env_C-C^star(G)-2}
Y_p=\Env_{\mathcal C}{\mathcal C}^\star(G:L)=\eqref{stroenie-nepr-Env-Lie-Moore}=\prod_{t\in T}{\mathcal C}(N_i)
\eeq
Отсюда
\begin{multline*}
X_i\odot Y_p=\Env_{\mathcal C}{\mathcal C}^\star(G:K)\odot\Env_{\mathcal C}{\mathcal C}^\star(G:L) =
\prod_{s\in S}{\mathcal C}(M_s)\odot \prod_{t\in T}{\mathcal C}(N_t)=\\=\cite[(4.151)]{Akbarov-De-Gruyter-I}=
\prod_{s\in S, \ t\in T}{\mathcal C}(M_s)\odot {\mathcal C}(N_t)=\cite[(4.274)]{Akbarov-De-Gruyter-I}=
\prod_{s\in S, \ t\in T}{\mathcal C}(M_s\times N_t),
\end{multline*}
--- и эти пространства, очевидно, полны и насыщены.

Теперь мы, наконец, применяем теорему \ref{TH:o-nasl-ster-odot}, и у нас получается, что пространство 
$\Env_{\mathcal C}{\mathcal C}^\star(G)\odot \Env_{\mathcal C}{\mathcal C}^\star(G)$ совпадает с проективным пределом в категории $\LCS$ пространств $\Env_{\mathcal C}{\mathcal C}^\star(G:K)\odot \Env_{\mathcal C}{\mathcal C}^\star(G:L)$:
\begin{multline*}
\Env_{\mathcal C}{\mathcal C}^\star(G)\odot \Env_{\mathcal C}{\mathcal C}^\star(G)=
\edif^{\Ste}\varPhi\odot\edif^{\Ste}\varPhi=\edif^{\Ste}(\varPhi\odot\varPhi)=\\=
\Ste\text{-}\kern-15pt\projlim_{(K,L)\in \lambda(G)\times\lambda(G)}\kern-15pt \Env_{\mathcal C}{\mathcal C}^\star(G:K)\odot \Env_{\mathcal C}{\mathcal C}^\star(G:L)=
\LCS\text{-}\kern-15pt\projlim_{(K,L)\in \lambda(G)\times\lambda(G)}\kern-15pt \Env_{\mathcal C}{\mathcal C}^\star(G:K)\odot \Env_{\mathcal C}{\mathcal C}^\star(G:L)
\end{multline*}
Отсюда следует, что для того, чтобы доказать, что $\Env_{\mathcal C}{\mathcal C}^\star(G)\circledast \Env_{\mathcal C}{\mathcal C}^\star(G)$ плотно в $\Env_{\mathcal C}{\mathcal C}^\star(G)\odot \Env_{\mathcal C}{\mathcal C}^\star(G)$, достаточно убедиться, что каждая компонента $\Env_{\mathcal C}{\mathcal C}^\star(G:K)\circledast \Env_{\mathcal C}{\mathcal C}^\star(G:L)$ плотна в $\Env_{\mathcal C}{\mathcal C}^\star(G:K)\odot \Env_{\mathcal C}{\mathcal C}^\star(G:L)$.
$$
\overline{@\Big(\Env_{\mathcal C}{\mathcal C}^\star(G:K)\circledast \Env_{\mathcal C}{\mathcal C}^\star(G:L)\Big)}=
\Env_{\mathcal C}{\mathcal C}^\star(G:K)\odot \Env_{\mathcal C}{\mathcal C}^\star(G:L)
$$
Разложения \eqref{PROOF:@:Env_C-C^star(G)-circledast-Env_C-C^star(G)->Env_C-C^star(G)-odot-Env_C-C^star(G)-1} и
\eqref{PROOF:@:Env_C-C^star(G)-circledast-Env_C-C^star(G)->Env_C-C^star(G)-odot-Env_C-C^star(G)-2} позволяют свести это к проверке плотности ${\mathcal C}(M_s)\circledast {\mathcal C}(N_t)$ в ${\mathcal C}(M_s)\odot {\mathcal C}(N_t)$:
$$
\overline{@\Big({\mathcal C}(M_s)\circledast {\mathcal C}(N_t)\Big)}=
{\mathcal C}(M_s)\odot {\mathcal C}(N_t)=\cite[(4.274)]{Akbarov-De-Gruyter-I}=
{\mathcal C}(M_s\times N_t)
$$
Но это  очевидно, потому что цилиндрические функции $u\boxdot v$, $u\in {\mathcal C}(M_s)$, $v\in {\mathcal C}(N_t)$, полны в ${\mathcal C}(M_s\times N_t)$. 
\epr

\subsection{$\Env_{\mathcal C} {\mathcal C}^\star(G)$ и ${\mathcal K}(G)$ как алгебры Хопфа}

\paragraph{Теорема о структуре алгебр Хопфа на $\Env_{\mathcal C} {\mathcal C}^\star(G)$ и ${\mathcal K}(G)$.}

Для групп Мура удается доказать, что алгебры $\Env_{\mathcal C}{\mathcal C}^\star(G)$ и ${\mathcal K}(G)$ являются инволютивными алгебрами Хопфа.

\btm\label{TH:Env_C-C^*(G)-odot-Hopf}
 Если $G$ -- группа Мура, то
\bit{
\item[(i)] непрерывная оболочка $\Env_{\mathcal C} {\mathcal C}^\star(G)$ ее групповой алгебры ${\mathcal C}^\star(G)$ является инволютивной алгеброй Хопфа в категории стереотипных пространств $(\tt{Ste},\odot)$,

\item[(ii)] двойственная ей алгебра ${\mathcal K}(G)$ является инволютивной алгеброй Хопфа в категории  стереотипных пространств $(\tt{Ste},\circledast)$.
}\eit
\etm

Для доказательства нам понадобятся некоторые предварительные утверждения.

\paragraph{Предварительные леммы.}

Напомним, что $C^*$-алгебра $A$ называется {\it строгой}, если ее умножение (непрерывно) продолжается до оператора $\mu:A\tilde{\otimes}_{\e}A\to A$ на инъективное локально выпуклое тензорное произведение. В \cite{Akbarov-De-Gruyter-I} мы отмечали, что это условие эквивалентно цепочке равенств
$$
A\underset{\min}{\otimes}A=A\underset{\max}{\otimes}A=A\tilde{\otimes}_{\e}A=A\odot A
$$
и это, в частности, означает, что умножение в $A$ (непрерывно) продолжается до оператора $\mu:A\odot A\to A$ на слабое стереотипное тензорное произведение.

\blm\label{LM:E(H)-dlya-Moore} Пусть $G$ -- группа Мура. Тогда
 для всякой полунормы $p\in{\mathcal P}({\mathcal C}^\star(G))$ фактор-алгебра ${\mathcal C}^\star(G)/p$ является строгой $C^*$-алгеброй.
\elm
\bpr
Если $\tau:{\mathcal C}^\star(G)/p\to{\mathcal B}(X)$ -- какое-нибудь унитарное неприводимое представление, то его композиция с фактор-отображением $\pi_p:{\mathcal C}^\star(G)\to {\mathcal C}^\star(G)/p$ будет унитарным неприводимым представлением инволютивной алгебры ${\mathcal C}^\star(G)$. При этом группа $G$, вкладываясь в ${\mathcal C}^\star(G)$ дельта-функционалами, полна в ${\mathcal C}^\star(G)$. Отсюда следует, что композиция
$$
\xymatrix @R=2.pc @C=4.0pc 
{
G\ar[r]^{\delta}\ar@/_4ex/[rrr]_{\pi}& {\mathcal C}^\star(G)\ar[r]^{\pi_p}& {\mathcal C}^\star(G)/p\ar[r]^{\tau}& {\mathcal B}(X)
}
$$
является унитарным неприводимым представлением группы Мура $G$, и поэтому оно должно быть конечномерно: $\dim X<\infty$. Ограниченность этих чисел (при фиксированном $p$ и меняющихся $\tau$) доказывается в несколько этапов (мы используем идеи \cite[Lemma 5.8]{Kuznetsova}).

1. Сначала рассмотрим случай, когда $G$ -- компактная группа. Тогда в силу предложения \ref{PROP:nepr-obol-komp-gruppy} непрерывная оболочка алгебры мер ${\mathcal C}^\star(G)$ представляет собой алгебру
$\prod_{\sigma\in\widehat{G}}{\mathcal B}(X_\sigma)$. Поскольку фактор-алгебра ${\mathcal C}^\star(G)/p$ представляет собой $C^*$-алгебру, фактор-отображение $\pi_p$ можно продолжить на непрерывную оболочку $\prod_{\sigma\in\widehat{G}}{\mathcal B}(X_\sigma)$:
$$
\xymatrix @R=2.pc @C=1.7pc 
{
{\mathcal C}^\star(G)\ar@/_3ex/[rd]_{\pi_p}\ar[rr]^{\env_{\mathcal C}{\mathcal C}^\star(G)}& & \prod_{\sigma\in\widehat{G}}{\mathcal B}(X_\sigma)\ar@/^3ex/@{-->}[ld]^{\widetilde{\pi_p}} \\
& {\mathcal C}^\star(G)/p &
}
$$
Заметим, что $\widetilde{\pi_p}$ -- эпиморфизм локально выпуклых пространств (поскольку композиция $\widetilde{\pi_p}\circ\env_{\mathcal C}{\mathcal C}^\star(G)=\pi_p$ является эпиморфизмом локально выпуклых пространств). С другой стороны, поскольку все алгебры ${\mathcal B}(X_\sigma)$ конечномерны, алгебра $\prod_{\sigma\in\widehat{G}}{\mathcal B}(X_\sigma)$, как локально выпуклое пространство, представляет собой декартову степень $\C^{\mathfrak m}$ поля $\C$ (где ${\mathfrak m}$ -- некоторое кардинальное число).

Отсюда следует, что линейное непрерывное отображение $\widetilde{\pi_p}$ должно иметь ядро конечной коразмерности (потому что оно действует в банахово пространство ${\mathcal C}^\star(G)/p$). Это в свою очередь означает, что в семействе $\{{\mathcal B}(X_\sigma);\ \sigma\in\widehat{G}\}$ имеется конечный набор алгебр ${\mathcal B}(X_1),...,{\mathcal B}(X_n)$ такой, что $\widetilde{\pi_p}$ представляет собой проекцию на их произведение:
\beq\label{C^star(G)/p-cong-prod_1^n-B(X_i)}
\widetilde{\pi_p}: \prod_{\sigma\in\widehat{G}}{\mathcal B}(X_\sigma)\to \prod_{i=1}^n{\mathcal B}(X_i)\cong {\mathcal C}^\star(G)/p
\eeq
Понятно, что ${\mathcal C}^\star(G)/p$ в таком случае является строгой $C^*$-алгеброй.

2. Пусть далее $G$ -- компактная надстройка абелевой группы, то есть $G=Z\cdot K$, где $Z$ -- абелева, $K$ -- компактная, и они коммутируют (см. определение на с.\pageref{DEF:komp-nadstr-abelevoi-gruppy}).
Рассмотрим ограничения  $\rho=\pi\big|_K$ и $\sigma=\pi\big|_Z$ и обозначим символами $C_\pi$, $C_\rho$ и $C_\sigma$ соответственно $C^*$-подалгебры в ${\mathcal B}(X)$, порожденные образами $\pi$, $\rho$ и $\sigma$. Поскольку $C_\rho$ и $C_\sigma$ коммутируют, $C_\pi$ является непрерывным образом максимального тензорного произведения $C_\rho\underset{\max}{\otimes}C_\sigma$. Заметим далее, что из теорем \ref{TH:nepr-predstavleniya} и \ref{TH:nepr-po-norme-predstavleniya} следует, что $C_\rho$ и $C_\sigma$ являются $C^*$-фактор-алгебрами алгебр ${\mathcal C}^\star(K)$ и ${\mathcal C}^\star(Z)$. При этом из \eqref{C^star(G)/p-cong-prod_1^n-B(X_i)} следует, что $C_\rho=\prod_{i=1}^n{\mathcal B}(X_{\pi_i})$, где $\pi_i\in\widehat{K}$ -- конечный набор унитарных неприводимых представлений $K$, а из предложения \ref{PROP:Env_C-C^star(G)=C(widehat(G))} -- что $C_\sigma={\mathcal C}(T)$ для некоторого компакта $T\subseteq\widehat{Z}$. Отсюда
$$
C_\rho\underset{\max}{\otimes}C_\sigma=\left(\prod_{i=1}^n{\mathcal B}(X_{\pi_i})\right)\underset{\max}{\otimes}{\mathcal C}(T)=\prod_{i=1}^n\left({\mathcal B}(X_{\pi_i})\underset{\max}{\otimes}{\mathcal C}(T)\right)=\prod_{i=1}^n{\mathcal C}\left(T,{\mathcal B}(X_{\pi_i})\right)={\mathcal C}\left(\bigsqcup_{i=1}^n T_i,{\mathcal B}(X_{\pi_i})\right),
$$
где $T_i$ -- копии компакта $T$. Теперь из \cite[10.4.4]{Dixmier} следует, что всякое унитарное неприводимое представление последней алгебры изоморфно какому-нибудь $X_{\pi_i}$, и поэтому имеет размерность, не превосходящую $\max_{i=1,...,n}X_{\pi_i}$. То же самое верно и для унитарных неприводимых представлений алгебры $C_\pi$, потому что они будучи перенесены на $C_\rho\underset{\max}{\otimes}C_\sigma$ тоже становятся унитарными и неприводимыми представлениями.

3. Пусть далее $G$ -- группа Ли-Мура. По теореме \ref{TH:Lie-Moore} $G$ является конечным расширением некоторой компактной надстройки $H=Z\cdot K$ некоторой абелевой группы $Z$. Пусть $m=\card G/H$ -- индекс $H$ в $G$, и пусть $p\in{\mathcal P}({\mathcal C}^\star(G))$ и $\tau$ -- унитарное неприводимое представление алгебры ${\mathcal C}^\star(G)/p$. В силу \cite[Theorem 1]{Clifford}, ограничение $\tau$ на $H$ раскладывается в сумму не более чем $m$ унитарных неприводимых представлений группы $H$, и значит, алгебры ${\mathcal C}^\star(H)/p$. Но в предыдущем пункте мы уже доказали, что размерность унитарных неприводимых представлений алгебры ${\mathcal C}^\star(H)/p$ ограничена некоторым числом. Если его обозначить $n$ ($n\in\N$), то мы получим, что размерность представления $\tau$ не превосходит $m\cdot n$.

4. Пусть, наконец, $G$ -- произвольная группа Мура и $p\in{\mathcal P}({\mathcal C}^\star(G))$. Гомоморфизм $G\to {\mathcal C}^\star(G)/p$ непрерывен по норме, поэтому в силу \cite[Theorem 1]{Shtern}, он пропускается через некоторое фактор-отображение $G\to G/H$, в котором $G/H$ -- группа Ли. По теореме  \ref{TH:Moore->quotient} $G/H$ будет группой Мура. В результате все сводится к случаю 3.
\epr

\blm\label{COR:nepr-Env-SIN-gruppy=LCS-lim} Пусть $G=Z\cdot K$ --- надстройка абелевой локально компактной группы $Z$ с помощью компактной группы $K$. Тогда непрерывная оболочка групповой алгебры ${\mathcal C}^\star(G)$ совпадает с проективным пределом ее $C^*$-фактор-алгебр в категории локально выпуклых пространств, и в категории $\Ste$ стереотипных пространств\footnote{Содержательный смысл леммы \ref{COR:nepr-Env-SIN-gruppy=LCS-lim} состоит в том, что, во-первых, топология предела $\LCS\text{-}\kern-3pt\projlim_{p\in{\mathcal P}({\mathcal C}^\star(G))}{\mathcal C}^\star(G)/p$ псевдонасыщена, и, во-вторых, $\Env_{\mathcal C}{\mathcal C}^\star(G)$ как множество совпадает с $\LCS\text{-}\kern-3pt\projlim_{p\in{\mathcal P}({\mathcal C}^\star(G))}{\mathcal C}^\star(G)/p$.}:
\beq\label{E(C*(G))=LCS-leftlim}
\Env_{\mathcal C}{\mathcal C}^\star(G)=\LCS\text{-}\kern-3pt\projlim_{p\in{\mathcal P}({\mathcal C}^\star(G))}{\mathcal C}^\star(G)/p=\Ste\text{-}\kern-3pt\projlim_{p\in{\mathcal P}({\mathcal C}^\star(G))}{\mathcal C}^\star(G)/p
\eeq
\elm
\bpr
Первое равенство сразу видно из формулы \eqref{Env_C-C*(Z-cdot-K)}, а из него следует второе.
\epr

\paragraph{Случай компактной надстройки абелевой группы.}

Теорему \ref{TH:Env_C-C^*(G)-odot-Hopf} удобно сначала доказать для случая компактной надстройки абелевой группы:

\blm\label{TH:Env_C-C^*(Z-cdot-K)-odot-Hopf} Пусть $G=Z\cdot K$ --- надстройка абелевой локально компактной группы $Z$ с помощью компактной группы $K$. Тогда непрерывная оболочка $\Env_{\mathcal C} {\mathcal C}^\star(G)$ ее групповой алгебры ${\mathcal C}^\star(G)$ является инволютивной алгеброй Хопфа в категории стереотипных пространств $(\tt{Ste},\odot)$.
\elm
\bpr

1. Прежде всего, ${\mathcal C}^\star(G)$ является стереотипной алгеброй (то есть алгеброй в категории $(\tt{Ste},\circledast)$), и ее оболочка $\Env_{\mathcal C}{\mathcal C}^\star(G)$ также является стереотипной алгеброй  (то есть алгеброй в категории  $(\tt{Ste},\circledast)$). Покажем, что умножение в $\Env_{\mathcal C}{\mathcal C}^\star(G)$ продолжается до оператора на $\Env_{\mathcal C}{\mathcal C}^\star(G)\odot \Env_{\mathcal C}{\mathcal C}^\star(G)$.

По лемме \ref{LM:E(H)-dlya-Moore}, ${\mathcal C}^\star(G)/p$ -- строгая $C^*$-алгебра. Поэтому умножение $\mu_p$ в ней продолжается до оператора $\mu'_p:{\mathcal C}^\star(G)/p\odot {\mathcal C}^\star(G)/p\to {\mathcal C}^\star(G)/p$. Рассмотрим диаграмму
$$
\xymatrix @R=2.pc @C=4.0pc 
{
{\mathcal C}^\star(G)\circledast{\mathcal C}^\star(G)\ar[d]^{\mu}\ar[r]^{\pi_p\circledast\pi_p}& {\mathcal C}^\star(G)/p\circledast{\mathcal C}^\star(G)/p\ar[d]^{\mu_p}\ar[r]^{@}& {\mathcal C}^\star(G)/p\odot{\mathcal C}^\star(G)/p\ar[d]^{\mu'_p} \\
{\mathcal C}^\star(G)\ar[r]^{\pi_p}& {\mathcal C}^\star(G)/p\ar[r]^{\id}& {\mathcal C}^\star(G)/p
}
$$
Выбросив средний столбец и перейдя к проективному пределу в $\tt{Ste}$, получим:
$$
\xymatrix @R=2.pc @C=2.0pc 
{
{\mathcal C}^\star(G)\circledast{\mathcal C}^\star(G)\ar[d]^{\mu}\ar@{-->}[r]& \kern-13pt\leftlim_{p\in{\mathcal P({\mathcal C}^\star(G))}}\kern-13pt\Big({\mathcal C}^\star(G)/p\odot{\mathcal C}^\star(G)/p\Big)\ar@{-->}[d]\ar@{=}[r]^{\cite[(4.153)]{Akbarov-De-Gruyter-I}}&
\kern-13pt
\leftlim_{p\in{\mathcal P({\mathcal C}^\star(G))}}\kern-13pt{\mathcal C}^\star(G)/p\odot\kern-13pt \leftlim_{p\in{\mathcal P({\mathcal C}^\star(G))}}\kern-13pt{\mathcal C}^\star(G)/p\ar@{=}[r]^{\eqref{E(C*(G))=LCS-leftlim}} &
\Env_{\mathcal C}{\mathcal C}^\star(G)\odot \Env_{\mathcal C}{\mathcal C}^\star(G)\ar@{-->}[d]^{\mu'}\\
{\mathcal C}^\star(G)\ar@{-->}[r]& \kern-13pt\leftlim_{p\in{\mathcal P({\mathcal C}^\star(G))}}\kern-13pt{\mathcal C}^\star(G)/p\ar@{=}[rr]& & \Env_{\mathcal C}{\mathcal C}^\star(G)
}
$$
Оператор $\mu'$ и будет искомым продолжением умножения на инъективный тензорный квадрат.

2. Далее, ${\mathcal C}^\star(G)$ является биалгеброй в категории $(\tt{Ste},\circledast)$, поэтому по теореме \ref{TH:C-obolochka-sohranyaet-Hopfov} $\Env_{\mathcal C}{\mathcal C}^\star(G)$ является коалгеброй в категории
$(\tt{Ste},\odot)$. Таким образом, $\Env_{\mathcal C}{\mathcal C}^\star(G)$ -- и алгебра, и коалгебра в $(\tt{Ste},\odot)$. Рассмотрим диаграмму

{\scriptsize
 $$\dgARROWLENGTH=-6em
\begin{diagram}
  \node[3]{H}\arrow[2]{se,l}{\varkappa}\arrow[4]{s,r,-}{\env_{\mathcal C}H}
 \\
 \\
 \node{H\circledast H}\arrow[3]{s,t}{@}
 \arrow{sse,l}{\varkappa\circledast\varkappa}
 \arrow[2]{ne,l}{\mu}
 \node[4]{H\circledast H}\arrow[3]{s,r}{@} \\
  \\
 \node[2]{(H\circledast H)\circledast (H\circledast H)}
 \arrow[2]{e,l,3}{\theta}
 \arrow[2]{s,l,1}{@}
 \node{}\arrow[3]{s}
 \node{(H\circledast H)\circledast (H\circledast H)}
 \arrow{nne,l}{\mu\circledast\mu}
 \arrow[2]{s,r,1}{@}
 \\
  \node{H\odot H}\arrow[4]{s,b}{\env_{\mathcal C}H\odot \env_{\mathcal C}H}
 \node[4]{H\odot H}\arrow[4]{s,l}{\env_{\mathcal C}H\odot \env_{\mathcal C}H}
 \\
 \node[2]{(H\odot H)\odot (H\odot H)}
 \arrow[5]{s,r}{(\env_{\mathcal C}H\odot \env_{\mathcal C}H)\odot(\env_{\mathcal C}H\odot \env_{\mathcal C}H)}
 \node[2]{(H\odot H)\odot (H\odot H)}
 \arrow[5]{s,l,3}{(\env_{\mathcal C}H\odot \env_{\mathcal C}H)\odot(\env_{\mathcal C}H\odot \env_{\mathcal C}H)}
 \\
  \node[3]{\Env_{\mathcal C}H}\arrow{se,-}
 \\
 \node[2]{}\arrow{ne,l,3}{\mu'}
 \node[2]{}\arrow{se,t,3}{\varkappa'}
 \\
 \node{\Env_{\mathcal C}H\odot \Env_{\mathcal C}H}
 \arrow{sse,l}{\varkappa'\odot\varkappa'}
 \arrow{ne,-}
 \node[4]{\Env_{\mathcal C}H\odot \Env_{\mathcal C}H}
  \\ \\
 \node[2]{\big(\Env_{\mathcal C}H\odot \Env_{\mathcal C}H\big)\odot \big(\Env_{\mathcal C}H \odot \Env_{\mathcal C}H\big)}
 \arrow[2]{e,l}{\theta'}
 \node[2]{\big(\Env_{\mathcal C}H\odot \Env_{\mathcal C}H\big)\odot \big(\Env_{\mathcal C}H \odot \Env_{\mathcal C}H\big)}
 \arrow{nne,l}{\mu'\odot\mu'}
\end{diagram}
 $$
 }
в которой $H={\mathcal C}^\star(G)$, а смысл остальных обозначений очевиден. Здесь верхняя грань коммутативна, потому что $H$ -- алгебра Хопфа в $(\tt{Ste},\circledast)$, а коммутативность боковых граней проверяется прямым вычислением. Вдобавок, по теореме \ref{TH:@:Env_C-C^star(G)-circledast-Env_C-C^star(G)->Env_C-C^star(G)-odot-Env_C-C^star(G)}  самый левый морфизм 
$$
\env_{\mathcal C}H\odot \env_{\mathcal C}H\circ @=@\circ \env_{\mathcal C}H\circledast \env_{\mathcal C}H
$$
 --- является эпиморфизмом в $\tt{Ste}$. Отсюда следует, что нижняя грань призмы также коммутативна.

Тем же манером доказывается коммутативность остальных диаграмм в определении алгебры Хопфа. Например, диаграмма, связывающая коумножение с единицей, проверяется достраиванием до следующей призмы:
 $$\dgARROWLENGTH=1em
\begin{diagram}
  \node[3]{\C\circledast\C}\arrow{se,l}{\iota\circledast\iota}\arrow[2]{s,-}
 \\
   \node[4]{H\circledast H}\arrow[2]{s,l}{@}
 \\
 \node{\C}\arrow[2]{ne,l}{l_{\C}^{-1}}\arrow{se,l}{\iota}\arrow[4]{s,l}{1_{\C}}
 \node[2]{}\arrow[2]{s,l}{@}
 \\
 \node[2]{H}\arrow[2]{ne,r,3}{\varkappa}\arrow[4]{s,r,3}{\env_{\mathcal C}H}\node[2]{H\odot H}\arrow[2]{s,r}{\env_{\mathcal C}H\odot \env_{\mathcal C}H}
 \\
  \node[3]{\C\odot\C}\arrow{se,l}{\iota\circledast\iota}
 \\
  \node[2]{}\arrow{ne,l}{l_{\C}^{-1}} \node[2]{\Env_{\mathcal C}H\odot \Env_{\mathcal C}H}
 \\
 \node{\C}\arrow{ne,-}\arrow{se,l}{\iota}
 \\
 \node[2]{\Env_{\mathcal C}H}\arrow[2]{ne,r,3}{\varkappa}
 \\
\end{diagram}
 $$
В ней верхнее основание коммутативно, потому что $H$ -- алгебра Хопфа в $(\tt{Ste},\circledast)$, а коммутативность боковых граней следует из свойств функтора $\Env_{\mathcal C}$. Поэтому нижнее основание также должно быть коммутативно.

3. Теперь покажем, что инволюция $\bullet$ в алгебре Хопфа ${\mathcal C}^\star(G)$ порождает инволюцию $\bullet'$ в алгебре Хопфа $\Env_{\mathcal C}{\mathcal C}^\star(G)$. Для всякой полунормы $p\in{\mathcal P}({\mathcal C}^\star(G))$ рассмотрим естественную проекцию $\pi_p:\Env_{\mathcal C}{\mathcal C}^\star(G)\to {\mathcal C}^\star(G)/p$. Пусть $\bullet/p$ -- инволюция на ${\mathcal C}^\star(G)/p$, порожденная инволюцией $\bullet$ на ${\mathcal C}^\star(G)$. Обозначим $\bullet_p=\bullet/p\ \circ\pi_p$:
$$
\xymatrix @R=2.pc @C=4.0pc 
{
\Env_{\mathcal C}{\mathcal C}^\star(G)\ar[r]^{\pi_p}\ar@/_2ex/[rd]_{\bullet_p} & {\mathcal C}^\star(G)/p\ar[d]^{\bullet/p} \\
& {\mathcal C}^\star(G)/p
}
$$
Для любых полунорм $p,q\in{\mathcal P}({\mathcal C}^\star(G))$, $p\le q$, в категории стереотипных пространств над вещественным полем $\R$ будет коммутативна диаграмма
$$
\xymatrix @R=2.pc @C=2.0pc 
{
& \Env_{\mathcal C}{\mathcal C}^\star(G)\ar@/_2ex/[dl]_{\bullet_q}\ar@/^2ex/[dr]^{\bullet_p} &  \\
{\mathcal C}^\star(G)/q\ar[rr]_{\pi^q_p} & & {\mathcal C}^\star(G)/p
}
$$
в которой $\pi^q_p$ -- естественная проекция. Это означает, что семейство морфизмов $\bullet_p$ является проективным конусом для системы $\pi^q_p$, и значит существует единственный морфизм $\bullet'$, замыкающий все диаграммы
$$
\xymatrix @R=2.pc @C=2.0pc 
{
& & \Env_{\mathcal C}{\mathcal C}^\star(G)\ar@/_2ex/[dl]_{\bullet'}\ar@/^2ex/[dr]^{\bullet_p} &  \\
\Env_{\mathcal C}{\mathcal C}^\star(G)\ar@{=}[r]& \leftlim_{q}{\mathcal C}^\star(G)/q\ar[rr]_{\pi_p} & & {\mathcal C}^\star(G)/p
}
$$
Покажем, что морфизм $\bullet'$ и будет нужной инволюцией в $\Env_{\mathcal C}{\mathcal C}^\star(G)$. Прежде всего, он связан с исходной инволюцией $\bullet$ коммутативной диаграммой:
$$
\xymatrix @R=3.pc @C=6.0pc 
{
{\mathcal C}^\star(G)\ar[r]^{\env_{\mathcal C}{{\mathcal C}^\star(G)}}\ar[d]_{\bullet} & \Env_{\mathcal C}{\mathcal C}^\star(G)\ar[d]^{\bullet'} \\
{\mathcal C}^\star(G)\ar[r]^{\env_{\mathcal C}{{\mathcal C}^\star(G)}} & \Env_{\mathcal C}{\mathcal C}^\star(G) \\
}
$$
Из нее следует равенство
$$
\bullet'\circ \bullet'\circ \env_{\mathcal C}{{\mathcal C}^\star(G)}=\env_{\mathcal C}{{\mathcal C}^\star(G)}\circ \bullet\circ \bullet=\env_{\mathcal C}{{\mathcal C}^\star(G)}=
\id_{\Env_{\mathcal C}{\mathcal C}^\star(G)}\circ \env_{\mathcal C}{{\mathcal C}^\star(G)}
$$
которое, в силу эпиморфности $\env_{\mathcal C}{{\mathcal C}^\star(G)}$, дает
$$
\bullet'\circ \bullet'=\id_{\Env_{\mathcal C}{\mathcal C}^\star(G)}.
$$

Чтобы доказать согласованность умножения с инволюцией \cite[(5.134)]{Akbarov-De-Gruyter-I}, рассмотрим диаграмму (в категории стереотипных пространств над $\R$), похожую на ту, что была выше в пункте 2 (здесь $H={\mathcal C}^\star(G)$, а $\br:x\otimes y\mapsto y\otimes x$ -- морфизм заузливания в моноидальной категории):
 $$\dgARROWLENGTH=-2em
\begin{diagram}
  \node[3]{H}\arrow[2]{se,l}{\bullet}\arrow[4]{s,r,-}{\env_{\mathcal C}H}
 \\
 \\
 \node{H\circledast H}\arrow[3]{s,t}{@}
 \arrow{sse,l}{\br}
 \arrow[2]{ne,l}{\mu}
 \node[4]{H}\arrow[7]{s,r}{\env_{\mathcal C}H} \\
  \\
 \node[2]{H\circledast H}
 \arrow[2]{e,l,3}{\bullet\circledast\bullet}
 \arrow[2]{s,l,1}{@}
 \node{}\arrow[3]{s}
 \node{H\circledast H}
 \arrow{nne,l}{\mu}
 \arrow[2]{s,r,1}{@}
 \\
  \node{H\odot H}\arrow[4]{s,t}{\env_{\mathcal C}H\odot \env_{\mathcal C}H}
 \\
 \node[2]{H\odot H}
 \arrow[5]{s,r,3}{\env_{\mathcal C}H\odot \env_{\mathcal C}H}
 \node[2]{H\odot H}
 \arrow[5]{s,l,3}{\env_{\mathcal C}H\odot \env_{\mathcal C}H}
 \\
  \node[3]{\Env_{\mathcal C}H}\arrow{se,-}
 \\
 \node[2]{}\arrow{ne,r,3}{\mu'}
 \node[2]{}\arrow{se,l,3}{\bullet'}
 \\
 \node{\Env_{\mathcal C}H\odot \Env_{\mathcal C}H}
 \arrow{sse,r}{\br}
 \arrow{ne,-}
 \node[4]{\Env_{\mathcal C}H}
  \\ \\
 \node[2]{\Env_{\mathcal C}H\odot \Env_{\mathcal C}H}
 \arrow[2]{e,l}{\bullet'\odot\bullet'}
 \node[2]{\Env_{\mathcal C}H\odot \Env_{\mathcal C}H}
 \arrow{nne,r}{\mu'}
\end{diagram}
 $$
В этой призме коммутативны верхнее основание (как \cite[(5.134)]{Akbarov-De-Gruyter-I} для $H$), и кроме того, в ней коммутативны боковые грани, а по теореме \ref{TH:@:Env_C-C^star(G)-circledast-Env_C-C^star(G)->Env_C-C^star(G)-odot-Env_C-C^star(G)} левое ребро 
$$
\env_{\mathcal C}H\odot \env_{\mathcal C}H\circ @=@\circ \env_{\mathcal C}H\circledast \env_{\mathcal C}H
$$
является эпиморфизмом. Отсюда следует, что нижнее основание тоже коммутативно:
$$
\bullet'\circ \mu=\mu\circ \bullet'\odot\bullet'\circ\br.
$$

Для доказательства согласованности коумножения с инволбцией \cite[(5.135)]{Akbarov-De-Gruyter-I}, рассмотрим диаграмму (в категории стереотипных пространств над $\R$):
 $$\dgARROWLENGTH=1em
\begin{diagram}
  \node[3]{H}\arrow{se,l}{\varkappa}\arrow[2]{s,-}
 \\
   \node[4]{H\circledast H}\arrow{s,r}{@}
 \\
 \node{H}\arrow[2]{ne,l}{\bullet}\arrow{se,l}{\varkappa}\arrow[4]{s,l}{\env_{\mathcal C}H}
 \node[2]{}\arrow[2]{s,l}{\env_{\mathcal C}H}\node{H\odot H}\arrow[3]{s,r,3}{\env_{\mathcal C}H\odot \env_{\mathcal C}H}
 \\
 \node[2]{H\circledast H}\arrow[2]{ne,l,1}{\bullet\circledast\bullet}\arrow{s,r}{@}
 \\
 \node[2]{H\odot H}\arrow[3]{s,r,3}{\env_{\mathcal C}H\odot \env_{\mathcal C}H} \node{\Env_{\mathcal C}H}\arrow{se,l}{\varkappa'}
 \\
  \node[2]{}\arrow{ne,r}{\bullet'} \node[2]{\Env_{\mathcal C}H\odot \Env_{\mathcal C}H}
 \\
 \node{\Env_{\mathcal C}H}\arrow{ne,-}\arrow{se,l}{\varkappa'}
 \\
 \node[2]{\Env_{\mathcal C}H\odot \Env_{\mathcal C}H}\arrow[2]{ne,r,3}{\bullet'\odot\bullet'}
 \\
\end{diagram}
 $$
Здесь верхнее основание коммутативно (как \cite[(5.135)]{Akbarov-De-Gruyter-I} для $H$), и кроме того, коммутативны боковые грани, а самый левый морфизм $\env_{\mathcal C}H$ является эпиморфизмом. Отсюда следует, что справедливо аналогичное равенство в $\Env_{\mathcal C}H$:
$$
\varkappa'\circ\bullet'=\bullet'\odot\bullet'\circ\varkappa'.
$$
\epr

\paragraph{Доказательство теоремы \ref{TH:Env_C-C^*(G)-odot-Hopf}.}

\bpr[Доказательство теоремы \ref{TH:Env_C-C^*(G)-odot-Hopf}.]
1. Рассмотрим сначала случай, когда $G$ -- группа Ли---Мура. Тогда по теореме \ref{TH:Lie-Moore} $G$ --- конечное расширение некоторой компактной надстройки абелевой группы (Ли):
$$
1\to Z\cdot K=N\to G\to F\to 1
$$
($Z$ -- абелева группа Ли, $K$ -- компактная группа Ли, $F$ -- конечная группа). Эта цепочка порождает цепочку гомоморфизмов групповых алгебр
$$
\C\to {\mathcal C}^\star(Z\cdot K)={\mathcal C}^\star(N)\to {\mathcal C}^\star(G)\to {\mathcal C}^\star(F)=\C_F\to \C.
$$
(второе равенство в этой цепочке следует из примера \ref{EX:Env_C-C_F=C_F}) и цепочку гомоморфизмов непрерывных оболочек
$$
\C\to \Env_{\mathcal C} {\mathcal C}^\star(Z\cdot K)=\Env_{\mathcal C} {\mathcal C}^\star(N)\to \Env_{\mathcal C} {\mathcal C}^\star(G)\to \Env_{\mathcal C} {\mathcal C}^\star(F)=\eqref{Env_C-C_F=C_F}=\C_F\to \C.
$$

Вспомним пространства $E_L$, определенные в \eqref{DEF:E_L=x_L-cdot-E_N}. По формуле \eqref{Env_C-C*(konech-rasshirenie-ZK)} пространство $\Env_{\mathcal C} {\mathcal C}^\star(G)$ представляется как прямая сумма подпространств $E_L$:
\beq\label{Env_C-C*(konech-rasshirenie-ZK)-0}
\Env_{\mathcal C}  {\mathcal C}^\star(G)=\bigoplus_{L\in G/N} E_L.
\eeq
Наша задача -- показать, что эта алгебра слабая (то есть является алгеброй относительно слабого тензорного произведения $\odot$). Прежде всего, заметим, что по лемме \ref{TH:Env_C-C^*(Z-cdot-K)-odot-Hopf}, алгебра $E_N=\Env_{\mathcal C}  {\mathcal C}^\star(N)$ слабая. Пусть
$$
m_N: E_N\odot E_N\to E_N
$$
-- непрерывное продолжение ее оператора умножения. Нам нужно достроить его до оператора
\beq\label{DEF:umnozhenie-na-Env-E*(G)-odot-Env-E*(G)}
m: \Env_{\mathcal C}  {\mathcal C}^\star(G)\odot \Env_{\mathcal C}  {\mathcal C}^\star(G)\to \Env_{\mathcal C}  {\mathcal C}^\star(G),
\eeq
который будет продолжением оператора умножения в $\Env_{\mathcal C}  {\mathcal C}^\star(G)$.

Это можно сделать так. Каждому классу смежности $L\in G/N$ поставим в соответствие какой-нибудь элемент $g_L\in L$. Рассмотрим операторы сдвига
$$
T^l_L:E_N\to E_L\quad\Big|\quad T^l_Lx=\env_{\mathcal C} \delta^{g_L}\cdot x,\qquad x\in E_N
$$
$$
T^r_L:E_N\to E_L\quad\Big|\quad T^r_Lx=x\cdot \env_{\mathcal C} \delta^{g_L},\qquad x\in E_N
$$
и обратного сдвига
$$
(T^l_L)^{-1}:E_N\gets E_L\quad\Big|\quad (T^l_L)^{-1}y=\env_{\mathcal C} \delta^{g_L^{-1}}\cdot y,\qquad x\in E_L
$$
$$
(T^r_L)^{-1}:E_N\gets E_L\quad\Big|\quad (T^r_L)^{-1}y=y\cdot\env_{\mathcal C} \delta^{g_L^{-1}},\qquad x\in E_L.
$$
Разложим инъективный тензорный квадрат пространства $\Env_{\mathcal C}  {\mathcal C}^\star(G)$ в произведение его компонент:
$$
\Env_{\mathcal C}  {\mathcal C}^\star(G)\odot \Env_{\mathcal C}  {\mathcal C}^\star(G)\cong
\bigoplus_{L,M\in G/N} E_L\odot E_M.
$$
Тогда оператор \eqref{DEF:umnozhenie-na-Env-E*(G)-odot-Env-E*(G)} можно определить на каждой компоненте $E_L\odot E_M$ формулой
\beq\label{PROOF:TH:gladk-obolochka-komp-porozhd-Lie-Moore-1}
m=T^l_L\circ T^r_M\circ m_N\circ \big((T^l_L)^{-1}\odot (T^r_M)^{-1}\big).
\eeq
Чтобы убедиться, что получаемый оператор продолжает умножение в $\Env_{\mathcal C}  {\mathcal C}^\star(G)$ достаточно проверить это на дельта-функционалах (которые полны в ${\mathcal C}^\star(G)$ и значит в $\Env_{\mathcal C}  {\mathcal C}^\star(G)$ тоже). Зафиксируем два элемента $a,b\in G$ и подберем классы смежности $L,M\in G/N$ такие что $a\in L$ и $b\in M$. Тогда
\begin{multline*}
m(\env_{\mathcal C} \delta^a\odot\env_{\mathcal C} \delta^b)=\Big(T^l_L\circ T^r_M\circ m_N\circ \big((T^l_L)^{-1}\odot (T^r_M)^{-1}\big)\Big)(\env_{\mathcal C} \delta^a\odot\env_{\mathcal C} \delta^b)=\\=
\Big(T^l_L\circ T^r_M\circ m_N\Big)\big((T^l_L)^{-1}(\env_{\mathcal C} \delta^a)\odot (T^r_M)^{-1}(\env_{\mathcal C} \delta^b)\big)=\\=
\Big(T^l_L\circ T^r_M\circ m_N\Big)\big((\env_{\mathcal C} \delta^{g_L^{-1}}\cdot\env_{\mathcal C} \delta^a)\odot (\env_{\mathcal C} \delta^b\cdot\env_{\mathcal C} \delta^{g_M^{-1}})\big)=\\=
\Big(T^l_L\circ T^r_M\circ m_N\Big)\big(\env_{\mathcal C} (\delta^{g_L^{-1}}*\delta^a)\odot \env_{\mathcal C} (\delta^b*\delta^{g_M^{-1}})\big)=
\Big(T^l_L\circ T^r_M\Big)\big(\env_{\mathcal C} (\delta^{g_L^{-1}}*\delta^a)\cdot \env_{\mathcal C} (\delta^b*\delta^{g_M^{-1}})\big)=\\=
\Big(T^l_L\circ T^r_M\Big)\big(\env_{\mathcal C} (\delta^{g_L^{-1}}*\delta^a*\delta^b*\delta^{g_M^{-1}})\big)=
\env_{\mathcal C} (\delta^{g_L})\cdot\env_{\mathcal C} (\delta^{g_L^{-1}}*\delta^a*\delta^b*\delta^{g_M^{-1}})
\cdot\env_{\mathcal C} (\delta^{g_M})=\\=
\env_{\mathcal C} (\delta^{g_L}*\delta^{g_L^{-1}}*\delta^a*\delta^b*\delta^{g_M^{-1}}*\delta^{g_M})=
\env_{\mathcal C} (\delta^a*\delta^b)=
\env_{\mathcal C} (\delta^a)\cdot\env_{\mathcal C} (\delta^b)
\end{multline*}

2. Теперь пусть $G$ -- произвольная группа Мура. Тогда по теореме \ref{TH:Moore=lim-Lie-Moore} ее можно представить как проективный предел своих факторгрупп Ли---Мура:
$$
G=\projlim_{K\in\lambda(G)} G/K
$$
При этом по теореме \ref{TM:Env_C-C*(G)=projlim-Env_C-C*(G/K)}, пространство $\Env_{\mathcal C} {\mathcal C}^\star(G)$ есть проективный предел пространств $\Env_{\mathcal C} {\mathcal C}^\star(G)$:
\beq\label{PROOF:nepr-obolochka-Moore-2}
\Env_{\mathcal C} {\mathcal C}^\star(G)=\projlim_{K\in\lambda(G)}
\Env_{\mathcal C} {\mathcal C}^\star (G/K).
\eeq
Поскольку $G_i$ есть группы Ли---Мура, по уже доказанному алгебры $\Env_{\mathcal C} {\mathcal C}^\star (G_i)$ будут инволютивными алгебрами Хопфа в категории $(\tt{Ste},\odot)$. Отсюда по следствию \ref{COR:proj-lim-Hopf-algebr-odot} мы получаем, что их проективный предел $\Env_{\mathcal C} {\mathcal C}^\star (G)$ также является инволютивной алгеброй Хопфа в категории $(\tt{Ste},\odot)$.
\epr

\section{Непрерывная двойственность}

\subsection{Непрерывная рефлексивность.}

Применим теперь конструкцию описанную на с.\pageref{DEF:algebra-Hopfa-reflex-otn-obolochki} к оболочке $\Env_{\mathcal C}$. Условимся говорить, что инволютивная стереотипная алгебра Хопфа $H$ в категории $({\tt Ste},\circledast)$ {\it непрерывно рефлексивна}\label{DEF:reflexiv-otn-obolochki}, если на ее непрерывной оболочке $\Env_{\mathcal C} H$ определена структура инволютивной алгебры Хопфа в категории $({\tt Ste},\odot)$ так, что выполняются следующие два условия:
\bit{
\item[(i)] морфизм непрерывной оболочки $\env_{\mathcal C} H:H\to \Env_{\mathcal C} H$ является гомоморфизмом алгебр Хопфа в том смысле, что коммутативны следующие диаграммы:
\beq\label{DIAG:reflex-otn-obolochki-1}
 \xymatrix @R=2.pc @C=2.pc
{
& H\odot H\ar[dr]^{\quad \env_{\mathcal C} H\odot \env_{\mathcal C} H} & \\
H\circledast H\ar[ur]^{@}\ar[dr]^{\quad \env_{\mathcal C} H\, \circledast\, \env_{\mathcal C} H}\ar[dd]_{\mu} & & \Env_{\mathcal C} H\odot \Env_{\mathcal C} H\ar[dd]_{\mu_E} \\
& \Env_{\mathcal C} H\circledast \Env_{\mathcal C} H\ar[ur]^{@} & \\
H\ar[rr]^{\env_{\mathcal C} H} && \Env_{\mathcal C} H
}
\eeq
\beq\label{DIAG:reflex-otn-obolochki-2}
 \xymatrix @R=2.pc @C=2.pc
{
& H\odot H\ar[dr]^{\quad\env_{\mathcal C} H\odot \env_{\mathcal C} H} & \\
H\circledast H\ar[ur]^{@}\ar[dr]^{\quad \env_{\mathcal C} H\circledast \env_{\mathcal C} H} & & \Env_{\mathcal C} H\odot \Env_{\mathcal C} H \\
& \Env_{\mathcal C} H\circledast \Env_{\mathcal C} H\ar[ur]^{@} & \\
H\ar[rr]^{\env_{\mathcal C} H}\ar[uu]^{\varkappa} && \Env_{\mathcal C} H\ar[uu]^{\varkappa_E}
}
\eeq
\beq\label{DIAG:reflex-otn-obolochki-3}
 \xymatrix @R=2.pc @C=2.pc
{
H\ar[rr]^{\env_{\mathcal C} H} & & \Env_{\mathcal C} H\\
& \C\ar[ul]^{\iota}\ar[ur]_{\iota_E} &
}\qquad
 \xymatrix @R=2.pc @C=2.pc
{
H\ar[rr]^{\env_{\mathcal C} H}\ar[dr]_{\e} & & \Env_{\mathcal C} H\ar[dl]^{\e_E} \\
& \C &
}
\eeq
\beq\label{DIAG:reflex-otn-obolochki-4}
 \xymatrix @R=3.pc @C=4.pc
{
H\ar[r]^{\env_{\mathcal C} H}\ar[d]_{\sigma} & \Env_{\mathcal C} H\ar[d]^{\sigma_E} \\
H\ar[r]^{\env_{\mathcal C} H} & \Env_{\mathcal C} H
}\qquad
 \xymatrix @R=3.pc @C=4.pc
{
H\ar[r]^{\env_{\mathcal C} H}\ar[d]_{\bullet} & \Env_{\mathcal C} H\ar[d]^{\bullet_E} \\
H\ar[r]^{\env_{\mathcal C} H} & \Env_{\mathcal C} H
}
\eeq
-- здесь $@$ -- преобразование Гротендика\footnote{См.\cite{Ak03} или \cite{Akbarov-De-Gruyter-I}}, $\mu$, $\iota$, $\varkappa$, $\e$, $\sigma$, $\bullet$ -- структурные морфизмы (умножение, единица, коумножение, коединица, антипод, инволюция) в $H$, а $\mu_E$, $\iota_E$, $\varkappa_E$, $\e_E$, $\sigma_E$, $\bullet_E$ -- структурные морфизмы в $\Env_{\mathcal C} H$.

\item[(ii)]\label{(env-H)^star:H^star-gets-(Env-H)^star} отображение $(\env_{\mathcal C} H)^\star:H^\star\gets (\Env_{\mathcal C} H)^\star$, сопряженное к морфизму непрерывной оболочки $\env_{\mathcal C} H:H\to\Env_{\mathcal C} H$, --- тоже непрерывная оболочка:
$$
(\env_{\mathcal C} H)^\star=\env_{\mathcal C} (\Env_{\mathcal C} H)^\star
$$

\item[(iii)]\label{@:H-circledast-H->H-odot-H-Mono-C} преобразование Гротендика $@:H\circledast H\to H\odot H$ является мономорфизмом стереотипных пространств (то есть инъективным отображением),

\item[(iv)]\label{@:E(H)-circledast-E(H)->E(H)-odot-E(H)-Epi-C}  преобразование Гротендика $@:E(H)\circledast E(H)\to E(H)\odot E(H)$ является эпиморфизмом стереотипных пространств (то есть его образ плотен в области значений).

}\eit

Условия (i) и (ii) удобно изображать в виде диаграммы
 \beq\label{obshaya-diagramma-refleksivnosti}
 \xymatrix @R=1.pc @C=1.pc
 {
 H
 & \ar@{|->}[r]^{\env_{\mathcal C}} & &
\Env_{\mathcal C} H
 \\
 & & &
 \ar@{|->}[d]^{\star}
 \\
 \ar@{|->}[u]^{\star}
 & & &
 \\
 H^\star
 & &
 \ar@{|->}[l]_{\env_{\mathcal C}}
 &
 (\Env_{\mathcal C} H)^\star
 }
 \eeq
которую мы называем {\it диаграммой рефлексивности}, и в которую вкладываем следующий смысл:
 \bit{
\item[1)] в углах квадрата стоят инволютивные алгебры Хопфа, причем $H$ -- алгебра Хопфа в $({\tt Ste},\circledast)$, затем следует алгебра Хопфа $\Env_{\mathcal C} H$ в $({\tt Ste},\odot)$, и далее категории $({\tt Ste},\circledast)$ и $({\tt Ste},\odot)$ чередуются,

\item[2)] чередование операций $\env_{\mathcal C}$ и $\star$ (с какого места ни начинай) на четвертом шаге возвращает к исходной алгебре Хопфа (с точностью до изоморфизма функторов).
 }\eit

Смысл термина ``рефлексивность'' здесь состоит в следующем. Обозначим однократное последовательное применение операций $\env_{\mathcal C}$ и $\star$ каким-нибудь символом, например $\widehat{\ }\,\,$,
$$
\widehat H:=(\Env_{\mathcal C} H)^\star
$$
Поскольку на $\Env_{\mathcal C} H$ определена структура инволютивной алгебры Хопфа относительно $\odot$, на сопряженном пространстве $\widehat H=(\Env_{\mathcal C} H)^\star$ определена структура инволютивной алгебры Хопфа относительно $\circledast$. Более того, $\widehat H=(\Env_{\mathcal C} H)^\star$ будет алгеброй Хопфа, рефлексивной относительно $\Env_{\mathcal C}$, потому что применение опрерации $\star$ к диаграммам \eqref{DIAG:reflex-otn-obolochki-1}-\eqref{DIAG:reflex-otn-obolochki-4} дает те же самые диаграммы, только с заменой $H$ на $\widehat H=(\Env_{\mathcal C} H)^\star$ (здесь нужно воспользоваться условием (ii) на с.\pageref{(env-H)^star:H^star-gets-(Env-H)^star}).

Объект $\widehat H=(\Env_{\mathcal C} H)^\star$ называется {\it алгеброй Хопфа, двойственной к
$H$ относительно оболочки $\Env_{\mathcal C}$}\index{алгебра!Хопфа!двойственная относительно оболочки}. Диаграмма \eqref{obshaya-diagramma-refleksivnosti} означает, что $H$ будет естественно изоморфна своей второй двойственной в этом смысле алгебре Хопфа:
 \beq\label{H-cong-(H^*)^*}
H\cong \widehat{\widehat H}
 \eeq

\subsection{Непрерывная двойственность для групп Мура.} 
Из теорем \ref{TH:E(K(G))=C(G)}, \ref{TH:@:C^star(G)-circledast-C^star(G)->C^star(G)-odot-C^star(G)}, \ref{TH:@:Env_C-C^star(G)-circledast-Env_C-C^star(G)->Env_C-C^star(G)-odot-Env_C-C^star(G)}, \ref{TH:Env_C-C^*(G)-odot-Hopf} следует

 \btm\label{TH:nepr-dvoistvennost}
Если $G$ -- группа Мура, то алгебры ${\mathcal C}^\star(G)$ и ${\mathcal
K}(G)$ непрерывно рефлексивны, а диаграмма рефлексивности для них
принимает вид \eqref{chetyrehugolnik-C-C*-0}:
 \beq\label{chetyrehugolnik-C-C*}
 \xymatrix @R=1.pc @C=2.pc
 {
 {\mathcal C}^\star(G)
 & \ar@{|->}[r]^{\Env_{\mathcal C} } & &
 \Env_{\mathcal C} {\mathcal C}^\star(G)
 \\
 & & &
 \ar@{|->}[d]^{\star}
 \\
 \ar@{|->}[u]^{\star}
 & & &
 \\
 {\mathcal C}(G)
 & &
 \ar@{|->}[l]_{\Env_{\mathcal C} }
 &
 {\mathcal K}(G)
 }
 \eeq
 \etm
\bpr
Двигаясь по цепочке \eqref{chetyrehugolnik-C-C*} от левого нижнего угла, ${\mathcal C}(G)$, мы придем к алгебре ${\mathcal K}(G)$, которая по теореме \ref{TH:E(K(G))=C(G)}, превращается в алгебру ${\mathcal C}(G)$ под действием оболочки $\Env_{\mathcal C}$, и таким образом цепочка \eqref{chetyrehugolnik-C-C*} замыкается. С другой стороны, в этой диаграмме элементы ${\mathcal C}(G)$ и $\Env_{\mathcal C}{\mathcal C}^\star(G)$ являются $\odot$-алгебрами Хопфа (${\mathcal C}(G)$ --- в силу \cite[Theorem 5.2.3]{Akbarov-De-Gruyter-I}, а $\Env_{\mathcal C}{\mathcal C}^\star(G)$ --- по теореме \ref{TH:Env_C-C^*(G)-odot-Hopf}), а элементы ${\mathcal C}^\star(G)$ и ${\mathcal K}(G)$ --- $\circledast$-алгебрами Хопфа (${\mathcal C}^\star(G)$ --- опять в силу \cite[Theorem 5.2.3]{Akbarov-De-Gruyter-I}, а ${\mathcal K}(G)$ --- опять по теореме \ref{TH:Env_C-C^*(G)-odot-Hopf}).
\epr

\noindent\rule{160mm}{0.1pt}\begin{multicols}{2}

\bex Из предложения \ref{PROP:Env_C-C^star(G)=C(widehat(G))} следует, что для абелевых локально компактных групп $C$ диаграмма рефлексивности имеет вид
 \beq\label{chetyrehugolnik-C-C*-F}
 \xymatrix @R=1.pc @C=2.pc
 {
 {\mathcal C}^\star(C)
 & \ar@{|->}[r]^{{\mathcal F}_C} & &
 {\mathcal C}(\widehat{C})
 \\
 & & &
 \ar@{|->}[d]^{\star}
 \\
 \ar@{|->}[u]^{\star}
 & & &
 \\
 {\mathcal C}(C)
 & &
 \ar@{|->}[l]_{{\mathcal F}_{\widehat{C}}}
 &
 {\mathcal C}^\star(\widehat{C})
  }
\eeq
(здесь $\widehat{C}$ -- двойственная по Понтрягину группа к $C$, а $\mathcal F$ -- преобразование Фурье, определенное выше формулой \eqref{Fourier-transform}).
\eex

\bex Из предложения \ref{PROP:nepr-obol-diskr-gruppy} мы получаем, что для дискретных групп Мура $D$ диаграмма рефлексивности выглядит так:
$$
 \xymatrix @R=1.pc @C=2.pc
 {
 \C_D
 & \ar@{|->}[r]^{\Env_{\mathcal C}} & &
 C^*(D)
 \\
 & & &
 \ar@{|->}[d]^{\star}
 \\
 \ar@{|->}[u]^{\star}
 & & &
 \\
 \C^D
 & &
 \ar@{|->}[l]_{\Env_{\mathcal C}}
 &
 {\mathcal K}(D)
 }
$$
Здесь $C^*(G)$ -- обычная $C^*$-алгебра группы $G$ \cite{Dixmier}, а алгебра ${\mathcal K}(D)$ по множеству элементов (и алгебраическим операциям) совпадает с алгеброй Фурье-Стилтьеса $B(G)$ группы $G$ \cite{Eymard}, отличаясь только более слабой топологией (мы это уже отмечали в доказательстве леммы \ref{LM:Spec-K(G)-dlya-diskretnoi-gruppy}).
\eex

\bex Из предложения \ref{PROP:nepr-obol-komp-gruppy} следует, что для компактных групп $K$ диаграмма рефлексивности становится такой:
$$
 \xymatrix @R=1.pc @C=2.pc
 {
 {\mathcal C}^\star(K)
 & \ar@{|->}[r]^{\Env_{\mathcal C}} & &
\prod_{\pi\in\widehat{K}}{\mathcal B}(X_\pi)
 \\
 & & &
 \ar@{|->}[d]^{\star}
 \\
 \ar@{|->}[u]^{\star}
 & & &
 \\
 {\mathcal C}(K)
 & &
 \ar@{|->}[l]_{\Env_{\mathcal C}}
 &
 \Trig(K)
 }
$$
\eex

\end{multicols}\noindent\rule[10pt]{160mm}{0.1pt}

\subsection{Группы, различаемые $C^*$-алгебрами.}
Условимся говорить, что локально компактная группа $G$ {\it различается $C^*$-алгебрами}\label{DEF:gruppy-razlich-C^*-algebrami}, если (непрерывные инволютивные) гомоморфизмы ее алгебры мер ${\mathcal C}^\star(G)\to B$  во всевозможные $C^*$-алгебры $B$ различают элементы $G$ (при вложении $G$ в ${\mathcal C}^\star(G)$ дельта-функциями). Понятно, что если групповая алгебра мер ${\mathcal C}^\star(G)$ непрерывно рефлексивна, то группа $G$ различается $C^*$-алгебрами, поэтому этот класс групп представляет интерес для оценки того, насколько широко обобщается теорема \ref{TH:nepr-dvoistvennost}. Из самой теоремы \ref{TH:nepr-dvoistvennost} следует, что все группы Мура различаются $C^*$-алгебрами. В работе Ю.~Н.~Кузнецовой \cite{Kuznetsova} показано, что все SIN-группы различаются $C^*$-алгебрами. Следующий результат дает представление о том, как устроены связные группы Ли, элементы которых различаются $C^*$-алгебрами:

\btm[Д.Люмине, А.Валетт \cite{Luminet-Valette}]\label{TH:Luminet-Valette} Если связная вещественная группа Ли $G$ различается $C^*$-алгебрами, то $G$ -- линейная группа (то есть она вкладывается как замкнутая подгруппа в некоторую полную линейную группу $GL_n(\R)$).
\etm

\chapter{ДВОЙСТВЕННОСТЬ В ДИФФЕРЕНЦИАЛЬНОЙ ГЕОМЕТРИИ}

\section{Гладкие оболочки}

\subsection{Присоединенные самосопряженные нильпотенты и системы частных производных}

\paragraph{Мультииндексы.}

Пусть $d\in\N$ -- натуральное число\footnote{Всюду под натуральными числами $\N$ мы понимаем целые неотрицательные числа: $\N=\{d\in\Z:\ d\ge 0\}$.}. Условимся {\it мультииндексом} длины $d$ называть произвольную конечную последовательность длины $d$ натуральных чисел
$$
k=(k_1,...,k_d),\qquad k_i\in\N.
$$
Для двух мультииндексов $k,l\in\N^d$ неравенство $l\le k$ определяется покоординатно:
$$
l\le k\qquad\Longleftrightarrow\qquad \forall i=1,...,d\quad l_i\le k_i.
$$
{\it Суммой} двух мультииндексов $k,l\in\N^d$ мы будем называть мультииндекс
$$
k+l=(k_1+l_1,...,k_d+l_d).
$$
Если же $l\le k$, то их {\it разность} определяется равенством
$$
k-l=(k_1-l_1,...,k_d-l_d).
$$
Кроме того, {\it порядок} и {\it факториал} мультииндекса $k\in\N^d$ определяются равенствами
$$
\abs{k}=k_1+...+k_d,\qquad k!=k_1!\cdot...\cdot k_d!.
$$
В соответствии с последней формулой, биномиальный коэффициент представляет собой число
$$
\begin{pmatrix}k \\ l\end{pmatrix}=\frac{k!}{l!\cdot (k-l)!}
$$

\paragraph{Алгебры степенных рядов с коэффициентами в заданной алгебре.}
Пусть $A$ -- произвольная инволютивная стереотипная алгебра. Рассмотрим алгебру
$$
A[[d]]=A^{\N^d},
$$
состоящую из всевозможных отображений $x:\N^d\to A$, или, что то же самое,
семейств $x=\{x_k;\ k\in\N^d\}$ элементов из $A$, индексированных
мультииндексами длины $d$. На $A[[d]]$ задается топология покоординатной
сходимости
$$
x^i\overset{A[[d]]}{\underset{i\to\infty}{\longrightarrow}}x\qquad\Longleftrightarrow\qquad
\forall k\in\N^d\quad
x_k^i\overset{A}{\underset{i\to\infty}{\longrightarrow}}x_k,
$$
а алгебраические операции в $A[[d]]$ -- инволюция, сумма, умножение на скаляр и
произведение -- определяются формулами
 \begin{align}
& (x^\bullet)_k=(x_k)^\bullet, && x\in A[[d]],\ k\in \N^d \label{involutsiya-v-B[[d]]}\\
& (x+y)_k=x_k+y_k, && x,y\in A[[d]],\ k\in \N^d \label{summa-v-B[[d]]}\\
& (\lambda\cdot x)_k=\lambda\cdot x_k, && \lambda\in\C, \ x\in A[[d]],\  k\in \N^d  \label{umnopzh-na-skalar-v-B[[d]]}\\
& (x\cdot y)_k=\sum_{0\le l\le k}x_{k-l}\cdot y_l, &&  x,y\in A[[d]],\ k\in \N^d \label{proizv-v-B[[d]]}
 \end{align}
Единицей в $A[[d]]$ будет, как легко понять, семейство
\beq\label{1_(B[[d]])}
\underset{\scriptsize\begin{matrix}\text{\rotatebox{90}{$\owns$}} \\ A[[d]]\end{matrix}}{1}\kern-5pt{}_k=\begin{cases}
1,& k=0 \\ 0, & k\ne 0\end{cases}
\eeq

Элементы $A[[d]]$ удобно представлять себе как степенные ряды от $d$ переменных $\tau_1,...,\tau_d$:
$$
x=\sum_{k\in\N^d}x_k\cdot\tau^k,
$$
где под $\tau^k$ понимается формальное произведение
$$
\tau^k=\tau^{k_1}\cdot...\cdot\tau^{k_d},
$$
и предполагается выполнение тождеств
$$
(\tau^k)^\bullet=\tau^k,\qquad a\cdot\tau^k=\tau^k\cdot a,\qquad \tau^k\cdot\tau^l=\tau^{k+l},\qquad a\in A,\quad \tau\in\N^d.
$$
Тогда сумма, умножение на скаляр и произведение в $A[[d]]$ описываются формулами для степенных рядов:
 \begin{align}\label{alg-operatsii-kak-ryady}
& x+y=\sum_{k\in\N^d}(x_k+y_k)\cdot\tau^k,
&& \lambda\cdot x=\sum_{k\in\N^d}(\lambda\cdot x_k)\cdot\tau^k,
&& x\cdot y=\sum_{k\in\N^d}\l\sum_{0\le l\le k}x_{k-l}\cdot y_l\r\cdot\tau^k.
 \end{align}

Другой способ -- представлять элементы $A[[d]]$ как ряды Тейлора от переменных $\tau_1,...,\tau_d$:
\beq\label{predstavlenie-ryadom-Teilora}
x=\sum_{k\in\N^d}\frac{x^{(k)}}{k!}\cdot\tau^k,
\eeq
При таком представлении коэффициенты $x^{(k)}$ ряда связаны с обычными коэффициентами $x_k$ формулами
\beq\label{predstavlenie-ryadom-Teilora-1}
x^{(k)}=k!\cdot x_k,
\eeq
а формулы для алгебраических операций \eqref{involutsiya-v-B[[d]]}-\eqref{proizv-v-B[[d]]} принимают вид
 \begin{align}\label{(xy)^(k)}
&(x^\bullet)^{(k)}=(x^{(k)})^\bullet &&
(\lambda\cdot x)^{(k)}=\lambda\cdot x^{(k)} &&
(x+y)^{(k)}=x^{(k)}+y^{(k)} &&
(x\cdot y)^{(k)}=\sum_{0\le l\le k}\begin{pmatrix}k \\ l\end{pmatrix}\cdot x^{(k-l)}\cdot y^{(l)}
 \end{align}

\paragraph{Алгебры с присоединенными самосопряженными нильпотентами.}

 \bit{\label{DEF:C^*[m]}
\item[$\bullet$] Пусть по-прежнему $A$ -- инволютивная стереотипная алгебра, $d\in\N$ и $m\in\N^d$. Обозначим через
$I_m$ замкнутый идеал в алгебре $A[[d]]$ степенных рядов с коэффициентами в
$A$, состоящий из рядов, у которых коэффициенты с индексами $k\le m$
обнуляются:
$$
I_m=\{x\in A[[d]]:\ \forall k\in\N^d\quad k\le m \quad\Longrightarrow\quad x_k=0\}.
$$
Фактор-алгебра
$$
A[m]:=A[[d]]/I_m
$$
называется {\it алгеброй $A$ с присоединенными самосопряженными нильпотентами (порядка $m$)}.
 }\eit

Обозначим символом $\N[m]$ множество мультииндексов, не превосходящих мультииндекса $m$
$$
\N[m]=\{k\in\N^d:\ k\le m\}.
$$
Тогда алгебру $A[m]$ можно представлять себе как пространство всевозможных семейств $x=\{x_k;\ k\in \N[m]\}$ элементов из $A$, индексированных мультииндексами $k\in\N[m]$. Инволюция, сумма, умножение на скаляр и произведение в $A[m]$ определяются формулами
 \begin{align}
& (x^\bullet)_k=(x_k)^\bullet, && x\in A[m],\ k\in \N[m] \label{involutsiya-v-B[m]}\\
& (x+y)_k=x_k+y_k, && x,y\in A[m],\ k\in \N[m] \label{summa-v-B[m]}\\
& (\lambda\cdot x)_k=\lambda\cdot x_k, && \lambda\in\C, \ x\in A[m],\  k\in \N[m]  \label{umnopzh-na-skalar-v-B[m]}\\
& (x\cdot y)_k=\sum_{0\le l\le k}x_{k-l}\cdot y_l, &&  x,y\in A[m],\ k\in
\N[m] \label{proizv-v-B[m]}
 \end{align}
Единицей в $A[m]$ будет, понятное дело, семейство
\beq\label{1_(B[[d]])}
\underset{\scriptsize\begin{matrix}\text{\rotatebox{90}{$\owns$}} \\
A[m]\end{matrix}}{1_k}\kern-5pt=\begin{cases} 1,& k=0 \\ 0, & k\ne
0\end{cases}
\eeq

Элементы $A[m]$ удобно представлять себе как многочлены степени $n$ от $d$
переменных $\tau_1,...,\tau_d$:
\beq\label{predstavl-B[m]-ryadami}
x=\sum_{k\in\N[m]}x_k\cdot\tau^k=\sum_{k\in\N[m]}\frac{x^{(k)}}{k!}\cdot\tau^k,
\eeq
где под $\tau^k$ понимается формальное произведение
$$
\tau^k=\tau^{k_1}\cdot...\cdot\tau^{k_d},
$$
причем считается, что
\beq\label{tozhdestva-dlya-prisoedinennyh-elementov}
(\tau_i)^\bullet=\tau_i, \qquad a\cdot\tau_i=\tau_i\cdot a,\qquad \tau_i\cdot\tau_j=\tau_j\cdot\tau_i,
\qquad\tau_i^{m_i+1}=0,\qquad a\in A,\quad i=1,...,d.
\eeq
-- при таком представлении переменные $\tau_1,...,\tau_d$ как раз можно понимать, как присоединенные к алгебре $A$ элементы, подчиненные условиям \eqref{tozhdestva-dlya-prisoedinennyh-elementov} (и это оправдывает данное выше название алгебры $A[m]$). Тогда сумма, умножение на скаляр и произведение в $A[m]$ описываются теми же формулами \eqref{alg-operatsii-kak-ryady} что и для $A[[d]]$. В частности, при представлении многочлена рядом Тейлора (второе равенство в \eqref{predstavl-B[m]-ryadami}) сохраняются формулы для алгебраических операций \eqref{(xy)^(k)}.

\btm\label{TH:norm(x)=sum_(k-le-n)norm(x_k)_B} Пусть $B$ -- инволютивная банахова алгебра с инволютивной субмультипликативной нормой $\norm{\cdot}_B$. Тогда для любого мультииндекса $n\in\N^d$ алгебра $B[n]$ является инволютивной банаховой с инволютивной субмультипликативной нормой
\beq\label{norm(x)=sum_(k-le-n)norm(x_k)_B}
\norm{x}=\sum_{k\le n}\norm{x_k}_B,\qquad x\in B[n].
\eeq
\etm
\bpr
Обозначим $M=\{(k,l)\subseteq\N[n]^2; \ l\le k\}$ и заметим, что отображение $(k,l)\in M\mapsto(k-l,l)\in\N[n]^2$ инъективно:
$$
(k,l),(k',l')\in M\quad\&\quad (k-l,l)=(k'-l',l')\quad\Longrightarrow\quad l=l'\quad\&\quad k=k'.
$$
Отсюда следует третий знак неравенства в цепочке
\begin{multline*}
\norm{x\cdot y}=\sum_{k\le n}\norm{(x\cdot y)_k}_B=\sum_{k\le n}\norm{\sum_{l\le k}x_{k-l}\cdot y_l}_B\le
\sum_{k\le n}\sum_{l\le k}\norm{x_{k-l}\cdot y_l}_B\le
\sum_{k\le n}\sum_{l\le k}\norm{x_{k-l}}_B\cdot \norm{y_l}_B=\\=\sum_{(k,l)\in M}\norm{x_{k-l}}_B\cdot \norm{y_l}_B
\le \sum_{(m,l)\in\N[n]^2}\norm{x_m}_B\cdot \norm{y_l}_B=\sum_{m\in\N[n]}\norm{x_m}_B\cdot\sum_{l\in\N[n]} \norm{y_l}_B=\norm{x}\cdot\norm{y}.
\end{multline*}
\epr

В дальнейшем нас будет интересовать исключительно случай, когда алгебра $A$, содержащая коэффициенты многочленов, является $C^*$-алгеброй. Следующий пример показывает, что свойство быть $C^*$-алгеброй не наследуется при переходе от $A$ к $A[m]$.

\noindent\rule{160mm}{0.1pt}\begin{multicols}{2}

\bex При ненулевых $m$ и $d$ алгебра $A[m]$ не может быть $C^*$-алгеброй.
\eex
\bpr
Предположим, что топология $A[m]$ порождена какой-то $C^*$-нормой. Рассмотрим какой-нибудь присоединенный элемент $\tau_i$, $i\in\{1,...,n\}$. Пусть $B$ -- замкнутая (содержащая единицу) подалгебра в $A[m]$, порожденная элементом $\tau_i$. Тогда $B$ -- коммутативная $C^*$-алгебра, и поэтому она изоморфна некоторой алгебре ${\mathcal C}(K)$ функций на компакте. С другой стороны, последнее условие в \eqref{tozhdestva-dlya-prisoedinennyh-elementov} означает, что $\tau_i$ должен быть нильпотентным элементом:
$$
\tau_i^{m_i+1}=0.
$$
Но такого не может быть, потому что в алгебрах вида $B={\mathcal C}(K)$ не бывает ненулевых нильпотентных элементов.
\epr

\end{multicols}\noindent\rule[10pt]{160mm}{0.1pt}

Для любых двух мультииндексов $m\in\N^d$ и $n\in\N^{d'}$ (необязательно, чтобы длина $d$ мультииндекса $m$ совпадала с длиной $d'$ мультииндекса $n$) определим их {\it прямую сумму} $m\oplus n$ как мультииндекс длины $d+d'$, то есть элемент пространства $\N^{d+d'}$, формулой
\beq\label{(k-oplus-l)_i}
(m\oplus n)_i=\begin{cases}m_i, & 1\le i\le d\\ n_{i-d}, & d<i\le d' \end{cases}.
\eeq
Тогда автоматически будет справедлива формула
\beq\label{N[m]-times-N[n]=N[m-oplus-n]}
\N[m]\times\N[n]=\N[m\oplus n]
\eeq

\bprop\label{PROP:A[m]-circledast-B}
Для любых двух инволютивных стереотипных алгебр $A$ и $B$ и любых мультииндексов $m\in\N^d$ и $n\in\N^{d'}$ справедливы следующие естественные изоморфизмы инволютивных стереотипных алгебр:
\beq\label{A[m][n]=A[m-oplus-n]}
\big(A[m]\big)[n]\cong A[m\oplus n]\cong A[n\oplus m]\cong \big(A[n]\big)[m]
\eeq
\beq\label{A[m]-circledast-B}
A[m]\circledast B\cong (A\circledast B)[m]\cong A\circledast \big(B[m]\big)
\eeq
Если $M$ -- паракомпактное локально компактное пространство, то
\beq\label{C(M,B)[n]=C(M,B[n])}
{\mathcal C}(M,B)[n]\cong {\mathcal C}(M,B[n]).
\eeq
Существует также естественный гомоморфизм стереотипных инволютивных алгебр
\beq\label{A[m]-oplus-B[n]-to-(A-oplus-B)[m-oplus-n]}
\ph:A[m]\oplus B[n]\to (A\oplus B)[m\oplus n]
\eeq
причем если $A$ и $B$ -- банаховы алгебры, то норма $\ph$ оценивается снизу единицей, а сверху двойкой:
\beq\label{1-le-norm(ph)-le-2}
1\le \norm{\ph}\le 2,
\eeq
если считать, что $A$ и $B$ наделены субмультипликативными нормами, и всякий раз при переходе к прямой сумме норма определяется как максимум норм слагаемых,
$$
\norm{x\oplus y}=\max\{\norm{x},\norm{y}\}
$$
а при переходе к алгебре с присоединенными самосопряженными нильпотентами -- как сумма норм компонент (по формуле \eqref{norm(x)=sum_(k-le-n)norm(x_k)_B}).
\eprop
\bpr 1. Формула \eqref{A[m][n]=A[m-oplus-n]} доказывается цепочкой
$$
\big(A[m]\big)[n]=\big(A^{\N[m]}\big)^{\N[n]}=A^{\N[m]\times\N[n]}=\eqref{N[m]-times-N[n]=N[m-oplus-n]}=A^{\N[m\oplus n]}=A[m\oplus n]
$$

2. Для доказательства \eqref{A[m]-circledast-B} определим отображение
$$
\gamma: A[m]\circledast B\to (A\circledast B)[m]
$$
формулой
\beq\label{DEF:A[m]-circledast-B->(A-circledast-B)[m]}
\gamma(x\circledast b)_k=x_k\circledast b,\qquad x\in A[m],\quad b\in B.
\eeq
(каждому семейству $x=\{x_k:\ k\in\N[m]\}$ элементов из $A$ и любому элементу $b\in B$ отображение $\gamma$ ставит в соответствие семейство $\gamma(x\circledast b)=\{\gamma(x\circledast b)_k;\ k\in\N[m]\}$ элементов из $A\circledast B$, определенное равенством \eqref{DEF:A[m]-circledast-B->(A-circledast-B)[m]}). Это отображение будет изоморфизмом стереотипных пространств, потому что тензорное произведение $\circledast$ дистрибутивно с операцией взятия прямой суммы \cite[(2.52)]{Ak16}. Проверим, что оно сохраняет умножение: для любых семейств $x=\{x_k:\ k\in\N[m]\}$ и $x'=\{x'_k:\ k\in\N[m]\}$ из $A$ и любых элементов $b,b'\in B$ мы получаем:
\begin{multline*}
\gamma\Big((x\circledast b)\cdot(x'\circledast b')\Big)_k=\gamma\big((x\cdot x')\circledast(b\cdot b')\big)_k=\eqref{DEF:A[m]-circledast-B->(A-circledast-B)[m]}=(x\cdot x')_k\circledast(b\cdot b')=\eqref{proizv-v-B[m]}=\\=\Big(\sum_{0\le l\le k}x_{k-l}\cdot x'_l\Big)\circledast (b\cdot b')=
\sum_{0\le l\le k}(x_{k-l}\cdot x'_l)\circledast (b\cdot b')=\sum_{0\le l\le k}(x_{k-l}\circledast b)\cdot (x'_l\circledast b')=\sum_{0\le l\le k}(x_{k-l}\circledast b)\cdot (x'_l\circledast b')=\\=\eqref{DEF:A[m]-circledast-B->(A-circledast-B)[m]}=\sum_{0\le l\le k}\gamma(x\circledast b)_{k-l}\cdot \gamma(x'\circledast b')_l=\eqref{proizv-v-B[m]}=\Big(\gamma(x\circledast b)\cdot\gamma(x'\circledast b')\Big)_k
\end{multline*}
Это доказывает первое равенство в \eqref{A[m]-circledast-B}. Точно так же доказывается второе.

3. Тождество \eqref{C(M,B)[n]=C(M,B[n])} очевидно.

4. Пусть $\sigma_1,...,\sigma_d$ -- последовательность присоединенных нильпотентов к алгебре $A$ в $A[m]$, а $\tau_1,...,\tau_{d'}$  -- последовательность присоединенных нильпотентов к алгебре $B$ в $B[n]$. Рассмотрим последовательность присоединенных нильпотентов к алгебре $A\oplus B$ в $(A\oplus B)[m\oplus n]$, и присвоим ее элементам обозначения $\overline{\sigma_i}$ и $\overline{\tau_j}$, расположив их так, чтобы сначала в порядке возрастания индексов шли $\overline{\sigma_i}$, а после них $\overline{\tau_j}$:
$$
\overline{\sigma_1},...,\overline{\sigma_d},\overline{\tau_1},...,\overline{\tau_{d'}}.
$$
Пусть
$$
\overline{\sigma}^k=\overline{\sigma_1}^{k_1}\cdot...\cdot \overline{\sigma_d}^{k_d},\qquad
\overline{\tau}^l=\overline{\tau_1}^{l_1}\cdot...\cdot \overline{\tau_{d'}}^{l_{d'}},\qquad k\in\N[m],\quad l\in\N[n].
$$
Тогда гомоморфизм $\ph$ в формуле \eqref{A[m]-oplus-B[n]-to-(A-oplus-B)[m-oplus-n]} можно определить правилом
$$
\ph(1_{A[m]}\oplus 1_{B[n]})=1_{(A\oplus B)[m\oplus n]},\quad \ph(\sigma_i\oplus 0_{B[n]})=(1_A\oplus 0_B)\cdot\overline{\sigma_i},\quad \ph(0_{A[m]}\oplus\tau_j)=(0_A\oplus 1_B)\cdot\overline{\tau_j}.
$$
или, что эквивалентно, правилом
$$
\ph\Big(\underbrace{\sum_{0\le k\le m}x_k\cdot\sigma^k}_{\scriptsize\begin{matrix}\text{\rotatebox{90}{$\owns$}}\\ A[m]\end{matrix}}\oplus\underbrace{\sum_{0\le l\le n} y_l\cdot\tau^l}_{\scriptsize\begin{matrix}\text{\rotatebox{90}{$\owns$}}\\ B[n]\end{matrix}}\Big)=\sum_{0\le k\le m}(\underset{\scriptsize\begin{matrix}\text{\rotatebox{90}{$\owns$}}\\ A\end{matrix}}{x_k}\oplus \underset{\scriptsize\begin{matrix}\text{\rotatebox{90}{$\owns$}}\\ B\end{matrix}}{0})\cdot\overline{\sigma}^k+ \sum_{0\le l\le n} (\underset{\scriptsize\begin{matrix}\text{\rotatebox{90}{$\owns$}}\\ A\end{matrix}}{0}\oplus \underset{\scriptsize\begin{matrix}\text{\rotatebox{90}{$\owns$}}\\ B\end{matrix}}{y_l})\cdot\overline{\tau}^l
$$
Это отображение сохраняет единицу и мультипликативно:
\begin{multline*}
\ph\Big((x\oplus y)\cdot (x'\oplus y')\Big)=\ph\left(\Big(\sum_{0\le k\le m}x_k\cdot \sigma^k\oplus\sum_{0\le l\le n}y_l\cdot\tau^l\Big) \cdot \Big(\sum_{0\le k'\le m}x'_{k'}\cdot \sigma^{k'}\oplus\sum_{0\le l'\le n}y'_{l'}\cdot\tau^{l'}\Big)\right)=\\=
\ph\left(\Big(\sum_{0\le k\le m}x_k\cdot \sigma^k \cdot \sum_{0\le k'\le m}x'_{k'}\cdot \sigma^{k'}\Big)\oplus\Big(\sum_{0\le l\le n}y_l\cdot\tau^l\cdot\sum_{0\le l'\le n}y'_{l'}\cdot\tau^{l'}\Big)\right)=
\\=
\ph\left(\sum_{0\le k\le m}\Big(\sum_{0\le p\le k}x_{k-p}\cdot x'_{p}\Big)\cdot \sigma^k\oplus\sum_{0\le l\le n}\Big(\sum_{0\le q\le l}y_{l-q}\cdot y'_q\Big)\cdot\tau^l\right)=\\=
\sum_{0\le k\le m}\Big(\sum_{0\le p\le k}x_{k-p}\cdot x'_{p}\oplus 0\Big)\cdot \overline{\sigma}^k+ \sum_{0\le l\le n}\Big(0\oplus\sum_{0\le q\le l}y_{l-q}\cdot y'_q\Big)\cdot\overline{\tau}^l=\\=
\sum_{0\le k\le m}(x_k\oplus 0)\cdot \overline{\sigma}^k \cdot \sum_{0\le k'\le m}(x'_{k'}\oplus 0)\cdot \overline{\sigma}^{k'}+\sum_{0\le l\le n}(0\oplus y_l)\cdot\overline{\tau}^l\cdot\sum_{0\le l'\le n}(0\oplus y'_{l'})\cdot\overline{\tau}^{l'}=\\=
\sum_{0\le k\le m}(x_k\oplus 0)\cdot \overline{\sigma}^k \cdot \sum_{0\le k'\le m}(x'_{k'}\oplus 0)\cdot \overline{\sigma}^{k'}+\sum_{0\le l\le n}(0\oplus y_l)\cdot\overline{\tau}^l\cdot\sum_{0\le l'\le n}(0\oplus y'_{l'})\cdot\overline{\tau}^{l'}+\\+
\underbrace{\sum_{0\le k\le m}(x_k\oplus 0)\cdot \overline{\sigma}^k \cdot\sum_{0\le l'\le n}(0\oplus y'_{l'})\cdot\overline{\tau}^{l'}}_{\scriptsize\begin{matrix}\| \\ 0\end{matrix}}+\underbrace{\sum_{0\le l\le n}(0\oplus y_l)\cdot\overline{\tau}^l\cdot
\sum_{0\le k'\le m}(x'_{k'}\oplus 0)\cdot \overline{\sigma}^{k'}}_{\scriptsize\begin{matrix}\| \\ 0\end{matrix}}=\\=
\left(
\sum_{0\le k\le m}(x_k\oplus 0)\cdot \overline{\sigma}^k \oplus \sum_{0\le l\le n}(0\oplus y_l)\cdot\overline{\tau}^l
\right)\cdot \left(
\sum_{0\le k'\le m}(x'_{k'}\oplus 0)\cdot \overline{\sigma}^{k'} \oplus \sum_{0\le l'\le n}(0\oplus y'_{l'})\cdot\overline{\tau}^{l'}
\right)=\\=
\ph\Big(\sum_{0\le k\le m}x_k\cdot \sigma^k\oplus\sum_{0\le l\le n}y_l\cdot\tau^l\Big) \cdot\ph\Big(\sum_{0\le k'\le m}x'_{k'}\cdot \sigma^{k'}\oplus\sum_{0\le l'\le n}y'_{l'}\cdot\tau^{l'}\Big)=
\ph(x\oplus y)\cdot \ph(x'\oplus y')
\end{multline*}
В оценке \eqref{1-le-norm(ph)-le-2} первое неравенство доказывается цепочкой
\begin{multline*}
\norm{\ph\Big(\sum_{0\le k\le m}x_k\cdot \sigma^k\oplus\sum_{0\le l\le n}y_l\cdot\tau^l\Big)}_{(A\oplus B)[m\oplus n]}=\norm{\sum_{0\le k\le m}(x_k\oplus 0)\cdot \overline{\sigma}^k \oplus \sum_{0\le l\le n}(0\oplus y_l)\cdot\overline{\tau}^l}_{(A\oplus B)[m\oplus n]}=\\=
\sum_{0\le k\le m}\norm{x_k\oplus 0}_{A\oplus B}+\sum_{0\le l\le n}\norm{0\oplus y_l}_{A\oplus B}=
\sum_{0\le k\le m}\norm{x_k}_A+\sum_{0\le l\le n}\norm{y_l}_B=\\=\norm{\sum_{0\le k\le m}x_k\cdot \sigma^k}_{A[m]}+\norm{\sum_{0\le l\le n}y_l\cdot\tau^l}_{B[n]}\ge
\max\left\{\norm{\sum_{0\le k\le m}x_k\cdot \sigma^k}_{A[m]},\norm{\sum_{0\le l\le n}y_l\cdot\tau^l}_{B[n]}\right\}=\\=
\norm{\sum_{0\le k\le m}x_k\cdot \sigma^k\oplus\sum_{0\le l\le n}y_l\cdot\tau^l}_{A[m]\oplus B[n]}
\end{multline*}
а второе -- цепочкой
\begin{multline*}
\norm{\ph\Big(\sum_{0\le k\le m}x_k\cdot \sigma^k\oplus\sum_{0\le l\le n}y_l\cdot\tau^l\Big)}_{(A\oplus B)[m\oplus n]}=\norm{\sum_{0\le k\le m}(x_k\oplus 0)\cdot \overline{\sigma}^k \oplus \sum_{0\le l\le n}(0\oplus y_l)\cdot\overline{\tau}^l}_{(A\oplus B)[m\oplus n]}=\\=
\sum_{0\le k\le m}\norm{x_k\oplus 0}_{A\oplus B}+\sum_{0\le l\le n}\norm{0\oplus y_l}_{A\oplus B}=
\sum_{0\le k\le m}\norm{x_k}_A+\sum_{0\le l\le n}\norm{y_l}_B=\\=\norm{\sum_{0\le k\le m}x_k\cdot \sigma^k}_{A[m]}+\norm{\sum_{0\le l\le n}y_l\cdot\tau^l}_{B[n]}\le
2\max\left\{\norm{\sum_{0\le k\le m}x_k\cdot \sigma^k}_{A[m]},\norm{\sum_{0\le l\le n}y_l\cdot\tau^l}_{B[n]}\right\}=\\=
2\norm{\sum_{0\le k\le m}x_k\cdot \sigma^k\oplus\sum_{0\le l\le n}y_l\cdot\tau^l}_{A[m]\oplus B[n]}
\end{multline*}
\epr

\paragraph[Системы частных производных.]{Системы частных производных и морфизмы со значениями в алгебрах степенных рядов.}

Пусть $A$ -- инволютивная стереотипная алгебра, а $B$ -- $C^*$-алгебра. Система операторов $D_k:A\to B$, $k\in\N[m]$, называется {\it системой частных производных} на алгебре $A$ с коэффициентами в алгебре $B$, если она удовлетворяет условиям:
\begin{align}
& D_k(a^\bullet)=D_k(a)^\bullet, && k\in\N[m], \label{DEF:sist-chastn-proizv-*} \\
& D_k(1)=\begin{cases}1, & k=0
\\
0, & k\ne 0
\end{cases}, && k\in\N[m],\ 1\in A, \label{DEF:sist-chastn-proizv-1} \\
& D_k(a\cdot b)=\sum_{0\le l\le k}\begin{pmatrix}k \\ l\end{pmatrix}\cdot D_{k-l}(a)\cdot D_l(b), && k\in\N[m],\ a,b\in A. \label{DEF:sist-chastn-proizv-2}
\end{align}
В частности, это означает, что оператор $D_0:A\to B$ должен быть инволютивным гомоморфизмом алгебр,
\beq\label{DEF:sist-chastn-proizv-0}
D_0(a^\bullet)=D_0(a)^\bullet,\qquad D_0(1)=1,\qquad  D_0(a\cdot b)=D_0(a)\cdot
D_0(b),\qquad a,b\in A.
\eeq
А в случае $\abs{k}=1$ операторы $D_k:A\to B$ должны быть дифференцированиями относительно гомоморфизма $D_0$:
$$
 D_k(a\cdot b)=D_k(a)\cdot D_0(b)+D_0(a)\cdot D_k(b),\qquad a,b\in A.
$$

\btm\label{TH:homomorphizm-v-B[m]<->chast-proizv-v-B} Для любой инволютивной
стереотипной алгебры $A$ и любой $C^*$-алгебры $B$ формула\footnote{Здесь $D(a)^{(k)}$ означает коэффициент в разложении Тейлора \eqref{predstavlenie-ryadom-Teilora} элемента $D(a)\in B[m]$.}
 \beq\label{homomorphizm-v-B[m]<->chast-proizv-v-B}
D_k(a)=D(a)^{(k)},\qquad k\in\N[m],\quad a\in A,
 \eeq
или, что эквивалентно, формула
 \beq\label{D(a)=sum-D_k(a)/k!-tau^k}
D(a)=\sum_{k\in\N[m]}\frac{D_k(a)}{k!}\cdot\tau^k,\qquad a\in A,
 \eeq
устанавливают взаимно однозначное соответствие между гомоморфизмами
инволютивных стереотипных алгебр $D:A\to B[m]$ и системами частных
производных $\{D_k;\ k\in\N[m]\}$ на $A$ с коэффициентами в $B$. \etm

\bpr 1. Если $D:A\to B[m]$ -- гомоморфизм инволютивных стереотипных алгебр, то определив отображения
$\{D_k;\ k\in\N[m]\}:A\to B$ по формуле
\eqref{homomorphizm-v-B[m]<->chast-proizv-v-B}, мы получим,  во-первых, для
любого $a\in A$
$$
D_k(a^\bullet)=D(a^\bullet)^{(k)}=\big(D(a)^\bullet\big)^{(k)}=\eqref{(xy)^(k)}=\big(D(a)^{(k)}\big)^\bullet=\big(D_k(a)\big)^\bullet,
$$
во-вторых,
$$
D_k(1)=D(1)^{(k)}=1^{(k)}=1_k=\eqref{1_(B[[d]])}=\begin{cases}1, & k=0
\\
0, & k\ne 0
\end{cases},
$$
и, в-третьих, для любых $a,b\in A$
\begin{multline*}
D_k(a\cdot b)=\eqref{homomorphizm-v-B[m]<->chast-proizv-v-B}=D(a\cdot b)^{(k)}=\big(D(a)\cdot D(b)\big)^{(k)}=\eqref{(xy)^(k)}=\\=
\sum_{0\le l\le k}\begin{pmatrix}k \\ l\end{pmatrix}\cdot  D(a)^{(k-l)}\cdot D(b)^{(l)}=\eqref{homomorphizm-v-B[m]<->chast-proizv-v-B}=\sum_{0\le l\le k}\begin{pmatrix}k \\ l\end{pmatrix}\cdot  D_{k-l}(a)\cdot D_l(b).
\end{multline*}
То есть выполняются тождества \eqref{DEF:sist-chastn-proizv-*}, \eqref{DEF:sist-chastn-proizv-1} и \eqref{DEF:sist-chastn-proizv-2}, и это значит, что семейство $\{D_k;\ k\in\N[m]\}$ есть система частных производных на $A$ с коэффициентами в $B$.

2. Наоборот, если $\{D_k;\ k\in\N[m]\}$ -- система частных производных на $A$ с коэффициентами в $B$, то, определив отображение $D:A\to B[m]$ формулой \eqref{D(a)=sum-D_k(a)/k!-tau^k}, мы получим, во-первых, для любого $a\in A$
$$
D(a^\bullet)^{(k)}=D_k(a^\bullet)=\eqref{DEF:sist-chastn-proizv-*}=D_k(a)^\bullet=\big(D(a)^{(k)}\big)^\bullet,
$$
во-вторых,
$$
D(1)^{(k)}=D_k(1)=\eqref{DEF:sist-chastn-proizv-1}=\begin{cases}1, & k=0
\\
0, & k\ne 0
\end{cases}=\eqref{1_(B[[d]])}=1_k\qquad\Longrightarrow\qquad D(1)=1,
$$
и, в-третьих, для любых $a,b\in A$
\begin{multline*}
D(a\cdot b)^{(k)}=\eqref{homomorphizm-v-B[m]<->chast-proizv-v-B}=D_k(a\cdot b)=\eqref{DEF:sist-chastn-proizv-2}=\sum_{0\le l\le k}\begin{pmatrix}k \\ l\end{pmatrix}\cdot D_{k-l}(a)\cdot D_l(b)=\eqref{homomorphizm-v-B[m]<->chast-proizv-v-B}=\\=\sum_{0\le l\le k}\begin{pmatrix}k \\ l\end{pmatrix}\cdot  D(a)^{(k-l)}\cdot D(b)^{(l)}=\eqref{(xy)^(k)}=\big(D(a)\cdot D(b)\big)^{(k)} \qquad\Longrightarrow\qquad D(a\cdot b)=D(a)\cdot D(b).
\end{multline*}
То есть $D:A\to B[m]$ является инволютивным гомоморфизмом.
\epr

\noindent\rule{160mm}{0.1pt}\begin{multicols}{2}

\bex\label{EX:sist-chast-proizv-na-C-infty(M)} Пусть $M$ -- гладкое многообразие  размерности $d\in\N$, $\ph:U\to V$ -- локальная карта, где $U\subseteq M$, $V\subseteq\R^d$ и $K\subseteq U$ -- компакт. Тогда система операторов
\beq\label{sist-chast-proizv-na-C-infty(M)}
D_k:{\mathcal E}(M)\to\mathcal{C}(K)\quad\Big|\quad
D_k(a)=\frac{\partial^{\abs{k}}(a\circ\ph^{-1})}{\partial t_1^{k_1}...\partial
t_d^{k_d}}\circ\ph
\eeq
является системой частных производных на ${\mathcal E}(M)$.
Соответствующий ей гомоморфизм алгебр
$D:{\mathcal E}(M)\to\mathcal{C}(K)[m]$ представляет собой систему
ограничений на компакт $K$ частных производных данной функции:
$$
(D(a))_k=D_k(a)\big|_K
$$
\eex

\end{multicols}\noindent\rule[10pt]{160mm}{0.1pt}

\paragraph{Частные производные как дифференциальные операторы.}

Пусть, по-прежнему, $A$ -- стереотипная алгебра с инволюцией, а $B$ -- $C^*$-алгебра, и $\{D_k;\ k\in\N[m]\}$ -- система частных производных на $A$ с коэффициентами в $B$. Тогда гомоморфизм $\ph=D_0:A\to B$ превращает $B$ в модуль над $A$, и поэтому для всякого оператора $P:A\to B$ и любого элемента $a\in A$ определен коммутатор $[P,a]:A\to B$.

\bprop Для всякой системы частных производных $\{D_k;\ k\in\N[m]\}$ на $A$ с коэффициентами в $B$ коммутатор операторов $D_k$ с произвольным элементом $a\in A$ относительно гомоморфизма $D_0:A\to B$ действует по формуле
\beq\label{[D_k,a]}
[D_k,a]=\sum_{0\le l< k}\begin{pmatrix}k \\ l\end{pmatrix}\cdot D_{k-l}(a)\cdot D_l
\eeq
причем, при $k=0$ эта формула приобретает вид
\beq\label{[D_0,a]}
[D_0,a]=0.
\eeq
\eprop
\bpr Тождество \eqref{[D_0,a]} выполняется тривиально, поскольку $D_0$ -- гомоморфизм. А \eqref{[D_k,a]} доказывается цепочкой
$$
[D_k,a](x)=D_k(a\cdot x)-\ph(a)\cdot D_k(x)=\sum_{0\le l\le k}\begin{pmatrix}k \\ l\end{pmatrix}\cdot D_{k-l}(a)\cdot D_l(x)-D_0(a)\cdot D_k(x)=
\sum_{0\le l< k}\begin{pmatrix}k \\ l\end{pmatrix}\cdot D_{k-l}(a)\cdot D_l(x).
$$
\epr

\btm\label{TH:differentsialnye-chast-proizv} Для любой системы частных производных $\{D_k;\ k\in\N[m]\}$ из инволютивной стереотипной алгебры $A$ в $C^*$-алгебру $B$ следующие условия эквивалентны:
 \bit{
 \item[(i)] операторы $\{D_k;\ k\in\N[m]\}$ являются дифференциальными операторами\footnote{Понятие дифференциального оператора было определено на с.\pageref{DEF:diff-oper}.} из $A$ в $B$ относительно гомоморфизма $D_0:A\to B$ с порядками, не превосходящими модулей своих индексов:
     \beq\label{D_k-in-D^|k|(D_0)}
     D_k\in \Diff^{\abs{k}}(D_0),
     \eeq

 \item[(ii)] для любого мультииндекса $k>0$ значения оператора $D_k$ лежат в пространстве $Z^{|k|}(D_0)$:
     \beq\label{D_k(A)-subseteq-Z^|k|(D_0)}
     D_k(A)\subseteq Z^{|k|}(D_0),\qquad k>0.
     \eeq

\item[(iii)] для любого мультииндекса $k>0$ значения оператора $D_k$ лежат в пространстве $Z^1(D_0)=D_0(A)^!$ (коммутанте образа оператора $D_0$):
     \beq\label{D_k(A)-subseteq-Z^1(D_0)}
     D_k(A)\subseteq Z^1(D_0)=D_0(A)^!,\qquad k>0.
     \eeq
  }\eit
\etm

 \bit{
\item[$\bullet$]\label{DEF:differentsialnye-chast-proizv} Систему частных
производных $\{D_k;\ k\in\N[m]\}$ на $A$ с коэффициентами в $B$ мы будем называть {\it
дифференциальной}, если она удовлетворяет эквивалентным условиям (i)-(iii)
теоремы \ref{TH:differentsialnye-chast-proizv}.

\item[$\bullet$] Гомоморфизм инволютивных стереотипных алгебр $D:A\to B[m]$ называется {\it
дифференциальным}, если определяемая им по формуле \eqref{homomorphizm-v-B[m]<->chast-proizv-v-B} система частных производных $\{D_k;\ k\in\N[m]\}$ является дифференциальной. Класс всех дифференциальных гомоморфизмов мы обозначаем $\DiffMor$.
 }\eit

\bpr Заметим сразу, что равносильность (ii) и (iii) есть следствие формулы \eqref{Z^n(ph)=Z^(n+1)(ph)}. Поэтому нужно проверить равносильность (i) и (ii).

1. (i)$\Longrightarrow$(ii). Пусть выполнено (i). Мы докажем (ii) по индукции.

1) Пусть $|k|=1$. Тогда для любых $a,a_1\in A$ мы получим:
$$
[D_k,a]=\eqref{[D_k,a]}=\sum_{0\le l< k}\begin{pmatrix}k \\ l\end{pmatrix}\cdot D_{k-l}(a)\cdot D_l=D_k(a)\cdot D_0
$$
$$
\Downarrow
$$
$$
0=[[D_k,a],a_1]=[D_k(a)\cdot D_0,a_1]=\eqref{[b-cdot-P, a]}=D_k(a)\cdot \kern-6pt\underbrace{[D_0,a_1]}_{\tiny\begin{matrix}
\phantom{\eqref{[D_0,a]}}
\ \text{\rotatebox{90}{$=$}}\ \eqref{[D_0,a]} \\ 0\end{matrix}}\kern-6pt +[D_k(a),D_0(a_1)]\cdot D_0
$$
$$
\Downarrow
$$
$$
[D_k(a),D_0(a_1)]=0
$$
Последнее верно для любого $a_1$, поэтому $D_k(a)\in Z^1(D_0)$. Это в свою очередь верно для любого $a\in A$, поэтому $D_k(A)\subseteq Z^1(D_0)$.

2) Предположим, что мы уже доказали \eqref{D_k(A)-subseteq-Z^|k|(D_0)} для всех $k$ таких, что $|k|\le n$:
     \beq\label{D_k(A)-subseteq-Z^|k|(D_0),n}
     D_k(A)\subseteq Z^{|k|}(D_0),\qquad 0<|k|\le n.
     \eeq
Тогда для $|k|=n+1$ мы получим: при любом $a\in A$
$$
\underbrace{[D_k,a]}_{\tiny\begin{matrix}
\phantom{\eqref{[D^n,A]-subseteq-D^(n-1)}}
\ \text{\rotatebox{90}{$\owns$}}\ \eqref{[D^n,A]-subseteq-D^(n-1)} \\ \Diff^{\abs{k}-1}\\ \text{\rotatebox{90}{$=$}} \\ \Diff^n \end{matrix}}\kern-5pt=\eqref{[D_k,a]}=\sum_{0\le l< k}\begin{pmatrix}k \\ l\end{pmatrix}\cdot D_{k-l}(a)\cdot D_l=D_k(a)\cdot D_0+
\sum_{0<l< k}\begin{pmatrix}k \\ l\end{pmatrix}\cdot\underbrace{\kern-5pt\overbrace{D_{k-l}(a)}^{\tiny\begin{matrix} Z^{|k-l|} \\
\eqref{D_k(A)-subseteq-Z^|k|(D_0),n}
\ \text{\rotatebox{90}{$\in$}}\ \phantom{\eqref{D_k(A)-subseteq-Z^|k|(D_0),n}} \end{matrix}}\kern-5pt\cdot\kern-15pt \overbrace{D_l}^{\tiny\begin{matrix} \Diff^{\abs{l}} \\
\phantom{\eqref{D_k-in-D^|k|(D_0)}}
\ \text{\rotatebox{90}{$\in$}}\ \eqref{D_k-in-D^|k|(D_0)} \end{matrix}}\kern-15pt}_{\tiny\begin{matrix}
\phantom{\eqref{Z^q.D^p->D^(q+p-1)}}
\ \text{\rotatebox{90}{$\owns$}}\ \eqref{Z^q.D^p->D^(q+p-1)} \\ \Diff^{\abs{k-l}+|l|-1}\\ \text{\rotatebox{90}{$=$}} \\ \Diff^n \end{matrix}}
$$
\vglue-40pt
$$
\Downarrow
$$
$$
D_k(a)\cdot D_0\in \Diff^n
$$
$$
\phantom{\tiny \eqref{b.ph-in-D^n(ph)<=>b-in-Z^(n+1)(ph)}}\quad \Downarrow \quad {\tiny \eqref{b.ph-in-D^n(ph)<=>b-in-Z^(n+1)(ph)}}
$$
$$
D_k(a)\in Z^{n+1}=Z^{|k|}
$$

2. (i)$\Longleftarrow$(ii). Пусть наоборот, выполнено (ii), тогда (i) доказывается также по индукции.

0) Для $k=0$ утверждение \eqref{D_k-in-D^|k|(D_0)} выполняется вообще всегда, потому что гомоморфизм $\ph=D_0$ является дифференциальным оператором нулевого порядка в силу формулы \eqref{[ph,a]=0}.

1) Предположим, что мы доказали \eqref{D_k-in-D^|k|(D_0)} для всех $k$ таких, что $|k|\le n$:
     \beq\label{D_k-in-D^|k|(D_0),n}
     D_k\in \Diff^{\abs{k}}(D_0),\qquad |k|\le n.
     \eeq
Тогда для $|k|=n+1$ мы получим:
$$
\forall  a\in A\qquad
[D_k,a]=\eqref{[D_k,a]}=
\sum_{0\le l< k}\begin{pmatrix}k \\ l\end{pmatrix}\cdot\underbrace{\kern-5pt\overbrace{D_{k-l}(a)}^{\tiny\begin{matrix} Z^{|k-l|} \\
\eqref{D_k(A)-subseteq-Z^|k|(D_0)}
\ \text{\rotatebox{90}{$\in$}}\ \phantom{\eqref{D_k(A)-subseteq-Z^|k|(D_0)}} \end{matrix}}\kern-5pt\cdot\kern-15pt \overbrace{D_l}^{\tiny\begin{matrix} \Diff^{\abs{l}} \\
\phantom{\eqref{D_k-in-D^|k|(D_0),n}}
\ \text{\rotatebox{90}{$\in$}}\ \eqref{D_k-in-D^|k|(D_0),n} \end{matrix}}\kern-15pt}_{\tiny\begin{matrix}
\phantom{\eqref{Z^q.D^p->D^(q+p-1)}}
\ \text{\rotatebox{90}{$\owns$}}\ \eqref{Z^q.D^p->D^(q+p-1)} \\ \Diff^{\abs{k-l}+|l|-1}\\ \text{\rotatebox{90}{$=$}} \\ \Diff^n \end{matrix}},
$$
$$
\Downarrow
$$
$$
\forall  a\in A\qquad [D_k,a]\in \Diff^n,
$$
$$
\Downarrow
$$
$$
D_k\in \Diff^{n+1}.
$$
\epr

\noindent\rule{160mm}{0.1pt}\begin{multicols}{2}

\bex\label{EX:ph-diff-morfizm}
Всякий гомоморфизм $D:A\to B$ инволютивной стереотипной алгебры $A$ в $C^*$-алгебру $B$ является дифференциальным, потому что в формуле \eqref{homomorphizm-v-B[m]<->chast-proizv-v-B} $k=0$, и поэтому система частных производных $\{D_k;\ k\in\N[m]\}$ состоит только из морфизма $D_0=D$, который является дифференциальным оператором порядка 0 относительно самого себя в силу примера \ref{EX:ph-diff-oper-por-0}.
\eex

\end{multicols}\noindent\rule[10pt]{160mm}{0.1pt}

\medskip
\centerline{\bf Свойства дифференциальных гомоморфизмов:}

\bit{\it

\item[$1^\circ$.]\label{1^0-diff-morf} Класс $\DiffMor$ дифференциальных морфизмов выходит\footnote{См. определение на с.\pageref{DEF:goes-from}.} из категории $\InvSteAlg$: какой бы ни была инволютивная стереотипная алгебра $A$, всегда найдется выходящий из нее дифференциальный морфизм $D:A\to B[m]$.
    
\item[$2^\circ$.]\label{2^0-diff-morf} Класс $\DiffMor$ дифференциальных морфизмов является левым идеалом в $\InvSteAlg$: если $\sigma:A\to B$ --- гомоморфизм инволютивных стереотипных алгебр, а $D:B\to C[m]$ --- дифференциальный гомоморфизм, то композиция $D\circ\sigma:A\to C[m]$ --- тоже дифференциальный гомоморфизм.

\item[$3^\circ$.]\label{3^0-diff-morf} Класс $\DEpi$ плотных эпиморфизмов подталкивает\footnote{См. определение на с.\pageref{Omega-podderzhivaet-Phi}.} класс $\DiffMor$ дифференциальных морфизмов: если $\sigma:A\to B$ --- плотный эпиморфизм инволютивных стереотипных алгебр, а $D:B\to C[m]$ --- гомоморфизм в $C^*$-алгебру с присоединенными нильпотентными элементами, причем композиция $D\circ\sigma:A\to C[m]$ является дифференциальным гомоморфизмом, то $D:B\to C[m]$ --- тоже дифференциальный гомоморфизм.  
}\eit

\bpr
1. В утверждении $1^\circ$ можно в качестве $B$ взять нулевую алгебру $B=0$ и положить $m=0$ и $D=0$.

2. Пусть $D\in\DiffMor$, $D:B\to C[m]$, и $\sigma\in\Mor$, $\sigma:A\to B$. Тогда $D\circ\sigma:A\to C[m]$. Для всякого $k>0$ мы получаем цепочку:
\begin{multline*}
(D\circ\sigma)_k(A)=D(\sigma(A))^{(k)}\subseteq D(B)^{(k)}=D_k(B)\overset{\eqref{D_k(A)-subseteq-Z^1(D_0)}}{\subseteq} D_0(B)^!\kern-30pt\overset{\scriptsize\begin{matrix}B\supseteq\sigma(A)\\ \Downarrow \\ D_0(B)\supseteq D_0(\sigma(A))\\ \Downarrow\end{matrix}}{\subseteq}\kern-30pt D_0(\sigma(A))^!=\\=
\Big(D(\sigma(A))^{(0)}\Big)^!=\Big((D\circ\sigma)(A)^{(0)}\Big)^!=
(D\circ\sigma)_0(A)^!
\end{multline*}
Таким образом, морфизм $D\circ\sigma:A\to C[m]$ удовлетворяет условию (iii) теоремы \ref{TH:differentsialnye-chast-proizv}, то есть $D\circ\sigma$ -- дифференциальный морфизм.

3. Пусть морфизмы $D:B\to C[m]$, и $\sigma:A\to B$ такие, что $\sigma\in\DEpi$ и $D\circ\sigma\in\DiffMor$. Тогда для всякого $k>0$ мы получим:
$$
D_k(\sigma(A))=(D\circ\sigma)_k(A)\subseteq (D\circ\sigma)_0(A)^!=\Big((D\circ\sigma)(A)^{(0)}\Big)^!=\Big(D(\sigma(A))^{(0)}\Big)^! \kern-50pt\overset{\scriptsize\begin{matrix}\overline{\sigma(A)}=B\\ \Downarrow \\ D_0(\sigma(A))^!= D_0\Big(\overline{\sigma(A)}\Big)^!=D_0(B)^!\\ \Downarrow\end{matrix}}{\subseteq}\kern-50pt D_0(B)^!
$$
$$
\Downarrow
$$
$$
D_k(B)=D_k\Big(\overline{\sigma(A)}\Big)\subseteq\overline{D_k(\sigma(A))}\subseteq D_0(B)^!.
$$
Таким образом, морфизм $D:B\to C[m]$ удовлетворяет условию (iii) теоремы \ref{TH:differentsialnye-chast-proizv}, то есть $D$ -- дифференциальный морфизм.

\epr

\subsection{Гладкие оболочки}

\paragraph{Определение гладкой оболочки и функториальность.}

\bit{

\item[$\bullet$] {\it Гладкой оболочкой}\label{DEF:C^infty-obolochka} $\env_{\mathcal E} A:A\to\Env_{\mathcal E} A$ инволютивной стереотипной алгебры $A$ называется ее оболочка в классе $\DEpi$ плотных эпиморфизмов категории $\InvSteAlg$ инволютивных стереотипных алгебр относительно класса $\DiffMor$ дифференциальных гомоморфизмов во всевозможные $C^*$-алгебры $B[m]$ с присоединенными самосопряженными нильпотентными элементами:
  $$
  \Env_{\mathcal E} A=\Env_{\DiffMor}^{\DEpi}A
  $$
}\eit

Более подробно, под {\it гладким расширением} инволютивной стереотипной алгебры $A$ понимается плотный эпиморфизм $\sigma:A\to A'$ инволютивных стереотипных алгебр такой, что для любой $C^*$-алгебры $B$, любого мультииндекса $m\in\N^d$ и любого дифференциального инволютивного гомоморфизма $\ph:A\to B[m]$ найдется единственный гомоморфизм инволютивных стереотипных алгебр $\ph':A'\to B[m]$, замыкающий диаграмму
$$
 \xymatrix @R=2pc @C=1.2pc
 {
  A\ar[rr]^{\sigma}\ar[dr]_{\ph} & & A'\ar@{-->}[dl]^{\ph'} \\
  & B[m] &
 }
$$
А под {\it гладкой оболочкой} алгебры $A$ понимается такое гладкое расширение $\rho:A\to \Env_{\mathcal E} A$, что для любого гладкого расширения $\sigma:A\to A'$ найдется единственный гомоморфизм инволютивных стереотипных алгебр $\upsilon:A'\to \Env_{\mathcal E} A$, замыкающий диаграмму
$$
 \xymatrix @R=2pc @C=1.2pc
 {
  & A\ar[ld]_{\sigma}\ar[rd]^{\rho} &   \\
  A'\ar@{-->}[rr]_{\upsilon} &  & \Env_{\mathcal E} A
 }
$$

\btm\label{TH:Env_C^infty-reg-obolochka} Гладкая оболочка $\Env_{\mathcal E}$ регулярна и согласована с проективным тензорным произведением\footnote{В смысле определений на страницах \pageref{DEF:reg-obolochka} и \pageref{DEF:obolochka-soglasovana-s-tenz-proizv}.} $\circledast$ в $\InvSteAlg$.
\etm
\bpr Здесь используются те же приемы, что и в теореме \ref{TH:Env_C-reg-obolochka}, причем для второй части утверждения -- согласованности с тензорным произведением -- доказательство в точности повторяется, а для первой части отличия появляются только в последних двух пунктах: условие RE.4 уже зафиксировано в свойствах $1^\circ$ и $2^\circ$ на с.\pageref{1^0-diff-morf}, а условие RE.5 --- в свойстве $3^\circ$ на с.\pageref{3^0-diff-morf}.
\epr

\bcor\label{COR:Env_E-idemp-funktor}
Гладкую оболочку можно определить как идемпотентный ковариантный функтор из $\InvSteAlg$ в $\InvSteAlg$: существуют
\bit{
\item[1)] отображение $A\mapsto (\Env_{\mathcal E}A,\env_{\mathcal E}A)$, сопоставляющее каждой инволютивной стереотипной алгебре $A$ инволютивную стереотипную алгебру $\Env_{\mathcal E}A$ и морфизм инволютивных стереотипных алгебр $\env_{\mathcal E}A:A\to\Env_{\mathcal E}A$, являющийся гладкой оболочкой алгебры $A$, и

\item[2)] отображение $\ph\mapsto \Env_{\mathcal E}(\ph)$, сопоставляющее каждому морфизму инволютивных стереотипных алгебр $\ph:A\to B$ морфизм инволютивных стереотипных алгебр $\Env_{\mathcal E}(\ph):\Env_{\mathcal E}A\to \Env_{\mathcal E}B$, замыкающий диаграмму
\beq\label{DIAGR:funktorialnost-env_E}
\xymatrix @R=2.pc @C=5.0pc 
{
A\ar[d]^{\ph}\ar[r]^{\env_{\mathcal E}A} & \Env_{\mathcal E}A\ar@{-->}[d]^{\Env_{\mathcal E}(\ph)} \\
B\ar[r]^{\env_{\mathcal E}B} & \Env_{\mathcal E}B \\
}
\eeq
}\eit
причем выполняются тождества
\beq\label{funktorialnost-env_E-v-Ste^circledast}
\Env_{\mathcal E}(1_A)=1_{\Env_{\mathcal E}A},\quad \Env_{\mathcal E}(\beta\circ\alpha)=\Env_{\mathcal E}(\beta)\circ \Env_{\mathcal E}(\alpha),
\eeq
\beq\label{funktorialnost-env_E-v-Ste^circledast-2}
\Env_{\mathcal E}(\Env_{\mathcal E}A)=\Env_{\mathcal E}A,\quad \env_{\mathcal E}\Env_{\mathcal E}A=1_{\Env_{\mathcal E}A},
\eeq
\beq\label{Env_E(C)=C}
\Env_{\mathcal E}\C=\C
\eeq
\ecor

\btm[связь с непрерывной оболочкой]\label{7:Env_infty->Env} Гладкая оболочка вкладывается в непрерывную оболочку: для всякой инволютивной стереотипной алгебры $A$ существует единственный морфизм инволютивных стереотипных алгебр $\zeta_A:\Env_{\mathcal E}A\to\Env_{\mathcal C} A$, замыкающий диаграмму
\beq\label{Env_infty->Env}
\begin{diagram}
\node[2]{A} \arrow{sw,t}{\env_{\mathcal E} A} \arrow{se,t}{\env_{\mathcal C} A}\\
\node{\Env_{\mathcal E} A}\arrow[2]{e,b,--}{\zeta_A}   \node[2]{\Env_{\mathcal C} A}
\end{diagram}
\eeq
Этот морфизм $\zeta_A$ всегда является плотным эпиморфизмом.
\etm
\bpr
Это следует из \cite[(3.10)]{Ak16}.
\epr

\paragraph{Сеть дифференциальных фактор-отображений в категории $\InvSteAlg$.}
Конструкцию гладкой оболочки можно описать несколько более наглядно следующим образом. Условимся окрестность нуля $U$ в алгебре $A$ называть {\it дифференциальной}, если она является прообразом единичного шара при некотором дифференциальном гомоморфизме $D:A\to B[m]$ в $C^*$-алгебру $B$ с присоединенными нильпотентными элементами:
$$
U=\{x\in A:\ \norm{D(x)}\le 1\}
$$
(норма $\norm{\cdot}$ определена формулой \eqref{norm(x)=sum_(k-le-n)norm(x_k)_B}). У каждой дифференциальной окрестности нуля $U$ в $A$ ядро
$$
\Ker U=\bigcap_{\lambda>0}\lambda\cdot U
$$
совпадает с ядром гомоморфизма $D$ и поэтому является замкнутым идеалом в $A$. Рассмотрим фактор-алгебру $A/\Ker U$ и наделим ее нормой, в которой класс $U+\Ker U$ является единичным шаром. Эту алгебру $A/\Ker U$ можно считать подалгеброй в $B[m]$ с индуцированной из этого пространства нормой. Ее пополнение
\beq\label{A/U=(A/Ker U)^blacktriangledown}
A/U=(A/\Ker U)^\blacktriangledown
\eeq
является банаховой алгеброй (и ее можно считать замкнутой подалгеброй в $B[m]$). Условимся называть $A/U$ {\it фактор-алгеброй алгебры $A$ по дифференциальной окрестности нуля $U$}, а соответствующее отображение
$$
\pi_U:A\to A/U
$$
-- {\it фактор-отображением} алгебры $A$ по дифференциальной окрестности нуля $U$, или {\it дифференциальным фактор-отображением} алгебры $A$.

Множество всех дифференциальных окрестностей нуля в алгебре $A$ мы будем обозначать символом ${\mathcal U}_{C^\infty}^A$,
а множество всех дифференциальных фактор-отображений алгебры $A$ мы будем обозначать символом ${\mathcal N}_{C^\infty}^A$.

\medskip
\centerline{\bf Свойства дифференциальных окрестностей нуля в категории $\InvSteAlg$:}

\bit{\it

\item[$1^\circ$.]\label{LM:diff-okr-nulya-predstavlyayut-diff-morfizmy}
Для всякого дифференциального гомоморфизма $\ph:A\to B[n]$ найдутся дифференциальная окрестность нуля $U\subseteq A$ и гомоморфизм $\ph_U:A/U\to B[n]$, замыкающий диаграмму
\beq\label{ph=ph_U-circ-pi_U}
 \xymatrix @R=2pc @C=1.2pc
 {
  A\ar[rr]^{\ph}\ar[dr]_{\pi_U} & & B[n]  \\
  & A/U\ar@{-->}[ur]_{\ph_U} &
  }
\eeq

\item[$2^\circ$.]\label{LM:diff-okr-nulya-svyazany}
Если $U$ и $U'$ -- две дифференциальные окрестности нуля в $A$, причем $U\supseteq U'$, то найдется единственный гомоморфизм $\varkappa^{U'}_U:A/U\gets A/U'$, замыкающий диаграмму
\beq\label{varkappa^U'_U}
 \xymatrix @R=2pc @C=1.2pc
 {
  & A\ar[ld]_{\pi_U}\ar[rd]^{\pi_{U'}} &   \\
  A/U &  & A/ U'\ar@{-->}[ll]^{\varkappa^{U'}_U}
 }
\eeq

\item[$3^\circ$.]\label{LM:diff-okr-nulya-uporyadocheny} Любые две дифференциальные окрестности нуля $U$ и $U'$ в $A$ содержат некую третью дифференциальную окрестность нуля $U''$:
$$
U\cap U'\supseteq U''.
$$

}\eit

\bpr Здесь,  как и в случае со свойствами на с.\pageref{LM:ph=ph_U-circ-pi_U-0}, не вполне очевидно свойство $3^\circ$. Докажем его. Пусть $D:A\to B[m]$ и $D':A\to B'[m']$ -- дифференциальные гомоморфизмы, порождающие $U$ и $U'$.
$$
U=\{x\in A:\ \norm{D(x)}\le 1\},\qquad U'=\{x\in A:\ \norm{D'(x)}\le 1\}.
$$
Рассмотрим гомоморфизм
$$
\ph:B[m]\oplus B'[m']\to (B\oplus B')[m\oplus m'],
$$
описываемый в \eqref{A[m]-oplus-B[n]-to-(A-oplus-B)[m-oplus-n]}. Отображение
$$
D'':A\to (B\oplus B')[m\oplus m']\quad\Big|\quad D''(x)=\ph(D(x)\oplus D'(x)), \quad x\in A,
$$
будет гомоморфизмом, и если положить
$$
U''=\{x\in A:\ \norm{D''(x)}\le 1\},
$$
то для всякого $x\in U''$ мы получим
$$
1\ge \norm{D''(x)}=\norm{\ph(D(x)\oplus D'(x))}\ge \eqref{1-le-norm(ph)-le-2}\ge \norm{D(x)\oplus D'(x)}=
\max\{\norm{D(x)},\norm{D'(x)}\}
$$
то есть $x\in U$ и $x\in U'$.
\epr

Из свойств $1^\circ$---$3^\circ$ следует

\btm\label{TH:A/U-set-epimorfizmov-v-InvSteAlg-diff}
Система $\pi_U:A\to A/U$ дифференциальных фактор-отображений образует сеть эпиморфизмов в категории $\InvSteAlg$ инволютивных стереотипных алгебр (в смысле \cite{Ak16}), то есть обладает следующими свойствами:
  \bit{
\item[(a)] у всякой алгебры $A$ есть хотя бы одна дифференциальная окрестность нуля $U$, и множество всех дифференциальных окрестностей нуля в $A$ направлено относительно предпорядка
$$
U\le U'\quad\Longleftrightarrow\quad U\supseteq U',
$$

\item[(b)] для всякой алгебры $A$ система морфизмов $\varkappa_U^{U'}$ из \eqref{varkappa^U'_U}
ковариантна, то есть для любых трех дифференциальных окрестностей нуля $U\supseteq U'\supseteq U''$ коммутативна диаграмма
$$
 \xymatrix @R=2pc @C=1.2pc
 {
  A/U &  & A/ U''\ar[ll]_{\varkappa^{U''}_U}\ar[dl]^{\varkappa^{U''}_{U'}}\\
  & A/U \ar[ul]^{\varkappa^{U'}_U} &
 }
$$
и эта система $\varkappa_U^{U'}$ обладает проективным пределом в $\InvSteAlg$;

\item[(c)] для всякого морфизма $\alpha:A\gets A'$ в $\InvSteAlg$ и любой дифференциальной окрестности нуля $U$ в $A$ найдется дифференциальная окрестность нуля $U'$ в $A'$ и морфизм $\alpha_U^{U'}:A/U\gets A'/U'$ такие, что коммутативна
диаграмма
 \beq\label{DIAGR:set-differ} \xymatrix @R=2.5pc @C=4.0pc {
 A\ar[d]_{\pi_U} & A'\ar@{-->}[d]^{\pi_{U'}}\ar[l]_{\alpha} \\
A/U & A'/U'\ar@{-->}[l]^{\alpha_U^{U'}}
 } \eeq
 }\eit
\etm

В соответствии с пунктом (b) этой теоремы, существует проективный предел $\projlim_{U'\in {\mathcal U}_{C^\infty}^A} A/ U'$ системы $\varkappa_U^{U'}$. Как следствие, существует единственная стрелка $\pi:A\to \projlim_{U'\in {\mathcal U}_{C^\infty}^A} A/ U'$ в  $\InvSteAlg$, замыкающая все диаграммы
\beq\label{A-to-leftlim-A/U-diff}
 \xymatrix @R=2pc @C=1.2pc
 {
  & A\ar[ld]_{\pi_U}\ar@{-->}[rd]^{\pi_{U'}} &   \\
  A/U &  & \projlim_{U'\in {\mathcal U}_{C^\infty}^A} A/ U'\ar[ll]^{\varkappa_U}
 }
\eeq
Образ $\pi(A)$ отображения $\pi$ является (инволютивной подалгеброй и) подпространством в стереотипном пространстве $\projlim_{U'\in {\mathcal U}_{C^\infty}^A} A/ U'$. Поэтому оно порождает некое непосредственное подпространство в $\projlim_{U'\in {\mathcal U}_{C^\infty}^A} A/ U'$, или оболочку $\Env\pi(A)$ \cite[(4.68)]{Akbarov-De-Gruyter-I}, то есть наибольшее стереотипное пространство содержащееся в  $\projlim_{U'\in {\mathcal U}_{C^\infty}^A} A/ U'$ и имеющее $\pi(A)$ плотным подпространством. 
Иными словами, $\Env\pi(A)$ --- узловой образ морфизма $\pi$:
$$
\Env\pi(A)=\Im_\infty \pi
$$
Обозначим через $\rho:A\to \Env\pi(A)$ поднятие морфизма $\pi$ в $\Env\pi(A)$:
\beq\label{rho:A-to-Env(pi(A))-diff}
 \xymatrix @R=2pc @C=1.2pc
 {
  A\ar[rr]^{\pi}\ar@{-->}[rd]_{\rho} & & \projlim_{U'\in {\mathcal U}_{C^\infty}^A} A/ U'   \\
   & \Env\pi(A)=\Im_\infty \pi\ar[ru]_{\iota} & 
 }
\eeq
(здесь $\iota$ --- естественное погружение подпространства $\Env\pi(A)$ в объемлющее пространство $\projlim_{U'\in {\mathcal U}_{C^\infty}^A} A/ U'$).

\btm\label{TH:opisanie-gladkoi-obolochki} Морфизм $\rho:A\to \Env\pi(A)$ является гладкой оболочкой алгебры $A$:
  \beq\label{Env_EA=Env-pi(A)-diff}
\Env_{\mathcal E} A=\Env\pi(A)=\Im_\infty \pi.
  \eeq
\etm

\bpr
В силу \eqref{ph=ph_U-circ-pi_U-*}, система дифференциальных фактор-отображений $\pi_U:A\to A/U$ порождает класс $\varPhi$ дифференциальных гомоморфизмов изнутри. С другой стороны, по теореме \ref{TH:SMono-circledcirc-DEpi=InvSteAlg}, класс $\DEpi$ всех плотных эпиморфизмов, в которых ищется оболочка, мономорфно дополняем в категории $\InvSteAlg$. Отсюда по теореме \ref{TH:funktorialnost-pri-seti-Epi-i-dolonyaemosti}, оболочки алгебры $A$ относительно классов $\varPhi$ и ${\mathcal N}$ совпадают между собой и представляют собой морфизм $\rho$.
\epr

\brem
Как и в случае с непрерывной оболочкой, гладкая оболочка $\rho:A\to\Env_{\mathcal E}A$ представляет собой композицию элементов $\red_\infty$ и $\coim_\infty$ узлового разложения морфизма $\pi:A\to\projlim_{U'\in {\mathcal U}_{C^\infty}^A} A/ U'$ в категории $\tt Ste$ стереотипных пространств (не алгебр!):
  \beq\label{C^infty-envelope=im_infty-lim N_X}
\env_{\mathcal E} A=\red_\infty\pi\circ\coim_\infty\pi.
  \eeq
Наглядно это изображается диаграммой
  \beq\label{DIAGR:E-envelope=im_infty-lim N_X}
\xymatrix @R=3.pc @C=4.0pc 
{
A\ar[d]_{\coim_\infty\pi}\ar@{-->}[drrr]^{\env_{\mathcal E} A}\ar[rrr]^{\pi=\projlim_{U'\in {\mathcal U}_{C^\infty}^A}\pi_{U'}} &&& \projlim_{U'\in {\mathcal U}_{C^\infty}^A} A/U' &   \\
\Coim_\infty\pi\ar[rrr]_{\red_\infty\pi} &&&  \Im_\infty \pi \ar[u]_{\im_\infty\pi}\ar@{=}[r] & \Env_{\mathcal E}A
}
 \eeq
\erem

\bcor\label{COR:0-ne-x-in-Env_E(A)}
Для всякой инволютивной стереотипной алгебры $A$ и любого ненулевого элемента $x\in \Env_{\mathcal E} A$ ее гладкой оболочки
$$
0\ne x\in \Env_{\mathcal E} A
$$
найдется  $C^*$-алгебра $B$ с присоединенными нильпотентными элементами и дифференциальный гомоморфизм $D:A\to B[m]$ такой что
$$
D(x)\ne 0.
$$
\ecor
\bpr
По теореме \ref{TH:opisanie-gladkoi-obolochki} оболочка $\Env_{\mathcal E} A$ представима как непосредстввенное подпространство в проективном пределе фактор-алгебр $A/U$ по дифференциальным окрестностям нуля:
$$
\Env_{\mathcal E} A\osubarr \projlim_{U\to 0} A/ U
$$
Отсюда следует, что для элемента $x\in \Env_{\mathcal E} A$, $x\ne 0$, должна найтись окрестность нуля $U$ такая, что фактор-отображение $\pi_U:A\to A/U$ будет ненулевым на $x$:
$$
\pi_U(x)\ne 0.
$$
По определению дифференциальной окрестности нуля, $U$ есть прообраз единичного шара при дифференциальном гомоморфизме $D:A\to B[m]$
в некоторую $C^*$-алгебру $B$ с присоединенными нильпотентными элементами. Это значит, что при гомоморфизме $D$ элемент $x$ тоже не обнуляется:
$$
D(x)\ne 0.
$$
\epr

\subsection{Гладкие алгебры}

Инволютивную стереотипную алгебру $A$ мы называем {\it гладкой алгеброй}, если она является полным объектом в категории $\InvSteAlg$ инволютивных стереотипных алгебр относительно гладкой оболочки $\Env_{\mathcal E}$ в смысле определения на с.\pageref{DEF:polnye-objekty}, то есть удовлетворяет следующим равносильным условиям\footnote{Условие (iv) в этом списке --- действие предложения \ref{PROP:polnota-pri-podtalkivanii}.}:
\bit{
\item[(i)] всякое гладкое расширение $\sigma:A\to A'$ алгебры $A$ является изоморфизмом в $\InvSteAlg$;

\item[(ii)] локальная единица $1_A:A\to A$ является гладкой оболочкой $A$;

\item[(iii)] гладкая оболочка алгебры $A$ является изоморфизмом в $\InvSteAlg$: $\env_{\mathcal E} A\in\Iso$,

\item[(iv)] алгебра $A$ изоморфна в $\InvSteAlg$ гладкой оболочке какой-нибудь алгебры $B$: 
$$
A\cong\Env_{\mathcal E}B.
$$
}\eit\noindent 
Класс всех гладких алгебр мы обозначаем ${\mathcal E}\text{-}{\tt Alg}$. Он образует полную подкатегорию в категории $\InvSteAlg$.

\noindent\rule{160mm}{0.1pt}\begin{multicols}{2}

\bex\label{EX:C(M)-nepr-alg}
Пусть $M$ --- гладкое многообразие. Алгебра ${\mathcal E}(M)$ гладких функций $u:M\to \C$ (с поточечными алгебраическими операциями) наделяется стандартной топологией, которую легче всего описать с помощью следующей системы определений.
\biter{
\item[$\bullet$] Говорят, что две функции $u,v\in {\mathcal E}(M)$ имеют одинаковый {\it росток} в точке $t\in M$, и обозначают это записью
    \beq\label{DEF:rostok}
    u\equiv v \ (\kern-10pt\mod t),
    \eeq
    если эти функции совпадают в некоторой окрестности $U$ точки $t$:
    \beq\label{DEF:rostok-1}
    u(s)=v(s),\quad s\in U.
    \eeq
    
\item[$\bullet$] Линейное отображение $D: {\mathcal E}(M)\to {\mathcal E}(M)$ называется {\it дифференциальным оператором}, если оно сохраняет ростки функций:
    \begin{multline}\label{DEF:diff-oper}
\forall t\in M\quad u,v\in {\mathcal E}(M) \\
    u\equiv v \ (\kern-10pt\mod t)\quad\Rightarrow\quad Du\equiv Dv \ (\kern-10pt\mod t)
    \end{multline}
     Это условие эквивалентно тому, что $D$ не увеличивает {\it носитель} функции
    \beq\label{DEF:diff-oper}
\forall u\in {\mathcal E}(M)\qquad  \supp Du\subseteq\supp u
    \eeq
    где носитель функции $u$ определяется как множество точек, в которых росток функции $u$ ненулевой:  
\beq\label{DEF:nositel}
\supp u=\{t\in M:\ u\nequiv 0 \ (\kern-10pt\mod t) \}
\eeq

\item[$\bullet$]  Считается, что направленность функций $u_i\in {\mathcal E}(M)$ сходится к функции $u\in {\mathcal E}(M)$ в пространстве ${\mathcal E}(M)$, если для всякого дифференциального оператора $D: {\mathcal E}(M)\to {\mathcal E}(M)$ направленность функций $Du_i$ сходится к функции $Du$ в пространстве ${\mathcal C}(M)$, то есть равномерно на каждом компакте $T\subseteq M$:
    $$
    \max_{t\in T}\abs{Du_i(t)-Du(t)}\underset{i\to\infty}{\longrightarrow}0.
    $$
}\eiter\noindent
Важное для нас наблюдение состоит в том, что {\it алгебра ${\mathcal E}(M)$ является гладкой}.
\eex
\bpr
Мы сошлемся здесь на доказательство следующего примера \ref{EX:E(M,B(X))-glad-alg}, в котором описывается более общий случай. 
\epr

\brem\label{EX:E(M,B(X))-glad-alg}
В примере \ref{EX:C(M,B(X))-nepr-alg} мы доказывали непрерывность алгебры ${\mathcal C}(M,{\mathcal B}(X))$ непрерывных функций на паракомпактном локально компактном пространстве $M$ со значениями а алгебре ${\mathcal B}(X)$ операторов на конечномерном гильбертовом пространстве $X$. В теории гладких оболочек аналогичное утверждение, видимо, неверно: алгебра ${\mathcal E}(M,{\mathcal B}(X))$ гладких функций на гладком многообразии $M$ со значениями а алгебре ${\mathcal B}(X)$ не будет гладкой. Я думаю, этот эффект есть результат выбора слишком жесткого определения гладкой оболочки. Если в определении на с.\pageref{DEF:C^infty-obolochka} отказаться от требования, чтобы тестовые морфизмы $\ph\in\varPhi$ были дифференциальными гомоморфизмами, а просто считать, что это должны быть гомоморфизмы во всевозможные $C^*$-алгебры $B[m]$ с присоединенными самосопряженными нильпотентными элементами, то эта асимметрия между теориями, я предполагаю, должна исчезнуть. Но, как всегда, когда дело касается изменений в фундаменте теории, масштабы переделок заставляют меня отложить эту работу на какое-то время.   
\erem

\end{multicols}\noindent\rule[10pt]{160mm}{0.1pt}

Из теоремы \ref{TH:proizvedenie-polnyh-obyektov} следует

\btm\label{TH:proizvedenie-gladkih-algebr}\phantom{.}
  \bit{
\item[(i)] Произведение $\prod_{i\in I}A_i$ любого семейства $\{A_i;\ i\in I\}$ гладких алгебр является гладкой алгеброй.

\item[(ii)] Проективный предел $\projlim_{i\in I}A_i$ любой ковариантной (контравариантной) системы $\{A_i,\iota^j_i;\ i\in I\}$ гладких алгебр является гладкой алгеброй.
 }\eit
 \etm

\paragraph{Гладкое тензорное произведение инволютивных стереотипных алгебр.}

Пусть $\Env_{\mathcal E}$ -- функтор гладкой оболочки, определенный в следствии \ref{COR:Env_E-idemp-funktor}.
Для любых двух инволютивных стереотипных алгебр $A$ и $B$ их {\it гладким тензорным произведением} условимся называть алгебру
\beq\label{E/circledast}
A\overset{\mathcal E}{\circledast} B=\Env_{\mathcal E}(A\circledast B)
\eeq

\btm\label{TH:C^infty-circledast->odot} Для любых двух инволютивных стереотипных алгебр $A$ и $B$ существует единственное линейное непрерывное отображение $\eta^\infty_{A,B}:A\overset{\mathcal E}{\circledast}B\to A\odot B$, замыкающее диаграмму
    \beq\label{DIAGR:eta^infty_(A,B)}
\xymatrix @R=2.pc @C=5.0pc 
{
A\circledast B\ar[dr]_{\env_{\mathcal E}{A\circledast B}\quad}\ar[rr]^{@_{A,B}} & & A\odot B \\
& A\overset{\mathcal E}{\circledast}B\ar[ur]_{\eta^\infty_{A,B}} & \\
}
\eeq
а система отображений  $\eta^\infty_{A,B}:A\overset{\mathcal E}{\circledast}B\to A\odot B$ является естественным преобразованием функтора $(A,B)\in{\mathcal E}\text{-}{\tt Alg}^2\mapsto A\overset{\mathcal E}{\circledast}B\in{\mathcal E}\text{-}{\tt Alg}$ в функтор $(A,B)\in{\mathcal E}\text{-}{\tt Alg}^2\mapsto A\odot B\in {\tt Ste}$.
\etm
\bpr Рассмотрим диаграмму
$$
\xymatrix @R=3.pc @C=8.0pc 
{
A\circledast B\ar[d]_{\env_{\mathcal E}{A\circledast B}}\ar[dr]_{\env_{\mathcal C}{A\circledast B}\quad}\ar[r]^{@_{A,B}} &  A\odot B \\
A\overset{\mathcal E}{\circledast}B\ar[r]_{\zeta_{A\circledast B}} & A\overset{\mathcal C}{\circledast}B\ar[u]_{\eta_{A,B}}  \\
}
$$
В ней левый нижний треугольник -- диаграмма \eqref{Env_infty->Env} для алгебры $A\circledast B$, а правый верхний треугольник -- диаграмма \eqref{DIAGR:eta_(A,B)}. Положив
$$
\eta^\infty_{A,B}=\eta_{A,B}\circ\alpha_{A,B}
$$
мы получим нужный морфизм.
 \epr

\paragraph{Гладкое тензорное произведение гладких алгебр.}

Из теорем \ref{TH:Env_C^infty-reg-obolochka} и \ref{TH:sushestvovanie-tenz-proizv-v-L} следует

\btm\label{TH:E-obolochka=monoidalnyi-funktor} Формула \eqref{C/circledast} определяет в ${\mathcal E}\text{-}{\tt Alg}$ тензорное произведение, превращающее ${\mathcal E}\text{-}{\tt Alg}$ в моноидальную категорию, а функтор гладкой оболочки $\Env_{\mathcal E}$ является моноидальным функтором из моноидальной категории $(\InvSteAlg,\circledast)$ инволютивных стереотипных алгебр в моноидальную категорию $({\mathcal E}\text{-}{\tt Alg},\overset{\mathcal E}{\circledast})$ гладких алгебр. Соответствующий морфизм бифункторов
 $$
\Big((A,B)\mapsto \Env_{\mathcal E}(A)\overset{\mathcal E}{\circledast} \Env_{\mathcal E}(B)\Big)\overset{E^{\circledast}}{\rightarrowtail} \Big((A,B)\mapsto \Env_{\mathcal E}(A\circledast B)\Big)
 $$
определяется формулой
$$
E^\circledast_{A,B}=\Env_{\mathcal E}(\env_{\mathcal E}A\circledast \env_{\mathcal E}B)^{-1}:\Env_{\mathcal E}(A)\overset{\mathcal E}{\circledast} \Env_{\mathcal E}(B)=\Env_{\mathcal E}(\Env_{\mathcal E}(A)\circledast \Env_{\mathcal E}(B))\to \Env_{\mathcal E}(A\circledast B),
$$
а морфизмом $E^{\C}$ в ${\mathcal C}\text{-}{\tt Alg}$, переводящий единичный объект $\C$ категории ${\mathcal C}\text{-}{\tt Alg}$ в образ $\Env_{\mathcal E}(\C)$ единичного объекта $\C$ категории $\InvSteAlg$, будет локальная единица:
$$
E^{\C}=1_{\C}:\C\to\C=\Env_{\mathcal E}(\C).
$$
\etm

Каждой паре элементов $a\in A$, $b\in B$ можно поставить в соответствие элементарный тензор
\beq\label{DEF:a-overset-C^infty-circledast-b}
a\overset{\mathcal E}{\circledast}b=\env_{\mathcal E}(a\circledast b)
\eeq

\blm\label{LM:polnota-a-overset-C^infty-circledast-b}
Элементарные тензоры $a\overset{\mathcal E}{\circledast}b$, $a\in A$, $b\in B$, полны в $A\overset{\mathcal E}{\circledast} B$ и при отображении $\eta_{A,B}$ переходят в элементарные тензоры $a\odot b$:
\beq\label{eta(a-overset-C^infty-circledast-b)=a-odot-b}
\eta_{A,B}(a\overset{\mathcal E}{\circledast}b)=a\odot b.
\eeq
\elm
\bpr
Тензоры $a\circledast b$ полны в $A\circledast B$, а образ $\env_{\mathcal E}$ плотен в $A\overset{\mathcal E}{\circledast} B$. Тождество \eqref{eta(a-overset-C^infty-circledast-b)=a-odot-b} следует из диаграммы \eqref{DIAGR:eta^infty_(A,B)}.
\epr

\paragraph{Действие гладкой оболочки на биалгебры.}

Следующие три утверждения аналогичны теоремам \ref{LM:koalg-v-CAlg->koalg-v-odot}, \ref{TH:C-obolochka-sohranyaet-Hopfov} и \ref{TH:C-obolochka-sohranyaet-inv-Hopfov}, и доказываются так же.

\blm\label{LM:koalg-v-C^inftyAlg->koalg-v-odot}
Если $A$ -- коалгебра в моноидальной категории $({\mathcal E}\text{-}{\tt Alg},\overset{\mathcal E}{\circledast})$  гладких алгебр со структурными морфизмами
$$
\varkappa:A\to A\overset{\mathcal E}{\circledast} A,\qquad \e:A\to\C,
$$
то $A$ является коалгеброй в моноидальной категории $({\tt Ste},\odot)$ стереотипных пространств со структурными морфизмами
$$
\lambda=\eta_{A,A}\circ\varkappa:A\to A\odot A,\qquad \e:A\to\C.
$$
\elm

\btm\label{TH:C^infty-obolochka-sohranyaet-Hopfov}
Пусть $H$ -- биалгебра в категории $({\tt Ste},\circledast)$ стереотипных пространств, или, что эквивалентно, коалгебра в категории ${\tt Ste}^{\circledast}$ стереотипных алгебр с коумножением $\varkappa$ и коединицей $\e$. Тогда
 \bit{
\item[(i)] гладкая оболочка $\Env_{\mathcal E}H$ является коалгеброй в моноидальной категории $({\mathcal C}\text{-}{\tt Alg},\overset{\mathcal E}{\circledast})$ гладких алгебр с коумножением и коединицей
    \beq\label{varkappa^E^infty,e^E^infty}
    \varkappa_{\Env_{\mathcal E}}=\Env_{\mathcal E}(\env_{\mathcal E}H\circledast \env_{\mathcal E}H)\circ \Env_{\mathcal E}(\varkappa),\qquad \e_{\Env_{\mathcal E}}=\Env_{\mathcal E}(\e),
    \eeq

\item[(ii)] гладкая оболочка $\Env_{\mathcal E}H$ является коалгеброй в моноидальной категории $({\tt Ste},\odot)$ стереотипных пространств с коумножением и коединицей
    \beq\label{varkappa^odot,e^odot-infty}
    \varkappa_\odot=\eta_{\Env_{\mathcal E}H,\Env_{\mathcal E}H}\circ \Env_{\mathcal E}(\env_{\mathcal E}H\circledast \env_{\mathcal E}H)\circ \Env_{\mathcal E}(\varkappa)=
    \eta_{\Env_{\mathcal E}H,\Env_{\mathcal E}H}\circ \varkappa_{\Env_{\mathcal E}},\qquad \e_\odot=\Env_{\mathcal E}(\e),
    \eeq

\item[(iii)] морфизм $(\env_{\mathcal E}H)^\star:H^\star\gets \Env_{\mathcal E}H^\star$, сопряженный к морфизму оболочки $\env_{\mathcal E}H:H\to \Env_{\mathcal E}H$, является морфизмом стереотипных алгебр, если $\Env_{\mathcal E}H^\star$ рассматривается как алгебра с умножением и единицей, сопряженными к \eqref{varkappa^odot,e^odot-infty}, а $H^\star$ -- как алгебра с умножением и единицей
    $$
    \varkappa^\star\circ @_{H^\star,H^\star},\qquad \e^\star.
    $$

 }\eit
\etm

\btm\label{TH:C^infty-obolochka-sohranyaet-inv-Hopfov}
Пусть $H$ -- инволютивная алгебра Хопфа в категории $({\tt Ste},\circledast)$ стереотипных пространств. Тогда
 \bit{
\item[(i)] гладкая оболочка $\Env_{\mathcal E}H$, как коалгебра в моноидальных категориях $({\mathcal E}\text{-}{\tt Alg},\overset{\mathcal E}{\circledast})$ и $({\tt Ste},\odot)$, обладает согласованными между собой антиподом $\Env_{\mathcal E}(\sigma)$ и инволюцией $\Env_{\mathcal E}(\bullet)$, однозначно определяемыми диаграммами в категории $\tt{Ste}$
    \beq\label{C^infty-obolochka-sohranyaet-inv-Hopfov}
    \xymatrix @R=2.pc @C=4.0pc 
{
H\ar[d]_{\sigma}\ar[r]^{\env_{\mathcal E}H} & \Env_{\mathcal E}H\ar@{-->}[d]^{\Env_{\mathcal E}(\sigma)} \\
H\ar[r]^{\env_{\mathcal E}H} & \Env_{\mathcal E}H
}
\qquad
    \xymatrix @R=2.pc @C=4.0pc 
{
H\ar[d]_{\bullet}\ar[r]^{\env_{\mathcal E}H} & \Env_{\mathcal E}H\ar@{-->}[d]^{\Env_{\mathcal E}(\bullet)} \\
H\ar[r]^{\env_{\mathcal E}H} & \Env_{\mathcal E}H
}
\eeq

\item[(ii)] морфизм $(\env_{\mathcal E}H)^\star:H^\star\gets \Env_{\mathcal E}H^\star$, сопряженный к морфизму оболочки $\env_{\mathcal E}H:H\to \Env_{\mathcal E}H$, является инволютивным гомоморфизмом стереотипных алгебр над $\circledast$, если $H^\star$ и $\Env_{\mathcal E}H^\star$ наделяются структурой сопряженных инволютивных алгебр к инволютивным коалгебрам с антиподом $H$ и $\Env_{\mathcal E}H$ по \cite[p.457, $4^\circ$]{Akbarov-De-Gruyter-I}.
 }\eit
\etm

\paragraph{Гладкое тензорное произведение с ${\mathcal E}(M)$.}

Пусть $X$ -- стереотипное пространство, и $M$ -- гладкое (локально евклидово) многообразие. Рассмотрим алгебру ${\mathcal E}(M)$ гладких функций на $M$ и пространство ${\mathcal E}(M,X)$ гладких отображений на $M$ со значениями в $X$. Мы наделяем ${\mathcal E}(M)$ и ${\mathcal E}(M,X)$ стандартной топологией равномерной сходимости на компактах по каждой частной производной
$$
u_i\overset{{\mathcal E}(M,X)}{\longrightarrow} 0\quad\Longleftrightarrow\quad \forall U\subseteq M\quad \forall k\in\N^d \qquad u_i^{(k)}|_U\overset{{\mathcal C}(U,X)}{\longrightarrow} 0,\qquad u\in {\mathcal E}(M,X),\quad t\in M.
$$
(где $u^{(k)}$ -- частная производная вдоль локальной карты на открытом множестве $U\subseteq M$) и поточечными алгебраическими операциями:
$$
(\lambda\cdot u)(t)=\lambda\cdot u(t)\qquad (u+v)(t)=u(t)+v(t),\qquad u,v\in {\mathcal C}(M,X),\quad \lambda\in\C,\quad t\in M.
$$

Из \cite[Theorem 8.9]{Ak03} следует

\bprop Справедливо тождество
\beq\label{E(M,X)-cong-E(M)-odot-X}
{\mathcal E}(M,X)\cong {\mathcal E}(M)\odot X
\eeq
\eprop

В дальнейшем нас будет интересовать случай, когда $A$ -- гладкая (и поэтому стереотипная) алгебра. Пространство ${\mathcal E}(M,A)$ при этом мы также наделяем структурой стереотипной алгебры с поточечными операциями
$$
(u\cdot v)(t)=u(t)\cdot v(t),\qquad u,v\in {\mathcal C}(M,A),\quad t\in M.
$$
Из \eqref{E(M,X)-cong-E(M)-odot-X} следует, что ${\mathcal E}(M,A)$ является стереотипным $A$-модулем.

\btm\label{TH:C^infty(M)-circledast-A->C(M,A)} Для всякой гладкой алгебры $A$ и любого гладкого многообразия $M$ естественное отображение
\beq\label{C^infty(M)-circledast-A->C(M,A)}
\iota:{\mathcal E}(M)\circledast A\to {\mathcal E}(M,A) \quad\Big|\quad \iota(u\circledast a)(t)=u(t)\cdot a,\quad u\in {\mathcal E}(M),\ a\in A,\ t\in M,
\eeq
является гладкой оболочкой и порождает изоморфизм стереотипных алгебр:
\beq\label{E(M)-circledast-A=E(M,A)}
{\mathcal E}(M)\overset{\mathcal E}{\circledast} A\cong{\mathcal E}(M,A).
\eeq
\etm

Мы разобьем доказательство на 5 лемм.

\blm\label{LM:iota-in-DEpi-E(M)}
Отображение $\iota:{\mathcal E}(M)\circledast A\to {\mathcal E}(M,A)$ является плотным эпиморфизмом.
\elm
\bpr
Это доказывается аналогично лемме \ref{LM:iota-in-DEpi-C(M)}.
\epr

\blm\label{LM:J^n_E(M)E(M)-circledast-A-cong-J^n_E(M)E(M,A)}
Модули ${\mathcal E}(M)\circledast A$ и ${\mathcal E}(M,A)$ над алгеброй ${\mathcal E}(M)$ имеют изоморфные расслоения струй:
\beq\label{J^n_E(M)E(M)-circledast-A-cong-J^n_E(M)E(M,A)}
\Jet^n_{{\mathcal E}(M)}{\mathcal E}(M)\circledast A\cong \Jet^n_{{\mathcal E}(M)}{\mathcal E}(M,A),\qquad n\in\N
\eeq
\elm
\bpr
Для каждой точки $t\in M$ идеал $I_t^{n+1}$ имеет конечную коразмерность в ${\mathcal E}(M)$, поэтому можно воспользоваться леммой \ref{PROP:[(X-circledast-Z)/(Y-circledast-Z)]^triangledown-cong-[(X-odot-Z)/(Y-odot-Z)]^triangledown}:
\begin{multline*}
\Jet^n_{{\mathcal E}(M)}{\mathcal E}(M)\circledast A=[({\mathcal E}(M)\circledast A)/ (I_t^{n+1}\circledast A)]^\vartriangle=\eqref{[(X-circledast-Z)/(Y-circledast-Z)]^triangledown-cong-[(X-odot-Z)/(Y-odot-Z)]^triangledown}=\\=
[({\mathcal E}(M)\odot A)/ (I_t^{n+1}\odot A)]^\vartriangle=\Jet^n_{{\mathcal E}(M)}{\mathcal E}(M)\odot A=\eqref{E(M,X)-cong-E(M)-odot-X}=
\Jet^n_{{\mathcal E}(M)}{\mathcal E}(M,A)
\end{multline*}
\epr

\blm\label{TH:diff-oper<-morfizm-rassl-struuj} Пусть $M$ -- гладкое многообразие и $\ph:{\mathcal E}(M)\to F$ --
гомоморфизм инволютивных стереотипных алгебр, причем $F$ --
$C^*$-алгебра, и $\ph({\mathcal E}(M))$ лежит в центре $F$:
$$
\ph({\mathcal E}(M))\subseteq Z(F).
$$
Тогда для любого стереотипного пространства $X$ всякий морфизм расслоений струй $\nu:\Jet_{{\mathcal E}(M)}^n({\mathcal E}(M,X))\to \Jet_{{\mathcal E}(M)}^0(F)$ определяет единственный дифференциальный оператор порядка $n$ между стереотипными ${\mathcal E}(M)$-модулями $D:{\mathcal E}(M)\to F$, удовлетворяющий
тождеству
\beq\label{j^0(Dx)=mu-circ-j^n(x)}
 \xymatrix @R=2pc @C=1.2pc
 {
 \Jet_{{\mathcal E}(M)}^0[{\mathcal E}(M,X)]\ar[dd]_{\nu} & \\
 & M \ar[ul]_{\jet^n(u)}\ar[dl]^{\jet^n(Du)} \\
 \Jet_{{\mathcal E}(M)}^0[F] &
 }
\qquad \jet^0(Du)=\nu\circ \jet^n(u),\qquad
u\in  {\mathcal E}(M,X).
 \eeq
 то есть такой, что $\nu$ является морфизмом раслоений струй, определяемым дифференциальным оператором $D$ по теореме \ref{TH:diff-oper->rassl-struuj}:
 $$
 \nu=\jet_n[D].
 $$
 \elm
\bpr По теореме \ref{TH:B-cong-Sec(val_AB)}, отображение
$v:F\to\Sec(\pi^0_{{{\mathcal E}(M)},F})$, переводящее $F$ в алгебру непрерывных сечений расслоения
значений $\pi^0_{{{\mathcal E}(M)},F}: \Jet^0_{{\mathcal E}(M)}F\to\Spec({{\mathcal E}(M)})$ над алгеброй ${{\mathcal E}(M)}$, является
изоморфизмом $C^*$-алгебр:
$$
F\cong\Sec(\pi^0_{{{\mathcal E}(M)},F}).
$$
Рассмотрим обратный изоморфизм $v^{-1}:\Sec(\pi^0_{{{\mathcal E}(M)},F})\to F$:
\beq\label{lambda(j^0(b))=b} v^{-1}(\jet^0(b))=b,\qquad b\in F.
\eeq
Тогда всякому морфизму расслоений струй $\nu:\Jet_{{\mathcal E}(M)}^n({{\mathcal E}(M,X)})\to \Jet_{{\mathcal E}(M)}^0(F)$ можно поставить в соответствие оператор $D:{{\mathcal E}(M,X)}\to F$ по формуле
\beq\label{Da=lambda(mu-circ-j^n(a))} Du=v^{-1}\Big(\nu\circ \jet^n(u)\Big),\qquad
u\in {{\mathcal E}(M,X)}.
\eeq
Он, очевидно, будет удовлетворять тождеству
\eqref{j^0(Dx)=mu-circ-j^n(x)}. С другой стороны, в пространстве гладких функций на $M$ со значениями в произвольном стереотипном пространстве $X$, как и в обычном пространстве функций со значениями в $\C$, справедлива формула Ньютона-Лейбница, и поэтому, при заданной локальной карте, справедливы лемма Адамара \cite{Petrovsky} и разложение Тейлора с остаточным членом в модуле $\overline{I_t^{n+1}\cdot X}$ с подходящим значением $n\in\N$. Поэтому, поскольку действие $D$ на элемент $x$ пропускается через струю $\jet^n(x)$, оно линейно выражается через коэффициенты Тейлора разложения элемента $x$ в окрестности данной точки $t\in M$ (при выборе локальной карты). Эти коэффициенты Тейлора являются дифференциальными операторами над ${\mathcal E}(M)$, и, как следствие, $D$ тоже должен быть дифференциальным оператором над ${\mathcal E}(M)$.
\epr

\blm Отображение $\iota:{\mathcal E}(M)\circledast A\to {\mathcal E}(M,A)$ является гладким расширением.
\elm
\bpr Пусть $D:{\mathcal E}(M)\circledast A\to B[m]$ -- морфизм в $C^*$-алгебру $B$ с присоединенными самосопряженными нильпотентами. Представим $D$ как семейство частных производных $D_k:{\mathcal E}(M)\circledast A\to B$, и положим
$$
\eta_k(u)=D_k(u\circledast 1),\qquad \alpha_k(a)=D_k(1\circledast a),\qquad u\in {\mathcal E}(M), \ a\in A.
$$
Тогда $\eta:{\mathcal E}(M)\to B[m]$, $\alpha:A\to B[m]$ будут морфизмами инволютивных стереотипных алгебр, и по лемме \ref{LM:ph:A-circledast-B->C},
\beq\label{D(u-circledast-a)=eta(u)-cdot-alpha(a)}
D(u\circledast a)=\eta(u)\cdot\alpha(a)=\alpha(a)\cdot\eta(u),\qquad u\in {\mathcal E}(M), \ a\in A,
\eeq
В частности,
\beq\label{D_0(u-circledast-a)=eta(u)-cdot-alpha(a)}
D_0(u\circledast a)=\eta_0(u)\cdot\alpha_0(a)=\alpha_0(a)\cdot\eta_0(u),\qquad u\in {\mathcal E}(M), \ a\in A.
\eeq

Рассмотрим оператор $\eta_0$ и обозначим буквой $C$ его образ в $B$:
$$
C=\overline{\eta_0({\mathcal E}(M))}.
$$
Пусть $F$ -- коммутант алгебры $C$ в $B$:
$$
F=C^!=\{x\in B:\quad \forall c\in C\quad x\cdot c=c\cdot x\}.
$$
Поскольку алгебра $C$ коммутативна, она также лежит в $F$, и более того, в центре $F$:
$$
C\subseteq Z(F).
$$
Заметим еще, что образы всех операторов $D_k$ лежат в $F$:
\beq\label{Im-D_k-subseteq-F}
D_k({\mathcal E}(M)\circledast A)\subseteq F.
\eeq
Для $k=0$ это можно доказать напрямую:
\begin{multline*}
D_0(v\circledast a)\cdot \eta_0(u)=D_0(v\circledast a)\cdot D_0(u\circledast 1)=D_0\big((v\circledast a)\cdot (u\circledast 1)\big)=D_0\big((v\cdot u)\circledast 1\big)=D_0\big((u\cdot v)\circledast 1\big)=\\=D_0\big((u\circledast 1)\cdot (v\circledast a)\big)=D_0(u\circledast 1)\cdot D_0(v\circledast a)=\eta_0(u)\cdot D_0(v\circledast a)
\end{multline*}
А для $k>0$ нужно применить соотношение \eqref{D_k(A)-subseteq-Z^1(D_0)}: поскольку при $k>0$ значения операторов $D_k$ и $D_0$ коммутируют, мы получаем
$$
D_k(v\circledast a)\cdot\eta_0(u)=D_k(v\circledast a)\cdot D_0(u\circledast 1)=D_0(u\circledast 1)\cdot D_k(v\circledast a)=\eta_0(u)\cdot D_k(v\circledast a).
$$

Чтобы убедиться, что $\iota$ -- гладкое расширение, нам надо показать, что существует система дифференциальных частных производных $\{D_k';\ k\in\N^d\}$ (над ${\mathcal E}(M,A)$!), продолжающих операторы $D_k$ с ${\mathcal E}(M)\circledast A$ на ${\mathcal E}(M,A)$ и принимающих значения в $F$:
 \beq\label{prodolzhenie-D_k-na-E(M,A)}
 \xymatrix @R=2pc @C=1.2pc
 {
 {\mathcal E}(M)\circledast A\ar[rr]^{\iota}\ar[dr]_{D_k} & & {\mathcal E}(M,A)\ar@{-->}[dl]^{D_k'}\\
  & F &
 }
 \eeq

Всякому дифференциальному оператору $D_k:{\mathcal E}(M)\circledast A\to F$ по теореме
\ref{TH:diff-oper->rassl-struuj} соответствует некий морфизм расслоений струй
$\jet_n[D_k]:\Jet_{{\mathcal E}(M)}^n[{\mathcal E}(M)\circledast A]\to \Jet_{{\mathcal E}(M)}^0(F)=\pi^0_A F$, где $n=|k|$, удовлетворяющий
тождеству
$$
\jet^0(D_k x)=\jet_n[D_k]\circ \jet^n(x),\qquad x\in{\mathcal E}(M)\circledast A.
$$
По лемме \ref{LM:J^n_E(M)E(M)-circledast-A-cong-J^n_E(M)E(M,A)}, расслоения струй алгебр ${\mathcal E}(M)\circledast A$ и ${\mathcal E}(M,A)$ изоморфны. Обозначим этот изоморфизм $\mu:\Jet_{{\mathcal E}(M)}^n[{\mathcal E}(M)\circledast A]\gets \Jet^n_{{\mathcal E}(M)}[{\mathcal E}(M,A)]$. Рассмотрим композицию $\nu=\jet_n[D_k]\circ\mu:\Jet_{{\mathcal E}(M)}[{\mathcal E}(M,A)]\to \Jet_{{\mathcal E}(M)}^0(F)=\pi^0_A F$:
$$
 \xymatrix @R=2pc @C=1.2pc
 {
 \Jet_{{\mathcal E}(M)}^n[{\mathcal E}(M)\circledast A]\ar[dr]_{\jet_n[D_k]} & & \Jet^n_{{\mathcal E}(M)}[{\mathcal E}(M,A)]\ar@{-->}[dl]^{\quad\nu=\jet_n[D_k]\circ\mu}\ar[ll]_{\mu}\\
  & \Jet^0_{{\mathcal E}(M)}[F] &
 }
$$
По лемме \ref{TH:diff-oper<-morfizm-rassl-struuj} этой пунктирной стрелке $\nu$ соответствует дифференциальный оператор $D_k':{\mathcal E}(M,A)\to F$ (над алгеброй ${\mathcal E}(M)$), удовлетворяющий тождеству
$$
\jet^0(D_k'f)=\jet_n[D_k]\circ \jet^n(f),\qquad  f\in {\mathcal E}(M,A).
$$
Для всякого $x\in {\mathcal E}(M)\circledast A$ мы получим
\beq\label{PROOF:TH:E(M,A)-1}
\jet^0\big(D_k'\iota(x)\big)=\jet_n[D_k]\circ \jet^n(\iota(x))=\jet_n[D_k]\circ
\jet^n(x)=\jet^0(D_kx).
\eeq
Заметим далее, что по теореме
\ref{TH:B-cong-Sec(val_AB)}, отображение
$\jet^0=v:F\to\Sec(\pi^0_{{\mathcal E}(M)}F)=\Sec(\Jet^0_{{\mathcal E}(M)}F)$, переводящее $F$ в алгебру непрерывных сечений
расслоения значений $\pi^0_{{\mathcal E}(M)}F: \Jet^0_{{\mathcal E}(M)}F\to\Spec({{\mathcal E}(M)})$ над алгеброй ${{\mathcal E}(M)}$, является
изоморфизмом $C^*$-алгебр:
$$
F\cong\Sec(\pi^0_{{\mathcal E}(M)}F).
$$
Поэтому к \eqref{PROOF:TH:E(M,A)-1} можно применить оператор, обратный $\jet^0$, и мы получим равенство
$$
D_k'\iota(x)=D_kx.
$$
То есть $D_k'$ продолжает $D_k$ в диаграмме
\eqref{prodolzhenie-D_k-na-E(M,A)}. Кроме того, из того, что $\iota$ отображает
${\mathcal E}(M)\circledast A$ плотно в ${\mathcal E}(M,A)$ следует, что условия
\eqref{DEF:sist-chastn-proizv-*}-\eqref{DEF:sist-chastn-proizv-2} переносятся с
оператора $D_k$ на оператор $D_k'$.

По построению, операторы $D_k'$ будут дифференциальными относительно алгебры ${\mathcal E}(M)$, однако нам этого недостаточно: нужно чтобы каждый $D_k'$ был дифференциальным оператором порядка $|k|$ относительно алгебры ${\mathcal E}(M,A)$, и чтобы операторы $D_k'$ образовывали систему частных производных на ${\mathcal E}(M,A)$.

И то и другое следует из того, что операторы $D_k$ образуют дифференциальную систему частных производных на ${\mathcal E}(M)\circledast A$. Во-первых, каждый $D_k$ является дифференциальным оператором порядка $|k|$ относительно ${\mathcal E}(M)\circledast A$, поэтому для любых $u_0,u_1,...,u_{|k|}\in {\mathcal E}(M)$ выполняется равенство
$$
[...[[D_k,u_0\circledast a_0],u_1\circledast a_1],... u_{|k|}\circledast a_{|k|}]=0.
$$
Из него следует
$$
[...[[D_k',\iota(u_0\circledast a_0)],\iota(u_1\circledast a_1)],... \iota(u_{|k|}\circledast a_{|k|})]=0,
$$
и, поскольку элементы вида $\iota(u\circledast a)$ полны в ${\mathcal E}(M,A)$, отсюда следует, что их можно заменить произвольными векторами из ${\mathcal E}(M,A)$, и мы получаем, что $D_k'$ -- дифференциальный оператор порядка $|k|$ над ${\mathcal E}(M,A)$.
Во-вторых, формулы \eqref{DEF:sist-chastn-proizv-*}-\eqref{DEF:sist-chastn-proizv-2} точно так же переносятся с $D_k$ на $D_k'$. Например, из \eqref{DEF:sist-chastn-proizv-2} для операторов $D_k$,
$$
D_k(x\cdot y)=\sum_{0\le l\le k}\begin{pmatrix}k \\ l\end{pmatrix}\cdot D_{k-l}(x)\cdot D_l(y),\qquad x,y\in {\mathcal E}(M)\circledast A,
$$
следует
$$
D_k'(\iota(x)\cdot\iota(y))=D_k'(\iota(x\cdot y))=D_k(x\cdot y)=\sum_{0\le l\le k}\begin{pmatrix}k \\ l\end{pmatrix}\cdot D_{k-l}(x)\cdot D_l(y)=\sum_{0\le l\le k}\begin{pmatrix}k \\ l\end{pmatrix}\cdot D_{k-l}(\iota(x))\cdot D_l(\iota(y)),
$$
Это верно для любых $x,y\in {\mathcal E}(M)\circledast A$. Поскольку образ $\iota$ плотен в ${\mathcal E}(M,A)$ (лемма \ref{LM:iota-in-DEpi-E(M)}), мы получаем, что $\iota(x)$ и $\iota(y)$ можно заменить на произвольные векторы из ${\mathcal E}(M,A)$, то есть \eqref{DEF:sist-chastn-proizv-2} справедливо и для операторов $D_k'$.
\epr

\blm Отображение $\iota:{\mathcal E}(M)\circledast A\to {\mathcal E}(M,A)$ является гладкой оболочкой.
\elm
\bpr
Пусть $\sigma:{\mathcal E}(M)\circledast A\to C$ -- какое-то другое гладкое расширение. Нам нужно убедиться, что существует морфизм $\upsilon$, замыкающий диаграмму
$$
 \xymatrix @R=2pc @C=1.2pc
 {
{\mathcal E}(M)\circledast A\ar[rr]^{\sigma}\ar[dr]_{\iota} & & C \ar@{-->}[dl]^{\upsilon}\\
  & {\mathcal E}(M,A) &
 }
$$

Зафиксируем локальную карту $\ph:U\to V$, где $U\subseteq M$, $V\subseteq\R^d$, и пусть и $K\subseteq U$ -- компакт, совпадающий с замыканием своей внутренности:
$\overline{\Int(K)}=K$. Операторы из примера \ref{EX:sist-chast-proizv-na-C-infty(M)}
$$
\varPhi_k(u)=\frac{\partial^{\abs{k}}(u\circ\ph^{-1})}{\partial t_1^{k_1}...\partial
t_d^{k_d}}\circ\ph,\qquad u\in {\mathcal E}(M),\quad k\in\N^d,
$$
образуют систему частных производных из ${\mathcal E}(M)$ со значениями в ${\mathcal C}(K)$.
Пусть $\varPhi:{\mathcal E}(M)\to \mathcal{C}(K)[m]$ -- дифференциальный гомоморфизм, соответствующий этой системе $\{\varPhi_k\}$.

Зафиксируем далее произвольный гомоморфизм $\eta:A\to B[n]$ в какую-нибудь $C^*$-алгебру с присоединенными нильпотентными элементами $B[n]$ и положим
$$
D(u\circledast a)=\varPhi(u)\circledast \eta(a),\qquad u\in {\mathcal E}(M),\quad a\in A.
$$
Отображение $D$ будет гомоморфизмом из ${\mathcal E}(M)\circledast A$ в алгебру ${\mathcal C}(K)[m]\circledast B[n]$, которая в силу \eqref{A[m]-circledast-B}, изоморфна $({\mathcal C}(K)\circledast B)[m\oplus n]$. Стереотипное тензорное произведение ${\mathcal C}(K)\circledast B$ естественно отображается в максимальное тензорное произведение $C^*$-алгебр ${\mathcal C}(K)\underset{\max}{\otimes}B$, которое в свою очередь, изоморфно ${\mathcal C}(K)\odot B$ и ${\mathcal C}(K,B)$:
$$
{\mathcal C}(K)\circledast B\to {\mathcal C}(K)\underset{\max}{\otimes}B\cong \cite[(5.25)]{Akbarov-De-Gruyter-I}
\cong
{\mathcal C}(K)\odot B\cong \eqref{C(M,X)-cong-C(M)-odot-X} \cong {\mathcal C}(K,B).
$$

Поэтому мы можем считать $D$ морфизмом в $C^*$-алгебру с присоединенными нильпотентными элементами ${\mathcal C}(K,B)[m\oplus n]$,
\begin{multline*}
D:{\mathcal E}(M)\circledast A\to {\mathcal C}(K)[m]\circledast B[n]\cong\eqref{A[m]-circledast-B},\eqref{A[m][n]=A[m-oplus-n]}\cong ({\mathcal C}(K)\circledast B)[m\oplus n]\to   ({\mathcal C}(K)\underset{\max}{\otimes}B)[m\oplus n]\cong\\ \cong\cite[(5.25)]{Akbarov-De-Gruyter-I}\cong({\mathcal C}(K)\odot B)[m\oplus n]\cong\eqref{C(M,X)-cong-C(M)-odot-X}\cong {\mathcal C}(K,B)[m\oplus n]
\end{multline*}
Поскольку $\sigma:{\mathcal E}(M)\circledast A\to C$ есть гладкое расширение, гомоморфизм $D:{\mathcal E}(M)\circledast A\to {\mathcal C}(K,B)[m\oplus n]$ должен однозначно продолжаться до некоторого
гомоморфизма $D':C\to {\mathcal C}(K,B)[m\oplus n]$:
 \beq\label{DIAGR:env-E(M)-circledast-A-0}
 \xymatrix @R=2pc @C=1.2pc
 {
{\mathcal E}(M)\circledast A\ar[rr]^{\sigma}\ar[dr]_{D} & & C \ar@{-->}[dl]^{D'}\\
  & {\mathcal C}(K,B)[m\oplus n] &
 }
 \eeq
Заметим, что ${\mathcal C}(K,B)[m\oplus n]$ изоморфно ${\mathcal C}\big(K, B[n]\big)[m]$,
$$
{\mathcal C}(K,B)[m\oplus n]\cong\eqref{A[m][n]=A[m-oplus-n]}\cong {\mathcal C}(K,B)[n][m]\cong\eqref{C(M,B)[n]=C(M,B[n])}\cong {\mathcal C}\big(K, B[n]\big)[m],
$$
поэтому диаграмму \eqref{DIAGR:env-E(M)-circledast-A-0} можно поправить так:
 \beq\label{DIAGR:env-E(M)-circledast-A-1}
 \xymatrix @R=2pc @C=1.2pc
 {
{\mathcal E}(M)\circledast A\ar[rr]^{\sigma}\ar[dr]_{D} & & C \ar@{-->}[dl]^{D'}\\
  & {\mathcal C}(K,B[n])[m] &
 }
 \eeq
Теперь мы можем вернуться к системе частных производных $D_k$, и для всякого
индекса $k\in\N[m]$ мы получим диаграмму
 $$
  \xymatrix @R=2pc @C=1.2pc
 {
 {\mathcal E}(M)\circledast A\ar[rr]^{\sigma}\ar[dr]_{D_k} & & C \ar@{-->}[dl]^{D_k'}\\
  & {\mathcal C}(K,B[n]) &
 }
 $$

Зафиксируем какой-нибудь элемент $c\in C$. Поскольку $\sigma:{\mathcal E}(M)\circledast A\to C$ -- плотный
эпиморфизм, найдется направленность элементов $x_i\in {\mathcal E}(M)\circledast A$ такая, что
$$
\sigma(x_i)\overset{C}{\underset{i\to\infty}{\longrightarrow}} c.
$$
Для всякого индекса $k\in\N[m]$ мы получим
 \beq\label{D_k(alpha_i)->D_k'(c)}
D_k(x_i)=D_k'(\sigma(x_i))\overset{\mathcal{C}(K)}{\underset{i\to\infty}{\longrightarrow}}
D_k'(c).
 \eeq
Теперь рассмотрим какую-нибудь гладкую кривую в $K$, точнее, гладкое
отображение $\gamma:[0,1]\to K$. Пусть для всякого индекса $k\in\N[m]$ порядка
$|k|=1$ и любой точки $t\in[0,1]$ символ $\gamma^k(t)$ обозначает $k$-ю
компоненту производной $\gamma'(t)$ в разложении по локальным координатам на
$U$. Для всякой функции $u\in{\mathcal E}(M)$ мы по теореме Ньютона-Лейбница получим
$$
\varPhi_0(u)(\gamma(1))-\varPhi_0(u)(\gamma(0))=\sum_{|k|=1}\int_0^1\gamma^k(t)\cdot
\varPhi_k(u)(\gamma(t))\d t.
$$
Помножив $u$ на произвольный элемент $a\in A$, мы получим
\begin{multline*}
D_0(u\circledast a)(\gamma(1))-D_0(u\circledast a)(\gamma(0))=\varPhi_0(u)(\gamma(1))\circledast\eta(a) -\varPhi_0(u)(\gamma(0))\circledast\eta(a)=\\=\sum_{|k|=1}\int_0^1\gamma^k(t)\cdot
\varPhi_k(u)(\gamma(t))\circledast\eta(a)\d t=\sum_{|k|=1}\int_0^1\gamma^k(t)\cdot
D_k(u\circledast a)(\gamma(t))\d t.
\end{multline*}
Поскольку элементы вида $u\circledast a$ полны в ${\mathcal E}(M)\circledast A$, мы можем заменить их в этом равенстве на произвольный элемент $x\in{\mathcal E}(M)\circledast A$.

Вместе с \eqref{D_k(alpha_i)->D_k'(c)} это дает
 \begin{multline*}
D_0'(c)(\gamma(1))-D_0'(c)(\gamma(0)) \underset{\infty\gets i}{\longleftarrow}
D_0(x_i)(\gamma(1))-D_0(x_i)(\gamma(0))=\\
=\sum_{|k|=1}\int_0^1\gamma^k(t)\cdot
D_k(x_i)(\gamma(t))\d t \underset{i\to\infty}{\longrightarrow}
\sum_{|k|=1}\int_0^1\gamma^k(t)\cdot D_k'(c)(\gamma(t))\d t
 \end{multline*}
и значит
$$
D_0'(c)(\gamma(1))-D_0'(c)(\gamma(0))=\sum_{|k|=1}\int_0^1\gamma^k(t)\cdot
D_k'(c)(\gamma(t))\d t
$$
Эта связь между функцией $D_0'(c)\in {\mathcal C}(K,B)$ и функциями $D_k'(c)\in
{\mathcal C}(K,B)$, $|k|=1$, означает, что $D_0'(c)$ непрерывно дифференцируема
на $K$, причем ее частными производными в выбранных нами локальных координатах
будут функции $D_k'(c)$, $|k|=1$.

Выбрав после этого какую-нибудь из производных $D_k'(c)$, $|k|=1$, и рассмотрев
индексы порядка 2, мы точно тем же приемом получим, что $D_k'(c)$ также
непрерывно дифференцируема. И вообще, организовав индукцию по индексам, мы
сможем показать, что все функции $D_k'(c)$ бесконечно дифференцируемы, и
связаны между собой как частные производные функции $D_0'(c)$ (относительно
выбранных нами локальных координат).
Это означает, что должна быть коммутативна диаграмма
\eqref{DIAGR:env-E(M)-circledast-A-1}:
 \beq\label{DIAGR:env-E(M)-circledast-A-*-K}
 \xymatrix @R=2.5pc @C=2pc
 {
 {\mathcal E}(M)\circledast A\ar[rr]^{\sigma}\ar@/_5ex/[ddr]_{D}\ar[dr]_{\iota_{K,B}} & & C \ar@/^5ex/[ddl]^{D'}\ar@{-->}[dl]^{\iota_{K,B}'}\\
  & {\mathcal E}(K,B[n])\ar[d]_{\varPhi_B} & \\
  & \mathcal{C}(K,B[n])[m] &
 }
 \eeq
в которой
$$
\iota_{K,B}(u\circledast a)=u(t)\cdot\eta(a),\qquad u\in {\mathcal E}(M),\quad a\in A,
$$
$$
(\varPhi_B)_k(f)=\frac{\partial^{\abs{k}}(f\circ\ph^{-1})}{\partial t_1^{k_1}...\partial
t_d^{k_d}}\circ\ph,\qquad f\in {\mathcal E}(M,B[n]),\quad k\in\N^d,
$$
Из этой диаграммы следует, что $\iota_{K,B}'$ должно быть
непрерывно, потому что если $c_i\to c$, то это условие сохраняется под
действием каждого оператора $D_k'$, то есть $D_k'(c_i)\to D_k'(c)$, а это как
раз и есть сходимость в пространстве ${\mathcal E}(K,B[n])$.

Если теперь менять компакт $K\subset U$ и открытое множество $U\subseteq M$, то
возникающие при этом гладкие функции $D_0'(c)$ на $K$ будут согласованы между
собой тем, что на пересечении своих областей определения они совпадают. Поэтому
определена некая общая гладкая функция $\iota_B'(c):M\to B[n]$, обладающая тем
свойством, что ее ограничение на каждый компакт $K$ будет совпадать с
соответствующей функцией $D_0'(c)$:
$$
\iota'(c)\big|_K=D_0'(c),\qquad K\subset U\subseteq M.
$$
а частные производные при выбранной системе локальных координат совпадают с
действием операторов $D_k'$ на $c$. Иными словами, определено некое отображение
$\iota'_B:C\to{\mathcal E}(M,B[n])$ (по построению это будет гомоморфизм алгебр),
для которого будет коммутативна следующая диаграмма, уточняющая
\eqref{DIAGR:env-E(M)-circledast-A-*-K}:
 \beq\label{DIAGR:env-E(M)-circledast-A-*-M-0}
 \xymatrix @R=2.5pc @C=2pc
 {
 {\mathcal E}(M)\circledast A\ar[rr]^{\sigma}\ar@/_5ex/[ddr]_{\iota_{K,B}}\ar[dr]_{\iota_B} & & C \ar@/^5ex/[ddl]^{\iota'_{K,B}}\ar@{-->}[dl]^{\iota'_B}\\
  & {\mathcal E}(M,B[n])\ar[d]_{\rho_K} & \\
  & {\mathcal E}(K,B[n]) &
 }
 \eeq
(здесь $\rho_K$ -- отображение ограничения на компакт $K$).

Пусть теперь $U$ -- дифференциальная окрестность нуля в $A$, соответствующая гомоморфизму $\eta:A\to B[n]$. Из того, что $\sigma$ -- плотный эпиморфизм, следует, что верхний внутренний треугольник в \eqref{DIAGR:env-E(M)-circledast-A-*-M-0} можно достроить до диаграммы
 \beq\label{DIAGR:env-E(M)-circledast-A-*-M}
 \xymatrix @R=2.5pc @C=2pc
 {
 {\mathcal E}(M)\circledast A\ar[rr]^{\sigma}\ar@/_5ex/[ddr]_{\iota_B}\ar[dr]_{\vartheta_U} & & C \ar@/^5ex/[ddl]^{\iota'_B}\ar@{-->}[dl]^{\vartheta'_U}\\
  & {\mathcal E}(M,A/U)\ar[d]_{\eta_U\oslash 1_M} & \\
  & {\mathcal E}(M,B[n]) &
 }
 \eeq
где $\eta_U:A/U\to B[n]$ -- морфизм из \eqref{ph=ph_U-circ-pi_U}, и
$$
\vartheta_U(u\circledast a)(t)=u(t)\cdot\pi_U(a),\qquad u\in{\mathcal E}(M),\quad a\in A,\quad t\in M,
$$
$$
(\eta_U\oslash 1_M)(h)(t)=\eta_U(h(t)),\qquad h\in {\mathcal E}(M,A/U),\quad t\in M.
$$
Из определения $\vartheta_U$ сразу следует, что если $U'\subseteq U$ -- какая-то другая дифференциальная окрестность нуля, то
\beq\label{vartheta_U=(varkappa^U'_U-oslash-1_M)-cdot-vartheta_U'}
\vartheta_U=(\varkappa^{U'}_U\oslash 1_M)\cdot\vartheta_{U'},\qquad U\supseteq U',
\eeq
где $\varkappa^{U'}_U$ -- морфизм из \eqref{varkappa^U'_U}, и
$$
(\varkappa^{U'}_U\oslash 1_M)(h)(t)=\varkappa^{U'}_U(h(t)),\qquad h\in {\mathcal E}(M,A/U'),\quad t\in M.
$$
Равенство \eqref{vartheta_U=(varkappa^U'_U-oslash-1_M)-cdot-vartheta_U'} будет левым нижним внутренним треугольником в диаграмме $$
 \xymatrix @R=2.5pc @C=3pc
 {
{\mathcal E}(M)\circledast A\ar[rr]^{\sigma}\ar[dr]_(.6){\vartheta_{U'}}\ar@/_5ex/[ddr]_{\vartheta_U} & & C \ar[dl]^(.6){\vartheta'_{U'}} \ar@/^5ex/[ddl]^{\vartheta'_U}\\
 & {\mathcal E}(M,A/U')\ar@{-->}[d]^{\varkappa^{U'}_U\oslash 1_M} & \\
  & {\mathcal E}(M,A/U) &
 }
$$
При этом периметр и верхний внутренний треугольник здесь будут вариантами верхнего внутреннего треугольника в \eqref{DIAGR:env-E(M)-circledast-A-*-M}, и вдобавок $\sigma$ -- эпиморфизм. Как следствие, оставшийся правый нижний внутренний треугольник тоже должен быть коммутативен.

Это означает, что морфизмы $\vartheta'_U:C\to {\mathcal E}(M,A/U)$ образуют проективный конус системы $\varkappa^{U'}_U\oslash 1_M$, и поэтому существует морфизм $\vartheta'$ в проективный предел:
$$
 \xymatrix @R=2.5pc @C=3pc
 {
{\mathcal E}(M)\circledast A\ar[rr]^{\sigma}\ar[dr]_(.6){\vartheta}\ar@/_5ex/[ddr]_{\vartheta_U} & & C \ar@{-->}[dl]^(.6){\vartheta'} \ar@/^5ex/[ddl]^{\vartheta'_U}\\
 & \projlim\limits_{U'\to 0}{\mathcal E}(M,A/U')\ar[d]^{\varkappa_U\oslash 1_M} & \\
  & {\mathcal E}(M,A/U) &
 }
$$
Теперь заметим цепочку
$$
\projlim\limits_{U'\to 0}{\mathcal E}(M,A/U')=\projlim\limits_{U'\to 0}({\mathcal E}(M)\odot A/U')=\cite[(4.153)]{Akbarov-De-Gruyter-I}=
{\mathcal E}(M)\odot \projlim\limits_{U'\to 0}A/U'={\mathcal E}(M,\projlim\limits_{U'\to 0}A/U')
$$
и подставим последнее пространство в нашу диаграмму:
$$
 \xymatrix @R=2.5pc @C=3pc
 {
{\mathcal E}(M)\circledast A\ar[rr]^{\sigma}\ar[dr]_(.6){\vartheta}\ar@/_5ex/[ddr]_{\vartheta_U} & & C \ar@{-->}[dl]^(.6){\vartheta'} \ar@/^5ex/[ddl]^{\vartheta'_U}\\
 & {\mathcal E}(M,\projlim\limits_{U'\to 0}A/U')\ar[d]^{\varkappa_U\oslash 1_M} & \\
  & {\mathcal E}(M,A/U) &
 }
$$

Еще раз вспомним, что $\sigma$ -- плотный эпиморфизм. Из этого следует, что стрелка $\vartheta'$ поднимается до некоторой стрелки $\upsilon$ со значениями в пространстве ${\mathcal E}(M,\Im\pi)$ функций, принимающих значения в образе отображения $\pi:A\to\projlim\limits_{U'\to 0}A/U'$, или, что то же самое, в непосредственном подпространстве, порожденном множеством значений отображения $\pi$, а это пространство как раз совпадает с $A$, поскольку $A$ -- гладкая алгебра:
$$
\Im\pi\cong\Env_{\mathcal E}A\cong A
$$
Мы получаем диаграмму
$$
 \xymatrix @R=2.5pc @C=3pc
 {
{\mathcal E}(M)\circledast A\ar[rr]^{\sigma}\ar[dr]_(.6){\iota}\ar@/_5ex/[ddr]_{\vartheta} & & C \ar@{-->}[dl]^(.6){\upsilon} \ar@/^5ex/[ddl]^{\vartheta'}\\
 & {\mathcal E}(M,A)\ar[d]^{\im\pi\oslash 1_M} & \\
  & {\mathcal E}(M,\projlim\limits_{U'\to 0}A/U') &
 }
$$
где $\pi$ -- морфизм из \eqref{A-to-leftlim-A/U}.
\epr

\subsection{Гладкие оболочки коммутативных алгебр}

\paragraph{${\mathcal E}(M)$, как гладкая оболочка своих подалгебр.}

Пусть $M$ -- гладкое многообразие и пусть, как и раньше ${\mathcal E}(M)$ -- алгебра гладких функций на $M$. Пусть кроме того, $T_s(M)$, $T_s^\star(M)$ и $\Jet^n_s(M)$ обозначают обычные касательное пространство, кокасательное пространство и алгебру струй в точке $s\in M$. С введенными нами на страницах \pageref{DEF:kasatelnyj-vektor}, \pageref{DEF:T_s^(C-star)[A]} и \pageref{DEF:prostr-struj} объектами эти пространства связаны равенствами
$$
T_s(M)=T_s[{\mathcal E}(M)],\quad T_s^\star(M)=T_s^\star[{\mathcal E}(M)],\quad \Jet_s^n(M)=\Jet_s^n[{\mathcal E}(M)].
$$

\btm\label{C^infty-obolochka-podalgebry-v-C^infty}
Пусть $A$ -- инволютивная стереотипная подалгебра в алгебре $\mathcal{E}(M)$ гладких функций на гладком (локально веклидовом) многообразии $M$, то есть задан (непрерывный и сохраняющий единицу) мономорфизм инволютивных стереотипных алгебр
$$
\iota:A\to \mathcal{E}(M).
$$
Для того, чтобы гладкая оболочка алгебры $A$ совпадала с
алгеброй ${\mathcal E}(M)$
\beq\label{Env_C^infty_A=C^infty(M)}
\Env_{\mathcal E} A={\mathcal E}(M)
\eeq
(то есть чтобы $\iota$ был гладкой оболочкой $A$), необходимо и достаточно выполнение следующих двух условий:
 \bit{
\item[(i)] сопряженное отображение спектров
$$
\Spec(A)\gets M
$$
является точным наложением (в смысле определения на с.\pageref{DEF:nalozhenie});

\item[(ii)] для всякой точки $s\in M$ естественное отображение касательных
пространств
 $$
T_s[A]\gets T_s(M)
 $$
является изоморфизмом (конечномерных векторных пространств).
 }\eit
\etm

Доказательство мы разобьем на несколько лемм.

\blm\label{LM:neobhodimost-v-TH:C^infty-obolochka-podalgebry-v-C^infty} Условия (i) и (ii) необходимы для того, чтобы вложение $A\subseteq{\mathcal E}(M)$ было гладким расширением алгебры $A$.
\elm
 \bpr
Обозначим через $\iota:A\to {\mathcal E}(M)$ естественное вложение $A$ в ${\mathcal E}(M)$. Предположим, что $\iota$ является гладким расширением, то есть расширением в классе $\DEpi$ инволютивных плотных эпиморфизмов относительно класса дифференциальных инволютивных гомоморфизмов в $C^*$-алгебры с присоединенными самосопряженными нильпотентами. Тогда $\iota$ является и расширением в $\DEpi$ и относительно класса гомоморфизмов в $C^*$-алгебры, поскольку $C^*$-алгебра $B$ может считаться $C^*$-алгеброй с присоединенным пустым множеством нильпотентов: $B=B[0]$. То есть $\iota$ является непрерывным расширением. Тем же приемом, что и в теореме \eqref{C-obolochka-podalgebry-v-C(M)} доказывается, что отображение спектров $\iota^{\Spec}:\Spec(A)\gets M$
является точным наложением.

Поэтому неочевидным здесь будет только выполнение условия (ii). Рассмотрим алгебру $\C_1[[1]]$ многочленов
степени 1 от одной переменной. Как векторное пространство она изоморфна прямой
сумме $\C\oplus\C$, а морфизмы $A\to \C_1[[1]]$ представляют собой пары
$(s,\sigma)$, где $s\in\Spec(A)$ -- точка спектра, а $\sigma\in T_s[A]$ --
касательный вектор в этой точке. Для точки $t\in M$, о которой речь идет
в формулировке, мы получаем коммутативную диаграмму
$$
 \xymatrix @R=2pc @C=1.2pc
 {
 A\ar[rr]^{\iota}\ar[dr]_{(t\circ\iota,\tau\circ\iota)} & & {\mathcal E}(M)\ar[dl]^{(t,\tau)}\\
  & \C_1[[1]] &
 }
$$
Поскольку $\iota:A\to {\mathcal E}(M)$ есть гладкое расширение, всякой стрелке
$A\to\C_1[[1]]$ соответствует единственная стрелка ${\mathcal E}(M)\to\C_1[[1]]$, и это означает,
что отображение касательных пространств
$\tau\mapsto\iota^\star(\tau)=\tau\circ\iota$ является биекцией и изоморфизмом векторных пространств. \epr

\blm\label{LM:Ker-rho_t^(T^star)=0} Пусть $\iota:A\to B$ -- гомоморфизм
инволютивных стереотипных алгебр, и для точки $t\in\Spec(B)$ соответствующий
морфизм кокасательных пространств $\iota_t^{\C
T^\star}:{\C}T_{t\circ\iota}^\star[A]\to{\C}T_t^\star[B]$ инъективен:
 \beq\label{Ker-rho_t^(T^star)=0}
\Ker\iota_t^{T^\star}=0
 \eeq
Тогда
 \beq\label{overline(I_s^2)[A]=rho^(-1)(overline(I_t^2)[B])}
\overline{I_{t\circ\iota}^2}[A]=\iota^{-1}\Big(\overline{I_t^2}[B]\Big).
 \eeq
\elm \bpr Здесь прямое вложение следует сразу из гомоморфности (и
непрерывности) $\iota$ (мы пользуемся обозначениями на
с.\pageref{DEF:M-cdot-N}):
 \begin{multline*}
\iota\Big(I_{t\circ\iota}[A]\Big)\subseteq I_t[B]\quad\Longrightarrow\\
\Longrightarrow\quad
\iota\Big(I_{t\circ\iota}^2[A]\Big)=\iota\Big(I_{t\circ\iota}[A]\cdot
I_{t\circ\iota}[A]\Big)=
\iota\Big(I_{t\circ\iota}[A]\Big)\cdot\iota\Big(I_{t\circ\iota}[A]\Big) \subseteq
I_s[B]\cdot I_s[B]=
I_s^2[B]\subseteq \overline{I_t^2}[B] \quad\Longrightarrow\\
\Longrightarrow\quad
I_{t\circ\iota}^2[A]\subseteq\iota^{-1}\Big(\overline{I_t^2}[B]\Big)
\quad\Longrightarrow\quad
\overline{I_{t\circ\iota}^2}[A]\subseteq\iota^{-1}\Big(\overline{I_t^2}[B]\Big).
 \end{multline*}
Полученное вложение можно переписать в виде
$$
\iota\Big(\overline{I_{t\circ\iota}^2}[A]\Big)\subseteq \overline{I_t^2}[B],
$$
и мы можем сделать вывод, что справедлива диаграмма
$$
 \xymatrix @R=2pc @C=3pc
 {
 & \overline{I_{t\circ\iota}^2}[A]\ar[r]^{\iota}\ar[d]_{\sigma_A} & \overline{I_t^2}[B]\ar[d]^{\sigma_B} &\\
 & I_{t\circ\iota}[A]\ar[r]^{\iota}\ar[d]_{\pi_A} & I_t[B]\ar[d]^{\pi_B} & \\
 {\C}T_{t\circ\iota}[A]\ar@{=}[r] & \Big(I_{t\circ\iota}[A]\Big/\overline{I_{t\circ\iota}^2}[A]\Big)^\triangledown\ar[r]^{\iota_t^{T^\star}} &
 \Big(I_t[B]\Big/\overline{I_t^2}[B]\Big)^\triangledown \ar@{=}[r] &  {\C}T_t[B]
 }
$$
Из нее следует обратная цепочка, необходимая для доказательства
\eqref{overline(I_s^2)[A]=rho^(-1)(overline(I_t^2)[B])}:
 \begin{multline*}
a\in\iota^{-1}\Big(\overline{I_t^2}[B]\Big) \quad\Longrightarrow\quad \iota(a)\in
\overline{I_t^2}[B] \quad\Longrightarrow\quad
\iota_t^{T^\star}(\pi_A(a))=\pi_B(\iota(a))=0
\quad\Longrightarrow\\
\Longrightarrow\quad
\pi_A(a)\in\Ker\iota_t^{T^\star}=\eqref{Ker-rho_t^(T^star)=0}=0
\quad\Longrightarrow\quad a\in\Ker(\pi_A)=I_{t\circ\iota}[A].
 \end{multline*}
\epr

\blm\label{LM-3-0-dlya-TH:C^infty-obolochka-podalgebry-v-C^infty}
При выполнении условий (i) и (ii) теоремы \ref{C^infty-obolochka-podalgebry-v-C^infty} естественный морфизм расслоений
$$
\Jet^n_{A}[{\mathcal E}(M)]\longleftarrow \Jet^n_{{\mathcal E}(M)}[{\mathcal E}(M)]
$$
является послойным изоморфизмом.
\elm
\bpr
По теореме Нахбина \ref{TH:Nachbin} алгебра $A$ плотна в алгебре ${\mathcal E}(M)$. Поэтому по лемме \ref{LM:o-plotnom-ideale}, идеал $I_t(A)=\{a\in A: a(t)=0\}$ плотен в идеале $I_t({\mathcal E}(M))=\{u\in {\mathcal E}(M): u(t)=0\}$. Отсюда
$$
\overline{I_t^n(A)}=\overline{I_t^n({\mathcal E}(M))}=I_t^n({\mathcal E}(M))
$$
и поэтому
$$
{\mathcal E}(M)/\overline{I_t^n(A)}={\mathcal E}(M))/\overline{I_t^n({\mathcal E}(M))}={\mathcal E}(M)/I_t^n({\mathcal E}(M)).
$$
\epr

\blm\label{LM-3-dlya-TH:C^infty-obolochka-podalgebry-v-C^infty}
Если $M$ -- гладкое многообразие и $\iota:A\to{\mathcal E}(M)$ -- мономорфизм инволютивных стереотипных алгебр, то (независимо от выполнения (i)) условие (ii) теоремы \ref{C^infty-obolochka-podalgebry-v-C^infty} эквивалентно условиям
\bit{

\item[(iii)] для всякой точки $s\in M$ естественное отображение кокасательных
пространств
 $$
T_s^\star[A]\to T_s^\star(M)
 $$
является изоморфизмом (конечномерных векторных пространств),

\item[(iv)] для всякой точки $s\in M$ и любого числа $n\in\N$ естественное
отображение пространств струй
 $$
\Jet_s^n[A]\to \Jet_s^n(M)
 $$
является изоморфизмом (конечномерных векторных пространств).
 }\eit
и влечет за собой условие
\bit{

\item[(v)] для всякого $n\in\N$ существует непрерывное отображение $\mu:\Jet^n_A[A]\gets\Jet^n_{{\mathcal E}(M)}[{\mathcal E}(M)]$, которое в паре с отображением спектров $\iota^{\Spec}:\Spec(A)\gets M$ образует морфизм расслоенией,
$$
 \xymatrix @R=2pc @C=4pc
 {
\Jet^n_A[A]\ar[d]^{\pi^n_{A,A}} &  \Jet^n_{{\mathcal E}(M)}[{\mathcal E}(M)]\ar[d]^{\pi^n_{{\mathcal E}(M),{\mathcal E}(M)}}\ar[l]_{\mu} \\
\Spec(A)  & M\ar[l]_{\iota^{\Spec}}
 }
$$
    удовлетворяющий тождеству
    \beq
    \mu(\jet^n(a\circ\iota^{\Spec})(t))=\jet^n(a)(\iota^{\Spec}(t)),\qquad a\in A,\quad t\in M.
    \eeq
 }\eit\noindent
Если вдобавок выполнено условие (i) теоремы \ref{C^infty-obolochka-podalgebry-v-C^infty}, то морфизм расслоений $\mu$ является биекцией.

\elm
\bpr
 1. Покажем сначала, что условия (ii)-(iv) эквивалентны. Эквивалентность (ii) и (iii) следует сразу
из \cite[Theorem 5.1.27]{Akbarov-De-Gruyter-I}: если $T_s[A]$ и $T_s(M)$ изоморфны как
векторные пространства, то из-за конечномерности $T_s(M)$ они изоморфны и как
стереотипные пространства, значит их сопряженные пространства $T_s^\star[A]$ и
$T_s^\star(M)$ также изоморфны (как стереотипные пространства и как векторные
пространства). И наоборот.

С другой стороны, из условия (iv) следует условие (iii), потому что при
изоморфизме
$$
\Jet_s^1[A]\cong \Jet_s^1(M)
$$
идеал $I_s[A]\cap \Jet_s^1[A]$ переходит в идеал $I_s(M)\cap \Jet_s^1(M)$, а это как
раз означает изоморфизм
 \begin{multline*}
T_s^\star[A]=\Big(I_s[A]\Big/\overline{I_s^2}[A]\Big)^\triangledown=I_s[A]\Big/\overline{I_s^2}[A]
=I_s[A]\cap \Jet_s^1[A]\cong I_s(M)\cap \Jet_s^1(M)=\\=
I_s(M)\Big/\overline{I_s^2}(M)=
\Big(I_s(M)\Big/\overline{I_s^2}(M)\Big)^\triangledown=T_s^\star(M)
 \end{multline*}

Таким образом в эквивалентностях $(ii)\Longleftrightarrow
(iii)\Longleftrightarrow (iv)$ неочевидным звеном остается только импликация
$(iii)\Longrightarrow (iv)$. Докажем ее. Пусть выполняется (iii). Тогда, прежде
всего, выполняется равенство \eqref{Ker-rho_t^(T^star)=0}, которое
применительно к данному случаю можно записать так:
 \beq\label{overline(I_s^2)[A]=A-cap-I_s^2(M)}
\overline{I_s^2}[A]=A\cap I_s^2(M)
 \eeq
Далее, поскольку пространство $T_s^\star[A]=\Big(
I_s[A]\Big/\overline{I_s^2}[A]\Big)^\triangledown$ конечномерно, оно совпадает
с пространством $I_s[A]\Big/\overline{I_s^2}[A]$, которое также конечномерно, и
поэтому  у него можно выбрать конечный базис. То есть должна существовать
конечная последовательность векторов $e^1,...,e^d\in I_s[A]$, для которых
классы эквивалентности $e^1+\overline{I_s^2}[A],...,e^d+\overline{I_s^2}[A]$
образуют базис в $I_s[A]\Big/\overline{I_s^2}[A]$.
Выберем окрестность $U$ точки $s$, в которой функции $e^1,...,e^d$ образуют локальную карту многообразия $M$. После этого для любой точки $t\in U$ и любого мультииндекса
$k\in\N^d$ обозначим
 \beq\label{e^k=(e^1)^(k_1)...(e^d)^(k_d)}
e_t^k=(e^1-e^1(t))^{k_1}\cdot...\cdot (e^d-e^d(t))^{k_d}.
 \eeq
Пусть для любой функции $f\in{\mathcal E}(M)$ и любого номера $n\in\N$
символ $E_s^n[f]$ обозначает линейную комбинацию функций $e^k$, $\abs{k}\le n$,
$$
E_s^n[f]=\sum_{|k|\le n}\lambda_k\cdot e_s^k.
$$
имеющую в точке $s$ ту же струю порядка $n$, что и функция $f$:
 \beq\label{f-stackrel(I_s^(n+1)(M))(equiv)E_s^n[f]}
f\stackrel{I_s^{n+1}(M)}{\equiv} E_s^n[f]
 \eeq
(такая функция существует и единственна, потому что функции $e_s^1,...,e_s^d$
образуют локальную карту в некоторой окрестности $V_s$ точки $s$, которая диффеоморфно отображает $V_s$ на некоторую окрестность нуля в $\R^d$, и при этом
диффеоморфизме функции $E_s^n[f]$ превращаются в точности в многочлены Тейлора
функции $f$ в точке 0).

2. Заметим, что операция $f\mapsto E_s^n[f]$ мультипликативна
 \beq\label{E_s^n[f-cdot-g]=E_s^n[f]-cdot-E_s^n[g]}
E_s^n[f\cdot g]=E_s^n[f]\cdot E_s^n[g],\qquad f,g\in{\mathcal E}(M)
 \eeq
обладает вариантом свойства идемпотентности
 \beq\label{idempotentnost-E_s^n[f]}
E_s^q\big[E_s^p[f]\big]=E_s^p[f]=E_s^p\big[E_s^q[f]\big],\qquad p\le
q\in\N,\quad f\in{\mathcal E}(M),
 \eeq
и непрерывна в $A$:
 \beq\label{nepreryvnost-E_s^n}
a_i\overset{A}{\underset{i\to\infty}{\longrightarrow}}a
\quad\Longrightarrow\quad
E_s^n[a_i]\overset{A}{\underset{i\to\infty}{\longrightarrow}} E_s^n[a].
 \eeq
Первые два свойства следуют из представления $E_s^n[f]$ в виде многочленов
Тейлора при диффеоморфизме, образованном локальной картой $e_1,...,e_d$, а
третье -- из конечномерности алгебры многочленов с фиксированной степенью: если
$a_i\overset{A}{\underset{i\to\infty}{\longrightarrow}}a$, то
$a_i\overset{{\mathcal E}(M)}{\underset{i\to\infty}{\longrightarrow}}a$,
поэтому
$E_s^n[a_i]\overset{{\mathcal E}(M)}{\underset{i\to\infty}{\longrightarrow}}E_s^n[a]$,
и поскольку левая и правая часть принадлежит конечномерному пространству
$P=\sp\{e^k;\ |k|\le n\}$, сходимость в ${\mathcal E}(M)$ здесь можно
заменить сходимостью в $P$, а затем и сходимостью в $A$.

3. Докажем теперь формулу
 \beq\label{a-stackrel(overline(I_s^(n+1))[A])(equiv)E_s^n[a]}
a\stackrel{\overline{I_s^{n+1}}[A]}{\equiv} E_s^n[a],\qquad a\in
\overline{I_s^n}[A].
 \eeq
Это делается по индукции. Прежде всего, эта формула верна для $n=0$,
 \beq\label{a-stackrel(I_s[A])(equiv)E_s^0[a]}
a\stackrel{I_s[A]}{\equiv} E_s^0[a],\qquad a\in A,
 \eeq
потому что
$$
a-E_s^0[a]=a-a(s)\cdot 1\in I_s[A].
$$
И при $n=1$,
 \beq\label{a-stackrel(overline(I_s^2)[A])(equiv)E_s^2[a]}
a\stackrel{\overline{I_s^2}[A]}{\equiv} E_s^1[a],\qquad a\in\overline{I_s}[A],
 \eeq
потому что из \eqref{f-stackrel(I_s^(n+1)(M))(equiv)E_s^n[f]} мы получаем
$a-E_s^1[a]\in I_s^2(M)$, а с другой стороны, $a-E_s^1[a]\in A$, и вместе это
означает, что $a-E_s^1[a]\in A\cap
I_s^2(M)=\eqref{overline(I_s^2)[A]=A-cap-I_s^2(M)}=\overline{I_s^2}[A]$.

Предположим далее, что формула
\eqref{a-stackrel(overline(I_s^(n+1))[A])(equiv)E_s^n[a]} верна для какого-то
числа $n-1$:
 \beq\label{a-stackrel(overline(I_s^n)[A])(equiv)E_s^(n-1)[a]}
a\stackrel{\overline{I_s^n}[A]}{\equiv} E_s^{n-1}[a],\qquad a\in
\overline{I_s^{n-1}}[A].
 \eeq
Пусть $a\in\overline{I_s^n}[A]$. Это означает, что найдется некая
направленность $a_i\in I_s^n[A]$, стремящаяся к $a$ в $A$:
$$
a_i\overset{A}{\underset{i\to\infty}{\longrightarrow}}a.
$$
Каждый элемент $a_i$ лежит в $I_s^n[A]$, поэтому представим в виде
$$
a_i=\sum_{j=1}^p b_i^j\cdot c_i^j,\qquad b_i^j\in I_s^{n-1}[A],\quad c_i^j\in
I_s[A].
$$
Из \eqref{a-stackrel(overline(I_s^n)[A])(equiv)E_s^(n-1)[a]} и
\eqref{a-stackrel(overline(I_s^2)[A])(equiv)E_s^2[a]} мы получаем цепочку:
$$
b_i^j\stackrel{\overline{I_s^n}[A]}{\equiv} E_s^{n-1}[b_i^j],\quad
c_i^j\stackrel{\overline{I_s^2}[A]}{\equiv} E_s^1[c_i^j],
$$
$$
\Downarrow
$$
$$
b_i^j\cdot c_i^j\stackrel{\overline{I_s^{n+1}}[A]}{\equiv}
E_s^{n-1}[b_i^j]\cdot E_s^1[c_i^j]=E_s^n[b_i^j\cdot c_i^j],
$$
$$
\Downarrow
$$
$$
a \overset{A}{\underset{\infty\gets i}{\longleftarrow}} a_i=\sum_{j=1}^p
b_i^j\cdot c_i^j\stackrel{\overline{I_s^{n+1}}[A]}{\equiv} \sum_{j=1}^p
E_s^n[b_i^j\cdot
c_i^j]=E_s^n[a_i]\overset{\eqref{nepreryvnost-E_s^n}}{\overset{A}{\underset{i\to\infty}{\longrightarrow}}}E_s^n[a],
$$
$$
\Downarrow
$$
$$
a\stackrel{\overline{I_s^{n+1}}[A]}{\equiv}E_s^n[a].
$$

После того, как формула
\eqref{a-stackrel(overline(I_s^(n+1))[A])(equiv)E_s^n[a]} доказана, в ней можно
заменить $a\in \overline{I_s^{n-1}}[A]$ на $a\in A$:
 \beq\label{a-stackrel(overline(I_s^(n+1))[A])(equiv)E_s^n[a],a-in-A}
a\stackrel{\overline{I_s^{n+1}}[A]}{\equiv} E_s^n[a],\qquad a\in A.
 \eeq
Действительно, при $n=0$ это превращается в уже отмечавшуюся формулу
\eqref{a-stackrel(I_s[A])(equiv)E_s^0[a]}. Далее, если
\eqref{a-stackrel(overline(I_s^(n+1))[A])(equiv)E_s^n[a],a-in-A} доказана для
какого-то $n-1$,
$$
a\stackrel{\overline{I_s^n}[A]}{\equiv} E_s^{n-1}[a],\qquad a\in A,
$$
то, записав это в виде
$$
a-E_s^{n-1}[a]\in \overline{I_s^n}[A],
$$
мы в силу \eqref{a-stackrel(overline(I_s^(n+1))[A])(equiv)E_s^n[a]} получим:
$$
a-E_s^{n-1}[a]\stackrel{\overline{I_s^{n+1}}[A]}{\equiv}
E_s^n\big[a-E_s^{n-1}[a]\big]=\eqref{idempotentnost-E_s^n[f]}=E_s^n[a]-E_s^{n-1}[a],
$$
откуда
$$
a\stackrel{\overline{I_s^{n+1}}[A]}{\equiv} E_s^n[a].
$$

4. Наконец, из формулы
\eqref{a-stackrel(overline(I_s^(n+1))[A])(equiv)E_s^n[a],a-in-A} следует, что
фактор-отображение $a\mapsto a+\overline{I_s^{n+1}}[A]$ пропускается через
отображение $a\mapsto E_s^n[a]$, и поэтому фактор-алгебра $\Jet_s^n[A]$ изоморфна
образу $E_s^n[A]$ отображения $a\mapsto E_s^n[a]$:
$$
\Jet_s^n[A]\cong E_s^n[A].
$$
С другой стороны, алгебра струй $\Jet_s^n(M)$ на многообразии $M$ также изоморфна
$E_s^n[A]$:
$$
\Jet_s^n(M)\cong E_s^n[A].
$$
Вместе то и другое означает выполнение условия (iv).

5. Докажем теперь, что условия (ii), (iii), (iv) влекут (v). Мы поняли выше, что в окрестности $U$ точки $t\in M$ расслоения $\Jet_A^n[A]$ и $\Jet_{C^\infty(M)}^n[C^\infty(M)]$ совпадают как множества. Пусть теперь
$$
\xi_i\overset{\Jet_{C^\infty(M)}^n[C^\infty(M)]}{\underset{i\to\infty}{\longrightarrow}}\xi,
$$
Обознчим $s_i=\pi(\xi_i)$, $s=\pi(\xi)$ и для всякого $i$ рассмотрим линейные комбинации
$$
h_i=\sum_{|k|\le n}\lambda_k\cdot e_{s_i}^k: \qquad h_i=E_{s_i}[h_i],
$$
и
$$
h=\sum_{|k|\le n}\lambda_k\cdot e_s^k: \qquad h_i=E_s[h].
$$
Отображение спектров $\Spec(A)\gets\Spec(C^\infty(M))=M$ непрерывно, поэтому мы получаем цепочку:
$$
s_i=\pi(\xi_i)\overset{M}{\underset{i\to\infty}{\longrightarrow}}\pi(\xi)=s\in U
$$
$$
\Downarrow
$$
$$
s_i=\pi(\xi_i)\overset{\Spec(A)}{\underset{i\to\infty}{\longrightarrow}}\pi(\xi)=s
$$
$$
\Downarrow
$$
$$
e_{s_i}^k=(e^1-e^1(s_i))^{k_1}\cdot...\cdot (e^d-e^d(s_i))^{k_d}\overset{A}{\underset{i\to\infty}{\longrightarrow}}
(e^1-e^1(s))^{k_1}\cdot...\cdot (e^d-e^d(s))^{k_d}=e_s^k
$$
$$
\Downarrow
$$
$$
h_i=\sum_{|k|\le n}\lambda_k\cdot e_{s_i}^k\overset{A}{\underset{i\to\infty}{\longrightarrow}}h=\sum_{|k|\le n}\lambda_k\cdot e_{s}^k
$$
$$
\Downarrow
$$
$$
\xi_i=\jet^n(h_i)\overset{\Jet_A^n[A]}{\underset{i\to\infty}{\longrightarrow}}\jet^n(h)=\xi.
$$

6. Если выполнено условие (i) теоремы \ref{C^infty-obolochka-podalgebry-v-C^infty}, то есть отображение спектров $\iota^{\Spec}:M\to\Spec(A)$ является биекцией, то морфизм расслоений $\mu$ является биекцией, потому что оно будет биекцией между слоями и биекцией на каждом слое.
\epr

\blm\label{LM:dostatochnost-v-TH:C^infty-obolochka-podalgebry-v-C^infty} Условия (i) и (ii) достаточны для того, чтобы морфизм $\iota:A\to{\mathcal E}(M)$ был гладкой оболочкой алгебры $A$.
\elm
\bpr 1. Сначала проверим, что при выполнении (i) и (ii) морфизм $\iota:A\to
{\mathcal E}(M)$ является гладким расширением.  Из теоремы Нахбина \ref{TH:Nachbin} следует, что при выполнении (i) и (ii) естественное вложение $\iota:A\to {\mathcal E}(M)$ становится плотным
эпиморфизмом. Пусть далее $B$ -- $C^*$-алгебра и $\ph:A\to B[m]$ -- дифференциальный гомоморфизм, то есть такой, для которого соответствующая система частных производных
$\{D_k;\ k\in\N[m]\}$ состоит из дифференциальных операторов $D_k:A\to B$ с порядками $\ord D_k\le |k|$.
Существование дифференциального гомоморфизма $\ph'$, замыкающего диаграмму
$$
 \xymatrix @R=2pc @C=1.2pc
 {
 A\ar[rr]^{\iota}\ar[dr]_{\ph} & & {\mathcal E}(M)\ar@{-->}[dl]^{\ph'}\\
  & B[m] &
 }
$$
эквивалентно существованию системы дифференциальных частных производных
$\{D_k';\ k\in\N[m]\}$, продолжающих операторы $D_k$ с $A$ на
${\mathcal E}(M)$ и принимающих значения в $B$:
$$
 \xymatrix @R=2pc @C=1.2pc
 {
 A\ar[rr]^{\iota}\ar[dr]_{D_k} & & {\mathcal E}(M)\ar@{-->}[dl]^{D_k'}\\
  & B &
 }
$$
Обозначим
$$
C=\overline{D_0(A)},\qquad
F=Z^1(D_0)=\{b\in B:\ \forall a\in A\quad [b,D_0(a)]=0\}=C^!.
$$
Поскольку алгебра $A$ коммутативна, $C$ тоже коммутативна, отсюда следует, что
$$
C\subseteq C^!=F.
$$
Это значит, что $D_0$ принимает значение в $F$. С другой стороны, по теореме \ref{TH:differentsialnye-chast-proizv}, все операторы $\{D_k;\ k\in\N[m], \ k>0\}$ также принимают значения в алгебре $F$. Поэтому нам достаточно доказать
существование системы дифференциальных частных производных $\{D_k';\
k\in\N[m]\}$, продолжающих операторы $D_k$ с $A$ на ${\mathcal E}(M)$ и
принимающих значения в $F$:
 \beq\label{prodolzhenie-D_k-na-C^infty}
 \xymatrix @R=2pc @C=1.2pc
 {
 A\ar[rr]^{\iota}\ar[dr]_{D_k} & & {\mathcal E}(M)\ar@{-->}[dl]^{D_k'}\\
  & F &
 }
 \eeq

2. Это нужно сначала доказать для индекса $k=0$:
 \beq\label{prodolzhenie-D_0-na-C^infty}
 \xymatrix @R=2pc @C=1.2pc
 {
 A\ar[rr]^{\iota}\ar[dr]_{D_0} & & {\mathcal E}(M)\ar@{-->}[dl]^{D_0'}\\
  & F &
 }
 \eeq
Алгебру ${\mathcal E}(M)$ можно вложить в алгебру ${\mathcal C}(M)$, тогда $\iota$ можно будет считать гомоморфизмом из $A$ в ${\mathcal C}(M)$, и условие (i) будет означать, что $\iota$ удовлетворяет условию теоремы  \ref{C-obolochka-podalgebry-v-C(M)}. Поэтому $D_0$ непрерывно продолжается сначала на ${\mathcal C}(M)$, а потом его можно ограничить на ${\mathcal E}(M)$, и мы получим $D_0'$.

3. После того, как построен гомоморфизм $D_0':{\mathcal E}(M)\to F$, алгебра $F$ становится модулем над ${\mathcal E}(M)$. Если обозначить через $I_t^A$ и $I_t^{\mathcal E}$ идеалы в $A$ и ${\mathcal E}(M)$, состоящие из функций, равных нулю в точке $t\in M$, то они связаны вложением
$$
I_t^A\subseteq I_t^{\mathcal E}
$$
причем по лемме \ref{LM:o-plotnom-ideale} $I_t^A$ будет плотно в $I_t^{\mathcal E}$. Поэтому
$$
\overline{I_t^A\cdot F}=\overline{I_t^{\mathcal E}\cdot F}.
$$
Обозначим $T=\Spec(C)$, тогда можно считать, что $C={\mathcal C}(T)$. Рассмотрим два случая.
\bit{
\item[---] Пусть $t\notin T$. Тогда функции из $I_t^{\mathcal E}$ аппроксимируют характеристическую функцию $\chi_T$ множества $T$, и поэтому то же справедливо для функций из $I_t^A$:
$$
\chi_T\in\overline{I_t^A}=\overline{I_t^{\mathcal E}}\subseteq {\mathcal E}(M).
$$
Отсюда следует, что
$$
\overline{I_t^A\cdot F}=\overline{I_t^{\mathcal E}\cdot F}=\chi_T\cdot F=1\cdot F=F,\qquad t\notin T.
$$
Это означает, что расслоения значений модуля $F$ над алгебрами $A$ и ${\mathcal E}(M)$ в точках вне компакта $T$ вырождены:
\beq\label{F/IF-t-notin-T}
F/\overline{I_t^A\cdot F}=F/\overline{I_t^{\mathcal E}\cdot F}=0,\qquad t\notin T=\Spec(C).
\eeq

\item[---] Пусть $t\in T$. Тогда, поскольку ${\mathcal E}(M)$ (будучи ограничено на $T$) плотно в $C={\mathcal C}(T)$, опять по лемме \ref{LM:o-plotnom-ideale} мы получаем, что идеал $I_t^{\mathcal E}$ плотен в идеале $I_t^C$ функций из $C={\mathcal C}(T)$, равных нулю в точке $t$:
$$
\overline{I_t^{\mathcal E}}=I_t^C.
$$
Поэтому
$$
\overline{I_t^A\cdot F}=\overline{I_t^{\mathcal E}\cdot F}=\overline{I_t^C\cdot F},\qquad t\in T,
$$
и это означает, что над каждой точкой $t\in T=\Spec(C)$ слои значений модуля $F$ над алгебрами $A$, ${\mathcal E}(M)$ и $C$ одинаковы:
\beq\label{F/IF-t-in-T}
F/\overline{I_t^A\cdot F}=F/\overline{I_t^{\mathcal E}\cdot F}=F/\overline{I_t^C\cdot F},\qquad t\in T=\Spec(C).
\eeq
}\eit
Помимо этого, по условию (i) компакт $T$ непрерывно (инъективно и с сохранением топологии) вкладывается в $\Spec(A)$ и в $M$. Вместе все это дает описание расслоений значений $\Jet^0_A[F]$ и $\Jet^0_{{\mathcal E}(M)}[F]$: они тривиальны вне компакта $T$, а на нем совпадают с расслоением значений $\Jet^0_C[F]$. Мы можем сделать вывод, что естественное отображение расслоений значений
$$
\Jet^0_A[F]=\bigsqcup_{t\in M} (F/\overline{I_t^A\cdot F})^\vartriangle\overset{\lambda}{\longleftarrow} \bigsqcup_{t\in M} (F/\overline{I_t^{\mathcal E}\cdot F})^\vartriangle=\Jet^0_{{\mathcal E}(M)}[F]
$$
является изоморфизмом.

4. Чтобы построить $D'_k$ для $k>0$, вспомним, что всякому дифференциальному оператору $D_k$ по теореме
\ref{TH:diff-oper->rassl-struuj} соответствует некий морфизм расслоений струй
$\jet_n[D_k]:\Jet_A^n[A]\to \Jet_A^0(F)=\pi^0_A F$, где $n=|k|$, удовлетворяющий
тождеству
$$
\jet^0(D_ka)=\jet_n[D_k]\circ \jet^n(a),\qquad a\in A.
$$
По условию (v) леммы \ref{LM-3-dlya-TH:C^infty-obolochka-podalgebry-v-C^infty}, расслоения струй алгебр $A$ и ${\mathcal E}(M)$ связаны между собой естественным морфизмом $\mu:\Jet_A^n[A]\gets \Jet^n_{{\mathcal E}(M)}[{\mathcal E}(M)]$, который вдобавок является биекцией. Мы получаем диаграмму
$$
 \xymatrix @R=2pc @C=4pc
 {
 \Jet^n_A[A]\ar[dd]_{\jet_n[D_k]}\ar[dr]^{\pi^n_{A,A}} & & & \Jet^n_{{\mathcal E}(M)}[{\mathcal E}(M)]\ar@{-->}[dd]^{\nu=\jet_n[D_k']}\ar[lll]_{\mu}\ar[dl]_{\pi^n_{{\mathcal E}(M),{\mathcal E}(M)}}  \\
  &  \Spec(A) &  M\ar[l]_{\iota^{\Spec}} &  \\
 \Jet^0_A[F]\ar[rrr]^{\lambda^{-1}}\ar[ur]_{\pi^0_{A,F}} & & & \Jet^0_{{\mathcal E}(M)}[F]  \ar[ul]^{\pi^0_{{\mathcal E}(M),F}}
 }
$$
Рассмотрим композицию (пунктирную стрелку)
$$
\nu=\lambda^{-1}\circ\jet_n[D_k]\circ\mu
$$
Она является морфизмом расслоений, поэтому по лемме \ref{TH:diff-oper<-morfizm-rassl-struuj} ей соответствует дифференциальный оператор $D_k':{\mathcal E}(M)\to F$ (над алгеброй ${\mathcal E}(M)$), удовлетворяющий тождеству
$$
\jet^0(D_k'u)=\jet_n[D_k]\circ \jet^n(u),\qquad u\in {\mathcal E}(M).
$$
Для всякого $a\in A$ мы получим \beq\label{PROOF:TH:C^infty-1}
\jet^0\big(D_k'\iota(a)\big)=\jet_n[D_k]\circ \jet^n(\iota(a))=\jet_n[D_k]\circ
\jet^n(a)=\jet^0(D_ka). \eeq Заметим далее, что по теореме
\ref{TH:B-cong-Sec(val_AB)}, отображение
$\jet^0=v:F\to\Sec(\pi^0_AF)=\Sec(\Jet^0_AF)$, переводящее $F$ в алгебру непрерывных сечений
расслоения значений $\pi^0_AF: \Jet^0_AF\to\Spec(A)$ над алгеброй $A$, является
изоморфизмом $C^*$-алгебр:
$$
F\cong\Sec(\pi^0_AF).
$$
Поэтому к \eqref{PROOF:TH:C^infty-1} можно применить оператор, обратный $\jet^0$, и мы получим равенство
$$
D_k'\iota(a)=D_ka.
$$
То есть $D_k'$ продолжает $D_k$ в диаграмме
\eqref{prodolzhenie-D_k-na-C^infty}. Кроме того, из того, что $\iota$ отображает
$A$ плотно в ${\mathcal E}(M)$ следует, что условия
\eqref{DEF:sist-chastn-proizv-*}-\eqref{DEF:sist-chastn-proizv-2} переносятся с
операторов $D_k$ на операторы $D_k'$. Поэтому семейство $\{D_k',\ k\in\N[m]\}$
представляет собой систему частных производных на ${\mathcal E}(M)$.

5. Теперь убедимся, что при выполнении (i) и (ii) расширение $\iota:A\to
{\mathcal E}(M)$ является гладкой оболочкой. Пусть $\sigma:A\to C$ -- какое-то
другое гладкое расширение. Зафиксируем локальную карту $\ph:U\to V$, где $U\subseteq M$, $V\subseteq\R^d$, и пусть и $K\subseteq U$ -- компакт, совпадающий с замыканием своей внутренности:
$\overline{\Int(K)}=K$. Тогда операторы из примера \ref{EX:sist-chast-proizv-na-C-infty(M)}
$$
D_k:{\mathcal E}(M)\to\mathcal{C}(K)\quad\Big|\quad
D_k(a)=\frac{\partial^{\abs{k}}(a\circ\ph^{-1})}{\partial t_1^{k_1}...\partial
t_d^{k_d}}\circ\ph,\qquad a\in A,\quad k\in\N[m],
$$
образуют систему частных производных из $A$ в $\mathcal{C}(K)$. Пусть $D:A\to
\mathcal{C}(K)[m]$ -- дифференциальный гомоморфизм, соответствующий этой системе
$\{D_k\}$. Поскольку $\sigma:A\to C$ есть гладкое расширение,
гомоморфизм $D:A\to C(K)[[d]]$ должен однозначно продолжаться до некоторого
гомоморфизма $D':C\to C(K)[[d]]$:
 \beq\label{DIAGR:env-A=E(M)}
 \xymatrix @R=2pc @C=1.2pc
 {
 A\ar[rr]^{\sigma}\ar[dr]_{D} & & C \ar@{-->}[dl]^{D'}\\
  & \mathcal{C}(K)[m] &
 }
 \eeq
Вернувшись назад к системе частных производных $D_k$, мы получим для всякого
индекса $k\in\N[m]$ диаграмму
 $$
  \xymatrix @R=2pc @C=1.2pc
 {
 A\ar[rr]^{\sigma}\ar[dr]_{D_k} & & C \ar@{-->}[dl]^{D_k'}\\
  & \mathcal{C}(K) &
 }
 $$

Зафиксируем какой-нибудь элемент $c\in C$. Поскольку $\sigma:A\to C$ -- плотный
эпиморфизм, найдется направленность элементов $a_i\in A$ такая, что
$$
\sigma(a_i)\overset{C}{\underset{i\to\infty}{\longrightarrow}} c.
$$
Для всякого индекса $k\in\N[m]$ мы получим
 \beq\label{D_k(a_i)->D_k'(b)}
D_k(a_i)=D_k'(\sigma(a_i))\overset{\mathcal{C}(K)}{\underset{i\to\infty}{\longrightarrow}}
D_k'(c).
 \eeq
Теперь рассмотрим какую-нибудь гладкую кривую в $K$, точнее, гладкое
отображение $\gamma:[0,1]\to K$. Пусть для всякого индекса $k\in\N[m]$ порядка
$|k|=1$ и любой точки $t\in[0,1]$ символ $\gamma^k(t)$ обозначает $k$-ю
компоненту производной $\gamma'(t)$ в разложении по локальным координатам на
$U$. Для всякой функции $a\in A$ мы по теореме Ньютона-Лейбница получим
$$
D_0(a)(\gamma(1))-D_0(a)(\gamma(0))=\sum_{|k|=1}\int_0^1\gamma^k(t)\cdot
D_k(a)(\gamma(t))\d t.
$$
Вместе с \eqref{D_k(a_i)->D_k'(b)} это дает
 \begin{multline*}
D_0'(c)(\gamma(1))-D_0'(c)(\gamma(0)) \underset{\infty\gets i}{\longleftarrow}
D_0(a_i)(\gamma(1))-D_0(a_i)(\gamma(0))=\\
=\sum_{|k|=1}\int_0^1\gamma^k(t)\cdot
D_k(a_i)(\gamma(t))\d t \underset{i\to\infty}{\longrightarrow}
\sum_{|k|=1}\int_0^1\gamma^k(t)\cdot D_k'(c)(\gamma(t))\d t
 \end{multline*}
и значит
$$
D_0'(c)(\gamma(1))-D_0'(c)(\gamma(0))=\sum_{|k|=1}\int_0^1\gamma^k(t)\cdot
D_k'(c)(\gamma(t))\d t
$$
Эта связь между функцией $D_0'(c)\in {\mathcal C}(K)$ и функциями $D_k'(c)\in
{\mathcal C}(K)$, $|k|=1$, означает, что $D_0'(c)$ непрерывно дифференцируема
на $K$, причем ее частными производными в выбранных нами локальных координатах
будут функции $D_k'(c)$, $|k|=1$.

Выбрав после этого какую-нибудь из производных $D_k'(c)$, $|k|=1$, и рассмотрев
индексы порядка 2, мы точно тем же приемом получим, что $D_k'(c)$ также
непрерывно дифференцируема. И вообще, организовав индукцию по индексам, мы
сможем показать, что все функции $D_k'(c)$ бесконечно дифференцируемы, и
связаны между собой как частные производные функции $D_0'(c)$ (относительно
выбранных нами локальных координат).

Если теперь менять компакт $K\subset U$ и открытое множество $U\subseteq M$, то
возникающие при этом гладкие функции $D_0'(c)$ на $K$ будут согласованы между
собой тем, что на пересечении своих областей определения они совпадают. Поэтому
определена некая общая гладкая функция $\iota'(c):M\to\C$, обладающая тем
свойством, что ее ограничение на каждый компакт $K$ будет совпадать с
соответствующей функцией $D_0'(c)$:
$$
\iota'(c)\big|_K=D_0'(c),\qquad K\subset U\subseteq M.
$$
а частные производные при выбранной системе локальных координат совпадают с
действием операторов $D_k'$ на $c$. Иными словами, определено некое отображение
$\iota':C\to{\mathcal E}(M)$ (по построению это будет гомоморфизм алгебр),
для которого будет коммутативна следующая диаграмма, уточняющая
\eqref{DIAGR:env-A=E(M)}:
 \beq\label{DIAGR:env-A=E(M)-*}
 \xymatrix @R=2pc @C=1.2pc
 {
 A\ar[rr]^{\sigma}\ar@/_4ex/[ddr]_{D}\ar[dr]_{\iota} & & C \ar@/^4ex/[ddl]^{D'}\ar@{-->}[dl]^{\iota'}\\
  & {\mathcal E}(M)\ar[d]_{\lambda} & \\
  & \mathcal{C}(K)[m] &
 }
 \eeq
(здесь отображение $\lambda$ есть разложение гладкой функции в многочлен Тейлора на компакте $K\subset U\subseteq M$ при выбранных локальных координатах на $U$). Из этой диаграммы следует, что $\iota'$ должно быть
непрерывно, потому что если $c_i\to c$, то это условие сохраняется под
действием каждого оператора $D_k'$, то есть $D_k'(c_i)\to D_k'(c)$, а это как
раз и есть сходимость в пространстве ${\mathcal E}(M)$.

Верхний внутренний треугольник в \eqref{DIAGR:env-A=E(M)-*} будет как раз той
диаграммой, которая нам нужна для того, чтобы убедиться, что расширение $\iota$
является оболочкой.
 \epr

\noindent\rule{160mm}{0.1pt}\begin{multicols}{2}

\paragraph{Контрпримеры.}

Следующий пример показывает, что ослабление условия (i) в теореме \ref{C^infty-obolochka-podalgebry-v-C^infty} делает это утверждение неверным.

\bex\label{EX:E(8)} {\it
Существует плотная инволютивная подалгебра $A$ в ${\mathcal E}(\R)$, у которой
\bit{
\item[(i)] сопряженное отображение спектров $\iota^{\Spec}:\Spec(A)\gets M$ является биекцией (но не наложением),

\item[(ii)] для всякой точки $s\in\R$ естественное отображение касательных пространств
 $T_s[A]\gets T_s(\R)$ является изоморфизмом (векторных пространств),

\item[(iii)] гладкая оболочка $A$ не изоморфна ${\mathcal E}(\R)$:
$$
\Env_{\mathcal E}A\not\cong {\mathcal E}(\R)
$$
}\eit
}
\eex
\bpr Это модификация примера \ref{EX:C(8)}. Рассмотрим окружность $\T=\R/\Z$ и на ней алгебру ${\mathcal E}(\T)$ гладких функций. В декартовом квадрате ${\mathcal E}(\T)^2$ рассмотрим подалгебру ${\mathcal E}(8)$, состоящую из пар функций $(u,v)\in{\mathcal E}(\T)^2$, у которых все производные в точке 0 совпадают:
\begin{multline*}
(u,v)\in{\mathcal E}(8)\quad\Longleftrightarrow\\ \Longleftrightarrow\quad u,v\in{\mathcal E}(\T)\quad\&\quad \forall k\ge 0\quad u^{(k)}(0)=v^{(k)}(0).
\end{multline*}
Это будет замкнутая подалгебра в ${\mathcal E}(\T)^2$, и ее спектром будет пространство, которое в примере \ref{EX:C(8)} мы обозначили символом $8$, и которое можно понимать как результат склеивания двух копий окружности $\T$ в точке 0:
\beq\label{Spec(E(8))}
\Spec{\mathcal E}(8)=8.
\eeq
Всякая точка $s\in 8$ должна лежать либо в первой копии окружности $\T$, которую мы будем обозначать $\T_1$, либо во второй ее копии, $\T_2$ (либо, если $s=0$, то мы считаем, что $s$ принадлежит обеим окружностям $\T_1$ и $\T_2$). Покажем, что для всякой точки $s\in 8$ касательное пространство к ${\mathcal E}(8)$ изоморфно касательному пространству к ${\mathcal E}(\T)$ в этой точке:
\beq\label{T_s(E(8))}
T_s({\mathcal E}(8))\cong\begin{cases} T_s({\mathcal E}(\T_1)), & s\in\T_1\\  T_s({\mathcal E}(\T_2)), & s\in\T_2 \end{cases},\qquad s\in 8.
\eeq

Зафиксируем для этого какую-нибудь функцию $\ph\in{\mathcal E}(\T)$, которая в некоторой окрестности $U$ точки 0 в $\R$ является обратным отображением для проекции $\pi:\R\to\R/\Z=\T$:
$$
\ph(\pi(t))=t,\qquad t\in U.
$$
Пусть $I_s(\T)$ идеал в ${\mathcal E}(\T)$, состоящий из функций, равных нулю в точке $s$. Тогда по лемме Адамара \cite{Petrovsky} всякая функция $u\in{\mathcal E}(\T)$ будет отличаться от функции $u(s)+u'(s)\ph$ на функцию из $I_s^2(\T)$:
\beq\label{Hadamard-dlya-T}
u\overset{I_0^2(\T)}{\equiv} u(0)+u'(0)\ph.
\eeq
Рассмотрим два случая.
\bit{

\item[1)] Зафиксируем точку $s\in 8\setminus\{0\}$. Как элемент пространства 8, $s$ должна лежать либо в  $\T_1$, либо в $\T_2$. Пусть $\ph_s$ обозначает сдвиг функции $\ph$ на элемент $s$ в $\T$:
$$
\ph_s(t)=\ph(t-s),\quad t\in\T,
$$
и пусть $\psi_s$ -- какая-нибудь функция с нулевым ростком в точке 0, и одинаковым с $\ph_s$ ростком в точке $s$:
$$
\psi_s(t)\equiv\begin{cases} 0, & \mod 0\\ \ph_s, & \mod s\end{cases}.
$$
Из \eqref{Hadamard-dlya-T} следует
$$
(u,v)\overset{I_s^2(\T)}{\equiv} \begin{cases}u(s)+u'(s)\cdot(\psi_s,0),& s\in\T_1\\ v(s)+v'(s)\cdot(0,\psi_s),& s\in\T_2 \end{cases},
$$
Как следствие, действие  касательного вектора $\tau\in T_s({\mathcal E}(8))$ на элемент $(u,v)$ описывается формулой
$$
\tau(u,v)=\begin{cases}u(s)+u'(s)\cdot\tau(\psi_s,0),& s\in\T_1\\ v(s)+v'(s)\cdot\tau(0,\psi_s),& s\in\T_2 \end{cases}
$$
Отсюда следует, что (при $s\ne 0$) касательное пространство $T_s({\mathcal E}(8))$ изоморфно $\R$ и $T_s(\T)$.
$$
T_s({\mathcal E}(8))\cong\R\cong T_s({\mathcal E}(\T)).
$$

\item[2)] Рассмотрим точку $s=0$.  Тогда из \eqref{Hadamard-dlya-T} мы получаем
\begin{multline*}
(u,v)\overset{I_0^2(\T)}{\equiv} u(s)\cdot(1,1)+u'(s)\cdot(\ph,\ph)=\\=v(s)\cdot(1,1)+v'(s)\cdot(\ph,\ph),
\end{multline*}
поэтому действие касательного вектора $\tau\in T_0({\mathcal E}(8))$ на этот элемент описывается формулой
$$
\tau(u,v)=u'(s)\cdot\tau(\ph,\ph)=v'(s)\cdot\tau(\ph,\ph).
$$
Отсюда следует, что (в точке $s=0$) касательное пространство $T_0({\mathcal E}(8))$ изоморфно $\R$ и $T_0(\T)$.
$$
T_0({\mathcal E}(8))\cong\R\cong T_0({\mathcal E}(\T)).
$$

}\eit

После того, как доказаны \eqref{Spec(E(8))} и \eqref{T_s(E(8))} зафиксируем какое-нибудь гладкое биективное отображение $\omega:\R\to 8$ (такое отображение всегда есть, и оно не единственно). Оно порождает гомоморфизм алгебр ${\mathcal E}(\R)\gets{\mathcal E}(8)$. Если через $A$ обозначить образ этого гомоморфизма (с индуцированными из ${\mathcal E}(8)$ алгебраическими операциями и топологией), то мы получим как раз алгебру с аннонсированными выше свойствами (i) и (ii), и значит, также и свойством (iii).
\epr

Следующие два примера показывают независимость друг от друга условий (i) и (ii) в теореме \ref{C^infty-obolochka-podalgebry-v-C^infty}.

\bex
Существует замкнутая подалгебра $A$ в алгебре ${\mathcal E}(\R)$ со следующими свойствами:
 \bit{
\item[(i)] спектр $A$ не совпадает с $\R$ как множество,
$$
\Spec(A)\ne \R,
$$

\item[(ii)] для всякой точки $s\in\R$ естественное отображение касательных
пространств
 $$
T_s[A]\gets T_s(\R)
 $$
является изоморфизмом (векторных пространств),
 }\eit
\eex \bpr Такими свойствами обладает алгебра периодических гладких функций на
$\R$ с каким-нибудь фиксированным периодом, например, $1$. \epr

\bex\label{EX:kontrprimer-podalgebra-E(R)}
Существует замкнутая подалгебра $A$ в алгебре ${\mathcal E}(\R)$ со следующими свойствами:
 \bit{
\item[(i)] спектр $A$ совпадает с $\R$,
$$
\Spec(A)=\R,
$$

\item[(ii)] в точке $s=0$ касательное пространство к $A$ вырождено:
$$
T_0[A]=0
$$
 }\eit
\eex
Для доказательства нам понадобится

\blm\label{LM:ploskaya-v-nule-funktsija} Пусть $x$ -- функция из ${\mathcal E}(\R)$, у которой все производные (включая саму функцию $x$ как свою производную нулевого порядка) равны нулю в точке 0,
$$
\forall n\ge 0\qquad x^n(0)=0.
$$
Тогда ее можно приблизить в ${\mathcal E}(\R)$ функциями с нулевым ростком в точке $0$.
\elm
\bpr
Зафиксируем $n\in\N$, $T>1$, $\e>0$ и подберем $\delta>0$ так, чтобы
$$
\forall k\in\{0,...,n\}\qquad \sup_{|t|\le\delta}\abs{x^{(k)}(t)}\le\e.
$$
Выберем функцию $\eta\in{\mathcal E}(\R)$ со свойствами:
$$
0\le\eta(s)\le 1,\qquad \eta(s)=\begin{cases}0,& |s|\le\frac{\delta}{2}\\ 1,& |s|\ge\delta \end{cases}
$$
и положим
\begin{multline*}
y_n(t)=\int_0^t \eta(s)\cdot x^{(n+1)}(s)\d s,\\
 y_{n-1}(t)=\int_0^t y_n(s)\d s,\\ ...,\\ y_0(t)=\int_0^t y_1(s)\d s.
\end{multline*}
Тогда
\begin{multline*}
\sup_{|t|\le T}\abs{x^n(t)-y_n(t)}=\\=
\sup_{|t|\le T}\abs{\int_0^t x^{(n+1)}(s)\d s-\int_0^t \eta(s)\cdot x^{(n+1)}(s)\d s}\le\\
\le\sup_{|t|\le T}\int_0^t (1-\eta(s))\cdot \abs{x^n(s)}\d s\le\\ \le \sup_{|s|\le
T}\abs{x^n(s)}\cdot \int_0^t (1-\eta(s))\d s\le\\ \le \sup_{|s|\le
T}\abs{x^{(n+1)}(s)}\cdot\delta.
\end{multline*}
$$
\Downarrow
$$
\begin{multline*}
\sup_{|t|\le T}\abs{x^{(n-1)}(t)-y_{n-1}(t)}=\\=
\sup_{|t|\le T}\abs{\int_0^t
x^n(s)\d s-\int_0^t y_n(s)\d s}\le\\ \le\sup_{|t|\le T}\int_0^t
\abs{x^n(s)-y_{n}(s)}\d s\le \sup_{|s|\le T}\abs{x^{(n+1)}(s)}\cdot
 \delta\cdot T.
\end{multline*}
$$
\Downarrow
$$
$$
...
$$
$$
\Downarrow
$$
\begin{multline*}
\sup_{|t|\le T}\abs{x^{(0)}(t)-y_0(t)}=\\=
\sup_{|t|\le T}\abs{\int_0^t x^1(s)\d
s-\int_0^t y_1(s)\d s}\le\\ \le\sup_{|t|\le T}\int_0^t \abs{x^1(s)-y_{1}(s)}\d
s\le\\ \le \sup_{|s|\le T}\abs{x^{(n+1)}(s)}\cdot
 \delta\cdot T^n.
\end{multline*}
Теперь видно, что при фиксированных $n\in\N$, $T>1$, $\e>0$ число $\delta>0$ можно подобрать так, чтобы
$$
\delta <\frac{\e}{T^n\cdot \sup_{|s|\le T}\abs{x^{(n+1)}(s)}},
$$
и тогда построенная функция $y_0$ будет отличаться от $x$ меньше, чем на $\e$ равномерно по производным порядка не больше $n$ на отрезке $[-T,T]$:
\begin{multline*}
\max_{0\le k\le n}\sup_{|t|\le T}\abs{x^{(k)}(t)-y^{(k)}(t)}=\\=
\max_{0\le k\le n}\sup_{|t|\le T}\abs{x^{(k)}(t)-y_k(t)}\le\\ \le T^n\cdot\sup_{|s|\le T}\abs{x^{(n+1)}(s)}<\e.
\end{multline*}
\epr

\bpr[Доказательство примера \ref{EX:kontrprimer-podalgebra-E(R)}] Такими свойствами обладает алгебра гладких функций, имеющих в точке
$s=0$ нулевые производные любого положительного порядка:
$$
A=\{x\in {\mathcal E}(M):\quad \forall n\ge 1\quad x^n(0)=0\}.
$$
1. Покажем сначала, что $\Spec(A)=\R$. Пусть $s\in\Spec(A)$, то есть $s$ -- инволютивный, непрерывный, линейный, мультипликативный и сохраняющий единицу функционал на $A$:
\begin{multline*}
s(x^\bullet)=\overline{s(x)},\quad s(\lambda\cdot x+y)=\lambda\cdot s(x)+s(y),\\ s(x\cdot y)=s(x)\cdot s(y),\quad s(1)=1.
\end{multline*}
Рассмотрим его ядро $\Ker s=\{x\in A:\ s(x)=0\}$ и покажем, что найдется точка $t\in\R$ такая, что если ее рассматривать, как функционал на $A$, то его ядро содержит $\Ker(s)$:
\beq\label{Ker s-subseteq-Ker t}
\Ker s\subseteq\Ker t=\{x\in A:\ x(t)=0\}.
\eeq
Предположим, что это не так, то есть
\beq\label{forall-t-in-R-exists-x_t-in-Ker(s)-x_t(t)-ne-0}
\forall t\in\R\quad \exists x_t\in\Ker s\quad x_t(t)\ne 0.
\eeq
Рассмотрим семейство множеств
$$
U_t=\{ r\in\R: \ x_t(r)\ne 0\}.
$$
Они открыты в $\R$, а условие \eqref{forall-t-in-R-exists-x_t-in-Ker(s)-x_t(t)-ne-0} означает, что они образуют покрытие прямой $\R$. Поэтому в это семейство можно вписать гладкое локально конечное разбиение единицы:
$$
\eta_t\in {\mathcal E}(M),\quad \supp\eta_t\subseteq U_t,\quad \eta_t\ge 0,\quad \sum_{t\in\R}\eta_t=1.
$$
Это разбиение можно выбрать так, чтобы оно удовлетворяло следующим дополнительным условиям:
$$
\forall n\ge 1\quad \eta_0^n(0)=0\quad\&\quad \forall t\ne 0\quad
0\notin\supp\eta_t.
$$
Из них автоматически будет следовать, что все функции $\eta_t$ лежат в $A$, а ряд
$\sum_{t\in\R}\eta_t\cdot x_t\cdot x_t^\bullet$
сходится в $A$ (относительно топологии, индуцированной из ${\mathcal E}(\R)$). Обозначим его сумму буквой $y$,
$$
y=\sum_{t\in\R}\eta_t\cdot x_t\cdot x_t^\bullet,
$$
и заметим, что, поскольку все $x_t$ лежат в замкнутом идеале $\Ker s$, элемент $y$ тоже должен лежать в нем:
$$
y\in\Ker s.
$$
С другой стороны, функция $y$ нигде не обращается в нуль. Действительно, для всякой точки $r\in\R$ найдется какая-то функция $\eta_{t_r}$, не равная нулю в этой точке, и поэтому $\eta_{t_r}(r)>0$. При этом, из условия $\supp\eta_{t_r}\subseteq U_{t_r}=\{q\in\R:\ x_{t_r}(q)\ne 0\}$ следует, что  $\abs{x_{t_r}(r)}>0$, и мы получаем
$$
y(r)=\sum_{t\in\R}\eta_t(r)\cdot \abs{x_t(r)}^2\ge \underbrace{\eta_{t_r}(r)}_{\tiny\begin{matrix}
\text{\rotatebox{90}{$<$}}\\ 0 \end{matrix}}\cdot \underbrace{\abs{x_{t_r}(r)}^2}_{\tiny\begin{matrix}
\text{\rotatebox{90}{$<$}}\\ 0 \end{matrix}}>0.
$$
То есть $y$ обратима, как элемент алгебры ${\mathcal E}(\R)$. Вдобавок, по построению, все производные функции $y$ в точке $0$ равны нулю, поэтому у ее обратной функции $\frac{1}{y}$ все производные в точке $0$ тоже равны нулю. Значит,
$$
\frac{1}{y}\in A.
$$
То есть функция $y$ обратима, как элемент алгебры $A$.

Итак, функция $y$ с одной стороны лежит в идеале $\Ker s$ алгебры $A$, а с другой -- обратима в $A$. Значит, $\Ker s=A$, а это противоречит условию $s(1)=1$. То есть наше предположение \eqref{forall-t-in-R-exists-x_t-in-Ker(s)-x_t(t)-ne-0} оказалось неверно. Таким образом, для некоторого $t\in\R$ должно выполняться \eqref{Ker s-subseteq-Ker t}. Мы получаем, что на $A$ заданы два ненулевых линейных функционала, $s$ и $t$, причем, во-первых, ядро $s$ содержится в ядре $t$, и, во-вторых, на элементе $1\in A$ их значения совпадают: $s(1)=t(1)=1$. Такое возможно только если они совпадают всюду на $A$: $s=t$.

2. Теперь докажем условие (ii). Для этого покажем, что у характера $x\mapsto
x(0)$ ядро $I_0[A]=\{x\in A:\ x(0)=0\}$ обладает следующим свойством:
\beq\label{overline(I_0^2)[A]=I_0[A]} 
\overline{I_0^2}[A]=I_0[A]. 
\eeq 
--- тогда
из формулы \cite[(5.82)]{Akbarov-De-Gruyter-I} будет следовать
$$
T_0[A]=\Real\big(I_0[A]/\overline{I_0^2}[A]\big)^\triangledown=0.
$$
Действительно, пусть $x\in I_0[A]$, то есть $x$ -- функция из ${\mathcal E}(\R)$, у которой не только производные, но и значение в точке 0 равно нулю,
$$
\forall n\ge 0\qquad x^n(0)=0.
$$
По лемме \ref{LM:ploskaya-v-nule-funktsija} ее можно приблизить в ${\mathcal E}(\R)$ функциями с нулевым ростком в точке $0$. С другой стороны, всякую функцию $y$ с нулевым ростком в точке $0$ можно представить как элемент идеала $I_0^2[A]$ -- для этого ее достаточно просто домножить на другую функцию с нулевым ростком в $0$, но равную 1 на носителе $y$. Мы получаем, что $x$ приближается в ${\mathcal E}(\R)$ (а значит и в $A$) функциями из $I_0^2[A]$.
\epr

Следующий пример, подсказанный автору М.Бахтольдом, показывает, что в лемме \ref{LM:Ker-rho_t^(T^star)=0} инъективность отображения $\iota_t^{\C T^\star}$ существенна.

\bex
Если у гомоморфизма инволютивных стереотипных алгебр $\iota:A\to B$ в какой-то точке $t\in\Spec(B)$ соответствующий
морфизм кокасательных пространств $\iota_t^{\C T^\star}:{\C}T_{t\circ\iota}^\star[A]\to{\C}T_t^\star[B]$ не инъективен,
$$
\Ker\iota_t^{T^\star}\ne 0,
$$
то формула \eqref{overline(I_s^2)[A]=rho^(-1)(overline(I_t^2)[B])} не обязана выполняться:
$$
\overline{I_{t\circ\iota}^2}[A]\ne \iota^{-1}\Big(\overline{I_t^2}[B]\Big).
$$
 \eex
\bpr
Рассмотрим алгебру ${\mathcal E}(\R)$ гладких функций на прямой $\R$ и в ней подалгебру $A$, (алгебраически) порожденную функциями $x^2$ и $x\cdot e^x$. Таким образом, $A$ можно представлять себе как линейную оболочку в ${\mathcal E}(\R)$ функций вида
$$
x^{2p+q}\cdot e^{qx},\qquad p,q\in\N=\{0,1,2,...\}
$$
(эти функции образуют алгебраический базис в $A$). Алгебра $A$ разделяет точки прямой $\R$
$$
s\ne t\in\R\qquad\Longrightarrow\qquad \exists a\in A\quad a(s)\ne a(t),
$$
и касательные векторы алгебры ${\mathcal E}(\R)$
$$
\forall s\in\R\quad\forall\tau\in T_s[{\mathcal E}(\R)]\quad \exists a\in A\quad \tau(a)\ne 0.
$$
Отсюда следует по теореме Нахбина \ref{TH:Nachbin}, что $A$ плотна в ${\mathcal E}(\R)$. Наделив $A$ сильнейшей локально выпуклой топологией, мы превратим $A$ в стереотипную алгебру (в силу \cite[Example 5.1.4]{Akbarov-De-Gruyter-I}), а естественное вложение в ${\mathcal E}(\R)$ станет эпиморфизмом стереотипных алгебр $\iota:A\to {\mathcal E}(\R)$. Рассмотрим идеал $I_0[A]$ функций в $A$, обращающихся в нуль в точке $0\in\R$. Он будет линейной оболочкой функций
$$
x^{2p+q}\cdot e^{qx},\quad p,q\in\N=\{0,1,2,...\},\quad p+q\ne 0.
$$
Его квадрат $I_0^2[A]$ будет линейной оболочкой функций
$$
x^{2p+q}\cdot e^{qx},\quad p,q\in\N=\{0,1,2,...\},\quad p>0,\quad q>0.
$$
В этом списке функций ровно на 2 меньше, чем в списке для $I_0[A]$ (в нем недостает $x^2$ и $x\cdot e^x$). Поэтому
фактор-пространство $\C T_0^\star[A]=I_0[A]/I_0^2[A]$ имеет комплексную размерность 2. Значит его вещественная часть имеет вещественную рамерность 2:
$$
\dim_{\R}T_0^\star[A]=2.
$$
Поэтому отображение кокасательных пространств $\iota_0^{T^\star}:T_0^\star[A]\to T_0^\star[{\mathcal E}(\R)]$ не может быть инъективно. При этом функция $x^2$, лежащая в $A$ и в $I_0^2[{\mathcal E}(\R)]$ (и значит, лежащая в $\iota^{-1}\Big(\overline{I_t^2}[B]\Big)$), не лежит в $I_0^2[A]$.
\epr

\end{multicols}\noindent\rule[10pt]{160mm}{0.1pt}

\paragraph{Гладкая оболочка алгебры $k(G)$ на компактной группе Ли $G$.}

Вспомним алгебру $k(G)$ непрерывных по норме матричных элементов, определенную на с.\pageref{DEF:k(G)}.

\btm Гладкая оболочка алгебры $k(G)$ матричных элементов на компактной группе Ли $G$ совпадает с алгеброй ${\mathcal E}(G)$ гладких функций на $G$:
\beq
\Env_{\mathcal E}k(G)={\mathcal E}(G).
\eeq
\etm
\bpr
Алгебра $k(G)$ вкладывается в ${\mathcal E}(G)$, причем их спектры совпадают, а по следствию \ref{COR:kasat-prostr-k-Trig(G)}, в каждой точке спектра касательные пространства тоже совпадают. Дальше работает теорема \ref{C^infty-obolochka-podalgebry-v-C^infty}.
\epr

\subsection{Предварительные результаты}

\paragraph{Гладкие оболочки аугментированных стереотипных алгебр.}

\btm\label{LM:A-augste=>Env_E-A-augste}
Пусть $(A,\e)$ --- инволютивная стереотипная алгебра с аугментацией. Тогда
 \bit{

 \item[(i)] гладкая оболочка $\Env_{\mathcal E} \e:\Env_{\mathcal E} A\to\Env_{\mathcal E} \C=\C$ аугментации $\e$ на $A$ является аугментацией на гладкой оболочке $\Env_{\mathcal E} A$ алгебры $A$;

 \item[(ii)] оболочка  $\env_{\mathcal E} A:A\to\Env_{\mathcal E} A$ является морфизмом аугментированых инволютивных стереотипных алгебр.
    \beq\label{env_E-e}
     \xymatrix @R=2.pc @C=4.pc
 {
 A\ar[r]^{\env_{\mathcal E} A}\ar[d]_{\e} &  \Env_{\mathcal E} A \ar@{-->}[d]^{\Env_{\mathcal E} \e} \\
 \C\ar[r]_{\id_{\C}=\env_{\mathcal E} \C} & \C
 }
    \eeq
  }\eit
\etm
\bpr
Это сразу следует из диаграммы \eqref{env_E-e}, которая коммутативна в силу \eqref{DIAGR:funktorialnost-env_E} и \eqref{Env_E(C)=C}.
\epr

\btm\label{TH:Env_E-functor-v-AugSteAlg}
Если $\ph:(A,\e_A)\to(B,\e_B)$ -- морфизм в категории аугментированных инволютивных стереотипных алгебр, то его гладкая оболочка $\Env_{\mathcal E} \ph:(\Env_{\mathcal E} A,\Env_{\mathcal E} \e_A)\to(\Env_{\mathcal E} B,\Env_{\mathcal E} \e_B)$ -- тоже морфизм в категории аугментированных инволютивных стереотипных алгебр в силу коммутативности диаграммы
    \beq\label{Env_E-functor-v-AugSteAlg}
     \xymatrix @R=2.pc @C=4.pc
 {
 A\ar[rr]^{\env_{\mathcal E} A}\ar[dd]_{\ph}\ar[dr]_{\e_A} & &  \Env_{\mathcal E} A \ar@{-->}[dd]^{\Env_{\mathcal E} \ph}\ar@{-->}[dl]^{\ \Env_{\mathcal E} \e_A} \\
 & \C & \\
 B\ar[rr]_{\env_{\mathcal E} B}\ar[ur]^{\e_B} & & \Env_{\mathcal E} B\ar@{-->}[ul]_{\ \Env_{\mathcal E} \e_B}
 }
    \eeq
\etm
\bpr
Здесь коммутативны периметр и все внутренние треугольники, кроме правого (у которого все стороны пунктирны), и его коммутативность мы должны доказать. Но он будет коммутативен из-за того, что $\env_{\mathcal E} A$ --- эпиморфизм.
\epr

\bprop\label{LM:exten_C=>exten_C^Aug-E}
Если $A$ --- инволютивная стереотипная алгебра с аугментацией $\e:A\to\C$, и $\sigma:A\to A'$ --- ее гладкое расширение, то алгебра $A'$ обладает единственной аугментацией $\e':A'\to\C$, для которой морфизм $\sigma$ становится морфизмом аугментированных стереотипных алгебр, и более того, расширением в категории $\AugInvSteAlg$ аугментированных инволютивных стереотипных алгебр в классе $\DEpi$ плотных эпиморфизмов относительно класса $\AugC*$ аугментированных $C^*$-алгебр.
\eprop
\bpr
1. Рассмотрим диаграмму:
    \beq\label{exten_C=>exten_C^Aug-E}
     \xymatrix 
 {
 A\ar[rr]^{\sigma}\ar[rd]_{\e} & & A' \ar@{-->}[dl]^{\e'} \\
 & \C &
 }
    \eeq
Поскольку $\C$ --- $C^*$-алгебра, морфизм $\e'$ существует и однозначно определен. Коммутативность этой диаграммы означает, что $\sigma$ будет морфизмом аугментированных стереотипных алгебр $\sigma:(A,\e)\to(A',\e')$.

2. Покажем, что этот морфизм $\sigma:(A,\e)\to(A',\e')$ является расширением в $\DEpi$ относительно $\AugC*$. Пусть $\ph:(A,\e)\to(B,\delta)$ --- морфизм в аугментированную $C^*$-алгебру. Рассмотрим диаграмму в категории стеретипных алгебр:
$$
     \xymatrix 
 {
 A\ar@/_5ex/[rdd]_{\ph}\ar[rr]^{\sigma}\ar[rd]_{\e} & & A' \ar[dl]^{\e'}\ar@{-->}@/^5ex/[ldd]^{\ph'} \\
 & \C & \\
 & B\ar[u]^{\delta} &
 }
$$
В ней пунктирная стрелка, $\ph'$, существует, единственна и замыкает периметр, поскольку $\ph$ --- дифференциальный гомоморфизм  в $C^*$-алгебру в силу примера \ref{EX:ph-diff-morfizm}, а $\sigma$ --- гладкая оболочка. Одновременно верхний внутренний треугольник коммутативен, поскольку это просто диаграмма \eqref{exten_C=>exten_C^Aug-E}, а левый внутренний треугольник коммутативен, поскольку $\ph:(A,\e)\to(B,\delta)$ --- морфизм аугментированных алгебр. Вдобавок $\sigma$ --- эпиморфизм, поэтому правый внутренний треугольник тоже должен быть коммутативен:
$$
\delta\circ\ph'=\e'.
$$
Это означает, что $\ph'$ --- морфизм в категории $\AugInvSteAlg$, и поскольку он единственен, $\sigma$ --- расширение в $\AugInvSteAlg$ (в классе $\DEpi$ относительно класса $\AugC*$ аугментированных $C^*$-алгебр).
\epr

\bcor
Для всякой стереотипной алгебры $A$ с аугментацией $\e:A\to\C$ существует единственный морфизм стереотипных алгебр с аугментацией $\upsilon:(A,\e)\to\Env_{\AugC*}^{\DEpi}(A,\e)$, замыкающий диаграмму
    \beq\label{DEF:Env_C-e-E}
     \xymatrix 
 {
 & (A,\e)\ar[dl]_{\env_{\mathcal E} A}\ar[dr]^{\env_{\AugC*}^{\DEpi}(A,\e)} & \\
 (\Env_{\mathcal E} A,\Env_{\mathcal E} \e) \ar@{-->}[rr]_{\upsilon}  & & \Env_{\AugC*}^{\DEpi}(A,\e)
 }
    \eeq
\ecor
\bpr
В обозначениях предложения \ref{LM:exten_C=>exten_C^Aug-E}, здесь $\sigma$ --- не просто расширение, а оболочка $\env_{\mathcal E} A$. И в этом случае аугментация $\e'$ на $A'=\Env_{\mathcal E} A$ --- это морфизм $\Env_{\mathcal E} \e$ из \eqref{env_E-e}. Поскольку по предложению \ref{LM:exten_C=>exten_C^Aug-E} $\sigma=\env_{\mathcal E} A$ --- расширение в категории аугментированных стереотипных алгебр, должен существовать и быть единственным морфизм $\upsilon$ в диаграмме \eqref{DEF:Env_C-e-E}.
\epr

\section{Гладкие оболочки групповых алгебр}

\subsection{Алгебра ${\mathcal K}_\infty(G)$}

Для всякой группы Ли $G$ ее групповая алгебра распределений ${\mathcal E}^\star(G)$ является инволютивной алгеброй Хопфа относительно проективного стереотипного тензорного произведения $\circledast$. Поэтому, по теоремам  \ref{TH:C^infty-obolochka-sohranyaet-Hopfov} и \ref{TH:C^infty-obolochka-sohranyaet-inv-Hopfov}, ее гладкая оболочка $\Env_{\mathcal E}({\mathcal E}^\star(G))$ должна быть коалгеброй с согласованными антиподом и инволюцией в категориях  ${\tt E}\text{-}{\tt Alg}$ гладких алгебр и $({\tt Ste},\odot)$ стереотипных пространств. Обозначим символом ${\mathcal K}_\infty(G)$ пространство, стереотипно сопряженное к $\Env_{\mathcal E}{\mathcal E}^\star(G)$:
 \beq\label{DEF:K_infty(G)}
{\mathcal K}_\infty(G):=\Big(\Env_{\mathcal E}{\mathcal E}^\star(G)\Big)^\star.
 \eeq
Это сопряженное пространство к коалгебре в $({\tt Ste},\odot)$ с согласованными антиподом и инволюцией, поэтому в силу \cite[p.657, $4^\circ$]{Akbarov-De-Gruyter-I}, справедлива

\btm Для всякой группы Ли $G$ пространство ${\mathcal K}_\infty(G)$ является алгеброй в категории $({\tt Ste},\circledast)$ (то есть стереотипной алгеброй) с согласованными антиподом и инволюцией.
\etm

По теореме \ref{TH:C^infty-obolochka-sohranyaet-inv-Hopfov}(ii) морфизм
$$
\big(\env_{\mathcal E}{{\mathcal E}^\star(G)}\big)^\star:{\mathcal K}_\infty(G)=\Big(\Env_{\mathcal E}{\mathcal E}^\star(G)\Big)^\star\to {\mathcal E}^\star(G)^\star={\mathcal E}(G),
$$
сопряженный к морфизму оболочки, является инволютивным гомоморфизмом алгебр. Те же соображения, что применялись на с.\pageref{Ker-e_C-C^star(G)^star=0} для алгебры ${\mathcal K}(G)$, показывают, что морфизм $\big(\env_{\mathcal E}{{\mathcal E}^\star(G)}\big)^\star$ имеет  нулевое ядро.

Как следствие, алгебру ${\mathcal K}_\infty(G)$ можно понимать как некую инволютивную подалгебру в ${\mathcal E}(G)$:

\btm\label{TH:K_infty(G)->E(G)} Отображение $u\mapsto u\circ \env_{\mathcal E}{{\mathcal E}^\star(G)}\circ\delta$
совпадает с отображением $\big(\env_{\mathcal E}{{\mathcal E}^\star(G)}\big)^\star$, сопряженным к $\env_{\mathcal E}{{\mathcal E}^\star(G)}$:
\beq\label{K_infty(G)->E(G)}
\env_{\mathcal E}{{\mathcal E}^\star(G)}^\star(u)=u\circ \env_{\mathcal E}{{\mathcal E}^\star(G)}\circ\delta
\eeq
и инъективно вкладывает ${\mathcal K}_\infty(G)$ в ${\mathcal E}(G)$ в качестве инволютивной подалгебры (и поэтому операции сложения, умножения и инволюции в ${\mathcal K}_\infty(G)$ являются поточечными).
\etm

Следующее утверждение аналогично теореме \ref{TH:K(G)=lim}:

\btm\label{TH:K_infty(G)=lim}
Алгебра ${\mathcal K}_\infty(G)$ как стереотипное пространство представима в виде узлового кообраза (в категории стереотипных пространств)
 \beq\label{K_infty(G)=lim}
{\mathcal K}_\infty(G)=\Coim_{\infty}\ph^\star
 \eeq
отображения $\ph^\star$, сопряженного к естественному морфизму стереотипных пространств
$$
\ph:{\mathcal E}^\star(G)\to \projlim_{U}{\mathcal E}^\star(G)/U.
$$
где $U$ пробегает систему дифференциальных окрестностей нуля в ${\mathcal E}^\star(G)$.
\etm

Вспомним диаграмму \eqref{Env_E-C^star(G)-Env_E-E^star(G)} и достроим ее до диаграммы
\beq\label{Env_E-C^star(G)-Env_E-E^star(G)-NEW}
 \xymatrix @R=4pc @C=4pc
 {
{\mathcal C}^\star(G)\ar[r]^{\lambda^\star}\ar@{..>}@/_8ex/[dd]_{\env_{\mathcal C}{\mathcal C}^\star(G)}\ar@{..>}[d]^{\env_{\mathcal E}{\mathcal C}^\star(G)} & {\mathcal E}^\star(G)\ar@{..>}[d]_{\env_{\mathcal E}{\mathcal E}^\star(G)}\ar@{..>}@/^8ex/[dd]^{\env_{\mathcal C}{\mathcal E}^\star(G)}
\\
\Env_{\mathcal E}{\mathcal C}^\star(G)\ar@{..>}[d]^{\zeta_{{\mathcal C}^\star(G)}}\ar@{=>}[r]^{\Env_{\mathcal C}(\lambda^\star)} & \Env_{\mathcal E}{\mathcal E}^\star(G)\ar@{..>}[d]_{\zeta_{{\mathcal E}^\star(G)}}
\\
\Env_{\mathcal C}{\mathcal C}^\star(G)\ar@{=>}[r]^{\Env_{\mathcal C}(\lambda^\star)} & \Env_{\mathcal C}{\mathcal E}^\star(G)
 }
\eeq
Здесь периметр -- диаграмма \eqref{Env_C-C^star(G)-Env_C-E^star(G)}, а треугольники по бокам -- диаграммы \eqref{Env_infty->Env}. При этом верхняя ("сплошная") стрелка $\lambda^\star$ является плотной инъекцией, ниже
две горизонтальные ("двойные") стрелки являются изоморфизмами силу \eqref{Env_E-C^star(G)=Env_E-E^star(G)} и \eqref{Env_C-C^star(G)=Env_C-E^star(G)}, а остальные ("пунктирные") стрелки являются плотными морфизмами в силу определения оболочек.

Двойственная диаграмма будет выглядеть так:
\beq\label{Env_E-C^star(G)-Env_E-E^star(G)-NEW-*}
 \xymatrix @R=3pc @C=3pc
 {
& {\mathcal C}(G) & {\mathcal E}(G)\ar[l]^{\lambda} &
\\
& \Big(\Env_{\mathcal E}{\mathcal C}^\star(G)\Big)^\star\ar@{..>}[u] & \Big(\Env_{\mathcal E}{\mathcal E}^\star(G)\Big)^\star\ar@{=>}[l]\ar@{=}[r]\ar@{..>}[u] & {\mathcal K}_\infty(G)
\\
{\mathcal K}(G)\ar@{=}[r] & \Big(\Env_{\mathcal C}{\mathcal C}^\star(G)\Big)^\star\ar@{..>}[u] & \Big(\Env_{\mathcal C}{\mathcal E}^\star(G)\Big)^\star \ar@{=>}[l]\ar@{..>}[u] &
 }
\eeq
и при этом верхняя стрелка $\lambda$ является плотной инъекцией, две горизонтальные ("двойные") стрелки -- изомофризмы, а остальные ("пунктирные") стрелки -- инъекции. Мы можем заключить поэтому, что справедлива цепочка  инъекций
$$
{\mathcal K}(G)\subseteq {\mathcal K}_\infty(G)\subseteq {\mathcal E}(G)\subseteq {\mathcal C}(G).
$$
Ее можно дополнить цепочкой \eqref{V-k-K-C}, и мы получим следующее утверждение

\btm\label{TH:V-k-K-K_infty-E-C}
Для всякой вещественной группы Ли справедлива цепочка теоретико-множественных включений
 \beq\label{V-k-K-K_infty-E-C}
\Trig(G)\subseteq k(G)\subseteq {\mathcal K}(G)\subseteq {\mathcal K}_\infty(G)\subseteq {\mathcal E}(G)\subseteq {\mathcal C}(G),
 \eeq
причем
\bit{

\item[(i)] всегда
$$
\overline{\Trig(G)}={\mathcal K}(G),\qquad \overline{{\mathcal E}(G)}={\mathcal C}(G)
$$

\item[(ii)] если $G=C\times K$, где $C$ -- абелева компактно порожденная группа Ли, а $K$ -- компактная группа Ли, то
\beq\label{overline-K(G)=K_infty(G)}
\overline{{\mathcal K}(G)}={\mathcal K}_\infty(G),
\eeq

\item[(iii)] если $G$ -- SIN-группа, то
\beq\label{overline-K(G)=E(G)}
\overline{{\mathcal K}(G)}={\mathcal E}(G).
\eeq
}\eit
\etm

Нам понадобится

\blm\label{LM:T_a[K(G)]}
Если $G$ -- SIN-группа, являющаяся одновременно группой Ли, то касательное пространство к алгебре ${\mathcal K}(G)$ в любой точке $a\in G$ совпадает с касательным пространством к $G$ в этой точке:
\beq\label{T_a[K(G)]=T_a(G)}
T_a[{\mathcal K}(G)]=T_a(G).
\eeq
\elm
\bpr
1. Пусть сначала $G$ -- абелева группа. Поскольку по теореме \ref{LM:sdvig-v-K(G)}, ${\mathcal K}(G)$ инвариантна относительно сдвигов, $a$ можно взять равным нулю: $a=0\in G$. Тогда
\begin{multline*}
T_a[{\mathcal K}(G)]=T_0\Big[\big(\Env_{\mathcal C}{\mathcal C}^\star(G)\big)^\star\Big]=\eqref{Env_C-C^star(G)=C(widehat(G))}=
T_0\Big[\big({\mathcal C}(\widehat{G})\big)^\star\Big]=\\=T_0\Big[{\mathcal C}^\star(\widehat{G})\Big]=\Hom(\widehat{G},\R)=\Hom(\R,\widehat{\widehat{G}})=\Hom(\R,G)=T_0(G)=T_a(G).
\end{multline*}

2. Пусть далее $G$ -- компактная группа. Тогда
$$
T_a[{\mathcal K}(G)]=\eqref{Trig=K}=T_a[k(G)]=\eqref{T_a[k(G)]=T_a(G)}=T_a(G).
$$

3. Пусть $G=\R^n\times K$. Для любых $s\in\R^n$ и $t\in K$ мы по уже доказанному получим:
\begin{multline*}
T_{s,t}[{\mathcal K}(\R^n\times K)]=\eqref{K(A-times-K)-cong-K(A)-circledast-K(K)}=T_{s,t}\Big[{\mathcal K}(\R^n)\circledast {\mathcal K}(K)\Big]=\cite[(5.88)]{Akbarov-De-Gruyter-I}=\\=
T_s[{\mathcal K}(\R^n)]\oplus T_t[{\mathcal K}(K)]=
T_s(\R^n)\times T_t(K)=T_{s,t}(\R^n\times K).
\end{multline*}

4. Пусть, наконец, $G$ -- произвольная SIN-группа (и одновременно группа Ли). Представим $G$ как дискретное расширение \eqref{SIN-kak-rasshirenie} некоторой группы $N=\R^n\times K$. Зафиксируем точку $a\in G$, выберем класс смежности $L$ относительно подгруппы $N$, содержащий $a$, и рассмотрим алгебру ${\mathcal K}_L(G)$, определенную формулой \eqref{DEF:K_L(G)}. Тогда, если $e$ -- единица группы $G$, то
\begin{multline*}
T_a[{\mathcal K}(G)]=T_a[{\mathcal K}_L(G)]=(\text{теорема \ref{LM:sdvig-v-K(G)}})=
T_e[{\mathcal K}_N(G)]=(\text{лемма \ref{LM:K(N)=K_N(G)}})=\\=T_e[{\mathcal K}(N)]=(\text{уже доказано})=T_e(N)=T_e(G)=T_a(G)
\end{multline*}
(равенства означают изоморфизмы при очевидных преобразованиях).
\epr

\bpr[Доказательство теоремы \ref{TH:V-k-K-K_infty-E-C}.]
1. Первая формула в (i) уже доказана в теореме \ref{TH:V-k-K-C}, а вторая -- стандартное соотношение между пространствами гладких и непрерывных функций на гладком многообразии.

2. Если $G=C\times K$, то с одной стороны,
$$
\Env_{\mathcal C}{\mathcal C}^\star(G)=\Env_{\mathcal C}{\mathcal C}^\star(C\times K)=
\eqref{env_C^star(R^n-times-K)}={\mathcal C}\Big(\widehat{C},\prod_{\sigma\in\widehat{K}}{\mathcal B}(X_\sigma)\Big)
$$
а с другой,
$$
\Env_{\mathcal E}{\mathcal E}^\star(G)=\Env_{\mathcal E}{\mathcal E}^\star(C\times K)=\eqref{env_E-C^star(R^n-times-K)}={\mathcal E}\Big(\widehat{C},\prod_{\sigma\in\widehat{K}}{\mathcal B}(X_\sigma)\Big)
$$
и поэтому второе пространство инъективно вкладывается в первое. Отсюда следует, что сопряженные пространства плотно отображаются одно на другое:
$$
\overline{{\mathcal K}(G)}={\mathcal K}_\infty(G).
$$

3. Пусть $G$ -- SIN-группа. Представим ее как дискретное расширение \eqref{SIN-kak-rasshirenie} некоторой группы $N=\R^n\times K$. По лемме \ref{LM:K(N)=K_N(G)}, ограничение пространства ${\mathcal K}(G)$ на $N$ изоморфно ${\mathcal K}(N)$, а по лемме \ref{LM:Spec(R^n-times-K)}, спектр алгебры ${\mathcal K}(N)$ совпадает с $N$:
$$
\Spec{\mathcal K}(N)=N.
$$
Мы можем сделать вывод, что алгебра ${\mathcal K}(G)$ разделяет точки $N$. По теореме \ref{LM:sdvig-v-K(G)}, сдвиги являются изоморфизмами ${\mathcal K}(G)$, поэтому ${\mathcal K}(G)$ разделяет точки каждого класса смежности $L\in G/N$. Кроме того, по лемме \ref{LM:1_L-in-K(G)}, характеристическая функция $1_L$ каждого такого класса $L$ принадлежит ${\mathcal K}(G)$, и отсюда следует, что ${\mathcal K}(G)$ разделяет точки не только внутри каждого $L$, но и между любыми двумя классами $L,M\in G/N$. В итоге мы получаем, что ${\mathcal K}(G)$ разделяет точки на всей группе $G$.

С другой стороны, справедливо равенство \eqref{T_a[K(G)]=T_a(G)}. Вместе это означает, что ${\mathcal K}(G)$ удовлетворяет условиям теоремы Нахбина \ref{TH:Nachbin}, поэтому алгебра ${\mathcal K}(G)$ плотна в ${\mathcal E}(G)$.
\epr

\paragraph{Отображение ${\mathcal K}_\infty(G)\circledast{\mathcal K}_\infty(H)\to {\mathcal K}_\infty(G\times H)$.}

Пусть $G$ и $H$ -- две группы Ли. По аналогии с отображением $\omega_{G,H}:{\mathcal K}(G)\circledast{\mathcal K}(H)\to {\mathcal K}(G\times H)$ (определенным выше на с.\pageref{omega_G,H}) вводится отображение
${\mathcal K}_\infty(G)\circledast{\mathcal K}_\infty(H)\to {\mathcal K}_\infty(G\times H)$.
И аналогично теореме \ref{TH:K(A-times-K)-cong-K(A)-circledast-K(K)} доказывается

\btm\label{TH:K_infty(A-times-K)-cong-K_infty(A)-circledast-K_infty(K)}
Если $C$ -- абелева компактно порожденная локально компактная группа, а $K$ -- компактная группа, то отображение ${\mathcal K}_\infty(C)\circledast {\mathcal K}_\infty(K)\to {\mathcal K}_\infty(C\times K)$ является изоморфизмом:
\beq\label{K_infty(A-times-K)-cong-K_infty(A)-circledast-K_infty(K)}
{\mathcal K}_\infty(C\times K)\cong{\mathcal K}_\infty(C)\circledast {\mathcal K}_\infty(K)
\cong{\mathcal K}_\infty(C)\circledast {\mathcal K}(K)
\eeq
\etm
\bpr
Это следует из \eqref{env_E-C^star(R^n-times-K)}:
\begin{multline*}
{\mathcal K}_\infty(C\times K)=\Big(\Env_{\mathcal E}{\mathcal E}^\star(C\times K)\Big)^\star\cong \eqref{env_E-C^star(R^n-times-K)}\cong \Big(\Env_{\mathcal E}{\mathcal E}^\star(C)\odot\Env_{\mathcal E}{\mathcal E}^\star(K) \Big)^\star\cong\\ \cong \Big(\Env_{\mathcal E}{\mathcal E}^\star(C)\Big)^\star\circledast\Big(\Env_{\mathcal E}{\mathcal E}^\star(K)\Big)^\star\cong{\mathcal K}_\infty(C)\circledast {\mathcal K}_\infty(K)
\end{multline*}
\epr

\subsection{Гладкие оболочки групповых алгебр}

\paragraph{Совпадение гладких оболочек алгебр ${\mathcal C}^\star(G)$ и ${\mathcal E}^\star(G)$.}

Выше на с.\pageref{Env_C-C^star(G)-Env_C-E^star(G)} мы уже отмечали, что у вещественной группы Ли $G$ можно рассматривать две групповые алгебры: алгебру ${\mathcal C}^\star(G)$ мер с компактным носителем и алгебру ${\mathcal E}^\star(G)$ распределений с компактным носителем (эти конструкции подробно описаны в \cite{Ak03}). Если обозначить естественное включение ${\mathcal E}(G)\subseteq{\mathcal C}(G)$ символом $\lambda$, то мы получим диаграмму
\beq\label{Env_E-C^star(G)-Env_E-E^star(G)}
 \xymatrix @R=2.5pc @C=4pc
 {
{\mathcal C}^\star(G)\ar[r]^{\lambda^\star}\ar[d]_{\env_{\mathcal E}{\mathcal C}^\star(G)} & {\mathcal E}^\star(G)\ar[d]^{\env_{\mathcal E}{\mathcal E}^\star(G)}
\\
\Env_{\mathcal E}{\mathcal C}^\star(G)\ar[r]^{\Env_{\mathcal C}(\lambda^\star)} & \Env_{\mathcal E}{\mathcal E}^\star(G)
 }
\eeq

По аналогии с теоремой \ref{TH:Env_C-C^star(G)=Env_C-E^star(G)} доказывается

\btm\label{TH:Env_E-C^star(G)=Env_E-E^star(G)}
Для всякой вещественной группы Ли $G$ гладкие оболочки групповых алгебр ${\mathcal C}^\star(G)$ и ${\mathcal E}^\star(G)$ совпадают:
\beq\label{Env_E-C^star(G)=Env_E-E^star(G)}
\Env_{\mathcal E}{\mathcal C}^\star(G)=\Env_{\mathcal E}{\mathcal E}^\star(G).
\eeq
\etm

\paragraph{Преобразование Фурье на коммутативной группе Ли.}
Вспомним преобразование Фурье на коммутативной локально компактной группе $G$, определенное выше формулой \eqref{Fourier-transform}. Если $G$ -- компактно порожденная коммутативная группа Ли, то двойственная группа $\widehat{G}$ также будет (компактно порожденной) группой Ли. Поэтому в этом случае определены алгебры гладких функций ${\mathcal E}(G)$, ${\mathcal E}(\widehat{G})$ и распределений ${\mathcal E}^\star(G)$, ${\mathcal E}^\star(\widehat{G})$. Та же самая формула \eqref{Fourier-transform} будет определять отображение
$$
{\mathcal F}_G:{\mathcal E}^\star(G)\to {\mathcal E}(\widehat{G}).
$$
которое мы также будем называть {\it преобразованием Фурье} на группе $G$.

\blm Для компактно порожденной коммутативной группы Ли $G$ преобразование Фурье непрерывно переводит алгебру мер ${\mathcal C}^\star(G)$ и алгебру распределений ${\mathcal E}^\star(G)$ в алгебру ${\mathcal E}(\widehat{G})$ гладких функций на  $\widehat{G}$.
$$
{\mathcal F}_G:{\mathcal C}^\star(G)\to {\mathcal E}(\widehat{G}),\qquad {\mathcal F}_G:{\mathcal E}^\star(G)\to {\mathcal E}(\widehat{G}).
$$
\elm
\bpr
Пусть $\alpha\in{\mathcal E}^\star(G)$. Всякий вещественный характер (непрерывный гомоморфизм) $r:G\to\R$ определяет однопараметрическую подгруппу  (непрерывный гомоморфизм)
$$
h:\R\to\widehat{G}\quad\Big|\quad h(t)=e^{itr},
$$
и наоборот, всякая однопараметрическая подгруппа $h:\R\to\widehat{G}$ имеет такой вид \cite[(24.43)]{Hewitt-Ross}. Поэтому
\begin{multline*}
\frac{{\mathcal F}_G(\alpha)(\chi\cdot h(t))-{\mathcal F}_G(\alpha)(\chi)}{t}=
\frac{\alpha(\chi\cdot e^{itr})-\alpha(\chi)}{t}=\alpha\left(\frac{\chi\cdot e^{itr}-\chi}{t}\right)=
\alpha\left(\chi\cdot \frac{e^{itr}-1}{t}\right)\underset{t\to 0}{\longrightarrow}
\alpha(\chi\cdot ir)
\end{multline*}
и из этого соотношения видно, что ${\mathcal F}_G(\alpha)$ непрерывно дифференцируема вдоль любой однопараметрической подгруппы в $\widehat{G}$. Точно так же проверяется непрерывная дифференцируемость любого порядка, и мы получаем ${\mathcal F}_G(\alpha)\in{\mathcal E}(\widehat{G})$, причем частные производные будут вычисляться по формулам
$$
\partial_{h_k}...\partial_{h_1}{\mathcal F}_G(\alpha)(\chi)=\alpha(\chi\cdot (ir_1)\cdot...\cdot (ir_k)).
$$
Далее, если $\alpha_\nu\to 0$ в ${\mathcal C}^\star(G)$, то для всякого компакта $K\subseteq\widehat{G}$ и любых $r_1,...,r_k$ множество $\{\chi\cdot (ir_1)\cdot...\cdot (ir_k);\ \chi\in K \}$ образует компакт в ${\mathcal C}(\widehat{G})$, поэтому
$$
\partial_{h_k}...\partial_{h_1}{\mathcal F}_G(\alpha_\nu)(\chi)=\alpha(\chi\cdot (ir_1)\cdot (ir_k))\underset{\nu\to \infty}{\longrightarrow}0
$$
равномерно по $\chi\in K$. Это значит, что $\partial_{h_k}...\partial_{h_1}{\mathcal F}_G(\alpha_\nu)\to 0$ в пространстве ${\mathcal C}(\widehat{G})$. Поскольку это верно для любых $h_1,...,h_k$, мы получаем, что ${\mathcal F}_G(\alpha_\nu)\to 0$ в пространстве ${\mathcal E}(\widehat{G})$. Точно так же рассматривается случай $\alpha\in{\mathcal C}^\star(G)$.
\epr

\btm\label{TH:Fourier-Lie=obolochka} Для компактно порожденной коммутативной группы Ли $G$ оба преобразования Фурье
$$
{\mathcal F}_G:{\mathcal C}^\star(G)\to {\mathcal E}(\widehat{G}),\qquad {\mathcal F}_G:{\mathcal E}^\star(G)\to {\mathcal E}(\widehat{G}).
$$
являются гладкими оболочками групповых алгебр:
\beq\label{Env_E-C^star(G)-commut}
\Env_{\mathcal E} {\mathcal C}^\star(G)=\Env_{\mathcal E} {\mathcal E}^\star(G)={\mathcal E}(\widehat{G}).
\eeq
 \etm
\bpr В силу теоремы \ref{TH:Env_E-C^star(G)=Env_E-E^star(G)}, достаточно рассмотреть случай ${\mathcal E}^\star(G)$. А для него надо просто проверить условия (i) и (ii) теоремы \ref{C^infty-obolochka-podalgebry-v-C^infty}.  Пусть $\delta:G\to {\mathcal E}^\star(G)$ -- вложение группы $G$ в групповую алгебру ${\mathcal E}^\star(G)$
в виде дельта-функций:
 $$
\delta_a(u)=u(a),\qquad u\in {\mathcal E}(G).
 $$

Тогда, во-первых, в силу \cite[Theorem 10.12]{Ak03}, формула
$$
\chi=\ph\circ\delta
$$
задает биекцию между характерами $\ph:{\mathcal E}^\star(G)\to\C$ на алгебре
${\mathcal E}^\star(G)$ и комплексными характерами $\chi:G\to\C^\times$ на
группе $G$. При этом инволютивным характерам $\ph:{\mathcal E}^\star(G)\to\C$
будут соответствовать обычные характеры $\chi:G\to\T$ (со значениями в
окружности $\T$). Это соответствие $\ph\leftrightarrow\chi$ непрерывно в обе
стороны, поэтому задает изоморфизм
$$
\Spec{\mathcal E}^\star(G)\cong \widehat{G}.
$$
Это условие (i) теоремы \ref{C^infty-obolochka-podalgebry-v-C^infty}.

Во-вторых, пусть $\tau:{\mathcal E}^\star(G)\to\C$ -- касательный вектор в
точке спектра $\e\in\Spec\Big({\mathcal E}^\star(G)\Big)$, соответствующей
единичному характеру $1(a)=1$, $a\in G$. Тождество Лейбница
\eqref{DEF:kasatelnyj-vektor} для него принимает вид
$$
\tau(\alpha *
\beta)=\tau(\alpha)\cdot\e(\beta)+\e(\alpha)\cdot\tau(\beta),\qquad
\alpha,\beta\in {\mathcal E}^\star(G),
$$
и если заменить $\alpha$ на $\delta_a$, а $\beta$ на $\delta_b$, то мы получаем
$$
\tau(\delta_{a\cdot b})= \tau(\delta_a *\delta_b)
=\tau(\delta_a)\cdot\e(\delta_b)+\e(\delta_a)\cdot\tau(\delta_b)=
\tau(\delta_a)+\tau(\delta_b),\qquad a,b\in G.
$$
Это значит, что отображение $a\mapsto\tau(\delta_a)$ является гомоморфизмом из
$G$ в аддитивную группу $\C$. Если вдобавок потребовать, чтобы $\tau$ был
инволютивным касательным вектором, то числа $\tau(\delta_a)$ станут
вещественными, поэтому отображение $a\mapsto\tau(\delta_a)$ превратится в
(непрерывный) гомоморфизм групп $G\to\R$, то есть станет вещественным
характером. Это устанавливает биекцию между касательным пространством
$T_{\e}\Big({\mathcal E}^\star(G)\Big)$ к алгебре ${\mathcal E}^\star(G)$ и
группой $\Hom(G,\R)$ вещественных характеров на $G$. Но $\Hom(G,\R)$ изоморфна
группе однопараметрических подгрупп в $\widehat{G}$, то есть группе
$\Hom(\R,\widehat{G})$ (непрерывных) гомоморфизмов $\R\to \widehat{G}$ (см. напр.
\cite[(24.42)]{Hewitt-Ross}). Как следствие, $\Hom(G,\R)$ изоморфна
касательному пространству к группе $\widehat{G}$ в точке $1\in \widehat{G}$, и мы
получаем цепочку изоморфизмов
$$
T_{\e}\Big({\mathcal E}^\star(G)\Big)\cong
\Hom(G,\R)\cong\Hom(\R,\widehat{G})\cong T_1(\widehat{G}).
$$
Если опустить промежуточные равенства, то точку 1 можно заменить на любую
другую точку $\chi\in \widehat{G}$, поэтому мы получаем изоморфизм
$$
T_{\chi}\Big({\mathcal E}^\star(G)\Big)\cong T_{\chi}(\widehat{G})\cong
T_{\chi}\Big({\mathcal E}(\widehat{G})\Big).
$$
Это как раз условие (ii) теоремы \ref{C^infty-obolochka-podalgebry-v-C^infty}.
 \epr

\paragraph{Гладкая оболочка групповой алгебры компактной группы.}

\btm\label{PROP:glad-obol-komp-gruppy}
Для компактной группы $K$ гладкой оболочкой ее групповой алгебры мер ${\mathcal C}^\star(K)$ является декартово произведение алгебр ${\mathcal B}(X_\pi)$, где $\pi$ пробегает двойственный объект $\widehat{K}$, а $X_\pi$ представляет собой пространство представления $\pi$:
\beq\label{glad-obol-komp-gruppy}
\Env_{\mathcal E}{\mathcal C}^\star(K)=\prod_{\pi\in\widehat{K}}{\mathcal B}(X_\pi).
\eeq
Если вдобавок $K$ -- группа Ли, то та же алгебра будет гладкой оболочкой групповой алгебры распределений ${\mathcal E}^\star(K)$:
\beq\label{glad-obol-komp-gruppy-E^star(G)}
\Env_{\mathcal E}{\mathcal E}^\star(K)=\Env_{\mathcal E}{\mathcal C}^\star(K)=\prod_{\pi\in\widehat{K}}{\mathcal B}(X_\pi).
\eeq
\etm

\blm\label{LM:D_k|_K=0}
Если $K$ -- компактная группа, а $B$ -- $C^*$-алгебра, то в любой системе частных производных $D_k:{\mathcal C}^\star(K)\to B$, $k\in\N[m]$, только оператор $D_0$ может быть отличен от нуля:
\beq\label{D_k|_K=0}
\forall k>0\qquad D_k=0.
\eeq
\elm

\bpr Рассмотрим отображения
$$
\ph_k=D_k\circ\delta:K\to B.
$$
Нам нужно убедиться, что
\beq\label{ph_k=0,k-ne-0}
\forall k\ne 0\qquad \ph_k=0
\eeq
Для всякого $a\in K$ мы получаем
$$
\ph_0(a)^{-1}=\ph_0(a^{-1})=D_0(\delta^{a^{-1}})=D_0((\delta^a)^\bullet)=D_0(\delta^a)^\bullet=\ph_0(a)^\bullet,
$$
то есть $\ph_0(a)$ является унитарным элементом в $B$. Как следствие,
\beq\label{||ph_0(a)||=1}
||\ph_0(a)||=1,\qquad a\in K.
\eeq
Пусть теперь $k$ -- мультииндекс порядка 1. Тогда во-первых,
\beq\label{norm(ph_k(x))=C<infty}
\sup_{x\in K}\norm{\ph_k(x)}=C<\infty
\eeq
(потому что $\ph_k:K\to B$ -- непрерывное отображение на компакте $K$).
И, во-вторых,
\beq\label{ph_k(a-cdot-b)=ph_0(a)-cdot-ph_k(b)+ph_k(a)-cdot-ph_0(b)}
\ph_0(a)\cdot\ph_k(b)=\ph_k(b)\cdot\ph_0(a),\qquad
\ph_k(a\cdot b)=\ph_0(a)\cdot\ph_k(b)+\ph_k(a)\cdot\ph_0(b),\qquad a,b\in K
\eeq
$$
\Downarrow
$$
$$
\ph_k(a^p)=p\cdot \ph_0(a)^{p-1}\cdot\ph_k(a),\qquad a\in K,\quad p\in\N.
$$
$$
\Downarrow
$$
$$
\ph_k(a)=\frac{1}{p}\cdot \ph_0(a^{-1})^{p-1}\cdot\ph_k(a^p),\qquad a\in K,\quad p\in\N.
$$
$$
\Downarrow
$$
$$
\norm{\ph_k(a)}=\frac{1}{p}\cdot\norm{\ph_0(a^{-1})^{p-1}\cdot\ph_k(a^p)}\le
\frac{1}{p}\cdot\underbrace{\norm{\ph_0(a^{-1})}^{p-1}}_{\scriptsize \begin{matrix}\eqref{||ph_0(a)||=1}\ \|\ \phantom{\eqref{||ph_0(a)||=1}} \\ 1\end{matrix}}\cdot\kern-10pt\underbrace{\norm{\ph_k(a^p)}}_{\scriptsize \begin{matrix} \eqref{norm(ph_k(x))=C<infty}\ \text{\rotatebox{90}{$\ge$}} \ \phantom{\eqref{norm(ph_k(x))=C<infty}} \\ C\end{matrix}}\kern-10pt
 \le \frac{1}{p}\cdot 1\cdot C\underset{p\to\infty}{\longrightarrow}0
$$
$$
\Downarrow
$$
$$
\ph_k(a)=0,\qquad a\in K.
$$

Мы доказали \eqref{ph_k=0,k-ne-0} для мультииндексов $k$ порядка 1. Если теперь $|k|=2$, то мы получим
\begin{multline*}
\ph_k(a\cdot b)=\sum_{0\le l\le k}\begin{pmatrix}k \\ l\end{pmatrix}\cdot \ph_{k-l}(a)\cdot \ph_l(b)
=\\=\ph_k(a)\cdot \ph_0(b)+\sum_{|l|=1}\begin{pmatrix}k \\ l\end{pmatrix}\cdot \ph_{k-l}(a)\cdot \underbrace{\ph_l(b)}_{\scriptsize\begin{matrix}\| \\ 0\end{matrix}}+
\ph_0(a)\cdot \ph_k(b)=\ph_k(a)\cdot \ph_0(b)+\ph_0(a)\cdot \ph_k(b)
\end{multline*}
То есть $\ph_k$ удовлетворяет правому тождеству в \eqref{ph_k(a-cdot-b)=ph_0(a)-cdot-ph_k(b)+ph_k(a)-cdot-ph_0(b)}
(а значит и обоим тождествам), и по тем же самым причинам $\ph_k=0$. В общем случае нужно провести индукцию по $k$.
\epr

\bpr[Доказательство теоремы \ref{PROP:glad-obol-komp-gruppy}]
Здесь достаточно доказать \eqref{glad-obol-komp-gruppy}, потому что \eqref{glad-obol-komp-gruppy-E^star(G)} будет следовать из теоремы \ref{TH:Env_E-C^star(G)=Env_E-E^star(G)}.
Соотношение \eqref{D_k|_K=0} означает, что при взятии гладкой оболочки у алгебры ${\mathcal C}^\star(G)$ класс тестовых морфизмов совпадает с классом морфизмов в $C^*$-алгебры. Это означает, что гладкая оболочка алгебры ${\mathcal C}^\star(G)$ совпадает с ее непрерывной оболочкой, и значит, по предложению \ref{PROP:nepr-obol-komp-gruppy},
$$
\Env_{\mathcal E}{\mathcal C}^\star(K)=\prod_{\pi\in\widehat{K}}{\mathcal B}(X_\pi).
$$
\epr

\paragraph{Гладкая оболочка групповой алгебры группы $Z\times K$.}

\btm\label{TH:env_E-C(R^n-times-K)}
Пусть $Z$ -- абелева компактно порожденная группа Ли, $K$ -- компактная группа. Тогда гладкой оболочкой групповой алгебры мер ${\mathcal C}^\star(Z\times K)$ является алгебра ${\mathcal E}(\widehat{Z},\prod_{\sigma\in\widehat{K}}{\mathcal B}(X_\sigma))$ гладких отображений из двойственной по Понтрягину группы $\widehat{Z}$ в декартово произведение алгебр ${\mathcal B}(X_\sigma)$, где $\sigma$ пробегает двойственный объект $\widehat{K}$, а $X_\sigma$  представляет собой пространство представления $\sigma$:
\beq\label{env_E-C(R^n-times-K)}
\Env_{\mathcal E}{\mathcal C}^\star(Z\times K)={\mathcal E}\Big(\widehat{Z},\prod_{\sigma\in\widehat{K}}{\mathcal B}(X_\sigma)\Big)=
\prod_{\sigma\in\widehat{K}}{\mathcal E}\big(\widehat{Z},{\mathcal B}(X_\sigma)\big)={\mathcal E}(\widehat{Z})\odot \prod_{\sigma\in\widehat{K}}{\mathcal B}(X_\sigma)=
\Env_{\mathcal E}{\mathcal C}^\star(Z)\odot \Env_{\mathcal E}{\mathcal C}^\star(K).
\eeq
Если вдобавок $K$ -- группа Ли, то та же алгебра будет гладкой оболочкой групповой алгебры распределений ${\mathcal E}^\star(Z\times K)$:
\begin{multline}\label{env_E-C^star(R^n-times-K)}
\Env_{\mathcal E}{\mathcal E}^\star(Z\times K)=\Env_{\mathcal C}{\mathcal E}^\star(Z\times K)={\mathcal E}\Big(\widehat{Z},\prod_{\sigma\in\widehat{K}}{\mathcal B}(X_\sigma)\Big)=\\=
\prod_{\sigma\in\widehat{K}}{\mathcal E}\big(\widehat{Z},{\mathcal B}(X_\sigma)\big)={\mathcal E}(\widehat{Z})\odot \prod_{\sigma\in\widehat{K}}{\mathcal B}(X_\sigma)=
\Env_{\mathcal E}{\mathcal E}^\star(Z)\odot \Env_{\mathcal E}{\mathcal E}^\star(K).
\end{multline}
\etm
\bpr Оба утверждения следуют из теорем \ref{TH:Fourier-Lie=obolochka} и \ref{PROP:glad-obol-komp-gruppy}. Например, для ${\mathcal E}^\star(Z\times K)$ цепочка рассуждений выглядит так:
\begin{multline*}
\Env_{\mathcal E}{\mathcal E}^\star(Z\times K)=\Env_{\mathcal E}\Big({\mathcal E}^\star(Z)\circledast{\mathcal E}^\star(K)\Big)=
\cite[(1.129)]{Ak16}=
\Env_{\mathcal E}\Big(\Env_{\mathcal E}{\mathcal E}^\star(Z)\circledast\Env_{\mathcal E}{\mathcal E}^\star(K)\Big)=\\=\eqref{Env_E-C^star(G)-commut},\eqref{glad-obol-komp-gruppy}=
\Env_{\mathcal E}\Big({\mathcal E}(\widehat{Z})\circledast\prod_{\sigma\in\widehat{K}}{\mathcal B}(X_\sigma)\Big)=\eqref{E/circledast}=
{\mathcal E}(\widehat{Z})\overset{\mathcal E}{\circledast}\prod_{\sigma\in\widehat{K}}{\mathcal B}(X_\sigma)=\eqref{E(M)-circledast-A=E(M,A)}=
{\mathcal E}\Big(\widehat{Z},\prod_{\sigma\in\widehat{K}}{\mathcal B}(X_\sigma)\Big)
\end{multline*}
Это доказывает второе равенство в \eqref{env_E-C^star(R^n-times-K)}. Точно так же доказывается первое равенство в \eqref{env_E-C(R^n-times-K)}, и вместе эти равенства дают первое равенство в \eqref{env_E-C^star(R^n-times-K)}:
$$
\Env_{\mathcal E}{\mathcal C}^\star(Z\times K)={\mathcal E}\Big(\widehat{Z},\prod_{\sigma\in\widehat{K}}{\mathcal B}(X_\sigma)\Big)=\Env_{\mathcal E}{\mathcal E}^\star(Z\times K).
$$
Третье равенство в \eqref{env_E-C^star(R^n-times-K)} очевидно. Четвертое следует из \cite[Theorem 8.9]{Ak03}. Наконец, последнее равенство в \eqref{env_E-C^star(R^n-times-K)} следует из \eqref{Env_E-C^star(G)-commut} и \eqref{glad-obol-komp-gruppy}.
\epr

\subsection{Гладкая оболочка алгебры ${\mathcal K}_\infty(Z\times K)$.}

\btm\label{TH:E_E(K_infty(G))=E(G)}
Пусть $Z$ -- абелева группа Ли, $K$ -- компактная группа Ли, и $G=Z\times K$. Тогда
\bit{

\item[(i)] спектр алгебры ${\mathcal K}_\infty(G)$ топологически изоморфен $G$:
\beq\label{Spec-K_infty(G)=G}
\Spec{\mathcal K}_\infty(G)=G
\eeq

\item[(ii)] в любой точке $a\in G$ касательные пространства у ${\mathcal K}_\infty(G)$ и $G$ изоморфны:
\beq\label{T_a[K_infty(G)]=T_a(G)}
T_a[{\mathcal K}_\infty(G)]=T_a(G)
\eeq

\item[(iii)] гладкая оболочка алгебры ${\mathcal K}_\infty(G)$ совпадает с ${\mathcal E}(G)$:
\beq\label{Env_E-K_infty(G)=E(G)}
\Env_{\mathcal E}{\mathcal K}_\infty(G)={\mathcal E}(G).
\eeq

}\eit
\etm
\bpr
1. Из цепочки \eqref{V-k-K-K_infty-E-C} извлечем фрагмент
\beq\label{K(G)->K_infty(G)->E(G)}
{\mathcal K}(G)\subseteq {\mathcal K}_\infty(G)\subseteq {\mathcal E}(G).
\eeq
Из \eqref{overline-K(G)=K_infty(G)} и \eqref{overline-K(G)=E(G)} следует, что это цепочка (непрерывных и) плотных инъекций. Поэтому переходя к спектрам мы получим цепочку (непрерывных) инъекций
$$
G=\Spec{\mathcal K}(G)\gets\Spec{\mathcal K}_\infty(G)\gets\Spec{\mathcal E}(G)=G
$$
(здесь первое равенство доказывается в теореме \ref{TH:Spec-K(G)=G}, а последнее очевидно). Понятно, что это должны быть гомеоморфизмы.

2. Из того, что инъекции в \eqref{K(G)->K_infty(G)->E(G)} плотны, следует также, что цепочка отображений касательных пространств состоит из инъекций
$$
T_a(G)=\eqref{T_a[K(G)]=T_a(G)}=T_a[{\mathcal K}(G)]\gets T_a[{\mathcal K}_\infty(G)]\gets T_a[{\mathcal E}(G)]=T_a(G),
$$
Это доказывает \eqref{T_a[K_infty(G)]=T_a(G)}.

3. Мы доказали (i) и (ii), и по теореме \ref{C^infty-obolochka-podalgebry-v-C^infty} из этого следует \eqref{Env_E-K_infty(G)=E(G)}.
\epr

\subsection{Структура алгебр Хопфа на $\Env_{\mathcal E} {\mathcal E}^\star(G)$ и ${\mathcal K}_\infty(G)$ в случае $G=Z\times K$.}

\btm\label{TH:Env_E(E*(G))-Hopf-dlya-R^n-times-K} Пусть $Z$ -- абелева компактно порожденная группа Ли, $K$ -- компактная группа Ли, и $G=Z\times K$. Тогда
\bit{
\item[(i)] гладкая оболочка $\Env_{\mathcal E}{\mathcal E}^\star(G)$ групповой алгебры ${\mathcal E}^\star(G)$ является инволютивной алгеброй Хопфа в категории стереотипных пространств $(\tt{Ste},\odot)$.

\item[(ii)] двойственная ей алгебра ${\mathcal K}_\infty(G)$ является инволютивной алгеброй Хопфа в категории  стереотипных пространств $(\tt{Ste},\circledast)$
}\eit
\etm
\bpr Здесь достаточно доказать (i). Прежде всего,
$$
\Env_{\mathcal E}{\mathcal E}^\star(Z)=\eqref{Env_E-C^star(G)-commut}={\mathcal E}(\widehat{Z}),
$$
и это будет алгебра Хопфа в категории $(\tt{Ste},\odot)$, например, в силу \cite[Example 10.25]{Ak03}. С другой стороны,
$$
\Env_{\mathcal E}{\mathcal E}^\star(K)=\eqref{glad-obol-komp-gruppy-E^star(G)}=\prod_{\pi\in\widehat{K}}{\mathcal B}(X_\pi)=\eqref{nepr-obol-komp-gruppy}=\Env_{\mathcal C}{\mathcal C}^\star(K),
$$
и это будет алгебра Хопфа в категории $(\tt{Ste},\odot)$ по теореме \ref{TH:Env_C-C^*(G)-odot-Hopf}. Отсюда
пространство
\beq\label{Env_E(E*(G))-Hopf-dlya-R^n-times-K}
\Env_{\mathcal E}{\mathcal E}^\star(Z\times K)=\eqref{env_E-C^star(R^n-times-K)}=
\Env_{\mathcal E}{\mathcal E}^\star(Z)\odot\Env_{\mathcal E}{\mathcal E}^\star(K).
\eeq
должно быть алгеброй Хопфа в $(\tt{Ste},\odot)$ как тензорное произведение алгебр Хопфа.
\epr

\section{Гладкая двойственность}

\paragraph{Гладко рефлексивные алгебры Хопфа.}
Пусть $H$ -- инволютивная стереотипная алгебра Хопфа относительно тензорного произведения $\circledast$. Мы говорим, что $H$ {\it гладко рефлексивна}, если она рефлексивна относительно гладкой оболочки $\Env_{\mathcal E}$ (в смысле определения на с.\pageref{DEF:reflexiv-otn-obolochki}).

Из теорем \ref{TH:E_E(K_infty(G))=E(G)} и \ref{TH:Env_E(E*(G))-Hopf-dlya-R^n-times-K} следует

 \btm\label{TH:glad-dvoistvennost}
Пусть $Z$ -- абелева компактно порожденная группа Ли, $K$ -- компактная группа Ли, и $G=Z\times K$. Тогда
алгебры ${\mathcal E}^\star(G)$ и ${\mathcal K}_\infty(G)$ гладко рефлексивны, а диаграмма рефлексивности для них принимает вид:
 \beq\label{chetyrehugolnik-E-E*}
 \xymatrix @R=1.pc @C=2.pc
 {
 {\mathcal E}^\star(G)
 & \ar@{|->}[r]^{\Env_{\mathcal E}} & &
 \Env_{\mathcal E}{\mathcal E}^\star(G)
 \\
 & & &
 \ar@{|->}[d]^{\star}
 \\
 \ar@{|->}[u]^{\star}
 & & &
 \\
 {\mathcal E}(G)
 & &
 \ar@{|->}[l]_{\Env_{\mathcal E}}
 &
 {\mathcal K}_\infty(G)
 }
 \eeq
 \etm

\bex Из теоремы \ref{TH:Fourier-Lie=obolochka} следует, что для компактно порожденных абелевых групп Ли $Z$ диаграмма рефлексивности принимает вид
 \beq\label{chetyrehugolnik-E-E*-F}
  \xymatrix @R=1.pc @C=2.pc
 {
 {\mathcal E}^\star(Z)
 & \ar@{|->}[r]^{{\mathcal F}_Z} & &
 {\mathcal E}(\widehat{Z})
 \\
 & & &
 \ar@{|->}[d]^{\star}
 \\
 \ar@{|->}[u]^{\star}
 & & &
 \\
 {\mathcal E}(Z)
 & &
 \ar@{|->}[l]_{{\mathcal F}_{\widehat{Z}}}
 &
 {\mathcal E}^\star(\widehat{Z})
  }
\eeq
(здесь $\widehat{Z}$ -- двойственная по Понтрягину группа к $Z$, а ${\mathcal F}_Z$ -- преобразование Фурье, определенное формулой \eqref{Fourier-transform}).
\eex

\bex Из теоремы \ref{PROP:glad-obol-komp-gruppy} следует, что для компактных групп Ли $K$ диаграмма рефлексивности становится такой:
$$
 \xymatrix @R=1.pc @C=2.pc
 {
 {\mathcal E}^\star(K)
 & \ar@{|->}[r]^{\Env_{\mathcal E}} & &
\prod_{\pi\in\widehat{K}}{\mathcal B}(X_\pi)
 \\
 & & &
 \ar@{|->}[d]^{\star}
 \\
 \ar@{|->}[u]^{\star}
 & & &
 \\
 {\mathcal E}(K)
 & &
 \ar@{|->}[l]_{\Env_{\mathcal E}}
 &
 \Trig(K)
 }
$$
\eex

\paragraph{Группы, различаемые $C^*$-алгебрами с присоединенными самосопряженными нильпотентами.} По аналогии с определением на с.\pageref{DEF:gruppy-razlich-C^*-algebrami}, будем говорить, что локально компактная группа $G$ {\it различается  $C^*$-алгебрами с присоединенными самосопряженными нильпотентами}\label{DEF:gruppy-razlich-C^*-algebrami-s-nilp},
если (непрерывные инволютивные) гомоморфизмы ее алгебры мер ${\mathcal C}^\star(G)\to B$  во всевозможные алгебры вида $B[m]$, где $B$ -- $C^*$-алгебра, а $m\in\N[n]$, различают элементы $G$ (при вложении $G$ в ${\mathcal C}^\star(G)$ дельта-функциями).

\btm\label{TH:gruppy-razlich-C*-s-pris-nilp}
Если вещественная группа Ли $G$ различается $C^*$-алгебрами с присоединенными самосопряженными нильпотентами, то $G$ различается $C^*$-алгебрами \etm

Вместе с теоремой Люмине---Валетта \ref{TH:Luminet-Valette} это дает

\bcor
Если связная вещественная группа Ли $G$ различается $C^*$-алгебрами с присоединенными самосопряженными нильпотентами, то $G$  является линейной группой.
\ecor

Нам понадобится

\blm\label{LM:x->e^x}
Для всякой $C^*$-алгебры экспоненциальное отображение
$$
x\mapsto e^x=\sum_{n=0}^\infty\frac{x^n}{n!}\quad : A\to A
$$
инъективно на самосопряженных элементах:
$$
x\ne y\in\Real A\quad\Longrightarrow\quad e^x\ne e^y.
$$
\elm
\bpr
Пусть $\Real_+ A$ обозначает множество положительных самосопряженных элементов в $A$.
Для всякого $z\in\Real_+ A$ и любого $n\in\N$ однозначно определен корень $\sqrt[n]{z}\in\Real_+ A$ (это следует из спектральной теоремы \cite{Kad-Ring}). С другой стороны, для всякого $x\in\Real A$ его экспонента $e^x$ лежит в $\Real_+A$, потому что
$$
e^x=e^{\frac{x}{2}}\cdot e^{\frac{x}{2}}=e^{\frac{x}{2}}\cdot\left( e^{\frac{x}{2}}\right)^\bullet\ge 0.
$$
Значит, однозначно определены корни $\sqrt[n]{e^x}=e^{\frac{x}{n}}$, и мы получаем
$$
x=\lim_{n\to\infty}n\left(e^{\frac{x}{n}}-1\right)=\lim_{n\to\infty}n\left(\sqrt[n]{e^x}-1\right),
$$
то есть $x$ однозначно восстанавливается по $e^x$.
\epr

\bpr[Доказательство теоремы \ref{TH:gruppy-razlich-C*-s-pris-nilp}.]
Предположим, что $G$ не различается $C^*$-алгебрами, то есть у всех гомоморфизмов $\ph:G\to B$ в различные $C^*$ алгебры $B$ имеется нетривиальное ядро $N$, $\{e\}\ne N\subseteq G$. Рассмотрим какой-нибудь гомоморфизм $D:G\to B[m]$. Для всякого мультииндекса $k\in\N[m]$ первого порядка, $|k|=1$, и любых $x,y\in N$ мы получим:
$$
D_k(x\cdot y)=D_k(x)\cdot\underbrace{D_0(y)}_{\scriptsize\begin{matrix}\|\\ 1\end{matrix}}+\underbrace{D_0(x)}_{\scriptsize\begin{matrix}\|\\ 1\end{matrix}}\cdot D_k(y)=D_k(x)+D_k(y).
$$
То есть $D_k:N\to B$ -- логарифм (отображение, переводящее умножение в $N$ в сложение в $B$). Кроме того,
$$
D_k(x)^\bullet=D_k(x^\bullet)=D_k(x^{-1})=-D_k(x),
$$
откуда следует, что всякий элемент $iD_k(x)$, $x\in N$, самосопряжен:
$$
(iD_k(x))^\bullet=iD_k(x).
$$
Рассмотрим отображение
$$
\ph(x)=e^{iD_k(x)}=\sum_{l=0}^\infty\frac{i^{|l|}}{l!}\cdot D_k(x)^l.
$$
Оно является инволютивным гомоморфизмом, то есть представлением группы $N$ в $C^*$-алгебре $B$. Можно рассмотреть его индуцированное представление $\psi:G\to{\mathcal B}(X)$ в каком-то гильбертовом пространстве $X$. Это будет гомоморфизм $G$ в $C^*$-алгебру ${\mathcal B}(X)$, поэтому на подгруппе $N$ отображение $\psi$ должно быть тривиально:
$$
\psi(x)=1,\qquad x\in N.
$$
Из этого следует, что исходный гомоморфизм $\ph:N\to B$ тоже должен быть тривиальным:
$$
\ph(x)=e^{iD_k(x)}=1_B.
$$
Поскольку всякий элемент $iD_k(x)$ самосопряжен, по лемме \ref{LM:x->e^x} это означает, что
$$
D_k(x)=0,\qquad x\in N.
$$
Мы получили, что на подгруппе $N$ все частные производные порядка 1 нулевые. Если теперь $k$ -- мультииндекс порядка 2, то для $x,y\in N$ мы получим
\begin{multline*}
D_k(x\cdot y)=\sum_{0\le l\le k}\begin{pmatrix}k \\ l\end{pmatrix}\cdot D_{k-l}(x)\cdot D_l(y)=\\=
\overbrace{D_k(x)\cdot \underbrace{D_0(y)}_{\scriptsize\begin{matrix}\|\\ 1\end{matrix}}}^{\scriptsize |l|=0}+\sum_{|l|=1}\begin{pmatrix}k \\ l\end{pmatrix}\cdot D_{k-l}(x)\cdot \underbrace{D_l(y)}_{\scriptsize\begin{matrix}\|\\ 0\end{matrix}}+\overbrace{\underbrace{D_0(x)}_{\scriptsize\begin{matrix}\|\\ 1\end{matrix}}\cdot D_k(y)}^{\scriptsize |l|=2}=D_k(x)+D_k(y)
\end{multline*}
То есть $D_k$ в этом случае тоже является логарифмом. По тем же самым причинам, что и раньше, $D_k=0$ на $N$. И так далее.
\epr

\chapter{ДВОЙСТВЕННОСТЬ В КОМПЛЕКСНОЙ ГЕОМЕТРИИ}

\section{Многообразия Штейна: прямоугольники в $\mathcal O(M)$ и ромбы в
$\mathcal{O}^\star(M)$}\label{rectangles-buses}

В этом параграфе мы обсудим некоторые свойства пространств голоморфных функций
на комплексных многообразиях. Для наглядности нам будет полезно условие, чтобы
глобальные функции разделяли точки многообразия, поэтому мы формулируем
результаты для многообразий Штейна. Мы используем терминологию учебников
\cite{Shabat,Taylor,Grauert-Remmert}.

\subsection{Многообразия Штейна}

Пусть $M$ -- комплексное многообразие. Символом $\mathcal O(M)$ мы обозначаем
алгебру всех голоморфных функций на $M$ (с обычными поточечными операциями и
топологией равномерной сходимости на компактах в $M$). Хорошо известно
\cite{Taylor}, что, как топологическое векторное пространство, $\mathcal O(M)$
является пространством Монтеля.

Многообразие $M$ называется {\it многообразием Штейна}\index{многообразие
Штейна} \cite{Shabat}, если оно удовлетворяет следующим трем условиям:
 \bit{
\item[1)] {\it голоморфная отделимость:}\label{holom-separability} для любых
двух точек $x\ne y\in M$ найдется функция $u\in\mathcal O(M)$ такая, что
$$
u(x)\ne u(y)
$$
\item[2)] {\it голоморфная униформизация:} для любой точки $x\in M$ найдутся
функции $u_1,...,u_n\in\mathcal O(M)$, образующие локальные координаты
многообразия $M$ в окрестности $x$;

\item[3)] {\it голоморфная выпуклость:} для любого компакта $K\subseteq M$ его
{\it голоморфно выпуклая оболочка}, то есть множество
$$
\widehat{K}=\{x\in M: \ \forall u\in\mathcal O(M)\ |u(x)|\le\max_{y\in K}|u(y)| \}
$$
является компактом в $M$.
 }\eit

\noindent\rule{160mm}{0.1pt}\begin{multicols}{2}

\bex
Комплексное пространство $\C^n$ является многообразием Штейна.
\eex

\bex
Многообразиями Штейна будут всевозможные {\it области голоморфности} в $\C^n$, то есть открытые множества $\varOmega\subseteq \C^n$, удовлетворяющие следующим эквивалентным условиям:
\biter{

\item[(i)] голоморфная выпуклость:   для любого компакта $K\subseteq \varOmega$ его
{\it голоморфно выпуклая оболочка}, то есть множество
$$
\widehat{K}=\{x\in \varOmega: \ \forall u\in{\mathcal O}(\varOmega)\ |u(x)|\le\max_{y\in K}|u(y)| \}
$$
является компактом в $\varOmega$,

\item[(ii)] в $\varOmega$ не существует открытого подмножества $U\subseteq\varOmega$ и связного открытого множества $V\subseteq\C^n$ со свойствами 
    $$
    U\subseteq \varOmega\cap V,\qquad V\nsubseteq \varOmega,
    $$
    и таких, что всякая функция $u\in {\mathcal O}(\varOmega)$ голомофрно продолжается с $U$ на $V$:
    $$
    \exists v\in{\mathcal O}(V)\quad u\big|_U=v\big|_U
    $$

\item[(iii)] всякая точка на границе $x\in\partial\varOmega$ обладает окрестностью $U$ такой что какая-нибудь функция $u\in {\mathcal O}(U\cap\varOmega)$ не продолжается до голоморфной функции на какую-нибудь окрестность $x$.     

}\eiter
\eex

\bex
 Простейший пример комплексного многообразия, не
являющегося многообразием Штейна -- {\it комплексный тор},\index{комплексный
тор} то есть факторгруппа комплексной плоскости $\C$ по решетке $\Z+i\Z$. Структурная функциональная алгебра такого многообразия состоит из постоянных функций:
\beq
{\mathcal O}\l\C/(\Z+i\Z)\r=\C.
\eeq
\eex

\bex
Cумма $\coprod_{i\in I}M_i$ (как топологических пространств) многообразий Штейна $M_i$ также будет многообразием Штейна. Структурная функциональная алгебра такого многообразия будет прямым проихзведением алгебр для $M_i$:
\beq
{\mathcal O}\l\coprod_{i\in I}M_i\r=\prod_{i\in I}{\mathcal O}(M_i)
\eeq
\eex

\bex\label{EX:discr-prostr-Stein}
Любое дискретное топологическое пространство $M$ можно считать примером многообразия Штейна (нулевой размерности) со структурной функциональной алгеброй $\C^M$:
\beq\label{O(M)-dlya-discr-prostr}
{\mathcal O}(M)=\C^M
\eeq
\eex

\end{multicols}\noindent\rule[10pt]{160mm}{0.1pt}

\subsection{Внешние огибающие на $M$ и прямоугольники в $\mathcal O(M)$}
\label{SUBSEC:vnesh-ogib}

\paragraph{Операции $\protect\BSQ$ и $\protect\SQ$.}

В этом пункте нас будут интересовать вещественные функции $f$ на многообразии
$M$, ограниченные снизу единицей,
$$
f\ge 1,
$$
то есть принимающие значения в множестве $[1,+\infty)$. Естественно, мы будем
использовать запись $f:M\to[1;+\infty)$. Функцию $f:M\to[1;+\infty)$ мы, как
обычно, называем локально ограниченной, если для всякой точки $x\in M$ можно
подобрать окрестность $U\owns x$ такую, что
$$
\sup_{y\in U} |f(y)|<\infty
$$
Поскольку $f$ ограничена снизу единицей, это условие равносильно условию
$$
\sup_{y\in U} f(y)<\infty
$$

 \bprop\label{PROP:f^blacksquare}
Для каждой локально ограниченной функции $f:M\to[1;+\infty)$ формула
 \beq\label{f^blacksquare}
f^\text{\BSQ} :=\{u\in \mathcal O(M): \quad \forall x\in M\quad |u(x)|\le
f(x)\}
 \eeq
определяет абсолютно выпуклое компактное множество функций
$f^\text{\BSQ}\subseteq\mathcal O(M)$, содержащее тождественную единицу:
$$
1\in f^\text{\BSQ}.
$$
\eprop \bpr Множество $f^\text{\BSQ}$ будет компактом, потому что оно замкнуто
и ограничено в пространстве Монтеля $\mathcal O(M)$. \epr

 \bprop\label{PROP-nepr-vnesh-ogib}
Для каждого ограниченного множества функций $D\subseteq \mathcal O(M)$,
содержащего тождественную единицу,
$$
1\in D,
$$
формула
 \beq\label{fi^square}
 D^\text{\SQ} (x):=\sup_{u\in D} |u(x)|, \qquad x\in M
 \eeq
определяет непрерывную вещественную функцию $D^\text{\SQ}:M\to\R$, ограниченную
снизу единицей:
$$
 D^\text{\SQ}\ge 1.
$$
\eprop \bpr Заметим с самого начала, что $D$ можно считать компактом. Для этого
рассмотрим замыкание $\overline{D}$ множества $D$. Поскольку $D$ ограничено в
пространстве Монтеля $\mathcal O(M)$, $\overline{D}$ будет компактом в
$\mathcal O(M)$. При этом, поскольку при непрерывном отображении $u\mapsto
\delta^x(u)=u(x)$ образ замыкания $\delta^x(\overline{D})$ содержится в
замыкании образа $\overline{\delta^x( D)}$
$$
\delta^x(\overline{D})\subseteq \overline{\delta^x( D)},
$$
выполняется цепочка неравенств
$$
\sup_{\lambda\in \delta^x( D)} |\lambda|\le \sup_{\lambda\in
\delta^x(\overline{D})} |\lambda|\le \sup_{\lambda\in \overline{\delta^x( D)}}
|\lambda|=\sup_{\lambda\in \delta^x( D)} |\lambda|
$$
откуда
$$
\sup_{\lambda\in \delta^x(\overline{D})} |\lambda|=\sup_{\lambda\in \delta^x(
D)} |\lambda|
$$
то есть функции $\overline{D}^\text{\SQ}$ и $D^\text{\SQ}$ совпадают:
$$
\overline{D}^\text{\SQ} (x)=\sup_{u\in \overline{D}}
|\delta^x(u)|=\sup_{\lambda\in \delta^x(\overline{D})}
|\lambda|=\sup_{\lambda\in \delta^x( D)} |\lambda|=\sup_{u\in D} |\delta^x(u)|=
D^\text{\SQ} (x)
$$

Итак, $D$ можно считать компактом. Зафиксируем произвольный компакт $K\subseteq
M$ и рассмотрим пространство $C(K)$ непрерывных функций на $K$ (с обычной
топологией равномерной сходимости на $K$). Отображение ограничения $u\in
\mathcal O(M)\mapsto u|_K\in C(K)$ будет непрерывным отображением пространства
Фреше $\mathcal O(M)$ в банахово пространство $C(K)$. Если $D$ -- компакт в
$\mathcal O(M)$ , то его образ $D|_K$ -- компакт в $C(K)$. Отсюда, по теореме
Арцела, $D|_K$ поточечно ограничено и равностепенно непрерывно на $K$.
Следовательно, функция
$$
 D^\text{\SQ} (x):=\sup_{u\in D} |u(x)|, \qquad x\in K
$$
непрерывна на $K$. Поскольку это верно для любого компакта $K$ в $M$, мы
получаем, что эта функция непрерывна на $M$. \epr

\medskip

\centerline{\bf Свойства операций $\text{\BSQ}$ и $\text{\SQ}$:}
 \begin{align}\label{sv-0}
f&\le g\quad\Longrightarrow\quad f^\text{\BSQ}\subseteq g^\text{\BSQ}, &
 D &\subseteq E\quad\Longrightarrow\quad D^\text{\SQ}\le
 E^\text{\SQ}
 \end{align}
 \begin{align}\label{sv-1}
(f^\text{\BSQ})^\text{\SQ} &\le f, & D &\subseteq ( D^\text{\SQ})^\text{\BSQ}
 \end{align}
 \begin{align}\label{sv-2}
((f^\text{\BSQ})^\text{\SQ})^\text{\BSQ} &=f^\text{\BSQ}, & ((
D^\text{\SQ})^\text{\BSQ})^\text{\SQ} &= D^\text{\SQ}
 \end{align}
\medskip

\bpr Свойства \eqref{sv-0} и \eqref{sv-1} очевидны, а \eqref{sv-2} следуют из
них:
$$
\left\{\begin{matrix}(f^\text{\BSQ})^\text{\SQ} \le f
\quad\Longrightarrow\quad\text{(применяем операцию
$\text{\BSQ}$)}\quad\Longrightarrow\quad
((f^\text{\BSQ})^\text{\SQ})^\text{\BSQ} \subseteq f^\text{\BSQ} \\ D \subseteq
( D^\text{\SQ})^\text{\BSQ} \quad\Longrightarrow\quad \text{(подстановка:
$D=f^\text{\BSQ}$)}\quad\Longrightarrow\quad
f^\text{\BSQ}\subseteq((f^\text{\BSQ})^\text{\SQ})^\text{\BSQ}\end{matrix}\right\}
\quad\Longrightarrow\quad ((f^\text{\BSQ})^\text{\SQ})^\text{\BSQ} =
f^\text{\BSQ}
$$
$$
\left\{\begin{matrix} D \subseteq ( D^\text{\SQ})^\text{\BSQ}
\quad\Longrightarrow\quad\text{(применяем операцию
$\text{\SQ}$)}\quad\Longrightarrow\quad D^\text{\SQ}\le
(( D^\text{\SQ})^\text{\BSQ})^\text{\SQ} \\
(f^\text{\BSQ})^\text{\SQ} \le f\quad\Longrightarrow\quad \text{(подстановка:
$f= D^\text{\SQ}$)}\quad\Longrightarrow\quad ((
D^\text{\SQ})^\text{\BSQ})^\text{\SQ} \le
 D^\text{\SQ}
\end{matrix}\right\}\quad\Longrightarrow\quad
(( D^\text{\SQ})^\text{\BSQ})^\text{\SQ} = D^\text{\SQ}
$$
\epr

\paragraph{Внешние огибающие на $M$.}

Если ввести обозначения
 \begin{align}\label{sv-3}
f^{\text{\BSQ}\text{\SQ}} &:=(f^\text{\BSQ})^\text{\SQ} &
 D^{\text{\SQ}\text{\BSQ}} &:=( D^\text{\SQ})^\text{\BSQ}
 \end{align}
то из \eqref{sv-0}, \eqref{sv-1}, \eqref{sv-2} будет следовать
 \begin{align}\label{sv-4}
f^{\text{\BSQ}\text{\SQ}} &\le f, & D &\subseteq D^{\text{\SQ}\text{\BSQ}}
 \end{align}
 \begin{align}\label{sv-0-0}
f&\le g\quad\Longrightarrow\quad f^{\text{\BSQ}\text{\SQ}}\le
g^{\text{\BSQ}\text{\SQ}}, & D &\subseteq E\quad\Longrightarrow\quad
 D^{\text{\SQ}\text{\BSQ}}\subseteq E^{\text{\SQ}\text{\BSQ}}
 \end{align}
 \begin{align}\label{sv-6}
(f^{\text{\BSQ}\text{\SQ}})^{\text{\BSQ}\text{\SQ}} &=f^{\text{\BSQ}\text{\SQ}}
& ( D^{\text{\SQ}\text{\BSQ}})^{\text{\SQ}\text{\BSQ}} &=
D^{\text{\SQ}\text{\BSQ}}
 \end{align}

Назовем локально ограниченную функцию $g:M\to[1;+\infty)$,
 \bit{
\item[---] {\it внешней огибающей для ограниченного множества}\index{внешняя
огибающая} $D\subseteq \mathcal O(M)$, $1\in D$, если
$$
g= D^\text{\SQ}
$$

\item[---] {\it внешней огибающей для локально ограниченной функции}
$f:M\to[1;+\infty)$ если
$$
g=f^{\text{\BSQ}\text{\SQ}}
$$

\item[---] {\it внешней огибающей на $M$}, если она удовлетворяет следующим
равносильным условиям:
 \bit{
\item[(i)] $g$ является внешней огибающей для некоторого ограниченного
множества $D\subseteq \mathcal O(M)$, $1\in D$,
$$
g= D^\text{\SQ}
$$

\item[(ii)] $g$ является внешней огибающей для некоторой локально ограниченной
функции $f:M\to[1;+\infty)$,
$$
g=f^{\text{\BSQ}\text{\SQ}}
$$

\item[(iii)] $g$ является внешней огибающей для самой себя:
$$
g^{\text{\BSQ}\text{\SQ}} = g
$$
 }\eit
 }\eit
\bpr Равносильность условий (i), (ii), (iii) требует некоторых пояснений.

$(i)\Longrightarrow (ii)$. Если $g$ -- внешняя огибающая для какого-то
множества $D$, то есть $g= D^\text{\SQ}$, то $g= D^\text{\SQ}=\eqref{sv-2}=((
D^\text{\SQ})^\text{\BSQ})^\text{\SQ}=( D^\text{\SQ})^{\text{\BSQ}\text{\SQ}}$,
то есть $g$ -- внешняя огибающая для функции $f= D^\text{\SQ}$.

$(ii)\Longrightarrow (iii)$. Если $g$ -- внешняя огибающая для некоторой
функции $f$, то есть $g=f^{\text{\BSQ}\text{\SQ}}$, то
$g^{\text{\BSQ}\text{\SQ}}=(f^{\text{\BSQ}\text{\SQ}})^{\text{\BSQ}\text{\SQ}}=\eqref{sv-6}=f^{\text{\BSQ}\text{\SQ}}=g$,
то есть $g$ -- внешняя огибающая для самой себя.

$(iii)\Longrightarrow (i)$. Если $g$ -- внешняя огибающая для самой себя, то
есть $g=g^{\text{\BSQ}\text{\SQ}}=(g^\text{\BSQ})^\text{\SQ}$, то, положив
$D=g^\text{\BSQ}$, мы получим, $g= D^\text{\SQ}$, то есть $g$ -- внешняя
огибающая для множества $D$. \epr

\bigskip

\centerline{\bf Свойства внешних огибающих:}
 {\it
 \bit{
\item[$1^{\text{\SQ}}$.]\label{1^SQ} Всякая внешняя огибающая $g$ на $M$ является непрерывной функцией на
$M$.

\item[$2^{\text{\SQ}}$.]\label{2^SQ} Для любой локально ограниченной функции $f:M\to[1;+\infty)$ ее
внешняя огибающая $f^{\text{\BSQ}\text{\SQ}}$ будет наибольшей внешней
огибающей на $M$, мажорируемой $f$:
 \bit{
\item[(a)] $f^{\text{\BSQ}\text{\SQ}}$ -- внешняя огибающая на $M$,
мажорируемая $f$:
\beq\label{2^SQ-1}
f^{\text{\BSQ}\text{\SQ}} \le f
\eeq
\item[(b)] если $g$ -- другая внешняя огибающая на $M$, мажорируемая $f$,
\beq\label{2^SQ-2}
g\le f
\eeq
то $g$ мажорируется и $f^{\text{\BSQ}\text{\SQ}}$:
\beq\label{2^SQ-3}
g\le f^{\text{\BSQ}\text{\SQ}}
\eeq
 }\eit
}\eit
 }

\bigskip
\bpr
Свойство $1^{\text{\SQ}}$ следует из предложения \ref{PROP-nepr-vnesh-ogib}. Докажем $2^{\text{\SQ}}$. Пусть $f:M\to[1;+\infty)$ --- локально ограниченная функция. Тогда, во-первых, \eqref{2^SQ-1} следует из \eqref{sv-4}. А, во-вторых, если $g$  --- внешняя огибающая, то есть $g^{\text{\BSQ}\text{\SQ}} = g$, то мы получаем импликацию  
$$
g\le f\quad\stackrel{\eqref{sv-0-0}}{\Rightarrow}\quad \underbrace{g^{\text{\BSQ}\text{\SQ}}}_{\scriptsize\begin{matrix}\|\\ g\end{matrix}}\le f^{\text{\BSQ}\text{\SQ}} \quad\Rightarrow\quad g\le f^{\text{\BSQ}\text{\SQ}}
$$
\epr

\paragraph{Прямоугольники в $\mathcal O(M)$.}

Назовем множество $ E\subseteq \mathcal O(M)$, $1\in E$,
 \bit{
\item[---] {\it прямоугольником, порожденным локально ограниченной
функцией}\index{прямоугольник} $f:M\to[1;+\infty)$ если
$$
 E=f^\text{\BSQ}
$$

\item[---] {\it прямоугольником, порожденным ограниченным множеством}
$D\subseteq \mathcal O(M)$, $1\in D$, если
$$
 E= D^{\text{\SQ}\text{\BSQ}}
$$

\item[---] {\it прямоугольником} в $\mathcal O(M)$, если выполняются следующие
равносильные условия:
 \bit{
\item[(i)] $ E$ является прямоугольником, порожденным некоторой локально
ограниченной функцией $f:M\to[1;+\infty)$
$$
 E=f^\text{\BSQ}
$$

\item[(ii)] $ E$ является прямоугольником, порожденным некоторым ограниченным
множеством $D\subseteq \mathcal O(M)$, $1\in D$,
$$
 E= D^{\text{\SQ}\text{\BSQ}}
$$

\item[(iii)] $ E$ является прямоугольником, порожденным самим собой:
$$
 E= E^{\text{\SQ}\text{\BSQ}}
$$
 }\eit
 }\eit
\bpr Равносильность условий (i), (ii), (iii) доказывается так же как в случае с
внешними огибающими. \epr

\bigskip

\centerline{\bf Свойства прямоугольников:}
 {\it
 \bit{
\item[$1^{\text{\BSQ}}$.]\label{1^BSQ} Всякий прямоугольник $E$ в $\mathcal O(M)$ является абсолютно
выпуклым компактом в $\mathcal O(M)$.

\item[$2^{\text{\BSQ}}$.]\label{2^BSQ} Для любого ограниченного множества $D\subseteq \mathcal O(M)$
порожденный им прямоугольник $D^{\text{\SQ}\text{\BSQ}}$ будет наименьшим
прямоугольником в $\mathcal O(M)$, содержащим $D$:
 \bit{
\item[(a)] $D^{\text{\SQ}\text{\BSQ}}$ -- прямоугольник в $\mathcal O(M)$,
содержащий $D$:
$$
 D\subseteq D^{\text{\SQ}\text{\BSQ}}
$$
\item[(b)] если $ E$ -- другой прямоугольник в $\mathcal O(M)$, содержащий $D$,
$$
 D\subseteq E,
$$
то $ E$ содержит и $D^{\text{\SQ}\text{\BSQ}}$:
$$
 D^{\text{\SQ}\text{\BSQ}}\subseteq E
$$
 }\eit

\item[$3^{\text{\BSQ}}$.]\label{3^BSQ} Прямоугольники в $\mathcal O(M)$ образуют фундаментальную систему
компактов в $\mathcal O(M)$: всякий компакт $D$ в $\mathcal O(M)$ содержится в
некотором прямоугольнике.

 }\eit
 }
\bigskip
\bpr
Здесь свойство $1^{\text{\BSQ}}$ следует из предложения \ref{PROP:f^blacksquare}.
\epr

\btm\label{TH:f<->D} Формулы
 \begin{align}
 D &=f^\text{\BSQ}, & f &= D^\text{\SQ}
 \end{align}
устанавливают биекцию между внешними огибающими $f$ на $M$ и прямоугольниками
$D$ в $\mathcal O(M)$. \etm
 \bpr
По определению внешних огибающих и прямоугольников, операции $f\mapsto
f^\text{\BSQ}$ и $D\mapsto D^\text{\SQ}$ переводят внешние огибающие в
прямоугольники, а прямоугольники во внешние огибающие. Более того, на этих двух
классах эти операции будут взаимно обратны. Например, если $f$ -- внешняя
огибающая, то $f^{\text{\BSQ}\text{\SQ}}=f$, то есть композиция операций
$\text{\BSQ}$ и $\text{\SQ}$ возвращает к исходной функции:
$$
f\mapsto f^\text{\BSQ}\mapsto f^{\text{\BSQ}\text{\SQ}}=f
$$
Точно также, композиция операций $\text{\SQ}$ и $\text{\BSQ}$ возвращает к
исходному множеству $D$, если оно изначально выбиралось как прямоугольник:
$$
D\mapsto D^\text{\SQ}\mapsto D^{\text{\SQ}\text{\BSQ}}=D
$$
 \epr

\subsection{Лемма о полярах}\label{SUBSEC:lemma-o-polyarah}

Напомним, что {\it полярой}\index{поляра} множества $A$ в локально выпуклом
пространстве $X$ называется множество $A^\circ$ линейных непрерывных
функционалов $f:X\to\C$, ограниченных на $A$ единицей:
$$
A^\circ=\{f\in X^\star:\quad\sup_{x\in A}|f(x)|\le 1\}
$$
Если $X$ -- стереотипное пространство и $A$ -- подмножество в сопряженном
пространстве $X^\star$, то из-за равенства $(X^\star)^\star=X$ поляру
$A^\circ\subseteq (X^\star)^\star$ удобно считать подмножеством в $X$ и
обозначать ее $^\circ\kern-2pt A$:
$$
{^\circ\kern-2pt A}=\{x\in X:\quad\sup_{f\in A}|f(x)|\le 1\}
$$
Важное для нас наблюдение здесь состоит в том, что если $A\subseteq X^\star$, и
мы сначала берем поляру $^\circ\kern-2pt A$, а затем поляру от поляры
$(^\circ\kern-2pt A)^\circ$ -- это множество называется {\it
биполярой}\index{биполяра} множества $A$ -- то оказывается, что
$(^\circ\kern-2pt A)^\circ$ представляет собой замкнутую абсолютно выпуклую
оболочку (замыкание множества линейных комбинаций вида
$\sum_{i=1}^n\lambda_i\cdot a_i$, где $a_i\in A$, $\sum_{i=1}^n|\lambda_i|\le
1$) множества $A$ в $X^\star$:
 \beq\label{bipolar}
(^\circ\kern-2pt A)^\circ=\cabsconv A
 \eeq
В этом заключается содержание классической {\it теоремы о
биполяре}\index{теорема о биполяре} применительно к стереотипным пространствам.

В частном случае, когда $A= D$ -- множество в $\mathcal O(M)$, его полярой
$D^\circ$ будет множество аналитических функционалов $\alpha\in \mathcal
O^\star(M)$, ограниченных на $D$ единицей:
$$
 D^\circ=\{\alpha\in \mathcal O^\star(M):\quad
\sup_{u\in D}|\alpha(u)|\le 1\}
$$
Наоборот, если $A$ -- какое-то множество аналитических функционалов,
$A\subseteq\mathcal O^\star(M)$, то его полярой в $\mathcal O(M)$ будет
множество, обозначаемое $^\circ\kern-2pt A$ и состоящее из функций
$u\in\mathcal O(M)$, на которых все функционалы $\alpha\in A$ ограничены
единицей:
$$
{^\circ\kern-2pt A}=\{u\in \mathcal{O}(M):\quad \sup_{\alpha\in
A}|\alpha(u)|\le 1\}
$$

\blm[\bf о полярах]\label{LM-o-polyarah} Операции перехода к поляре
удовлетворяют следующим тождествам.
 \bit{
 \item[(a)] Для ограниченного множества $D$ в $\mathcal O(M)$, содержащего
единицу, его внешняя огибающая $D^\text{\SQ}$ связана с его полярой $D^\circ$
тождеством
 \beq\label{frac_1_varPhi^square_x}
\frac{1}{D^\text{\SQ}(x)} =\max \{\lambda>0:\quad \lambda\cdot\delta^x\in
D^\circ\}
 \eeq

 \item[(b)] Для локально ограниченной функции $f:M\to[1;+\infty)$ порождаемый ей
прямоугольник $f^\text{\BSQ}$ есть поляра системы функционалов вида
$\frac{1}{f(x)}\cdot\delta^x$:
 \beq
f^\text{\BSQ} = {^\circ\kern-2pt\left\{\frac{1}{f(x)}\cdot\delta^x;\quad x\in
M\right\}}
 \eeq

 \item[(c)] Поляра прямоугольника $f^\text{\BSQ}$ есть абсолютная выпуклая оболочка системы
функционалов вида $\frac{1}{f(x)}\cdot\delta^x$:
 \beq\label{f^blacksquare^circ}
(f^\text{\BSQ})^\circ = \cabsconv \left\{\frac{1}{f(x)}\cdot\delta^x;\quad x\in
M\right\}
 \eeq
 }\eit
 \elm
\bpr (a) Для $\lambda>0$ получаем:
$$
\lambda\cdot\delta^x\in D^\circ\quad\Longleftrightarrow\quad \sup_{u\in
D}|\lambda\cdot\delta^x(u)|\le 1 \quad\Longleftrightarrow\quad
 D^\text{\SQ}(x)=\sup_{u\in D}|\delta^x(u)|\le\frac{1}{\lambda}
\quad\Longleftrightarrow\quad \lambda\le \frac{1}{D^\text{\SQ}(x)}
$$
(b) есть переформулировка определения $f^\text{\BSQ}$:
 \begin{multline*}
u\in
f^\text{\BSQ}\quad\overset{\eqref{f^blacksquare}}{\Longleftrightarrow}\quad
\forall x\in M\quad |u(x)|=|\delta^x(u)|\le f(x)
\quad\Longleftrightarrow \\
\Longleftrightarrow\quad\sup_{x\in
M}\left|\frac{1}{f(x)}\cdot\delta^x(u)\right|\le 1\quad\Longleftrightarrow\quad
u\in {^\circ\left\{\frac{1}{f(x)}\cdot\delta^x;\quad x\in M\right\}}
 \end{multline*}
(c) есть следствие (b) и теоремы о биполяре:
 $$
(f^\text{\BSQ})^\circ = \left({^\circ\left\{\frac{1}{f(x)}\cdot\delta^x;\quad
x\in M\right\}}\right)^\circ=\eqref{bipolar}=\cabsconv
\left\{\frac{1}{f(x)}\cdot\delta^x;\quad x\in M\right\}
 $$ \epr

\subsection{Внутренние огибающие на $M$ и ромбы в $\mathcal O^\star(M)$}
\label{SUBSEC:vnut-ogib}

В этом пункте нас будут интересовать замкнутые абсолютно выпуклые окрестности
нуля $\varDelta$ в $\mathcal O^\star(M)$, удовлетворяющие следующим двум
равносильным условиям:
 {\it
 \bit{
\item[(A)]\label{opredelenie-Delta} поляра ${^\circ\kern-2pt \varDelta}$
множества $\varDelta$ содержит тождественную единицу $1\in\mathcal{O}(M)$:
 \beq\label{uslovie-na-A-2}
1\in {^\circ\kern-2pt \varDelta}
 \eeq
\item[(B)] на тождественной единице $1\in\mathcal{O}(M)$ значение функционалов
$\alpha\in \varDelta$ не превосходит единицы:
 \beq\label{uslovie-na-A-1}
\forall \alpha\in \varDelta\quad |\alpha(1)|\le 1
 \eeq
 }\eit }\noindent
Эти два условия влекут еще одно (неэквивалентное (A) и (B)):
 {\it
 \bit{
\item[(C)] функционал вида $\lambda\cdot\delta^x$, где $\lambda>0$, может
принадлежать $\varDelta$ только если $\lambda\le 1$:
 \beq\label{uslovie-na-A}
\forall x\in M\quad\forall\lambda>0 \quad \Big(\lambda\cdot\delta^x\in
\varDelta\;\Longrightarrow\; \lambda\le 1\Big).
 \eeq
 }\eit }
\bpr Если $\lambda\cdot\delta^x\in \varDelta$, то $\forall u\in
{^\circ\varDelta}$ $|\lambda\cdot\delta^x(u)|\le 1$. В частности, при $u=1$ мы
получаем $|\lambda\cdot\delta^x(1)|=\lambda\cdot 1 \le 1$, то есть $\lambda\le
1$.
 \epr

Функцию $\ph:M\to(0,+\infty)$ мы называем {\it локально отделенной от нуля},
если для всякой точки $x\in M$ можно подобрать окрестность $U\owns x$ такую,
что
$$
\inf_{y\in U} \ph(y)>0
$$

\paragraph{Операции $\protect\BLZ$ и $\protect\LZ$.}

 \bprop
Для каждой локально отделенной от нуля функции $\ph:M\to(0,1]$ формула
 \beq\label{f^blacklozenge}
\ph^{\BLZ} :=\cabsconv\{\ph(x)\cdot\delta^x;\; x\in M\}
 \eeq
определяет абсолютно выпуклую окрестность нуля в пространстве аналитических
функционалов $\mathcal O^\star(M)$, удовлетворяющую условиям
\eqref{uslovie-na-A-2}-\eqref{uslovie-na-A}.

\eprop \bpr Функция $f(x)=\frac{1}{\ph(x)}$ будет локально ограниченной и
принимающей значения в интервале $[1;+\infty)$. Значит по лемме
\ref{LM-o-polyarah},
$$
(f^\text{\BSQ})^\circ=\eqref{f^blacksquare^circ}=\cabsconv
\left\{\frac{1}{f(x)}\cdot\delta^x;\quad x\in M\right\}=\cabsconv
\left\{\ph(x)\cdot\delta^x;\quad x\in M\right\}
$$
Это множество будет замкнутой абсолютно выпуклой окрестностью нуля в $\mathcal
O^\star(M)$, поскольку оно является полярой компакта $f^\text{\BSQ}$ в
$\mathcal O(M)$. Кроме того, поскольку $f\ge 1$, компакт $f^\text{\BSQ}$
содержит единицу, поэтому его поляра $\cabsconv\left\{\ph(x)\cdot\delta^x;\quad
x\in M\right\}$ удовлетворяет условиям (A),(B),(C) с.
\pageref{opredelenie-Delta}.
 \epr

 \bprop
Для каждой замкнутой абсолютно выпуклой окрестности нуля $\varDelta$ в
пространстве аналитических функционалов $\mathcal O^\star(M)$, удовлетворяющей
условиям \eqref{uslovie-na-A-2}-\eqref{uslovie-na-A-1} формула
 \beq\label{fi^lozenge}
\varDelta^{\LZ} (x):=\sup\{\lambda>0: \; \lambda\cdot\delta^x\in \varDelta\}
 \eeq
определяет непрерывную локально отделенную от нуля функцию
$\varDelta^{\LZ}:M\to(0;1]$.
 \eprop
\bpr Поскольку $\varDelta$ -- окрестность нуля в $\mathcal O^\star(M)$, ее
поляра $D={^\circ\kern-2pt \varDelta}$ будет компактом в $\mathcal O^\star(M)$,
причем $1\in D$. Внешняя огибающая $D^\text{\SQ}$ этого компакта будет
непрерывна и локально ограничена по предложению \ref{PROP-nepr-vnesh-ogib}.
Значит, функция
$$
\varDelta^{\LZ} (x):=\sup\{\lambda>0: \; \lambda\cdot\delta^x\in \varDelta=
D^\circ\}=\eqref{frac_1_varPhi^square_x}=\frac{1}{D^\text{\SQ}(x)}
$$
непрерывна и локально отделена от нуля.
 \epr
Следующие свойства доказываются так же, как \eqref{sv-0}, \eqref{sv-1} и
\eqref{sv-2}.

\medskip

\centerline{\bf Свойства операций ${\BLZ}$ и ${\LZ}$:}
 \begin{align}\label{sv-0-int}
\ph&\le \psi\quad\Longrightarrow\quad \ph^{\BLZ}\subseteq \psi^{\BLZ}, &
\varDelta &\subseteq \varGamma\quad\Longrightarrow\quad
 \varDelta^{\LZ}\le \varGamma^{\LZ}
 \end{align}
 \begin{align}\label{sv-1-int}
\ph& \le (\ph^{\BLZ})^{\LZ}, & ( \varDelta^{\LZ})^{\BLZ}&\subseteq \varDelta
 \end{align}
 \begin{align}\label{sv-2-int}
((\ph^{\BLZ})^{\LZ})^{\BLZ} &=\ph^{\BLZ}, & (( \varDelta^{\LZ})^{\BLZ})^{\LZ}
&= \varDelta^{\LZ}
 \end{align}

\medskip

\paragraph{Внутренние огибающие на $M$.}

Если ввести обозначения
 \begin{align}\label{sv-3-1}
\ph^{{\BLZ}\kern-0.5pt{\LZ}} &:=(\ph^{\BLZ})^{\LZ} &
 \varDelta^{{\LZ}\kern-0.5pt{\BLZ}} &:=( \varDelta^{\LZ})^{\BLZ}
 \end{align}
то из \eqref{sv-0-int}, \eqref{sv-1-int}, \eqref{sv-2-int} будет следовать
 \begin{align}\label{sv-0-0-1}
\ph&\le \psi\quad\Longrightarrow\quad \ph^{{\BLZ}\kern-0.5pt{\LZ}}\le
\psi^{{\BLZ}\kern-0.5pt{\LZ}}, & \varDelta &\subseteq
\varGamma\quad\Longrightarrow\quad
 \varDelta^{{\LZ}\kern-0.5pt{\BLZ}}\subseteq \varGamma^{{\LZ}\kern-0.5pt{\BLZ}}
 \end{align}
 \begin{align}\label{sv-4-1}
\ph &\le \ph^{{\BLZ}\kern-0.5pt{\LZ}}, & \varDelta^{{\LZ}\kern-0.5pt{\BLZ}}
&\subseteq
 \varDelta
 \end{align}
 \begin{align}\label{sv-6-1}
(\ph^{{\BLZ}\kern-0.5pt{\LZ}})^{{\BLZ}\kern-0.5pt{\LZ}}
&=\ph^{{\BLZ}\kern-0.5pt{\LZ}} & (
\varDelta^{{\LZ}\kern-0.5pt{\BLZ}})^{{\LZ}\kern-0.5pt{\BLZ}} &=
\varDelta^{{\LZ}\kern-0.5pt{\BLZ}}
 \end{align}

Назовем локально отделенную от нуля функцию $\psi:M\to(0;1]$
 \bit{
\item[---] {\it внутренней огибающей\index{внутренняя огибающая} для абсолютно
выпуклой окрестности нуля $\varDelta$} в $\mathcal O^\star(M)$,
$1\in{^\circ\varDelta}$, если
$$
\psi=\varDelta^{\LZ}
$$

\item[---] {\it внутренней огибающей для локально отделенной от нуля функции}
$\ph:M\to(0;1]$, если
$$
\psi=\ph^{{\BLZ}\kern-0.5pt{\LZ}}
$$

\item[---] {\it внутренней огибающей на $M$}, если она удовлетворяет следующим
равносильным условиям:
 \bit{
\item[(i)] $\psi$ является внутренней огибающей для некоторой абсолютно
выпуклой окрестности нуля $\varDelta$ в $\mathcal O^\star(M)$,
$$
\psi=\varDelta^{\LZ}
$$

\item[(ii)] $\psi$ является внутренней огибающей для некоторой локально
отделенной от нуля функции $\ph:M\to(0;1]$.
$$
\psi=\ph^{{\BLZ}\kern-0.5pt{\LZ}}
$$

\item[(iii)] $\psi$ является внутренней огибающей для самой себя:
$$
\psi^{{\BLZ}\kern-0.5pt{\LZ}} = \psi
$$
 }\eit
 }\eit
Равносильность условий (i), (ii), (iii) доказывается так же, как для внешних
огибающих.

\bigskip

\centerline{\bf Свойства внутренних огибающих:}
 {\it
 \bit{
\item[$1^{\LZ}$.]\label{1^LZ} Всякая внутренняя огибающая $\psi$ на $M$ является непрерывной
функцией на $M$.

\item[$2^{\LZ}$.]\label{2^LZ} Для любой локально отделенной от нуля функции $\ph:M\to(0;1]$ ее
внутренняя огибающая $\ph^{{\BLZ}\kern-0.5pt{\LZ}}$ будет наименьшей внутренней
огибающей на $M$, мажорирующей $\ph$:
 \bit{
\item[(a)] $\ph^{{\BLZ}\kern-0.5pt{\LZ}}$ -- внутренняя огибающая на $M$,
мажорирующая $\ph$:
$$
\ph \le \ph^{{\BLZ}\kern-0.5pt{\LZ}}
$$
\item[(b)] если $\psi$ -- другая внутренняя огибающая на $M$, мажорирующая
$\ph$,
$$
\ph\le \psi
$$
то $\psi$ мажорирует и $\ph^{{\BLZ}\kern-0.5pt{\LZ}}$:
$$
\ph^{{\BLZ}\kern-0.5pt{\LZ}}\le \psi
$$
 }\eit
}\eit
 }
\bigskip

\paragraph{Ромбы в $\mathcal O^\star(M)$.}

Назовем множество $\varGamma\subseteq \mathcal O^\star(M)$,
 \bit{
\item[---] {\it ромбом, порожденным локально отделенной от нуля
функцией}\index{ромб} $\ph:M\to(0;1]$, если
$$
\varGamma=\ph^{\BLZ}
$$

\item[---] {\it ромбом, порожденным абсолютно выпуклой окрестностью нуля}
$\varDelta\subseteq \mathcal O^\star(M)$, $1\in{^\circ\varDelta}$, если
$$
\varGamma=\varDelta^{{\LZ}\kern-0.5pt{\BLZ}}
$$

\item[---] {\it ромбом} в $\mathcal O^\star(M)$, если оно удовлетворяет
следующим равносильным условиям:
 \bit{
\item[(i)] $\varGamma$ является ромбом, порожденным некоторой локально
отделенной от нуля функцией $\ph:M\to(0;1]$,
$$
\varGamma=\ph^{\BLZ}
$$

\item[(ii)] $\varGamma$ является ромбом, порожденным некоторой абсолютно
выпуклой окрестностью нуля $\varDelta\subseteq \mathcal O^\star(M)$,
$1\in{^\circ\varDelta}$,
$$
\varGamma=\varDelta^{{\LZ}\kern-0.5pt{\BLZ}}
$$

\item[(iii)] ромб, порожденный $\varGamma$, совпадает с $\varGamma$:
$$
\varGamma=\varGamma^{{\LZ}\kern-0.5pt{\BLZ}}
$$
 }\eit
 }\eit
\bpr Равносильность условий (i), (ii), (iii) доказывается так же как в случае с
внешними огибающими. \epr

\bigskip

\centerline{\bf Свойства ромбов:}
 {\it
 \bit{
\item[$1^{\BLZ}$.]\label{1^BLZ} Всякий ромб $\varDelta$ в $\mathcal O^\star(M)$ является замкнутой
абсолютно выпуклой окрестностью нуля в $\mathcal O^\star(M)$.

\item[$2^{\BLZ}$.]\label{2^BLZ} Для любой окрестности нуля $\varDelta\subseteq \mathcal O^\star(M)$
порождаемый ею ромб $\varDelta^{{\LZ}\kern-0.5pt{\BLZ}}$ будет наибольшим
ромбом в $\mathcal O^\star(M)$, содержащимся в $\varDelta$:
 \bit{
\item[(a)] $\varDelta^{{\LZ}\kern-0.5pt{\BLZ}}$ -- ромб в $\mathcal
O^\star(M)$, содержащийся в $\varDelta$:
\beq\label{varDelta^LZ-subseteq-varDelta}
\varDelta^{{\LZ}\kern-0.5pt{\BLZ}}\subseteq \varDelta
\eeq
\item[(b)] если $\varGamma$ -- другой ромб в $\mathcal O^\star(M)$,
содержащийся в $\varDelta$,
$$
\varGamma\subseteq \varDelta,
$$
то $\varGamma$ содержится и в $\varDelta^{{\LZ}\kern-0.5pt{\BLZ}}$:
$$
\varGamma\subseteq \varDelta^{{\LZ}\kern-0.5pt{\BLZ}}
$$
 }\eit

\item[$3^{\BLZ}$.]\label{3^BLZ} Ромбы в $\mathcal O^\star(M)$ образуют фундаментальную систему
окрестностей нуля в $\mathcal O^\star(M)$: всякая окрестность нуля $\varGamma$
в $\mathcal O^\star(M)$ содержит некоторый ромб.
 }\eit
 }
\bigskip

По аналогии с теоремой \ref{TH:f<->D} доказывается

\btm Формулы
 \begin{align}
\varDelta &=\ph^{\BLZ}, & \ph &=\varDelta^{\LZ}
 \end{align}
устанавливают биекцию между внутренними огибающими $\ph$ на $M$ и ромбами
$\varDelta$ в $\mathcal O^\star(M)$. \etm

\subsection{Двойственность между прямоугольниками и ромбами}
\label{SUBSEC:pryamoug<->romb}

Из леммы о полярах \ref{LM-o-polyarah} следует

\btm\label{TH:formuly-dlya-pryamoug-i-rombov} Справедливы следующие формулы
\begin{align}
 \label{f^blacksquare^circ=frac-1-f^blacklozenge}
(f^\text{\BSQ})^\circ &=\left(\frac{1}{f}\right)^{\BLZ}, &
{^\circ\kern-2pt(\ph^{\BLZ})} &=\left(\frac{1}{\ph}\right)^\text{\BSQ},
\\
 \label{D^circ^lozenge=frac-1-D^square}
(D^\circ)^{\LZ} &=\frac{1}{D^\text{\SQ}}\;\;, &
\left({^\circ\kern-2pt\varDelta}\right)^\text{\SQ}
&=\frac{1}{\varDelta^{\LZ}}\;\;,
\\
\label{1/f^blacksquare^square=(1/f)^blacklozenge^lozenge}
 \frac{1}{f^{\text{\BSQ}\text{\SQ}}}
&=\left(\frac{1}{f}\right)^{{\BLZ}\kern-0.5pt{\LZ}}, &
\frac{1}{\ph^{{\BLZ}\kern-0.5pt{\LZ}}}
&=\left(\frac{1}{\ph}\right)^{\text{\BSQ}\text{\SQ}}, \\
 \label{D^SQ^BSQ^circ=D^circ^lozenge^blacklozenge}
\Big(D^{\text{\SQ}\text{\BSQ}}\Big)^\circ
&=\Big(D^\circ\Big)^{{\LZ}\kern-0.5pt{\BLZ}}, &
{^\circ\kern-0.5pt\Big(\varDelta^{{\LZ}\kern-0.5pt{\BLZ}}\Big)}
&=\Big({^\circ\kern-2pt\varDelta}\Big)^{\text{\SQ}\text{\BSQ}},
\end{align}
-- где $f:M\to[1;+\infty)$ -- произвольная локально ограниченная функция,
$\ph:M\to(0;1]$ -- произвольная функция, локально отделенная от нуля, $D$ --
произвольный абсолютно выпуклый компакт в $\mathcal O(M)$, $\varDelta$ --
произвольная замкнутая абсолютно выпуклая окрестность нуля в $\mathcal
O^\star(M)$. \etm

\bpr 1. Первая формула в \eqref{f^blacksquare^circ=frac-1-f^blacklozenge}
следует из \eqref{f^blacksquare^circ}:
$$
(f^\text{\BSQ})^\circ =\eqref{f^blacksquare^circ}= \cabsconv
\left\{\frac{1}{f(x)}\cdot\delta^x;\quad x\in M\right\}=\eqref{f^blacklozenge}=
\left(\frac{1}{f}\right)^{\BLZ}
$$
Из нее подстановкой $f=\frac{1}{\ph}$ следует вторая:
$$
\left(\frac{1}{\ph}\right)^{\text{\BSQ}\ \circ}
=\ph^{\BLZ}\qquad\Longrightarrow\qquad
\left(\frac{1}{\ph}\right)^{\text{\BSQ}}={^\circ\left(\left(\frac{1}{\ph}\right)^{\text{\BSQ}\
\circ}\right)} ={^\circ\left(\ph^{\BLZ}\right)}
$$

2. Далее, первая формула в \eqref{D^circ^lozenge=frac-1-D^square} следует из
\eqref{frac_1_varPhi^square_x}:
$$
\frac{1}{D^\text{\SQ}(x)}=\eqref{frac_1_varPhi^square_x} =\max \{\lambda>0:\
\lambda\cdot\delta^x\in D^\circ\}=\eqref{fi^lozenge}=(D^\circ)^{\LZ}(x)
 \qquad\Longrightarrow\qquad
(D^\circ)^{\LZ} =\frac{1}{D^\text{\SQ}}
$$
А из нее подстановкой $D={^\circ\varDelta}$ получается вторая:
$$
D^\text{\SQ}=\frac{1}{(D^\circ)^{\LZ}}
 \qquad\Longrightarrow\qquad
({^\circ\varDelta})^\text{\SQ}=\frac{1}{(({^\circ\varDelta})^\circ)^{\LZ}}=\frac{1}{\varDelta^{\LZ}}
$$

3. Теперь первая формула в
\eqref{1/f^blacksquare^square=(1/f)^blacklozenge^lozenge} получается из первой
формулы в \eqref{D^circ^lozenge=frac-1-D^square} и первой формулы в
\eqref{f^blacksquare^circ=frac-1-f^blacklozenge} подстановкой
$D=f^{\text{\BSQ}}$:
$$
\frac{1}{D^\text{\SQ}}=(D^\circ)^{\LZ}
 \qquad\Longrightarrow\qquad
\frac{1}{f^\text{\BSQ\SQ}}=(f^{\text{\BSQ}\
\circ})^{\LZ}=\eqref{f^blacksquare^circ=frac-1-f^blacklozenge}=\left(\frac{1}{f}\right)^{{\BLZ}{\LZ}}
$$
Вторая формула в \eqref{1/f^blacksquare^square=(1/f)^blacklozenge^lozenge}
получается из второй формулы в \eqref{D^circ^lozenge=frac-1-D^square} и второй
формулы в \eqref{f^blacksquare^circ=frac-1-f^blacklozenge} подстановкой
$\varDelta=\ph^{\BLZ}$:
$$
\frac{1}{\varDelta^{\LZ}}=\left({^\circ\kern-2pt\varDelta}\right)^\text{\SQ}
 \qquad\Longrightarrow\qquad
\frac{1}{\ph^{{\BLZ}{\LZ}}}=\left({^\circ\kern-2pt(\ph^{\BLZ})}\right)^\text{\SQ}=
\eqref{f^blacksquare^circ=frac-1-f^blacklozenge}=\left(\frac{1}{\ph}\right)^{\text{\BSQ}\text{\SQ}}
$$

4. Первая формула в \eqref{D^SQ^BSQ^circ=D^circ^lozenge^blacklozenge} следует
из первой формулы в \eqref{D^circ^lozenge=frac-1-D^square} и первой формулы в
\eqref{f^blacksquare^circ=frac-1-f^blacklozenge}:
$$
(D^\circ)^{\LZ} =\frac{1}{D^\text{\SQ}}
 \qquad\Longrightarrow\qquad
(D^\circ)^{{\LZ}{\BLZ}}=\left(\frac{1}{D^\text{\SQ}}\right)^{\BLZ}
=\eqref{f^blacksquare^circ=frac-1-f^blacklozenge}=(D^\text{\SQ\BSQ})^\circ
$$
Наконец, вторая формула в \eqref{D^SQ^BSQ^circ=D^circ^lozenge^blacklozenge}
следует из второй формулы в \eqref{D^circ^lozenge=frac-1-D^square} и второй
формулы в \eqref{f^blacksquare^circ=frac-1-f^blacklozenge}:
$$
\left({^\circ\kern-2pt\varDelta}\right)^\text{\SQ}=\frac{1}{\varDelta^{\LZ}}
 \qquad\Longrightarrow\qquad
\left({^\circ\kern-2pt\varDelta}\right)^\text{\SQ\BSQ}=\left(\frac{1}{\varDelta^{\LZ}}\right)^\text{\BSQ}
=\eqref{f^blacksquare^circ=frac-1-f^blacklozenge}={^\circ\kern-0.5pt\Big(\varDelta^{{\LZ}\kern-0.5pt{\BLZ}}\Big)}
$$
 \epr

Теорема \ref{TH:formuly-dlya-pryamoug-i-rombov} влечет за собой еще два важных
утверждения.

\btm\label{TH-biject-vnesh-i-vnutr-ogib} Операция перехода к обратной функции
\begin{align}
f&=\frac{1}{\ph}, & \ph &=\frac{1}{f}
\end{align}
устанавливает биекцию между внешними огибающими $f$ и внутренними огибающими
$\ph$ на $M$.
 \etm
\bpr Если $\ph$ -- внутренняя огибающая, то $\ph^{{\BLZ}{\LZ}}=\ph$, и поэтому
$\left(\frac{1}{\ph}\right)^\text{\BSQ\SQ}=\eqref{1/f^blacksquare^square=(1/f)^blacklozenge^lozenge}=
\frac{1}{\ph^{{\BLZ}\kern-0.5pt{\LZ}}}=\frac{1}{\ph}$, то есть $\frac{1}{\ph}$
-- внешняя огибающая. И наоборот, если $f$ -- внешняя огибающая, то
$f^\text{\BSQ\SQ}=f$, и поэтому
$\left(\frac{1}{f}\right)^{{\BLZ}\kern-0.5pt{\LZ}}=\eqref{1/f^blacksquare^square=(1/f)^blacklozenge^lozenge}=
\frac{1}{f^\text{\BSQ\SQ}}=\frac{1}{f}$, то есть $\frac{1}{f}$ -- внутренняя
огибающая.
 \epr

\btm\label{TH-biject-pryam-i-romb} Операции перехода к полярам
\begin{align}
D&={^\circ\kern-2pt\varDelta}, & \varDelta &=D^\circ
\end{align}
устанавливают биекцию между прямоугольниками $D$ в $\mathcal O(M)$, и ромбами
$\varDelta$ в $\mathcal O^\star(M)$.\etm
 \bpr
Если $\varDelta$ -- ромб, то $\varDelta^{{\LZ}\kern-0.5pt{\BLZ}}=\varDelta$, и
поэтому
$({^\circ\kern-2pt\varDelta})^\text{\SQ\BSQ}=\eqref{D^SQ^BSQ^circ=D^circ^lozenge^blacklozenge}=
{^\circ\kern-2pt(\varDelta^{{\LZ}\kern-0.5pt{\BLZ}})}={^\circ\kern-2pt\varDelta}$,
то есть ${^\circ\kern-2pt\varDelta}$ -- прямоугольник. Наоборот, если $D$ --
прямоугольник, то $D^\text{\BSQ\SQ}=D$, поэтому
$(D^\circ)^{{\LZ}\kern-0.5pt{\BLZ}}=\eqref{D^SQ^BSQ^circ=D^circ^lozenge^blacklozenge}=
(D^\text{\BSQ\SQ})^\circ=D^\circ$, то есть $D^\circ$ -- ромб. Поскольку переход
к поляре является биекцией между замкнутыми абсолютно выпуклыми множествами, он
будет биекцией и между прямоугольниками и ромбами.
 \epr

Из теорем \ref{TH-biject-vnesh-i-vnutr-ogib}, \ref{TH-biject-pryam-i-romb} и
\ref{TH:formuly-dlya-pryamoug-i-rombov} следует теперь

\btm Справедливы следующие формулы:
\begin{align}
\left( (f^\text{\BSQ})^\circ\right)^{\LZ} &=\frac{1}{f}, &
\left({^\circ\kern-2pt(\ph^{\BLZ})}\right)^\text{\SQ} &=\frac{1}{\ph} \\
\left(\frac{1}{(D^\circ)^{\LZ}}\right)^\text{\BSQ} &=D, &
\left(\frac{1}{\left({^\circ\kern-2pt\varDelta}\right)^\text{\SQ}}\right)^{\BLZ}
&=\varDelta,
\end{align}
где $f:M\to[1;+\infty)$ -- произвольная внешняя огибающая, $\ph:M\to(0;1]$ --
произвольная внутренняя огибающая, $D$ -- произвольный прямоугольник в
$\mathcal O(M)$, $\varDelta$ -- произвольный ромб в $\mathcal O^\star(M)$. \etm
 \bpr
Если $f$ -- внешняя огибающая, то по теореме
\ref{TH-biject-vnesh-i-vnutr-ogib}, $\frac{1}{f}$ -- внутренняя огибающая,
поэтому
$$
\left(
(f^\text{\BSQ})^\circ\right)^{\LZ}=\eqref{f^blacksquare^circ=frac-1-f^blacklozenge}=
\left(\frac{1}{f}\right)^{{\BLZ}\kern-0.5pt{\LZ}}=\frac{1}{f}
$$
Остальные формулы доказываются по аналогии.
 \epr

\section{Группы Штейна и связанные с ними алгебры Хопфа}
\label{SEC:stein-groups}

\subsection{Некоторые классы групп Штейна}\label{SUBSEC-lin-groups}

Комплексная группа Ли $G$ называется {\it группой Штейна}\index{группа!Штейна},
если $G$ является многообразием Штейна \cite{Grauert-Remmert}. По теореме
Мацусимы-Моримото \cite[XIII.5.9]{Neeb}, для комплексных групп это равносильно
условию голоморфной отделимости:
$$
\forall x\ne y\in G\quad \exists u\in \mathcal{O}(G)\quad u(x)\ne u(y)
$$
{\it Размерностью} группы Штейна называется ее размерность, как комплексного
многообразия.

Группа Штейна $G$ называется {\it компактно порожденной}\index{группа!компактно
порожденная}, если в ней существует порождающий компакт, то есть компакт
$K\subseteq G$ со свойством
$$
G=\bigcup_{n\in\N}K^n,\qquad K^n=\underbrace{K\cdot...\cdot K}_{\text{$n$
множителей}}
$$

\noindent\rule{160mm}{0.1pt}\begin{multicols}{2}

\bex
Частными случаями групп Штейна являются {\it линейные комплексные
группы}\index{линейная комплексная группа}, определяемые как комплексные группы
Ли, представимые в виде замкнутой подгруппы в полной линейной группе
$\GL(n,\C)$. Иными словами, комплексная группа $G$ считается линейной, если она
изоморфна некоторой замкнутой комплексной подгруппе $H$ в $\GL(n,\C)$ (то есть
существует изоморфизм групп $\ph:G\to H$ являющийся биголоморфным
отображением).

\eex

\bex
Еще более узкий класс групп -- {\it комплексные аффинные алгебраические
группы}\index{группа!алгебраическая}. Это подгруппы $H$ в $\GL(n,\C)$,
являющиеся одновременно алгебраическими подмногообразиями. Это означает, что
группа $H$ должна быть общим множеством нулей некоторого набора многочленов
$u_1,...,u_k$ на $\GL(n,\C)$ (под многочленом здесь можно понимать многочлен от
матричных элементов):
$$
H=\{x\in \GL(n,\C):\ u_1(x)=...=u_k(x)=0\}
$$
Если комплексная группа $G$ изоморфна какой-то алгебраической группе $H$ (то
есть существует изоморфизм групп $\ph:G\to H$ являющийся биголоморфным
отображением), то $G$ также считается алгебраической группой, потому что
алгебраические операции на $G$ будут регулярными отображениями относительно
индуцированной из $H$ структуры алгебраического многообразия.
\eex

\bex
В частности, сама полная линейная группа $\GL(n,\C)$, то есть группа невырожденных линейных
преобразований пространства $\C^n$ является комплексной алгебраической группой
(размерности $n^2$). Понятно, что $\GL(n,\C)$ компактно порождена. Как
следствие, {\it любая линейная группа компактно порождена}.
\eex

\bex Уже упоминавшийся нами в \ref{rectangles-buses} {\it комплексный тор}, то
есть фактор-группа $\C/(\Z+i\Z)$ является примером комплексной группы, не
являющейся группой Штейна. \eex

\bex\label{EX:discr-gruppa-Stein} Любая {\it дискретная группа} $G$ является группой Штейна (нулевой
размерности). Компактами в $G$ будут только конечные подмножества, поэтому $G$
будет компактно порождена в том и только в том случае, если она конечно
порождена. Поэтому, скажем, свободную группу с бесконечным числом образующих
можно считать примером не компактно порожденной группы Штейна.
 \eex

\bex Дискретная группа $G$ будет алгебраической тогда и только тогда, когда она
{\it конечна}. В этом случае ее можно рассматривать как группу преобразований
пространства $\C^n$, где $n=\card G$ -- число элементов $G$. Для этого $\C^n$
реализуют в виде пространства отображений $G$ в $\C$:
$$
x\in\C^G\quad\Longleftrightarrow\quad x:G\to\C,
$$
тогда вложение $G$ в $\GL(\C^G)$ (группу невырожденных линейных преобразований
пространства $\C^G$) описывается формулой
\begin{multline*}
\ph:G\to \GL(\C^G):\quad \ph(g)(x)(h)=x(h\cdot g),\\ g,h\in G,\quad x\in
\C^G
\end{multline*}
\eex

\bex
Если у комплексной группы Ли $G$ связная компонента единицы $G_e$ является группой Штейна, то $G$ тоже является группой Штейна.
\eex

\bex Аддитивная группа комплексных чисел $\C$ является комплексной
алгебраической группой (размерности 1), потому что она вложима как
алгебраическая подгруппа в $\GL(2,\C)$ по формуле
$$
\ph(\lambda)=\begin{pmatrix}1 & \lambda \\ 0 & 1 \end{pmatrix},\qquad
\lambda\in\C
$$
\eex

\bex\label{EX:Z} Аддитивная группа $\Z$ целых чисел является комплексной
линейной группой (нулевой размерности), потому что она вложима в группу
$\GL(2,\C)$ по формуле
$$
\ph(n)=\begin{pmatrix}1 & n \\ 0 & 1 \end{pmatrix},\qquad n\in\Z
$$
Однако $\Z$ не будет алгебраической группой, потому что ни при этом, ни при
каком-либо другом вложении в $\GL(n,\C)$ она не является общим множеством нулей
какой-либо системы многочленов (потому что дискретна и бесконечна). \eex

\bex\label{EX:C^x} Мультипликативную группу $\C^\times$ ненулевых комплексных
чисел
$$
\C^\times=\C\setminus\{0\}
$$
можно представлять себе как полную линейную группу на $\C$,
$$
\C^\times\cong\GL(1,\C),
$$
поэтому она будет линейной комплексной группой (размерности 1). Мы условимся
называть эту группу {\it комплексной окружностью}.\index{комплексная
окружность} \eex

\bex Рассмотрим действие группы $\C$ на себе самой экспонентами:
$$
\ph:\C\to\Aut(\C),\qquad \ph(a)(x)=x\cdot e^a
$$
Полупрямое произведение $\C$ на $\C$ относительно такого действия, то есть
группа $\C\ltimes\C$, совпадающая как множество с декартовым произведением
$\C\times\C$, но наделенная более сложным умножением
$$
(a,x)\cdot(b,y):=(a+b,x\cdot e^b+y)
$$
будет (связной) линейной комплексной группой, потому что она вложима в
$\GL_3(\C)$ гомоморфизмом
$$
(x,a)\mapsto \begin{pmatrix} e^a & 0 & 0 \\ x & 1 & 0 \\ 0 & 0 &
e^{ia}\end{pmatrix}
$$
Однако $\C\ltimes\C$ не будет комплексной алгебраической группой, потому что ее
центр
$$
Z(\C\ltimes\C)=\{(2\pi i n,0);\; n\in\Z\}
$$
представляет собой бесконечную дискретную подгруппу (чего у алгебраических
групп не бывает).
 \eex

\end{multicols}\noindent\rule[10pt]{160mm}{0.1pt}

\subsection{Алгебры Хопфа ${\mathcal O}(G)$, ${\mathcal O}^\star(G)$, ${\mathcal P}(G)$, ${\mathcal P}^\star(G)$}\label{Examples-of-Stein-groups}

В \cite{Akbarov-De-Gruyter-I} мы ввели следующие обозначения для важных алгебр Хопфа в комплексном анализе и алгебраической геометрии: 
\bit{
\item[---] ${\mathcal O}(G)$ --- алгебра голоморфных функций на группе Штейна $G$,

\item[---] ${\mathcal O}^\star(G)$ --- сопряженная ей алгебра аналитических функцционалов на группе Штейна $G$,

\item[---] ${\mathcal P}(G)$ --- алгебра многочленов (регулярных функций) на аффинной алгебраической группе $G$,

\item[---] ${\mathcal P}^\star(G)$ --- сопряженная ей алгебра потоков  на аффинной алгебраической группе $G$.
}\eit
Сейчас мы обсудим строение этих алгебр в простейших случаях.

\paragraph{Алгебры $\mathcal{O}(\Z)$ и $\mathcal{O}^\star(\Z)$.}

Мы уже отмечали в примере \ref{EX:Z}, что группу $\Z$ целых чисел можно
рассматривать как комплексную группу (нулевой размерности). Поскольку $\Z$
дискретна, любая функция на ней автоматически будет голоморфна. Поэтому алгебра
$\mathcal{O}(\Z)$ формально будет совпадать с алгеброй $\C^{\Z}$, а алгебра
$\mathcal{O}^\star(\Z)$ -- с алгеброй $\C_{\Z}$:
 $$
\mathcal{O}(\Z)=\C^{\Z},\qquad \mathcal{O}^\star(\Z)=\C_{\Z}
 $$
Как следствие, строение этих алгебр описывается формулами для алгебр вида $\C^M$:
характеристические функции одноточечных множеств
 \beq\label{DEF:1_n,n-in-Z}
1_n(m)=\begin{cases}0,& m=n\\ 1,& m=n\end{cases},\qquad m\in\Z,\; n\in\Z
 \eeq
образуют базис в стереотипном пространстве $\mathcal{O}(\Z)=\C^{\Z}$, а
дельта-функционалы
$$
\delta^k(u)=u(k),\qquad u\in\mathcal{O}(\Z)
$$
-- сопряженный ему алгебраический базис в $\mathcal{O}^\star(\Z)=\C_{\Z}$:
 $$
\langle 1_n,\delta^k\rangle=\begin{cases}0,& n\ne k\\ 1,& n=k\end{cases},
 $$
Элементы $\mathcal{O}(\Z)$ и $\mathcal{O}^\star (\Z)$ удобно представляются в
виде (сходящихся в этих пространствах) рядов
 \begin{align}
 \label{razlozhenie-u-v-Z}
u &\in \mathcal{O}(\Z)=\C^\Z & & \Longleftrightarrow & u&=\sum_{n\in\Z}
u(n)\cdot 1_n, & u(n) & =\delta^n(u),
 \\
 \label{razlozhenie-alpha-v-Z}
\alpha &\in \mathcal{O}^\star(\Z)=\C_\Z & & \Longleftrightarrow & \alpha
&=\sum_{n\in\Z} \alpha_n\cdot \delta^n, & \alpha_n & =\alpha(1_n),
 \end{align}
причем действие $\alpha$ на $u$ описывается формулой
 $$
\langle u,\alpha\rangle=\sum_{n\in\Z} u(n)\cdot\alpha_n
 $$
Операции умножения в $\mathcal{O}(\Z)=\C^{\Z}$ и $\mathcal{O}^\star
(\Z)=\C_{\Z}$ записываются рядами:
 \beq\label{umnozhenie-v-O(Z)-i-O^star(Z)}
u\cdot v=\sum_{n\in\Z} u(n)\cdot v(n)\cdot 1_n, \qquad \alpha
* \beta=\sum_{k\in\Z} \left(\sum_{i\in\Z}\alpha_i\cdot \beta_{k-i}\right)\cdot \delta^k,
 \eeq
(в первом случае это покоординатное умножение, а во втором -- умножение
степенных рядов).

\bprop Алгебра $\mathcal{O}(\Z)=\C^{\Z}$ функций на $\Z$ является ядерной
алгеброй Хопфа-Фреше с топологией, порожденной полунорами
 \beq\label{polunormy-v-H(Z)}
||u||_N=\sum_{|n|\le N} |u(n)|,\qquad N\in\N.
 \eeq
и алгебраическим операциями, определяемыми своими значениями на базисных
элементах $1_k$ по формулам
 \begin{align}
\label{umn-v-C^Z} & 1_m\cdot 1_n=\begin{cases}1_m,& m=n
\\ 0,& m\ne n
\end{cases} && 1_{\mathcal{P}^\star(\C^\times)}=\sum_{n\in\Z} 1_n \\
\label{koumn-v-C^Z}  &\varkappa(1_n)=\sum_{m\in\Z} 1_m\odot 1_{n-m} & &
\e(1_n)=\begin{cases}1,& n=0\\ 0,& n\ne 0\end{cases} \\
\label{antipod-v-C^Z} & \sigma(1_n)=1_{-n} &&
 \end{align}
 \eprop

\bprop Алгебра $\mathcal{O}^\star(\Z)=\C_{\Z}$ точечных зарядов на группе $\Z$
является ядерной алгеброй Хопфа-Браунера с топологией, порожденной полунорами
 \beq\label{|alpha|_r-v-Z}
|||\alpha|||_r=\sum_{n\in\Z} r_n\cdot |\alpha_n|\qquad (r_n\ge 0)
 \eeq
и алгебраическим операциями, определяемыми своими значениями на мономах
$\delta^k$ по формулам
 \begin{align}
\label{umn-v-C_Z} & \delta^k*\delta^l=\delta^{k+l} && 1_{\mathcal{P}(\C^\times)}=\delta^0 \\
\label{koumn-v-C_Z}  &\varkappa(\delta^k)= \delta^k\circledast \delta^k & &
\e(z^k)=1 \\
\label{antipod-v-C_C_Z} & \sigma(\delta^k)=\delta^{-k} &&
 \end{align}
 \eprop
\bpr Здесь может быть неочевидно, что полунормы \eqref{|alpha|_r-v-Z}
действительно определяют топологию в $\mathcal{O}^\star(\Z)=\C_{\Z}$.
Соответствующее утверждение удобно сформулировать отдельно, поскольку ниже оно
нам понадобится в лемме \ref{LM-||alpha||_C-v-Z}: \epr

\blm\label{LM-polunormy-v-H*(Z)} Если $p$ -- непрерывная полунорма на
$\mathcal{O}^\star(\Z)$, и $r_n=p(\delta^n)$, то $p$ мажорируется полунормой
\eqref{|alpha|_r-v-Z}:
 \beq\label{p(alpha)-le-|||alpha|||_r-O*(Z)}
p(\alpha)\le |||\alpha|||_r
 \eeq
 \elm
\bpr
 \beq\label{tsep-p-le-|||.|||_r}
p(\alpha)=p\left(\sum_{n\in\Z} \alpha_n\cdot\delta^n\right)\le \sum_{n\in\Z}
|\alpha_n|\cdot p(\delta^n)=\sum_{n\in\Z} |\alpha_n|\cdot r_n=|||\alpha|||_r
 \eeq
 \epr

\paragraph{Алгебры $\mathcal{P}(\C^\times)$, $\mathcal{P}^\star(\C^\times)$, $\mathcal{O}(\C^\times)$, $\mathcal{O}^\star(\C^\times)$.}\label{SUBSEC:algebry-na-C^x}

В примере \ref{EX:C^x} мы условились обозначать символом $\C^\times$
мультипликативную группу ненулевых комплексных чисел:
$$
\C^\times:=\C\setminus\{0\}.
$$
(умножение в $\C^\times$ -- обычное умножение комплексных чисел) и назвать ее
{\it комплексной окружностью}.

Алгебра $\mathcal{P}(\C^\times)$ многочленов на $\C^\times$ состоит из
многочленов Лорана, то есть функций вида
$$
u=\sum_{n\in\Z} u_n\cdot z^n
$$
где $z^n$ -- мономы на $\C^\times$:
 \beq\label{DEF:z^n-in-R(C^x)}
z^n(x):=x^n,\qquad x\in \C^\times,\qquad n\in\Z
 \eeq
и почти все $u_n\in\C$ равны нулю:
$$
\card\{n\in\Z:\ u_n\ne 0\}<\infty.
$$

А алгебра $\mathcal{P}^\star(\C^\times)$ потоков на $\C^\times$ состоит из
функционалов
$$
\alpha=\sum_{k\in\Z} \alpha_k\cdot \zeta^k
$$
где $\zeta_k$ -- функционал вычисления коэффициента степени $k$ ряда Лорана:
 \beq\label{DEF:zeta_n}
\zeta_k(u)=\frac{1}{2\pi i}\int_{|z|=1} \frac{u(z)}{z^{k+1}}\d z=\int_0^1
e^{-2\pi i k t}u(e^{2\pi i t})\d t,\qquad u\in\mathcal{P}(\C^\times)
 \eeq
и $\alpha_k\in\C$ -- произвольная последовательность.

Мономы $\zeta_k$ и $z^n$ действуют друг на друга правилом
 \beq\label{zeta^k(z^n)}
\langle z^n,\zeta_k\rangle=\begin{cases}1,& n=k \\ 0,& n\ne k \end{cases},
 \eeq
поэтому действие потока $\alpha$ на многочлен $u$ будет описываться формулой
 $$
\langle u,\alpha\rangle=\sum_{n\in\Z} u_n\cdot\alpha_n
 $$
а операции умножения в $\mathcal{P}(\C^\times)$ и $\mathcal{P}^\star
(\C^\times)$ записываются рядами:
 \beq\label{umnozhenie-v-R(C-x)-i-R^star(C-x)}
u\cdot v=\sum_{n\in\Z} \left(\sum_{i\in\Z} u_i\cdot v_{n-i}\right)\cdot z^n,
\qquad \alpha * \beta=\sum_{k\in\Z} \alpha_k\cdot \beta_k\cdot \zeta_k,
 \eeq
(в первом случае это умножение степенных рядов, а во втором -- покоординатное
умножение).

 \bprop Отображение
$$
u\in\mathcal{P}(\C^\times)\mapsto \{u_k;\;k\in\Z\}\in\C_{\Z}
$$
является изоморфизмом ядерных алгебр Хопфа-Браунера
  \beq\label{R(C^times)=C_Z}
\mathcal{P}(\C^\times)\cong \C_{\Z}
  \eeq
 \eprop

\bprop Алгебра $\mathcal{P}(\C^\times)$ многочленов на комплексной окружности
$\C^\times$ является ядерной алгеброй Хопфа-Браунера с топологией, порожденной
полунорами
 \beq\label{|||u|||_r-v-R(C^x)}
|||u|||_r=\sum_{n\in\Z} r_n\cdot |u_n|\qquad (r_n\ge 0)
 \eeq
и алгебраическим операциями, определяемыми своими значениями на мономах $z^k$
по формулам
 \begin{align}
\label{umn-v-R(C^times)} & z^k\cdot z^l=z^{k+l} && 1_{\mathcal{P}(\C^\times)}=z^0 \\
\label{koumn-v-R(C^times)}  &\varkappa(z^k)= z^k\odot z^k & &
\e(z^k)=1 \\
\label{antipode-in-R(C^times)} & \sigma(z^k)=z^{-k} &&
 \end{align}
 \eprop

 \bprop Отображение
$$
\alpha\in\mathcal{P}^\star(\C^\times)\mapsto \{\alpha_k;\;k\in\Z\}\in\C^{\Z}
$$
является изоморфизмом ядерных алгебр Хопфа-Фреше
  \beq\label{R^star(C^times)=C^Z}
\mathcal{P}^\star(\C^\times)\cong \C^{\Z}
  \eeq
 \eprop

\bprop Алгебра $\mathcal{P}^\star(\C^\times)$ потоков на комплексной
окружности $\C^\times$ является ядерной алгеброй Хопфа-Фреше с топологией,
порожденной полунорами
 \beq\label{polunormy-v-R*(C^x)}
||\alpha||_N=\sum_{|n|\le N} |\alpha_n|\qquad (N\in\N)
 \eeq
и алгебраическим операциями, определяемыми своими значениями на мономах
$\zeta_k$ по формулам
 \begin{align}
\label{umn-v-R^star(C^times)} & \zeta_k*\zeta_l=\begin{cases}\zeta_l,& k=l
\\ 0,& k\ne l
\end{cases} && 1_{\mathcal{P}^\star(\C^\times)}=\sum_{n\in\Z} \zeta_n \\
\label{koumn-v-R^star(C^times)}  &\varkappa(\zeta_k)=\sum_{l\in\Z}
\zeta_l\circledast \zeta_{k-l} & &
\e(\zeta_k)=\begin{cases}1,& k=0\\ 0,& k\ne 0\end{cases} \\
\label{antipode-in-R^star(C^times)} & \sigma(\zeta_k)=\zeta_{-k} &&
 \end{align}
 \eprop

Символом $\mathcal{O}(\C^\times)$ мы, как обычно, обозначаем алгебру
голоморфных функций на комплексной окружности $\C^\times$ (с топологией
равномерной сходимости на компактах в $\C^\times$), а символом $\mathcal{O}^\star(\C^\times)$ -- сопряженную к ней алгебру аналитических функционалов на
$\C^\times$.

Элементы алгебр $\mathcal{O}(\C^\times)$ и $\mathcal{O}^\star (\C^\times)$
удобно представлять себе в виде формальных рядов
 \begin{align}
 \label{razlozhenie-u-v-C-x}
 u &\in \mathcal{O}(\C^\times) & &
\Longleftrightarrow & u&=\sum_{n\in\Z} u_n\cdot z^n, & u_n\in\C:\quad & \forall
C>0\quad \sum_{n\in\Z}|u_n|\cdot C^{|n|}<\infty & \Big(u_n & =\zeta_n(u)\Big)
 \\
 \label{razlozhenie-alpha-v-C-x}
\alpha &\in \mathcal{O}^\star(\C^\times) & & \Longleftrightarrow & \alpha
&=\sum_{n\in\Z} \alpha_n\cdot \zeta_n, & \alpha_n\in\C:\quad & \exists C>0\quad
\forall n\in\N\quad |\alpha_n|\le C^{|n|} & \Big(\alpha_n & =\alpha(z^n)\Big)
 \end{align}
Как и в случае $\mathcal{P}(\C^\times)$ и $\mathcal{P}^\star (\C^\times)$,
действие $\alpha$ на $u$ описывается формулой
 $$
\langle u,\alpha\rangle=\sum_{n\in\Z} u_n\cdot\alpha_n
 $$
а операции умножения в $\mathcal{O}(\C^\times)$ и $\mathcal{O}^\star
(\C^\times)$ записываются, в первом случае, как обычное умножение степенных
рядов, а во втором -- как покоординатное умножение:
 \beq\label{umnozhenie-v-O(C-x)-i-O^star(C-x)}
u\cdot v=\sum_{n\in\Z} \left(\sum_{i\in\Z} u_i\cdot v_{n-i}\right)\cdot z^n,
\qquad \alpha * \beta=\sum_{n\in\Z} \alpha_n\cdot \beta_n\cdot \zeta_n,
 \eeq

\bprop Алгебра $\mathcal{O}(\C^\times)$ голоморфных функций на комплексной
окружности $\C^\times$ является ядерной алгеброй Хопфа-Фреше с топологией,
порожденной полунормами
 \beq\label{polunormy-v-H(C-x)}
||u||_C=\sum_{n\in\Z} |u_n|\cdot C^{|n|},\qquad C\ge 1.
 \eeq
и алгебраическими операциями, определяемыми своими значениями на мономах $z^k$
по тем же формулам \eqref{umn-v-R(C^times)}-\eqref{antipode-in-R(C^times)}, что
и в случае $\mathcal{P}(\C^\times)$:
 \begin{align*}
 & z^k\cdot z^l=z^{k+l} && 1_{\mathcal{O}(\C^\times)}=z^0 \\
  &\varkappa(z^k)= z^k\odot z^k & &
\e(z^k)=1 \\
 & \sigma(z^k)=z^{-k} &&
 \end{align*}
 \eprop
\bpr Любая обычная полунорма $|u|_K=\max_{x\in K}|u(x)|$, где $K$ -- компакт в
$\C^\times$, мажорируется некоторой полунормой $||u||_C$, а именно для этого
нужно взять $C=\max_{x\in K}\max\{|x|,\frac{1}{|x|}\}$:
$$
|u|_K=\max_{x\in K}|u(x)|=\max_{x\in K}\left|\sum_{n\in\Z} u_n\cdot
x^n\right|\le\max_{x\in K}\sum_{n\in\Z} |u_n|\cdot |x^n|\le\sum_{n\in\Z}
|u_n|\cdot C^{|n|}=||u||_C
$$
Наоборот, для любого числа $C\ge 1$ мы можем выбрать компакт
$K=\{t\in\C:\;\frac{1}{C+1}\le |t|\le C+1\}$, и тогда из формул Коши для
коэффициентов ряда Лорана
$$
|u_n|\le
|u|_K\cdot\min\left\{(C+1)^n;\frac{1}{(C+1)^n}\right\}=|u|_K\cdot(C+1)^{-|n|}=
\frac{|u|_K}{(C+1)^{|n|}}
$$
следует, что $||u||_C$ подчинена $|u|_K$ (с положительным множителем
$(C+1)^2$):
$$
||u||_C=\sum_{n\in\Z} |u_n|\cdot C^{|n|}\le \sum_{n\in\Z}
\frac{|u|_K}{(C+1)^{|n|}}\cdot C^{|n|}\le |u|_K\cdot\sum_{n\in\Z}
\left(\frac{C}{C+1}\right)^{|n|}=(C+1)^2\cdot|u|_K
$$
\epr

\bprop\label{PROP:O*(C^x)-Hopf} Алгебра $\mathcal{O}^\star(\C^\times)$
аналитических функционалов на комплексной окружности $\C^\times$ является
ядерной алгеброй Хопфа-Браунера с топологией, порожденной полкнормами
 \beq\label{|alpha|_r-v-C-x}
|||\alpha|||_r=\sup_{u\in E_r}|\alpha(u)|= \sum_{n\in\Z} r_n\cdot
|\alpha_n|\qquad \left(r_n\ge 0: \forall C>0\qquad \sum_{n\in\Z} r_n\cdot
C^{|n|}<\infty\right)
 \eeq
и алгебраическим операциями, определяемыми своими значениями на мономах
$\zeta_k$ по тем же формулам
\eqref{umn-v-R^star(C^times)}-\eqref{antipode-in-R^star(C^times)}, что и в
случае $\mathcal{P}^\star(\C^\times)$:
 \begin{align*}
 & \zeta_k*\zeta_l=\begin{cases}\zeta_l,& k=l
\\ 0,& k\ne l
\end{cases} && 1_{\mathcal{O}^\star(\C^\times)}=\sum_{n\in\Z} \zeta_n \\
  &\varkappa(\zeta_k)=\sum_{l\in\Z} \zeta_l\circledast
\zeta_{k-l} & &
\e(\zeta_k)=\begin{cases}1,& k=0\\ 0,& k\ne 0\end{cases} \\
 & \sigma(\zeta_k)=\zeta_{-k} &&
 \end{align*}
 \eprop
\bpr Для всякой последовательности неотрицательных чисел, $r_n\ge 0$,
удовлетворяющих условию
 \beq\label{forall-R>0-sum_r_n-R^n<infty-x}
\forall C>0\qquad \sum_{n\in\Z} r_n\cdot C^{|n|}<\infty
 \eeq
множество
 \beq\label{E_r-v-H(C-x)}
E_r=\{u\in \mathcal{O}(\C^\times):\; \forall n\in\Z\quad |u_n|\le r_n\}
 \eeq
является компактом в $\mathcal{O}(\C^\times)$ потому что оно замкнуто и
содержится в прямоугольнике $f^\text{\BSQ}$, где
$$
f(t)=\sum_{n\in\Z} r_n\cdot |t|^n,\qquad t\in\C^\times
$$
Значит, $E_r$ порождает непрерывную полунорму $\alpha\mapsto \max_{u\in
E_r}|\langle u,\alpha\rangle|$ на $\mathcal{O}^\star(\C^\times)$. Но это как
раз полунорма \eqref{|alpha|_r-v-C-x}:
$$
\max_{u\in E_r}|\langle u,\alpha\rangle|=\max_{|u|\le r_n} \left|\sum_{n\in\Z}
u_n\cdot\alpha_n \right|=\sum_{n\in\Z} r_n\cdot|\alpha_n|=|||\alpha|||_r
$$
Остается убедиться, что полунормы $|||\cdot|||_r$ действительно порождают
топологию пространства $\mathcal{O}^\star(\C^\times)$. Это вытекает из
следующей леммы: \epr

\blm\label{LM-polunormy-v-H*(C-x)-1} Если $p$ -- непрерывная полунорма на
$\mathcal{O}^\star(\C^\times)$, то семейство чисел
$$
r_n=p(\zeta_n)
$$
удовлетворяет условию \eqref{forall-R>0-sum_r_n-R^n<infty-x}, а $p$
мажорируется полунормой $|||\cdot|||_r$:
 \beq\label{p(alpha)-le-|||alpha|||_r-O^(C^x)}
p(\alpha)\le |||\alpha|||_r
 \eeq
 \elm
\bpr Множество
$$
D=\{u\in\mathcal{O}(\C^\times):\; \sup_{\alpha\in \mathcal{O}^\star(\C^\times):\; p(\alpha)\le 1}|\alpha(u)|\le 1\}
$$
будет компактом в $\mathcal{O}(\C^\times)$, порождающим полунорму $p$:
$$
p(\alpha)=\sup_{u\in D}|\alpha(u)|
$$
Поэтому
$$
r_n=p(\zeta_n)=\sup_{u\in D} |\zeta_n(u)|=\sup_{u\in D} |u_n|
$$
Для любого $C>0$ получаем:
$$
\infty>\sup_{u\in D}||u||_{C+1}=\sup_{u\in D}\sum_{n\in\Z} |u_n|\cdot
(C+1)^{|n|}\ge \sup_{u\in D}\sup_n |u_n|\cdot (C+1)^{|n|}=\sup_n \Big(r_n\cdot
(C+1)^{|n|}\Big)
$$
$$
\Downarrow
$$
$$
\exists M>0\quad \forall n\in\Z\qquad r_n\le \frac{M}{(C+1)^{|n|}}
$$
$$
\Downarrow
$$
$$
\sum_{n\in\Z} r_n\cdot C^{|n|} \le \sum_{n\in\Z} \frac{M}{(C+1)^{|n|}}\cdot
C^{|n|}<\infty
$$
То есть $r_n$ действительно удовлетворяет условию
\eqref{forall-R>0-sum_r_n-R^n<infty-x}. Формула
\eqref{p(alpha)-le-|||alpha|||_r-O^(C^x)} доказывается цепочкой неравенств,
аналогичной \eqref{tsep-p-le-|||.|||_r}. \epr

\paragraph{Цепочка $\mathcal{P}(\C)\subset\mathcal{O}(\C)\subset\mathcal{O}^\star(\C)\subset\mathcal{P}^\star(\C)$.}

Символом $\mathcal{P}(\C)$ мы обозначаем обычную алгебру многочленов на
комплексной плоскости $\C$. Пусть $t^k$ обозначает моном степени $k\in\N$ на
$\C$:
 \beq\label{def-t^k}
t^k(x):=x^k,\qquad x\in\C,\; k\in\N
 \eeq
Поскольку всякий многочлен $u\in \mathcal{P}(\C)$ однозначно представляется
рядом (с конечным числом ненулевых слагаемых)
 \beq\label{razlozhenie-u-v-R(C)}
u=\sum_{k\in\N} u_k\cdot t^k,\qquad u_k\in\C: \quad \card\{k\in\N: \; u_k\ne
0\}<\infty
 \eeq
мономы $t^k$ образуют алгебраический базис в пространстве $\mathcal{P}(\C)$.
Умножение в $\mathcal{P}(\C)$ есть обычное умножение многочленов
$$
u\cdot v=\sum_{k\in\N} \left(\sum_{i=0}^n u_{k-i}\cdot v_i\right)\cdot t^k
$$
а топология в $\mathcal{P}(\C)$ определяется как сильнейшая локально выпуклая
топология. Отсюда следует,

 \bprop Отображение
$$
u\in\mathcal{P}(\C)\mapsto \{u_k;\;k\in\N\}\in\C_{\N}
$$
является изоморфизмом топологических векторных пространств
  \beq\label{R(C)=C_N}
\mathcal{P}(\C)\cong \C_{\N}
  \eeq
 \eprop

\brem Формулу \eqref{R(C)=C_N} можно также понимать, как изоморфизм алгебр,
если умножение в $\C_{\N}$ задать той же формулой, что и для алгебры $\C_G$ точечных
зарядов на группе $G$ (определенной в
\cite[5.1.5]{Akbarov-De-Gruyter-I}), 
$$
\delta^x*\delta^y=\delta^{x\cdot y},
$$
но только при этом нужно помнить,
что множество $\N$ не будет группой, а лишь моноидом относительно используемой
на нем операции сложения.

С другой стороны, формула \eqref{R(C)=C_N} не будет изоморфизмом алгебр Хопфа,
например, потому что $\C_{\N}$ вообще не будет алгеброй Хопфа относительно операций,
определенных в \cite[5.2.1]{Akbarov-De-Gruyter-I}: для
$x\in\N$ обратный элемент $x^{-1}=-x$ не существует, если $x\ne 0$, и значит
антипод здесь не может определяться равенством 
$$
\sigma(\delta^x)=\delta^{x^{-1}}
$$
\erem

\bprop Алгебра $\mathcal{P}(\C)$ многочленов на комплексной плоскости $\C$
является ядерной алгеброй Хопфа-Браунера с топологией, порожденной полунормами
 \beq
\norm{u}_r=\sum_{k\in\N}r_k\cdot |u_k|,\qquad (r_k\ge 0)
 \eeq
и алгебраическим операциями, определяемыми своими значениями на мономах $t^k$
по формулам
 \begin{align}
\label{umn-v-R(C)} & t^k\cdot t^l=t^{k+l} && 1_{\mathcal{P}(\C)}= t^0 \\
\label{koumn-v-R(C)}  &\varkappa(t^k)=\sum_{i=0}^k\begin{pmatrix}k
\\ i\end{pmatrix}\cdot t^{k-i}\odot t^i & &
\e(t^k)=\begin{cases}1,& k=0\\ 0,& k>0\end{cases} \\
\label{antipode-in-R(C)} & \sigma(t^k)=(-1)^k\cdot t^k &&
 \end{align}
 \eprop
\bpr Алгебраические операции, которые мы еще не нашли --- коумножение, коединица
и антипод --- вычисляются по формулам \cite[(5.117)-(5.121)]{Akbarov-De-Gruyter-I} (тем же, что и для ${\mathcal C}(G)$). Например, коумножение:
 $$
\widetilde{\varkappa}(t^k)(x,y)=t^k(x+y)=(x+y)^k=\sum_{i=0}^k\begin{pmatrix}k
\\ i\end{pmatrix}\cdot x^{k-i}\cdot y^i=\sum_{i=0}^k\begin{pmatrix}k
\\ i\end{pmatrix}\cdot t^{k-i}\boxdot t^i(x,y)
 $$
 $$
 \Downarrow
 $$
 $$
 \widetilde{\varkappa}(t^k)= \sum_{i=0}^k\begin{pmatrix}k
\\ i\end{pmatrix}\cdot t^{k-i}\boxdot t^i
 $$
 $$
 \Downarrow
 $$
 $$
\varkappa(t^k)=\rho_{\C,\C}(\widetilde{\varkappa}(t^k))=
\sum_{i=0}^k\begin{pmatrix}k
\\ i\end{pmatrix}\cdot t^{k-i}\odot t^i
 $$
\epr

Следуя терминологии \cite{Ak03}, мы называем {\it потоком степени
$0$}\index{поток}, или просто {\it потоком} на комплексной плоскости $\C$
произвольный линейный (и тогда он автоматически будет непрерывным) функционал
$\alpha:\mathcal{P}(\C)\to\C$ на пространстве многочленов $\mathcal{P}(\C)$.
Из формулы \eqref{R(C)=C_N} следует, что всякий компакт в пространстве
многочленов $\mathcal{P}(\C)$ конечномерен, и значит, содержится в выпуклой
оболочке конечного набора базисных мономов $t^k$. Поэтому топологию
пространства $\mathcal{P}^\star(\C)$ потоков на $\C$ (стандартно определяемую
как топологию равномерной сходимости на компактах в $\mathcal{P}(\C)$) можно
описать как топологию сходимости на мономах $t^k$, то есть как топологию,
порождаемую полунормами
$$
\norm{\alpha}_N=\sum_{k\in N} \abs{\alpha(t^k)},
$$
где $N$ -- произвольное конечное множество в $\N$.

Типичный пример потока -- функционал вычисления производной степени $k$ в точке
$0$:
 \beq\label{def-tau^k}
\tau^k(u)=\left(\frac{\d^k}{\d x^k}u(x)\right)\Bigg|_{x=0},\qquad u\in\mathcal{P}(\C)
 \eeq
По теореме Тейлора, эти функционалы связаны с коэффициентами $u_k$ в разложении
\eqref{razlozhenie-u-v-R(C)} многочлена $u\in\mathcal{P}(\C)$ формулой
 \beq
u_k=\frac{1}{k!}\cdot\tau^k(u)
 \eeq
Отсюда следует, что действие произвольного потока $\alpha\in \mathcal{P}^\star(\C)$ на многочлен $u\in\mathcal{P}(\C)$ можно записать формулой
 \beq \label{deistvie-alpha-na-u-v-R(C)}
\alpha(u)=\alpha\left(\sum_{k\in\N}u_k\cdot t^k\right)=\sum_{k\in\N}u_k\cdot
\alpha(t^k)=\sum_{k\in\N}\frac{1}{k!}\cdot\tau^k(u)\cdot
\alpha(t^k)=\left(\sum_{k\in\N}\frac{1}{k!}\cdot
\alpha(t^k)\cdot\tau^k\right)(u)
 \eeq
Это значит, что $\alpha$ раскладывается в ряд по $\tau^k$ следующим образом:
 \beq \label{razlozhenie-alpha-v-R(C)}
\alpha=\sum_{k\in\N} \alpha_k\cdot \tau^k,\qquad
\alpha_k=\frac{1}{k!}\cdot\alpha(t^k)
 \eeq
Из этого можно сделать вывод, что потоки $\tau^k$ образуют базис в
топологическим векторном пространстве $\mathcal{P}^\star(\C)$: всякий
функционал $\alpha\in\mathcal{P}^\star(\C)$ однозначно представим в виде
сходящегося в $\mathcal{P}^\star(\C)$ ряда \eqref{razlozhenie-alpha-v-R(C)},
коэффициенты которого $\alpha_k\in\C$ непрерывно зависят от $\alpha\in\mathcal{P}^\star(\C)$.

Из \eqref{deistvie-alpha-na-u-v-R(C)} следует, что действие потока $\alpha$ на
многочлен $u$ описывается формулой
$$
\alpha(u)=\sum_{k\in\N}\underbrace{\frac{1}{k!}\cdot\tau^k(u)}_{u_k}\cdot
\underbrace{\frac{1}{k!}\cdot \alpha(t^k)}_{\alpha_k}\cdot k!=\sum_{k\in\N}
u_k\cdot\alpha_k\cdot k!
$$
а базисные потоки $\tau^k$ действуют на мономы $t^k$ по формуле
$$
\tau^k(t^n)=\begin{cases}0,& n\ne k\\ n!,& n=k\end{cases},
$$
Это означает, что система $\tau^k$ не является сопряженным базисом к $t^k$: она
отличается от сопряженного базиса скалярными множителями $k!$. Тем не менее,
поскольку, с одной стороны, $\mathcal{P}^\star(\C)\cong \C^{\N}$ (это следует
из формулы \eqref{R(C)=C_N}) и, с другой, в $\C^{\N}$ любые два базиса
изоморфны \cite[Theorem 4.5.11]{Akbarov-De-Gruyter-I}, справедливо

\bprop Отображение
$$
\alpha\in\mathcal{P}^\star(\C)\mapsto \{\alpha_k;\; k\in\N\}\in\C^{\N}
$$
является изоморфизмом топологических векторных пространств:
 \beq\label{R*(C)=C^N}
\mathcal{P}^\star(\C)\cong \C^{\N}
 \eeq
 \eprop
Этот изоморфизм не будет, однако, изоморфизмом алгебр, потому что в $\mathcal{P}^\star(\C)$ элементы умножаются друг на друга не покомпонентно, как можно было бы определить
в $\C^{\N}$ по аналогии со случаем $\C^{M}$, а как степенные ряды:
 \beq\label{umnozhenie-v-R^star(C)}
\alpha * \beta=\sum_{k\in\N} \left(\sum_{i=0}^k \alpha_{k-i}\cdot
\beta_i\right)\cdot \tau^k,
 \eeq
Это следует из формулы \eqref{umn-v-R*(C)} ниже:

\bprop Алгебра $\mathcal{P}^\star(\C)$ потоков нулевой степени на комплексной
плоскости $\C$ является ядерной алгеброй Хопфа-Фреше с топологией, порожденной
полунормами
 \beq
\norm{\alpha}_K=\sum_{k=0}^K |\alpha_k|,\qquad (K\in\N)
 \eeq
и алгебраическим операциями, определенными на базисных элементах $\tau^k$
формулами
 \begin{align}
\label{umn-v-R*(C)} & \tau^k*\tau^l=\tau^{k+l} && 1_{\mathcal{P}^\star(\C)}=\tau^0 \\
\label{koumn-v-R*(C)}  &\varkappa(\tau^k)=\sum_{i=0}^k\begin{pmatrix}k
\\ i\end{pmatrix}\cdot \tau^{k-i}\circledast\tau^i & &
\e(\tau^k)=\begin{cases}1,& k=0\\ 0,& k>0\end{cases} \\
\label{antipode-in-R*(C)} & \sigma(\tau^k)=(-1)^k\cdot \tau^k &&
 \end{align}
 \eprop
\bpr Доказывать эти формулы можно по-разному, например, можно воспользоваться
формулой \eqref{koumn-v-R(C)} для коумножения в $\mathcal{P}(\C)$:
 \begin{multline*}
(\tau^k*\tau^l)(t^m)=(\tau^k\circledast\tau^l)(\varkappa(t^m))=\eqref{koumn-v-R(C)}=
(\tau^k\circledast\tau^l)\left(\sum_{i=0}^m \begin{pmatrix}m \\ i \end{pmatrix}
\cdot t^{m-i}\odot t^i \right)=\\=\sum_{i=0}^m \begin{pmatrix}m \\ i
\end{pmatrix}\cdot \tau^k(t^{m-i})\cdot\tau^l(t^i)=\left\{\begin{matrix} \begin{pmatrix}m \\
i \end{pmatrix}\cdot (m-l)!\cdot l!,& m=k+l \\ 0,& m\ne
k+l\end{matrix}\right\}=\left\{\begin{matrix} m!,& m=k+l \\ 0,& m\ne
k+l\end{matrix}\right\}=\tau^{k+l}(t^m).
 \end{multline*}
(на любом мономе $t^m$ действие функционалов $\tau^k*\tau^l$ и $\tau^{k+l}$
совпадает, значит они сами совпадают). \epr

Рассмотрим теперь алгебру $\mathcal{O}(\C)$ целых функций и сопряженную ей
алгебру $\mathcal{O}^\star(\C)$ аналитических функционалов на комплексной
плоскости $\C$. Пусть, по-прежнему, $t^k$ -- моном степени $k\in\N$ на $\C$, а
$\tau^k$ -- функционал вычисления производной степени $k\in\N$:
$$
t^k(z)=z^k\qquad \tau^k(u)=\left(\frac{\d^k}{\d
z^k}u(z)\right)\Bigg|_{z=0}\qquad \Big(x\in\C,\qquad u\in\mathcal{O}(\C)\Big)
$$
Тогда элементы $\mathcal{O}(\C)$ и $\mathcal{O}^\star (\C)$ удобно
представлять себе в виде (сходящихся в этих пространствах) рядов
 \begin{align}
\label{razlozhenie-u-v-C} u &\in \mathcal{O}(\C) & & \Longleftrightarrow &
u&=\sum_{n=0}^\infty u_n\cdot t^n, & u_n=\frac{1}{n!}\tau^n(u)\in\C:\quad &
\forall C>0\quad \sum_{n=0}^\infty|u_n|\cdot C^n<\infty
 \\
\label{razlozhenie-alpha-v-C} \alpha &\in \mathcal{O}^\star(\C) & &
\Longleftrightarrow & \alpha &=\sum_{n=0}^\infty \alpha_n\cdot \tau^n,  &
\alpha_n=\frac{1}{n!}\alpha(t^n)\in\C:\quad & \exists M,C>0\quad \forall
n\in\N\quad |\alpha_n|\le M\cdot \frac{C^n}{n!}
 \end{align}
Действие аналитического функционала $\alpha$ на целую функцию $u$ будет
описываться формулой
$$
\langle u,\alpha\rangle=\sum_{n=0}^\infty u_n\cdot\alpha_n\cdot n!
$$
а операции умножения в $\mathcal{O}(\C)$ и $\mathcal{O}^\star (\C)$ задаются
теми же формулами, что для обычных степенных рядов:
 \beq\label{umnozhenie-v-O(C)-i-O^star(C)}
u\cdot v=\sum_{n=0}^\infty \left(\sum_{i=0}^n u_i\cdot v_{n-i}\right)\cdot t^n,
\qquad \alpha * \beta=\sum_{k\in\N} \left(\sum_{i=0}^k \alpha_i\cdot
\beta_{k-i}\right)\cdot \tau^k,
 \eeq

\bprop Алгебра $\mathcal{O}(\C)$ целых функций на комплексной плоскости $\C$
является ядерной алгеброй Хопфа-Фреше с топологией, порожденной полунормами
 \beq\label{polunormy-v-H(C)}
\norm{u}_C=\sum_{k\in\N} |u_k|\cdot C^k\qquad (C\ge 1)
 \eeq
и алгебраическим операциями, определяемыми своими значениями на мономах $t^k$
по тем же формулам \eqref{umn-v-R(C)}-\eqref{antipode-in-R(C)}, что и в случае
$\mathcal{P}(\C)$:
 \begin{align*}
& \mu(t^k\odot t^l)=t^{k+l} && 1_{\mathcal{P}(\C)}= t^0 \\
&\varkappa(t^k)=\sum_{i=0}^k\begin{pmatrix}k
\\ i\end{pmatrix}\cdot t^{k-i}\odot t^i & &
\e(t^k)=\begin{cases}1,& k=0\\ 0,& k>0\end{cases} \\
& \sigma(t^k)=(-1)^k\cdot t^k &&
 \end{align*}
 \eprop
\bpr Формулы для алгебраических операций доказываются так же, как и у
$\mathcal{P}(\C)$, поэтому нам нужно лишь объяснить, почему топология задается
полунормами \eqref{polunormy-v-H(C)}. Это тоже предмет математического
фольклора (см.напр.\cite{Pi2008}): ясно, что любая обычная полунорма
$|u|_K=\max_{z\in K}|u(z)|$, где $K$ -- компакт в $\C$, мажорируется некоторой
полунормой $||u||_C$, а именно для этого нужно взять $C=\max_{z\in K}|z|$:
$$
|u|_K=\max_{z\in K}|u(z)|=\max_{z\in K}\left|\sum_{n=0}^\infty u_n\cdot
z^n\right|\le\max_{z\in K}\sum_{n=0}^\infty |u_n|\cdot
|z^n|\le\sum_{n=0}^\infty |u_n|\cdot C^n=||u||_C
$$
Наоборот, для любого числа $C>0$ мы можем выбрать компакт $K=\{z\in\C:\; |z|\le
C+1\}$, и тогда из формул Коши для коэффициентов
$$
|u_n|\le \frac{|u|_K}{(C+1)^n}
$$
следует, что $||u||_C$ подчинена $|u|_K$ (с положительным множителем $C+1$):
$$
||u||_C=\sum_{n=0}^\infty |u_n|\cdot C^n\le \sum_{n=0}^\infty
\frac{|u|_K}{(C+1)^n}\cdot C^n\le |u|_K\cdot\sum_{n=0}^\infty
\left(\frac{C}{C+1}\right)^n=|u|_K\cdot
\frac{1}{1-\frac{C}{C+1}}=(C+1)\cdot|u|_K
$$
\epr

\bprop\label{PROP:O*(C)-Hopf} Алгебра $\mathcal{O}^\star(\C)$ аналитических
функционалов на комплексной плоскости $\C$ является ядерной алгеброй
Хопфа-Браунера с топологией, порожденной полунормами
 \beq\label{|alpha|_r-v-C}
|||\alpha|||_r=\sum_{k\in\N} r_k\cdot |\alpha_k|\cdot k!\qquad\left(r_k\ge 0:\
\forall C>0\ \sum_{k\in\N} r_k\cdot C^k<\infty \right)
 \eeq
и алгебраическим операциями, определяемыми своими значениями на базисных
элементах $\tau^k$ по тем же формулам
\eqref{umn-v-R*(C)}-\eqref{antipode-in-R*(C)}, что и в случае $\mathcal{P}^\star(\C)$:
 \begin{align*}
& \mu(\tau^k\circledast\tau^l)=\tau^{k+l} && 1_{\mathcal{P}^\star(\C)}=\tau^0 \\
&\varkappa(\tau^k)=\sum_{i=0}^k\begin{pmatrix}k
\\ i\end{pmatrix}\cdot \tau^{k-i}\circledast\tau^i & &
\e(\tau^k)=\begin{cases}1,& k=0\\ 0,& k>0\end{cases} \\
& \sigma(\tau^k)=(-1)^k\cdot \tau^k &&
 \end{align*}
 \eprop
 \bpr
Здесь тоже формулы для алгебраических операций доказываются так же, как и у
$\mathcal{P}^\star(\C)$, поэтому нам нужно лишь объяснить, почему топология
задается полунормами \eqref{|alpha|_r-v-C}. Заметим сначала, что для всякой
последовательности неотрицательных чисел, $r_k\ge 0$, удовлетворяющих условию
 \beq\label{forall-R>0-sum_r_n-R^n<infty}
\forall C>0\qquad \sum_{k\in\N} r_k\cdot C^k<\infty
 \eeq
множество
 \beq\label{E_r-v-H(C)}
E_r=\{u\in \mathcal{O}(\C):\; \forall k\in\N\quad |u_k|\le r_k\}
 \eeq
является компактом в $\mathcal{O}(\C)$, потому что замкнуто и содержится в
прямоугольнике $f^\text{\BSQ}$, где
$$
f(x)=\sum_{k\in\N} r_k\cdot |x|^k,\qquad x\in\C
$$
Значит, $E_r$ порождает непрерывную полунорму $\alpha\mapsto \max_{u\in
E_r}|\langle u,\alpha\rangle|$ на $\mathcal{O}^\star(\C)$. Но это как раз
полунорма \eqref{|alpha|_r-v-C}:
$$
\max_{u\in E_r}|\langle u,\alpha\rangle| =\max_{|u_k|\le r_k}
\left|\sum_{k\in\N} u_k\cdot\alpha_k\cdot k! \right|=\sum_{k\in\N}
r_k\cdot|\alpha_k|\cdot k!=|||\alpha|||_r
$$
Остается убедиться, что полунормы $|||\cdot|||_r$ действительно порождают
топологию пространства $\mathcal{O}^\star(\C)$. Это вытекает из следующей
леммы: \epr

\blm\label{LM-polunormy-v-H*(C)} Если $p$ -- непрерывная полунорма на
$\mathcal{O}^\star(\C)$, то числа $r_k=\frac{1}{k!}p(\tau^k)$ удовлетворяют
условию \eqref{forall-R>0-sum_r_n-R^n<infty}, а $p$ мажорируется полунормой
$|||\cdot|||_r$:
 \beq\label{p(alpha)-le-|||alpha|||_r}
p(\alpha)\le |||\alpha|||_r
 \eeq
\elm
 \bpr Множество
$$
D=\{u\in\mathcal{O}(\C):\; \sup_{\alpha\in \mathcal{O}^\star(\C):\;
p(\alpha)\le 1}|\alpha(u)|\le 1\}
$$
будет компактом в $\mathcal{O}(\C)$, порождающим полунорму $p$:
$$
p(\alpha)=\sup_{u\in D}|\alpha(u)|
$$
Поэтому
$$
r_k=\frac{1}{k!}p(\tau^k)=\frac{1}{k!}\sup_{u\in D} |\tau^k(u)|=\sup_{u\in D}
|u_k|
$$
Для любого $C>0$ получаем:
$$
\infty>\sup_{u\in D}||u||_{C+1}=\sup_{u\in D}\sum_{k\in\N} |u_k|\cdot
(C+1)^k\ge \sup_{u\in D}\sup_k |u_k|\cdot (C+1)^k=\sup_k \Big(r_k\cdot
(C+1)^k\Big)
$$
$$
\Downarrow
$$
$$
\exists M>0\quad \forall k\in\N\qquad r_k\le \frac{M}{(C+1)^k}
$$
$$
\Downarrow
$$
$$
\sum_{k\in\N} r_k\cdot C^k \le \sum_{k\in\N} \frac{M}{(C+1)^k}\cdot C^k<\infty
$$
То есть $r_k$ действительно удовлетворяет условию
\eqref{forall-R>0-sum_r_n-R^n<infty}. Формула \eqref{p(alpha)-le-|||alpha|||_r}
доказывается цепочкой неравенств, аналогичной \eqref{tsep-p-le-|||.|||_r}.
 \epr

\bprop Отображения
$$
t^k\mapsto t^k\mapsto \tau^k\mapsto\tau^k \qquad (k\in\N)
$$
определяют цепочку гомоморфизмов жестких алгебр Хопфа:
 \beq\label{R->O->O*->R*}
\mathcal{P}(\C)\to \mathcal{O}(\C)\to \mathcal{O}^\star(\C)\to \mathcal{P}^\star(\C)
 \eeq
\eprop

\subsection{Двойственность Понтрягина для абелевых компактно порожденных групп Штейна}
\label{SUBSEC:pont-duality-for-Stein}

Комплексная окружность $\C^\times$, о которой мы говорили в
\ref{SEC:stein-groups}\ref{SUBSEC-lin-groups}, занимает среди всех абелевых
компактно порожденных групп Штейна то же место, что и обычная <<вещественная>>
окружность $\T=\R / \Z$ среди всех локально компактных абелевых групп (или
среди компактно порожденных вещественных абелевых групп Ли), потому что для
абелевых компактно порожденных групп Штейна справедлив следующий вариант теории
двойственности Понтрягина.

Пусть $G$ --  абелева компактно порожденная группа Штейна. Назовем {\it
голоморфным характером} на $G$ произвольный голоморфный гомоморфизм группы $G$
в комплексную окружность:
\beq\label{G^bullet-v-kompl-geom}
\chi\in G^\bullet\quad\Longleftrightarrow\quad \chi:G\to\C^\times
\eeq
Множество $G^\bullet$ всех голоморфных характеров на $G$ является
топологической группой относительно поточечной операции умножения и топологии
равномерной сходимости на компактах. Следующая теорема показывает, что операция
$G\mapsto G^\bullet$ аналогична понтрягинской операции перехода к двойственной
локально компактной абелевой группе:

\btm\label{PONT-AB} Если $G$ -- абелева компактно порожденная группа Штейна, то
ее двойственная группа $G^\bullet$ -- тоже абелева компактно порожденная группа
Штейна, а отображение
$$
\i_G:G\to G^{\bullet\bullet},\quad \i_G(x)(\chi)=\chi(x),\qquad x\in G,\;
\chi\in G^\bullet
$$
является естественным изоморфизмом функторов $G\mapsto G$ и $G\mapsto
G^{\bullet\bullet}$:
$$
G^{\bullet\bullet}\cong G
$$
\etm

Ввиду этой теоремы мы называем $G^\bullet$ {\it двойственной комплексной
группой к группе}\index{двойственная комплексная группа} $G$.

 \bpr
Сначала нужно заметить, что это верно для частных случаев $G=\C, \;\C^\times,\;
\Z$ и случая конечной абелевой группы $G=F$. Это следует из очевидных формул:
$$
\C^\bullet\cong \C,\quad (\C^\times)^\bullet\cong \Z,\quad \Z^\bullet\cong
\C^\times,\quad F^\bullet\cong F
$$
После этого останется заметить, что всякая абелева компактно порожденная группа
Штейна имеет вид
$$
G\cong \C^l\times (\C^\times)^m\times \Z^n\times F\qquad (l,m,n\in\Z_+)
$$
и поэтому ее двойственная группа будет иметь вид
$$
G^\bullet\cong \C^l\times \Z^m\times (\C^\times)^n\times F\qquad (l,m,n\in\Z_+)
$$
то есть тоже будет абелевой компактно порожденной группой Штейна, а вторая
двойственная группа $G^{\bullet\bullet}$ получается изоморфной $G$:
$$
G^{\bullet\bullet}\cong \C^l\times (\C^\times)^m\times \Z^n\times F\cong G
$$
 \epr

\btm Конструкция \eqref{G^bullet-v-kompl-geom} определяет обобщение теории двойственности
Понтрягина с категории конечных абелевых групп на категорию
компактно порожденных абелевых групп Штейна: {\sf
 $$
 \xymatrix
 {
 \boxed{\begin{matrix}
 \text{абелевы}\\
 \text{компактно порожденные группы Штейна}
 \end{matrix}}
 \ar[rr]^{G\mapsto G^\bullet} & &
 \boxed{\begin{matrix}
 \text{абелевы}\\
 \text{компактно порожденные группы Штейна}
 \end{matrix}}
 \\ & & \\
 \boxed{\begin{matrix}
 \text{абелевы}\\
 \text{конечные группы}
 \end{matrix}} \ar[uu]^{\mathfrak{e}} \ar[rr]^{G\mapsto G^\circ} & &
  \boxed{\begin{matrix}
 \text{абелевы}\\
 \text{конечные группы}
\end{matrix}}\ar[uu]_{\mathfrak{e}}
 }
 $$ }
и коммутативность этой диаграммы обеспечивается просто совпадением функторов на категории конечных абелевых групп:
$$
G^\circ = G^\bullet.
$$
 \etm

\section{Функции экспоненциального типа на комплексной группе Ли} \label{SEC-O_exp(G)}

\subsection{Полухарактеры и обратные полухарактеры}

Пусть $G$ -- комплексная группа Ли. Тогда
 \bit{
\item[---] локально ограниченная функция $f:G\to[1,+\infty)$ называется {\it
полухарактером},\index{полухарактер} если она удовлетворяет неравенству
субмультипликативности:
 \beq
f(x\cdot y)\le f(x)\cdot f(y),\qquad x,y\in G
 \eeq

\item[---] локально отделенная от нуля функция $\ph:G\to(0;1]$ называется {\it
обратным полухарактером},\index{полухарактер!обратный} если она удовлетворяет
обратному неравенству:
 \beq
\ph(x)\cdot \ph(y)\le\ph(x\cdot y),\qquad x,y\in G
 \eeq
 }\eit
Ясно, что если $f:G\to[1,+\infty)$ -- полухарактер, то обратная функция
 \beq\label{poluharakter-obratnyi-poluharakter}
\ph(x)=\frac{1}{f(x)}
 \eeq
будет обратным полухарактером, и наоборот.

\bigskip

\centerline{\bf Операции над полухарактерами и обратными
полухарактерами:}\label{PR-f+g}

 \bit{
\item[(i)] Множество всех полухарактеров на $G$ устойчиво относительно
следующих операций:
 \bit{
\item[---] умножение на достаточно большую константу: $C\cdot f$  $\quad (C\ge
1)$,

\item[---] умножение: $ f\cdot g,$

\item[---] сложение: $ f+g,$

\item[---] взятие максимума: $\max\{f,g\}.$
 }\eit

\item[(ii)] Множество всех обратных полухарактеров на $G$ устойчиво
относительно следующих операций:
 \bit{
\item[---] умножение на достаточно малую константу: $C\cdot \ph$ $\quad (C\le
1)$,

\item[---] умножение: $\ph\cdot\psi$,

\item[---] взятие половинки от среднего гармонического:
$\frac{\ph\cdot\psi}{\ph+\psi}$

\item[---] взятие минимума: $ \min\{\ph,\psi\}.$
 }\eit
 }\eit
 \bpr В силу естественной связи между полухарактерами и обратными полухарактерами,
обеспечиваемой формулой \eqref{poluharakter-obratnyi-poluharakter}, нам
достаточно рассмотреть только случай полухарактеров.

Если $f$ -- полухарактер на $G$ и $C\ge 1$, то
$$
C\cdot f(x\cdot y)\le C\cdot f(x)\cdot f(y)\le \big(C\cdot f(x)\big)\cdot
\big(C\cdot f(y)\big)
$$
Если $f$ и $g$ -- полухарактеры на $G$, то для произведения получаем:
 \begin{multline*}
(f\cdot g)(x\cdot y)=f(x\cdot y)\cdot g(x\cdot y)\le f(x)\cdot f(y)\cdot
g(x)\cdot g(y)=\\=f(x)\cdot g(x)\cdot f(y)\cdot g(y) =(f\cdot g)(x)\cdot
(f\cdot g)(y)
 \end{multline*}
Для суммы:
 \begin{multline*}
(f+g)(x\cdot y)=f(x\cdot y)+g(x\cdot y)\le f(x)\cdot f(y)+ g(x)\cdot g(y)\le\\
\le f(x)\cdot f(y)+g(x)\cdot f(y)+f(x)\cdot g(y)+ g(x)\cdot
g(y)=\Big(f(x)+g(x)\Big)\cdot \Big(f(y)+g(y)\Big) =(f+g)(x)\cdot (f+g)(y)
 \end{multline*}
А для $\max\{f,g\}$ доказательство основывается на следующем очевидном
неравенстве:
 \beq\label{max}
\max\{a\cdot b, c\cdot d\}\le \max\{a, c\}\cdot \max\{b, d\}\qquad (a,b,c,d>0)
 \eeq
Из него мы получаем:
 \begin{multline*}
\max\{f,g\}(x\cdot y)=\max\{f(x\cdot y), g(x\cdot y)\}\le \max\{f(x)\cdot f(y),
g(x)\cdot g(y)\}\le\eqref{max}\le\\ \le \max\{f(x), g(x)\}\cdot
\max\{f(y),g(y)\}= \max\{f, g\}(x)\cdot \max\{f,g\}(y)
 \end{multline*}

 \epr

Пусть, для произвольных множеств $G$ и $H$ и функций $g:G\to\C$ и $h:H\to\C$ символ $g\boxdot h$ обозначает функцию на декартовом произведении $G\times H$, определенную равенством
 \beq\label{g-h-na-S-T}
(g\boxdot h)(s,t):=g(s)\cdot h(t),\qquad s\in G,\quad t\in H.
 \eeq

\medskip
\centerline{\bf Свойства полухарактеров:}

\bit{\it

\item[$1^\circ$.] Если $f$ --- полухарактер на $G$, то для всякой константы $C\ge 1$ функция $C\cdot f$ --- тоже полухарактер на $G$;

\item[$2^\circ$.]\label{f-in-SC=>f^sigma-in-SC} Если $f$ --- полухарактер на $G$, то функция
$$
f^\sigma(t)=f(t^{-1}),\qquad t\in G,
$$
 --- тоже полухарактер на $G$;

\item[$3^\circ$.] Если $f$ и $g$ --- полухарактеры на $G$, то функции $f+g$, $f\cdot g$, $\max\{f,g\}$ --- тоже полухарактеры на $G$;

\item[$4^\circ$.]\label{f,g-in-SC=>f-boxdot-g-in-SC} Если $f$ и $g$ --- полухарактеры на группах $G$ и $H$, то функция $f\boxdot g$ --- полухарактер на группе $G\times H$;

\item[$5^\circ$.]\label{f-in-SC=>f^Delta-in-SC} Если $f$ --- полухарактер на декартовом квадрате  $G\times G$ группы $G$, то функция
$$
f^\varDelta(t)=f(t,t),\qquad t\in G,
$$
 --- полухарактер на $G$.

}\eit

\bex Пусть $\|\cdot\|$ --- субмультипликативная матричная норма на $\M_n(\C)$  \cite{Horn-Johnson}, например
 \beq\label{matrichnye-normy}
\|x\|=\sum_{i,j=1}^n|x_{i,j}|,\quad \|x\|=\sqrt{\sum_{i,j=1}^n|x_{i,j}|^2},\quad
\|x\|=\sup_{\xi\in\C^n: \|\xi\|\le 1}\|x(\xi)\|
 \eeq
(где $\|\xi\|=\sqrt{\sum_{i=1}^n |\xi_i|^2}$). Тогда для любых $C\ge 1$ и $N\in\N$ отображение $r^N_C:\GL_n(\C)\to\R$
 \beq\label{r_N(x)=||x||^N-cdot-||x^-1||^N}
r^N_C(x)=C\cdot\max\{\|x\|; \|x^{-1}\|\}^N,\qquad C\ge 1,\; N\in\N
 \eeq
является полухарактером на $\GL_n(\C)$.
\eex
\bpr
Здесь нужно только проверить неравенство
\beq\label{max|x|;|x^-1|-ge-1}
\max\{\|x\|; \|x^{-1}\|\}\ge 1.
\eeq
Оно следует из невырожденности и мультипликативности нормы $\|\cdot\|$. Во-первых, из невырожденности нормы $\|\cdot\|$ мы получаем $\|1\|>0$, откуда
$$
\|1\|=\|1\cdot 1\|\le \|1\|\cdot \|1\|\quad\Rightarrow\quad \|1\|\ge 1
$$
А, во-вторых,
$$
1\le \|1\|=\|x\cdot x^{-1}\|\le \|x\|\cdot \|x^{-1}\|\quad\Rightarrow\quad \eqref{max|x|;|x^-1|-ge-1}.
$$
\epr

\bprop\label{PROP-stroenie-poluharakterov-na-GL_n} Для любой
субмультипликативной матричной нормы $||\cdot||$ на $\GL_n(\C)$ полухарактеры
вида \eqref{r_N(x)=||x||^N-cdot-||x^-1||^N} мажорируют все остальные
полухарактеры на $\GL_n(\C)$. \eprop

\bpr Заметим сразу, что нам достаточно рассмотреть случай, когда $||\cdot||$ --
евклидова норма в алгебре $\M_n(\C)$ всех комплексных матриц $n\times n$:
$$
||x||=\sup_{\xi\in\C^n: ||\xi||\le 1}||x(\xi)||,\qquad\text{где}\quad
||\xi||=\sqrt{\sum_{i=1}^n |\xi_i|^2},
$$
-- поскольку любая другая норма на $\M_n(\C)$ мажорирует евклидову норму с
точностью до константы, это будет доказывать наше утверждение.

Рассмотрим множество
$$
K:=\big\{x\in\GL_n(\C):\; \max\{||x||; ||x^{-1}||\}\le 2 \big\}=
\big\{x\in\GL_n(\C):\;\forall\xi\in\C^n\quad
\frac{1}{2}\cdot||\xi||\le||x(\xi)||\le 2\cdot||\xi|| \big\}
$$
Оно замкнуто и ограничено в алгебре матриц $\M_n(\C)$, и поэтому компактно. Это
будет порождающий компакт в $\GL_n(\C)$
 \beq\label{GL_n(C)=bigcup}
\GL_n(\C)=\bigcup_{m=1}^\infty K^m,\qquad K^m=\underbrace{K\cdot...\cdot
K}_{\text{$m$ множителей}},
 \eeq
обладающий к тому же следующим свойством:
 \beq\label{K^m:=max-||x||^m;||x^-1||^m-le-2^m}
K^m=\big\{x\in\GL_n(\C):\; \max\{||x||; ||x^{-1}||\}\le 2^m \big\}
 \eeq
Действительно, если $x\in K^m$, то $x=x_1\cdot...\cdot x_m$, где $||x_i||\le 2$
и $||x_i^{-1}||\le 2$, и поэтому
$$
||x||\le \prod_{i=1}^m ||x_i||\le 2^m,\qquad ||x^{-1}||\le \prod_{i=1}^m
||x_i^{-1}||\le 2^m
$$
Наоборот, пусть $\max\{||x||; ||x^{-1}||\}\le 2^m$. Рассмотрим полярное
разложение: $x=r\cdot u$, где $r$ -- положительно определенная эрмитова
матрица, а $u$ -- унитарная. Разложим $r$ в произведение $r=v\cdot d\cdot
v^{-1}$, где $v$ -- унитарная, а $d$ -- диагональная матрицы:
$$
d=\begin{pmatrix}d_1 & 0 & ... & 0 \\ 0 & d_2 & ... & 0 \\ ... & ... & ... &
... \\ 0 & 0 & ... & d_n \\
 \end{pmatrix},\qquad d_i\ge 0
$$
Корень степени $m$
$$
\root{m}\of{d}=\begin{pmatrix}\root{m}\of{d_1} & 0 & ... & 0 \\ 0 &
\root{m}\of{d_2} & ... & 0
\\ ... & ... & ... &
... \\ 0 & 0 & ... & \root{m}\of{d_n} \\
 \end{pmatrix}
$$
обладает такими свойствами:
$$
\max_{i=1,...,n} d_i=||d||=||r||=||x||\le 2^m\quad\Longrightarrow\quad
||\root{m}\of{d}||=\max_{i=1,...,n} \root{m}\of{d_i}\le 2,
$$
$$
\max_{i=1,...,n} d_i^{-1}=||d^{-1}||=||r^{-1}||=||x^{-1}||\le
2^m\quad\Longrightarrow\quad ||\root{m}\of{d^{-1}}||=\max_{i=1,...,n}
\root{m}\of{d_i^{-1}}\le 2
$$
Как следствие, матрица
$$
y=v\cdot \root{m}\of{d}\cdot v^{-1}
$$
лежит в $K$:
$$
\max\{||y||; ||y^{-1}||\}=\max\{||\root{m}\of{d}||;
||\root{m}\of{d}^{-1}||\}\le 2 \quad\Longrightarrow\quad y\in K
$$
Поэтому матрица $y\cdot u$ также лежит в $K$, и мы получаем
$$
x=v\cdot d\cdot v^{-1}\cdot u=v\cdot (\root{m}\of{d})^{m-1}\cdot v^{-1}\cdot
v\cdot (\root{m}\of{d})\cdot v^{-1}\cdot u= \underbrace{y^{m-1}}_{\scriptsize
\begin{matrix}\text{\rotatebox{-90}{$\in$}} \\ K^{m-1}\end{matrix}}\cdot \underbrace{y\cdot v}_{\scriptsize
\begin{matrix}\text{\rotatebox{-90}{$\in$}} \\ K\end{matrix}}\in
K^m
$$
Мы доказали формулу \eqref{K^m:=max-||x||^m;||x^-1||^m-le-2^m}. Теперь пусть
$f$ -- произвольный полухарактер на $\GL_n(\C)$. Положим
$$
C=\sup_{x\in K}f(x),\qquad N\ge \log_2 C
$$
и покажем, что
 \beq
\forall x\in\GL_n(\C)\qquad f(x)\le r^N_C(x)
 \eeq
Зафиксируем $x\in\GL_n(\C)$. Из \eqref{GL_n(C)=bigcup} следует, что $x\in
K^m\setminus K^{m-1}$ для некоторого $m\in\N$. По формуле
\eqref{K^m:=max-||x||^m;||x^-1||^m-le-2^m} получим:
$$
2^{m-1}<\max\{||x||; ||x^{-1}||\}\le 2^m
$$
$$
\Downarrow
$$
$$
m-1<\log_2 \Big(\max\{||x||; ||x^{-1}||\}\Big)
$$
$$
\Downarrow
$$
 \begin{multline*}
f(x)\le\sup_{y\in K}f(y)^m\le C^m\le C\cdot C^{m-1}<C\cdot C^{\log_2
\Big(\max\{||x||; ||x^{-1}||\}\Big)}=\\ =C\cdot \Big(\max\{||x||;
||x^{-1}||\}\Big)^{\log_2 C} \le C\cdot \Big(\max\{||x||;
||x^{-1}||\}\Big)^N=r^N_C(x)
 \end{multline*}
\epr

В частном случае, когда $n=1$ мы получаем:

\bcor На комплексной окружности $\C^\times$ полухарактеры вида
$$
r_C^N(t)=C\cdot \max\{|t|;|t|^{-1}\}^N,\qquad C\ge 1,\; N\in\N
$$
мажорируют все остальные полухарактеры. \ecor

\bprop\label{PROP-poluh-porozhd-K} Если $G$ -- компактно порожденная группа
Штейна и $K$ -- какая-нибудь компактная окрестность единицы порождающая $G$:
 $$
G=\bigcup_{n=1}^\infty K^n,\qquad K^n=\underbrace{K\cdot...\cdot K}_{\text{$n$
множителей}},
 $$
то для всякого $C\ge 1$ правило
 \beq\label{poluh-opred-K}
h_C(x)=C^n\quad\Longleftrightarrow\quad x\in K^n\setminus K^{n-1}
 \eeq
определяет полухарактер $h_C$ на $G$. Такие полухарактеры образуют
фундаментальную систему среди всех полухарактеров на $G$ в том смысле, что
всякий полухарактер $f$ на $G$ будет подчинен некоторому полухарактеру $h_C$
$$
f(x)\le h_C(x),\qquad x\in G,
$$
-- для этого константу $C$ достаточно выбрать условием
 \beq\label{C-ge-max-t-in-K-f(t)}
C\ge\sup_{t\in K}f(t)
 \eeq
 \eprop
\bpr Локальная ограниченность $h_C$ очевидна, нужно проверить
субмультипликативность. Пусть $x,y\in G$, подберем $k,l\in\N$ такие что
$$
x\in K^k\setminus K^{k-1},\qquad y\in K^l\setminus K^{l-1}.
$$
Тогда $x\cdot y\in K^{k+l}$, и поэтому
$$
h_C(x\cdot y)\le C^{k+l}=C^k\cdot C^l=h_C(x)\cdot h_C(y)
$$

Если теперь $f$ -- произвольный полухарактер, и $C$ подчинена условию
\eqref{C-ge-max-t-in-K-f(t)}, то
$$
x\in K^n\setminus K^{n-1} \quad\Longrightarrow\quad f(x)\le \Big(\sup_{t\in K}
f(t)\Big)^n\le C^n=h_C(x)
$$
\epr

\subsection{Функции длины.}

{\it Функцией длины} на группе $G$ называется произвольное локально ограниченное отображение $\ell:G\to \R_+$, удовлетворяющее условию
\beq\label{DEF:func-dliny}
\ell(x\cdot y)\le \ell(x)+\ell(y).
\eeq
Пусть $G$ --- группа. Зафиксируем какое-нибудь множество $S\subseteq G$, порождающее $G$ как полугруппу:
$$
\forall x\in G\quad \exists a_{k_1},...,a_{k_j}\in S\qquad x=a_{k_1}\cdot ...\cdot a_{k_j}.
$$
Для всякой функции
$$
F:S\to\R_+
$$
положим
\beq\label{DEF:ell_F}
\ell_F(x)=\inf\left\{ \sum_{i=1}^m F(a_i);\ a_1,...,a_m\in S: \ x=a_1\cdot ...\cdot a_m \right\}.
\eeq

\bprop
Отображение $\ell_F:G\to\R_+$ является функцией длины:
\beq\label{ell_F(xy)-le-ell_F(x)+ell_F(y)}
\ell_F(x\cdot y)\le \ell_F(x)+\ell_F(y).
\eeq
\eprop
\bpr
Если
$x=a_1\cdot ...\cdot a_m$, $a_i\in S$, и $y=b_1\cdot ...\cdot b_n$, $b_j\in S$, то
$x\cdot y=a_1\cdot ...\cdot a_m\cdot b_1\cdot ...\cdot b_n$, поэтому
$$
\ell_F(x\cdot y)\le \sum_{i=1}^m F(a_i)+\sum_{i=1}^n F(b_n).
$$
Это верно для любых представлений $x=a_1\cdot ...\cdot a_m$, $y=b_1\cdot ...\cdot b_n$, с $a_i,b_j\in S$. Значит,
$$
\ell_F(x\cdot y)\le \ell_F(x)+\ell_F(y).
$$
\epr

Из \eqref{ell_F(xy)-le-ell_F(x)+ell_F(y)} следует, что для всякой функции $F:S\to\R_+$ функция
$$
f_F(x)=e^{\ell_F(x)}
$$
является полухарактером на $G$.

\bprop\label{PROP:poluhar-mazhor-f_F^C}
Любой полухарактер $f:G\to[1,+\infty)$ мажорируется некоторым полухарактером $f_F:G\to[1,+\infty)$:
$$
f(x)\le f_F(x)=e^{\ell_F(x)},\quad x\in G.
$$
\eprop
\bpr
Рассмотрим функцию
$$
F(a)=\log f(a),\qquad a\in S.
$$
Если $x=a_1\cdot...\cdot a_m$, с $a_i\in S$, то
$$
\log f(x)=\log f(a_1\cdot...\cdot a_m)\le \log \left(f(a_1)\cdot...\cdot f(a_m)\right)=\log f(a_1)+...+ \ln f(a_m)=F(a_1)+...+F(a_m).
$$
Это верно для всякого представления $x=a_1\cdot...\cdot a_m$, с $a_i\in S$, поэтому
$$
\log f(x)\le \ell_F(x).
$$
\epr

\brem\label{REM:neprodolzhaemyi-poluharakter}
Полухарактер не всегда продолжается с подгруппы на группу и не всегда даже можно найти продолжаемый полухарактер, мажорируюший данный: если $h:H\to[1,\infty)$ -- полухарактер на подгруппе $H$ группы $G$, то необязательно существует полухарактер $g:G\to[1,\infty)$ такой, что
$$
h(x)\le g(x),\qquad x\in H.
$$

Контрпримером служит стандартная конструкция из метрической геометрии групп \cite[Remark 4.B.12]{Cornulier}: пусть $G$ --- дискретная группа Гейзенберга, то есть группа целочисленных матриц вида
$$
X=\begin{pmatrix}1& a& c \\ 0& 1& b\\ 0& 0& 1\end{pmatrix},\qquad a,b,c\in\Z.
$$
Отображение
$$
\ph:\Z\to G\quad\Big|\quad \ph(n)=\begin{pmatrix}1& 0& n \\ 0& 1& 0\\ 0& 0& 1\end{pmatrix}
$$
является инъективным гомоморфизмом (вложением) групп. Мы утверждаем, что при таком вложении полухарактер на $\Z$
$$
h(n)=2^{\abs{n}}
$$
не мажорируется никаким полухарактером на $G$: неравенство
\beq\label{neprodolzhaemyi-poluharakter-0}
h(n)\le g(\ph(n)),\qquad n\in\Z
\eeq
не выпоняется ни для какого $g\in\sf{SC}(G)$.

Чтобы это понять, нужно заметить, что матрицы
$$
A=\begin{pmatrix}1& 1& 0 \\ 0& 1& 0\\ 0& 0& 1\end{pmatrix},\qquad
A^{-1}=\begin{pmatrix}1& -1& 0 \\ 0& 1& 0\\ 0& 0& 1\end{pmatrix},\qquad
B=\begin{pmatrix}1& 0& 0 \\ 0& 1& 1\\ 0& 0& 1\end{pmatrix},\qquad
B^{-1}=\begin{pmatrix}1& 0& 0 \\ 0& 1& -1\\ 0& 0& 1\end{pmatrix}
$$
порождают $G$ как полугруппу. В частности, они порождают матрицу с единицей в правом верхнем углу по формуле
\beq
\begin{pmatrix}1& 0& 1 \\ 0& 1& 0\\ 0& 0& 1\end{pmatrix}=
\begin{pmatrix}1& 1& 0 \\ 0& 1& 0\\ 0& 0& 1\end{pmatrix}\cdot \begin{pmatrix}1& 0& 0 \\ 0& 1& 1\\ 0& 0& 1\end{pmatrix}=A\cdot B.
\eeq
Из того, что множество $S=\{A,A^{-1},B,B^{-1}\}$ порождает $G$ как полугруппу, по предложению \ref{PROP:poluhar-mazhor-f_F^C} следует, что всякий полухарактер $g$ на $G$ подчинен некоторому полухарактеру вида $f_F$.
\beq\label{g(X)-le-f_F(X)}
g(X)\le f_F(X),\qquad X\in G.
\eeq
С другой стороны, поскольку система образующих $S$ конечна, функцию $F$ можно считать постоянной
$$
F(s)=C,\qquad s\in S.
$$
Тогда
\beq\label{neprodolzhaemyi-poluharakter-1}
f_F(X)=e^{\ell_F(X)}=e^{\inf\left\{ C\cdot m;\ A_1,...,A_m\in S: \ X=A_1\cdot ...\cdot A_m \right\}}
=e^{ C\cdot \inf\left\{m;\ A_1,...,A_m\in S: \ X=A_1\cdot ...\cdot A_m \right\}}
\eeq
Теперь заметим тождество
\begin{multline}\label{neprodolzhaemyi-poluharakter-2}
\ph(n^2)=\begin{pmatrix}1& 0& n^2 \\ 0& 1& 0\\ 0& 0& 1\end{pmatrix}=
\begin{pmatrix}1& 0& 0 \\ 0& 1& -n\\ 0& 0& 1\end{pmatrix}\cdot
\begin{pmatrix}1& n& 0 \\ 0& 1& 0\\ 0& 0& 1\end{pmatrix}\cdot \begin{pmatrix}1& 0& 0 \\ 0& 1& n\\ 0& 0& 1\end{pmatrix}\cdot \begin{pmatrix}1& -n& 0 \\ 0& 1& 0\\ 0& 0& 1\end{pmatrix}=\\=
\begin{pmatrix}1& 0& 0 \\ 0& 1& -1\\ 0& 0& 1\end{pmatrix}^n\cdot
\begin{pmatrix}1& 1& 0 \\ 0& 1& 0\\ 0& 0& 1\end{pmatrix}^n\cdot \begin{pmatrix}1& 0& 0 \\ 0& 1& 1\\ 0& 0& 1\end{pmatrix}^n\cdot \begin{pmatrix}1& -1& 0 \\ 0& 1& 0\\ 0& 0& 1\end{pmatrix}^n=B^{-n}\cdot A^n\cdot B^n\cdot A^{-n}.
\end{multline}
Оно дает цепочку
\begin{multline*}
2^{n^2}=h(n^2)\le\eqref{neprodolzhaemyi-poluharakter-0}\le g(\ph(n^2))\le
\eqref{g(X)-le-f_F(X)}\le f_F(\ph(n^2))
= \\ =
\eqref{neprodolzhaemyi-poluharakter-2}= f_F(B^{-n}\cdot A^n\cdot B^n\cdot A^{-n}) \le\eqref{neprodolzhaemyi-poluharakter-1}\le e^{4n\cdot C},\qquad n\in\N.
\end{multline*}
Понятно, что ни при каком $C$ такое невозможно.
\erem

\subsection{Субмультипликативные ромбы и дуально субмультипликативные прямоугольники}
\label{SUBSEC:submultiplik-romb}

Введем следующие определения:

 \bit{
\item[---] замкнутую абсолютно выпуклую окрестность нуля $\varDelta$ в
$\mathcal O^\star(G)$ мы будем называть {\it
субмультипликативной},\index{субмультипликативная окрестность нуля} если для
любых функционалов $\alpha,\beta$ из $\varDelta$ их свертка $\alpha*\beta$ тоже
лежит в $\varDelta$:
$$
\forall\alpha,\beta\in\varDelta\qquad \alpha*\beta\in\varDelta
$$
коротко это изображается вложением:
$$
\varDelta * \varDelta\subseteq \varDelta
$$

\item[---] абсолютно выпуклый компакт $D\subseteq \mathcal O(G)$ мы будем
называть {\it дуально субмультипликативным},\index{дуально субмультипликативный
компакт} если его поляра $D^\circ$ является субмультипликативной окрестностью
нуля:
$$
 D^\circ * D^\circ\subseteq D^\circ
$$
}\eit

\blm \label{f-obr-poluh->F-dual-subm}$\phantom{.}$
 \bit{
\item[(a)] Если $\ph:G\to(0;1]$ -- обратный полухарактер на $G$, то порождаемый
им ромб $\ph^{\BLZ}$ является (замкнутой, абсолютно выпуклой и)
субмультипликативной окрестностью нуля в $\mathcal O^\star(G)$.

\item[(b)] Если $\varDelta\subseteq \mathcal O^\star(G)$ -- замкнутая абсолютно
выпуклая субмультипликативная окрестность нуля в $\mathcal O^\star(G)$, то ее
внутренняя огибающая $\varDelta^{\LZ}$ -- обратный полухарактер на $G$. }\eit
\elm
 \bpr
1. Обозначим
$$
\e_x=\ph(x)\cdot\delta^x\in \mathcal{O}^\star(G),
$$
тогда
 \beq
\ph^{\BLZ} = \cabsconv \left\{\ph(x)\cdot\delta^x;\quad x\in
M\right\}=\cabsconv \left\{\e_x;\quad x\in M\right\}
 \eeq
и вложение
$$
\ph^{\BLZ} * \ph^{\BLZ}\subseteq \ph^{\BLZ}
$$
проверяется в три этапа. Сначала нужно заметить, что
$$
 \forall x,y\in G\qquad \e_x*\e_y\in \ph^{\BLZ}
$$
Действительно,
 $$
 \e_x*\e_y=\ph(x)\cdot\delta^x * \ph(y)\cdot\delta^y=
 \ph(x)\cdot \ph(y)\cdot\delta^{x\cdot y}=
 \underbrace{\frac{\ph(x)\cdot \ph(y)}{\ph(x\cdot y)}}_{\scriptsize \begin{matrix}\text{\rotatebox{-90}{$\le$}} \\ 1\end{matrix}}
 \cdot\e_{x\cdot y}\in
 \cabsconv \left\{\e_z;\quad z\in M\right\}=\ph^{\BLZ}
 $$
Затем берутся конечные абсолютно выпуклые комбинации функционалов $\e_x$:
 \beq\label{abs-vyp-komb}
\alpha=\sum_{i=1}^k \lambda_i\cdot \e_{x_i},\quad \beta=\sum_{j=1}^l \mu_j\cdot
\e_{y_j}\qquad \left( \sum_{i=1}^k |\lambda_i|\le 1,\quad \sum_{j=1}^l
|\mu_j|\le 1 \right)
 \eeq
для которых получается
 $$
\alpha * \beta=\sum_{i=1}^k \lambda_i\cdot \e_{x_i} * \sum_{j=1}^l \mu_j\cdot
\e_{y_j}=\sum_{1\le i\le k,1\le j\le l}
 \underbrace{\lambda_i\cdot\mu_j}_{\scriptsize
\begin{matrix}\uparrow \\ \sum_{1\le i\le k,1\le j\le l} |\lambda_i\cdot\mu_j| \\ \text{\rotatebox{-90}{$=$}}
\\ \sum_{i=1}^k |\lambda_i|\cdot \sum_{j=1}^l |\mu_j| \\
\text{\rotatebox{-90}{$\le$}}\\ 1
 \end{matrix}}\cdot
\underbrace{\e_{x_i} *\e_{y_j}}_{\scriptsize
\begin{matrix}\text{\rotatebox{-90}{$\in$}} \\ \\ \ph^{\BLZ} \end{matrix}}\in
 \underbrace{\ph^{\BLZ}}_{\scriptsize\begin{matrix} \text{абсолютно} \\
\text{выпуклое} \\
\text{множество}\end{matrix}}
 $$
И, наконец, для произвольных $\alpha,\beta\in \ph^{\BLZ}$ включение
$\alpha*\beta\in \ph^{\BLZ}$ получается как следствие плотности функционалов
вида \eqref{abs-vyp-komb} в множестве $\ph^{\BLZ} = \cabsconv \left\{\e_x;\quad
x\in M\right\}$ и непрерывности операции свертки.

2. Здесь используется формула \eqref{fi^lozenge}:
 \begin{multline*}
\varDelta * \varDelta\subseteq \varDelta\quad\Longrightarrow\quad \forall
x,y\in G\qquad \underbrace{\varDelta^{\LZ}(x)\cdot\delta^x}_{\scriptsize
\begin{matrix}\text{\rotatebox{-90}{$\in$}} \\ \varDelta, \\ \text{в силу \eqref{fi^lozenge}}\end{matrix}} *
\underbrace{\varDelta^{\LZ}(y)\cdot\delta^y}_{\scriptsize
\begin{matrix}\text{\rotatebox{-90}{$\in$}} \\ \varDelta, \\ \text{в силу \eqref{fi^lozenge}}\end{matrix}}
=\varDelta^{\LZ}(x)\cdot\varDelta^{\LZ}(y)\cdot\delta^{x\cdot y} \in
 \varDelta \quad\Longrightarrow\\ \Longrightarrow\quad
\varDelta^{\LZ}(x)\cdot\varDelta^{\LZ}(y)\le \max \{\lambda>0:\;
\lambda\cdot\delta^{x\cdot y}\in \varDelta
\}=\eqref{fi^lozenge}=\varDelta^{\LZ}(x\cdot y)
\quad\Longrightarrow\quad \varDelta^{\LZ}(x)\cdot\varDelta^{\LZ}(y)\le
\varDelta^{\LZ}(x\cdot y)
 \end{multline*}
 \epr

\blm \label{f-poluh->F-subm}$\phantom{.}$
 \bit{
\item[(a)] Если $f:G\to[1;\infty)$ -- полухарактер на $G$, то порождаемый им
прямоугольник $f^\text{\BSQ}$ дуально субмультипликативен.

\item[(b)] Если $D\subseteq \mathcal O(G)$ -- дуально субмультипликативный абсолютно выпуклый компакт, то его внешняя огибающая $D^\text{\SQ}$ -- полухарактер на $G$. }\eit \elm
 \bpr
1. Если $f:G\to[1;\infty)]$ -- полухарактер, то $\frac{1}{f}$ -- обратный
полухарактер, поэтому по лемме \ref{f-obr-poluh->F-dual-subm} (a) ромб
$$
\left(\frac{1}{f}\right)^{\BLZ}=\eqref{f^blacksquare^circ=frac-1-f^blacklozenge}=
(f^\text{\BSQ})^\circ
$$
будет субмультипликативной окрестностью нуля. Значит $f^\text{\BSQ}$ -- дуально
субмультипликативный прямоугольник.

2. Если $D\subseteq \mathcal O(G)$ -- дуально субмультипликативный абсолютно выпуклый компакт, то его поляра $D^\circ\subseteq \mathcal O^\star(G)$ --
субмультипликативная окрестность нуля в $\mathcal O^\star(G)$, поэтому по лемме
\ref{f-obr-poluh->F-dual-subm} (b) внутренняя огибающая
$$
(D^\circ)^{\LZ}=\eqref{D^circ^lozenge=frac-1-D^square}=\frac{1}{D^\text{\SQ}}
$$
-- обратный полухарактер. Значит $D^\text{\SQ}$ -- полухарактер.
 \epr

Леммы \ref{f-obr-poluh->F-dual-subm} и \ref{f-poluh->F-subm} вместе с формулами
$\varDelta=\varDelta^{{\LZ}{\BLZ}}$ и $D=D^\text{\SQ\BSQ}$ для ромбов и
прямоугольников дают следующую теорему:

 \btm\label{TH-prayam-poluharakter}$\phantom{.}$
 \bit{
\item[(a)] Ромб $\varDelta\in\mathcal O^\star(G)$ субмультипликативен тогда и
только тогда, когда его внутренняя огибающая $\varDelta^{\LZ}$ -- обратный
полухарактер на $G$.

\item[(b)] Прямоугольник $D\in\mathcal O(G)$ дуально субмультипликативен тогда
и только тогда, когда его внешняя огибающая $D^\text{\SQ}$ -- полухарактер на
$G$.

}\eit \etm

Следующий результат показывает, что субмультипликативные ромбы образуют
фундаментальную систему среди всех субмультипликативных замкнутых абсолютно
выпуклых окрестностей нуля в $\mathcal O^\star(G)$:

\btm\label{subm-in-trubk} $\phantom{.}$
 \bit{
\item[(a)] Всякая замкнутая абсолютно выпуклая субмультипликативная окрестность
нуля $\varDelta$ в $\mathcal O^\star(G)$ содержит некоторый
субмультипликативный ромб, а именно ромб $\varDelta^{{\LZ}\kern-0.5pt{\BLZ}}$.

\item[(b)] Всякое дуально субмультипликативное ограниченное множество
$D\subseteq \mathcal O(G)$ содержится в некотором дуально субмультипликативном
прямоугольнике, а именно, в прямоугольнике $D^{\text{\SQ}\text{\BSQ}}$.
 }\eit
 \etm
 \bpr
1. Если $\varDelta$ -- замкнутая абсолютно выпуклая субмультипликативная
окрестность нуля в $\mathcal O^\star(G)$, то по лемме
\ref{f-obr-poluh->F-dual-subm}(b), $\varDelta^{\LZ}$ -- обратный полухарактер,
и значит, по лемме \ref{f-obr-poluh->F-dual-subm}(a), ромб
$\varDelta^{{\LZ}\kern-0.5pt{\BLZ}}$ будет субмультипликативным.

2. Если $D$ -- дуально субмультипликативное ограниченное множество в $\mathcal
O(G)$, то его поляра $D^\circ$ -- замкнутая абсолютно выпуклая
субмультипликативная окрестность нуля в $\mathcal O^\star(G)$, и по уже
доказанному $(D^\circ)^{{\LZ}\kern-0.5pt{\BLZ}}$ -- субмультипликативный ромб.
То есть множество
$(D^{\text{\SQ}\text{\BSQ}})^\circ=(D^\circ)^{{\LZ}\kern-0.5pt{\BLZ}}$
субмультипликативно. Значит, прямоугольник $D^{\text{\SQ}\text{\BSQ}}$ дуально
субмультипликативен.
 \epr

В соответствии с определениями \ref{rectangles-buses}, мы называем функцию $f$
на $G$ {\it огибающим полухарактером},\index{полухарактер!огибающий} если она
является внешней огибающей и одновременно полухарактером на $G$.

\btm\label{TH:f_N} Если $G$ -- компактно порожденная группа Штейна, то системы
всех полухарактеров, всех огибающих полухарактеров, всех дуально
субмультипликативных прямоугольников на $G$ и всех замкнутых абсолютно выпуклых субмультипликативных окрестностей нуля в $\mathcal O^\star(G)$ содержат счетные конфинальные
подсистемы:
 \bit{
\item[(i)] существует последовательность полухарактеров $h_N$ на $G$ такая, что
всякий полухарактер $g$ мажорируется некоторым полухарактером $h_N$:
$$
g(x)\le h_N(x),\quad x\in G
$$

\item[(ii)] существует последовательность огибающих полухарактеров $f_N$ на $G$
такая, что всякий огибающий полухарактер $g$ мажорируется некоторым
полухарактером $f_N$:
$$
g(x)\le f_N(x),\quad x\in G
$$

\item[(iii)] существует последовательность $E_N$ дуально субмультипликативных
прямоугольников на $G$ такая, что всякий дуально субмультипликативный
прямоугольник $D$ в $G$ содержится в некотором $ E_N$:
 $$
 D\subseteq E_N
 $$

\item[(iv)] существует последовательность $\Delta_N$  замкнутых абсолютно выпуклых субмультипликативных окрестностей нуля в $\mathcal O^\star(G)$ такая, что всякая замкнутая абсолютно выпуклая субмультипликативная окрестность нуля $U$ в $\mathcal O^\star(G)$  содержит некоторую окрестность $\Delta_N$:
 $$
 U\supseteq \Delta_N
 $$
 }\eit
\etm \bpr Это следует из предложения \ref{PROP-poluh-porozhd-K}:
последовательность полухарактеров $h_N$, $N\in\N$, определенных формулой
\eqref{poluh-opred-K}, будет искомой последовательностью полухарактеров на $G$.

А $E_N$, $f_N$ и $\varDelta_N$ определяются формулами
 \beq\label{f_N^black}
 E_N=\{u\in\mathcal O(G):\; \max_{x\in K^n}|u(x)|\le N^n\}=(h_N)^\text{\BSQ}
 \eeq
 \beq\label{f_N-f_N}
f_N=(E_N)^\text{\SQ}=(h_N)^{\text{\BSQ}\text{\SQ}}
 \eeq
 \beq\label{f_N-Delta_N}
\varDelta_N=\left(\frac{1}{h_N}\right)^{\BLZ}=\eqref{f^blacksquare^circ=frac-1-f^blacklozenge}=
(h_N^\text{\BSQ})^\circ
 \eeq
($E_N$ будут дуально субмультипликативными прямоугольниками по лемме
\ref{f-poluh->F-subm}(a), $f_N$ -- огибающими полухарактерами по лемме
\ref{f-poluh->F-subm}(b), а $\varDelta_N$ --- субмультипликативными окрестностями нуля по лемме \ref{f-obr-poluh->F-dual-subm}(a)).

Если теперь $D$ -- дуально субмультипликативный прямоугольник, то по лемме
\ref{f-poluh->F-subm}(b), его внешняя огибающая $D^\text{\SQ}$ будет
полухарактером, значит, по предложению \ref{PROP-poluh-porozhd-K}, найдется
$N\in\N$ такое, что
$$
D^\text{\SQ}\le h_N
$$
Отсюда получаем
$$
D=(D^\text{\SQ})^\text{\BSQ}\subseteq (h_N)^\text{\BSQ}=E_N
$$
то есть $D$ обязан содержаться в некотором $E_N$. Из этого в свою очередь
следует
$$
D^\text{\SQ}\le (E_N)^\text{\SQ}=f_N,
$$
и это можно понимать так, что всякий огибающий полухарактер (всегда имеющий вид
$D^\text{\SQ}$ по теореме \ref{TH-prayam-poluharakter}) мажорируется некоторым
$f_N$.

Наконец, если $U$ --- замкнутая абсолютно выпуклая субмультипликативная окрестность нуля в ${\mathcal O}^\star(G)$, то по теореме \ref{subm-in-trubk}(a) она должна содержать некоторый субмультипликативный ромб $\ph^{\BLZ}$
$$
U\supseteq \ph^{\BLZ},
$$
где $\ph=U^{\LZ}$ --- некий обратный полухарактер по лемме \ref{f-obr-poluh->F-dual-subm}(b). Функция $\frac{1}{\ph}$ является полухарактером на $G$. И по уже доказанному, она должна мажорироваться неким полухарактером $h_N$:
$$
\frac{1}{\ph(x)}\le h_N(x),\qquad x\in G.
$$
Отсюда мы получаем 
$$
\ph(x)\ge\frac{1}{h_N(x)}
$$
и, как следствие,
$$
U\supseteq\eqref{sv-0-int}\supseteq \ph^{\BLZ}\supseteq\l\frac{1}{h_N}\r^{\BLZ}=\eqref{f_N-Delta_N}=\varDelta_N
$$

 \epr

\subsection{Голоморфные функции экспоненциального типа}\label{SUBSEC:algebra-O_exp(G)}

\paragraph{Алгебра $\mathcal O_{\exp} (G)$ голоморфных функций
экспоненциального типа.}

Голоморфную функцию $u\in\mathcal O(G)$ на компактно порожденной группе Штейна
$G$ мы называем {\it функцией экспоненциального типа},\index{функция
экспоненциального типа} если она ограничивается некоторым полухарактером:
$$
|u(x)|\le f(x),\qquad x\in G\quad \Big(f(x\cdot y)\le f(x)\cdot f(y)\Big)
$$
Множество всех функций экспоненциального типа на $G$ мы обозначаем $\mathcal{O}_{\exp}(G)$. Это подпространство в $\mathcal{O}(G)$, и по теореме
\ref{subm-in-trubk}, его можно рассматривать как объединение дуально
субмультипликативных прямоугольников в $\mathcal{O}(G)$:
$$
\mathcal{O}_{\exp}(G)=\bigcup_{\scriptsize\begin{matrix}\text{$D$ --
дуально}\\ \text{субмультипликативный}\\ \text{прямоугольник в $\mathcal
O(G)$}\end{matrix}}
 D=\bigcup_{\text{$f$ -- полухарактер на $G$}} f^\text{\BSQ}
$$
или, что то же самое, как объединение подпространств вида $\C D$, где $D$ --
дуально субмультипликативный прямоугольник в $\mathcal O(G)$:
 \beq\label{O-exp}
\mathcal{O}_{\exp}(G)=\bigcup_{\scriptsize\begin{matrix}\text{$D$ --
дуально}\\ \text{субмультипликативный}\\ \text{прямоугольник в $\mathcal
O(G)$}\end{matrix}} \C D
 \eeq
Поэтому $\mathcal{O}_{\exp}(G)$ естественно наделяется топологией инъективного
(локально выпуклого) предела подпространств Смит $\C D$:
 \beq\label{O-exp-top}
\mathcal{O}_{\exp}(G)=\underset{\scriptsize\begin{matrix}\text{$D$ --
дуально}\\ \text{субмультипликативный}\\ \text{прямоугольник в $\mathcal
O(G)$}\end{matrix}}{\LCS\text{-}\injlim} \C D
 \eeq
Из теоремы \ref{TH:f_N} следует, что в этом пределе систему всех дуально
субмультипликативных прямоугольников можно заменить на счетную подсистему:
 \beq\label{O-exp-top-E_N}
\mathcal{O}_{\exp}(G)=\underset{\scriptsize N\to \infty}{\LCS\text{-}\injlim} \C E_N
 \eeq
Как следствие, справедлива

\btm\label{TH:O_exp--Brauner} Пространство $\mathcal{O}_{\exp}(G)$ функций экспоненциального типа на компактно порожденной группы Штейна $G$ является пространством Браунера, а локально выпуклый  инъективный предел в ее определении можно заменить на инъективный предел в категории стереотипных пространств $\Ste$:
 \beq\label{O-exp-top-Ste}
\mathcal{O}_{\exp}(G)=\underset{\scriptsize\begin{matrix}\text{$D$ --
дуально}\\ \text{субмультипликативный}\\ \text{прямоугольник в $\mathcal
O(G)$}\end{matrix}}{\LCS\text{-}\injlim} \C D
= \underset{\scriptsize\begin{matrix}\text{$D$ --
дуально}\\ \text{субмультипликативный}\\ \text{прямоугольник в $\mathcal
O(G)$}\end{matrix}}{\Ste\text{-}\injlim} \C D
 \eeq

\etm

\bcor\label{O_exp(G):comp-subseteq-f^BSQ} Если $G$ -- компактно порожденная
группа Штейна, то всякое ограниченное множество $D$ в $\mathcal{O}_{\exp}(G)$
содержится в некотором прямоугольнике вида $f^\text{\BSQ}$, где $f$ --
некоторый полухарактер на $G$:
$$
D\subseteq f^\text{\BSQ}
$$
 \ecor
\bpr В силу \cite[Theorem 4.1.22]{Akbarov-De-Gruyter-I}, $D$ должен содержаться в каком-то
компакте $E_N$, который, будучи дуально субмультипликативным прямоугольником,
должен по теореме \ref{TH-prayam-poluharakter} иметь вид $f_N^\text{\BSQ}$ для
некоторого полухарактера $f_N$, а именно для $f_N=E_N^\text{\SQ}$. \epr

\btm  Пространство $\mathcal{O}_{\exp}(G)$ функций экспоненциального типа на
компактно порожденной группе Штейна $G$ является проективной стереотипной
алгеброй относительно обычного поточечного умножения функций. \etm

\bpr Заметим, что если функции $u,v\in\mathcal{O}(G)$ подчинены полухарактерам
$f$ и $g$, то их произведение $u\cdot v$ подчинено полухарактеру $f\cdot g$:
$$
u\in f^\text{\BSQ},\ v\in g^\text{\BSQ}\quad\Longrightarrow\quad u\cdot v\in
(f\cdot g)^\text{\BSQ}
$$
Это можно интерпретировать так, что операция умножения $(u,v)\mapsto u\cdot v$
в пространстве $\mathcal{O}(G)$ переводит всякий компакт вида
$f^\text{\BSQ}\times g^\text{\BSQ}$ (где $f$ и $g$ -- полухарактеры) в компакт
$(f\cdot g)^\text{\BSQ}$.
$$
(u,v)\in f^\text{\BSQ}\times g^\text{\BSQ}\quad\mapsto\quad u\cdot v\in (f\cdot
g)^\text{\BSQ}
$$
Причем эта операция в пространстве $\mathcal{O}(G)$ непрерывна, значит она
непрерывно переводит $f^\text{\BSQ}\times g^\text{\BSQ}$ в $(f\cdot
g)^\text{\BSQ}$. С другой стороны, $(f\cdot g)^\text{\BSQ}$, как дуально
субмультипликативный прямоугольник непрерывно вкладывается в пространство
$\mathcal{O}_{\exp}(G)$. Значит мы получаем непрерывное отображение
$$
(u,v)\in f^\text{\BSQ}\times g^\text{\BSQ}\quad\mapsto\quad u\cdot v\in
\mathcal{O}_{\exp}(G)
$$
Если теперь $D$ и $E$ -- произвольные компакты в пространстве Браунера
$\mathcal{O}_{\exp}(G)$, то  по следствию \ref{O_exp(G):comp-subseteq-f^BSQ},
они должны содержаться в компактах вида $f^\text{\BSQ}$ и $g^\text{\BSQ}$, где
$f$ и $g$ -- полухарактеры. Поэтому
$$
D\times E\subseteq f^\text{\BSQ}\times g^\text{\BSQ}
$$
Отсюда можно заключить, что операция умножения должна быть непрерывной из
$D\times E$ в $\mathcal{O}_{\exp}(G)$:
$$
(u,v)\in D\times E\subseteq f^\text{\BSQ}\times g^\text{\BSQ}\quad\mapsto\quad
u\cdot v\in \mathcal{O}_{\exp}(G)
$$
Это верно для произвольных компактов $D$ и $E$ в $\mathcal{O}_{\exp}(G)$.
Значит, мы можем применить \cite[Theorem 3.6.5]{Akbarov-De-Gruyter-I}: поскольку $\mathcal{O}_{\exp}(G)$ как пространство Браунера, кополно (то есть его стереотипное сопряженное пространство полно), оно в силу \cite[Theorem 3.4.16]{Akbarov-De-Gruyter-I}
насыщено, значит любая билинейная форма на нем, непрерывная на компактах вида
$D\times E$, должна быть непрерывна в смысле \cite[3.6.1]{Akbarov-De-Gruyter-I}. То есть операция умножения
$(u,v)\mapsto u\cdot v$ должна быть непрерывна в каком нам нужно смысле, и значит $\mathcal{O}_{\exp}(G)$
--- сильная стереотипная алгебра. \epr

Отметим следующий очевидный факт:

\btm\label{TH-ogranichenie} Если $H$ -- замкнутая подгруппа в группе Штейна
$G$, то ограничение всякой голоморфной функции экспоненциального типа с $G$ на
$H$ также является голоморфной функцией экспоненциального типа:
 \beq
u\in \mathcal{O}_{\exp}(G)\quad\Longrightarrow\quad u|_H\in \mathcal{O}_{\exp}(H)
 \eeq
 \etm

\paragraph{Алгебра $\mathcal O_{\exp}^\star(G)$ экспоненциальных аналитических функционалов.}

Элементы сопряженного стереотипного пространства $\mathcal{O}_{\exp}^\star(G)$, то есть линейные непрерывные функционалы на $\mathcal{O}_{\exp}(G)$, мы будем называть {\it экспоненциальными аналитическими
функционалами}\index{экспоненциальный аналитический функционал} на группе $G$,
а само пространство $\mathcal{O}_{\exp}^\star(G)$ -- {\it пространством
экспоненциальных функционалов} на $G$.

Из теоремы \ref{TH:O_exp--Brauner} следует

\btm\label{TH:O-exp^star-top-Ste} Пространство $\mathcal{O}_{\exp}^\star(G)$ экспоненциальных функционалов
на компактно порожденной группе Штейна $G$ является пространством Фреше, 
и его можно отождествить с проективным пределом в категориях $\LCS$ и $\Ste$ пространств, двойственных к пространствам Смит $\C D$, где $D$ --- дуально субмультипликативные прямоугольники в $\mathcal
O(G)$:
 \beq\label{O-exp^star-top-Ste}
\mathcal{O}_{\exp}^\star(G)=\underset{\scriptsize\begin{matrix}\text{$D$ --
дуально}\\ \text{субмультипликативный}\\ \text{прямоугольник в $\mathcal
O(G)$}\end{matrix}}{\LCS\text{-}\projlim} (\C D)^\star
= \underset{\scriptsize\begin{matrix}\text{$D$ --
дуально}\\ \text{субмультипликативный}\\ \text{прямоугольник в $\mathcal
O(G)$}\end{matrix}}{\Ste\text{-}\projlim} (\C D)^\star
 \eeq
\etm

Отметим следующее предложение, которое можно считать очевидным:

\bprop\label{PROP:delta-func-polny-v-O_exp(G)} В пространстве $\mathcal{O}_{\exp}^\star(G)$ экспоненциальных функционалов
на произвольной группе Штейна $G$ дельта функции полны, то есть всякий функционал $\beta\in \mathcal{O}_{\exp}^\star(G)$ можно приблизить  в $\mathcal{O}_{\exp}^\star(G)$ направленностью функционалов
$\{\beta_i;
i\to\infty\}\subset \mathcal{O}_{\exp}^\star (H)$, являющихся линейными
комбинациями дельта-функционалов:
$$
\beta_i=\sum_{k}\lambda^k_i\cdot\delta^{a^k_i},\qquad \beta_i\overset{\mathcal{O}_{\exp}^\star (H)}{\underset{i\to\infty}{\longrightarrow}}\beta
$$
\eprop
\bpr Это сразу следует из того, что общий аннулятор дельта-функционалов состоит из одного нуля:
$$
\{\delta_x;\ x\in G\}^\perp=\{0\}.
$$
\epr

Формулируемая за ним теорема, наоборот, требует некоторых усилий:

\btm\label{TH:O_exp^*(G)-algebra} Пространство $\mathcal{O}_{\exp}^\star(G)$
экспоненциальных функционалов на компактно порожденной группе Штейна $G$
является сильной стереотипной алгеброй относительно обычной операции
свертки функционалов, определенной формулами
\cite[p.673]{Akbarov-De-Gruyter-I}:
$$
(\alpha,\beta)\quad\mapsto\quad \alpha*\beta
$$
\etm

\bpr 1. Зафиксируем функцию $u\in \mathcal{O}_{\exp}(G)$ и заметим, что для
любой точки $s\in G$ ее сдвиг $u\cdot s$ (определенный формулой
\eqref{DEF:u-cdot-a}) снова является функцией из $\mathcal{O}_{\exp}(G)$:
 $$
\forall s\in G\quad \forall u\in\mathcal{O}_{\exp}(G)\quad u\cdot s\in
\mathcal{O}_{\exp}(G)
 $$
Действительно, подобрав полухарактер $f$ мажорирующий $u$, мы получим:
 $$
\forall t\in G \quad |(u\cdot s)(t)|=|u(s\cdot t)|\le f(s\cdot t)\le f(s)\cdot
f(t)
 $$
То есть,
 \beq\label{u-in-praym=>us-in-s-pryam}
u\in f^\text{\BSQ}\quad\Longrightarrow\quad u\cdot s\in f(s)\cdot f^\text{\BSQ}
 \eeq

2. Обозначим получающееся отображение $s\mapsto u\cdot s$ каким-нибудь
символом, например, $\widehat{u}$:
$$
\widehat{u}:G\to \mathcal{O}_{\exp}(G)\quad\Big|\quad \widehat{u}(s):=u\cdot
s,\qquad s\in G
$$
и покажем, что оно непрерывно. Пусть $s_i$ -- последовательность точек в $G$,
сходящаяся к точке $s$:
$$
s_i\overset{G}{\underset{i\to\infty}{\longrightarrow}} s
$$
Тогда последовательность голоморфных функций $\widehat{u}(s_i)\in \mathcal{O}(G)$ стремится к голоморфной функции $\widehat{u}(s)\in \mathcal{O}(G)$
равномерно на каждом компакте $K\subseteq G$, то есть в пространстве $\mathcal{O}(G)$:
$$
\widehat{u}(s_i)\overset{\mathcal{O}(G)}
{\underset{i\to\infty}{\longrightarrow}} \widehat{u}(s)
$$
При этом из \eqref{u-in-praym=>us-in-s-pryam} следует, что все эти функции
подчинены полухарактеру $C\cdot f$, где $C=\max\{\sup_{i}f(s_i), f(s)\}$ (эта
величина конечна, потому что сходящаяся последовательность $s_i$ образует
вместе со своим пределом $s$ компакт), и поэтому лежат в прямоугольнике,
порожденном полухарактером $C\cdot f$:
$$
\widehat{u}(s_i)=u\cdot s_i\in C\cdot f^\text{\BSQ},\qquad
\widehat{u}(s)=u\cdot s\in C\cdot f^\text{\BSQ}
$$
Иными словами, $\widehat{u}(s_i)$ сходится к $\widehat{u}(s)$ в компакте
$C\cdot f^\text{\BSQ}$
$$
\widehat{u}(s_i)=u\cdot s_i\overset{C\cdot f^\text{\BSQ}}
{\underset{i\to\infty}{\longrightarrow}} u\cdot s=\widehat{u}(s)
$$
Значит, $\widehat{u}(s_i)$ сходится к $\widehat{u}(s)$ в объемлющем этот
компакт пространстве $\mathcal{O}_{\exp}(G)$:
$$
\widehat{u}(s_i)\overset{\mathcal{O}_{\exp}(G)}
{\underset{i\to\infty}{\longrightarrow}} \widehat{u}(s)
$$

3. Из непрерывности отображения $\widehat{u}:G\to\mathcal{O}_{\exp}(G)$
следует, что для всякого функционала $\beta\in \mathcal{O}_{\exp}^\star (G)$
функция $\beta\circ \widehat{u}:G\to\C$ голоморфна. Для доказательства можно
воспользоваться теоремой Мореры: рассмотрим замкнутую ориентированную
гиперповерхность $\varGamma$ в $G$, имеющую достаточно малый диаметр так, чтобы
интеграл по ней любой голоморфной фунции равнялся нулю, и покажем, что интеграл
по ней от $\beta\circ \widehat{u}$ тоже равен нулю:
 \beq\label{morera-1}
\int_{\varGamma} (\beta\circ \widehat{u})(s) \d s=0
 \eeq
Действительно, подберем направленность функционалов $\{\beta_i;
i\to\infty\}\subset \mathcal{O}_{\exp}^\star (G)$, являющихся линейными
комбинациями дельта-функционалов, и аппроксимирующих $\beta$ в $\mathcal{O}_{\exp}^\star (G)$:
$$
\beta_i=\sum_{k}\lambda^k_i\cdot\delta^{a^k_i},\qquad \beta_i\overset{\mathcal{O}_{\exp}^\star (G)}{\underset{i\to\infty}{\longrightarrow}}\beta
$$
Тогда получим: поскольку $\widehat{u}:G\to\mathcal{O}_{\exp}(G)$ непрерывно,
$$
\beta\circ \widehat{u}\overset{\mathcal{C}(G)} {\underset{\infty\gets
i}{\longleftarrow}}\beta_i\circ \widehat{u}
$$
Отсюда следует, что для всякой меры Радона $\alpha\in \mathcal{C}(G)$
$$
\alpha(\beta\circ \widehat{u}) \underset{\infty\gets i}{\longleftarrow}
\alpha(\beta_i\circ \widehat{u})
$$
В частности, для функционала интегрирования по выбранной нами гиперповерхности
$\varGamma$ получаем
 \begin{multline*}
\int_{\varGamma} (\beta\circ \widehat{u})(s) \d s \underset{\infty\gets
i}{\longleftarrow} \int_{\varGamma} (\beta_i\circ \widehat{u})(s)\d s=
\int_{\varGamma} \left(\sum_{k}\lambda^k_i\cdot\delta^{a^k_i}\circ
\widehat{u}\right)(s)\d s=\\= \sum_{k}\lambda^k_i\cdot\int_{\varGamma}
\left(\delta^{a^k_i}\circ \widehat{u}\right)(s)\d s=
\sum_{k}\lambda^k_i\cdot\int_{\varGamma} \delta^{a^k_i}(\widehat{u}(s))\d s=\\=
 \sum_{k}\lambda^k_i\cdot\int_{\varGamma}
\widehat{u}(s)(a^k_i)\d s= \sum_{k}\lambda^k_i\cdot\underbrace{\int_{\varGamma}
u(s\cdot a^k_i)\d s}_{\scriptsize\begin{matrix}\| \\ 0, \\ \text{поскольку $u$
голоморфна} \end{matrix}}=0
 \end{multline*}
То есть действительно справедливо \eqref{morera-1}.

4. Мы показали, что для всякого функционала $\beta\in \mathcal{O}_{\exp}^\star
(G)$ функция $\beta\circ \widehat{u}:G\to\C$ голоморфна. Покажем теперь, что
она имеет экспоненциальный тип:
 \beq\label{beta-circ-widetilde-u-in-mathcal O-exp-G}
\forall u\in \mathcal{O}_{\exp}(G)\quad \forall \beta\in \mathcal{O}_{\exp}^\star (G)\qquad \beta\circ \widehat{u}\in \mathcal{O}_{\exp}(G)
 \eeq
Действительно, поскольку функционал $\beta\in \mathcal{O}_{\exp}^\star (G)$
ограничен на компакте $f^\text{\BSQ}\subseteq \mathcal{O}_{\exp}(G)$, он
является ограниченным функционалом на банаховом отпечатке пространства Смит
$\C f^\text{\BSQ}$, то есть выполняется неравенство
 \beq
\forall v\in\C f^\text{\BSQ}\qquad |\beta(v)|\le M\cdot
\norm{v}_{f^\text{\BSQ}}
 \eeq
где
$$
M=\norm{\beta}_{(f^\text{\BSQ})^\circ}:=\max_{v\in
f^\text{\BSQ}}|\beta(v)|,\qquad \norm{v}_{f^\text{\BSQ}}:=\inf\{\lambda>0:\
v\in \lambda\cdot f^\text{\BSQ}\}
$$
Теперь из формулы \eqref{u-in-praym=>us-in-s-pryam} получаем:
 \begin{multline*}
\forall s\in G\quad u\cdot s\in f(s)\cdot f^\text{\BSQ}
\quad\Longrightarrow\quad \norm{u\cdot s}_{f^\text{\BSQ}}:=\inf\{\lambda>0:\
u\cdot s\in \lambda\cdot f^\text{\BSQ}\}\le f(s) \quad\Longrightarrow\\
\Longrightarrow\quad |(\beta\circ \widehat{u})(s)|=|\beta(u\cdot s)|\le
M\cdot\norm{u\cdot s}_{f^\text{\BSQ}}\le M\cdot f(s)
 \end{multline*}
То есть функция $\beta\circ \widehat{u}$ подчинена полухарактеру $M\cdot
f=\norm{\beta}_{(f^\text{\BSQ})^\circ}\cdot f$:
 \beq\label{beta-circ-widehat-u-in-norm-beta-cdot-f^BSQ}
u\in f^\text{\BSQ}\quad\Longrightarrow\quad \beta\circ \widehat{u}\in
\norm{\beta}_{(f^\text{\BSQ})^\circ}\cdot f^\text{\BSQ}
 \eeq

5. Заметим теперь, что
$$
(\beta\circ\widehat{u})(s)=\beta(u\cdot
s)=\cite[(5.55)]{Akbarov-De-Gruyter-I}=(u*\widetilde{\beta})(s),\qquad s\in G
$$
Получается, что мы доказали, что функция
$u*\widetilde{\beta}=\beta\circ\widehat{u}$ принадлежит пространству $\mathcal{O}_{\exp}(G)$, поэтому для любого функционала $\alpha\in \mathcal{O}_{\exp}^\star(G)$ определена свертка
$$
\alpha * \beta (u)=\eqref{DEF:svertka}=\alpha(u* \widetilde{\beta})
$$
Остается убедиться, что операция $(\alpha,\beta)\mapsto \alpha*\beta$ является
непрерывной билинейной формой (в смысле \cite[3.6.1]{Akbarov-De-Gruyter-I}).

Пусть $\alpha_i$ -- направленность, стремящаяся к нулю в $\mathcal{O}_{\exp}(G)$, а $B$ -- компакт в $\mathcal{O}_{\exp}^\star(G)$. Рассмотрим
произвольный компакт $D$ в $\mathcal{O}_{\exp}(G)$. По следствию
\ref{O_exp(G):comp-subseteq-f^BSQ}, он должен содержаться в каком-то
прямоугольнике $f^\text{\BSQ}$, где $f$ -- полухарактер:
$$
D\subseteq f^\text{\BSQ}
$$
С другой стороны, на пространстве Смит $\C f^\text{\BSQ}$ нормы функционалов
$\beta\in B$ должны быть ограничены:
 $$
\sup_{\beta\in B}\norm{v}_{f^\text{\BSQ}}=M<\infty
 $$
Поэтому в силу \eqref{beta-circ-widehat-u-in-norm-beta-cdot-f^BSQ} мы получаем:
$$
\forall u\in D\quad \forall \beta\in B\qquad u*\widetilde{\beta}=\beta\circ
\widehat{u}\in \norm{\beta}_{(f^\text{\BSQ})^\circ}\cdot f^\text{\BSQ}\subseteq
M\cdot f^\text{\BSQ}
$$
Таким образом, множество $\{u*\widetilde{\beta};\ u\in D,\ \beta\in B\}$
содержится в компакте $M\cdot f^\text{\BSQ}$ пространства $\mathcal{O}_{\exp}^\star(G)$. Значит, направленность функционалов $\alpha_i$ должна
равномерно на нем стремиться к нулю в пространстве $\C$:
 $$
(\alpha_i*\beta)(u)=\alpha_i(u*\widetilde{\beta})\overset{\C}{\underset{i\to\infty}{\rightrightarrows}}
0\qquad \Big(u\in D,\ \beta\in B\Big)
 $$
Это верно для любого компакта $D$, значит направленность функционалов
$\alpha_i*\beta$ стремится к нулю в пространстве $\mathcal{O}_{\exp}^\star(G)$
равномерно по $\beta\in B$:
 $$
\alpha_i*\beta\overset{\mathcal{O}_{\exp}^\star(G)}{\underset{i\to\infty}{\rightrightarrows}} 0\qquad
\Big(\beta\in B\Big)
 $$
Случай противоположного порядка множителей можно не рассматривать из-за
тождества
$$
\widetilde{\alpha*\beta}=\widetilde{\beta}*\widetilde{\alpha}
$$
 \epr

\subsection{Примеры}
\label{SUBSEC:exmples-of-groups}

\paragraph{Конечные группы.}

\label{EX-finite-groups} Как мы отмечали в
\ref{SEC:stein-groups}\ref{SUBSEC-lin-groups}, всякая конечная группа $G$,
рассматриваемая как нульмерное комплексное многообразие, является линейной
комплексной группой Ли. При этом голоморфной функцией на $G$ должна считаться
просто произвольная функция $u:G\to\C$. Поскольку любая такая функция
обязательно ограничена (ее множество значений конечно), мы получаем, что она
должна быть функцией экспоненциального типа. Поэтому алгебры $\mathcal{O}_{\exp}(G)$ и $\mathcal{O}(G)$ в этом случае совпадают и равны алгебре
$\C^G$ всех функций на $G$ с поточечными алгебраическими операциями и
топологией поточечной сходимости:
$$\mathcal{O}_{\exp}(G)=\mathcal{O}(G)=\C^G
$$

\paragraph{Группы $\C^n$.}

Для случая комплексной плоскости $\C$ наше определение, разумеется, совпадает с
обычным -- функциями экспоненциального типа на $\C$ будут целые функции,
растущие не быстрее экспоненты:
 $$
\exists A>0:\quad |u(x)|\le A^{|x|}\qquad (x\in\C).
 $$
В соответствии с классическими теоремами о росте целых функций
\cite{Levin,Shabat-1} это условие эквивалентно тому, что производные $u$ в
фиксированной точке, например, в нуле, растут не быстрее экспоненты:
 $$
\exists B>0:\quad |u^{(k)}(0)|\le B^k\qquad (k\in\N).
 $$
То же самое справедливо и для многих переменных: функциями экспоненциального
типа на $\C^n$ согласно нашему определению будут в точности функции
$u\in\mathcal{O}(\C^n)$, удовлетворяющие условию
 \beq\label{exp-type-for-C}
\exists A>0:\quad |u(x)|\le A^{|x|}\qquad (x\in\C^n),
 \eeq
которое оказывается эквивалентно условию
 \beq\label{exp-type-for-C-u'}
\exists B>0:\quad |u^{(k)}(0)|\le B^{|k|}\qquad (k\in\N^n),
 \eeq
где факториал $k!$ и модуль $|k|$ мультииндекса $k=(k_1,...,k_n)\in\N^n$
определяются равенствами
$$
k!:=k_1\cdot ...\cdot k_n,\qquad |k|:=k_1+...+k_n
$$
Эквивалентность \eqref{exp-type-for-C} нашему определению следует из
предложения \ref{PROP-poluh-porozhd-K}, а равносильность \eqref{exp-type-for-C}
и \eqref{exp-type-for-C-u'} доказывается так же, как и в случае одной
переменной: импликация \eqref{exp-type-for-C-u'} $\Rightarrow$
\eqref{exp-type-for-C} очевидна, а обратная \eqref{exp-type-for-C}
$\Rightarrow$ \eqref{exp-type-for-C-u'} является следствием неравенств Коши
(см. \cite{Shabat}) для коэффициентов $c_k=\frac{u^{(k)}(0)}{k!}$ ряда Тейлора
функции $u$:
$$
\forall R>0\qquad |c_k|\le \frac{\max\limits_{|x|\le
R}|u(x)|}{R^{|k|}}\le\eqref{exp-type-for-C}\le \frac{A^R}{R^{|k|}}
$$
$$
\Downarrow
$$
$$
|c_k|\le
\min_{R>0}\frac{A^R}{R^{|k|}}=\frac{A^R}{R^{|k|}}\Bigg|_{R=\frac{|k|}{\ln
A}}=\frac{A^{\frac{|k|}{\ln A}}}{|k|^{|k|}}\cdot (\ln
A)^{|k|}=\frac{1}{|k|^{|k|}}\cdot B^{|k|},\qquad B=A^{\frac{1}{\ln A}}\cdot \ln
A
$$
$$
\Downarrow
$$
$$
|u^{(k)}(0)|=k!\cdot |c_k|\le \frac{k!}{|k|^{|k|}}\cdot B^{|k|}\le B^{|k|}
$$

\paragraph{Группы $\GL_n(\C)$.}

\label{EX-mnogochleny-na-GL=f-exp-rosta} На группе $\GL_n(\C)$ функциями
экспоненциального типа будут в точности многочлены, то есть функция вида
 \beq\label{mnogochleny-na-GL}
u(x)=\frac{P(x)}{(\det x)^m},\qquad x\in \GL(n,\C),\qquad m\in\Z_+
 \eeq
где $P$ -- обычный многочлен от матричных элементов матрицы $x$, а $\det x$ --
определитель матрицы $x$. Таким образом, справедливо равенство:
 \beq\label{O(GL)=H_exp(GL)}
\mathcal{O}_{\exp}(\GL_n(\C))=\mathcal{P}(\GL_n(\C))
 \eeq
\bpr 1. Сначала докажем вложение $\mathcal{P}(\GL_n(\C))\subseteq \mathcal{O}_{\exp}(\GL_n(\C))$. Заметим для этого, что любой матричный элемент
$$
x\mapsto x_{k,l}
$$
является функцией экспоненциального типа, потому что он подчинен, например,
первой из матричных норм \eqref{matrichnye-normy}:
$$
|x_{k,l}|\le \sum_{i,j=1}^n|x_{i,j}|=||x||
$$
Как следствие, всякий многочлен $x\mapsto P(x)$ от матричных элементов также
является функцией экспоненциального типа (потому что функции экспоненциального
типа образуют алгебру).

Далее, всякая степень определителя, $x\mapsto (\det x)^{-m}$, мультипликативна
$$
(\det (x\cdot y))^{-m}=(\det x)^{-m}\cdot (\det y)^{-m}
$$
и поэтому ее модуль также мультипликативен:
$$
\left|(\det (x\cdot y))^{-m}\right|=\left|(\det x)^{-m}\right| \cdot
\left|(\det y)^{-m}\right|
$$
Значит, функция $x\mapsto \left|(\det x)^{-m}\right|$ является полухарактером
на $\GL(n,\C)$, а $x\mapsto (\det x)^{-m}$ -- функцией экспоненциального типа
на $\GL(n,\C)$.

Умножая теперь две функции экспоненциального типа $x\mapsto P(x)$ и $x\mapsto
(\det x)^{-m}$, мы получим функцию экспоненциального типа $x\mapsto P(x)/ (\det
x)^m$.

2. Теперь докажем обратное вложение: $\mathcal{O}_{\exp}(\GL_n(\C))\subseteq
\mathcal{P}(\GL_n(\C))$. Если $u$ -- голоморфная функция экспоненциального
типа на $\GL_n(\C)$, то, по предложению
\ref{PROP-stroenie-poluharakterov-na-GL_n}, она должна быть подчинена
некоторому полухарактеру вида \eqref{r_N(x)=||x||^N-cdot-||x^-1||^N}. В
частности, для некоторых $C\ge 1$ и $N\in\N$
 \beq\label{|u(t)|-le-max-||x||^N-||x^-1||^N}
|u(x)| \le C\cdot\max\{||x||; ||x^{-1}||\}^N, \qquad
||x||=\sum_{i,j=1}^n|x_{i,j}|
 \eeq
Элементы обратной матрицы $x^{-1}$ получаются из $x$ как дополнительные миноры,
поделенные на определитель, поэтому их можно представлять себе, как многочлены
(степени $n-1$) от $x_{i,j}$ и $(\det x)^{-1}$. Это означает, что можно оценить
правую часть \eqref{|u(t)|-le-max-||x||^N-||x^-1||^N} многочленом (степени
$N(n-1)$) от $|x_{i,j}|$ и $|\det x|^{-1}$ с неотрицательными коэффициентами:
$$
|u(x)| \le C\cdot\max\{||x||; ||x^{-1}||\}^N\le C\cdot
P\Big(\{|x_{i,j}|\}_{1\le i,j\le n},\; |\det x|^{-1}\Big)
$$
Поэтому если домножить $u$ на функцию $(\det x)^{N(n-1)}$, то мы получим
голоморфную функцию на $\GL_n(\C)$, ограниченную многочленом от $|x_{i,j}|$:
$$
\left|u(x)\cdot (\det x)^{N(n-1)}\right| \le C\cdot Q\Big(\{|x_{i,j}|\}_{1\le
i,j\le n}\Big)
$$
Поскольку такая функция локально ограничена в точках аналитического множества
$\det x=0$, по теореме Римана о продолжении \cite{Shabat}, она продолжается до
голоморфной функции на многообразие $\M_n(\C)$ всех комплексных матриц. Значит,
$u(x)\cdot (\det x)^{N(n-1)}$ можно считать голоморфной функцией на
$\M_n(\C)=\C^{n^2}$. Поскольку она имеет полиномальный рост, она представляет
собой некий многочлен $q$ от матричных элементов $x_{i,j}$:
$$
u(x)\cdot (\det x)^{N(n-1)}=q(x)
$$
Отсюда
$$
u(x)=\frac{q(x)}{(\det x)^{N(n-1)}},
$$
то есть $u\in \mathcal{P}(\GL_n(\C))$.
 \epr

\bcor На комплексной окружности $\C^\times$ функциями экспоненциального типа
будут многочлены Лорана (и только они):
$$
u(t)=\sum_{|n|\le N} u_n\cdot t^n,\qquad N\in\N
$$
 \ecor

\subsection{Инъекция $\imath_G:\mathcal{O}_{\exp}(G)\to \mathcal{O}(G)$}\label{SUBSEC:O_exp->O}

Для вложения алгебры $\mathcal{O}_{\exp}(G)$ в алгебру $\mathcal{O}(G)$ нам
будет необходим какой-то символ, мы условимся использовать в этих целях
$\imath_G$:
 \beq\label{H_exp->H}
\imath_G:\mathcal{O}_{\exp}(G)\to \mathcal{O}(G)
 \eeq
Это отображение всегда инъективно, является гомоморфизмом алгебр, и, по
определению топологии в $\mathcal{O}_{\exp}(G)$, непрерывно.

Равенство \eqref{O(GL)=H_exp(GL)} вместе с теоремой \ref{TH-ogranichenie}
означают, что если $G$ -- произвольная линейная группа (с фиксированным
представлением в качестве замкнутой подгруппы в $\GL_n(\C)$), то всякая функция
на $G$, продолжающаяся до многочлена на объемлющую группу $\GL_n(\C)$, является
функцией экспоненциального типа. Таким образом, справедлива цепочка вложений:
 \beq\label{O-subset-H_exp-subset-H}
\mathcal{P}(G)\subseteq \mathcal{O}_{\exp}(G)\subseteq \mathcal{O}(G)
 \eeq
(в которой $\mathcal{P}(G)$ обозначает функции, продолжающиеся до многочленов
на $\GL_n(\C)$ -- это уточнение необходимо потому что $G$ не обязана быть
алгебраической группой).

\btm\label{TH-H_exp-plotno-v-H} Если $G$ -- линейная комплексная группа, то
алгебра $\mathcal{O}_{\exp}(G)$ голоморфных функций экспоненциального типа на
$G$ плотна в алгебре $\mathcal{O}(G)$ всех вообще голоморфных функций на $G$.
\etm
 \bpr
Пусть $\ph:G\to\GL(n,\C)$ -- голоморфное вложение в виде замкнутой подгруппы.
По следствию из теорем Картана \cite[11.5.2]{Taylor}, всякую голоморфную
функцию $v\in \mathcal{O}(G)$ можно продолжить до голоморфной функции $u\in
\mathcal{O}(\GL_n(\C))$. Если затем приблизить $u$ равномерно на компактах
многочленами $u_i\in \mathcal{P}(\GL_n(\C))$ , то, поскольку, в силу
\eqref{O(GL)=H_exp(GL)}, полиномы $u_i$ -- функции экспоненциального типа на
$\GL(n,\C)$, их ограничения $u_i|_G$ будут функциями экспоненциального типа на
$G$, приближающими $v$ равномерно на компактах в $G$.
 \epr

\subsection{Ядерность пространств $\mathcal{O}_{\exp}(G)$ и $\mathcal{O}^\star_{\exp}(G)$}
\label{SUBSEC:nuclearity-of-O_exo(G)}

\btm\label{TH-nuclear} Для всякой компактно порожденной группы Штейна $G$
пространство $\mathcal{O}_{\exp}(G)$ является ядерным пространством Браунера,
а его сопряженное пространство $\mathcal{O}^\star_{\exp}(G)$ -- ядерным
пространством Фреше. \etm

Доказательству этого факта мы предпошлем две леммы. Первая из них верна для
произвольного комплексного многообразия $M$ и доказывается теми же приемами,
что применяются при доказательстве ядерности пространства $\mathcal{O}(\C)$
(см. \cite[6.4.2]{Pietsch}):

\blm\label{LM-ver-mera} Если $M$ -- комплексное многообразие и $K$ и $L$ -- два
компакта в $M$, причем $K$ строго содержится в $L$,
$$
K\subseteq \Int L
$$
то существуют константа $C\ge 0$ и вероятностная мера $\mu$ на $L$ такие, что
для любого $u\in\mathcal{O}(G)$ выполняется неравенство
 \beq\label{ver-mera}
|u|_K\le C\cdot\int_L |\alpha(u)|\ \mu(\d \alpha)
 \eeq
Как следствие, оператор $u|_L\mapsto u|_K$ ограничения с более широкого
компакта на компакт поменьше будет абсолютно суммирующим, причем его квазинорма
абсолютного суммирования оценивается сверху величиной $C$: для любых
$u_1,...,u_l\in\mathcal{O}(M)$
 $$
\sum_{i=1}^l |u_i|_K\le C\cdot \int_L \sum_{i=1}^l |\alpha(u_i)|\; \mu(\d
\alpha)\le C\cdot \sup_{\alpha\in\cabsconv\left(\delta^L\right)} \sum_{i=1}^l
|\alpha(u_i)|
 $$
\elm

\brem Здесь $\cabsconv\left(\delta^L\right)$ обозначает универсальный компакт в
пространстве Смит $\mathcal C^\star(L)$, сопряженном к банахову пространству
$\mathcal C(L)$ непрерывных функций на $L$, или, что то же самое, единичный шар
в банаховом пространстве $\mathcal M(L)$ мер Радона на $L$. Мы используем такое
обозначение, потому что символом $\delta^L=\{\delta^x;\ x\in L\}$ удобно
обозначать множество дельта-функционалов на $\mathcal C(L)$ -- тогда поляра
$B^\circ$ единичного шара $B$ в $\mathcal C(L)$ совпадает с абсолютно выпуклой
замкнутой оболочкой множества $\delta^L$:
$$
B^\circ=\cabsconv\left(\delta^L\right)=\{\alpha\in \mathcal M(L):\
||\alpha||\le 1\}
$$
(замыкание относительно топологии пространства $\mathcal C^\star(L)$). \erem

\blm Для всякой порождающей компактной окрестности единицы $K$ в $G$
$$
G=\bigcup_{n=1}^\infty K^n
$$
найдутся константы $C\ge 0$, $\lambda\ge 0$ такие, что для любого $l\in\N$ и
произвольных $u_1,...,u_l\in\mathcal{O}(G)$, $n\in\N$ выполняется неравенство
 \beq\label{quasinor-abs-sum}
\sum_{i=1}^l |u_i|_{K^n}\le
C\cdot\lambda^{n-1}\cdot\sup_{\alpha\in\cabsconv\left(\delta^{K^{2n+1}}\right)}
\sum_{i=1}^l |\alpha(u_i)|
 \eeq
\elm
 \bpr 1. Множество $U=\Int K$ есть открытая окрестность единицы в $G$, поэтому
система сдвигов $\{x\cdot U;\; x\in K^2\}$ будет открытым покрытием компакта
$K^2$. Выберем из него конечное подпокрытие, то есть конечное множество
$F\subseteq K^2$ такое, что
$$
K^2\subseteq \bigcup_{x\in F} x\cdot U= F\cdot U
$$
Тогда мы получим:
 \beq\label{K^n-subseteq-F^n-1-cdot K}
K^n\subseteq F^{n-1}\cdot K\subseteq F^{n-1}\cdot K^2\subseteq K^{2n+1},\qquad
n\in\N
 \eeq
Это доказывается так: сначала нужно заметить, что
$$
F\subseteq K^2\quad \Longrightarrow\quad F^n\subseteq K^{2n},\qquad n\in\N
$$
После этого в цепочке \eqref{K^n-subseteq-F^n-1-cdot K} достаточно будет
проверить только первое включение:
 \beq\label{K^n-subseteq-F^n-1-cdot K-1}
K^n\subseteq F^{n-1}\cdot K,\qquad n\in\N
 \eeq
Это делается индукцией: при $n=2$ получаем
$$
K^2\subseteq F\cdot U\subseteq F\cdot K
$$
и, если \eqref{K^n-subseteq-F^n-1-cdot K-1} верно для какого-то $n$, то для
$n+1$ получаем:
$$
K^{n+1}=K^n\cdot K\subseteq F^{n-1}\cdot K\cdot K\subseteq F^{n-1}\cdot
K^2\subseteq F^{n-1}\cdot F\cdot K=F^n\cdot K
$$

2. Поскольку компакт $K^2$ строго содержит компакт $K$, по лемме
\ref{LM-ver-mera} существуют $C$, $\mu$ такие, что
 \beq\label{ver-mera-1}
|u|_K\le C\cdot\int_{K^2} |\alpha(u)|\mu(\d \alpha)
 \eeq
При сдвиге компакта $K$ на произвольный элемент группы $x\in G$ эта формула
принимает вид
 \beq\label{ver-mera-2}
|u|_{x\cdot K}\le C\cdot\int_{x\cdot K^2} |\alpha(u)|(x\cdot \mu)(\d \alpha)
 \eeq
где $x\cdot \mu$ -- сдвиг меры $\mu$:
 $$
(x\cdot \mu)(u)=\mu(u\cdot x),\qquad (u\cdot x)(t)=u(x\cdot t)
 $$
Отсюда следует, что если $E$ -- произвольное конечное множество в $G$, то
положив
 $$
\nu=\frac{1}{\card E}\sum_{x\in E}x\cdot \mu
 $$
мы получим вероятностную меру на множестве $E\cdot K^2$, со свойством
 \begin{multline*}
|u|_{E\cdot K}\le \sum_{x\in E} |u|_{x\cdot K}\le C\cdot\sum_{x\in E}
\int_{x\cdot K^2} |\alpha(u)|(x\cdot \mu)(\d \alpha)=\\= C\cdot\sum_{x\in E}
\int_{E\cdot K^2} |\alpha(u)|(x\cdot \mu)(\d \alpha)= C\cdot \int_{E\cdot K^2}
|\alpha(u)| \left(\sum_{x\in E} x\cdot \mu\right)(\d \alpha)=\\=
C\cdot\card(E)\cdot \int_{E\cdot K^2} |\alpha(u)|\; \nu(\d \alpha)
 \end{multline*}
Отсюда в свою очередь можно сделать вывод, что оператор $u|_{E\cdot K^2}\mapsto
u|_{E\cdot K}$ ограничения с более широкого компакта на компакт поменьше будет
абсолютно суммирующим, причем его квазинорма абсолютного суммирования
оценивается сверху величиной $C\cdot\card(E)$: для любых
$u_1,...,u_l\in\mathcal{O}(G)$
 $$
\sum_{i=1}^l |u_i|_{E\cdot K}\le C\cdot\card(E)\cdot \int_{E\cdot K^2}
\sum_{i=1}^l |\alpha(u_i)|\; \nu(\d \alpha)\le C\cdot\card(E)\cdot
\sup_{\alpha\in\cabsconv\left(\delta^{E\cdot K^2}\right)} \sum_{i=1}^l
|\alpha(u_i)|
 $$
 Теперь из формул \eqref{K^n-subseteq-F^n-1-cdot K} получаем:
 \begin{multline*}
\sum_{i=1}^l |u_i|_{K^n}\le \sum_{i=1}^l |u_i|_{F^{n-1}\cdot K}\le
C\cdot\card(F^{n-1})\cdot \sup_{\alpha\in\cabsconv\left(\delta^{F^{n-1}\cdot
K^2}\right)} \sum_{i=1}^l |\alpha(u_i)|\le\\ \le C\cdot\card(F^{n-1})\cdot
\sup_{\alpha\in\cabsconv\left(\delta^{K^{2n+1}}\right)} \sum_{i=1}^l
|\alpha(u_i)|\le C\cdot(\card(F))^{n-1}\cdot
\sup_{\alpha\in\cabsconv\left(\delta^{K^{2n+1}}\right)} \sum_{i=1}^l
|\alpha(u_i)|
 \end{multline*}
и, обозначив $\lambda=\card F$, мы как раз получим формулу
\eqref{quasinor-abs-sum}.
 \epr

 \bpr[Доказательство теоремы \ref{TH-nuclear}]
Рассмотрим последовательность прямоугольников как в теореме \ref{TH:f_N}:
$$
E_N=(f_N)^\text{\BSQ}
$$
и представим $\mathcal{O}_{\exp}(G)$ как инъективный предел последовательности
пространств Смит $\C E_N$ по формуле \eqref{O-exp-top-E_N}:
 $$
\mathcal{O}_{\exp}(G)=\underset{\scriptsize N\to \infty}{\lim} \C E_N
 $$
Поскольку $\mathcal{O}_{\exp}(G)$ -- пространство Браунера, а $\mathcal{O}^\star_{\exp}(G)$ -- пространство Фреше, нам, чтобы показать, что оба
пространства ядерные, достаточно убедиться, что $\mathcal{O}_{\exp}(G)$
коядерно. Рассмотрим $\C E_N$ как банаховы пространства с единичными шарами $
E_N$. В силу \eqref{f_N^black}, полунорма на $\C E_N$ определяется равенством
 $$
p_N(u)=\sup_{n\in\N}\frac{1}{N^n}|u|_{K^n}
 $$
Ее единичным шаром будет множество $ E_N$:
 \begin{multline*}
 E_N=\{u\in\mathcal{O}(G):\;\forall n\in\N\quad |u|_{K^n}\le
N^n\}=\left\{u\in\mathcal{O}(G):\;\forall n\in\N\quad
\frac{1}{N^n}|u|_{K^n}\le 1\right\}=\\=\left\{u\in\mathcal{O}(G):\;p_N(u)=\sup_{n\in\N}\frac{1}{N^n}|u|_{K^n}\le 1\right\}
 \end{multline*}
Чтобы доказать, что пространство $\mathcal{O}_{\exp}(G)=\underset{N\to\infty}{\underset{\rightarrow}{\lim}} \C E_N$
коядерно, нам достаточно убедиться, что для всякого $N\in\N$ найдется $M\in\N$,
$M>N$, такое что отображение вложения $\C E_N\to\C E_M$ будет абсолютно
суммирующим, то есть, таким, что для некоторой константы $L>0$, любого $l\in\N$
и любых $u_1,...,u_l\in \C E_N$ выполняется неравенство:
 \beq\label{psi_N-psi_M-abs-summ}
\sum_{i=1}^l p_M(u_i)\le L\cdot\sup_{\alpha\in(
E_N)^\circ}\sum_{i=1}^l|\alpha(u_i)|
 \eeq
Это доказывается так. Сначала нужно заметить, что для любых $n,N\in\N$
выполняется вложение
 \beq\label{varPsi_N-in-N^n-B_n}
 E_N=\{u\in\mathcal{O}(G):\;\forall n\in\N\quad |u|_{K^n}\le
N^n\}\subseteq \{u\in\mathcal{O}(G):\; |u|_{K^n}\le N^n\}=N^n\cdot
 {^\circ{\delta^{K^n}}}
 \eeq
Из него мы получаем цепочку:
 $$ E_N\subseteq N^n\cdot
 {^\circ\kern-2pt\left(\delta^{K^n}\right)}
 $$
 $$
 \Downarrow
 $$
 $$
( E_N)^\circ\supseteq (N^n\cdot
 {^\circ\kern-2pt\left(\delta^{K^n}\right)})^\circ=\frac{1}{N^n}\cdot
( {^\circ\kern-2pt\left(\delta^{K^n}\right)})^\circ=\frac{1}{N^n}\cdot
\cabsconv (\delta^{K^n})
 $$
 $$
 \Downarrow
 $$
 $$
 \cabsconv (\delta^{K^n})\subseteq N^n\cdot( E_N)^\circ
 $$
 $$
 \Downarrow
 $$
 $$
 \cabsconv (\delta^{K^{2n+1}})\subseteq N^{2n+1}\cdot( E_N)^\circ
 $$
 $$
 \Downarrow
 $$
 \begin{multline*}
\sum_{i=1}^l |u_i|_{K^n}\le\eqref{quasinor-abs-sum}\le C\cdot\lambda^{n-1}\cdot
\sup_{\alpha\in\cabsconv\left(\delta^{K^{2n+1}}\right)}\sum_{i=1}^l
|\alpha(u_i)|\le C\cdot\lambda^{n-1}\cdot\sup_{\alpha\in N^{2n+1}\cdot(
E_N)^\circ}\sum_{i=1}^l |\alpha(u_i)|=\\=C\cdot\lambda^{n-1}\cdot\sup_{\beta\in
( E_N)^\circ} \sum_{i=1}^l |N^{2n+1}\cdot\beta(u_i)|= C\cdot\lambda^{n-1}\cdot
N^{2n+1}\cdot\sup_{\beta\in ( E_N)^\circ}\sum_{i=1}^l |\beta(u_i)|
 \end{multline*}
 $$
 \Downarrow
 $$
 $$
 \forall M>0\qquad
\sum_{i=1}^l\frac{1}{M^n}|u_i|_{K^n}\le C\cdot\frac{\lambda^{n-1}\cdot
N^{2n+1}}{M^n}\cdot\sup_{\beta\in ( E_N)^\circ}\sum_{i=1}^l |\beta(u_i)|
 $$
 $$
 \Downarrow
 $$
 \begin{multline*}
\sum_{i=1}^l p_M(u_i)=\sum_{i=1}^l\sup_{n\in\N}\frac{1}{M^n}|u_i|_{K^n}\le
\sum_{i=1}^l\left(\sum_{n=1}^\infty\frac{1}{M^n}|u_i|_{K^n}\right)=
\sum_{n=1}^\infty\left(\sum_{i=1}^l\frac{1}{M^n}|u_i|_{K^n}\right)\le\\
\le \sum_{n=1}^\infty\left(C\cdot\frac{\lambda^{n-1}\cdot
N^{2n+1}}{M^n}\cdot\sup_{\beta\in ( E_N)^\circ}\sum_{i=1}^l |\beta(u_i)|\right)
= C\cdot\left(\sum_{n=1}^\infty\frac{\lambda^{n-1}\cdot N^{2n+1}}{M^n}\right)
\cdot\sup_{\beta\in ( E_N)^\circ}\sum_{i=1}^l |\beta(u_i)|
 \end{multline*}
Отсюда видно, что если подобрать $M\in\N$ достаточно большим так, чтобы
 $$
C\cdot \sum_{n=1}^\infty\frac{\lambda^{n-1}\cdot N^{2n+1}}{M^n} \le 1
 $$
то константу $L$ в \eqref{psi_N-psi_M-abs-summ} можно будет взять равной $1$:
 $$
 \sum_{i=1}^l p_M(u_i)\le
\underbrace{C\cdot \sum_{n=1}^\infty\frac{\lambda^{n-1}\cdot
N^{2n+1}}{M^n}}_{\scriptsize \begin{matrix} \text{\rotatebox{-90}{$\le$}} \\ 1
\end{matrix}}\cdot\sup_{\beta\in ( E_N)^\circ}\sum_{i=1}^l |\beta(u_i)|
 $$
\epr

\subsection{Голоморфные отображения экспоненциального типа и тензорные
произведения пространств $\mathcal{O}_{\exp}(G)$ и $\mathcal{O}^\star_{\exp}(G)$} \label{SUBSEC:O_exp(G-times-H)}

\btm\label{mathcal-O-exp-G-times-H-cong-mathcal-O-exp-G-mathcal-O-exp-H} Пусть
$G$ и $H$ -- две компактно порожденные группы Штейна. Формула
 \beq\label{rho_G,H}
\rho_{G,H}(u\boxdot v)=u\odot v,\qquad u\in\mathcal{O}_{\exp}(G),\quad
v\in\mathcal{O}_{\exp}(H)
 \eeq
(где функция $u\boxdot v$ определена равенством \eqref{g-h-na-S-T}) задает
линейное непрерывное отображение
$$
\rho_{G,H}:\mathcal{O}_{\exp}(G\times H)\to \mathcal{O}_{\exp}(G)\oslash
\mathcal{O}_{\exp}^\star(H)=\mathcal{O}_{\exp}(G)\odot \mathcal{O}_{\exp}(H),
$$
Это отображение является изоморфизмом стереотипных пространств и естественно по
$G$ и $H$, то есть представляет собой изоморфизм бифункторов из категории
$\mathfrak{SG}$ групп Штейна в категорию ${\tt Ste}$ стереотипных
пространств:
$$
\Big((G;H)\mapsto \mathcal{O}(G\times H)\Big)\rightarrowtail \Big((G;H)\mapsto
\mathcal{O}(G)\odot \mathcal{O}(H)\Big)
$$
Эквивалентно это отображение определяется формулой
 \beq\label{widetilde-w-s-t-star}
\rho_{G,H}(w)(\beta)=\beta\circ\widehat{w},\qquad w\in \mathcal{O}_{\exp}(G\times H),\quad \beta\in \mathcal{O}_{\exp}^\star (H)
 \eeq
где $\widehat{w}:G\to \mathcal{O}_{\exp}(H)$ -- отображение, заданное
равенством
 \beq\label{widetilde-w-s-t=w-s-t}
\widehat{w}(s)(t)=w(s,t),\qquad s\in G,\; t\in H
 \eeq
\etm

 \bcor
Справедливы следующие изоморфизмы функторов:
 \beq\label{tenz-pr-O-exp}
\mathcal{O}_{\exp}(G\times H)\cong \mathcal{O}_{\exp}(G)\odot \mathcal{O}_{\exp}(H)\cong \mathcal{O}_{\exp}(G)\circledast \mathcal{O}_{\exp}(H)
 \eeq
 \beq\label{tenz-pr-O-exp-star}
\mathcal{O}_{\exp}^\star(G\times H)\cong \mathcal{O}_{\exp}^\star(G)\odot
\mathcal{O}_{\exp}^\star(H)\cong \mathcal{O}_{\exp}^\star(G)\circledast
\mathcal{O}_{\exp}^\star(H)
 \eeq
 \ecor

Доказательство теоремы
\ref{mathcal-O-exp-G-times-H-cong-mathcal-O-exp-G-mathcal-O-exp-H} опирается на
следующую ниже лемму
\ref{(f-otimes-g)^blacksquare=f^blacksquare-odot-g^blacksquare}, в формулировке
которой используется инъективное тензорное произведение $A\odot B$ множеств $A$
и $B$ в стереотипных пространствах $X$ и $Y$. Напомним, что согласно
обозначениям \cite[(4.155)]{Akbarov-De-Gruyter-I}, $A\odot B$ определяется как подмножество в
пространстве $X\odot Y=X\oslash Y^\star$ операторов $\varphi:Y^\star\to X$,
состоящее из тех операторов, которые удовлетворяют условию
$\varphi(B^\circ)\subseteq A$. В этом заключается смысл используемой нами ниже
формулы
$$
A\odot B=A\oslash B^\circ.
$$

\blm\label{(f-otimes-g)^blacksquare=f^blacksquare-odot-g^blacksquare} Если
$g:G\to\R_+$ и $h:H\to\R_+$ -- два полухарактера на $G$ и $H$, то $g\boxdot h$
-- полухарактер на $G\times H$, а отображение $\rho_{G,H}$, определяемое
формулами \eqref{widetilde-w-s-t-star} - \eqref{widetilde-w-s-t=w-s-t},
является гомеоморфизмом между компактами $(g\boxdot h)^\text{\BSQ}\subseteq
\mathcal{O}_{\exp}(G\times H)$ и $g^\text{\BSQ}\odot h^\text{\BSQ}\subseteq
\mathcal{O}_{\exp}(G)\odot \mathcal{O}_{\exp}(H)$:
 \beq
 (g\boxdot h)^\text{\BSQ}\cong g^\text{\BSQ}\odot h^\text{\BSQ}
 \eeq
 \elm

\bpr Здесь на первых шагах мы будем повторять рассуждения, употреблявшиеся при
доказательстве теоремы \ref{TH:O_exp^*(G)-algebra}.

1. Заметим сначала, что для любой функции $w\in(g\boxdot
h)^\text{\BSQ}\subseteq \mathcal{O}_{\exp}(G\times H)$ формула
\eqref{widetilde-w-s-t=w-s-t} определяет некое отображение
 $$
\widehat{w}:G\to\mathcal{O}_{\exp}(H)
 $$
Действительно, поскольку $w$ голоморфна на $G\times H$, она голоморфна по обеим
переменным, поэтому при фиксированном $s\in G$ функция $\widehat{w}(s):H\to\C$
тоже будет голоморфна. При этом она подчинена полухарактеру $g(s)\cdot h$:
 \beq\label{widetilde-w-(s)-in-g(s)-cdot-h^blacksquare}
\forall s\in G\qquad \widehat{w}(s)\in g(s)\cdot h^\text{\BSQ}
 \eeq
потому что
 $$
w\in(g\boxdot h)^\text{\BSQ}\quad\Longrightarrow\quad
|\widehat{w}(s)(t)|=|w(s,t)|\le g(s)\cdot h(t)\quad\Longrightarrow\quad
\widehat{w}(s)\in g(s)\cdot h^\text{\BSQ}
 $$
Таким образом, $\widehat{w}(s)$ -- всегда голоморфная функция экспоненциального
типа на $H$, то есть $\widehat{w}((s)\in \mathcal{O}_{\exp}(H)$

2. Покажем, что отображение $\widehat{w}:G\to\mathcal{O}_{\exp}(H)$
непрерывно. Пусть $s_i$ -- последовательность точек в $G$, сходящаяся к точке
$s$: $$ s_i\overset{G}{\underset{i\to\infty}{\longrightarrow}} s
$$ Тогда последовательность голоморфных функций
$\widehat{w}(s_i)\in \mathcal{O}(H)$ стремится к голоморфной функции
$\widehat{w}(s)\in \mathcal{O}(H)$ равномерно на каждом компакте $K\subseteq
H$, то есть в пространстве $\mathcal{O}(H)$:
$$
\widehat{w}(s_i)\overset{\mathcal{O}(H)}
{\underset{i\to\infty}{\longrightarrow}} \widehat{w}(s)
$$
При этом все эти функции подчинены полухарактеру $C\cdot h$, где
$C=\max\{\sup_{i}g(s_i), g(s)\}$ -- конечная величина, потому что сходящаяся
последовательность $s_i$ вместе со своим пределом образует компакт: $$
|\widehat{w}(s_i)(t)|\le g(s_i)\cdot h(t)\le C\cdot h(t) $$ Таким образом,
функции $\widehat{w}(s_i)$ и $\widehat{w}(s)$ лежат в прямоугольнике,
порожденном полухарактером $C\cdot h$: $$
\{\widehat{w}(s_i);\widehat{w}(s)\}\subseteq(C\cdot h)^\text{\BSQ} $$ Иными
словами, $\widehat{w}(s_i)$ сходится к $\widehat{w}(s)$ в компакте $(C\cdot
h)^\text{\BSQ}$ $$ \widehat{w}(s_i)\overset{(C\cdot h)^\text{\BSQ}}
{\underset{i\to\infty}{\longrightarrow}} \widehat{w}(s) $$ поэтому
$\widehat{w}(s_i)$ сходится к $\widehat{w}(s)$ в объемлющем этот компакт
пространстве $\mathcal{O}_{\exp}(H)$: $$ \widehat{w}(s_i)\overset{\mathcal{O}_{\exp}(H)} {\underset{i\to\infty}{\longrightarrow}} \widehat{w}(s)
$$

3. Из непрерывности отображения $\widehat{w}:G\to\mathcal{O}_{\exp}(H)$
следует, что для всякого функционала $\beta\in \mathcal{O}_{\exp}^\star (H)$
функция $\beta\circ \widehat{w}:G\to\C$ голоморфна. Для доказательства можно
воспользоваться теоремой Мореры: рассмотрим замкнутую ориентированную
гиперповерхность $\varGamma$ в $G$ размерности $\dim G-1$, имеющую достаточно малый
диаметр, и покажем, что интеграл по ней равен нулю:
 \beq\label{morera}
\int_{\varGamma} (\beta\circ \widehat{w})(s) \d s=0
 \eeq
Действительно, подберем (с помощью предложения \ref{PROP:delta-func-polny-v-O_exp(G)}) направленность функционалов $\{\beta_i;
i\to\infty\}\subset \mathcal{O}_{\exp}^\star (H)$, являющихся линейными
комбинациями дельта-функционалов, и аппроксимирующих $\beta$ в $\mathcal{O}_{\exp}^\star (H)$:
$$
\beta_i=\sum_{k}\lambda^k_i\cdot\delta^{a^k_i},\qquad \beta_i\overset{\mathcal{O}_{\exp}^\star (H)}{\underset{i\to\infty}{\longrightarrow}}\beta
$$
Тогда получим: поскольку $\widehat{w}:G\to\mathcal{O}_{\exp}(H)$ непрерывно,
$$
\beta\circ \widehat{w}\overset{\mathcal{C}(G)} {\underset{\infty\gets
i}{\longleftarrow}}\beta_i\circ \widehat{w}
$$
Отсюда следует, что для всякой меры Радона $\alpha\in \mathcal{C}(G)$
$$
\alpha(\beta\circ \widehat{w}) \underset{\infty\gets i}{\longleftarrow}
\alpha(\beta_i\circ \widehat{w})
$$
В частности, для функционала интегрирования по выбранной нами гиперповерхности
$\varGamma$ получаем
 \begin{multline*}
\int_{\varGamma} (\beta\circ \widehat{w})(s) \d s \underset{\infty\gets
i}{\longleftarrow} \int_{\varGamma} (\beta_i\circ \widehat{w})(s)\d s=
\int_{\varGamma} \left(\sum_{k}\lambda^k_i\cdot\delta^{a^k_i}\circ
\widehat{w}\right)(s)\d s=\\= \sum_{k}\lambda^k_i\cdot\int_{\varGamma}
\left(\delta^{a^k_i}\circ \widehat{w}\right)(s)\d s=
\sum_{k}\lambda^k_i\cdot\int_{\varGamma} \delta^{a^k_i}(\widehat{w}(s))\d s=\\=
 \sum_{k}\lambda^k_i\cdot\int_{\varGamma}
\widehat{w}(s)(a^k_i)\d s= \sum_{k}\lambda^k_i\cdot\underbrace{\int_{\varGamma}
w(s,a^k_i)\d
s}_{\scriptsize\begin{matrix}\| \\ 0, \\ \text{поскольку $w$ голоморфна} \\
\text{по первому переменному}\end{matrix}}=0
 \end{multline*}
То есть действительно справедливо \eqref{morera}.

4. Мы показали, что для всякого функционала $\beta\in \mathcal{O}_{\exp}^\star
(H)$ функция $\beta\circ \widehat{w}:G\to\C$ голоморфна. Покажем теперь, что
она имеет экспоненциальный тип:
 \beq\label{beta-circ-widetilde-w-in-mathcal O-exp-G}
\forall w\in \mathcal{O}_{\exp}(G\times H)\quad \forall \beta\in \mathcal{O}_{\exp}^\star (H)\qquad \beta\circ \widehat{w}\in \mathcal{O}_{\exp}(G)
 \eeq
Действительно, поскольку функционал $\beta\in \mathcal{O}_{\exp}^\star (H)$
ограничен на компакте $h^\text{\BSQ}\subseteq \mathcal{O}_{\exp}(H)$, он
является ограниченным функционалом на банаховом отпечатке пространства Смит
$\C h^\text{\BSQ}$, то есть выполняется неравенство
 \beq
\forall v\in\C h^\text{\BSQ}\qquad |\beta(v)|\le M\cdot
\norm{v}_{h^\text{\BSQ}}
 \eeq
где
$$
M=\norm{\beta}_{(h^\text{\BSQ})^\circ}:=\max_{v\in
h^\text{\BSQ}}|\beta(v)|,\qquad \norm{v}_{h^\text{\BSQ}}:=\inf\{\lambda>0:\
v\in \lambda\cdot h^\text{\BSQ}\}
$$
Поэтому из формулы \eqref{widetilde-w-(s)-in-g(s)-cdot-h^blacksquare} получаем:
 \begin{multline*}
\widehat{w}(s)\in g(s)\cdot h^\text{\BSQ} \quad\Longrightarrow\quad
\norm{\widehat{w}(s)}_{h^\text{\BSQ}}:=\inf\{\lambda>0:\ \widehat{w}(s)\in
\lambda\cdot h^\text{\BSQ}\}\le g(s) \quad\Longrightarrow\\
\Longrightarrow\quad |\beta(\widehat{w}(s))|\le M\cdot g(s)
 \end{multline*}
То есть функция $\beta\circ \widehat{w}$ подчинена полухарактеру $M\cdot g$ :
 \beq\label{|beta-widetilde-w-s|-le-B-cdot-g-s}
\beta\circ \widehat{w}\in M\cdot g^\text{\BSQ}
 \eeq

5. Мы доказали \eqref{beta-circ-widetilde-w-in-mathcal O-exp-G}. Теперь
покажем, что при фиксированном $w\in(g\boxdot h)^\text{\BSQ}\subseteq \mathcal{O}_{\exp}(G\times H)$ отображение
 \beq\label{beta-mapsto-beta-circ-widetilde-w-in-mathcal O-exp-G}
\beta\in \mathcal{O}_{\exp}^\star (H)\mapsto \rho_{G,H}(w)(\beta)=\beta\circ
\widehat{w}\in \mathcal{O}_{\exp}(G)
 \eeq
непрерывно, то есть
 \beq
 \rho_{G,H}(w)\in \mathcal{O}_{\exp}(G)\oslash \mathcal{O}_{\exp}^\star (H)
 \eeq
Это следует из \eqref{|beta-widetilde-w-s|-le-B-cdot-g-s}: если $\beta_i$ --
направленность, стремящаяся к нулю в $\mathcal{O}_{\exp}^\star (H)$, то
$$
\beta_i\circ\widehat{w}\in M_i\cdot g^\text{\BSQ},\qquad M_i=\max_{v\in
h^\text{\BSQ}}|\beta_i(v)|\underset{i\to\infty}{\longrightarrow} 0
$$
$$
\Downarrow
$$
$$
\beta_i\circ\widehat{w}\overset{g^\text{\BSQ}}{\underset{i\to\infty}{\longrightarrow}}
0
$$
$$ \Downarrow $$ $$ \beta_i\circ\widehat{w}\overset{\mathcal{O}_{\exp}(G)}{\underset{i\to\infty}{\longrightarrow}} 0
$$

6. Теперь нужно убедиться, что
 \beq
\rho_{G,H}(w)\in g^\text{\BSQ}\odot h^\text{\BSQ}=g^\text{\BSQ}\oslash
(h^\text{\BSQ})^\circ\subseteq \mathcal{O}_{\exp}(G)\oslash \mathcal{O}_{\exp}^\star (H)
 \eeq
Или, иными словами,
 $$
\rho_{G,H}(w)\left((h^\text{\BSQ})^\circ\right)\subseteq g^\text{\BSQ}
 $$
Это следует из \eqref{widetilde-w-(s)-in-g(s)-cdot-h^blacksquare}:
 $$
 \forall s\in G\qquad\widehat{w}(s)\in g(s)\cdot h^\text{\BSQ}
 $$
 $$
 \Downarrow
 $$
 $$
 \forall s\in G\qquad \frac{1}{g(s)}\widehat{w}(s)\in h^\text{\BSQ}
 $$
 $$
 \Downarrow
 $$
 $$
\forall s\in G\qquad\forall \beta\in (h^\text{\BSQ})^\circ\qquad  1\ge
\left|\beta\left(\frac{1}{g(s)}\widehat{w}(s)\right)\right|=
\frac{1}{g(s)}|(\beta\circ\widehat{w})(s)|=
\frac{1}{g(s)}|(\rho_{G,H}(w)(\beta))(s)|
 $$
 $$
 \Downarrow
 $$
 $$
\forall s\in G\qquad\forall \beta\in (h^\text{\BSQ})^\circ\qquad
|\rho_{G,H}(w)(\beta)(s)|\le g(s)
 $$
$$
 \Downarrow
 $$
 $$
\forall \beta\in (h^\text{\BSQ})^\circ\qquad \rho_{G,H}(w)(\beta)\in
g^\text{\BSQ}
 $$
 $$
 \Downarrow
 $$
$$
\rho_{G,H}(w)\left((h^\text{\BSQ})^\circ\right)\subseteq g^\text{\BSQ}
 $$

7. Покажем, что отображение
 \beq
w\in (g\boxdot h)^\text{\BSQ} \quad\mapsto\quad \rho_{G,H}(w)\in
g^\text{\BSQ}\odot h^\text{\BSQ}
 \eeq
инъективно. Для этого рассмотрим функционалы вида
 \beq
\delta^{s,t}:\mathcal{O}_{\exp}(G)\oslash \mathcal{O}_{\exp}^\star(H)\to\C\quad\Big|\quad \delta^{s,t}(\ph)=\ph(\delta^t)(s)
 \eeq
Теперь получаем: если $w\ne 0$, то для некоторых $s\in G$, $t\in H$ имеем
$w(s,t)\ne 0$, поэтому
$$
\delta^{s,t}(\rho_{G,H}(w))=\rho_{G,H}(w)(\delta^t)(s)=
(\delta^t\circ\widehat{w})(s)=
\delta^t\Big(\widehat{w}(s)\Big)=\widehat{w}(s)(t)=w(s,t)\ne 0
$$
и значит, $\rho_{G,H}(w)\ne 0$.

8. Точно так же обнаруживается, что отображение
 \beq
w\in (g\boxdot h)^\text{\BSQ} \quad\mapsto\quad \rho_{G,H}(w)\in
g^\text{\BSQ}\odot h^\text{\BSQ}
 \eeq
сюръективно. Для любого $\ph\in g^\text{\BSQ}\odot
h^\text{\BSQ}=g^\text{\BSQ}\oslash (h^\text{\BSQ})^\circ\subset\mathcal{O}_{\exp}(G)\oslash \mathcal{O}_{\exp}^\star(H)$ полагаем
$$
w(s,t)=\delta^s\circledast\delta^t(\ph)=\ph(\delta^t)(s),\qquad s\in G,\ t\in H
$$
и тогда, во-первых, $w$ будет голоморфной функцией на $G\times H$, потому что
она голоморфна по каждой переменной: при фиксированном $t\in H$ объект
$\ph(\delta^t)$ есть элемент пространства $\mathcal{O}_{\exp}(G)$, то есть
голоморфная функция (экспоненциального типа) на $G$, поэтому $w(\cdot,t)$
голоморфна по первой переменной, а при фиксированном $s\in G$ отображение
$\beta\in \mathcal{O}_{\exp}^\star(H) \mapsto (\delta^s\circ\ph)(\beta)$ есть
непрерывный функционал на пространстве $\mathcal{O}_{\exp}^\star(H)$, то есть,
в силу стереотипной двойственности, элемент пространства $\mathcal{O}_{\exp}(H)$:
 $$
(\delta^s\circ\ph)(\beta)=\beta(v),\qquad v\in\mathcal{O}_{\exp}(H)
 $$
Отсюда
 $$
w(s,t)=\ph(\delta^t)(s)=(\delta^s\circ\ph)(\delta^t)=\delta^t(v)=v(t)
 $$
то есть, функция $w(s,\cdot)$ голоморфна по второй переменной.

Далее, разворачивая цепочку пункта 6 в обратном направлении, получаем:
 $$
\ph\in g^\text{\BSQ}\odot h^\text{\BSQ}
 $$
 $$
 \Downarrow
 $$
 $$
\ph\left((h^\text{\BSQ})^\circ\right)\subseteq g^\text{\BSQ}
 $$
 $$
 \Downarrow
 $$
 $$
\forall \beta\in (h^\text{\BSQ})^\circ\qquad \ph(\beta)\in g^\text{\BSQ}
 $$
 $$
 \Downarrow
 $$
 $$
\forall t\in H\quad \frac{1}{h(t)}\delta^t\in
(h^\text{\BSQ})^\circ\quad\Longrightarrow\quad
\ph\left(\frac{1}{h(t)}\delta^t\right)\in g^\text{\BSQ}
 $$
 $$
 \Downarrow
 $$
 $$
\forall t\in H\qquad  \ph(\delta^t)\in h(t)\cdot g^\text{\BSQ}
 $$
 $$
 \Downarrow
 $$
 $$
\forall s\in G\quad \forall t\in H\qquad  |w(s,t)|=|\ph(\delta^t)(s)|\le
h(t)\cdot g(s)=|(g\boxdot h)(s,t)|
 $$
 $$
 \Downarrow
 $$
 $$
w\in (g\boxdot h)^\text{\BSQ}
 $$
И, наконец, остается заметить, что $w$ является прообразом $\ph$ при
отображении $w\mapsto \rho_{G,H}(w)$:
$$
\forall s,t\qquad \rho_{G,H}(w)(\delta^t)(s)=(\delta^t\circ\widehat{w})(s)=
\delta^t(\widehat{w}(s))=\widehat{w}(s)(t)=w(s,t)=\ph(\delta^t)(s)
$$
 $$
 \Downarrow
 $$
$$
\rho_{G,H}(w)=\ph
$$

9. Мы получили, что отображение
 \beq
w\in (g\boxdot h)^\text{\BSQ} \quad\mapsto\quad \rho_{G,H}(w)\in
g^\text{\BSQ}\odot h^\text{\BSQ}
 \eeq
биективно, и теперь остается показать, что оно непрерывно в обе стороны. Это
следует из того, что $(g\boxdot h)^\text{\BSQ}$ и $g^\text{\BSQ}\odot
h^\text{\BSQ}$ -- компакты. Поскольку функционалы $\delta^s\circledast\delta^t$
разделяют точки компакта $g^\text{\BSQ}\odot h^\text{\BSQ}$, порождаемая ими
(отделимая) топология на $g^\text{\BSQ}\odot h^\text{\BSQ}$ совпадает с
топологией компакта $g^\text{\BSQ}\odot h^\text{\BSQ}$: если
$$
\delta^s\circledast\delta^t(\ph_i)\underset{i\to\infty}{\longrightarrow}\delta^s\circledast\delta^t(\ph)
$$
для любых $s\in G$, $t\in H$, то
$$
\ph_i\overset{g^\text{\BSQ}\odot
h^\text{\BSQ}}{\underset{i\to\infty}{\longrightarrow}}\ph
$$
Отсюда получаем, что отображение $w\mapsto \rho_{G,H}(w)$ непрерывно в прямую
сторону:
$$
w_i\overset{(g\boxdot h)^\text{\BSQ}}{\underset{i\to\infty}{\longrightarrow}} w
$$
 $$
 \Downarrow
 $$
$$
\forall (s,t)\in G\times H\qquad
\delta^s\circledast\delta^t(\rho_{G,H}(w_i))=w_i(s,t)\underset{i\to\infty}{\longrightarrow}
w(s,t)=\delta^s\circledast\delta^t(\rho_{G,H}(w))
$$
 $$
 \Downarrow
 $$
$$
\rho_{G,H}(w_i)\overset{g^\text{\BSQ}\odot
h^\text{\BSQ}}{\underset{i\to\infty}{\longrightarrow}}\rho_{G,H}(w)
$$
Итак, операция $w\mapsto \rho_{G,H}(w)$ есть непрерывное биективное отображение
компакта $(g\boxdot h)^\text{\BSQ}$ в компакт $g^\text{\BSQ}\odot
h^\text{\BSQ}$. Значит, $\rho_{G,H}$ есть гомеоморфизм между $(g\boxdot
h)^\text{\BSQ}$ и $g^\text{\BSQ}\odot h^\text{\BSQ}$.
 \epr

\bpr[Доказательство теоремы
\ref{mathcal-O-exp-G-times-H-cong-mathcal-O-exp-G-mathcal-O-exp-H}]  Отметим с
самого начала, что если  $f$ -- полухарактер на $G\times H$, то функции
 $$
g(s)=f(s,1_H),\qquad h(t)=f(1_G,t),\qquad s\in G,\; t\in H
 $$
тоже должны быть полухарактерами (как ограничения $f$ на подгруппы), а функция
$g\boxdot h$ является полухарактером на $G\times H$, мажорирующим $f$:
 \beq\label{f-le-g-otimes-h}
f\le g\boxdot h
 \eeq
Действительно,
$$
f(s,t)=f( (s,1_H)\cdot (1_G,t) )\le f(s,1_H)\cdot f(1_G,t)=g(s)\cdot
h(t)=(g\boxdot h)(s,t)
$$

Отсюда следует, что всякая функция $w\in \mathcal{O}_{\exp}(G\times H)$
содержится в некотором компакте вида $(g\boxdot h)^\text{\BSQ}$ (поскольку $w$
всегда содержится в компакте вида $f^\text{\BSQ}$). При этом объект
$\rho_{G,H}(w)$ является элементом множества $g^\text{\BSQ}\odot
h^\text{\BSQ}$, то есть элементом содержащего его пространства $\mathcal{O}_{\exp}(G)\oslash \mathcal{O}_{\exp}^\star(H)$.

Таким образом, формулы \eqref{widetilde-w-s-t=w-s-t} и
\eqref{widetilde-w-s-t-star} корректно определяют отображение
 \beq
w\in \mathcal{O}_{\exp}(G\times H) \quad\mapsto\quad \rho_{G,H}(w)\in
\mathcal{O}_{\exp}(G)\oslash \mathcal{O}_{\exp}^\star(H)
 \eeq
и нам нужно лишь проверить его биективность и непрерывность в обе стороны.

1. Инъективность этого отображения следует из его инъективности на компактах
$(g\boxdot h)^\text{\BSQ}$ (и того факта, что компакты $(g\boxdot
h)^\text{\BSQ}$ и $g^\text{\BSQ}\odot h^\text{\BSQ}$ инъективно вкладываются в
пространства $\mathcal{O}_{\exp}(G\times H)$ и $\mathcal{O}_{\exp}(G)\oslash
\mathcal{O}_{\exp}^\star(H)$).

2. Сюръективность этого отображения следует из того, что оно сюръективно
отображает компакты $(g\boxdot h)^\text{\BSQ}$ на компакты $g^\text{\BSQ}\odot
h^\text{\BSQ}$ (и из того факта, что компакты $g^\text{\BSQ}\odot
h^\text{\BSQ}$ покрывают все пространство $\mathcal{O}_{\exp}(G)\oslash
\mathcal{O}_{\exp}^\star(H)$).

3. Непрерывность следует из того, что оно непрерывно отображает всякий компакт
$(g\boxdot h)^\text{\BSQ}$ в компакт $g^\text{\BSQ}\odot h^\text{\BSQ}$, и
поэтому в пространство $\mathcal{O}_{\exp}(G)\oslash \mathcal{O}_{\exp}^\star(H)$. Это значит, оно непрерывно на каждом компакте $K$ в
пространстве Браунера $\mathcal{O}_{\exp}(G\times H)$, и поэтому непрерывно на
всем пространстве $\mathcal{O}_{\exp}(G\times H)$.

4. Непрерывность в обратную сторону доказывается так же: поскольку обратное
отображение непрерывно переводит всякий компакт $g^\text{\BSQ}\odot
h^\text{\BSQ}$ в компакт $(g\boxdot h)^\text{\BSQ}$, оно непрерывно на каждом
компакте в пространстве Браунера $\mathcal{O}_{\exp}(G)\oslash \mathcal{O}_{\exp}^\star(H)$. Значит, оно непрерывно на всем пространстве $\mathcal{O}_{\exp}(G)\oslash \mathcal{O}_{\exp}^\star(H)$.

5. Мы доказали, что отображение, определенное формулами
\eqref{widetilde-w-s-t-star} - \eqref{widetilde-w-s-t=w-s-t} является
изоморфизмом стереотипных пространств: $\mathcal{O}_{\exp}(G\times H)\cong
\mathcal{O}_{\exp}(G)\oslash \mathcal{O}_{\exp}^\star(H)=\mathcal{O}_{\exp}(G)\odot \mathcal{O}_{\exp}(H)$ (и значит справедливы тождества
\eqref{tenz-pr-O-exp}). Покажем, что это отображение удовлетворяет тождеству
\eqref{rho_G,H}. Если $u\in\mathcal{O}_{\exp}(G)$, $v\in\mathcal{O}_{\exp}(H)$, то для отображения $\widehat{u\boxdot v}:G\to\mathcal{O}_{\exp}(H)$, определенного формулой \eqref{widetilde-w-s-t=w-s-t}, получаем
логическую цепочку:
$$
\widehat{u\boxdot v}(s)(t)=(u\boxdot v)(s,t)=u(s)\cdot v(t)
$$
$$
\Downarrow
$$
$$
\widehat{u\boxdot v}(s)=u(s)\cdot v
$$
$$
\Downarrow
$$
$$
\forall\beta\in\mathcal{O}_{\exp}^\star(H)\qquad \beta(\widehat{u\boxdot
v}(s))=u(s)\cdot \beta(v)
$$
$$
\Downarrow
$$
$$
\forall\beta\in\mathcal{O}_{\exp}^\star(H)\qquad \rho_{G,H}(u\boxdot
v)(\beta)=\beta\circ\widehat{u\boxdot v}=\beta(v)\cdot
u=\cite[(4.131)]{Akbarov-De-Gruyter-I}=(u\odot v)(\beta)
$$
$$
\Downarrow
$$
$$
\rho_{G,H}(u\boxdot v)=u\odot v
$$
Теперь остается заметить, что, поскольку элементы вида $u\circledast v$
порождают плотное подпространство в $\mathcal{O}_{\exp}(G)\circledast
\mathcal{O}_{\exp}(H)$, в силу доказанной уже ядерности пространств $\mathcal{O}_{\exp}$ соответствующие им элементы вида $u\odot v$ должны порождать плотное
подпространство в $\mathcal{O}_{\exp}(G)\odot \mathcal{O}_{\exp}(H)$, а
элементы $u\boxdot v$ -- плотное подпространство в $\mathcal{O}_{\exp}(G\times
H)$. Отсюда следует, что свойство \eqref{rho_G,H} однозначно определяет
отображение $\rho_{G,H}$.
 \epr

\subsection{Структура алгебр Хопфа на $\mathcal{O}_{\exp}(G)$ и $\mathcal{O}^\star_{\exp}(G)$}\label{proverka-Hopfa-dlya-H(G)}

В \cite{Akbarov-De-Gruyter-I} мы говорили о стандартном приеме, с помощью
которого доказывается, что функциональные алгебры данного класса на группах
являются алгебрами Хопфа -- достаточным условием для этого является
естественный изоморфизм между функциональной алгеброй на декартовом
произведении групп $\times$ и тензорным произведением функциональных алгебр.
Теорема \ref{mathcal-O-exp-G-times-H-cong-mathcal-O-exp-G-mathcal-O-exp-H},
устанавливающая естественный изоморфизм
 $$
\mathcal{O}_{\exp}(G\times H)\stackrel{\rho_{G,H}}{\cong} \mathcal{O}_{\exp}(G)\odot \mathcal{O}_{\exp}(H)
 $$
позволяет теперь доказать то же самое для алгебр $\mathcal{O}_{\exp}(G)$:

\btm\label{O-exp-(G)-algebra-Hopfa} Для всякой компактно порожденной группы
Штейна $G$
 \bit{
\item[--] пространство $\mathcal{O}_{\exp}(G)$ голоморфных функций
экспоненциального типа на $G$ является ядерной алгеброй Хопфа-Браунера
относительно алгебраических операций, определенных формулами, аналогичными
\cite[(5.117)-(5.121)]{Akbarov-De-Gruyter-I};

\item[--] его сопряженное пространство $\mathcal{O}_{\exp}^\star(G)$ является
ядерной алгеброй Хопфа-Фреше относительно сопряженных алгебраических операций.
 }\eit
\etm

\section{Оболочка Аренса-Майкла и рефлексивность относительно нее}\label{SEC:Arens-Michael}

В этом параграфе мы переносимся из категории $\SteAlg$ в категорию $\TopAlg$ ассоциативных локально выпуклых топологических алгебр над $\C$ с единицей и раздельно непрерывным умножением, и морфизмами в виде линейных непрерывных мультипликативных и сохраняющих единицу отображений. Под топологической алгеброй мы понимаем именно объект категории $\TopAlg$.

\subsection{Оболочка Аренса-Майкла}
\label{SUBSEC:Arens-Michael}

\bit{

\item[$\bullet$] Эпиморфизм топологических алгебр $\sigma:A\to A'$ называется {\it
расширением Аренса-Майкла}, если для любой банаховой алгебры $B$ и любого
морфизма $\ph:A\to B$ найдется единственный морфизм $\ph':A'\to B$, замыкающий
диаграмму:
\beq\label{DEF:rasshir-AM}
\begin{diagram}
\node{A} \arrow[2]{e,t}{\sigma} \arrow{se,b}{\forall\ph} \node[2]{A'}
\arrow{sw,b,--}{\exists!\ph'}
\\
\node[2]{B}
\end{diagram}
\eeq

\item[$\bullet$] Расширение Аренса---Майкла $\rho:A\to E$ называется {\it
оболочкой Аренса-Майкла}, если для любого расширения Аренса-Майкла $\sigma:A\to
A'$ найдется единственный морфизм $\upsilon:A'\to E$, замыкающий диаграмму:
\beq\label{DEF:obol-AM}
\begin{diagram}
\node[2]{A} \arrow{sw,t}{\forall\sigma} \arrow{se,t}{\rho}\\
\node{A'}  \arrow[2]{e,b,--}{\exists!\upsilon} \node[2]{E}
\end{diagram}
\eeq
}\eit

Из этого определения ясно, что если $\rho:A\to E$ и $\sigma:A\to F$ -- две
оболочки Аренса-Майкла алгебры $A$, то возникающий гомоморфизм $\upsilon:F\to
E$ (из-за своей единственности) будет изоморфизмом (топологических алгебр).
Поэтому оболочка Аренса-Майкла алгебры $A$ определяется однозначно с точностью
до изоморфизма, и, как следствие, для нее можно ввести специальное обозначение:
$$
\heartsuit_A:A\to A^\heartsuit
$$
Его нужно понимать так: если нам дан какой-то гомоморфизм $\rho:A\to E$, то
запись $\rho=\heartsuit_A$ означает, что $\rho:A\to E$ является оболочкой
Аренса-Майкла алгебры $A$; если же нам дана алгебра $E$, то запись
$E=A^\heartsuit$ означает, что существует гомоморфизм $\rho:A\to E$, являющийся
оболочкой Аренса-Майкла алгебры $A$ --- в этом случае алгебру $E$ также принято
называть {\it оболочкой Аренса-Майкла алгебры} $A$.

Выражаясь терминами главы \ref{CH:obolochki}, оболочка Аренса---Майкла --- это оболочка в категории $\TopAlg$ топологических алгебр в классе $\Epi$
всех эпиморфизмов относительно класса $\Ban$ морфизмов в банаховы алгебры:
$$
A^\heartsuit=\Env_{\Ban}^{\Epi}A
$$
Ниже в следствии \ref{COR:A-plotna-v-A^heartsuit} мы увидим, что оболочка Аренса---Майкла --- не просто эпиморфизм, а плотный эпиморфизм в категории $\TopAlg$:
$$
\heartsuit_A\in\DEpi(\TopAlg).
$$

\paragraph{Сеть банаховых фактор-отображений в категории $\TopAlg$.}

Объясним, как такое определение связано с обычным определением оболочки Аренса-Майкла топологической алгебры.

Абсолютно выпуклая замкнутая окрестность нуля $U$ в топологической алгебре $A$
называется {\it субмультипликативной}, если выполняется вложение:
\beq\label{DEF:submult-okr}
U\cdot U\subseteq U
\eeq
Это эквивалентно тому, что порождаемая $U$ полунорма (функционал Минковского)
на $A$
$$
p:A\to\R_+,\qquad p(u)=\inf\{\lambda>0:\quad u\in\lambda\cdot U\}
$$
является {\it субмультипликативной}\index{субмультипликативная полунорма}, то
есть удовлетворяет условию
$$
p(u\cdot v)\le p(u)\cdot p(v),\qquad u,v\in A
$$

Всякой субмультипликативной абсолютно выпуклой окрестности нуля $U$ в $A$ можно
поставить в соответствие замкнутый идеал $\Ker U$ в $A$, определяемый
равенством
$$
\Ker U=\bigcap_{\varepsilon>0}\varepsilon\cdot U
$$
и фактор-алгебру
$$
A/\Ker U
$$
наделенную топологией нормированного пространства с единичным шаром $U+\Ker U$.
Тогда пополнение $(A/\Ker U)^\blacktriangledown$ будет банаховой алгеброй, и  мы обозначаем эту алгебру
символом $A/U$:
\beq\label{DEF:A/U-compl}
A/U:=(A/\Ker U)^\blacktriangledown
\eeq
Мы называем эту алгебру {\it банаховой фактор-алгеброй} алгебры $A$, а отображение
\beq\label{DEF:pi_U-complex}
\pi_U=\pi_p:A\to A/U=A/p
\eeq
--- {\it банаховым фактор-отображением} алгебры $A$.

Множество всех субмультипликативных окрестностей нуля в алгебре $A$ мы будем обозначать символом ${\mathcal U}_{\Ban}^A$:
\beq\label{DEF:U_C^*^A-complex}
U\in {\mathcal U}_{\Ban}^A\quad\Leftrightarrow\quad U=\{x\in A:\ p(x)\le 1\},\qquad p(x\cdot y)\le p(x)\cdot p(y),
\eeq
а множество всех банаховых фактор-отображений алгебры $A$ мы будем обозначать символом ${\mathcal N}_{\Ban}^A$:
\beq\label{DEF:N_C^*^A-complex}
{\mathcal N}_{\Ban}^A=\{\pi_U:A\to A/U;\ U\in {\mathcal U}_{\Ban}^A\}
\eeq

\medskip
\centerline{\bf Свойства субмультипликативных окрестностей нуля в категории $\TopAlg$:}

\bit{\it

\item[$1^\circ$.]\label{LM:ph=ph_U-circ-pi_U-0-compl-Ste}
Для всякого гомоморфизма $\ph:A\to B$ в банахову алгебру $B$ найдется субмультипликативная окрестность нуля $U\subseteq A$ и гомоморфизм $\ph_U:A/U\to B$, замыкающий диаграмму
\beq\label{ph=ph_U-circ-pi_U-0-compl}
 \xymatrix @R=2pc @C=1.2pc
 {
  A\ar[rr]^{\ph}\ar[dr]_{\pi_U} & & B  \\
  & A/U\ar@{-->}[ur]_{\ph_U} &
  }
\eeq

\item[$2^\circ$.]\label{LM:varkappa^U'_U-0-compl-Ste}
Если $U$ и $U'$ -- две субмультипликативные окрестности нуля в $A$, причем $U\supseteq U'$, то найдется единственный морфизм $\varkappa^{U'}_U:A/U\gets A/U'$, замыкающий диаграмму
\beq\label{varkappa^U'_U-0-compl}
 \xymatrix @R=2pc @C=1.2pc
 {
  & A\ar[ld]_{\pi_U}\ar[rd]^{\pi_{U'}} &   \\
  A/U &  & A/ U'\ar@{-->}[ll]^{\varkappa^{U'}_U}
 }
\eeq

\item[$3^\circ$.]\label{LM:diff-okr-nulya-uporyadocheny-compl-Ste} Пересечение $U\cap U'$ любых двух субмультипликативных окрестностей нуля $U$ и $U'$ в $A$ является субмультипликативной крестностью нуля.

}\eit

\brem\label{REM:ph=ph_U-circ-pi_U-0-compl-*}
Свойство $1^\circ$  в этом списке означает, что система банаховых фактор-отображений $\pi_U:A\to A/U$ порождает класс $\varPhi$ морфизмов в банаховы алгебры изнутри:
\beq\label{ph=ph_U-circ-pi_U-0-compl-*}
    {\mathcal N}_{\Ban}\subseteq\varPhi\subseteq\Mor(\TopAlg)\circ {\mathcal N}_{\Ban}.
\eeq
\erem

\bpr Здесь мы докажем первое свойство. Пусть $\ph:A\to B$ --- гомоморфизм в банахову алгебру $B$. Если $W$ --- единичный шар в $B$, то множество $U=\ph^{-1}(W)$ будет окрестностью нуля в $A$, причем из условия $W\cdot W\subseteq W$ следует условие $U\cdot U\subseteq U$:
 $$
x,y\in U\quad\Rightarrow\quad \ph(x),\ph(y)\in W\quad
\Rightarrow\quad \ph(x\cdot y)=\ph(x)\cdot\ph(y)\in W\quad\Rightarrow\quad
x\cdot y\in U=\ph^{-1}(W)
 $$
Отсюда следует, что морфизм $\ph:A\to B$ пропускается через морфизм
$\pi_U:A\to A/U$
$$
\begin{diagram}
\node{A} \arrow[2]{e,t}{\pi_U} \arrow{se,b}{\ph} \node[2]{A/U}
\arrow{sw,b,--}{\ph_U}
\\
\node[2]{B}
\end{diagram}
$$
\epr

\btm\label{LM:A/U-set-epimorfizmov-v-InvSteAlg-complex}
Система ${\mathcal N}_{\Ban}^A$ банаховых фактор-отображений образует сеть эпиморфизмов\footnote{См. определение на с.\pageref{DEF:set-epimorf}.} в категории $\TopAlg$ топологических алгебр, то есть обладает следующими свойствами:
  \bit{
\item[(a)] у всякой алгебры $A$ есть хотя бы одна субмультипликативная окрестность нуля $U$, и множество всех субмультипликативных окрестностей нуля в $A$ направлено относительно предпорядка
$$
U\le U'\quad\Longleftrightarrow\quad U\supseteq U',
$$

\item[(b)] для всякой алгебры $A$ система морфизмов $\varkappa_U^{U'}$ из \eqref{varkappa^U'_U-0-compl} ковариантна, то есть для любых трех субмультипликативных окрестностей нуля $U\supseteq U'\supseteq U''$ коммутативна диаграмма
$$
 \xymatrix @R=2pc @C=1.2pc
 {
  A/U &  & A/ U''\ar[ll]_{\varkappa^{U''}_U}\ar[dl]^{\varkappa^{U''}_{U'}}\\
  & A/U \ar[ul]^{\varkappa^{U'}_U} &
 }
$$
и эта система $\varkappa_U^{U'}$ обладает проективным пределом в $\TopAlg$;

\item[(c)] для всякого гомоморфизма $\alpha:A\gets A'$ в $\TopAlg$ и любой субмультипликативной окрестности нуля $U$ в $A$ найдется субмультипликативная окрестность нуля $U'$ в $A'$ и гомоморфизм $\alpha_U^{U'}:A/U\gets A'/U'$ такие, что коммутативна
диаграмма
 \beq\label{DIAGR:set-Ban} \xymatrix @R=2.5pc @C=4.0pc {
 A\ar[d]_{\pi_U} & A'\ar@{-->}[d]^{\pi_{U'}}\ar[l]_{\alpha} \\
A/U & A'/U'\ar@{-->}[l]^{\alpha_U^{U'}}
 } \eeq
 }\eit
\etm

В соответствии с пунктом (b) этой теоремы, существует проективный предел $\projlim_{U'\in {\mathcal U}_{\Ban}^A} A/ U'$ ковариантной системы $\varkappa_U^{U'}$ в  $\TopAlg$. Как следствие, существует единственная стрелка $\rho:A\to \projlim_{U'\in {\mathcal U}_{\Ban}^A} A/ U'$ в  $\TopAlg$, замыкающая все диаграммы
\beq\label{A-to-leftlim-A/U-0-complex}
 \xymatrix @R=2pc @C=1.2pc
 {
   A\ar[rd]_{\pi_U}\ar@{-->}[rr]^(.3){\rho} & & \projlim_{U'\in {\mathcal U}_{\Ban}^A} A/ U'\ar[ld]^{\varkappa_U}   \\
   & A/U &  
 },\qquad U\in {\mathcal U}_{\Ban}^A
\eeq

\btm\label{TH:harakt-obol-A-M} В категории $\TopAlg$ гомоморфизм
 \beq\label{vlozhenie-A-v-proj-limit}
\rho:A\to\projlim_{V\in {\mathcal U}_{\Ban}^A}A/V
 \eeq
является оболочкой Аренса-Майкла алгебры $A$:
 \beq\label{A^heartsuit=proj-limit}
A^\heartsuit=\TopAlg\text{-}\projlim_{V\in {\mathcal U}_{\Ban}^A}A/V=\LCS\text{-}\projlim_{V\in {\mathcal U}_{\Ban}^A}A/V
 \eeq
\etm
 \bpr
1. Сначала нужно заметить, что каждый морфизм $\pi_U:A\to A/U$ является эпиморфизмом локально выпуклых пространств. Поэтому морфизм $\rho:A\to\projlim_{V\in {\mathcal U}_{\Ban}^A}A/V$ тоже должен быть эпиморфизмом  локально выпуклых пространств\footnote{Это особенность категории $\TopAlg$, в категории $\SteAlg$ это правило не работает!}. Действительно, пусть $U\subseteq A$ --- субмультипликативная окрестность нуля и пусть $c\in \LCS\text{-}\projlim_{V\in {\mathcal U}_{\Ban}^A}A/V$. Положим 
$$
b=\varkappa_U(c).
$$
Поскольку $\pi_U:A\to A/U$ --- плотный эпиморфизм, найдется элемент $a\in A$ такой что
$$
\pi_U(a)-b\in U.
$$
Теперь мы получаем
$$
\pi_U(a)-b=\varkappa_U(\rho(a))-\varkappa_U(c)=\varkappa_U(\rho(a)-c)\in U\quad\Rightarrow\quad \rho(a)-c\in \varkappa_U^{-1}(U).
$$
Поскольку окрестности $\varkappa_U^{-1}(U)$ образуют локальную базу окрестностей нуля в $\LCS\text{-}\projlim_{V\in {\mathcal U}_{\Ban}^A}A/V$, мы получаем, что $\rho(A)$ плотно в $\LCS\text{-}\projlim_{V\in {\mathcal U}_{\Ban}^A}A/V$.

  Как следствие, морфизм $\rho$ является эпиморфизмом и в категории $\TopAlg$ топологических алгебр.

2. Далее покажем, что $\rho$ является расширением Аренса-Майкла алгебры $A$.
Пусть $\ph:A\to B$ -- морфизм в какую-нибудь банахову алгебру $B$. 
По лемме \ref{LM:ph=ph_U-circ-pi_U-0-compl}
найдется субмультипликативная окрестность нуля $U$ в $A$ такая, что морфизм $\ph$ пропускается через морфизм
$\pi_U:A\to A/U$:
\eqref{A-to-leftlim-A/U-0-complex}:
\beq\label{TH:harakt-obol-A-M-0}
\begin{diagram}
\node{A} \arrow[2]{e,t}{\pi_U} \arrow{se,b}{\ph} \node[2]{A/U}
\arrow{sw,b,--}{\ph_U}
\\
\node[2]{B}
\end{diagram}
\eeq
Теперь положив $\ph'=\ph_U\circ\varkappa_U$, мы получим диаграмму
$$
\begin{diagram}\dgARROWLENGTH=4em
\node{A} \arrow[4]{e,t}{\rho} \arrow{see,t}{\pi_U}\arrow{sssee,b}{\ph}
\node[4]{\projlim_{V\in {\mathcal U}_{\Ban}^A}A/V} \arrow{sww,t,--}{\varkappa_U} \arrow{sssww,b,--}{\ph'}
\\
\node[3]{A/U}\arrow[2]{s,b}{\ph_U}\\ \\ \node[3]{B}
\end{diagram}
$$
В ней верхний внутренний треугольник есть диаграмма \eqref{A-to-leftlim-A/U-0-complex}, левый внутренний треугольник --- диаграмма \eqref{TH:harakt-obol-A-M-0}, а правый внутренний треугольник --- определение морфизма $\ph'$. Как следствие, периметр этой диаграммы коммутативен, а это и есть нужная нам диаграмма \eqref{DEF:rasshir-AM}.
Стрелка $\ph'$ здесь будет единственной, потому что, как мы уже поняли в пункте 1, $\rho$ --- эпиморфизм.

3. Покажем наконец, что расширение $\rho$ является оболочкой Аренса-Майкла. Пусть
$\sigma:A\to A'$ -- какое-нибудь другое расширение. Тогда всякий морфизм
$\pi_U:A\to A/U$ порождает единственный морфизм $(\pi_U)'$, для которого будет коммутативна диаграмма \eqref{A-to-leftlim-A/U-0-complex}:
\beq\label{TH:harakt-obol-A-M-1}
\begin{diagram}
\node{A} \arrow[2]{e,t}{\sigma} \arrow{se,b}{\pi_U} \node[2]{A'}
\arrow{sw,b,--}{(\pi_U)'}
\\
\node[2]{A/U}
\end{diagram}
\eeq
Зафиксируем теперь еще одну субмультипликативную окрестность $V\subseteq U$ и рассмотрим диаграмму
$$
\begin{diagram}\dgARROWLENGTH=6em
\node[3]{A}\arrow[2]{s,b,3}{\sigma} \arrow{sssee,t}{\pi_U} \arrow{sssww,t}{\pi_V} \\ \\
\node[3]{A'} \arrow{see,b,--}{(\pi_U)'}\arrow{sww,b,--}{(\pi_V)'}
 \\
\node{A/V} \arrow[4]{e,b,--,3}{\varkappa_U^V}\node[4]{A/U}
\end{diagram}
$$
Здесь левый и правый внутренние треугольники коммутативны, потому что это варианты диаграммы \eqref{TH:harakt-obol-A-M-1}, а периметр коммутативен, потому что это диаграмма \eqref{varkappa^U'_U-0-compl}. 
Поскольку вдобавок стрелка $\sigma$ --- эпиморфизм, это значит, что нижний внутренний треугольник здесь тоже должен быть коммутативен:
$$
 \xymatrix @R=2pc @C=1.2pc
 {
 & A' \ar[dr]^{(\pi_U)'}\ar[dl]_{(\pi_V)'} &
 \\
A/V \ar[rr]_{\varkappa_U^V} & & A/U
}
$$
Это значит, что система морфизмов $(\pi_U)':A'\to A/U$ --- проективный конус ковариантной системы $\varkappa^V_U:A/V\to A/U$.

Как следствие существует единственный морфизм $\upsilon$ из $A'$ в проективный предел $\projlim_{V\in {\mathcal U}_{\Ban}^A}A/V$, замыкающий все диаграммы
\beq\label{TH:harakt-obol-A-M-2}
 \xymatrix @R=2pc @C=1.2pc
 {
 & A' \ar[dr]^{(\pi_U)'}\ar@{-->}[dl]_{\upsilon} &
 \\
\projlim_{V\in {\mathcal U}_{\Ban}^A}A/V \ar[rr]_{\varkappa_U^V} & & A/U
}
\eeq
Рассмотрим теперь такую диаграмму:
$$
 \xymatrix @R=3pc @C=3pc
 {
   A\ar[rd]_{\pi_U}\ar[r]^{\sigma} & A'\ar[d]^{(\pi_U)'}\ar[r]^(.3){\upsilon} & \projlim_{U'\in {\mathcal U}_{\Ban}^A} A/ U'\ar[ld]^{\varkappa_U}   \\
   & A/U &  
 }
$$
В ней левый внутренний треугольник коммутативен, потому что это диаграмма 
\eqref{TH:harakt-obol-A-M-1}, а правый внутренний треугольник коммутативен, потому что это диаграмма \eqref{TH:harakt-obol-A-M-2}. 
Значит, периметр тоже должен быть коммутативен:
\beq\label{TH:harakt-obol-A-M-3}
 \xymatrix @R=3pc @C=3pc
 {
   A\ar[rd]_{\pi_U}\ar@{-->}[rr]^(.4){\upsilon\circ\sigma} &  & \projlim_{U'\in {\mathcal U}_{\Ban}^A} A/ U'\ar[ld]^{\varkappa_U}   \\
   & A/U &  
 }
\eeq
И это верно для всякой окрестности $U\in {\mathcal U}_{\Ban}^A$. Это можно понимать так, что морфизм
$\upsilon\circ\sigma$ --- морфизм проективного конуса $\pi_U:A\to A/U$ в проективный конус $\varkappa_U:\projlim_{U'\in {\mathcal U}_{\Ban}^A} A/ U'\to A/U$. Но второй из этих конусов --- это проективный предел. Поэтому такой морфизм конусов единственен, и значит морфизм $\upsilon\circ\sigma$ просто совпадает с морфизмом $\rho$:
$$
\upsilon\circ\sigma=\rho.
$$
Мы получили коммутативную диаграмму \eqref{DEF:obol-AM}:
$$
\begin{diagram}
\node[2]{A} \arrow{sw,t}{\sigma} \arrow{se,t}{\rho}\\
\node{A'}  \arrow[2]{e,b,--}{\upsilon} \node[2]{\projlim_{U'\in {\mathcal U}_{\Ban}^A} A/ U'}
\end{diagram}
$$
В ней морфизм $\upsilon$ будет единственным, потому что $\sigma$, как всякое расширение Аренса---Майкла, эпиморфизм.
\epr

\bcor\label{COR:A-plotna-v-A^heartsuit} Алгебра $A$ всегда плотна в своей оболочке Аренса-Майкла $A^\heartsuit$ (то есть морфизм $\heartsuit_A:A\to A^\heartsuit$ всегда является плотныи эпиморфизмом).
\ecor
\bpr
В пункте 1 доказательства теоремы \ref{TH:harakt-obol-A-M} мы показали, что $\rho(A)$ плотно в $\LCS\text{-}\projlim_{V\in {\mathcal U}_{\Ban}^A}A/V$.
\epr

\bcor\label{COR:Arens-Michael-v-Ste=v-TopAlg} Если в контраварианной системе банаховых фактор-алгебр $A/U$ имеется счетная конфинальная подсистема, то оболочка Аренса-Майкла $A^\heartsuit$ алгебры $A$ является алгеброй Фреше (и поэтому стереотипной алгеброй).
\ecor
\bpr
Проективный предел счетной системы банаховых алгебр --- это алгебра Фреше.
\epr

\paragraph{Алгебры Аренса-Майкла}

 \bit{
\item[$\bullet$] Топологическая алгебра $A$ называется {\it алгеброй
Аренса-Майкла}, если ее оболочка Аренса-Майкла $\heartsuit:A\to A^\heartsuit$
является изоморфизмом топологических алгебр.
 }\eit

Теорема \ref{TH:harakt-obol-A-M} позволяет дать следующую характеризацию
алгебрам Аренса-Майкла:

 \bprop
Топологическая алгебра $A$ является алгеброй Аренса-Майкла если она полна (как
топологическое векторное пространство) и удовлетворяет следующим равносильным
условиям:
 \bit{
\item[(i)] топология $A$ задается системой субмультипликативных полунорм;

\item[(ii)] топология $A$ является псевдонасыщением топологии, в которой  субмультипликативные замкнутые
абсолютно выпуклые окрестности нуля образуют локальную базу.
 }\eit
 \eprop

\noindent\rule{160mm}{0.1pt}\begin{multicols}{2}

\bex\label{EX-polunormy-v-O(Z)} Полунормы \eqref{polunormy-v-H(Z)}, задающие
топологию на $\mathcal{O}(\Z)$, субмультипликативны:
 $$
||u||_N=\sum_{|n|\le N} |u(n)|,\qquad N\in\N,
 $$
и поэтому $\mathcal{O}(\Z)$ -- алгебра Аренса-Майкла.
 \eex \bpr
Действительно,
 \begin{multline*}
\norm{u\cdot v}_N=\sum_{|n|\le N}|(u\cdot
v)(n)|=\eqref{umnozhenie-v-O(Z)-i-O^star(Z)}=\\=\sum_{|n|\le N}|u(n)\cdot v(n)|\le\\ \le
\left(\sum_{|n|\le N}|u(n)|\right)\cdot\left(\sum_{|n|\le
N}|v(n)|\right)=\norm{u}_N\cdot\norm{v}_N
\end{multline*}
 \epr

\bex\label{EX-polunormy-v-O(C^x)} Полунормы \eqref{polunormy-v-H(C-x)},
задающие топологию на $\mathcal{O}(\C^\times)$, субмультипликативны:
 $$
||u||_C=\sum_{n\in\Z} |u_n|\cdot C^{|n|},\qquad C\ge 1,
 $$
и поэтому $\mathcal{O}(\C^\times)$ -- алгебра Аренса-Майкла.
 \eex \bpr
Действительно,
 \begin{multline*}
||u\cdot v||_C=\sum_{n\in\Z} |(u\cdot v)_n|\cdot C^{|n|}
=\eqref{umnozhenie-v-O(C-x)-i-O^star(C-x)}=\\=\sum_{n\in\Z} \left|\sum_{i\in\Z}
u_i\cdot v_{n-i}\right|\cdot C^{|n|}\le\\ \le \sum_{n\in\Z} \sum_{i\in\Z}
|u_i|\cdot |v_{n-i}|\cdot C^{|i|}\cdot C^{|n-i|}=\\=\left(\sum_{k\in\Z} |u_k|\cdot
C^{|k|}\right)\cdot \left(\sum_{l\in\Z} |u_l|\cdot C^{|l|}\right)=||u||_C\cdot
||v||_C
 \end{multline*}\epr

\bex\label{EX-polunormy-v-H(C)} Полунормы \eqref{polunormy-v-H(C)}, задающие
топологию на $\mathcal{O}(\C)$, субмультипликативны:
 $$
||u||_C=\sum_{n=0}^\infty |u_n|\cdot C^n,\qquad C\ge 0,
 $$
и поэтому $\mathcal{O}(\C)$ --- алгебра Аренса-Майкла.
 \eex \bpr
Действительно,
 \begin{multline*}
||u\cdot v||_C=\sum_{n=0}^\infty |(u\cdot v)_n|\cdot C^n
=\eqref{umnozhenie-v-O(C)-i-O^star(C)}=\\=\sum_{n=0}^\infty \left|\sum_{i=0}^n
u_i\cdot v_{n-i}\right|\cdot C^n\le\\ \le \sum_{n=0}^\infty \sum_{i=0}^n |u_i|\cdot
|v_{n-i}|\cdot C^n=\\=\left(\sum_{k\in\N} |u_k|\cdot C^k\right)\cdot
\left(\sum_{l=0}^\infty |u_l|\cdot C^l\right)=\\=||u||_C\cdot ||v||_C
 \end{multline*}
 \epr

\end{multicols}\noindent\rule[10pt]{160mm}{0.1pt}

\bprop\label{PROP:O(M)-AM} Алгебра $\mathcal{O}(M)$ голоморфных функций на комплексном
многообразии $M$ является алгеброй Аренса-Майкла. \eprop
\bpr
Каждому компакту $T\subseteq M$ поставим в соответствие полунорму
$$
p_T(u)=\max_{t\in T}\abs{u(t)},\quad u\in {\mathcal O}(M).
$$
Это субмультипликативная полунорма на $\mathcal{O}(M)$, поэтому ее единичный шар
$$
U_K=\{u\in {\mathcal O}(M): \ p_T(u)\le 1\}
$$
--- субмультипликативная окрестность нуля в $\mathcal{O}(M)$. Ясно, что алгебра $\mathcal{O}(M)$ совпадает с локально выпуклым проективным пределом банаховых фактор-алгебр по таким окрестностям нуля:
$$
{\mathcal O}(M)=\projlim_{T\subseteq M}{\mathcal O}(M)/U_T.
$$ 
С другой стороны, всякая непрерывная субмультипликативная полунорма $p$ на $\mathcal{O}(M)$ мажорируется некоторой полунормой $p_T$.
$$
p(u)\le C\cdot p_T(u), \quad C>0.
$$
потому что полунормы $p_T$ задают топологию на $\mathcal{O}(M)$. Как следствие, всякая субмультипликативная окрестность нуля $U$ в $\mathcal{O}(M)$ содержит некоторую гомотетию некоторой окрестности $U_T$:
$$
\e\cdot U_T\subseteq U, \quad \e>0.
$$
Из этого следует, что окрестности вида $U_T$ образуют конфинальную подсистему в системе всех субмультипликативных окрестностей в $\mathcal{O}(M)$. Теперь мы получаем цепочку
$$
\mathcal{O}(M)=\projlim_{T\subseteq M}{\mathcal O}(M)/U_T=\projlim_{U\in {\mathcal U}_{\Ban}^A}{\mathcal O}(M)/U=\mathcal{O}(M)^\heartsuit.
$$
\epr

\bprop Алгебра $\mathcal{O}_{\exp}^\star(G)$ экспоненциальных функционалов на
любой группе Штейна является алгеброй Аренса-Майкла. \eprop

\bprop Алгебра $\mathcal{P}(M)$ многочленов на комплексном алгебрическом
многообразии $M$, снабженная сильнейшей локально выпуклой топологией является
алгеброй Аренса-Майкла, если и только если многообразие $M$ конечно. \eprop

Из \eqref{ph=ph_U-circ-pi_U-0-compl-*} и теоремы \ref{TH:kriterij-obolochki-v-term-polnyh-objektov} следует

\btm\label{TH:kriterij-obolochki-AM-v-term-algebr-AM}
Для всякого морфизма стереотипных алгебр\footnote{Здесь не предполагается, что $\eta$ лежит в $\DEpi$.} $\eta:A\to S$ в произвольную алгебру Аренса---Майкла $S$ следующие условия эквивалентны:
\bit{

\item[(i)] морфизм $\eta:A\to S$ является оболочкой Аренса---Майкла;

\item[(ii)] для любого морфизма  стереотипных алгебр\footnote{Здесь не предполагается, что $\theta$ лежит в $\DEpi$.} $\theta:A\to E$ в произвольную алгебру Аренса---Майкла $E$ найдется единственный морфизм стереотипных алгебр $\upsilon:S\to E$, замыкающий диаграмму
 \beq\label{kriterij-obolochki-AM-v-term-algebr-AM}
\xymatrix @R=2.pc @C=2.0pc 
{
& A \ar[dl]_{\eta} \ar[dr]^{\theta} &\\
S \ar@{-->}[rr]^{\upsilon} & & E
}
 \eeq
}\eit
\etm

\bprop\label{kriterij-A-M} Морфизм топологических алгебр $\pi:A\to B$ является
оболочкой Аренса-Майкла для топологической алгебры $A$, если и только если
 \bit{
\item[\rm (i)] $B$ -- алгебра Аренса-Майкла,

\item[\rm (ii)]  образ $\pi(A)$ алгебры $A$ под действием $\pi$ плотен в $B$,

\item[\rm (iii)] для любой непрерывной субмультипликативной полунормы
$p:A\to\R_+$ найдется непрерывная субмультипликативная полунорма
$\widetilde{p}:B\to\R_+$, такая что полунорма $\widetilde{p}\circ\pi:A\to\R_+$
мажорирует полунорму $p$:
$$
p(a)\le\widetilde{p}(\pi(a)),\qquad a\in A
$$
 }\eit\noindent
Последнее условие в этой системе можно заменить условием
 \bit{
\item[\rm (iii)'] любая непрерывная субмультипликативная полунорма $p:A\to\R_+$
продолжается до непрерывной субмультипликативной полунормы
$\widetilde{p}:B\to\R_+$:
$$
p(a)=\widetilde{p}(\pi(a)),\qquad a\in A
$$
 }\eit
 \eprop

Следующие утверждения показывают, что операция взятия оболочки Аренса-Майкла
коммутирует с операциями перехода к прямой сумме и фактор-алгебре.

\bprop Оболочка Аренса-Майкла прямой суммы $A_1\oplus ...\oplus A_n$ конечного
набора топологических алгебр $A_1,...,A_n$ совпадает с прямой суммой оболочек
Аренса-Майкла этих алгебр:
$$
(A_1\oplus ...\oplus A_n)^\heartsuit\cong A_1^\heartsuit\oplus ...\oplus
A_n^\heartsuit
$$
 \eprop

\bprop Пусть $\pi:A\to A^\heartsuit$ -- оболочка Аренса-Майкла алгебры $A$, и
пусть $I$ -- замкнутый идеал в $A$. Тогда оболочка Аренса-Майкла фактор-алгебры
$A/I$ совпадает с пополнением фактор-алгебры $A^\heartsuit/\overline{\pi(I)}$
по замыканию $\overline{\pi(I)}$ в $A^\heartsuit$ образа идеала $I$ под
действием отображения $\pi$:
$$
(A/I)^\heartsuit\cong
\big(A^\heartsuit/\overline{\pi(I)}\big)^\blacktriangledown
$$
 \eprop

\paragraph{Теорема Пирковского.}

Важный пример оболочки Аренса-Майкла был построен А.Ю.Пирковским:

\btm[А.Ю.Пирковский, \cite{Pi2008}]\label{TH-Pirkovskii} Оболочка
Аренса-Майкла алгебры $\mathcal{P}(M)$ многочленов на аффинном алгебраическом
многообразии $M$ совпадает с алгеброй $\mathcal{O}(M)$ голоморфных функций на
$M$:
 \beq\label{R(M)^heartsuit=O(M)}
\Big(\mathcal{P}(M)\Big)^\heartsuit\cong \mathcal{O}(M)
 \eeq
\etm

\subsection{Отображение $\imath_G^\star:\mathcal{O}^\star(G)\to \mathcal{O}_{\exp}^\star(G)$
как оболочка Аренса-Майкла} \label{SUBSEC:O*(G)->O*_exp(G)}

\btm\label{TH-AM-O-star} Для любой компактно порожденной группы Штейна $G$ отображение
$$
\imath_G^\star:\mathcal{O}^\star(G)\to \mathcal{O}_{\exp}^\star(G)
$$
является оболочкой Аренса-Майкла алгебры $\mathcal{O}^\star(G)$:
 \beq\label{AM-O-star}
\Big(\mathcal{O}^\star(G)\Big)^\heartsuit\cong \mathcal{O}^\star_{\exp}(G)
 \eeq
\etm

\bpr 
 \begin{multline*}
\mathcal{O}^\star_{\exp}(G)=\eqref{O-exp-top}= \Bigg(\underset{\scriptsize\begin{matrix}\text{$D$
-- дуально}\\ \text{субмультипликативный}\\ \text{прямоугольник в $\mathcal
O(G)$}\end{matrix}}{\LCS\text{-}\injlim}\kern-15pt \C
D\Bigg)^\star=\eqref{O-exp-top-Ste}=
\Bigg(\underset{\scriptsize\begin{matrix}\text{$D$
-- дуально}\\ \text{субмультипликативный}\\ \text{прямоугольник в $\mathcal
O(G)$}\end{matrix}}{\Ste\text{-}\injlim}\kern-15pt \C
D\Bigg)^\star=\\=
 \underset{\scriptsize\begin{matrix}\text{$D$ -- дуально}\\
\text{субмультипликативный}\\ \text{прямоугольник в $\mathcal
O(G)$}\end{matrix}}{\Ste\text{-}\projlim } (\C
D)^\star=\eqref{O-exp^star-top-Ste}=
\underset{\scriptsize\begin{matrix}\text{$D$ -- дуально}\\
\text{субмультипликативный}\\ \text{прямоугольник в $\mathcal
O(G)$}\end{matrix}}{\LCS\text{-}\projlim }  (\C
D)^\star=\cite[(4.35)]{Akbarov-De-Gruyter-I}=\\= \underset{\scriptsize\begin{matrix}\text{$D$ -- дуально}\\
\text{субмультипликативный}\\ \text{прямоугольник в $\mathcal
O(G)$}\end{matrix}}{\LCS\text{-}\projlim } \mathcal{O}^\star(G)/ D^\circ=
\underset{\scriptsize\begin{matrix}\text{$\varDelta$ -- субмультипликативный}\\
\text{ромб в $\mathcal
O^\star(G)$}\end{matrix}}{\LCS\text{-}\projlim }\mathcal{O}^\star(G)/ \varDelta =(\text{теорема \ref{subm-in-trubk}(a)})=\\=
\underset{\scriptsize\begin{matrix}\text{$U$ -- субмультипликативная}\\
\text{окрестность нуля в $\mathcal
O^\star(G)$}\end{matrix}}{\LCS\text{-}\projlim } \mathcal{O}^\star(G)/ U=\eqref{A^heartsuit=proj-limit} =\Big(\mathcal{O}^\star(G)\Big)^\heartsuit
 \end{multline*}
 \epr

\subsection{Отображение $\imath_G:\mathcal{O}_{\exp}(G)\to \mathcal{O}(G)$ как оболочка
Аренса-Майкла}\label{SUBSEC:O_exp(G)->O(G)}

\paragraph{Случай конечно порожденной дискретной группы.}

\btm\label{TH-AM-O-diskr-gr} Пусть $G$ -- конечно порожденная дискретная группа. Тогда
отображение
$$
\imath_G:{\mathcal O}_{\exp}(G)\to {\mathcal O}(G)
$$
является оболочкой Аренса-Майкла алгебры ${\mathcal O}_{\exp}(G)$:
 \beq\label{AM-O-diskr-gr}
\Big({\mathcal O}_{\exp}(G)\Big)^\heartsuit\cong {\mathcal O}(G)
 \eeq
 \etm

Это утверждение мы докажем в несколько этапов. Напомним, что в \eqref{1_x} мы
условились обозначать символом $1_x$ характеристические функции одноточечных
подмножеств $\{x\}$ в $G$:
 \beq\label{1_x}
1_x(y)=\begin{cases}1, & y=x \\ 0& y\ne x\end{cases}
 \eeq
(из-за дискретности $G$, функцию $1_x$ можно считать элементом обеих алгебр
$\mathcal O(G)$ и ${\mathcal O}_{\exp}(G)$).

\blm\label{LM_1_x} Функции $\{1_x;\; x\in G\}$ образуют базис в топологических
векторных пространствах ${\mathcal O}(G)$ и ${\mathcal O}_{\exp}(G)$: для
всякой функции $u\in {\mathcal O}(G)$ ($u\in {\mathcal O}_{\exp}(G)$)
справедливо равенство
 \beq\label{series-in-H}
u=\sum_{x\in G} u(x)\cdot 1_x
 \eeq
где ряд справа сходится в ${\mathcal O}(G)$ (${\mathcal O}_{\exp}(G)$), а его
коэффициенты непрерывно зависят от $u\in {\mathcal O}(G)$ ($u\in {\mathcal
O}_{\exp}(G)$).\elm
 \bpr Для пространства ${\mathcal O}(G)$ это очевидно, потому что в случае дискретной группы $G$
оно совпадает с пространством $\C^G$ всех функций на $G$. Докажем это для
${\mathcal O}_{\exp}(G)$: если $u\in {\mathcal O}_{\exp}(G)$, то, подобрав
мажорирующий полухарактер $f:G\to\R_+$,
$$
|u(x)|\le f(x),\qquad x\in G
$$
мы получим, что частичные суммы ряда \eqref{series-in-H} содержатся в
прямоугольнике $f^\text{\BSQ}$, поэтому ряд \eqref{series-in-H} сходится (не
только в ${\mathcal O}(G)$, но и) в ${\mathcal O}_{\exp}(G)$. С другой стороны,
каждый коэффициент $u(x)$ непрерывно зависит от $u$, если $u$ бегает по
прямоугольнику вида $f^\text{\BSQ}$. По определению топологии в ${\mathcal
O}_{\exp}(G)$, это означает, что $u(x)$ непрерывно зависит от $u$, когда $u$
бегает по ${\mathcal O}_{\exp}(G)$.
 \epr

\blm\label{sum-q(1_x)-f(x)<infty} Если $G$ -- дискретная конечно порожденная
группа, то для всякой непрерывной полунормы $q:{\mathcal O}_{\exp}(G)\to\R_+$ и
любого полухарактера $f:G\to[1;+\infty)$ числовое семейство $\{f(x)\cdot
q(1_x);\; x\in G\}$ суммируемо:
$$
\sum_{x\in G} f(x)\cdot q(1_x)<\infty
$$
\elm
 \bpr
Пусть $T$ -- абсолютно выпуклый компакт в ${\mathcal O}_{\exp}^\star(G)$,
соответствующий полунорме $q$:
 $$
q(u)=\sup_{\alpha\in T}|\alpha(u)|
 $$
Всякий прямоугольник $f^\text{\BSQ}$ является компактом в ${\mathcal
O}_{\exp}^\star(G)$, поэтому
 \begin{multline*}
\infty>\sup_{u\in f^\text{\BSQ}}\sup_{\alpha\in T}
|\alpha(u)|=\eqref{series-in-H}=\sup_{u\in f^\text{\BSQ}}\sup_{\alpha\in T}
\Big|\alpha\Big(\sum_{x\in G} u(x)\cdot 1_x\Big)\Big|=\sup_{u\in
f^\text{\BSQ}}\sup_{\alpha\in T}  \Big|\sum_{x\in G} u(x)\cdot
\alpha(1_x)\Big|\ge\\ \ge \sup_{\alpha\in T} \Big|\sum_{x\in G}
\underbrace{\frac{f(x)\cdot\overline{\alpha(1_x)}}{|\alpha(1_x)|}}_{\scriptsize\begin{matrix}\uparrow\\
\text{одно из значений}\\ \text{$u\in f^\text{\BSQ}$}\end{matrix}}\cdot
\alpha(1_x)\Big|= \sup_{\alpha\in T} \sum_{x\in G} f(x)\cdot |\alpha(1_x)|\ge
\sup_{\alpha\in T} \sup_{x\in G} f(x)\cdot |\alpha(1_x)|=\\=\sup_{x\in G}
f(x)\cdot \sup_{\alpha\in T} |\alpha(1_x)|=\sup_{x\in G} f(x)\cdot q(1_x)
  \end{multline*}
Мы получили, что для всякого полухарактера $f:G\to[1;+\infty)$
$$
\sup_{x\in G} f(x)\cdot q(1_x)<\infty
$$
Если теперь взять какое-нибудь конечное множество $K$, порождающее $G$,
$$
\bigcup_{n=1}^\infty K^n=G
$$
и определить полухарактер $g:G\to[1;+\infty)$ формулой
$$
g(x)=R^n\quad\Longleftrightarrow\quad x\in K^n\setminus K^{n-1}
$$
где $R$ -- какое-нибудь число, большее мощности множества $K$:
$$
R>\card K,
$$
то, поскольку произведение $g\cdot f$ также будет полухарактером, получаем:
$$
\sup_{x\in G} \Big[ g(x)\cdot f(x)\cdot q(1_x)\Big]<\infty
$$
$$
\Downarrow
$$
$$
\exists C>0\qquad \forall x\in G \qquad f(x)\cdot q(1_x)\le \frac{C}{g(x)}
$$
$$
\Downarrow
$$
 \begin{multline*}
\sum_{x\in G} f(x)\cdot q(1_x)\le \sum_{x\in G} \frac{C}{g(x)}=
\sum_{n=1}^\infty \sum_{x\in K^n\setminus K^{n-1}} \frac{C}{g(x)}=
\sum_{n=1}^\infty \sum_{x\in K^n\setminus K^{n-1}} \frac{C}{R^n}\le
\sum_{n=1}^\infty  \frac{C\cdot \card (K^n)}{R^n}\le\\ \le \sum_{n=1}^\infty
\frac{C\cdot (\card K)^n}{R^n}=C\cdot \sum_{n=1}^\infty \left(\frac{\card
K}{R}\right)^n<\infty
  \end{multline*}
 \epr

Если $q:{\mathcal O}_{\exp}(G)\to\R_+$ -- непрерывная полунорма на ${\mathcal
O}_{\exp}(G)$, то ее {\it носителем} условимся называть множество
 \beq
\supp(q)=\{x\in G:\; q(1_x)\ne 0\}
 \eeq

\blm\label{LM-supp-q} Если $G$ -- дискретная конечно порожденная группа, то для
всякой субмультипликативной непрерывной полунормы $q:{\mathcal O}_{\exp}(G)\to\R_+$
 \bit{
\item[(a)] носитель $\supp(q)$ является конечным множеством:
\beq\label{supp-q-1}
\card \supp(q)<\infty,
\eeq
\item[(b)] для любой точки $x\in\supp(q)$ значение полунормы $q$ на функции
$1_x$ не меньше единицы:
\beq\label{supp-q-2}
q(1_x)\ge 1.
\eeq
 }\eit
\elm \bpr Сначала докажем (b). При $x\in\supp(q)$, то есть при $q(1_x)>0$,
получаем логическую цепочку:
$$
1_x=1_x^2\quad\Longrightarrow\quad q(1_x)=q(1_x^2)\le
q(1_x)^2\quad\Longrightarrow\quad 1\le q(1_x)
$$
Теперь (a). Поскольку тождественная единица $f(x)=1$ является полухарактером на
$G$, по лемме \ref{sum-q(1_x)-f(x)<infty} числовое семейство $\{q(1_x);\; x\in
G\}$ должно быть суммируемо:
$$
\sum_{x\in G} q(1_x)<\infty
$$
С другой стороны, по уже доказанному условию (b), все ненулевые слагаемые в
этом ряду оцениваются снизу единицей. Значит, их должно быть конечное число.
 \epr

\blm\label{LM:u|_supp-q=0=>q(u)=0}
Если функция $u\in{\mathcal O}_{\exp}(G)$ обнуляется на носителе непрерывной субмультипликативной полунормы $q:{\mathcal O}_{\exp}(G)\to\R_+$,  то $q$ обнуляется на $u$:
\beq\label{u|_supp-q=0=>q(u)=0}
u\Big|_{\supp q}=0 \quad\Rightarrow\quad q(u)=0.
\eeq
\elm
\bpr
Пусть $K=\supp q$ (по лемме \ref{LM-supp-q} это конечное множество). Если $u\in{\mathcal O}_{\exp}(G)$ обнуляется на $K$, то в разложении \eqref{series-in-H} все слагаемые с $x\in K$ нулевые:
$$
u=\sum_{x\in G} u(x)\cdot 1_x=\sum_{x\notin K} u(x)\cdot 1_x.
$$
Это значит, что
$$
u=\lim_{N\to\infty, \ N\cap K=\varnothing}\sum_{x\in N} u(x)\cdot 1_x.
$$
Отсюда
\begin{multline*}
q(u)=q\left(\lim_{N\to\infty, \ N\cap K=\varnothing}\sum_{x\in N} u(x)\cdot 1_x\right)=
\lim_{N\to\infty, \ N\cap K=\varnothing}q\left(\sum_{x\in N} u(x)\cdot 1_x\right)\le \\ \le
\lim_{N\to\infty, \ N\cap K=\varnothing}\sum_{x\in N} u(x)\cdot q(1_x)=0.
\end{multline*}
\epr

\bpr[Доказательство теоремы \ref{TH-AM-O-diskr-gr}] Пусть $G$ -- произвольная конечно порожденная дискретная группа и $q:{\mathcal O}_{\exp}(G)\to\R_+$ -- субмультипликативная непрерывная
полунорма. По лемме \ref{LM-supp-q} ее носитель $K=\supp q$ конечен.
Положим
$$
C=\sum_{x\in K} q(1_x)
$$
(в силу \eqref{supp-q-1}, сумма конечна, а в силу \eqref{supp-q-2}, $C\ge 1$) и обозначим
$$
p(u)=C\cdot \max_{x\in K}\abs{u(x)},\qquad u\in{\mathcal O}(G).
$$
Поскольку $C\ge 1$, это будет субмультипликативная полунорма на ${\mathcal O}(G)$. Она будет мажорировать $q$ на ${\mathcal O}_{\exp}(G)$, потому что для всякого $u\in {\mathcal O}_{\exp}(G)$ мы получим
\begin{multline*}
q(u)=\eqref{series-in-H}=q\left(\sum_{x\in G} u(x)\cdot 1_x\right)\le
q\left(\sum_{x\in K} u(x)\cdot 1_x\right)+\underbrace{q\left(\sum_{x\notin K} u(x)\cdot 1_x\right)}_{\scriptsize\begin{matrix}\| \put(2,0){\eqref{u|_supp-q=0=>q(u)=0}}\\ 0\end{matrix}} =\\=
q\left(\sum_{x\in K} u(x)\cdot 1_x\right)=\sum_{x\in G} \abs{u(x)}\cdot q(1_x)=\sum_{x\in K} \abs{u(x)}\cdot q(1_x)\le
C\cdot \max_{x\in K}\abs{u(x)}=p(u).
\end{multline*}
По предложению \ref{kriterij-A-M}, мы получаем, что отображение $\imath_G:{\mathcal O}_{\exp}(G)\to {\mathcal O}(G)$ --- оболочка Аренса---Майкла.
 \epr

\paragraph{Случай линейной группы с конечным числом связных компонент}
\label{SUBSEC:O_exp(G)->O(G)}

Следующее утверждение было доказано в \cite[Theorem 3.12]{ArHRC}:

\btm\label{TH-AM-O-kon-rassh-svyaz-lin-gruppy} Пусть $G$ --- линейная комплексная группа с конечным числом связных компонент. Тогда отображение
$$
\imath_G:\mathcal{O}_{\exp}(G)\to \mathcal{O}(G)
$$
является оболочкой Аренса-Майкла алгебры $\mathcal{O}_{\exp}(G)$:
 \beq\label{AM-O}
\Big(\mathcal{O}_{\exp}(G)\Big)^\heartsuit\cong \mathcal{O}(G)
 \eeq
 \etm

\subsection{Преобразование Фурье как оболочка Аренса---Майкла}
\label{SUBSEC:Fourier}

Если $G$ -- абелева компактно порожденная группа Штейна, то всякий ее
голоморфный характер $\chi:G\to\C^\times$ будет, как нетрудно сообразить,
голоморфной функцией на $G$. То есть двойственную комплексную группу
$G^\bullet$ можно представлять себе как подгруппу в группе обратимых элементов
алгебры $\mathcal{O}(G)$ голоморфных функций на $G$:
$$
G^\bullet\subset\mathcal{O}(G)
$$
Переходя по теореме \ref{PONT-AB} к двойственным объектам, мы получаем, что
сама группа $G$ вкладывается преобразованием $\i_G$ в группу обратимых
элементов алгебры $\mathcal{O}(G^\bullet)$ голоморфных функций на $G^\bullet$:
$$
\i_G:G\to G^{\bullet\bullet}\subset\mathcal{O}(G^\bullet).
$$
С другой стороны, $G$, очевидно, вкладывается (в виде дельта-функционалов) в
алгебру $\mathcal{O}^\star(G)$:
$$
\delta: G\to \mathcal{O}^\star(G)\qquad (x\mapsto\delta^x).
$$
В силу \cite[Theorem 5.1.22]{Akbarov-De-Gruyter-I} отсюда следует, что существует единственный
гомоморфизм стереотипных алгебр
$$
{\mathscr F}_G:\mathcal{O}^\star(G)\to \mathcal{O}(G^\bullet),
$$
замыкающий диаграмму
$$
\begin{diagram}
\node[2]{G} \arrow{se,t}{\i_G}\arrow{sw,t}{\delta}
\\
\node{\mathcal{O}^\star(G)}\arrow[2]{e,b,--}{{\mathscr F}_G} \node[2]{\mathcal{O}(G^\bullet)}
\end{diagram}
$$
(в этом состоит свойство $\mathcal{O}^\star(G)$ быть групповой алгеброй).
Гомоморфизм ${\mathscr F}_G:\mathcal{O}^\star(G)\to \mathcal{O}(G^\bullet)$
естественно называть (обратным) {\it преобразованием
Фурье}\index{преобразование!Фурье} на группе Штейна $G$, потому что он явно
задается формулой, которой определяется (обратное) преобразование Фурье для мер
и обобщенных функций \cite[31.2]{Hewitt-Ross-2}:
 \beq\label{alpha-chi=w-chi}
\overbrace{{\mathscr F}_G(\alpha)(\chi)}^{\scriptsize \begin{matrix}
\text{значение функции ${\mathscr F}_G(\alpha)\in \C^{G^\bullet}$}\\
\text{в точке $\chi\in G^\bullet$} \\ \downarrow \end{matrix}}\kern-35pt=\kern-50pt\underbrace{\alpha(\chi)}_{\scriptsize \begin{matrix}\uparrow \\
\text{действие функционала $\alpha\in\mathcal{O}^\star(G)$}\\
\text{на функцию $\chi\in G^\bullet\subseteq \mathcal{O}(G)$ }\end{matrix}}
\kern-50pt ,\qquad (\chi\in G^\bullet,\quad \alpha\in \mathcal{O}^\star(G))
 \eeq

\btm\label{O-exp^star=O} Для всякой абелевой компактно порожденной группы
Штейна $G$ ее преобразование Фурье
$$
{\mathscr F}_G:\mathcal{O}^\star(G)\to \mathcal{O}(G^\bullet),
$$
является:
 \bit{
\item[(a)] гомоморфизмом жестких алгебр Хопфа и

\item[(b)] оболочкой Аренса-Майкла алгебры $\mathcal{O}^\star(G)$.
 }\eit
Как следствие, справедливы изоморфизмы алгебр Хопфа-Фреше:
 \beq\label{O_exp^*_G=O-G^bullet}
\mathcal{O}_{\exp}^\star(G)\cong \Big(\mathcal{O}^\star(G)\Big)^\heartsuit
\cong \mathcal{O}(G^\bullet)
 \eeq
а диаграмма рефлексивности для $G$ принимает вид
 \beq\label{chetyrehugolnik-O-O*-Abel}
 \xymatrix @R=1.pc @C=2.pc
 {
 \mathcal{O}^\star(G)
 &\ar@{|->}[rrr]^{{\mathscr F}_G}_{\eqref{O_exp^*_G=O-G^bullet}} & &
 &&
 \mathcal{O}(G^\bullet)
 \\
 & & & &&
 \ar@{|->}[d]^{\star}
 \\
 \ar@{|->}[u]^{\star}
 & & & & &
 \\
 \mathcal{O}(G)
 & & & &
 \ar@{|->}[lll]_{{\mathscr F}_{G^\bullet}}^{\eqref{O_exp^*_G=O-G^bullet}} &
  \mathcal{O}^\star(G^\bullet)
 }
 \eeq

 \etm

Как и теорема \ref{PONT-AB}, это утверждение доказывается последовательным
рассмотрением случаев $G=\C, \;\C^\times,\; \Z$ и случая конечной абелевой
группы $G=F$. В оставшейся части этого параграфа до ``диаграммы вложения'' мы
займемся этим.

\paragraph{Конечная абелева группа.}\label{ex-finite}

Как уже говорилось в \ref{SEC:stein-groups}\ref{SUBSEC-lin-groups}, всякую
конечную группу $G$ можно считать линейной комплексной группой Ли (нулевой
размерности), на которой любая функция считается голоморфной, причем в первом
примере \ref{SEC-O_exp(G)}\ref{SUBSEC:algebra-O_exp(G)} отмечалось, что более
того, любая функция на конечной группе является голоморфной функцией
экспоненциального типа, и поэтому алгебры $\mathcal{O}_{\exp}(G)$, $\mathcal{O}(G)$ и $\C^G$ в этом случае совпадают:
$$
\mathcal{O}_{\exp}(G)=\mathcal{O}(G)=\C^G
$$
Если $G$ вдобавок коммутативна, то иллюстрируемая нами здесь теорема
\ref{O-exp^star=O} превращается в формально более сильное утверждение:

\bprop Если $G$ -- конечная абелева группа, то формула \eqref{alpha-chi=w-chi}
устанавливает изоморфизм алгебр Хопфа
 \beq\label{finite-groups}
\mathcal{O}_{\exp}^\star(G)=\mathcal{O}^\star(G)=
\C_G\cong\C^{G^\bullet}=\mathcal{O}(G^\bullet)= \mathcal{O}_{\exp}(G^\bullet)
 \eeq
 \eprop

\paragraph{Комплексная плоскость $\C$.}\label{ex-C} Пусть для всякого
$\lambda\in\C$ символ $\chi_\lambda$ обозначает характер на группе $\C$,
заданный формулой
$$
\chi_\lambda(t)=e^{\lambda\cdot t}
$$
Отображение $\lambda\in\C\mapsto\chi_\lambda\in\C^\bullet$ является
изоморфизмом комплексных групп
$$
\C\cong\C^\bullet
$$
а формула \eqref{alpha-chi=w-chi} при таком изоморфизме приобретает вид
 \beq\label{alpha-chi-lambda=w-lambda}
{\mathscr F}_{\C}(\alpha)(\lambda)=\alpha(\chi_\lambda),\qquad \lambda\in \C\qquad
(\alpha\in \mathcal{O}_{\exp}^\star(\C),\quad w\in \mathcal{O}(\C))
 \eeq
(мы обозначаем это отображение тем же символом ${\mathscr F}_{\C}$, хотя формально оно
представляет собой композицию отображения \eqref{alpha-chi=w-chi} с
отображением $\lambda\mapsto\chi_\lambda$). В результате теорема
\ref{O-exp^star=O}, применительно к группе $G=\C$ превращается в

\bprop\label{PROP-H*(C)->H(C)} Формула \eqref{alpha-chi-lambda=w-lambda}
определяет гомоморфизм жестких стереотипных алгебр Хопфа
$$
{\mathscr F}_{\C}:\mathcal{O}^\star(\C)\to \mathcal{O}(\C)
$$
являющийся оболочкой Аренса-Майкла алгебры $\mathcal{O}^\star(\C)$, и поэтому
устанавливающий изоморфизм алгебр Хопфа-Фреше
 \beq\label{O_exp*(C)=O(C)}
\mathcal{O}_{\exp}^\star(\C)\cong \mathcal{O}(\C)
 \eeq
\eprop

Для доказательства нам понадобится

\blm\label{LM-||alpha||_C-v-C} Полунормы вида
 \beq\label{||alpha||_C-v-C}
||\alpha||_C=\sum_{k\in\N} |\alpha_k|\cdot C^k,\qquad C\ge 0,
 \eeq
(частный случай полунорм \eqref{|alpha|_r-v-C}, когда $r_k=\frac{C^k}{k!}$)
образуют фундаментальную систему в множестве всех субмультипликативных
непрерывных полунорм на $\mathcal{O}^\star(\C)$. \elm

\bpr Как мы уже отмечали в
\ref{SEC:stein-groups}\ref{Examples-of-Stein-groups}, операции умножения в
$\mathcal{O}(\C)$ и $\mathcal{O}^\star (\C)$ задаются одинаковыми формулами
на рядах \eqref{umnozhenie-v-O(C)-i-O^star(C)}. Поэтому субмультипликативность
полунорм \eqref{||alpha||_C-v-C} можно считать доказанной, после того, как в
примере \ref{EX-polunormy-v-H(C)} тот же факт мы проверили для полунорм
\eqref{polunormy-v-H(C)}, задаваемых той же формулой на рядах.

Покажем, что полунормы \eqref{||alpha||_C-v-C} образуют фундаментальную систему
среди всех субмультипликативных полунорм на $\mathcal{O}^\star(\C)$. Это
делается так же, как в предложении \ref{PROP:O*(C)-Hopf}. Пусть $p$ --
субмультипликативная непрерывная полунорма:
$$
p(\alpha*\beta)\le p(\alpha)\cdot p(\beta)
$$
Положим
$$
r_k=\frac{1}{k!}p(\tau_k)
$$
Тогда
$$
(k+l)!\cdot r_{k+l}=p(\tau_{k+l})=p(\tau_k*\tau_l)\le p(\tau_k)\cdot
p(\tau_l)=(k!\cdot r_k)\cdot (l!\cdot r_l)
$$
Если обозначить $A_k=r_k\cdot k!$, то для этой последовательности получается
рекуррентное неравенство $A_{k+1}\le A_k\cdot A_1$, из которого следует $A_k\le
C^k$, для $C=A_1$. Это в свою очередь влечет неравенства
 $$
r_k\le\frac{C^k}{k!}
 $$
Теперь используя те же рассуждения, что и в доказательстве предложения
\ref{PROP:O*(C)-Hopf}, получаем:
$$
p(\alpha)\le\eqref{p(alpha)-le-|||alpha|||_r}\le |||\alpha|||_r=\sum_{k\in\N}
r_k\cdot |\alpha_k|\cdot k!\le \sum_{k\in\N} |\alpha_k|\cdot C^k=||\alpha||_C
$$ \epr

\bpr[Доказательство предложения \ref{PROP-H*(C)->H(C)}] Заметим сразу, что
отображение ${\mathscr F}_{\C}:\mathcal{O}^\star(\C)\to \mathcal{O}(\C)$,
определенное формулой \eqref{alpha-chi-lambda=w-lambda} непрерывно: в силу
непрерывности отображения $\lambda\in\C\mapsto \chi_\lambda\in \mathcal{O}(\C)$, всякий компакт $T$ в $\C$ превращается в компакт
$\{\chi_\lambda;\;\lambda\in T\}$ в $\mathcal{O}(\C)$, поэтому если
направленность функционалов $\alpha_i$ сходится к нулю в $\mathcal{O}^\star(\C)$, то для всякого компакта $T$ в $\C$ получаем
$$
{\mathscr F}_{\C}(\alpha_i)(\lambda)=\alpha_i(\chi_\lambda)\underset{\lambda\in
T}{\rightrightarrows} 0,\qquad i\to\infty
$$
То есть функции ${\mathscr F}_{\C}(\alpha_i)$ стремятся к нулю в $\mathcal{O}(\C)$.

Далее заметим, что отображение ${\mathscr F}_{\C}$ переводит функционалы $\tau_k$
(определенные формулой \eqref{def-tau^k}) в функции $t^k$ (определенные
формулой \eqref{def-t^k}):
$$
{\mathscr F}_{\C}(\tau_k)(\lambda)=\tau_k(\chi_\lambda)=\left(\frac{\d^k}{\d
x^k}e^{\lambda x}\right)\Bigg|_{t=0}=\lambda^k=t^k(\lambda)
$$
Отсюда и из непрерывности ${\mathscr F}_{\C}$ следует, что это отображение действует
на функционалы $\alpha$ заменой в разложении \eqref{razlozhenie-alpha-v-C}
мономов $\tau_k$ на мономы $t^k$:
$$
{\mathscr F}_{\C}(\alpha)={\mathscr F}_{\C}\left(\sum_{k=0}^\infty \alpha_k\cdot
\tau_k\right)=\sum_{k=0}^\infty \alpha_k\cdot
{\mathscr F}_{\C}(\tau_k)=\sum_{k=0}^\infty \alpha_k\cdot t^k
$$
Отсюда сразу следует все остальное.

1. Во-первых, отображение ${\mathscr F}_{\C}:\mathcal{O}^\star(\C)\to \mathcal{O}(\C)$ будет оболочкой Аренса-Майкла, потому что по лемме
\ref{LM-||alpha||_C-v-C}, всякая субмультипликативная непрерывная полунорма на
$\mathcal{O}^\star(\C)$ мажорируется полунормой вида \eqref{||alpha||_C-v-C},
которая в свою очередь продолжается отображением ${\mathscr F}_{\C}$ до полунормы
\eqref{polunormy-v-H(C)} на $\mathcal{O}(\C)$.

2. Во-вторых, отображение ${\mathscr F}_{\C}:\mathcal{O}^\star(\C)\to \mathcal{O}(\C)$ будет гомоморфизмом алгебр, потому что эти алгебры мы можем
представлять себе по формулам
\eqref{razlozhenie-u-v-C}-\eqref{razlozhenie-alpha-v-C} алгебрами степенных
рядов, в которых умножение задается обычными для степенных рядов правилами
\eqref{umnozhenie-v-O(C)-i-O^star(C)}, и ${\mathscr F}_{\C}$ тогда будет просто
вложением одной алгебры в другую, более широкую.

3. Чтобы доказать, что отображение ${\mathscr F}_{\C}:\mathcal{O}^\star(\C)\to
\mathcal{O}(\C)$ -- изоморфизм коалгебр, заметим, что сопряженное отображение
$$
({\mathscr F}_{\C})^\star:\mathcal{O}^\star(\C)\to(\mathcal{O}^\star(\C))^\star=\mathcal{O}(\C)^{\star\star}
$$
с точностью до изоморфизма $\i_{\mathcal{O}(\C)}:\mathcal{O}(\C)\cong
\mathcal{O}(\C)^{\star\star}$ совпадает с ${\mathscr F}_{\C}$:
 \beq\label{sharp_C^star=i_H(C)-circ-sharp_C}
\begin{diagram}
\node{\mathcal{O}^\star(\C)}
\arrow[2]{e,t}{({\mathscr F}_{\C})^\star}\arrow{se,b}{{\mathscr F}_{\C}} \node[2]{\mathcal{O}(\C)^{\star\star}}
\\
\node[2]{\mathcal{O}(\C)}\arrow{ne,b}{\i_{\mathcal{O}(\C)}}
\end{diagram}
 \eeq
Это следует из формулы
 \beq\label{sharp_C-delta_x}
{\mathscr F}_{\C}(\delta^x)=\chi_x
 \eeq
Действительно,
$$
{\mathscr F}_{\C}(\delta^x)(\lambda)=\delta^x(\chi_\lambda)=\chi_\lambda(x)=e^{\lambda
x}=\chi_x(\lambda)
$$
Теперь получаем:
$$
{\mathscr F}_{\C}^\star(\delta^a)(\delta^b)=\delta^a({\mathscr F}_{\C}(\delta^b))=\eqref{sharp_C-delta_x}
=\delta^a(\chi_b)={\mathscr F}_{\C}(\delta^a)(b)=\delta^b({\mathscr F}_{\C}(\delta^a))=
=\i_{\mathcal{O}(\C)}({\mathscr F}_{\C}(\delta^a))(\delta^b)=(\i_{\mathcal{O}(\C)}\circ\; {\mathscr F}_{\C})(\delta^a)(\delta^b)
$$
Это верно для любых $a,b\in\C$, а дельта-функционалы полны в $\mathcal{O}^\star$, поэтому
$$
({\mathscr F}_{\C})^\star=\i_{\mathcal{O}(\C)}\circ\;{\mathscr F}_{\C}
$$
то есть диаграмма \eqref{sharp_C^star=i_H(C)-circ-sharp_C} коммутативна. Теперь
мы можем заметить, что про ${\mathscr F}_{\C}$ мы уже доказали, что это гомоморфизм
алгебр, а для $\i_{\mathcal{O}(\C)}$ это очевидно, значит мы получаем, что
$({\mathscr F}_{\C})^\star$ -- также гомоморфизм алгебр, и это означает, что
${\mathscr F}_{\C}$ -- гомоморфизм коалгебр.

4. Теперь остается проверить, что ${\mathscr F}_{\C}$ сохраняет антипод:
$$
\sigma_{\mathcal{O}(\C)}(\chi_\lambda)(x)=\chi_\lambda(-x)=e^{-\lambda
x}=(e^{\lambda x})^{-1}=\chi_\lambda(x)^{-1}
$$
$$
\Downarrow
$$
$$
\sigma_{\mathcal{O}(\C)}(\chi_\lambda)=\chi_\lambda^{-1}
$$
$$
\Downarrow
$$
$$
{\mathscr F}_{\C}(\sigma_{\mathcal{O}^\star(\C)}(\alpha))(\lambda)= (\sigma_{\mathcal{O}^\star(\C)}(\alpha))(\chi_\lambda)= (\alpha\circ\sigma_{\mathcal{O}(\C)})(\chi_\lambda)= \alpha(\sigma_{\mathcal{O}(\C)}(\chi_\lambda))=
\alpha(\chi_\lambda^{-1})={\mathscr F}_{\C}(\alpha)(-\lambda)=\sigma_{\mathcal{O}(\C)}({\mathscr F}_{\C}(\alpha))(\lambda)
$$
$$
\Downarrow
$$
$$
{\mathscr F}_{\C}(\sigma_{\mathcal{O}^\star(\C)}(\alpha))=\sigma_{\mathcal{O}(\C)}({\mathscr F}_{\C}(\alpha))
$$
\epr

\paragraph{Комплексная окружность $\C^\times$.}\label{ex-C-x} Пусть для всякого
$n\in\Z$ символ $z^n$ обозначает характер на группе $\C^\times$, заданный
формулой
 $$
z^n(t)=t^n
 $$
 Отображение
$n\in\Z\mapsto z^n\in(\C^\times)^\bullet$ является изоморфизмом комплексных
групп
 $$ \Z\cong(\C^\times)^\bullet
 $$
а формула \eqref{alpha-chi=w-chi} при таком изоморфизме приобретает вид
 \beq\label{alpha-chi-lambda=w-lambda-C-x}
{\mathscr F}_{\C^\times}(\alpha)(n)=\alpha(z^n),\qquad n\in\Z\qquad (\alpha\in \mathcal{O}_{\exp}^\star(\C^\times))
 \eeq
(как и в предыдущем примере мы обозначаем получаемое отображение тем же
символом ${\mathscr F}_{\C^\times}$, хотя формально оно представляет собой композицию отображения
\eqref{alpha-chi=w-chi} с отображением $n\mapsto z^n$). В результате теорема
\ref{O-exp^star=O}, применительно к группе $G=\C^\times$ превращается в

\bprop\label{PROP-H*(C-x)->H(Z)} Формула \eqref{alpha-chi-lambda=w-lambda-C-x}
определяет гомоморфизм жестких стереотипных алгебр Хопфа
$$
{\mathscr F}_{\C^\times}:\mathcal{O}^\star(\C^\times)\to \mathcal{O}(\Z)
$$
являющийся оболочкой Аренса-Майкла алгебры $\mathcal{O}^\star(\C^\times)$, и
поэтому устанавливающий изоморфизм алгебр Хопфа-Фреше
 \beq\label{O_exp*(C^x)=O(Z)}
\mathcal{O}_{\exp}^\star(\C^\times)\cong \mathcal{O}(\Z)=\C^\Z
 \eeq
\eprop

Нам понадобится

\blm\label{LM-||alpha||_C-v-C-x} Полунормы вида
 \beq\label{||alpha||_C-v-C-x}
||\alpha||_N=\sum_{|n|\le N}|\alpha_n|,\qquad N\in\N
 \eeq
-- частный случай полунорм \eqref{|alpha|_r-v-C-x}, когда
$$
r_n=\begin{cases}1,& |n|\le N \\ 0,& |n|>N \end{cases}
$$
-- образуют фундаментальную систему в множестве всех субмультипликативных
непрерывных полунорм на $\mathcal{O}^\star(\C^\times)$.
 \elm
\bpr Субмультипликативность полунорм \eqref{||alpha||_C-v-C-x} следует из
формулы для операции умножения в $\mathcal{O}^\star (\C^\times)$:
$$
||\alpha*\beta||_N=\sum_{|n|\le
N}|(\alpha*\beta)_n|=\eqref{umnozhenie-v-O(C-x)-i-O^star(C-x)}=\sum_{|n|\le
N}|\alpha_n\cdot\beta_n|\le \left(\sum_{|n|\le
N}|\alpha_n|\right)\cdot\left(\sum_{|n|\le
N}|\beta_n|\right)=||\alpha||_N\cdot||\beta||_N
$$
Покажем, что полунормы \eqref{|alpha|_r-v-C-x} образуют фундаментальную систему
среди всех субмультипликативных непрерывных полунорм на $\mathcal{O}^\star(\C^\times)$. Пусть $p$ -- субмультипликативная непрерывная полунорма:
$$
p(\alpha*\beta)\le p(\alpha)\cdot p(\beta)
$$
Положим
$$
r_n=p(\zeta_n)
$$
Тогда
$$
r_n=p(\zeta_n)=p(\zeta_n*\zeta_n)\le p(\zeta_n)\cdot p(\zeta_n)=r_n^2
$$
То есть $0\le r_n\le r_n^2$, а это возможно только если $r_n\ge 1$ или $r_n=0$.
Но по лемме \ref{LM-polunormy-v-H*(C-x)-1}, числа $r_n$ должны удовлетворять
условию \eqref{forall-R>0-sum_r_n-R^n<infty-x}, из которого следует в
частности, что $r_n\to 0$. Такое возможно только если все они, кроме конечного
набора, равны нулю:
$$
\exists N\in\N\quad \forall n\in\Z\quad |n|>N\quad\Longrightarrow\quad r_n=0
$$
Положим $M=\max_{n}r_n$, тогда по лемме \ref{LM-polunormy-v-H*(C-x)-1}
получаем:
$$
p(\alpha)\le |||\alpha|||_r=\sum_{n\in\Z} r_n\cdot |\alpha_n|= \sum_{|n|\le N}
r_n\cdot |\alpha_n|\le \sum_{|n|\le N} M\cdot|\alpha_n|=M\cdot ||\alpha||_N
$$ \epr

\bpr[Начало доказательство предложения \ref{PROP-H*(C-x)->H(Z)}] Заметим, что
отображение ${\mathscr F}_{\C^\times}:\mathcal{O}^\star(\C^\times)\to \mathcal{O}(\Z)$, определенное формулой \eqref{alpha-chi-lambda=w-lambda-C-x}
непрерывно: если направленность функционалов $\alpha_i$ сходится к нулю в
$\mathcal{O}^\star(\C^\times)$, то для всякого $n\in\Z$ получаем
$$
{\mathscr F}_{\C^\times}(\alpha_i)(n)=\alpha_i(z^n)\longrightarrow 0,\qquad i\to\infty
$$
Это означает, что ${\mathscr F}_{\C^\times}(\alpha_i)$ сходится к нулю в $\mathcal{O}(\Z)=\C^{\Z}$.

Далее заметим, что отображение ${\mathscr F}_{\C^\times}$ переводит функционалы
$\zeta_k$ в характеристические функции одноэлементных множеств в $\Z$:
$$
{\mathscr F}_{\C^\times}(\zeta_k)(n)=\zeta_k(z^n)=\eqref{zeta^k(z^n)}=\left\{\begin{matrix}1,& n=k \\
0,& n\ne k\end{matrix}\right\}=\eqref{DEF:1_n,n-in-Z}=1_k(n)
$$
Отсюда и из непрерывности ${\mathscr F}_{\C^\times}$ следует, что это отображение
действует на функционалы $\alpha$ заменой в разложении
\eqref{razlozhenie-alpha-v-C} мономов $\zeta_n$ на мономы $1_n$:
$$
{\mathscr F}_{\C^\times}(\alpha)={\mathscr F}_{\C^\times}\left(\sum_{n\in\Z} \alpha_n\cdot
\zeta_n\right)=\sum_{n\in\Z} \alpha_n\cdot
{\mathscr F}_{\C^\times}(\zeta_n)=\sum_{n\in\Z} \alpha_n\cdot 1_n
$$
Это в свою очередь влечет за собой большую часть предложения
\ref{PROP-H*(C-x)->H(Z)}.

1. Во-первых, отображение ${\mathscr F}_{\C^\times}:\mathcal{O}^\star(\C^\times)\to
\mathcal{O}(\Z)$ будет оболочкой Аренса-Майкла, потому что по лемме
\ref{LM-||alpha||_C-v-C-x}, всякая субмультипликативная непрерывная полунорма
на $\mathcal{O}^\star(\C^\times)$ мажорируется полунормой вида
\eqref{||alpha||_C-v-C-x}, которая в свою очередь продолжается отображением
${\mathscr F}_{\C}$ до непрерывной полунормы на $\mathcal{O}(\Z)=\C^{\Z}$.

2. Во-вторых, отображение ${\mathscr F}_{\C^\times}:\mathcal{O}^\star(\C^\times)\to
\mathcal{O}(\Z)$ будет гомоморфизмом алгебр
\label{homomorphism-sharp-C^times}, потому что в силу формулы умножения
\eqref{umnozhenie-v-O(C-x)-i-O^star(C-x)} в $\mathcal{O}^\star(\C^\times)$,
это пространство можно представлять себе как пространство двусторонних
последовательностей $\alpha_n$ с покоординатным умножением, которое
отображением ${\mathscr F}_{\C^\times}$ тождественно вкладывается в более широкое
пространство $\mathcal{O}(\Z)=\C^{\Z}$ всех двусторонних последовательностей с
покоординатным умножением.

3. Доказательство того, что отображение ${\mathscr F}_{\C^\times}:\mathcal{O}^\star(\C^\times)\to \mathcal{O}(\Z)$ является изоморфизмом коалгебр нам
придется отложить до следующего примера
(с.\pageref{end-of-proof-PROP-H*(C-x)->H(Z)}).

4. Проверим, что ${\mathscr F}_{\C^\times}$ сохраняет антипод:
$$
\sigma_{\mathcal{O}(\C^\times)}(z^n)(x)=z^n(x^{-1})=x^{-n}=z^{-n}(x)
$$
$$
\Downarrow
$$
$$
\sigma_{\mathcal{O}(\C^\times)}(z^n)=z^{-n}
$$
$$
\Downarrow
$$
\begin{multline*}
{\mathscr F}_{\C^\times}(\sigma_{\mathcal{O}^\star(\C^\times)}(\alpha))(n)= (\sigma_{\mathcal{O}^\star(\C^\times)}(\alpha))(z^n)= (\alpha\circ\sigma_{\mathcal{O}(\C^\times)})(z^n)=\\=
\alpha(\sigma_{\C^\times}(z^n))=\alpha(z^{-n})={\mathscr F}_{\C^\times}(\alpha)(-n)=\sigma_{\mathcal{O}(\Z)}({\mathscr F}_{\C^\times}(\alpha))(n)
\end{multline*}
$$
\Downarrow
$$
$$
{\mathscr F}_{\C^\times}(\sigma_{\mathcal{O}^\star(\C^\times)}(\alpha)) =\sigma_{\mathcal{O}(\Z)}({\mathscr F}_{\C^\times}(\alpha))
$$

\epr

\paragraph{Группа целых чисел $\Z$.}\label{ex-Z} Пусть для всякого $t\in\C^\times$
символ $\chi_t$ обозначает характер на группе $\Z$, заданный формулой
 \beq\label{chi_t(n)=t^n}
\chi_t(n)=t^n
 \eeq
Отображение $t\in\C^\times\mapsto\chi_t\in\Z^\bullet$ является изоморфизмом
комплексных групп
 $$
\C^\times\cong\Z^\bullet
 $$
а формула \eqref{alpha-chi=w-chi} при таком изоморфизме приобретает вид
 \beq\label{alpha-chi-lambda=w-lambda-Z}
{\mathscr F}_{\Z}(\alpha)(t)=\alpha(\chi_t),\qquad t\in\C^\times\qquad (\alpha\in \mathcal{O}_{\exp}^\star(\Z))
 \eeq
(мы обозначаем получаемое отображение тем же символом ${\mathscr F}_{\Z}$, хотя формально
оно представляет собой композицию отображения \eqref{alpha-chi=w-chi} с
отображением $t\mapsto \chi_t$). В результате теорема \ref{O-exp^star=O},
применительно к группе $G=\Z$ превращается в

\bprop\label{PROP-H*(Z)->H(C-x)} Формула \eqref{alpha-chi-lambda=w-lambda-Z}
определяет гомоморфизм жестких стереотипных алгебр Хопфа
$$
{\mathscr F}_{\Z}:\mathcal{O}^\star(\Z)\to \mathcal{O}(\C^\times)
$$
являющийся оболочкой Аренса-Майкла алгебры $\mathcal{O}^\star(\Z)$, и поэтому
устанавливающий изоморфизм алгебр Хопфа-Фреше
 \beq\label{O_exp*(Z)=O(C^x)}
\mathcal{O}_{\exp}^\star(\Z)\cong \mathcal{O}(\C^\times)
 \eeq
\eprop

Нам понадобится

\blm\label{LM-||alpha||_C-v-Z} Полунормы вида
 \beq\label{||alpha||_C-v-Z}
||\alpha||_C=\sum_{n\in\Z} |\alpha_n|\cdot C^{|n|},\qquad C\ge 1
 \eeq
(частный случай полунорм \eqref{|alpha|_r-v-Z}, когда $r_n=C^{|n|}$) образуют
фундаментальную систему в множестве всех субмультипликативных непрерывных
полунорм на $\mathcal{O}^\star(\Z)$.
 \elm
\bpr Субмультипликативность полунорм вида \eqref{||alpha||_C-v-Z} мы уже
отмечали в примере \ref{EX-polunormy-v-O(C^x)}. Покажем, что они образуют
фундаментальную систему среди всех субмультипликативных полунорм на $\mathcal{O}^\star(\Z)$. Пусть $p$ -- субмультипликативная непрерывная полунорма:
$$
p(\alpha*\beta)\le p(\alpha)\cdot p(\beta)
$$
Положим $r_n=p(\delta^n)$. Тогда
$$
r_{k+l}=p(\delta^{k+l})=p(\delta^k*\delta^l)\le p(\delta^k)\cdot
p(\delta^l)=r_k\cdot r_l
$$
Из этого рекуррентного соотношения следует, что
$$
r_n\le M\cdot C^{|n|}
$$
где $M=r_0$, $C=\max\{r_1;r_{-1}\}$, и теперь по лемме
\ref{LM-polunormy-v-H*(Z)} получаем:
$$
p(\alpha)\le\eqref{p(alpha)-le-|||alpha|||_r-O*(Z)}\le
|||\alpha|||_r=\sum_{n\in\Z} r_n\cdot |\alpha_n|\le \sum_{n\in\Z} M\cdot
C^{|n|}\cdot |\alpha_n|=M\cdot ||\alpha||_C
$$ \epr

\bpr[Доказательство предложения \ref{PROP-H*(Z)->H(C-x)}] Прежде всего,
заметим, что отображение ${\mathscr F}_{\Z}:\mathcal{O}^\star(\Z)\to \mathcal{O}(\C^\times)$, определенное формулой \eqref{alpha-chi-lambda=w-lambda-Z}
непрерывно: в силу непрерывности отображения $t\in\C^\times\mapsto \chi_t\in
\mathcal{O}(\Z)$, всякий компакт $T$ в $\C^\times$ превращается в компакт
$\{\chi_t;\; t\in T\}$ в $\mathcal{O}(\Z)$, поэтому если направленность
функционалов $\alpha_i$ сходится к нулю в $\mathcal{O}^\star(\Z)$, то для
всякого компакта $T$ в $\C^\times$ получаем
$$
{\mathscr F}_{\Z}(\alpha_i)(t)=\alpha_i(\chi_t)\underset{\lambda\in T}{\rightrightarrows}
0,\qquad i\to\infty
$$
Это означает, что ${\mathscr F}_{\Z}(\alpha_i)$ сходится к нулю в $\mathcal{O}(\C^\times)$.

Далее заметим, что отображение ${\mathscr F}_{\Z}$ переводит функционалы $\delta^n$ в
мономы $z^n$:
$$
{\mathscr F}_{\Z}(\delta^n)(t)=\delta^n(\chi_t)=\chi_t(n)=t^n=z^n(t)
$$
Отсюда и из непрерывности ${\mathscr F}_{\Z}$ следует, что это отображение действует
на функционалы $\alpha$ заменой в разложении \eqref{razlozhenie-alpha-v-Z}
мономов $\delta^n$ на мономы $z^n$:
$$
{\mathscr F}_{\Z}(\alpha)={\mathscr F}_{\Z}\left(\sum_{n\in\Z} \alpha_n\cdot
\delta^n\right)=\sum_{n\in\Z} \alpha_n\cdot
{\mathscr F}_{\Z}(\delta^n)=\sum_{n\in\Z} \alpha_n\cdot z^n
$$
Отсюда следует остальное.

1. Во-первых, отображение ${\mathscr F}_{\Z}:\mathcal{O}^\star(\Z)\to \mathcal{O}(\C^\times)$ будет оболочкой Аренса-Майкла, потому что по лемме
\ref{LM-||alpha||_C-v-Z}, всякая субмультипликативная полунорма на $\mathcal{O}^\star(\Z)$ мажорируется полунормой вида \eqref{||alpha||_C-v-Z}, которая в
свою очередь продолжается отображением ${\mathscr F}_{\Z}$ до полунормы
\eqref{polunormy-v-H(C-x)} на $\mathcal{O}(\C^\times)$.

2. Во-вторых, отображение ${\mathscr F}_{\Z}:\mathcal{O}^\star(\Z)\to \mathcal{O}(\C^\times)$ будет гомоморфизмом алгебр, потому что эти алгебры мы можем
представлять себе в соответствии с формулами
\eqref{umnozhenie-v-O(Z)-i-O^star(Z)}-\eqref{umnozhenie-v-O(C-x)-i-O^star(C-x)}
алгебрами степенных рядов, в которых умножение задается обычными для степенных
рядов правилами, и ${\mathscr F}_{\C}$ тогда будет просто вложением одной алгебры в
другую, более широкую.

3. Чтобы доказать, что отображение ${\mathscr F}_{\Z}:\mathcal{O}^\star(\Z)\to
\mathcal{O}(\C^\times)$ -- изоморфизм коалгебр, заметим, что сопряженное
отображение
$$
({\mathscr F}_{\Z})^\star:\mathcal{O}^\star(\C^\times)\to(\mathcal{O}^\star(\Z))^\star=\mathcal{O}(\Z)^{\star\star}
$$
с точностью до изоморфизма $\i_{\mathcal{O}(\Z)}:\mathcal{O}(\Z)\cong
\mathcal{O}(\Z)^{\star\star}$ совпадает с ${\mathscr F}_{\C^\times}$:
 \beq\label{sharp_Z^star=i_H(Z)-circ-sharp_C-x}
\begin{diagram}
\node{\mathcal{O}^\star(\C^\times)}
\arrow[2]{e,t}{({\mathscr F}_{\Z})^\star}\arrow{se,b}{{\mathscr F}_{\C^\times}}
\node[2]{\mathcal{O}(\Z)^{\star\star}}
\\
\node[2]{\mathcal{O}(\Z)}\arrow{ne,b}{\i_{\mathcal{O}(\Z)}}
\end{diagram}
 \eeq
Это следует из формулы
 \beq\label{sharp_Z-C-x-delta_x}
\delta^t({\mathscr F}_{\Z}(\delta^n))=t^n=\delta^n({\mathscr F}_{\C^\times}(\delta^t)),\qquad
t\in\C^\times,\quad n\in\Z
 \eeq
Действительно,
$$
\delta^t({\mathscr F}_{\Z}(\delta^n))={\mathscr F}_{\Z}(\delta^n)(t)=\delta^n(\chi_t)=\chi_t(n)=t^n
$$
и
$$
\delta^n({\mathscr F}_{\C^\times}(\delta^t))={\mathscr F}_{\C^\times}(\delta^t)(n)=\delta^t(z^n)=z^n(t)=t^n
$$
Теперь получаем: для $t\in\C^\times$ и $n\in\Z$
$$
({\mathscr F}_{\Z})^\star(\delta^t)(\delta^n)=\delta^t({\mathscr F}_{\Z}(\delta^n))=\eqref{sharp_Z-C-x-delta_x}
=\delta^n({\mathscr F}_{\C^\times}(\delta^t))=\i_{\mathcal{O}(\Z)}({\mathscr F}_{\C^\times}(\delta^t))(\delta^n)=(\i_{\mathcal{O}(\Z)}\circ\;
{\mathscr F}_{\C^\times})(\delta^t)(\delta^n)
$$
Это верно для любых $t\in\C^\times$ и $n\in\Z$, а дельта-функционалы плотны в
$\mathcal{O}^\star$, поэтому
$$
({\mathscr F}_{\Z})^\star=\i_{\mathcal{O}(\Z)}\circ\;{\mathscr F}_{\C^\times}
$$
то есть диаграмма \eqref{sharp_Z^star=i_H(Z)-circ-sharp_C-x} коммутативна.
Теперь мы можем заметить, что в первой части доказательства предложения
\ref{PROP-H*(C-x)->H(Z)} (c.\pageref{homomorphism-sharp-C^times}) мы уже
убедились, что ${\mathscr F}_{\C^\times}$ является гомоморфизмом алгебр. А для
$\i_{\mathcal{O}(\Z)}$ это очевидно, значит мы получаем, что
$({\mathscr F}_{\Z})^\star$ -- также гомоморфизм алгебр, и это означает, что
${\mathscr F}_{\Z}$ -- гомоморфизм коалгебр.

4. Теперь остается проверить, что ${\mathscr F}_{\Z}$ сохраняет антипод:
$$
\sigma_{\mathcal{O}(\Z)}(\chi_t)(n)=\chi_t(-n)=t^{-n}=(t^{-1})^n=\chi_{t^{-1}}(n)
$$
$$
\Downarrow
$$
$$
\sigma_{\mathcal{O}(\Z)}(\chi_t)=\chi_{t^{-1}}
$$
$$
\Downarrow
$$
\begin{multline*}
{\mathscr F}_{\Z}(\sigma_{\mathcal{O}^\star(\Z)}(\alpha))(t)= (\sigma_{\mathcal{O}^\star(\Z)}(\alpha))(\chi_t)= (\alpha\circ\sigma_{\mathcal{O}(\Z)})(\chi_t)=\\=
\alpha(\sigma_{\mathcal{O}(\Z)}(\chi_t))=
\alpha(\chi_{t^{-1}})={\mathscr F}_{\Z}(\alpha)(t^{-1})=\sigma_{\mathcal{O}(\C^\times)}({\mathscr F}_{\Z}(\alpha))(t)
\end{multline*}
$$
\Downarrow
$$
$$
{\mathscr F}_{\Z}(\sigma_{\mathcal{O}^\star(\Z)}(\alpha))=\sigma_{\mathcal{O}(\C^\times)}({\mathscr F}_{\Z}(\alpha))
$$
\epr

\bpr[Окончание доказательства предложения
\ref{PROP-H*(C-x)->H(Z)}]\label{end-of-proof-PROP-H*(C-x)->H(Z)} В предложении
\ref{PROP-H*(C-x)->H(Z)} нам осталось доказать, что отображение
${\mathscr F}_{\C^\times}:\mathcal{O}^\star(\C^\times)\to \mathcal{O}(\Z)$ является
изоморфизмом коалгебр. Заметим, что сопряженное отображение
$$
({\mathscr F}_{\C^\times})^\star:\mathcal{O}^\star(\Z)\to(\mathcal{O}^\star(\C^\times))^\star=\mathcal{O}(\C^\times)^{\star\star}
$$
с точностью до изоморфизма $\i_{\mathcal{O}(\C^\times)}:\mathcal{O}(\C^\times)\cong \mathcal{O}(\C^\times)^{\star\star}$ совпадает с
${\mathscr F}_{\Z}$:
 \beq\label{sharp_C-x^star=i_H(C-x)-circ-sharp_Z}
\begin{diagram}
\node{\mathcal{O}^\star(\Z)}
\arrow[2]{e,t}{({\mathscr F}_{\C^\times})^\star}\arrow{se,b}{{\mathscr F}_{\Z}}
\node[2]{\mathcal{O}(\C^\times)^{\star\star}}
\\
\node[2]{\mathcal{O}(\C^\times)}\arrow{ne,b}{\i_{\mathcal{O}(\C^\times)}}
\end{diagram}
 \eeq
Это также следует из формулы \eqref{sharp_Z-C-x-delta_x}: для $t\in\C^\times$ и
$n\in\Z$
$$
({\mathscr F}_{\C^\times})^\star(\delta^n)(\delta^t)=\delta^n({\mathscr F}_{\C^\times}(\delta^t))=\eqref{sharp_Z-C-x-delta_x}
=\delta^t({\mathscr F}_{\Z}(\delta^n))=\i_{\mathcal{O}(\C^\times)}({\mathscr F}_{\Z}(\delta^n))(\delta^t)=(\i_{\mathcal{O}(\Z)}\circ\;
{\mathscr F}_{\C^\times})(\delta^n)(\delta^t)
$$
Это верно для любых $t\in\C^\times$ и $n\in\Z$, а дельта-функционалы плотны в
$\mathcal{O}^\star$, поэтому
$$
({\mathscr F}_{\C^\times})^\star=\i_{\mathcal{O}(\C^\times)}\circ\;{\mathscr F}_{\Z}
$$
то есть диаграмма \eqref{sharp_C-x^star=i_H(C-x)-circ-sharp_Z} коммутативна.
Теперь мы можем заметить, что в доказательстве предложения
\ref{PROP-H*(Z)->H(C-x)} мы уже убедились, что ${\mathscr F}_{\Z}$ является
гомоморфизмом алгебр. А для $\i_{\mathcal{O}(\C^\times)}$ это очевидно, значит
мы получаем, что $({\mathscr F}_{\C^\times})^\star$ -- также гомоморфизм алгебр, и
это означает, что ${\mathscr F}_{\C^\times}$ -- гомоморфизм коалгебр. \epr

\paragraph{Доказательство теоремы \ref{O-exp^star=O}} После того, как формулами
\eqref{finite-groups}, \eqref{O_exp*(C)=O(C)}, \eqref{O_exp*(C^x)=O(Z)},
\eqref{O_exp*(Z)=O(C^x)} мы установили изоморфизм \eqref{O_exp^*_G=O-G^bullet}
для случаев $G=\C,\; \C^\times,\; \Z$ и для случая конечной группы $G$, нам
остается просто применить формулы \eqref{tenz-pr-O-exp} и
\eqref{tenz-pr-O-exp-star}: разложив произвольную абелеву компактно порожденную
группу Штейна $G$ в произведение
$$
G=\C^l\times(\C^\times)^m\times\Z^n\times F
$$
где $F$ -- конечная группа, мы получим:
 \begin{multline*}
\mathcal{O}_{\exp}^\star(G)=\mathcal{O}_{\exp}^\star(\C^l\times(\C^\times)^m\times\Z^n\times
F)=\eqref{tenz-pr-O-exp}=\\= \mathcal{O}_{\exp}^\star(\C)^{\odot l}\odot
\mathcal{O}_{\exp}^\star(\C^\times)^{\odot m}\odot \mathcal{O}_{\exp}^\star(\Z)^{\odot n}\odot \mathcal{O}_{\exp}^\star(F)=\\= \mathcal{O}(\C^\bullet)^{\odot l}\odot \mathcal{O}((\C^\times)^\bullet)^{\odot m}\odot
\mathcal{O}(\Z^\bullet)^{\odot n}\odot \mathcal{O}(F^\bullet)=\\=\eqref{tenz-pr-O-exp}= \mathcal{O}((\C^\bullet)^l\times((\C^\times)^\bullet)^m\times(\Z^\bullet)^n\times
F^\bullet)= \mathcal{O} (G^\bullet)
 \end{multline*}

\subsection{Рефлексивность относительно оболочки Аренса---Майкла}

Условимся стереотипное пространство $H$ называть {\it жестким}, если 
преобразования Гротендика для пары $(H;H)$, для тройки $(H;H;H)$ и для четверки
$(H;H;H;H)$ являются изоморфизмами стереотипных пространств:
 \begin{align*}
& @_{H,H}:H\circledast H\cong H\odot H,\\ & @_{H,H,H}:H\circledast H\circledast
H\cong H\odot H\odot H\\ & @_{H,H,H,H}:H\circledast H\circledast H\circledast
H\cong H\odot H\odot H\odot H
 \end{align*}
(это всегда так, если $H$ -- ядерное пространство Фреше или ядерное
пространство Браунера). Тогда, очевидно, задание структуры сильной 
стереотипной алгебры Хопфа\footnote{То есть алгебры Хопфа в категории $(\Ste,\circledast)$.} на $H$ эквивалентно заданию структуры слабой 
стереотипной алгебры Хопфа\footnote{То есть алгебры Хопфа в категории $(\Ste,\odot)$.} на $H$: структурные элементы алгебр Хопфа в
$(\mathfrak{Ste},\circledast)$ и $(\mathfrak{Ste},\odot)$ (мы отличаем их
индексами $\circledast$ и $\odot$) будут или совпадать
$$
\iota_{\circledast}=\iota_{\odot},\qquad \e_{\circledast}=\e_{\odot},\qquad
\sigma_{\circledast}=\sigma_{\odot}
$$
или будут связаны диаграммами
$$
\begin{diagram}
\node{H\circledast H} \arrow{se,r}{\mu_{\circledast}} \arrow[2]{e,l}{@_{H,H}}
\node[2]{H\odot H} \arrow{sw,r}{\mu_{\odot}}
\\
\node[2]{H}
\end{diagram}
\qquad
\begin{diagram}
\node{H\circledast H}  \arrow[2]{e,l}{@_{H,H}} \node[2]{H\odot H}
\\
\node[2]{H}\arrow{nw,r}{\varkappa_{\circledast}}
\arrow{ne,r}{\varkappa_{\odot}}
\end{diagram}
$$
Такие (одновременно сильные и слабые) алгебры Хопфа на жестких стереотипных пространствах $H$ мы будем называть
{\it жесткими стереотипными алгебрами Хопфа}.\index{алгебра!Хопфа!жесткая}

Далее, пусть $H$ --- жесткая стереотипная алгебра Хопфа, и пусть ее оболочка
Аренса-Майкла $H^\heartsuit$ --- тоже жесткое стереотипное пространство, и на ней тоже определена структура (жесткой) стереотипной
алгебры Хопфа, причем выполняются следующие условия:
 \bit{
\item[(i)] естественный гомоморфизм алгебр
$$
\heartsuit_H:H\to H^\heartsuit
$$
является гомоморфизмом жестких алгебр Хопфа, и

\item[(ii)] сопряженное отображение
$$
(\heartsuit_H)^\star:(H^\heartsuit)^\star\to H^\star
$$
является оболочкой Аренса-Майкла алгебры $(H^\heartsuit)^\star$:
$$
(\heartsuit_H)^\star=\heartsuit_{(H^\heartsuit)^\star}
$$
 }\eit
Заметим в связи с этим вот что:

\bprop\label{PROP:edinstvennost-Hopfa-na-H^heartsuit} Для произвольной жесткой
стереотипной алгебры Хопфа $H$ структура жесткой алгебры Хопфа на ее оболочке
Аренса-Майкла $H^\heartsuit$, удовлетворяющая условиям (i) и (ii), если она
существует, определяется однозначно. \eprop
 \bpr
Заметим вначале, что из (i) и (ii) сразу следует условие
 \bit{
\item[(iii)] отображения $\heartsuit_H$ и $(\heartsuit_H)^\star$ являются
биморфизмами стеретипных пространств (то есть инъективны и имеют плотный образ
в области значений).
 }\eit
Действительно, отображения $\heartsuit_H:H\to H^\heartsuit$ и
$(\heartsuit_H)^\star:(H^\heartsuit)^\star\to H^\star$ являются эпиморфизмами
(имеют плотные образы), потому что это оболочки Аренса-Майкла. Поскольку они
сопряжены друг другу, они являются и мономорфизмами (то есть инъективны).

Отсюда следует все остальное. Прежде всего, умножение и единица на
$H^\heartsuit$ определяются однозначно условием, что $\heartsuit_H:H\to
H^\heartsuit$ -- оболочка Аренса-Майкла алгебры $H$. Рассмотрим теперь
сопряженное отображение $(\heartsuit_H)^\star:(H^\heartsuit)^\star\to H^\star$.
Как и $\heartsuit_H$, оно должно быть гомоморфизмом алгебр Хопфа. Значит, оно
является гомоморфизмом алгебр, причем инъективным, в силу доказанного свойства
(iii). Отсюда следует, что умножение и единица на $(H^\heartsuit)^\star$ также
определяются однозначно, поскольку они индуцируются из $H^\star$.

Таким образом, условия (i) и (ii) накладывают жесткие условия на умножение,
единицу, коумножение и коединицу в $H^\heartsuit$, позволяя определить на этом
пространстве не более одной структуры биалгебры. С другой стороны, мы знаем,
что антипод у биалгебры, если он существует, определяется однозначно, поэтому
структура алгебры Хопфа на $H^\heartsuit$ также единственна.
 \epr

Условия (i) и (ii) удобно изображать в виде диаграммы
 \beq\label{diagramma-refleksivnosti-AM}
 \xymatrix @R=1.pc @C=1.pc
 {
 H
 & \ar@{|->}[r]^{\heartsuit} & &
 H^{\heartsuit}
 \\
 & & &
 \ar@{|->}[d]^{\star}
 \\
 \ar@{|->}[u]^{\star}
 & & &
 \\
 H^\star
 & &
 \ar@{|->}[l]_{\heartsuit}
 &
 (H^{\heartsuit})^\star
 }
 \eeq
в которую мы предлагаем вкладывать вот какой смысл: во-первых, в углах квадрата
стоят жесткие стереотипные алгебры Хопфа, причем горизонтальные стрелки
(операции Аренса-Майкла $\heartsuit$) являются их гомоморфизмами, и, во-вторых,
чередование операций $\heartsuit$ и $\star$ (с какого места ни начинай) на
четвертом шаге возвращает к исходной алгебре Хопфа (конечно, с точностью до
изоморфизма функторов).

Жесткие стереотипные алгебры Хопфа $H$, удовлетворяющие условиям (i) и (ii), мы
будем называть {\it голоморфно рефлексивными}\index{алгебра!Хопфа!голоморфно
рефлексивная}, а диаграмму \eqref{diagramma-refleksivnosti-AM} для таких алгебр --
{\it диаграммой рефлексивности}\index{диаграмма рефлексивности}. Смысл термина
``рефлексивность'' в этом случае состоит в том, что если однократное
последовательное применение операций $\heartsuit$ и $\star$ обозначить
каким-нибудь символом, например $\dagger$,
$$
H^\dagger:=(H^\heartsuit)^\star
$$
и называть такой объект {\it алгеброй Хопфа, двойственной к
$H$ относительно оболочки Аренса---Майкла}\index{алгебра!Хопфа!двойственная относительно оболочки Аренса---Майкла}, то $H$ будет естественно
изоморфна своей второй двойственной в этом смысле алгебре Хопфа:
 \beq\label{H-cong-(H^*)^*}
H\cong (H^\dagger)^\dagger
 \eeq
Это можно считать следствием предложения
\eqref{PROP:edinstvennost-Hopfa-na-H^heartsuit}: поскольку для голоморфно
рефлексивных алгебр Хопфа переход $H\mapsto H^\heartsuit$ однозначно определяет
структуру алгебры Хопфа на $H^\heartsuit$, изоморфизм алгебр
$$
((H^\heartsuit)^\star)^\heartsuit\cong H^\star,
$$
постулируемый в аксиоме (ii), автоматически должен быть изоморфизмом алгебр
Хопфа. Переходя к сопряженным алгебрам Хопфа, мы как раз получаем
\eqref{H-cong-(H^*)^*}.

\paragraph{Рефлексивность относительно оболочки Аренса---Майкла для компактно порожденных абелевых групп Штейна.}
\label{SUBSEC:diagramma-vlozhenija}

 \btm\label{TH:reflex-AM-abelev-Stein}
Если $G$ --- компактно порожденная абелева группа Штейна, то алгебры $\mathcal{O}^\star(G)$ и $\mathcal{O}_{\exp}(G)$ голоморфно рефлексивны, а диаграмма рефлексивности для них принимает вид:
 \beq\label{chetyrehugolnik-O-O*-Abel}
 \xymatrix @R=1.pc @C=2.pc
 {
 \mathcal{O}^\star(G)
 & \ar@{|->}[r]^{\heartsuit}_{\eqref{AM-O-star}} & &
 \mathcal{O}_{\exp}^\star(G)\cong \mathcal{O}(G^\bullet)
 \\
 & & &
 \ar@{|->}[d]^{\star}
 \\
 \ar@{|->}[u]^{\star}
 & & &
 \\
 \mathcal{O}(G)
 & &
 \ar@{|->}[l]_{\heartsuit}^{\eqref{AM-O}}
 &
 \mathcal{O}_{\exp}(G)\cong \mathcal{O}^\star(G^\bullet)
 }
 \eeq
(цифры под горизонтальными стрелками -- ссылки на формулы в тексте).
 \etm

\paragraph{Рефлексивность Аренса---Майкла для конечно порожденных дискретных групп.}

Из теорем \ref{TH-AM-O-star} и \ref{TH-AM-O-diskr-gr} следует:

 \btm\label{TH:reflex-AM-kon-porozhd-diskr-gruppy-1}
Если $G$ --- конечно порожденная дискретная группа, то алгебры ${\mathcal O}^\star(G)$ и ${\mathcal
O}_{\exp}(G)$ рефлексивны относительно оболочки Аренса---Майкла с диаграммой рефлексивности
 \beq\label{chetyrehugolnik-O-O*}
 \xymatrix @R=1.pc @C=2.pc
 {
 {\mathcal O}^\star(G)
 & \ar@{|->}[r]^{\heartsuit}_{\eqref{AM-O-star}} & &
 {\mathcal O}_{\exp}^\star(G)
 \\
 & & &
 \ar@{|->}[d]^{\star}
 \\
 \ar@{|->}[u]^{\star}
 & & &
 \\
 {\mathcal O}(G)
 & &
 \ar@{|->}[l]_{\heartsuit}^{\eqref{AM-O}}
 &
 {\mathcal O}_{\exp}(G)
 }
 \eeq
 \etm

\btm\label{TH:reflex-AM-kon-porozhd-diskr-gruppy-2} 
Рефлексивность относительно оболочки Аренса---Майкла продолжает классическую двойственность Понтрягина для конечных абелевых групп на класс конечно порожденных дискретных групп:
{\sf
 \beq\label{diagramma-kategorij-dlya-konechno-porozhd-diskr-grupp}
 \xymatrix @R=3.pc @C=2.pc
 {
 \boxed{\begin{matrix}
  \text{алгебры Хопфа в $(\Ste,\circledast)$,}\\
  \text{рефлексивные относительно}\\
  \text{оболочки Аренса---Майкла}
 \end{matrix}}
 \ar[rr]^{H\mapsto (H^\heartsuit)^\star} & &
 \boxed{\begin{matrix}
  \text{алгебры Хопфа в $(\Ste,\circledast)$,}\\
  \text{рефлексивные относительно}\\
  \text{оболочки Аренса---Майкла}
 \end{matrix}}
 \\
 \boxed{\begin{matrix}
  \text{конечно порожденныые}\\
  \text{дискретные группы}
 \end{matrix}} \ar[u]^(0.45){\scriptsize\begin{matrix} \mathcal{O}^\star(G)\\
 \text{\rotatebox{90}{$\mapsto$}} \\ G\end{matrix}} & &
 \boxed{\begin{matrix}
  \text{конечно порожденныые}\\
  \text{дискретные группы}
 \end{matrix}} \ar[u]_(0.45){\scriptsize\begin{matrix} \mathcal{O}^\star(G) \\
 \text{\rotatebox{90}{$\mapsto$}} \\ G\end{matrix}} \\
 \boxed{\begin{matrix}
 \text{абелевы конечные группы}
 \end{matrix}} \ar[u]^{e} \ar[rr]^{G\mapsto \widehat{G}} & &
  \boxed{\begin{matrix}
 \text{абелевы конечные группы}
  \end{matrix}}\ar[u]_{e}
 }
 \eeq }
\etm

\paragraph{Рефлексивность относительно оболочки Аренса---Майкла для конечных расширений связных линейных групп.}

Из теорем \ref{TH-AM-O-star} и \ref{TH-AM-O-kon-rassh-svyaz-lin-gruppy} следует:

\btm\label{TH:reflex-AM-kon-rassh-svyaznyh-lin-grupp-1}
Если $G$ --- конечное расширение связной линейной группы, то алгебры $\mathcal{O}^\star(G)$ и $\mathcal{O}_{\exp}(G)$ голоморфно рефлексивны, а диаграмма рефлексивности для них принимает вид:
 \beq\label{chetyrehugolnik-O-O*}
 \xymatrix @R=1.pc @C=2.pc
 {
 \mathcal{O}^\star(G)
 & \ar@{|->}[r]^{\heartsuit}_{\eqref{AM-O-star}} & &
 \mathcal{O}_{\exp}^\star(G)
 \\
 & & &
 \ar@{|->}[d]^{\star}
 \\
 \ar@{|->}[u]^{\star}
 & & &
 \\
 \mathcal{O}(G)
 & &
 \ar@{|->}[l]_{\heartsuit}^{\eqref{AM-O}}
 &
 \mathcal{O}_{\exp}(G)
 }
 \eeq
(цифры под горизонтальными стрелками -- ссылки на формулы в тексте).
 \etm

\bex Для группы $\GL_n(\C)$ диаграмма рефлексивности
\eqref{chetyrehugolnik-O-O*} принимает вид
 \beq\label{chetyrehugolnik-O(GL)-O*(GL)}
 \xymatrix @R=1.pc @C=2.pc
 {
 \mathcal{O}^\star(\GL_n(\C))
 & \ar@{|->}[r]^{\heartsuit}_{\eqref{O(GL)=H_exp(GL)}} & &
 \mathcal{P}^\star(\GL_n(\C))
 \\
 & & &
 \ar@{|->}[d]^{\star}
 \\
 \ar@{|->}[u]^{\star}
 & & &
 \\
 \mathcal{O}(\GL_n(\C))
 & &
 \ar@{|->}[l]_{\heartsuit}^{\eqref{R(M)^heartsuit=O(M)}}
 &
 \mathcal{P}(\GL_n(\C))
 }
 \eeq
\eex

\btm\label{TH:reflex-AM-kon-rassh-svyaznyh-lin-grupp-2}
Рефлексивность относительно оболочки Аренса---Майкла продолжает классическую двойственность Понтрягина для конечных абелевых групп на класс конечных расширений связных линейных комплексных групп:
{\sf
 \beq\label{diagramma-kategorij-dlya-konechnyh-lineinyh-grupp}
 \xymatrix @R=3.pc @C=2.pc
 {
 \boxed{\begin{matrix}
  \text{алгебры Хопфа в $(\Ste,\circledast)$,}\\
  \text{рефлексивные относительно}\\
  \text{оболочки Аренса---Майкла}
 \end{matrix}}
 \ar[rr]^{H\mapsto (H^\heartsuit)^\star} & &
 \boxed{\begin{matrix}
  \text{алгебры Хопфа в $(\Ste,\circledast)$,}\\
  \text{рефлексивные относительно}\\
  \text{оболочки Аренса---Майкла}
 \end{matrix}}
 \\
 \boxed{\begin{matrix}
  \text{конечные расширения}\\
  \text{связных линейных групп}
 \end{matrix}} \ar[u]^(0.45){\scriptsize\begin{matrix} \mathcal{O}^\star(G)\\
 \text{\rotatebox{90}{$\mapsto$}} \\ G\end{matrix}} & &
 \boxed{\begin{matrix}
  \text{конечные расширения}\\
  \text{связных линейных групп}
 \end{matrix}} \ar[u]_(0.45){\scriptsize\begin{matrix} \mathcal{O}^\star(G) \\
 \text{\rotatebox{90}{$\mapsto$}} \\ G\end{matrix}} \\
 \boxed{\begin{matrix}
 \text{абелевы конечные группы}
 \end{matrix}} \ar[u]^{e} \ar[rr]^{G\mapsto \widehat{G}} & &
  \boxed{\begin{matrix}
 \text{абелевы конечные группы}
  \end{matrix}}\ar[u]_{e}
 }
 \eeq }
\etm

\chapter{ДОПОЛНЕНИЯ}

\section{Некоторые факты из теории категорий}

\subsection{Категории мономорфизмов и эпиморфизмов.}

\paragraph{Категория мономорфизмов $\varGamma_X$ и системы подобъектов.}

Пусть в категории ${\tt K}$ задан некий класс мономорфизмов $\varGamma$, содержащий все локальные единицы
$$
\{1_X;\ X\in\Ob({\tt K})\}\subseteq\varGamma\subseteq\Mono({\tt K})
$$
(руководящими примерами для нас являются классы  $\varGamma=\Mono$ и $\varGamma=\SMono$). Для всякого объекта $X$ в ${\tt K}$ обозначим через
$\varGamma_X$ класс всех морфизмов из $\varGamma$ с областью значений $X$:
 \beq\label{DEF:varGamma(X)}
 \varGamma_X=\{\sigma\in\varGamma:\quad \Ran\sigma=X\}.
 \eeq
Этот класс образует категорию, в которой морфизмом объекта $\rho\in\varGamma_X$ в объект
$\sigma\in\varGamma_X$, то есть мономорфизма $\rho:A\to X$ в мономорфизм
$\sigma:B\to X$, считается любой морфизм $\varkappa:A\to B$ в ${\tt K}$, замыкающий
диаграмму \beq\label{morphism-v-Mono(X)} \xymatrix @R=1pc @C=2pc {
A\ar[rd]^{\rho}\ar[dd]_{\varkappa} &   \\
  & X \\
B\ar[ru]_{\sigma} & } \eeq Саму эту диаграмму в исходной категории ${\tt K}$
можно представлять себе, как морфизм
$\rho\overset{\varkappa}{\longrightarrow}\sigma$ в категории $\varGamma_X$.
Композицией двух таких морфизмов
$\rho\overset{\varkappa}{\longrightarrow}\sigma$ и
$\sigma\overset{\lambda}{\longrightarrow}\tau$, то есть диаграмм
$$
\xymatrix @R=1pc @C=2pc
{
A\ar[rd]^{\rho}\ar[dd]_{\varkappa} &   \\
  & X \\
B\ar[ru]_{\sigma} & } \qquad \xymatrix @R=1pc @C=2pc {
B\ar[rd]^{\sigma}\ar[dd]_{\lambda} &   \\
  & X \\
C\ar[ru]_{\tau} & }
$$
будет морфизм $\rho\overset{\lambda\circ\varkappa}{\longrightarrow}\tau$, или,
что то же самое, диаграмма,
$$
\xymatrix @R=1pc @C=2pc
{
A\ar[rd]^{\rho}\ar[dd]_{\lambda\circ\varkappa} &   \\
  & X \\
C\ar[ru]_{\tau} & }
$$
которую следует представлять себе, как результат склеивания исходных диаграмм
по общему ребру $\sigma$, подрисовывания возникающей стрелки-композиции
$\varkappa\circ\lambda$, а затем выбрасывания промежуточной вершины $B$ и всех
входящих в нее и выходящих из нее ребер:
$$
\xymatrix 
{
A\ar@/^2ex/[rrd]^{\rho}\ar[dd]_{\lambda\circ\varkappa}\ar@{-->}[rd]_{\varkappa} & &  \\
  & B\ar@{-->}[r]_{\sigma}\ar@{-->}[ld]_{\lambda} &  X \\
  C\ar@/_2ex/[rru]_{\tau} & &
}
$$
Локальными единицами в $\varGamma_X$, понятное дело, считаются диаграммы вида
$$
\xymatrix @R=1pc @C=2pc
{
A\ar[rd]^{\rho}\ar[dd]_{1_A} &   \\
  & X \\
A\ar[ru]_{\rho} & }
$$

\brem
Композицию морфизмов в $\varGamma_X$ можно определять иначе, мы ее определили так, чтобы ей соответствовала в точности композиция морфизмов в ${\tt K}$:
$$
\lambda\kern-3pt\underset{\varGamma_X}{\circ}\kern-3pt\varkappa=\lambda\underset{\tt K}{\circ}\varkappa.
$$
\erem

\btm\label{card-Mono(a,b)-le-1} Для всякого объекта $X$ категория $\varGamma_X$
является графом. \etm \bpr Нужно убедиться, что для любых двух объектов
$\rho:A\to X$ и $\sigma:B\to X$ существует не более одного морфизма
$\rho\overset{\varkappa}{\longrightarrow}\sigma$. Действительно, морфизм
$\varkappa$ в диаграмме \eqref{morphism-v-Mono(X)} будет единственным, потому
что мономорфность $\sigma$ влечет импликацию:
$\sigma\circ\varkappa=\rho=\sigma\circ\varkappa'$ $\Longrightarrow$
$\varkappa=\varkappa'$. \epr

\brem\label{REM:svoistva-Mono(X)} В силу \cite[Example 1.3.2]{Akbarov-De-Gruyter-I} это означает, что {\it в категории
$\varGamma_X$ все морфизмы являются биморфизмами}. Связь между свойствами морфизма
$\rho\overset{\varkappa}{\longrightarrow}\sigma$ в категории $\varGamma_X$ и его
же, рассматриваемого, как морфизм $\varkappa:A\to B$ в исходной категории ${\tt
K}$ (замыкающий диаграмму \eqref{morphism-v-Mono(X)}), выражается следующими
наблюдениями:

\bit{\it

\item[---] \label{PROP:mono-v-Mono(X)} всякий морфизм
$\rho\overset{\varkappa}{\longrightarrow}\sigma$ в категории $\varGamma_X$
является мономорфизмом в исходной категории ${\tt K}$,

\item[---] \label{PROP:iso-v-Mono(X)} морфизм
$\rho\overset{\varkappa}{\longrightarrow}\sigma$ в категории $\varGamma_X$
является изоморфизмом в $\varGamma_X$ $\Longleftrightarrow$ $\varkappa$ является
изоморфизмом в исходной категории ${\tt K}$. }\eit \erem \bpr 1. Морфизм
$\varkappa$, замыкающий диаграмму \eqref{morphism-v-Mono(X)}, должен быть
мономорфизмом, поскольку $\sigma\circ\varkappa$
является мономорфизмом \cite[p.80 $1^\circ$]{Akbarov-De-Gruyter-I}.

2. Если морфизм $\varkappa:A\to B$, замыкающий диаграмму \eqref{morphism-v-Mono(X)}, является изоморфизмом в ${\tt K}$, то, положив $\lambda=\varkappa^{-1}:A\gets B$, мы получим диаграммы
\beq\label{iso-v-Mono(X)}
\xymatrix 
{
A\ar[dd]_{1_A}\ar@/^2ex/[rrd]^{\rho}\ar@{-->}[rd]_{\varkappa} & &  \\
  & B\ar@{-->}[r]_{\sigma}\ar@{-->}[ld]_{\lambda} &  X \\
  A\ar@/_2ex/[rru]_{\rho} & &
}
\qquad
\xymatrix 
{
B\ar[dd]_{1_B}\ar@/^2ex/[rrd]^{\sigma}\ar@{-->}[rd]_{\lambda} & &  \\
  & A\ar@{-->}[r]_{\rho}\ar@{-->}[ld]_{\varkappa} &  X \\
B\ar@/_2ex/[rru]_{\sigma} & & } \eeq которые будут коммутативны, потому что
$\rho$ и $\sigma$ в них -- мономорфизмы. Эти диаграммы означают, что морфизмы
$\rho\overset{\varkappa}{\longrightarrow}\sigma$ и
$\sigma\overset{\lambda}{\longrightarrow}\rho$ в категории $\varGamma_X$ взаимно
обратны. Наоборот, если даны взаимно обратные морфизмы
$\rho\overset{\varkappa}{\longrightarrow}\sigma$ и
$\sigma\overset{\lambda}{\longrightarrow}\rho$ в категории $\varGamma_X$, то это
означает, что должны быть коммутативны диаграммы \eqref{iso-v-Mono(X)}, и
значит $\varkappa$ и $\lambda$ -- взаимно обратные морфизмы в категории ${\tt
K}$, и поэтому $\varkappa$ должен быть изоморфизмом. \epr

Для отношения предпорядка в $\varGamma_X$ удобно ввести специальное обозначение,
$\to$, определяемое правилом
 \beq\label{DEF:to-in-Mono(X)}
\rho\to\sigma\Longleftrightarrow\quad \exists \varkappa\in\Mor({\tt K})\quad
\rho=\sigma\circ\varkappa.
 \eeq
Здесь морфизм $\varkappa$, если он существует, должен быть единственным, и, кроме того, он будет мономорфизмом (это
следует из того, что $\sigma$ -- мономорфизм). Поэтому определена операция,
которая каждой паре морфизмов $\rho,\sigma\in\varGamma_X$, удовлетворяющей условию
$\rho\to\sigma$, ставит в соответствие морфизм
$\varkappa=\varkappa^\sigma_\rho$ в \eqref{DEF:to-in-Mono(X)}:
 \beq\label{DEF:to-in-Mono(X)-*}
\rho=\sigma\circ\varkappa^\sigma_\rho.
 \eeq
При этом, если $\rho\to\sigma\to\tau$, то из цепочки
$$
\tau\circ\varkappa^\tau_\rho=\rho=\sigma\circ\varkappa^\sigma_\rho=
\tau\circ\varkappa^\tau_\sigma\circ\varkappa^\sigma_\rho,
$$
в силу мономорфности $\tau$, следует равенство
\beq\label{varkappa^gamma_alpha=varkappa^gamma_beta-circ-varkappa^beta_alpha}
\varkappa^\tau_\rho=\varkappa^\tau_\sigma\circ\varkappa^\sigma_\rho.
 \eeq

\bit{

\item[$\bullet$] {\it Системой подобъектов класса $\varGamma$} в объекте $X$ категории ${\tt K}$ называется любой скелет $S$ категории $\varGamma_X$, содержащий морфизм $1_X$ в качестве одного из элементов. Иными словами, подкласс $S$ в классе $\varGamma_X$ называется системой подобъектов в $X$, если
\bit{

\item[(a)]\label{DEF:toch-sist-podobj-a} локальная единица объекта $X$ принадлежит классу $S$:
$$
1_X\in S,
$$

\item[(b)]\label{DEF:toch-sist-podobj-b} для всякого мономорфизма $\mu\in\varGamma_X$ в классе $S$ найдется изоморфный ему мономорфизм $\sigma$:
$$
\forall\mu\in\varGamma_X\qquad\exists\sigma\in S\qquad \mu\cong\sigma.
$$

\item[(c)]\label{DEF:toch-sist-podobj-c} в классе $S$ изоморфизм (в смысле категории $\varGamma_X$) эквивалентен равенству:
$$
\forall\sigma,\tau\in S\qquad \Big( \sigma\cong\tau\quad\Longleftrightarrow\quad \sigma=\tau \Big)
$$
}\eit
По теореме о существовании скелета \cite[Proposition 1.2.34]{Akbarov-De-Gruyter-I}, такой класс $S$ всегда существует. Его элементы называются {\it подобъектами объекта $X$ класса $\varGamma$}. Сам класс $S$ наделяется структурой полной подкатегории в $\varGamma_X$.
}\eit

\btm\label{S-chast-upor-klass} Система подобъектов $S$ объекта $X$ всегда
является частично упорядоченным классом. \etm \bpr Пусть у подобъектов $\rho\in
S$ и $\sigma\in S$ существуют два взаимно обратных морфизма $\varkappa:A\gets
B$ и $\lambda:A\to B$, то есть
$$
\rho=\sigma\circ\varkappa,\qquad \sigma=\rho\circ\lambda.
$$
Тогда
$$
\rho\circ\lambda\circ\varkappa=\rho=\rho\circ 1_A,\qquad
\sigma\circ\varkappa\circ\lambda=\sigma=\sigma\circ 1_B,
$$
и, поскольку $\rho$ и $\sigma$ -- мономорфизмы в ${\tt K}$, на них можно
сократить,
$$
\lambda\circ\varkappa=1_A,\qquad \varkappa\circ\lambda=1_B,
$$
то есть $\varkappa$ и $\lambda$ -- изоморфизмы. Мы получаем, что
$\rho\cong\sigma$, поэтому в силу свойства (c), $\rho=\sigma$. \epr

\btm\label{PROP:dost-mnozh-v-klasse-podobjektov} Если $S$ -- система подобъектов на объекте $X$, то для всякого подобъекта $\sigma\in S$, $\sigma:Y\to X$, класс мономорфизмов
$$
A=\{\alpha\in\varGamma_Y:\ \sigma\circ\alpha\in S\}
$$
является системой подобъектов на объекте $Y$. Если вдобавок $S$ -- множество, то $A$ -- тоже множество.
\etm
\bpr
1. Условие (a) очевидно: поскольку $\sigma\circ 1_Y=\sigma\in S$, мы получаем, что $1_Y\in A$.

2. Условие (b). Пусть $\beta:B\to Y$ -- какой-нибудь мономорфизм. Рассмотрим композицию $\sigma\circ\beta:B\to X$. Это мономорфизм из $\varGamma_X$, поэтому, поскольку $S$ -- система подобъектов на $X$, должен существовать $\tau\in S$ такой, что
$$
\tau\cong\sigma\circ\beta
$$
Это значит, что
$$
\tau=\sigma\circ\beta\circ\iota
$$
для некоторого изоморфизма $\iota$. Теперь мы получаем, что мономорфизм $\alpha=\beta\circ\iota$ изоморфен $\beta$
$$
\alpha\cong\beta
$$
и при этом лежит в $A$, потому что $\sigma\circ\alpha=\tau\in S$.

3. Условие (c). Пусть $\alpha,\beta\in A$ -- какие-то два изоморфных мономорфизма, то есть
$$
\alpha=\beta\circ\iota
$$
для некоторого изоморфизма $\iota$. Тогда, во-первых, морфизмы $\sigma\circ\alpha$ и $\sigma\circ\beta$ тоже будут изоморфны, потому что
$$
\sigma\circ\alpha=\sigma\circ\beta\circ\iota
$$
И, во-вторых, они будут лежать в $S$, потому что $\alpha$ и $\beta$ лежат в $A$. Поскольку $S$ удовлетворяет условию (c),
морфизмы $\sigma\circ\alpha$ и $\sigma\circ\beta$ должны совпадать:
$$
\sigma\circ\alpha=\sigma\circ\beta
$$
Поскольку вдобавок $\sigma$ -- мономорфизм, мы получаем $\alpha=\beta$.

4. Остается проверить, что если $S$ -- множество, то $A$ -- тоже множество. Для этого достаточно убедиться, что отображение $\alpha\in A\mapsto \sigma\circ\alpha\in S$ является инъекцией. Действительно, если для каких-то $\alpha,\alpha'\in A$ выполняется
$$
\sigma\circ\alpha=\sigma\circ\alpha',
$$
то, поскольку $\sigma$ -- мономорфизм, это означает, что $\alpha=\alpha'$.
\epr

 \bit{
\item[$\bullet$] Категорию ${\tt K}$ мы называем  {\it локально малой в подобъектах класса $\varGamma$}, если в ней всякий объект $X$ обладает системой подобъектов $S$, образующей множество (а не просто класс объектов); это эквивалентно тому, что всякая категория $\varGamma_X$ является скелетно малым графом.
 }\eit

\bex  Часто используемые в качестве примеров категории множеств, групп,
векторных пространств, алгебр над данным полем, топологических пространств,
топологических векторных пространств, топологических алгебр, и т.д., очевидно,
локально малы в подобъектах класса $\Mono$. \eex

 \btm\label{TH:o-lok-malosti-v-podobjektah}
Если категория ${\tt K}$ локально мала в подобъектах класса $\varGamma$, то существует отображение $X\mapsto S_X$ которое каждому объекту $X$ в ${\tt K}$ ставит в соответствие его систему подобъектов $S_X$ класса $\varGamma$, являющуюся множеством.
 \etm
 \bpr
Здесь используется теорема о полном упорядочении класса всех множеств \cite[Theorem 1.1.1]{Akbarov-De-Gruyter-I}: каждому $X$ можно поставить в соответствие систему подобъектов $S$, наименьшую относительно этого упорядочения.
\epr

\paragraph{Категория эпиморфизмов $\varOmega^X$ и системы фактор-объектов.}

Пусть в категории ${\tt K}$ задан некий класс эпиморфизмов $\varOmega$, содержащий все локальные единицы
$$
\{1_X;\ X\in\Ob({\tt K})\}\subseteq\varOmega\subseteq\Epi({\tt K})
$$
(руководящими примерами для нас являются классы  $\varOmega=\Epi$ и $\varOmega=\SEpi$). Для всякого объекта $X$ в ${\tt K}$ обозначим через $\varOmega^X$ класс всех морфизмов из $\varOmega$ с областью определения $X$:
 \beq\label{DEF:varOmega(X)}
 \varOmega^X=\{\sigma\in\varOmega:\quad \Dom\sigma=X\}.
 \eeq
Этот класс образует категорию, в которой морфизмом объекта $\rho\in\varOmega^X$ в объект
$\sigma\in\varOmega^X$, то есть эпиморфизма $\rho:X\to A$ в эпиморфизм $\sigma:X\to
B$, считается любой морфизм $\varkappa:A\to B$, замыкающий диаграмму
\beq\label{morphism-v-Epi(X)} \xymatrix @R=1pc @C=2pc {
 & A\ar[dd]^{\varkappa}    \\
X\ar[ru]^{\rho}\ar[rd]_{\sigma}  &  \\
 & B
} \eeq Саму эту диаграмму в исходной категории ${\tt K}$ можно представлять
себе, как морфизм $\rho\overset{\varkappa}{\longrightarrow}\sigma$ в категории
$\varOmega^X$. Композицией двух таких морфизмов
$\rho\overset{\varkappa}{\longrightarrow}\sigma$ и
$\sigma\overset{\lambda}{\longrightarrow}\tau$, то есть диаграмм
$$
\xymatrix @R=1pc @C=2pc
{
 & A\ar[dd]^{\varkappa}    \\
X\ar[ru]^{\rho}\ar[rd]_{\sigma}  &  \\
 & B
}
\qquad
\xymatrix @R=1pc @C=2pc
{
 & B\ar[dd]^{\lambda}    \\
X\ar[ru]^{\sigma}\ar[rd]_{\tau}  &  \\
 & C
}
$$
будет морфизм $\rho\overset{\lambda\circ\varkappa}{\longrightarrow}\tau$, или,
что то же самое, диаграмма,
$$
\xymatrix @R=1pc @C=2pc
{
 & A\ar[dd]^{\lambda\circ\varkappa}    \\
X\ar[ru]^{\rho}\ar[rd]_{\tau}  &  \\
 & C
}
$$
которую следует представлять себе, как результат склеивания исходных диаграмм
по общему ребру $\sigma$, подрисовывания возникающей стрелки-композиции
$\lambda\circ\varkappa$, а затем выбрасывания промежуточной вершины $B$ и всех
входящих в нее и выходящих из нее ребер:
$$
\xymatrix 
{
 & & A\ar[dd]^{\lambda\circ\varkappa}\ar@{-->}[ld]^{\varkappa}    \\
X\ar@/^2ex/[rru]^{\rho}\ar@{-->}[r]_{\sigma}\ar@/_2ex/[rrd]_{\tau}  & B\ar@{-->}[rd]^{\lambda} &  \\
 & & C
}
$$
Локальными единицами в $\varOmega^X$, понятное дело, считаются диаграммы вида
$$
\xymatrix @R=1pc @C=2pc
{
 & A\ar[dd]^{1_A}    \\
X\ar[ru]^{\rho}\ar[rd]_{\sigma}  &  \\
 & A
}
$$

\brem
Композицию морфизмов в $\varOmega^X$ можно определять иначе, мы ее определили так, чтобы ей соответствовала в точности композиция морфизмов в ${\tt K}$:
$$
\lambda\kern-3pt\underset{\varOmega^X}{\circ}\kern-3pt\varkappa=\lambda\underset{\tt K}{\circ}\varkappa.
$$
\erem

По аналогии со свойствами $\varGamma_X$ доказываются свойства категории $\varOmega^X$.

\btm\label{card-Epi(a,b)-le-1}
Для всякого объекта $X$ категория $\varOmega^X$ является графом.
\etm

\brem\label{REM:svoistva-Epi(X)} В силу \cite[Example 1.3.2]{Akbarov-De-Gruyter-I} это означает, что {\it в категории $\varOmega^X$
все морфизмы являются биморфизмами}. Связь между свойствами морфизма
$\rho\overset{\varkappa}{\longrightarrow}\sigma$ в категории $\varOmega^X$ и его
же, рассматриваемого, как морфизм $\varkappa:A\to B$ в исходной категории ${\tt
K}$ (замыкающий диаграмму \eqref{morphism-v-Epi(X)}), выражается следующими
наблюдениями:

\bit{\it

\item[---] \label{PROP:epi-v-Epi(X)} всякий морфизм
$\rho\overset{\varkappa}{\longrightarrow}\sigma$ в категории $\varOmega^X$ является
эпиморфизмом в исходной категории ${\tt K}$,

\item[---] \label{PROP:iso-v-Epi(X)} морфизм
$\rho\overset{\varkappa}{\longrightarrow}\sigma$ в категории $\varOmega^X$ является
изоморфизмом в $\varOmega^X$ $\Longleftrightarrow$ $\varkappa$ является
изоморфизмом в исходной категории ${\tt K}$. }\eit \erem

Для отношения предпорядка в $\varOmega^X$ удобно ввести специальное обозначение,
$\to$, определяемое правилом
 \beq\label{DEF:le-in-F_X}
\rho\to\sigma\Longleftrightarrow\quad \exists \iota\in\Mor({\tt K})\quad
\sigma=\iota\circ\rho.
 \eeq
Здесь морфизм $\iota$, если он существует, должен быть единственным и, кроме того, он будет эпиморфизмом (это
следует из того, что $\rho$ и $\sigma$ -- эпиморфизмы). Поэтому определена операция,
которая каждой паре морфизмов $\rho,\sigma\in\varOmega^X$, удовлетворяющей условию
$\rho\to\sigma$, ставит в соответствие морфизм $\iota=\iota^\sigma_\rho$ в
\eqref{DEF:le-in-F_X}:
 \beq\label{DEF:le-in-F_X-*}
\sigma=\iota^\sigma_\rho\circ\rho.
 \eeq
При этом, если $\pi\to\rho\to\sigma$, то из цепочки
$$
\iota^\sigma_\pi\circ\pi=\sigma=\iota^\sigma_\rho\circ\rho=\iota^\sigma_\rho\circ\iota^\rho_\pi\circ\pi,
$$
опять в силу эпиморфности $\pi$, следует равенство
\beq\label{iota_rho^tau=iota_rho^sigma-circ-iota_sigma^tau}
\iota^\sigma_\pi=\iota^\sigma_\rho\circ\iota^\rho_\pi.
 \eeq

\bit{

\item[$\bullet$] {\it Системой фактор-объектов класса $\varOmega$} на объекте $X$ категории ${\tt K}$ называется любой скелет $Q$ категории $\varOmega^X$, содержащий морфизм $1_X$ в качестве одного из элементов. Иными словами, подкласс $Q$ в классе $\varOmega^X$ называется системой фактор-объектов объекта $X$, если
\bit{

\item[(a)] локальная единица объекта $X$ принадлежит классу $Q$:
$$
1_X\in Q,
$$

\item[(b)] для всякого эпиморфизма $\e\in\varOmega^X$ в классе $Q$ найдется изоморфный ему эпиморфизм $\pi$:
$$
\forall\e\in\varOmega^X\qquad\exists\pi\in Q\qquad \e\cong\pi,
$$

\item[(c)] в классе $Q$ изоморфизм (в смысле категории $\varOmega^X$) эквивалентен равенству:
$$
\forall\pi,\rho\in Q\qquad \Big( \pi\cong\rho\quad\Longleftrightarrow\quad \pi=\rho \Big)
$$
}\eit
По теореме о существовании скелета \cite[Proposition 1.2.34]{Akbarov-De-Gruyter-I}, такой класс $Q$ всегда существует. Его элементы  называются {\it фактор-объектами объекта $X$ класса $\varOmega$}. Сам класс $Q$ наделяется структурой полной подкатегории в $\varOmega^X$.

}\eit

По аналогии с теоремами \ref{S-chast-upor-klass} и \ref{PROP:dost-mnozh-v-klasse-podobjektov} доказываются

\btm\label{Q-chast-upor-klass} Система $Q$ фактор-объектов объекта $X$ всегда является частично упорядоченным классом.
\etm

\btm\label{PROP:dost-mnozh-v-klasse-faktor-objektov} Если $Q$ -- система фактор-объектов объекта $X$, то для всякого фактор-объекта $\pi\in Q$, $\pi:X\to Y$, класс эпиморфизмов
$$
A=\{\alpha\in\varOmega^Y:\ \alpha\circ\pi\in Q\}
$$
является системой фактор-объектов объекта $Y$. Если вдобавок $Q$ -- множество, то $A$ -- тоже множество.
\etm

 \bit{
\item[$\bullet$] Категорию ${\tt K}$ мы называем  {\it локально малой в фактор-объектах класса $\varOmega$}, если в ней всякий объект $X$ обладает системой фактор-объектов $Q$ класса $\varOmega$, образующей множество (а не просто класс объектов); это эквивалентно тому, что всякая категория $\varOmega^X$ является скелетно малым графом.
 }\eit

\bex Среди обычно приводимых в качестве примеров категорий множеств, групп,
векторных пространств, алгебр над данным полем, топологических пространств,
топологических векторных пространств, топологических алгебр, и т.д., некоторые, как категория векторных
пространств, локально малы и в фактор-объектах класса $\Epi$. Однако для
некоторых, как для категории алгебр, локальная малость в фактор-объектах
доказывается необычайно сложно (см.\cite{Adamek-Rosicky}). При этом, более
слабое условие локальной малости в фактор-объектах класса $\SEpi$ проверяется намного
проще. \eex

По аналогии с теоремой \ref{TH:o-lok-malosti-v-podobjektah} доказывается

 \btm\label{TH:o-lok-malosti-v-faktor-objektah}
Если категория ${\tt K}$ локально мала в фактор-объектах класса $\varOmega$, то существует отображение $X\mapsto Q_X$ которое каждому объекту $X$ в ${\tt K}$ ставит в соответствие его систему фактор-объектов $Q_X$ класса $\varOmega$, являющуюся множеством.
 \etm

\subsection{Алгебры Хопфа в моноидальной категории.}

Нам понадобится несколько определений из теории категорий. Пусть $(I,\le)$ -- частично упорядоченное множество, и пусть
 \bit{
\item[--] каждому индексу $i\in I$ поставлен в соответствие объект $X^i$
категории ${\tt K}$, а

\item[--] каждой паре индексов $i\le j$ поставлен в соответствие морфизм
$X^i\overset{\iota_i^j}{\to} X^j$,
 }\eit
причем
 \bit{
\item[1)] морфизмы $\iota_i^i$ (с совпадающими нижним и верхним индексами) есть просто единичные морфизмы:
$$
\iota_i^i=1_{X^i}
$$

\item[2)] для любых трех индексов $i\le j\le k$ коммутативна диаграмма
$$
\begin{diagram}
\node[2]{X^j}\arrow{se,t}{\iota_j^k}  \\
\node{X^i}\arrow[2]{e,b}{\iota_i^k}\arrow{ne,t}{\iota_i^j}
\node[2]{X^k}
\end{diagram}
$$
 }\eit\noindent
Тогда семейство $\{X^i;\iota_i^j\}$ называется {\it ковариантной системой} в
категории ${\tt K}$ над частично упорядоченным множеством $(I,\le)$.

Если $\{X_i,\iota_i^j\}$ и $\{Y_i,\varkappa_j^i\}$ -- ковариантные системы над частично упорядоченным множеством $I$, то {\it морфизмом} $\ph:\{X_i,\iota_i^j\}\to \{Y_i,\varkappa_j^i\}$ называется система морфизмов $\ph_i:X_i\to Y_i$, замыкающая все диаграммы
$$
 \xymatrix 
 {
 X_i\ar[r]^{\iota_i^j}\ar[d]_{\ph_i}& X_j\ar[d]^{\ph_j} \\
 Y_i\ar[r]_{\varkappa_j^i} & Y_j
 }
$$

 Пусть $\{X^i;\iota_i^j\}$  -- ковариантная система над частично упорядоченным множеством $I$ в категории ${\tt K}$.
 \bit{

\item[---] Для всякого объекта $X$ в ${\tt K}$ {\it проективным конусом} ковариантной системы $\{X^i;\iota_i^j\}$ с вершиной $X$ называется система морфизмов $\pi^i:X \to X^i$ такая, что для любых индексов $i\le j$ коммутативна диаграмма:
$$
\begin{diagram}
\node[2]{X}\arrow{sw,t}{\pi^i}\arrow{se,t}{\pi^j}  \\
\node{X^i}\arrow[2]{e,b}{\iota_i^j} \node[2]{X^j}
\end{diagram}
$$

\item[---] Проективный конус $\{X,\pi^i\}$ ковариантной системы $\{X^i;\iota_i^j\}$
называется {\it проективным пределом}\label{DEF:proj-limit} этой системы, если для любого другого ее проективного
конуса $\{Y,\rho^i\}$ найдется единственный морфизм
$\tau:Y\to X$, такой, что для любого индекса $i$ будет коммутативна
диаграмма:
\beq\label{DEF:proj-limit-diagr}
\begin{diagram}
\node{Y}\arrow{se,b}{\rho^i}\arrow[2]{e,t,--}{\tau}\node[2]{X}\arrow{sw,b}{\pi^i}
\\
\node[2]{X^i}
\end{diagram}
\eeq
В этом случае для объекта $X$ и морфизмов $\pi^i$ и $\tau$ используются обозначения:
$$
X=\lim_{\infty\gets j}X^j,\qquad \pi^i=\lim_{\infty\gets j}\iota_i^j,\qquad \tau=\lim_{\infty\gets j}\rho^j.
$$
 }\eit
Двойственным образом определяются инъективные пределы.

Условимся говорить, что в моноидальной категории ${\sf M}$ {\it тензорное произведение $\otimes$ коммутирует с проективными пределами}, если в ней выполняются естественные тождества
$$
\projlim_{(i,j)\to\infty}(X_i\otimes Y_j)=(\projlim_{i\to\infty}X_i)\otimes (\lim_{\infty\gets j}Y_j).
$$

\bex\label{EX:odot-perestanov-s-proj-lim}
В категории $({\sf Ste},\odot)$ тензорное произведение $\odot$ коммутирует с проективными пределами \cite[(4.153)]{Akbarov-De-Gruyter-I}.
\eex

Точно так же мы говорим, что в моноидальной категории ${\sf M}$ {\it тензорное произведение $\otimes$ коммутирует с инъективными пределами}, если в ней выполняются естественные тождества
$$
\lim_{(i,j)\to\infty}(X_i\otimes Y_j)=(\injlim_{i\to\infty}X_i)\otimes (\lim_{j\to\infty}Y_j).
$$

\bex\label{EX:circledast-perestanov-s-inj-lim}
В категории $({\sf Ste},\circledast)$ тензорное произведение $\circledast$ коммутирует с инъективными пределами \cite[(4.153)]{Akbarov-De-Gruyter-I}.
\eex

Напомним, что частично упорядоченное множество $(I,\le)$ называется {\it направленным по возрастанию} (соответственно, {\it по убыванию}), если для любых $i,j\in I$ найдется $k\in I$ такой что
$$
i\le k\ \&\ j\le k
$$
(соответственно, $k\le i\ \&\ k\le j$).

\btm\label{TH:nasledovanie-Hopfa-predelami}
Пусть в моноидальной категории ${\sf M}$ тензорное произведение $\otimes$ коммутирует с проективными (соответственно, инъективными) пределами, и пусть $\{H_i;\iota_i^j\}$ --- контравариантная (соответственно, ковариантная) система в категории ${\sf Hopf}_{\sf M}$ алгебр Хопфа в ${\sf M}$ над направленным по возрастанию множеством $I$. Тогда если система $\{H_i;\iota_i^j\}$ обладает проективным (соответственно, инъективным) пределом в категории ${\sf M}$, то этот предел обладает струтурой алгебры Хопфа в ${\sf M}$, превращающей его в проективный (соответственно, инъективный) предел системы $\{H_i;\iota_i^j\}$ в категории ${\sf Hopf}_{\sf M}$:
$$
\overset{{\sf M}}{\projlim_{i\to\infty}}H_i=\overset{{\sf Hopf}_{\sf M}}{\projlim_{i\to\infty}}H_i \qquad (\overset{{\sf M}}{\injlim_{i\to\infty}}H_i=\overset{{\sf Hopf}_{\sf M}}{\injlim_{i\to\infty}}H_i.
)
$$
\etm

\brem
Теорема \ref{TH:nasledovanie-Hopfa-predelami} остается справедливой если в ней заменить алгебры Хопфа на биалгебры, алгебры и коалгебры. Мы ее формулируем для алгебр Хопфа потому что именнно в таком виде она понадобится нам далее в следствии \ref{COR:proj-lim-Hopf-algebr-odot}, которое в свою очередь будет затем использоваться в теореме \ref{TH:Env_C-C^*(G)-odot-Hopf}.
\erem

Для доказательства нам понадобятся две леммы.

\blm\label{ph:lim-X_i->lim-Y_i}
Пусть $I$ -- частично упорядоченное множество и $\{X_i,\iota_i^j\}$ и $\{Y_i,\varkappa_j^i\}$ --- две контравариантные системы над ним, имеющие проективные пределы
$$
X=\projlim_{i\to\infty}X_i,\qquad Y=\projlim_{i\to\infty}Y_i.
$$
Для всякого морфизма этих контравариантных систем $\ph:\{X_i,\iota_i^j\}\to \{Y_i,\varkappa_j^i\}$ существует морфизм между проективными пределами $\ph:X\to Y$, замыкающий все диаграммы
$$
 \xymatrix 
 {
  & X\ar[dl]_{\pi_i}\ar[dr]^{\pi_j}\ar@{--}[d] & \\
 X_i\ar[dd]_{\ph_i}& \ar@{-->}[d]_{\ph} & X_j\ar[dd]^{\ph_j}\ar[ll]^(.7){\iota_i^j} \\
  & Y\ar[dl]_{\rho_i}\ar[dr]_{\rho_j} & \\
 Y_i & & Y_j\ar[ll]^{\varkappa_j^i}
 }
$$
\elm
\bpr
Здесь нужно просто заметить, что система морфизмов $\{\ph_i\circ\pi_i\}$ является проективным конусом для ковариантной системы $\{Y_i,\varkappa_j^i\}$, поэтому должен существовать морфизм $\ph$, для которого будут коммутативны все дальние боковые грани призмы:
$$
\ph_i\circ\pi_i=\rho_i\circ\ph.
$$
А все остальные грани коммутативны в силу определений конуса и морфизма ковариантных систем.
\epr

\blm\label{LM:lim-X_(i,i)=lim-X_(i,j)}
Если $I$ -- направленное по возрастанию множество и $\{X_{i,j}\}$ -- контравариантная система над его декартовым квадратом $I\times I$, имеющая проективный предел, то диагональ этой системы
$\{X_{i,i}\}$ также имеет проективный предел, и эти пределы совпадают:
\beq\label{lim-X_(i,i)=lim-X_(i,j)}
\projlim_{i\to\infty}X_{i,i}=\projlim_{(i,j)\to\infty}X_{i,j}.
\eeq
\elm
\bpr
Это следует из определения направленного множества.
\epr

\bpr[Доказательство теоремы \ref{TH:nasledovanie-Hopfa-predelami}.] Здесь нужно понять, как определяются структурные морфизмы на проективном пределе
$$
H=\overset{{\sf M}}{\projlim_{i\to\infty}}H_i,
$$
превращающие его в алгебру Хопфа. Для этого нужно просто воспользоваться формулой, вытекающей из \eqref{lim-X_(i,i)=lim-X_(i,j)}:
$$
\projlim_{i\to\infty}X_i\otimes X_i=\projlim_{(i,j)\to\infty}X_i\otimes X_j.
$$
Из этой формулы следует, что, например, умножение в $H$ можно определить, перейдя от умножений $\mu_i:H_i\otimes H_i\to H_i$ к их проективному пределу по лемме \ref{ph:lim-X_i->lim-Y_i}
$$
\mu:\projlim_{i\to\infty}(H_i\otimes H_i)\to \projlim_{i\to\infty} H_i=H,
$$
затем по лемме \ref{LM:lim-X_(i,i)=lim-X_(i,j)} заменив $\projlim_{i\to\infty}(H_i\otimes H_i)$ на $\projlim_{(i,j)\to\infty}(H_i\otimes H_j)$, а после этого воспользовавшись перестановочностью $\otimes$ с проективными пределами, нужно будет заменить $\projlim_{(i,j)\to\infty}(H_i\otimes H_j)$ на $\projlim_{i\to\infty}H_i\otimes\lim_{\infty\gets j}H_j$.
$$
\mu:H\otimes H=\projlim_{i\to\infty}H_i\otimes\lim_{\infty\gets j}H_j\cong\projlim_{(i,j)\to\infty}(H_i\otimes H_j)\cong\projlim_{i\to\infty}(H_i\otimes H_i)\to \projlim_{i\to\infty} H_i=H,
$$
Таким же образом определяются остальные структурные морфизмы.
\epr

Пусть далее ${\sf Hopf}_\odot$ обозначает класс алгебр Хопфа в моноидальной категории стереотипных пространств $(\Ste,\odot,\C)$, а ${\sf Hopf}_\circledast$ --- класс  алгебр Хопфа в категории $(\Ste,\circledast,\C)$. Из полноты категории {\sf Ste} \cite[Theorem 4.21]{Ak03} и примеров \ref{EX:odot-perestanov-s-proj-lim} и \ref{EX:circledast-perestanov-s-inj-lim} мы получаем два важных следствия:

\bcor\label{COR:proj-lim-Hopf-algebr-odot}
Всякая контравариантная система $\{H_i;\iota_i^j\}$ в категории ${\sf Hopf}_\odot$  стереотипных алгебр Хопфа над направленным по возрастанию множеством $I$ обладает проективным пределом, причем им будет проективный предел системы $\{H_i;\iota_i^j\}$ в категории ${\sf Ste}$ с подходящей структурой инъективной алгебры Хопфа:
$$
\overset{{\sf Ste}}{\projlim_{i\to\infty}}H_i=\overset{{\sf Hopf}_\odot}{\projlim_{i\to\infty}}H_i.
$$
\ecor

\bcor\label{COR:inj-lim-Hopf-algebr-circledast}
Всякая ковариантная система $\{H_i;\iota_i^j\}$ в категории ${\sf Hopf}_\circledast$  стереотипных алгебр Хопфа над направленным по возрастанию множеством $I$ обладает инъективным пределом, причем им будет инъективный предел системы $\{H_i;\iota_i^j\}$ в категории ${\sf Ste}$ с подходящей структурой проективной алгебры Хопфа:
$$
\overset{{\sf Ste}}{\injlim_{i\to\infty}}H_i=\overset{{\sf Hopf}_\circledast}{\injlim_{i\to\infty}}H_i.
$$
\ecor

\subsection{Теоремы о повторных пределах}

\paragraph{Теорема о повторном проективном пределе.}

Пусть $(I,\le)$ --- частично упорядоченное множество, и пусть ${\mathscr J}\subseteq 2^I$ --- система его подмножеств со следующими свойствами:
\bit{

\item[(i)] система ${\mathscr J}$ покрывает множество $I$:
\beq\label{DEF:pokrytie-indeksov-1}
\bigcup_{J\in{\mathscr J}} J=I,
\eeq

\item[(ii)] система ${\mathscr J}$ направлена по включению:
\beq\label{DEF:pokrytie-indeksov-2}
\forall J,K\in {\mathscr J}\qquad \exists L\in {\mathscr J}\qquad J\cup K\subseteq L.
\eeq

}\eit\noindent
Тогда мы будем говорить, что {\it система множеств ${\mathscr J}$ является покрытием\label{DEF:pokrytie-POSet} частично упорядоченного множества $(I,\le)$}.

\btm[о повторном проективном пределе]\label{TH:o-povtornom-proj-predele}
Пусть даны:
\bit{

\item[(a)] проективно полная категория ${\sf K}$;

\item[(b)] контравариантная система $\{X_i,\pi_i^j;\ i\le j,\ i,j\in I\}$ в категории ${\sf K}$,

\item[(c)] покрытие ${\mathscr J}$ частично упорядоченного множества $(I,\le)$.

}\eit\noindent
Пусть  
\beq\label{TH:o-povtornom-proj-predele-5-2}
X=\projlim_{i\in I} X_i
\eeq  
--- проективный предел контравариантной системы $\{X_i,\pi_i^j;\ i\le j,\ i,j\in I\}$ со структурными проекциями 
\beq\label{TH:o-povtornom-proj-predele-5-3}
\pi_i:X\to X_i.
\eeq
Тогда 
\bit{

\item[(i)] для всякого множества $J\in{\mathscr J}$ семейство $\{X_i,\pi_i^j;\ i\le j,\ i,j\in J\}$ образует контравариантную систему в категории ${\sf K}$, у которой имеется проективный предел
\beq\label{TH:o-povtornom-proj-predele-1}
X_J:=\projlim_{j\in J}X_j
\eeq
с неким структурным семейством проекций 
\beq\label{TH:o-povtornom-proj-predele-1-1}
\pi_j^J:X_J\to X_j;\ j\in J,
\eeq

\item[(ii)] для  любых множеств $J,K\in{\mathscr J}$, связанных включением $J\subseteq K$ существует единственный морфизм
\beq\label{TH:o-povtornom-proj-predele-2}
\pi_J^K:X_J\gets X_K,
\eeq
замыкающий все диаграммы
\beq\label{TH:o-povtornom-proj-predele-2-1}
 \xymatrix 
 {
 X_J\ar[dr]_{\pi_i^J} & & X_K\ar@{-->}[ll]_{\pi_J^K}\ar[dl]^{\pi_i^K} \\
  & X_i &  
 },\qquad i\in J
\eeq

\item[(iii)] семейство $\{X_J,\pi_J^K;\ J,K\in {\mathscr J}\}$ образует контравариантную систему в категории ${\sf K}$, 
\beq\label{TH:o-povtornom-proj-predele-3-2}
 \xymatrix 
 {
 X_J & & X_L\ar[ll]_{\pi_J^L}\ar[dl]^{\pi_K^L} \\
  & X_K\ar[ul]^{\pi_J^K} &  
 },\qquad J\subseteq K\subseteq L, \ J,K,L\in{\mathscr J}
\eeq
и поэтому у нее имеется проективный предел  
\beq\label{TH:o-povtornom-proj-predele-3}
X_I=\projlim_{J\in{\mathscr J}} X_J
\eeq
с неким структурным семейством проекций 
\beq\label{TH:o-povtornom-proj-predele-3-1}
\pi_J^I:X_J\gets X_I;
\eeq

\item[(iv)] для любого $i\in I$ и любого $J\in {\mathscr J}$ со свойством $i\in J$ композиция
\beq\label{TH:o-povtornom-proj-predele-4}
\pi_i^I:=\pi_i^J\circ \pi_J^I:X_i\gets X_I,
\qquad \xymatrix 
 {
 X_J\ar[dr]_{\pi_i^J} & & X_I\ar[ll]_{\pi_J^I}\ar@{-->}[dl]^{\pi_i^I} \\
  & X_i & 
 },\qquad i\in J\in {\mathscr J}
\eeq   
не зависит от выбора $J\in {\mathscr J}$ (со свойством $i\in J$), 

\item[(v)] семейство морфизмов $\{\pi_i^I;\ i\in I\}$ из \eqref{TH:o-povtornom-proj-predele-4} представляет собой проективный конус над контравариантной системой  $\{\pi_i^j;\ i,j\in I\}$: 
\beq\label{TH:o-povtornom-proj-predele-4-1}
 \xymatrix 
 {
 & X_I\ar[ld]_{\pi_i^I}\ar[rd]^{\pi_j^I} \\
   X_i & & X_j \ar[ll]_{\pi_i^j}
 },\qquad i\le j\in I
\eeq 
поэтому существует единственный морфизм в проективный предел
\beq\label{TH:o-povtornom-proj-predele-5}
\pi^I:X_I\to X=\projlim_{i\in I} X_i,
\eeq
замыкающий все диаграммы
\beq\label{TH:o-povtornom-proj-predele-5-1}
 \xymatrix 
 {
 X\ar[dr]_{\pi_i} & & X_I\ar@{-->}[ll]_{\pi^I}\ar[dl]^{\pi_i^I} \\
  & X_i &  
 },\qquad i\in I
\eeq

\item[(vi)]  морфизм $\pi^I$ из \eqref{TH:o-povtornom-proj-predele-5-1} является изоморфизмом.

}\eit
\etm

\bit{

\item[$\bullet$] Проективный предел 
\beq
X_I=\projlim_{J\in{\mathscr J}} X_J=\projlim_{J\in{\mathscr J}}\Big(\projlim_{j\in J} X_j\Big)
\eeq  
называется {\it повторным проективным пределом контравариантной системы $\{X_i,\pi_i^j;\ i\le j\in I\}$ относительно покрытия ${\mathscr J}$ индексного множества $I$}. В этой ситуации проективный предел этой ковариантной системы $X=\projlim_{i\in I} X_i$  
желательно называть каким-то новым термином, мы условимся использовать для этого словосочетание {\it исходный проективный предел}.
Теорема \ref{TH:o-povtornom-proj-predele} таким образом утверждает, что повторный проективный предел контравариантной системы всегда совпадает с ее исходным проективным пределом:
\beq\label{DEF:povtornyj-proj-limit}
\projlim_{J\in{\mathscr J}}\Big(\projlim_{j\in J} X_j\Big) =\projlim_{i\in I} X_i.
\eeq  
}\eit

\bpr[Доказательство теоремы \ref{TH:o-povtornom-proj-predele}.] 
1. Утверждение (i) очевидно.

2. Утверждение (ii) следует из того, что для любого $K\supseteq J$ система проекций \eqref{TH:o-povtornom-proj-predele-1-1}
$$
\pi_j^K:X_K\to X_j;\ j\in J,
$$
будет проективным конусом над ковариантной системой  $\{X_i,\pi_i^j;\ i\le j,\ i,j\in J\}$, и поэтому существует единственный морфизм \eqref{TH:o-povtornom-proj-predele-2} в проективный конус $X_J=\projlim_{j\in J}X_j$.

3.  Докажем утверждение (iii). Пусть $J\subseteq K\subseteq L$, $J,K,L\in{\mathscr J}$. Рассмотрим диаграмму \eqref{TH:o-povtornom-proj-predele-3-2} и, зафиксировав $i\in J$, дополним ее до диаграммы
\beq\label{PROOF:o-povtornom-proj-predele-0-1}
 \xymatrix 
 {
 X_J\ar@/_4ex/[ddr]_{\pi_i^J} & & X_L\ar@{-->}[ll]_{\pi_J^L}\ar@{-->}[dl]^{\pi_K^L}\ar@/^4ex/[ddl]^{\pi_i^L} \\
  & X_K\ar@{-->}[ul]^{\pi_J^K}\ar[d]_{\pi_i^K} &  \\
  & X_i & 
 },\qquad i\in J\subseteq K\subseteq L, \ J,K,L\in{\mathscr J}
\eeq
Здесь периметр и нижние внутренние треугольники (то есть те треугольники, у которых по одной пунктирной стрелке) коммутативны, потому что  варианты диаграммы \eqref{TH:o-povtornom-proj-predele-2-1} (а про верхний внутренний треугольник, который состоит из одних пунктирных стрелок, мы пока не знаем, коммутативен он или нет). Если из диаграммы выбросить стрелку $\pi_i^K$, то мы получим диаграмму  
\beq\label{PROOF:o-povtornom-proj-predele-0-2}
 \xymatrix 
 {
 X_J\ar@/_4ex/[ddr]_{\pi_i^J} & & X_L\ar@{-->}[ll]_{\pi_J^L}\ar@{-->}[dl]^{\pi_K^L}\ar@/^4ex/[ddl]^{\pi_i^L} \\
  & X_K\ar@{-->}[ul]^{\pi_J^K} &  \\
  & X_i & 
 },\qquad i\in J\subseteq K\subseteq L, \ J,K,L\in{\mathscr J}
\eeq
в которой периметр и нижний внутренний четырехугольник коммутативны. Эту диаграмму можно перерисовать так:
\beq\label{PROOF:o-povtornom-proj-predele-0-3}
 \xymatrix @R=2.5pc @C=2.5pc
 {
 \projlim_{j\in J}X_j\ar@{=}[r] & X_J\ar[dr]_{\pi_i^J} & & X_L\ar@{-->}@/_2ex/[ll]_{\pi_J^L}\ar@{-->}@/^2ex/[ll]^{\pi_J^K\circ\pi_K^L}\ar[dl]^{\pi_i^L} \\
   & & X_i & 
 },\qquad i\in J\subseteq K\subseteq L, \ J,K,L\in{\mathscr J}
\eeq
Здесь периметр и нижний внутренний треугольник коммутативны. Это миожно интерпретировать так: морфизмы $\pi_J^L$ и $\pi_J^K\circ\pi_K^L$ являются морфизмами проективного конуса $\pi_i^L:X_L\to X_i$ над контравариантной системой $\{X_i,\pi_i^j;\ i\le j,\ i,j\in J\}$ в проективный предел $X_J=\projlim_{j\in J}X_j$ этой контравариантной системы. Поскольку такой морфизм единственен, эти морфизмы должны совпадать:
\beq\label{PROOF:o-povtornom-proj-predele-0-3}
\pi_J^L=\pi_J^K\circ\pi_K^L.
\eeq
То есть коммутативна диаграмма \eqref{TH:o-povtornom-proj-predele-3-2}.

4.  Докажем утверждение (iv). Зафиксируем $i\in I$ и два множества $J,K\in{\mathscr J}$, содержащие $i$:
$$
i\in J,\quad i\in K.
$$ 
Нам нужно доказать коммутативность диаграммы
\beq\label{PROOF:o-povtornom-proj-predele-00-1}
 \xymatrix 
 {
 & X_I\ar[dl]_{\pi_J^I}\ar[dr]^{\pi_K^I} & \\
 X_J\ar[dr]_{\pi_i^J} && X_K\ar[dl]^{\pi_i^K} \\
  & X_i &  
 }
\eeq
Воспользуемся условием \eqref{DEF:pokrytie-indeksov-2} и подберем множество $L\in {\mathscr J}$ со свойством
$$
J\cup K\subseteq L
$$
Тогда диаграмму \eqref{PROOF:o-povtornom-proj-predele-00-1} можно будет дополнить до диаграммы 
\beq\label{PROOF:o-povtornom-proj-predele-00-2}
 \xymatrix  @R=2.5pc @C=4pc
 {
 & X_I\ar[dl]_{\pi_J^I}\ar[d]^{\pi_L^I}\ar[dr]^{\pi_K^I} & \\
 X_J\ar[dr]_{\pi_i^J} & X_L\ar[d]^{\pi_i^L}\ar[l]_{\pi_J^L}\ar[r]^{\pi_K^L} & X_K\ar[dl]^{\pi_i^K} \\
  & X_i &  
 }
\eeq
В ней верхние внутренние треугольники коммутативны, потому что $X_I$, как проективный предел контравариантной системы $\{X_J\}$, является ее проективным конусом. А нижние внутренние треугольники коммутативны, потому что это варианты диаграммы \eqref{TH:o-povtornom-proj-predele-2-1}. Из этого следует, что периметр, то есть диаграмма \eqref{PROOF:o-povtornom-proj-predele-00-1}, тоже коммутативен.  

5. Теперь докажем коммутативность диаграммы \eqref{TH:o-povtornom-proj-predele-4-1}. Для этого подберем множество $L\in {\mathscr J}$ так, чтобы
$$
i,j\in L
$$
(это всегда можно сделать, воспользовавшись условием \eqref{DEF:pokrytie-indeksov-2}). Тогда диаграмму \eqref{TH:o-povtornom-proj-predele-4-1} можно дополнить до диаграммы
\beq\label{TH:o-povtornom-proj-predele-00-3}
 \xymatrix 
 {
  & X_I\ar[d]^{\pi_L^I}\ar@/_4ex/[ddl]_{\pi_j^I}\ar@/^4ex/[ddr]^{\pi_j^I} &  \\
 & X_L\ar[dl]_{\pi_i^L}\ar[dr]^{\pi_j^L} &  \\
   X_i & & X_j \ar[ll]_{\pi_i^j}
 },\qquad i\le j\in I
\eeq
В ней верхние внутренние треугольники коммутативны, потому что это варианты диаграммы \eqref{TH:o-povtornom-proj-predele-4}. А нижний внутренний трегугольник коммутативен, потому то проективный предел $X_L=\projlim_{j\in L}X_i$ --- проективный конус над контравариантной системой $\{X_i,\pi_i^j;\ i\le j,\ i,j\in L\}$. Из этого следует, что периметр --- то есть диаграмма \eqref{TH:o-povtornom-proj-predele-4-1} --- тоже коммутативен.

6. Нам остается доказать утверждение (vi). 

6.1. Для всякого $J\in{\mathscr J}$ семейство $\{\pi_j: X\to X_j;\ j\in J\}$ является проективным конусом над конравариантной системой $\{X_i,\pi_i^j;\ i\le j,\ i,j\in J\}$. Поэтому существует единственный морфизм в ее проективный предел
\beq\label{PROOF:o-povtornom-proj-predele-1}
\eta_J:X\to X_J
\eeq
замыкающий все диаграммы
\beq\label{PROOF:o-povtornom-proj-predele-1-1}
 \xymatrix 
 {
 X\ar[dr]_{\pi_i}\ar@{-->}[rr]^{\eta_J} & & X_J\ar[dl]^{\pi_i^J} \\
  & X_i &  
 },\qquad i\in J
\eeq

6.2. Семейство $\{\pi_J: X\to X_J;\ J\in {\mathscr J}\}$ является проективным конусом над конравариантной системой $\{X_J,\pi_J^K;\ J\subseteq K,\ J,K\in {\mathscr J}\}$:
\beq\label{PROOF:o-povtornom-proj-predele-2}
 \xymatrix 
 {
 & X\ar[dl]_{\pi_J}\ar[dr]^{\pi_K} & \\
 X_J & & X_K\ar[ll]^{\pi_J^K} 
 },\qquad J\subseteq K,\ J,K\in {\mathscr J}
\eeq

6.3. Из пункта 6.2 следует, что существует единственный морфизм в проективный предел $X_I=\projlim_{J\in{\mathscr J}}X_J$ конравариантной системы $\{X_J,\pi_J^K;\ J\subseteq K,\ J,K\in {\mathscr J}\}$,
\beq\label{PROOF:o-povtornom-proj-predele-3}
\eta_I:X\to X_I
\eeq
замыкающий все диаграммы
\beq\label{PROOF:o-povtornom-proj-predele-4}
 \xymatrix 
 {
 X\ar[dr]_{\eta_J}\ar@{-->}[rr]^{\eta_I} & & X_I\ar[dl]^{\pi_J^I} \\
  & X_J &  
 },\qquad J\in{\mathscr J}
\eeq

6.4. Для всякого $i\in J$ диаграмму \eqref{PROOF:o-povtornom-proj-predele-4} можно дополнить до диаграммы
\beq\label{PROOF:o-povtornom-proj-predele-5}
 \xymatrix 
 {
 X\ar[dr]_{\eta_J}\ar@/_4ex/[ddr]_{\pi_i}\ar[rr]^{\eta_I} & & X_I\ar[dl]^{\pi_J^I}\ar@/^4ex/[ddl]^{\pi_i^I} \\
  & X_J\ar[d]^{\pi_i^J} &  \\
  & X_i &
 },\qquad i\in J\in{\mathscr J}
\eeq
в которой левый внутренний треугольник --- это диаграмма \eqref{PROOF:o-povtornom-proj-predele-2}, а правый внутренний треугольник --- это диаграмма \eqref{TH:o-povtornom-proj-predele-4}.

6.5. Если из диаграммы \eqref{PROOF:o-povtornom-proj-predele-5} выбросить внутреннюю вершину, то мы получим диаграмму
\beq\label{PROOF:o-povtornom-proj-predele-6}
 \xymatrix 
 {
 X\ar[dr]_{\pi_i}\ar[rr]^{\eta_I} & & X_I\ar[dl]^{\pi_i^I} \\
    & X_i &
 },\qquad i\in I.
\eeq
(коммутативную для любого $i\in I$). Ее можно дполнить до диаграммы
\beq\label{PROOF:o-povtornom-proj-predele-7}
 \xymatrix 
 {
 X\ar[dr]_{\pi_i}\ar[r]^{\eta_I} &  X_I\ar[d]^{\pi_i^I}\ar[r]^{\pi^I} & X \ar[dl]^{\pi_i} \\
    & X_i &
 },\qquad i\in I.
\eeq
в которой правый внутренний треугольник --- это диаграмма \eqref{TH:o-povtornom-proj-predele-5-1}. Теперь если выкинуть верхнюю вершину посередине, то мы получим диаграмму (коммутативную для любого $i\in I$):
\beq\label{PROOF:o-povtornom-proj-predele-8}
 \xymatrix 
 {
 X\ar[dr]_{\pi_i}\ar[rr]^{\pi^I\circ\eta_I} &   & X \ar[dl]^{\pi_i}\ar@{=}[r] & \projlim_{i\in I} X_i \\
    & X_i &
 },\qquad i\in I.
\eeq
С другой стороны, очевидно, для всех $i\in I$ коммутативна диаграмма 
\beq\label{PROOF:o-povtornom-proj-predele-8-1}
 \xymatrix 
 {
 X\ar[dr]_{\pi_i}\ar[rr]^{1_X} &   & X \ar[dl]^{\pi_i}\ar@{=}[r] & \projlim_{i\in I} X_i \\
    & X_i &
 },\qquad i\in I.
\eeq
Поскольку $X=\projlim_{i\in I} X_i$ --- проективный предел, морфизмы $\pi^I\circ\eta_I$ и $1_X$ должны совпадать:
\beq\label{PROOF:o-povtornom-proj-predele-8-2}
\pi^I\circ\eta_I=1_X
\eeq

6.6. Рассмотрим диаграмму \eqref{TH:o-povtornom-proj-predele-5-1}, поменяв в ней местами верхние вершины и 
заменив условие $i\in I$ на условие $i\in J\in {\mathscr J}$:
\beq\label{TH:o-povtornom-proj-predele-9-1}
 \xymatrix 
 {
 X_I\ar@{-->}[rr]^{\pi^I}\ar[dr]_{\pi_i^I} & & X\ar[dl]^{\pi_i}   \\
  & X_i &  
 },\qquad i\in J\in {\mathscr J}.
\eeq
Ее можно дополнить до диаграммы  
\beq\label{TH:o-povtornom-proj-predele-9-2}
 \xymatrix 
 {
 X_I\ar[dr]_{\pi_J^I}\ar@{-->}[rr]^{\pi^I}\ar@/_4ex/[ddr]_{\pi_i^I} & & X\ar[dl]^{\eta_J}\ar@/^4ex/[ddl]^{\pi_i}   \\
  & X_J\ar[d]^{\pi_i^J} &  \\
    & X_i &  
 },\qquad i\in J\in {\mathscr J}.
\eeq
Здесь левый нижний внутренний треугольник --- это диаграмма \eqref{TH:o-povtornom-proj-predele-4}, а правый нижний внутренний треугольник --- это диаграмма \eqref{PROOF:o-povtornom-proj-predele-1-1}. Из этой диаграммы мы теперь можем выкинуть правую нижнюю стрелку, и у нас получится диаграмма 
\beq\label{TH:o-povtornom-proj-predele-9-3}
 \xymatrix 
 {
 X_I\ar[dr]_{\pi_J^I}\ar@{-->}[rr]^{\pi^I}\ar@/_4ex/[ddr]_{\pi_i^I} & & X\ar[dl]^{\eta_J}   \\
  & X_J\ar[d]^{\pi_i^J} &  \\
    & X_i &  
 },\qquad i\in J\in {\mathscr J}.
\eeq
в которой периметр и нижний внутренний треугольник коммутативны. Это можно интерпретировать так: у нас имеется два параллельных морфизма, $\eta_J\cdot\pi^I$ и $\pi_J^I$, каждый из которых является морфизмом проективного конуса $\pi_i^I:X_I\to X_i$ над контравариантной системой $\{X_i,\pi_i^j;\ i\le j,\ i,j\in J\}$ в конус проективного предела $\pi_i^J:X_J\to X_i$ над этой системой:
\beq\label{TH:o-povtornom-proj-predele-9-4}
 \xymatrix  @R=2.5pc @C=2.5pc
 {
 X_I\ar@{-->}@/_2ex/[rr]_{\pi_J^I}\ar@{-->}@/^2ex/[rr]^{\eta_J\cdot\pi^I}\ar[dr]_{\pi_i^I} & & X_J\ar[dl]^{\pi_i^J}\ar@{=}[r] & \projlim_{j\in J}X_j \\
    & X_i &  
 },\qquad i\in J\in {\mathscr J}.
\eeq
Поскольку такой морфизм конусов единственен, мы получаем, что эти пунктирные стрелки должны совпадать:
\beq\label{TH:o-povtornom-proj-predele-9-5}
\eta_J\circ\pi^I=\pi_J^I,\qquad J\in {\mathscr J}.
\eeq
То есть должна быть коммутативна диаграмма 
\beq\label{TH:o-povtornom-proj-predele-9-6}
 \xymatrix 
 {
 X_I\ar[dr]_{\pi_J^I}\ar@{-->}[rr]^{\pi^I} & & X\ar[dl]^{\eta_J}   \\
  & X_J &  \\
     },\qquad J\in {\mathscr J}.
\eeq
Дополним ее до диаграммы
\beq\label{PROOF:o-povtornom-proj-predele-9}
 \xymatrix 
 {
 X_I\ar[dr]_{\pi_J^I}\ar[r]^{\pi^I} &  X\ar[d]^{\eta_J}\ar[r]^{\eta_I} & X_I \ar[dl]^{\pi_J^I} \\
    & X_J &
 },\qquad J\in {\mathscr J},
\eeq
в которой правый внутренний треугольник --- это диаграмма \eqref{PROOF:o-povtornom-proj-predele-4}.
Теперь если выкинуть верхнюю вершину посередине, то мы получим диаграмму (коммутативную для любого $J\in {\mathscr J}$):
\beq\label{PROOF:o-povtornom-proj-predele-9-7}
 \xymatrix 
 {
 X_I\ar[dr]_{\pi_J^I}\ar[rr]^{\eta_I\circ\pi^I} &   & X_I \ar[dl]^{\pi_J^I}\ar@{=}[r] & \projlim_{K\in {\mathscr J}} X_K \\
    & X_J &
 },\qquad J\in {\mathscr J}.
\eeq
С другой стороны, очевидно, для всех $i\in I$ коммутативна диаграмма 
\beq\label{PROOF:o-povtornom-proj-predele-9-8}
 \xymatrix 
 {
 X_I\ar[dr]_{\pi_J^I}\ar[rr]^{1_{X_I}} &   & X_I \ar[dl]^{\pi_J^I}\ar@{=}[r] & \projlim_{K\in {\mathscr J}} X_K \\
    & X_J &
 },\qquad J\in {\mathscr J}.
\eeq
Поскольку $X_I=\projlim_{K\in {\mathscr J}} X_K$ --- проективный предел, морфизмы $\eta_I\circ\pi^I$ и $1_{X_I}$ должны совпадать:
\beq\label{PROOF:o-povtornom-proj-predele-9-9}
\eta_I\circ\pi^I=1_{X_I}
\eeq

6.7. Равенства \eqref{PROOF:o-povtornom-proj-predele-8-2} и \eqref{PROOF:o-povtornom-proj-predele-9-9} вместе означают, что $\pi^I$ --- изоморфизм.
\epr

\paragraph{Теорема об узловом образе в повторном проективном пределе.}

\bit{
\item[$\bullet$] Напомним, что категория ${\sf K}$ называется категорией с {\it узловым разложением} \cite[1.5]{Akbarov-De-Gruyter-I}, если всякий морфизм $\ph:X\to Y$ раскладывается в композицию
\beq\label{uzlovoe-razlozhenie}
\ph=\im_\infty\ph\circ \red_\infty\ph\circ\coim_\infty\ph
\eeq
в которой
\bit{
\item[---] $\coim_\infty\ph$ --- строгий эпиморфизм,

\item[---] $\red_\infty\ph$ --- биморфизм,

\item[---] $\im_\infty\ph$ --- строгий мономорфизм.

}\eit 
Ниже мы почти всегда будем рассматривать композицию $\red_\infty\ph\circ\coim_\infty\ph$ как единое целое, не выделяя в нем множители, и для этого желательно иметь какое-то обозначение для этого образования. Мы будем пользоваться обозначением
\beq\label{e(ph)-uzlovoe-razlozhenie}
\e(\ph)=\red_\infty\ph\circ\coim_\infty\ph,
\eeq
и в свойствах этого морфизма для нас будет важно только то, что он является эпиморфизмом. Таким образом, в категории с узловым разложением всякий морфизм 
$\ph:X\to Y$ раскладывается в композицию
\beq\label{uzlovoe-razlozhenie-1}
\ph=\im_\infty\ph\circ \e(\ph).
\eeq
в которой
\bit{
\item[---] $\im_\infty\ph$ --- строгий мономорфизм,

\item[---] $\e(\ph)$ --- эпиморфизм.

}\eit 
 
}\eit

\btm[об узловом образе в повторном проективном пределе]\label{TH:ob-uzlov-obraze-v-povtornom-proj-predele}
Пусть даны:
\bit{

\item[(a)] проективно полная категория ${\sf K}$ с узловым разложением;

\item[(b)] контравариантная система $\{X_i,\pi_i^j;\ i\le j\in I\}$ в категории ${\sf K}$.

\item[(c)] покрытие ${\mathscr J}$ частично упорядоченного множества $(I,\le)$;

\item[(d)] проективный конус $\{\rho_i:Y\to X_i; \ i\in I\}$ над контравариантной системой $\{X_i,\pi_i^j;\ i\le j\in I\}$.
}\eit
Пусть  
\beq\label{TH:ob-uzlov-obraze-v-povtornom-proj-predele-*}
X=\projlim_{i\in I} X_i
\eeq  
--- проективный предел контравариантной системы $\{X_i,\pi_i^j;\ i\le j,\ i,j\in I\}$ со структурными проекциями 
\beq\label{TH:ob-uzlov-obraze-v-povtornom-proj-predele-**}
\pi_i:X\to X_i.
\eeq
Тогда 
\bit{

\item[(i)] для каждого множества $J\in{\mathscr J}$ семейство $\{\rho_i:Y\to X_i; \ i\in J\}$ является проективным конусом над контравариантной системой $\{X_i,\pi_i^j;\ i\le j\in J\}$, 
\beq\label{TH:ob-uzlov-obraze-v-povtornom-proj-predele-1-2}
 \xymatrix 
 {
  & Y\ar[dl]_{\rho_i}\ar[dr]^{\rho_j} & \\ 
 X_i & & X_j\ar[ll]_{\pi_i^j} \\
 },\qquad i\le j,\ i,j\in J
\eeq
и поэтому существует единственный морфизм
\beq\label{TH:ob-uzlov-obraze-v-povtornom-proj-predele-1}
\rho_J:X_J\gets Y,
\eeq
в проективный предел $X_J:=\projlim_{j\in J}X_j$ из \eqref{TH:o-povtornom-proj-predele-1}, замыкающий все диаграммы
\beq\label{TH:ob-uzlov-obraze-v-povtornom-proj-predele-1-1}
 \xymatrix 
 {
 X_J\ar[dr]_{\pi_i^J} & & Y\ar@{-->}[ll]_{\rho_J}\ar[dl]^{\rho_i} \\
  & X_i &  
 },\qquad i\in J
\eeq

\item[(ii)] семейство $\{\rho_J:Y\to X_J; \ J\in{\mathscr J}\}$ является проективным конусом контравариантной системы $\{X_J,\pi_J^K;\ J,K\in {\mathscr J}\}$:
\beq\label{TH:ob-uzlov-obraze-v-povtornom-proj-predele-2}
 \xymatrix 
 {
 X_J & & Y\ar[ll]_{\rho_J}\ar[dl]^{\rho_K} \\
  & X_K\ar[ul]^{\pi_J^K} &  
 },\qquad J\subseteq K,\ J,K\in {\mathscr J}
\eeq

\item[(iii)] каждый морфизм $\rho_J:X_J\gets Y$ из \eqref{TH:ob-uzlov-obraze-v-povtornom-proj-predele-1} однозначно определяет объект
\beq\label{TH:ob-uzlov-obraze-v-povtornom-proj-predele-3}
Y_J=\Im_\infty\rho_J,
\eeq
эпиморфизм
\beq\label{TH:ob-uzlov-obraze-v-povtornom-proj-predele-3-1}
\tau_J=\e(\rho_J),
\eeq
и строгий мономорфизм
\beq\label{TH:ob-uzlov-obraze-v-povtornom-proj-predele-3-2}
\sigma^J=\im_\infty\rho_J,
\eeq
для которых будет коммутативна диаграмма
\beq\label{TH:ob-uzlov-obraze-v-povtornom-proj-predele-3-4}
 \xymatrix 
 {
 X_J & & Y\ar[ll]_{\rho_J}\ar@{-->}[dl]^{\tau_J} \\
  & Y_J\ar@{-->}[ul]^{\sigma^J} &  
 }
\eeq

\item[(iv)] для  любых множеств $J,K\in{\mathscr J}$, связанных включением $J\subseteq K$, коммутативна диаграмма
\beq\label{TH:ob-uzlov-obraze-v-povtornom-proj-predele-4}
 \xymatrix  @R=2.5pc @C=4pc
 {
 & Y\ar@/_2ex/[dl]_{\tau_J}\ar@/^2ex/[dr]^{\tau_K} & \\
 Y_J\ar[d]_{\sigma_J} &  & Y_K\ar[d]^{\sigma_K} \\
 X_J &  & X_K\ar[ll]_{\pi_J^K} \\
 },\qquad J\subseteq K
\eeq
и существует единственный морфизм
\beq\label{TH:ob-uzlov-obraze-v-povtornom-proj-predele-4-1}
\sigma_J^K:Y_J\gets Y_K,
\eeq
замыкающий диаграмму
\beq\label{TH:ob-uzlov-obraze-v-povtornom-proj-predele-4-2}
 \xymatrix  @R=2.5pc @C=4pc
 {
& Y\ar@/_2ex/[dl]_{\tau_J}\ar@/^2ex/[dr]^{\tau_K} & \\
 Y_J\ar[d]_{\sigma_J} &  & Y_K\ar@{-->}[ll]_{\sigma_J^K}\ar[d]^{\sigma_K} \\
 X_J &  & X_K\ar[ll]_{\pi_J^K} \\
 },\qquad J\subseteq K
\eeq

\item[(v)] семейство $\{Y_J,\sigma_J^K;\ J,K\in {\mathscr J}\}$ образует контравариантную систему в категории ${\sf K}$, 
\beq\label{TH:ob-uzlov-obraze-v-povtornom-proj-predele-5-5}
 \xymatrix 
 {
 Y_J & & Y_L\ar[ll]_{\sigma_J^L}\ar[dl]^{\sigma_K^L} \\
  & Y_K\ar[ul]^{\sigma_J^K} &  
 },\qquad J\subseteq K\subseteq L, \ J,K,L\in{\mathscr J}
\eeq
над  которой система морфизмов $\tau_J:Y\to Y_J$ будет проективным конусом,  
\beq\label{TH:ob-uzlov-obraze-v-povtornom-proj-predele-5-1}
 \xymatrix 
 {
 & Y\ar[ld]_{\tau_J}\ar[rd]^{\tau_K} \\
   Y_J & & Y_K \ar[ll]_{\sigma_J^K}
 },\qquad J\subseteq K\in {\mathscr J}
\eeq
и по этой причине если обозначить через
\beq\label{TH:ob-uzlov-obraze-v-povtornom-proj-predele-5-2}
Y_I=\projlim_{J\in{\mathscr J}} Y_J
\eeq
проективный предел этой контравариантной системы, а через
\beq\label{TH:ob-uzlov-obraze-v-povtornom-proj-predele-5-3}
\sigma_J^I:Y_J\gets Y_I;
\eeq
семейство его структурных проекций, то найдется единственный морфизм 
\beq\label{TH:ob-uzlov-obraze-v-povtornom-proj-predele-5}
\tau_I:Y_I\gets Y
\eeq
замыкающий все диаграммы
\beq\label{TH:ob-uzlov-obraze-v-povtornom-proj-predele-5-4}
 \xymatrix 
 {
 Y_I\ar[dr]_{\sigma_J^I} & & Y\ar@{-->}[ll]_{\tau_I}\ar[ld]^{\tau_J} \\
  & Y_J & 
 },\qquad J\in {\mathscr J}
\eeq

\item[(vi)] для любого $i\in I$ и любого $J\in {\mathscr J}$ со свойством $i\in J$ композиция
\beq\label{TH:ob-uzlov-obraze-v-povtornom-proj-predele-6}
\sigma_i^I:=\pi_i^J\circ \sigma_J \circ \sigma_J^I:Y_I\to Y_i
\qquad
 \xymatrix 
 {
Y_J\ar[d]_{\sigma_J} & Y_I\ar[l]_{\sigma_J^I}\ar@{-->}[d]^{\sigma_i^I} \\
   X_J\ar[r]_{\pi_i^J} & X_i 
 },\qquad i\in J\in{\mathscr J}
\eeq
в которой $\sigma_J^I$ --- морфизм из \eqref{TH:ob-uzlov-obraze-v-povtornom-proj-predele-4-1}, 
$\sigma_J$ --- морфизм из \eqref{TH:ob-uzlov-obraze-v-povtornom-proj-predele-3-2}, а $\pi_i^J$ --- морфизм из \eqref{TH:o-povtornom-proj-predele-1-1}, не зависит от выбора $J\in {\mathscr J}$ (со свойством $i\in J$),

\item[(vii)] для всякого $i\in I$ коммутативна диаграмма
\beq\label{TH:ob-uzlov-obraze-v-povtornom-proj-predele-6*}
 \xymatrix  
 {
 Y_I\ar[dr]_{\sigma_i^I} & & Y \ar[ll]_{\tau_I}\ar[dl]^{\rho_i} \\
  & X_i & 
 }
\eeq
в которой $\tau_I$ --- морфизм из \eqref{TH:ob-uzlov-obraze-v-povtornom-proj-predele-5}. а $\sigma_i^I$ --- морфизм из  
\eqref{TH:ob-uzlov-obraze-v-povtornom-proj-predele-6},

\item[(viii)] семейство морфизмов $\{\sigma_i^I;\ i\in I\}$ из \eqref{TH:ob-uzlov-obraze-v-povtornom-proj-predele-6} представляет собой проективный конус над контравариантной системой $\{\pi_i^j;\ i,j\in I\}$: 
\beq\label{TH:ob-uzlov-obraze-v-povtornom-proj-predele-7-2}
 \xymatrix 
 {
 & Y_I\ar[ld]_{\sigma_i^I}\ar[rd]^{\sigma_j^I} \\
   X_i & & X_j \ar[ll]_{\pi_i^j}
 },\qquad i\le j\in I
\eeq
поэтому существует единственный морфизм\footnote{Здесь как было условлено в \eqref{TH:ob-uzlov-obraze-v-povtornom-proj-predele-*}, $X=\projlim_{i\in I} X_i$ --- проективный предел контравариантной системы $\{X_i,\pi_i^j;\ i\le j,\ i,j\in I\}$, а $\pi_i$ --- его структурные  проекции.} 
\beq\label{TH:ob-uzlov-obraze-v-povtornom-proj-predele-7}
\sigma^I:X\gets Y_I,
\eeq
замыкающий все диаграммы
\beq\label{TH:ob-uzlov-obraze-v-povtornom-proj-predele-7-3}
 \xymatrix 
 {
 X\ar[dr]_{\pi_i} & & Y_I\ar@{-->}[ll]_{\sigma^I}\ar[dl]^{\sigma_i^I} \\
  & X_i &  
 },\qquad i\in I
\eeq

\item[(ix)] для морфизма $\sigma^I$ из \eqref{TH:ob-uzlov-obraze-v-povtornom-proj-predele-7} коммутативна также диаграмма
\beq\label{TH:ob-uzlov-obraze-v-povtornom-proj-predele-8}
 \xymatrix  @R=2pc @C=4pc
 {
 & Y \ar@/_2ex/[ld]_{\rho_I}\ar[dd]^{\tau_I} \\
 X & \\
 &  Y_I\ar@/^2ex/@{-->}[lu]^{\sigma^I} \\
 }
\eeq
в которой $\tau_I$ --- морфизм из \eqref{TH:ob-uzlov-obraze-v-povtornom-proj-predele-5}, а $\rho_I$ --- естественный морфизм проективного конуса $\{\rho_i:Y\to X_i; \ i\in I\}$ в проективный предел $X$ контравариантной системы $\{X_i,\pi_i^j;\ i\le j\in I\}$, то есть такой, для которого коммутативны все диаграммы
\beq\label{TH:ob-uzlov-obraze-v-povtornom-proj-predele-8-1}
 \xymatrix @R=2.5pc @C=4pc
 {
 X\ar[dr]_{\pi_i} & & Y\ar@{-->}[ll]_{\rho_I}\ar[dl]^{\rho_i} \\
  & X_i &  
 },\qquad i\in I,
\eeq

\item[(x)] существует единственный морфизм, связывающий узловой образ морфизма $\tau_I$ из \eqref{TH:o-povtornom-proj-predele-5} с узловым образом морфизма $\rho_I$,
\beq\label{TH:ob-uzlov-obraze-v-povtornom-proj-predele-10}
\zeta:\Im_\infty\rho_I\gets\Im_\infty\tau_I,
\eeq
замыкающий диаграмму
\beq\label{TH:ob-uzlov-obraze-v-povtornom-proj-predele-10-1}
 \xymatrix  @R=2.5pc @C=8pc
 {
 & & Y \ar@/_8ex/[lldd]_{\rho_I}\ar@/^12ex/[dddd]^{\tau_I}\ar[dl]^{\e(\rho_I)}
 \ar[dd]_{\e(\tau_I)} \\
 & \Im_\infty\rho_I\ar[dl]^{\quad\im_\infty\rho_I} & \\
 X & & \Im_\infty\tau_I\ar@{-->}[ul]_{\zeta} \ar[dd]_{\im_\infty\tau_I} \\
 &  & \\
 &  & Y_I\ar@/^8ex/[lluu]_{\sigma^I} \\
 }
\eeq

\item[(xi)] морфизм $\zeta$ из \eqref{TH:ob-uzlov-obraze-v-povtornom-proj-predele-10} является изоморфизмом.

}\eit
\etm

\bit{

\item[$\bullet$] Если условиться обозначать символом 
\beq
\projlim_{i\in I}\rho_i
\eeq
морфизм, связывающий проективный конус $\{\rho_i:Y\to X_i; \ i\in I\}$ над контравариантной системой $\{X_i,\pi_i^j;\ i\le j\in I\}$ с проективным пределом этой конравариантной системы, то есть такой, для которого будут коммутативны все диаграммы
$$
 \xymatrix @R=2.5pc @C=4pc
 {
 \projlim_{i\in I} X_i\ar[dr]_{\pi_i} & & Y\ar@{-->}[ll]_{\projlim_{i\in I}\rho_i}\ar[dl]^{\rho_i} \\
  & X_i &  
 },\qquad i\in I,
$$
то из \eqref{DEF:povtornyj-proj-limit} будет следовать равенство
\beq
\projlim_{J\in{\mathscr J}}\Big(\projlim_{j\in J}\rho_j\Big)=\projlim_{i\in I}\rho_i
\eeq
а содержание теоремы \ref{TH:ob-uzlov-obraze-v-povtornom-proj-predele} можно будет выразить формулой
\beq
\im_\infty\left(
\projlim_{J\in{\mathscr J}}\Big(\projlim_{j\in J}\rho_j\Big)\right)=\im_\infty\projlim_{J\in{\mathscr J}}\Big(\im_\infty\projlim_{j\in J} \rho_j\Big)
\eeq  
или формулой
\beq
\im_\infty\left(
\projlim_{J\in{\mathscr J}}\circ\projlim_{j\in J}\right)=\Big(\im_\infty\projlim_{J\in{\mathscr J}}\Big)\circ \Big(\im_\infty\projlim_{j\in J}\Big)
\eeq  
которую можно интерпретировать как утверждение, что операция взятия узлового образа представляет собой в некотором смысле ``гомоморфизм для операции взятия проективного предела''.

}\eit

\bpr[Доказательство теоремы \ref{TH:ob-uzlov-obraze-v-povtornom-proj-predele}.]
1. Утверждение (i) очевидно.

2. Докажем утверждение (ii). Пусть $J\subseteq K$, $J,K\in{\mathscr J}$. Рассмотрим диаграмму \eqref{TH:ob-uzlov-obraze-v-povtornom-proj-predele-2} и, зафиксировав $i\in J$, дополним ее до диаграммы
\beq\label{PROOF:ob-uzlov-obraze-v-povtornom-proj-predele-2-1}
 \xymatrix 
 {
 X_J\ar@/_4ex/[ddr]_{\pi_i^J} & & Y\ar@{-->}[ll]_{\rho_J}\ar@{-->}[dl]^{\rho_K}\ar@/^4ex/[ddl]^{\rho_i} \\
  & X_K\ar@{-->}[ul]^{\pi_J^K}\ar[d]_{\pi_i^K} &  \\
  & X_i & 
 },\qquad i\in J\subseteq K, \ J,K\in{\mathscr J}
\eeq
Здесь периметр и правый внутренний треугольник коммутативны, потому что это варианты диаграммы \eqref{TH:ob-uzlov-obraze-v-povtornom-proj-predele-1-1}, а левый внутренний треугольник коммутативен, потому что это диаграмма \eqref{TH:o-povtornom-proj-predele-2-1} (а про верхний внутренний треугольник, который состоит из одних пунктирных стрелок, мы пока не знаем, коммутативен он, или нет). Если из диаграммы выбросить стрелку $\pi_i^K$, то мы получим диаграмму  
\beq\label{PROOF:ob-uzlov-obraze-v-povtornom-proj-predele-2-2}
 \xymatrix 
 {
 X_J\ar@/_4ex/[ddr]_{\pi_i^J} & & Y\ar@{-->}[ll]_{\rho_J}\ar@{-->}[dl]^{\rho_K}\ar@/^4ex/[ddl]^{\rho_i} \\
  & X_K\ar@{-->}[ul]^{\pi_J^K} &  \\
  & X_i & 
 },\qquad i\in J\subseteq K, \ J,K\in{\mathscr J}
\eeq
в которой периметр и нижний внутренний четырехугольник коммутативны. Эту диаграмму можно перерисовать так:
\beq\label{PROOF:ob-uzlov-obraze-v-povtornom-proj-predele-2-3}
 \xymatrix @R=2.5pc @C=2.5pc
 {
 \projlim_{j\in J}X_j\ar@{=}[r] & X_J\ar[dr]_{\pi_i^J} & & Y\ar@{-->}@/_2ex/[ll]_{\rho_J}\ar@{-->}@/^2ex/[ll]^{\pi_J^K\circ\rho_K}\ar[dl]^{\rho_i} \\
   & & X_i & 
 },\qquad i\in J\subseteq K\subseteq L, \ J,K,L\in{\mathscr J}
\eeq
Здесь периметр и нижний внутренний треугольник коммутативны и это можно интерпретировать так: морфизмы $\rho_J$ и $\pi_J^K\circ\rho_K$ являются морфизмами проективного конуса $\rho_i:Y\to X_i$ над контравариантной системой $\{X_i,\pi_i^j;\ i\le j,\ i,j\in J\}$ в проективный предел $X_J=\projlim_{j\in J}X_j$ этой контравариантной системы. Поскольку такой морфизм единственен, эти морфизмы должны совпадать:
\beq\label{PROOF:ob-uzlov-obraze-v-povtornom-proj-predele-2-4}
\rho_J=\pi_J^K\circ\rho_K.
\eeq
То есть коммутативна диаграмма \eqref{TH:ob-uzlov-obraze-v-povtornom-proj-predele-2}.

3. Утверждение (iii) следует из условия (a), по которому $\sf K$ --- категория с узловым разложением.

4. Диаграмму \eqref{TH:ob-uzlov-obraze-v-povtornom-proj-predele-4} можно дополнить  до диаграммы
\beq\label{PROOF:ob-uzlov-obraze-v-povtornom-proj-predele-4-1}
 \xymatrix  @R=2.5pc @C=4pc
 {
 & Y\ar@/_2ex/[dl]_{\tau_J}\ar@/^2ex/[dr]^{\tau_K}\ar@{-->}[ddl]^{\rho_J}\ar@{-->}[ddr]_{\rho_K} & \\
 Y_J\ar[d]_{\sigma_J} &  & Y_K\ar[d]^{\sigma_K} \\
 X_J &  & X_K\ar[ll]_{\pi_J^K} \\
 },\qquad J\subseteq K
\eeq
в которой верхние внутренние треугольники будут коммутативны, потому что это диаграммы \eqref{TH:ob-uzlov-obraze-v-povtornom-proj-predele-3-4}, а нижний внутренний треугольник коммутативен, потому что это диаграмма   
\eqref{TH:ob-uzlov-obraze-v-povtornom-proj-predele-2}. Поскольку все внутренние треугольники коммутативны, периметр, то есть диаграмма \eqref{TH:ob-uzlov-obraze-v-povtornom-proj-predele-4}, тоже должен быть коммутативен.

Теперь если \eqref{TH:ob-uzlov-obraze-v-povtornom-proj-predele-4} перерисовать в виде
\beq\label{PROOF:ob-uzlov-obraze-v-povtornom-proj-predele-4-2}
 \xymatrix  @R=2.5pc @C=4pc
 {
 & Y\ar@/_2ex/[dl]_{\tau_J}\ar@/^2ex/[dr]^{\tau_K} & \\
 Y_J\ar[d]_{\sigma_J} &  & Y_K\ar@/^2ex/[dll]^{\pi_J^K\circ\sigma_K} \\
 X_J &  &  \\
 },\qquad J\subseteq K
\eeq
и заметить, что в силу (iii), здесь $\tau_K$ --- эпиморфизм, а $\sigma_J$ --- строгий моономорфизм, то станет понятно, что должна  существовать (и быть единственной) диагональ $\sigma_J^K$ этого четырехугольника:
\beq\label{PROOF:ob-uzlov-obraze-v-povtornom-proj-predele-4-2}
 \xymatrix  @R=2.5pc @C=4pc
 {
 & Y\ar@/_2ex/[dl]_{\tau_J}\ar@/^2ex/[dr]^{\tau_K} & \\
 Y_J\ar[d]_{\sigma_J} &  & Y_K\ar@/^2ex/[dll]^{\pi_J^K\circ\sigma_K}\ar@{-->}[ll]_{\sigma_J^K} \\
 X_J &  &  \\
 },\qquad J\subseteq K
\eeq
Это тот самый морфизм, который делает коммутативной диаграмму \eqref{TH:ob-uzlov-obraze-v-povtornom-proj-predele-4-2}.

5. В пункте (v) можно сразу заметить, что диаграмма \eqref{TH:ob-uzlov-obraze-v-povtornom-proj-predele-5-1} коммутативна просто потому, что это верхний внутренний треугольник в диаграмме \eqref{TH:ob-uzlov-obraze-v-povtornom-proj-predele-4-2}. Поэтому здесь важно только доказать коммутативность диаграммы \eqref{TH:ob-uzlov-obraze-v-povtornom-proj-predele-5-5}. Впишем ее в диаграмму
\beq\label{PROOF:ob-uzlov-obraze-v-povtornom-proj-predele-5-5}
 \xymatrix 
 {
 Y_J\ar[dd]_{\sigma_J} & & Y_L\ar[dd]^{\sigma_L}\ar[ll]_{\sigma_J^L}\ar[dl]_{\sigma_K^L} \\
  & Y_K\ar[ul]_{\sigma_J^K}\ar[dd]^(.3){\sigma_K} &  \\
   X_J & & X_L\ar@{-->}[ll]_(.7){\pi_J^L}\ar[dl]^{\pi_K^L} \\
  & X_K\ar[ul]^{\pi_J^K} &  
 },\qquad J\subseteq K\subseteq L, \ J,K,L\in{\mathscr J}
\eeq
Здесь боковые стороны коммутативны,  потому что это внутренний прямоугольник в  \eqref{TH:ob-uzlov-obraze-v-povtornom-proj-predele-4-2}, а нижнее основание --- потому что это диаграмма \eqref{TH:o-povtornom-proj-predele-4-1}. Из этого следует, что если двигаться из начала этой диаграмма (вершины $Y_L$) в ее конец (вершину $X_J$), то выбираемые маршруты должны давать совпадающие морфизмы. В частности, мы получим равенство
$$
\sigma_J\circ \sigma_J^L=\sigma_J\circ \sigma_J^K\circ \sigma_K^L
$$
Поскольку $\sigma_J$ --- мономорфизм, мы можем на него сократить:
$$
\sigma_J^L=\sigma_J^K\circ \sigma_K^L
$$
Это как раз означает коммутативность диаграммы \eqref{TH:ob-uzlov-obraze-v-povtornom-proj-predele-5-5}.

6. Докажем утверждение (vi). Зафиксируем $i\in I$ и два множества $J,K\in{\mathscr J}$, содержащие $i$:
$$
i\in J,\quad i\in K.
$$ 
Нам нужно доказать коммутативность диаграммы
\beq\label{PROOF:ob-uzlov-obraze-v-povtornom-proj-predele-00-1}
 \xymatrix 
 {
 & Y_I\ar[dl]_{\sigma_J^I}\ar[dr]^{\sigma_K^I} & \\
 Y_J\ar[d]_{\sigma_J} && Y_K\ar[d]^{\sigma_K} \\ 
 X_J\ar[dr]_{\pi_i^J} && X_K\ar[dl]^{\pi_i^K} \\
  & X_i &  
 }
\eeq
Воспользуемся условием \eqref{DEF:pokrytie-indeksov-2} и подберем множество $L\in {\mathscr J}$ со свойством
$$
J\cup K\subseteq L
$$
Тогда диаграмму \eqref{PROOF:ob-uzlov-obraze-v-povtornom-proj-predele-00-1} можно будет дополнить до диаграммы 
\beq\label{PROOF:ob-uzlov-obraze-v-povtornom-proj-predele-00-2}
 \xymatrix  @R=2.5pc @C=4pc
 {
 & Y_I\ar[dl]_{\sigma_J^I}\ar[d]^{\sigma_L^I}\ar[dr]^{\sigma_K^I} & \\
 Y_J\ar[d]_{\sigma_J} & Y_L\ar[d]^{\sigma_L}\ar[l]^{\sigma_J^L}\ar[r]_{\sigma_K^L} & Y_K\ar[d]^{\sigma_K} \\ 
 X_J\ar[dr]_{\pi_i^J} & X_L\ar[d]^{\pi_i^L}\ar[l]_{\pi_J^L}\ar[r]^{\pi_K^L} & X_K\ar[dl]^{\pi_i^K} \\
  & X_i &  
 }
\eeq
В ней верхние внутренние треугольники коммутативны, потому что $Y_I$, как проективный предел контравариантной системы $\{Y_J\}$, является ее проективным конусом. Внутренние прямоугольники коммутативны, потому что это прмоугольники из диаграмм вида \eqref{TH:ob-uzlov-obraze-v-povtornom-proj-predele-4-2}. А нижние внутренние треугольники коммутативны, потому что это варианты диаграммы \eqref{TH:o-povtornom-proj-predele-2-1}. Из этого следует, что периметр, то есть диаграмма \eqref{PROOF:ob-uzlov-obraze-v-povtornom-proj-predele-00-1}, тоже коммутативен.

7. Выберем произвольное множество $J$ со свойством $i\in J$ и представим диаграмму \eqref{TH:ob-uzlov-obraze-v-povtornom-proj-predele-6*} в виде
\beq\label{PROOF:ob-uzlov-obraze-v-povtornom-proj-predele-6*-1}
 \xymatrix  
 {
 Y_I\ar[d]_{\sigma_J^I} & & Y \ar[ll]_{\tau_I}\ar[d]^{\rho_i} \\
 Y_J\ar[r]_{\sigma_J} & X_J\ar[r]_{\pi_i^J} & X_i
 }
\eeq
Ее можно дополнить до диаграммы
\beq\label{PROOF:ob-uzlov-obraze-v-povtornom-proj-predele-6*-2}
 \xymatrix  
 {
 Y_I\ar[d]_{\sigma_J^I} & & Y \ar[ll]_{\tau_I}\ar[d]^{\rho_i}\ar[lld]_{\tau_J}\ar[ld]^{\rho_J} \\
 Y_J\ar[r]_{\sigma_J} & X_J\ar[r]_{\pi_i^J} & X_i
 }
\eeq
Здесь все внутренние треугольники коммутативны, потому что это (слева направо) диаграммы \eqref{TH:ob-uzlov-obraze-v-povtornom-proj-predele-5-4}, \eqref{TH:ob-uzlov-obraze-v-povtornom-proj-predele-3-4} и \eqref{TH:ob-uzlov-obraze-v-povtornom-proj-predele-1-1}. Поэтому периметр тоже коммутативен, а это диаграмма \eqref{PROOF:ob-uzlov-obraze-v-povtornom-proj-predele-6*-1}.

8. Теперь докажем коммутативность диаграммы \eqref{TH:ob-uzlov-obraze-v-povtornom-proj-predele-7-2}. Для этого подберем множество $L\in {\mathscr J}$ так, чтобы
$$
i,j\in L
$$
(это всегда можно сделать, воспользовавшись условием \eqref{DEF:pokrytie-indeksov-2}). Тогда диаграмму \eqref{TH:o-povtornom-proj-predele-4-1} можно дополнить до диаграммы
\beq\label{TH:o-povtornom-proj-predele-00-3}
 \xymatrix 
 {
  & Y_I\ar[d]^{\sigma_L^I}\ar@/_4ex/[dddl]_{\sigma_j^I}\ar@/^4ex/[dddr]^{\sigma_j^I} &  \\
 & Y_L\ar[d]_{\sigma_L} &  \\
 & X_L\ar[dl]_{\pi_i^L}\ar[dr]^{\pi_j^L} &  \\
   X_i & & X_j \ar[ll]_{\pi_i^j}
 },\qquad i\le j\in I
\eeq
В ней верхние внутренние треугольники коммутативны, потому что это варианты диаграммы \eqref{TH:ob-uzlov-obraze-v-povtornom-proj-predele-6}. А нижний внутренний трегугольник коммутативен, потому то проективный предел $X_L=\projlim_{j\in L}X_i$ --- проективный конус над контравариантной системой $\{X_i,\pi_i^j;\ i\le j,\ i,j\in L\}$. Из этого следует, что периметр --- то есть диаграмма \eqref{TH:o-povtornom-proj-predele-4-1} --- тоже коммутативен.

9. Представим диаграмму \eqref{TH:ob-uzlov-obraze-v-povtornom-proj-predele-8} в виде
\beq\label{PROOF:ob-uzlov-obraze-v-povtornom-proj-predele-8-1}
 \xymatrix  @R=2pc @C=4pc
 {
 X & & Y \ar[ll]_{\rho_I}\ar[dl]^{\tau_I} \\
 &  Y_I\ar[lu]^{\sigma^I} & \\
 }
\eeq
и впишем ее в диаграмму
\beq\label{PROOF:ob-uzlov-obraze-v-povtornom-proj-predele-8-2}
 \xymatrix  @R=2pc @C=4pc
 {
 X\ar@/_4ex/[ddr]_{\pi_i} & & Y \ar[ll]_{\rho_I}\ar[dl]^{\tau_I}\ar@/^4ex/[ddl]^{\rho_i} \\
 &  Y_I\ar[lu]^{\sigma^I}\ar[d]^{\sigma_i^I} & \\
 & X_i & 
 }
\eeq
Здесь периметр коммутативен, потому что это диаграмма \eqref{TH:ob-uzlov-obraze-v-povtornom-proj-predele-8-1}, левый нижний внутренний треугольник коммутативен, потому что это диаграмма \eqref{TH:ob-uzlov-obraze-v-povtornom-proj-predele-7-3}, 
а правый нижний внутренний треугольник --- потому что это диаграмма \eqref{TH:ob-uzlov-obraze-v-povtornom-proj-predele-6*}.

Если из \eqref{PROOF:ob-uzlov-obraze-v-povtornom-proj-predele-8-2} выбросить стррелку $\sigma_i^I$, то мы получим диаграмму
\beq\label{PROOF:ob-uzlov-obraze-v-povtornom-proj-predele-8-3}
 \xymatrix  @R=2pc @C=4pc
 {
 X\ar@/_4ex/[ddr]_{\pi_i} & & Y \ar[ll]_{\rho_I}\ar[dl]^{\tau_I}\ar@/^4ex/[ddl]^{\rho_i} \\
 &  Y_I\ar[lu]^{\sigma^I} & \\
 & X_i & 
 }
\eeq
в которой периметр и нижний внутренний четрехугольник коммутативны. Эту диаграмму можно перерисовать так:
\beq\label{PROOF:ob-uzlov-obraze-v-povtornom-proj-predele-8-4}
 \xymatrix @R=2.5pc @C=2.5pc
 {
 \projlim_{i\in I}X_i\ar@{=}[r] & X\ar[dr]_{\pi_i} & & Y\ar@{-->}@/_2ex/[ll]_{\rho_I}\ar@{-->}@/^2ex/[ll]^{\sigma^I\circ\tau_I}\ar[dl]^{\rho_i} \\
   & & X_i & 
 },\qquad i\in J\subseteq K\subseteq L, \ J,K,L\in{\mathscr J}
\eeq
Здесь периметр и нижний внутренний треугольник коммутативны и это можно интерпретировать так: морфизмы $\rho_I$ и $\sigma^I\circ\tau_I$ являются морфизмами проективного конуса $\rho_i:Y\to X_i$ над контравариантной системой $\{X_i,\pi_i^j;\ i\le j,\ i,j\in I\}$ в проективный предел $X=\projlim_{i\in I}X_i$ этой контравариантной системы. Поскольку такой морфизм единственен, эти морфизмы должны совпадать:
\beq\label{PROOF:ob-uzlov-obraze-v-povtornom-proj-predele-8-5}
\rho_I=\sigma^I\circ\tau_I.
\eeq
То есть коммутативна диаграмма \eqref{PROOF:ob-uzlov-obraze-v-povtornom-proj-predele-8-1}.

10. Рассмортрим диаграмму
\beq\label{PROOF:ob-uzlov-obraze-v-povtornom-proj-predele-10}
 \xymatrix  @R=2.5pc @C=8pc
 {
 & & Y \ar@/_8ex/[lldd]_{\rho_I}\ar@/^12ex/[dddd]^{\tau_I}\ar[dl]^{\e(\rho_I)}
 \ar[dd]_{\e(\tau_I)} \\
 & \Im_\infty\rho_I\ar[dl]^{\quad\im_\infty\rho_I} & \\
 X & & \Im_\infty\tau_I \ar[dd]_{\im_\infty\tau_I} \\
 &  & \\
 &  & Y_I\ar@/^8ex/[lluu]_{\sigma^I} \\
 }
\eeq
Во внутреннем четырехугольнике
\beq\label{PROOF:ob-uzlov-obraze-v-povtornom-proj-predele-10-1}
 \xymatrix  @R=2.5pc @C=8pc
 {
 & & Y \ar[dl]^{\e(\rho_I)}
 \ar[dd]_{\e(\tau_I)} \\
 & \Im_\infty\rho_I\ar[dl]^{\quad\im_\infty\rho_I} & \\
 X & & \Im_\infty\tau_I\ar@{-->}[ul]_{\zeta} \ar[dd]_{\im_\infty\tau_I} \\
 &  & \\
 &  & Y_I\ar@/^8ex/[lluu]_{\sigma^I} \\
 }
\eeq
$\e(\tau_I)$ --- эпиморфизм, а $\im_\infty\rho_I$ --- строгий мономорфизм. Поэтому существует единственный морфизм $\zeta:\Im_\infty\tau_I\to \Im_\infty\rho_I$, для которой коммутативна диаграмма
\beq\label{PROOF:ob-uzlov-obraze-v-povtornom-proj-predele-10-2}
 \xymatrix  @R=2.5pc @C=8pc
 {
 & & Y \ar[dl]^{\e(\rho_I)}
 \ar[dd]_{\e(\tau_I)} \\
 & \Im_\infty\rho_I\ar[dl]^{\quad\im_\infty\rho_I} & \\
 X & & \Im_\infty\tau_I\ar@{-->}[ul]_{\zeta} \ar[dd]_{\im_\infty\tau_I} \\
 &  & \\
 &  & Y_I\ar@/^8ex/[lluu]_{\sigma^I} \\
 }
\eeq

11. Нам остается убедиться, что морфизм $\zeta$ является изоморфизмом. 

11.1. Сначала докажем коммутативность диаграммы
\beq\label{PROOF:ob-uzlov-obraze-v-povtornom-proj-predele-11-1}
 \xymatrix  @R=2.5pc @C=2.5pc
 {
X\ar[d]_{\eta_J} & & Y\ar[ll]_{\rho_I}\ar[d]^{\tau_J} \\
X_J & & Y_J\ar[ll]^{\sigma_J}
 }
\eeq
в которой $\eta_J:X\to X_J$ --- морфизм из \eqref{PROOF:o-povtornom-proj-predele-1}.

11.2. Перерисуем диаграмму \eqref{PROOF:ob-uzlov-obraze-v-povtornom-proj-predele-11-1} в виде
\beq\label{PROOF:ob-uzlov-obraze-v-povtornom-proj-predele-11-2}
 \xymatrix  @R=2.5pc @C=8pc
 {
X\ar[d]_{\eta_J} & \Im_\infty\rho_I\ar[l]_{\im_\infty\rho_I} & Y\ar@/_8ex/[ll]_{\rho_I}\ar[l]_{\e(\rho_I)} \ar[d]^{\tau_J} \\
X_J & & Y_J\ar[ll]^{\sigma_J}
 }
\eeq
Здесь $\e(\rho_I)$ --- эпиморфизм, а $\sigma_J$ --- строгий мономорфизм. Поэтому существует и единственна стрелка $\xi_J:\Im_\infty\rho_I\to Y_J$, для которой коммутативна диаграмма
\beq\label{PROOF:ob-uzlov-obraze-v-povtornom-proj-predele-11-2-1}
 \xymatrix  @R=2.5pc @C=8pc
 {
X\ar[d]_{\eta_J} & \Im_\infty\rho_I\ar[l]_{\im_\infty\rho_I}\ar@{-->}[dr]_{\xi_J} & Y\ar@/_8ex/[ll]_{\rho_I}\ar[l]_{\e(\rho_I)} \ar[d]^{\tau_J} \\
X_J & & Y_J\ar[ll]^{\sigma_J}
 }
\eeq

11.3. Покажем, что система морфизмов $\xi_J:\Im_\infty\rho_I\to Y_J$ --- проективный конус над контравариантной системой $\{\sigma_J^K;\ J,K\in{\mathscr J}\}$:
\beq\label{PROOF:ob-uzlov-obraze-v-povtornom-proj-predele-11-3}
 \xymatrix  @R=2.5pc @C=2.5pc
 {
 & \Im_\infty\rho_I \ar@/_2ex/[dl]_{\xi_J}\ar@/^2ex/[dr]^{\xi_K} &  \\
Y_J & & Y_K\ar[ll]^{\sigma_J^K}
 },\qquad J\subseteq K, \ J,K\in{\mathscr J}.
\eeq

11.4. Из \eqref{PROOF:ob-uzlov-obraze-v-povtornom-proj-predele-11-3} следует, что существует единственный морфизм
\beq\label{PROOF:ob-uzlov-obraze-v-povtornom-proj-predele-11-4}
\xi_I:\Im_\infty\rho_I\to Y_I=\projlim_{K\in{\mathscr J}}Y_K
\eeq
замыкающий все диаграммы
\beq\label{PROOF:ob-uzlov-obraze-v-povtornom-proj-predele-11-4-1}
 \xymatrix  @R=2.5pc @C=2.5pc
 {
  \Im_\infty\rho_I \ar[dr]_{\xi_J}\ar@{-->}[rr]^{\xi_I} & & Y_I\ar[dl]^{\sigma_J^I}\ar@{=}[r] & \projlim_{K\in{\mathscr J}}Y_K \\
  & Y_J & & 
 },\qquad  J\in{\mathscr J}.
\eeq

11.5. Морфизм $\xi_I$ из \eqref{PROOF:ob-uzlov-obraze-v-povtornom-proj-predele-11-4} замыкает диаграмму
\beq\label{PROOF:ob-uzlov-obraze-v-povtornom-proj-predele-11-5}
 \xymatrix  @R=2.5pc @C=8pc
 {
X\ar[d]_{\eta_I} & \Im_\infty\rho_I\ar[l]_{\im_\infty\rho_I}\ar@{-->}[dr]_{\xi_I} & Y\ar@/_8ex/[ll]_{\rho_I}\ar[l]_{\e(\rho_I)} \ar[d]^{\tau_I} \\
X_I & & Y_I\ar[ll]^{\sigma_I}
 }
\eeq
где $\eta_I$ --- морфизм из \eqref{PROOF:o-povtornom-proj-predele-3}, а $\sigma_I:X_I=\projlim_{K\in{\mathscr J}}X_K\gets \projlim_{K\in{\mathscr J}}Y_K=Y_I$ --- проективный предел морфизмов 
$\sigma_J:X_J\gets Y_J$.

11.6. Выделим в \eqref{PROOF:ob-uzlov-obraze-v-povtornom-proj-predele-11-5} фрагмент
\beq\label{PROOF:ob-uzlov-obraze-v-povtornom-proj-predele-11-6}
 \xymatrix  @R=2.5pc @C=8pc
 {
 \Im_\infty\rho_I\ar[dr]_{\xi_I} & Y\ar[l]_{\e(\rho_I)} \ar[d]^{\tau_I} \\
 & Y_I
 }
\eeq
и перерисуем его в виде
\beq\label{PROOF:ob-uzlov-obraze-v-povtornom-proj-predele-11-6-1}
 \xymatrix  @R=2.5pc @C=8pc
 {
 \Im_\infty\rho_I\ar@/_4ex/[ddr]_{\xi_I} & Y\ar[l]_{\e(\rho_I)}\ar[d]_{\e(\tau_I)} \ar@/^12ex/[dd]^{\tau_I} \\
 & \Im_\infty\tau_I\ar[d]_{\im_\infty\tau_I} \\
 & Y_I
 }
\eeq
Здесь во внутреннем четырехугольнике $\e(\rho_I)$ --- эпиморфизм, а $\im_\infty\tau_I$ --- строгий мономорфизм. Поэтому существует и единственна стрелка $\xi:\Im_\infty\rho_I\to \Im_\infty\tau_I$, для которой будет коммутативна диаграмма 
\beq\label{PROOF:ob-uzlov-obraze-v-povtornom-proj-predele-11-6-2}
 \xymatrix  @R=2.5pc @C=8pc
 {
 \Im_\infty\rho_I\ar@/_4ex/[ddr]_{\xi_I}\ar@{-->}@/_2ex/[dr]_{\xi} & Y\ar[l]_{\e(\rho_I)}\ar[d]_{\e(\tau_I)} \ar@/^12ex/[dd]^{\tau_I} \\
 & \Im_\infty\tau_I\ar[d]_{\im_\infty\tau_I} \\
 & Y_I
 }
\eeq

11.7. Нам остается убедиться, что $\zeta$ и $\xi$ --- взаимно обратные морфизмы. Во-первых,
$$
\zeta\circ\overbrace{\xi\circ\underbrace{\e(\rho_I)}_{\scriptsize \begin{matrix} \text{\rotatebox{90}{$\owns$}} \\ \Epi \end{matrix}}}^{\scriptsize \begin{matrix}\e(\tau_I)\\ \| \end{matrix}}=\zeta\circ\e(\tau_I)=\e(\rho_I)=1_{\Im_\infty\rho_I}\circ\underbrace{\e(\rho_I)}_{\scriptsize \begin{matrix} \text{\rotatebox{90}{$\owns$}} \\ \Epi \end{matrix}}
\qquad\Rightarrow\qquad \zeta\circ\xi=1_{\Im_\infty\rho_I}
$$
И, во-вторых,
$$
\xi\circ\overbrace{\zeta\circ\underbrace{\e(\tau_I)}_{\scriptsize \begin{matrix} \text{\rotatebox{90}{$\owns$}} \\ \Epi \end{matrix}}}^{\scriptsize \begin{matrix}\e(\rho_I)\\ \| \end{matrix}}=\xi\circ\e(\rho_I)=\e(\tau_I)=1_{\Im_\infty\tau_I}\circ\underbrace{\e(\tau_I)}_{\scriptsize \begin{matrix} \text{\rotatebox{90}{$\owns$}} \\ \Epi \end{matrix}}
\qquad\Rightarrow\qquad \xi\circ\zeta=1_{\Im_\infty\tau_I}
$$
\epr

\section{Некоторые факты из теории локально выпуклых и стереотипных пространств}

\subsection{Свойства поляр}

Полярой множества $A$ в локально выпуклом пространстве $X$ называется множество
\beq\label{DEF:A^circ}
A^\circ=\{f\in X^\star: \ \sup_{x\in A}\abs{f(x)}\le 1\}.
\eeq
Обратной полярой множества $F\subseteq X^\star$ называется множество
\beq\label{DEF:^circ-F}
^\circ F=\{x\in X: \ \sup_{f\in F}\abs{f(x)}\le 1\}.
\eeq

Пусть $X$ ---- локально выпуклое пространство. Тогда для любых множеств $A,B\subseteq X$
\beq\label{(cabsconv A)^circ=A^circ}
(\cabsconv A)^\circ=A^\circ,
\eeq 
\beq\label{(C-cdot-A)^circ=1/C-cdot-A^circ}
(\lambda\cdot A)^\circ=\frac{1}{\lambda}\cdot A^\circ,\quad \lambda\in\C\setminus\{0\}
\eeq 
\beq\label{(A-cup-B)^circ=A^circ-cap-B^circ}
(A\cup B)^\circ=A^\circ\cap B^\circ,
\eeq 
а если $A$ и $B$ замкнуты и абсолютно выпуклы в $X$, то 
\beq\label{A^circ-circ=A}
^\circ(A^\circ)=A
\eeq 
\beq\label{(A-cap-B)^circ=cabsconv(A^circ-cup-B^circ)}
(A\cap B)^\circ=\cabsconv(A^\circ\cup B^\circ),
\eeq 
\bpr
Здесь последняя формула доказана в \cite[(3.30)]{Akbarov-De-Gruyter-I}, предпоследняя -- вариант теоремы о биполяре \cite[Theorem 3.1.11]{Akbarov-De-Gruyter-I} а первые три очевидны, например, третья доказывается цепочкой
\begin{multline*}
f\in (A\cap B)^\circ \quad\Leftrightarrow\quad \sup_{x\in A\cup B}\abs{f(x)}\le 1
\quad\Leftrightarrow\quad \sup_{x\in A}\abs{f(x)}\le 1\ \& \ \sup_{x\in B}\abs{f(x)}\le 1
\quad\Leftrightarrow\quad f\in A^\circ \ \& \ f\in B^\circ.
\end{multline*}
\epr

\subsection{Кополные пространства}
Мы, как обычно, называем локально выпуклое пространство $X$ {\it полным}, если в нем всякая направленность Коши сходится. Двойственным к этому свойству в теории стереотипных пространств является следующее: локально выпуклое пространство $X$ называется {\it кополным}\label{DEF:kopolno}, если всякий линейный функционал $f:X\to\C$, непрерывный на любом вполне ограниченном множестве $K\subseteq X$, непрерывен на всем пространстве $X$. Двойственность между этими понятиями в классе стереотипных пространств описывается эквивалентностью \cite[Section 2.2]{Ak03}:
$$
\text{$X$ кополно}\qquad\Longleftrightarrow\qquad \text{$X^\star$ полно}.
$$

Ниже мы будем пользоваться следующими двумя определениями из работы \cite{Ak16} ($\mathcal{U}(X)$ обозначает систему окрестностей нуля в пространстве $X$).
\bit{
\item[$\bullet$]
Линейное отображение локально выпуклых пространств $\ph:X\to Y$ мы называем {\it открытым}\label{DEF:open-map}, если образ $\ph(U)$ всякой окрестности нуля $U\subseteq X$ является окрестностью нуля в подпространстве $\ph(X)$ ЛВП $Y$ (с индуцированной из $Y$ топологией):
$$
\forall
U\in \mathcal{U}(X) \quad \exists V\in \mathcal{U}(Y) \quad \ph(U)\supseteq
\ph(X)\cap V.
$$

\item[$\bullet$]
Линейное непрерывное отображение локально выпуклых пространств $\ph:X\to Y$ мы называем
{\it замкнутым}, если для всякого вполне ограниченного множества $T\subseteq \overline{\ph(X)}\subseteq Y$ найдется вполне ограниченное множество $S\subseteq X$ такое, что $\ph(S)\supseteq T$. (Понятно, что это, в частности, означает, что множество значений $\ph(X)$ отображения $\ph$ должно быть замкнуто в $Y$.)
}\eit
В классе стереотипных пространств открытость и замкнутость линейного непрерывного отображения (то есть морфизма в этой категории) -- двойственные свойства \cite[Теорема 2.10]{Ak16}:
$$
\text{$\ph:X\to Y$ замкнуто}\qquad\Longleftrightarrow\qquad \text{$\ph^\star:Y^\star\to X^\star$ открыто}.
$$

\btm\label{TH:zamk-biektsija-v-kopolnoe-prostrancsvo}
Пусть $\ph:X\to Y$ -- замкнутое биективное линейное непрерывное отображение стереотипных пространств, причем $Y$ кополно. Тогда $\ph:X\to Y$ -- изоморфизм стереотипных пространств.
\etm
\bpr
Пространство $X$ можно понимать как новую, более тонкую стереотипную топологизацию стереотипного пространства $Y$, сохраняющую систему вполне ограниченных множеств и топологию на каждом вполне ограниченном множестве. Если перейти к сопряженным пространствам, то $Y^\star$ будет подпространством в $X^\star$ с топологией, индуцированной из $X^\star$, причем $Y^\star$ будет плотно в $X^\star$ (потому что сопряженное отображение $\ph:X\to Y$ -- инъекция), и одновременно $Y^\star$ будет полно (потому что $Y$ кополно). Вместе это означает, что $Y^\star$ и $X^\star$ совпадают как локально выпуклые пространства, и значит, то же самое верно для $Y$ и $X$.
\epr

\subsection{Лемма о факторе тензорных произведений}

Следующее утверждение понадобится нам в лемме \ref{LM:J^0_C(M)C(M)-circledast-A-cong-J^0_C(M)C(M,A)}:

\blm\label{PROP:[(X-circledast-Z)/(Y-circledast-Z)]^triangledown-cong-[(X-odot-Z)/(Y-odot-Z)]^triangledown}
Пусть $X$, $Y$ и $Z$ -- стереотипные пространства, причем $Y$ -- дополняемое замкнутое подпространство в $X$. Тогда
\beq\label{[X/Y]^triangledown-circledast-Z}
[(X\circledast Z)/(Y\circledast Z)]^\triangledown\cong[X/Y]^\triangledown\circledast Z
\eeq
и
\beq\label{[X/Y]^triangledown-odot-Z}
[(X\odot Z)/(Y\odot Z)]^\triangledown\cong[X/Y]^\triangledown\odot Z
\eeq
В частности, если $Y$ -- замкнутое подпространство конечной коразмерности в $X$, то
\beq\label{[(X-circledast-Z)/(Y-circledast-Z)]^triangledown-cong-[(X-odot-Z)/(Y-odot-Z)]^triangledown}
[(X\circledast Z)/(Y\circledast Z)]^\triangledown\cong [X/Y]\circledast Z\cong [X/Y]\odot Z\cong
[(X\odot Z)/(Y\odot Z)]^\triangledown
\eeq
\elm
\bpr
Пусть $P$ -- дополнение к $Y$ в $X$:
$$
X=Y\oplus P.
$$
Тогда $X/Y\cong P$ и фактор-отображение $\ph:X\to P$ можно представлять себе, как проекцию на вторую компоненту. Тензорное произведение при этом принимет вид
$$
X\circledast Z=(Y\oplus P)\circledast Z=\cite[(4.151)]{Akbarov-De-Gruyter-I}=(Y\circledast Z)\oplus (P\circledast Z),
$$
и $\ph\circledast\id_Z:X\circledast Z\to P\circledast Z$ тоже становится проекцией на вторую компоненту. Мы получаем
\eqref{[X/Y]^triangledown-circledast-Z}:
$$
[(X\circledast Z)/(Y\circledast Z)]^\triangledown=\Coim(\ph\circledast\id_Z)=P\circledast Z=[X/Y]^\triangledown\circledast Z.
$$
Аналогично доказывается \eqref{[X/Y]^triangledown-odot-Z}. После этого \eqref{[(X-circledast-Z)/(Y-circledast-Z)]^triangledown-cong-[(X-odot-Z)/(Y-odot-Z)]^triangledown} становится следствием того факта, что тензорные произведения $\circledast$ и $\odot$ совпадают, если один из множителей конечномерен:
$$
P\circledast Z\cong P\odot Z,\qquad \dim P<\infty.
$$
\epr

\subsection{Лемма об эпиморфизме}

Напомним, что в работе \cite{Ak16} было введено понятие узлового разложения морфизма $\ph:X\to Y$. Так называется представление $\ph$ в виде композиции
$$
\ph=\sigma\circ\beta\circ\pi,
$$
в которой $\sigma$ --- строгий мономорфизм, $\beta$ --- биморфизм, а $\pi$ --- строгий эпиморфизм. Если такое представление существует, то оно единственно с точностью до изоморфизма компонент, поэтому этим компонентам можно присвоить обозначения:
$$
\sigma=\im_\infty\ph,\quad \beta=\red_\infty\ph,\quad \pi=\coim_\infty\ph.
$$
При этом область определения морфизма $\sigma$ обозначается $\Im_\infty\ph$, а область значений морфизма $\pi$ --- $\Coim_\infty\ph$. Таким образом, морфизм $\ph$ распадается в композицию
\beq\label{DEF:oboznacheniya-dlya-uzlov-razlozh}
\begin{diagram}
\node{X}\arrow{s,l}{\coim_\infty\ph}\arrow{e,t}{\ph}\node{Y} \\
\node{\Coim_\infty\ph}\arrow{e,t}{\red_\infty\ph}\node{\Im_\infty\ph}\arrow{n,r}{\im_\infty\ph}
\end{diagram}
\eeq
Известно \cite[Theorem 4.100]{Ak16}, что в категории $\Ste$ стереотипных пространств любой морфизм обладает узловым разложением.

\blm\label{LM:epimorphism}
Пусть $\e:X\to Y$ и $\ph:Y\to Z$ -- морфизмы стереотипных пространств, и $\psi=\ph\circ\e$ -- их композиция. Тогда если $\e$ -- эпиморфизм, то существует единственный морфизм $\ph_\infty:Y\to\Im_\infty\psi$, замыкающий диаграмму
\beq\label{diagr:epimorphism}
 \xymatrix 
 {
 X\ar[dr]_{\psi_\infty}\ar[rr]^{\e}\ar@/_5ex/[ddr]_{\psi} & & Y\ar@{-->}[dl]^{\ph_\infty}
 \ar@/^5ex/[ddl]^{\ph} \\
  & \Im_\infty\psi \ar[d]_{\im_\infty\psi} & \\
   & Z  &
 }
\eeq
где $\psi_\infty=\red_\infty\psi\circ\coim_\infty\psi$.
\elm
\bpr
Здесь нужно воспользоваться формулой \cite[(4.85)]{Ak16}, по которой узловой образ отображения $\psi$ совпадает с оболочкой $\Env^Z\psi(X)$ (в смысле \cite[$\S$ 2, (e)]{Ak16}) образа $\psi(X)$ отображения $\psi$:
$$
\im_\infty\psi=\Env^Z\psi(X).
$$
Морфизм $\ph_\infty$ строится трансфинитной индукцией.

Опишем подробно ее нулевой шаг. Зафиксируем точку $y\in Y$. Поскольку $\e$ -- эпиморфизм, найдется направленность $\{x_i\}\subseteq X$ такая, что
$$
\e(x_i)\overset{Y}{\underset{i\to\infty}{\longrightarrow}}y.
$$
Отсюда
$$
\psi(x_i)=\ph(\e(x_i))\overset{Z}{\underset{i\to\infty}{\longrightarrow}}\ph(y).
$$
и поэтому $\ph(y)\in\overline{\psi(X)}^Z$. Это верно для любого $y\in Y$, значит,
$$
\ph(Y)\subseteq\overline{\psi(X)}^Z
$$
Теперь можно поглядеть на $\ph$, как на непрерывное отображение стереотипного пространства $Y$ в (необязательно стереотипное, но псевдополное) пространство $\overline{\psi(X)}^Z$ (с топологией, индуцированной из $Z$):
$$
\ph: Y\to \overline{\psi(X)}^Z.
$$
Поскольку пространство $Y$ псевдонасыщено, из \cite[(1.26)]{Ak03} мы получаем, что если считать $\ph$ отображением, действующим из $Y$ в псевдонасыщение $\Big(\overline{\psi(X)}^Z\Big)^\vartriangle$ пространства $\overline{\psi(X)}^Z$, то $\ph$ также будет непрерывно:
$$
\ph: Y\to\Big(\overline{\psi(X)}^Z\Big)^\vartriangle.
$$
Обозначим $E_0=\Big(\overline{\psi(X)}^Z\Big)^\vartriangle$, и пусть $\ph_0$ и $\psi_0$ -- отображения $\ph$ и $\psi$, рассмотренные как действующие со значениями в $E_0$. Тогда мы получим диаграмму, представляющую собой нулевое приближение к \eqref{diagr:epimorphism}:
$$
 \xymatrix 
 {
 X\ar[dr]_{\psi_0}\ar[rr]^{\e}\ar@/_5ex/[ddr]_{\psi} & & Y\ar@{-->}[dl]^{\ph_0}
 \ar@/^5ex/[ddl]^{\ph} \\
  & E_0 \ar[d] & \\
   & Z  &
 }
$$

Далее для всякого ординала $k$ мы определяем $E_k$, $\ph_k$ и $\psi_k$ следующим образом:  \bit{

\item[---] если $k$ -- изолированный ординал, то есть $k=j+1$ для некоторого ординала $j$, то мы применяем тот же самый прием, что и при $k=0$: полагаем $E_k=\Big(\overline{\psi(X)}^{E_j}\Big)^\vartriangle$ и получаем диаграмму
\beq\label{diagr:epimorphism-1}
 \xymatrix 
 {
 X\ar[dr]_{\psi_{j+1}}\ar[rr]^{\e}\ar@/_5ex/[ddr]_{\psi_j}\ar@/_10ex/[dddr]_{\psi} & & Y\ar@{-->}[dl]^{\ph_{j+1}}
 \ar@/^5ex/[ddl]^{\ph_j} \ar@/^10ex/[dddl]^{\ph} \\
  & E_{j+1} \ar[d] & \\
  & E_j \ar[d] & \\
   & Z  &
 }
\eeq

\item[---] если $k$ -- предельный ординал, то есть не существует такого $j$, что $k=j+1$, то мы полагаем $E_k=\lim_{k\gets j} E_j$ (проективный предел в категории стеретипных пространств), и получаем для этого случая диаграмму
$$
 \xymatrix 
 {
 X\ar[dr]_{\lim\limits_{k\gets j}\psi_j}\ar[rr]^{\e}\ar@/_8ex/[ddr]_{\psi} & & Y\ar@{-->}[dl]^{\lim\limits_{k\gets j}\ph_j}
 \ar@/^8ex/[ddl]^{\ph} \\
  & \lim\limits_{k\gets j} E_j \ar[d] & \\
   & Z  &
 }
$$

}\eit

В результате мы получим трансфинитную последовательность стереотипных пространств $E_k$ и диаграмм \eqref{diagr:epimorphism-1}. В силу \cite[(4.58)]{Ak16}, эта последовательность стабилизируется, и ее пределом будет как раз диаграмма \eqref{diagr:epimorphism}.

\epr

\bcor\label{COR:Im_infty-ph-circ-e=Im_infty-ph}
Пусть $\e:X\to Y$ и $\ph:Y\to Z$ -- морфизмы стереотипных пространств, причем $\e$ -- эпиморфизм. Тогда
\beq\label{Im_infty-ph-circ-e=Im_infty-ph}
\Im_\infty(\ph\circ\e)=\Im_\infty\ph
\eeq
\ecor
\bpr
С одной стороны, справедлива очевидная импликация
$$
(\ph\circ\e)(X)=\ph(\e(X))\subseteq \ph(X)\quad\Longrightarrow\quad
\Im_\infty(\ph\circ\e)\subseteq\Im_\infty\ph.
$$
С другой стороны, из леммы \ref{LM:epimorphism} следует вложение
$$
\ph(Y)\subseteq \Im_\infty(\ph\circ\e),
$$
которое в свою очередь влечет вложение
$$
\Im_\infty\ph(Y)\subseteq \Im_\infty(\ph\circ\e).
$$
\epr

\subsection{Строгие пределы стереотипных пространств}

\bit{

\item[$\bullet$] Пусть $\{X_i;\ i\in\N\}$ --- последовательность локально выпуклых пространств, и пусть нам  задана последовательность их морфизмов (то есть линейных непрерывных отображений)
$$
\iota_i:X_i\to X_{i+1}
$$
со следующими свойствами:
\bit{

\item[1)] все $\iota_i$ --- вложения (то есть инъективны и такие, что топология $X_i$  наследуется из $X_{i+1}$),  и

\item[2)] все $\iota_i$ --- слабо замкнуты (то есть образ $\iota_i(X_i)$ замкнут в $X_{i+1}$).

}\eit
Тогда инъективный предел такой последовательности в категории $\LCS$
$$
\LCS\text{-}\injlim_{i\in\N}X_i
$$
--- мы будем называть ее {\it строгим инъективным пределом}.

}\eit

\btm\label{TH:strogij-injlim-ogr}
В строгом инъективном пределе последовательности локально выпуклых пространств
$$
\LCS\text{-}\injlim_{i\in\N}X_i,
$$
множество $B$ ограничено тогда и только тогда, когда оно лежит в некотором пространстве $X_i$ и ограничено в нем.
\etm
\bpr
Этот классический результат доказывается, например,  в \cite[VII, 1.4]{Robertson-Robertson}.

\epr

\btm\label{TH:strogij-injlim}
Строгий инъективный предел последовательности стереотипных пространств совпадает с их инъективным пределом в категории $\Ste$ 
$$
\Ste\text{-}\injlim_{i\in\N}X_i=\LCS\text{-}\injlim_{i\in\N}X_i,
$$
а если вдобавок все пространства $X_i$ полны, то этот предел --- тоже полное пространство.
\etm
\bpr
1. По теореме \ref{TH:strogij-injlim-ogr}, всякое ограниченное множество в $\LCS\text{-}\injlim_{i\in\N}X_i$ содержится в некотором $X_i$ и ограничено в нем. Отсюда следует, что пространство  $\LCS\text{-}\injlim_{i\in\N}X_i$ псевдополно, и поэтому
$$
\Ste\text{-}\injlim_{i\to\infty}X_i^\star=(\LCS\text{-}\injlim_{i\to\infty}X_i^\star)^\triangledown= \LCS\text{-}\injlim_{i\to\infty}X_i^\star.
$$

2. Если вдобавок все пространства $X_i$ полны, то по другому свойству индуктивных пределов \cite[VII, 1.3]{Robertson-Robertson}, пространство  $\LCS\text{-}\injlim_{i\in\N}X_i$ полно, и поэтому $\Ste\text{-}\injlim_{i\in\N}X_i$ тоже полно.
\epr

\bit{

\item[$\bullet$] Пусть $\{X_i;\ i\in\N\}$ --- последовательность локально выпуклых пространств и пусть нам задана последовательность их морфизмов (то есть линейных непрерывных отображений)
$$
\varkappa_i:X_{i+1}\to X_i
$$
со следующими свойствами:
\bit{

\item[1)] все $\varkappa_i$ --- наложения (то есть каждое вполне ограниченное множество $S\subseteq X_i$ является образом некоторого вполне ограниченного множества $T\subseteq X_{i+1}$), и 
    
\item[2)] все $\varkappa_i$ слабо открыты (то есть всякий функционал $f\in X_{i+1}^\star$, обнуляющийся на ядре $\Ker\varkappa_i$, продолжается до функционала $g\in X_i^\star$), 
}\eit    
Тогда проективный предел такой последовательности в категории $\LCS$
$$
\LCS\text{-}\projlim_{i\in\N}X_i
$$
--- мы будем называть ее {\it строгим проективным пределом}.

}\eit

\btm\label{TH:strogij-projlim}
Строгий проективный предел последовательности стереотипных пространств совпадает с их проективным пределом в категории $\Ste$ 
$$
\Ste\text{-}\projlim_{i\in\N}X_i=\LCS\text{-}\projlim_{i\in\N}X_i,
$$
а если вдобавок все пространства $X_i$ насыщены, то этот предел --- тоже насыщенное пространство.
\etm
\bpr
1. Морфизмы 
$$
\varkappa_i:X_{i+1}\to X_i
$$ 
являются слабо открытыми наложениями, поэтому их сопряженные морфизмы 
$$
(\varkappa_i)^\star:X_i^\star\to X_{i+1}^\star
$$ 
являются слабо замкнутыми вложениями. По теореме \ref{TH:strogij-injlim}, инъективный предел такой системы в $\Ste$ совпадает с ее инъективным пределом в $\LCS$:
$$
\Ste\text{-}\injlim_{i\to\infty}X_i^\star= \LCS\text{-}\injlim_{i\to\infty}X_i^\star
$$
Отсюда мы можем сделать вывод, что проективный предел в $\Ste$ пространств $X_i$, который по двойственности, совпадает с сопряженным пространством к $\Ste\text{-}\injlim_{i\to\infty}X_i^\star$, представляет собой пространство, сопряженное к $\LCS\text{-}\injlim_{i\to\infty}X_i^\star$:
$$
\Ste\text{-}\projlim_{i\to\infty}X_i=(\Ste\text{-}\injlim_{i\to\infty}X_i^\star)^\star=
(\LCS\text{-}\injlim_{i\to\infty}X_i^\star)^\star
$$
С другой стороны, $\Ste\text{-}\projlim_{i\to\infty}X_i$ есть псевдонасыщение локально выпуклого проективного предела, и мы получаем цепочку
$$
(\LCS\text{-}\projlim_{i\to\infty}X_i)^\vartriangle=\Ste\text{-}\projlim_{i\to\infty}X_i= (\LCS\text{-}\injlim_{i\to\infty}X_i^\star)^\star
$$
Это можно понимать так: $\LCS\text{-}\projlim_{i\to\infty}X_i$ есть, как множество, пространство функционалов на $\LCS\text{-}\injlim_{i\to\infty}X_i^\star$, а топология у него, --- обозначим ее $\sigma$, --- если ее псевдонасытить, превращается в топологию, --- обозначим ее $\tau$, --- равномерной сходимости на вполне ограниченных множествах в $\LCS\text{-}\injlim_{i\to\infty}X_i^\star$. 
$$
\sigma^\vartriangle=\tau.
$$
Но по свойству строгого индуктивного предела \cite[VII, 1.3]{Robertson-Robertson}, вполне ограниченные множества из $\LCS\text{-}\injlim_{i\to\infty}X_i^\star$ -- это те, которые содержатся в каком-то $X_i^\star$. Значит топология $\tau$ --- это топология равномерной сходимости на вполне ограниченных множествах из $\LCS\text{-}\injlim_{i\to\infty}X_i^\star$,  которые содержатся в каком-то $X_i^\star$. А это при таком отождествлении --- в точности топология пространства  $\LCS\text{-}\projlim_{i\to\infty}X_i$, то есть $\sigma$:
$$
\tau=\sigma.
$$
Мы получаем, что $\sigma^\vartriangle=\sigma$, то есть топологии у  $\LCS\text{-}\projlim_{i\to\infty}X_i$ и у $(\LCS\text{-}\projlim_{i\to\infty}X_i)^\vartriangle$ совпадают:
$$
\LCS\text{-}\projlim_{i\to\infty}X_i=(\LCS\text{-}\projlim_{i\to\infty}X_i)^\vartriangle=
\Ste\text{-}\projlim_{i\to\infty}X_i.
$$

2. Пусть далее все пространства $X_i$ насыщены. Тогда все пространства $X_i^\star$ полны и по теореме  \ref{TH:strogij-injlim}, пространство $\Ste\text{-}\injlim_{i\to\infty}X_i^\star$ полно. Значит, сопряженное ему пространство  
$\Ste\text{-}\projlim_{i\to\infty}X_i$ насыщено.
\epr

\subsection{Непосредственные подпространства и фактор-пространства}

\paragraph{Непосредственные подпространства.}

Напомним определение непосредственного подпространства в стереотипном пространстве \cite{Akbarov-De-Gruyter-I}.

\bit{
\item[$\bullet$]
Пусть $Y$ --- подмножество в стереотипном пространстве $X$, наделенное структурой стереотипного пространства так, чтобы теоретико-множественное включение $Y\subseteq X$ было морфизмом стереотипных пространств (то есть линейным и непрерывным отображением). Тогда стереотипное пространство $Y$ называется {\it подпространством}\index{подпространство!стереотипного пространства} стереотипного пространства $X$, а теоретико-множественное включение $\sigma:Y\subseteq X$ --- его {\it представляющим мономорфизмом}\index{мономорфизм!представляющий}. Запись\index{$\subarr$}
$$
Y\subarr X
$$
или запись
$$
X\suparr Y
$$
будет означать, что $Y$ является подпространством стереотипного пространства $X$.

 \item[$\bullet$]
Пусть нам дана последовательность из двух вложенных друг в друга подпространств
$$
Z\subarr Y\subarr X,
$$
причем включение $Z\subarr Y$ является биморфизмом стереотипных пространств, то есть, проще говоря, $Z$, помимо того, что линейно и непрерывно отображается в $Y$, еще и плотно в нем (относительно топологии $Y$):
$$
\overline{Z}^Y=Y.
$$
Тогда мы будем говорить, что подпространство $Y$ является {\it посредником}\index{посредник!для подпространства} для подпространства $Z$ в пространстве $X$.

 \item[$\bullet$]\label{DEF:neposr-podprostr}
Подпространство $Z$ стереотипного пространства  $X$  мы называем {\it непосредственным подпространством}\index{подпространство!непосредственное} в $X$, если у него нет посредников, то есть
 для любого его посредника $Y$ в $X$ соответствующее включение $Z\subarr Y$ является изоморфизмом. В этом случае мы пользуемся записью $Z\osubarr X$:\index{$\osubarr$}
 $$
 Z\osubarr X\qquad\Leftrightarrow\qquad \forall Y\quad \bigg( \Big(Z\subarr Y\subarr X \quad\&\quad \overline{Z}^Y=Y\Big)\quad\Rightarrow\quad Z=Y\bigg).
 $$
Понятно, что это просто означает, что представляющий мономорфизм $\upsilon:Y\to X$ должен быть непосредственным мономорфизмом\footnote{Понятие непосредственного мономорфизма было определено на с.\pageref{DEF:neposr-mono}.} в категории $\Ste$.

}\eit

\btm\label{TH:Y=Env^X Y<=>neposr-podpr}
Для подпространства $Y\subarr X$ стереотипного пространства $X$ следующие условия эквивалентны:
\bit{

\item[(i)] $Y$ --- непосредственное подпространство  в $X$:
$$
Y\osubarr X;
$$

\item[(ii)] $Y$ совпадает со своей оболочкой\footnote{Понятие оболочки $\Env^X M$ множества $M$ в стереотипном пространстве $X$ определено в \cite[p.534]{Akbarov-De-Gruyter-I}.} в $X$:
$$
Y=\Env^X Y
$$
(то есть естественное включение $Y\to\Env^X Y$ является изоморфизмом стереотипных пространств). 

}\eit
\etm

\btm\label{TH:Z-subarr-X,Y-osubarr-X,Z-subseteq-Y}
Пусть $Z$ и $Y$ --- подпространства стереотипного пространства $X$, причем $Y$ --- непосредственное подпространство:
$$
Z\subarr X, \qquad Y\osubarr X.
$$
Тогда если $Z$ --- подмножество в $Y$
$$
Z\subseteq Y,
$$
то $Z$ --- подпространство в $Y$:
$$
Z\subarr Y
$$
(то есть включение $Z\subseteq Y$ является непрерывным отображением).
\etm

\paragraph{Непосредственные фактор-пространства.}

 \bit{
\item[$\bullet$] Пусть нам дано стереотипное пространство $X$, и пусть
 \bit{

 \item[1)] в $X$, как в локально выпуклом пространстве, задано некое замкнутое подпространство $E$,

 \item[2)] на фактор-пространстве $X/E$ задана некая топология $\tau$, мажорируемая естественной  топологией фактор-пространства, и превращающее  $X/E$ в локально выпуклое пространство,

 \item[3)] в пополнении $(X/E)^\blacktriangledown$ локально выпуклого пространства $X/E$ относительно топологии $\tau$ задано некое подпространство $Y$, содержащее $X/E$ и являющееся стереотипным пространством относительно топологии, индуцированной из $(X/E)^\blacktriangledown$.
 }\eit
Тогда стереотипное пространство $Y$ мы называем {\it фактор-пространством стереотипного пространства}\label{DEF:faktor-prostr-Ste}\index{фактор-пространство} $X$, а композиция $\upsilon=\sigma\circ\pi:X\to Y$ фактор-отображения $\pi:X\to X/E$ и естественного включения $\sigma:X/E\to Y$ называется {\it представляющим эпиморфизмом} фактор-пространства $Y$. Запись\index{$\quarr$}
$$
Y\quarr X
$$
или запись
$$
X\qparr Y
$$
будет означать, что $Y$ является фактор-пространством стереотипного пространства $X$ (разумеется, $Y$ не определется пространствами $X$ и $E$).

 \item[$\bullet$]
Пусть $Y$ и $Z$ --- два фактор-пространства для $X$, связанных неравенством
$$
Z\le Y,
$$
причем эпиморфизм $\varkappa:Z\gets Y$ из диаграммы \eqref{DEF:le-dlya-faktor-prostranstv} является мономорфизмом (и значит, биморфизмом) стереотипных пространств. Тогда мы будем говорить, что фактор-пространство $Y$ является {\it посредником}\index{посредник!для фактор-пространства} для фактор-пространства $Z$ пространства $X$. Можно заметить, что в этом случае $Y$ просто является подмножеством в $Z$, поэтому мы будем писать $Z\supseteq Y$.

 \item[$\bullet$] Пусть $Y\quarr X$ и $Z\quarr X$. Условимся говорить, что фактор-пространство $Y$ {\it подчиняет} фактор-пространство $Z$ и изображать это записью $Z\le Y$, если существует морфизм $\varkappa:Y\to Z$ такой, что
 \beq\label{DEF:le-dlya-faktor-prostranstv}
\xymatrix @R=1pc @C=2pc
{
Y\ar[dd]_{\varkappa} &   \\
  & X\ar[lu]_{\upsilon_Y}\ar[ld]^{\upsilon_Z} \\
Z &
}
\eeq
(здесь $\upsilon_Y$ и $\upsilon_Z$ --- представляющие эпиморфизмы для $Y$ и $Z$). Морфизм $\varkappa$, если он существует, будет, во-первых, единственным и, во-вторых, эпиморфизмом.

 \item[$\bullet$]
Пусть $Y$ и $Z$ --- два фактор-пространства для $X$, связанных неравенством
$$
Z\le Y,
$$
причем эпиморфизм $\varkappa:Z\gets Y$ из диаграммы \eqref{DEF:le-dlya-faktor-prostranstv} является мономорфизмом (и значит, биморфизмом) стереотипных пространств. Тогда мы будем говорить, что фактор-пространство $Y$ является {\it посредником}\index{посредник!для фактор-пространства} для фактор-пространства $Z$ пространства $X$. Можно заметить, что в этом случае $Y$ просто является подмножеством в $Z$, поэтому мы будем писать $Z\supseteq Y$.

 \item[$\bullet$]\label{DEF:neposr-faktor-prost}
Фактор-пространство $Z$ стереотипного пространства  $X$  мы называем {\it непосредственным фактор-пространством}\index{фактор-пространство!непосредственное} в $X$, если у него нет посредников, то есть для любого его посредника $Y$ в $X$ соответствующий эпиморфизм $Z\quarr Y$ является изоморфизмом. В этом случае мы пользуемся записью $Z\oquarr X$:\index{$\oquarr$}
 $$
 Z\oquarr X\qquad\Leftrightarrow\qquad \forall Y\quad \bigg( \Big(Z\le Y\quad \&\quad Y\quarr X \quad\&\quad Z\supseteq Y\Big)\quad\Rightarrow\quad Z=Y\bigg).
 $$
Понятно, что это просто означает, что представляющий эпиморфизм $\upsilon:X\to Y$ должен быть непосредственным эпиморфизмом\footnote{Понятие непосредственного эпиморфизма было определено на с.\pageref{DEF:neposr-epi}.} в категории $\Ste$.

}\eit

\btm\label{TH:Y=Ref^X<=>neposr-faktor}
Для фактор-пространства $Y\quarr X$ стереотипного пространства $X$ с представляющим эпиморфизмом $\upsilon:X\to Y$ следующие условия эквивалентны:
\bit{

\item[(i)] $Y$ --- непосредственное фактор-пространство пространства $X$:
$$
Y\oquarr X;
$$

\item[(ii)] $Y$ совпадает с детализацией\footnote{Понятие детализации $\Rf^X F$ множества $F$ функционалов на стереотипном пространстве $X$ определено в \cite[p.554]{Akbarov-De-Gruyter-I}.} множества функционалов $Y^\star\circ\upsilon=\{f\circ\upsilon;\ f\in Y^\star\}$ на $X$:
$$
Y=\Rf^X (Y^\star\circ\upsilon)
$$
(то есть естественное отображение $\Rf^X (Y^\star\circ\upsilon)\to Y$ является изоморфизмом стереотипных пространств). 

}\eit
\etm

\btm\label{TH:Z-quarr-X,Y-oquarr-X,Ker}
Пусть $Z$ и $Y$ --- фактор-пространства стереотипного пространства $X$ с представляющими эпиморфизмами $\upsilon:X\to Y$ и $\zeta:X\to Z$, причем $Y$ --- непосредственное фактор-пространство:
$$
Z\quarr X, \qquad Y\oquarr X.
$$
Тогда если 
$$
\Ker\zeta\supseteq\Ker\upsilon,
$$
то $Z$ --- фактор-пространство пространства $Y$:
$$
Z\quarr Y
$$
(то есть естественное отображение $Y\to Z$ является эпиморфизмом в категории $\Ste$).
\etm

\paragraph{Непосредственные подпространства и фактор-пространства в пространствах операторов и тензорных произведениях}

\btm\label{TH:oslash-sohranyaet-oquarr-i-osubarr}
Пусть $\pi:X\to A$ --- непосредственный эпиморфизм стереотипных пространств, а $\sigma:B\to Y$  --- непосредственный мономорфизм  стереотипных пространств. Тогда $(\sigma\oslash\pi): B\oslash A\to Y\oslash X$  --- непосредственный мономорфизм стереотипных пространств.
\etm

\brem\label{REM:oslash-sohranyaet-oquarr-i-osubarr}
В эквивалентной формулировке это звучит так. {\it Пусть $A$ --- непосредственное фактор-пространство стереотипного пространства $X$, а $B$  --- непосредственное подпространство в стереотипном пространстве $Y$:
$$
A\oquarr X,\qquad B\osubarr Y.
$$
Тогда $B\oslash A$  --- непосредственное подпространство в стереотипном пространстве $Y\oslash X$:
$$
B\oslash A\osubarr Y\oslash X.
$$
}
\erem

Доказательству этого утверждения мы предпошлем две леммы.

\blm\label{LM:sigma-oslash-1_X-in-SMono}
Если $\sigma:B\to Y$ --- непосредственный мономорфизм стереотипных пространств, то для любого стереотипного пространства $X$ морфизм
$\sigma\oslash 1_X: B\oslash X\to Y\oslash X$ --- тоже непосредственный мономорфизм стереотипных пространств.
\elm
\bpr
1. Сначала рассмотрим случай, когда $B$ представляет собой замкнутое непосредственное подпространство в $Y$, то есть имеется некое замкнутое подпространство $H$ в $Y$, как в локально выпуклом пространстве, причем топология $H$ индуцируется из $Y$, и если обозначить вложение $H$ в $Y$ буквой $\chi$, 
$$
\chi:H\to Y,
$$
то псевдонасыщение $H$ будет $B$, а псевдонасыщение $\chi$ будет $\sigma$:
$$
\sigma=\chi^\vartriangle:B=H^\vartriangle\to Y^\vartriangle=Y.
$$
Рассмотрим морфизм пространств операторов:
$$
(\chi:1_X):(H:X)\to (Y:X).
$$ 
то есть отображение, действующее по формуле
$$
(\chi:1_X)(\ph)=\chi\circ\ph, \qquad \ph\in H:X.
$$
Из того, что $H$ --- замкнутое подпространство в $Y$ с индуцированной из $Y$ топологией, следует, что $(\chi:1_X):(H:X)\to (Y:X)$ --- вложение локально выпуклых пространств (то есть оно непрерывно и открыто). Отсюда в свою очередь следует \cite[Example 4.3.11]{Akbarov-De-Gruyter-I}, что псевдонасыщение этого морфизма 
$$
\underbrace{(\chi:1_X)^\vartriangle}_{\scriptsize \begin{matrix}\|  \\ \cite[(3.203)]{Akbarov-De-Gruyter-I}\\ \| \\ (\chi^\vartriangle:1_X^\triangledown)^\vartriangle \\ \| \\ (\sigma:1_X)^\vartriangle \\ \| \\ \sigma\oslash 1_X \end{matrix}}: \underbrace{(H:X)^\vartriangle}_{\scriptsize \begin{matrix}\|  \\ \cite[(3.203)]{Akbarov-De-Gruyter-I}\\ \| \\ (H^\vartriangle:X^\triangledown)^\vartriangle \\ \| \\ (B:X)^\vartriangle \\ \| \\ B\oslash X \end{matrix}}\to \underbrace{(Y:X)^\vartriangle}_{\scriptsize \begin{matrix}\|  \\ \cite[(4.106)]{Akbarov-De-Gruyter-I}\\ \| \\ Y\oslash X \end{matrix}}.
$$
--- является непосредственным мономорфизмом.
 
2. Теперь рассмотрим общий случай, когда $B$ --- произвольное непосредственное подпространство в $Y$. Тогда оно  совпадает со своей оболочкой в $Y$ \cite[p.534]{Akbarov-De-Gruyter-I}:
$$
B=\Env^Y B.
$$
При этом оболочка $\Env^Y B$ строится как проективный предел подпространств
$$
E_0=(\overline{B}^Y)^\vartriangle, \quad E_j=\begin{cases}(\overline{B}^{E_i})^\vartriangle, & j=i+1 \\ \Ste\text{-}\projlim_{i\to j} E_i, & \text{$j$ --- предельный ординал}\end{cases}
$$  
Проследим, как это отражается на отображениях $E_i\oslash X\to Y\oslash X$.

а) Отображение $E_0\oslash X=(\overline{B}^Y)^\vartriangle\oslash X\to Y\oslash X$ является непосредственным мономорфизмом, потому что это тип отображения рассмотренный в пункте 1. 

б) Предположим, мы доказали, что для данного ординала $i$ отображение  
$$
E_i\oslash X\to Y\oslash X
$$
является непосредственным мономорфизмом. Заметим тогда, что отображение $E_{i+1}\oslash X=(\overline{B}^{E_{i+1}})^\vartriangle\oslash X\to E_i\oslash X$ является непосредственным мономорфизмом, потому что это тип отображения рассмотренный в пункте 1. Отсюда следует, что это отображение является строгим мономорфизмом \cite[Theorem 4.3.55]{Akbarov-De-Gruyter-I}. Поэтому его композиция с отображением $E_i\oslash X\to Y\oslash X$ (которое тоже, будучи непосредственным мономорфизмом является строгим мономорфизмом), 
$$
E_{i+1}\oslash X\to E_i\oslash X\to Y\oslash X,
$$
--- тоже является строгим мономорфизмом \cite[Theorem 1.3.16]{Akbarov-De-Gruyter-I}. Значит, это отображение --- непосредственный мономорфизм.

в) Если мы доказали, что для данного предельного ординала $j$ система отображений
$$
E_i\oslash X\to Y\oslash X,\qquad i<j
$$
состоит из непосредственных  (=строгих) мономорфизмов, то их проективный предел  
$$
E_j\oslash X=(\Ste\text{-}\projlim_{i\to j} E_i)\oslash X=\Ste\text{-}\projlim_{i\to j} (E_i\oslash X)\to Y\oslash X,
$$
--- тоже является непосредственным (=строгим) мономорфизмом \cite[Theorem 1.4.29]{Akbarov-De-Gruyter-I}.

\epr

\blm\label{LM:1_Y-oslash-pi-in-SMono}
Если $\pi:X\to A$ --- непосредственный эпиморфизм стереотипных пространств, то для любого стереотипного пространства $Y$ морфизм
$1_Y\oslash\pi: Y\oslash A\to Y\oslash X$ --- тоже непосредственный мономорфизм стереотипных пространств.
\elm
\bpr
Змаетим, что сопряженное отображение $\pi^\star:A^\star\to X^\star$ является непоссредственным мономорфизмом. Поэтому по лемме \ref{LM:sigma-oslash-1_X-in-SMono} морфизм
$$
\pi^\star\oslash 1_{Y^\star}:A^\star\oslash Y^\star\to X^\star\oslash Y^\star
$$
является непосредственным мономорфизмом. Значит (здесь мы снова применяем \cite[Theorem 4.3.55]{Akbarov-De-Gruyter-I}) он является строгим мономорфизмом. Это значит, что в любой диаграмме
$$
 \xymatrix  @R=2.5pc @C=8pc
 {
 P\ar[d]_{\alpha}\ar[r]^{\e} & Q\ar[d]_{\beta}\\
 A^\star\oslash Y^\star \ar[r]^{\pi^\star\oslash 1_{Y^\star}} & X^\star\oslash Y^\star
 }.
$$
у которой $\e$ --- эпиморфизм, существует диагональ, то есть морфизм $\delta$> для которого будет коммутативна диаграмма
$$
 \xymatrix  @R=2.5pc @C=8pc
 {
 P\ar[d]_{\alpha}\ar[r]^{\e} & Q\ar[d]_{\beta}\ar@{-->}[dl]_{\delta}\\
 A^\star\oslash Y^\star \ar[r]^{\pi^\star\oslash 1_{Y^\star}} & X^\star\oslash Y^\star
 }.
$$
Если теперь достроить ее до диаграммы
$$
 \xymatrix  @R=2.5pc @C=8pc
 {
 P\ar[d]_{\alpha}\ar[r]^{\e} & Q\ar[d]_{\beta}\ar@{-->}[dl]_{\delta}\\
 A^\star\oslash Y^\star \ar[r]^{\pi^\star\oslash 1_{Y^\star}} & X^\star\oslash Y^\star \\
 Y\oslash A \ar[r]^{1_Y\oslash \pi}\ar@{=>}[u]^{\ph\mapsto \ph^\star} & Y\oslash X \ar@{=>}[u]_{\psi\mapsto \psi^\star} 
 }
$$
(двойные стрелки $\Rightarrow$ одначают изоморфизмы), то станет видно, что любая диаграмма 
$$
 \xymatrix  @R=2.5pc @C=8pc
 {
 P\ar[d]_{\alpha}\ar[r]^{\e} & Q\ar[d]_{\beta}\\
 Y\oslash A \ar[r]^{1_Y\oslash \pi} & Y\oslash X  
 }.
$$
(в которой $\e$ --- по-прежнему, эпиморфизм) дополняется до диаграммы
$$
 \xymatrix  @R=2.5pc @C=8pc
 {
 P\ar[d]_{\alpha}\ar[r]^{\e} & Q\ar[d]_{\beta}\ar@{-->}[dl]_{\delta}\\
 Y\oslash A \ar[r]^{1_Y\oslash \pi} & Y\oslash X  
 }.
$$
Это значит, что $1_Y\oslash \pi$ --- строгий мономорфизм (и, как следсвие, непосредственный мономорфизм).
\epr

\bpr[Доказательство теоремы \ref{TH:oslash-sohranyaet-oquarr-i-osubarr}.]
Пусть $\pi:X\to A$ --- непосредственный эпиморфизм стереотипных пространств, а $\sigma:B\to Y$  --- непосредственный мономорфизм  стереотипных пространств. Тогда по лемме \ref{LM:sigma-oslash-1_X-in-SMono} морфизм 
$\sigma\oslash 1_X: B\oslash X\to Y\oslash X$ --- непосредственный мономорфизм стереотипных пространств. А по лемме \ref{LM:1_Y-oslash-pi-in-SMono} морфизм $1_Y\oslash\pi: Y\oslash A\to Y\oslash X$ --- тоже непосредственный мономорфизм стереотипных пространств. Это значит, что оба этих морфизма --- строгие мономорфизмы, и поэтому их композиция
$$
(\sigma\oslash 1_X)\circ (1_Y\oslash\pi)=\sigma\oslash\pi
$$
--- тоже строгий (и поэтому непосредственный) мономорфим стереотипных пространств. 
\epr

Из теоремы \ref{TH:oslash-sohranyaet-oquarr-i-osubarr} сразу следуют два результата о сохранении непосредственных морфизмов тензорными произведениями.

\btm\label{TH:odot-sohranyaet-osubarr}
Пусть $\alpha:A\to X$ и $\beta:B\to Y$  --- непосредственные мономорфизмы стереотипных пространств. Тогда $\alpha\odot\beta: A\odot B\to X\odot Y$  --- непосредственный мономорфизм стереотипных пространств.
\etm

\brem Иначе это можно сформулировать так. {\it 
Пусть $A$ --- непосредственное подпространство в стереотипном пространстве $X$, а $B$  --- непосредственное подпространство в стереотипном пространстве $Y$:
$$
A\osubarr X,\qquad B\osubarr Y.
$$
Тогда $A\odot B$  --- непосредственное подпространство в стереотипном пространстве $X\odot Y$:
$$
A\odot B\osubarr X\odot Y.
$$
}
\erem

\btm\label{TH:odot-sohranyaet-oquarr}
Пусть $\alpha:X\to A$ и $\beta:Y\to B$  --- непосредственные эпиморфизмы стереотипных пространств. Тогда $\alpha\circledast\beta: X\circledast Y\to A\circledast B$  --- непосредственный эпиморфизм стереотипных пространств.
\etm

\brem Иначе это можно сформулировать так. {\it 
Пусть $A$ --- непосредственное фактор-пространство стереотипного пространства $X$, а $B$  --- непосредственное фактор-пространство стереотипного пространства $Y$:
$$
A\oquarr X,\qquad B\oquarr Y.
$$
Тогда $A\circledast B$  --- непосредственное фактор-пространство стереотипного пространства $X\circledast Y$:
$$
A\circledast B\oquarr X\circledast Y.
$$
}
\erem

\section{Некоторые факты из теории топологических и стереотипных алгебр}

\subsection{Теорема Нахбина}

Пусть  $\mathcal{E}(M)$ -- алгебра комплекснозначных гладких функций на гладком многообразии $M$, наделенная обычной топологией равномерной сходимости на компактах по каждой производной. Следующий вариант теоремы Нахбина о подалгебрах в алгебре $\mathcal{E}(M)$ играет в дифференциальной геометрии ту же роль, какую теорема Стоуна-Вейерштрасса играет в топологии (\cite{Nachbin}, доказательство можно найти в \cite{Llavona}).

\btm\label{TH:Nachbin} Пусть $A$ -- инволютивная стереотипная подалгебра в алгебре $\mathcal{E}(M)$ гладких функций на многообразии $M$, причем выполняется условие
 \bit{
 \item[(i)] $A$ разделяет точки $M$:
 $$
 \forall s\ne t\in M\qquad \exists a\in A\qquad a(s)\ne a(t).
 $$
}\eit\noindent
Тогда следующие условия для $A$ эквивалентны:

\bit{
 \item[(ii)] для всякой точки $s\in M$ и любого ненулевого касательного вектора $\tau\in T_s(M)$ найдется функция $a\in A$, дифференциал которой в точке $s$ на векторе $\tau$ не равен нулю:
 $$
 \forall \tau\in T_s(M)\setminus\{0\}\qquad \exists a\in A\qquad \d a(s)(\tau)\ne 0.
 $$

 \item[(iii)] в каждой точке $s\in M$ естественное отображение кокасательных пространств
 $$
 T_s^*[A]\to T_s^*(M)
 $$
 является изоморфизмом векторных пространств;

\item[(iv)] для любого числа $n\in\N$ и любой точки $s\in M$ естественное отображение кокасательных алгебр струй
 $$
 J_s^n[A]\to J_s^n(M)
 $$
 является изоморфизмом векторных пространств, или, что эквивалентно, алгебр;

\item[(v)] алгебра $A$ плотна в $\mathcal{E}(M)$.

 }\eit
\etm

\subsection{Лемма о плотном идеале}

\blm[о плотном идеале]\label{LM:o-plotnom-ideale}
Пусть $A$ --- плотная (и содержащая единицу) подалгебра в функциональной алгебре $B=\mathcal{C}(M)$ или $B=\mathcal{E}(M)$ (в частности, $A$ содержит константы). Тогда для всякой точки $t\in M$ идеал
$$
I_t^A=\{a\in A:\quad a(t)=0\}
$$
плотен в идеале
$$
I_t^B=\{u\in \mathcal{C}(M):\quad u(t)=0\}.
$$
\elm
\bpr
Пусть $u\in I_t^{B}$, то есть $u(t)=0$. Поскольку алгебра $A$ плотна в $\mathcal{C}(M)$,
найдется направленность $a_i\in A$
$$
a_i\overset{B}{\underset{i\to\infty}{\longrightarrow}}u.
$$
Положим $b_i=a_i-a_i(t)$. Тогда $b_i(t)=0$, поэтому $b_i\in I_t^A$, а с другой стороны,
$$
b_i=a_i-a_i(t)\overset{B}{\underset{i\to\infty}{\longrightarrow}}u-\underbrace{u(t)}_{\scriptsize\begin{matrix}\|\\ 0\end{matrix}}=u.
$$
\epr

\subsection{Лемма о тензорном произведении стереотипных алгебр}

\blm\label{LM:ph:A-circledast-B->C}
Пусть $\alpha:A\to C$ и $\beta:B\to C$ -- морфизмы стереотипных алгебр с коммутирующими образами:
\beq\label{alpha-betta-kommut}
\alpha(a)\cdot\beta(b)=\beta(b)\cdot\alpha(a),\qquad a\in A,\ b\in B.
\eeq
Тогда отображение $\ph:A\circledast B\to C$, определенное тождеством
\beq\label{ph:A-circledast-B->C}
\ph(a\circledast b)=\alpha(a)\cdot\beta(b),\qquad a\in A,\ b\in B,
\eeq
является морфизмом стереотипных алгебр. И наоборот, всякий морфизм стереотипных алгебр $\ph:A\circledast B\to C$ единственным образом представим в виде \eqref{ph:A-circledast-B->C}, где $\alpha:A\to C$ и $\beta:B\to C$ -- морфизмы стереотипных алгебр, удовлетворяющие \eqref{alpha-betta-kommut}. Если $A$, $B$, $C$ -- инволютивные алгебры и $\ph$ -- морфизм инволютивных алгебр, то $\alpha$ и $\beta$ -- тоже морфизмы инволютивных алгебр.
\elm
\bpr Здесь не вполне очевидна только вторая часть утверждения. Пусть $\ph:A\circledast B\to C$ -- морфизм алгебр.
Положим
$$
\alpha(a)=\ph(a\circledast 1_B),\qquad \beta(b)=\ph(1_A\circledast b),\qquad a\in A,\ b\in B.
$$
Тогда, во-первых,
$$
\ph(a\circledast b)=\ph\Big((a\circledast 1_B)\cdot(1_A\circledast b)\Big)=\ph(a\circledast 1_B)\cdot\ph(1_A\circledast b)=\alpha(a)\cdot\beta(b),
$$
и, во-вторых,
$$
\alpha(a)\cdot\beta(b)=\ph(a\circledast b)=\ph\Big((1_A\circledast b)\cdot(a\circledast 1_B)\Big)=\ph(1_A\circledast b)\cdot\ph(a\circledast 1_B)=\beta(b)\cdot\alpha(a).
$$
\epr

\subsection{Действие гомоморфизмов в $C^*$-алгебры на спектры}

\btm\label{LM:overline(ph(I_t)-cdot-B)=B} Пусть $\ph :A\to B$ --- инволютивный гомоморфизм стереотипных алгебр, причем $A$ коммутативна, $B$ --- $C^*$-алгебра и $\ph(A)$ плотно в $B$. Тогда

\bit{

\item[(i)] для всякой точки $t\in\Spec(B)$ выполняется равенство
\beq\label{overline(ph(Ker(t-circ-ph)))=Ker_t}
\overline{\ph\big(\Ker (t\circ\ph)\big)}=\Ker t
\eeq

\item[(ii)] для всякой точки
$s\in\Spec(A)\setminus\big(\Spec(B)\circ\ph \big)$ выполняется равенство
\beq\label{overline(ph(Ker_s))=B}
\overline{\ph(\Ker s)}=B
\eeq
}\eit
\etm

Для доказательства нам понадобится следующее обозначение. Пусть $X$ -- левый стереотипный модуль над стереотипной алгеброй $A$. Если $M$
-- подпространство в $A$ и $N$ -- подпространство в $X$, то мы обозначаем
\beq\label{DEF:M-cdot-N} M\cdot N=\left\{ \sum_{i=1}^k m_i\cdot n_i ;\ m_i\in
M, n_i\in N, k\in \Bbb N\right\} \eeq

Очевидно, для любых $L,M,N\subseteq A$
\beq\label{(LM)N=L(MN)}
(L\cdot M)\cdot N=L\cdot (M\cdot N)
\eeq
и
\beq\label{(M^N)^=(MN)^}
\overline{M\cdot N}=\overline{\overline{M} \cdot N}= \overline{M\cdot
\overline{N}}=\overline{\overline{M}\cdot \overline{N}}
\eeq

\blm\label{LM:I-ideal=>ph(I)-ideal}
Пусть  $\ph :A\to B$ -- морфизм стереотипных алгебр. Тогда
\bit{
\item[(i)] для любых подмножеств $M,N\subseteq A$ выполняется включение
$$
\ph(\overline{M\cdot N})\subseteq \overline{\ph (M)\cdot \ph (N)}
$$
\item[(ii)] для всякого левого (правого) идеала $I$ в $A$ множество $\overline{\ph(I)}$ есть левый (правый) идеал в $\overline{\ph(A)}$.
}\eit
\elm
\bpr Часть (i) отмечалась в \cite[Lemma 13.10]{Ak03}, поэтому остается доказать (ii).
Пусть $b\in \overline{\ph(A)}$ и $y\in\overline{\ph(I)}$. Тогда
$$
b\underset{\infty\gets i}{\longleftarrow} \ph(a_i),\qquad
y\underset{\infty\gets j}{\longleftarrow} \ph(x_j),
$$
для некоторых направленностей $a_i\in A$ и $x_j\in I$. Отсюда
$$
b\cdot y\underset{\infty\gets i}{\longleftarrow} \ph(a_i)\cdot y\underset{\infty\gets j}{\longleftarrow} \ph(a_i)\cdot \ph(x_j)=\ph(a_i\cdot x_j).
$$
и поэтому $b\cdot y\in \overline{\ph(I)}$.
\epr

\bpr[Доказетельство теоремы \ref{LM:overline(ph(I_t)-cdot-B)=B}]
1. Пусть $t\in\Spec(B)$. Прежде всего,
$$
a\in \Ker (t\circ\ph)\quad\Rightarrow\quad t(\ph(a))=(t\circ\ph)(a)=0\quad\Rightarrow\quad \ph(a)\in\Ker t.
$$
Отсюда $\ph\big(\Ker (t\circ\ph)\big)\subseteq\Ker t$, и поэтому $\overline{\ph\big(\Ker (t\circ\ph)\big)}\subseteq\Ker t$.

Пусть наоборот, $b\in \Ker t$. Поскольку $\ph(A)$ плотно в $B$, для всякого $\e>0$ можно подобрать элемент $a_{\e}\in A$ такой, что
$$
\norm{b-\ph(a_{\e})}<\e.
$$
Положим
$$
a'_{\e}=a_{\e}-t\big(\ph(a_{\e})\big)\cdot 1_A.
$$
Тогда, во-первых,
$$
t\big(\ph(a'_{\e})\big)=t\big(\ph(a_{\e})\big)-t\big(\ph(a_{\e})\big)\cdot t\big(\ph(1_A)\big)=0.
$$
То есть $a'_{\e}\in\Ker (t\circ\ph)$. А, во-вторых,
$$
\norm{\ph(a_{\e})-\ph(a'_{\e})}=\norm{\ph(a_{\e})-\ph(a_{\e})+t\big(\ph(a_{\e})\big)\cdot\ph(1_A)}=\abs{t\big(\ph(a_{\e})\big)}=
\Big|\underbrace{t(b)}_{\scriptsize\begin{matrix}\text{\rotatebox{90}{$=$}}\\ 0\end{matrix}}-t\big(\ph(a_{\e})\big)\Big|\le \norm{t}\cdot \norm{b-\ph(a_{\e})}<\e.
$$
Отсюда
$$
\norm{b-\ph(a'_{\e})}\le \norm{b-\ph(a_{\e})}+\norm{\ph(a_{\e})-\ph(a'_{\e})}<2\e.
$$
Мы получили, что для любых $b\in \Ker t$ и $\e>0$ найдется элемент $a'_{\e}\in\Ker(t\circ\ph)$ такой, что $\norm{b-\ph(a'_{\e})}<2\e$. Это доказывает включение $\overline{\ph\big(\Ker (t\circ\ph)\big)}\supseteq\Ker t$.

2. Зафиксируем $s\in\Spec(A)\setminus\big(\Spec(B)\circ\ph \big)$. Для всякой точки $t\in \Spec(B)$ найдется элемент $x_t\in A$, отделяющий $s$ от $t\circ\ph$:
$$
s(x_t)\ne t(\ph(x_t)).
$$
Как следствие, элемент $a_t=x_t-s(x_t)\cdot 1_A$ должен обладать свойством
\beq\label{s(a_t)=0}
s(a_t)=s(x_t)-s(x_t)\cdot s(1_A)=s(x_t)-s(x_t)=0
\eeq
и свойством
\beq\label{t(ph(a_t))-ne-0}
t(\ph(a_t))=t(\ph(x_t))-s(x_t)\cdot t(\ph(1_A))=t(\ph(x_t))-s(x_t)\ne 0.
\eeq
Из \eqref{t(ph(a_t))-ne-0} следует, что множества
$$
U_t=\{r\in\Spec(A):\ r(a_t)\ne 0 \}
$$
покрывают компакт $\Spec(B)\circ\ph$. Значит, среди них имеется конечное подпокрытие $U_{t_1},...,U_{t_n}$:
$$
\bigcup_{i=1}^n U_{t_i}\supseteq \Spec(B)\circ\ph.
$$
Положим
$$
a=\sum_{i=1}^n a_{t_i}\cdot a_{t_i}^\bullet.
$$
Этот элемент в точке $s$ равен нулю, потому что в силу \eqref{s(a_t)=0}, все $a_{t_i}$ в ней равны нулю,
$$
s(a)=\sum_{i=1}^n s(a_{t_i})\cdot \overline{s(a_{t_i})}=0.
$$
Как следствие, $a\in \Ker s$, и значит,
$$
\ph(a)\in \overline{\ph(\Ker s)}.
$$
А, с другой стороны, элемент $\ph(a)$ отличен от нуля всюду на $\Spec(B)$, потому что на каждой точке $t\in\Spec(B)$ какой-нибудь элемент $\ph(a_{t_i})$ отличен от нуля:
$$
t(\ph(a))=\sum_{i=1}^n t(\ph(a_{t_i}))\cdot \overline{t(\ph(a_{t_i}))}=\sum_{i=1}^n \abs{t(\ph(a_{t_i}))}^2>0.
$$
То есть $\ph(a)$ отличен от нуля всюду на спектре $\Spec(B)$ коммутативной $C^*$-алгебры $B\cong C(\Spec(B))$, и как следствие, он обратим в $B$:
$$
b\cdot\ph(a)=1_B
$$
для некоторого $b\in B$. Поскольку здесь $\ph(a)$ лежит в левом идеале $\overline{\ph(\Ker s)}$ алгебры $B=\overline{\ph(A)}$ (по лемме \ref{LM:I-ideal=>ph(I)-ideal}(ii)), мы получаем, что $1_B$ тоже лежит в $\overline{\ph(\Ker s)}$. Поэтому
$$
\overline{\ph(\Ker s)}\supseteq B\cdot\overline{\ph(\Ker s)}\supseteq B\cdot 1_B=B.
$$
\epr

\bcor\label{COR:overline(ph(I_t)-cdot-B)=B} Пусть $\ph :A\to B$ --- инволютивный гомоморфизм инволютивных стереотипных алгебр, причем $A$ коммутативна, а $B$ --- $C^*$-алгебра. Тогда для всякой точки
$$
s\in\Spec(A)\setminus\big(\Spec\big(\overline{\ph(A)}\big)\circ\ph\big)
$$
выполняется равенство
\beq
\overline{\ph(\Ker s)\cdot B}=B
\eeq
\ecor
\bpr По теореме \ref{LM:overline(ph(I_t)-cdot-B)=B}(ii), множество $\overline{\ph(\Ker s)}$ содержит единицу алгебры $\overline{\ph(A)}$, которая совпадает с единицей алгебры $B$:
$$
\overline{\ph(\Ker s)}\owns 1_B.
$$
Поэтому
$$
\overline{\ph(\Ker s)\cdot B}\supseteq \overline{\ph(\Ker s)}\cdot B\supseteq 1_B\cdot B=B.
$$
\epr

\subsection{Идеалы в алгебрах ${\mathcal C}(M)$ и ${\mathcal E}(M)$}

\paragraph{Идеалы в алгебре ${\mathcal C}(M)$.}

\bit{

\item[$\bullet$] Для всякого идеала $I\subseteq {\mathcal C}(M)$ обозначим символом ${\mathcal Z}(I)$ множество его нулей, то есть тех точек в $M$, в которых все функции $u\in I$ обнуляются:
    \beq\label{DEF:Z(F)}
    {\mathcal Z}(I)=\{t\in M: \quad \forall u\in I \quad u(t)=0\}
    \eeq 
    Понятно, что ${\mathcal Z}(I)$ --- всегда замкнутое подмножество в $M$.

\item[$\bullet$] Для всякого множества $N\subseteq M$ обозначим символом ${\mathcal I}(N)$ множество функций на $M$, обнуляющихся на $N$:
    \beq\label{DEF:I(N)}
    {\mathcal I}(N)=\{u\in {\mathcal C}(M): \quad u\big|_N=0\}
    \eeq 
    Понятно, что ${\mathcal I}(N)$ --- всегда замкнутый идеал в ${\mathcal C}(M)$.

\item[$\bullet$] Введем еще такие два дополниительных обозначения. Пусть $X$ и $Y$ --- подпространства в стереотипнной алгебре $A$. Тогда
    
    \bit{
\item[---] символом $X+Y$ мы обозначаем подпространство в $A$, состоящее из всевозможных сумм вида $x+y$, где $x\in X$, $y\in Y$:
\beq\label{DEF:X+Y}
X+Y=\{x+y; \ x\in X, \ y\in Y\},
\eeq

\item[---] символом $X\cdot Y$ мы обозначаем подпространство в $A$, состоящее из всевозможных (конечных) сумм вида $\sum_{i=1}^n x_i\cdot y_i$, где $x_i\in X$, $y_i\in Y$:
\beq\label{DEF:X-cdot-Y}
X\cdot Y=\Big\{\sum_{i=1}^n x_i\cdot y_i; \ x_i\in X, \ y_i\in Y,\ n\in\N \Big\}
\eeq
(понятно, что если $X$ и $Y$ --- идеалы в $A$, то $X\cdot Y$ --- тоже идеал в $A$).
}\eit
}\eit

\medskip
\centerline{\bf Свойства операций ${\mathcal Z}$ и ${\mathcal I}$ в ${\mathcal C}(M)$:}

\bit{\it

\item[$1^\circ$.] Антимонотонность ${\mathcal Z}$:
\beq\label{antimonot-Z(I)}
I\subseteq J\quad\Rightarrow\quad {\mathcal Z}(I)\supseteq {\mathcal Z}(J)
\eeq

\item[$2^\circ$.] Экстремальные значения ${\mathcal Z}$: если $I$ --- замкнутый идеал в ${\mathcal C}(M)$, то 
\beq\label{extremumy-Z(I)-1}
{\mathcal Z}(I)=M \quad\Leftrightarrow\quad I=0.
\eeq
\beq\label{extremumy-Z(I)-0}
{\mathcal Z}(I)=\varnothing \quad\Leftrightarrow\quad I={\mathcal C}(M).
\eeq

\item[$3^\circ$.] Антимонотонность ${\mathcal I}$:
\beq\label{antimonot-I(N)}
L\subseteq N\quad\Rightarrow\quad {\mathcal I}(L)\supseteq {\mathcal I}(N)
\eeq

\item[$4^\circ$.] Экстремальные значения ${\mathcal I}$: если $N$ --- замкнутое подмножество в $M$, то 
\beq\label{extremumy-I(N)-1}
{\mathcal I}(N)={\mathcal C}(M) \quad\Leftrightarrow\quad N=\varnothing.
\eeq
\beq\label{extremumy-I(N)-0}
{\mathcal I}(N)=0 \quad\Leftrightarrow\quad N=M.
\eeq

\item[$5^\circ$.] Теоремы о нулях: 

\bit{\it

\item[(a)] если $N$ --- подмножество в $M$, то последовательное применение к нему операций $\mathcal I$ и $\mathcal Z$ дает замыкание $N$ в $M$:
\beq\label{Z(I(N))=N}
{\mathcal Z}({\mathcal I}(N))=\overline{N}
\eeq

\item[(b)] если $I$ --- идеал в ${\mathcal C}(M)$, то последовательное применение к нему операций $\mathcal Z$ и $\mathcal I$ дает замыкание $I$ в ${\mathcal C}(M)$: 
\beq\label{I(Z(I))=I}
{\mathcal I}({\mathcal Z}(I))=\overline{I}
\eeq

}\eit

\item[$6^\circ$.] Теоремы о сумме идеалов: 
\bit{\it

\item[(a)] если $L$ и $N$ --- замкнутые подмножества в $M$, то\footnote{Здесь сумма ${\mathcal I}(L)+{\mathcal I}(N)$ понимается в смысле \eqref{DEF:X+Y}.} 
\beq\label{Z(I(L)+I(N))=L-cap-N}
{\mathcal Z}({\mathcal I}(L)+{\mathcal I}(N))=L\cap N
\eeq

\item[(b)] если $I$ и $J$ --- идеалы в ${\mathcal C}(M)$, то\footnote{Здесь сумма $I+J$ понимается в смысле \eqref{DEF:X+Y}.} 
\beq\label{I(Z(I)-cap-Z(J))=I+J}
{\mathcal I}({\mathcal Z}(I)\cap {\mathcal Z}(J))=\overline{I+J}
\eeq
}\eit

\item[$7^\circ$.] Теоремы о произведении идеалов: 
\bit{\it

\item[(a)] если $L$ и $N$ --- подмножества в $M$, то\footnote{Здесь произведение ${\mathcal I}(L)\cdot {\mathcal I}(N)$ понимается в смысле \eqref{DEF:X-cdot-Y}.} 
\beq\label{Z(I(L)-cdot-I(N))=L-cup-N}
{\mathcal Z}({\mathcal I}(L)\cdot {\mathcal I}(N))=\overline{L\cup N},
\eeq

\item[(b)] если $I$ и $J$ --- идеалы в ${\mathcal C}(M)$, то\footnote{Здесь произведение $I\cdot J$ понимается в смысле \eqref{DEF:X-cdot-Y}.}  
\beq\label{I(Z(I)-cup-Z(J))=I-cdot-J}
{\mathcal I}({\mathcal Z}(I)\cup {\mathcal Z}(J))=\overline{I\cdot J}
\eeq

}\eit

}\eit

\blm\label{LM:f_T|_T=1}
Пусть $I$ --- идеал в ${\mathcal C}(M)$. Тогда для всякого компакта $T\subseteq M\setminus {\mathcal Z}(I)$ найдется функция $f\in I$ такая, что 
\beq\label{sup|f(t)|=1}
\sup_{t\in M}\abs{f(t)}=1 
\eeq
и
\beq\label{f|_T=1}
f\big|_T=1. 
\eeq
\elm
\bpr
Пусть $t\notin {\mathcal Z}(I)$. Тогда найдется функция $u\in I$, которая в этой точке не равна нулю:
$$
u(t)\ne 0.
$$ 
Умножив эту функцию на сопряженную ей функцию $\overline{u}$, мы получим
$$
0\le u(t)\cdot\overline{u(t)}\ne 0,
$$ 
то есть 
$$
u(t)\cdot\overline{u(t)}>0,
$$ 
и, поскольку $I$ --- идеал, мы получаем, что  в $I$ имеется некая вещественная функция $v=u\cdot\overline{u}\in I$, положительная в точке $t$:
$$
v(t)>0.
$$   
Умножив ее на ``шапочку'' с достаточно малым носителем $\eta\in {\mathcal C}(M)$, мы получим функцию $w=v\cdot\eta \in I$, такую что
$$
w\ge 0,\qquad w(t)>0.
$$ 
Отсюда в свою очередь следует, что для всякого компакта $T\subseteq M\setminus {\mathcal Z}(I)$ можно найти некую функцию $\omega\in I$ со свойствами
\beq\label{PROOF:I(Z(I))=I-1}
\omega\ge 0,\qquad \forall t\in T\qquad \omega(t)>0.
\eeq
Подберем теперь окрестность $V$ компакта $T$ так, чтобы
$$
\forall s\in V\quad \omega(s)>0.
$$
Поскольку $M$ --- паракомпактное пространство, для $T$ и $V$ можно выбрать функцию $\theta\in{\mathcal C}(M)$ со свойствами
$$
\theta\ge 0,\qquad \theta\big|_T=0,\qquad \theta\big|_{M\setminus V}=1. 
$$
Функция $\theta+\omega$ будет обладать следующими свойствами:
$$
\theta+\omega>0,\qquad (\theta+\omega)\big|_T=\omega\big|_T. 
$$
Поэтому функция $f=\frac{\omega}{\theta+\omega}$ будет обладать свойствами \eqref{sup|f(t)|=1} и \eqref{f|_T=1}. С другой стороны, $f\in I$, потому что она получена умножением функции $\omega\in I$ на функцию $\frac{1}{\theta+\omega}\in{\mathcal C}(M)$.
\epr

\bpr[Доказательство свойств $1^\circ-7^\circ$ на с.\pageref{antimonot-Z(I)}.]
1. Свойство $1^\circ$ очевидно.

2. Пусть $I$ --- замкнутый идеал в ${\mathcal C}(M)$. Свойство \eqref{extremumy-Z(I)-1} очевидно, а в \eqref{extremumy-Z(I)-0} нужно доказать импликацию
$$
{\mathcal Z}(I)=\varnothing \quad\Rightarrow\quad I={\mathcal C}(M).
$$
Пусть ${\mathcal Z}(I)=\varnothing$. Тогда по лемме \ref{LM:f_T|_T=1} для всякого компакта $K\subseteq M=M\setminus {\mathcal Z}(I)$ найдется функция $f\in I$ со свойством 
$$
f\big|_T=1. 
$$
Это значит, что тождественная единица 1 приближается элементами идеала $I$. А, поскольку $I$ замкнут, мы получаем, что
$1\in I$. Отсюда следует, что $I={\mathcal C}(M)$.

3. Свойство $3^\circ$ очевидно.

4. Пусть $N$ --- замкнутое подмножество в $M$. Тогда свойство \eqref{extremumy-I(N)-1} очевидно, а в свойстве \eqref{extremumy-I(N)-0} нужно доказать импликацию
$$
{\mathcal I}(N)=0 \quad\Rightarrow\quad N=M.
$$
Это легче сделать, превратив ее в импликацию
$$
N\ne M  \quad\Rightarrow\quad {\mathcal I}(N)\ne 0.
$$
Действительно, если $N\ne M$, то можно взять точку $t\in M\setminus N$ и подобрать для нее функцию  $u\in {\mathcal C}(M)$ так, чтобы
$$
u\Big|_N=0,\qquad u(t)=1.
$$
Это всегда можно сделать, поскольку $M$ --- паракомпактное пространство. Тогда мы получим, что $0\ne u\in {\mathcal I}(N)$, и значит ${\mathcal I}(N)\ne 0$. 

5. Докажем сначала \eqref{Z(I(N))=N}. Очевидно, 
$$
N\subseteq {\mathcal Z}({\mathcal I}(N))
$$
Поскольку ${\mathcal Z}({\mathcal I}(N))$ замкнуто, мы отсюда получаем
$$
\overline{N}\subseteq {\mathcal Z}({\mathcal I}(N)).
$$
Наоборот, если $t\notin \overline{N}$, то, поскольку $M$ --- паракомпактное пространство, найдется функция $u\in {\mathcal C}(M)$ такая, что
$$
u\Big|_{\overline{N}}=0,\qquad u(t)=1.
$$
Здесь первое условие влечет за собой $u\in {\mathcal I}(N)$, и вместе со вторым условием, это означает, что $t\notin {\mathcal Z}({\mathcal I}(N))$.  

Теперь докажем \eqref{I(Z(I))=I}. Пусть $I$ --- идеал в ${\mathcal C}(M)$. Очевидно, 
$$
I\subseteq {\mathcal I}({\mathcal Z}(I))
$$
Поскольку ${\mathcal I}({\mathcal Z}(I))$ --- замкнутое множество, мы отсюда получаем
$$
\overline{I}\subseteq {\mathcal I}({\mathcal Z}(I)).
$$
Нам теперь нужно доказать обратное включение 
\beq\label{PROOF:I(Z(I))=I-0}
{\mathcal I}({\mathcal Z}(I))\subseteq \overline{I}.
\eeq
Зафиксируем функцию $u\in {\mathcal I}({\mathcal Z}(I))$, компакт $K\subseteq M$ и число $\e>0$. Пусть 
\beq\label{PROOF:I(Z(I))=I-1}
U=\Big\{t\in M: \ \abs{u(t)}<\frac{\e}{2}\Big\}.
\eeq
Условие $u\in {\mathcal I}({\mathcal Z}(I))$ означает, что 
$$
u\big|_{{\mathcal Z}(I)}=0
$$  
и поэтому
$$
{\mathcal Z}(I)\subseteq U.
$$
Как следствие, компакт 
$$
T=K\setminus U
$$
не пересекается с ${\mathcal Z}(I)$:
$$
T\subseteq M\setminus {\mathcal Z}(I).
$$
Поэтому по лемме \ref{LM:f_T|_T=1} найдется функция $f\in I$ со свойствами 
\beq\label{PROOF:I(Z(I))=I-2}
\sup_{t\in M}\abs{f(t)}=1,\qquad f\big|_T=1. 
\eeq
Теперь рассмотрим функцию
$$
v=u\cdot f.
$$
С одной стороны, она лежит в идеале $I$, 
\beq\label{PROOF:I(Z(I))=I-3}
v\in I,
\eeq
потому что $u\in I$. А, с другой, на компакте $K\subseteq M$ она отстоит от функции $u$ меньше чем на $\e$, 
\beq\label{PROOF:I(Z(I))=I-4}
\sup_{t\in K}\abs{v(t)-u(t)}<\e
\eeq
потому что, с одной стороны,
\begin{multline*}
\sup_{t\in K\cap U}\abs{v(t)-u(t)}=
\sup_{t\in K\cap U}\abs{f(t)\cdot u(t)-u(t)}=\\=\sup_{t\in K\cap U}\abs{(f(t)-1)\cdot u(t)}=
\underbrace{\sup_{t\in K\cap U}\abs{f(t)-1}}_{\scriptsize\begin{matrix}
\text{\rotatebox{90}{$\ge$}} \\
\sup_{t\in K\cap U}\abs{f(t)}+1
\\
\text{\rotatebox{90}{$\ge$}}\put(2,0){\eqref{PROOF:I(Z(I))=I-2}}
\\ 2 \end{matrix}}\cdot 
\underbrace{\sup_{t\in K\cap U}\abs{u(t)}}_{\scriptsize\begin{matrix}
\text{\rotatebox{90}{$\ge$}}\put(2,0){\eqref{PROOF:I(Z(I))=I-1}}
\\ \frac{\e}{2} \end{matrix}}
\le 2\cdot \frac{\e}{2}= \e
\end{multline*}
а, с другой,
\begin{multline*}
\sup_{t\in K\setminus U}\abs{v(t)-u(t)}=
\sup_{t\in T}\abs{v(t)-u(t)}=
\sup_{t\in T}\abs{f(t)\cdot u(t)-u(t)}=\\=
\sup_{t\in T}\abs{(f(t)-1)\cdot u(t)}=
\underbrace{\sup_{t\in T}\abs{f(t)-1}}_{\scriptsize\begin{matrix}
\| \put(2,0){\eqref{PROOF:I(Z(I))=I-2}}
\\ 0 \end{matrix}}\cdot 
\sup_{t\in T}\abs{u(t)}=0
\end{multline*}

Условия \eqref{PROOF:I(Z(I))=I-3} и \eqref{PROOF:I(Z(I))=I-4} вместе означают, что функция $u\in {\mathcal I}({\mathcal Z}(I))$ приближается в ${\mathcal C}(M)$ функциями $v\in I$. Это доказывает \eqref{PROOF:I(Z(I))=I-0}.

6. Для доказательства \eqref{Z(I(L)+I(N))=L-cap-N}, зафиксируем замкнутые подмножества $L,N\subseteq M$ и покажем сначала, что 
\beq\label{PROOF:Z(I(L)+I(N))=L-cap-N-1}
{\mathcal Z}({\mathcal I}(L)+{\mathcal I}(N))\subseteq L\cap N
\eeq
Действительно,  с одной стороны,
$$
{\mathcal I}(L)+{\mathcal I}(N)\supseteq {\mathcal I}(L)
$$
$$
\Downarrow
$$
\beq\label{PROOF:Z(I(L)+I(N))=L-cap-N-2}
{\mathcal Z}({\mathcal I}(L)+{\mathcal I}(N))\subseteq {\mathcal Z}({\mathcal I}(L))=\eqref{Z(I(N))=N}=L
\eeq
А, с другой стороны,
$$
{\mathcal I}(L)+{\mathcal I}(N)\supseteq {\mathcal I}(N)
$$
$$
\Downarrow
$$
\beq\label{PROOF:Z(I(L)+I(N))=L-cap-N-3}
{\mathcal Z}({\mathcal I}(L)+{\mathcal I}(N))\subseteq {\mathcal Z}({\mathcal I}(N))=\eqref{Z(I(N))=N}=N
\eeq
И вместе \eqref{PROOF:Z(I(L)+I(N))=L-cap-N-2} и \eqref{PROOF:Z(I(L)+I(N))=L-cap-N-3} дают \eqref{PROOF:Z(I(L)+I(N))=L-cap-N-1}.

Теперь докажем обратное включение:
\beq\label{PROOF:Z(I(L)+I(N))=L-cap-N-4}
{\mathcal Z}({\mathcal I}(L)+{\mathcal I}(N))\supseteq L\cap N
\eeq
Действительно, 
$$
t\in L\cap N
$$
$$
\Downarrow
$$
$$
\underbrace{\underbrace{t\in L}_{\scriptsize \begin{matrix}\Downarrow \\ \forall u\in {\mathcal I}(L)\quad u(t)=0  \end{matrix}} 
\quad \& \quad \underbrace{t\in N}_{\scriptsize \begin{matrix}\Downarrow \\ \forall v\in {\mathcal I}(N)\quad v(t)=0  \end{matrix}} 
}_{\scriptsize \begin{matrix}\Downarrow \\ \forall u\in {\mathcal I}(L)\quad \forall v\in {\mathcal I}(N)\quad u(t)+v(t)=0  \end{matrix}}
$$
$$
\Downarrow
$$
$$
t\in {\mathcal Z}({\mathcal I}(L)+{\mathcal I}(N))
$$
 
Теперь докажем \eqref{I(Z(I)-cap-Z(J))=I+J}. Пусть $I$ и $J$ --- идеалы в ${\mathcal C}(M)$. 
Обозначим
$$
L={\mathcal Z}(I),\qquad N={\mathcal Z}(J).
$$
Тогда 
\beq\label{PROOF:I(Z(I)-cap-Z(J))=I+J-1}
{\mathcal I}(L)={\mathcal I}({\mathcal Z}(I))=\eqref{I(Z(I))=I}=\overline{I}
,\qquad 
{\mathcal I}(N)={\mathcal I}({\mathcal Z}(J))=\eqref{I(Z(I))=I}=\overline{J}
\eeq
и поэтому
$$
{\mathcal Z}(I)\cap {\mathcal Z}(J)=L\cap N=\eqref{Z(I(L)+I(N))=L-cap-N}={\mathcal Z}({\mathcal I}(L)+{\mathcal I}(N))=
\eqref{PROOF:I(Z(I)-cap-Z(J))=I+J-1}= {\mathcal Z}(\overline{I}+\overline{J}).
$$
Отсюда 
$$
{\mathcal I}({\mathcal Z}(I)\cap {\mathcal Z}(J))={\mathcal I}({\mathcal Z}(\overline{I}+\overline{J}))=\eqref{I(Z(I))=I}=
\overline{\overline{I}+\overline{J}}=\overline{I+J}.
$$

7. Чтобы доказать \eqref{Z(I(L)-cdot-I(N))=L-cup-N}, зафиксируем подмножества $L,N\subseteq M$ и покажем сначала, что 
\beq\label{PROOF:Z(I(L)-cdot-I(N))=L-cup-N-1}
{\mathcal Z}({\mathcal I}(L)\cdot {\mathcal I}(N))\supseteq \overline{L\cup N}.
\eeq
Действительно, с одной стороны,
$$
\forall u\in {\mathcal I}(L)\quad u\big|_L=0
$$
$$
\Downarrow
$$
$$
\forall u_i\in {\mathcal I}(L)\quad \forall v_i\in {\mathcal I}(N) \quad \Big(\sum_{i=1}^n u_i\cdot v_i\Big)\Big|_L=0
$$
$$
\Downarrow
$$
$$
\forall u_i\in {\mathcal I}(L)\quad \forall v_i\in {\mathcal I}(N) \quad \sum_{i=1}^n u_i\cdot v_i\in {\mathcal I}(L)
$$
$$
\Downarrow
$$
$$
{\mathcal I}(L)\cdot {\mathcal I}(N) \subseteq {\mathcal I}(L)
$$
$$
\Downarrow
$$
$$
{\mathcal Z}({\mathcal I}(L)\cdot {\mathcal I}(N)) \supseteq {\mathcal Z}({\mathcal I}(L))=\eqref{Z(I(N))=N}=\overline{L}
$$
А, с другой, по тем же причинам,
$$
{\mathcal Z}({\mathcal I}(L)\cdot {\mathcal I}(N)) \supseteq {\mathcal Z}({\mathcal I}(N))=\eqref{Z(I(N))=N}=\overline{N}
$$
И вместе мы получаем
$$
{\mathcal Z}({\mathcal I}(L)\cdot {\mathcal I}(N)) \supseteq \overline{L}\cup \overline{N}=\overline{L\cup N}
$$

Теперь докажем обратное включение:
\beq\label{PROOF:Z(I(L)-cdot-I(N))=L-cup-N-2}
{\mathcal Z}({\mathcal I}(L)\cdot {\mathcal I}(N))\subseteq \overline{L\cup N}.
\eeq
Это удобно сделать так:
$$
t\notin \overline{L\cup N}=\overline{L}\cup \overline{N}=\eqref{Z(I(N))=N}={\mathcal Z}({\mathcal I}(L))\cup {\mathcal Z}({\mathcal I}(N))
$$
$$
\Downarrow
$$
$$
t\notin {\mathcal Z}({\mathcal I}(L))\quad \& \quad t\notin {\mathcal Z}({\mathcal I}(N))
$$
$$
\Downarrow
$$
$$
\exists u\in {\mathcal I}(L)\quad u(t)\ne 0\quad\&\quad \exists v\in {\mathcal I}(N)\quad v(t)\ne 0
$$
$$
\Downarrow
$$
$$
\exists u\in {\mathcal I}(L)\quad \exists v\in {\mathcal I}(N)\quad u(t)\cdot v(t)\ne 0
$$
$$
\Downarrow
$$
$$
t\notin {\mathcal Z}({\mathcal I}(L)\cdot {\mathcal I}(N))
$$

Наконец, нам остается доказать \eqref{I(Z(I)-cup-Z(J))=I-cdot-J}. Пусть $I$ и $J$ --- идеалы в ${\mathcal C}(M)$. Обозначим
$$
L={\mathcal Z}(I), \qquad N={\mathcal Z}(J).
$$
Тогда в силу \eqref{I(Z(I))=I},
\beq\label{PROOF:I(Z(I)-cup-Z(J))=I-cdot-J-0}
\overline{I}={\mathcal I}(L), \qquad \overline{J}={\mathcal I}(N).
\eeq
и поэтому 
\beq\label{PROOF:I(Z(I)-cup-Z(J))=I-cdot-J-1}
{\mathcal Z}(I\cdot J)={\mathcal Z}(\overline{I}\cdot \overline{J})=\eqref{PROOF:I(Z(I)-cup-Z(J))=I-cdot-J-0}={\mathcal Z}({\mathcal I}(L)\cdot {\mathcal I}(N))=\overline{L\cup N}=\overline{L}\cup \overline{N}=L\cup N={\mathcal Z}(I)\cup {\mathcal Z}(J).
\eeq
Отсюда 
$$
\overline{I\cdot J}=\eqref{I(Z(I))=I}={\mathcal I}({\mathcal Z}(I\cdot J))=\eqref{PROOF:I(Z(I)-cup-Z(J))=I-cdot-J-1}={\mathcal I}({\mathcal Z}(I)\cup {\mathcal Z}(J)).
$$
\epr

\paragraph{Идеалы в алгебре ${\mathcal E}(M)$.}

\bit{

\item[$\bullet$] Для всякого идеала $I\subseteq {\mathcal E}(M)$ обозначим символом ${\mathcal Z}(I)$ множество его нулей, то есть тех точек в $M$, в которых все функции $u\in I$ обнуляются:
    \beq\label{DEF:Z(F)-E}
    {\mathcal Z}(I)=\{t\in M: \quad \forall u\in I \quad u(t)=0\}
    \eeq 
    Понятно, что ${\mathcal Z}(I)$ --- всегда замкнутое подмножество в $M$.

\item[$\bullet$] Для всякого множества $N\subseteq M$ обозначим символом ${\mathcal I}_0(N)$ множество функций на $M$, обнуляющихся на $N$:
    \beq\label{DEF:I_0(N)-E}
    {\mathcal I}_0(N)=\{u\in {\mathcal E}(M): \quad u\big|_N=0\}
    \eeq 
Кроме того, для любого $n\in\N$ символом ${\mathcal I}_n(N)$ мы обозначим множество функций на $M$, обнуляющихся на $N$ вместе со всеми своими дифференциалами до поряка $n$ включительно:
    \beq\label{DEF:I_n(N)-E}
    {\mathcal I}_n(N)=\{u\in {\mathcal E}(M): \quad \forall k=0,...,n\quad \d^k u\big|_N=0\}
    \eeq 
И. наконец, символом ${\mathcal I}_\infty(N)$ мы обозначим множество функций на $M$, обнуляющихся на $N$ вместе со всеми своими дифференциалами\footnote{В формуле \eqref{DEF:I_0(N)-E} считается, что $0\in\N$.}:
    \beq\label{DEF:I_0(N)-E}
    {\mathcal I}_\infty(N)=\{u\in {\mathcal E}(M): \quad \forall k\in\N\quad \d^k u\big|_N=0\}
    \eeq 
    Понятно, что ${\mathcal I}_0(N)$, ${\mathcal I}_n(N)$, ${\mathcal I}_\infty(N)$ --- всегда замкнутые идеалы в ${\mathcal E}(M)$, причем при $0\le m\le n$ справедливо включение:
    \beq\label{DEF:I_infty(N)-I_0(N)}
    {\mathcal I}_\infty(N)\subseteq {\mathcal I}_n(N)\subseteq {\mathcal I}_m(N)\subseteq{\mathcal I}_0(N)
    \eeq

}\eit

\medskip
\centerline{\bf Свойства операций ${\mathcal Z}$, ${\mathcal I}_0$ и ${\mathcal I}_\infty$ в ${\mathcal E}(M)$:}

\bit{\it

\item[$1^\circ$.] Антимонотонность ${\mathcal Z}$:
\beq\label{E:antimonot-Z(I)}
I\subseteq J\quad\Rightarrow\quad {\mathcal Z}(I)\supseteq {\mathcal Z}(J)
\eeq

\item[$2^\circ$.] Экстремальные значения ${\mathcal Z}$: если $I$ --- замкнутый идеал в ${\mathcal E}(M)$, то 
\beq\label{E:extremumy-Z(I)-1}
{\mathcal Z}(I)=M \quad\Leftrightarrow\quad I=0.
\eeq
\beq\label{E:extremumy-Z(I)-0}
{\mathcal Z}(I)=\varnothing \quad\Leftrightarrow\quad I={\mathcal E}(M).
\eeq

\item[$3^\circ$.] Антимонотонность ${\mathcal I}_0$ и ${\mathcal I}_\infty$:
\beq\label{E:antimonot-I(N)}
L\subseteq N\quad\Rightarrow\quad {\mathcal I}_\infty(L)\supseteq {\mathcal I}_\infty(N)\quad \& \quad {\mathcal I}_0(L)\supseteq {\mathcal I}_0(N)
\eeq

\item[$4^\circ$.] Экстремальные значения ${\mathcal I}_0$ и ${\mathcal I}_\infty$: если $N$ --- замкнутое подмножество в $M$, то 
\beq\label{E:extremumy-I(N)-1}
{\mathcal I}_\infty(N)={\mathcal E}(M) \quad\Leftrightarrow\quad {\mathcal I}_0(N)={\mathcal E}(M) \quad\Leftrightarrow\quad N=\varnothing.
\eeq
\beq\label{E:extremumy-I(N)-0}
{\mathcal I}_\infty(N)=0 \quad\Leftrightarrow\quad {\mathcal I}_0(N)=0 \quad\Leftrightarrow\quad N=M.
\eeq

\item[$5^\circ$.] Теоремы о нулях: 

\bit{\it

\item[(a)] если $N$ --- подмножество в $M$, то последовательное применение к нему операций $\mathcal I_0$ и $\mathcal Z$, и операций $\mathcal I_\infty$ и $\mathcal Z$, дает замыкание $N$ в $M$:
\beq\label{E:Z(I(N))=N}
{\mathcal Z}({\mathcal I}_\infty(N))=\overline{N}={\mathcal Z}({\mathcal I}_0(N))
\eeq

\item[(b)] если $I$ --- идеал в ${\mathcal E}(M)$, то последовательное применение к нему операций $\mathcal Z$ и $\mathcal I_0$, и операций $\mathcal Z$ и $\mathcal I_\infty$, дает оценки для замыкания $I$ в ${\mathcal E}(M)$: 
\beq\label{E:I(Z(I))=I}
{\mathcal I}_\infty({\mathcal Z}(I))\subseteq\overline{I}\subseteq {\mathcal I}_0({\mathcal Z}(I))
\eeq

}\eit

\item[$6^\circ$.] Теоремы о сумме идеалов: 
\bit{\it

\item[(a)] если $L$ и $N$ --- замкнутые подмножества в $M$, то\footnote{Здесь суммы ${\mathcal I}_\infty(L)+{\mathcal I}_\infty(N)$ и ${\mathcal I}_0(L)+{\mathcal I}_0(N)$ понимаются в смысле \eqref{DEF:X+Y}.} 
\beq\label{E:Z(I(L)+I(N))=L-cap-N}
{\mathcal Z}({\mathcal I}_\infty(L)+{\mathcal I}_\infty(N))=L\cap N={\mathcal Z}({\mathcal I}_0(L)+{\mathcal I}_0(N))
\eeq

\item[(b)] если $I$ и $J$ --- идеалы в ${\mathcal E}(M)$, то\footnote{Здесь сумма $I+J$ понимается в смысле \eqref{DEF:X+Y}.} 
\beq\label{E:I(Z(I)-cap-Z(J))=I+J}
{\mathcal I}_\infty({\mathcal Z}(I)\cap {\mathcal Z}(J))\subseteq\overline{I+J}\subseteq {\mathcal I}_0({\mathcal Z}(I)\cap {\mathcal Z}(J))
\eeq
}\eit

\item[$7^\circ$.] Теоремы о произведении идеалов: 
\bit{\it

\item[(a)] если $L$ и $N$ --- подмножества в $M$, то\footnote{Здесь произведения ${\mathcal I}_\infty(L)\cdot {\mathcal I}_\infty(N)$ и ${\mathcal I}_0(L)\cdot {\mathcal I}_0(N)$ понимаются в смысле \eqref{DEF:X-cdot-Y}.} 
\beq\label{E:Z(I(L)-cdot-I(N))=L-cup-N}
{\mathcal Z}({\mathcal I}_\infty(L)\cdot {\mathcal I}_\infty(N))=\overline{L\cup N}={\mathcal Z}({\mathcal I}_0(L)\cdot {\mathcal I}_0(N)),
\eeq

\item[(b)] если $I$ и $J$ --- идеалы в ${\mathcal E}(M)$, то\footnote{Здесь произведение $I\cdot J$ понимается в смысле \eqref{DEF:X-cdot-Y}.}  
\beq\label{E:I(Z(I)-cup-Z(J))=I-cdot-J}
{\mathcal I}_\infty({\mathcal Z}(I)\cup {\mathcal Z}(J))\subseteq\overline{I\cdot J}\subseteq {\mathcal I}_0({\mathcal Z}(I)\cup {\mathcal Z}(J))
\eeq

}\eit

}\eit

\subsection{Идеалы в алгебрах ${\mathcal C}(M,{\mathcal B}(X))$}

\blm\label{LM:idealy-v-C(M,B(X))}\footnote{Лемма \ref{LM:idealy-v-C(M,B(X))} подсказана автору Робертом Исраэлем.}
Пусть $M$ --- паракомпактное локально компактное топологическое пространство, и $X$ -- конечномерное векторное пространство над $\C$. Всякий замкнутый двусторонний идеал $I$ в алгебре ${\mathcal C}\big(M,{\mathcal B}(X)\big)$ имеет вид
$$
I=\{f\in {\mathcal C}\big(M,{\mathcal B}(X)\big):\ \forall t\in N\ f(t)=0\}
$$
где $N$ -- некоторое замкнутое подмножество в $M$, а соответствующая фактор-алгебра при этом имеет вид
$$
{\mathcal C}\big(M,{\mathcal B}(X)\big)/I={\mathcal C}\big(N,{\mathcal B}(X)\big).
$$
\elm
\bpr
Для всякой точки $t\in M$ обозначим 
$$
I(t)=\{u(t);\ u\in I\}\qquad \big(I(t)\subseteq {\mathcal B}(X)\big)
$$ 
и положим
$$
N=\{t\in M:\ I(t)=0\},\qquad I(N)=\{f\in {\mathcal C}\big(M,{\mathcal B}(X)\big):\ f\big|_N=0\}.
$$
Очевидно, $I(N)$ --- замкнутый двусторонний идеал в ${\mathcal C}\big(M,{\mathcal B}(X)\big)$, и $I\subseteq I(N)$. Покажем, что обратное включение тоже верно. Пусть $f\in I(N)$. Для всякой точки $t\notin N$, $I(t)\ne 0$, то есть $I(t)$ --- ненулевой двусторонний идеал в ${\mathcal B}(X)$, и поскольку ${\mathcal B}(X)$ --- простая алгебра, $I(t)={\mathcal B}(X)$. Поэтому
$$
\forall t\notin N\quad \exists u_t\in I\quad u_t(t)=f(t).
$$
Теперь для всякого компакта $T\subseteq M$ и любого $\e>0$ мы можем подобрать разбиение единицы $\{\eta_t;\ t\in T\}$ на $T$ так, чтобы
$$
\sup_{s\in T}\norm{f(s)-\sum_{t\in T}\eta_t(s)\cdot u_t(s)}<\e.
$$
Суммы $\sum_{t\in T}\eta_t\cdot u_t$ лежат в $I$, и поскольку $I$ --- замкнутый идеал, мы получаем, что $f\in I$.
\epr

\subsection{Идеалы в алгебрах $\prod_{\sigma\in S}A_\sigma$}

\blm\label{LM:prod-A/prod-I}
Пусть $\{A_\sigma;\sigma\in S\}$ -- семейство стереотипных алгебр,
$$
A=\prod_{\sigma\in S}A_\sigma,
$$
--- их произведение, $I$ --- замкнутый двусторонний идеал в $A$, и
$$
I_\sigma=\{x_\sigma;\ x\in I\}
$$
--- его проекции на компоненты $A_\sigma$. Тогда
\bit{

\item[(i)] каждый $I_\sigma$ --- замкнутый двусторонний идеал в $A_\sigma$,

\item[(ii)] $I$ восстанавливается по $I_\sigma$ по формуле
\beq\label{I=prod-I_s}
I=\prod_{\sigma\in S}I_\sigma,
\eeq

\item[(iii)] существует естественный изоморфизм
$$
A/I\cong\prod_{\sigma\in S}(A_\sigma/I_\sigma).
$$
}\eit
\elm
\bpr
Для всякого $\sigma\in S$ определим вложение $\iota_\sigma:A_\sigma\to A$
$$
\iota_\sigma(p)_\tau=\begin{cases}p, & \tau=\sigma\\ 0, & \tau\ne\sigma \end{cases}.
$$
Оно будет мультипликативно, но не унитально, потому что единицу $1_\sigma\in A_\sigma$ оно переводит не в единицу $1\in A$, а в семейство
$$
\iota_\sigma(1_\sigma)_\tau=\begin{cases}1_\sigma, & \tau=\sigma\\ 0, & \tau\ne\sigma \end{cases}.
$$
Заметим, что
\beq\label{I_sigma=iota_sigma^(-1)(I)}
I_\sigma=\iota_\sigma^{-1}(I).
\eeq
Действительно,
\begin{multline*}
p\in I_\sigma\quad\Longleftrightarrow\quad \exists x\in I\quad p=x_\sigma
\quad\Longleftrightarrow\quad \exists x\in I\quad  \iota_\sigma(p)=x\cdot\iota_\sigma(1_\sigma)
\quad\Longleftrightarrow\\ \Longleftrightarrow\quad \iota_\sigma(p)\in I\quad\Longleftrightarrow\quad p\in \iota_\sigma^{-1}(I).
\end{multline*}

1. Из \eqref{I_sigma=iota_sigma^(-1)(I)} сразу следует, что $I_\sigma$ --- замкнутое множество в $A_\sigma$. С другой стороны, для любых $p\in I_\sigma$ и $q\in A_\sigma$ мы получим
$$
\iota_\sigma(p\cdot q)=\underbrace{\iota_\sigma(p)}_{\scriptsize\begin{matrix}\text{\rotatebox{90}{$\owns$}} \\ I\end{matrix}}\cdot \underbrace{\iota_\sigma(q)}_{\scriptsize\begin{matrix}\text{\rotatebox{90}{$\owns$}} \\ A\end{matrix}}\in I
\quad\overset{\eqref{I_sigma=iota_sigma^(-1)(I)}}{\Longrightarrow}\quad p\cdot q\in I_\sigma.
$$
И точно так же $q\cdot p\in I_\sigma$. Это доказывает (i).

2. Формула \eqref{I=prod-I_s} очевидна.

3. Формула
$$
\varPhi\{a_\sigma+I_\sigma;\ \sigma\in S\}=\{a_\sigma;\ \sigma\in S\}+I
$$
корректно определяет отображение
$$
\varPhi:\prod_{\sigma\in S}(A_\sigma/I_\sigma)\to A/I
$$
которое, как нетрудно проверить, является изоморфизмом локально выпуклых пространств.
\epr

\subsection{Аугментированные алгебры}

\paragraph{Категория $\AugInvSteAlg$ аугментированных стереотипных алгебр.}

{\it Аугментация} на инволютивной стереотипной алгебре $A$ определяется как произвольный морфизм $\e:A\to\C$ (то есть непрерывное линейное мультипликативное и сохраняющее единицу и инволюцию отображение). {\it Аугментированной инволютивной стереотипной алгеброй}\label{DEF:augmented-ster-alghebras} называется пара $(A,\e)$, в которой $A$ --- стереотипная алгебра над $\C$, а $\e:A\to\C$ --- аугментация на ней. Задать аугментацию на инволютивной стереотипной алгебре $A$ --- то же, что задать двусторонний замкнутый и инвариантный относительно инволюции идеал $I_A$ в $A$ такой, что $A$ является прямой суммой векторных пространств над $\C$
$$
A=I_A\oplus \C\cdot 1_A,
$$
где $1_A$ -- единица $A$.

Класс всех аугментированных инволютивных стереотипных алгебр мы будем обозначать $\AugInvSteAlg$. Он
образует категорию, в которой морфизмами $\ph :(A,\e_A)\to (B,\e_B)$ являются морфизмы стереотипных алгебр $\ph :A\to B$, сохраняющие инволюцию и аугментацию: 
\beq\label{e_A=e_B-circ-ph}
\e_A=\e_B\circ\ph
\eeq

В категории $\AugInvSteAlg$, алгебра $\C$ с тождественным отображением $\id_{\C}:\C\to\C$ в качестве аугментации является нулем. Как следствие, нулевым морфизмом в $\AugInvSteAlg$ будет всякий морфизм $\ph:A\to B$, представимый в виде композиции $\ph=\iota_B\circ\e_A$, где $\e_A:A\to\C$ --- аугментация на $A$, а $\iota_B:\C\to B$ --- (единственный) морфизм $\C$ в $B$. Таким образом,
\beq\label{0=iota_B-circ-e_A}
0_{A,B}=\iota_B\circ\e_A
\eeq

\bprop\label{PROP:0:A->B}
Для морфизма аугментированных инволютитвных стереотипных алгебр $\ph :(A,\e_A)\to (B,\e_B)$ следующие условия экивалентны:
\bit{
\item[(i)] $\ph:A\to B$ --- нулевой морфизм: 
$$
\ph=0_{A,B};
$$

\item[(ii)] $\ph:A\to B$ пропускается через нулевой объект $\C$:
\beq\label{kappa-circ-lambda=0}
 \xymatrix   @R=2.pc @C=4.pc
{
A\ar[rr]^{\ph}\ar@/_2ex/@{-->}[dr]_{\e_A} & & B \\
& \C\ar@/_2ex/@{-->}[ur]_{\iota_B} &
}
\eeq
(здесь $\e_A$ --- аугментация в $A$, а $\iota_B$ --- вложение $\C$ в $B$ в виде подалгебры $\C\cdot 1_B$, порожденной единицей $1_B$), или, иными словами, выполняется тождество
\beq\label{0:A->B-1}
\ph(x)=\e_A(x)\cdot 1_B, \qquad x\in A;
\eeq

\item[(iii)] множество значений $\ph$ содержится в подалгебре, порожденной единицей $1_B$ алгебры $B$:
\beq\label{0:A->B}
\ph(A)\subseteq\C\cdot 1_B.
\eeq
}\eit
\eprop
\bpr
1. Условия (i) и (ii) эквивалентны просто по определению нулевого морфизма \cite[p.159]{Akbarov-De-Gruyter-I}.

2. Импликация (ii)$\Rightarrow$(iii) очевидна. Докажем обратную: (ii)$\Leftarrow$(iii). Пусть выполненео \eqref{0:A->B}. Тогда для всякого $x\in A$ мы получим:
\beq\label{PROOF:0:A->B-1}
\ph(x)=\lambda\cdot 1_B.
\eeq
Применив к этому равенству аугментацию $\e_B$, мы получим: 
$$
\e_A(x)=\eqref{e_A=e_B-circ-ph}=\e_B(\ph(x))=\lambda\cdot\e_B(1_B)=\lambda.
$$
Теперь \eqref{PROOF:0:A->B-1} превращается в тождество \eqref{0:A->B-1}.
\epr

\bex\label{EX:functional-algebras} {\sf Функциональные алгебры на группах \cite{Akbarov-De-Gruyter-I}.}
Типичный пример аугментированной инволютивной стереотипной алгебры --- алгебра ${\mathcal C}(G)$ непрерывных функций на локально компактной группе $G$ с поточечными алгебраическими операциями и топологией равномерной сходимости на компактах. Что она является стереотипной алгеброй доказано в \cite[Example 10.3]{Ak03}, а инволюция и аугментация в ней описываются правилами
\beq\label{u^bullet(t)=overline{u(t)}}
u^\bullet(t)=\overline{u(t)},\quad \e(u)=u(1_G),\quad u\in {\mathcal C}(G).
\eeq
Этот пример можно дополнить другой классической функциональной алгеброй на группах,  алгеброй ${\mathcal E}(G)$ гладких функций на группе Ли $G$ (с обычной топологией равномерной сходимости на компактах по каждой частной производной \cite[Example 10.4]{Ak03}). Инволюция и аугментация на ${\mathcal E}(G)$  определяются теми же формулами \eqref{u^bullet(t)=overline{u(t)}}.
\eex

\bex\label{EX:group-algebras} {\sf Групповые алгебры \cite{Akbarov-De-Gruyter-I}.}
Переходя к стереотипно сопряженным пространствам мы получим два других примера --- алгебры относительно свертки:
\bit{

\item[---] алгебру ${\mathcal C}^\star(G)$ мер с компактным носителем на локально компактной группе $G$ \cite[Example 10.7]{Ak03}, и

\item[---] алгебру ${\mathcal E}^\star(G)$ распределений с компактными носителем на группе Ли $G$ (с обычной топологией равномерной сходимости на компактах в ${\mathcal E}(G)$ \cite[Example 10.8]{Ak03}).
}\eit
Инволюция на этих алгебрах определяется формулой
 \beq\label{involution-in-C*(G)}
\alpha^\bullet(u)=\overline{\alpha(\widetilde{u}^\bullet)}, \qquad \alpha\in {\mathcal C}^\star(G)\quad ({\mathcal E}^\star(G)),
 \eeq
где $u^\bullet$ --- инволюция функции, определенная равенством \eqref{u^bullet(t)=overline{u(t)}}, а
$$
\widetilde{u}(t)=u(t^{-1})
$$
--- антипод функции. А аугментация в сверточных алгебрах определяется равенством
 \beq\label{augmentation-in-C*(G)}
\e(\alpha)=\alpha(1),\qquad \alpha\in {\mathcal C}^\star(G)\quad ({\mathcal E}^\star(G)),
 \eeq
($1$ --- функция, тождественно равная единице на $G$).
\eex

\paragraph{Ядро и коядро в категории $\AugInvSteAlg$.}

Ниже нам понадобятся следующие две конструкции.
\bit{
\item[$\bullet$] Пусть $B$ --- стереотипная алгебра и $A$ --- замкнутая подалгебра в $B$ (то есть $A$ является унитальной подалгеброй в $B$ в чисто алгебраическом смысле и одновременно замкнутым подпространством в локально выпуклом пространстве $B$). Наделим $A$ топологией, индуцированной из $B$. Тогда псевдонасыщение $A^\vartriangle$ пространства  $A$ является стереотипной алгеброй, называемой {\it замкнутой непосредственной подалгеброй в стереотипной алгебре $B$, порожденной подалгеброй $A$} \cite[Theorem 5.14]{Ak16}.

\item[$\bullet$] Снова пусть $B$ --- стереотипная алгебра и $I$ --- замкнутый двусторонний идеал в $B$ (то есть $I$ является двусторонним идеалом в $B$ в чисто алгебраическом смысле и одновременно замкнутым подпространством в локально выпуклом пространстве $B$). Рассмотрим фактор-пространство $B/I$. Оно является алгеброй в чисто алгебраическом смысле, как фактор-алгебра алгебры $B$ по идеалу $I$, и одновременно $B/I$ будет локально выпуклым пространством с обычной фактор-топологией, унаследованной из $B$. Псевдопополнение $(B/I)^\triangledown$ пространства $B/I$ является стереотипной алгеброй в $B$, называемой {\it открытой непосредственной фактор-алгеброй стереотипной алгебры $B$ по идеалу $I$} \cite[Theorem 5.17]{Ak16}.
}\eit

\btm\label{PROP:AugSteAlg-imeet-ker}
Для морфизма $\ph:(A,\e_A)\to (B,\e_B)$ в категории $\AugInvSteAlg$ аугментированных инволютивных стереотипных алгебр над $\C$
\bit{

\item[---] ядром является замкнутая непосредственная подалгебра в $A$, порожденная прообразом подалгебры $\C\cdot 1_B$ в $B$ при отображении $\ph$:
\beq\label{ker-v-AugSteAlg}
\Ker\ph=\Big(\ph^{-1}(\C\cdot 1_B)\Big)^\vartriangle
\eeq

\item[---] коядром является открытая непосредственная фактор-алгебра алгебры $B$ по замкнутому двустороннему идеалу $I$ в $B$, порожденному образом $\ph(\e_A^{-1}(0))$ идеала $\e_A^{-1}(0)$ в $A$:
\beq\label{coker-v-AugSteAlg}
\Coker\ph=\Big(B/I\Big)^\triangledown
\eeq
}\eit
Как следствие, категория $\AugInvSteAlg$ обладает ядрами и коядрами.
\etm
\bpr
Пусть $\ph:(A,\e_A)\to (B,\e_B)$ -- морфизм аугментированных стереотипных алгебр, и пусть идеалы
$$
I_A=\e_A^{-1}(0),\qquad I_B=\e_B^{-1}(0)
$$
наделены структурой непосредственных подпространств в $A$ и $B$ (в данном случае топологией, являющейся псевдонасыщением топологии, индуцированной из $A$ и $B$, см. детали в \cite{Ak16}).

1. Докажем \eqref{ker-v-AugSteAlg}. Положим $K=\ph^{-1}(\C\cdot 1_B)$, причем снова наделим $K$ структурой непосредственного подпространства в $A$. Это будет стереотипная алгебра и подалгебра в $A$. Определим $\varkappa:K\to A$ как теоретико-множественное вложение, 
\beq\label{varkappa(x)=x}
\varkappa(x)=x,\qquad x\in K,
\eeq
а $\e_K:K\to\C$ --- как ограничение $\e_A$ на $K$. 
\beq\label{e_K=e_A-circ-varkappa}
\e_K=\e_A\circ\varkappa.
\eeq
Из \eqref{e_K=e_A-circ-varkappa} будет следовать, что композиция $\ph\circ\varkappa$ --- морфизм в $\AugInvSteAlg$.

a) Заметим сначала, что $\varkappa$ обнуляет $\ph$ изнутри, то есть что $\ph\circ\varkappa$ --- нулевой морфизм в $\AugInvSteAlg$: 
\beq\label{PROOF:AugSteAlg-imeet-ker-0}
\ph\circ\varkappa=0_{K,B}
\eeq
По предложению \ref{PROP:0:A->B}, это эквивалентно условию
\beq\label{PROOF:AugSteAlg-imeet-ker-1}
(\ph\circ\varkappa)(K)\subseteq \C\cdot 1_B,
\eeq
которое очевидно, потому что если $x\in K=\ph^{-1}(\C\cdot 1_B)$, то $\ph(\varkappa(x))=\ph(x)\in \C\cdot 1_B$.

b) Пусть далее $\nu:N\to A$ --- какой-нибудь другой морфизм, обнуляющий $\ph$ изнутри:
\beq\label{PROOF:AugSteAlg-imeet-ker-2}
\ph\circ\nu=0_{K,B}
\eeq
то есть по предложению \ref{PROP:0:A->B}, удовлетворяющий условию
\beq\label{PROOF:AugSteAlg-imeet-ker-3}
(\ph\circ\nu)(N)\subseteq \C\cdot 1_B.
\eeq
Это значит, что 
\beq\label{PROOF:AugSteAlg-imeet-ker-4}
\nu(N)\subseteq \ph^{-1}(\C\cdot 1_B)=K.
\eeq
И поскольку $K$ --- непосредственное подпространство, $\nu$ поднимается до линейного отображения $\nu':N\to K$:
\beq\label{PROOF:AugSteAlg-imeet-ker-5}
 \xymatrix   @R=2.pc @C=4.pc
{
K\ar[r]^{\varkappa} & A\ar[r]^{\ph} & B \\
& N\ar[ur]_{0_{N,B}}\ar[u]^{\nu}\ar@{-->}[ul]^{\nu'} &
}
\eeq
Это отображение автоматически мультипликативно и сохраняет единицу, потому что $\nu$ обладает этими свойствами. Кроме того, оно сохраняет аугментацию, потому что
$$
\e_K\circ\nu'=\eqref{e_K=e_A-circ-varkappa}=\e_A\circ\varkappa\circ\nu'=\e_A\circ\nu=\e_N.
$$

2. Докажем \eqref{coker-v-AugSteAlg}. Пусть $I$ --- замкнутый непосредственный двусторонний идеал в $B$, порожденный множеством $\ph(I_A)$ (то есть $I$, как множество, совпадает с замыканием в $B$ двустороннего идеала, порожденного множеством $\ph(I_A)$, и наделяется топологией псевдонасыщения топологии, индуцированной из $B$). Из равенства \eqref{e_A=e_B-circ-ph}
$$
\e_B\circ\ph=\e_A
$$
следует
$$
\ph(I_A)\subseteq I_B,
$$
поэтому
$$
\ph(I_A)\subseteq I\subseteq I_B.
$$
Положим $C=B/I=(B/I)^\triangledown$, и пусть $\gamma:B\to B/I=C$ -- фактор-отображение. Из условия $I\subseteq I_B=\e_B^{-1}(0)$ следует, что функционал $\e_B$ можно продолжить до некоторого функционала $\e_C$ на фактор-пространство $(B/I)^\triangledown=B/I=C$:
$$
 \xymatrix @R=2.pc @C=1.pc %
 {
 B\ar[rr]^{\gamma}\ar[dr]_{\e_B}& & C\ar[dl]^{\e_C} \\
 & \C &
 }
$$
Поэтому $\gamma$ является морфизмом аугментированных стереотипных алгебр $(B,\e_B)\to(C,\e_C)$. Покажем, что это и будет коядром морфизма $\ph$.

а) Покажем, что $\gamma$ обнуляет $\ph$ снаружи:
$$
\gamma\circ\ph=0_{A,C},
$$
то есть, по предложению \ref{PROP:0:A->B}, 
\beq\label{PROOF:AugSteAlg-imeet-ker-6}
(\gamma\circ\ph)(A)\subseteq \C\cdot 1_C.
\eeq
Действительно, если взять $a\in A$, то положив $x=a-\e_A(a)\cdot 1_A$ мы получим $x\in I_A$, откуда $\ph(x)\in I$, и значит   $\gamma(\ph(x))=0$, и поэтому,
$$
(\gamma\circ\ph)(a)=(\gamma\circ\ph)(\e_A(a)\cdot 1_A+x)=\e_A(a)\cdot \gamma(\ph(1_A))+\gamma(\ph(x))=\e_A(a)\cdot 1_C.
$$

б) Пусть теперь $\delta:(B,\e_B)\to(D,\e_D)$ --- другой морфизм, обнуляющий $\ph$ снаружи, 
$$
\delta\circ\ph=0_{A,D},
$$
то есть, по предложению \ref{PROP:0:A->B}, 
\beq\label{PROOF:AugSteAlg-imeet-ker-7}
\delta\circ\ph=\iota_D\circ\e_A
\eeq
Тогда мы получим
$$
\forall x\in I_A=\e_A^{-1}(0) \qquad \delta(\ph(x))=\iota_D(\e_A(x))=0
$$
$$
\Downarrow
$$
$$
\delta(\ph(I_A))=0
$$
$$
\Downarrow
$$
$$
\ph(I_A)\subseteq \delta^{-1}(0)
$$
$$
\Downarrow
$$
$$
I\subseteq \delta^{-1}(0)
$$
Последнее означает, что морфизм $\delta$ пропускается через $\gamma$: для некоторого морфизма $\lambda$ коммутативна диаграмма
 \beq\label{DIAGR:coker}
 \xymatrix @R=2.pc @C=4.pc
 {
 A\ar[r]^{\ph}& B \ar[r]^{\gamma}\ar[d]_{\delta}& C\ar@{-->}[dl]^{\lambda} \\
 & D &
 }
 \eeq
Морфизм $\lambda$ сохраняет аугментацию, потому что 
$$
\e_D\circ\lambda\circ\gamma=\e_D\circ\delta=\e_B=\e_C\circ\gamma
$$ 
и, поскольку $\gamma$ --- эпиморфизм, его можно отбросить: 
$$
\e_D\circ\lambda=\e_C.
$$ 
 
\epr

\section{О перестановочности сумм и произведений с тензорными произведениями $\circledast$ и $\odot$}

\bit{

\item[$\bullet$] Условимся говорить, что стереотипное пространство $X$ {\it не имеет бесконечных произведений в качестве подпространств}, если не существует мономорфизмов вида 
    \beq
    \sigma:\C^{\N}\to X.
    \eeq 
    Это эквивалентно тому что в $X$ не может существовать последовательность ненулевых векторов $0\ne x_n\in X$ такая, что для любой последовательности чисел $\lambda_n\in \C$ ряд
    $$
    \sum_{n=1}^\infty \lambda_n\cdot x_n
    $$
    сходится в $X$ безусловно.
 
\item[$\bullet$] Условимся говорить, что стереотипное пространство $X$ {\it не имеет бесконечных сумм в качестве фактор-пространств}, если не существует эпиморфизмов вида 
    \beq
    \pi:X\to\C_{\N}.
    \eeq 

}\eit

\blm\label{LM:Brauner-ne-soderzhit-besk-proizv}
Никакое пространство Фреше $X$ не может иметь бесконечную сумму в качестве фактор-пространства. Двойственным образом, никакое  пространство Браунера $Y$ не может иметь бесконечное произведение в качестве подпространства.
\elm
\bpr
Здесь достаточно доказать второе утверждение. Пусть $Y$ --- пространство Браунера с фундаментальной системой компактов $\{K_n;\ n\in\N\}$. Предположим, что в $Y$ имеется последовательность ненулевых векторов $0\ne y_n\in Y$ такая, что для любой последовательности чисел $\lambda_n\in \C$ ряд
\beq\label{Brauner-ne-soderzhit-besk-proizv-0}
    \sum_{n=1}^\infty \lambda_n\cdot y_n
\eeq
сходится в $Y$ безусловно. Построим последовательность скаляров $\lambda_n\in \C$ следуюшим образом. 

1) Выберем число $\lambda_1\in\C$ так, чтобы
$$
\lambda_1\cdot y_1\notin K_1.
$$   
Это всегда можно сделать, потому что $K_1$ --- компакт, и поэтому он не может содержать одномерное подпространство $\C\cdot y_1$.

2) Выберем число $\lambda_2\in\C$ так, чтобы
$$
\lambda_1\cdot y_1+\lambda_2\cdot y_2\notin K_2.
$$   
Это всегда можно сделать, потому что $K_2-\lambda_1\cdot y_1$ --- компакт, и поэтому он не может содержать одномерное подпространство $\C\cdot y_2$.

n) На $n$-м шаге выберем число $\lambda_n\in\C$ так, чтобы
$$
\sum_{i=1}^{n-1}\lambda_i\cdot y_i+\lambda_n\cdot y_n\notin K_n.
$$   
Это всегда можно сделать, потому что $K_n-\sum_{i=1}^{n-1}\lambda_i\cdot y_i$ --- компакт, и поэтому он не может содержать одномерное подпространство $\C\cdot y_n$.

В результате мы получим последовательность чисел $\lambda_n\in \C$ со свойством
\beq\label{Brauner-ne-soderzhit-besk-proizv-1}
\sum_{i=1}^n\lambda_i\cdot y_i\notin K_n,\qquad n\in\N.
\eeq
Но с другой стороны, ряд \eqref{Brauner-ne-soderzhit-besk-proizv-0} должен сходиться в $Y$, и потому его частичные суммы 
$$
\sum_{i=1}^n\lambda_i\cdot y_i
$$
должны образовывать вполне ограниченную последовательность в $Y$. Значит, они должны содержаться в некотором компакте $K_N$:
\beq\label{Brauner-ne-soderzhit-besk-proizv-2}
\sum_{i=1}^n\lambda_i\cdot y_i\notin K_N,\qquad n\in\N.
\eeq
Условия \eqref{Brauner-ne-soderzhit-besk-proizv-1} и \eqref{Brauner-ne-soderzhit-besk-proizv-2} противоречат друг другу.
\epr

\blm\label{LM:Ker-ph-supseteq-prod}
Если $\{X_i;\ i\in I\}$ --- семейство пространств Фреше, а $Y$ --- пространство Браунера, то всякий морфизм $\ph:\prod_{i\in I}X_i\to Y$ пропускается через конечное произведение пространств $X_i$:
\beq\label{Ker-ph-supseteq-prod}
 \xymatrix @R=2.pc @C=4.pc
 {
 \prod_{i\in I}X_i\ar[rd]_{\pi} \ar[rr]^{\ph}& & Y  \\
 & \prod_{k=1}^l X_{i_k}\ar@{-->}[ru]_{\psi}  &
 }
\eeq
(здесь $l\in\N$, а $\pi$ --- естественная проекция).
\elm
\bpr
Здесь нужно убедиться, что почти все пространства $X_i$ (то есть все, кроме конечного набора) переводятся отображением $\ph$  в нуль. Предположим, что это не так, то есть что существует последовательность неповторяющихся индексов $i_k$, $k\in\N$, такая, что
$$
\forall k\in\N\quad \ph(X_{i_k})\ne 0.
$$ 
Выберем векторы $x_k\in X_{i_k}$ так, чтобы 
$$
\forall k\in\N\quad \ph(x_k)\ne 0.
$$ 
Полскольку индексы $i_k$ не повторяются, векторы $x_k$ лежат в разных пространствах $X_{i_k}$:
$$
\nexists k,l\in\N \qquad x_k\in X_{i_k}\ \&\ x_l\in X_{i_k}
$$
Отсюда следует, что для любой последовательности чисел $\lambda_k\in\C$ в пространстве $\prod_{i\in I}X_i$ сходится ряд
$$
\sum_{k=1}^\infty \lambda_k\cdot x_k\in \prod_{i\in I}X_i.
$$
Как следствие ряд
$$
\sum_{k=1}^\infty \lambda_k\cdot \ph(x_k)
$$
тоже сходится в пространстве $Y$ для любой последовательности чисел $\lambda_k\in\C$. Но это означает, что пространство Браунера $Y$ имеет бесконечное произведение в качестве подпространства, что невозможно по лемме \ref{LM:Brauner-ne-soderzhit-besk-proizv}.
\epr

Двойственное к лемме \ref{LM:Ker-ph-supseteq-prod} утверждение выглядит так:

\blm\label{LM:Im-ph-subseteq-bigoplus}
Если $X$ --- пространство Фреше, а $\{Y_j;\ j\in J\}$ --- семейство пространств Браунера, то всякий морфизм $\ph:X\to\bigoplus_{j\in J}Y_j$ пропускается через конечную сумму пространств $Y_j$:
\beq\label{Im-ph-subseteq-bigoplus}
 \xymatrix @R=2.pc @C=4.pc
 {
 X\ar@{-->}[dr]_{\psi}\ar[rr]^{\ph}& & \bigoplus_{j\in J}Y_j \\
 & \bigoplus_{m=1}^n Y_{j_m}\ar[ru]_{\sigma} &
 }
 \eeq
(здесь $n\in\N$, а $\sigma$ --- естественная инъекция).
\elm

Леммы \ref{LM:Ker-ph-supseteq-prod} и \ref{LM:Im-ph-subseteq-bigoplus} используются в доказательстве следующего важного утверждения:

\btm\label{TM:prod_X_i->bigoplus_Y_j}
Если $\{X_i;\ i\in I\}$ --- семейство пространств Фреше,  а $\{Y_j;\ j\in J\}$ --- семейство пространств Браунера, то 
\bit{

\item[(i)] всякий морфизм $\ph:\prod_{i\in I}X_i\to\bigoplus_{j\in J}Y_j$ пропускается через некоторое конечное произведение пространств $X_i$ и некоторую конечную сумму пространств $Y_j$:
\beq\label{prod_X_i->bigoplus_Y_j}
 \xymatrix @R=2.pc @C=4.pc
 {
\prod_{i\in I}X_i\ar[d]_{\pi} \ar[r]^{\ph}&  \bigoplus_{j\in J}Y_j \\
\prod_{k=1}^l X_{i_k} \ar@{-->}[r]_{\psi} & \bigoplus_{m=1}^n Y_{j_m}\ar[u]_{\sigma} 
 }
 \eeq
(здесь $l,n\in\N$, $\pi$ --- естественная проекция, а $\sigma$ --- естественная инъекция),

\item[(ii)] справедлив естественный изоморфизм
\beq\label{bigoplus_Y_j-oslash-prod_X_i}
\Big(\bigoplus_{j\in J}Y_j\Big)\oslash\Big(\prod_{i\in I}X_i\Big)\cong \bigoplus_{j\in J,i\in I}(Y_j\oslash X_i).
\eeq

}\eit

\etm
\bpr
i. Обозначим
$$
X=\prod_{i\in I}X_i
$$
и рассмотрим множество $J_0$ индексов $j\in J$, у которых пространства $Y_j$ имеют нетривиальное пересечение с множеством значений $\ph(X)$:  
\beq\label{PROOF:prod_X_i->bigoplus_Y_j-0}
J_0=\{j\in J:\ \ph(X)\cap Y_j\ne 0\}.
\eeq

1. Покажем, что это множество конечно. Для этого предположим, что оно наоборот, бесконечно:
\beq\label{PROOF:prod_X_i->bigoplus_Y_j-1}
\card J_0=\infty
\eeq
Организуем следующую индукцию.

1.1. Зафиксируем произвольный индекс $j_1\in J_0$. Из неравенства $\ph(X)\cap Y_{j_1}\ne 0$ следует, что для какого-то $i_1\in I$ должно выполняться неравенство
$$
\ph(X_{i_1})\cap Y_{j_1}\ne 0
$$
(потому что иначе мы бы получили импликацию  $(\forall i\quad \ph(X_{i})\cap Y_{j_1}= 0)\ \Rightarrow\ \ph(X)\cap Y_{j_1}=0$). Запомним эти индексы $j_1\in J_0$ и $i_1\in I$.

1.2. Из леммы \ref{LM:Im-ph-subseteq-bigoplus} следует, что ограничение 
$$
\ph\Big|_{X_{i_1}}:X_{i_1}\to\bigoplus_{j\in J}Y_j
$$
имеет образом подпространство в $\bigoplus_{j\in J}Y_j$, нетривиально пересекающееся только с конечным набором пространств $Y_j$:
$$
\card \{j\in J:\ \ph(X_{i_1})\cap Y_j\ne 0\}<\infty.
$$ 
Поэтому в (бесконечном) множестве $J_0\setminus\{j_1\}$ найдется индекс $j_2\in J_0\setminus\{j_1\}$ такой, что 
\beq\label{PROOF:prod_X_i->bigoplus_Y_j-2}
\ph(X_{i_1})\cap Y_{j_2}=0.
\eeq
При этом из неравенства $\ph(X)\cap Y_{j_2}\ne 0$ будет следовать, что найдется индекс $i_2\in I$ такой, что 
\beq\label{PROOF:prod_X_i->bigoplus_Y_j-3}
\ph(X_{i_2})\cap Y_{j_2}\ne 0.
\eeq
А вместе с \eqref{PROOF:prod_X_i->bigoplus_Y_j-2} это будет означать, что
\beq\label{PROOF:prod_X_i->bigoplus_Y_j-4}
i_2\ne i_1.
\eeq
Запомним эти индексы $j_2\in J_0$ и $i_2\in I$.

1.3. Из леммы \ref{LM:Im-ph-subseteq-bigoplus} следует, что ограничение 
$$
\ph\Big|_{X_{i_1}\times X_{i_2}}:X_{i_1}\times X_{i_2}\to\bigoplus_{j\in J}Y_j
$$
имеет образом подпространство в $\bigoplus_{j\in J}Y_j$, нетривиально пересекающееся только с конечным набором пространств $Y_j$:
$$
\card \{j\in J:\ \ph(X_{i_1}\times X_{i_2})\cap Y_j\ne 0\}<\infty.
$$ 
Поэтому в (бесконечном) множестве $J_0\setminus\{j_1,j_2\}$ найдется индекс $j_3\in J_0\setminus\{j_1,j_2\}$ такой, что 
\beq\label{PROOF:prod_X_i->bigoplus_Y_j-5}
\ph(X_{i_1}\times X_{i_2})\cap Y_{j_3}=0.
\eeq
При этом из неравенства $\ph(X)\cap Y_{j_3}\ne 0$ будет следовать, что найдется индекс $i_3\in I$ такой, что 
\beq\label{PROOF:prod_X_i->bigoplus_Y_j-6}
\ph(X_{i_3})\cap Y_{j_3}\ne 0.
\eeq
А вместе с \eqref{PROOF:prod_X_i->bigoplus_Y_j-5} это будет означать, что
\beq\label{PROOF:prod_X_i->bigoplus_Y_j-7}
i_3\notin \{i_1,i_2\}.
\eeq
Запомним эти индексы $j_3\in J_0$ и $i_3\in I$.

1.$k+1$. На $k+1$-м шаге мы замечаем, что из леммы \ref{LM:Im-ph-subseteq-bigoplus} следует, что ограничение 
$$
\ph\Big|_{X_{i_1}\times...\times X_{i_k}}:X_{i_1}\times...\times X_{i_k}\to\bigoplus_{j\in J}Y_j
$$
имеет образом подпространство в $\bigoplus_{j\in J}Y_j$, нетривиально пересекающееся только с конечным набором пространств $Y_j$:
$$
\card \{j\in J:\ \ph(X_{i_1}\times ...\times X_{i_k})\cap Y_j\ne 0\}<\infty.
$$ 
Поэтому в (бесконечном) множестве $J_0\setminus\{j_1,...,j_k\}$ найдется индекс $j_{k+1}\in J_0\setminus\{j_1,...,j_k\}$ такой, что 
\beq\label{PROOF:prod_X_i->bigoplus_Y_j-8}
\ph(X_{i_1}\times ...\times X_{i_k})\cap Y_{j_{k+1}}=0.
\eeq
При этом из неравенства $\ph(X)\cap Y_{j_{k+1}}\ne 0$ будет следовать, что найдется индекс $i_{k+1}\in I$ такой, что 
\beq\label{PROOF:prod_X_i->bigoplus_Y_j-9}
\ph(X_{i_{k+1}})\cap Y_{j_{k+1}}\ne 0.
\eeq
А вместе с \eqref{PROOF:prod_X_i->bigoplus_Y_j-8} это будет означать, что
\beq\label{PROOF:prod_X_i->bigoplus_Y_j-10}
i_{k+1}\notin \{i_1,...,i_k\}.
\eeq
Запомним эти индексы $j_{k+1}\in J_0$ и $i_{k+1}\in I$.

Продолжая эту процедуру, мы получим последовательности индексов $i_k\in I$ и $j_k\in J$ со следующими свойствами:
\bit{

\item[a)] они инъективны: 
\beq\label{PROOF:prod_X_i->bigoplus_Y_j-11}
k\ne l \quad \Rightarrow\quad ( i_k\ne i_l \quad \&\quad j_k\ne j_l ),
\eeq

\item[b)] образ $X_{i_k}$ нетривиально пересекается с $Y_{j_k}$:
\beq\label{PROOF:prod_X_i->bigoplus_Y_j-12}
\ph(X_{i_k})\cap Y_{j_k}\ne 0.
\eeq

}\eit

Из свойства \eqref{PROOF:prod_X_i->bigoplus_Y_j-12} следует, что можно выбрать последовательность
\beq\label{PROOF:prod_X_i->bigoplus_Y_j-13}
x_k\in X_{i_k}
\eeq 
так, чтобы  
\beq\label{PROOF:prod_X_i->bigoplus_Y_j-14}
0\ne \ph(x_k)\in Y_{j_k}.
\eeq
При этом из свойства \eqref{PROOF:prod_X_i->bigoplus_Y_j-11} будет следовать, что для любой числовой последовательности $\lambda_k\in\C$ в пространстве $\prod_{i\in I}X_i$ сходится ряд
\beq\label{PROOF:prod_X_i->bigoplus_Y_j-15}
\sum_{k=1}^\infty \lambda_k\cdot x_k
\eeq
и поэтому в пространстве $\bigoplus_{j\in J}Y_j$ должен сходиться ряд 
\beq\label{PROOF:prod_X_i->bigoplus_Y_j-16}
\sum_{k=1}^\infty \lambda_k\cdot \ph(x_k).
\eeq
И, поскольку $\ph(x_k)\in Y_{j_k}$, мы получаем, что этот ряд сходится в пространстве $\bigoplus_{k=1}^\infty Y_{j_k}$.
И при этом $\ph(x_k)\ne 0$. Такое невозможно, потому что элемент прямой суммы $\bigoplus_{k=1}^\infty Y_{j_k}$ может иметь ненулевую проекцию только на конечное число пространств $Y_{j_k}$.  

2. Мы доказали, что множество $J_0$ из \eqref{PROOF:prod_X_i->bigoplus_Y_j-0} конечно. Из этого следует, что 
$\ph:\prod_{i\in I}X_i\to\bigoplus_{j\in J}Y_j$ пропускается через конечную сумму пространств $Y_j$:
\beq\label{PROOF:prod_X_i->bigoplus_Y_j-17}
 \xymatrix @R=2.pc @C=4.pc
 {
\prod_{i\in I}X_i\ar@{-->}[rd]_{\chi} \ar[r]^{\ph}&  \bigoplus_{j\in J}Y_j \\
  & \bigoplus_{m=1}^n Y_{j_m}\ar[u]_{\sigma} 
 }
 \eeq
(здесь $n\in\N$, а $\sigma$ --- естественная инъекция).

3. Теперь земетим, что пространство $\bigoplus_{m=1}^n Y_{j_m}$ --- пространство Браунера. 
А $\{X_i;\ i\in I\}$ --- по-прежнему, семейство пространств Фреше. Поэтому по лемме \ref{LM:Ker-ph-supseteq-prod}, морфизм $\chi: \prod_{i\in I}X_i\to \bigoplus_{m=1}^n Y_{j_m}$ пропускается через  конечное произведение пространств $X_i$:
\beq\label{PROOF:prod_X_i->bigoplus_Y_j-18}
 \xymatrix @R=2.pc @C=4.pc
 {
\prod_{i\in I}X_i\ar[d]_{\pi}\ar[rd]^{\chi} &   \\
\prod_{k=1}^l X_{i_k} \ar@{-->}[r]_{\psi} & \bigoplus_{m=1}^n Y_{j_m}
 }
 \eeq
(здесь $l\in\N$, а $\pi$ --- естественная проекция).
Диаграмма \eqref{PROOF:prod_X_i->bigoplus_Y_j-17} и \eqref{PROOF:prod_X_i->bigoplus_Y_j-18} вместе дают \eqref{prod_X_i->bigoplus_Y_j}.

ii. Теперь докажем формулу \eqref{bigoplus_Y_j-oslash-prod_X_i}. Сначала заметим, что она тривиально верна, когда индексные множества $I$ и $J$ конечны. В этом случае
\begin{multline}\label{PROOF:prod_X_i->bigoplus_Y_j-19}
\Big(\bigoplus_{j\in J_0}Y_j\Big)\oslash\Big(\prod_{i\in I_0}X_i\Big)
\cong \Big(\prod_{j\in J_0}Y_j\Big)\oslash\Big(\bigoplus_{i\in I_0}X_i\Big)
\cong \prod_{j\in J_0,i\in I_0}(Y_j\oslash X_i)
\cong \bigoplus_{j\in J_0,i\in I_0}(Y_j\oslash X_i),\\ \qquad \card I_0<\infty, \ \card J_0<\infty.
\end{multline}
Далее зафиксируем произвольные $I$ и $J$ и выберем произвольные конечные подмножества $I_0\subseteq I$ и $J_0\subseteq J$. 
Пусть 
$$
\pi_{I_0}:\prod_{i\in I}X_i\to \prod_{i\in I_0}X_i,\qquad \sigma_{J_0}:\bigoplus_{j\in J_0}Y_j\to \bigoplus_{j\in J}Y_j
$$
--- соответствующие проекция и вложение. Они порождают естественный морфизм пространств
\beq\label{PROOF:prod_X_i->bigoplus_Y_j-20}
\sigma_{J_0}\oslash\pi_{I_0}:\bigoplus_{j\in J_0}Y_j\oslash\prod_{i\in I_0}X_i\to \bigoplus_{j\in J}Y_j\oslash\prod_{i\in I_0}X_i
\eeq
действующий по формуле
$$
(\sigma_{J_0}\oslash\pi_{I_0})(\psi)=\sigma_{J_0}\circ\psi\circ\pi_{I_0},\qquad \psi\in \bigoplus_{j\in J_0}Y_j\oslash\prod_{i\in I}X_i
$$
или, что то же самое, его действие описывается диаграммой
\beq\label{PROOF:prod_X_i->bigoplus_Y_j-21}
 \xymatrix @R=2.pc @C=8.pc
 {
\prod_{i\in I}X_i\ar[d]_{\pi_{I_0}} \ar@{-->}[r]^{(\sigma_{J_0}\oslash\pi_{I_0})(\psi)}&  \bigoplus_{j\in J}Y_j \\
\prod_{i\in I_0}X_i \ar[r]_{\psi} & \bigoplus_{j\in J_0}Y_j\ar[u]_{\sigma_{J_0}} 
 }
 \eeq
 Если в левой части \eqref{PROOF:prod_X_i->bigoplus_Y_j-20} перейти к инъективному пределу, то мы получим морфизм
\beq\label{PROOF:prod_X_i->bigoplus_Y_j-22}
\omega=\injlim_{I_0\to I, \ J_0\to J}\sigma_{J_0}\oslash\pi_{I_0}: \injlim_{I_0\to I, \ J_0\to J}\bigoplus_{j\in J_0}Y_j\oslash\prod_{i\in I_0}X_i\to \bigoplus_{j\in J}Y_j\oslash\prod_{i\in I}X_i
\eeq
Заметим, что 
\beq\label{PROOF:prod_X_i->bigoplus_Y_j-24}
\injlim_{I_0\to I, \ J_0\to J}\bigoplus_{j\in J_0}Y_j\oslash\prod_{i\in I_0}X_i
\cong \eqref{PROOF:prod_X_i->bigoplus_Y_j-19} \cong 
\injlim_{I_0\to I, \ J_0\to J}\bigoplus_{j\in J_0,\ i\in I_0}\Big(Y_j\oslash X_i\Big) 
\cong 
\bigoplus_{j\in J,\ i\in I}\Big(Y_j\oslash X_i\Big)
\eeq
Поэтому морфизм \eqref{PROOF:prod_X_i->bigoplus_Y_j-22} можно понимать как морфизм 
\beq\label{PROOF:prod_X_i->bigoplus_Y_j-25}
\omega=\injlim_{I_0\to I, \ J_0\to J}\sigma_{J_0}\oslash\pi_{I_0}: \bigoplus_{j\in J,\ i\in I}\Big(Y_j\oslash X_i\Big) 
\to \bigoplus_{j\in J}Y_j\oslash\prod_{i\in I}X_i
\eeq

Далее, из уже доказанного свойства (i) следует, что морфизм $\omega$  является сюръекцией. 

Покажем, этот морфизм является инъекцией. Действительно, морфизмы, образующие ковариантную  систему инъективного предела
\beq\label{PROOF:prod_X_i->bigoplus_Y_j-26}
\bigoplus_{j\in J_0}Y_j\oslash\prod_{i\in I_0}X_i\to \bigoplus_{j\in J_1}Y_j\oslash\prod_{i\in I_1}X_i,
\qquad I_0\subseteq I_1\subseteq I,\ J_0\subseteq J_1\subseteq J,
\eeq
по формуле \eqref{PROOF:prod_X_i->bigoplus_Y_j-19} представляют собой естестввенные вложения прямых сумм
\beq\label{PROOF:prod_X_i->bigoplus_Y_j-27}
\bigoplus_{j\in J_0,\ i\in I_0}Y_j\oslash X_i\to \bigoplus_{j\in J_1,\ i\in I_1}Y_j\oslash X_i,
\qquad I_0\subseteq I_1\subseteq I,\ J_0\subseteq J_1\subseteq J,
\eeq
и поэтому они являются инъекциями. Отсюда следует, что ненулевой элемент в 
$$
0\ne \psi\in \injlim_{I_0\to I, \ J_0\to J}\bigoplus_{j\in J_0}Y_j\oslash\prod_{i\in I_0}X_i=\bigoplus_{j\in J,\ i\in I}\Big(Y_j\oslash X_i\Big)
$$
--- это просто ненулевой морфизм из какого-то пространства $\bigoplus_{j\in J_0}Y_j\oslash\prod_{i\in I_0}X_i$:
$$
0\ne \psi\in \bigoplus_{j\in J_0}Y_j\oslash\prod_{i\in I_0}X_i
$$
и он будет оставаться ненулевым во всех более широких пространствах $\bigoplus_{j\in J_1}Y_j\oslash\prod_{i\in I_1}X_i$, где $I_0\subseteq I_1$ и $J_0\subseteq J_1$. Из диаграммы \eqref{PROOF:prod_X_i->bigoplus_Y_j-21} видно, что под действием отображения 
$\sigma_{J_0}\oslash\pi_{I_0}$ он превратится в ненулевой элемент
$$
(\sigma_{J_0}\oslash\pi_{I_0})(\psi)\ne 0.
$$
и это значение не будет меняться если мы будем менять индексы  $I_0$ и $J_0$ на более широкие индексы $I_1\supseteq I_0$ и $J_1\supseteq J_0$. Это значит, что значение морфизма \eqref{PROOF:prod_X_i->bigoplus_Y_j-22} на элементе $\psi$ будет ненулевым:
$$
\omega(\psi)\ne 0.
$$

Остается проверить открытость отображения $\omega$. Пусть $U$ --- базисная окрестность нуля в 
$$
\injlim_{I_0\to I, \ J_0\to J}\bigoplus_{j\in J_0}Y_j\oslash\prod_{i\in I_0}X_i \cong 
\bigoplus_{j\in J,\ i\in I}\Big(Y_j\oslash X_i\Big)
$$
Это значит, что $U$ представляет собой прямую сумму
$$
U=\bigoplus_{j\in J,\ i\in I} V_{i,j}
$$
некоторого семейства окрестностей нуля 
$$
V_{i,j}\subseteq Y_j\oslash X_i
$$
Поскольку здесь $X_i$ --- пространство Фреше, а $Y_j$ --- пространство Браунера, мы можем считать окрестности $V_{i,j}$ имеющими вид дробей:
$$
V_{i,j}=W_j\oslash K_i
$$
где $W_j$ --- окрестность нуля в $Y_j$, а $K_i$ --- компакт в $X_i$. То есть
$$
U=\bigoplus_{j\in J,\ i\in I} W_j\oslash K_i.
$$
Понятно, что под действием морфизма $\omega$ такая окрестность првратится в окрестность нуля
$$
\omega(U)=\bigoplus_{j\in J} W_j\oslash \prod_{i\in I}K_i.
$$
\epr

\btm\label{TH:prod_X_i-circledast-prod_Y_j}
Пусть $X_i$ и $Y_j$ --- пространства Фреше. Тогда
\beq\label{prod_X_i-circledast-prod_Y_j}
\l\prod_{i\in I}X_i\r\circledast \l\prod_{j\in J}Y_j\r\cong\prod_{i\in I,\ j\in J} \Big(X_i\circledast Y_j\Big)
\eeq
\etm
\bpr
\begin{multline*}
\l\prod_{i\in I}X_i\r\circledast \l\prod_{j\in J}Y_j\r=\l\Big(\prod_{i\in I}X_i\Big)^\star\oslash\Big(\prod_{j\in J}Y_j\Big)\r^\star=\l\Big(\bigoplus_{i\in I}X_i^\star\Big)\oslash\Big(\prod_{j\in J}Y_j\Big)\r^\star\cong \eqref{bigoplus_Y_j-oslash-prod_X_i}\cong \\ \cong
\l\bigoplus_{i\in I,j\in J}\Big( X_i^\star\oslash Y_j\Big)\r^\star\cong 
\prod_{i\in I,j\in J}\Big( X_i^\star\oslash Y_j\Big)^\star\cong 
\prod_{i\in I,j\in J}\Big( X_i\circledast Y_j\Big)
\end{multline*}
\epr

Двойственное утверждение выглядит так:

\btm\label{TH:bigoplus_X_i-odot-bigoplus_Y_j}
Пусть $X_i$ и $Y_j$ --- пространства Браунера. Тогда
\beq\label{bigoplus_X_i-odot-bigoplus_Y_j}
\l\bigoplus_{i\in I}X_i\r\odot \l\bigoplus_{j\in J}Y_j\r\cong\bigoplus_{i\in I,\ j\in J} \Big(X_i\odot Y_j\Big)
\eeq
\etm

\section{Исправления в книге ``Стереотипные пространства и алгебры''}

После выхода моей предыдущей монографии \cite{Akbarov-De-Gruyter-I} я заметил в ней некоторое количество ошибок и пробелов. Я думаю, что будет правильно исправить их здесь, потому что неизвестно когда у меня еще будет эта возможность.

\paragraph{Мелкие исправления.}

\begin{center}
\begin{tabular}{|c|p{2.5in}|p{2.5in}|} \hline
Page & Fragment & Must be replaced by \\ \hline
2 & No other mathematical discipline is separated by such a powerful barrier from the rest of mathematics and does not require & No other mathematical discipline is separated by such a powerful barrier from the rest of mathematics and does require \\ \hline
7 & These constructions automatically give a solution to problem 0.1.1. & These constructions automatically give a solution to Exercise 0.1.1. \\ \hline
8 & consists of the fact that on the dual space & consists of the convention that on the dual space \\ \hline
10 & with their own ideas about right and wrong & with their own understanding of right and wrong \\ \hline
10 &	leaved the university (Berkeley) & left the university (Berkeley) \\ \hline
11 & all Cauchy nets in these spaces converge & all totally bounded Cauchy nets in these spaces converge \\ \hline
26 & also have stereotype approximations as well & 	also have stereotype approximation  \\ \hline
26 & ${\mathcal L}(X)$ also has a stereotype approximation &	${\mathcal L}(X)$ also has stereotype approximation  \\ \hline
37 & [13, Theorems 6.4 7.2] & [21, Theorem 4.3] \\ \hline
46 & geometric discipline is a category	a geometric & discipline is a category \\ \hline
46 & (respectively, in $\Env^\varOmega_\varPhi\SteAlg$)	& (respectively, in $\Env^\varOmega_\varPhi\InvSteAlg$) \\ \hline
57 & if $\ph,\chi,\psi\in{\tt{K}}$,  $\psi\circ\chi$ and $\chi\circ\ph$ are defined, and	& 
if $\ph,\chi,\psi\in{\tt{K}}$,  and $\psi\circ\chi$, $\chi\circ\ph$ are defined,  \\ \hline
60 & It is clear that this condition is equivalent to specifying the structure of a directed graph & It is clear that this condition is equivalent to specifying the structure of a reflexive transitive directed graph  \\ \hline
60 & and therefore there must be $1_A$	& and therefore this must be $1_A$  \\ \hline
60 & and therefore there must be $1_B$ & and therefore this must be $1_B$  \\ \hline
78 & For each $i\in{\Ord}$ we choose a subobject $G(i)\in S$ such that & For each $i\in{\Ord}$ we choose an object $G(i)\in S$ such that  \\ \hline
78 & the condition $F(k_i,k_j)\notin\Iso$ follows from $G( i)\ne G(j)$ & the condition $F(k_i,k_j)\notin\Iso$ implies that $G( i)\ne G(j)$  \\ \hline
\end{tabular}
\end{center}

\begin{center}
\begin{tabular}{|c|p{2.5in}|p{2.5in}|} \hline
78 &	
then we choose the limit $k_i$ as $k_j$ & then we choose $k_j$ as the limit of $k_i$ 
\\ \hline
78 &	
We obtain a mapping $G:{\Ord}\to S$ that satisfies the injectivity condition, because for any $i\ne j$ the condition $F(k_i,k_j)\notin\Iso$ follows from $G( i)\ne G(j)$. 
&
We obtain a mapping $G:{\Ord}\to S$ that satisfies the injectivity condition, because for any $i\ne j$ the condition $F(k_i,k_j)\notin\Iso$ implies $G( i)\ne G(j)$.
\\ \hline
87	& By immediate mono- and epimorphisms we mean morphisms, which in some manuals are called extremal. The reasons why we use this term are explained below.	& \\ \hline
88 & $1^\circ$. If $\ph\circ\mu$ is an immediate monomorphism, then $\mu$  is an immediate monomorphism as well. &
$1^\circ$. If $\ph$ is a monomorphism, and $\ph\circ\mu$ is an immediate monomorphism, then $\mu$  is an immediate monomorphism as well. \\ \hline
88 & $3^\circ$. If $\e\circ\ph$ is an immediate epimorphism, then $\e$  is an immediate epimorphism as well. & 
$3^\circ$. If $\ph$ is a epimorphism, and $\e\circ\ph$ is an immediate epimorphism, then $\e$  is an immediate epimorphism as well.
\\ \hline
144	& Fiber product and fiber sums	& Fiber products and fiber sums \\ \hline
263	& $\chi\oslash\psi=\bullet_{A,B,D}\circ \bullet_{B,C,D}\otimes 1_{ B\oslash A}
\circ \chi\otimes 1_{ C\oslash B}\otimes\psi\circ \triangleleft_{ D\oslash C}\otimes 1_{ C\oslash B}\otimes\triangleright_{\frac{B}{C}}$ & 
$\chi\oslash\psi=\bullet_{A,B,D}\circ \bullet_{B,C,D}\otimes 1_{ B\oslash A}
\circ \chi\otimes 1_{ C\oslash B}\otimes\psi\circ (\triangleright_{I\otimes  C\oslash B})^{-1}\circ
(\triangleleft_{ C\oslash B})^{-1}$ 
\\ \hline  	 
334	& However, this is true, since by the definition & This is true, since by the definition \\ \hline
352	& let $\tau$ be the locally convex topology	& let $\tau$ be a locally convex topology \\ \hline
592 & that maps the elementary tensors & that turns the elementary tensors \\ \hline
592 & where the family of these operators & and the family of these operators \\ \hline
592 & For any triple of stereotype spaces $X$, $Y$, $Z$, there exists an operator &	For each stereotype spaces $X$, $Y$, $Z$, there exists a linear continuous mapping  \\ \hline
460 & Theorem 3.6.15. Let $X$ be pseudosaturated and $Y$ be complete. & Theorem 3.6.15. Let $X$ be pseudosaturated and $Y$ be pseudocomplete.  \\ \hline
713 & $\d(a\cdot  t) = \d t	$ & $I(\d(a\cdot  t)) = I(\d t)$  \\ \hline
713 & $\d(t\cdot  a) = \d t$ & $I(\d(t\cdot a)) = I(\d t)$ \\ \hline
713 & $\d(t^{-1}) = \d t$ & $I(\d(t^{-1})) = I(\d t)$ \\ \hline
713 & $\d I(a^{-1}\cdot  s) = \d I(s)$ & $I(\d (a^{-1}\cdot  s)) = I(\d s) $ \\ \hline
\end{tabular}
\end{center}

\paragraph{На с.91.}
Некоторые исправления не умещающаются в таблицу, и их поэтому придется выписать отдельно. В частности, на с.91 доказательство свойства $1^\circ$ должно быть заменено следующим:

1. Let morphisms $C\overset{\mu}{\longrightarrow}D$ and $D\overset{\ph}{\longrightarrow}E$ be given, and it is known that $\ph\circ\mu$ is a strong monomorphism. Since $\ph\circ\mu$ is a monomorphism, it follows from the property $1^\circ$ on p.80  that $\mu$ is also a monomorphism, so it remains to prove that in the diagram (1.41) the commutative square splits into commutative triangles if $\e$ is an epimorphism. For this we first complete the square
$$
 
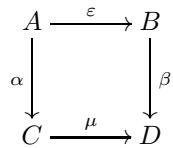
\xymatrix @R=2.5pc @C=2.5pc
 {
 A\ar[d]_{\alpha}\ar[r]^{\e} & B\ar[d]^{\beta}\\
 C\ar[r]^{\mu} & D
 }
$$
to the diagram (in which, by definition, $\beta'=\ph\circ\beta$)
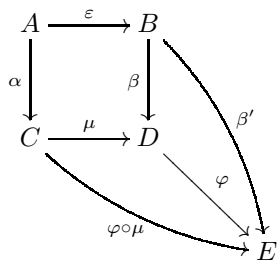
$$
 \xymatrix @R=2.5pc @C=2.5pc
 {
 A\ar[d]_{\alpha}\ar[r]^{\e} & B\ar[d]_{\beta}\ar@/^2ex/[rdd]^{\beta'} & \\
 C\ar[r]^{\mu}\ar@/_2ex/[rrd]_{\ph\circ\mu} & D\ar[dr]^{\ph} & \\
                        &&E
 }
$$
then we discard the inner vertex, and we get a new square
$$
 \xymatrix @R=2.5pc @C=2.5pc
 {
 A\ar[d]_{\alpha}\ar[r]^{\e} & B\ar[d]^{\beta'}\\
 C\ar[r]^{\ph\circ\mu} & E
 }
$$
which must be divided into triangles, since by condition $\ph\circ\mu$ is a strong monomorphism:
$$
\begin{diagram}
\node{A}\arrow{s,l}{\alpha}\arrow{e,t}{\e}\node{B}\arrow{s,r}{\beta'}\arrow{sw,t,--}{\delta}
\\
\node{C}\arrow{e,b}{\ph\circ\mu}\node{E}
\end{diagram}
$$
What is important for us, is the equality $\alpha=\delta\circ\e$. Now we consider the diagram  
\beq\label{PROOF:ph-mu=Smono-=>-mu=Smono}
\begin{diagram}
\node{A}\arrow{s,l}{\alpha}\arrow{e,t}{\e}\node{B}\arrow{s,r}{\beta}\arrow{sw,t,--}{\delta}
\\
\node{C}\arrow{e,b}{\mu}\node{D}
\end{diagram}
\eeq
Here the equality $\alpha=\delta\circ\e$ imply by (1.42) the equality $\beta=\mu\circ\delta$. I.e. the whole diagram (1.41) is commutatuve. This means that $\mu$ is a strong monomorphism.

\paragraph{На с.263.}

На с.263 диаграмма 
$$
\xymatrix @R=2.pc @C=6.0pc 
{
 C\oslash B\ar@{-->}[rr]^{\chi\oslash\psi}\ar[d]_(.4){\triangleleft_{ D\oslash C}\otimes 1_{ C\oslash B}\otimes\ \triangleright_{\frac{B}{C}}} & &  D\oslash A \\
I\otimes  C\oslash B\otimes I\ar[r]_(.4){\chi\otimes 1_{ C\oslash B}\otimes\ \psi} &
 D\oslash C\otimes  C\oslash B\otimes  B\oslash A\ar[r]_(.6){\bullet_{B,C,D}\otimes 1_{ B\oslash A}} &
D\oslash B\otimes  B\oslash A\ar[u]_{\bullet_{A,B,D}}
\\
}
$$ 
--- должна быть заменена на диаграмму
$$
\xymatrix @R=2.pc @C=5.0pc 
{
 C\oslash B\ar@{-->}[r]^{\chi\oslash\psi}\ar[d]_(.4){(\triangleleft_{ C\oslash B})^{-1}} &   D\oslash A \\
I\otimes  C\oslash B\ar[d]_(.4){(\triangleright_{I\otimes  C\oslash B})^{-1}} & D\oslash B\otimes  B\oslash A\ar[u]_{\bullet_{A,B,D}}\\
I\otimes  C\oslash B\otimes I\ar[r]_(.4){\chi\otimes 1_{ C\oslash B}\otimes\ \psi} &
 D\oslash C\otimes  C\oslash B\otimes  B\oslash A\ar[u]_(.6){\bullet_{B,C,D}\otimes 1_{ B\oslash A}}
}
$$

\paragraph{На с.719.}

На с.719 в Examples 5.2.26-5.2.28 всюду инволюция должна быть без подкрутки!

\paragraph{Тождества $\C^M\circledast\C^N\cong \C^{M\times N}$ и $\C_M\odot\C_N\cong \C_{M\times N}$}

В \cite{Akbarov-De-Gruyter-I} на с.625 в теореме 4.5.13 приводится группа тождеств
\beq\label{C^M-odot-C^N-cong-C^(M-times-N)}
\C^M\odot\C^N\cong \C^{M\times N},\qquad \C_M\circledast\C_N\cong \C_{M\times N}
\eeq 
\beq\label{C^M-circledast-C^N-cong-C^(M-times-N)}
\C^M\circledast\C^N\cong \C^{M\times N},\qquad \C_M\odot\C_N\cong \C_{M\times N}
\eeq 
Из них только первые два тождества --- \eqref{C^M-odot-C^N-cong-C^(M-times-N)} --- доказаны в теореме 4.5.13, а вторые --- \eqref{C^M-circledast-C^N-cong-C^(M-times-N)} --- приводятся без доказательств (и из доказанного в тексте не следуют). 
В нашем тексте эти тождества можно считать частными случаями уже доказанных тождеств \eqref{prod_X_i-circledast-prod_Y_j} и
\eqref{bigoplus_X_i-odot-bigoplus_Y_j}, но у них есть более простое доказательство, которое мы здесь тоже приведем. Оно опирается на следующее наблюдение.

\blm\label{LM:operatory-C^M->C_N}
Всякий оператор $\ph:\C^M\to\C_N$ конечномерен.
\elm
\bpr
Рассмотрим базис $\{1_x;\ x\in M\}$ в $\C^M$ из характеристических функций одноточечных множеств:
\beq\label{DEF:1_x}
1_x(y)=\begin{cases}1, & y=x\\ 0, & y\ne x\end{cases}.
\eeq
Множество $\{1_x;\ x\in M\}$ ограничено (и даже вполне ограничено) в $\C^M$, поэтому под действием непрерывного линейного отображения $\ph$ оно должно превратиться в ограниченное множество $\{\ph(1_x);\ x\in M\}$ в $\C_N$. Но в $\C_N$ всякое ограниченное множество содержится в некотором конечномерном подпространстве. Значит, $\{\ph(1_x);\ x\in M\}$ тоже содержится в некотором конечномерном подпространстве $Y\subseteq \C_N$. Отсюда следует, что вообще весь образ $\ph$ содержится в $Y$,
$$
\ph(\C^M)\subseteq Y,\qquad \dim Y<\infty.
$$
То есть $\ph$ --- конечномерный оператор.
\epr

Из леммы \ref{LM:operatory-C^M->C_N} следует двойное тождество
\beq\label{C_N-oslash-C^M-cong-(C_M)_N}
\C_N:\C^M\cong \C_N\oslash\C^M\cong (\C_M)_N
\eeq
Действительно, поскольку всякий элемент $\ph\in \C_N:\C^M$, то есть оператор $\ph:\C^M\to\C_N$, конечномерен, он представляет собой конечную линейную комбинацию вида
\beq\label{C_N-oslash-C^M-cong-(C_M)_N-1}
\ph(u)=\sum_{i=1}^n \alpha_{y_i}(u)\cdot \delta^{y_i},\qquad u\in\C^M, 
\eeq 
где $y_i\in N$ --- некая конечная последовательность в $N$, $\delta^y$ --- базис в $\C_N$, состоящий из характеристических функций одноточечных множеств (в \cite{Akbarov-De-Gruyter-I} мы условились обозначать его иначе, чем \eqref{DEF:1_x}, чтобы было понятно, что это базис в пространстве другого типа), а $\alpha_y$ --- линейные непрерывные функционалы на $\C^M$, определяемые равенством
\beq\label{C_N-oslash-C^M-cong-(C_M)_N-2}
\alpha_y(u)=\ph(u)_y,\qquad u\in\C^M, \ y\in N
\eeq 
(по определению, $\ph(u)$, как элемент пространства $\C_N$, представляет собой семейство $\{\ph(u)_y;\ y\in N\}$).

Из \eqref{C_N-oslash-C^M-cong-(C_M)_N-1} и \eqref{C_N-oslash-C^M-cong-(C_M)_N-2} следует, что $\ph$ можно понимать как семейство функционалов 
$$
\{\alpha_y;\ y\in N\}\subseteq (\C^M)^\star
$$
из которых почти все равны нулю, и, поскольку $(\C^M)^\star\cong\C_M$, это эквивалентно представлению $\ph$ как семейство элементов пространства $\C_M$
$$
\{\alpha_y;\ y\in N\}\subseteq \C_M
$$
из которых почти все равны нулю. То есть $\ph$ представляется как элемент пространства $(\C_M)_N$:
$$
\{\alpha_y;\ y\in N\}\in (\C_M)_N.
$$
Теперь нужно заметить, что непрерывному изменению $\ph\in \C_N:\C^M$ будет соответствовать непрерывное изменение 
$\{\alpha_y;\ y\in N\}\in (\C_M)_N$, и это будет означать изоморфизм локально выпуклых пространств
$$
\C_N:\C^M\cong (\C_M)_N.
$$
Отсюда, в свою очередь можно будет заключить, что пространство $\C_N:\C^M$ псевдонасыщено, и поэтому оно совпадает с пространством $\C_N\oslash\C^M$:
$$
\C_N:\C^M=\C_N\oslash\C^M.
$$ 
И вместе это доказывает \eqref{C_N-oslash-C^M-cong-(C_M)_N}.

Из \eqref{C_N-oslash-C^M-cong-(C_M)_N} мы теперь получаем правое тождество в \eqref{C^M-circledast-C^N-cong-C^(M-times-N)},
$$
\C_M\odot\C_N\cong \C_N\oslash(\C_M)^\star\cong \C_N\oslash\C^M\cong\eqref{C_N-oslash-C^M-cong-(C_M)_N}\cong (\C_M)_N
\cong \C_{M\times N}
$$
а из него будет следовать левое:
$$
\C^M\circledast\C^N\cong (\C_M\odot\C_N)^\star\cong (\C_{M\times N})^\star\cong \C^{M\times N}.
$$

\renewcommand{\thesection}{\Alph{section}.}
\renewcommand{\thesubsection}{(\alph{subsection})}


\tableofcontents

\end{document}